\documentclass[11pt]{amsart}
\usepackage{amsfonts}
\usepackage{fullpage, amsmath, amsthm,amsfonts,amssymb,stmaryrd, mathrsfs,amscd}\usepackage{graphicx}
\usepackage{hyperref}
\usepackage{cleveref}
\makeatletter
\usepackage[pdftex]{color}
\usepackage[all, 2cell]{xy}
\UseAllTwocells

\usepackage{tikz}
\usepackage{tikz-cd}
\usetikzlibrary{matrix,arrows,decorations.pathmorphing}
\usepackage{faktor}

\usepackage[shortlabels]{enumitem}

\numberwithin{equation}{section}

\theoremstyle{definition}
\newtheorem{theorem}{Theorem}
\numberwithin{theorem}{section}
\newtheorem{corollary}[theorem]{Corollary}

\newtheorem{lemma}[theorem]{Lemma}
\newtheorem{definition}[theorem]{Definition}
\newtheorem{proposition}[theorem]{Proposition}
\newtheorem{remark}[theorem]{Remark}
\newtheorem{example}[theorem]{Example}

\newtheorem{assumptions}[theorem]{Assumptions}
\newtheorem{assumption}[theorem]{Assumption}
\newtheorem{notation}[theorem]{Notation}

\newcommand{\quash}[1]{}  

\newcommand{\fraka}{{\mathfrak a}}
\newcommand{\frakb}{{\mathfrak b}}
\newcommand{\frakc}{{\mathfrak c}}

\newcommand{\frakg}{{\mathfrak g}}
\newcommand{\frakh}{{\mathfrak h}}

\newcommand{\frakm}{{\mathfrak m}}

\newcommand{\frakp}{{\mathfrak p}}

\newcommand{\fraks}{{\mathfrak s}}
\newcommand{\frakt}{{\mathfrak t}}
\newcommand{\fraku}{{\mathfrak u}}

\newcommand{\frakA}{{\mathfrak A}}

\newcommand{\frakC}{{\mathfrak C}}

\newcommand{\frakM}{{\mathfrak M}}

\newcommand{\frakX}{{\mathfrak X}}

\newcommand{\bA}{{\mathbb A}}
\newcommand{\bB}{{\mathbb B}}
\newcommand{\bC}{{\mathbb C}}
\newcommand{\bD}{{\mathbb D}}

\newcommand{\bF}{{\mathbb F}}
\newcommand{\bG}{{\mathbb G}}

\newcommand{\bL}{{\mathbb L}}

\newcommand{\bN}{{\mathbb N}}
\newcommand{\bO}{{\mathbb O}}
\newcommand{\bP}{{\mathbb P}}
\newcommand{\bQ}{{\mathbb Q}}
\newcommand{\bR}{{\mathbb R}}

\newcommand{\bT}{{\mathbb T}}

\newcommand{\bX}{{\mathbb X}}

\newcommand{\bZ}{{\mathbb Z}}

\newcommand{\bfA}{{\mathbf A}}
\newcommand{\bfB}{{\mathbf B}}
\newcommand{\bfC}{{\mathbf C}}
\newcommand{\bfD}{{\mathbf D}}
\newcommand{\bfE}{{\mathbf E}}
\newcommand{\bfF}{{\mathbf F}}
\newcommand{\bfG}{{\mathbf G}}
\newcommand{\bfH}{{\mathbf H}}
\newcommand{\bfI}{{\mathbf I}}

\newcommand{\bfK}{{\mathbf K}}
\newcommand{\bfL}{{\mathbf L}}
\newcommand{\bfM}{{\mathbf M}}
\newcommand{\bfN}{{\mathbf N}}

\newcommand{\bfP}{{\mathbf P}}

\newcommand{\bfR}{{\mathbf R}}
\newcommand{\bfS}{{\mathbf S}}
\newcommand{\bfT}{{\mathbf T}}

\newcommand{\bfa}{{\mathbf a}}

\newcommand{\bff}{{\mathbf f}}

\newcommand{\mA}{{\mathcal A}}
\newcommand{\mB}{{\mathcal B}}

\newcommand{\mE}{{\mathcal E}}
\newcommand{\mF}{{\mathcal F}}
\newcommand{\mG}{{\mathcal G}}
\newcommand{\mH}{{\mathcal H}}
\newcommand{\mI}{{\mathcal I}}
\newcommand{\mJ}{{\mathcal J}}
\newcommand{\mK}{{\mathcal K}}
\newcommand{\mL}{{\mathcal L}}
\newcommand{\mM}{{\mathcal M}}
\newcommand{\mN}{{\mathcal N}}
\newcommand{\mO}{{\mathcal O}}
\newcommand{\mP}{{\mathcal P}}
\newcommand{\mQ}{{\mathcal Q}}

\newcommand{\mS}{{\mathcal S}}
\newcommand{\mT}{{\mathcal T}}
\newcommand{\mU}{{\mathcal U}}
\newcommand{\mV}{{\mathcal V}}
\newcommand{\mW}{{\mathcal W}}
\newcommand{\mX}{{\mathcal X}}

\newcommand{\mZ}{{\mathcal Z}}

\newcommand{\scrA}{{\mathscr A}}
\newcommand{\scrB}{{\mathscr B}}

\newcommand{\scrS}{{\mathscr S}}

\newcommand{\al}{{\alpha}}

\newcommand{\ga}{{\gamma}}
\newcommand{\Ga}{{\Gamma}}
\newcommand{\la}{{\lambda}}
\newcommand{\La}{{\Lambda}}

\newcommand{\barcons}{\mathrm{Bar}}
\newcommand{\cycl}{|| \cdot ||}
\newcommand{\Sym}{\mathrm{Sym}}
\newcommand{\ext}{\mathrm{ext}}
\newcommand{\id}{\mathrm{id}}
\newcommand{\pr}{\mathrm{pr}}
\newcommand{\pt}{\mathrm{pt}}

\newcommand{\bs}{\backslash}
\newcommand{\Quad}{\mathrm{Quad}}
\newcommand{\sw}{\mathrm{sw}}

\newcommand{\spc}{\mathrm{Ani}}
\newcommand{\cat}{\widehat{\mathcal{C}\mathrm{at}}_\infty}

\newcommand{\lincat}{\mathrm{Lincat}}
\newcommand{\cptcat}{\lincat^{\mathrm{cg}}}
\newcommand{\catid}{\mathrm{Lincat}^{\mathrm{Perf}}}

\newcommand{\prl}{\mathcal{P}\mathrm{r}^{\mathrm{L}}}

\newcommand{\colim}{\mathrm{colim}}
\newcommand{\op}{\mathrm{op}}
\newcommand{\rev}{\mathrm{rev}}
\newcommand{\fun}{\mathrm{Fun}}
\newcommand{\Fun}{\mathrm{Fun}}
\newcommand{\map}{\mathrm{Map}}
\newcommand{\Map}{\mathrm{Map}}

\newcommand{\Hom}{\mathrm{Hom}}
\newcommand{\rhom}{\underline{\mathrm{Hom}}}

\newcommand{\End}{\mathrm{End}}
\newcommand{\Ext}{\mathrm{Ext}}
\newcommand{\aut}{\mathrm{Aut}}
\newcommand{\Aut}{\mathrm{Aut}}

\newcommand{\ren}{\mathrm{ren}}
\newcommand{\ind}{\mathrm{Ind}}
\newcommand{\Ind}{\mathrm{Ind}}
\newcommand{\res}{\mathrm{res}}

\newcommand{\tot}{\mathrm{Tot}}

\newcommand{\cpt}{\omega}

\newcommand{\alg}{\mathrm{Alg}}
\newcommand{\calg}{\mathrm{CAlg}}
\newcommand{\modu}{\mathrm{Mod}}
\newcommand{\Mod}{\mathrm{Mod}}
\newcommand{\lmodu}{\mathrm{LMod}}

\newcommand{\rmodu}{\mathrm{RMod}}

\newcommand{\bmodu}{\mathrm{BMod}}
\newcommand{\BMod}{\mathrm{BMod}}

\newcommand{\spec}{\operatorname{spec}}
\newcommand{\Spec}{\operatorname{Spec}}
\newcommand{\Spf}{\operatorname{Spf}}

\newcommand{\Spd}{\operatorname{Spd}}

\newcommand{\aff}{\mathrm{CAlg}}
\newcommand{\Aff}{\mathrm{Aff}}
\newcommand{\Sch}{\mathrm{Sch}}
\newcommand{\algsp}{\mathrm{AlgSp}}

\newcommand{\ArStk}{\mathrm{AlgStk}}
\newcommand{\prestk}{\mathrm{PreStk}}

\newcommand{\qc}{\mathrm{qc}}
\newcommand{\qs}{\mathrm{qs}}
\newcommand{\qcqs}{\mathrm{qcqs}}
\newcommand{\aft}{\mathrm{afp}}
\newcommand{\laft}{\mathrm{l}\aft}
\newcommand{\sep}{\mathrm{sep}}
\newcommand{\pfp}{\mathrm{Pfp}}

\newcommand{\coh}{\mathrm{coh}}
\newcommand{\Coh}{\mathrm{Coh}}
\newcommand{\qcoh}{\mathrm{QCoh}}
\newcommand{\Qcoh}{\mathrm{QCoh}}
\newcommand{\QC}{\mathrm{QC}^!}
\newcommand{\indcoh}{\mathrm{IndCoh}}

\newcommand{\et}{\Acute{e}\mathrm{t}}
\newcommand{\proet}{\mathrm{pro\Acute{e}t}}
\newcommand{\fpqc}{\mathrm{fpqc}}

\newcommand{\corr}{\mathrm{Corr}}
\newcommand{\verti}{\mathrm{V}}
\newcommand{\horiz}{\mathrm{H}}
\newcommand{\all}{\mathrm{All}}
\newcommand{\isom}{\mathrm{Iso}}
\newcommand{\adm}{\mathrm{Adm}}
\newcommand{\verd}{\mathbb{D}}

\newcommand{\rg}{\mathrm{R}\Gamma}
\newcommand{\cons}{\mathrm{cons}}
\newcommand{\ctf}{\mathrm{ctf}}
\newcommand{\shv}{\mathrm{Shv}}
\newcommand{\Shv}{\mathrm{Shv}}
\newcommand{\fgshv}{\mathrm{Shv}_{\mathrm{f.g.}}}

\newcommand{\sshv}{\shv^{*}}
\newcommand{\cshv}{\shv_{\mathrm{c}}}
\newcommand{\scshv}{\shv_{\mathrm{c}}^{*}}
\newcommand{\rshv}{\ind\fgshv}

\newcommand{\os}{\otimes^!}
\newcommand{\cohdual}{\omega}
\newcommand{\consdual}{\omega}
\newcommand{\frobdual}{\omega}
\newcommand{\red}{\mathrm{red}}

\newcommand{\rstar}{^{\sharp}}
\newcommand{\horizl}{^{\star}}
\newcommand{\horizr}{_{\star}}

\newcommand{\vertl}{_{\dagger}} 
\newcommand{\vertr}{^{\dagger}} 
\newcommand{\renflat}{\blacktriangle} 

\newcommand{\tateshift}[1]{\langle{#1}\rangle}

\newcommand{\mon}{\mathrm{mon}}

\newcommand{\ftor}{\mathrm{ftor}}

\newcommand{\hchi}{\hat{\chi}}
\newcommand{\pialg}{\pi^c_{1}}
\newcommand{\Vogan}{V}

\newcommand{\tr}{\mathrm{Tr}}
\newcommand{\geo}{\mathrm{geo}}

\newcommand{\trg}{\mathrm{Tr}_{\mathrm{geo}}}

\newcommand{\der}{\mathcal{D}}

\newcommand{\cl}{\mathrm{cl}}
\newcommand{\proj}{\mathcal{P}}
\newcommand{\HH}{\mathrm{HH}}

\newcommand{\sht}{\mathrm{Sht}}
\newcommand{\Sht}{\mathrm{Sht}}
\newcommand{\glob}{\mathrm{glob}}
\newcommand{\loc}{\mathrm{loc}}
\newcommand{\Loc}{\mathrm{Loc}}
\newcommand{\kot}{\mathrm{Isoc}}

\newcommand{\plstk}{\mathrm{Stk}^{\mathrm{pl}}}

\newcommand{\vplstk}{\mathrm{Stk}^{\mathrm{vpl}}}
\newcommand{\qplstk}{\mathrm{Stk}^{\mathrm{qpl}}}
\newcommand{\qcqplstk}{\mathrm{Stk}^{\qc.\mathrm{qpl}}}

\newcommand{\plsp}{\mathrm{AlgSp}^{\mathrm{spl}}}

\newcommand{\perf}{\mathrm{CAlg}^{\mathrm{perf}}}
\newcommand{\Perf}{\mathrm{Perf}}
\newcommand{\pf}{\mathrm{perf}}
\newcommand{\psch}{\mathrm{Sch}^{\mathrm{perf}}}
\newcommand{\pfpsch}{\mathrm{Sch}^{\mathrm{pfp}}}
\newcommand{\qcsp}{\mathrm{AlgSp}^{\mathrm{perf}}}
\newcommand{\pfpsp}{\mathrm{AlgSp}^{\mathrm{pfp}}}

\newcommand{\indsp}{\mathrm{IndAlgSp}}
\newcommand{\indarstk}{\mathrm{IndArStk}}
\newcommand{\indsch}{\mathrm{IndSch}}
\newcommand{\indplstk}{\mathrm{IndStk}^{\mathrm{pl}}}
\newcommand{\qcindplstk}{\mathrm{IndStk}^{\qc.\mathrm{pl}}}
\newcommand{\indqplstk}{\mathrm{IndStk}^{\mathrm{qpl}}}
\newcommand{\qcindqplstk}{\mathrm{IndStk}^{\qc.\mathrm{qpl}}}
\newcommand{\indvplstk}{\mathrm{IndStk}^{\mathrm{vpl}}}
\newcommand{\sfplstk}{\mathrm{sIndStk}^{\mathrm{pl}}}
\newcommand{\qcsfplstk}{\mathrm{sIndStk}^{\qc.\mathrm{pl}}}
\newcommand{\sfvplstk}{\mathrm{sIndStk}^{\mathrm{vpl}}}

\newcommand{\sing}{\mathrm{Sing}}
\newcommand{\Sing}{\mathrm{Sing}}

\newcommand{\can}{\mathrm{can}}

\newcommand{\rep}{\mathrm{Rep}}
\newcommand{\Rep}{\mathrm{Rep}}
\newcommand{\crep}{\rep_{\mathrm{c}}}
\newcommand{\cRep}{\mathrm{Rep}_{\mathrm{c}}}
\newcommand{\fgrep}{\rep_{\mathrm{f.g.}}}

\newcommand{\rrep}{\ind\fgrep}

\newcommand{\iw}{\mathrm{Iw}}
\newcommand{\sph}{\mathrm{K}}

\newcommand{\hk}{\mathrm{Hk}}
\newcommand{\Hk}{\mathrm{Hk}}
\newcommand{\locsys}{\mathrm{Loc}}

\newcommand{\unip}{\mathrm{unip}}
\newcommand{\tame}{\mathrm{tame}}
\newcommand{\unr}{\mathrm{unr}}
\newcommand{\cohspr}{\mathrm{CohSpr}}
\newcommand{\cont}{\mathrm{cont}}
\newcommand{\Ch}{\mathrm{Ch}}
\newcommand{\ad}{\mathrm{ad}}
\newcommand{\Ad}{\mathrm{Ad}}
\newcommand{\af}{\mathrm{aff}}
\newcommand{\cind}{c\text{-}\mathrm{ind}}

\newcommand{\Bun}{\mathrm{Bun}}

\newcommand{\pFr}{\mathrm{pFr}}

\newcommand{\Nt}{\mathrm{Nt}}
\newcommand{\av}{\mathrm{Av}}

\newcommand{\Gr}{\mathrm{Gr}}
\newcommand{\Fl}{\mathrm{Fl}}
\newcommand{\xch}{\bX^\bullet}
\newcommand{\xcoch}{\bX_\bullet}

\newcommand{\iso}{\mathrm{iso}}
\newcommand{\GL}{\mathrm{GL}}
\newcommand{\SL}{\mathrm{SL}}
\newcommand{\fg}{\mathrm{f.g.}}

\newcommand{\coWhit}{\mathrm{coWhit}}
\newcommand{\IW}{\mathrm{IW}}
\newcommand{\asp}{\mathrm{asp}}
\newcommand{\GG}{\mathrm{GG}}

\newcommand{\Til}{\mathrm{Til}}

\newcommand{\Sh}{\mathrm{Sh}}
\newcommand{\bfSh}{\mathbf{Sh}}
\newcommand{\Igs}{\mathrm{Igs}}

\newcommand{\Igss}{\mI gs}

\title{Tame categorical local Langlands correspondence}
\author{Xinwen Zhu}

\begin{document}
\maketitle
\begin{abstract}
   In one of our previous articles, we outlined the formulation of a version of the categorical arithmetic local Langlands conjecture. The aims of this article are threefold. First, we provide a detailed account of one component of this conjecture: the local Langlands category. Second, we aim to prove this conjecture in the tame case for quasi-split unramified reductive groups. Finally, we will explore the first applications of such categorical equivalence.
\end{abstract}

\tableofcontents

\section{Introduction}\label{sec: introduction}

\subsection{Backgrounds and motivations}

In \cite{zhu2020coherent}, we sketched the formulation of a version of the categorical arithmetic local Langlands conjecture. 
The aims of this article are threefold. First, we provide a detailed account of one component of this conjecture: the local Langlands category. Second, we aim to prove this conjecture in the tame case for quasi-split unramified reductive groups. Finally, we will explore the first applications of such categorical equivalence.

Let us start with some motivations for the categorial arithmetic local Langlands conjecture.
Let $F$ be a non-archimedean local field, i.e., a finite extension of $\bQ_p$ or of $\bF_p(\!(\varpi)\!)$, and let 
$W_F\subset \Ga_F$ be its Weil group and the Galois group. Let $G$ be connected reductive group over $F$, and let ${}^LG=\hat{G}\rtimes\Ga_{\tilde F/F}$ be its Langlands dual group.

Recall that the \textbf{classical local Langlands correspondence} roughly predicts a natural bijection: 
\begin{multline*}
\Bigl\{\mbox{Smooth irreducible representations of } G(F) \Bigr\}\leftrightarrow \\
 \Bigl\{\mbox{Langlands parameters } \varphi\colon W_F \rightarrow {}^LG \mbox{ up to } \hat{G} \mbox{ conjugation by } G\Bigr\}. 
\end{multline*}

For $\GL_n$, ``naturality" can be made precise and the local Langlands correspondence is a theorem, proved by Laumon-Rapoport-Stuhler \cite{Laumon.Rapoport.Stuhler}  when $F$ is of positive characteristic,
and by Harris-Taylor \cite{Harris.Taylor}, and independently by Henniart \cite{Henniart} when $F$ is of characteristic zero.

For a general reductive group $G$, however, ``naturality" is hard to formulate. In fact, the set of Langlands parameters needs to be enhanced. 
For example, Kazhdan-Lusztig \cite{Kazhdan.Lusztig.Deligne.Langlands} constructed (for $G$ split)  an injective map
\begin{multline*}
\Bigl\{\mbox{Smooth  irreducible representations  of } G(F) \mbox{ with Iwahori fixed vectors} \Bigr\}\hookrightarrow \\
 \Bigl\{(\varphi,r) \mid \varphi \colon  W_F \rightarrow \hat{G}, r\in \rep(C_{\hat{G}}(\varphi))\Bigr\}/ \hat{G},
\end{multline*}
Here $\varphi$ is a Langlands parameter as described above, and $r$ is a representation of the stabilizer $C_{\hat{G}}(\varphi)$ of $\varphi$ under the conjugation action of $\hat{G}$.
The appearance of $r$ suggests that there are stacks involved in the story. Namely, such $r$ can be interpreted as a coherent sheaf on the stack 
\[
\{\varphi\}/C_{\hat{G}}(\varphi)\cong \big\{\hat{G}\mbox{-orbit of } \varphi: W_F\to {}^LG\big\}/\hat{G}.
\]

The geometric Langlands program suggests that the local Langlands correspondence can--and probably needs to--be lifted to an equivalence of categories. 
Namely, instead of considering the set of isomorphism classes of pairs $(\varphi,r)$, one should consider the category of coherent sheaves on $\locsys_{{}^cG,F}$, where $\locsys_{{}^cG,F}$ is the stack of local Langlands parameters, classifying continuous (in appropriate sense) $\ell$-adic representations of $W_F$ with values in the $C$-group ${}^cG$ of $G$ (which is a slight variant the usual Langlands dual group ${}^L\!G$ of $G$).  
Such a stack $\locsys_{{}^cG,F}$ indeed exists, see \cite[\textsection{3.1}]{zhu2020coherent}, and also \cite{dat2020moduli} and \cite[Chapter VIII]{Fargues.Scholze.geometrization}. It is a classical algebro-geometric object, specifically the disjoint union of affine schemes of finite type (over $\bZ_\ell$) modulo the action of $\hat{G}$. Therefore, the category $\Coh(\locsys_{{}^cG,F})$ of coherent sheaves on $\locsys_{{}^cG,F}$ makes sense and serves as the replacement for the set of Langlands parameters in the categorical local Langlands conjecture. 

The categorification of the representation-theoretic side turns out to be much more involved. 
Naively, one might guess that we could replace the set of smooth irreducible representations of $G(F)$ by the (derived) category $\Rep(G(F))$ of smooth representations. However, this is not quite sufficient. 
As has long been observed, to obtain a good parameterization of representations in terms of Langlands parameters, it is better to consider not only the representations of the $p$-adic group $G(F)$ itself, but also the representations of its various (extended, pure) inner forms.
However, there is considerable evidence suggesting that we should study the representation theory of  $G(F)$ alongside a collection of groups $\{G_b(F)\}_{b\in B(G)}$, indexed by a certain set $B(G)$. Each $G_b(F)$ is a(n extended) inner form of a Levi subgroup of $G$ (say $G$ is quasi-split). In addition, the categories $\{\Rep(G_b(F))\}_b$ can be glued  together as the category of sheaves on certain geometric objects. Indeed, the set $B(G)$ was first introduced by Kottwitz (and is now referred to as the Kottwitz set) in the study of mod $p$ points of Shimura varieties.

There are two ways to make this idea precise. One is developed by Fargues-Scholze in their monumental
document \cite{Fargues.Scholze.geometrization}. 
In this approach, the set $B(G)$ is regarded as the set of points of the $v$-stack  $\Bun_G$ of $G$-bundles on the Fargues-Fontaine curve, and the glued category is defined as the category of appropriately defined $\ell$-adic sheaves on $\Bun_G$. This definition is quite sophisticated, relying on recent progress in $p$-adic geometry and condensed mathematics. 

In this work, we take a different approach to introduce another category $\Shv(\kot_G)$, which can be regarded as an alternative candidate on the representation-theoretic side of the categorical local Langlands conjecture. This approach, although still involved, remains within the realm of traditional $\ell$-adic formalism in algebraic geometry. This category is implicitly considered in \cite{xiao2017cycles}, and its definition is outlined in \cite{zhu2020coherent}. See also \cite{gaitsgory2016geometric} for an informal account. We will let $\La$ be a certain $\bZ_\ell$-algebra (e.g. $\La=\bF_\ell, \bQ_\ell, \bZ_\ell$ or finite extensions of such), which serves as the coefficient ring for our sheaf theory in the sequel.

To introduce $\kot_G$, let us first recall the definition of the Kottwitz set $B(G)$.
Let $k$ be an algebraic closure of the residue field $k_F$ of $F$. We write $q=\sharp k_F$.
Let $\breve F$ be the completion of the maximal unramified extension of $F$, and let $\sigma\in\Aut(\breve F/F)$ be the automorphism that lifts the $q$-Frobenius automorphism of $k$. Then $B(G)$ is defined as the isomorphism classes of $F$-isocrystals with $G$-structures $(\mE,\psi)$, which consist of a $G$-torsor $\mE$ over $\Spec \breve F$ equipped with a $G$-torsor isomorphism $\psi: \sigma^*\mE\simeq \mE$. When $G=\GL_n$, these can be further explicitly described as pairs $(V,\psi)$, consisting of an $n$-dimensional $\breve F$-vector space $V$ equipped with a $\sigma$-semilinear bijection. Since any $G$-torsor over $\Spec \breve F$ is trivial, the set $B(G)$ can be identified as the quotient set $G(\breve F)/\sim$, where, and $\sim$ is the equivalence relation given by $g_1\sim g_2$ if $g_1=h^{-1}g_2\sigma(h)$ for some $h\in G(\breve F)$. This is naturally an infinite poset. Minimal elements are called basic elements.

Recall that $F$-isocrystals with $G$-structure appears as the ``crystalline realization" of motives with $G$-structures over $k$. For example, giving an abelian variety $A$ over $k$, its rational Dieudonn\'e module is an $F$-isocrystal. Since abelian varieties (with additional structures) over $k$ form moduli spaces (known as mod $p$ fibers of Shimura varieties), it is natural to expect that $F$-isocrystals with $G$-structures over $k$ also form a moduli space, whose $k$-points are classified by $B(G)$. In addition, by sending an abelian variety over $k$ to its rational Dieudonn\'e module, there should exist morphisms from the mod $p$ Shimura varieties to such moduli spaces of $F$-isocrystals (with additional structures).

This is indeed the case, although the resulting moduli space is  not a familiar geometric object in classical algebraic geometry. To describe it, let $LG$ denote the loop group of $G$, which is a (perfect) ind-group scheme over $k_F$ such that its $k_F$-points are  $G(F)$ and its $k$-points are $G(\breve F)$.
Being an ind-scheme over $k_F$, it admits a $\sharp k_F$-Frobenius endomorphism, denoted by $\sigma$. Then
we consider the (\'etale) quotient stack\footnote{Using $h$-sheafification instead of \'etale sheafification give another version of $\kot_G$. See \Cref{K-points-of-kot(G)} for a discussion.} 
\[
\kot_G:=\frac{LG}{\Ad_\sigma LG},
\] 
where $\Ad_\sigma$ denotes the Frobenius twisted conjugation given by 
\[
\Ad_\sigma: LG\times LG\to LG,\quad  (h,g)\mapsto hg\sigma(h)^{-1}.
\] 
Therefore, $\kot_G$ is a quotient of an infinite dimensional space by an infinite dimensional group, which is a wild object in classical algebraic geometry.
However, it still has many geometric structures. In particular, the category of $\ell$-adic sheaves over $\kot_G$ has nice properties, as we shall see shortly.

But before that, let us mention that the space $\kot_G$ arises naturally from another perspective. This viewpoint also clarifies that why we should consider the category of $\ell$-adic sheaves on $\kot_G$.
To explain this, let us temporarily switch the setting and let $H$ be a reductive group over a finite field $\kappa$. 
Let $\Rep(H(\kappa),\La)$ denote the (derived) category of representations of the finite group $H(\kappa)$ with $\La$-coefficients, where $\La$ is a certain $\bZ_\ell$-algebra as above (e.g., $\La=\bZ_\ell, \overline\bQ_\ell, \overline\bF_\ell$).
On the other hand, we can regard the finite group $H(\kappa)$ as an affine algebraic group over $k=\overline{\kappa}$. Then the classifying stack $\bB H(\kappa)$ of $H(\kappa)$ makes sense as an algebraic stack. Let $ \shv(\bB H(\kappa),\La)$ denote the (derived) category of $\La$-sheaves on $\bB H(\kappa)$.
The starting point of the Deligne-Lusztig theory is following two observations:
\begin{itemize}
\item There is a canonical equivalence of categories $\Rep(H(\kappa),\La)\cong \shv( \bB H(\kappa),\La)$.
\item There is a natural isomorphism of algebraic stacks $\bB H(\kappa)\cong H/\Ad_\sigma H$. Here as above $\Ad_\sigma$ denotes $\sigma$-conjugation, i.e.
$\Ad_\sigma(h)(g)=h^{-1} g \sigma(h),\ g,h\in H$.
\end{itemize}

If we choose a (rational) Borel subgroup $B_H\subset H$. Then the (unipotent part of the) Deligne-Lusztig theory can be regarded as a construction of representations of $H(\kappa)$ via the correspondence
\[
B_H\backslash H/B_H\xleftarrow{\delta} H/\Ad_\sigma B_H\xrightarrow{\Nt} H/\Ad_\sigma H.
\]
Namely, for every complex of $\ell$-adic constructible sheaf $\mF$ on $B_H\backslash H/B_H$, we let
\[
\Ch_{H,\phi}^{\unip}(\mF):=\Nt_*(\delta^!\mF),
\] 
which is a complex of $\ell$-adic constructible sheaf on $H/\Ad_\sigma H$ and can therefore be viewed as a representation of $H(\kappa)$.
For example, if we apply this construction to the $*$-pushforward of the constant sheaf along the locally closed embedding $B_H\bs B_H w B_H/ B_H\subset B_H\bs H/ B_H$, where  $w$ is an element in the (absolute) Weyl group of $H$, we obtain the famous Deligne-Lusztig representation of $H(\kappa)$ on the cohomology of Deligne-Lusztig variety 
\[
X_w=\{ gB_H\in H/B_H\mid g^{-1}\sigma(g)\in B_HwB_H\}.
\] 

From this perspective, $\kot_G$ is clearly an analogue of $\bB H(\kappa)$ when $\kappa$ is replaced by a local field $F$. In addition, the category of $\ell$-adic sheaves on $\kot_G$, if it makes sense, would be the analogue of the category of representations of $H(\kappa)$.
However, there is a significant difference. Namely, unlike $\bB H(\kappa)$, the underlying set of points of $\kot_G$ is no longer a singleton. Indeed, the underlying set of $\kot_G(k)$ is just the Frobenius conjugacy classes in $G(\breve F)$, and therefore it is identified with the Kottwitz set. Additionally, for $b\in G(\breve F)/\sim$, regarded as an object in the groupoid $\kot_G(k)$, its automorphism group
\[
G_b(F)= \{h\in G(\breve F)\mid h^{-1}b\sigma(h)=b\}
\]
is in general not the $p$-adic group
group $G(F)$ itself, but rather  the set of $F$-points of an inner form of a Levi subgroup of $G$.
Only when $b=1$ do we have $G_b(F)=G(F)$. Therefore, the category of  $\ell$-adic sheaves on $\kot_G$, even if it makes sense, will not simply be the category of smooth representations of $G(F)$, but rather a collection of categories of smooth representations of all these groups $G_b(F)$, glued together in an intricate way.

We note that the classical local Langlands correspondence primarily focuses on the smooth representations of $G$ and its (extended, pure) inner forms. Traditionally, there exists another formulation of the local Langlands conjecture (mostly advanced by Vogan), also of a categorical nature, that relates the representations of $G$ and its (extended pure) inner forms in terms of constructible sheaves on some other version of the spaces of Langlands parameters. This raises a question: Do the representations of $G_b(F)$ for non-basic $b$ (or genuine $\ell$-adic sheaves on $\kot_G$) in our story merely serve an artificial extension that could make our categorical conjecture potentially valid, or do they possess substaintial significance within the classical Langlands correspondence? We present an additional motivation for introducing our story: from this perspective, the existence of representations of $G_b(F)$ for non-basic $b$ is not a drawback, but rather an essential feature.

This motivation is rooted in global considerations and applications to arithmetic geometry (see \cite{zhu2020coherent, zhu2022icm} for some surveys), which originally inspired our desire to develop the categorical local Langlands correspondence. In the classical global Langlands correspondence, one studies not just the space of automorphic forms, but also various cohomology groups associated with Shimura varieties or more general locally symmetric spaces in the number field case, and the cohomology of moduli spaces of Shtukas in the function field scenario. As explained in \cite[\textsection{4.7}]{zhu2020coherent}, there exists a conjectural formula for computing such cohomology groups in terms of the coherent cohomology of certain (ind-)coherent sheaves on the stack of global Langlands parameters. The input for this formula--the (ind-)coherent sheaf to compute--is provided by the categorical local Langlands correspondence. Crucially, under the categorical local Langlands correspondence as we are going to develop, these coherent sheaves should correspond to $\ell$-adic sheaves on $\kot_G$ spreaded out over different points of $\kot_G$. In other words, genuine sheaves on $\kot_G$ (rather than merely representations of specific $G_b(F)$) naturally emerge in the study of the global Langlands correspondence.

\subsection{Main results}
  Now we will discuss some of our main results. Along the way, we will provide additional background and motivations.
  
 \subsubsection{Local Langlands category} 
  We start with some geometry of the stack $\kot_G$. For an element $b\in B(G)$, we consider  substacks
    \[
    i_b: \kot_{G,b}\stackrel{j_{b}}{\hookrightarrow}\kot_{G,\leq b}\stackrel{i_{\leq b}}{\hookrightarrow} \kot_G,
    \] 
    where $\kot_{G,\leq b}$ and $\kot_G$ are defined as
    \[
    \kot_{G,\leq b}(R) = \left\{(\mE,\psi)\in \kot_G(R) \big|\ b_x:=(\mE_x,\psi_x)\leq b, \  x\in \Spec R \right\},
    \]  
    \[
    \kot_{G,b} =  \kot_{G,\leq b}\setminus \cup_{b'< b}\kot_{G,\leq b'}.
    \]

Although the above definition may seem bizarre from the perspective of classical algebraic geometry, what we have defined is, in fact, quite reasonable. The following result related to the geometry of 
$\kot_G$ is essentially known before. However, we will provide a new proof of these results in \Cref{SS: Newton stratification}.

\begin{theorem}\label{intro:thm: geometry of isoc}
We have
\begin{enumerate}    
\item\label{intro:thm: geometry of isoc-1} $\kot_{G,b}\cong \bB_{\proet} G_b(F)$;
\item\label{intro:thm: geometry of isoc-2}  $i_{\leq b}$ is a (perfectly) finitely presented closed embedding;
\item\label{intro:thm: geometry of isoc-3}  $j_b$ is a (perfectly) finitely presented affine open embedding and $\kot_{G,\leq b}$ is the closure of $\kot_{G,b}$;
\item\label{intro:thm: geometry of isoc-4} $\pi_0(\kot_G)=\pi_1(G)_{\Ga_F}$.
\end{enumerate}
\end{theorem}
Here, we regard the locally profinite group $G_b(F)$ as a group ind-scheme over $k$ (see the beginning of \Cref{SS: rep of loc profinite group} for detailed discussions) and let $\bB_{\proet} G_b(F)$ denote its classifying stack in the pro-\'etale topology.
We note that although we only consider the quotient of $LG/\Ad_\sigma LG$ in the \'etale topology, the pro-\'etale topology appears naturally.
Additionally, we note that when $b\in B(G)$ is basic,  $\kot_{G,b} = \kot_{G,\leq b}$ is closed in $\kot_G$.

In  \Cref{sec:pspl-stacks}, we will carefully develop a theory of $\ell$-adic (co)sheaves on a very general class of geometric objects called prestacks, which includes usual algebraic stacks, as well as $\bB_{\proet} G_b(F)$ and $\kot_G$ as examples.
Thus, for a coefficient ring $\La$ as mentioned above, the (stable $\infty$-)categories of $\ell$-adic sheaves $\shv(\bB_{\proet} G_b(F),\La)$ on $\bB_{\proet} G_b(F)$ and $\shv(\kot_G,\La)$ on $\kot_G$ are well-defined.

However, in this formalism, only $!$-pullback functors are defined for general maps between (pre)stacks. The $*$- and $!$-pushforward functors, as well as the $*$-pullback functors, are only defined for certain classes of maps. The above theorem provides the necessary geometric ingredients to guarantee the existence of all the functors in the following theorem, which will be proved in \Cref{SS: definition of LL category}, \Cref{SS: coh. duality. Kot. stack} and \Cref{SSS: t-structure on shvkot}.

\begin{theorem}\label{intro:thm-LL category}
\begin{enumerate}
    \item\label{intro:thm-LL category-1} For every $b\in B(G)$ choosing a geometric point of $\kot_{G,b}$ induces a natural equivalence 
    \[
    \shv(\kot_{G,b},\La)\cong \Rep(G_b(F),\La).
    \]
        \item The category $\shv(\kot_G,\La)$ is compactly generated, and the subcategory $\shv(\kot_G,\La)^\cpt$ of compact objects consist of those $\mF$ such that $(i_b)^!\mF\in \shv(\kot_{G,b},\La)\cong \rep(G_b(F),\La)$ is a compact object and is zero for almost all $b$'s. 
    \item There are adjoint functors
    \begin{equation}\label{eq: intro-open-closed gluing}
  \xymatrix{
    \shv(\kot_{G,b})\ar@/^/[rr]^{(j_b)_!}\ar@/_/[rr]_{(j_b)_*} && \ar[ll]|{(j_b)^!}\shv(\kot_{G,\leq b}) \ar@/^/[rr]^{(i_{<b})^*}\ar@/_/[rr]_{(i_{<b})^!} && \ar[ll]|{(i_{<b})_*}\shv(\kot_{G,<b}),
    }
    \end{equation}
    inducing a semi-orthogonal decomposition of $\shv(\kot_G,\La)$ in terms of $\bigl\{(i_b)_*(\rep(G_b(F),\La))\bigr\}_b$, as well as in terms of $\bigl\{(i_b)_!(\rep(G_b(F),\La))\bigr\}_b$.
All categories in the diagram are compactly generated and all functors preserve subcategories of compact objects.
    \item There is a canonical self-duality $\shv(\kot_G,\La)$
    \[
    (\verd_{\kot_G}^{\can})^{\cpt}\colon (\shv(\kot_G,\La)^\omega)^{\op}\simeq \shv(\kot_G,\La)^{\omega}
    \]
    such that for every $b\in B(G)$, there are canonical isomorphisms of functors
   \[
(\verd_{\kot_G}^{\can})^\cpt\circ (i_{b})_*\cong (i_{b})_!\circ (\verd^{\can}_{G_b(F)})^\cpt[-2\langle2\rho, \nu_b\rangle](-\langle2\rho, \nu_b\rangle),
\] 
\[
(i_{b})^*\circ (\verd_{\kot_G}^{\can})^\cpt\cong (\verd^{\can}_{G_b(F)})^\cpt\circ (i_{b})^! [-2\langle2\rho, \nu_b\rangle](-\langle2\rho, \nu_b\rangle).
\]
Here $\nu_b$ is the Newton cocharacter associated to $b$, and $(\verd^{\can}_{G_b(F)})^\cpt$ denotes the cohomological duality (or known as the Bernstein-Zelevinsky duality) of the category of smooth representations of $G_b(F)$.
   \item Let $\shv(\kot_G)^{2\rho\mbox{-}p,\leq 0}\subset \shv(\kot_G)$ be the full subcategory generated under small colimits and extensions by objects of the form 
\[
(i_b)_!\cind_K^{G_b(F)}\La[n-\langle 2\rho,\nu_b\rangle], \quad b\in B(G),\ n\geq 0, \ K\subset G_b(F) \mbox{ prop-}p \mbox{ open compact}.
\]
Then $\shv(\kot_G)^{2\rho\mbox{-}p,\leq 0}$ form a connective part of an admissible $t$-structure on $\shv(\kot_G)$. The coconnective part can be described as
\[
\shv(\kot_G)^{2\rho\mbox{-}p,\geq 0}=\bigl\{\mF\in\shv(\kot_G)\mid (i_b)^!\mF\in\rep(G_b(F))^{\geq \langle \chi,\nu_b\rangle}\bigr\}.
\]

  \end{enumerate} 
\end{theorem}

This theorem provides the construction of the local Langlands category, with the promised properties that it glues various categories $\{\rep(G_b(F),\La)\}_{b\in B(G)}$.
On the other hand, recall that the classical local Langlands correspondence aims to classify irreducible smooth representations  of $p$-adic groups. A natural abelian category containing all irreducible representations is the category of admissible representations. It turns out that (the derived version of) this notion has a purely categorical interpretation and we have the full subcategory 
\[
\shv(\kot_G,\La)^{\adm}\subset \shv(\kot_G,\La)
\]
of admissible objects in $\shv(\kot_G,\La)$. We will introduce and study the notion of admissible objects in dualizable categories in details in \Cref{SS: admissible objects}. But as a first approximation, the notion of admissible objects is dual to the notion of compact objects.
Namely, recall that an object $c$ in a (presentable $\La$-linear stable $\infty$-)category $\bfC$ can be regarded as a $\La$-linear functor $F_c$ from the (stable $\infty$-)category $\Mod_\La$ of $\La$-modules to $\bfC$. The object $c$ is called compact if $F_c$ admits a $\La$-linear right adjoint functor. Dually, we call an object admissible if $F_c$ admits a $\La$-linear left adjoint functor $F_c^L$. One can check that admissible objects in $\bfC=\rep(G_b(F),\overline\bQ_\ell)$ are precisely the (derived) admissible representations of $G_b(F)$.
The following statement, in some sense, is dual to \Cref{intro:thm-LL category} and will be proved in \Cref{SS: coh. duality. Kot. stack} and \Cref{SSS: t-structure on shvkot}.
We will let  $(i_b)_\flat$ denote the right adjoint of $(i_b)^!$ and let $(i_b)^\sharp$ denote the right adjoint of $(i_b)_*$. Thanks to \Cref{intro:thm-LL category}, both $(i_b)_\flat$ and $(i_b)^\sharp$ are $\La$-linear continuous functors and, by general nonsense, preserve the subcategory of admissible objects.

\begin{theorem}\label{intro:thm-LL category-adm}
\begin{enumerate}
\item An object $\mF\in\shv(\kot_G,\La)$ is admissible if and only if $(i_b)^!\mF\in \rep(G_b(F),\La)$ is admissible for every $b\in B(G)$, if and only if $(i_b)^\sharp\mF\in \rep(G_b(F),\La)$ is admissible for every $b\in B(G)$. 
\item The canonical duality $(\verd_{\kot_G}^{\can})^{\cpt}$ in \Cref{intro:thm-LL category} induces a duality 
\[
(\verd_{\kot_G}^{\can})^{\adm}\colon (\shv(\kot_G,\La)^{\adm})^{\op}\simeq \shv(\kot_G,\La)^{\adm}
\]
such that for every $b\in B(G)$, we have
\[
(\verd_{\kot_G}^\can)^\adm\circ (i_b)_* \cong (i_b)_\flat \circ (\verd^\can_{G_b(F)})^{\adm}[-2\langle 2\rho,\nu_b\rangle ](-\langle 2\rho,\nu_b\rangle),
\]
\[
 (i_b)\rstar\circ (\verd_{\kot_G}^\can)^\adm[2\langle 2\rho,\nu_b\rangle ](\langle 2\rho,\nu_b\rangle) \cong (\verd^\can_{G_b(F)})^{\adm}\circ (i_b)^!.
\]
Here $(\verd^\can_{G_b(F)})^{\adm}$ is the (derived version of the) usual smooth duality for the category of admissible representations.

\item The following pair of subcategories of $\shv(\kot_G,\La)$
   \[
    \shv(\kot_G,\La)^{2\rho\mbox{-}e,\leq 0}=\bigr\{\mF\in \shv(\kot_G,\La)\mid (i_b)^!\mF\in\rep(G_b(F))^{\leq \langle 2\rho,\nu_b\rangle} \mbox{ for all } b\in B(G) \bigl\}
    \]
   \[
    \shv(\kot_G,\La)^{2\rho\mbox{-}e,\geq 0}=\bigr\{\mF\in \shv(\kot_G,\La)\mid (i_b)^\sharp\mF\in\rep(G_b(F))^{\geq \langle 2\rho,\nu_b\rangle} \mbox{ for all } b\in B(G) \bigl\}
    \]
define an accessible $t$-structure on $\shv(\kot_G,\La)$, which further restricts to a $t$-structure on $\shv(\kot_G,\La)^\adm$. When $\La$ is a field, the abelian category 
\[
\shv(\kot_G,\La)^{2\rho\mbox{-}e,\heartsuit}\cap \shv(\kot_G,\La)^{\adm}
\] 
is stable under the duality $(\verd^\can_{\kot_G})^{\adm}$.
   \end{enumerate}
\end{theorem}

With $\shv(\kot_G)$ defined and its basic properties discussed, we can thus formulate the categorical arithmetic local Langlands correspondence (when $\La=\overline\bQ_\ell$) as a canonical equivalence
\[
\bL_G: \shv(\kot_G,\overline\bQ_\ell)^\cpt\cong \Coh(\locsys_{{}^cG,F}\otimes\overline\bQ_\ell),
\]
which should satisfy a set of compatibility conditions. We shall not discuss these compatibility conditions in the introduction. 

The precise formulation of the conjecture for more general coefficients $\La$ (e.g. $\overline\bF_\ell$) is more subtle. In general, we only expect a natural fully faithful embedding
\[
\shv(\kot_G,\La)^\cpt\hookrightarrow \Coh(\locsys_{{}^cG,F}\otimes \La).
\]
This can be easily seen even when $G=\bG_m$.
To obtain an equivalence, one needs either to replace $\Coh(\locsys_{{}^cG,F}\otimes\La)$ with a smaller subcategory or to enlarge $\shv(\kot_G,\La)^\cpt$. In \cite[Conjecture 4.6.4]{zhu2020coherent}, we explained the first formulation. See also \cite{Fargues.Scholze.geometrization} for the corresponding formulation in their set-up. The second formulation was also indicated in \cite[Remark 4.6.7]{zhu2020coherent}. Here we discuss this second formulation, as it seems to be more convenient for arithmetic applications (as in \Cref{SS: coh of shimura varieties} and also in \cite{yangIhara, yangzhutorsion}). 

For this purpose, we need to introduce a variant of $\shv(\kot_G,\La)$. For each $b\in B(G)$, let
\[
\fgshv(\bB_{\proet}G_b(F),\La)\subset \shv(\bB_{\proet}G_b(F),\La)\cong \rep(G_b(F),\La)
\]
be the smallest full stable subcategory generated by objects $\cind_{K}^{G_b(F)}\La$ under finite colimits and retracts, where $K\subset G_b(F)$ is an open compact subgroup. Let 
\[
\fgshv(\kot_G,\La)\subset \shv(\kot_G,\La)
\] 
be the smallest full stable subcategory generated by objects $(i_b)_*\pi$ under finite colimits and retracts, where $b\in B(G)$ and $\pi\in\fgshv(\bB_{\proet}G_b(F),\La)\}$. If $\La=\overline\bQ_\ell$, then every $(i_b)_*\cind_{K}^{G_b(F)}\La$ is compact and therefore we have $\fgshv(\kot_G,\La)=\shv(\kot_G,\La)^{\cpt}$. However, in general, we only have $\shv(\kot_G,\La)^\cpt\subset \fgshv(\kot_G,\La)$. We can then formulate the categorical local Langlands correspondence (now for general coefficients $\La$) as a canonical equivalence\footnote{When $\La=\overline\bF_\ell$ and $\ell$ is very small (e.g. $\ell$ is not good for $\hat{G}$), we do not have much evidence of the conjecture and the statement might need modifications.}
\[
\bL_G: \fgshv(\kot_G,\La)\cong \Coh(\locsys_{{}^cG,F}\otimes\La),
\]
which again should satisfy a set of compatibility conditions.

Before moving to the next topic, let us make some comments regarding the category $\fgshv(\kot_G,\La)$. First, the actually definition of $\fgshv(\kot_G,\La)$ given in the main context is different. In fact, in \Cref{sec:pspl-stacks}, we will construct another sheaf theory $\fgshv$ for a very general class of stacks $X$ including $\bB G_b(F)$ and $\kot_G$, which can be thought as a theory of constructible sheaves on these geometric objects. Indeed, there is always a functor $\fgshv(X,\La)\to \shv(X)$, which identifies $\fgshv(X,\La)$ with the subcategory of constructible sheaves for familiar geometric objects such as quasi-compact schemes or algebraic stacks. However, the functor $\fgshv(X,\La)\to \shv(X)$ may not be fully faithful in general. It is a non-trivial fact, which will be proved in \Cref{SS: f.g. representations for p-adic groups} and \Cref{SSS: fgshv for isoc}, that for $X=\bB G_b(F)$ and $\kot_G$ the corresponding functors are indeed fully faithful, and the essential images can be described explicitly as above. We shall also mention that various results statement in \Cref{intro:thm-LL category} have counterparts for the theory $\fgshv$, as will be discussed in \Cref{SSS: fgshv for isoc}.
 
This concludes our general discussion of the local Langlands category $\shv(\kot_G,\La)$ and its variants, and the formulation of the categorical local Langlands conjecture.
Next we turn to certain subcategories of both sides, for which we can establish the desired equivalence.

Recall that the stack $\locsys_{{}^cG,F}$, which classifies continuous representations of the Weil group $\varphi: W_F\to {}^cG$, breaks into connected components according to the ``ramification" of $\varphi$. In particular, when $G$ is tamely ramified, there is a well-defined open and closed substack 
\[
\locsys_{{}^cG,F}^\tame\subset \locsys_{{}^cG,F}
\] 
classifying those parameters $\varphi$ that factor through $W_F\to W_F/P_F\to {}^cG$, where $P_F\subset W_F$ denotes the wild inertia. If $G$ additionally splits over an unramified extension, there is also the substack 
\[
\locsys_{{}^cG,F}^{\widehat\unip}\subset  \locsys_{{}^cG,F}^\tame
\] 
of unipotent Langlands parameters, roughly speaking classifying those $\varphi: W_F/P_F\to {}^cG$ sending a generator of the tame inertia to a unipotent element.\footnote{There are actually different versions of the stack of unipotent Langlands parameters. We refer to \Cref{rem: three versions of unipotent stack} for such subtleties.} 
When $\La$ is a field, then  $\locsys_{{}^cG,F}^{\widehat\unip}\otimes\La$ is a connected component of $ \locsys_{{}^cG,F}^\tame\otimes\La$. On the Galois side, we thus have the corresponding subcategories
\[
\Coh(\locsys_{{}^cG,F}^{\widehat\unip}\otimes\La)\subset\Coh(\locsys_{{}^cG,F}^\tame\otimes\La)\subset \Coh(\locsys_{{}^cG,F}\otimes\La).
\]

On the representation theoretic side,
recall that there is a notion of ``depth" for representations of $p$-adic groups. In particular, when $G$ splits over a tamely ramified extension, there is a decomposition
\[
\rep(G(F),\La)=\rep^{\tame}(G(F),\La) \oplus \rep^{>0}(G(F),\La),
\]
where $\rep^{\tame}(G(F),\La)$ denotes the subcategory of depth zero representations and $ \rep^{>0}(G(F),\La)$ denotes the subcategory of representation of $G(F)$ of positive depths.
We let 
\[
\shv^{\tame}(\kot_G,\La)\subset \shv(\kot_G,\La),\quad (\mbox{resp.} \ \shv^{>0}(\kot_G,\La)\subset \shv(\kot_G,\La))
\] 
be the full subcategory consisting of those $\mF$ such that $(i_b)^!\mF\in \rep^{\tame}(G_b(F),\La)$ (resp. $(i_b)^!\mF\in \rep^{>0}(G_b(F),\La)$) for every $b\in B(G)$. For $?$ being $<$ or $\leq$, we denote 
\[
\shv^{\tame}(\kot_{G,? b},\La)=\shv^\tame(\kot_{G},\La)\cap \shv(\kot_{G,?b},\La).
\] 

\begin{theorem}\label{intro:thm-tame-LL category}
\begin{enumerate}
\item\label{intro:thm-tame-LL category-1} The category $\shv^{\tame}(\kot_G,\La)$ is compactly generated by compact objects of the form $(i_b)_*\pi$ with $b\in B(G)$ and $\pi\in \rep^{\tame}(G_b(F),\La)^\cpt$. The pair 
\[
(\shv^{\tame}(\kot_G,\La), \shv^{>0}(\kot_G,\La))
\] 
form a semi-orthogonal decomposition of $\shv(\kot_G,\La)$. Let $\proj^{\tame}$ denote the right adjoint of the inclusion $\shv^{\tame}(\kot_G,\La)\subset \shv(\kot_G,\La)$.
\item Diagram \eqref{eq: intro-open-closed gluing} restricts to a diagram with ``$\tame$" added everywhere, which also induces corresponding semi-orthogonal decompositions of $\shv^{\tame}(\kot_G,\La)$. 
\item The canonical self-duality $(\verd_{\kot_G}^{\can})^\cpt$ restricts to a self-duality $(\verd_{\kot_G}^{\tame,\can})^\cpt$ of $\shv^{\tame}(\kot_G,\La)^\cpt$.
\item The category $\shv^{\tame}(\kot_G,\La)\cap \shv(\kot_G,\La)^\adm$ coincides with the category $\shv^{\tame}(\kot_G,\La)^{\adm}$ of admissible objects of $\shv^{\tame}(\kot_G,\La)$.
\item\label{intro:thm-tame-LL category-adm-dual} The duality $(\verd_{\kot_G}^{\tame,\can})^\cpt$ induces a duality $(\verd_{\kot_G}^{\tame,\can})^\adm$ of $\shv^{\tame}(\kot_G,\La)^{\adm}$. In addition we have 
\[
(\verd_{\kot_G}^{\tame,\can})^\adm= \proj^{\tame}\circ (\verd_{\kot_G}^{\can})^\adm.
\] 
\item\label{intro:thm-tame-LL category-t-structure} The following pair of subcategories of $\shv^{\tame}(\kot_G,\La)$
\begin{align*}
    \shv^{\tame}(\kot_G,\La)^{2\rho\mbox{-}e,\leq 0} &=\bigr\{\mF\in \shv^{\tame}(\kot_G,\La)\mid (i_b)^!\mF\in\rep^{\tame}(G_b(F))^{\leq \langle 2\rho,\nu_b\rangle} \mbox{ for all } b\in B(G) \bigl\}\\
    \shv^{\tame}(\kot_G,\La)^{2\rho\mbox{-}e,\geq 0}&=\bigr\{\mF\in \shv^{\tame}(\kot_G,\La)\mid \proj^\tame((i_b)^\sharp\mF)\in\rep^{\tame}(G_b(F))^{\geq \langle 2\rho,\nu_b\rangle} \mbox{ for all } b\in B(G) \bigl\}
    \end{align*}
define an accessible $t$-structure on $\shv^{\tame}(\kot_G,\La)$, which restricts to a $t$-structure on $\shv^{\tame}(\kot_G,\La)^\adm$. When $\La$ is a field, the abelian category 
\[
\shv^{\tame}(\kot_G,\La)^{2\rho\mbox{-}e,\heartsuit}\cap \shv^{\tame}(\kot_G,\La)^{\adm}
\] 
is stable under the duality $(\verd^{\tame,\can}_{\kot_G})^{\adm}$.
\end{enumerate}
\end{theorem}

\begin{remark}
We expect that Part \eqref{intro:thm-tame-LL category-1} of the above theorem can be strengthened. Namely, the pair $(\shv^{\tame}(\kot_G,\La), \shv^{>0}(\kot_G,\La))$ should form an orthogonal decomposition of $\shv(\kot_G,\La)$. If this is the case, then the further projection $\proj^{\tame}$ in Parts \eqref{intro:thm-tame-LL category-adm-dual} and \eqref{intro:thm-tame-LL category-t-structure} are not necessary.
\end{remark}

There is also a notion of unipotent representations of $p$-adic groups. When $\La=\overline\bQ_\ell$, this was defined by Lusztig in \cite{Lusztig.classification.unipotent.IMRN}. For general coefficients $\La$, see \Cref{SS: 1st approach tame Langlands category}. When $\La$ is a field, unipotent representations also form a subcategory $\rep^{\widehat\unip}(G(F),\La)$, which is in fact a direct summand of $\rep^\tame(G(F),\La)$. Then one can similarly define $\shv^{\widehat\unip}(\kot_G,\La)$. \Cref{intro:thm-tame-LL category} has an analogue in the unipotent case.

\subsubsection{Tame and unipotent categorical Langlands correspondence}
Having the category $\shv(\kot_G,\La)$ and its tame and unipotent parts precisely defined, let us state the one of the main results of the article, which verifies the tame part of the categorical arithmetic local Langlands conjecture under some mild assumptions on the reductive groups.

We assume that $G$ is an unramified reductive group over $F$, equipped with a pinning $(B,T,e)$ defined over $\mO_F$. Such data determine a standard hyperspecial integral model $\underline G$ and an Iwahori integral model $\mI$ of $G$ over $\mO_F$. Let $\iw=L^+\mI$ be the positive loop group of $\mI$, and let $\iw^u\subset \iw$ be the pro-unipotent radical of $\iw$. We let $\sph=\underline G(\mO_F)\subset G(F)$ be the corresponding hyperspecial subgroup,
 let $I=\mI(\mO_F)=\iw(k_F)\subset G(F)$ be the corresponding Iwahori subgroup of $G(F)$, and let $I^u\subset I$ be the pro-$p$-radical of $I$.
For an open compact subgroup $Q\subset G(F)$, we let
\[
\delta_Q:=\cind_{Q}^{G(F)}\La
\]
denote the compact induction of the trivial representation of $Q$. We let $1$ denote the element in $B(G)$ given by $1\in G(F)$.

In the sequel, we will fix a non-trivial additive character $\psi: k_F\to \La^\times$.  We let 
\[
\IW=\cind_{I^u}^{G(F)}\psi_e
\] 
be the compact induction of the character $\psi_e: I^u\to U(k_F)\xrightarrow{e} k_F\xrightarrow{\psi}\La^\times$. This $G(F)$-representation is sometimes called the Iwahori-Whittaker module.

Here is the tame part of the categorical local Langlands correspondence.

\begin{theorem}\label{intro: thm-langlands-categorical-tame}
Let $G$ be a connected unramified reductive group equipped with a pinning (defined over $\mO_F$). Suppose $\La=\overline\bQ_\ell$. 
\begin{enumerate}
\item Then there is a canonical equivalence of categories
\[
\mathbb{L}^{\tame}_{G}\colon  \shv^{\tame}(\kot_G,\La)^\cpt\cong \Coh(\locsys_{^{c}G}^{\tame}\otimes\La),
\]
which restricts to an equivalence
\[
\mathbb{L}^{\widehat\unip}_{G}\colon  \shv^{\widehat\unip}(\kot_G,\La)^\cpt\cong \Coh(\locsys_{^{c}G}^{\widehat\unip}\otimes\La).
\]
\item The equivalence intertwines the canonical duality of $\shv^{\tame}(\kot_G,\La)^\cpt$ (as in \Cref{intro:thm-tame-LL category}) and the twisted Grothendieck-Serre duality of $\Coh(\locsys_{{}^cG}^{\tame}\otimes\La)$.
\item The equivalence is compatible with the natural $\pi_1(G)_{\Ga_F}\cong \xch(Z_{\hat{G}}^{\Ga_F})$-gradings on both sides.
\item We have the following matching of objects under the above equivalences
\begin{align*}
    &\mathbb{L}^{\tame}_{G}((i_{1})_{*}\delta_{I^u}) \cong \cohspr_{^{c}G,F},                                      &\mathbb{L}^{\tame}_{G}((i_{1})_{*}\delta_{I}) \cong \cohspr_{^{c}G,F}^{\unip},\\
    &\mathbb{L}^{\tame}_G((i_{1})_{*}\delta_{\sph})\cong \mO_{\locsys_{{}^cG,F}^{\mathrm{unr}}},    & \\
    &\mathbb{L}^{\tame}_{G}((i_{1})_{*}\IW)\cong \mathcal{O}_{\locsys_{^{c}G,F}}^{\tame},                & \mathbb{L}^{\tame}_{G}((i_{1})_{*}\IW^{\unip})\cong \mathcal{O}_{\locsys_{^{c}G,F}}^{\widehat\unip}.  
\end{align*}
\end{enumerate}
\end{theorem}

We briefly explain some notations and terminology in the theorem. By the twisted Grothedieck-Serre duality, we mean the composition of the usual Grothendieck-Serre duality with an automorphism of $\locsys_{{}^cG}^{\tame}$ induced by the Cartan involution of $\hat{G}$ (see \eqref{eq:Cartan involution of loc}). The $\pi_1(G)_{\Ga_F}$-grading of $\shv(\kot_G,\La)$ is induced by the decomposition of $\kot_G$ into connected components (see \Cref{intro:thm: geometry of isoc} \eqref{intro:thm: geometry of isoc-4}), and the $\xch(Z_{\hat{G}}^{\Ga_F})$-grading of $\Coh(\locsys_{^{c}G}^{\tame}\otimes\La)$ is induced from a $Z_{\hat{G}}^{\Ga_F}$-gerbe structure on $\locsys_{{}^cG,F}$. 
The stack $\locsys_{{}^cG,F}^{\mathrm{unr}}\subset \locsys_{{}^cG,F}^{\tame}$ classifies unramified Langlands parameters.
The coherent sheaf $\cohspr_{^{c}G,F}$ (resp. $\cohspr_{^{c}G,F}^{\unip}$) is called the tame coherent Springer sheaf (resp. unipotent coherent Springer sheaf), which is defined as the $*$-pushforward of the dualizing sheaf of $\locsys_{{}^cB, F}^{\tame}$ (resp. $\locsys_{{}^cB, F}^{\unip}$) to $\locsys_{{}^cG, F}^{\tame}$. Here $\locsys_{{}^cB, F}^{\tame}$ (resp. $\locsys_{{}^cB, F}^{\unip}$) classifies ${}^cB$-valued continuous $\ell$-adic representations of the tame Weil group, where  ${}^cB\subset {}^cG$ is the Borel subgroup of ${}^cG$. See \cite[\textsection{4.4}]{zhu2020coherent} and \Cref{ex: spectral DL induction w equal 1}.

When $\La=\overline\bF_\ell$, we can only prove a weaker version, which is sufficient for some arithmetic applications. First, there is certain subcategory
\[
\fgshv^{\unip}(\kot_G,\La)\subset \fgshv(\kot_G,\La)\cap \shv^{\widehat\unip}(\kot_G,\La).
\]
It contains $(i_1)_*\delta_I$. Under some mild assumption on the characteristic $\ell$ (which will be satisfied in the following theorem), it also contains $(i_1)_*\delta_P$ for every parahoric subgroup $P$ of $G(F)$.

\begin{theorem}\label{intro: thm-langlands-categorical-unip}
Suppose $\La=\overline\bF_\ell$ with $\ell$ bigger than the Coxeter number of any simple factor of $G$, and $\ell\neq 19$ (resp. $\ell\neq 31$) if $G$ has a simple factor of type $E_7$ (resp. $E_8$). Then there is a fully faithful embedding
\[
\mathbb{L}_{G}^{\unip}\colon  \rshv^{\unip}(\kot_G,\La)\hookrightarrow \indcoh(\locsys_{^{c}G,F}^{\widehat\unip}\otimes \La),
\]
with the essential image stable under the action of $\ind\Perf(\locsys_{^{c}G,F}^{\widehat\unip}\otimes \La)$. We have
\[
\mathbb{L}^{\unip}_{G}((i_{1})_{*}\delta_{I}) \cong \cohspr_{^{c}G,F}^{\unip}.
\] 
If $Z_G$ is connected,  then essential image contains the category $\ind\Perf(\locsys_{^{c}G,F}^{\widehat\unip}\otimes \La)$. 
\end{theorem} 

\begin{remark}\label{intro: rmk-langlands-categorical-unip-improvement}
We mention that the restrictions of the characteristic are largely due to the current restriction of the characteristic in the modular local geometric Langlands as established in \cite{Bez.Riche.modular}. We expect the theorem holds under a much milder restriction of the characteristic. We also expect that the functor $\mathbb{L}_{G}^{\unip}$ will send $(i_{1})_{*}\IW^{\unip}$ to $\mathcal{O}_{\locsys_{^{c}G,F}^{\widehat\unip}}$. This again would follow if certain result in the modular local geometric Langlands is established.
\end{remark}
As a corollary, we obtain the following result. The functor $\End$ below is the derived endomorphism.
\begin{corollary}\label{intro: thm: end of coh spr}
There are natural isomorphisms 
\begin{enumerate}
\item For $\La=\overline\bQ_\ell$ or $\overline\bF_\ell$ (with $\ell$ satisfying condition as in \Cref{intro: thm-langlands-categorical-unip}), we have 
\[
\End_{\locsys_{{}^c G,F}\otimes \La}(\cohspr^{\unip}_{{}^c G,F}\otimes\La)\cong C_c(I\backslash G(F)/I,\La),
\] 
where $C_c(I\backslash G(F)/I,\La)$ is the \emph{derived} Iwahori-Hecke algebra (which is non-derived if $\La=\overline\bQ_\ell$ or $\La=\overline\bF_\ell$ if $\ell$ is banal).
\item Suppose $\La=\overline\bQ_\ell$. Then we have 
\begin{align*}
&\rg(\locsys^{\tame}_{{}^c G,F}, \mO)\cong C_c((I^u,\psi)\backslash G(F)/(I^u,\psi)),  \\
&\rg(\locsys^{\tame}_{{}^c G,F}, \cohspr^{\tame}_{{}^c G,F})\cong C_c((I^u,\psi)\backslash G(F)/I^u),\\
&\rg(\locsys^{\widehat\unip}_{{}^c G,F}, \cohspr^{\unip}_{{}^c G,F})\cong C_c((I^u,\psi)\backslash G(F)/I).
\end{align*}
\end{enumerate}
\end{corollary}
Again we expect the last isomorphism still holds when $\La=\overline\bF_\ell$, by virtue of \Cref{intro: rmk-langlands-categorical-unip-improvement}.

We can also prove the following result.
\begin{theorem}\label{thm: compact induction in the heart}
Suppose $\La=\overline\bQ_\ell$. Then for every basic element $b\in B(G)$, and for every pair $(P,\varrho)$, where $P\subset G_b(F)$ is a parahoric subgroup and $\varrho$ is a finite dimensional representation of $P$ obtained by inflation of a representation of the Levi quotient $L_P$ of $P$, the object
\[
\mathbb{L}^{\tame}_{G}((i_b)_*\cind_P^G\varrho)
\] 
is in the abelian category $\Coh(\locsys^{\tame}_{{}^cG,F})^{\heartsuit}$, and is a maximal Cohen-Macaulay coherent sheaf.
\end{theorem}
\begin{remark}\label{rem: compact induction in the heart}
For $\La=\overline\bF_\ell$, we do not expect the same statement holds for arbitrary $\varrho$. However, we expect it remains to hold if $\varrho$ is a projective object in $\rep(L_P,\La)^{\heartsuit}$. In fact, given \Cref{intro: DL for tilting} below, this will be the case if the last expectation of \Cref{intro: rmk-langlands-categorical-unip-improvement} holds.
\end{remark}

\subsubsection{Some applications to the classical Langlands program}
Now we discuss the relation between the categorical local Langlands correspondence and the classical local Langlands correspondence.
We assume that $\La=\overline\bQ_\ell$.

As the category $\shv^{\tame}(\kot_G)$ is equivalent to $\indcoh(\locsys_{{}^cG,F}^{\tame})$, every object in $\shv^{\tame}(\kot_G)$ is acted by the tame spectral Bernstein center 
\[
Z^{\tame}_{{}^cG,F}:=H^0\rg(\locsys_{{}^cG,F}^{\tame},\mO).
\] 
In particular, if $\pi$ is a depth zero irreducible representation of $G_b(F)$ for some basic $b$ (or more generally $\pi$ is a representation of $G_b(F)$ such that $H^0\End(\pi)$ is a local artinian $\La$-algebra), then  $Z^{\tame}_{{}^cG,F}$ acts on $(i_b)_*\pi$ through a local artinian quotient, which determines a unique maximal ideal of  $Z^{\tame}_{{}^cG,F}$. Since closed points of $\Spec Z^{\mathrm{spec},\tame}_{{}^cG,F}$ are in bijection to continuous semisimple representations $W_F$ up to $\hat{G}$-conjugacy, we obtain the following.

\begin{theorem}\label{intro: thm ss parameter}
One can attach  to every irreducible depth zero representation $\pi\in \rep^{\tame}(G_b(F))$  a tame semisimple Langlands parameter $\varphi_\pi^{ss}$. 
\end{theorem}

\begin{remark}\label{intro: rem FS parameter}
When $F$ is an equal characteristic local field, Genestier-Lafforgue's paramterization attaches to every (not necessarily depth zero) irreducible representation $\pi$ a semisimple Langlands parameter $\varphi_\pi^{ss}$. It is not difficult to show that  our parameterization given in the above theorem is the restriction of Genestier-Lafforgue's to depth zero representations. We will discuss this in another place.
On the other hand, for $F$ being a general local field, Fargues-Scholze also associate to every (not necessarily depth zero) irreducible representation $\pi$ a semisimple Langlands parameter $\varphi_\pi^{ss}$. It is known that Fargues-Scholze's  and Genestier-Lafforgue's parameterizations coincide when $F$ is of equal characteristic by \cite{Li-Huerta}. We expect that, when $F$ is a $p$-adic field, our parameterization will also be the restriction of Fargues-Scholze's  to the depth zero representations. 
\end{remark}

To lift semisimple Langlands parameters attached to $\pi$ to a true parameter $\varphi_\pi$ is more subtle, even with the categorical equivalence at hand. Here we only discuss such liftings for supercuspidal representations.

We assume that  $\La=\overline{\bQ}_\ell$. For simplicity,  we assume that $G$ is semisimple in the introduction. (We allow general $G$ in the main body of article.)  
Recall a parameter $\varphi: W_F\to {}^LG$ is called discrete if $C_{\hat{G}}(\varphi)$ is finite. This is equivalent to saying that $\{\varphi\}/C_{\hat{G}}(\varphi)$ is an open point of $\Loc_{{}^cG,F}\otimes \overline{\bQ}_\ell$. 
One can show that its closure, denoted by $\overline{\{\varphi\}}$ for simplicity, is a smooth irreducible component of $\Loc_{{}^cG,F}\otimes\overline{\bQ}_\ell$. In fact, it is always the quotient of a prehomogeneous space by a reductive group.  See \Cref{prop: geometry of discrete component}.

\begin{theorem}\label{intro: parameter of ss}
Let $\pi$ be a depth zero supercuspidal irreducible representation of $G_b$, for $b$ basic. Then $\bL_G((i_b)_*\pi)$ is a vector bundle on $\overline{\{\varphi_\pi\}}$, for some discrete tame parameter $\varphi_\pi$. 
If $\pi$ is generic (with respect to our choice of Whittaker datum), then such vector bundle is just the structure sheaf of $\overline{\{\varphi_\pi\}}$. Consequently, the semisimple parameter $\varphi_\pi^{ss}$ attached to $\pi$ as from \Cref{intro: thm ss parameter} can be lifted to an enhanced Langlands parameter $(\varphi_\pi, r_\pi)$ attached to $\pi$, consisting of a discrete Langlands parameter $\varphi_\pi: W_F\to {}^cG(\La)$ whose semisimplification is $\varphi_\pi^{ss}$ and a finite dimensional representation $r_\pi$ of $C_{\hat{G}}(\varphi_\pi)$.  If $\pi$ admits a Whittaker model (with respect to our choice of Whittaker datum), then $r$ is the trivial representation of $C_{\hat{G}}(\varphi_\pi)$.
\end{theorem}

The above assignment
\[
\pi\leadsto (\varphi_\pi, r_\pi)
\]
is a candidate of the Langlands parameterization of depth zero supercuspidal representations. 
To the best of our knowledge, this is the first construction of the Langlands parameterization for \emph{all} depth zero supercuspidal representations; previously, only specific cases had been associated with enhanced Langlands parameters.
In these instances, it would be intriguing to compare our parameterization with those found in the existing literature. In \Cref{SSS: regular supercuspidal}, we study this question in the simplest case. Namely, we will demonstrate that when 
$\pi$ is as in the work of DeBecker-Reeder \cite{Debacker.Reeder} and Kazhdan-Varshavsky, $\varphi_\pi$ coincides with the attached local Langlands parameter by \emph{loc. cit.}
On the other hand, we expect that in the case when $\pi$ is a unipotent supercuspidal representation of $G_b$, $\varphi_\pi$ coincides with the local Langlands parameter attached to $\pi$ by Lusztig \cite{Lusztig.classification.unipotent.IMRN} and Morris \cite{Morris.tame.supercuspidal}. We hope to address this question in another occasion.

Let us also mention that it is well-known that given a semisimple parameter $h: W_F\to {}^L\!G$, there is at most one discrete parameter $\varphi: W_F\to {}^L\!G$ such that $h=\varphi^{ss}$ (up to $\hat{G}$-conjugation). In other words, for $\pi$ being supercuspidal, if $\varphi_{\pi}^{ss}$ can be lifted to $\varphi_\pi$, then such lifting is unique. However, to assign the additional representation $r_\pi$ of $C_{\hat{G}}(\varphi_\pi)$ is much more subtle.  We will study properties of this parameterization $\pi\leadsto (\varphi_\pi,r_\pi)$ in another place.

In another direction, we can attach  an admissible representation of the $p$-adic group to certain Langlands parameters. Naively, one may expect the following recipe as indicated before.
Let $(\varphi, r)$ be an enhanced Langlands parameter. I.e. $\varphi: W_F\to {}^cG$ a Langlands parameter, and $r\in C_{\hat{G}}(\varphi)$. Then we may regard $\varphi$ as a stacky point $\{\varphi\}/C_{\hat{G}}(\varphi)$ of $\locsys_{{}^cG,F}$ and $r$ as a vector bundle $\mV_r$ on $\{\varphi\}/C_{\hat{G}}(\varphi)$. Then under the equivalence $\bL_G$, $\bL^{-1}_G(\mV_r)$ should give corresponds to the representation attached to the parameter $(\varphi,r)$.
This idea works in many cases as follows. (But it fails in general.)

\begin{theorem}\label{intro: parameter to representation}
Let $\varphi: W_F\to {}^cG$ be a parameter such that 
\begin{itemize}
\item $H^2(W_F, \Ad^0)=0$, where $\Ad^0$ denotes the adjoint representation of $W_F$ on $\hat\frakg$ via the representation $\varphi$;
\item $C_{\hat{G}}(\varphi)$ is reductive.
\end{itemize}
Let $r$ be an irreducible representation of the $C_{\hat{G}}(\varphi)$. Let  $r_0$ be its restriction to $Z_{\hat{G}}^{\Ga_{\widetilde F/F}}$, which corresponds to an element $\al_r\in \pi_1(G)_{\Ga_F}=\xch(Z_{\hat{G}}^{\Ga_{\widetilde F/F}})$. Let $b\in B(G)$ be the unique basic element which maps to $\al_r$ under the Kottwitz map.
Then 
\[
\bL_G^{-1}(\mV_{\varphi,r})=:\mF_{(\varphi,r)}\in \shv^{\tame}(\kot_G)
\]
is an admissible, supported on the connected component of $\kot_G$ corresponding to $\al_r$ (see \Cref{intro:thm: geometry of isoc}), and is in the heart of the $t$-structure of $\Shv^{\tame}(\kot_G)^{\adm}$ as constructed in \Cref{intro:thm-tame-LL category}.
In particular, the $!$-fiber of $\mF_{(\varphi,r)}$ at $b$ is an admissible representation of $G_b(F)$.
\end{theorem}

The assignment 
\[
(\varphi, r)\leadsto (i_b)^!\mF_{(\varphi,r)}\in \rep(G_b(F))^{\adm}\cap \rep(G_b(F))^{\heartsuit}
\] 
thus can be regarded as a candidate of the construction of the $L$-packets for certain depth zero Langlands parameters. 
Unfortunately, currently we can say very little about $(i_b)^!\mF_{(\varphi,r)}$. We do not even know when it is non-zero, and if it is non-zero, when it is irreducible. The only exception is that when $r=\mathbf{1}$ is the trivial representation, then we know that $(i_1)^!\mF_{\varphi, \mathbf{1}}\neq 0$, and admits a Whittaker model. We shall also mention that if the parameter $\varphi$ is not a smooth point in $\locsys_{{}^cG,F}^{\tame}$, the above result needs to be modified. 

\subsubsection{Cohomology of Shimura varieties via coherent sheaves}
On of the main motivations of the categorical local Langlands correspondence is to understand the cohomology of Shimura varieties via the local-global compatibility. See \cite[\textsection{4.7}]{zhu2020coherent} for some discussions and speculations. 
We state a result in this direction.
Let $(G,X)$ be a Shimura datum of Hodge type. Let $p$ be a prime such that $G_{\bQ_p}$ is unramified.
Let $K\subset G(\bA_f)$ be an open compact subgroup written as $K=K_pK^p$ where $K_p=I\subset G(\bQ_p)$ is an Iwahori subgroup and $K^p\subset G(\bA_f^p)$ is a prime-to-$p$ level. Let $d=\dim \bfSh_K(G,X)$. 
Let $\bfSh_K(G,X)$ be the corresponding Shimura variety defined over the reflex field $E=E(G,X)\subset \bC$. We shall fix an embedding $\iota: E\subset \overline\bQ_p$, determining a $p$-adic place $v$ of $E$ over $p$. Let $E_v$ be the completion of $E$.

Let $\La$ be either $\overline\bF_\ell$ or $\overline\bQ_\ell$. We will be interested in the \'etale cohomology $C(\bfSh_K(G,X)_{\overline\bQ_p}, \La[d])$ of the Shimura variety $\bfSh_K(G,X)$ base changed to $\overline\bQ_p$, equipped with an action of the Hecke algebra $H_K:=H^0C_c(K\bs G(\bA_f)/K,\La)$, as well as the action of the Galois group $\Ga_{E_v}=\mathrm{Gal}(\overline\bQ_p/E_v)$. We shall write $H_K=H_I\otimes_\La H_{K^p}$, where $H_I$ is the Iwahori-Hecke algebra and $H_{K^p}$ is the prime-to-$p$ Hecke algebra.

The Shimura datum gives a conjugacy class of minuscule cocharacters $\{\mu\}$ of $G_{\bQ_p}$ with field of definition $E_v$. Let $V_\mu$ be the associated highest weight irreducible representation of $\hat{G}\otimes E_v$ with coefficient in $\La$.
 As before, we let $\locsys_{{}^cG,\bQ_p}^{\widehat\unip}$ denote the stack of unipotent Langlands parameters and we use the same notation to denote its base change to $\La$. Then $V_\mu$ gives an ``evaluation" bundle $\widetilde{V_\mu}$ on $\locsys_{{}^cG,\bQ_p}^{\widehat\unip}$, equipped with an action of $W_{E_v}$.

We have the following theorem, which is a special case of \Cref{thm: etale coh Shimura variety vai coh sheaf}.
\begin{theorem}\label{intro: coh of Shimura}
Assume that either $\La=\overline\bQ_\ell$ or $\La=\overline\bF_\ell$ with $\ell$  bigger than the Coxeter number of any simple factor of $G$.
There is an object 
\[
\Igss_{K^p}^{\spec,\unip}\in \indcoh(\locsys_{{}^cG,\bQ_p}^{\widehat\unip})^{\adm},
\] 
equipped with an action of $H_{K^p}$, such that there is an $H_{K^p}\times W_{E_v}$-equivariant isomorphism
\[
C(\bfSh_K(G,X)_{\overline\bQ_p}, \La[d])\cong \Hom_{\indcoh(\locsys_{{}^cG,\bQ_p}^{\tame})}(\cohspr^{\unip}_{{}^cG,\bQ_p}\otimes\widetilde{V_\mu}, \Igss_{K^p}^{\spec,\unip}).
\]
Here on the right hand side $H_{K^p}$ acts on $\Igss_{K^p}^{\spec,\unip}$, and $W_{E_v}$ acts $\widetilde{V_\mu}$.
\end{theorem}

We refer to \cite{yangIhara, yangzhutorsion} for some applications of this formula. We also mention that the isomorphism is compatible with the $H_I$-action on both sides, where $H_I$ acts on the right hand side via the action of $\cohspr^{\unip}_{{}^cG,\bQ_p}$ through \Cref{intro: thm: end of coh spr}. This will be proved in \cite{yangzhutorsion}.

\subsection{Ideas of proof and some other results}

Now we briefly discuss the main ideas behind the proof of our results.

\subsubsection{Categorical trace}
As mentioned before, the Deligne-Lusztig theory provides a way to construct representations of finite groups of Lie type from
the category $\shv(B_H\backslash H/B_H)$. The category $\shv(B_H\backslash H/B_H)$ with a natural monoidal structure is usually called the (finite) Hecke category\footnote{There are actually different versions of Hecke categories, see \Cref{SS: Aff Hecke cat} for a discussion.}, and has been extensively studied in geometric representation theory. 
In recent years, it has been realized that the Deligne-Lusztig induction functor can be regarded as a Frobenius-twisted categorical trace construction, and induces an equivalence from the Frobenius-twisted categorical trace of the monoidal category $\shv(B_H\bs H/B_H)$ to (the unipotent part of) the category of representations of $H(\kappa)$.
See 
\cite{Lusztig2015unipotent, Lusztig2017, eteve.monodromic, eteve.DL} for various versions of this ideas.

We will apply similar ideas in the affine setting. Namely, we shall look at the correspondence
\[
\iw\backslash LG/\iw\xleftarrow{\delta} LG/ \Ad_\sigma \iw \xrightarrow{\Nt} \kot_G.
\]
Here $\iw\subset LG$ is an Iwahori subgroup of $LG$, defined over $k_F$. 
The stack $\iw\backslash LG/\iw$ is usually called the Hecke stack and the stack 
\[
\Sht^{\loc}= LG/ \Ad_\sigma \iw
\] 
is sometimes called the stack of local Shtukas.
Then we can construct objects in $\shv(\kot_G)$ via the pull-push of sheaves on $\shv(\iw\bs LG/\iw)$.
The category $\shv(\iw\backslash LG/\iw)$ with a natural monoidal structure is usually called the affine Hecke category. Then we can similarly define the affine Deligne-Lusztig induction, which instead of producing representations of $G(F)$ now produces sheaves on $\kot_G$. 
Similarly, the affine Deligne-Lusztig induction should induce an equivalence from the Frobenius-twisted categorical trace of the monoidal category $\shv(\iw\bs LG/\iw)$ to (the unipotent part of) the category $\shv(\kot_G)$. As explained above, the category $\shv(\kot_G)$ is obtained by gluing categories of representations of various $p$-adic groups related to $G$. Therefore, we produce representations of $p$-adic groups via the affine Deligne-Lusztig induction.

\medskip

\begin{center}
\begin{tabular}{ c | c  }
\hline 
$H$ over $\kappa$                    & $G$ over $F$  \\ \hline 
$\bB H(\kappa)$                        & $\kot_G$ \\ 
$\rep(H(\kappa))$                      & $\shv(\kot_G)$ \\
$H/\Ad_\sigma B_H$     & $\Sht^\loc$ \\ 
$\shv(B_H\backslash H/B_H)$ & $\shv(\iw\backslash LG/\iw)$\\
\end{tabular}
\end{center}

\medskip

Although this idea has been in the air for sometime (e.g. see \cite{gaitsgory2016geometric, Zhu2016} for some informal accounts), to make it really work for representation theory of $p$-adic groups is non-trivial, as we need to work in a highly infinite dimensional set-up and to work with some exotic (from the traditional point of view) geometric object such as $\kot_G$.
In some sense, a considerable portion of the second part of this article is to review and further develop necessary foundational materials to make sure such procedure is valid.

While making the above construction work in the affine setting is challenging, there is a reward.
The affine Hecke category $\fgshv(\iw\backslash LG/\iw)$ admits another realization via the coherent sheaves on certain algebraic stack $S_{{}^cG,\breve F}^{\unip}$ constructed from the Langlands dual group.
This is a celebrated result of Bezrukavnikov see \cite{{bezrukavnikov2016two}}.  (As far as we know, there is no such coherent description of finite Hecke category.)
One can then similarly taking the twisted categorical cocenter of the category of $\Coh(S_{{}^cG,\breve F}^{\unip})$, which can be realized via what we call (in \cite{zhu2020coherent}) the spectral Deligne-Lusztig induction
\[
S_{{}^cG,\breve F}^{\unip}\xleftarrow{\delta^{\unip}} \widetilde{\locsys}_{{}^cG,F}^{\unip}\xrightarrow{\tilde\pi^{\tame}} \locsys_{{}^cG,F}.
\]
Therefore, the category of coherent sheaves on the stack of unipotent Langlands parameters appears naturally.

To summarize, we will deduce \Cref{intro: thm-langlands-categorical-unip} from taking the Frobenius-twisted categorical trace of the tame local geometric Langlands correspondence as proved in \cite{arkhipov2009perverse}, \cite{bezrukavnikov2016two} \cite{Bez.Riche.modular} and \cite{DYYZ2}. We shall, however, emphasize that even with the local geometric Langlands correspondence at hand and with the general formalism of taking categorical traces being developed, there are additional challenges to obtain \Cref{intro: thm-langlands-categorical-unip}. 
We explain these additional difficulties in the unipotent case. 

The general formalism developed in the second  part of this article will imply that there are fully faithful embeddings
\[
\tr(\rshv(\iw\bs LG/\iw,\La),\phi)\hookrightarrow \rshv(\kot_G,\La),
\]
and
\[
 \tr(\indcoh(S_{{}^cG,\breve F}^{\unip}\otimes\La), \phi)\hookrightarrow \indcoh(\locsys_{{}^cG,F}\otimes\La).
\] 
Here $\tr(-,\phi)$ denotes the Frobenius-twisted categorical trace of the corresponding affine Hecke categories in representation theoretic side and in spectral side. 
To obtain \Cref{intro: thm-langlands-categorical-unip}, we need to identify essential images of these functors. 

In the representation theory side, we need to show that $\ind\fgshv^{\unip}(\kot_G,\La)$ is generated by the essential image of the unipotent affine Deligne-Lusztig induction. 
While in the finite-dimensional case this is simply the definition of unipotent representations (of finite group of Lie type), this is not the case in the affine setting.
We deduce the essential surjectivity by analyzing the geometry of the map $\Nt: \Sht^{\loc}\to \kot_G$, making use of some beautiful results of He and Nie-He (\cite{he2014geometric, He.Nie.minimal.length}) regarding the combinatorics of the Iwahori-Weyl group. 

\quash{
The analogous statement in the classical Deligne-Lusztig theory amounts to saying that all representations of $H(\kappa)$ can be constructed by the Deligne-Lusztig induction. This is a quite non-trivial fact, due to Deligne-Lusztig (\cite[Corollary 7.7]{Deligne.Lusztig}) when $\La=\overline\bQ_\ell$ and due to Bonnaf\'e-Rouquier (\cite[\textsection{9}, Theorem A]{BR.modular}) when $\La=\overline\bF_\ell$.  
We will deduce the essential surjectivity in the affine case from the finite case, by analyzing some geometry of the map $\Nt: \Sht^{\loc}\to \kot_G$, making use of some deep results of He and Nie-He on some combinatorics of Iwahori-Weyl group. 
}

In the spectral side, if $\La$ is a field of characteristic zero, then the general theory of singular support of coherent sheaves developed by Arinkin-Gaitsgory in \cite{arinkin2015singular} together with a computation of pull-push singular supports is enough to show that $ \tr(\indcoh(S_{{}^cG,\breve F}^{\unip}\otimes\La), \phi)\to \indcoh(\locsys^{\widehat\unip}_{{}^cG,F}\otimes\La)$ is essential surjective. In fact, such computations have been essentially done by Ben-Zvi-Nadler-Pregyel \cite{ben2017spectral}.
 However, when $\La$ is a field of positive characteristic, the theory of coherent sheaves on the stack over $\La$ is very subtle and many arguments in characteristic zero fail. We must analyze the geometry of the spaces involved in the spectral Deligne-Lusztig induction more carefully.

\subsubsection{Whittaker coefficient}
Next we now discuss the main idea behinds the proof of \Cref{thm: compact induction in the heart}.
We assume that $\La=\overline\bQ_\ell$, although the same strategy should work for $\La=\overline\bF_\ell$ once certain result in the local geometric Langlands correspondence is established.
 
Since $\locsys_{{}^cG,F}=\locsys_{{}^cG, F}^{\Box}/\hat{G}$, it is enough to show that for all finite dimensional representations $V$ of $\hat{G}$, giving the ``evaluation" vector bundle $\widetilde{V}$ on $\locsys_{{}^cG,F}$, we have
\begin{equation*}\label{intro: eq vanishing coherent}
H^i\rg(\locsys_{{}^cG,F}, \widetilde{V}\otimes \mathbb{L}^{\tame}_{G}((i_b)_*\cind_{P}^{G_b}\varrho))=0, \quad \mbox{ for } i\neq 0.
\end{equation*}

Via the equivalence $\bL^{\tame}_G$ we may translate this question back to show that the Whittaker model of the cohomology of certain sheaves on affine Deligne-Lusztig varieties concentrate in middle degree. More precisely, we will show that
\begin{equation}\label{intro: eq vanishing whittaker}
H^i\Hom_{\shv(\kot_G)}(\Ch_{LG,\phi}^{\tame}(\mZ^{\mon}(V)\star^u \widetilde{\mathrm{Til}}_{\dot{w}}^{\mon}), (i_{1})_*\IW_{\psi_1})=0,\quad i>0,
\end{equation}
Here $\widetilde{\mathrm{Til}}_{\dot{w}}^{\mon}$ is a monodromic version of the tilting sheaf on $\iw^u\backslash LG/\iw^u$. 

The above formula can be regarded as a (correct) generalization of a result by Dudas (\cite{dudas.DL.restriction.Gelfand-Graev}) on the Gelfand-Graev model of the compactly supported cohomology $C_c(Y_{\dot{w}},\La)$ of the classical Deligne-Lusztig variety $Y_{\dot{w}}$. 
But Dudas' method does not seem to generalize in the affine setting. 
Note that our argument is applicable even in the classical Deligne-Lusztig setting, giving a simpler proof of  Dudas'  result. See \Cref{prop:multiplicity one GG rep}.
In the process, we also discovered class of projective generators of the category of representations of finite group of Lie type coming from the Deligne-Lusztig\footnote{This class of representations is also discovered by Eteve \cite{eteve.tilting} independently.}.

\begin{theorem}\label{intro: DL for tilting}
Let $H$ be a connected reductive group over a finite field $\kappa$.
For each $u\in W_H$, there is a representation $\widetilde{R}^T_{\dot{u}}$ of $H(\kappa)\times T_H^{u\sigma}$ on a finite projective $\La$-module. When regarded as a representation of $H(\kappa)$, it is a projective object.
In addition, for every representation $\pi$ of $H(\kappa)$, there is some $u\in W_H$ and a non-zero map $\widetilde{R}^T_{\dot{u}}\to \pi$.
\end{theorem}

The representation $\widetilde{R}^T_{\dot{u}}$ in the above theorem arises as the Deligne-Lusztig induction of tilting sheaves  $\Ch_{H,\phi}(\widetilde\Til^{\mon}_{\dot{u}})$. See \Cref{prop: exactness of DL induction on tilting}.

Now using the geometry of $\Nt: \sht^{\loc}\to \kot_G$, one deduces from \eqref{intro: eq vanishing whittaker} that when $b$ is basic,
\[
H^i\rg(\locsys_{{}^cG,F}, \widetilde{V}\otimes \mathbb{L}^{\tame}_{G}((i_b)_*\cind_{P}^{G_b}(\widetilde{R}^T_{\dot{u}})))=0, \quad \mbox{ for } i\neq 0,
\]
where $\widetilde{R}^T_{\dot{u}}$ range over those representations of the Levi quotient $L_{P}$ of $P$ from \Cref{intro: DL for tilting}.

When $\La=\overline\bQ_\ell$, every irreducible irreducible of $L_{P}$ is a direct summand of $\widetilde{R}^T_{\dot{u}}$ for some $w_f$. This gives \Cref{thm: compact induction in the heart}.
We also notice that as mentioned in \Cref{rem: compact induction in the heart}, for $\La=\overline\bF_\ell$ this type of argument should work for projective representation of $L_{P}$. (For $\La=\overline\bF_\ell$, the current missing ingredient to translate \Cref{thm: compact induction in the heart} to the vanishing result of Whittaker model.)

\quash{

A similar ideal lead to the following result.

\begin{theorem}
Suppose $G$ is a split reductive group with simply-connected derived subgroup. Then $C_c(G(F)/G(\mO))$ is a flat $C_c(G(\mO)\backslash G(F)/G(\mO))$-module.  
\end{theorem}
When $G=\GL_n$, this result is proved in \cite{} (originally conjectured in \cite{}). 
In fact, the assumption on $G$ can be weakened. See \Cref{prop: flatness of universal spherical rep} for the precise statement.
Nevertheless to say, our proof is completely different.

Finally we mention that the crucial ingredient to compare our construction of Langlands parameters and other approaches is the so-called $S=T$ type statement.
}

\subsubsection{Supercusipdal representations}
Having \Cref{thm: compact induction in the heart} at hand, we explain how to deduce \Cref{intro: parameter of ss}. For simplicity, we assume that $G$ is simply-connected. Then every supercuspidal representation of $G$ is of the form $\cind_{P}^{G(F)}\varrho$ for $P$ a maximal parahoric subgroup of $G(F)$ and $\varrho$ a cuspidal irreducible representation of the Levi quotient $L_P$ of $P$. Then $\bL_G^{\tame}((i_1)_*\pi)$ is a maximal Cohen-Macaulay sheaf on $\locsys_{{}^cG,F}^{\tame}$, and therefore is supported on the union of several  irreducible components of $\locsys_{{}^cG,F}^{\tame}$. As $\End((i_1)_*\pi)$ is $\La=\overline\bQ_\ell$, we see that the tame spectral Bernstein center $Z^{\mathrm{spec},\tame}_{{}^cG,F}$ acts on $\bL_G^{\tame}(\pi)$ via scalar. Then it follows from analysis of the geometry of the stack $\locsys_{{}^cG,F}^{\tame}$ that $\bL_G^{\tame}(\pi)$ must be scheme-theoretically supported on one irreducible component of $\locsys_{{}^cG,F}^{\tame}$. In addition, this component must be smooth and contains an open point which then must be a discrete parameter. The Cohen-Macaulayness of $\bL_G^{\tame}((i_1)_*\pi)$ then also implies that it must be a vector bundle on this irreducible component, giving the desired claim.

\quash{
\subsubsection{Cohomology of Shimura varieties} Once the unipotent categorical local Langlands correspondence is established, one can deduce \Cref{intro: coh of Shimura} easily using the existence of the following Cartesian diagram of perfect stacks over $k$
}

\subsection{Origin of ideas, some history, and relations to other works}
\subsubsection{} 
We briefly discuss the origin of ideas of this work and some history of this work.
In \cite{xiao2017cycles}, together we Liang Xiao, we applied the geometric Satake (in mixed characteristic) to construct correspondences between mod $p$ fibers of (different) Shimura varieties with hyperspecial level, which realizes certain cases of the Jacquet-Langlands correspondence in a geometric way (i.e. via cohomology of Shimura varieties). It was soon realized that the local theory of \emph{loc. cit.} is the application of the categorical trace construction in a very simple situation. See \cite{Zhu2016} for a survey.
However, in many applications in number theory (e.g. see \cite{LTXZZ}), it is desirable to generalize the constructions of \cite{xiao2017cycles} to Shimura varieties with the Iwahori level (or general parahoric level). 
This is the main motivation of the current work, although we will not really discuss such generalizations in this article. The current work can be regarded as a generalization of the local part of \cite{xiao2017cycles}.
It turns out that while at the hyperspecial level, we could work within the abelian category of perverse sheaves and could realize the categorical trace construction ``by hand", at the Iwahori level one must deal with the whole derived categories of $\ell$-adic sheaves and make use of machinery of higher categories to rigorously make sense of the categorical trace construction. It makes the whole story significantly more complicated.

This project begain with a collaboration with Tamir Hemo in 2019. In fact, the unipotent part of the categorical equivalence for $\overline\bQ_\ell$-coefficient (namely, the unipotent part of \Cref{intro: thm-langlands-categorical-tame}) was already established with Hemo at that time (see \cite{zhu2020coherent, zhu2022icm} for the announcement of some the results.) Along the way we have established some foundational results about $\shv(\kot_G)$ (such as \Cref{intro:thm-LL category}).
Since then, several new developments in the local geometric Langlands correspondence (see \cite{Bez.Riche.modular,DYYZ}) have enabled us to significantly extend our results. Specifically, we now have the categorical local Langlands correspondence at the tame level, and we also allow for modular coefficients. These generalizations have important applications (e.g. see \cite{yangIhara, yangzhutorsion} for applications of the modular coefficient categorical local Langlands). But achieving them required a major revision and generalization of the previous results obtained jointly with Hemo.
We apologize for the long delay in releasing the article. 

Ultimately, Hemo decided to let us retain the article in its entirety without being listed as a coauthor. Some of the key ideas in the article, such as the consideration of a geometric version of the categorical trace in the context of an abstract setting (\Cref{sec:trace-geometric-trace-main-general}), belong to him (see also his thesis \cite{Hemo}). This approach allows one to bypass the integral transform found in the works of Ben-Zvi-Nadler (e.g. see \cite{benzvi2009character,benintegralcoh, ben2017spectral}),  which is crucial in application, as such integral transform results typically do not hold in the $\ell$-adic setting. The development of a general theory of $\ell$-adic sheaves in \Cref{sec:pspl-stacks} is also largely joint with Hemo. In particular, the terminology of sind-placid stack is suggested by him.

\subsubsection{} Let us also briefly discuss the relation between our work and related works in this subject. As already mentioned, Fargues-Scholze \cite{Fargues.Scholze.geometrization} proposed another version of the categorical local Langlands conjecture, in which instead of the category $\shv(\kot_G,\La)$, they use the category $\der_{\mathrm{lis}}(\Bun_G,\La)$ of lisse sheaves on $\Bun_G$, whose definition is quite different from $\shv(\kot_G,\La)$ given in this article.
The main achievement of \cite{Fargues.Scholze.geometrization} is the construction of the so-called spectral action on $\der_{\mathrm{lis}}(\Bun_G)$, from which they extracted semi-simple Langlands parameters for every irreducible representation of the $p$-adic group, as mentioned in \Cref{intro: rem FS parameter}. However, \cite{Fargues.Scholze.geometrization} did not prove any equivalence of categories. They did not construct any functor from one side to another. A candidate of the local Langlands functor in Fargues-Scholze's approach was constructed later on by Hansen \cite{Hansen.Beijing}.

So besides the formal analogy of the categorical local Langlands conjecture, there is no direction relation between our work and the work of Fargues-Scholze. In other words, our work is independent of Fargues-Scholze's work. 
Nevertheless, one expects that the category $\shv(\kot_G,\La)$ considered in this work and  $\der_{\mathrm{lis}}(\Bun_G)$ considered in \cite{Fargues.Scholze.geometrization} are canonically equivalent. In addition, one expects that there is also the spectral action on $\shv(\kot_G,\La)$ and the equivalence is compatible with the spectral actions.
There are notable advances towards such expectations. Indeed, by a work in preparation by Gleason, Hamann, Ivanov, Louren\c{c}o and Zou \cite{GHILZ}, there is a canonical defined equivalence $\shv(\kot_G,\La)\cong \der_{\mathrm{lis}}(\Bun_G,\La)$, at least when $\La$ is a torsion ring. On the other hand, very recently Eteve, Gaitsgory, Genestier, Lafforgue have announced a construction of a spectral action on $\shv(\kot_G,\La)$ when $F$ is a field of positive characteristic.
Anyway, if the above expectations hold in general, our categorical conjecture then would agree with the categorical Langlands conjecture in \cite{Fargues.Scholze.geometrization}. Such expectation also leads us to discover an exotic $t$-structure on $\shv(\kot_G,\La)$ in \Cref{prop: t-structures on llc-2}. Some applications to the cohomology of Shimura varieties are also inspired by advances in Fargues-Scholze's program, although the actually proofs are quite different.

As mentioned earlier, the idea of studying the classical local Langlands correspondence via taking the cateogrical trace of the local geometric Langlands correspondence has been in the air for sometime. E.g. see \cite{gaitsgory2016geometric, Zhu2016} for some general discussions/speculations. An important work towards this direction is the work by Ben-Zvi-Chen-Helm-Nadler \cite{benzvi2020CohSpr} (built on \cite{ben2017spectral}), which constructed a fully faithful embedding of the Iwahori block $\rep(G)^{[I]}$ of the category of smooth representations of $G(F)$ into the category of (ind)coherent sheaves on the stack of unipotent Langlands parameters when $G$ is a split reductive group, and when the coefficient $\La$ is a characteristic zero field. (Partial results in this direction were also obtained earlier by Hellmann \cite{Hellmann} via a more down-to-earth approach.)  
Although both \cite{benzvi2020CohSpr} and our work use categorical trace construction, these two works are using this construction in different ways.
For example,  \cite{benzvi2020CohSpr} constructed the fully faithful embedding $\rep(G)^{[I]}\to \indcoh(\locsys_{{}^cG,F}^{\widehat\unip})$ as a consequence of the identification of the endomorphisms of the coherent Springer sheaf with the (extend) affine Hecke algebra of $G(F)$. This amounts our \Cref{intro: thm: end of coh spr} for split group $G$ and characteristic zero coefficient field $\La$. 
However, we deduce \Cref{intro: thm: end of coh spr} as a consequence of our categorical equivalence (so the logic is reversed). 
In addition, Ben-Zvi-Chen-Helm-Nadler did not define $\shv(\kot_G)$ or anything similar. As a result, they did not have equivalence of categories. In fact, they did not say anything about $\rep(G(F))$ when the group $G(F)$ is not split. Let us also mention that under the same assumption of $G$ and $\La$ as in \cite{benzvi2020CohSpr}, Propp \cite{Propp} also proved that the unipotent coherent Springer sheaf is an honest coherent sheaf (rather than a complex), by a different method of ours. As far as I understand, he did not deal with any other coherent object corresponding to compact inductions, as we do in \Cref{thm: compact induction in the heart}.


\subsection{Organization, notations and conventions}
\subsubsection{Organization}
The article consists of two parts. The first is the main part, which deals with the categorical local Langlands correspondence and some of its consequences. 

In \Cref{sec: L-parameter}, we review and further study the stack of local Langlands parameters. The main results include: the study of geometry of the stack around the (essentially) discrete Langlands parameters, the study of the tame and unipotent part of the stack of local Langlands parameters, in particular the tame and unipotent the spectral Deligne-Lusztig induction. We also explain how to put such construction into the framework of categorical trace construction.

In \Cref{sec:kot-stack}, we define and study the local Langlands category $\shv(\kot_G)$. We prove the basic categorical properties of  $\shv(\kot_G)$, such as compact generation, canonical self-duality, semi-orthogonal decomposition, $t$-structure on the subcategory of admissible objects. As a warm-up, we explain how to relate the category of smooth representations of a $p$-adic group to the category of $\ell$-adic sheaves on the classifying stack of the $p$-adic group. 
Along the way, we also revisit the geometry of $\kot_G$, giving new proofs of some known results about the geometry of $\kot_G$. 

In \Cref{sec: unipotent and tame LL category}, we restrict our attention to the tame and the unipotent part of $\shv(\kot_G)$. Main results include: developing a general theory of monodromic sheaves on stacks with group action (\Cref{SS: monodromic sheaves}) which might be of independent interests, developing an affine Deligne-Lusztig theory parallel to the classical Deligne-Lusztig theory and put it into the framework of categorical trace construction. Along the way, we also discover a class of projective objects in the category of representations of finite group of Lie type.

In \Cref{sec:unip-cat-langlands}, we review input from the local geometric Langlands correspondence, and put everything together to prove our main theorems. We establish the categorical equivalence and prove a few additional properties of such equivalences. We give some first applications. In particular, we attach every depth zero supercuspidal representation an enhanced Langlands parameter.

In \Cref{sec: local-global-comp}, we express the \'etale cohomology of Shimura varieties of Hodge type over a $p$-adic field in terms of the coherent cohomology on the stack of local Langlands parameters. Besides the unipotent categorical local Langlands, another ingredient is the Igusa stack as constructed in \cite{DHKZ}. However, for our purpose, we just need perfect Igusa stack, for which we give a direct construction in \Cref{prop: perfect Igusa stack}.

In the very long second part, we assembly various general sense in category theory, and the basic facts about coherent sheaves and constructible sheaves.  

In \Cref{sec:categorical-preliminaries} we review and further develop the general formalism of trace construction in (higher) categories. As mentioned before, we also introduce the notion of admissible objects in general dualizable categories, which might be of independent interest.

In \Cref{sec:trace} we review and further develop the general sheaf theory. We also review and further develop some methods computing categorical traces arising from the convolution pattern in geometry.

In \Cref{S: theory of coherent sheaves} we review and further develop the theory of coherent sheaves in the derived algebraic geometry. Notably, we discuss the theory of coherent sheaves for algebraic stacks over fields of positive characteristic. As is well-known, the theory is much more subtle than the theory for stacks in characteristic zero. Many crucial facts in characteristic zero simply fail in positive characteristic. The theory of singular supports of coherent sheaves in positive characteristic also need some extra care (even for schemes).

In \Cref{sec:pspl-stacks} we carefully develop the theory of $\ell$-adic sheaves, for a general coefficient ring $\La$ (which is a $\bZ_\ell$-algebra satisfying certain conditions). We will first assemble various ingredients in literature to write down a six functor
formalism for ind-constructible sheaves on prestacks, making use of the full strength of extension of sheaf
theories as developed in \Cref{sec:symmetric-monoidal-and-projection-for-corr}. Then we restrict our attend to a large class of infinite dimensional stacks (which we call sind placid stacks), where the theory has better properties. Such class of stacks include classifying stack of locally profinite groups, as well as $\kot_G$. 
The materials developed in this section should be useful in other context (in particular in geometric representation theory).

\medskip

\subsubsection{Notations and conventions}We will make use of the following notations and conventions throughout the article.

\begin{itemize}
\item For a Galois extension $E/F$ of fields, let $\Ga_{E/F}$ denote the Galois group. For a field $F$, let $\overline{F}$ denote a fixed separable closure, and let $\Ga_F=\Ga_{\overline{F}/F}$.

\item We refer to the beginning of \Cref{sec: L-parameter} for our notations and conventions related to Galois groups for non-archimedean local fields.

\item Let $A\to B$ be a homomorphism of commutative rings. For an $A$-module $M$, let $M_B:=M\otimes_AB$ denote its base change to $B$. Similarly, if $X$ is a scheme (or a more general geometric object such as a stack) over $\Spec A$, we write $X_B=X\times_{\Spec A}\Spec B$.

\item Let $H$ be an algebraic group over a field. Let $H^\circ$ denote the neutral connected component of $H$. More generally, if $\mH$ is an affine smooth group scheme over a base commutative ring $B$, let $\mH^\circ\subset \mH$ denote the open group subscheme that is fiberwise connected.

\item For a positive integer $n$, let $\mu_n$ denote the finite group scheme (over a base scheme) of $n$th roots of unity.

\item If $A$ is a group of multiplicative type over a field $k$, we let 
\[
\xch(A)=\Hom(A_{\overline k}, \bG_{m,\overline k}),\quad \xcoch(A)=\Hom(\bG_{m,\overline k}, A_{\overline k}),
\] 
regarded as $\Ga_k$-modules. If $A$ is a split torus over a base scheme, we also write 
\[
\xch(A)=\Hom(A,\bG_m),\quad \xcoch(A)=\Hom(\bG_m,A),
\] 
and call them the weight lattice and coweight lattice of $A$.

\item Let $G$ be a connected reductive group over a field $E$. Let $Z_G$ denote the center of $G$.
 Let  $G_{\mathrm{der}}$ denote its derived group, which is a connected semisimple group. Let $G_{\mathrm{sc}}$ be the simply-connected cover of $G_{\mathrm{der}}$, and $G_\ad$ the adjoint quotient of $G$. Let $G_{\mathrm{ab}}=G/G_{\mathrm{der}}$ be the abelianization of $G$. Let $\pi_1(G)$ be the algebraic fundemantal group of $G$, regarded as a $\Ga_E$-module.
For further notations and conventions related to reductive groups over local fields, we refer to \Cref{SS: Iwahori-Weyl-group}.

\item Let $\hat{G}$ be the dual group of $G$, regarded as a reductive group scheme over $\bZ$, equipped with a pinning $(\hat{B},\hat{T},\hat{e})$, where $\hat{B}$ is a Borel subgroup of $\hat{G}$ with $\hat{U}\subset \hat{B}$ its unipotent radical, where $\hat{T}\subset\hat{B}$ is a maximal torus, and where $\hat{e}: \hat{U}\to \bG_a$ is a homomorphism such that its restriction to every simple root subgroup is an isomorphism. 
Let 
\[
2\rho: \bG_m\to \hat{G}
\]
be the cocharacter given by the sum of positive coroots of $\hat{G}$ (with respect to $(\hat{B},\hat{T})$). 
Let $\hat{G}_{\ad}$ be the adjoint group of $\hat{G}$ and let
\[
\rho_{\ad}: \bG_m\to\hat{G}_{\ad}
\] 
be the cocharacter given by the half sum of positive coroots of $\hat{G}$ (with respect to $(\hat{B},\hat{T})$). 

\item There is an action of $\Ga_E$ on $\hat{G}$ via the homomorphism 
\[
\xi: \Ga_E\to \mathrm{Out}(G)\cong \mathrm{Out}(\hat{G}) \cong \Aut(\hat{G},\hat{B},\hat{T},\hat{e}).
\] 
Let $\pr: \Ga_E\to\Ga_{\widetilde E/E}$ be the finite quotient of $\Ga_E$ by $\ker\xi$.  Let 
\[
{}^cG=\hat{G}\rtimes(\bG_m\times\Ga_{\widetilde E/E})
\] 
be the $C$-group of $G$, regarded as a group scheme over $\bZ$, where $\bG_m$ acts on $\hat{G}$ via the homomorphism $\bG_m\xrightarrow{\rho_{\ad}}\hat{G}_{\ad}\subset\Aut(\hat{G})$, and $\Ga_{\widetilde E/E}$ acts via $\xi$.

\item In this article, we will extensively use the language of higher categories. For our notations and conventions, we refer to \Cref{Subsec: recollection of category}.
\item Our notations and conventions regarding derived algebraic geometry can be found in \Cref{SS: derived algebraic geometry}.
\item Our notations and conventions regarding $\ell$-adic sheaves can be found in \Cref{sec:pspl-stacks}.
\end{itemize}

\noindent \bf Acknowledgement. \rm 
As mentioned above, the project started with a collaboration with Tamir Hemo. We wholeheartedly thank him for the collaboration at early stage of the project.
We thank Linus Hamann for teaching us Fargues-Scholze's program, Xuhua He for teaching us beautiful combinatorics of affine Weyl groups, and Weizhe Zheng for teaching us abstract formalism of sheaf theory.
We thank David Ben-Zvi, Roman Bezrukavnikov, Gurbir Dhillon, Dennis Gaitsgory, Yau-Wing Li, Simon Riche, Yakov Varshavsky, Xiangqian Yang, Zhiwei Yun for various very useful discussions.

The work is supported by NSF grant under DMS-2200940.

\newpage

\part{Main Content}

\section{The stack of local Langlands parameters}\label{sec: L-parameter}
In this section, we study the spectral side of the categorical local Langlands correspondence. That is, the category of coherent sheaves on the stack of local Langlands parameters.  We make use of the following notations throughout this section.
\begin{itemize}
\item We fix a non-archimedean local field $F$, and a separable closure $\overline{F}$ of $F$. Let $F^u\subset F^t\subset \overline{F}$ be the maximal unramified and tamely ramified extension of $F$ in $\overline{F}$. Let $\breve{F}$ be the completion of $F^u$. We also fix a separable closure $\overline{\breve F}$ of $\breve F$ and embedding $\overline{F}\subset\overline{\breve F}$. Let $\breve F^t=F^t\breve F\subset\overline{\breve F}$.
\item Let $\Ga_F$ be the absolute Galois group of $F$, and let $W_F\subset \Ga_F$ be the Weil group of $F$.
We write 
\[
P_F=\Ga_{\overline{F}/F^t}\cong \Ga_{\breve F^t} \subset I_F=\Ga_{\overline{F}/F^u}\cong \Ga_{\breve F}\subset W_F
\] 
for the inertia and wild inertia subgroups of $F$. Let
$W_F^t=W_F/P_F$ for the tame Weil group. Write $I_F^t=I_F/P_F=\Ga_{F^t/F^u}\cong \Ga_{\breve F^t/\breve F}$ for the tame inertia. 
\item Let 
\begin{equation}\label{eq: cycl character}
\cycl: W_F\to W_F/I_F\cong \bZ\to q^{\bZ}\subset \widehat{\bZ}^\times
\end{equation}
be the cyclotomic character, which sends the arithmetic Frobenius to $q:=\sharp k_F$. Let $\Ga_q$ be the $q$-tame group with two generators $\tau,\sigma$ satisfying the relation $\sigma\tau\sigma^{-1}=\tau^q$. 
\item Let
\begin{equation}\label{eq:iota-vs-tau}
t: I_F\to I^t_F\cong \widehat{\bZ}^p(1):=\lim_{(n,p)=1}\mu_n(\breve F)
\end{equation}
be the homomorphism obtained as follows: for each $n$ coprime to $p$, let $\varpi^{1/n}$ be a uniformizer of the unique degree $n$ extension of $\breve F$ in $\breve F^t$. Then $\tau(\varpi^{1/n})=a_n\varpi^{1/n}$ for some $a_n\in \mu_n(k)$ which is in fact independent of the choice of $\varpi^{1/n}$. Then $t$ sends $\tau$ to the compatible system $\{a_n\}_n$ of roots of unity. For $\ell\neq p$, let $t_\ell: I_t^F\to \bZ_\ell(1)$ be the projection of $t$ to the pro-$\ell$-part.
\end{itemize}

Let $G$ be a reductive group over $F$.  We write $\hat{G}$ be the dual group of $G$ and ${}^cG$ be the $C$-group of $G$. Let
\begin{equation}\label{eq: d and tildepr}
d: {}^cG\to \bG_m\times\Ga_{\widetilde F/F},\quad\quad \widetilde{\pr}=(\cycl^{-1},\pr): W_F\to \bZ[1/p]^\times\times\Ga_{\widetilde F/F},
\end{equation}
where the first map is the natural projection.

\subsection{Some geometry of the stack of local Langlands parameters}

\subsubsection{Space of continuous representations}\label{SS: Space of continuous representations}

Recall that there is the stack of local Langlands parameters $\locsys_{{}^cG, F}$ over $\bZ_\ell$, which classifies continuous homomorphisms $\rho: W_F\to {}^cG$ such that $d\circ\rho=\widetilde{\pr}$ up to $\hat{G}$-conjugation. 
We recall the construction following \cite{zhu2020coherent}. 

Let $\Ga$ be a locally profinite group and let $H$ be a flat affine group scheme of finite type over $\bZ_\ell$. Then there is the moduli space $(R_{\Ga,H})_\cl$ of \emph{strongly continuous} homomorphisms from $\Ga$ to $H$. 
By definition,
\[
(R_{\Ga,H})_\cl: \aff_{\bZ_\ell}^{\heartsuit}\to \spc,\quad R\mapsto \Hom_{cts}(\Ga, H(R)),
\]
where $\Hom_{cts}(\Ga, H(R))$ consist of  homomorphisms $\varphi:\Ga\to H(R)$ such that for one (and therefore every) faithful representation $H\to \GL(M)$ on a finite free $\bZ_\ell$-module $M$, for every $m\in M\otimes R$, and for every open compact subgroup $\Ga_0$ of $\Ga$, the $\bZ_\ell$-module in $N\subset M\otimes \bR$ spanned by $\varphi(\Ga_0)m$ is finite and the resulting representation of $\varphi(\Ga_0)$ on $N$ is continuous (in the usual sense). For our purpose, we also need to recall how to extend $(R_{\Ga,H})_\cl$ as functor $R_{\Ga,H}$ from the category $\aff_{\bZ_\ell}$ of animated $\bZ_\ell$-algebras to $\spc$. 

We work with the ordinary category of ind-profinite sets, and write $C_{cts}(-,-)$ for the hom set in this category. We may regard a $\bZ_\ell$-module as an ind-profinite set (by writing a $\bZ_\ell$-module as an inductive limit of finitely generated ones, which can be regarded as profinite sets), and then regard a $\bZ_\ell$-algebra as an ind-profinite set by regarding it as a $\bZ_\ell$-module.
If $S$ is a profinite set, then we may regard $C_{cts}(S,-)$ as a functor from $\aff_{\bZ_\ell}^{\heartsuit}$ to itself, which preserves sifted colimits. Taking animation gives $C_{cts}(S,-): \aff_{\bZ_\ell}\to \aff_{\bZ_\ell}$. If $S=\sqcup_j S_j$ is a disjoint union of profinite set, we let $C_{cts}(S,-)=\prod_j C_{cts}(S_j,-)$. Now, we consider the simplicial set $\Ga^\bullet$ given by the group structure of $\Ga$. Then we have $C_{cts}(\Ga^\bullet,-): \aff_{\bZ_\ell}\to \aff_{\bZ_\ell}^{\Delta}$, where $\aff_{\bZ_\ell}^\Delta$ denotes the category of cosimplicial animated $\bZ_\ell$-algebras. Then we define
\[
R_{\Ga,H}: \aff_{\bZ_\ell}\to \spc,\quad R\mapsto \mathrm{Map}_{\aff_{\bZ_\ell}^\Delta}(\bZ_\ell[H^\bullet], C_{cts}(\Ga^\bullet, R)).
\]
(The space $R_{\Ga,H}$ was denoted by $R_{\Ga,H}^{sc}$ in \cite[\textsection{2.4}]{zhu2020coherent}.)
One checks without difficulty that if $R$ is an ordinary $\bZ_\ell$-algebra, $R_{\Ga, H}(R)=\Hom_{cts}(\Ga, H(R))$.

The conjugation action of $H$ on itself induces a conjugation action of $H$ on $R_{\Ga,H}$. We let $\mX_{\Ga,H}:=R_{\Ga,H}/H$ denote the quotient stack.
If $H$ is smooth over $\bZ_\ell$, one can show that the tangent complex of the quotient stack $\mX_{\Ga,H}$ exists and at a classical point $\varphi$ is given by $C_{cts}(\Ga^\bullet, \Ad_\varphi)[1]$, where 
\[
\Ad_\varphi\colon \Ga\xrightarrow{\varphi} H\xrightarrow{\Ad} \GL(\frakh)
\] 
denote the induced representation of $\Ga$ on the Lie algebra $\frakh$ of $H$, and $C_{cts}(\Ga^\bullet, \Ad_\varphi)$ is considered as a chain complex via the Dold-Kan correspondence. In particular, the degree $i$ term of the tangent complex at $\varphi$ is given by $H^{i+1}_{cts}(\Ga, \Ad_\varphi)$, the $(i+1)$th continuous cohomology of $\Ga$ with coefficient $\Ad_\varphi$.

\begin{remark}\label{rem: repn space of discrete group}
We note that when $\Ga$ is a discrete group, the moduli space $R_{\Ga,H}$ makes sense over $\La$ for any commutative ring $\La$, as soon as $H$ is an affine smooth group scheme defined over $\La$. In addition, it is easy to see that in this case $R_{\Ga,H}$ is represented by a (possibly derived) affine scheme.
\end{remark}

\begin{example}\label{ex: continuous representation of Zhat}
Suppose the neutral connected component $H^\circ$ of $H$ is reductive and $H/H^\circ$ is finite \'etale. Let $H/\!\!/H$ be the GIT quotient of $H$ by adjoint action. Let $(H/\!\!/H)^\wedge\subset H/\!\!/H$ be the union of closed subschemes that are finite over $\bZ_\ell$. Then it is easy to see that $R_{\widehat{\bZ}, H}\cong H\times_{H/\!\!/ H}(H/\!\!/H)^\wedge$.
In particular, $R_{\widehat{\bZ},H}$ is represented by an ind-affine scheme ind-of finite type $\bZ_\ell$. Note that the map $R_{\widehat{\bZ},H}\to H$ induces isomorphisms of tangent spaces. In particular, $R_{\widehat{\bZ},H}$ is formally smooth over $\bZ_\ell$.

We let $(H/\!\!/H)^{\wedge,p}\subset (H/\!\!/H)^{\wedge}$ be the union of those subschemes in $Z\subset H/\!\!/H$ that are finite over $\bZ_\ell$, such that $Z(\overline{\bF}_\ell)$ lift to points in $H(\overline\bF_\ell)$ of order prime-to-$p$.
Let $\widehat{\bZ}^p=\prod_{\ell\neq p}\bZ_\ell$ be the maximal pro-$p$-quotient of $\widehat{\bZ}$. 
Then the map $R_{\widehat{\bZ}^p, H}\to R_{\widehat{\bZ},H}$ induces the isomorphism 
\[
R_{\widehat{\bZ}^p, H}\cong H\times_{H/\!\!/ H}(H/\!\!/H)^{\wedge,p}.
\] 
Again the map $R_{\widehat{\bZ}^p, H}\to H$ induces isomorphisms of tangent spaces. Therefore, $R_{\widehat{\bZ}^p, H}$ is also an ind-affine scheme, ind-of finite type and formally smooth over $\bZ_\ell$.
\end{example}

\subsubsection{Space of Langlands parameters}\label{SS: Space of Langlands parameters}
Now we let the space of $L$-parameters as
\begin{equation}\label{eq: moduli of L-parameters}
\locsys_{{}^cG, F}=\locsys_{{}^cG, F}^{\Box}/\hat{G},\quad\quad \locsys_{{}^cG, F}^{\Box}= R_{W_F,{}^cG}\times_{R_{W_F, \bG_m\times \Ga_{\widetilde F/F}}} \{\widetilde{\pr}\},
\end{equation}
where we regard $\widetilde{\pr}: W_F\to \bG_m\times\Ga_{\widetilde F/F}$ as a $\bZ_\ell$-point of $R_{W_F, \bG_m\times \Ga_{\widetilde F/F}}$.

If $L$ is a Galois extensions of $F$ (in $\overline{F}$) that is finite over $F^t\widetilde{F}$, let $\Ga_L\subset \Ga_F$ be the Galois group of $L$. Then we can define $\locsys_{{}^cG, L/F}^{\Box}$ as above, with $W_F$ replaced by $W_F/\Ga_L$. 

We recall the following basic facts about  $\locsys_{{}^cG, F}$ (see \cite[\textsection{3.1}]{zhu2020coherent}, and also \cite{dat2020moduli} and \cite[Chapter VIII]{Fargues.Scholze.geometrization}).
\begin{theorem}\label{thm: geometry of stack of L-parameters}
The moduli space $\locsys_{{}^cG, F}^{\Box}$ is represented by a classical scheme over $\bZ_\ell$, which is a union 
\[
\locsys_{{}^cG, F}^{\Box}=\colim_{L} \locsys_{{}^cG, L/F}^{\Box},
\]
where $L$ ranges over all Galois extensions of $F$ (in $\overline{F}$) that are finite over $F^t\widetilde{F}$. Each $ \locsys_{{}^cG, L/F}^{\Box}$ is represented by a reduced
affine scheme flat and of finite type over $\bZ_\ell$, is equidimensional of dimension$=\dim\hat{G}$, and is a local complete intersection. If $L'/L$ is finite, then the inclusion $\locsys_{{}^cG, L/F}^{\Box}\subset \locsys_{{}^cG, L'/F}^{\Box}$ is open and closed.
It follows that
\begin{equation}\label{eq: ind-presentation of locsys}
\locsys_{{}^cG,F}=\colim_{L} \locsys_{{}^cG, L/F},
\end{equation}
where $\locsys_{{}^cG, L/F}=\locsys_{{}^cG, L/F}^{\Box}/\hat{G}$ is a classical algebraic stack of relative dimension zero over $\bZ_\ell$. 
\end{theorem}

Let $Z_{{}^cG,L/F}:= H^0(\locsys_{{}^cG,L/F},\mO)$ be the ring of regular functions on $\locsys_{{}^cG,L/F}$.
We regard $\locsys_{{}^cG,F}$ as ind-algebraic stack via the presentation \eqref{eq: ind-presentation of locsys}. Then 
we let
\begin{equation}\label{eq: stable Bernstein center}
Z_{{}^cG,F}=H^0(\locsys_{{}^cG,F},\mO):=\lim_{L} H^0(\locsys_{{}^cG,L/F},\mO)
\end{equation}
be the ring of regular functions on $\locsys_{{}^cG, F}$, which then is regarded as a pro-algebra. Let
\[
\Spf Z_{{}^cG,F}:=\colim_L \Spec Z_{{}^cG,L/F},
\]
which can be regarded as the coarse moduli space of $\locsys_{{}^cG,F}$.

Let $P$ be a rational parabolic subgroup of $G$ with Levi quotient $M$. Let $\hat{P}$ and $\hat{M}$ be the corresponding dual. The action of $\bG_m\times\Ga_{\widetilde F/F}$ on $\hat{G}$ preserves $\hat{P}$ and $\hat{M}$, so we can form ${}^cP$ and ${}^cM$ respectively and similarly define 
\begin{equation}\label{eq:locP and locM}
\locsys_{{}^cP, F}=\locsys_{{}^cP, F}^{\Box}/\hat{P}\to  \locsys_{{}^cM, F}=\locsys_{{}^cM, F}^{\Box}/\hat{M}.
\end{equation}
It turns out that in general $\locsys_{{}^c\!P, F}$ has non-trivial derived structure. But it is still quasi-smooth. 

There is the following commutative diagram over $\bZ_\ell$
\begin{equation}\label{E:LocMtoG}
\xymatrix{
&\ar@/_/[dl]_-{r}\locsys_{{}^cP,F}\ar^-{\pi}[dr]&  \\
\ar_{\varpi_{{}^cM}}[d]\locsys_{{}^cM,F}\ar@/_/[ur]_-{i}&& \locsys_{{}^cG,F}\ar^{\varpi_{{}^cG}}[d]  \\
\Spec Z_{{}^cM,F}\ar[rr] && \Spec Z_{{}^cG,F}.
}
\end{equation}
where $\pi,r,i$ are induced by the corresponding morphisms between $\hat{G},\hat{P},\hat{M}$, and where the bottom map is induced by $\pi\circ i: \locsys_{{}^cM,F}\to \locsys_{{}^cG,F}$.
The morphism $\pi$ is schematic and is proper, while $r$ is quasi-smooth.

After a choice of a homomorphism $\iota: \Ga_q\to W_F^t$ sending $\tau$ to a generator of the tame inertia and $\sigma$ to a lifting of the Frobenius, there is also an algebraic stack $\locsys_{{}^cG,F,\iota}$ over $\bZ[1/p]$, together with a canonical isomorphism $\locsys_{{}^cG,F,\iota}\otimes_{\bZ[1/p]}\bZ_\ell\cong \locsys_{{}^cG,F}$. Namely, let 
\begin{equation}\label{eq: GaFiota}
\Ga_{F,\iota}:=W_F\times_{W_F^t,\iota}\Ga_q.
\end{equation}
Then $\Ga_{F,\iota}$ is an extension of $\Ga_q$ by $P_F$.
Similarly, for a Galois extensions $L/F$ that is finite over $F^t\widetilde{F}$, let $\Ga_{L/F,\iota}:= W_F/\Ga_L\times_{W_F^t,\iota}\Ga_q$, which is an extension of $\Ga_q$ by a finite $p$-group $Q_L=\mathrm{Gal}(L/F^t)$. 
The map $\widetilde{\pr}$ from \eqref{eq: d and tildepr} induces a homomorphism $\Ga_{L/F,\iota}\to q^\bZ\times\Ga_{\widetilde F/F}$ still denoted by $\widetilde{\pr}$.
Let
\[
\locsys_{{}^cG,F,\iota}^\Box=\colim_L R_{\Ga_{L/F,\iota},{}^cG}\times_{R_{\Ga_{L/F,\iota}, \bG_m\times \Ga_{\widetilde{F}/F}}}\{\widetilde{\pr}\},\quad \locsys_{{}^cG,F,\iota}=\locsys_{{}^cG,F,\iota}^\Box/\hat{G}.
\]
As $\Ga_{L/F,\iota}$ now is a discrete group, by \Cref{rem: repn space of discrete group} the above spaces make sense over $\bZ[1/p]$\footnote{They make sense even over $\bZ$ but we shall only consider them over $\bZ[1/p]$.}.  
In addition, by \cite[Proposition 3.1.6]{zhu2020coherent} each space $R_{\Ga_{L/F,\iota},{}^cG}\times_{R_{\Ga_{L/F,\iota}, \bG_m\times \Ga_{\widetilde{F}/F}}}\{\widetilde{\pr}\}$ as above is represented by an affine scheme flat, local complete intersection, of finite type over $\bZ[1/p]$. If $L'/L$ is a finite, then $R_{\Ga_{L/F,\iota},{}^cG}\times_{R_{\Ga_{L/F,\iota}, \bG_m\times \Ga_{\widetilde{F}/F}}}\{\widetilde{\pr}\}\subset R_{\Ga_{L'/F,\iota},{}^cG}\times_{R_{\Ga_{L'/F,\iota}, \bG_m\times \Ga_{\widetilde{F}/F}}}\{\widetilde{\pr}\}$ is open and closed.

For different choice of $\iota,\iota': \Ga_q\to W_F^t$, the resulting spaces $\locsys_{{}^cG,F,\iota}$ and $\locsys_{{}^cG,F,\iota'}$ (over $\bZ[1/p]$) are in general different. However, by \cite[Lemma 3.1.8, Corollary 3.1.12]{zhu2020coherent}, we have
\begin{itemize}
\item the natural inclusion $\Ga_{F,\iota}\subset W_F$ induces a canonical isomorphism
\[
\locsys_{{}^cG,F,\iota}^\Box\otimes\bZ_\ell\cong \locsys_{{}^cG,F}^\Box, \quad \locsys_{{}^cG,F,\iota}\otimes\bZ_\ell\cong \locsys_{{}^cG,F};
\]
\item the ring of regular functions $H^0(\locsys_{{}^cG,F,\iota},\mO)$ is independent of the choice of $\iota$ and gives a canonical extension of $Z_{{}^cG,F}$ as a (pro-)algebra over $\bZ[1/p]$.
\end{itemize}
Note that the first point above (together with the geometry of $\locsys_{{}^cG,F,\iota}$) in particular implies \Cref{thm: geometry of stack of L-parameters}.

\begin{remark}\label{rem: cG vs LG}
Let ${}^L\!G=\hat{G}\rtimes\Ga_{\widetilde{F}/F}$ be the usual full Langlands dual group of $G$. One can define a version of moduli $\locsys_{{}^L\!G,F}$ of $L$-parameters by replacing ${}^cG$ everywhere in the above discussions by ${}^L\!G$ (and replacing the requirement $d\circ \varphi=\widetilde{\pr}$ by the requirement $d\circ \varphi=\pr$, where $d: {}^L\!G\to \Ga_{\widetilde{F}/F}$ is the projection). 

If we fix $\sqrt{q}$, then the cyclotomic character $\cycl$ admits a square root $\cycl^{\frac{1}{2}}$, which induces a homomorphism over $\bZ_\ell[\sqrt{q}^{\pm}]$
\[
\locsys_{{}^L\!G,F}\otimes \bZ_\ell[\sqrt{q}^\pm]\cong \locsys_{{}^cG,F}\otimes \bZ_\ell[\sqrt{q}^\pm],\quad  \varphi\mapsto \tilde\varphi, 
\]
where if we write $\varphi(\ga)=(\varphi_0(\ga), \pr(\ga)))\in \hat{G}\rtimes \Ga_{\widetilde{F}/F}$, then 
\[
\tilde\varphi(\ga)=(\varphi_0(\ga)2\rho(||\ga||^{\frac{1}{2}}),\widetilde{\pr}(\ga)))\in \hat{G}\times(\bG_m\times\Ga_{\widetilde{F}/F}),
\] 
and  where $2\rho$ denotes the sum of all positive coroots of $\hat{G}$.
\end{remark}

We recall some symmetries of $\locsys_{{}^cG,F}$. 

\begin{enumerate}
\item Let $\theta$ be an automorphism of the pinned group $(\hat{G},\hat{B},\hat{T}, \hat{e})$ that sends $\mu\in \xch(\hat{T})$ to $-w_0(\mu)$. This is usually called Cartan involution of $\hat{G}$, which commutes with any pinned automorphism of $\hat{G}$, as well as the conjugation action by $\rho_\ad(\bG_m)$. Therefore, the Cartan involution induces an automorphism of ${}^cG$ (and ${}^L\!G$), and therefore an automorphism
\begin{equation}\label{eq:Cartan involution of loc}
\theta: \locsys_{{}^cG,F}\to \locsys_{{}^cG,F}.
\end{equation}

\item
The $\Ga_{\widetilde F/F}$-fixed point subscheme $Z_{\hat{G}}^{\Ga_{\widetilde F/F}}$ of the center $Z_{\hat{G}}$ of $\hat{G}$ is a flat group scheme of multiplicative type over $\bZ_\ell$ (and is smooth if $\La$ is a field of characteristic zero). Let 
\begin{equation}\label{eq: connected center of dual group}
C_{{}^cG}\subset Z_{\hat{G}}^{\Ga_{\widetilde F/F}}
\end{equation}  
be the maximal subtorus.
Let  $G'$ be the intersection of all kernels of rational characters of $G$. We note that $G/G'$ is a split torus over $F$. If we let $Z_G^s\subset Z_G$ denote the maximally $F$-split torus in the center of $G$, then the composed map
\begin{equation}\label{eq: split center to split abelianizable}
Z_G^s\to G\to G/G'
\end{equation}
is an isogeny. We note that $C_{{}^cG}$ is identified as the dual torus of $G/G'$ and ${}^cG'={}^cG/C_{{}^cG}$.

Note that for every $L$ as above, $R_{W_F/\Ga_L, C_{{}^cG}}$ has a natural structure as a group scheme over $\bZ_\ell$, and there is a free action 
\begin{equation}\label{eq: center action on loc}
R_{W_F/\Ga_L, C_{{}^cG}}\times \locsys_{{}^cG,L/F}\to  \locsys_{{}^cG,L/F},\quad (\psi: W_F\to C_{{}^cG},\ \varphi: W_F\to {}^cG)\mapsto \psi\varphi.
\end{equation}
This induces an isomorphism
\begin{equation}\label{eq: get rid of center of the dual group}
\locsys_{{}^cG, L/F}/(R_{W_F/\Ga_L,C_{{}^cG}}/C_{{}^cG})=\locsys_{{}^cG',L/F},
\end{equation}
where we consider the trivial action of $C_{{}^cG}$ on $R_{W_F/\Ga_L,C_{{}^cG}}$.
It follows that we get a free action of
$R_{W_F/\Ga_L,C_{{}^cG}}$ on $\Spec Z_{{}^cG,L/F}$, and $\Spec Z_{{}^cG,L/F}/R_{W_F/\Ga_L,C_{{}^cG}}=\Spec Z_{{}^cG',L/F}$.
\end{enumerate}

\subsubsection{$\phi$-fixed point construction}\label{SS: phi-fixed point construction spectral side}
In spirit of the trace construction, we would like to express $\locsys_{{}^cG,F}$ as a $\phi$-fixed point subscheme. Recall that there is a general $\phi$-fixed points construction (as from \eqref{eq:tau-fixed-point}). Namely, if $X$ is an object equipped with an automorphism $\phi$ in a category $\bfC$ (admitting
finite products), then we let
\[
\mL_\phi(X):=X\times_{\id\times\phi, X\times X,\Delta} X. 
\]
Now if $\phi_1$ and $\phi_2$ are two automorphisms of $X$ and $\al: \phi_1\simeq \phi_2$ is an isomorphism, there $\al$ induces an isomorphism 
\begin{equation}\label{eq: isomorphism of loop spaces}
\mL_\al:\mL_{\phi_1}(X)\simeq \mL_{\phi_2}(X).
\end{equation}

We specialize to the case where the category $\bfC$ is the category of ind-algebraic stacks (as defined in \Cref{def-ind-scheme-derived}) over $\La$. Let $V$ be a(n ind-)scheme equipped with an action $\mathrm{act}: V\times H\to V$ by an affine flat group scheme $H$ of finite type over $\La$. Suppose $V$ and $H$ are equipped with automorphisms $\phi_V$ and $\phi_H$ compatible with the action map. Then the quotient stack $X=V/H$ is equipped with an automorphism $\phi$. In this case
\[
\mL_\phi(X)\cong  \bigl(V \times_{\id\times \phi_V,V\times V, \pr_1\times \mathrm{act}}(V\times H)\bigr)/H,
\]
Here, in the formulation of the quotient, $H$ acts on $V$ via the action map $\mathrm{act}$ and on $H$ via the $\phi_H$-twisted conjugation $\Ad_{\phi_H}$, i.e. $h\in H$ acts on $H$ by sending $h'\mapsto h^{-1}h'\phi_H(h)$. Therefore, (algebraically closed field valued) points of $\mL_\phi(X)$ can be identified with pairs $(v,h)\in V\times H$ satisfying $vh=\phi_V(v)$, up to $H$-conjugacy.

If we replace $\phi_V(-)$ by $ \phi_V(-)h_0$ for some $h_0\in H$, and replace $\phi_H$ by $h_0^{-1}\phi_Hh_0$, then we obtain a new automorphism of $X$, denoted by $\phi_{h_0}$. We have a canonical isomorphism
\begin{equation}\label{eq: isomorphisms between fixed points construction}
\mL_{\phi_{h_0}}(X)\cong \mL_{\phi}(X), 
\end{equation}
induced by the map $V\times H\to V\times H, \ (w,h)\mapsto (w, hh_0)$.

Now we apply the above considerations to the study of the stack of local Langlands parameters.
Recall that we can identify $I_F\cong \mathrm{Gal}(\overline{\breve F}/\breve F)$. Let $\Ga_{\widetilde{\breve F}/\breve F}$ denote the image of $I_F$ in $\Ga_{\widetilde F/F}$, and let ${}^L\!G_{\breve F}:=\hat{G}\rtimes \Ga_{\widetilde{\breve F}/\breve F}$. This is the Langlands dual group of $G_{\breve F}$.
Then we consider a moduli space as the same definition of $\locsys_{{}^cG,F}$ but with $W_F$ replaced by $I_F$. Explicitly, 
\[
\locsys_{{}^cG,\breve F}=\locsys_{{}^cG,\breve F}^{\Box}/\hat{G},\quad \locsys_{{}^cG,\breve F}^{\Box}=R_{I_F,{}^cG}\times_{R_{I_F,\bG_m\times\Ga_{\widetilde{F}/F}}}\{\widetilde{\pr}\},
\] 
which classifies all strongly continuous homomorphisms $\breve\varphi: I_F\to {}^cG$ such that $d\circ \breve\varphi= \widetilde{\pr}$.
 If $L/\breve F$ is a Galois extension (in $\overline{\breve F}$) finite over $\widetilde{F}\breve F^t$, we also have $\locsys_{{}^cG,L/\breve F}=\locsys_{{}^cG,L/\breve F}^{\Box}/\hat{G}$ as above, with $I_F$ replaced by $\mathrm{Gal}(L/\breve F)$ in the definition.  We note that such $\breve\varphi$ necessarily sends $I_F$ to  so one can replace ${}^cG$ by ${}^LG_{\breve F}$ in the definition, and write $\locsys_{{}^L\!G_{\breve F},\breve F}$ instead of $\locsys_{{}^cG,\breve F}$.

The difference now is that $\locsys_{{}^cG,L/\breve F}^{\Box}$ is no longer represented by
an affine scheme, but rather by an ind-affine scheme. More precisely, we have the following. 

\begin{proposition}\label{prop: geometric Langlands parameters as ind-stack}
We have 
\[
\locsys_{{}^cG, \breve F}^{\Box}=\colim_{L} \locsys_{{}^cG, L/\breve F}^{\Box},
\]
where $L$ ranges over all Galois extensions of $\breve F$ (in $\overline{\breve F}$) that are finite over $\widetilde{F}\breve F^t$. Each $ \locsys_{{}^cG, L/\breve F}^{\Box}$ is represented by an ind-affine scheme, ind-of finite type and formally smooth over $\bZ_\ell$. If $L'/L$ is finite, then $\locsys_{{}^cG, L/\breve F}^{\Box}\subset \locsys_{{}^cG, L'/\breve F}^{\Box}$ is open and closed.
\end{proposition}
\begin{proof}
We use the same argument as in \cite[Proposition 2.3.9]{zhu2020coherent}, and reduce to show that if $H$ is an affine smooth group scheme over $\mO$ (a finite extension of $\bZ_\ell$), with its neutral connected component $H^\circ$ reductive over $\mO$ and $H/H^\circ$ (finite) \'etale, then $R_{I_F^t, H}$ is represented by an ind-affine scheme, formally smooth over $\mO$. 

We choose a topological generator $\tau$ of $I_F^t$, given an isomorphism $\widehat{\bZ}^p\cong I_F^t$. This induces an isomorphism $R_{I_F^t, H}\cong H\times_{H/\!\!/H}(H/\!\!/H)^{\wedge,p}$ (using \Cref{ex: continuous representation of Zhat}).  The proposition then follows.
\end{proof}

Consider the morphism
\begin{equation}\label{eq: res from arith parameter to geom parameter}
\res: \locsys_{{}^cG,F}\to \locsys_{{}^cG,\breve F}
\end{equation}
obtained by restriction along $I_F\subset W_F$. 

By abuse of notations, we will use $\sigma$ to denote a lifting of the arithmetic Frobenius to $W_F$. Let 
\begin{equation}\label{eq: barsigma}
\bar\sigma=\widetilde{\pr}(\sigma)\in \bG_m(\bZ_\ell)\times \Ga_{\widetilde F/F}\subset {}^cG(\bZ_\ell).
\end{equation} 

Then the conjugation action of $\sigma$ on $I_F$ and the action of $\bar\sigma$ on ${}^cG$ by conjugation together induce an automorphism  
\begin{equation}\label{eq: phi automorphism of loc}
\phi: \locsys^\Box_{{}^cG,\breve F}\to \locsys^\Box_{{}^cG,\breve F}, \quad \breve\varphi\mapsto (\phi(\breve\varphi): \ga\mapsto \bar\sigma(\breve\varphi(\sigma^{-1}\ga\sigma)), \quad \ga\in I_F).
\end{equation}
We still use $\phi$ to denote the induced automorphism of $\locsys_{{}^cG,\breve F}$.

\begin{lemma}\label{lem: arith parameter as fix point of geometric parameter}
We have a canonical isomorphism $\locsys_{{}^cG,F}\cong \mL_\phi (\locsys_{{}^cG,\breve F})$.
\end{lemma}
\begin{proof}
Note that the map \eqref{eq: res from arith parameter to geom parameter} fits into the following commutative diagram
\begin{equation}\label{eq: arith parameter as fix point of geometric parameter}
\xymatrix{
\locsys_{{}^cG,F}\ar^{\res_\phi}[r]\ar_{\res}[d] & \locsys_{{}^cG,\breve F} \ar^{\Delta}[d]\\
\locsys_{{}^cG,\breve F}\ar^-{\id\times\phi}[r] & \locsys_{{}^cG,\breve F}\times_{\bZ_\ell} \locsys_{{}^cG,\breve F},
}
\end{equation}
which induces a map $\locsys_{{}^cG,F}\to \mL_\phi (\locsys_{{}^cG,\breve F})$. Indeed, as all the moduli spaces in the above diagram are classicial, to check its commutativity, it is enough to check the commutativity when evaluated at classical $\bZ_\ell$-algebras. In this case, it follows that giving a point $\varphi$ of $\locsys_{{}^cG,F}^{\Box}$ is the same as giving a point $\breve\varphi$ of $\locsys_{{}^cG,\breve F}^{\Box}$ and an element $g\in \hat{G}$ such that $\breve\varphi=g\phi(\breve\varphi)g^{-1}$. Namely, given $\breve\varphi$ and $g\in\hat{G}$ we can define $\varphi$ such that $\varphi|_{I_F}=\breve\varphi$ and $\varphi(\sigma)=g\bar\sigma$, and vice versa. 

This in fact already implies that the map $\locsys_{{}^cG,F}=(\locsys_{{}^cG,F})_\cl\to (\mL_\phi (\locsys_{{}^cG,\breve F}))_\cl$ is an isomorphism.
To check that it is an isomorphism at the derived level, it is enough to check that the map induces an isomorphism of tangent spaces at classical points. 
Now the tangent space of the left hand side at $\varphi$ is given by $C_{cts}((W_F)^\bullet, \Ad_\varphi^0)$, where $\Ad_\varphi^0$ denotes the representation of $W_F$ on $\hat\frakg$ via $W_F\xrightarrow{\varphi} {}^cG\xrightarrow{\Ad} \hat\frakg$, while the tangent space of the right hand side at $\varphi$ is the fiber of $C_{cts}((I_F)^{\bullet}, \Ad_{\varphi}^0)\xrightarrow{1-\phi} C_{cts}((I_F)^\bullet, \Ad_{\varphi}^0)$. Now the desired isomorphism follows from the fiber sequence 
\[
C_{cts}((W_F)^\bullet, \Ad_\varphi^0)\to C_{cts}((I_F)^{\bullet}, \Ad_{\varphi}^0)\xrightarrow{1-\phi} C_{cts}((I_F)^\bullet, \Ad_{\varphi}^0).
\]
\end{proof}

\begin{remark}\label{rem: fix point construction depending on outer automorphism}
We fix a lifting $\sigma$.
For every automorphism $a:  {}^L\!G_{\breve F}\to {}^L\!G_{\breve F}$ such that the induced automorphism of $\Ga_{\widetilde{\breve F}/\breve F}$ coincides with the automorphism induced by conjugation by $\sigma$ on $I_F$, one can similarly define an automorphism $\phi_a$ of $\locsys_{{}^cG, \breve F}$ sending $\breve\varphi$ to $\phi_a(\breve\varphi)$ where $\phi_a(\breve\varphi)(\ga)=a(\breve\varphi(\sigma^{-1}\ga\sigma))$. Then we have
the space $\mL_{\phi_a}(\locsys_{{}^cG, \breve F})$. If $b(-)=\delta^{-1}a(-)\delta$ for some $\delta\in \hat{G}$, then by \eqref{eq: isomorphisms between fixed points construction} we have an natural isomorphism 
\[
\mL_{\phi_a}(\locsys_{{}^cG, \breve F})\stackrel{\cong}{\to} \mL_{\phi_{b}}(\locsys_{{}^cG, \breve F}), \quad (\breve\varphi, g)\mapsto (\breve\varphi, g\delta).
\]
Therefore, up to isomorphism the space $\mL_{\phi_a}(\locsys_{{}^cG, \breve F})$ depends only on the image of $a$ in $\Aut({}^LG_{\breve F})/\hat{G}$.

We apply the above discussion to the following situations. 
\begin{enumerate}
\item Let $\sigma'$ be another lifting of $\sigma$, giving another automorphism $\phi'$ of $\locsys_{{}^cG,\breve F}^{\Box}$. As $\sigma'= \sigma\delta$ for some $\delta\in I_F$, we see that $\phi'(-)= \bar\sigma(\delta)^{-1}\phi(-)\bar\sigma(\delta)$. Then we have 
\[
\mL_\delta: \mL_{\phi}(\locsys_{{}^cG,\breve F})\cong \mL_{\phi'}(\locsys_{{}^cG,\breve F})
\] 
sending $(\breve\varphi, g)$ to $\breve\varphi, g\bar\sigma(\delta)$.  It is easy to see that $\mL_\delta$ is compatible with the isomorphism in \Cref{lem: arith parameter as fix point of geometric parameter}.

\item Let $a=\bar\sigma$, we have $a(-)=\delta^{-1}a(-)\delta$ for every $\delta\in C_{{}^cG}$. Thus, every $\delta\in C_{{}^cG}$ gives rise to an automorphism of $\mL_\phi(\locsys_{{}^cG,\breve F})$. On the other hand, we may regard $\delta$ as an element in $R_{W_F, {}^cG}$ which sends $I_F$ to $1$ and $\sigma$ to $\delta$. Therefore, \eqref{eq: center action on loc} provides another automorphism of $\mL_\phi(\locsys_{{}^cG,\breve F})$. Clearly, these two automorphisms match
each other under the isomorphism from \Cref{lem: arith parameter as fix point of geometric parameter}.

\item We apply the above consideration to $a=\bar\sigma$ and $b(-)=2\rho(\sqrt{q})a(-) 2\rho(\sqrt{q}^{-1})$, we recover the isomorphism
\Cref{rem: cG vs LG} between the two versions of Langlands parameters over $\bZ_\ell[\sqrt{q}^{\pm 1}]$.
\end{enumerate}
\end{remark}

\begin{notation}\label{notation: base change of loc}
Let $Z\to \locsys_{{}^cG,\breve F}$ be a morphism. In the sequel, we write
\[
\locsys_{{}^cG,F}^{Z}:=Z\times_{\locsys_{{}^cG,\breve F}}\locsys_{{}^cG,F}.
\]
\end{notation}

The same proof of \Cref{lem: arith parameter as fix point of geometric parameter} gives the following.
\begin{lemma}\label{lem: phi fixed point for Z}
Let $Z\subset  \locsys_{{}^cG,\breve F}$ be a $\phi$-stable (finitely presented) locally closed embedding, 
and let $\widehat{Z}$ be its formal completion in $\locsys_{{}^cG,\breve F}$. Then
we have a natural isomorphism $\mL_{\phi}(\widehat{Z})=\locsys_{{}^cG,F}^{\widehat{Z}}$.
\end{lemma}

The presentation of $\locsys_{{}^cG,F}$ as $\phi$-fixed points of $\locsys_{{}^cG,\breve F}$ leads a decomposition  $\locsys_{{}^cG,F}$ into open and closed substacks refining 
\eqref{eq: ind-presentation of locsys}. It also leads a parameterization of  irreducible components of $\locsys_{{}^cG,F}$. We start with the discussion of the former.

Similar to \eqref{eq: stable Bernstein center}, we define a pro-algebra $Z_{{}^cG,L/\breve F}= H^0\Gamma(\locsys_{{}^cG,L/\breve F},\mO)$ for finite extension $L/\breve F$ as in \Cref{prop: geometric Langlands parameters as ind-stack} and let
\[
Z_{{}^cG,\breve F}=H^0\Gamma(\locsys_{{}^cG,\breve F},\mO):= \lim_L Z_{{}^cG,L/\breve F}
\]
As explained in \cite[Remark 2.2.20]{zhu2020coherent}, $\La$-points of $\Spf Z_{{}^cG, \breve F}$ are the same as (continuous) pseudorepresentations of $I_F$. 
Recall that by  \cite[Proposition 2.3.25]{zhu2020coherent} (see also discussions around displayed equation (2.34) in \emph{loc. cit.}),
each $\Spf Z_{{}^cG,L/\breve F}$ is a formal scheme formally of finite type over $\bZ_\ell$, with reduced subscheme finite over $\bZ_\ell$.

The automorphism $\phi$ of $\locsys_{{}^cG,\breve F}$ induces an automorphism of $Z_{{}^cG,\breve F}$, still denoted by $\phi$. We let 
\[
(\Spf Z_{{}^cG,\breve F})^\phi=\colim_{L} (\Spf Z_{{}^cG,L/\breve F})^\phi
\] 
be the \emph{classical} $\phi$-fixed point subscheme of $\Spf Z_{{}^cG,\breve F}$.
Note that \eqref{eq: res from arith parameter to geom parameter} induces a morphism $\Spf Z_{{}^cG,F}\to\Spf Z_{{}^cG, \breve F}$, which clearly factors as $\Spf Z_{{}^cG,F}\to (\Spf Z_{{}^cG, \breve F})^{\phi}\subset \Spf Z_{{}^cG, \breve F}$.

\begin{lemma}\label{rem: 3rd def of inertia type}
Every connected component of  $(\Spf Z_{{}^cG, \breve F})^{\phi}$ is a scheme finite over $\bZ_\ell$. 
\end{lemma}
\begin{proof}
We follow the argument of \cite[Lemma 3.4.3]{zhu2020coherent} (with slightly different notations). Let $Z_{{}^cG,\breve F}\to A^{\Theta}$ be a surjective homomorphism, corresponding to a connected component of $\Spf Z_{{}^cG,\breve F}$, and let
$B^{\Theta}:=A^{\Theta}/(1-\phi)A^{\Theta}$, which is a complete noetherian $\bZ_\ell$-algebra. We need to show that it is finite over $\bZ_\ell$. It is enough to show that $B^{\Theta}/\ell$ is artinian over $\bF_\ell$. Let $x: B^{\Theta}\to \kappa[[t]]$ be a homomorphism, where $\kappa$ is finite field extension of $\bF_\ell$. As argued in \cite[Lemma 3.4.3]{zhu2020coherent} (which relies on \cite[Lemma 2.4.14]{zhu2020coherent}), there is some $\varphi \in \locsys_{{}^cG,F}^{\Box}(\Spf \mO_K)$ for some finite extension $K/\kappa((t))$ (such $\varphi$ corresponds to a continuous representation $\varphi: W_F\to {}^cG(\mO_K)$,  where $\mO_K$ is equipped with $t$-adic topology), such that $\varphi|_{I_F}\in \locsys_{{}^cG, \breve F}^\Box(\Spf \mO_K)$ is over $x\in \Spf Z_{{}^cG,F}(\Spf \kappa[[t]])$. As $\locsys_{{}^cG, F}$ is an algebraic stack locally of finite presentation, $\varphi$ comes from a $\Spec \mO_K$-point of $ \locsys_{{}^cG,F}^{\Box}$. I.e., $\varphi$ is continuous now $\mO_K$ is equipped with the discrete topology.
It follows that $\varphi(I_F)$ has finite image. This will imply that the image of the map $B^\Theta/\ell\to \kappa[[t]]$ is contained in $\kappa$. The lemma follows.
\end{proof}

\begin{remark}
As is clear from the above argument, the key ingredient is the algebraicity of $\locsys_{{}^cG,F}$, which implies that the image of $\varphi(I_F)$ is finite, for every continuous representation $W_F\to {}^cG(\kappa((t)))$ where $\kappa((t))$ is equipped with the $t$-adic topology. The analogous statement when $F$ is a global function field is known as de Jong's conjecture, which is much deeper and was proved by Gaitsgory (via the global Langlands correspondence). In fact, in \cite{zhu2020coherent}, de Jong's conjecture was the key input to prove that the analogous stack $\locsys_{{}^cG,F}$ of global Langlands parameters (for global function field $F$)  is algebraic.
\end{remark}

\begin{definition}\label{def: inertia type}
Let $\La$ be an algebraically closed field. 
An inertia type $\zeta$ of ${}^cG$ over $\La$ is a $\La$-point of $(\Spf Z_{{}^cG, \breve F})^\phi$. 
\end{definition}

Note that by \Cref{rem: 3rd def of inertia type}, every inertia type is defined over an algebraic extension of $\bF_\ell$ or $\bQ_\ell$.
Here is yet another equivalent definition. 

\begin{lemma}\label{lem: second definition of inertia type}
Let $\La$ be an algebraically closed field.
There is a bijection between inertia types $\zeta$ over $\La$ and $\hat{G}(\La)$-conjugacy class of completely reducible representations $\breve\varphi^{ss}: I_F\to {}^L\!G(\La)$ with finite image
that can be extended to a homomorphism $W_F\to {}^cG(\La)$ giving a $\La$-valued point of $\locsys_{{}^cG,F}$.
\end{lemma}
Here a representation $\breve\varphi^{ss}:I_F\to {}^L\!G(\La)$ is called completely reducible if  the image $\breve\varphi^{ss}(I_F)$ in ${}^L\!G$ is completely reducible. I.e. if $\breve\varphi^{ss}(I_F)$ is contained in an $R$-parabolic subgroup of ${}^L\!G$, then it is contained in an $R$-Levi subgroup of this parabolic subgroup. (We refer to \cite[\textsection{6}]{BMR} for the notions of $R$-parabolic and $R$-Levi in a possibly disconnected reductive group.)
\begin{proof}
Let $\zeta$ be an inertia type.
By definition, there is a finite extension $L/\widetilde{F}\breve F^t$ such that the inertia type $\zeta$ comes from a $\La$-point of $\Spf H^0\Gamma(\locsys_{{}^cG, L/\breve F},\mO)$. Then by \cite[11.7]{Lafforgue2018Chtoucas} (and \cite[Remark 2.2.20]{zhu2020coherent}), $\zeta$ can be lifted to a unique $\La$-point of $\locsys_{{}^cG, L/\breve F}$, corresponding to a completely reducible continuous representation $\breve\varphi^{ss}: \mathrm{Gal}(L/\breve F) \to {}^L\!G(\La)$ up to $\hat{G}$-conjugacy. 
As $\phi(\breve\varphi)$ is still completely reducible, giving $\phi(\zeta)$ in the coarse moduli space, and
as $\zeta$ is $\phi$-fixed, there is some $g\in \hat{G}$ such that $g\phi(\breve\varphi)g^{-1}=\breve\varphi$.  The argument of \Cref{lem: arith parameter as fix point of geometric parameter} implies that $\breve\varphi$ extends to a $W_F$-representation.  It remains to prove that $\breve\varphi$ is of finite image.

If $\La$ is of characteristic $\ell$, then $\breve\varphi^{ss}(I_F)$ is finite as $\breve\varphi^{ss}$ is continuous. So we assume that $\La$ is of characteristic zero.
It follows from the standard argument that for any topological generator $\tau$ of $I_F^t$, lifted to $\mathrm{Gal}(L/\breve F)$,  the semisimple part $\breve\varphi^{ss}(\tau)_s$ of $\breve\varphi(\tau)$ is of finite order. We write $\breve\varphi^{ss}(\tau)=\breve\varphi(\tau)_s \cdot \breve\varphi^{ss}(\tau)_u$ for the Jordan decomposition. 
We claim that $\breve\varphi^{ss}(\tau)_u=1$. Indeed, $\breve\varphi^{ss}$ induces a map $\overline{\breve\varphi^{ss}}: I_F^t\to N_{{}^L\!G}(\breve\varphi(P_F))/\breve\varphi(P_F)$ which is still semisimple. Therefore $\overline{\breve\varphi^{ss}}(\tau)_u=1$. It follows that $\breve\varphi^{ss}(\tau)_u$ belongs to $\breve\varphi^{ss}(P_F)$, which is a finite $p$-group. Therefore, we must have $\breve\varphi^{ss}(\tau)_u=1$. So in any case $\breve\varphi^{ss}(\tau)$ is of finite order. The lemma is proved.
\end{proof}

\begin{remark}\label{rem: complete reducible vs semisimple}
Let $\breve\varphi^{ss}:I_F\to {}^L\!G(\La)$ be a completely reducible representation associated to an inertia type as above. When $\La$ is of characteristic zero, then $\breve\varphi^{ss}(\ga)$ is always a semisimple element of ${}^L\!G$.
This, however, may not be the case when $\La$ is a field over $\bF_\ell$. Indeed, when $\ell$ divides the order of $\Ga_{\widetilde F/F}$, then the homomorphism $I_F\to \Ga_{\widetilde{\breve F}/\breve F}\xrightarrow{\ga\mapsto (1,\ga)}{}^L\!G$ gives an example of $\breve\varphi^{ss}$ that associates to an inertia type. But the image of this map contains non semisimple elements.
\end{remark}

In the sequel, for an inertia type $\zeta$ over $\La$, we let $\hat\zeta$ denote the formal completion of $\Spf Z_{{}^cG,\breve F}\otimes \La$ at $\zeta$. 
Note that if $\La=\overline\bF_\ell$, 
$\Spf Z_{{}^cG,\breve F}\otimes \overline\bF_\ell$ is formal at $\zeta$ so  $\hat\zeta$ is the connected component of $\Spf Z_{{}^cG, \breve F}\otimes\overline\bF_\ell$ that contains $\zeta$ as the unique closed point. We also let
\[
\locsys_{{}^cG,F}^{\hat\zeta}\to \locsys_{{}^cG,\breve F}^{\hat\zeta}
\] 
denote the preimages of $\hat\zeta$ under the maps $\locsys_{{}^cG,F}\to \locsys_{{}^cG,\breve F}\to \Spf Z_{{}^cG,\breve F}$. 
As $\zeta$ is $\phi$-fixed, the $\phi$-action on $\locsys_{{}^cG, \breve F}$ restricts to a $\phi$-action of $\locsys_{{}^cG, \breve F}^{\hat\zeta}$, and by \Cref{lem: phi fixed point for Z} we have
\[
\locsys_{{}^cG,F}^{\hat\zeta}\cong \mL_\phi(\locsys_{{}^cG,\breve F}^{\hat\zeta}).
\]

Note that a priori, $\locsys_{{}^cG,F}^{\hat\zeta}$ is a formal stack. But we have the following.

\begin{lemma}\label{lem: algebraicity of zeta component}
The formal stack $\locsys_{{}^cG,F}^{\hat\zeta}$ is a finite union of connected components $\locsys_{{}^cG,F}$, and therefore is an algebraic stack of finite presentation over $\La$.
\end{lemma} 
\begin{proof}
Note that $\locsys_{{}^cG,F}$ maps to $(\Spf Z_{{}^cG, \breve F})^{\phi}$, which is a disjoint union of schemes finite over $\bZ_\ell$. It follows that every connected component of $\locsys_{{}^cG,F}\otimes \La$ maps set-theoretically  to one point of $(\Spf Z_{{}^cG, \breve F})^{\phi}\otimes\La$. (If $\La$ is $\overline\bF_\ell$, see also \cite[Lemma 2.4.25]{zhu2020coherent} and \cite[Remark 3.1.2]{zhu2020coherent}.) Therefore, every connected component of $\locsys_{{}^cG,F}\otimes \La$ will map to some $\hat\zeta$. On the other hand, given an inertia type $\zeta$, there are only finitely many connected components of $\locsys_{{}^cG,F}$ that map to $\hat\zeta$ (as all of such components must be contained in $\locsys_{{}^cG,L/F}$ from \Cref{thm: geometry of stack of L-parameters}, for some $L$). The lemma is proved.
\end{proof}

\begin{remark}\label{rem: locsyszeta comment}
\begin{enumerate}
\item The stack $\locsys_{{}^cG,F}^{\hat\zeta}$ may still be disconnected (e.g. see \Cref{ex: regular supercuspidal} below). But in some important cases, it is connected (e.g. see \Cref{lem: connectedness of unipotent stack}).
\item We may regard $\zeta$ as the closed point of $\hat\zeta$. Then we have $\locsys_{{}^cG,F}^\zeta\to \locsys_{{}^cG,\breve F}^\zeta$. The inclusion $\locsys_{{}^cG,F}^\zeta\subset \locsys_{{}^cG,F}^{\hat\zeta}$ induces an isomorphism of the underlying reduced substacks. But $\locsys_{{}^cG,F}^\zeta$ usually has non-trivial derived structure.
\end{enumerate}
\end{remark}

In the above discussions we decompose $\locsys_{{}^cG, F}$ according to points of $\Spf Z_{{}^cG, \breve F}$ that are fixed by $\phi$.
Next we consider irreducible components of $\locsys_{{}^cG,F}$. Informally, the idea is to consider (finite type) points of $\locsys_{{}^cG, \breve F}$ that are fixed by $\phi$.
We assume that $\La$ is algebraically closed in the sequel, and base change everything to $\La$. To simplify notations, we omit $\La$ from the subscriptions. 

Let $\bO$ be a (finite type) point of $\locsys_{{}^cG, \breve F}$ over $\zeta$, regarded as a locally closed substack (more precisely as the residual gerbe at this point in the sense of \cite[\href{https://stacks.math.columbia.edu/tag/06ML}{Section 06ML}]{stacks-project}) of $\locsys_{{}^cG, \breve F}$.
Let $\bO^{\Box}$ be its preimage in $\locsys_{{}^cG, \breve F}^{\Box}$. So $\bO=\bO^{\Box}/\hat{G}$.

\begin{lemma}\label{lem: locO}
Suppose $\bO$ is $\phi$-stable, i.e., for $\breve\varphi\in \bO^\Box$, we have $\delta\phi(\breve\varphi)\delta^{-1}=\breve\varphi$ for some $\ga\in\hat{G}$ as in the proof of \Cref{lem: arith parameter as fix point of geometric parameter}. Then
$\locsys_{{}^cG,F}^{\bO}=\bO\times_{\locsys_{{}^cG, \breve F}}\locsys_{{}^cG,  F}$
is locally closed in $\locsys_{{}^cG,F}$ of dimension zero. Each connected component of $\locsys_{{}^cG,F}^{\bO}$ is irreducible.
\end{lemma} 
\begin{proof}
As mentioned above, we regard $\bO\subset \locsys_{{}^cG,\breve F}$ as a $\phi$-stable locally closed substack. Taking $\phi$-fixed points gives a morphism
\begin{equation}\label{eq: O-component of loc}
\mL_\phi(\bO)\to  \locsys_{{}^cG,F}^{\bO}\subset \mL_\phi(\locsys_{{}^cG,\breve F})=\locsys_{{}^cG,F}.
\end{equation}
It is enough to show that $\mL_\phi(\bO)$ is of dimension zero, whose connected components coincide with irreducible components, and
the first morphism induces an isomorphism of underlying classical stacks.

In the sequel of the proof, we will ignore the derived structure on the involved schemes/stacks. 
Let
\[
C_{\bO}=\left\{(\breve\varphi, g)\in \bO^\Box\times \hat{G}\mid g\breve\varphi g^{-1}=\breve\varphi\right\}.
\]
This is a flat group scheme over $\bO^\Box$, whose fiber over $\breve\varphi\in \bO^\Box$ is the centralizer $C_{\hat{G}}(\breve\varphi)$ of $\breve\varphi$ in $\hat{G}$. In particular, when $\La$ is a field of characteristic zero, this group scheme is smooth over $\bO^\Box$. 

If $g\in C_{\hat{G}}(\breve\varphi)$, then $\bar\sigma(g)$ in $C_{\hat{G}}(\phi(\breve\varphi))$. In addition, if $\bO$ is $\phi$-stable, so that there is $\delta\in \hat{G}$ such that $\delta\phi(\breve\varphi)\delta^{-1}=\breve\varphi$, then we obtain an automorphism 
\begin{equation}\label{eq: automorphism phidelta}
\phi_\delta: C_{\hat{G}}(\breve\varphi)\to C_{\hat{G}}(\breve\varphi),\quad h\mapsto  \delta\bar\sigma(h) \delta^{-1}.
\end{equation}
We let $\Ad_{\phi_\delta}$ be the $\phi_\delta$-twisted conjugation action of $C_{\hat{G}}(\breve\varphi)$ on itself. I.e. $g\in C_{\hat{G}}(\breve\varphi)$ acts on $C_{\hat{G}}(\breve\varphi)$ by sending $h\mapsto gh\phi_\delta(g)^{-1}$.

We can summarize the above discussions as saying that after choosing $\breve\varphi\in \bO^\Box$ and $\delta\in C_{\hat{G}}(\breve\varphi)$, we have
$\bO\cong \bB_{\mathrm{fppf}} C_{\hat{G}}(\breve\varphi)$, such that the $\phi$ action on $\bO$ is identified with the $\phi_\delta$ action on $C_{\hat{G}}(\breve\varphi)$. This implies that
\begin{equation}\label{eq: phi-fix point of bO}
\mL_\phi(\bO)\simeq C_{\hat{G}}(\breve\varphi)/\Ad_{\phi_\delta} C_{\hat{G}}(\breve\varphi).
\end{equation} 
So it is of dimension zero, with irreducible components and connected components coincide.

On the other hand, a choice of such $\delta$ amounts to an extension of $\breve\varphi$ to a Langlands parameter $\varphi$ by requiring $\varphi(\sigma)=\delta\bar\sigma$.
In this case, it is clear that $\res^{-1}(\bO^\Box)\subset \locsys_{{}^cG,F}^{\Box}$ is a (left) $C_\bO$-torsor. Namely,  an element $g\in C_{\hat{G}}(\breve\varphi)$ sends an extension  $\varphi_1:W_F\to {}^cG$ of $\varphi_0$ to another extension $\varphi_2$ with $\varphi_2(\sigma):=g\varphi_1(\sigma),\ \varphi_2|_{I_F}:=\varphi_1|_{I_F}=\breve\varphi$. There is another right $C_\bO$-torsor structure on $\res^{-1}(\bO)$, given by sending $(\varphi_1, g)$ to $\varphi_2$ with $\varphi_2|_{I_F}=\varphi_1|_{I_F}=\breve\varphi$ and $\varphi_2(\sigma)= \varphi_1(\sigma)\bar\sigma(g)\bar\sigma$.
Therefore, once we fix an extension $\varphi$ of $\breve\varphi$ to a Langlands parameter (equivalently an element $\ga\in \hat{G}$ such that $\ga\phi(\breve\varphi)\ga^{-1}=\breve\varphi$),
we have (at the level of classical stacks)
\begin{equation*}\label{eq: locsysO}
\locsys_{{}^cG,F}^{\bO}= \res^{-1}(\bO^\Box)/\hat{G}\simeq C_{\hat{G}}(\breve\varphi)/\Ad_{\phi_\delta} C_{\hat{G}}(\breve\varphi).
\end{equation*}
The lemma is proved.
\end{proof}

\Cref{lem: locO} implies that after ignoring possible derived and non-reduced structures, the closure of connected components of $\locsys_{{}^cG,F}^\bO$ inside $\locsys_{{}^cG,F}$ give irreducible components of $\locsys_{{}^cG,F}$.
We now would like to give a parameterization of $\pi_0\locsys_{{}^cG,F}^\bO$.

Let 
\[
A(\breve\varphi)=\pi_0 C_{\hat{G}}(\breve\varphi)
\] 
denote the group of connected components of $C_{\hat{G}}(\breve\varphi)$. 
The $\phi_\delta$-twisted conjugation $\Ad_{\phi_\delta}$ induces a $\phi_\delta$-twisted conjugation action of $A(\breve\varphi)$ on itself, still denoted by $\Ad_{\phi_\delta}$. 

Let $A(\breve\varphi)/\!\!/\Ad_{\phi_\delta}A(\breve\varphi)$ be the quotient set.
If we replace $\delta$ by $\delta'=g\delta$ for some $g\in C_{\hat{G}}(\breve\varphi)$, then $\phi_\ga$ is replaced by $\phi_{\delta'}=\Ad_g\phi_{\delta}$, and $A(\breve\varphi)/\!\!/\Ad_{\phi_\delta}A(\breve\varphi)$ is canonically identified with $A(\breve\varphi)/\!\!/\Ad_{\delta'}A(\breve\varphi)$ given by $x\mapsto x\bar{g}^{-1}$, where $\bar{g}$ is the image of $g$ in $A(\breve\varphi)$.
Therefore $\phi_\delta$ is well-defined up to inner automorphism of $C_{\hat{G}}(\breve\varphi)$, and $A(\breve\varphi)/\!\!/\Ad_{\phi_\delta}A(\breve\varphi)$ is independent of the choice of $\delta$ up to a canonical isomorphism. 

We will make the $\phi_\delta$-action on $A(\breve\varphi)$ more explicit when we restrict our attention to stack of unipotent Langlands parameters. But at the moment,
we arrive at the following statement. (See also \cite[Theorem 1.5]{dat2020moduli}.)

\begin{proposition}\label{prop: irr comp of loc}
Let $\La$ be an algebraically closed field. 
Irreducible components of $\locsys_{{}^cG, F}\otimes\La$ are indexed by $(\bO, x)$, where $\bO$ is a $\phi$-stable $\hat{G}$-orbit in $\locsys^{\Box}_{{}^cG,\breve F}$, and $x\in A(\breve\varphi)/\!\!/\Ad_{\phi_\delta}A(\breve\varphi)$. 
\end{proposition}

We also recall that there is the action of $R_{W_F,C_{{}^cG}}$ on $\locsys_{{}^cG,F}$ (see \eqref{eq: center action on loc}). Let $\psi: W_F\to C_{{}^cG}$,  and let $\breve\psi$ denote its restriction to $I_F$. Clearly, the action of $\psi$  on $\locsys_{{}^cG,F}$ will send $\locsys_{{}^cG,F}^{\bO}$ to $\locsys_{{}^cG,F}^{\breve\psi\bO}$. In particular, the torus $C_{{}^cG}$, regarded as the subspace of $R_{W_F,C_{{}^cG}}$ consisting of those $\psi$ such that $\breve\psi$ is trivial, will act freely on $\locsys_{{}^cG,F}^{\bO}$.
Let $\bO'$ be the image of $\bO$ under the map $\locsys_{{}^cG,\breve F}\to \locsys_{{}^cG', \breve F}$. Then \eqref{eq: get rid of center of the dual group} induces an isomorphism
\begin{equation}\label{eq: get rid of center of the dual group-2}
\locsys_{{}^cG, F}^{\bO}/(C_{{}^cG}/C_{{}^cG})=\locsys_{{}^cG',F}^{\bO'}.
\end{equation}

\subsubsection{Frobenius semisimplification and Weil-Deligne representations}
We assume that $\La=\overline\bQ_\ell$. We will fix $\sqrt{q}\in\overline\bQ_\ell$ and work of ${}^L\!G$-valued representations of $W_F$ (as in \Cref{rem: cG vs LG}). We let $\cycl^{\frac{1}{2}}$ be the square root of the cyclotomic character determined by $\sqrt{q}$.

We first recall the ``Jordan decomposition" of homomorphisms $\breve\varphi: I_F\to {}^L\!G$.

\begin{lemma}\label{lem: Jordan decomposition of representations}
Let $\varphi: W_F\to {}^L\!G(\La)$ be a point of $\locsys_{{}^L\!G,F}$, and let $\breve\varphi=\varphi|_{I_F}$. Then $\breve\varphi$ admits a unique ``Jordan decomposition" $\breve\varphi=\breve\varphi^{ss}\breve\varphi^u$, where
$\breve\varphi^{ss}: I_F\to {}^L\!G$ is a completely reducible representation with finite image associated to the inertia type $\zeta$ of $\varphi$ as in \Cref{lem: second definition of inertia type}, and where
$\breve\varphi^u: I_F\to C_{\hat{G}}(\breve\varphi^{ss})=:\hat{G}_\zeta$ is homomorphism which factors as $I_F \twoheadrightarrow I_F^t\twoheadrightarrow \bZ_\ell(1)\to \bG_a(\La)$ for some unipotent subgroup $\bG_a\subset C_{\hat{G}}(\breve\varphi^{ss})$. In addition, we have $C_{\hat{G}}(\breve\varphi)=C_{\hat{G}_\zeta}(\breve\varphi^u)$.
\end{lemma}
\begin{proof}
For every $\ga\in I_F$, we may write the Jordan decomposition $\breve\varphi(\ga)=\breve\varphi^{ss}(\ga)\breve\varphi^u(\ga)$ with $\breve\varphi^{ss}(\ga)$ semisimple and $\breve\varphi^{u}(\ga)$ unipotent. Note that $\breve\varphi^u(\ga)\in \hat{G}$ and $\breve\varphi^u(\ga)$ is trivial if $\ga\in P_F$. In addition, as argued in \Cref{lem: second definition of inertia type}, conjugation by $\breve\varphi(\ga)$ induces an automorphism of $\breve\varphi(P_F)$ which is a finite $p$-group. Therefore the unipotent part $\breve\varphi(\tau)_u$ of $\breve\varphi(\tau)$ acts trivially on $\breve\varphi(P_F)$. It follows that $\breve\varphi^{ss}: I_F\to {}^L\!G$ is a homomorphism with finite image. This is a completely reducible representation associated to the inertia type of $\varphi$.
In addition, $\breve\varphi^u$ is a continuous homomorphism from $I_F$ to $\hat{G}_\zeta$, trivial on $P_F$ with values in unipotent elements in $\hat{G}_\zeta$. By continuity, such homomorphism necessarily factors through $I_F^t\to \bZ_\ell(1)$.
The last statement is clear.
\end{proof}

\begin{remark}
We do not know whether the analogous statement holds when $\La$ is a field over $\bF_\ell$. When $\varphi$ is tame, i.e. $\varphi|_{P_F}$ is trivial, such a decomposition does exist. In fact, after fixing $\iota: \Ga_q\to W_F^t$ as before, this amounts to decomposing $g\bar\tau\in \hat{G}\bar\tau$ into $g\bar\tau= g_1g_2$ such that $g_1\in \hat{G}\bar\tau$ whose $\hat{G}$-orbit under conjugation is closed in $\hat{G}\bar\tau$, and $g_2\in C_{\hat{G}}(g_1)$ is unipotent. 
The existence of such decomposition follows from \cite[\textsection{5}]{xiao.zhu}. See more details in the proof of \Cref{prop: reducing to endoscopic group}. 
However, by virtual of \Cref{rem: complete reducible vs semisimple}, this decomposition may be different from the usual Jordan decomposition of $g\bar\tau$, regarded as an element in the non-connected algebraic group $\hat{G}\rtimes\langle\bar\tau\rangle$.
\end{remark}

Recall that  in the traditional formulation of the local Langlands correspondence, a local Langlands parameter is a continuous $\ell$-adic representation of $W_F$ (or a Weil-Deligne representation) such that the image of (a lifting of) the Frobenius is a semisimple element in ${}^L\!G$.

\begin{lemma}\label{lem: Frob semisimplification}
Let $\varphi: W_F\to {}^L\!G(\La)$ be a point on $\locsys_{{}^cG,F}$. Let $\varphi(\sigma)=\varphi(\sigma)^s\varphi(\sigma)^u$ be the Jordan decomposition of the $\varphi(\sigma)$.
We let $\varphi^{F\mbox{-}ss}: W_F\to {}^L\!G$ be the map sending $\ga\sigma^n$ to $\varphi(\ga)(\varphi(\sigma)^{s})^n$. Then $\varphi^{F\mbox{-}ss}$ also defines point on $\locsys_{{}^L\!G,F}$ which is independent of the choice of the lifting of the Frobenius $\sigma$. 
\end{lemma}
We call $\varphi^{F\mbox{-}ss}$ the Frobenius semisimplification of $\varphi$.
This fact is of course well-known, at least when $G=\GL_n$. We include a proof for completeness.
\begin{proof}
Let $\varphi: W_F\to {}^L\!G$ be parameter, and let $\breve\varphi$ be its restriction to $I_F$. Then $\breve\varphi=\breve\varphi^{ss}\breve\varphi^u$ admits a unique Jordan decomposition as in \Cref{lem: Jordan decomposition of representations}. We see that $\varphi(\sigma)$ normalizes both $\breve\varphi^{ss}$ and $\breve\varphi^{u}$. Then as argued in \Cref{lem: Jordan decomposition of representations}, $\varphi(\sigma)^u$ in fact centralizes $\breve\varphi^{ss}$. Let $\iota(\tau)$ be a tame generator of $I_F^t$. Then $\breve\varphi^u(\iota(\tau))=\exp(X)$ for some nilpotent element $X\in \hat\frakg$, which is an eigenvector of $\Ad_{\varphi(\sigma)}$ with eigenvalue $q$. It follows that $\varphi(\sigma)^u$ also centralizes $\breve\varphi^u$. This already implies that $\varphi^{F\mbox{-}ss}$ is well-defined. To see that $\varphi^{F\mbox{-}ss}$ is independent of the choice of the lifting $\sigma$, we need to show that $\varphi(\ga\sigma)^s=\varphi(\ga)\varphi(\sigma)^s$ for every $\ga\in I_F$. Since $\varphi(\ga\sigma)=\varphi(\ga)\varphi(\sigma)^s\varphi(\sigma)^u$, and since $\varphi(\sigma)^u$ is unipotent commuting with $\varphi(\ga)\varphi(\sigma)^s$, it is enough to show that $\varphi(\ga)\varphi(\sigma)^s=\breve\varphi^{ss}(\ga)\breve\varphi^u(\ga)\varphi(\sigma)^s$ is a semisimple element. Let $v$ be the unipotent element in $\hat{G}$ such that $v^{q-1}=\breve\varphi^u(\ga)$. Then $v$ commutes with $\breve\varphi^{ss}(\ga)$ and $\varphi(\sigma)^s v (\varphi(\sigma)^s)^{-1}=v^q$. Therefore, $v \varphi(\ga)\varphi(\sigma)^s v^{-1}=\breve\varphi^{ss}(\ga)\varphi(\sigma)^s$. So it remains to show that $\breve\varphi^{ss}(\ga)\varphi(\sigma)^s$ is semisimple. As $\varphi(\sigma)^s$ normalizes $\breve\varphi^{ss}$, we see that certain power of $\breve\varphi^{ss}(\ga)\varphi(\sigma)^s$ is a product of two commuting semisimple elements, and therefore is semisimple. This finishes the proof of the lemma.
\end{proof}

\begin{remark}\label{rem: non-frob semisimple parameters}
Here are some consequences of the argument.
Let $\varphi\in\locsys_{{}^L\!G,F}(\La)$ and let  $\varphi^{F\mbox{-}ss}$ be its Frobenius semisimplification. Then $v:=\varphi(\sigma)^u$ is independent of the choice of the lifting of the Frobenius $\sigma\in W_F$, and $v\in M:=C_{\hat{G}}(\varphi^{F\mbox{-}ss})$. Then we have $C_{\hat{G}}(\varphi)=C_M(v)$. It follows that there is a morphism $\mU_{M^\circ}/M$ from the adjoint quotient of the unipotent variety of $M$ to $\locsys_{{}^L\!G,F}^\bO$, where $\bO\in \locsys_{{}^L\!G,\breve F}$ is the point given by $\breve\varphi$, such that both $\varphi$ and $\varphi^{F\mbox{-}ss}$ are in the image of this map. All points in the image have the same Frobenius semisimplification. 

We also note that the above map $\mU_{M^\circ}/M\to \locsys_{{}^L\!G,F}^\bO$ in fact extends to a morphism $M/M\to \locsys_{{}^L\!G,F}^\bO$.
\end{remark}

As mentioned earlier, it is important not to impose Frobenius semisimplicity in the definition of $\locsys_{{}^cG, F}$, as this is not an algebraic condition when allowing Langlands parameters to vary in families. It turns out in each fiber of the map $\varpi_{{}^L\!G,F}$ (${}^L\!G$-version of the map in \eqref{E:LocMtoG}), Frobenius semisimple parameters do form a closed subspace. In fact, this space was originally introduced by Vogan. 

To explain this, it is convenient to recall the stack of Weil-Deligne parameters 
\[
\locsys^{\mathrm{WD}}_{{}^L\!G,F}=\locsys^{\mathrm{WD},\Box}_{{}^L\!G,F}/\hat{G}
\] 
over $\bQ(\sqrt{q})$. Here $\locsys^{\mathrm{WD},\Box}_{{}^L\!G,F}$ is a classical scheme classifying for every $\bQ(\sqrt{q})$-algebra $A$, the set of pairs $(h, X)$, where $h: W_F\to {}^L\!G(A)$ is a homomorphism and $X\in \mN_{\hat{G}}(A)$ is in the nilpotent cone of $\hat{G}$ such that 
\begin{itemize}
\item $d\circ h=\pr$;
\item  $\psi(I_F)$ has finite image; 
\item  $\Ad_{h(\ga)}X=||\ga||X$. 
\end{itemize}
Note that we consider the ${}^L\!G$-version of Weil-Deligne parameters here rather than the ${}^cG$-version considered in \cite[\textsection{3.1}]{zhu2020coherent}. But  since we fix $\sqrt{q}$, the two versions are equivalent (by a similar reasoning as in \Cref{rem: cG vs LG}). 

We shall also write
\[
\locsys^{\mathrm{W}}_{{}^L\!G,F}=\locsys^{\mathrm{W},\Box}_{{}^L\!G,F}/\hat{G}
\] 
for the stack of Weil parameters, where $\locsys^{\mathrm{W},\Box}_{{}^L\!G,F}$ just classifies those representations $h: W_F\to {}^L\!G(A)$ as in the definition of Weil-Deligne parameters. Clearly, we have a projection
\[
\la^{\mathrm{WD}}: \locsys^{\mathrm{WD}}_{{}^L\!G,F}\to \locsys^{\mathrm{W}}_{{}^L\!G,F},\quad (h,X)\mapsto h,
\] 
with a section sending $h$ to $(h,0)$. Note that
there is a $\bG_m$-action on $\locsys^{\mathrm{WD}}_{{}^L\!G,F}$ by scaling $X$. Then $\la$ is $\bG_m$-equivariant, with $ \locsys^{\mathrm{W}}_{{}^L\!G,F}$ equipped with the trivial $\bG_m$-action. In addition, the above section $ \locsys^{\mathrm{W}}_{{}^L\!G,F}\to  \locsys^{\mathrm{WD}}_{{}^L\!G,F}$ realizes $ \locsys^{\mathrm{W}}_{{}^L\!G,F}$ as fixed point loci of the $\bG_m$-action.

Recall that after choosing $\iota: \Ga_q\to W_F^t$ as before, there are isomorphisms of stacks over $\La=\overline\bQ_\ell$
\[
\locsys_{{}^L\!G,F}\simeq \locsys_{{}^L\!G,F,\iota}\simeq \locsys^{\mathrm{WD}}_{{}^L\!G,F}.
\]
We refer to \cite[\textsection{3.1}]{zhu2020coherent} for the chain of isomorphism. At the level of $\La=\overline\bQ_\ell$-points,
the isomorphism sends $\varphi$ to $(h,X)$ where 
\[
h(\ga)=\breve\varphi^{ss}(\ga) \mbox{ for }\ga\in I_F, \quad h(\iota(\sigma))=\varphi(\iota(\sigma)), \quad
X=\log(\breve\varphi^u(\iota(\tau))).
\]
This isomorphism induces isomorphisms of ring of functions
\[
H^0\rg(\locsys^{\mathrm{W}}_{{}^L\!G,F},\mO)=H^0\rg(\locsys^{\mathrm{WD}}_{{}^L\!G,F},\mO)=H^0\rg(\locsys_{{}^L\!G,F},\mO),
\]
which is independent of the choice of $\iota$. Thus the map $\varpi_{{}^L\!G, F}: \locsys_{{}^L\!G,F}\to \Spf Z_{{}^L\!G, F}$ factors as
\[
 \locsys_{{}^L\!G,F}\simeq \locsys_{{}^L\!G,F}^{\mathrm{WD}}\to  \locsys_{{}^L\!G,F}^{\mathrm{W}}\to \Spf Z_{{}^L\!G, F}.
\]
Note that the map $\locsys_{{}^L\!G,F}\to  \locsys_{{}^L\!G,F}^{\mathrm{W}},\ \varphi\mapsto h$ is independent of the choice of $\iota$. We denote this map by $\la$.

There is a similar notion of Frobenius semisimple Weil-Deligne (resp. Weil) parameters, and giving a Weil-Deligne (resp. Weil) parameter $(h,X)$ (resp. $h$), there is its Frobenius semisimplification $(h^{F\mbox{-}ss},X)$ (resp. $h^{F\mbox{-}ss}$).
Clearly if a Langlands parameter $\varphi$ matches a Weil-Deligne parameter $(h,X)$ under the above isomorphism, then $\varphi^{F\mbox{-}ss}$ matches $(h^{F\mbox{-}ss},X)$.

Now we fix a $\La$-point $z$ of $\Spf Z_{{}^L\!G,F}$, giving a strongly continuous completely reducible representation $h:W_F\to {}^L\!G(\La)$. This can be regarded as a $\La$-point of $\locsys_{{}^L\!G,F}$. But the corresponding Weil-Deligne representation is just $(h,0)$, and therefore its image in $\locsys_{{}^L\!G,F}^{\mathrm{W}}$ is given by the same representation $h$. 

The corresponding map $\{h\}/C_{\hat{G}}(h)\to \locsys_{{}^L\!G,F}^{\mathrm{W}}$ is a closed embedding. Let
$\Vogan^{\mathrm{WD}}_h=((\la^{\mathrm{WD}})^{-1}(h))_{\red}\subset  \locsys_{{}^L\!G,F}^{\mathrm{WD}}$ 
be the reduced fiber of $h$ for the map $\la^{\mathrm{WD}}$, which is a closed substack of $\locsys_{{}^L\!G,F}^{\mathrm{WD}}$. As closed substack in $\locsys_{{}^L\!G,F}^{\mathrm{WD}}$, we have
\[
\Vogan^{\mathrm{WD}}_h\cong \hat\frakg^{h(I_F)=1,h(\sigma)=q}/C_{\hat{G}}(h).
\]
(Note that elements in $\hat\frakg^{h(I_F)=1,h(\sigma)=q}$ are automatically nilpotent.) Similarly, we write
\begin{equation}\label{eq: vogan stack2}
\Vogan_h:=(\la^{-1}(h))_{\red}\subset \varpi_{{}^L\!G,F}^{-1}(z)\subset \locsys_{{}^L\!G,F},
\end{equation}
with both inclusions being closed embeddings. We have 
\[
\Vogan_h\cong (\hat{\frakg}\otimes \bZ_\ell(-1))^{h(W_F)}/C_{\hat{G}}(h).
\]

Using the fact that a Weil parameter $h:W_F\to {}^L\!G(\La)$ is completely reducible if and only if it is Frobenius semisimple, we obtain the following description of ($\La$-points of) $\Vogan_h$.
\begin{lemma}
The $\La$-points of the stack $\Vogan_h$ 
consist of Frobenius-semisimple representations $\varphi$ such that $\varpi_{{}^LG,F}(\varphi)=z$.
\end{lemma}

\begin{remark}
We note that the inclusion $\Vogan_h\subset \varpi_{{}^L\!G,F}^{-1}(z)$ is strict in general. This can be easily seen from \Cref{rem: non-frob semisimple parameters}. We thank Teruhisa Koshikawa for warning us this subtlety.
\end{remark}

\begin{remark}
Note that by definition $V_h$ is smooth. But the closed embedding $V_h\to \locsys_{{}^cG,F}$ is not of finite tor amplitude. It is not difficult to write down a derived enhancement of $V'_h$ so that  $V'_h\to  \locsys_{{}^cG,F}$ becomes quasi-smooth.
\end{remark}

\begin{remark}\label{rem: SL2 version of WD parameter}
In literature, people considers another form of Frobenius-semisimple Weil-Deligne parameters, which are  homomorphisms $\psi: W_F\times \SL_2\to {}^L\!G(\La)$ such that 
\begin{itemize}
\item $d\circ \psi|_{W_F}=\pr$;
\item $\psi(I_F)$ has finite image, $\psi(\sigma)$ is semisimple for one (and any) choice of lifting of the Frobenius; 
\item $\psi|_{\SL_2}:\SL_2\to \hat{G}$ is algebraic. 
\end{itemize}
It is well-known the two versions of Frobenius-semisimple Weil-Deligne parameters are equivalent. Namely, given $\psi$, one can construct $(h,X)$ as 
\begin{equation}\label{eq: WD SL2 parameter to WD parameter}
 h(\ga)= \psi(\ga, \begin{pmatrix} \|\ga\|^{\frac{1}{2}} & \\ & \|\ga\|^{-\frac{1}{2}}\end{pmatrix}) \ \mbox{for} \ \ga\in W_F, \quad X=d(\psi|_{\SL_2})(\begin{pmatrix}0& 1 \\  0 & 0 \end{pmatrix}).
\end{equation}
In fact, this construction $\psi\mapsto (h,X)$ does not make use of semisimplicity of $\psi(\sigma)$. However, when $(h,X)$ is Frobenius-semisimple,
this process can be reverse  by Jacobson-Morozov's lemma.

Clearly, we have $C_{\hat{G}}(\psi)\subset C_{\hat{G}}(h,X)$ but the inclusion might be strict. 
In fact, the neutral connected component of $C_{\hat{G}}(\psi)$ is always reductive but this may not be the case for $C_{\hat{G}}(h,X)$. But one knows that $\pi_0(C_{\hat{G}}(\psi))=\pi_0(C_{\hat{G}}(h,X))$ and $C_{\hat{G}}(\psi)^\circ \subset C_{\hat{G}}(h,X)^\circ$ is a Levi subgroup.
\end{remark}

\subsubsection{Discrete parameters}\label{SS: discrete parameter}
As an application of previous discussions, we study the geometry $\locsys_{{}^cG,F}$ around (essentially) discrete Langlands parameters. The materials here will be used in \Cref{SSS: cllc for supercuspidal} to study the categorical local Langlands correspondence for the supercuspidal representations.

Assume that $\La$ is an algebraically closed field (but not necessarily of characteristic zero at the moment).  

\begin{lemma}\label{lem:discrete parameter condition}
Let $\varphi: W_F\to {}^c G(\La)$ be a point in $\locsys_{{}^cG,F}(\La)$. The following are equivalent.
\begin{enumerate}
\item\label{lem:discrete-1} $H^0(W_F,\Ad_\varphi^0)=H^2(W_F,\Ad_\varphi^0)=0$.
\item\label{lem:discrete-2} The tangent complex of $\locsys_{{}^cG,F}$ at $\varphi$ is trivial. 
\item\label{lem:discrete-3} $\varphi$ is an open smooth point in $\locsys_{{}^cG,F}$. 
\end{enumerate}
When $\La$ is of characteristic zero, of characteristic $\ell$ with $\ell$ good for $\hat{G}$, these conditions are in addition equivalent to
\begin{enumerate}[resume]
\item\label{lem:discrete-4} The eigenvalues of the linear operator $\varphi(\sigma): \hat{\frakg}^{\varphi(I_F)}\to \hat{\frakg}^{\varphi(I_F)}$ does not contain  $1, q^{-1}$. (Here recall $\sigma$ is a lifting of arithmetic Frobenius.)
\end{enumerate}
\end{lemma}
\begin{proof}
The cohomology of  $\bT_\varphi\locsys_{{}^cG,F}$ at $\varphi$ are given by $H^i(W_F,\Ad_\varphi^0)$. So clearly \eqref{lem:discrete-2} implies \eqref{lem:discrete-1}. The converse follows from the fact that the Euler characteristic of $\bT_\varphi\locsys_{{}^cG,F}$ is zero. In addition, \eqref{lem:discrete-2} implies that $\varphi: \pt/Z_{\hat{G}}(\varphi)$ is smooth and the morphism $\varphi: \pt/Z_{\hat{G}}(\varphi)\to \locsys_{{}^cG,F}$ is an \'etale monomorphism, and therefore is an open embedding. Conversely, if $\pt/Z_{\hat{G}}(\varphi)\to \locsys_{{}^cG,F}$ is open and smooth, then $H^0(W_F,\Ad_\varphi^0)=H^2(W_F,\Ad_\varphi^0)=0$.

Finally for \eqref{lem:discrete-4},  clearly $H^0(W_F, \Ad_\varphi^0)=0$ is equivalent to the invertibility of $\varphi(\sigma)-1: \hat{\frakg}^{\varphi(I_F)}\to \hat{\frakg}^{\varphi(I_F)}$. On the other hand,
$H^2(W_F,\Ad_\varphi^0)=0$ is equivalent to $H^0(W_F,(\Ad_\varphi^0)^*(1))=0$, which in turn is equivalent to the invertibility of $q\varphi^*(\sigma)-1: (\hat{\frakg}^*)^{\varphi(I_F)}\to (\hat{\frakg}^*)^{\varphi(I_F)}$ .  
Now one uses the $\hat{G}\rtimes \mathrm{Out}(\hat{G})$-equivariant isomorphism $\hat\frakg\cong\hat\frakg^*$ to conclude.
\end{proof}

We call $\varphi$ a discrete parameter if the above equivalent conditions hold. Note that the space $(\hat{\frakg}^{\varphi(I_F)})^{\varphi(\sigma)=1}$ always contains the Lie algebra of $Z_{\hat{G}}^{\Ga_{\widetilde{F}/F}}$. Therefore a necessary condition for the existence of discrete parameter is that $Z_{\hat{G}}^{\Ga_{\widetilde{F}/F}}$ is finite \'etale over $\La$, which is restrictive. 
For this reason, we relax the condition.
\begin{definition}
A point $\varphi: W_F\to {}^cG(\La)$ is called an essentially discrete (local Langlands) parameter if $C_{\hat{G}}(\varphi)/Z_{\hat{G}}^{\Ga_{\widetilde{F}/F}}$ is finite.
\end{definition}

\begin{remark}
Note that $\varphi$ being essentially discrete is equivalent to requiring $C_{\hat{G}}(\varphi)/C_{{}^cG}$ is finite, where $C_{{}^cG}$ is the maximal subtorus of $Z_{\hat{G}}^{\Ga_{\widetilde F/F}}$ 
as in \eqref{eq: connected center of dual group}. 

If $\La$ is a field of characteristic zero, this is further equivalent to $H^0(W_F,\Ad_\varphi^0)=\mathrm{Lie}\ C_{{}^cG}$. In this case, we will see that $H^2(W_F,\Ad_\varphi^0)=0$ so $\varphi$ is a smooth point in $\locsys_{{}^cG,F}$.
\end{remark}

From now on we assume that $\La=\overline\bQ_\ell$ and fix $\sqrt{q}\in \overline\bQ_\ell$.
Our goal is to describe some geometry of irreducible components containing essentially discrete parameters. We start with some basic facts about these parameters.

Note that if we let $\varphi\mapsto \widetilde\varphi$ be the correspondence between points on $\locsys_{{}^L\!G,F}$ and on $\locsys_{{}^cG,F}$ as in \Cref{rem: cG vs LG}. Then $C_{\hat{G}}(\varphi)$ and $C_{\hat{G}}(\widetilde\varphi)$ are conjugate by $2\rho(\sqrt{q})$. Therefore, we will work with ${}^L\!G$ instead of ${}^cG$.
As mentioned before, homomorphisms $W_F\to {}^L\!G(\La)$ corresponding to points on $\locsys_{{}^L\!G,F}$ may not be Frobenius semisimple in general. But this is not a concern for essentially discrete parameters. (Of course the whole parameter may not be semisimple.)

\begin{lemma}
Let $\varphi: W_F\to {}^L\!G$ be an essentially discrete parameter. Then $\varphi=\varphi^{F\mbox{-}ss}$ is Frobenius semisimple.
\end{lemma}
\begin{proof}
By \Cref{lem: Frob semisimplification} and \Cref{rem: non-frob semisimple parameters}, we have a Frobenius semisimplification $\varphi^{F\mbox{-}ss}$ of $\varphi$ and a unipotent element $v:=\varphi(\sigma)^u$ commuting with $\varphi^{F\mbox{-}ss}$. Let $M=C_{\hat{G}}(\varphi^{F\mbox{-}ss})$. Then $v\in M$, and $C_{\hat{G}}(\varphi)=C_M(v)$. If $\dim M/Z_{\hat{G}}^{\Ga_{\widetilde F/F}}>0$, then
$\dim C_M(v)/Z_{\hat{G}}^{\Ga_{\widetilde F/F}}>0$. This shows that if $\varphi$ is essential discrete, then $M/Z_{\hat{G}}^{\Ga_{\widetilde F/F}}$ is finite so $v=1$. I.e. $\varphi=\varphi^{F\mbox{-}ss}$.
\end{proof}

To further study essentially discrete parameters, it is convenient to consider the associated Weil-Deligne representation.
 So in the sequel we will fix $\iota: \Ga_q\to W_F$.
Let $\varphi: W_F\to {}^L\!G$ be a $\La$-point of $\locsys_{{}^L\!G,F}$, let $(h,X)$ be the associated Weil-Deligne parameter. We suppose that $\varphi$ is Frobenius semisimple so $h$ is Frobenius semisimple. Then let
$\psi: W_F\times \SL_2\to {}^L\!G$ be the associated representation as in \Cref{rem: SL2 version of WD parameter}. 

\begin{lemma}\label{cor: purity of discrete parameters}
\begin{enumerate}
\item\label{cor: purity of discrete parameters-1} Suppose that $\varphi$ is essentially discrete. Then the inclusion $C_{\hat{G}}(\psi)\subset C_{\hat{G}}(h,X)$ is an isomorphism. In addition, the group $\psi(W_F)\subset {}^L\!G(\La)$ is finite modulo $Z_{\hat{G}}^{\Ga_{\widetilde F/F}}$. 
\item\label{cor: purity of discrete parameters-2} If $C_{\hat{G}}(\psi)/Z_{\hat{G}}^{\Ga_{\widetilde F/F}}$ is finite, then $C_{\hat{G}}(\psi)=C_{\hat{G}}(\varphi)$ and $\varphi$ is essentially discrete.
\end{enumerate}
\end{lemma}
It follows that there Weil-Deligne representation associated to an essentially discrete parameter is pure (in the sense of weights). 
\begin{proof}
For the first statement of Part \eqref{cor: purity of discrete parameters-1}, just notice that as $C_{\hat{G}}(\varphi)^\circ$ is reductive, we must have $C_{\hat{G}}(\psi)=C_{\hat{G}}(r,X)=C_{\hat{G}}(\varphi)$.

Now let $\sigma$ a lifting of the Frobenius of $W_F$.
Let $A$ be the algebraic group generated by $\psi(\sigma,1)\in {}^L\!G(\La)$. Then the neutral connected component $A^\circ\subset \hat{G}$ is a torus, normalizing $\psi(I_F)$. Therefore, $A^\circ$ is in the neutral connected component of the center of $C_{\hat{G}}(\psi)$.
Therefore $A^\circ\subset Z_{\hat{G}}^{\Ga_{\widetilde F/F}}$, and $\psi(W_F)$ is finite modulo $Z_{\hat{G}}^{\Ga_{\widetilde F/F}}$. This finishes the proof of Part \eqref{cor: purity of discrete parameters-1}.

For Part \eqref{cor: purity of discrete parameters-2}, we first recall the relation between $(h,X)$ and $\psi$ given in \eqref{eq: WD SL2 parameter to WD parameter}.
The argument of \Cref{cor: purity of discrete parameters} implies that $\psi|_{W_F}$ has finite image. Therefore $h(\ga)=\psi(1, \begin{pmatrix}||\ga||^{\frac{1}{2}} & \\ & ||\ga||^{-\frac{1}{2}}\end{pmatrix})$ for $\ga$ in a finite index subgroup of $W_F$. This implies that if $g\in C_{\hat{G}}(\varphi)$, then $g$ centralizes $\psi(\bG_m)$, where $\bG_m$ is the standard diagonal torus of $\SL_2$. It follows that $C_{\hat{G}}(\varphi)\subset C_{\hat{G}}(\psi)$ and the lemma follows.
\end{proof}

Recall the action of $R_{W_F, {}^cG}$ on $\locsys_{{}^cG,F}$ and on $\Spf Z_{{}^cG,F}$ from \eqref{eq: center action on loc}.
Note that the subspace of $R_{W_F, {}^cG}$ consisting of those $\psi$ that is trivial one $I_F$ is identified with $C_{{}^cG}$ by sending $\psi$ to $\psi(\sigma)$. It follows that we obtain a free action of $C_{{}^cG}$ on $\locsys_{{}^cG,F}$ and on $\Spf Z_{{}^cG, F}$.

\begin{proposition}\label{prop: geometry of discrete component}
Let $\varphi$ be an essentially discrete parameter. Let $(h,X)$ be the associated Weil-Deligne parameter. Then $\locsys_{{}^L\!G,F}$ is smooth at $\varphi$. In addition, the action map $j_{\varphi}: C_{{}^cG}\times \{\varphi\}/C_{\hat{G}}(\varphi)\to \locsys_{{}^L\!G,F}$ is an open embedding, which factors as
\[
\xymatrix{
&C_{{}^cG}\times \{\varphi\}/C_{\hat{G}}(\varphi)\ar@{^{(}->}[dl]\ar^-{j_\varphi}[dr]&\\
C_{{}^cG}\times V_h \ar@{^{(}->}[rr]&& \locsys_{{}^L\!G,F}
}
\]
with $C_{{}^cG}\times \{\varphi\}/C_{\hat{G}}(\varphi)\to C_{{}^cG}\times V_h$ is an open embedding, and $C_{{}^cG}\times V_h$ is the (unique) irreducible component of $\locsys_{{}^L\!G,F}$ containing $\varphi$. In particular, it is smooth.
\end{proposition}
\begin{proof}
\Cref{cor: purity of discrete parameters} implies that the eigenvalues of $\varphi(\sigma)$ on $\hat\frakg^{\varphi(I_F)}$ are of the form $\al \sqrt{q}^i$, where $\al$ is a root of unit and $i\geq 0$ is an integer. As argued in \Cref{lem:discrete parameter condition}, this implies that $H^2(W_F, \Ad_\varphi^0)=0$. In addition, we $\dim H^1(W_F, \Ad_\varphi^0)=\dim H^0(W_F,\Ad_\varphi^0)=\dim C_{{}^cG}$.
Clearly, the action map
$j_\varphi: C_{{}^cG}\times \{\varphi\}/C_{\hat{G}}(\varphi)\to \locsys_{{}^L\!G}$ is a monomorphism. In addition, it induces an isomorphism between tangent complexes. Therefore it is an open embedding.

Notice that $C_{{}^cG}$ also acts transitively on $\locsys^{\mathrm{W}}_{{}^L\!G,F}$ and $C_{{}^cG}\times \{h\}/C_{\hat{G}}(h)\to \locsys^{\mathrm{W}}_{{}^L\!G,F}$ is closed embedding. Then map $\la: \locsys_{{}^L\!G}\to \locsys^{\mathrm{W}}_{{}^L\!G}$ is $C_{{}^cG}$-equivariant so
\[
C_{{}^cG}\times V_h =\la^{-1}(C_{{}^cG}\times \{h\}/C_{\hat{G}}(h))_{\red}\to \locsys_{{}^L\!G,F}
\]
is a closed embedding. Clearly, $j_\varphi$ factors through $C_{{}^cG}\times \{\varphi\}/C_{\hat{G}}(\varphi)\to  C_{{}^cG}\times V_h$, and is open. 

Now it follows that $\dim (C_{{}^cG}\times V_h)=\dim \locsys_{{}^L\!G,F}$. Since $C_{{}^cG}\times V_h$ is smooth, and is closed in $\locsys_{{}^L\!G,F}$, and since $\locsys_{{}^L\!G,F}$ is reduced (over $\La=\overline\bQ_\ell$), we see that $C_{{}^cG}\times V_h$ is an irreducible component of $\locsys_{{}^L\!G,F}$. 
\end{proof}

\begin{corollary}\label{lem: extends ss parameter to discrete parameter}
Given a point $z\in \Spf Z_{{}^cG,F}(\overline\bQ_\ell)$, there at most one essentially discrete parameter $\varphi\in\locsys_{{}^L\!G,F}$ that maps to $z$.
\end{corollary}
\begin{proof}
Let $\varphi_i, \ i=1,2$ be two essentially discrete parameters that maps to $z$. Then their images in $\locsys_{{}^L\!G,F}^{\mathrm{W}}$ are the same, say $h$. Then $\{\varphi_i\}/\bB C_{\hat{G}}(\varphi_i)\subset V_h$ are two open subspaces. As $V_h$ is irreducible, we see that $\varphi_1$ and $\varphi_2$ give the same point in $\locsys_{{}^L\!G,F}$.
\end{proof}

We have seen that if an irreducible component of $\locsys_{{}^L\!G,F}$ contains an essential discrete parameter, then is image in $\Spf Z_{{}^L\!G,F}$ is a single  $C_{{}^cG}$-orbit. It turns out the converse is also true.
\begin{lemma}\label{lem: support of supercusp in one component}
Let $Z\subset \locsys_{{}^L\!G,F}$ be an irreducible component that maps to a single $C_{{}^cG}$-orbit in $\Spf Z_{{}^L\!G,F}$. Then $Z$ contains an essential discrete parameter $\varphi$.
\end{lemma}
\begin{proof}
Suppose $Z\subset \overline{\locsys_{{}^L\!G,F}^{\bO}}$ for some point $\bO$ in $\locsys_{{}^L\!G,\breve F}$. Then $Z\cap \locsys_{{}^L\!G,F}^{\bO}$  is dense in one connected component of
 $\locsys_{{}^L\!G,F}^{\bO}$. We choose $\varphi: W_F\to {}^L\!G$ representing a point of $Z\cap\locsys_{{}^cG,F}^{\bO}$. As explained in \Cref{rem: non-frob semisimple parameters}, $\varphi^{F\mbox{-}ss}$ is also a point on $Z\cap\locsys_{{}^L\!G,F}^{\bO}$. There, we have a point $\varphi\in Z\cap\locsys_{{}^L\!G,F}^{\bO}$ such that $\varphi(\sigma)$ is semisimple. 
 Then we have a morphism $M /M \to Z\cap \locsys^{\bO}_{{}^cG,F}$ sending $m$ to $\varphi_m\varphi$, where $\varphi_m:W_F\to \sigma^\bZ\to M$ is the homomorphism sending $\sigma$ to $m$. 
 Let $C_{\hat{G}}(\psi)\subset M$ as before.
If $C_{\hat{G}}(\psi)^\circ/C_{{}^cG}$ is non-trivial, then the image of the map $C_{\hat{G}}(\psi)^\circ/\!\!/C_{\hat{G}}(\psi)\to \Spf Z_{{}^cG,F}$ cannot be a single $C_{{}^cG}$-orbit. Contradiction. Thus $C_{\hat{G}}(\psi)^\circ/C_{{}^cG}$ is trivial. By \Cref{cor: purity of discrete parameters} \eqref{cor: purity of discrete parameters-2},
$\varphi$ is essential discrete.
\end{proof}

\subsection{The stack of tame and unipotent Langlands parameters}

\subsubsection{The stack of tame Langlands parameters}
We assume that $G$ is tamely ramified so $\widetilde F/F$ is a tame extension. We will allow $\La$ to be general.
Then we have an open and closed quasi-compact substack $\locsys^{\tame}_{{}^cG,F}\subset \locsys_{{}^cG,F}$ classifying those representations of $W_F$ that factors through the tame Weil group $W_{F}^{t}$.
We similarly have $\locsys^{\tame}_{{}^cG,\breve F}\subset \locsys_{{}^cG,\breve F}$. Then \eqref{eq: res from arith parameter to geom parameter} restricts to a morphism
\begin{equation}\label{eq: res from tame arith parameter to geom parameter}
\res^{\tame}: \locsys^{\tame}_{{}^cG,F}\to \locsys^{\tame}_{{}^cG,\breve F},
\end{equation}
and the isomorphism from \Cref{lem: arith parameter as fix point of geometric parameter} restricts to an isomorphism
\begin{equation}\label{eq: tame arith parameter as fixed point}
 \locsys^{\tame}_{{}^cG,F}\cong \mL_\phi(\locsys^{\tame}_{{}^cG,\breve F}).
\end{equation}

We fix an embedding $\iota: \Ga_q\to W_F^t$ such that $\iota(\tau)$ is a generator of the tame inertia and $\iota(\sigma)$ is a lifing of the Frobenius as before.
Then we similarly have open and closed substack $\locsys^{\tame}_{{}^cG,F,\iota}\subset \locsys_{{}^cG,F,\iota}$ over $\bZ[1/p]$, which can be described explicitly as follows:
Let $\bar\tau$ and and $\bar\sigma$ be the images of $\tau$ and $\sigma$ under the projection $\Ga_q\xrightarrow{\iota} W_F^t\xrightarrow{\widetilde{\pr}} \bG_m\times \Ga_{\widetilde F/F}$. (Note that $\bar\tau$ is trivial on the $\bG_m$-factor.) 
Then
\begin{equation}\label{eq:presentation-tame-stack}
\locsys^{\tame}_{{}^cG,F,\iota}\simeq \left\{(h,g)\in \hat{G}\bar\tau\times \hat{G}\bar\sigma \subset {}^cG\times {}^cG \mid ghg^{-1}=h^{q}\right\}/\hat{G}.
\end{equation}
We can similar define $\locsys^{\tame}_{{}^cG,\breve F,\iota}$, where we replace $\Ga_q$ by $\tau^{\bZ[1/p]}\subset \Ga_q$. If we let $\tau_i=\sigma^{-i}\tau\sigma^i\in \tau^{\bZ[1/p]}$, then
$\locsys^{\tame}_{{}^cG,\breve F,\iota}=\lim_i\hat{G}\bar\tau_i/\hat{G}$, with the transitioning maps given by $q$-power map, and $\iota$-version of \eqref{eq: res from tame arith parameter to geom parameter} is the map
\begin{equation}\label{eq:restr-to-inertia}
    \res: \locsys^{\tame}_{{}^cG,F,\iota}\to \locsys^{\tame}_{{}^cG,\breve F,\iota}=\lim_i\hat{G}\bar\tau_i/\hat{G}\to \hat{G}\bar\tau/\hat{G},
\end{equation}
which explicitly is the map sending $(h,g)$ to $h$ (up to $\hat{G}$ conjugacy). Similarly, there are $\iota$-version $Z^\tame_{{}^cG, F,\iota}$ and $Z^\tame_{{}^cG, \breve F,\iota}$, which are $\bZ[1/p]$-(pro)algebras.

The $\iota$-version of \eqref{eq: phi automorphism of loc} is explicitly given by
\begin{equation}\label{eq:endo-tau-of-hatG}
\phi: \lim_i\hat{G}\bar\tau_i/\hat{G}\to \lim_i\hat{G}\bar\tau_i/\hat{G},\quad  g_i\mapsto \bar\sigma(g_{i+1}),
\end{equation}
where we recall $\bar\sigma$ sends $\hat{G}\bar\tau_{i+1}$ to $\hat{G}\bar\tau_i$. 
Then
we have $\iota$-version of \eqref{eq: tame arith parameter as fixed point}
\[
\locsys^{\tame}_{{}^cG,F,\iota}\cong \mL_{\phi}(\locsys^{\tame}_{{}^cG,\breve F,\iota}).
\]
It is also convenient to consider the inverse map of \eqref{eq:endo-tau-of-hatG}, which induces an endomorphism of $\hat{G}\bar\tau$
\begin{equation}\label{eq: map [q]}
\hat{G}\bar\tau\to\hat{G}\bar\tau,\quad g\mapsto \bar\sigma^{-1}(g^q),
\end{equation} 
which induces a map $\hat{G}\bar\tau/\!\!/\hat{G}\to \hat{G}\bar\tau/\!\!/\hat{G}$ still denoted by $[q]$. Let $(\hat{G}\bar\tau/\!\!/\hat{G})^{[q]}$ be $[q]$-fixed point subscheme of $\hat{G}\bar\tau/\!\!/\hat{G}$.

Let $\hat{S}=\hat{T}/(1-\bar\tau)\hat{T}$ be the $\bar\tau$-coinvariants of $\hat{T}$. Then the action of $\bar\sigma$ on $\hat{T}$ induces an action on $\hat{S}$, still denoted by $\bar\sigma$. Note that the morphism $[q]$ then induces a morphism of $\hat{S}$, still denoted by $[q]$.
Let $W_0=W^{\bar\tau}$ be the $\bar\tau$-invariants of the absolute Weyl group $W$ of $\hat{G}$, which also acts on $\hat{S}$. Recall that we have the Chevalley restriction isomorphism $\hat{G}\bar\tau/\!\!/\hat{G}\cong \hat{S}/\!\!/W_0$ (e.g.  see \cite[Proposition 4.2.3]{xiao.zhu} in this generality). Then we the following identification of $\bZ[1/p]$-schemes.
\[
(\Spf Z^\tame_{{}^cG, \breve F,\iota})^\phi\cong(\hat{G}\bar\tau/\!\!/\hat{G})^{[q]}\cong (\hat{S}/\!\!/W_0)^{[q]}.
\]
Note that they are all finite over $\bZ[1/p]$, which is consistent with \Cref{rem: 3rd def of inertia type}.

The stack $\locsys^{\tame}_{{}^cG,F}$ (and $\locsys^{\tame}_{{}^cG,F,\iota}$) still breaks into connected components. Now we study some components over an algebraically closed field.

\begin{definition}
An inertia type $\zeta$ of ${}^cG$ is called tame if it is a $\La$-point of $(\Spf Z^\tame_{{}^cG, \breve F})^\phi$. Here, we denote $Z^\tame_{{}^cG, \breve F}=H^0\rg(\locsys_{{}^cG,\breve F}^\tame,\mO)$.
\end{definition}

So after fixing a choice of $\iota: \Ga_q\to W_F^t$, tame inertia types can be identified with the subset 
\[
(\hat{S}/\!\!/W_0)^{[q]}(\La)\subset (\hat{S}/\!\!/W_0)(\La)=\hat{S}(\La)/\!\!/W_0,
\] 
consisting of $W_0$-orbits of those $\chi\in\hat{S}(\La)$ such that there is some $w\in W_0$ such that $w(\bar\sigma(\chi))=\chi^q$.
But we can reinterpret such identification without referring a choice of $\iota$ as follows. (Note that $\hat{S}, W_0$ and the action of $\bar\sigma$ on $\hat{S}$ are canonically defined independent of the choice of $\iota$.)

\begin{lemma}\label{lem: classification tame inertia type}
There is a bijection between tame inertia type and finite order homomorphism $\chi: I_F^t\to \hat{S}(\La)$ up to $W_0$-conjugacy such that there is some $w\in W_0$ such that $w(\bar\sigma(\chi))=\chi^q$.
\end{lemma}

In the sequel, for a tame inertia type $\zeta$, we will $\Xi(\zeta)$ denote the set of finite order homomorphisms $\chi: I_F^t\to \hat{S}(\La)$ corresponding to $\zeta$. Then $W_0$ acts transitively on $\Xi(\zeta)$.

\subsubsection{The stack of unipotent Langlands parameters}
We look more closely into the part of the stack corresponding to the unipotent inertia type
\[
\zeta=\unip,
\] 
by which we mean $\bar\tau=1$\footnote{This requirement is not really necessary. We refer to \cite{zhu2020coherent} for some discussions when $\bar\tau\neq 1$. But when $\bar\tau\neq 1$, such stack does not really parameterize Langlands parameters with unipotent monodromy so it would be a little bit awkward to call it the stack of unipotent Langlands parameters. In addition, as we shall see soon in \Cref{prop: reducing to endoscopic group}, the study of tame inertia types can be more or less reduced to the study of case $\bar\tau=1$ and $\zeta=\unip$.}  and the corresponding homomorphism $\chi: I_F^t\to \hat{S}$ is trivial. In this case, we let $\locsys^{\widehat\unip}_{{}^cG, F}$ denote the corresponding stack. As mentioned in \Cref{rem: locsyszeta comment}. Namely, we can regard $\zeta=\unip$ as a $\La$-point of $\Spf Z_{{}^cG,\breve F}$ and let 
\[
\locsys^{\unip}_{{}^cG, F}=\locsys^{\tame}_{{}^cG, F}\times_{\Spf Z^{\tame}_{{}^cG,\breve F}} \{\unip\}= \locsys^{\widehat\unip}_{{}^cG, F}\times_{\widehat\unip}\{\unip\}.
\]
This is a closed substack of $\locsys^{\widehat\unip}_{{}^cG, F}$.

Note that implicitly in the definition, $\locsys^{\unip}_{{}^cG, F}\subset \locsys^{\widehat\unip}_{{}^cG, F}$ are stacks over an algebraically closed field $\La$ (due to our definition of inertia type). But in fact,
both $\locsys^{\unip}_{{}^cG, F}\subset \locsys^{\widehat\unip}_{{}^cG, F}$ are defined over $\bZ_\ell$, and even 
admit $\iota$-version defined over $\bZ[1/p]$. Let $\hat{G}\to \hat{G}/\!\!/\hat{G}$ the Chevalley map and let $\{1\}\in \hat{G}/\!\!/\hat{G}$ be the image of the unit of $\hat{G}$ in $\hat{G}/\!\!/\hat{G}$. Let $\widehat{\{1\}}$ be the formal completion of $\{1\}$ in $\hat{G}/\!\!/\hat{G}$. These (formal) schemes are defined over $\bZ[1/p]$ (in fact over $\bZ$). Therefore, if we fix $\iota: \Ga_q\to W_F^t$ as before, we can define stacks over $\bZ[1/p]$
\begin{equation}\label{eq: unipotent stack}
\locsys^{\unip}_{{}^cG,F,\iota}=\locsys^{\tame}_{{}^cG,F,\iota}\times_{\hat{G}/\!\!/\hat{G}}\{1\}, \quad \locsys^{\widehat\unip}_{{}^cG,F,\iota}=\locsys^{\tame}_{{}^cG,F,\iota}\times_{\hat{G}/\!\!/\hat{G}}\widehat{\{1\}},
\end{equation}
whose base change to $\bZ_\ell$ give promised stacks of unipotent Langlands parameters canonically defined independent of the choice of $\iota$.
Note that by definition, $\locsys^{\unip}_{{}^cG,F,\iota}$ is an algebraic stack of almost of finite presentation over $\bZ[1/p]$, which in general have non-trivial derived structure, while $\locsys^{\widehat\unip}_{{}^cG,F,\iota}$ is classical but is in general an ind-algebraic stack over $\bZ[1/p]$. But the base change of $\locsys^{\widehat\unip}_{{}^cG,F,\iota}$ to a field is always classical and algebraic (although may not be reduced) by \Cref{lem: algebraicity of zeta component}. 
In addition, by \Cref{lem: phi fixed point for Z} we have 
\[
\locsys^{\widehat\unip}_{{}^cG,F,\iota}\cong \mL_\phi(\hat{G}/\hat{G}\times_{\hat{G}/\!\!/\hat{G}}\widehat{\{1\}}).
\]

\begin{remark}\label{rem: three versions of unipotent stack}
Let $\mU_{\hat{G}}\subset \hat{G}$ be the variety of unipotent elements in $\hat{G}$. This is the reduced subscheme of the (possibly derived) scheme $\{1\}\times_{\hat{G}/\!\!/\hat{G}}\hat{G}$. Let $\widehat{\mU}_{\hat{G}}$ be its formal completion in $\hat{G}$. 
Then we have
\begin{equation}\label{eq: three versions of unipotent}
\mU_{\hat{G}}\subset \hat{G}\times_{\hat{G}/\!\!/\hat{G}}\{1\}\subset  \hat{G}\times_{\hat{G}/\!\!/\hat{G}}\widehat{\{1\}}=\widehat{\mU}_{\hat{G}}.
\end{equation}
Thus, there is a variant of $\locsys^{\unip}_{{}^cG,F,\iota}$ defined as $'\locsys^{\unip}_{{}^cG, F,\iota}=\locsys^{\tame}_{{}^cG,F,\iota}\times_{\hat{G}/\hat{G}}\mU_{\hat{G}}/\hat{G}$, which is a closed substack of $\locsys^{\unip}_{{}^cG,F,\iota}$. If the derived group of $\hat{G}$ is simply-connected, then 
the first inclusion in \eqref{eq: three versions of unipotent} is an isomorphism, and we have $'\locsys^{\unip}_{{}^cG,F,\iota}=\locsys^{\unip}_{{}^cG,F,\iota}$.
\end{remark}

To study $\locsys^{\widehat\unip}_{{}^cG, F}$ or its variants, we also need to recall some basic facts about unipotent (and nilpotent) elements and their centralizers in $\hat{G}$, when
$\La$ is an algebraically closed field (over $\bZ_\ell)$. Here are some ``standard hypothesis" on $\hat{G}$.
\begin{assumption}\label{ass: standard hypothesis for group in positive char}
Consider the following conditions for $\hat{G}$.
\begin{enumerate}
\item\label{ass: standard hypothesis for group in positive char-1} The characteristic $\ell$ is good for $\hat{G}$;
\item\label{ass: standard hypothesis for group in positive char-2} $\ell\nmid \sharp\pi_1(\hat{G}_{\mathrm{der}})$; 
\item\label{ass: standard hypothesis for group in positive char-3} There exists a $\hat{G}$-invariant non-degenerate bilinear form on $\hat\frakg$.
\end{enumerate}
\end{assumption}

It is known that under \Cref{ass: standard hypothesis for group in positive char} \eqref{ass: standard hypothesis for group in positive char-1},
there exists a $\hat{G}$-equivariant homeomorphism $\varepsilon\colon \mN_{\hat{G}}\to \mU_{\hat{G}}$ from the nilpotent cone of $\hat{g}$ to the unipotent variety of $\hat{G}$. If \Cref{ass: standard hypothesis for group in positive char} \eqref{ass: standard hypothesis for group in positive char-1} \eqref{ass: standard hypothesis for group in positive char-2} hold, then such $\varepsilon$ can be chosen to be an isomorphism. In any case,
we fix such $\varepsilon$. (Over a field of characteristic zero, $\varepsilon$ can be chosen to be the usual exponential map.) 
For $u\in \mU_{\hat{G}}$, let $X\in \mN_{\hat{G}}$ be the corresponding nilpotent element. It is known that under \Cref{ass: standard hypothesis for group in positive char} \eqref{ass: standard hypothesis for group in positive char-1}, $\ell\nmid\sharp\pi_0(C_{\hat{G}}(u))$ and under \Cref{ass: standard hypothesis for group in positive char} \eqref{ass: standard hypothesis for group in positive char-1}-\eqref{ass: standard hypothesis for group in positive char-3} $C_{\hat{G}}(u)$ is smooth (see \cite[5.10]{Jantzen}).

Recall that when $\La$ is a field of characteristic zero or characteristic $\ell$ large enough, the Jacobson-Morozov theorem implies that associated to $u$ there is a homomorphism $\SL_2\to \hat{G}$, unique up to conjugation by $C_{\hat{G}}(u)$, sending $\begin{pmatrix} 1 & 1 \\ & 1 \end{pmatrix}$ to $u$.
Under the generality we are considering, such $\SL_2$ may not exist.
However, there is a replacement. We say a cocharacter $\la: \bG_m\to \hat{G}$ is associated to $u$ if we write the grading induced by $\la$ as
\[
\hat\frakg=\bigoplus_i\hat\frakg_i,
\] 
then $X\in\frakg_2$ and the $C_{\hat{G}}(\la)$-orbit through $X$ is open dense in $\hat\frakg_2$. In particular $\Ad_{\la(t)}X=t^2X$, and $\la(\sqrt{q}) u\la(\sqrt{q})=u^q$. 
It is known that \Cref{ass: standard hypothesis for group in positive char} \eqref{ass: standard hypothesis for group in positive char-1}, such cocharacter exists and all such cocharacters are conjugate under $C_{\hat{G}}^\circ(u)$, the neutral connected component of the centralizer $C_{\hat{G}}(u)$ of $u$ in $\hat{G}$. In addition, the map from the set of unipotent conjugacy classes to the set of conjugacy classes of cocharacters of $\hat{G}$ is injective. In particular, there is a unique dominant cocharacter $\la$ (with respect to $(\hat{B},\hat{T})$) associated to the conjugacy class of $u$. 
Of course, if there is a homomorphism $\SL_2\to \hat{G}$ associated to $u$ as above, then the restriction of it to the standard diagonal torus $\bG_m\subset \SL_2$ is a cocharacter associated to $u$.

Next, let $\hat{P}_u$ be the attractor in $\hat{G}$ for the conjugation action of $\la(\bG_m)$ on $\hat{G}$. It is known that $C_{\hat{G}}(u)\subset\hat{P}_u$ (and stable under the conjugation action by $\la(\bG_m)$), and therefore $\hat{P}_u$ is independent of the choice of $\la$. It is called the canonical parabolic subgroup of $\hat{G}$ associated to $u$. Let $\hat{\frakp}_u$ and $\fraku_{\hat{\frakp}_u}$ be the Lie algebra of $\hat{P}_u$ and $U_{\hat{P}_u}$ respectively.
Then 
\[
\hat{\frakp}_u=\bigoplus_{i\geq 0}\hat{\frakg}_i,\quad \fraku_{\hat{\frakp}_u}=\bigoplus_{i> 0}\hat{\frakg}_i.
\]

Note that $C_{\hat{G}}(\la)\subset \hat{P}_u$ is a Levi subgroup. The Levi decomposition $\hat{P}_u= U_{\hat{P}_u}\rtimes C_{\hat{G}}(\la)$ induces a Levi decomposition
\begin{equation}\label{eq:Levi decomp of C(u)}
C_{\hat{G}}(u)= R_{\hat{G}}(u)\rtimes C_{\hat{G}}(\la,u),
\end{equation}
where $R_{\hat{G}}(u)=U_{\hat{P}_u}\cap C_{\hat{G}}(u)$ is the unipotent radical of $C_{\hat{G}}(u)$, and $C_{\hat{G}}(\la,u)=C_{\hat{G}}(\la)\cap C_{\hat{G}}(u)$ is isomorphic to the reductive quotient of $C_{\hat{G}}(u)$. (Here we assume that $C_{\hat{G}}(u)$ is smooth. Otherwise, one should replace $C_{\hat{G}}(u)$ by its reduced subgroup in the above discussions.)
In particular, 
\[
A(u):=\pi_0 C_{\hat{G}}(u)= \pi_0 C_{\hat{G}}(\la,u).
\]
In addition, it is known that every element in $A(u)$ can be lifted to a semisimple element in $C_{\hat{G}}(u)$.

We need some ``disconnected" version of the above discussions.

\begin{lemma}\label{lem: disconnected of C(u)}
Suppose $\widetilde{G}$ is an extension of a finite cyclic group $\langle c\rangle$ by $\hat{G}$. 
Let $u$ be a unipotent and let $\la$ be an associated cocharacter as above. Suppose that the conjugacy class of $u$ in $\hat{G}$ is stable under $\widetilde{G}$-conjugation (in which case we also say the conjugacy class of $u$ is $c$-stable). Then we have a short exact sequence
\[
1\to C_{\hat{G}}(\la,u)\to C_{\widetilde{G}}(\la,u)\to \langle c\rangle \to 1.
\]
In addition, we also have the decomposition
\begin{equation*}\label{eq:Levi decomp of C(u)-2}
C_{\widetilde{G}}(u)=R_{\hat{G}}(u)\rtimes C_{\widetilde{G}}(\la,u).
\end{equation*}
\end{lemma}
\begin{proof}
For the first statement, only surjectivity of $C_{\widetilde{G}}(\la,u)\to \langle c\rangle$ requires justification. 
By the assumption of $u$, we may choose $g_0\in \widetilde{G}$ such that $g_0$ maps to $c$ and $g_0 u g_0^{-1}=u$. Then $g_0\la g_0^{-1}$ is a cocharacter associated to $u$ as well. Then we may choose $h\in C_{\hat{G}}(u)$ such that $hg_0 \la (hg_0)^{-1}=\la$. Therefore, $hg_0\in C_{\widetilde{G}}(\la,u)$, which maps to $c$. 

The last statement follows from the first and \eqref{eq:Levi decomp of C(u)}.
\end{proof}

As a first application of the above facts, we can make the parameterization of irreducible components of $\locsys^{\widehat\unip}_{{}^cG,F}$ more explicit, when
$\La$ is an algebraically closed field (over $\bZ_\ell)$ of good characteristic (for $\hat{G}$). It is convenient to fix $\sqrt{q}$, and work with ${}^LG$ rather than ${}^cG$. (See \Cref{rem: cG vs LG}, and so in the sequel $\bar\sigma$ will denote the image of the arithmetic Frobenius in $\Ga_{\widetilde F/F}$ rather than in $\bG_m\times \Ga_{\widetilde F/F}$ as in \eqref{eq: barsigma}.)
We fix $\iota: \Ga_q\to W_F^t$, then $\locsys^{\widehat\unip}_{{}^cG, \breve{F}}\subset \hat{G}/\hat{G}$, so 
a point $\bO$ in  $\locsys^{\widehat\unip}_{{}^cG, \breve{F}}$ corresponds to a  unipotent conjugacy class of $\hat{G}$, and a choice $\breve\varphi$ amounts to choosing a unipotent element $u$ the conjugacy. In this case, $C_{\hat{G}}(\breve\varphi)=C_{\hat{G}}(u)$ and $A(\breve\varphi)=A(u)$. We shall also consider $C_{{}^L\!G}(u)$ and let $\widetilde{A(u)}:=\pi_0(C_{{}^L\!G}(u))=\pi_0(C_{{}^L\!G}(\la,u))$. That $\breve\varphi$ extends to a Langlands parameter means that this unipotent conjugacy class is $\bar\sigma$-stable.
Therefore by \Cref{lem: disconnected of C(u)} we have an exact sequence
\[
1\to A(u)\to \widetilde{A(u)}\to  \langle \bar\sigma\rangle\to 1.
\]
Let  $C_{{}^L\!G}(u)^1$ (resp. $\widetilde{A(u)}^1$) denote the preimage of $\bar\sigma$ in $C_{{}^L\!G}(u)$ (resp. in $\widetilde{A(u)}$).

\begin{proposition}\label{prop: irr. comp. of loc-unip and conj. class of Au}
Suppose $\La$ is an algebraically closed field of good characteristic for $\hat{G}$.
Then there is a canonical bijection between $\pi_0(\locsys^{\bO}_{{}^cG,F})$ and the quotient of $\widetilde{A(u)}^1$ by the conjugation action of $A(u)$. In particular, if $G$ is split, then  irreducible components of $\locsys_{{}^cG,F}^{\unip}$ are parameterized by pairs $(\bO, x)$ where $\bO$ is a unipotent conjugacy class in $\hat{G}$ and $x$ is a conjugacy class of $A(u)$ for some $u$ in $\bO$.
\end{proposition}
\begin{proof}
Let $\la:\bG_m\to\hat{G}$ be a cocharacter associated to $u$.
Recall that a lifting $\breve\varphi$ to a Langlands parameter $\varphi$ amounts to choose 
$g\bar\sigma\in \hat{G}\bar\sigma$ such that $g\bar\sigma(u)g^{-1}=u^q$. 
Then $x:=\la(\sqrt{q}^{-1})g\bar\sigma \in C_{{}^L\!G}(u)^1$. 
We thus obtain an isomorphism (of the underlying classical stacks)
\[
\locsys^{\bO}_{{}^cG, F}\cong \mL_\phi(\bO)\cong C_{{}^LG}(u)^1/\Ad_{\la(\sqrt{q})} C_{\hat{G}}(u), \quad (u,g\bar\sigma) \mapsto x=\la(\sqrt{q}^{-1})g\bar\sigma,
\]
where we recall $\la(\bG_m)$ acts on $C_{\hat{G}}(u)$ by conjugation (and the action of $\Ad_{\la(\sqrt{q})} C_{\hat{G}}(u)$ on $C_{{}^LG}(u)^1$ is given by $(h, x)\mapsto h x \la(\sqrt{q})h^{-1}\la(\sqrt{q}^{-1})$).

Clearly after taking $\pi_0$, the conjugation action of $\la(\bG_m)$ on $C_{\hat{G}}(u)$ becomes the trivial action on $A(u)$. It follows that
$\pi_0\locsys^{\bO}_{{}^cG,F}\cong \widetilde{A(u)}^1/A(u)$, as desired.
\end{proof}

\begin{example}\label{ex: unramified component}
Let $\bO$ denote the trivial unipotent conjugacy class of $\hat{G}$. More canonically, we denote it by
\[
\locsys_{{}^cG,\breve F}^{\mathrm{unr}}=\bB \hat{G},
\]
classifying the unramified (a.k.a. trivial) representation of $I_F^t$. Then 
\[
\locsys_{{}^cG,F}^{\mathrm{unr}}:=\mL_\phi\locsys_{{}^cG,\breve F}^{\mathrm{unr}}\cong \hat{G}\bar\sigma/\hat{G}
\]
is the substack of unramified parameters. Note that $\locsys_{{}^cG,F}^{\mathrm{unr}}$ exists over $\bZ_\ell$ (and there is an $\iota$-version over $\bZ[1/p]$).
It is smooth and is the reduced substack of an irreducible component of $\locsys^{\widehat\unip}_{{}^cG,F}$ (even integrally). 
\end{example}

\begin{proposition}\label{lem: connectedness of unipotent stack}
The stack $\locsys^{\widehat\unip}_{{}^cG,F,\iota}$ is connected over $\bZ[1/p]$. Base changed to $\bQ$, it is a (geometrically) connected component of $\locsys^{\tame}_{{}^cG,F,\iota}\otimes\bQ$. In particular, $\locsys^{\widehat\unip}_{{}^cG,F,\iota}\otimes \bQ$ is a reduced, local complete intersection over $\bQ$.
\end{proposition}
\begin{proof}
First we show that $\locsys^{\widehat\unip}_{{}^cG,F,\iota}$ is (geometrically) connected. If $G=T$ is a torus, this can be verified easily (e.g. see \cite[\textsection{3.2}]{zhu2020coherent}).
Then we consider the diagram \eqref{E:LocMtoG} for $\hat{P}=\hat{B}$. As explained in \cite[\textsection{3.3}]{zhu2020coherent}, there is a $\bG_m$-action on $\locsys_{{}^cB,F}$ which contracts $\locsys_{{}^cB,F}$ to $\locsys_{{}^cT,F}$. It follows that 
\[
\locsys_{{}^cB,F}^{\unip}:= \locsys_{{}^cT,F}^{\unip}\times_{\locsys_{{}^cT,F}} \locsys_{{}^cB,F}
\]
is connected. Now we note that the map $\locsys_{{}^cB,F}^{\unip}\to \locsys_{{}^cG,F}^{\widehat\unip}$ is surjective. 
Indeed, as this map is proper, it is enough to verify the surjectivity over $\bC$. Note that every pair $(g, u)\in \hat{G}\times \mU_{\hat{G}}$ satisfying $g\bar\sigma(u)g^{-1}=u^q$ is contained in Borel subgroup of $\hat{G}$. Indeed, we can fix a cocharacter $\la$ associated to $u$ and write $g=\la(\sqrt{q})x$ as in the proof of \Cref{prop: irr. comp. of loc-unip and conj. class of Au}. for some $x\in C_{\hat{G}}(u)$. In addition, we can write $x=x_0x_+\in   C_{\hat{G}}(\la, u)\ltimes R_{\hat{G}}(u)$. Choose a Borel $B'\subset C_{\hat{G}}(\la)$ containing $\la(\sqrt{q})x_0$. Then $B'\ltimes R_{\hat{G}}(u)$ is contained in a Borel $B''\subset \hat{G}$, which contains $(g,u)$.
This shows that $\locsys_{{}^cG,F}^{\widehat\unip}$ is connected. 

Next, consider the map $[q]: \hat{G}/\!\!/\hat{G}\to \hat{G}/\!\!/\hat{G}$ induced by \eqref{eq: map [q]}, and let $(\hat{G}/\!\!/\hat{G})^{[q]}$ be $[q]$-fixed point subscheme of $\hat{G}/\!\!/\hat{G}$ as before. Then the map $\locsys^{\tame}_{{}^cG,F,\iota}\to \hat{G}/\hat{G}\to \hat{G}/\!\!/\hat{G}$ factors through $(\hat{G}/\!\!/\hat{G})^{[q]}$. Over $\bQ$, $(\hat{G}/\!\!/\hat{G})^{[q]}$ is \'etale and $1$ is an isolated points. It follows that 
\[
\locsys^{\widehat\unip}_{{}^cG,F,\iota}\otimes \bQ=  \locsys^{\tame}_{{}^cG,F,\iota}\times_{(\hat{G}/\!\!/\hat{G})}\widehat{\{1\}} =\locsys^{\tame}_{{}^cG,F,\iota}\times_{(\hat{G}/\!\!/\hat{G})^{[q]}} \{1\}
\] 
is open and closed. 

Putting together, we see that $\locsys^{\widehat\unip}_{{}^cG,F,\iota}\otimes \bQ$ is a connected component of $\locsys^{\tame}_{{}^cG,F,\iota}\otimes\bQ$.
The last statement then follows as $\locsys_{{}^cG,F,\iota}$ is reduced l.c.i. over $\bZ[1/p]$. 
\end{proof}

\begin{remark}\label{rem: unipotent stack, banal }
Clearly, let $N$ be a finite product primes such that over $\bZ[\frac{1}{pN}]$, $\{1\}\subset (\hat{G}/\!\!/\hat{G})^{[q]}$ is a connected component. Then  $\locsys^{\widehat\unip}_{{}^cG,F,\iota}\otimes \bZ[\frac{1}{pN}]$ is a connected component of $\locsys^{\tame}_{{}^cG,F,\iota}\otimes \bZ[\frac{1}{pN}]$, and therefore is flat and l.c.i. over $\bZ[\frac{1}{pN}]$.
\end{remark}

\subsubsection{General tame inertia type}\label{SSS: general tame inertia type}
We assume that $\La$ is an algebraically closed field. 
Our next goal is to show that for a general tame inertia type $\zeta$, the geometry of the stack $\locsys_{{}^cG, F}^{\hat\zeta}$ is closely related to the stack of unipotent Langlands parameters of a smaller group. 
We will fix $\La$ an algebraic closure of either $\bF_\ell$ or $\bQ_\ell$ and base change everything to $\La$. But we omit $\La$ from the notations if no confusion will arise. We will also fix $\iota: \Ga_q\to W_F^t$ as before. Finally, we assume that the order of $\bar\tau\in \Ga_{\widetilde F/F}$ is invertible in $\La$.

For a tame inertia type $\zeta$ over $\La$, let $\breve\varphi^{ss}:I_F\to {}^LG$ be a semisimple representation associated to $\zeta$, and let $\delta\in \hat{G}$ such that $\delta\phi(\breve\varphi^{ss})\delta^{-1}=\breve\varphi^{ss}$, i.e. extend $\breve\varphi^{ss}$ to a Langlands parameter as in \Cref{lem: second definition of inertia type}.
Then we have the group $C_{\hat{G}}(\breve\varphi^{ss})$ equipped with an automorphism 
\[
\phi_\delta: C_{\hat{G}}(\breve\varphi^{ss})\to C_{\hat{G}}(\breve\varphi^{ss}),\quad g\mapsto \delta\bar\sigma(g)\delta^{-1}.
\]
As $\breve\varphi^{ss}$ is unique up to $\hat{G}$-conjugacy, the group $C_{\hat{G}}(\breve\varphi^{ss})$ is a well-defined conjugacy class of subgroups of $\hat{G}$, which (by abuse of notations) we denote by $\hat{G}_\zeta$. Similarly, $\phi_\delta$ is well-defined as an element in the group of outer automorphisms of $\hat{G}_\zeta$, denoted by $\phi_\zeta$.

The group $\hat{G}_\zeta$ may not be connected. We let $\hat{G}_\zeta^\circ$ be its neutral connected component. 
Under our assumption that the order of $\bar\tau$ is invertible in $\La$, $\hat{G}_\zeta^\circ$ is smooth and therefore is a connected reductive group over $\La$. Indeed, up to conjugacy, we may assume $\breve\varphi^{ss}(\iota(\tau))=t\bar\tau\in\hat{T}\bar\tau$. Then 
$\hat{G}_\zeta=\hat{G}^{t\bar\tau}$ is smooth. Let $\mU_{\hat{G}_\zeta}$ denote the unipotent variety of $\hat{G}_\zeta^\circ$, and let $\widehat{\mU}_{\hat{G}_\zeta}$ be the formal completion of  $\mU_{\hat{G}_\zeta}$ in $\hat{G}_\zeta^\circ$. 

By abuse of notations, we write 
\[
\locsys_{\hat{G}_\zeta,\breve F}=R_{I_F, \hat{G}_\zeta}/\hat{G}_\zeta,
\] 
equipped with an action of $\phi_\zeta:=\phi_{\sigma_\zeta}$ as considered in \Cref{rem: fix point construction depending on outer automorphism}. Explicitly, $\phi_\zeta$ sends a homomorphism $\breve\varphi: I_F\to \hat{G}_\zeta$ to the homomorphism $\delta\bar\sigma(\breve\varphi(\sigma^{-1}(-)\sigma))\delta^{-1}:I_F\to \hat{G}_\zeta$. There is a natural morphism
\begin{equation}\label{eq: endoscopic parameter}
\locsys_{\hat{G}_\zeta, \breve F}\to \locsys_{{}^cG, \breve F},\quad \breve\varphi\mapsto \breve\varphi \breve\varphi^{ss}.
\end{equation}
It is a direct computation to see that this morphism intertwines the action of $\phi_\zeta$ on the left hand and the action $\delta\phi(-)\delta^{-1}$ on the right hand side.

We write $\locsys^{\widehat\unip}_{\hat{G}_\zeta,\breve F}\subset \locsys_{\hat{G}_\zeta,\breve F}$ for those $\breve \varphi: I_F\to \hat{G}_\zeta$ that factors through $\mU^\wedge_{\hat{G}_\zeta}$.

\begin{proposition}\label{prop: reducing to endoscopic group}
Assume that the order of $\bar\tau$ is invertible in $\La$. Then restriction of \eqref{eq: endoscopic parameter} to $\locsys^{\widehat\unip}_{\hat{G}_\zeta,\breve F}$ induces an isomorphism
\[
\locsys^{\widehat\unip}_{\hat{G}_\zeta,\breve F}\cong \locsys_{{}^cG,\breve F}^{\hat\zeta},
\]
intertwining the $\phi_\zeta$-action on the left hand side and the $\phi$-action on the right hand side.
\end{proposition}
\begin{proof}
Recall that after choosing $\iota$, a tame inertia type $\zeta$ over $\La$ can be regarded as a $\La$-point of $\hat{G}\bar\tau/\!\!/\hat{G}$. Let $(\hat{G}\bar\tau/\!\!/\hat{G})_\zeta^\wedge$ be the formal completion of $\hat{G}\bar\tau/\!\!/\hat{G}$ at $\zeta$, and let $(\hat{G}\bar\tau/\hat{G})_\zeta^\wedge:= (\hat{G}\bar\tau/\hat{G}) \times_{\hat{G}\bar\tau/\!\!/\hat{G}}(\hat{G}\bar\tau/\!\!/\hat{G})_\zeta^\wedge$.
Then
by \Cref{ex: continuous representation of Zhat} we have the isomorphism
\[
\Spf Z^\zeta_{{}^cG, \breve F,\iota}\cong(\hat{G}\bar\tau/\!\!/\hat{G})_\zeta^\wedge,\quad \locsys_{{}^cG,\breve F,\iota}^\zeta\cong (\hat{G}\bar\tau/\hat{G})_\zeta^\wedge.
\]
Suppose $\breve\varphi^{ss}(\iota(\tau))=t\bar\tau\in \hat{T}\bar\tau$ as before.
We have the map 
\begin{equation}\label{eq: endoscopic embedding}
\hat{G}_\zeta/\hat{G}_\zeta\to \hat{G}\bar\tau/\hat{G},\quad g\mapsto gt\bar\tau
\end{equation} 
that induces a morphism $\hat{G}_\zeta/\!\!/\hat{G}_\zeta\to \hat{G}\bar\tau/\!\!/\hat{G}$, which sends $1\in \hat{G}_\zeta/\!\!/\hat{G}_\zeta$ to $\zeta\in \hat{G}\bar\tau/\!\!/\hat{G}$.
Then it induces a morphism of (formal) algebraic stacks
\begin{equation}\label{eq: endoscopic embedding-2}
\widehat{\mU}_{\hat{G}_\zeta}/ \hat{G}_\zeta\to (\hat{G}\bar\tau/\hat{G})_\zeta^\wedge.
\end{equation}
We claim that this is an isomorphism.
It is enough to show that for every algebraically closed field $K$ over $\La$, the map $\widehat{\mU}_{\hat{G}_\zeta}/ \hat{G}_\zeta)(K)\to (\hat{G}\bar\tau/\hat{G})_\zeta^\wedge(K)$ is an isomorphism (of groupoids), and \eqref{eq: endoscopic embedding-2} induces an isomorphism of tangent complexes at each (field valued) point.
Let $g\bar\tau\in \hat{G}\bar\tau$ be a $K$-point that maps to $\zeta$. As the ($\bar\tau$-twisted) Grothendieck-Springer map $\hat{B}\bar\tau/\hat{B}\to \hat{G}\bar\tau/\hat{G}$ is surjective\footnote{For  $\bar\tau\neq 1$, a situation which is possibly less familiar to many readers, see \cite[\textsection{5.3}]{xiao.zhu}. Note that the group is assumed to be semisimple and simply-connected in  \emph{loc. cit.} But such assumption is not needed for the surjectivity statement.}, we may assume that $g\bar\tau\in \hat{B}\bar\tau$ and is of the form $g=ut\bar\tau$ for some $u\in \hat{U}$. 
In addition, by \cite[Lemma 5.2.10]{xiao.zhu}, after conjugation we may assume $u\in \hat{U}^{t\bar\tau}$. This shows that $(\widehat{\mU}_{\hat{G}_\zeta}/ \hat{G}_\zeta)(K)\to (\hat{G}\bar\tau/\hat{G})_\zeta^\wedge(K)$ is surjective.
 
Let $u\in \mU_{\hat{G}_\zeta}(K)$. Note that $u\cdot t\bar\tau$ is the Jordan decomposition of this element in the disconnected group $\hat{G}\rtimes\langle\bar\tau\rangle$. It follows that
if $g\in \hat{G}(K)$ such that $gut\bar\tau g^{-1}=ut\bar\tau$, then $gug^{-1}=u$ and $gt\bar\tau g^{-1}=t\bar\tau$. Therefore, $g\in C_{\hat{G}_\zeta}(u)(K)$.
Using the same argument, we see that if $u_1t\bar\tau$ and $u_2t\bar\tau$ are conjugate in $(\hat{G}\bar\tau)(K)$ for $u_1,u_2\in \mU_{\hat{G}_\zeta}(K)$, then$u_1$ and $u_2$ are conjugate in by an element in $\hat{G}_\zeta(K)$.
Putting these facts together, we say that the map $(\widehat{\mU}_{\hat{G}_\zeta}/ \hat{G}_\zeta)(K)\to (\hat{G}\bar\tau/\hat{G})_\zeta^\wedge(K)$ is an isomorphism (of groupoids).

Next we compare the tangent complex of $\mU_{\hat{G}_\zeta}^\wedge/ \hat{G}_\zeta$ at $u\in \mU_{\hat{G}_\zeta}$ and the tangent complex of $(\hat{G}\bar\tau/\hat{G})_\zeta^\wedge$ at $ut\bar\tau$. The first is given by
$\hat\frakg^{t\bar\tau}\xrightarrow{1-\Ad_u} \hat\frakg^{t\bar\tau}$ and the second is given by $\hat\frakg\xrightarrow{1-\Ad_{ut}\bar\tau}\hat\frakg$. As the order of $\bar\tau$ is invertible in $\La$, the inclusion $\hat\frakg^{t\bar\tau}\subset\hat\frakg$ does induce a quasi-isomorphism between these two complexes, as desired.

Now, the morphism in the proposition is just the composed morphism
\[
\locsys^{\widehat\unip}_{\hat{G}_\zeta,\breve F}\cong  \mU^{\wedge}_{\hat{G}_\zeta}/ \hat{G}_\zeta\cong  (\hat{G}\bar\tau/\hat{G})_\zeta^\wedge\cong \locsys_{{}^cG,\breve F,\iota}^{\hat\zeta}\cong  \locsys_{{}^cG,\breve F}^{\hat\zeta},
\]
and therefore is an isomorphism, as desired.
\end{proof}

\begin{remark}\label{rem: component of centralizer of inertia type in hatG} 
We suppose $\breve\varphi^{ss}$ takes values in ${}^L\!T$, and write $\breve\varphi^{ss}(\iota(\tau))=t\bar\tau\in \hat{T}\bar\tau$ as before. Let $\chi\in \hat{S}$ be the image of $t\bar\tau$ as before. 
Then $\hat{G}_{\zeta}=\hat{G}^{t\bar\tau}$ for $t\in\hat{T}$. Let $\hat{G}'_\zeta=\hat{G}_\zeta^\circ\cdot \hat{T}^{\bar\tau}$. This is in fact a normal subgroup of $\hat{G}_\zeta$. 
so we have $\hat{G}_\zeta^\circ\subset\hat{G}'_\zeta\subset \hat{G}_\zeta$. Both inclusions could be straight. In fact,  it is proved in \cite{DYYZ2} that 
\[
\hat{G}'_\zeta/\hat{G}_\zeta^\circ=\hat{T}^{\bar\tau}/\hat{T}^{\bar\tau,\circ}, \quad \hat{G}_\zeta/\hat{G}'_\zeta\cong (W_0)_\chi/(W_0)_\chi^\circ,
\] 
where $(W_0)_\chi$ consists of $w\in W_0$ such that $w\chi=\chi$ while $(W_0)_\chi^\circ\subset (W_0)_\chi$ is the finite Weyl group of $\hat{G}_\zeta^\circ$. 
In addition, we may assume $\breve\varphi^{ss}$ extends to a Langlands parameter which sends $\sigma$ to $n\bar\sigma$ for some $n\in N_{\hat{G}}(\hat{T})$ such that $n\bar\sigma(n)^{-1}\in \hat{T}$. Then the image of $n$ in $W_0$ is $w$ as in \Cref{lem: classification tame inertia type}.
\end{remark}

\begin{example}\label{ex: simple case of endoscopic group}
Suppose $\hat{G}_\zeta=\hat{G}^{t\bar\tau}$ as above and is connected. Then $(\hat{B}^{t\bar\tau}, \hat{T}^{t\bar\tau})$ is a pair of Borel subgroup and maximal torus of $\hat{G}^{t\bar\tau}$. We may extend it to a pinning.
Then there is a unique automorphism of $\hat{G}_\zeta$ preserving the pinning and projecting to $\sigma_\zeta\in \mathrm{Out}(\hat{G}_\zeta)$. We still use $\sigma_\zeta$ to denote such element.
We thus obtain a unique unramified reductive group $G_\zeta$ over $F$, splits over $F_\zeta$, whose Langlands dual group is ${^L}\!G_\zeta=\hat{G}_\zeta\rtimes\Ga_{F_\zeta/F}$, where $\Ga_{F_\zeta/F}$ is generated by the (arithmetic) Frobenius, acting on $\hat{G}_\zeta$ by $\sigma_\zeta$. We can choose $\delta\in\hat{G}$ such that $\sigma_\zeta=\delta\bar\sigma(-)\delta^{-1}$ (beware that there are $Z_{\hat{G}_\zeta}^{\sigma_\zeta}$-family of choices).
Then we obtain an isomorphism (over $\La$)
\[
\locsys_{{}^cG,F}^{\hat\zeta}\simeq \locsys^{\widehat\unip}_{{}^L\!G_\zeta, F}.
\]
When $\hat{G}_\zeta$ is not connected, the situation is more complicated. 
\end{example}

\begin{example}\label{ex: regular supercuspidal}
Let $\zeta$ be a tame inertia type over $\La$, giving $\chi: I_F^t\to \hat{S}(\La)$ up to $W_0$-conjugacy as in \Cref{lem: classification tame inertia type}. We say $\zeta$ is 
\begin{itemize}
\item regular if $(W_0)_\chi=1$;
\item nonsingular if  $(W_0)_\chi^{\circ}=1$.
\end{itemize}

We lift $\chi$ to a homomorphism $\breve\varphi^{ss}:I_F^t\to {}^LT$, and write $\hat{G}_\zeta=C_{\hat{G}}(\breve\varphi^{ss})$ as before. Clearly $\hat{T}^{\bar\tau}\subset \hat{G}_\zeta$.
Then $\zeta$ is non-singular if and only if $\hat{G}_\zeta^\circ=\hat{T}^{\bar\tau,\circ}$, and $\zeta$ is regular if and only if $\hat{G}_\zeta=\hat{T}^{\bar\tau}$. We always have
\[
\locsys_{{}^cG,\breve F}^{\hat\zeta}\cong \widehat{\{1\}}/\hat{G}_\zeta,
\]
where $\widehat{\{1\}}$ denotes the formal completion of $1$ in $\hat{T}^{\bar\tau}$.

If $\zeta$ is regular, then the automorphism $\sigma_\zeta$ of $\hat{G}_\zeta$ is well-defined and is given by the natural action $w\bar\sigma$ for an element $w\in W_0$ on $\hat{T}^{\bar\tau}$. This is in fact the unique element in $W_0$ such that
$w(\bar\sigma(\chi))=\chi^q$. 
If $w\bar\sigma-q: \hat{\frakt}^{\bar\tau}\to \hat{\frakt}^{\bar\tau}$ is an isomorphism (e.g. this is the case when $\La$ is of characteristic zero), then there is an isomorphism
\[
\locsys_{{}^cG,F}^{\hat\zeta}\simeq \hat{T}^{\bar\tau}/(1-w\bar\sigma)\hat{T}^{\bar\tau}.
\]

Therefore, the geometry of connected components associated to regular inertia types are simple. However, this example also shows that $\locsys_{{}^cG,F}^{\hat\zeta}$ may not be connected. E.g. $G=\bG_m$ with $\bar\tau$ acting by inversion and $\bar\sigma$ acting trivially, then $\locsys_{{}^cG,F}^{\hat\zeta}$ has two connected components.
\end{example}

\begin{remark}
One can generalize \Cref{prop: reducing to endoscopic group} to relate $\locsys_{{}^cG,F}^{\hat\zeta}$ for any inertia type $\zeta$ to the stack of unipotent Langlands parameters to a subgroup of ${}^cG$, at least when $\La$ is a field of characteristic zero. We will discuss this in another occasion. 
\end{remark}

\subsubsection{Steinberg stack}
Let $\pi_{\breve F}: \locsys_{{}^cB, \breve F} \to  \locsys_{{}^cG, \breve F}$ be the morphism induced by ${}^cB\to {}^cG$. 
It restricts to a morphism
$\pi_{\breve F}^{\tame}: \locsys_{{}^cB, \breve F}^{\tame} \to  \locsys_{{}^cG, \breve F}^{\tame}$.
Let
\[
S^{\tame}_{{}^cG,\breve F}:=  \locsys_{{}^cB, \breve F}^{\tame}\times_{ \locsys_{{}^cG, \breve F}^{\tame}} \locsys_{{}^cB, \breve F}^{\tame}, 
\] 
which we call the (tame) Steinberg stack of ${}^cG$. 
We have the following result concerning the geometry of $S^{\tame}_{{}^cG,\breve F}$. In the course of proving it, we will also justify our choice of terminology for this stack.

\begin{proposition}\label{prop: derived geometry of Steinberg stack}
The stack $S^{\tame}_{{}^cG,\breve F}$ is a classical quasi-smooth formal algebraic stack, ind-almost of finite presentation over $\La$. The morphism $\pi_{\breve F}^{\tame}$ is quasi-smooth and proper.
\end{proposition}
\begin{proof}
We fix a topological generator $\iota(\tau) \in I_F^t$.
Let $\hat{S}=\hat{T}/(1-\bar\tau)\hat{T}$ be the $\bar\tau$-coinvariants of $\hat{T}$ as before.
Let $\hat{S}^{\wedge,p}\subset \hat{S}$ be the union of closed subschemes $Z\subset \hat{S}$ that are finite over $\bZ_\ell$ such that $Z(\overline\bF_\ell)$ are prime-to-$p$ order points in $\hat{S}(\overline\bF_\ell)$ (see \Cref{ex: continuous representation of Zhat}.)

Then as in \Cref{ex: continuous representation of Zhat}, we have
\[
\locsys_{{}^cB, \breve F, \iota}^{\tame} \cong (\hat{B}\bar\tau/\hat{B})\times_{\hat{S}}\hat{S}^{\wedge,p}, \quad
\locsys_{{}^cG, \breve F, \iota}^{\tame} \cong (\hat{G}\bar\tau/\hat{G})\times_{\hat{G}\bar\tau/\!\!/\hat{G} } (\hat{G}\bar\tau/\!\!/\hat{G})^{\wedge,p},
\]
and therefore there is a morphism
\begin{equation}\label{eq: embedding of restr Steinberg}
\hat\iota: S^{\tame}_{{}^cG,\breve F}\to \hat{B}\bar\tau/\hat{B}\times_{\hat{G}\bar\tau/\hat{G}}\hat{B}\bar\tau/\hat{B}=: S_{\hat{G}\bar\tau},
\end{equation}
which realizing $S^{\tame}_{{}^cG,\breve F}$ as the formal completion of $S_{\hat{G}\bar\tau}$ along certain closed subschemes. As $\hat{B}\bar\tau/\hat{B}\to\hat{G}\bar\tau/\hat{G}$ is quasi-smooth and proper, the second statement follows.

The stack $S_{\hat{G}\bar\tau}$ is usually called the (twisted) Steinberg stack of $\hat{G}$, and can be write as
\[
S_{\hat{G}\bar\tau}=S_{\hat{G}\bar\tau}^{\Box}/\hat{G},\quad S_{\hat{G}\bar\tau}^{\Box}=\widetilde{\hat{G}\bar\tau}\times_{\hat{G}\bar\tau}\widetilde{\hat{G}\bar\tau},
\]
where $\widetilde{\hat{G}\bar\tau}=\hat{G}\times^{\hat{B}}(\hat{B}\bar\tau)\to \hat{G}\bar\tau$ is usually called the (twisted) Grothendieck-Springer alteration of $\hat{G}$ (for the possibly less familiar twisted case, we refer to \cite[Section 5.3]{xiao.zhu}), and is a proper morphism of schemes. 
Points of $S^{\Box}_{\hat{G}\bar\tau}$ consist of $(g\bar\tau, g_1\hat{B}, g_2\hat{B})\in \hat{G}\bar\tau\times \hat{G}/\hat{B}\times\hat{G}/\hat{B}$ such that $g\in g_1\hat{B}\bar\tau(g_1)^{-1}\cap g_2\hat{B}\bar\tau(g_2)^{-1}$. Now the proposition reduces to the similar statements for $S_{\hat{G}\bar\tau}$, which are recalled in \Cref{lem: geometry of Steinberg stack}  below.
\end{proof}

We recall the following basic fact of $S_{\hat{G}\bar\tau}$.

\begin{lemma}\label{lem: geometry of Steinberg stack} 
The stack $S_{\hat{G}\bar\tau}$ is a classical local complete intersection. Its irreducible components are indexed by $W_0$. Its cotangent complex at $(g\bar\tau, g_1\hat{B}, g_2\hat{B})$ is given by the total complex of the following double complex (with the left upper corner in cohomological degree $-1$) 
\[\xymatrix{
\hat{\frakg}^*\ar[r]\ar_{\id-\Ad^*_{g}\bar\tau}[d] &  (\Ad_{g_1}\hat{\frakb})^*\oplus  (\Ad_{g_2}\hat\frakb)^* \ar^{\id-\Ad^*_{g}\bar\tau}[d] \\
\hat{\frakg}^*\ar[r]&  (\Ad_{g_1}\hat\frakb)^*\oplus  (\Ad_{g_2}\hat\frakb)^*.
}\] 
\end{lemma}
\begin{proof}
This is well-known when $\bar\tau=1$. The proof of the twisted version is the same. We include a proof for completeness.

Note that we may write $S_{\hat{G}\bar\tau}= \hat{G}\bar\tau/\hat{G}\times_{\Delta, \hat{G}\bar\tau/\hat{G}\times \hat{G}\bar\tau/\hat{G}} (\hat{B}\bar\tau/\hat{B}\times \hat{B}\bar\tau/\hat{B})$. So
the morphism $S_{\hat{G}\bar\tau}\to \hat{B}\bar\tau/\hat{B}\times \hat{B}\bar\tau/\hat{B}$ is quasi-smooth. It follows that $S_{\hat{G}\bar\tau}$ itself is quasi-smooth. To prove that it is a classical local complete intersection, it is enough to show that $\dim S_{\hat{G}\bar\tau}^\Box=\dim\hat{G}$. 

To prove this, consider the map $\widetilde{\hat{G}\bar\tau}=\hat{G}\times^{\hat{B}}(\hat{B}\bar\tau)\to \hat{G}/\hat{B}$, which induces a map $r: S_{\hat{G}\bar\tau}^\Box\to \hat{G}/\hat{B}\times \hat{G}/\hat{B}$. The fiber of the map over $(g_1\hat{B}, g_2\hat{B})$ is isomorphic to $g_1\hat{B}\tau(g_1)^{-1}\cap  g_2\hat{B}\tau(g_2)^{-1}$. It is easy to see that the intersection is nonempty only if $g:=g_1^{-1}g_2\in \hat{B}w\hat{B}$ for some $w\in W_0$, and in this case 
\[
\dim (g_1\hat{B}\tau(g_1)^{-1}\cap  g_2\hat{B}\tau(g_2)^{-1})= \dim (g^{-1}\hat{B}\tau(g)\cap \hat{B})=\dim (\Ad_w\hat{B}\cap \hat{B})= \dim \hat{B}- \ell(w),
\]
where $\ell(w)$ denotes the length of $w$ in $W$. As the variety $O(w)$ of Borels $(g_1\hat{B}, g_2\hat{B})\in \hat{G}/\hat{B}\times \hat{G}/\hat{B}$ such that $g_1^{-1}g_2\in \hat{B}w\hat{B}$ is of dimension $\dim\hat{G}-\dim\hat{B}+\ell(w)$, and that
\begin{equation}\label{eq: framed strate Steinberg}
r: \mathring{S}^{\Box}_{\hat{G}\bar\tau,w}:=r^{-1}(O(w))\to O(w)
\end{equation} 
is smooth (as $\hat{G}$ acts transitively along $O(w)$),
the desired dimension formula follows. In addition, we see that the irreducible components of $S_{\hat{G}\bar\tau}$ are indexed by $W_0$. Namely, for $w\in W_0$, there is a unique irreducible component of $S_{\hat{G}\bar\tau}$, denoted by $S_{\hat{G}\bar\tau,w}$, which contains $\mathring{S}_{\hat{G}\bar\tau,w}= \mathring{S}^{\Box}_{\hat{G}\bar\tau,w}/\hat{G}$.
For the last statement, we use the following well fact. 
\end{proof}
\begin{lemma}
Let $H$ be a smooth affine group with a finite order automorphism $\phi: H\to H$. Then the cotangent complex of $H/\Ad_\phi H$ at $h\in H$ is given by the two term complex $\frakh^*\xrightarrow{\id-\Ad_h^*\phi}\frakh^*$ in degree $[0,1]$.
\end{lemma}

\begin{remark}\label{rem: strata Steinberg stack}
Note that $\hat{G}$ acts transitively on $O(w)$ with the stabilizer at $(1, w)\in \hat{G}/\hat{B}\times \hat{G}/\hat{B}$ being $\hat{B}_w:=\Ad_w\hat{B}\cap \hat{B}$. As $w\in W_0$, the action of $\bar\tau$ on $\hat{B}$ restricts to an action of $\bar\tau$ on $\hat{B}_w$.
We see that \eqref{eq: framed strate Steinberg} descends to the map
$\frac{\hat{B}_w}{\Ad_{\bar\tau}\hat{B}_w}\to \bB\hat{B}_w$.
\end{remark}

\begin{remark}\label{rem: Normal CM of irr comp of Steinberg}
When $\bar\tau=1$, and either $\La$ is a field of characteristic zero of a field of characteristic $\ell$ bigger than the Coxeter number of $\hat{G}$, it is known that the irreducible component $S_{\hat{G}\bar\tau,w}$ as above is normal and Cohen-Macaulay (see \cite[Theorem 2.2.1]{Bez.Riche.braided}). 
\end{remark}

\begin{remark}\label{rem: Steinberg stack beyond tame level}
Our formulation of the (tame) Steinberg stack clearly suggests that there is a version of the Steinberg stack beyond the tame level, defined as
\[
S_{{}^cG,\breve F}:=  \locsys_{{}^cB, \breve F}\times_{\locsys_{{}^cG, \breve F}} \locsys_{{}^cB, \breve F}.
\] 
Then the statements in \Cref{prop: derived geometry of Steinberg stack} hold for $S_{{}^cG,\breve F}$. In addition, we similarly have \eqref{eq:trace-convolution-horocycle-diagram-spectral side} with $\tame$ removed from the subscripts. We will discuss these in details in another occassion. 
\end{remark}

For $w\in W_0$, we let 
\begin{equation}\label{eq: w part of Steinberg}
S^{\tame}_{{}^cG,\breve F,w}:=S^{\tame}_{{}^cG,\breve F}\times_{S_{\hat{G}\bar\tau}} S_{\hat{G}\bar\tau,w}.
\end{equation}
Here the map $S^{\tame}_{{}^cG,\breve F}\to S_{\hat{G}\bar\tau}$ is from the proof of \Cref{prop: derived geometry of Steinberg stack}, which depends on the choice of $\iota$. But $S^{\tame}_{{}^cG,\breve F,w}$ as a closed substack of $S^{\tame}_{{}^cG,\breve F}$ is clearly independent of the choice of $\iota$.

We specialize \eqref{eq:trace-convolution-horocycle-diagram-special} to the current setting, with $X\to Y$ being $\locsys_{{}^cB, \breve F}^{\tame}\to \locsys_{{}^cG,\breve F}^{\tame}$, equipped with the compatible $\phi$-action. It
gives rise to the following correspondence
\begin{equation}\label{eq:trace-convolution-horocycle-diagram-spectral side}
\xymatrix{
 \ar_-{\tilde{\pi}^{\tame}}[d] \widetilde{\locsys}_{{}^cG,F}^{\tame} \ar^-{\delta^{\tame}}[r] &  S^{\tame}_{{}^cG,\breve F}\\
 \locsys_{{}^cG, F}^{\tame}.   &
 }
\end{equation}
Here $\delta^{\tame}$ is the map induced by 
\begin{multline*}
\widetilde{\locsys}_{{}^cG,F}^{\tame}:= \locsys^{\tame}_{{}^cB, \breve F} \times_{\locsys_{{}^cG,\breve F}^{\tame}} \locsys_{{}^cG,F}^{\tame} \cong \locsys^{\tame}_{{}^cB,\breve F} \times_{\id\times\phi, (\locsys^{\tame}_{{}^cG,\breve F}\times \locsys^{\tame}_{{}^cG,\breve F})} \locsys^{\tame}_{{}^cG,\breve F} \\
\xrightarrow{\Delta_{\locsys^{\tame}_{{}^cB,\breve F}}\times \id }  (\locsys^{\tame}_{{}^cB,\breve F}\times \locsys^{\tame}_{{}^cB,\breve F}) \times_{\id\times\phi,(\locsys^{\tame}_{{}^cG,\breve F}\times \locsys^{\tame}_{{}^cG,\breve F})}  \locsys^{\tame}_{{}^cG,\breve F}  \cong S_{{}^cG,\breve F}^{\tame}.
\end{multline*}
Note that as mentioned in \Cref{rem: first projection vs second projection},  the map $\delta^{\tame}$ composed with the \emph{second} projection of $S_{{}^cG,\breve F}^{\tame}$ to $\locsys^{\tame}_{{}^cB, \breve F}$ is the natural projection of $\widetilde{\locsys}_{{}^cG,F}^{\tame}$ to $ \locsys^{\tame}_{{}^cB, \breve F}$.

For $w\in W_0$, we let 
\begin{equation}\label{eq: w part of tildeloc}
\widetilde{\locsys}_{{}^cG,F,w}^{\tame}:= \widetilde{\locsys}_{{}^cG,F}^{\tame}\times_{S^{\tame}_{{}^cG,\breve F}} S^{\tame}_{{}^cG,\breve F,w}.
\end{equation}
Following \cite{zhu2020coherent}, we call $\widetilde{\locsys}_{{}^cG,F,w}^{\tame}$ spectral Deligne-Lusztig stacks, which in general have non-trivial derived structures. In particular, when $w=1$ is the unit element in $W_0$, we have
\begin{equation}\label{eq: w=1 part of  tildeloc}
\widetilde{\locsys}_{{}^cG,F,1}^{\tame}\cong \mL_\phi(\locsys^{\tame}_{{}^cB,\breve F})\cong \locsys_{{}^cB,F}^{\tame}.
\end{equation}
We let
\[
\widetilde{\pi}^{\tame}_w: \widetilde{\locsys}_{{}^cG,F,w}^{\tame}\subset \widetilde{\locsys}_{{}^cG,F}^{\tame}\xrightarrow{\widetilde{\pi}^{\tame}} \locsys_{{}^cG,F}^{\tame}.
\]
Note that for $w=1$, we have $\widetilde{\pi}^{\tame}_1=\pi^{\tame}$ which is the restriction of the map $\pi$ in \eqref{E:LocMtoG} (for ${}^cP={}^cB$) to the tame part.

When $\bar\tau=1$, we also have unipotent version of the above discussions. Consider
\begin{equation}\label{eq:tildeloc}
\widetilde{\locsys}^{\tame}_{{}^cG,F}=\locsys^{\tame}_{{}^cB,\breve F}\times_{\locsys^{\tame}_{{}^cG,\breve F}}\locsys^{\tame}_{{}^cG,F} \to \locsys^{\tame}_{{}^cB,\breve F}\to \locsys^{\tame}_{{}^cT,\breve F}.
\end{equation}
Let 
\begin{equation}\label{eq:locBbreveunip}
\locsys^{\widehat\unip}_{{}^cB,\breve F}:=\locsys^{\tame}_{{}^cB,\breve F}\times_{\locsys^{\tame}_{{}^cT,\breve F}}\locsys^{\widehat\unip}_{{}^cT,\breve F},\quad \locsys^{\unip}_{{}^cB,\breve F}:=\locsys^{\tame}_{{}^cB,\breve F}\times_{\locsys^{\tame}_{{}^cT,\breve F}}\locsys^{\unip}_{{}^cT,\breve F}.
\end{equation} 
We similarly have 
\[
S_{{}^cG,\breve F}^{\widehat\unip}=\locsys_{{}^cB, \breve F}^{\widehat\unip}\times_{\locsys_{{}^cG,\breve F}^\tame} \locsys_{{}^cB, \breve F}^{\widehat\unip},\quad S_{{}^cG,\breve F}^{\unip}=\locsys_{{}^cB, \breve F}^{\unip}\times_{\locsys_{{}^cG,\breve F}^\tame} \locsys_{{}^cB, \breve F}^{\unip}.
\]
If we fix $\iota$, then $S_{{}^cG,\breve F}^{\unip}$ can be identified with
\[
S_{\hat{G}}^{\unip}:=\hat{U}/\hat{B}\times_{\hat{G}/\hat{G}}\hat{U}/\hat{B}\cong S^{\unip,\Box}_{\hat{G}}/\hat{G}, \quad \mbox{where} \quad S^{\unip,\Box}_{\hat{G}}= \widetilde{\mathcal{U}}_{\hat{G}}\times_{{\hat{G}}} \widetilde{\mathcal{U}}_{\hat{G}}.
\]
The scheme $S^{\unip,\Box}_{\hat{G}}$ is usually called the (multiplicative)  unipotent Steinberg variety. It has non-trivial derived structure, but is still quasi-smooth. 
Similar to \eqref{eq:trace-convolution-horocycle-diagram-spectral side}, we have
\begin{equation}\label{eq:trace-convolution-horocycle-diagram-spectral side-unip}
\xymatrix{
 \ar_-{\tilde{\pi}^{\widehat\unip}}[d] \widetilde{\locsys}_{{}^cG,F}^{\widehat\unip} \ar^-{\delta^{\widehat\unip}}[r] &  S^{\widehat\unip}_{{}^cG,\breve F}, & \ar_-{\tilde{\pi}^{\unip}}[d] \widetilde{\locsys}_{{}^cG,F}^{\unip} \ar^-{\delta^{\unip}}[r] &  S^{\unip}_{{}^cG,\breve F}\\
 \locsys_{{}^cG, F}^{\tame}.   & & \locsys_{{}^cG, F}^{\tame},   &
 }
\end{equation}
where
\[
\widetilde{\locsys}^{\widehat\unip}_{{}^cG,F}:=\widetilde{\locsys}^{\tame}_{{}^cG,F}\times_{\locsys^{\tame}_{{}^cT,\breve F}}\locsys^{\unip}_{{}^cT,\breve F}=\widetilde{\locsys}^{\tame}_{{}^cG,F}\times_{\locsys^{\tame}_{{}^cG,\breve F}}\locsys^{\widehat\unip}_{{}^cB,\breve F}.
\] 
and where $\widetilde{\locsys}^{\unip}_{{}^cG,F}$ is defined similarly.
As before, all these (ind-)stacks are in fact defined over $\bZ_\ell$, and once we fix $\iota$, they can be further extended to $\bZ[1/p]$.

For $w\in W_0$, let
\begin{equation}\label{eq: unipotent steinberg stack w} 
S^{\unip}_{{}^cG,\breve F, w}=(\locsys_{{}^cB,\breve F}^{\unip}\times_{\locsys_{{}^cG,\breve F}} \locsys^{\tame}_{{}^cB,\breve F})\cap S^{\tame}_{{}^cG,\breve F,w},
\end{equation}
where the intersection is taken in $S^{\tame}_{{}^cG,\breve F}$.
It is a classical stack, although it is not irreducible in general. In addition, the map $S^{\unip}_{{}^cG,\breve F, w}\to \locsys_{{}^cB,\breve F}^{\unip}\times_{\locsys_{{}^cG,\breve F}} \locsys_{{}^cB,\breve F}\xrightarrow{\pr_2} \locsys_{{}^cB,\breve F}$ factors through $\locsys_{{}^cB,\breve F}^{\unip}\subset \locsys_{{}^cB,\breve F}$ and therefore the natural closed embedding $S^{\unip}_{{}^cG,\breve F, w}\to \locsys_{{}^cB,\breve F}^{\unip}\times_{\locsys_{{}^cG,\breve F}} \locsys_{{}^cB,\breve F}$ factors as 
\[
S^{\unip}_{{}^cG,\breve F, w}\to  S_{{}^cG,\breve F}^{\unip}\to \locsys_{{}^cB,\breve F}^{\unip}\times_{\locsys_{{}^cG,\breve F}} \locsys_{{}^cB,\breve F}
\]
realizing $S^{\unip}_{{}^cG,\breve F, w}$ as a closed substack of $S_{{}^cG,\breve F}^{\unip}$. We denote  the composed morphism 
\[
\widetilde{\pi}_{w}^{\unip}:\widetilde{\locsys}_{{}^cG,F,w}^{\unip}\subset \widetilde{\locsys}_{{}^cG,F}^{\unip}\xrightarrow{\widetilde{\pi}^{\unip}} \locsys^{\tame}_{{}^cG,F}.
\]

Note that for $w=1$,
\[
\widetilde{\locsys}_{{}^cG,F,1}^{\unip}=\locsys^{\unip}_{{}^cB,F}.
\]
and we write $\widetilde{\pi}^{\unip}_{1}$ as $\pi^{\unip}: \locsys^{\unip}_{{}^cB,F}\to \locsys^{\tame}_{{}^cG,F}$.

\begin{remark}\label{rem: two different numeration of Sunip by W}
We note that we also have $S_{{}^cG,\breve F}^{\unip}\to \hat{G}\bs (\hat{G}/\hat{B}\times \hat{G}/\hat{B})$ so that the preimage of $\hat{G}\bs O(w)\subset \hat{G}\bs(\hat{G}/\hat{B}\times \hat{G}/\hat{B})$ defines a locally closed stack of $S_{{}^cG,\breve F}^{\unip}$, whose reduced substack will be denoted by $Z_w$. As in \Cref{rem: strata Steinberg stack}, we have
\[
Z_w=\frac{\hat{U}_w}{\Ad \hat{B}_w}.
\]
where $\hat{U}_w=\Ad_w\hat{U}\cap \hat{U}$ is the unipotent radical of $\hat{B}_w$.
We note that the closure of $Z_w$ in $S_{{}^cG,\breve F}^{\unip}$ is contained in but in general is not equal to (the reduced substack of) $S^{\unip}_{{}^cG,\breve F, w}$ as defined in \eqref{eq: unipotent steinberg stack w}. 
\end{remark}

Finally, we briefly discuss Steinberg stacks for general inertia types. We do not require $\bar\tau=1$.

\begin{notation}\label{not: Steinberg ass to monodromy}
For a map $Z\to \Spf Z^{\tame}_{{}^cT,\breve F}\cong R_{I_F^t,\hat{S}}$, we let $\locsys_{{}^cB,\breve F}^{Z}$ be its preimage under the map $\locsys_{{}^cB,\breve F}^{\tame}\to \locsys_{{}^cT,\breve F}^{\tame}\to \Spf Z_{{}^cT,\breve F}^{\tame}$. For a map $Z\times Z'\to \Spf Z^{\tame}_{{}^cT,\breve F}\times \Spf Z^{\tame}_{{}^cT,\breve F}$, let
\[
S_{{}^cG,\breve F}^{Z,Z'}=\locsys_{{}^cB,\breve F}^{Z}\times_{\locsys_{{}^cG,\breve F}^{\tame}}\locsys_{{}^cB,\breve F}^{Z'}.
\]
For each $w\in W_0$, there is similarly defined closed substack
\[
S_{{}^cG,\breve F,w}^{Z,Z'}=\locsys_{{}^cB,\breve F}^{Z}\times_{\locsys_{{}^cG,\breve F}^{\tame}}\locsys_{{}^cB,\breve F}^{Z'}=S_{{}^cG,\breve F}^{Z,Z'}\cap S_{{}^cG,\breve F}^{\tame}.
\]
\end{notation}

Now assume that $\La$ is an algebraically closed field.  We fix a tame inertia type $\zeta$ of ${}^cG$, and let $\{\chi: I_F^t\to \hat{S}\}_\chi$ be the $W_0$-orbit of homomorphisms corresponding to $\zeta$ as in \Cref{lem: classification tame inertia type}.
For each $\chi$, we let $\hchi\subset \Spf Z_{{}^cT,\breve F}^{\tame}\otimes\La$ be the formal completion at $\chi$.
Then \eqref{eq:trace-convolution-horocycle-diagram-spectral side} restricts to the following correspondence
\begin{equation}\label{eq:trace-convolution-horocycle-diagram-spectral side-zeta}
\xymatrix{
\widetilde{\locsys}^{\hat\zeta}_{{}^cG, F}:=\prod_{\chi}\locsys_{{}^cG,F}^{\tame}\times_{\locsys_{{}^cG,\breve F}^{\tame}} \locsys_{{}^cB, \breve F}^{\hchi}\ar[d]\ar^-{\delta^{\hat\zeta}}[r] & S_{{}^cG,\breve F}^{\hat\zeta}:= \prod_{\chi,\chi'} S_{{}^cG,\breve F}^{\hchi,\hchi'} \\
\locsys^{\hat\zeta}_{{}^cG, F}
}
\end{equation}

There is similarly defined correspondence  $\prod_{\chi,\chi'} S_{{}^cG,\breve F}^{\chi,\chi'}\leftarrow \widetilde{\locsys}^{\zeta}_{{}^cG, F}\to \widetilde{\locsys}^{\zeta}_{{}^cG, F}$.

\begin{example}
Assume that $\bar\tau=1$, and let $\zeta=\unip$ be the unipotent inertia type as before. 
In this case, $\chi: I_F^t\to \hat{S}=\hat{T}$ must be the trivial representation. Then we specialize to the unipotent Steinberg stacks as discussed above.
\end{example}

\begin{example}\label{ex: steinberg in regular case}
We continue \Cref{ex: regular supercuspidal}, but further assume that $\bar\tau=1$ and $\La=\overline\bQ_\ell$.

Let $\zeta$ be a regular inertia type, corresponding to a $W_0$-orbit of continuous homomorphisms $I_F^t\to \hat{T}(\La)$. The map $\locsys_{{}^c B, \breve F}\to \locsys_{{}^c G, \breve F}$ is finite \'etale $W_0$-cover when restricted to $\locsys_{{}^c G, \breve F}^{\hat\zeta}$. In this case, once we choose
$\chi: I_F^t\to \hat{T}$ lifting the inertia type and a lifting of the Frobenius $\sigma\in W_F$, the correspondence \eqref{eq:trace-convolution-horocycle-diagram-spectral side-zeta} can be identified with the following correspondence
\[
\xymatrix{
(W_0\times N_{\hat{G}}(\hat{T})_w\bar\sigma)/\hat{T} \ar[d]\ar[r]& \prod_{\chi_1,\chi_2} \hchi_1/\hat{T} \times_{\locsys_{{}^cG,\breve F}^{\hat\zeta}} \hchi_2/\hat{T}\\
N_{\hat{G}}(\hat{T})_w\bar\sigma /\hat{T}.&
}
\]
Here $w\in W_0$ is the element associated to $\chi$ as in  \Cref{ex: regular supercuspidal}, and $N_{\hat{G}}(\hat{T})_w\bar\sigma$ consist of those $n\in N_{\hat{G}}(\hat{T})\bar\sigma$ that maps to $w\bar\sigma$. The action of $\hat{T}$ on $N_{\hat{G}}(\hat{T})\bar\sigma$ is just the conjugation action. The vertical map is the natural projection, and the horizontal map sends $(w', n\bar\sigma)\in (W_0\times N_{\hat{G}}(\hat{T})_w\bar\sigma)/\hat{T}$ to $({w'}^{-1}(\chi), (ww')^{-1}(\chi))\in \prod_{\chi_1,\chi_2} \hchi_1/\hat{T} \times_{\locsys_{{}^cG,\breve F}^{\hat\zeta}} \hchi_2/\hat{T}$.

In addition, for every $w'\in W_0$, we have the stack
The stack $S_{{}^cG,\breve F,w'}^{\hchi,\hchi'}$, which is equal to
$S_{{}^cG,\breve F}^{\hchi,\hchi'}$ if $\chi=w\chi'$ and is empty otherwise.
\end{example}

As $\locsys_{{}^cG,F}$ is l.c.i. over $\bZ_\ell$, the stack $\Sing(\locsys_{{}^cG,F})$ of singularities of $\locsys_{{}^cG,F}$ is well-defined (see \Cref{def: Stack of singularieties}). Its points can be described as
\begin{equation}\label{eqn: Sing of loc}
\Sing(\locsys_{{}^cG,F})=\Bigl\{ (\varphi, \xi)\mid \varphi\in \locsys_{{}^cG,F},\ \xi\in H_2(W_F, \Ad^*_\varphi)=(\hat\frakg^*)^{\varphi(I_F)=1,\varphi(\sigma)=q^{-1}}\Bigr\}.
\end{equation}

Let $\hat\mN^*\subset\hat\frakg^*$ be the nilpotent cone of $\hat\frakg^*$. By definition, it is the (reduced) scheme theoretic image of the moment map $T^*(\hat{G}/\hat{B})\to \frakg^*$.
We define the ``global nilpotent cone" in this setting as
\begin{equation}\label{E: lgnlip}
\mN_{{}^cG,F}= \Bigl\{ (\varphi, \xi)\in\Sing(\Loc_{{}^cG,F}),\ \xi\in \hat\mN^* \Bigr\}.
\end{equation}
As explained in \cite[Lemma 3.3.3]{zhu2020coherent} and \cite[Proposition VIII.2.11]{Fargues.Scholze.geometrization}, if $\La$ is a field of characteristic zero, then $\mN_{{}^cG,F}= \Sing(\Loc_{{}^cG,F})$, but in general, it is a closed conic subset in $\Sing(\Loc_{{}^cG,F})$.  We similarly have $\mN_{{}^cG,F}^\tame\subset \Sing(\locsys_{{}^cG,F}^{\tame})$.

Later on we will compute the twisted categorical trace of the category of coherent sheaves on the Steinberg stack. For this purpose, we first compute pull-push of $\Sing(S_{{}^cG,\breve F}^{\tame})$ along the correspondence \eqref{eq:trace-convolution-horocycle-diagram-spectral side}. Here as $S_{{}^cG,\breve F}^{\tame}$ is a formal algebraic stack, its stack of singularities is defined as in \Cref{rem: singular support formal stack}, by choosing an embedding $S_{{}^cG,\breve F}^{\tame}\subset S_{\hat{G}\bar\tau}$ determined by $\iota(\tau)\in I_F^t$ as in \Cref{prop: derived geometry of Steinberg stack}.

The correspondence \eqref{eq:trace-convolution-horocycle-diagram-spectral side} induces the following correspondence
\[
\Sing(\locsys_{{}^cG,F}^{\tame})_{\widetilde{\locsys}_{{}^cG,F}^{\tame}}\xrightarrow{\Sing(\tilde\pi^{\tame})} \Sing(\widetilde{\locsys}_{{}^cG,F}^{\tame})\xleftarrow{\Sing(\delta^{\tame})} \Sing(S^{\tame}_{{}^cG,\breve F})_{\widetilde{\locsys}_{{}^cG,F}^{\tame} }
\]
See \eqref{eq: base change singular support} for the notation
$\Sing(\locsys_{{}^cG,F}^{\tame})_{\widetilde{\locsys}_{{}^cG,F}^{\tame}}$ and $\Sing(S^{\tame}_{{}^cG,F})_{\widetilde{\locsys}_{{}^cG,F}^{\tame} }$. There are similar correspondences induced by \eqref{eq:trace-convolution-horocycle-diagram-spectral side-unip}
\[
\Sing(\locsys_{{}^cG,F}^{\tame})_{\widetilde{\locsys}_{{}^cG,F}^{\tame}}\xrightarrow{\Sing(\tilde\pi^{\widehat\unip})} \Sing(\widetilde{\locsys}_{{}^cG,F}^{\widehat\unip})\xleftarrow{\Sing(\delta^{\widehat\unip})} \Sing(S^{\widehat\unip}_{{}^cG,\breve F})_{\widetilde{\locsys}_{{}^cG,F}^{\widehat\unip} },
\]
\[
\Sing(\locsys_{{}^cG,F}^{\tame})_{\widetilde{\locsys}_{{}^cG,F}^{\tame}}\xrightarrow{\Sing(\tilde\pi^{\unip})} \Sing(\widetilde{\locsys}_{{}^cG,F}^{\unip})\xleftarrow{\Sing(\delta^{\unip})} \Sing(S^{\unip}_{{}^cG,\breve F})_{\widetilde{\locsys}_{{}^cG,F}^{\unip} }.
\]

\begin{lemma}\label{lem: pull-push singular support}
We have 
\[
\Sing(\tilde\pi^{\tame})^{-1}(\Sing(\delta^{\tame})(\Sing(S_{{}^cG,\breve F})_{\widetilde{\locsys}_{{}^cG,F}^{\tame} }))=\mN_{{}^cG,F}^{\tame},
\]
\[
\Sing(\tilde\pi^{\widehat\unip})^{-1}(\Sing(\delta^{\widehat\unip})(\Sing(S^{\widehat\unip}_{{}^cG,\breve F})_{\widetilde{\locsys}_{{}^cG,F}^{\widehat\unip} }))=(\mN_{{}^cG,F}^{\tame})_{\locsys_{{}^cG,F}^{\widehat\unip}},
\]
\[
\Sing(\tilde\pi^{\unip})^{-1}(\Sing(\delta^{\unip})(\Sing(S^{\unip}_{{}^cG,\breve F})_{\widetilde{\locsys}_{{}^cG,F}^{\unip} }))=\Sing(\locsys_{{}^cG,F}^{\tame})_{\locsys_{{}^cG,F}^{\widehat\unip}}.
\]

\end{lemma}
\begin{proof}
We will fix $\iota: \Ga_q\to W_F^t$ as before.
By \Cref{lem: geometry of Steinberg stack}, we have
\[
\sing(S_{\hat{G}\bar\tau})=\Big\{(g\bar\tau, g_1\hat{B},g_2\hat{B}, \eta )\mid g\in  g_1\hat{B}\bar\tau(g_1)^{-1}\cap g_2\hat{B}\bar\tau(g_2)^{-1}, \ \eta\in (\hat\frakg^*)^{g\bar\tau=1} \cap (\hat\frakg/\Ad_{g_1}\hat\frakb)^*\cap  (\hat\frakg/\Ad_{g_2}\hat\frakb)^* \Big\}.
\]
We have
\[
\sing(\locsys^{\tame}_{^{c}G,F})=\Bigl\{ (\varphi, \xi)\mid \varphi\in \Loc^{\tame}_{{}^cG,F},\ \xi\in H_2(W_F, \Ad^*_\varphi)=(\hat\frakg^*)^{\varphi(\tau)=1,\varphi(\sigma)=q^{-1}}\Bigr\},
\]
and therefore
\[
\sing(\widetilde{\locsys}^{\tame}_{{}^cG,F})=\left\{ (\varphi, g\hat{B}g^{-1},\xi)\ \left| \ \begin{split} &\varphi\in \Loc^{\tame}_{{}^cG,F},\ \varphi(\tau)\in g\hat{B}g^{-1}, \\ & \xi\in(\hat\frakg^*)^{\varphi(\tau)=1}, \xi-q\Ad^*_{\varphi(\sigma)}(\xi)\in (\hat\frakg/\Ad_g\hat\frakb)^*\end{split} \right.\right\}.
\]
The map
\[
\Sing(\tilde\pi^{\tame}): \sing(\locsys_{{}^cG,F}^{\tame})_{\widetilde{\locsys}^{\tame}_{{}^cG,F}}\to \sing(\widetilde{\locsys}^{\tame}_{{}^cG,F})
\]
is given by the natural inclusion, and  the map 
\[
\Sing(\delta^{\tame}):
\sing(S_{\hat{G}\bar\tau})_{\widetilde{\locsys}^{\tame}_{{}^cG,F}}\to \sing(\widetilde{\locsys}^{\tame}_{{}^cG,F})
\] 
is given by $g\bar\tau=\varphi(\tau)$, $g_1\hat{B}g_1^{-1}=g\hat{B}g^{-1}$, $g_2\hat{B}g_2^{-1}=\varphi(\sigma)^{-1}g_1\hat{B}g_1^{-1}\varphi(\sigma)$, and $\xi=\eta$.
Now the lemma follows as for every nilpotent element $\xi\in\hat{\frakg}^*$, there is a Borel subgroup $g\hat{B}g^{-1}$ of $\hat{G}$ such that $\xi \in (\hat\frakg/\Ad_g\hat\frakb)^*$.

The second equality follows similarly. Using the variant that for every element $\xi\in\hat{\frakg}^*$, there is a Borel subgroup $g\hat{B}g^{-1}$ of $\hat{G}$ such that $\xi \in (\hat\frakg/\Ad_g\hat\fraku)^*$, the last equality also follows.
\end{proof}

\subsection{Coherent sheaves on the stack of local Langlands parameters}\label{sec:traces.spectral.aff.hecke}
In the sequel, we will assume that $\La$ is a Dedekind domain which is either an algebraic field extension of $\bQ_\ell$ or $\bF_\ell$, or a finite extension of $\bZ_\ell$.

\subsubsection{The category of coherent sheaves}
The main player in the spectral side of the categorical local Langlands correspondence is the category of coherent sheaves on $\locsys_{{}^cG,F}$ and its variants. First, we fix $L/F$ and consider $\locsys_{{}^cG,L/F}$, which is classical (and therefore eventually coconnective) almost of finite presentation over $\La$. We have the action of $\Perf(\locsys_{{}^cG,L/F})$ on $\Coh(\locsys_{{}^cG,L/F})$, inducing a fully faithful embedding $\Xi_L: \Perf(\locsys_{{}^cG,L/F})\to \Coh(\locsys_{{}^cG,L/F})$ (as in \eqref{eq: XiX functor}). If $\La$ is a field of characteristic zero, then $\ind\Perf(\locsys_{{}^cG,L/F})\cong \Qcoh(\locsys_{{}^cG,L/F})$.

Now we shall regard $\locsys_{{}^cG,F}=\colim_L \locsys_{{}^cG,L/F}$ as an ind-finite type algebraic stack over $\La$. Then we have category of quasi-coherent sheaves 
\[
\Qcoh(\locsys_{{}^cG,F})=\lim_L\Qcoh(\locsys_{{}^cG,L/F}),
\]
which is a $\La$-linear (presentable, stable $\infty$-)category. It contains a full subcategory of perfect complex
\[
\Perf(\locsys_{{}^cG,F})=\lim_{L} \Perf(\locsys_{{}^cG,L/F}).
\] 

\begin{example}\label{ex: evaluation bundle}
Let $V$ be a representation of $\hat{G}$ on finite projective $\La$-modules, which can be regarded as an object in $\Perf(\bB \hat{G})$. The pullback of $V$ along the natural morphism $\locsys_{{}^cG,F}\to \bB \hat{G}$ gives rise to an object in $\Perf(\Loc_{{}^cG,F})$ denoted by $\widetilde{V}$, and is usually called the ``evaluation bundle" or the "tautological bundle". For example when $V$ is the trivial representation (on $\La$), then $\widetilde{V}=\mO_{\locsys_{{}^cG,F}}$ is the structure sheaf. 

If $F\subset E\subset \widetilde{F}$ is a field such that $V$ extends to a representation of ${}^cG_E:=\hat{G}\rtimes (\bG_m\times \Ga_{\widetilde{F}/E})$, then $\widetilde{V}$ is canonically equipped with a tautological representation 
\begin{equation}\label{eq:univ Galois action}
\varphi_{V}^{\mathrm{univ}}: W_E\to \GL(\widetilde V)
\end{equation}
such that for every $\Spec A\to \locsys_{{}^cG,F}$ corresponding to $\varphi: W_F\to {}^cG(A)$, the pullback of $\varphi_{V}^{\mathrm{univ}}$ to $\spec A$ is the continuous representation $W_E\xrightarrow{\varphi|_{W_E}} {}^cG_E(A)\to \GL(V\otimes A)$. More generally, if $V$ is a representation of $\prod_i {}^cG_{E_i}$, for a finite collection of field extensions $F\subset E_i\subset \widetilde{F}$, then $\widetilde{V}$ is equipped with a representation of $\prod_i W_{E_i}$.

Let $\widetilde{V}^{\tame}$ (resp. $\widetilde{V}^{\unip}$) denote the restriction of $\widetilde{V}$ to $\locsys_{{}^cG,F}^{\tame}$ (resp.  $\locsys_{{}^cG,F}^{\widehat\unip}$). However, if the context is clear, we will sometimes simply denote them by $\widetilde{V}$.
\end{example}

\begin{remark}\label{rem: ev bundle on the stack of geometric parameters}
Let $V$ be a representation of $\hat{G}$ as above. Then we have $\breve{\widetilde{V}}$ the corresponding evaluation bundle over $\locsys_{{}^cG,\breve F}$. Its pullback along the morphism $\res: \locsys_{{}^cG, F}\to \locsys_{{}^cG,\breve F}$ is $\widetilde{V}$.
Now suppose $V\in ({}^cG)^I$. Regarding it as a $\hat{G}$-representation (via diagonal embedding). Then similar to \eqref{eq:univ Galois action}, we have a representation of $(I_F)^I$ on $\breve{\widetilde{V}}$ which extends to a representation of $(W_F)^I$ on $\widetilde{V}=\res^*\breve{\widetilde{V}}$ as above.

We fix a lifting of the Frobenius $\sigma\in W_F$. Then for each $i\in I$, there is a canonical isomorphism
\[
\Phi_{i}: \phi^*\breve{\widetilde{V}}\cong \breve{\widetilde{V}},
\]
where $\phi$ is the automorphism of $\locsys_{{}^cG,\breve F}$ defined before. Via pulling back to $\locsys_{{}^cG,F}$ via two maps (as in \eqref{eq: arith parameter as fix point of geometric parameter}), and under the canonical identification
\[
\res^*\breve{\widetilde{V}}\cong \widetilde{V}\cong (\res_\phi)^*\breve{\widetilde{V}},
\]
the isomorphism $\Phi_{i}$ is given by the tautological action of $\varphi^{\mathrm{univ}}(\sigma_i)$ on $\widetilde{V}$, where $\sigma_i\in (W_F)^I$ is equal to $\sigma$ in the $i$th component and is equal to $1$ otherwise.
\end{remark}

As we regard $\locsys_{{}^cG,F}$ as an ind-algebraic stack, the category $\Coh(\locsys_{{}^cG,F})$ of coherent sheaves on it satisfies
\[
\Coh(\locsys_{{}^cG,F})=\colim_{L} \Coh(\locsys_{{}^cG,L/F}).
\] 
By definition, an object in $\Coh(\locsys_{{}^cG,F})$ is supported on $\Coh(\locsys_{{}^cG,L/F})$ for some $L$. In particular, by our convention the structure sheaf of $\Loc_{{}^cG,F}$ itself is not regarded as a coherent sheaf. (But the structure sheaf of  $\Loc_{{}^cG,L/F}$ for each $L$ is a coherent sheaf.)
There is a natural action of $\Perf(\Loc_{{}^cG,F})$ on $\Coh(\locsys_{{}^cG,F})$.

As for a finite extension $L'/L$, $\locsys_{{}^cG,L/F}\subset \locsys_{{}^cG,L'/F}$ is open and closed, $*$-extension is also the left adjoint of $*$-restriction. Therefore, $\Coh(\locsys_{{}^cG,L/F})$ is a direct summand of $\Coh(\locsys_{{}^cG,L'/F})$. In particular, when $G$ splits over a tamely ramified extension, we write 
\[
\locsys_{{}^cG,F}=\locsys_{{}^cG,F}^{\tame}\sqcup \locsys_{{}^cG,F}^{>0}
\] 
as a disjoint union. 
If $\La$ is an algebraically closed field, we have a further decomposition 
\[
\locsys_{{}^cG,F}^{\tame}=\sqcup_{\zeta} \locsys^{\hat\zeta}_{{}^cG,F}
\]
according to the tame inertia types. In particular, if $G$ splits over an unramified extension, we have a connected component $\locsys_{{}^cG,F}^{\widehat\unip}$.
Then we have an orthogonal decomposition
\begin{eqnarray}\label{eq: depth decomposition spectral side}
\indcoh(\locsys_{{}^cG,F})&=&\indcoh(\locsys_{{}^cG,F}^{\tame})\bigoplus \indcoh(\locsys_{{}^cG,F}^{>0})\\
&=&(\bigoplus_{\zeta}\indcoh(\locsys_{{}^cG,F}^{\hat\zeta}))\bigoplus \indcoh(\locsys_{{}^cG,F}^{>0}). \nonumber
\end{eqnarray}

We let
\[
\Coh_{\mN_{{}^cG,F}}(\locsys_{{}^cG,F})\subset \Coh(\locsys_{{}^cG,F})
\]
be the full subcategory consisting of those coherent complexes $\mF$ such that $s.s.(\mF)\subset \mN_{{}^cG,F}$.  
Similarly, we have $\Coh_{\mN^\tame_{{}^cG,F}}(\locsys^\tame_{{}^cG,F})$.

Recall the automorphism $\theta$ of $\locsys_{{}^cG,F}$ from \eqref{eq:Cartan involution of loc}.
We let 
\begin{equation}\label{eq: modified GS duality}
\verd^{\Coh'}:=\theta_{*} \circ \verd^{\Coh}: \Coh(\locsys_{{}^cG,F})^{\op}\cong  \Coh(\locsys_{{}^cG,F}),
\end{equation}
where $\verd^{\Coh}=\verd^{\Coh}_{\locsys_{{}^cG,F}}$ is the usual Grothendieck-Serre duality functor \eqref{eq: GS duality, stack}. 
By \Cref{prop: GS duality and Singular support} this duality restricts to an anti-involution of $\Coh_{\mN_{{}^cG,F}}(\locsys_{{}^cG,F})$. It also restricts to an anti-involution of $\Coh_{\mN^\tame_{{}^cG,F}}(\locsys_{{}^cG,F}^{\tame})$. We will call $\verd^{\Coh'}$ (and its ind-completion $\verd^{\indcoh'}$) the twisted Grothendieck-Serre duality.

We have the following simple observation, essential due to the fact that $\locsys_{{}^cG,F}$ is of relative dimension zero over $\La$.
\begin{lemma}\label{lem:right t exact GS duality}
The Grothendieck-Serre duality $\verd^{\coh}: \Coh(\locsys_{{}^cG,F})^{\op}\to \Coh(\locsys_{{}^cG,F})$ is right $t$-exact, i.e. it sends $\Coh(\locsys_{{}^cG,F})^{\leq 0}=(\Coh(\locsys_{{}^cG,F})^{\op})^{\geq 0}$ to $(\Coh(\locsys_{{}^cG,F})^{\op})^{\geq 0}$.
The same statement holds for $\verd^{\Coh'}$.
\end{lemma}
\begin{proof}As $\theta_{*}$ is $t$-exact, the second statement follows from the first. For the first, let $\mF\in \Coh(\locsys_{{}^cG,F})^{\heartsuit}$. It is enough to show that $f^*\verd^{\Coh}(\mF)\in \Coh(\locsys_{{}^cG,F}^{\Box})^{\geq 0}$, where $f: \locsys_{{}^cG,F}^{\Box}\to \locsys_{{}^cG,F}$ is the natural smooth cover. By \Cref{lem: Grothendieck duality commute with proper push}, we have
\begin{multline*}
f^*\verd^{\Coh}(\mF)=\verd^{\Coh}(f^{\indcoh,!}\mF)=\underline\Hom(f^{\indcoh,!}\mF, \cohdual_{\locsys_{{}^cG,F}^{\Box}})\\
=\underline\Hom(f^*\mF, f^*\cohdual_{\locsys_{{}^cG,F}})=\underline\Hom(f^*\mF, \mO_{\locsys_{{}^cG,F}^{\Box}}),
\end{multline*}
which belongs to $\Coh(\locsys_{{}^cG,F}^{\Box})^{\geq 0}$, as desired.
\end{proof}

Recall the notion of admissible objects in a $\La$-linear dualizable category (see \Cref{def:adm vs compact}). By \Cref{lem: char of adm obj in cg cat}, the category $\indcoh(\locsys_{{}^cG,F})^{\adm}$ consist of objects in $\mG$ such that 
\[
\Hom_{\indcoh(\locsys_{{}^cG,F})}(\mF,\mG)\in \Perf_\La
\] 
for every $\mF\in \Coh(\locsys_{{}^cG,F})$. 
If $\mG$ is in addition coherent, such $\mG$ must support on $Z\times_{\Spf Z_{{}^cG,F}} \locsys_{{}^cG, F}$, where $Z$ is a closed subscheme of $\Spf Z_{{}^cG,F}$ finite over $\La$. We have
\begin{equation}\label{eq: coh and adm objects}
\indcoh(\locsys_{{}^cG,F})^{\cpt}\cap \indcoh(\locsys_{{}^cG,F})^{\adm}\subset \Coh(\locsys_{{}^cG,F}\times_{\Spf Z_{{}^cG,F}} \Spf Z_{{}^cG,F}^\wedge),
\end{equation}
Here $\Spf Z_{{}^cG,F}^\wedge$ denotes the formal completion of $\Spf Z_{{}^cG,F}$ at all closed points. The following lemma is easy.
\begin{lemma}
If $\La$ is a field of characteristic zero,  the inclusion \eqref{eq: coh and adm objects} is an equality.
\end{lemma}

We note that not all admissible objects in $\indcoh(\locsys_{{}^cG,F})$ are coherent. 
\begin{example}\label{ex: admissible ind-coherent sheaf on loc}
Suppose $\La$ is a field.
Let $\varphi$ be a parameter such that $H^2(W_F^t, \Ad^0_\varphi)=0$, or equivalently $q^{-1}$ is not an eigenvalue of the linear operator $\varphi(\sigma): \hat\frakg^{\varphi(I_F)}\to  \hat\frakg^{\varphi(I_F)}$.
By the proof of \Cref{lem:discrete parameter condition},  $\varphi$ gives rise to a smooth point of $\locsys_{{}^cG,F}$. Let $i_\varphi: \{\varphi\}/C_{\hat{G}}(\varphi)\to \locsys_{{}^cG,F}$ denote the locally closed embedding of residual gerbe, which is a schematic morphism of finite tor amplitude. Then 
\[
(i_\varphi)^{\indcoh}_*: \indcoh(\{\varphi\}/C_{\hat{G}}(\varphi))\cong \qcoh(\{\varphi\}/C_{\hat{G}}(\varphi))\cong \rep(C_{\hat{G}}(\varphi))\to \indcoh( \locsys_{{}^cG,F})
\]
sends admissible objects to admissible objects (as it admits an left adjoint given by $(i_\varphi)^{\indcoh,*}$). 

Suppose that $C_{\hat{G}}(\varphi)$ is smooth.
Then the regular representation $\mathrm{Reg}_{C_{\hat{G}}(\varphi)}$ of $C_{\hat{G}}(\varphi)$ (i.e. the ring of regular functions of $C_{\hat{G}}(\varphi)$ equipped with the action of $C_{\hat{G}}(\varphi)$ by left translation) is always admissible object in $\indcoh(\{\varphi\}/C_{\hat{G}}(\varphi))$, althought itself may not be coherent if $C_{\hat{G}}(\varphi)$ is not finite.  Therefore $(i_\varphi)_*^{\indcoh}(\mathrm{Reg}_{C_{\hat{G}}(\varphi)})$ is admissible in $ \indcoh( \locsys_{{}^cG,F})$.
In addition, if $\La$ is a field of characteristic zero, then every finite dimensional representation $V$ of $C_{\hat{G}}(\varphi)$, regarded as a vector bundle on $\{\varphi\}/C_{\hat{G}}(\varphi)$, is admissible. Therefore  $(i_\varphi)^{\indcoh}_*(V)$ is an admissible in $\indcoh( \locsys_{{}^cG,F})$.  But may not be coherent as $i_\varphi$ may not be a closed embedding. (E.g. when $\hat{G}$ is semisimple and $\varphi$ is a unipotent discrete parameter, then $i_\varphi$ is an open embedding.)
\end{example}

Now let 
\begin{equation}\label{eq: adm coh duality}
(\verd^{\indcoh'})^{\adm}: (\indcoh(\locsys_{{}^cG,F})^{\adm})^{\op}\to \indcoh(\locsys_{{}^cG,F})^{\adm}
\end{equation}
be the duality of admissible objects as from \eqref{eq: dual of admissible objects}. Note that $(\verd^{\indcoh'})^{\adm}=\theta_*\circ (\verd^{\indcoh})^{\adm}$.
We will make use of the following observation. 
\begin{lemma}\label{lem: adm coh duality discrete par}
Suppose $\La$ is a field of characteristic zero. Let $\varphi: W_F\to {}^cG(\La)$ be a parameter such that $H^2(W_F, \Ad_\varphi^0)=0$.
Let $V$ be a finite dimensional representation of $C_{\hat{G}}(\varphi)$, regarded as a coherent sheaf on $\{\varphi\}/C_{\hat{G}}(\varphi)$. 
Let $\fraku$ be the Lie algebra of the unipotent radical of $C_{\hat{G}}(u)$ and let $d=\dim_\La \fraku$
Then we have
\[
(\verd^{\indcoh'})^{\adm}((i_\varphi)^{\indcoh}_*V)\cong \theta^{\indcoh}_*((i_\varphi)^{\indcoh}_*(V^*\otimes (\wedge^{d} \fraku)[d]))).
\] 
\end{lemma}
As mentioned at the end of \Cref{rem: SL2 version of WD parameter}, $\fraku$ may not be trivial, even if $\varphi$ is Frobenius-semsimple. However, if $\varphi$ is essentially discrete, then $\fraku$ is trivial.
\begin{proof}
Note that $i_\varphi: \{\varphi\}/C_{\hat{G}}(\varphi)\to \locsys_{{}^cG,F}$ is a regular embedding of codimension $\dim C_{\hat{G}}(\varphi)$. 
We apply \Cref{lem: adm duality compatible with shrek pullback along closed embedding} to this setting. 
Note that $\cohdual_{\{\varphi\}/C_{\hat{G}}(\varphi)}$ is an invertible sheaf on $\{\varphi\}/C_{\hat{G}}(\varphi)$. We let  
\[
\verd_{\{\varphi\}/C_{\hat{G}}(\varphi)}^{\coh'}(-)=  \verd_{\{\varphi\}/C_{\hat{G}}(\varphi)}^{\coh}( (-)\otimes \cohdual_{\{\varphi\}/C_{\hat{G}}(\varphi)}^{-1})
\] 
be the modified Grothendieck Serre duality on $\{\varphi\}/C_{\hat{G}}(\varphi)$, which in fact is nothing but the naive duality sending a finite dimensional representation $V$ to its dual representation $V^*$ (regarded as coherent sheaves on $\{\varphi\}/C_{\hat{G}}(\varphi)$). We shall use $\verd_{\{\varphi\}/C_{\hat{G}}(\varphi)}^{\indcoh'}$ to denote its ind-extension. Then by \Cref{lem: adm duality compatible with shrek pullback along closed embedding} 
for any admissible object $\mV\in \indcoh(\{\varphi\}/C_{\hat{G}}(\varphi))$ we have 
\[
(\verd_{\locsys_{{}^cG,F}}^{\indcoh'})^{\adm}((i_\varphi)^{\indcoh}_*\mV)\cong \theta^{\indcoh}_*((i_\varphi)^{\indcoh}_*((\verd_{\{\varphi\}/C_{\hat{G}}(\varphi)}^{\indcoh'})^{\adm}(\mV)).
\]
Note that the category $\indcoh(\{\varphi\}/C_{\hat{G}}(\varphi))$ is a proper category over $\La$ so compact objects are admissible. 
Let $S$ be its Serre functor. Then by \Cref{ex: explicit Serre and admissible-2} for a finite dimensional representation $V$ of $C_{\hat{G}}(\varphi)$ regarded as a vector bundle on $\{\varphi\}/C_{\hat{G}}(\varphi)$, we have
\[
(\verd_{\{\varphi\}/C_{\hat{G}}(\varphi)}^{\indcoh'})^{\adm}(V)\cong S(V^*).
\]
Now the statement follows from discussions in \Cref{L: Serre functor coh classifying stack}.
\end{proof}

We let 
\begin{equation}\label{eq: spectral Bernstein center}
Z_{{}^cG,F}:=H^0\rg(\locsys_{{}^cG,F},\mO). .
\end{equation}
This is sometimes called the spectral Bernstein center. On the other hand, there is the $E_2$-center 
\[
Z(\indcoh(\locsys_{{}^cG,F})):=Z(\indcoh(\locsys_{{}^cG,F})/\Mod_\La)
\]
of the dualizable category $\indcoh(\locsys_{{}^cG,F})$, see \eqref{eq: def of cat center} and \eqref{eq: hor trace and center for id}. We note that there are natural functors
\[
\indcoh(\locsys_{{}^cG, F})\xrightarrow{\Delta_*^{\indcoh}} \indcoh(\locsys_{{}^cG, F}\times \locsys_{{}^cG,F})\to \End(\indcoh(\locsys_{{}^cG,F}))
\]
sending $\cohdual_{\locsys_{{}^cG, F}}$ to the identity functor. 
Here the last functor is given by the usual integral transform $\mF\mapsto  (\pr_2)^{\indcoh}_*((\pr_1)^{!,\indcoh}(-)\os \mF)$. 
We thus obtain a map
\begin{equation}\label{eq: stable center}
Z_{{}^cG,F}\to H^0Z(\indcoh(\locsys_{{}^cG,F})).
\end{equation}

\begin{remark}
If $\La=\bQ_\ell$, then
\[
\indcoh(\locsys_{{}^cG, F})\otimes_\La \indcoh(\locsys_{{}^cG,F})\to \indcoh(\locsys_{{}^cG, F}\times_\La \locsys_{{}^cG,F})
\]
is an equivalence, so $\indcoh(\locsys_{{}^cG, F}\times \locsys_{{}^cG,F})\to \End(\indcoh(\locsys_{{}^cG,F}))$ is an equivalence.
We believe this still holds over $\bZ_\ell$, excluding a few small $\ell$, although we have not checked this. 
\end{remark}

Now assume that $\La$ is a field of characteristic zero. We let
\[
\mL(\locsys_{{}^cG,F}):= \locsys_{{}^cG,F}\times_{\locsys_{{}^cG,F}\times \locsys_{{}^cG,F}}\locsys_{{}^cG,F}\cong \mL_\phi(\mL(\locsys_{{}^cG,\breve F})).
\]
It has a highly derived structure. Its underlying classical stack classifies
\[
\bigl\{(\varphi, \kappa)\mid \varphi: W_F\to {}^cG, \kappa\in \hat{G}, \kappa\varphi\kappa^{-1}=\varphi \bigr\}/\hat{G}.
\]
We apply the discussions in \Cref{rem: horizontal trace qcoh vs indcoh} to the current setting, giving a smooth map
\[
\mO_{\mL(\locsys_{{}^cG,F})}\to \cohdual_{\mL(\locsys_{{}^cG,F})}
\]
which induces
\[
\mathrm{tr}(\qcoh(\locsys_{{}^cG,F}))=\rg(\mL(\locsys_{{}^cG,F}), \mO)\to \mathrm{tr}(\indcoh(\locsys_{{}^cG,F}))=\rg(\mL(\locsys_{{}^cG,F}), \cohdual_{\mL(\locsys_{{}^cG,F})}).
\]
If $\varphi$ is a smooth point of $\locsys_{{}^cG,F}$, then $i_{\varphi}: \bB C_{\hat{G}}(\varphi)\to \locsys_{{}^cG,F}$ is a quasi-smooth locally closed embedding. The functor $(i_\varphi)^{\indcoh,*}$ admits a continuous right adjoint and therefore we have
\[
\mathrm{tr}(\qcoh(\locsys_{{}^cG,F}))\to \mathrm{tr}(\indcoh(\locsys_{{}^cG,F}))\to \mathrm{tr}(\indcoh(\bB C_{\hat{G}}(\varphi))),
\]
which is identified with the pullback 
\[
\rg(\mL(\locsys_{{}^cG,F}), \mO)\to \rg(C_{\hat{G}}(\varphi)/C_{\hat{G}}(\varphi),\mO)
\]
along the map $\mL(\bB C_{\hat{G}}(\varphi))\to \mL(\locsys_{{}^cG,F})$.

\subsubsection{Spectral affine Hecke categories}\label{SSS: spec affine Hecke}
From now on until the end of the section, we assume that $G$ is tamely ramified.
We first assume $\bar\tau=1$. 
Recall we have the proper morphism  $\pi_{\breve F}^{\unip}: \locsys_{{}^cB,\breve F}^{\unip}\to \locsys_{{}^cG,\breve F}^\tame$ and the unipotent Steinberg stack $S_{{}^cG,\breve F}^\unip$. Then $\indcoh(S_{{}^cG,\breve F}^\unip)$ also admits a monoidal structure by  \Cref{prop: trace of convolution category} \eqref{prop: trace of convolution category-1.5}, with the monoidal unit given by 
\[
(\Delta_{ \locsys_{{}^cB,\breve F}^\unip/ \locsys_{{}^cG,\breve F}^\tame})_*\cohdual_{ \locsys_{{}^cB,\breve F}^\unip}.
\]
We write the monoidal product as
\[
\indcoh(S_{{}^cG,\breve F}^\unip)\otimes_\La \indcoh(S_{{}^cG,\breve F}^\unip)\to \indcoh(S_{{}^cG,\breve F}^\unip),\quad (\mF,\mG)\mapsto \mF\star\mG,
\]
and call it the $!$-convolution product, or just convolution product for simplicity.

\begin{remark}\label{rem: shrek and star convolution of spectral hecke}
Note that one can apply \Cref{prop: trace of convolution category} \eqref{prop: trace of convolution category-1} to endow $\indcoh(S_{{}^cG,\breve F}^\unip)$ with another monoidal structure, which we call the $*$-convolution product.  The monoidal unit is
\[
(\Delta_{ \locsys_{{}^cB,\breve F}^\unip/ \locsys_{{}^cG,\breve F}^\tame})_*\mO_{ \locsys_{{}^cB,\breve F}^\unip}.
\]
As the $(\indcoh,!)$-pullback and $(\indcoh,*)$-pullback along $\locsys_{{}^cB,\breve F}^\unip\to  \locsys_{{}^cB,\breve F}^\unip\times  \locsys_{{}^cB,\breve F}^\unip$ defers by shifting by $\dim \hat{T}$, we see that $\mF\mapsto \mF[\dim \hat{T}]$ is a monoidal equivalence from $\indcoh(S_{{}^cG,\breve F}^\unip)$ equipped with the $!$-convolution product to $\indcoh(S_{{}^cG,\breve F}^\unip)$ equipped with the $*$-convolution product.

We shall mainly use the first monoidal structure.
\end{remark}

 We will need the following result.
 
 \begin{lemma}\label{lem: Coh Steinberg exterior tensor} 
 The exterior tensor product 
 \[
 \boxtimes: \indcoh(S_{{}^cG,\breve F}^\unip)\otimes_\La \indcoh(S_{{}^cG,\breve F}^\unip)\to \indcoh(S_{{}^cG,\breve F}^\unip\times S_{{}^cG,\breve F}^\unip)
 \]
 is an equivalence.
\end{lemma}
As explained in \Cref{SSS: star coh sheaf theory}, such type of result is subtle when $\La$ is a field of positive characteristic, or a more general base ring.
\begin{proof}
Recall from \Cref{rem: two different numeration of Sunip by W} that $S_{{}^cG,\breve F}^\unip$ admits a stratification with the underlying reduced substack of  each strata being $Z_w$, which is smooth.
Using \Cref{cor: tensor equiv coh}, it reduces to show that $\indcoh(Z_v)\otimes \indcoh(Z_w)\to \indcoh(Z_v\times Z_w)$ is an equivalence. Then by \Cref{lem: criterion of tensor product equivalence via diagonal} we may reduce to the case $w=v$. We claim that
$\indcoh(Z_w)$ is generated by the $*$-pullback of objects in $\indcoh(\bB \hat{B}_w)$ along the map $Z_w=\hat{U}_w/\Ad \hat{B}_w\to \Ad \hat{B}_w$. The same proof will also show a similar statement holds with $Z_w$ replaced by $Z_w\times Z_w$. We then reduce to show that $\indcoh(\bB \hat{B}_w)\otimes \indcoh(\bB \hat{B}_w)\to \indcoh(\bB \hat{B}_w\times \bB\hat{B}_w)$ is essential surjective, which is clear.

To prove the claim, by \Cref{lem: a criterion of generation under star pullback of perf}, it is enough to show that $(\Delta_{Z_w/\bB \hat{B}_w})_*\mO_{Z_w}\in \Perf((\hat{U}_w\times \hat{U}_w)/\Delta(\Ad\hat{B}_w))$ is contained in the idempotent complete subcategory generated by the essential image of $\Perf(\hat{B}_w)$ under $*$-pullback. We consider the map $\hat{U}_w\times \hat{U}_w\mapsto \hat{U}_w,\quad (z_1,z_2)\mapsto z_1^{-1}z_2$.
This morphism is $\hat{B}_w$-equivariant. 
 In addition, the diagonal is the preimage of $\{1\}\in \hat{U}_w$. Therefore, it is enough to show that
the $*$-pushforward of the structure sheaf $\mO_{\bB \hat{B}_w}$ along closed embedding $\bB \hat{B}_w=\{1\}/\Ad \hat{B}_w\to \hat{U}_w/\Ad \hat{B}_w$ 
belongs to the subcategory of $\Perf(Z_w)$  generated by $\Perf(\bB \hat{B}_w)$ (under pullbacks). Now we can filter $\hat{U}_w$ by normal subgroups $\hat{U}_w=Z_0\supset Z_1\supset Z_2\supset\cdots$, with each $Z_{i+1}$ codimension one in $Z_i$, and such that the ideal of definition of $Z_{i+1}$ inside $Z_i$ is generated by a function $f_i\in \mO(Z_i)$, on which $\hat{B}_w$ acts through a character. Then by induction on $i$, we see that each $\mO_{Z_i/\Ad \hat{B}_w}$, regarded as an object in $\Perf(\hat{U}_w/\Ad \hat{B}_w)$ via the $*$-pushforward along the closed embedding $Z_i\subset \hat{U}_w$, belongs to the subcategory of $\Perf(Z_w)$ generated by $\Perf(\bB \hat{B}_w)$ (under pullbacks). The claim and therefore the lemma is proved.
\end{proof}

\begin{corollary}\label{cor: Coh Steinberg exterior tensor-2}
For every smooth algebraic stack $Y\in \ArStk_\La^{\aft}$, the exterior tensor product functors
\[
\indcoh(\locsys_{{}^cB,\breve F}^{\unip})\otimes_\La \indcoh(Y)\to \indcoh(\locsys_{{}^cB,\breve F}^{\unip}\times Y)
\]
\[
\indcoh(S_{{}^cG,\breve F}^{\unip})\otimes_\La \indcoh(Y)\to \indcoh(S_{{}^cG,\breve F}^{\unip}\times Y)
\]
are equivalences.
\end{corollary}
\begin{proof}
Thanks to the proof of \Cref{lem: Coh Steinberg exterior tensor}, \Cref{lem: criterion of tensor product equivalence via diagonal} is applicable to the sheaf theory $\indcoh^*$ and $X=Z_w$. To deduce the second case, we apply \Cref{cor: tensor equiv coh} again.
\end{proof}
 
 Let us have more discussion of the Grothendieck-Serre duality of $S^{\unip}_{{}^cG,\breve F}$.

By \Cref{lem: Coh Steinberg exterior tensor} and \Cref{prop: trace of convolution category} \eqref{prop: trace of convolution category-2}, the monoidal category $\indcoh(S^{\unip}_{{}^cG,\breve F})$ is a rigid monoidal category, with a Frobenius structure given by
\[
\Hom(\Delta^{\indcoh}_* \cohdual_{\locsys_{{}^cB,\breve F}^{\unip}},-): \Coh(S^{\unip}_{{}^cG,\breve F})\to \Mod_\La.
\]
We let $\verd_{S^{\unip}_{{}^cG,\breve F}}^{\mathrm{sr}}$ denote the self-duality of $\indcoh(S^{\unip}_{{}^cG,\breve F})$ induced by this Frobenius algebra structure. See \Cref{ex: self-duality of cpt gen rigid monoidal cat}. 

On the other hand, the category $\indcoh(S^{\unip}_{{}^cG,\breve F})$ is also equipped with a symmetric monoidal product given by the $!$-tensor product $\os$ of coherent sheaves, which in fact also admits a Frobenius structure given by
\[
\rg^{\indcoh}(S^{\unip}_{{}^cG,\breve F},-): \indcoh(S^{\unip}_{{}^cG,\breve F})\to \Mod_\La.
\]
As explained in \Cref{SS: GS duality}, the induced self-duality of $\indcoh(S^{\unip}_{{}^cG,\breve F})$ is nothing but the Grothendieck-Serre duality $\verd^{\indcoh}_{S^{\unip}_{{}^cG,\breve F}}$ of $S^{\unip}_{{}^cG,\breve F}$.

Let $\sw: S^{\unip}_{{}^cG,\breve F}\to S^{\unip}_{{}^cG,\breve F}$ denote the involution induced by switching two factors $\locsys_{{}^cB,\breve F}^{\unip}\times \locsys_{{}^cB,\breve F}^{\unip}$. By abuse of notations, the induced involution $\sw^{\indcoh,!}$ of $\indcoh(S^{\unip}_{{}^cG,\breve F})$ is still denoted by $\sw$.

\begin{proposition}\label{lem: comparison of two duality for spectral affine Hecke}
We have $\verd^{\indcoh}_{S^{\unip}_{{}^cG,\breve F}}\cong \sw\circ \verd_{S^{\unip}_{{}^cG,\breve F}}^{\mathrm{sr}}  [\dim \hat{T}]$. Concretely, this means that for $\mF\in\Coh(S^{\unip}_{{}^cG,\breve F})$, we have
\[
\verd^{\Coh}_{S^{\unip}_{{}^cG,\breve F}}(\mF) \cong \sw(\mF^\vee)[\dim \hat{T}].
\]
Here $\mF^\vee$ is the (right) dual of $\mF$ with respect to the $!$-monoidal structure of $\indcoh(S_{{}^cG,\breve F}^{\unip})$.
\end{proposition}
\begin{proof}
As explained in \Cref{rem: pivotal structure}, have
\[
\Hom(\Delta_*^{\indcoh} \cohdual_{\locsys_{{}^cB,\breve F}^{\unip}},\mF_1\star\mF_2)=\Hom(\cohdual_{\locsys_{{}^cB,\breve F}^{\unip}}, (\pr_1)_{*}^{\indcoh}(\mF_1\os \sw(\mF_2))).
\]
Here $\pr_1$ denotes the first projection $S^{\unip}_{{}^cG,\breve F}=\locsys_{{}^cB,\breve F}^{\unip}\times_{\locsys_{{}^cG,\breve F}} \locsys_{{}^cB,\breve F}^{\unip}$.
We note that there is a canonical isomorphism
\begin{equation}\label{eq: dualizing vs structural springer}
\cohdual_{\locsys_{{}^cB,\breve F}^{\unip}}\cong \mO_{\locsys_{{}^cB,\breve F}^{\unip}}[-\dim \hat{T}].
\end{equation}
Therefore,
\[
\Hom(\cohdual_{\locsys_{{}^cB,\breve F}^{\unip}}, (\pr_1)_{*}^{\indcoh}(-))=\rg^{\indcoh}(S^{\unip}_{{}^cG,\breve F},-)[\dim \hat{T}].
\]
We thus obtain the first statement from \Cref{rem: pivotal structure}.

As explained in \Cref{ex: self-duality of cpt gen rigid monoidal cat}, the self-duality $\verd_{S^{\unip}_{{}^cG,\breve F}}^{\mathrm{sr}}$, when restricted to the subcategory of compact objects, just sends $\mF$ to its right dual $\mF^\vee$ (with respect to the convolution monoidal structure). The second statement follows.
\end{proof}

Now we drop the assumption that $\bar\tau=1$.
Similarly we have the tame version
$\pi_{\breve F}^\tame: \locsys_{{}^cB,\breve F}^\tame\to \locsys_{{}^cG,\breve F}^\tame$, and the
Steinberg stack $S_{{}^cG,\breve F}^\tame$. We apply \Cref{prop: trace of convolution category} \eqref{prop: trace of convolution category-1.5} to endow $\indcoh(S_{{}^cG,\breve F}^\tame)$ with a monoidal structure, with the monoidal unit given by 
\[
(\Delta_{ \locsys_{{}^cB,\breve F}^\tame/ \locsys_{{}^cG,\breve F}^\tame})_*^{\indcoh}\cohdual_{ \locsys_{{}^cB,\breve F}^\tame}.
\]

\begin{remark}\label{rem:big spectral affine Hecke}
As $S_{{}^cG,\breve F}^\tame$ is an ind-stack, sometimes it is convenient to consider the proper morphism $\hat{B}\bar\tau/\hat{B}\to \hat{G}\bar\tau/\hat{G}$ and $S_{\hat{G}\bar\tau}$. Then $\indcoh(S_{\hat{G}\bar\tau})$ also admits a monoidal structure by  \Cref{prop: trace of convolution category} \eqref{prop: trace of convolution category-1.5}. A choice of a tame generator induces an embedding $\hat\iota: S_{{}^cG,\breve F}^\tame\subset S_{\hat{G}\bar\tau}$ as in \eqref{eq: embedding of restr Steinberg}.
We have a pair of adjoint functors
\[
\hat\iota^{\indcoh}_*: \indcoh(S_{{}^cG,\breve F}^\tame)\rightleftharpoons \indcoh(S_{\hat{G}\bar\tau})\colon \hat\iota^{\indcoh,!}.
\]
One sees easily that $\hat\iota^{\indcoh}_*$ is non-unital monoidal and $\hat\iota^{\indcoh,!}$ is monoidal. 

Again one can apply \Cref{prop: trace of convolution category} \eqref{prop: trace of convolution category-1} to endow $\indcoh(S_{\hat{G}\bar\tau})$ with another monoidal structure, called the $*$-convolution. 
Note that as $(\indcoh,*)$-pullback and the $(\indcoh,!)$-pullback along the diagonal map $\hat{B}\bar\tau/\hat{B}\to \hat{B}\bar\tau/\hat{B}\times \hat{B}\bar\tau/\hat{B}$ coincide, we see that these two monoidal structures actually coincide.
\end{remark} 

We have analogues of \Cref{lem: Coh Steinberg exterior tensor} and \Cref{cor: Coh Steinberg exterior tensor-2}. We record them here.

\begin{lemma}\label{lem: Coh Steinberg exterior tensor-3}
 The exterior tensor product 
 \[
 \boxtimes: \indcoh(S_{{}^cG,\breve F}^\tame)\otimes_\La \indcoh(S_{{}^cG,\breve F}^\tame)\to \indcoh(S_{{}^cG,\breve F}^\tame\times S_{{}^cG,\breve F}^\tame)
 \]
 is an equivalence. For every smooth algebraic stack $Y\in \ArStk_\La^{\aft}$, the exterior tensor product functors
\[
\indcoh(\locsys_{{}^cB,\breve F}^{\tame})\otimes_\La \indcoh(Y)\to \indcoh(\locsys_{{}^cB,\breve F}^{\tame}\times Y)
\]
\[
\indcoh(S_{{}^cG,\breve F}^{\tame})\otimes_\La \indcoh(Y)\to \indcoh(S_{{}^cG,\breve F}^{\tame}\times Y)
\]
are equivalences.
\end{lemma}

We also have the following counterpart of \Cref{lem: comparison of two duality for spectral affine Hecke}.
\begin{proposition}\label{lem: comparison of two duality for spectral affine Hecke-tame}
We have $\verd^{\indcoh}_{S^{\tame}_{{}^cG,\breve F}}\cong \sw\circ \verd_{S^{\tame}_{{}^cG,\breve F}}^{\mathrm{sr}}$. Concretely, this means that for $\mF\in\Coh(S^{\tame}_{{}^cG,\breve F})$, we have
\[
\verd^{\Coh}_{S^{\tame}_{{}^cG,\breve F}}(\mF) \cong \sw(\mF^\vee).
\]
\end{proposition}
\begin{proof}
By \Cref{rem:big spectral affine Hecke}, we may replace $S^{\tame}_{{}^cG,\breve F}$ by $S_{\hat{G}\bar\tau}$. Then the argument of \Cref{lem: comparison of two duality for spectral affine Hecke} applies to $S_{\hat{G}\bar\tau}$.
\end{proof}

\subsubsection{Spectral Deligne-Lusztig induction}\label{SSS: spec DL induction}
We fix a lifting of the Frobenius $\sigma$ in $W_F$, which induces automorphisms of $\locsys_{{}^cB, \breve F}^{\tame}$ and $\locsys_{{}^cG, \breve F}^{\tame}$ respectively, denoted by $\phi$ as before. Clearly, $\pi_{\breve F}^{\tame}$ intertwines these two actions. If we let $\pi_{\breve F}^{\tame}$ as above be $f:X\to Y$ as in \Cref{subsec-trace.convolution.categories}, then
the diagram \eqref{eq:trace-convolution-horocycle-diagram-special} specializes to \eqref{eq:trace-convolution-horocycle-diagram-spectral side} mentioned before.
Note that all stacks belong to $\indarstk_\La^{\aft}$, the morphism
$\pi_{\breve F}^{\tame}: \locsys_{{}^cB, \breve F}^{\tame}\to \locsys_{{}^cG, \breve F}^{\tame}$ is proper, and $\delta^{\tame}$ is representable of finite tor-amplititude.
Therefore, we have the following well-defined functor
\begin{equation}\label{eq: spectral DL induction tame} 
\Ch^{\tame}_{{}^cG,\phi}:=(\tilde\pi^\tame)_*\circ (\delta^\tame)^{!}: \Coh(S^{\tame}_{{}^cG,\breve F})\to \Coh(\locsys_{{}^cG, F}^{\tame}),
\end{equation}
which we call the spectral Deligne-Lusztig induction.
Similarly, there is the unipotent version
\begin{equation}\label{eq: spectral DL induction unip}
\Ch^{\unip}_{{}^cG,\phi}:=(\tilde\pi^\unip)_*\circ (\delta^\unip)^{!}: \Coh(S^{\unip}_{{}^cG,\breve F})\to \Coh(\locsys_{{}^cG, F}^{\tame}).
\end{equation}
By abuse of notations, we will use the same notations for the ind-completion of these two functors.
We note the target of the functor $\Ch^{\unip}_{{}^cG,\phi}$ is still $\Coh(\locsys_{{}^cG, F}^{\tame})$.

\begin{remark}\label{rem: grading of coherent sheaves on loc}
Recall that since $Z_{\hat{G}}^{\Ga_F}\subset \hat{G}$ acts trivially on $\locsys_{{}^cG,F}^{\Box}$, the stack $\locsys_{{}^cG,F}$ is a $Z_{\hat{G}}^{\Ga_F}$-gerbe. It follows that there is a decomposition
\begin{equation}\label{eq: grading of coherent sheaves on loc}
\indcoh(\locsys_{{}^cG,F})=\bigoplus_{\beta\in \xch(Z_{\hat{G}}^{\Ga_F})} \indcoh^\beta(\locsys_{{}^cG,F}).
\end{equation}
See also \cite[\textsection{3.2}]{zhu2020coherent}. 
Such decomposition clearly restricts to a decomposition of $\indcoh(\locsys^{\tame}_{{}^cG,F})$.
Similarly, we have
\begin{equation}\label{eq: grading of coherent sheaves on St}
\indcoh(S_{{}^cG,F}^{\tame})=\bigoplus_{\beta\in \xch(Z_{\hat{G}}^{I_F})} \indcoh^\beta(S_{{}^cG,F}^{\tame}).
\end{equation}
Note that the whole correspondence \eqref{eq:trace-convolution-horocycle-diagram-spectral side} is relative over $\bB \hat{G}$, and the group $Z_{\hat{G}}^{\Ga_F}$ acts trivially on the base change of such correspondence along $\Spec \La\to \bB \hat{G}$.
Therefore, the functor $\Ch_{{}^cG,\phi}^{\tame}$ sends the direct summand $\indcoh^\beta(S_{{}^cG,F}^{\tame})$ to $\indcoh^{\bar\beta}(\locsys_{{}^cG,F})$, where $\bar\beta$ is the image of $\beta$ under the natural map $\xch(Z_{\hat{G}}^{I_F})\to \xch(Z_{\hat{G}}^{\Ga_F})$.
\end{remark}

\begin{lemma}
The spectral Deligne-Lusztig induction functor \eqref{eq: spectral DL induction tame} commute with Grothendieck-Serre duality. When $\bar\tau=1$, the functor \eqref{eq: spectral DL induction unip} commutes with Gorthendieck-Serre duality up to shift by the rank of $\hat{G}$.
\end{lemma}
\begin{proof}
Once we fix a topological generator of the tame inertia, the morphism $\delta^{\tame}$ is the base change along $\Delta_X: \hat{B}\bar\tau/\hat{B}\to \hat{B}\bar\tau/\hat{B}\times \hat{B}\bar\tau/\hat{B}$, which is a relative Gorenstein morphism with trivial relative dualizing sheaf. It follows that $(\delta^{\tame})^{\indcoh,!}=(\delta^{\tame})^{\indcoh,*}$ commutes with the Grothendieck-Serre duality. The morphism $\tilde\pi^{\tame}$ is the base change of $\pi_{\breve F}^{\tame}$ and therefore is proper. Therefore $(\tilde\pi^{\tame})_*$ also commutes with the Grothendieck-Serre duality. The unipotent case can be proved similarly.
\end{proof}

Combining with \Cref{lem: comparison of two duality for spectral affine Hecke}, we see that for $\mF\in \Coh(S_{{}^cG,\breve F}^{\unip})$, there is a canonical isomorphism
\begin{equation}\label{eq: Cohdual Ch vs Ch sw dual unip}
\verd^{\coh}_{\locsys_{{}^cG,F}^{\tame}}(\Ch_{{}^cG,\phi}^{\unip}(\mF))\cong \Ch_{{}^cG,\phi}^{\unip}( \sw(\mF^\vee)).
\end{equation}

We generalize the above isomorphism to the tame case.
\begin{proposition}\label{lem: Cohdual Ch vs Ch sw dual}
Let $\mF\in \indcoh(S_{{}^cG,\breve F}^{\tame})$ be a dualizable object with a right dual $\mF^\vee$ (with respect to the conovolution product). Suppose $\Ch_{{}^cG,\phi}^{\tame}(\mF)\in \Coh(\locsys_{{}^cG,F}^{\tame})$. Then there is a canonical isomorphism
\[
\verd^{\coh}_{\locsys_{{}^cG,F}^{\tame}}(\Ch_{{}^cG,\phi}^{\tame}(\mF))\cong \Ch_{{}^cG,\phi}^{\tame}( \sw(\mF^\vee)).
\]
\end{proposition}
\begin{proof}
Note that we have
\begin{eqnarray*}
\Hom(\Ch_{{}^cG,\phi}^{\tame}(\mF),\mG) & = &\Hom((\Delta_{\locsys_{{}^cB,\breve F}^{\tame}/\locsys_{{}^cG,\breve F}})_* \cohdual_{\locsys_{{}^cB,\breve F}^{\tame}},  (\delta^{\tame})_*^{\indcoh}((\tilde\pi^{\tame})^{\indcoh,!}\mG)\star \mF^\vee)\\
    & =& \Hom(\cohdual_{\locsys_{{}^cB,\breve F}^{\tame}}, (\pr_1)^{\indcoh}_*( (\delta^{\tame})_*^{\indcoh}((\tilde\pi^{\tame})^{\indcoh,!}\mG)\os\mathrm{sw}(\mF^\vee)))\\
    &\stackrel{(1)}{=}& \Hom(\cohdual_{S_{{}^cG,\breve F}^{\tame}}, (\delta^{\tame})_*^{\indcoh}((\tilde\pi^{\tame})^{\indcoh,!}\mG)\os\mathrm{sw}(\mF^\vee))\\
    &\stackrel{(2)}{=}& \Hom(\cohdual_{\widetilde{\locsys}^{\tame}_{{}^cG,F}},  (\tilde\pi^{\tame})^{\indcoh,!}\mG \os  (\delta^{\tame})^{\indcoh,!} \mathrm{sw}(\mF^\vee)) \\
    &\stackrel{(3)}{=}& \rg^{\indcoh}(\locsys_{{}^cG,F}^{\tame}, \mG\os\Ch_{{}^cG,\phi}^{\tame}(\mathrm{sw}(\mF^\vee))).
\end{eqnarray*}
Here 
\begin{itemize}
\item $(1)$ holds since $\pr_1$ is the base change of $\hat{B}\bar\tau/\hat{B}\to \hat{G}\bar\tau/\hat{G}$, which is quasi-smooth with trivial relative dualizing complex so $(\pr_1)^{\indcoh,!}=(\pr_1)^{\indcoh,*}$ is the left adjoint of $(\pr_1)_*^{\indcoh}$;
\item $(2)$ holds by projection formula and the fact $(\delta^{\tame})^{\indcoh,*}=(\delta^{\tame})^{\indcoh,!}$; and
\item $(3)$ holds by projection formula and the fact $\cohdual_{\widetilde{\locsys}^{\tame}_{{}^cG,F}}=\mO_{\widetilde{\locsys}^{\tame}_{{}^cG,F}}$.
\end{itemize}
The proposition is proved.
\end{proof}

For later purpose, we need to understand where some objects go under the functors. We start with introducing a few objects in  $\Coh(S_{{}^cG, \breve F}^{\tame})$.

Let $w\in W_0$. Recall the stack \eqref{eq: w part of Steinberg}.
Let $\cohdual_{S^{\tame}_{{}^cG,\breve F, w}}$ denote its dualizing sheaf, regarded as ind-coherent sheaves on $S^{\tame}_{{}^cG,\breve F}$ via $*$-pushforward along the closed embedding $S^{\tame}_{{}^cG,\breve F,w}\subset S^{\tame}_{{}^cG,\breve F}$. 

Similarly when $\bar\tau=1$, we have the stack \eqref{eq: unipotent steinberg stack w} and $\cohdual_{S^{\unip}_{{}^cG,\breve F, w}}$ and $\mO_{S^{\unip}_{{}^cG,\breve F, w}}$, as coherent sheaves on $S^{\unip}_{{}^cG,\breve F}$.

Next, we have the following commutative diagram
\begin{equation}\label{eq: relative diagonal over BhatG}
\xymatrix{
 \locsys_{{}^cB, \breve F}^{\tame}\ar^-{\Delta_{\locsys_{{}^cB, \breve F}^{\tame}/\locsys_{{}^cG, \breve F}^{\tame}}}[rrr]\ar[d]\ar^{s}[drrr] &&& \locsys_{{}^cB, \breve F}^{\tame}\times_{\locsys_{{}^cG, \breve F}^{\tame}}\locsys_{{}^cB, \breve F}^{\tame}=S_{{}^cG, \breve F}^{\tame}\ar^t[d]\\
  \bB \hat{B}\ar[rrr] &&& \bB \hat{G}.
}
\end{equation}
We let 
\begin{equation}\label{eq:spec central functor}
\mZ^{\spec,\tame}(-)=(\Delta_{\locsys_{{}^cB, \breve F}^{\tame}/\locsys_{{}^cG, \breve F}^{\tame}})^{\indcoh}_* (s^{\indcoh,!} (-))\colon \indcoh(\bB\hat{G})\to \indcoh(S_{{}^cG, \breve F}^{\tame}).
\end{equation}
Note that as $S_{{}^cG, \breve F}^{\tame}$ is just an ind-algebraic stack, the pullback $s^{\indcoh,!}$ does not preserve coherence.
We similarly have the unipotent version $\mZ^{\unip}(-)$, which sends $\Coh(\bB\hat{G})\to \Coh(S_{{}^cG, \breve F}^{\unip})$.
\begin{equation}\label{eq:spec central functor unip}
\mZ^{\spec,\unip}(-)=(\Delta_{\locsys_{{}^cB, \breve F}^{\unip}/\locsys_{{}^cG, \breve F}^{\tame}})^{\indcoh}_* (s^{\indcoh,!} (-))\colon \indcoh(\bB\hat{G})\to \indcoh(S_{{}^cG, \breve F}^{\unip}).
\end{equation}

\begin{notation}\label{notion: line bundle}
We make use of the following notations.
We identify $\la\in\xch(\hat{T})$ with the Picard group of line bundles on $\bB \hat{T}$. 
We let $\mO_Z(\la)\in\Perf(Z)$ denote its $*$-pullback along a map of stacks $Z\to \bB \hat{T}$, if such map is clear from the context.
If $\mF\in \Coh(Z)$, we will write $\mF(\la):=\mO_Z(\la)\otimes \mF$, where $\otimes$ denotes the action of $\Perf(Z)$ on $\Coh(Z)$ as in \Cref{rem: first remark of ind-coh} \eqref{rem: first remark of ind-coh-1}.
\end{notation}

\begin{remark}\label{rem: Borel Weil}
Consider $\hat{G}/\hat{B}\to \bB\hat{B}\to \bB\hat{T}$. Then for $\la\in\xch(\hat{T})$ we have the line bundle $\mO_{\hat{G}/\hat{B}}(\la)$. Note that according to our convention, if $\la$ is dominant (with respect to $\hat{B}$), then $\mO_{\hat{G}/\hat{B}}(w_0\la)$ is semi-ample whose global section is the Schur module (also called the dual of the Weyl module) $V_\la$ of highest weight $\la$. Here as usual $w_0$ denotes the longest element in the absolute Weyl group of $\hat{G}$.
If $\la$ is regular dominant, then $\mO(w_0\la)$ is ample.
\end{remark}

We will apply the above construction to the following set-up.
We consider the map
\[
S^{\tame}_{{}^cG,\breve F,w}\subset S^{\tame}_{{}^cG,\breve F}=\locsys_{{}^cB,\breve F}^{\tame}\times_{\locsys^{\tame}_{{}^cG,\breve F}} \locsys_{{}^cB,\breve F}^{\tame}\xrightarrow{\pr_1} \locsys_{{}^cB,\breve F}^{\tame}\to \locsys_{{}^cT,\breve F}^{\tame} \to \bB\hat{T}.
\]

In the following lemma, we regard $\cohdual_{\locsys_{{}^cB,\breve F}^{\tame}}(\la)$ as an ind-coherent sheaf on $S_{{}^cG,\breve F}^{\tame}$ via $*$-pushforward along the relative diagonal $\Delta_{\locsys_{{}^cB, \breve F}^{\tame}/\locsys_{{}^cG, \breve F}^{\tame}}$.

\begin{lemma}\label{lem:image-of-trace-of-some-objects-in-specHecke}
\begin{enumerate}
\item\label{lem:image-of-trace-of-some-objects-in-specHecke-1} We have
\[
\Ch_{{}^cG,\phi}^{\tame}( \cohdual_{S^{\tame}_{{}^cG,\breve F,w}}\star \cohdual_{\locsys_{{}^cB,\breve F}^{\tame}}(\la))\cong (\widetilde\pi^{\tame}_{w})_*\cohdual_{\widetilde{\locsys}^{\tame}_{{}^cG,F,w}}(\la),
\]

\item\label{lem:image-of-trace-of-some-objects-in-specHecke-2}
  Let $\mF\in \indcoh(S_{{}^cG,F}^{\tame})$. Then we have the canonical isomorphism
\[
\Ch_{{}^cG,\phi}^{\tame}(\mZ^{\spec,\tame}(V)\star \mF)\cong  \widetilde{V}\otimes \Ch_{{}^cG,\phi}^{\tame}(\mF).
\]

\item\label{lem:image-of-trace-of-some-objects-in-specHecke-3} When $\bar\tau=1$, there are similar statements for unipotent versions
\[
\Ch_{{}^cG,\phi}^{\unip}(\cohdual_{S^{\unip}_{{}^cG,\breve F,w}}\star \cohdual_{\locsys_{{}^cB,\breve F}^{\unip}}(\la))\cong (\widetilde\pi^{\unip}_{w})_*\cohdual_{\widetilde{\locsys}^{\unip}_{{}^cG,F,w}}(\la),
\]
\[
\Ch_{{}^cG,\phi}^{\unip}(\mO_{S^{\unip}_{{}^cG,\breve F,w}}[-\dim \hat{T}] \star \cohdual_{\locsys_{{}^cB,\breve F}^{\unip}}(\la))\cong (\widetilde\pi^{\unip}_{w})_*\mO_{\widetilde{\locsys}^{\unip}_{{}^cG,F,w}}(\la),
\]
\[
\Ch_{{}^cG,\phi}^{\unip}(\mZ^{\spec,\unip}(V)\star \mF)\cong  \widetilde{V}\otimes \Ch_{{}^cG,\phi}^{\unip}(\mF).
\]
\end{enumerate}
\end{lemma}
\begin{proof}
We recall that $(\delta^{\tame})^{\indcoh,!}=(\delta^{\tame})^{\indcoh,*}$.
For Part \eqref{lem:image-of-trace-of-some-objects-in-specHecke-1}, we consider the map
\[
S^{\tame}_{{}^cG,\breve F,w}\subset S^{\tame}_{{}^cG,\breve F}=\locsys_{{}^cB,\breve F}^{\tame}\times_{\locsys^{\tame}_{{}^cG,\breve F}} \locsys_{{}^cB,\breve F}^{\tame}\xrightarrow{\pr_2} \locsys_{{}^cB,\breve F}^{\tame}\to \locsys_{{}^cT,\breve F}^{\tame} \to \bB\hat{T}.
\]
(Note the projection is to the second factor. See \Cref{rem: first projection vs second projection}.) 
Using \Cref{notion: line bundle}, there are the canonical isomorphisms
\[
\cohdual_{S^{\tame}_{{}^cG,\breve F,w}}\star \cohdual_{\locsys_{{}^cB,\breve F}^{\tame}}(\la)\cong  \cohdual_{S^{\tame}_{{}^cG,\breve F,w}}(\la).
\] 
Now Part \eqref{lem:image-of-trace-of-some-objects-in-specHecke-1} follows from definitions.  

For Part \eqref{lem:image-of-trace-of-some-objects-in-specHecke-2}, notice that the whole correspondence \eqref{eq:trace-convolution-horocycle-diagram-spectral side} is over $\bB \hat{G}$ (in fact over $\locsys_{{}^cG,\breve F}^{\tame}$), we see that $\Ch_{{}^cG,\phi}^{\tame}$ are $\Perf(\bB\hat{G})$-linear. Note that the $\Perf(\bB\hat{G})$-linear structures on $\Coh(\locsys_{{}^cG,F}^{\tame})$ is given by 
\[
\Perf(\bB\hat{G})\otimes_\La \Coh(\locsys_{{}^cG,F}^{\tame})\to \Coh(\locsys_{{}^cG,F}^{\tame}),\quad (V, \mF)\mapsto \widetilde{V}\otimes\mF,
\]
and similarly on $\Coh(S_{{}^cG,F}^{\tame})$ is given by
\[
\Perf(\bB\hat{G})\otimes_\La \Coh(S_{{}^cG,F}^{\tame})\to \Coh(S_{{}^cG,F}^{\tame}), \quad (V,\mF)\mapsto t^* V\otimes \mF.
\]
But we have the canonical isomorphism 
\begin{equation}\label{eq: convolution by central vs tensoring}
t^* V\otimes \mF\cong \mZ^{\tame}(V)\star \mF,
\end{equation}
giving the desired isomorphism.

Part \eqref{lem:image-of-trace-of-some-objects-in-specHecke-3} is proved similarly.
\end{proof}

\begin{example}\label{ex: spectral DL induction w equal 1}
In particular, if $w=1$,  we see that
\[
\Ch_{{}^cG,\phi}^{\tame}(\cohdual_{\locsys_{{}^cB,\breve F}^{\tame}}(\la))\cong (\pi^{\tame})_*\cohdual_{\locsys^{\tame}_{{}^cB,F}}(\la),\quad \Ch_{{}^cG,\phi}^{\unip}(\cohdual_{\locsys_{{}^cB,\breve F}^{\unip}}(\la))\cong (\pi^{\unip})_*\cohdual_{\locsys^{\unip}_{{}^cB,F}}(\la).
\]
Recall that $\cohdual_{\locsys^{\tame}_{{}^cB,F}}\cong \mO_{\locsys^{\tame}_{{}^cB,F}}$ (e.g. see \cite[Proposition 2.3.7]{zhu2020coherent}), and
\begin{equation}\label{eq: def of coherent Springer sheaf tame}
\cohspr^{\tame}_{{}^cG,F}:=(\pi^{\tame})_*\cohdual_{\locsys^{\tame}_{{}^cB,F}}\cong(\pi^{\tame})_*\mO_{\locsys^{\tame}_{{}^cB,F}} , 
\end{equation}
is called the tame coherent Springer sheaf (\cite[\textsection{4.4}]{zhu2020coherent}). 
When $\bar\tau=1$, we also have $\cohdual_{\locsys^{\unip}_{{}^cB,F}}\cong \mO_{\locsys^{\unip}_{{}^cB,F}}$, and the unipotent version of the coherent Springer sheaf
\begin{equation}\label{eq: def of coherent Springer sheaf unip}
\cohspr^{\unip}_{{}^cG,F}:=(\pi^{\unip})_*\cohdual_{\locsys^{\unip}_{{}^cB,F}}\cong (\pi^{\unip})_*\mO_{\locsys^{\unip}_{{}^cB,F}}.
\end{equation}
In particular, we can say that the spectral Deligne-Lusztig induction sends the unit object of the spectral affine Hecke category to the coherent Springer sheaf.
\end{example}

Now,  suppose $\bar\tau=1$.
Let $\bfC\subset \indcoh(\locsys^{\tame}_{{}^cG,F})$ be the $\La$-linear presentable stable subcategory generated by the essential image of $\Ch^{\unip}_{{}^cG,\phi}$. It is known that $\indcoh(S_{{}^cG,\breve F}^{\unip})$ is generated as $\La$-linear presentable category by objects $\mO_{\locsys_{{}^cB,\breve F}^{\tame}}(\la)\star\mO_{S^{\tame}_{{}^cG,\breve F,w}}$ for $\la\in\xcoch(T)$ and $w\in W_0$ as in \Cref{lem:image-of-trace-of-some-objects-in-specHecke} \eqref{lem:image-of-trace-of-some-objects-in-specHecke-1}. (This for example follows from the Bezrukavnikov equivalence and a description of a set of generators of the affine Hecke category.) Therefore $\bfC$ is generated as $\La$-linear presentable stable category by objects $(\widetilde\pi^{\tame}_{w})_*\mO_{\widetilde{\locsys}^{\tame}_{{}^cG,F,w}}(\la)$. We expect that $\bfC= \indcoh(\locsys^{\widehat\unip}_{{}^cG,F})$, at least when the characteristic of $\La$ is not too small. Currently, the following weaker result is sufficient for many applications.

\begin{proposition}\label{lem: image of spectral DL induction}
Let $\La$ be a Dedekind domain, which is either an algebraic field extension of $\bQ_\ell$ or $\bF_\ell$, or a finite extension of $\bZ_\ell$.
Then $\bfC\subset \indcoh(\locsys^{\widehat\unip}_{{}^cG,F})$ is stable under the action of $\ind\Perf(\locsys^{\widehat\unip}_{{}^cG,F})$.
Assume that the derived group of $\hat{G}$ is simply-connected. Then we have
\[
\ind\Perf(\locsys^{\widehat\unip}_{{}^cG,F})\subset \bfC\subset \indcoh(\locsys^{\widehat\unip}_{{}^cG,F}).
\] 
In particular, $\cohdual_{\locsys^{\widehat\unip}_{{}^cG,F}}\otimes \widetilde{V}$ belongs to $\bfC$ for every $V\in \rep(\hat{G})$.
\end{proposition}

\begin{proof}
As the morphism $\widetilde{\locsys}^{\unip}_{{}^cG,F}\to \locsys^{\tame}_{{}^cG,F}$ factors through the connected component $\locsys^{\widehat\unip}_{{}^cG,F}\subset  \locsys^{\tame}_{{}^cG,F}$, we see that $\bfC\subset \indcoh(\locsys^{\widehat\unip}_{{}^cG,F})$.

Note that for every $\mF\in \indcoh(\widetilde\loc^{\unip}_{{}^cG,F})$ and $\mE\in \ind\Perf(\locsys_{{}^cG,F}^{\tame})$, we have
\[
\mE\otimes (\tilde\pi^{\unip})_*^{\indcoh}\mF\cong (\tilde\pi^{\unip})_*^{\indcoh}((\tilde\pi^{\unip})^*\mE\otimes \mF).
\]
Therefore, to show that $\bfC$ is stable under the action of $ \ind\Perf(\locsys_{{}^cG,F}^{\tame})$, it is enough to show that if $\mF$ is contained in the subcategory of $ \indcoh(\widetilde\loc^{\unip}_{{}^cG,F})$ generated by $(\delta^{\unip})^{\indcoh,*}(\indcoh(S^{\unip}_{{}^cG,\breve F}))$, so is $(\tilde\pi^{\unip})^*\mE\otimes \mF$. But this follows from  \Cref{lem: control support ofr deltaunip} below.

We next show that when  the derived group of $\hat{G}$ is simply-connected, then $\ind\Perf(\locsys^{\widehat\unip}_{{}^cG,F})\subset \bfC$.
Given that $\bfC$ is stable under the $\ind\Perf(\locsys^{\widehat\unip}_{{}^cG,F})$-action, it is enough to show $\cohdual_{\locsys^{\widehat\unip}_{{}^cG,F}}=\mO_{\locsys^{\widehat\unip}_{{}^cG,F}}$ belongs to $\bfC$.

As the derived group of $\hat{G}$ is simply-connected,  the Chevalley map $\hat{G}/\hat{G}\to \hat{G}/\!\!/\hat{G}$ is flat. Let
\[
\mU_{\hat{G}}=\hat{G}/\hat{G}\times_{\hat{G}/\!\!/\hat{G}}\{1\}.
\]
The base change of $\mU_{\hat{G}}$ to a field is the variety of unipotent elements of $\hat{G}$ as in \Cref{rem: three versions of unipotent stack}.
In addition, the Springer map $f$ factors as $\hat{U}/\hat{B}\to \mU_{\hat{G}}/\hat{G}\to \hat{G}/\hat{G}$ and the $*$-pushfoward of the $\mO_{\hat{U}/\hat{B}}$ along the first map is $\mO_{\mU_{\hat{G}}/\hat{G}}$. It follows that 
$f_*\cohdual_{\hat{U}/\hat{B}}=\cohdual_{\mU_{\hat{G}}/\hat{G}}$, where $\cohdual_{\mU_{\hat{G}}/\hat{G}}$ is regarded as a coherent sheaf on $\hat{G}/\hat{G}$ via the $*$-pushforward along the closed embedding $\mU_{\hat{G}}/\hat{G}\to \hat{G}/\hat{G}$.

Recall that once a topological generator $\tau$ of the tame inertia is chosen, the proper morphism $\tilde\pi^{\unip}: \widetilde{\locsys}^{\unip}_{{}^cG,F}\to \locsys^{\tame}_{{}^cG,F}$
is a base change of $f$.  It follows by base change that
\begin{equation}\label{eq: pushforward of omegaunip from tildelocunip}
(\tilde\pi^{\unip})_*\cohdual_{\widetilde{\locsys}^{\unip}_{{}^cG,F}}\cong i_*\cohdual_{\locsys^{\unip}_{{}^cG,F}}.
\end{equation}
Here we let $i: \locsys^{\unip}_{{}^cG,F}\to \locsys^{\tame}_{{}^cG,F}$ denote the closed embedding. 
Note that in $\indcoh( \locsys^{\tame}_{{}^cG,\breve F})$, $\cohdual_{\mU_{\hat{G}}^\wedge/\hat{G}}$ is in the $\La$-linear category generated by $\cohdual_{\mU_{\hat{G}}/\hat{G}}$, and the $!$-pullback of $\cohdual_{\mU_{\hat{G}}^\wedge/\hat{G}}$ along $\res: \locsys^{\tame}_{{}^cG,F}\to  \locsys^{\tame}_{{}^cG,\breve F}$ is $\cohdual_{\locsys^{\widehat\unip}_{{}^cG,F}}$, we see that $\cohdual_{\locsys^{\widehat\unip}_{{}^cG,F}}$ is contained in the $\La$-linear subcategory of $\indcoh(\locsys^{\widehat\unip}_{{}^cG,F})$ generated by $i_*\cohdual_{\locsys^{\unip}_{{}^cG,F}}$. Note that $\cohdual_{\locsys^{\widehat\unip}_{{}^cG,F}}$ itself is perfect, by \Cref{lem: compact object as a colimit} we see that it is in fact contained in the idempotent complete subcategory generated by $i_*\cohdual_{\locsys^{\unip}_{{}^cG,F}}$. The proposition is proved.
\end{proof}

\begin{lemma}\label{lem: control support ofr deltaunip}
The essential image of the functor $(\delta^{\unip})^{*}: \Coh(S_{{}^cG,\breve F}^{\unip})\to \Coh(\widetilde\locsys^{\unip}_{{}^cG,F})$
generates 
\[
\Coh_{\Sing(\delta^{\unip})(\Sing(S_{{}^cG,\breve F}^{\unip})_{\widetilde{\locsys_{{}^cG,F}^{\unip}}})}(\widetilde\locsys^{\unip}_{{}^cG,F})
\] 
as idempotent complete stable category.
\end{lemma}
\begin{proof}
When $\La$ is a field of characteristic zero, this follows from \Cref{prop: control functor image by singular support stacks}. (But notice that  \Cref{prop: control functor image by singular support stacks} fails in general by virtue of \Cref{rem: failure of generation by smooth pullback}.) The argument below works for more general base $\La$.

We write $X=\hat{U}/\hat{B}$ and $Y=\hat{G}/\hat{G}$.
As before, by choosing a generator of the tame inertia, we write $S_{{}^cG,F}^{\unip}=X\times_YX$ and $\widetilde\locsys^{\unip}_{{}^cG,F}=X\times_{Y\times Y}Y$. The map factors as
\[
X\times_{Y\times Y} Y\to (X\times_{\bB\hat{B}} X)\times_{Y\times Y} Y\to (X\times X)\times_{Y\times Y}Y=X\times_YX.
\]
We note that $(X\times_{\bB\hat{B}} X)\times_{Y\times Y} Y=(\hat{U}\times \hat{U})/\hat{B} \times_{\hat{G}/\hat{G}\times\hat{G}/\hat{G}}\hat{G}/\hat{G}$.
The first map is a quasi-smooth closed embedding between quasi-smooth algebraic stacks, induced by thel map $\hat{U}\xrightarrow{\id\times\phi} \hat{U}\times\hat{U}$. 
The second morphism is the base change of a smooth morphism $\bB \hat{B}\xrightarrow{\id\times\phi}\bB \hat{B}\times \bB\hat{B}$, and therefore is smooth.
Using \Cref{lem: pullback quasi-smooth closed embedding}, the desired statement follows from the control of the image of the $*$-pullback functor along the second smooth morphism, as given in the following lemma.
\end{proof}

\begin{lemma}
The essential image of the $*$-pullback functor
\[
\indcoh(X\times_YX)\to \indcoh( (X\times_{\bB\hat{B}} X)\times_{Y\times Y} Y)
\]
generate $\indcoh( (X\times_{\bB\hat{B}} X)\times_{Y\times Y} Y)$ as presentable $\La$-linear category. 
\end{lemma}

\begin{proof}

We have the natural maps $(X\times_{\bB\hat{B}} X)\times_{Y\times Y} Y\to X\times_YX=S_{{}^cG,F}^{\unip}\to \hat{B}\bs \hat{G}/\hat{B}$. For each $w$,  let $Z_w\subset S_{{}^cG,F}^{\unip}$ be the (reduced) locally substack of $S_{{}^cG,F}^{\unip}$ corresponding to $w$ as in \Cref{rem: two different numeration of Sunip by W}, and let $\widetilde{Z}_w$ be the preimage of $Z_w$ in $(X\times_{\bB\hat{B}} X)\times_{Y\times Y} Y$. 
Then we have
\[
Z_w\cong \frac{\hat{U}_w}{\Ad \hat{B}_w},\quad \widetilde{Z}_w\cong \frac{\hat{U}_w\times\hat{B}}{\Ad_w \hat{B}_w}.
\]
Here the action $\Ad_w$ of $\hat{B}_w$ on the first factor $\hat{U}_w$ is still the usual conjugation action but on the second factor $\hat{B}$ is given to $b\cdot b'= (\dot{w}^{-1}b\dot{w}) b' b^{-1}$.
Using \Cref{prop: open-closed gluing coh} (together with \Cref{lem: coh sheaf closed embedded}), the lemma is a consequence of the following statement.
\end{proof}

\begin{lemma}\label{lem: generation by O(lambda)}
The image of the $*$-pullback functor $\Perf(\bB \hat{B}_w)\to \Perf(\hat{U}_w/\Ad\hat{B}_w)$ generates the target as an idempotent complete category.
The image of the $*$-pullback functor $\Perf(\bB \hat{B}_w)\to \Perf((\hat{U}_w\times \hat{B})/\Ad_w\hat{B}_w)$ generates the target as an idempotent complete category.
\end{lemma}

\begin{proof}
The first statement has been proved in the course of proving \Cref{lem: Coh Steinberg exterior tensor}. The argument for 
the second statement is very similar. We only briefly explain needed modifications.

Let us write $\widetilde{Z}_w^{\Box}=\hat{U}_w\times \hat{B}$.
Note that it has a group structure. We consider the map
\[
\widetilde{Z}_w^{\Box}\times \widetilde{Z}_w^{\Box}\mapsto \widetilde{Z}_w^{\Box},\quad (z_1,z_2)\mapsto z_1^{-1}z_2.
\]
This morphism is $\hat{B}_w$-equivariant, where now $\hat{B}_w$ acts on the left diagonally as before, but on the target $\widetilde{Z}_w^{\Box}$ by usual conjugation. 
 In addition, the diagonal of $\widetilde{Z}_w^{\Box}\times \widetilde{Z}_w^{\Box}$ is the preimage of $\{1\}\in \widetilde{Z}_w^{\Box}$. Therefore, it is enough to show that
the $*$-pushforward of the structure sheaf $\mO_{\bB \hat{B}_w}$ along closed embedding 
\[
\bB \hat{B}_w=\frac{\{1\}}{\Ad \hat{B}}\to \frac{\widetilde{Z}_w^{\Box}}{\Ad \hat{B}_w}
\] 
belongs to the subcategory of $\Perf(\frac{\widetilde{Z}_w^{\Box}}{\mathrm{Ad} \hat{B}_w})$  generated by $\Perf(\bB \hat{B}_w)$ (under pullbacks). One then proceeds as in the proof of \Cref{lem: Coh Steinberg exterior tensor} by filtering $\widetilde{Z}_w^{\Box}$ as $\hat{B}_w$-conjugate invariant normal subgroups to conclude.
\end{proof}
We also consider the tame version of \Cref{lem: control support ofr deltaunip} 
We can drop the assumption $\bar\tau=1$.
First, the proof of \Cref{lem: control support ofr deltaunip} works in the tame setting with obvious modifications, giving the following. (To avoid working with ind-stacks, one can choose a generator of the tame inertia $\tau$ and work with $S_{\hat{G}\bar\tau}$ instead of $S_{{}^cG,F}^{\tame}$ as in \Cref{prop: derived geometry of Steinberg stack}.)

\begin{lemma}\label{lem: control support ofr deltatame}
The essential image of the functor $(\delta^{\tame})^{*}: \Coh(S_{{}^cG,F}^{\tame})\to \Coh(\widetilde\locsys^{\tame}_{{}^cG,F})$
generates 
\[
\Coh_{\Sing(\delta^{\tame})(\Sing(S_{{}^cG,\breve F}^{\tame})_{\widetilde{\locsys_{{}^cG,F}^{\tame}}})}(\widetilde\locsys^{\tame}_{{}^cG,F})
\] 
as idempotent complete stable category.
\end{lemma}

\subsubsection{Cateogrical trace computation}\label{SSS: cat trace spectral side}

Now we can state the outcome of the computation of the categorical trace of the spectral affine Hecke category. 

\begin{theorem}\label{prop:trace-spectral-affine-Hecke-category}
Assume that $\La=\overline\bQ_\ell$. 
\begin{enumerate}
\item\label{prop:trace-spectral-affine-Hecke-category-1} There is the following commutative diagram  with the bottom arrow an equivalence
\[
\xymatrix{
\indcoh(S^{\tame}_{{}^cG,\breve F})\ar[d]\ar[r] &\indcoh(\widetilde{\locsys}^{\tame}_{{}^cG,F})\ar[d]\\
\tr(\indcoh(S^{\tame}_{{}^cG, \breve F}),\phi)\ar^\cong[r] & \indcoh(\locsys^{\tame}_{{}^cG,F}).
}
\]
Suppose $\zeta$ is a tame inertia type. Then the above diagram restricts to commutative diagram with the bottom arrow an equivalence
\[
\xymatrix{
\indcoh(S^{\hat\zeta}_{{}^cG,\breve F})\ar[d]\ar[r] &\indcoh(\widetilde{\locsys}^{\hat\zeta}_{{}^cG,F})\ar[d]\\
\tr(\indcoh(S^{\hat\zeta}_{{}^cG, \breve F}),\phi)\ar^\cong[r] & \indcoh(\locsys^{\hat\zeta}_{^{c}G,F}).
}\]

\item Assume that $\bar\tau=1$. We also have a canonical equivalence
\[
\tr(\indcoh(S_{{}^cG,\breve F}^{\unip}),\phi) \cong \indcoh(\locsys^{\widehat\unip}_{^{c}G,F}),
\]
fitting into a commutative diagram as the one in Part \eqref{prop:trace-spectral-affine-Hecke-category-1}.
\end{enumerate}

Next assume that $\La=\overline\bF_\ell$.
\begin{enumerate}[resume]
\item\label{prop:trace-spectral-affine-Hecke-category-3} Assume that $\bar\tau=1$. Then there is a fully faithful embedding
\[
\tr(\indcoh(S_{{}^cG,\breve F}^{\unip}),\phi) \hookrightarrow \indcoh(\locsys^{\widehat\unip}_{^{c}G,F}).
\]
fitting into a commutative diagram as the one in Part \eqref{prop:trace-spectral-affine-Hecke-category-1}. The essential image is stable under the $\ind\Perf(\locsys^{\widehat\unip}_{^{c}G,F})$-action.
In addition, if  the derived group $\hat{G}$ is simply-connected, then $\ind\Perf(\locsys^{\widehat\unip}_{^{c}G,F})\subset \tr(\indcoh(S_{{}^cG,\breve F}^{\unip}),\phi)$.
\end{enumerate}
\end{theorem}

\begin{proof}
We apply \Cref{prop: trace of convolution category} by letting  $X\to Y$ be as in $\locsys_{{}^cB,\breve F}^{\tame}\to \locsys_{{}^cG,\breve F}^{\tame}$. We thus obtain the fully faithful embedding. The essential surjectivity follows from \Cref{prop: control functor image by singular support stacks} and the calculation made in \Cref{lem: pull-push singular support}.

More precisely, we will fix $\iota: \Ga_q\to W_F^t$ as before and consider the $\iota$-version of \eqref{eq:trace-convolution-horocycle-diagram-spectral side}
\begin{equation*}
\xymatrix{
 \ar_-{\tilde{\pi}}[d] \widetilde{\locsys}_{{}^cG,F,\iota}^{\tame}:=  \locsys_{{}^cG, F,\iota}^{\tame}\times_{\hat{G}\bar\tau/\hat{G}}\hat{B}\bar\tau/\hat{B}\ar^-{\delta}[r] &  S_{\hat{G}\bar\tau}=\hat{B}\bar\tau/\hat{B} \times_{\hat{G}\bar\tau/\hat{G}}\hat{B}\bar\tau/\hat{B}\\
 \locsys_{{}^cG, F,\iota}^{\tame}.   &
 }
\end{equation*}
We note that all the stacks in the above diagram are global complete intersection stack in the sense of \cite[\textsection{9.2}]{arinkin2015singular} so \Cref{prop: control functor image by singular support stacks} is indeed applicable. In addition, the map $\delta$ in the above diagram factors through $\indcoh(S_{{}^cG,\breve F}^{\tame})$ the essential image of $\Ch^{\tame}_{{}^cG,\phi}$ and its $\iota$-version coincide. This gives the equivalence $\tr(\indcoh(S^{\tame}_{{}^cG, \breve F}),\phi)\cong \indcoh(\locsys^{\tame}_{{}^cG,F})$. The rest equivalences in are similar.

Fully faithfulness of Part \eqref{prop:trace-spectral-affine-Hecke-category-3} still follows from \Cref{prop: trace of convolution category}. The rest statements follow from  \Cref{lem: image of spectral DL induction}.
\end{proof}

\begin{proposition}\label{prop: spectral side identifying duality}
Under the canonical equivalence from \Cref{prop:trace-spectral-affine-Hecke-category}, the self-duality of $\tr(\indcoh(S^{\tame}_{{}^cG,\breve F}),\phi)$ induced by the one on $\indcoh(S^{\tame}_{{}^cG,\breve F})$ is canonically identified with the modified Grothendieck-Serre duality of $\indcoh(\locsys_{{}^cG,F}^{\tame})$ from \eqref{eq: modified GS duality}.
\end{proposition}

On the other hand, recall that if $\bfM$ is a (left) dualizable $\indcoh(S_{{}^cG,\breve F}^{\tame})$-module, equipped with a left module functor $\phi_{\bfM}: \bfM\to {}^\phi\bfM$, then the map \eqref{eq:class-of-modules}
defines an object 
\[
[\bfM, \phi_\bfM]_{{}^\phi\indcoh(S_{{}^cG,\breve F}^{\tame})}\in \tr(\indcoh(S_{{}^cG,\breve F}^{\tame}),\phi)=\indcoh(\locsys^{\tame}_{{}^cG,F}).
\]
By abuse of notations, we will denote this object by $\Ch^{\tame}_{{}^cG,\phi}(\bfM,\phi_\bfM)$, although this is not really the spectral Deligne-Lusztig induction of a sheaf.

Similarly, if $\bfM$ is a left $\indcoh(S_{{}^cG,\breve F}^{\unip})$-module,  equipped with a left module functor $\phi_{\bfM}: \bfM\to {}^\phi\bfM$, then we write $[\bfM, \phi_\bfM]_{{}^\phi\indcoh(S_{{}^cG,\breve F}^{\unip})}$ by $\Ch^{\unip}_{{}^cG,\phi}(\bfM,\phi_\bfM)$.

The case we will be interested in will be 
\[
\bfM=\indcoh(\locsys_{{}^cB,\breve F}^{\tame}\times_{\locsys_{{}^cG,\breve F}^{\tame}}W),
\]
where $W$ is an ind-algebraic stack almost of finite presentation equipped with a map $g:W\to \locsys_{{}^cG,\breve F}^{\tame}$. 
In all the situations considered below, it is easy to see that both $g$ and the relative diagonal $W\to W\times_{\locsys_{{}^cG,\breve F}^{\tame}} W$ are proper.
We suppose $W$ is equipped with an automorphism $\phi=\phi_W$ and an isomorphism $g\circ \phi\cong \phi\circ g$.
Let 
\[
\mL_\phi g: \mL_\phi W\to \mL_\phi  \locsys_{{}^cG,\breve F}^{\tame}=\locsys_{{}^cG,F}^{\tame}
\] 
be the map induced by $g$, which is also ind-proper by \Cref{lem: loop of f in verti}.

\begin{proposition}\label{prop: class of diagonal}
Let $\bfM=\indcoh(\locsys_{{}^cB,\breve F}^{\tame})$ equipped with the natural $\phi$-structure. We regard $\bfM$ as a left $\indcoh(S_{{}^cG,\breve F}^{\tame})$-module by convolution.
Then
\[
\Ch^{\tame}_{{}^cG,\phi}(\bfM,\phi_\bfM)\cong \cohdual_{\locsys^{\tame}_{{}^cG,F}}.
\]
Similarly, for $\bfM=\indcoh(\locsys_{{}^cB,\breve F}^{\unip})$ equipped with the natural $\phi$-structure. We regard $\bfM$ as a left $\indcoh(S_{{}^cG,\breve F}^{\unip})$-module by convolution. Then
\[
\Ch^{\unip}_{{}^cG,\phi}(\bfM,\phi_\bfM)\cong \cohdual_{\locsys^{\widehat\unip}_{{}^cG,F}}.
\]
\end{proposition}
\begin{proof}
Note that thanks to \Cref{cor: Coh Steinberg exterior tensor-2} and \Cref{lem: Coh Steinberg exterior tensor-3},  \Cref{ex: geo phi-trace class} is applicable, giving the proposition.

To say a little bit more in the second case, we notice that $\tr(\indcoh(S_{{}^cG,\breve F}^{\unip}),\phi)\subset \indcoh(\locsys^{\tame}_{{}^cG,F})$, and $\cohdual_{\locsys^{\widehat\unip}_{{}^cG,F}}$ is contained in the essential image of $\Ch^{\unip}_{{}^cG,\phi}$. 
Then $\proj_{\trg}(\cohdual_{\locsys^{\tame}_{{}^cG,F}})=\cohdual_{\locsys^{\widehat\unip}_{{}^cG,F}}$.
\end{proof}

Next, let
\[
\bfM_{{}^cP}=\indcoh(\locsys_{{}^cB,\breve F}^{\tame}\times_{\locsys_{{}^cG,\breve F}^{\tame}} \locsys_{{}^cP,\breve F}^{\tame}\times_{\locsys_{{}^cM,\breve F}^{\tame}} \locsys_{{}^cM,\breve F}^{\unr}),
\]
with the natural $\phi$-structure.

\begin{proposition}\label{prop: class of parabolic spectral}
Assume that $\La=\overline\bQ_\ell$. Then 
\[
\Ch^{\tame}_{{}^cG,\phi}(\bfM_{{}^cP},\phi)\cong \pi_* (\cohdual_{\locsys_{{}^cP,F}\times_{\locsys_{{}^cM,F}}\locsys_{{}^cM,F}^{\unr}}),
\]
where the map $\pi$ is from \eqref{E:LocMtoG}. In particular, when ${}^cP={}^cG$, we have
\[
\Ch^{\tame}_{{}^cG,\phi}(\bfM,\phi_\bfM)\cong \cohdual_{\locsys^{\unr}_{{}^cG,F}}.
\]
\end{proposition}
\begin{proof}
This is again a consequence of \Cref{ex: geo phi-trace class}. 
\end{proof}

\subsubsection{Excursion algebra}
We recall the formulation of excursion algebra/$S$-operators \`a la Vincent Lafforgue \cite{Lafforgue2018Chtoucas} in the spectral side.  We follow the approach of \cite{zhu2020coherent}. We will fix $\iota: \Ga_q\to W_F^t$ and let $\Ga_{F,\iota}$ be defined as in \eqref{eq: GaFiota}.

Let $\mathbf{FFM}$ be the category of finitely generated free monoids. For a finite set $I$, let $\mathbf{FM}(I)$ be the free monoid generated by $I$. An $I$-uple $\ga^I\in (\Ga_{F,\iota})^I$ can be regarded as a homomorphism $\mathrm{FM}(I)\to \Ga_{F,\iota}$, inducing a map $\locsys_{{}^cG,F,\iota}\to ({}^cG)^I/\hat{G}$. Explicitly, these maps send a Langlands parameter $\varphi$ to $(\varphi(\ga_i))_{i\in I}\in {}^cG^I/\hat{G}$. They induce maps of ring of regular functions 
\[
 \chi_{(\ga_i)_i}: \La[({}^cG)^I]^{\hat{G}}\to Z_{{}^cG,F}=H^0\rg(\locsys_{{}^cG, F,\iota},\mO).
\]
Note that as the ring of regular functions of $\hat{G}$ as a $\hat{G}$-representation by conjugation action admits a good filtration, taking $\hat{G}$-invariants of $\La[({}^cG)^I]$ does not have higher cohomology.
These maps are compatible with homomorphisms $\mathrm{FM}(I)\to \mathrm{FM}(J)$. Therefore, they together induce a ring map
\[
H^0(\colim_{\mathbf{FFM}/\Ga_{F,\iota}} \La[{}^cG ^I]^{\hat{G}})\to Z_{{}^cG,F}.
\]
Here as the slice category $\mathbf{FFM}/\Ga_{F,\iota}$ is not filtered, the colimit on the left hand side might have derived structure. But we will only need its degree zero part. The algebra on the left hand side is usually called the excursion algebra. Its geometric points classify closed $\hat{G}$-orbits in $R_{\Ga_{F,\iota}, {}^cG}$ (see \Cref{rem: repn space of discrete group} for the space $R_{\Ga_{F,\iota}, {}^cG}$).

One also has the framed version: given $\mathrm{FM}(I)\to \Ga_{F,\iota}$, we have  $\locsys_{{}^cG,F,\iota}^{\Box}\to ({}^cG)^I$, induces
\[
 \colim_{\mathbf{FFM}/\Ga_{F,\iota}} \La[({}^cG)^I] \cong \La[R_{\Ga_{F,\iota}, {}^cG}] \twoheadrightarrow \La[\locsys_{{}^cG, F,\iota}^{\Box}].
\]
Here the first isomorphism follows from \cite[Proposition 2.2.3]{zhu2020coherent}.

Now, let $W$ be a representation of ${}^cG$ on a finite projective $\La$-module and let $W^*$ be its dual representation. Let $m_W: W^*\otimes W\to \La[{}^cG^I]$ be the matrix coefficient map. We let 
\[
\chi_{W,(\ga_i)_i}=\chi_{(\ga_i)_i}(m_W(u_W))\in Z_{{}^cG,F},
\] 
where  $u_W\in W^*\otimes W$ is the unit of the duality datum of $W$.

Now we restrict to tame and unipotent part. Let $V\in \rep(\hat{G})$, and let $\mZ^{\tame}(V)\in \indcoh(S_{{}^cG,\breve F}^{\tame})$ be as in  \eqref{eq:spec central functor}. Note that for every $\mF\in \indcoh(S_{{}^cG,\breve F}^{\tame})$, there are canonical isomorphisms 
\begin{equation}\label{eq:centrality of spec central sheaf}
\mF\star \mZ^{\tame}(V)\cong t^*V\otimes \mF\cong \mZ^{\tame}(V)\star \mF,
\end{equation}
where the map $t$ is as in \eqref{eq: relative diagonal over BhatG} (see \eqref{eq: convolution by central vs tensoring}).

Now suppose $V\in ({}^cG)^I$.  Note that the morphism $s: \locsys_{{}^cB,\breve F}^{\tame}\to \bB \hat{G}$ in \eqref{eq: relative diagonal over BhatG} factors through $\locsys_{{}^cB,\breve F}^{\tame}\to \locsys_{{}^cG,\breve F}^{\tame}$. It follows from the discussions in \Cref{rem: ev bundle on the stack of geometric parameters} that $\mZ^{\tame}(V)$ is equipped with an action 
\begin{equation}\label{eq: taut action of spectral central sheaf}
\varphi^{\mathrm{univ}}: (I_F^t)^I \to \End(\mZ^{\tame}(V))
\end{equation}
In addition, for each $i$ there is a canonical isomorphism
\begin{equation}\label{eq: phi equivariance of spectral central sheaf}
\Phi_i: \mZ^{\tame}(V)\cong \phi_*(\mZ^{\tame}(V)).
\end{equation}

 Now let $I=\{1,2\}$. For every $\ga\in \tau^{\bZ\frac[{1}{p}]}$, we define a map $\eta_{\ga}$ as in \eqref{eq:input for S operator} as
 \[
 \mF\star\mZ^{\tame}(V)\stackrel{\eqref{eq:centrality of spec central sheaf}}{\cong} \mZ^{\tame}(V)\star\mF\stackrel{\varphi^{\mathrm{univ}}(\ga,1)}{\cong}\mZ^{\tame}(V)\star\mF\stackrel{\Phi_2}{\cong} \phi_*(\mZ^{\tame}(V))\star\mF.
 \]
It follows from the abstract construction  \eqref{eq: abstract S-operator} that there is the $S$-operator  
\begin{equation}\label{eq: tame central S operator}
S_{(\mZ^{\mon}(V), \eta_\ga)}: \Ch_{{}^cG, \phi}^{\tame}(\mF)\to \Ch_{{}^cG, \phi}^{\tame}(\mF).
\end{equation}

\begin{lemma}\label{lem: spectral S operator vs abstract S operator}
Let $\mF, V$ be as above. 
Then the endomorphism of $\Ch_{{}^cG, \phi}^{\tame}(\mF)$ given \eqref{eq: tame central S operator} is the same as endomorphism by multiplying by $\chi_{V, (\ga, \sigma)}$.
\end{lemma}
\begin{proof}
We may write the multiplication by $\chi_V$ map as
\begin{equation}\label{eq: spectral S-operator}
\mF\xrightarrow{\id\otimes u_V} \mF\otimes\widetilde{V}\otimes \widetilde{V}^*\xrightarrow{\id\otimes(\ga,\sigma)\otimes \id} \mF\otimes\widetilde{V}\otimes \widetilde{V}^*\cong \mF\otimes\widetilde{V}^*\otimes\widetilde{V}\xrightarrow{\id\otimes e_V}\mF.
\end{equation}
It then follows from \Cref{lem:image-of-trace-of-some-objects-in-specHecke} \eqref{lem:image-of-trace-of-some-objects-in-specHecke-2} that this coincides the operator \eqref{eq: tame central S operator} as defined via \eqref{eq: abstract S-operator}.
\end{proof}

\begin{remark}\label{rem: tame central S operator}
Note that the abstract construction  \eqref{eq: abstract S-operator} of $S$-operators are only for the objects in the essential image of $\Ch_{{}^cG,\phi}^{\tame}$. Thanks to \Cref{lem: spectral S operator vs abstract S operator}, they are now defined on every object in $\indcoh(\locsys_{{}^cG,F}^{\tame}$. Namely, $S_{(\mZ^{\mon}(V), \eta_\ga)}$ is just given by multiplication by $\chi_{V, (\ga, \sigma)}\in Z_{{}^cG,F}^{\tame}$.  
\end{remark}

In the unipotent case, we can just to consider $I=\{1\}$. We have $\mZ^{\unip}(V)$ equipped with 
\begin{equation}\label{eq: taut action of spectral central sheaf unip} 
\varphi^{\mathrm{univ}}: I_F^t \to \End(\mZ^{\unip}(V))
\end{equation}
\begin{equation}\label{eq: phi equivariance of spectral central sheaf unip}
\Phi:\mZ^{\unip}(V)\cong \phi_*\mZ^{\unip}(V).
\end{equation} 
In addition, for $\mF\in \indcoh(S_{{}^cG,\breve F}^{\unip})$, we have $\mF\star \mZ^{\unip}(V)\cong t^*V\otimes \mF\cong \mZ^{\unip}(V)\star \mF$. Then
we have
 \[
 \mF\star\mZ^{\unip}(V)\cong \mZ^{\unip}(V)\star\mF\cong \mZ^{\unip}(V)\star\mF\stackrel{\Phi}{\cong} \phi_*(\mZ^{\unip}(V))\star\mF.
 \]
The corresponding $S$-operator will be denoted by 
\begin{equation}\label{eq: unip central S operator}
S_{V}: \Ch_{{}^cG, \phi}^{\unip}(\mF)\to \Ch_{{}^cG, \phi}^{\unip}(\mF),
\end{equation}
which as in \Cref{lem: spectral S operator vs abstract S operator} is isomorphic to multiplication by $\chi_V:=\chi_{V,(\sigma)}$.

\newpage

\section{The local Langlands category}\label{sec:kot-stack}

A general wisdom shared among various people is that in the local Langlands correspondence it is better not to just study representation theory of a single $p$-adic group $G$, but simultaneously to study representation theory of a collection of groups closely related to $G$. There are various ways to make this idea precise by appropriately choosing such collection, such as Vogan's pure inner forms, Kottwitz-Kaletha's extended pure inner forms, etc. It is the extended pure inner forms of $G$ that is most suitable for the geometric/categorical approach, as they arises naturally in the study of Shimura varieties and moduli of Shtukas.
It turns out one can go one step further to consider the representation theory of not just extended pure inner forms of $G$, but all extended pure inner forms of Levi subgroups of $G$ together. The representation categories of these groups glue nicely together to a category which is conjecturally equivalent to the category of (ind-)coherent sheaves on the stack of arithmetic local Langlands parameters, as we will explain at the end of this section.

In this section, we introduce one of the central players of this article, the stack of $G$-isocrystals and the category of $\ell$-adic sheaves on it. This framework realizes the idea of gluing the aforementioned representation categories together. The stack of $G$-isocrystals, along with several related objects considered here, may seem unconventional from a traditional algebraic geometry perspective, as it is the quotient of an ind-scheme by an ind-algebraic group. Nonetheless, we will demonstrate that the category of $\ell$-adic sheaves on this stack can still be understood within the framework developed in \Cref{sec:pspl-stacks}. We will show that the category possesses numerous favorable properties akin to those of the usual category of $\ell$-adic sheaves on (stratified) algebraic varieties. For instance, it is compactly generated, admits a canonical self-duality, and admits a natural $t$-structure.

It is worth noting that a very different approach for gluing these categories has been developed by Fargues-Scholze \cite{Fargues.Scholze.geometrization}. It is reasonable to conjecture that the two approaches lead to equivalent categories,  albeit through non-trivial means. For further discussion on this topic, see \Cref{rem: comparison between kotG and BunG}.

\subsection{The stack of local $G$-Shtukas}\label{sec:stack-Sht}
In this subsection, we review and further study some basic facts about the stack of local $G$-Shtukas.

\subsubsection{Iwahori-Weyl group and parahoric group schemes}\label{SS: Iwahori-Weyl-group}
We review a few facts about Iwahori-Weyl group and parahoric group schemes. We take the opportunity to also fix a few notations that will be used throughout this article.

Let $F$ be a non-archimedean local field with ring of integers $\mathcal{O}_{F}$ and finite residue field $k_{F}$ of $q=p^{[k_F:\bF_p]}$ elements. We fix a uniformizer $\varpi\in\mO_F$. We fix a separable closure $\overline F$ and let $\Ga_F$ be the Galois group of $F$ and $I_F\subset \Ga_F$ the inertia subgroup.
Let $\breve{F}$ be the completion of the maximal unramified extension of $F$ (in $\overline F$) and its ring of integers by $\mO_{\breve F}$ and its quotient field by $k$ (so that $k = \overline{k_{F}}$). Then $\aut(\breve F/F)$ contains a canonical element lifting the $q$-Frobenius element $\sigma$ in $\Ga_{k_F}$. By abuse of notations, we also use $\sigma$ to denote this element in $\aut(\breve F/F)$. 
Sometimes for simplicity we also write $\mO_F\subset \mO_{\breve F}$ simply as $\mO\subset \breve\mO$ if no confusion is to likely arise.

Let $G$ be a connected reductive group over $F$. 
Let $A$ be a maximally split torus of $G$ over $F$. Let $S\subset G$ be an $F$-rational torus containing $A$ such that $S_{\breve F}$ is a maximally split torus of $G_{\breve F}$. The pair $A\subset S$ is unique up to conjugation by an element in $G(F)$.
Let $T=Z_G(S)$, which is a maximal torus of $G$. Let $\mA\subset\mS\subset\mT$ be the Iwahori group schemes (over $\mO_F$) of $A\subset S\subset T$. Let 
\[
W_0=N_{G}(T)(\breve F)/T(\breve F),\quad \mbox{resp.} \ \widetilde W=N_{G}(T)(\breve F)/\mT(\mO_{\breve F})
\]
denote the relative finite Weyl group of $G_{\breve F}$, resp. the Iwahori-Weyl group of $G_{\breve F}$. They fit into the following short exact sequence
\begin{equation}\label{eq-Iwahori-Weyl-group-sequence}
1\to \xcoch(T)_{I_F}\to \widetilde W\to W_0\to 1.
\end{equation}
Elements of $\xcoch(T)_{I_F}\subset \widetilde{W}$ are usually called translation elements.
To avoid the confusion of notations, for $\la\in \xcoch(T)_{I_F}$ we will let $t_\la$ denote the corresponding translation element in $\widetilde{W}$.

We let $\scrB(G, \breve F)$ denote the (reduced) Bruhat-Tits building of $G$ over $\breve F$ and let $\scrA(G_{\breve F},S_{\breve F})\subset \scrB(G,\breve F)$ denote the apartment corresponding to $S_{\breve F}$. Let $\Phi\subset \xch(S_{\breve F})$ be the relative root system of $(G_{\breve F}, S_{\breve F})$, and let $\Phi_{\af}$ be the set of corresponding affine roots, regarded as affine functions on $\scrA(G_{\breve F},S_{\breve F})$. Let $\Phi_\af\to\Phi$ be the map sending an affine root $\al$ to its vector part $\dot{\al}$. For $\al\in \Phi_\af$, let $s_\al\in \widetilde{W}$ be the affine reflection corresponding to $\al$.
Let $W_\af\subset \widetilde{W}$ be the subgroup generated by affine reflections corresponds to affine roots. It is a normal subgroup, usually called the affine Weyl group of $G_{\breve F}$ (which can also be regarded as the Iwahori-Weyl group of the simply-connected cover of $G_{\breve F}$). It is known that 
\begin{equation}\label{eq-Omega-dual-group}
 \widetilde W/W_\af\cong \pi_1(G)_{I_F}.
\end{equation}
On the other hand, the group $\widetilde W$ is a quasi Coxeter group with a length function, once we fix  an alcove $\breve\bfa\subset \scrA(G_{\breve F}, S_{\breve F})$, or equivalently, an Iwahori group scheme $\breve\mI$ of $G$ (over $\mO_{\breve F}$) containing $\mT_{\mO_{\breve F}}$.
Let $\Omega_{\breve\bfa}\subset \widetilde W$ be the corresponding subgroup of length zero elements. Then 
\begin{equation}\label{eq-Iwahori-Weyl-semiproduct}
   \widetilde W=W_\af\rtimes\Omega_{\breve\bfa}. 
\end{equation}

Note that the $q$-Frobenius $\sigma$ acts on everything. In particular, $\scrA(G_{\breve F}, S_{\breve F})^\sigma=\scrA(G,A)$ is the apartment associated to $A$ in the building $\scrB(G,F)=\scrB(G,\breve F)^\sigma$. 
For an alcove $\breve\bfa\subset \scrA(G_{\breve F},S_{\breve F})$ such that  $\mathbf{a}=\breve\bfa\cap \scrA(G,A)$ is an alcove of $\scrA(G,A)$, 
the corresponding decomposition \eqref{eq-Iwahori-Weyl-semiproduct} is preserved under the action of $\sigma$. In addition, there is a canonical isomorphism
\[
   (\Omega_{\breve\bfa})_\sigma\cong \pi_1(G)_{\Ga_F}.
\]

We will occasionally also consider the extended building 
\begin{equation}\label{eq: ext BT building}
\scrB^{\ext}(G,\breve F)= \scrB(G, \breve F)\times \xcoch(Z_G)^{I_F}_\bR
\end{equation}
on which $G(\breve F)$-acts. 
If $D\subset \scrB(G, \breve F)$ is a subset, let $D^{\ext}=D\times  \xcoch(Z^1_G)^{I_F}_\bR \subset \scrB^{\ext}(G, \breve F)$. 
If $D$ is bounded in $\scrB(G,\breve F)$, we let $
\breve\mG_{D}$ denote the ``stabilizer" group scheme of $D^{\ext}$ as constructed by Bruhat-Tits. I.e. $\breve\mG_{D}$ is  
 the smooth affine group scheme over $\breve\mO$, with generic fiber $G_{\breve F}$, such that $\breve\mG_D(\breve \mO)$ consist of elements in $G(\breve F)$ that fix every point of $D^{\ext}$. If $D\subset \breve\bff$ is contained in a facet, then the neutral connected component 
 \[
 \breve\mP_{\breve \bff}:=\breve\mG_{D}^\circ
 \] 
 of $\breve\mG_{D}$ is the parahoric group scheme associated to $\breve\bff$. Sometimes, we simply denote it by $\breve\mP$.

If $\breve\bff\subset \scrA(G_{\breve F},S_{\breve F})$, or equivalently $\mT_{\mO_{\breve F}}\subset \breve\mP$, then $\breve\mP$ is called semi-standard. 
We let $L_{\breve\mP}$ or sometimes $L_{\breve\bff}$ denote the corresponding Levi ``quotient" (more precisely it is the Levi quotient of $\breve\mP_k$), which is a connected reductive group over $k$, containing $\mS_k$ as its maximal torus. 
Let $\Phi_{\breve\bff}\subset \Phi_\af$ be the subset consisting of affine roots that vanish on $\breve\bff$, and let $W_{\breve\bff}\subset \widetilde{W}$ be the subgroup generated by affine reflections corresponding to affine roots in $W_{\breve\bff}$. Then the map $\Phi_\af\to \Phi\subset \xch(S)=\xch(\mS_k)$ sends $\Phi_{\breve\bff}$ to  the root system of $(L_{\mP},\mS_k)$. Sometimes $(\Phi_{\breve\bff}, W_{\breve\bff})$ is also denoted as $(\Phi_{L_{\breve\mP}}, W_{L_{\breve\mP}})$ or $(\Phi_{\breve\mP}, W_{\breve\mP})$.
Once we fix an alcove $\breve\bfa\subset \scrA(G_{\breve F}, S_{\breve F})$, or equivalently an Iwahori group scheme $\breve\mI$ of $G$ (over $\mO_{\breve F}$) containing $\mT_{\mO_{\breve F}}$, we call a parahoric group scheme $\breve\mP$ standard if the corresponding facet $\breve\bff\subset\overline{\breve\bfa}$, or equivalently  $\breve\mP(\mO_{\breve F})\supset \breve\mI(\mO_{\breve F})$.

Similarly, a parabolic subgroup $\breve\bfP\subset G_{\breve F}$ is called semi-standard if $T_{\breve F}\subset \breve\bfP$. We write $\breve\bfP=\breve{M}_{\breve\bfP}U_{\breve\bfP}$, where $\breve{M}=\breve{M}_{\breve\bfP}$ is the unique Levi subgroup of $\breve\bfP$ containing $T_{\breve F}$. The root system of $(\breve{M},S_{\breve F})$ is denoted by $\Phi_{\breve{M}}\subset \Phi$, and the relative Weyl group is denoted as $W_{\breve{M}}$ (or sometimes by $W_{\breve\bfP}$), which is a subgroup of $W_0$. Let $\Phi_{U_{\breve\bfP}}\subset \Phi$ be the set of roots whose root groups are contained in $U_{\breve\bfP}$. Associated to $\breve{M}$, there is also the corresponding affine roots $\Phi_{\breve{M},\af}=\Phi_{\breve{M}}\times_{\Phi}\Phi_\af$, and the corresponding Iwahori-Weyl group $\widetilde{W}_{\breve{M}}=N_{\breve{M}}(T)(\breve F)/\mT(\mO_{\breve F})=W_{\breve{M}}\times_{W_0} \widetilde{W}$. 
Once we fix a Borel subgroup $\breve{B}\supset T_{\breve F}$, a parabolic subgroup $\breve\bfP\subset G_{\breve F}$ is called standard if $\breve\bfP\supset \breve{B}$.

Note that the inclusion $\breve\bfP\subset G_{\breve F}$ of a semi-standard parabolic induces a surjective map $\scrA(G_{\breve F},S_{\breve F})\to \scrA(\breve M,S_{\breve F})$.
 Given an alcove $\breve\bfa\subset \scrA(G_{\breve F}, S_{\breve F})$, its image in $\scrA(\breve M,S_{\breve F})$ is an alcove
 \begin{equation}\label{eq: alcove of levi}
 \breve\bfa_{\breve M}=\{v\in \scrA(\breve M,S_{\breve F})\mid \al(v)>0, \ \forall \al\in \Phi_{\breve M,\af}\cap \Phi_\af^+\},
 \end{equation}
where $ \Phi_\af^+\subset\Phi_\af$ is the set of positive affine roots determined by $\breve\bfa$. In particular, $\Phi_{\breve M,\af}^+=\Phi_{\breve M,\af}\cap\Phi_{\af}^+$.
Let $\Delta_\af\subset \Phi_\af^+$ and $\Delta_{\breve M,\af}\subset\Phi_{\breve M,\af}^+$ be the corresponding sets of simple affine roots.
Note that  $\Delta_{\af}\cap \Phi_{\breve M,\af}\subset\Delta_{\breve M,\af}$. In particular, the image of a facet $\breve\bff\subset\overline{\breve\bfa}$ under the map $\scrA(G_{\breve F}, S_{\breve F})\to \scrA(\breve M,S_{\breve F})$ is a facet 
\begin{equation}\label{eq: facet of levi}
\breve\bff_{\breve M}=(\bigcap_{\al\in \Phi_{\breve M,\af}\cap \Phi_{\breve\bff}}\{\al=0\})\bigcap \overline{\breve\bfa}_{\breve M}.
\end{equation}
We let $\ell_{\breve M}$ denote the length function on $\widetilde{W}_{\breve M}= \xcoch(T)_{I_F}\rtimes W_{\breve M}$ determined by $\breve\bfa_{\breve M}$. Note that $\ell_{\breve M}\neq \ell|_{\widetilde{W}_{\breve M}}$.

Now if $\breve\mQ$ is the parahoric group scheme of $G_{\breve F}$ corresponding to $\breve\bff\subset \breve\bfa$, let $\breve\mQ_{\breve M}$ denote the parahoric group scheme of $\breve M$ corresponding to $\breve\bff_{\breve M}$. Then 
\begin{equation}\label{eq: parahoric of levi}
\breve\mQ_{\breve M}(\breve \mO)=\breve\mQ(\breve\mO)\cap \breve M(\breve F).
\end{equation}
We also let $\breve\mQ_{U_{\breve\bfP}}$ be the (fiberwise connected) smooth affine group scheme over $\breve\mO$ such that $\breve\mQ_{U_{\breve\bfP}}(\breve\mO)=\breve\mQ(\breve\mO)\cap U_{\breve\bfP}(\breve F)$. Then $\breve\mQ_{\breve\bfP}=\breve\mQ_{\breve M}\breve\mQ_{U_{\breve\bfP}}$ is a smooth integral model of $\breve\bfP$.  Let $U_{\breve\bfP}^-$ be the unipotent radical of the opposite parabolic $\breve\bfP^-$. Then we similarly have $\breve\mQ_{\breve\bfP^-}=\breve\mQ_{\breve M_{\breve\bfP}}\breve \mQ_{U^-_{\breve\bfP}}$. On the other hand, the natural multiplication map $\breve\mQ_{U_{\breve\bfP}^-}\times \breve\mQ_{\breve M_{\breve\bfP}}\times \breve\mQ_{U_{\breve\bfP}}\to \breve\mQ$ is an open embedding. We need the following variant of \eqref{eq: parahoric of levi}, which will only be used in \Cref{eq: non-connected-local-shtuka-uw}.

\begin{lemma}\label{lem: parahoric intersect with Levi}
Let $D_G\subset \scrA(G_{\breve F},S_{\breve F})$ be a bounded subset and let $D_M\subset \scrA(\breve M,S_{\breve F})$ be its image. Then $\breve\mM_{D_M}(\breve \mO)=\breve\mG_{D_G}(\breve \mO)\cap M(\breve F)$.
\end{lemma}
\begin{proof}
We may choose an $\breve M(\breve F)$-equivariant embedding $\scrB^{\ext}(\breve M,\breve F)\to \scrB^{\ext}(G,\breve F)$ identifying $\scrA^{\ext}(\breve M,S_{\breve F})=\scrA^{\ext}(G_{\breve F},S_{\breve F})$. (Such an embedding is unique up to translation by $\xcoch(Z_{\breve M})^{I_F}_\bR$, and in particular the image is independent of the choice of the embedding.)
Then $D_M^{\ext}= \xcoch(Z_{\breve M})^{I_F}_\bR D_G^{\ext}$. The inclusion $\breve\mM_{D_M}(\breve \mO)\subset \breve\mG_{D_G}(\breve \mO)\cap M(\breve F)$ is obvious. On the other hand, if $g\in \breve\mG_{D_G}(\breve \mO)\cap M(\breve F)$, then $g=\la(\varpi)g\la(\varpi)^{-1}\in \breve\mG_{\la(\varpi)D_G^{\ext}}(\breve\mO)$ for every $\la\in \xcoch(Z_{\breve M})^{I_F}$. Therefore, $g\in \breve\mM_{D_M}(\breve \mO)$.
\end{proof}

Now we assume that $G$ is quasi-split. Then once we fix a pinning $(B,T,e)$ of $G$ over $F$, there is a natural choice of subtori of $G$ (over $F$) and Iwahori group scheme of $G$ (over $\mO_F$). Namely, we can take $A\subset S\subset T$, where $A$ is the maximal split subtorus of $T$, and $S$ is the maximal $\breve F$-split subtorus of $T$. Recall that the pinning determines an absolutely special vertex $v_0\in \scrA(G,A)$ (e.g. see \cite[\S 4.2]{zhu2020coherent}). 
We then obtain an identification 
\[
\xcoch(S/Z_G)^I_\bR\cong \scrA(G_{\breve F},S_{\breve F}).
\] 
Then there is the alcove $\breve\bfa\subset\xcoch(S/Z_G)^I_\bR$ that contains original and is contained in the dominant Weyl chamber determined by $B$. This alcove is $\sigma$-stable, and the corresponding Iwahori $\mI$ is defined over $\mO_F$, containing $\mT$. We equip $\widetilde{W}$ with the length function $\ell$ determined by $\breve\bfa$.  
We may also identify $W_0$ with the Weyl group of $v_0$.
Then
\begin{equation}\label{eq: Bernstein presentation of Iwahori-Weyl group}
\widetilde{W}=\xcoch(T)_{I_F}\rtimes W_0.
\end{equation}

\subsubsection{$\sigma$-conjugacy classes of the Iwahori-Weyl group}\label{sec:sigma-straight-element}
We assume that $G$ is quasi-split with a pinning $(B,T,e)$. It determines an alcove $\breve\bfa$ and an absolutely special vertex $v_0$ as above. Then we have the length function $\ell$ on $\widetilde{W}$ and a semi-direct product decomposition \eqref{eq: Bernstein presentation of Iwahori-Weyl group}.

We review some results from \cite{he2014geometric,he2016hecke,He.Nie.minimal.length}. First, recall from \cite[\textsection{1.7}]{he2014geometric} that there is a map
\begin{equation}
\widetilde W\to \xcoch(T)^{I_F}_\bQ\times \pi_1(G)_{I_F},\quad w\mapsto (\tilde\nu_w, \pi_0(w)).
\end{equation}
Namely, for $w\in \widetilde{W}$, $\pi_0(w)$ is just the image of $w$ in $\widetilde{W}/W_\af\cong \xch(Z_{\hat{G}}^{I_F})$. On the other hand,
choose $n$ such that $\sigma^n=1$ and $w\sigma(w)\cdots \sigma^{n-1}(w)\in \xcoch(T)_{I_F}$ (such $n$ always exists) and regard it as an element in $\xcoch(T)^{I_F}_\bQ$ under the natural map $\xcoch(T)_{I_F}\to (\xcoch(T)_{I_F})_\bQ\cong \xcoch(T)
^{I_F}_\bQ$. We can uniquely write this element as $n\tilde\nu_w$ with $\tilde\nu_w\in \xcoch(T)^{I_F}_{\bQ}$. Then $\tilde\nu_w$ is independent of the choice of $n$.

Let 
\[
B(\widetilde W)=\widetilde W/\sim,  \quad  w\sim vw' \sigma(v)^{-1}
\] 
denote the set of $\sigma$-conjugacy classes of $\widetilde{W}$. The above map induces a map
\begin{equation*}\label{eq-newton-kottwitz-for-BW}
B(\widetilde{W})\to \xcoch(T)^{+,\Ga_F}_\bQ\times \pi_1(G)_{\Ga_F},\quad w\mapsto (\nu_w, \kappa(w)).
\end{equation*}
Here, $\nu_w\in  \xcoch(T)_\bQ^{+}$ be corresponding dominant element (with respect to $B$) in the Weyl group orbit of $\tilde\nu_w$, called the Newton point of $w$. It is in fact $\Ga_F$-invariant and depends only on the $\sigma$-conjugacy class of $w$. On the other hand, $\kappa(w)$ is the image of $\pi_0(w)$ under the map $\pi_1(G)_{I_F}\to \pi_1(G)_{\Ga_F}$, called the Kottwitz point of $w$. It also only depends on the $\sigma$-conjugacy class of $w$.

Every $\sigma$-conjugacy class in $\widetilde{W}$ determines an $F$-rational Levi subgroup 
\[
M=Z_G(\nu_w).
\]
Then $\bfP=MB$ is a standard parabolic subgroup (with respect to $T\subset B$) defined over $F$.
Recall that we let $\breve\bfa_M$ denote the unique alcove in $\scrA(M_{\breve F}, S_{\breve F})$ such that $\breve\bfa_M^{\ext}$ contains $\breve\bfa^{\ext}$, and let $\ell_M$ denote the length function on $\widetilde{W}_M= \xcoch(T)_{I_F}\rtimes W_M$ determined by $\breve\bfa_M$.

Now let $w\in \widetilde{W}$. Let
\[
\breve M_w=Z_{G_{\breve F}}(\tilde\nu_w).
\]
This is a Levi of $G_{\breve F}$ defined over $\breve F$. It is related to the rational Levi $M$ attached to the $\sigma$-conjugacy class of $\widetilde{W}$ containing $w$ as follows: there is a unique element $y\in W_0$, of minimal length in $y W_M$, such that $y\nu_w=\tilde\nu_w$.
Then $\breve M_w= \dot{y} M_{\breve F}\dot{y}^{-1}$, where $\dot{y}$ is a lifting of $y$ to $N_G(T)(\breve F)$. Later on, we will consider
\begin{equation}\label{eq: translate w to x}
w^+:=y^{-1}w\sigma(y),
\end{equation}
which belongs to $\widetilde{W}_M$, and $\tilde\nu_{w^+}=\nu_w$. 

We note that for a general element $w\in \widetilde{W}$, we have 
\begin{equation}\label{eq: length vs Newton point}
\ell(w)\geq \langle 2\rho, \nu_w\rangle.
\end{equation}
If $\ell(w) = \langle 2\rho,\nu_{w} \rangle$, or equivalently
$\ell(w\sigma(w)\cdots \sigma^{n-1}(w)) = n\ell(w)$
for all $n>0$, then $w$ is called $\sigma$-straight.  A $\sigma$-conjugacy class of $\widetilde{W}$ is called $\sigma$-\textit{straight} if it contains a $\sigma$-straight element. Let $B(\widetilde W)_{\mathrm{str}}\subset B(\widetilde W)$ denote the set of straight $\sigma$-conjugacy classes.

We recall some remarkable combinatorics of $\sigma$-conjugacy classes in $\widetilde{W}$ due to He-Nie (\cite{He.Nie.minimal.length}). For the purpose, we need some notations and terminology. 
For $w, w', t\in \widetilde{W}$, we write $w\xrightarrow{t}_\sigma w'$ if $w'=tw\sigma(t)^{-1}$, $\ell(w')\leq \ell(w)$ and $\ell(t)\leq 1$. We write $w\rightarrow_\sigma w'$ if  there is a sequence of elements $w=w_0\xrightarrow{t_1}_\sigma w_1\xrightarrow{t_2}_\sigma \cdots \xrightarrow{t_n}_\sigma w_n=w'$. We write $w\leftrightarrow_\sigma w'$ if $w\rightarrow_\sigma w'$ and $w'\rightarrow_\sigma w$. In this case, we say $w$ and $w'$ are $\sigma$-conjugate by cyclic shift. By \cite[Lemma 1.6.4]{Deligne.Lusztig} (which works for affine Weyl group as well), it is easy to see that $w,w'$ are $\sigma$-conjugate by cyclic shift if and only if there is a sequence of elements $\{w'_0,w'_1,\ldots, w'_r\}\subset \widetilde{W}$ and for each $i$ there are $x_i,y_i\in\widetilde{W}$ such that $w'_0=w$ and $w'_r=w'$, $\ell(w'_i)=\ell(x_i)+\ell(y_i)=\ell(w'_i)$ and $w'_{i-1}=x_iy_i$, $w'_i=y_i\sigma(x_i)$.

\begin{theorem}\label{thm: reduction to min length elements}
Let $C\subset \widetilde{W}$ be a $\sigma$-conjugacy class, and let $C_{\mathrm{min}}\subset C$ be the subset of minimal length elements (with respect to the length function $\ell$ on $\widetilde{W}$).
\begin{enumerate}
\item\label{thm: reduction to min length elements-3} Suppose $C$ is a straight $\sigma$-conjugacy class. Then $C_{\mathrm{min}}$ is the set of $\sigma$-straight elements in $C$. 
Every two elements $w,w'\in C_{\mathrm{min}}$ are $\sigma$-conjugate by cyclic shift.
In addition, for $w\in C_{\mathrm{min}}$, the element $w^+$ from \eqref{eq: translate w to x} satisfies $\ell_M(w^+)=0$.

\item\label{thm: reduction to min length elements-1} For every $v\in C$, there is a sequence of elements $v=v_0\xrightarrow{s_1}_\sigma v_1\xrightarrow{s_2}_\sigma \cdots \xrightarrow{s_n}_\sigma v_n=v'$, with $s_i$ simple reflections, and a facet $\breve\bff\subset\breve\bfa$ such that $v_n\in C_{\mathrm{min}}$ is of the form $v_n=uw$ where $w$ is a $\sigma$-straight element and is of minimal length in $W_{\breve\bff} w$,  $w\sigma(W_{\breve\bff})w^{-1}=W_{\breve\bff}$, and $u\in W_{\breve\bff}$.
\end{enumerate}
\end{theorem}
\begin{proof}
Part \eqref{thm: reduction to min length elements-3} is \cite[Proposition 3.2]{He.Nie.minimal.length}. Part \eqref{thm: reduction to min length elements-1} is \cite[Theorem 2.9, Theorem 3.4]{He.Nie.minimal.length}.
\end{proof}

\begin{remark}\label{rem: refinement of B(W)}
Note that  \Cref{thm: reduction to min length elements} \eqref{thm: reduction to min length elements-1} in particular applies to $v$ that already is of minimal length in its $\sigma$-conjugacy class, in which case $v$ and $v_n$ are $\sigma$-conjugate by cyclic shift. But unlike $\sigma$-straight conjugacy classes, for general $C$, not every pair of elements in $C_{\mathrm{min}}$ are $\sigma$-conjugate by cyclic shift.
\end{remark}

\begin{remark}\label{lem: wsigma stable facet}
Let $\breve\bff\subset \overline{\breve\bfa}\subset \scrB(G, \breve F)$ be a facet, and let $\Phi_{\breve\bff}\subset \Phi_{\af}$ be the corresponding sub root system.
Let $\Phi_{\breve\bff}^+=\Phi_{\breve\bff}\cap\Phi_\af^+$. It is easy to check that the following are equivalent.
\begin{itemize}
\item $w$ is of minimal length in $W_{\breve\bff} w$ and $w\sigma(W_{\breve\bff})w^{-1}=W_{\breve\bff}$.
\item $w(\sigma(\Phi_{\breve\bff}^+))=\Phi_{\breve\bff}^+$.
\end{itemize}
In this case, if $\al\in\Phi_{\breve\bff}^+$ is a simple root, then $\al\in \widetilde{W}_{\breve M_w}$, and therefore is also simple in $\Phi_{\breve M_w, \af}$ (with respect to the alcove $\breve\bfa_{\breve M_w}$ as defined in \eqref{eq: alcove of levi}). In particular, $\Phi_{\breve\bff}\subset\Phi_{\breve M_w,\af}$. It follows that the length function $\ell$ of $\widetilde{W}$ and the length function $\ell_{\breve M_w}$ of $\widetilde{W}_{\breve M_w}$ coincide when restricted to $W_{\breve\bff}$.
\end{remark}

\begin{remark}\label{rem: u min in finite weyl}
Let $v_n=uw$ be as in \Cref{thm: reduction to min length elements} \eqref{thm: reduction to min length elements-1}.
Note that $u$ is minimal length in its $\Ad_w\sigma$-conjugacy class in $W_{\breve\bff}$. Otherwise, there would be some $t\in W_{\breve\bff}$ such that $\ell(tu w\sigma(t)^{-1}w^{-1})<\ell(u)$ so $\ell(tuw\sigma(t)^{-1})<\ell(uw)$, contradiction.
On the other hand, if $u'$ and $u$ are in the same $\Ad_w\sigma$-conjugacy class in $W_{\breve\bff}$ and if $\ell(u')=\ell(u)$, then $u'w$ and $uw$ are in the same $\sigma$-conjugacy class of $\widetilde{W}$ and $\ell(u'w)=\ell(uw)$. Therefore, $u'w$ is also of minimal length in $C$.
\end{remark}

\begin{remark}\label{rem: special form of minimal length element in sigma conjugacy class}
Let $\xcoch(T)_{I_F}^+=\{\la \mid (\la, a)\geq 0, \forall a\in \Phi^+\}$. Note that every $\sigma$-conjugacy class $C$ of $\widetilde{W}$, there is $w\in C_{\min}$ of the form $w=vt_\la$, with $v\in W_0$ and $\la\in \xcoch(T)_{I_F}^+$. To prove this, first 
recall that if $\la\in \xcoch(T)_{I_F}^+$, then $t_\la$ is of minimal length in the coset $t_\la W_0$ (e.g. see \cite[Lemma 9.2]{zhu2014coherence}). In fact, using \cite[Lemma 9.1]{zhu2014coherence}, one sees immediately that if $w=t_\la v$ with $t_\la \in\xcoch(T)_{I_F}, v\in W_0$ is of minimal length in $wW_0$, then $v^{-1}(\la)\in \xcoch(T)_{I_F}^+$. 
Now let $w'\in C_{\min}$. We write $w= w_1w_2$ with $w_2\in W_0$ and $w_1$ of minimal length in $w_1W_0$. Then $w_2\sigma(w_1)$ is also in $C_{\min}$. If we further write $w_1=t_\la v$, then $w_2\sigma(w_1)=w_2\sigma(v) t_{\sigma(v^{-1}(\la))}$ as desired.
\end{remark}

We will need to review a classification of $\sigma$-conjugacy classes of $\widetilde{W}$ given in \cite[Theorem 1.19]{he2016hecke}.
We may reinterpret standard quadruples of \emph{loc. cit.} as $(M, x, \breve\bff_M, \frakc)$, where 
\begin{itemize}
\item $M$ contains $T$ and is the Levi subgroup of a standard $F$-rational parahobic subgroup $P\subset G$;
\item $x\in \widetilde{W}_M\subset \widetilde{W}$ satisfying $\ell_M(x)=0$, $\tilde{\nu}_x$ is dominant (w.r.t. $B$), and $Z_G(\tilde{\nu}_x)=M$;
\item $\breve\bff_M\subset \overline{\breve\bfa_M}$ is a facet stable under the action of $x\sigma$, where $\breve\bfa_M$ is as in \eqref{eq: alcove of levi};
\item $\frakc\subset W_{\breve\bff_M}$ is an elliptic $\Ad_x\sigma$-conjugacy class.
\end{itemize} 
Note that our $M$ corresponds to the set $J$, and $\breve\bff_M$ corresponds to the set $K$ as in \emph{loc. cit.} Two standard quadruples $(M, x, \breve\bff_M, \frakc)$ and $(M', x', \breve\bff'_{M'}, \frakc')$ are called equivalent if $M=M'$, and there is an element $w\in \widetilde{W}_M$
such that $x'=w x \sigma(w)^{-1}$, $\breve\bff'_{M}= w(\breve\bff_M)$ and $\frakc'=w \frakc w^{-1}$.
Let $\Quad_\sigma$ be the set of equivalence classes of standard quadruples.

Then by \cite[Theorem 1.19]{he2016hecke}, we have a bijection
\begin{equation}\label{eq: standard 4uple to conjugacy class} 
\Quad_\sigma\to B(\widetilde{W}),\quad (M, x, \breve\bff_M, \frakc)\mapsto C:=\{ w\frakc x\sigma(w)^{-1}\mid w\in \widetilde{W}\}.
\end{equation}

We shall also need to recall the inverse map.
Let $C$ is a $\sigma$-conjugacy class in $\widetilde{W}$. Let $uw\in C_{\min}$ as in \Cref{thm: reduction to min length elements} \eqref{thm: reduction to min length elements-1}, and write $w=yw^+\sigma(y)^{-1}$ as in \eqref{eq: translate w to x}. 
Then let $M=Z_G(\nu_w)$ and $x=w^+$.
We may choose $\breve\bff$ in \Cref{thm: reduction to min length elements} \eqref{thm: reduction to min length elements-1} to be minimal. 
As $y$ is of minimal length in $yW_M$, $y(\Phi_{M,\af}^+)\subset \Phi_{\af}^+$.
Therefore if $\al\in\Phi_\af$ is a simple affine root vanishing on $\breve\bff$, then $y^{-1}(\al)$ is a simple affine root in $\Phi_{M,\af}$. We thus let $\breve\bff_M\subset \overline{\breve\bfa_M}$ be the zero locus of $y^{-1}(\al)$, for $\al$ simple affine root in $\Phi_\af$ that vanishes on $\breve\bff$. Finally, let $\frakc$ be the $\Ad_x\sigma$-conjugacy class containing $y^{-1}uy$. It is elliptic by the assumption of $\breve\bff$. 

\begin{remark}\label{rem: u min in finite weyl-2}
Note that if $W_{\breve\bff_M}\subset \widetilde{W}$ is equipped with the length function $\ell_M$ determined by $\breve\bfa_M$ (see \eqref{eq: alcove of levi}), then the isomorphism $W_{\breve\bff_M}\cong W_{\breve\bff},\ w\mapsto ywy^{-1}$ is compatible with the length functions, i.e. $\ell_M(w)=\ell(ywy^{-1})$.
It follows from \Cref{lem: wsigma stable facet} and \Cref{rem: u min in finite weyl} that $y^{-1}uy$ is of minimal length in its $\frakc$, and for every $u'\in \frakc_{\min}$, we have $yu'xy^{-1}=yu'y^{-1}w\in C_{\min}$.
\end{remark}

Finally, let us review a partial order on the set of straight $\sigma$-conjugacy classes. 
Let $C$ be a $\sigma$-straight conjugacy class of $\widetilde{W}$ and $w\in \widetilde{W}$.
Following \cite[\textsection{2}, \textsection{3}]{he2015kottwitzrapoport}, 
we write $C\preceq w$ if there is $v\in C_{\mathrm{min}}$ such that $v\leq w$ in the Bruhat order on $\widetilde{W}$. It is known that if $C\preceq w$ and $w'\rightarrow_\sigma w$, then $C\preceq w'$.
If $C_1, C_2$ are two $\sigma$-straight conjugacy class, we write $C_1\leq C_2$ if $C_1\preceq w$ for some (equivalently every) $w\in (C_2)_{\mathrm{min}}$.

\subsubsection{Loop groups}\label{sec:kot-actions-of-loop-groups}
We will introduce and study the algebro-geometric version of $B(G)$. First we introduce some notations.
For a perfect $k_F$-algebra $R$, let
\[
W_{\mO,n}(R)=W(R)\otimes_{W(k_F)}\mO_F/\varpi^n,\quad W_\mO(R):=\varprojlim_nW_{\mO,n}(R).
\]
We write $D_{F,R}=\Spec W_{\mO}(R)$ and $D_{F,R}^* = \Spec (W_{\mO}(R)[1/\varpi])$, or just $D_R$ and $D_R^*$ if $F$ is clear from the context. We write the automomrphism of $D_{F,R}$ and $D_{F,R}^*$ induced by the Frobenius automorphism of $R$ by $\sigma_R$.

Let $\mathcal{G}$ be an smooth affine model for $G$ over $\mathcal{O}_{F}$. Recall that the positive loop group $L^{+}\mathcal{G}$ and the loop group $LG$ associated to $G$ are defined as functors $\perf_{k_F} \rightarrow \mathrm{Set}$ by
\[
L^{+}\mathcal{G}(R) = \mG(W_{\mathcal{O}}(R)), \quad
LG(R) = G(W_{\mathcal{O}}(R)[1/\varpi]).
\]
We also define the $n$-th jet schemes $L^{n}\mathcal{G}$ by
\[
L^{n}\mathcal{G}(R) = \mathcal{G}(W_{\mathcal{O},n}(R)).
\]
For each $n\geq 0$ the functor $L^{n}\mathcal{G}$ is represented by a perfect affine scheme perfectly finite type over $k_{F}$. 
For $n' \geq n\geq 1$, the kernel of the surjective homomorphism $L^{n'}\mathcal{G}\rightarrow L^{n}\mathcal{G}$ is unipotent.
Then $L^{+}\mathcal{G} = \varprojlim_{n} L^{n}\mathcal{G}$ is represented by an affine scheme.  
We write $L^+\mG^{(n)}=\ker(L^+\mG\to L^n\mG)$, which is the $n$th congruence subgroup of $L^+\mG$. If $\mG$ is an Iwahori group scheme, sometimes we write $\iw=L^+\mG$, $\iw_n=L^n\mG$ and $\iw^{(n)}=L^+\mG^{(n)}$.

Given two $G$-torsors $\mE_1$ and $\mE_2$ over $D_R$, we write $\mE_1\dashrightarrow \mE_2$ for an isomorphism $\mE_1|_{D_R^*}\simeq \mE_2|_{D_R^*}$ and call it a modification of $G$-torsors.

Now we assume that $\mG=\mP$ is parahoric.
Let  $\mathrm{Gr}_{\mP} = LG/L^{+}\mP$ be the (partial) affine flag variety associated to $\mP$, which is an ind-projective scheme. 
It represents the moduli problem
\[
\mathrm{Gr}_{\mP}(R) = \bigg\{(\mathcal{E,\beta})
\Big| 
\begin{array}{ccc}
\mathcal{E}\ \text{is a}\ \mP\text{-torsor on } D_{R},\\
\beta: \mathcal{E} \dashrightarrow \mathcal{E}^{0}\text{ is a modification}
\end{array}
\bigg\}.
\]

Let $\mI$ be an Iwahori group scheme containing $\mT$, and write $\iw:=L^+\mI$. Recall that the $\iw_k$-orbits of $(\Gr_{\mI})_{k}$ are parameterized by $\widetilde W$. For a standard parahoric group scheme $\mP$,
the $(L^+\mP)_k$-orbits of $(\Gr_{\mP})_k$ are parametrized by $W_\mP\backslash \widetilde W/W_{\mP}$. For $w\in W_\mP\backslash \widetilde W/W_{\mP}$, let 
\begin{equation*}\label{eq:(locally-)closed-embedding-of-Schubert-variety}
    i_{\mP,\leq w}\colon \Gr_{\mP,\leq w}\subset (\Gr_{\mP})_k,\quad \mbox{resp. } i_{\mP,w}\colon \Gr_{\mP,w}\subset (\Gr_{\mP})_k
\end{equation*} 
be the corresponding Schubert variety (resp. Schubert cell). 

Recall that $\Gr_{\mP,\leq w}$ is a perfect projective irreducibe scheme over $k$ and $j_{\mP,w}\colon \Gr_{\mP,w}\to \Gr_{\mP,\leq w}$ is open. 
Let $LG_{\mP,\leq w}$ (resp. $LG_{\mP,w}$) be the pre-image of $\Gr_{\mP,\leq w}$ (resp. $\Gr_{\mP,w}$) under the projection $LG\to \Gr_\mP$. Then 
\[
LG_{\mP,\leq w}=\varprojlim_{n} \Gr_{\mP,\leq w}^{(n)}=\varprojlim_{n} LG_{\leq w}/L^+\mP^{(n)}_k, \quad LG_{\mP,w}=\varprojlim_{n} \Gr_{\mP,w}^{(n)}=\varprojlim_{n} LG_{\mP,w}/L^+\mP^{(n)}_k,
\]
with $\Gr_{\mP,\leq w}^{(n')}\to \Gr_{\mP,\leq w}^{(n)}$ (resp. $\Gr_{\mP,w}^{(n')}\to\Gr_{\mP,w}^{(n)}$) coh. unipotent if $n'\geq n\geq 1$.  Therefore, the morphisms
$LG_{\mP,w}\to \Spec k$ and $LG_{\mP,\leq w}\to \Spec k$ 
are ess. coh. pro-unipotent in the sense of \Cref{def-ess-coh-pro-unipotent}. In particular, $LG_{\mP,w}$ and $LG_{\mP,\leq w}$ are standard placid in the sense of \Cref{def: placid-algsp}. We still use $i_{\mP,\leq w}$ (resp. $i_{\mP,w}$) to denote the embedding $LG_{\mP,\leq w}\to LG_k$ (resp. $LG_{\mP,w}\to LG_k$).
Then $LG_k=\varinjlim_{w}  LG_{\mP,\leq w}$ is an ind-placid scheme over $k$. If we only take the colimit over those $w$ that are $\sigma$-invariant, then we see $LG$ is an ind-placid scheme over $k_F$.

\begin{remark}\label{rem: affine of LGw}
We note that $LG_{\mP,\leq w}$ and $LG_{\mI,w}$ are in fact affine. Indeed, $LG$ is ind-affine, and $LG_{\mP,\leq w}\subset LG_k$ is a closed embedding so $LG_{\mP,\leq w}$ is affine. On the other hand, $LG_{\mI,w}\subset LG_{\mI,\leq w}$ is an affine open embedding (as it is the base change of the affine open embedding $j_{\mI,w}: \Gr_{\mI,w}\subset \Gr_{\mI,\leq w}$). Therefore, $LG_{\mI,w}$ is also affine.
\end{remark}

In the sequel, when $\mP=\mI$, we usually omit $\mP$ from the subscripts in the above notations. We will sometimes also denote $\Gr_\mI$ by $\Fl$.

\begin{remark}
Of course, one can start with an integral model $\breve\mG$ of $G$ defined over $\breve\mO$ and all the discussions above (except those involving rationality) go through without change.
\end{remark}

\subsubsection{Moduli of local Shtukas}\label{sec.moduli.local.shtukas}
We introduce the stack
\[
\Sht_{\mP}^{\loc}:= \frac{LG}{\Ad_\sigma L^+\mP},
\]
where as before $\mP$ is a parahoric group scheme of $G$ over $\mO_F$, and the $\sigma$-conjugation action is given by
\[
\Ad_\sigma: L^+\mP\times LG\to LG,\quad (h,g)\mapsto hg\sigma(h)^{-1}.
\]
Recall from \cite[\textsection{5.3.2}]{xiao2017cycles} and \cite[\textsection{4}]{Zhu2016} that in this case the stack $LG/\Ad_{\sigma}L^{+}\mP$ can be identified with the moduli of local shtukas, and can be identified with the fiber product
\[
\begin{tikzcd}
\Sht_\mP^{\loc}\arrow[r,"\delta"]\arrow[d] & L^{+} \mP\backslash LG/ L^{+} \mP\arrow[d]\\
\bB L^{+} \mP \arrow[r,"\sigma \times \id"] &  \bB L^{+} \mP\times \bB L^{+} \mP
\end{tikzcd}
\]
Namely, for a perfect $k_F$-algebra $R$, 
\[
\sht_\mP^{\loc}(R)=\bigl\{(\mE,\varphi)\mid \mE \mbox{ is a } \mP\mbox{-torsors on } D_{R},\  \varphi: \sigma^{*}_{R}\mE\dashrightarrow \mE \bigr\}.
\]

\begin{remark}\label{rem: modification direction of Shtuka}
Note that our convention is different from \cite[\textsection{5}]{xiao2017cycles} and \cite[\textsection{4}]{Zhu2016}. In \emph{loc. cit.}, we defined $\sht_\mP^{\loc}$ such that its $R$-points classify $\mE\dashrightarrow \sigma^*_R\mE$.
We choose the convention here to be consistent with the $\sigma$-conjugation on $\widetilde W$ considered in \Cref{sec:sigma-straight-element}. However, the map 
\[
\delta: \Sht_\mP^{\loc}\to L^+\mP\bs LG/L^+\mP
\]
should send $ \sigma^{*}_{R}\mE\dashrightarrow \mE$ to the modifcation of $\mP$-bundles given by $\mE_1=\mE\dashrightarrow \mE_0= \sigma^{*}_{R}\mE$.  More explicitly, for $w\in  \widetilde{W}$, the map $\delta$ will send $\frac{LG_{\mP,w}}{\Ad_\sigma L^+\mP}$ to $L^+\mP\bs LG_{\mP,w^{-1}}/L^+\mP$. Later on, when we interpret $\delta$ as the horizontal map in \eqref{eq:trace-convolution-horocycle-diagram-special}, this convention is consistent with \Cref{rem: first projection vs second projection}.
\end{remark}

\begin{remark}\label{rem: allowable G-shtuka}
In fact, the definition of $\sht^\loc_{\breve\mG}$ makes sense even if $\breve\mG$ is just a smooth affine group scheme of $G$ over $\breve\mO$. Namely, we always have $\sigma: LG\to LG$, which may not necessarily send $L^+\breve\mG$ to itself. But the above quotient space still makes sense.
\end{remark}

\begin{remark}\label{rem: affine diagonal of Sht}
By \Cref{torsor.descent.gives.placid.stack}, $\Sht^{\loc}_\mP$ is ind-very placid. 
As mentioned in \Cref{rem: representability of diagonal}, the diagonal of $\Sht_\mP^{\loc}$ is affine.
\end{remark}

For every $n\geq 0$, we define the iterated $n$-th Hecke stack as the \'{e}tale quotient stack
\begin{equation}\label{eq: nth local Hecke stack for Sht}
\hk_{n}(\sht_\mP^{\loc})=L^{+}\mP^{n+1}\backslash LG^{n+1}
\end{equation}
with action given by
\[
(k_0,k_1,\dots,k_{n})\cdot (g_0,g_1,\dots,g_n) = (k_0 g_0 k_1^{-1}, k_1 g_1 k_2^{-1},\dots,k_{n-1}g_{n-1}k_{n}^{-1},k_{n}g_n\sigma(k_0)^{-1}).
\]
Similar to the case $n=0$, it represents the moduli problem:
\[
\hk_{n}(\sht_\mP^{\loc})(R) = \bigg\{\sigma^{*}_{R}\mathcal{E}_{0} \stackrel{g_n}{\dashrightarrow} \mathcal{E}_n \dashrightarrow \cdots \dashrightarrow \mathcal{E}_1 \stackrel{g_0}{\dashrightarrow} \mathcal{E}_{0} \ \bigg| \ 
\mE_i \mbox{ are } \mP\mbox{-torsors on } D_{R}.
\bigg\}.
\]
There is the important partial Frobenius endomorphism of $\hk_n(\sht_\mP^{\loc})$, which is induced by the endomorphism of $LG^{n+1}$ sending $(g_0,\ldots,g_n)$ to $(g_1,\ldots,g_n,\sigma(g_0))$.
At the level of the moduli problem, it can be described as
\begin{equation}\label{eq:partial-Frobenius}
    \pFr: \hk_n(\sht_\mP^{\loc})\to \hk_n(\sht_\mP^{\loc}), \quad ({}^\sigma\mathcal{E}_{0} \dashrightarrow \mathcal{E}_n \dashrightarrow \cdots \dashrightarrow  \mathcal{E}_{0})\mapsto ({}^\sigma\mathcal{E}_{1} \dashrightarrow {^{\sigma}\mathcal{E}_{0}} \dashrightarrow  \cdots  \dashrightarrow \mathcal{E}_{1}).
\end{equation}

One can organize $\{\Hk_n(\Sht^\loc)\}_n$ as a simplicial stack $\Hk_\bullet(\Sht^\loc)$ with the boundary maps are given by: 
\begin{equation}\label{eq: formula for face maps from iterated Sht to Sht}
d_i(g_0,\dots,g_{i},g_{i+1},\dots g_{n}) = \begin{cases}
(g_0,\dots g_{i}g_{i+1},\dots , g_n), & i\neq n\\
(g_1,g_2, \dots,g_{n}\sigma(g_0)), & i=n,
\end{cases}
\end{equation}
and degeneracy maps given by
\begin{equation}\label{eq: formula for deg maps from iterated Sht to Sht}
s_i(g_0,\dots,g_{i-1},g_{i},\dots g_{n}) = (g_0,\dots,g_{i-1},e,g_i,\dots,g_n).
\end{equation}
All the morphisms $\Hk_m(\Sht_\mP^{\loc})\to \Hk_n(\Sht_\mP^{\loc})$ are ind-pfp proper.
Note that $d_n=d_{n-1}\circ\pFr$.

\begin{example}\label{example-explicit-1-Hk}
We write $\hk_1(\sht_\mP^{\loc})$ by $\hk(\sht_\mP^{\loc})$ for simplicity. 
Then the two boundary maps $d_1,d_0\colon \hk(\sht_\mP^{\loc}) \rightarrow \sht_\mP^{\loc}$ are given by
\[
d_0(g_0,g_1) = g_0g_1, \quad d_1(g_0,g_1) = g_1 \sigma(g_0).
\]
Using the change of coordinates $g=g_0$ and $b_0=g_0g_1$ we can also identify
\begin{align}\label{eq:kot-hecke-ADLV-coordinates}
\hk(\sht_\mP^{\loc})\cong L^+\mP\backslash (\mathrm{Gr}_{\mP}\times LG),
\end{align}
where $L^{+}\mP$ acts on $\Gr_\mP$ by left translation and on $LG$ by $\sigma$-conjugation. Under this identification the boundary maps 
$d_0,d_1$
are given by sending a pair $(g,b) \in L^+\mP\backslash (\mathrm{Gr}_{\mP}\times LG)$ to
\[
d_0(g,b_0) = b_0, \quad d_1(g,b) = g^{-1} b_0 \sigma(g).
\]
Note that we also have the change of coordinates $g=g_0, b_1=g_1\sigma(g_0)$, giving
\begin{align}\label{eq:kot-hecke-ADLV-coordinates-2}
\hk(\sht_\mP^{\loc})\cong (L^+\mP\backslash LG)\times^{L^+\mP,\Ad_\sigma} LG, 
\end{align}
where $L^+\mP$ acts on $L^+\mP\backslash LG$ by right translation and on $LG$ by $\sigma$-conjugation.
Under this identification, we have
\[
d_0(g,b_1)=gb_1\sigma(g)^{-1},\quad d_1(g,b_1)=b_1.
\]
There is also the moduli theoretic interpretations of $\Hk(\sht_\mP^\loc)$. Namely, for a $k_F$-algebra $R$, $\Hk(\sht_\mP^\loc)(R)$ classify triples consisting of
\[
\Bigl\{(\mE_i,\varphi_i)\in \Sht_\mP^\loc(R), i=0,1, (\beta: \mE_1\dashrightarrow \mE_0)\in L^+\mP\backslash LG/L^+\mP\mid \varphi_0\circ \sigma_R(\beta)= \beta\circ \varphi_1\Bigr\}.
\]
I.e., $\Hk(\Sht_\mP^\loc)$ can be thought as the Hecke correspondence of $\Sht^\loc_\mP$. This justifies our notation. 
The relation between the moduli interpretation and previous discussions is encoded by the following commutative diagram
\begin{equation}\label{eq: hk-local-sht-diagram}
\xymatrix{
{}^\sigma \mE_0\ar^-{\varphi_0=b_0}@{-->}[rr]\ar^-{g_1}@{-->}[drr]&& \mE_0 \\
\ar^{\sigma(\beta)=\sigma(g_0)}@{-->}[u]{}^\sigma \mE_1\ar_-{\varphi_1=b_1}@{-->}[rr]&& \mE_1. \ar_{\beta=g_0}@{-->}[u]
}\end{equation}
\end{example}

To simplify notations, in the rest of this section, we base change all the geometric object to $k$ and omit $k$ from the subscripts, although some of them are defined over $k_F$ (or a finite extension of $k_F$). So $L^+\mP, LG$ will mean $L^+\mP_k, LG_k$ etc. in the sequel. 

Let $w\in W_\mP\backslash \widetilde W/W_\mP$. We let
\[
\Sht^{\loc}_{\mP,w}=\frac{LG_{\mP,w}}{\Ad_\sigma L^+\mP}\subset \Sht^{\loc}_{\mP,\leq w}= \frac{LG_{\mP,\leq w}}{\Ad_\sigma L^+\mP} \subset \Sht_\mP^{\loc}.
\]
As before, we write $i_{\mP,\leq w}: \Sht^{\loc}_{\mP,\leq w}\to \Sht_\mP^{\loc}$ (resp. $i_{\mP,w}: \Sht^{\loc}_{\mP,w}\to \Sht_\mP^{\loc}$) for the embedding. Each $ \Sht^{\loc}_{\mP,\leq w}$ can be ``approximated" by algebraic stacks perfectly of finite presentation. Namely, for each $n\geq 0$ and $m$ sufficiently large relative to $w$ and $n$, let
\begin{equation}\label{eq:restr local Sht} 
\Sht^{\loc(m,n)}_{\mP,\leq w}= \frac{L^{(n)}G\bs LG_{\mP,\leq w}}{\Ad_\sigma L^m\mP}
\end{equation}
Then we have $\Sht^{\loc(m',n')}_{\mP,\leq w}\to \Sht^{\loc(m',n)}_{\mP,\leq w}\rightarrow \Sht^{\loc(m',n')}_{\mP,\leq w}$ as soon as $m'\geq m$, $n'\geq n$ and $m'$ is sufficiently large relative to $w$ and $n'$. 
The first map is coh. unipotent as soon as $n\geq 1$). Namely, it is an $L^{n'}\mP/L^n\mP$-torsor. The second map is
weakly coh. pro-smooth but is not representable. In fact, it is a gerbe for the unipotent group $L^{m'}\mP/L^m\mP$ (as soon as $m\geq 1$). Then
\begin{equation}\label{eq: local Sht as lim of restr Sht}
\Sht_{\mP,\leq w}\cong \lim_{(m,n)} \Sht^{\loc(m,n)}_{\mP,\leq w}
\end{equation}

We have a similarly definied stacks $\Sht^{\loc(m,n)}_{\mP,w}$.

\begin{remark}
Note that the stack $\Sht^{\loc(m,n)}_{\mP,\leq w}$ defined as above and the one defined in \cite{xiao2017cycles} under the same notion (in the case when $\mP$ is hyperspecial) differ by a Frobenius. See \cite[Remark 4.1.9]{Zhu2016} for a discussion of this point. But the corresponding categories of  \'etale sheaves are canonically equivalent so we will ignore this difference.
\end{remark} 

For $w_0,\ldots,w_n$, one can similarly define 
\begin{equation}\label{eq:iterated-local-shtuka}
    \Sht^{\loc}_{\mP,w_0,\ldots,w_n}= \frac{LG_{\mP,w_0}\times^{L^+\mP}\cdots\times^{L^+\mP}LG_{\mP,w_n}}{\Ad_\sigma L^+\mP}\subset \hk_n(\Sht_\mP^{\loc}).    
\end{equation}
Note that the partial Frobenius \eqref{eq:partial-Frobenius} induces an isomorphism 
\begin{equation}\label{eq:pFr-iterated-local-shtuka}
\Sht^{\loc}_{\mP,w_0,\dots,w_n}\cong \Sht^{\loc}_{\mP,w_1,\ldots,w_n,\sigma(w_0)}.
\end{equation}
This are also similarly defined spaces with $LG_{\mP,w_i}$ replaced by $LG_{\mP,\leq w_i}$.

We will need the following simple observation.
\begin{lemma}\label{lem: invariant moduli of Shtuka}
Let $\mP=\mI$ be an Iwahori group scheme as in \Cref{sec:sigma-straight-element}.
If $w,w'\in\widetilde W$ are two elements such that there exist $x,y\in\widetilde{W}$ such that $w=xy$, $w'=y\sigma(x)$ and $\ell(w)=\ell(w')=\ell(x)+\ell(y)$, then $\Sht_{\mI,w}^{\loc}\cong \Sht^{\loc}_{\mI,w'}$. 
\end{lemma}
\begin{proof}
Recall that for $u_1,u_2\in \widetilde{W}$ with $\ell(u_1u_2)=\ell(u_1)+\ell(u_2)$, the multiplication map induces an isomorphism $LG_{u_1}\times^{\iw}LG_{u_2}\cong LG_{u_1u_2}$. 
Then the desired isomorphism for the first statement follows from 
\begin{equation}\label{eq:invariance-of-local-Shtuka-under-cyclic-shift}
\Sht_{\mI,w}^{\loc}\cong \Sht^{\loc}_{\mI,x,y}\cong \Sht^{\loc}_{\mI,y,\sigma(x)}\cong\Sht^{\loc}_{\mI,w'},
\end{equation}
where the isomorphism in the middle is induced by the partial Frobenius \eqref{eq:partial-Frobenius}. 
\end{proof}

\subsubsection{$\Sht^\loc_w$ for $\sigma$-straight element $w$}

Our next goal is to understand $\sht^{\loc}_w$ when $w$ is a $\sigma$-straight element. This will be the main tool for us to understand the geometry and category of sheaves on the stack of $G$-isocrystals studied later.
We fix a pinning $(B,T,e)$ of $G$ as before. Let $A\subset S\subset T$ be the corresponding tori, and let $\mI$ is the Iwahori group scheme of $G$ determined by the pinning as before. Let $\iw=L^+\mI$. As for before, we usually omit $\mP=\mI$ from the subscripts. E.g. we write $\Sht^\loc$ for $\Sht^\loc_\mI$.
Let $w\in\widetilde{W}$ be a $\sigma$-straight element.  

Let $\dot{w}$ be a lifting of $w$ to a $k$-point of $LG$. Let $\aut(\dot{w})$ be the stabilizer of $\dot{w}$ under the action of $\iw$ on $LG_{w}$ by $\sigma$-conjugation. 
Then $\aut(\dot{w})$ is an affine group scheme over $k$. Note that 
\begin{equation}\label{eq-group-Ib}
I_{\dot{w}}:=\aut(\dot{w})(k)=\{g\in \mI(\mO_{\breve F})\mid g^{-1}\dot{w}\sigma(g)=\dot{w}\},
\end{equation} 
is a profinite group. By abuse of notation, we also use it to denote the associated affine group scheme over $k$.

\begin{proposition}\label{prop: Sht-loc-w-straight}
We have $I_{\dot{w}}\cong \aut(\dot{w})$.
The morphism $\dot{w}\to \sht^{\loc}_{w}$ induces an isomorphism
\[
\bB_{\mathrm{profet}} I_{\dot{w}}\cong \Sht^\loc_{w}.
\]
\end{proposition}
Here we denote by $\bB_{\mathrm{profet}} \aut(\dot{w})$ the classifying stack of $\aut(\dot{w})$-torsors in pro-finite \'etale topology. (See \Cref{SSS:torsors convention} for our convention.)

\begin{proof} 
Let $\iw'=\iw\cap (\Ad_{\dot{w}}\sigma)^{-1}(\iw)\subset \iw$. We note that the map $\iw\to LG_w,\ g\mapsto g\cdot \dot{w}$ induces an isomorphism
\[
\frac{\iw}{\Ad_\sigma \iw'}\to \frac{LG_w}{\iw}=\sht^{\loc}_w,
\]
and  the proposition is equivalent to saying that
\[
L:\iw'\to \iw\cdot \dot{w},\quad g\mapsto g\dot{w}\sigma(g)^{-1}\dot{w}^{-1}\cdot \dot{w}
\] 
is an $\aut(\dot{w})$-torsor (in pro-finite \'etale topology). Indeed, as $\iw'=\dot{w}\times_{\sht^\loc_w}\iw$, the proposition clearly implies the above statement. Conversely, suppose $L$ is an $\aut(\dot{w})$-torsor
we show that for every morphism $x: S\to \sht^{\loc}_w$ (with $S$ qcqs), the base change $\spec k\times_{\sht^{\loc}_{w}}S\to S$ is a surjective pro-finite \'etale morphism. Now,
there is some \'etale cover $S'\to S$ such that $x$ lifts to $\tilde{x}: S'\to \iw\cdot\dot{w}$. Then the base change of $\spec k\times_{\Sht^{\loc}_w} S\to S$ along $S'\to S$ is identified with $\iw'\times_{\iw\cdot\dot{w}}S'\to S'$ which is a pro-finite \'etale morphism. Via \'etale descent of affine morphisms, we get the desired statement.

Note that a necessary condition that $L$ is a pro-finite \'etale torsor is that the map $L$ is surjective on $K$-points for any algebraically closed field $K$. 
\begin{lemma}\label{lem: single conjugacy class of straight schubert cell}
The map $L$ is surjective on $k$-points.
\end{lemma}
This follows from \cite[Theorem 3.3.1]{Gortz.He.Nie.nonempty} (generalizing \cite[Theorem 2.1.2]{gortz2010affine}). Note that in \emph{loc. cit.}, it is assumed that $G$ is tamely ramified but this assumption is not necessary.
 We sketch the arguments later.
Unfortunately, as $\iw$ and $LG_w$ are schemes that are not of perfectly finite type and $L$ is not (perfectly) finitely presented, surjectivity on $k$-points is insufficient to conclude that $L$ is surjective on $K$-points. Some extra cares are needed. The extra ingredient we need is the following.
\begin{lemma}\label{lem: principal congruence filtration of iwahori}
Let $\iw=\iw^{(0)}\supset \iw^{(1)}\supset \iw^{(2)}\supset\cdots$ be the filtration of $\iw$ by principal congruence subgroups, and let $\iw_n=\iw/\iw^{(n)}$ as before.
Let ${\iw'}^{(n)}:= \iw^{(n)}\cap (\Ad_{\dot{w}}\sigma)^{-1}(\iw^{(n)})\subset \iw^{(n)}$, and let $\iw'_n=\iw'/{\iw'}^{(n)}$. Then $\iw'={\iw'}^{(0)}\supset {\iw'}^{(1)}\supset {\iw'}^{(2)}\supset\cdots$ is a filtration of $\iw'$ by normal subgroups. In addition, $\dim \iw'_n=\dim \iw_n$.
\end{lemma}
Let we finish the proof of the proposition assuming these lemmas. The map $L$ induces a map
\[
L_n: \iw'_n\to \iw_n, \quad g\mapsto g^{-1}\dot{w}\sigma(g)\dot{w}^{-1}.
\]
This is the orbit map over $1\in \iw_n$ of the $\dot{w}\sigma$-twisted conjugation action of $\iw'_n$ on $\iw_n$. Let $\aut(\dot{w})_n\subset \iw'_n$ denote the stabilizer of $1\in \iw_n$, which is a perfect group scheme inside $\iw'_n$. 
By \cite[Proposition A.32]{zhu2017affine}, $L_n$ induces a locally closed embedding $\iw'_n/\aut(\dot{w})_n\to \iw_n$ of perfectly finitely presented algebraic spaces. As $L_n$ is surjective on $k$-points by \Cref{lem: single conjugacy class of straight schubert cell}, it is an isomorphism.
 Therefore, the morphism $L_n$ is an $\aut(\dot{w})_n$-torsor. In addition, by dimension reasons, $\aut(\dot{w})_n$ is finite. Therefore, for every algebraically closed field $K$, $L_n$ is surjective on $K$-points with finite fibers.
 
We have $\aut(\dot{w})=\varprojlim \aut(\dot{w})_n$.
Consider the following commutative diagram with the square Cartesian
\[
\xymatrix{
\iw'\ar[r]\ar_{L}[dr] &  \iw \times_{\iw_n}\iw'_n \ar[d]\ar[r] & \iw'_n\ar^{L_n}[d]\\
                    &  \iw \ar[r]              & \iw_n
}
\]
Then $\iw'=\varprojlim \iw'_n=\varprojlim \iw \times_{\iw_n}\iw'_n$ is a pseudo $\varprojlim \aut(\dot{w})_n=\aut(\dot{w})$-torsor over $\iw_w$ in pro-finite \'etale topology.  Since inverse limit of non-empty finite sets is non-empty, we see that after passing to the limit, $L$ is surjective on $K$-points. This shows that $L$ is indeed a $\aut(\dot{w})$-torsor.

It remains to prove \Cref{lem: single conjugacy class of straight schubert cell} and \Cref{lem: principal congruence filtration of iwahori}.
First, we can write $w=yw^+\sigma(y)^{-1}$ as in \eqref{eq: translate w to x}. See also \Cref{thm: reduction to min length elements} \eqref{thm: reduction to min length elements-3}. We fix a lifting $\dot{w}^+\in M(\breve F)$, $\dot{y}\in G(\breve F)$ and let $\dot{w}=\dot{y}\dot{w}^+\sigma(\dot{y})^{-1}$.

Recall that $M=Z_G(\nu_w)$ is the rational Levi associated to the $\sigma$-conjugacy class of $w$. Let $\bfP=MB$ be the standard parabolic. Let $\bfP=MU_\bfP$ be the Levi decomposition and let $U_\bfP^-$ be the unipotent radical of the opposite parabolic $\bfP^-$.
Let $\breve M_w:=Z_{G_{\breve F}}(\tilde{\nu}_w)= \dot{y}M_{\breve F}\dot{y}^{-1}$, $\breve\bfP_w=\dot{y}\bfP_{\breve F}\dot{y}^{-1}=\breve M_wU_{\breve\bfP_w}$ be the Levi decomposition. Similarly, we have $\breve\bfP_w^-=\dot{y}\bfP^-_{\breve F}\dot{y}^{-1}=\breve M_wU^-_{\breve\bfP_w}$.
Note that $\breve M_w, U_{\breve \bfP_w}, U_{\breve\bfP_w}^-$ may not be rational over $F$, but they
are invariant under $\Ad_{\dot{w}}\sigma$, and
we have the following  commutative diagrams (over $\breve F$)
\[
\xymatrix{
M\ar^-{\Ad_{\dot{y}}}[r] \ar_-{\Ad_{\dot{x}}\sigma}[d]& \breve M_w\ar^-{\Ad_{\dot{w}}\sigma}[d] && U_{\bfP}\ar^-{\Ad_{\dot{y}}}[r] \ar_-{\Ad_{\dot{x}}\sigma}[d]& U_{\breve\bfP_w}\ar^-{\Ad_{\dot{w}}\sigma}[d] && U^-_{\bfP}\ar^-{\Ad_{\dot{y}}}[r] \ar_-{\Ad_{\dot{x}}\sigma}[d]& U^-_{\breve\bfP_w}\ar^-{\Ad_{\dot{w}}\sigma}[d]\\
M\ar^-{\Ad_{\dot{y}}}[r]& \breve M_w, && U_{\bfP}\ar^-{\Ad_{\dot{y}}}[r]& U_{\breve\bfP_w},  && U^-_{\bfP}\ar^-{\Ad_{\dot{y}}}[r]& U^-_{\breve\bfP_w}.
}\]
Let $\iw_M=\iw\cap LM$, $\iw_{U_{\bfP}}=\iw\cap LU_{\bfP}$ and $\iw_{U_{\bfP}^-}=\iw\cap LU_{\mathbf{P}}^-$.
Let $\iw_{\breve M_w}=\iw\cap L\breve M_w$, and we similarly consider $\iw_{U_{\breve\bfP_w}}, \iw_{U_{\breve\bfP_w}^-}$, $\iw_{\breve\bfP_w}$ and $\iw_{\breve\bfP_w^-}$. We have the direct product decomposition (over $k$)
\begin{equation}\label{eq: triangular decomposition of iwahori}
\iw=\iw_{U_{\mathbf{P}}^-}\cdot \iw_M\cdot \iw_{U_{\mathbf{P}}}=\iw_{U_{\breve\bfP_w}^-}\cdot \iw_{\breve M_w}\cdot \iw_{U_{\breve\bfP_w}}. 
\end{equation}
We claim that $\Ad_{\dot{w}}\sigma: \breve M_w\to \breve M_w, \ \Ad_{\dot{w}}\sigma: U_{\breve\bfP_w}\to  U_{\breve\bfP_w}$ restrict to homomorphisms
\[
\Ad_{\dot{w}}\sigma: \iw_{\breve M_w}\to \iw_{\breve M_w},\quad \Ad_{\dot{w}}\sigma: \iw_{U_{\breve\bfP_w}}\to \iw_{U_{\breve\bfP_w}}.
\]
The first restriction holds as we have the isomorphism $\Ad_{\dot{y}}\colon \iw_M\cong \iw_{\breve M_w}$, which in turn follows from the fact that $y\in W_0$ is of minimal length in $yW_M$. For the second restriction,
we let $\al$ be a positive affine root $\al$ with $y^{-1}(\dot{\al})\in \Phi_{U_{\mathbf{P}}}$, where we recall $\dot{\al}$ denotes the vector part of $\al$ and $\Phi_{U_{\mathbf{P}}}$ denotes the set of finite roots whose root groups are contained in $U_{\mathbf{P}}$.
We need to show that $w\sigma(\al)$ remains to be positive affine root.
Note that for sufficiently large $n$, we have
$(w\sigma)^n=y (n \nu_w) y^{-1}$, which sends $\al$ to a positive affine root. As $w$ is $\sigma$-straight, it implies that $w\sigma(\al)$ is positive. 
We note that $\Ad_{\dot{w}}\sigma$, however, does not preserve $ \iw_{U_{\breve\bfP_w}^-}$. Rather, $(\Ad_{\dot{w}}\sigma)^{-1}$ preserves it, by the same reasoning. It follows that
\[
\iw'= \iw'_{U_{\breve\bfP_w}^-}\cdot \iw_{\breve M_w}\cdot \iw_{U_{\breve\bfP_w}},
\]
where $ \iw'_{U_{\breve\bfP_w}^-}=  (\Ad_{\dot{w}}\sigma)^{-1}(\iw_{U_{\breve\bfP_w}^-})$.

Now let $\iw=\iw^{(0)}\supset \iw^{(1)}\supset \iw^{(2)}\supset\cdots$ be the filtration of $\iw$ by principal congruence subgroups. 
The decomposition \eqref{eq: triangular decomposition of iwahori} implies that
 $\iw^{(n)}_{U_{\breve\bfP_w}^-}=\iw^{(n)}\cap \iw_{U_{\breve\bfP_w}^-}$ is the $n$th congruence subgroup of $\iw_{U_{\breve\bfP_w}^-}$, and similarly we have $\iw^{(n)}_{\breve M_w}$ and $\iw^{(n)}_{U_{\breve\bfP_w}}$. In addition,
we have the decomposition
\[
\iw^{(n)}= \iw^{(n)}_{U_{\breve\bfP_w}^-} \cdot \iw^{(n)}_{\breve M_w} \cdot  \iw^{(n)}_{U_{\breve\bfP_w}}.
\]
It follows that
\[
{\iw'}^{(n)}= \iw^{(n)}\cap (\Ad_{\dot{w}}\sigma)^{-1}(\iw^{(n)})= (\Ad_{\dot{w}}\sigma)^{-1}(\iw^{(n)}_{U_{\breve\bfP_w}^-})\cdot \iw^{(n)}_{\breve M_w}\cdot \iw^{(n)}_{U_{\breve\bfP_w}}.
\]
Then $\dim \iw'_n=\dim \iw_n$, as desired.

It remains to prove \Cref{lem: single conjugacy class of straight schubert cell}. For this, we follow \cite{gortz2010affine, Gortz.He.Nie.nonempty}: we can construct a filtration of $\iw$ by normal subgroups by refining the filtration of $\iw$ by principal congruence subgroups
\[
\iw=\iw[0]\supset \iw[1]\supset \iw[2]\supset \cdots
\]
such that $\iw\langle i\rangle:=\iw[i]/\iw[i+1]$ is one-dimensional (isomorphic to either $\bG_m$ or $\bG_a$ over $k$), and such that
\[
\iw[i]= (\iw[i]\cap \iw_{U_{\bfP_w}^-})\cdot (\iw[i]\cap \iw_{M_w})\cdot (\iw[i]\cap \iw_{U_{\bfP_w}}).
\] 
Then it is enough to show that for each $g_i\in  \iw_{\breve M_w} \iw[i] $ one can find $h_i\in \iw[i]$ such that
\[
h_i^{-1} g_i \dot{w}\sigma(h_i)\dot{w}^{-1} \in \iw_{\breve M_w} \iw[i+1].  
\] 
There are three cases. If $\iw_{M_w} \iw[i]=\iw_{\breve M_w} \iw[i+1]$, there is nothing to prove. If $\iw_{U_{\breve\bfP_w}} \iw[i]  = \iw_{U_{\breve\bfP_w}} \iw[i+1]$, then we write $g_i=u_i^-m_i u_i$, with $u_i^-, \iw[i+1]\cap \iw_{U_{\breve\bfP_w}^-}$, $m_i\in \iw_{\breve M_w}$ and $u_i\in \iw[i]\cap \iw_{U_{\breve\bfP_w}}$. Then by \cite[Lemma 3.4.1 (ii)]{Gortz.He.Nie.nonempty} (which is based on \cite[Lemma 5.1.1]{gortz2010affine})
there exists some $h_i\in \iw_{U_{\breve\bfP_w}}\cap \iw[i]$ such that $h_i^{-1}m_i u_i m_i^{-1} m_i\dot{w}\sigma(h_i)\dot{w}^{-1}m_i^{-1} \in  \iw_{U_{\breve\bfP_w}}\cap \iw[i+1]$. It follows that $h_i^{-1} g_i \dot{w}\sigma(h_i)\dot{w}^{-1}\in \iw_{\breve M_w} \iw[i+1]$. The case $\iw_{U_{\breve\bfP_w}^-} \iw[i]  = \iw_{U_{\breve\bfP_w}^-} \iw[i+1]$ is proved similarly using  \cite[Lemma 3.4.1 (i)]{Gortz.He.Nie.nonempty}.
\end{proof}

\begin{remark}\label{rem: combinatorics of Gb}
Continuing the notations of \Cref{prop: Sht-loc-w-straight}.
The automorphism 
\[
\sigma_{\dot{w}}:=\Ad_{\dot{w}}\sigma: \breve M_w\to \breve M_w
\] 
defines an $F$-rational structure on $\breve M_w$. We denote the corresponding $F$-group by $G_{\dot{w}}$. This coincides with the group $G_b$ introduced in \eqref{eq: def-of-Jb} below (for $b=\dot{w}$). 
We have 
\begin{equation}\label{eq: concrete model Gb} 
G_{\dot{w}}(F)=\{g\in G(\breve F)\mid g^{-1}\dot{w}\sigma(g)=\dot{w}\}=\{g\in \breve M_w(\breve F)\mid \sigma_{\dot{w}}(g)=g\}.
\end{equation}

The torus $S_{\breve F}\subset \breve M_w$ is stable under the Frobenius structure $\sigma_{\dot{w}}$ and therefore gives rise to a rational torus of $G_{\dot{w}}$, denoted by $S_{\dot{w}}$. 
As explained before, we have a surjective map $\scrA(G_{\breve F}, S_{\breve F})\to \scrA(\breve M_w, S_{\breve F})=\scrA((G_{\dot{w}})_{\breve F}, (S_{\dot{w}})_{\breve F})$ sending the alcove $\breve \bfa$ to an alcove $\breve \bfa_{\breve M_w}$ (see
\eqref{eq: alcove of levi}). The corresponding Iwahori subgroup is $\iw_{\breve M_w}$, equipped with the Frobenius structure given by $\sigma_{\dot{w}}$.

Passing to rational points, we see that 
\[
I_{\dot{w}}=\iw_{\breve M_w}(k)^{\sigma_{\dot{w}}}=(\iw(k)\cap L\breve M_w)^{\sigma_{\dot{w}}}
\] 
is an Iwahori subgroup of $G_{\dot{w}}(F)$. Here the second equality follows from \eqref{eq: parahoric of levi}. The Iwahori-Weyl group of $(G_{\dot{w}}(F), S_{\dot{w}}(F))$ is 
\[
\widetilde{W}^{\sigma_{w}}=\{v\in\widetilde{W}\mid w\sigma(v)w^{-1}=v\}.
\]

Now, let $\breve\bff\subset \overline{\breve\bfa}\subset \scrB(G, \breve F)$ be a facet as in \Cref{lem: wsigma stable facet}. Let $\breve\bff_{\breve M_w}$ be the corresponding facet in $\overline{\breve\bfa_{\breve M_w}}$ (see \eqref{eq: facet of levi}). Let $\breve\mP_{\breve\bff}$ (resp. $\breve\mP_{\breve\bff_{\breve M_w}}$) be the corresponding standard parahoric group schemes of $G_{\breve F}$ and $\breve M_w$. Note that  $\breve\mP_{\breve\bff_{\breve M_w}}$ is rational with respect to $\sigma_{\dot{w}}$ so
\begin{equation}\label{eq: concrete model Pb}
P_{\dot{w},{\breve\bff}}=\breve\mP_{\breve\bff_{\breve M_w}}(\breve\mO)^{\sigma_{\dot{w}}}= (\breve\mP_{\breve\bff}(\breve\mO)\cap \breve M_w(\breve F))^{\sigma_{\dot{w}}}=\{g\in \breve\mP_{\breve\bff}(\breve \mO)\mid g^{-1}\dot{w}\sigma(g)=\dot{w}\}
\end{equation}
is a standard parahoric subgroup of $G_{\dot{w}}$ (containing $I_{\dot{w}}$). 

We also notice that the Levi quotients $L_{\breve\bff}$ and $L_{\breve\bff_{\breve M_w}}$ are canonically identified, and therefore has a rational structure over $k_F$ given by $\sigma_{\dot{w}}$.
The image of $\iw$ in $L_{\breve\bff}$, denoted by $B_{L_{\breve\bff}}$, is a rational Borel of $L_{\breve\bff}$. The flag variety of $L_{\breve\bff}$ is identified with $L^+\breve\mP_{\breve\bff}/\iw$. 
\end{remark}

Let $LG_{W_{\breve\bff}w}=\cup_{w'\in W_{\breve\bff}w} LG_{w'}\subset LG$. Note that it is the inverse image of the Schubert cell $L^+\breve\mP_{\breve\bff}\backslash L^+\breve\mP_{\breve\bff} w \iw \subset L^+\breve\mP_{\breve\bff}\backslash LG$. Therefore, $LG_{W_{\breve\bff}w}$ is an affine scheme, and $LG_{W_{\breve\bff}w}\subset LG$ is a pfp locally closed embedding. In addition, notice that  since $w\sigma(W_{\breve\bff})=W_{\breve\bff}w$,
the $\sigma$-conjugation action of $L^+\breve\mP_{\breve\bff}$ on $LG$ preserves $LG_{W_{\breve\bff}w}$.

\begin{proposition}\label{prop-local-shtuka-uw} 
We have an isomorphism
\[
\bB_{\mathrm{profet}} P_{\dot{w},\mathbf{f}}\cong \frac{LG_{W_{\breve\bff}w}}{\Ad_\sigma L^+\breve\mP_{\breve\bff}}\subset \frac{LG}{\Ad_\sigma L^+\breve\mP_{\breve\bff}}.
\]
In addition, for every $uw\in W_{\mathbf{f}} w$, the fiber of $\sht^{\loc}_{\mI,uw}\to  \frac{LG_{W_{\breve\bff}w}}{\Ad_\sigma L^+\breve\mP_{\breve\bff}}$ is identified with a Deligne-Lusztig variety $X_{L_{\breve\bff},u}$ of $L_{\breve\bff}$ associated to $u\in W_{\breve\bff}$.
\end{proposition}
\begin{proof}
We follow the argument of \cite[Theorem 4.8]{he2014geometric}.
Using the Lang isogeny for $L_{\breve\bff}\to L_{\breve\bff},\ g\mapsto g^{-1}\sigma_{\dot{w}}(g)$, one sees that the map
\[
\sht^{\loc}_{\mI,w}\to \frac{LG_{W_{\breve\bff}w}}{\Ad_\sigma L^+\breve\mP_{\breve\bff}}
\] 
is finite \'etale. The first statement then follows easily from \Cref{prop: Sht-loc-w-straight}. For the second statement, we notice that the fiber of the map over $\dot{w}$ is identified with 
\[
\bigl\{g\iw \in L^+\breve\mP_{\breve\bff}/\iw \mid g^{-1}\dot{w}\sigma(g) \in \iw uw \iw\bigr\}\cong\bigl\{gB'_{L_{\breve\bff}}\in L_{\breve\bff}/B'_{L_{\breve\bff}}\mid g^{-1} \sigma_{\dot{w}}(g)\in B_{L_{\breve\bff}}uB_{L_{\breve\bff}}\bigr\}.
\] 
\end{proof}

We need the following invariant. 
Let $v\in (\breve \bff_{\breve M_w})^{\sigma_{\dot{w}}}$ be a point fixed by the $\sigma_{\dot{w}}$-action, and let  $K_{\dot{w},v}=(\breve\mM_w)_v(\breve\mO)^{\sigma_{\dot{w}}}$, which is an open compact subgroup of $G_{\dot{w}}(F)$, containing $P_{\dot{w},\breve\bff}$.
Let $v'\in\breve\bff$ be a lifting of $v$ under the map $\breve\bff\to \breve\bff_{\breve M_w}$ and let $\breve\mG_{v'}$ be its stabilizer group scheme. 

\begin{lemma}\label{eq: non-connected-local-shtuka-uw} 
We have a pfp locally closed  embedding
\begin{equation*}\label{eq: non-connected-local-shtuka-uw} 
\bB_{\mathrm{profet}} K_{\dot{w},v}\hookrightarrow \frac{LG}{\Ad_\sigma L^+\breve\mG_{v'}}.
\end{equation*}
\end{lemma}
\begin{proof}
First, by \Cref{lem: parahoric intersect with Levi}, 
\[
K_{\dot{w},v}=(\breve\mM_{w})_v(\breve\mO)^{\sigma_{\dot{w}}}=(\breve\mG_{v'}(\breve\mO)\cap \breve M_w(\breve F))^{\sigma_{\dot{w}}}=\{g\in\breve\mG_{v'}(\breve\mO))\mid g^{-1}\dot{w}\sigma(g)=\dot{w}\}.
\] 
We may write $L^+\breve\mG_{v'}=\sqcup_iL^+\breve\mP_{\breve\bff} \tau_i$, where $\tau_i$ are representatives of $\pi_0(L^+\breve\mG_{v'})$ in $L^+\breve\mG_{v'}$. In fact, we can choose $\tau_i$ to be liftings of elements in $\Omega_{\breve\bfa}$. Then $L^+\breve\mG_{v'} \cdot w\cdot \sigma(L^+\breve\mG_{v'})=\cup_{ij}L^+\breve\mP_{\breve\bff} \tau_iw\sigma(\tau_j)\sigma(L^+\breve\mP_{\breve\bff})$ is a union of connected components. It follows that
\[
\bB_{\mathrm{profet}} P'_{\dot{w},v}\cong \frac{\cup_{i} L^+\breve\mP_{\breve\bff}\tau_i w\sigma(\tau_i)^{-1}\sigma(L^+\breve\mP_{\breve\bff})}{L^+\breve\mG_{v'}}
\] 
is open and closed in $\frac{L^+\breve\mG_{v'} \cdot w\cdot \sigma(L^+\breve\mG_{v'})}{\Ad_\sigma L^+\breve\mG_{v'}}$.
\end{proof}

\subsection{The stack of $G$-isocrystals}\label{sec:kot-stack-geometry}

\subsubsection{The Kottwitz set $B(G)$}\label{SS: set BG}
In the study of mod $p$ points of Shimura varieties, Kottwitz introduced the set $B(G)$, which is the quotient of $G(\breve F)$ by itself under the $\sigma$-conjugation action
\[
\Ad_\sigma:  G(\breve F)\times G(\breve F)\to G(\breve F), \quad (b,g)\mapsto \Ad_\sigma(g)(b):=g^{-1}b\sigma(g).
\]
For our purpose, we will assume that $G$ is quasi-split over $F$, equipped with a pinning $(B,T,e)$ as before.
In this case, there an injective map
\[
B(G)\to \xcoch(T)^{+,\Ga_F}_{\bQ}\times \pi_1(G)_{\Ga_F}, \quad b\mapsto (\nu_b, \kappa_G(b)).
\]
The element $\nu_b\in \xcoch(T)^{+,\Ga_F}_{\bQ}$ is called the Newton point of $b$ and $\kappa_G(b)\in \pi_1(G)_{\Ga_F}$ is called the Kottwitz point of $b$.

Recall that
there is a partial order on $\xcoch(T)^+_{\bQ}$: we say $\nu_1\leq \nu_2$ if $\nu_2-\nu_1$ is a non-negative rational linear combination of positive (absolute) coroots of $G$ (with respect to $B,T$). The above map then induces a partial ordered on $B(G)$. We say 
\[
b\leq b'  \quad \mbox{if} \quad \kappa_G(b)=\kappa_G(b'), \mbox{ and } \nu_b\leq \nu_{b'}.
\] 
For each $b$, the set $\{b'\mid b'\leq b\}$ is finite (\cite[Proposition 2.4(iii)]{rapoport1996classification}). Minimal elements in $B(G)$ with respect to this partial order are called basic elements. The set of basic elements are denoted by $B(G)_{\mathrm{bsc}}$. The restriction of $\kappa_G$ to $B(G)_{\mathrm{bsc}}$ induces a bijection $\kappa_G: B(G)_{\mathrm{bsc}}\cong \pi_1(G)_{\Ga_F}$.

The inclusion $N_G(T)(\breve F) \rightarrow G(\breve F)$ induces maps
\begin{equation}\label{eq: str-vs-sigma-conj}
    B(\widetilde W)_{\mathrm{str}}\subset B(\widetilde{W})\to B(G).
\end{equation}
matching the Newton points and the Kottwitz points. In addition, by \cite[Theorem 3.7]{he2014geometric} the composed map is a bijection 
under which the partial order between $\sigma$-straight conjugacy classes matches the above mentioned partial order on $B(G)$ by  \cite[Theorem 3.1]{he2015kottwitzrapoport}. (Note that the article assumes that $G$ is semisimple and tamely ramified over $F$ of positive characteristic. But the identification of these two partial orders is purely a combinatoric problem related to the extended affine Weyl group $\widetilde{W}$ equipped with an action of $\sigma$, and holds without these assumptions.)
Using this bijection, we will also write $b\preceq w$ if $C\preceq w$ for the $\sigma$-straight conjugacy class $C$ corresponding to $b$.

For $b\in G(\breve F)$, there is an $F$-algebraic group defined by the functor sending an $F$-algebra $R$ to
\begin{equation}\label{eq: def-of-Jb}
G_b(R) = \{g\in G(\breve F\otimes_F R)|g^{-1}b\sigma(g) = b\}.    
\end{equation}
This group depends on $b$ up to $\sigma$-conjugation action. The set $\{G_b, b\in B(G)_{\mathrm{bsc}}\}$ is called the set of extended pure inner forms of $G$, since if $b$ is basic the group $G_b$ is naturally an inner form of $G$. The map
\[
B(G)_{\mathrm{bsc}}\cong \pi_1(G)_{\Ga_F}\to \pi_1(G_\ad)_{\Ga_F}\cong H^1(F,G_\ad)
\]
sends $b$ to the cohomology class given by $G_b$. In general, if $G$ is quasi-split, $G_b$ is naturally an extended form of $M=C_G(\nu_b)$.

\subsubsection{The stack of $G$-isocrystals}
Now we introduce the main geometric object of this work. 
Recall for our convention that for a group stack (in \'etale topology) acting on a stack $X$, the quotient stack $X/G$ is the \'etale sheafification of the prestack quotient.

\begin{definition}\label{def:kot-B(G)-stack}
For a smooth affine algebraic group $H$ over $F$, let $\kot_H$ be the prestack over $k_F$ defined as
\[
\kot_H = LH/\mathrm{Ad}_{\sigma}LH,
\]
i.e. the \'etale sheafification of the prestack quotient of $LH$ by the $\Ad_\sigma$-conjugation by $LH$. It is called the stack of $H$-isocrystals, or the stack of isocrystals with $H$-structure.
\end{definition}

\begin{lemma}\label{lem:kot-moduli}
There is a canonical isomorphism of prestacks $\kot_H\cong \mL_{\sigma}(\bB LH)$, where $\mathcal{L}_{\sigma}(\bB LH)$ is $\sigma$-fixed point prestack of $\bB LH$ (see \eqref{eq:tau-fixed-point}) defined by the pullback
\[
\begin{tikzcd}
\mathcal{L}_{\sigma}(\bB LH) \arrow[d]\arrow[r] & 
\bB LH \arrow[d,"\Delta_{\bB LH}"]\\
\bB LH \arrow[r,"\id\times \sigma"] & \bB LH\times \bB LH.
\end{tikzcd}
\]
In addition, $\kot_H$ is the moduli space assigning a perfect $k_F$-algebra $R$ the groupoid consisting of pairs $(\mE,\varphi)$, where $\mE$ is an $H$-torsor over $D_R^*$, which can be trivialized over $D_{R'}^*$ for some \'etale covering map $R\to R'$,  and $\varphi: \mE\simeq \sigma_R^*\mE$ is an isomorphism of $H$-torsors. 
\end{lemma}
\begin{proof}
The first claim is tautological.
By interpreting $\bB LH$ as the moduli functor assigning $R$ the groupoid of $H$-torsors on $D_R^*$ that can be trivialized over $D_{R'}^*$ for some \'etale cover $R\to R'$, the second statement also follows.
\end{proof}

\begin{remark}\label{rem: inperf-kot-stack}
It is possible to define an imperfect version of $\kot_H$. If $F=\bF_q[[\varpi]]$, we have for every (not necessarily perfect) $k_F$-algebra $R$, the disc $D_R=\Spec R[[\varpi]]$ and the punctured disc $D_R^*=\Spec R((\varpi))$. So the moduli problem makes sense as a prestack on $\aff_{k_F}$. However, in general, it is difficult to understand the geometry of these imperfect version (even in equal characteristic). 
\end{remark}

Now let $H=G$ be connected reductive.

\begin{remark}\label{rem: kot(G)-v-topology}
We do not know whether $G$-torsor over $D_R^*$ can be trivialized over $D_{R'}^*$ for an \'etale covering $R\to R'$. (In the non-perfect setting there exists a vector bundle on $R((\varpi))$ that cannot be trivialized \'etale locally on $R$, e.g.  see \cite[Example 5.1.24]{emerton2021schematic}. But this example does not pass in the perfect setting). On the other hand, by \cite[Lemma 11.1, Theorem 11.6]{Anschutz.extending.torsors}, every $G$-torsor on $D_R^*$ can be trivialized over $D_{R'}^*$ for an $h$-cover $R\to R'$. This suggests to further sheafify $\kot_G$ in $h$-topology to obtain a stack $\kot_G^h$ which then will represent the moduli functor sending $R$ to the groupoid of pairs $(\mE,\varphi)$ consisting of a $G$-torsor $\mE$ on $D_R^*$ and $\varphi$ is as in the above lemma. This is the stack of $G$-isocrystals considered in some literature,  e.g.\cite{Fargues.Scholze.geometrization} and \cite{Anschutz.extending.torsors}.

Our work mainly concerns the category of sheaves on the space rather than the space itself, and since $h$-sheafifcation will not change the category of sheaves by \Cref{prop-h-descent-shv}, either version of stacks of $G$-isocrystals works.
The advantage of \'etale sheafification is that it is easy to show that Newton map \eqref{eq:Newton-map-from-Sht-to-kot}  is ind-pfp proper by \Cref{lem-indproper-Nt} below so its category of sheaves can be studied via ind-proper descent. 
\end{remark}

\begin{remark}\label{rem: automorphism of sigma action}
Let $x\in G(\breve F)$. Then we can define an automorphism $\sigma_x: LG\to LG$ (over $k$) sending $g\mapsto x\sigma(g)x^{-1}$. It induces an automorphism of $\bB LG$ still denoted by $\sigma_x$. Note that $\sigma$ and $\sigma_x$ are canonical isomorphic as automorphisms of $\bB LG$. By \eqref{eq: isomorphism of loop spaces} we have a canonical isomorphism $\mL_\sigma (\bB LG)\cong \mL_{\sigma_x} (\bB LG)$ (over $k$). Explicitly, it is given by
\[
\frac{LG}{\Ad_\sigma LG}\cong \frac{LG}{\Ad_{\sigma_x} LG},\quad g\mapsto gx^{-1}.
\]
\end{remark}

\begin{proposition}\label{K-points-of-kot(G)}
For every separably closed field extensions $K_1\subset K_2$ over $k_F$, the natural functor $\kot_G(K_1)\to \kot_G(K_2)$ is an equivalence of groupoids. The
 set of isomorphism classes of $\kot_G(K_i)$ is identified with the Kottwitz set $B(G)$, and for every $b\in B(G)$, considered as a point of $\kot_G(K_i)$, its automorphism group is identified with $G_b(F)$.
\end{proposition}
\begin{proof}
If $K$ is separably closed, then $L=W_\mO(K)[1/\varpi]$ is a field of cohomological dimension one and therefore by Steinberg's theorem, every $G$-torsor over $D_K^*$ is trivial (since $G$ is connected). It follows that the groupoid $\kot_G(K)$ is given by the quotient of $G(L)$ by its $\sigma$-conjugation action, 
which is independent of $K$ by \cite[Proposition 1.16]{rapoport.zink}. (This was only stated in \emph{loc. cit.} when $F$ is a finite extension of $\bQ_p$ but the proof remains to work in general, e.g. see \cite[Lemma 1.3]{rapoport1996classification}). 
\end{proof}

As a $\mP$-torsor on $D_R$ can be trivialized over $D_{R'}$ for an \'etale covering (\cite[Lemma 1.3]{zhu2017affine}), there is a natural map 
\begin{equation}\label{eq:Newton-map-from-Sht-to-kot}
\Nt_\mP:\Sht_\mP^\loc\to\kot_G,
\end{equation}
which we call the Newton map.  
For $\Spec K\to \kot_G$, the fiber $\Nt^{-1}(x)$ is isomorphic to $\Gr_{\mP,K}$.  We note that the \v{C}ech nerve of the Newton map \eqref{eq:Newton-map-from-Sht-to-kot} is canonically identified with $\hk_{\bullet}(\sht_\mP^{\loc})$. 

\begin{lemma}\label{lem-indproper-Nt} 
The morphism $\Nt_\mP$ is ind-pfp proper in the sense of \Cref{def:ind-fp.and.proper.morphism}. 
\end{lemma}
This lemma is the main reason we define $\kot_G$ using \'etale sheafification rather than $h$-sheafification.
\begin{proof}
For every $\Spec R\to \kot_G$, there is some \'etale cover $\Spec R'$ of $\Spec R$ such that $\mE|_{D_{R'}^*}$ is trivial. Then $\sht_\mP^\loc|_{\Spec R'}\to \Spec R'$ is the affine Grassmannian of $\mP$ over $\Spec R'$. If follows that from
\Cref{lem-indfp-etale-local} that $\Sht_\mP^\loc|_{\Spec R}\to \Spec R$ is ind-pfp proper.
\end{proof}

Note that for a field valued point $b:\Spec K\to \kot_G$, 
\begin{equation}\label{eq: ADLV}
X_{\leq w}(b):= b\times_{\kot_G}\Sht^{\loc}_{\leq w}
\end{equation}
is pfp closed sub-ind-scheme $b\times_{\kot_G}\Sht^{\loc}\simeq (\Gr_{\mP})_K$, usually called the affine Deligne-Lusztig variety associated to $(b,w)$. It contains $X_w(b):=b\times_{\kot_G}\Sht^{\loc}_w$ as an open sub-ind-scheme. It is known that $X_{(\leq) w}(b)$ is in fact a scheme locally of perfectly finite presentation over $K$, and $\dim X_{(\leq) w}(b)<\infty$. It is in general a difficult question to determine for which pairs $(w,b)$, $X_{w}(b)$ is non-empty. For our purpose, we just need the following ``coarse estimate".

\begin{proposition}\label{prop: nonempty ADLV}
If $X_w(b)$ is non-empty, then $b\preceq w$.
\end{proposition}
\begin{proof}
This follows from \cite[Theorem 2.1]{he2015kottwitzrapoport}. We include a proof for completeness. 
First, notices and if $w$ and $w'$ are $\sigma$-conjugate by cyclic shift, then $X_w(b)\neq \emptyset \Leftrightarrow X_{w'}(b)\neq \emptyset$ (by  \Cref{lem: invariant moduli of Shtuka}) and $b\preceq w \Leftrightarrow b\preceq w'$.

Now we prove the proposition by induction on the length of $w$. If $w$ is of minimal length, then we may assume $w=ux$ as in  \Cref{thm: reduction to min length elements} \eqref{thm: reduction to min length elements-1}. It follows from \Cref{prop-local-shtuka-uw} that $\sht^\loc_w$ maps to the unique point $b_x$ given by the $\sigma$-straight element $x$. Therefore, $b=b_x$ and $b\preceq w$.

Now for general $w$, after $\sigma$-conjugation by cyclic shift, 
we write $w\xrightarrow{s}_\sigma w'$ for a simple reflection $s$ and $\ell(w)=\ell(w')+2$. Then 
\begin{equation}\label{eq: Sht simple reflection}
\Sht^{\loc}_w\cong \Sht^{\loc}_{s,w',\sigma(s)}\cong \Sht^{\loc}_{w',\sigma(s),\sigma(s)}\to \Sht^{\loc}_{w',\leq \sigma(s)},
\end{equation}
where the second isomorphism is given by the partial Frobenius  \eqref{eq:partial-Frobenius}.
Therefore, $X_w(b)\neq \emptyset$ implies that either $X_{w'}(b)\neq \emptyset$ or $X_{w'\sigma(s)}(b)\neq \emptyset$. As both $w'\leq w$ and $w'\sigma(s)\leq w$, the proposition follows.
\end{proof}

As a first application of ind-pfp properness of the Newton map $\Nt: \sht^{\loc}\to \kot_G$, we determine the connected components of $\kot_G$. We base change $\kot_G$ to $k$.

For every $\alpha\in \pi_1(G)_{\Ga_F}\cong \xch(Z_{\hat{G}}^{\Ga_F})$, let $\kot_{G}^{\alpha}\subset \kot_{G}$ be the subfunctor classifying those $(\mE,\varphi)$ such that at every $x\in\Spec R$, the Kottwitz point of the isomorphism class of $b_x:=(\mE_x,\varphi_x)$ is $\alpha$.

Let $|\kot_G|$ denote the topological space associated to $\kot_G$ (see \eqref{eq-top-space-prestack}).

\begin{proposition}\label{lem:local-constancy-Kot-map}
The stack $\kot_{G}^{\alpha}$ is connected and the inclusion $\kot_{G}^{\alpha}\subset\kot_G$ is open and closed.  
Therefore, there is a decomposition into connected components:
\[
\kot_G = \coprod_{\alpha\in  \pi_1(G)_{\Ga_F}} \kot_{G}^{\alpha}.
\]
\end{proposition}

\begin{proof}
The composed map $|LG|\to |\Sht^{\loc}|\to  |\kot_G|$ is a quotient map. Indeed, $LG\to \Sht^{\loc}$ is surjective strongly coh. pro-smooth (and ess. pro-unipotent), and therefore is open, see the paragraph before \Cref{composition.coh.pro-smooth}. In particular, $|LG|\to |\Sht^{\loc}|$ is a quotient map. The second map $\Sht^{\loc}\to \kot_G$ is surjective ind-pfp proper and therefore is also a submersion (see \Cref{rem-ind-proper-submersive}).
Then the claim follows immediately from the fact that the Kottwitz map induces an isomorphism $\pi_0(LG)\cong \pi_1(G)_{I_F}$, and each connected component of $LG$ is open and closed.
\end{proof}

\subsubsection{Newton stratification}\label{SS: Newton stratification}
The stack $\kot_G$ has underlying set of points given by $B(G)$. However, the stack itself is obtained by gluing these points in a non-trivial way. 

Recall that for $b\in B(G)$, there is an $F$-group $G_b$.
Then $G_{b}(F)$ is a locally profinite group, and by abuse of notation we still use it to denote the associated group ind-scheme over $k$. 
On the other hand,
we define 
\[
i_b:\kot_{G,b}\subset \kot_G,  \quad\mbox{resp. } i_{\leq b}: \kot_{G,\leq b}\subset \kot_G,\quad \mbox{resp. } i_{<b}: \kot_{G,<b}\subset\kot_G
\]
be the subfunctors consisting of those $(\mE,\varphi)$ such that for every point $x\in\Spec R$ and every geometric point $\bar x$ over $x$, the isomorphism class of $b_{\bar x}:=(\mE_{\bar x},\varphi_{\bar x})$ is equal to $b$, resp. is $\leq b$, resp. is $<b$ with respect to the partial order on $B(G)$. We factor $i_b$ as
\[
\kot_{G,b}\stackrel{j_b}{\hookrightarrow} \kot_{G,\leq b}\stackrel{i_{\leq b}}{\hookrightarrow} \kot_G.
\]

Although the above definitions look bizarre,
our goal is to prove the following result, which in particular says that $\{\kot_{G,b}\}_b$ form a stratification of $\kot_G$, called the Newton stratification.
\begin{theorem}\label{newton.stratification}
\begin{enumerate}
\item\label{newton.stratification-1} The morphism
$i_{\leq b}$ is a perfectly finitely presented closed embedding, and $j_b: \kot_{G,b}\subset\kot_{G,\leq b}$ is a quasi-compact open embedding. 
\item\label{newton.stratification-2} The closure of $\kot_{G,b}$ in $\kot_G$ is $\kot_{G,\leq b}$.
\item\label{newton.stratification-3} The morphism $j_b$ (and therefore $i_{\leq b}$) is affine.
\item\label{newton.stratification-4} We have $\kot_{G,b}\simeq \bB_{\proet} G_b(F)$.
\end{enumerate}
Here $\mathrm{proet}$ denote the pro-\'etale topology. 
\end{theorem}

The theorem is essentially known by combining various results from literature.  E.g. When $G=\GL_n$, and $F$ is in mixed characteristic, Part \eqref{newton.stratification-1} is a theorem of Grothendieck and Katz, Part \eqref{newton.stratification-2} is usually known as the (weak) Grothendieck conjecture, and Part \eqref{newton.stratification-3} is usually known as purity of Newton strata. 
We refer to \cite{Katz1979slopefil, rapoport1996classification, deJong.Oort2000, Oort2001, Vasiu2006crystalline, hartl2011newton, viehmann2013newton, Hamacher.Kim.Point.Counting} for (an incomplete list of) discussions of these results in various contents and generalities. 
Also see  \cite[Theorem 2.11]{Hamacher.Kim.Point.Counting} for Part \eqref{newton.stratification-4}. We note that in \eqref{newton.stratification-4}, one cannot replace $\proet$ by $\mathrm{profet}$ as in \Cref{prop: Sht-loc-w-straight}.\footnote{For locally profinite group $H$, $\bB_{\mathrm{profet}} H$ in general does not satisfy \'etale descent. (e.g. $H=\bZ$.)}

We here give a self-contained new proof, which provides some new information that will also be useful for later purpose. We shall mention the strategy for the proof of Part \eqref{newton.stratification-1} and \eqref{newton.stratification-2} are in fact borrowed from  \cite{HHZ}, where we prove an analogue of \Cref{newton.stratification} \eqref{newton.stratification-1} \eqref{newton.stratification-2} when $\sigma$-conjugation is replaced by the more general twisted conjugation (including the usual conjugation) of $LG$ on itself. We refer to \emph{loc. cit.} for details.
 
We first reformulate \Cref{newton.stratification}. For $b\in B(G)$, by abuse of notations, we also use it to denote the $\sigma$-conjugacy class in $G(L)$ giving by $b$, for every algebraically closed field $K$, and $L=W_\mO(K)[1/\varpi]$.  We let
\begin{equation}\label{eq: LGb in LG}
LG_{(\leq) b}:= LG\times_{\kot_G}\kot_{G,(\leq) b}.
\end{equation}
Then
\[
LG_{(\leq)b}(R)=\bigl\{g\in LG(R)\mid \forall x\in \Spec R, \ g_{\bar x} \in (b'\leq )b\bigr\},
\]
where $\bar x$ denotes a geometric point over $x$, and $g_{\bar x}$ denotes the restriction of $g$ to $\bar{x}$.
We may rephrase \Cref{newton.stratification} as follows.

\begin{theorem}\label{newton.stratification-equiv}
\begin{enumerate}
\item\label{newton.stratification-equiv-1} The morphism
$i_{\leq b}: LG_{\leq b}\to LG$ is a perfectly finitely presented closed embedding, and $j_b: LG_b\subset LG_{\leq b}$ is a quasi-compact open embedding. In particular, $LG_{(\leq)b}$ is an ind-placid scheme.
\item\label{newton.stratification-equiv-2} The closure of $LG_b$ in $LG$ is $LG_{\leq b}$.
\item\label{newton.stratification-equiv-3} The embedding $j_b$ is an affine morphism.
\item\label{newton.stratification-equiv-4} Fix $g_0\in LG(k)$ in the $\sigma$-conjugacy class $b$. Then the morphism $LG\to LG_b, \quad g\mapsto g^{-1}b\sigma(b)$ is a $G_b(F)$-torsor in pro-\'etale topology.
\end{enumerate}
\end{theorem}

Indeed, \Cref{newton.stratification} implies \Cref{newton.stratification-equiv} by base change. On the other hand, by definition $LG\to \kot_G$ is surjective in \'etale topology and as all the statements in \Cref{newton.stratification} can be checked \'etale locally, \Cref{newton.stratification-equiv} also implies \Cref{newton.stratification}.

We will prove \Cref{newton.stratification-equiv} by giving a different construction of $LG_{(\leq b)}$. To start with,
let $v,w\in \widetilde{W}$. Let $Z\subset \iw\backslash LG/\iw$ be a pfp closed embedding.
We consider the following locally closed substack $\Hk_{(\leq) v\mid (\leq) w}^{Z}(\Sht^\loc)\subset\Hk(\sht^\loc)$, classifying those as in \eqref{eq: hk-local-sht-diagram} such that $(\mE_0,\varphi_0)\in \Sht^\loc_{(\leq) v}$, $(\mE_1,\varphi_1)\in \Sht^\loc_{(\leq) w}$ and $(\beta: \mE_1\to \mE_0)\in Z$.
(Such correspondence was also considered in \cite{xiao2017cycles} at the hyperspecial level.) Let
\begin{equation*}
f_{(\leq) v,(\leq) w,Z}: \Hk_{(\leq) v\mid (\leq) w}^{Z}(\Sht^\loc)\to \shv^{\loc}_{(\leq) v}
\end{equation*}
be the morphism obtained by the restriction of $d_0: \Hk(\Sht^\loc)\to \Sht^\loc$.

\begin{lemma}
The morphism $f_{(\leq) v, \leq w, Z}$
is a representable pfp proper morphism, and  $f_{(\leq) v, w, Z}$ is representable pfp.
\end{lemma}
\begin{proof}
The second statement follows from the first as  $f_{(\leq) v, w, Z}$ is the composition of  $f_{(\leq) v, \leq w, Z}$ with a representable pfp open embedding.
Let $f_{(\leq) w,Z}: \Hk_{-\mid (\leq) w}^{Z}(\Sht^\loc)\to \sht^{\loc}$ be defined as $f_{(\leq) v,(\leq) w,Z}$, but without the requirement $(\mE_0,\varphi_0)\in \Sht^\loc_{(\leq) v}$. So $f_{(\leq) v,(\leq) w,Z}$ is given by base change along $\sht^\loc_{(\leq) v}\to \sht^\loc$ of $f_{(\leq) w,Z}$. We similarly have $f_Z: \Hk_{-\mid -}^{Z}(\Sht^\loc)\to \sht^{\loc}$.

Let $Z^{(n)}$ denote the preimage of $Z$ in $LG/\iw^{(n)}$, and let $Z^{(\infty)}$ be the preimage of $Z$ in $LG$. Note that $Z^{(0)}$ is quasi-compact and $Z^{(0)}\subset \Fl$ is a closed embedding. Therefore, $Z^{(0)}$ is pfp scheme proper over $k$.

We factorize $f_{\leq w,Z}=f_Z\circ i$, where $i: \Hk_{-\mid \leq w}^{Z}(\Sht^\loc)\to \Hk_{-\mid -}^{Z}(\Sht^\loc)$ is a pfp closed embedding (as it is the base change of $i_{\leq w}: \sht^\loc_{\leq w}\to \sht^\loc$).
Using \eqref{eq:kot-hecke-ADLV-coordinates} and \eqref{eq:kot-hecke-ADLV-coordinates-2}, we can identify this factorization as the maps in the first row of the following commutative diagram with Cartesian square. 
\begin{equation}\label{eq: adding Z and w to Hk}
\xymatrix{
 \iw\backslash Z^{(\infty)}\times^{\iw, \Ad_\sigma} LG_{\leq w} \ar^-{i}[rr]\ar[d]&& \iw\backslash(Z^{(0)}\times LG)\ar[d]\ar^-{f_Z(g,b_0)=b_0}[rr]  &&   \Sht^\loc\\
\iw\backslash LG\times^{\iw, \Ad_\sigma} LG_{\leq w}\ar^-{\id\times i_{\leq w}}[rr] &&  \iw\backslash LG\times^{\iw, \Ad_\sigma} LG\cong  \iw\backslash(\Gr\times LG) && 
}
\end{equation}
with the square Cartesian. Now as $f_Z$ is representable pfp proper and $i$ is pfp closed embedding, 
the morphism $f_{\leq w,Z}$, and therefore the morphism $f_{(\leq) v,\leq w,Z}$ is representable pfp proper. 
\end{proof}

Now let 
\[
\widetilde{f}_{(\leq )v,(\leq)w,Z}: \widetilde{\Hk}_{(\leq )v\mid (\leq)w}^{Z}(\Sht^\loc)\to LG_{\leq v}
\] 
be the base change of $f_{\leq v,(\leq)w,Z}$ along $LG_{\leq v}\to \sht^{\loc}_{\leq v}$, which in turn arises  as the base change of a representable pfp $f_{\leq v,(\leq)w,Z}^{(n)}: \Hk_{\leq v\mid (\leq)w}^{Z}(\Sht^\loc)^{(n)}\to \Gr^{(n)}_{\leq v}$ for some $n$ (depending on $v,w,Z$) by \Cref{prop:appr-fp-morphism} \eqref{prop:appr-fp-morphism-2}. We notice the following.

Notice that base change of the diagram \eqref{eq: adding Z and w to Hk} along $LG\to \Sht^\loc$ shows that 
\[
\widetilde{\Hk}_{-\mid (\leq )w}^{Z}(\Sht^\loc)\cong Z^{(\infty)}\times^{\iw,\Ad_\sigma} LG_{(\leq )w},
\] 
which is a qcqs scheme. In addition, 
the open embedding $\widetilde{\Hk}_{-\mid w}^{Z}(\Sht^\loc)\subset \widetilde{\Hk}_{-\mid \leq w}^{Z}(\Sht^\loc)$ has dense image.
 
\begin{lemma}\label{lem: affineness I}
Assume that $w$ is a $\sigma$-straight element. Then
$\widetilde{\Hk}_{\leq v\mid w}^{Z}(\Sht^\loc)$ is an affine scheme.
\end{lemma}
\begin{proof}
We already know that $\widetilde{\Hk}_{\leq v\mid w}^{Z}(\Sht^\loc)$ is a scheme. 
By \Cref{prop: Sht-loc-w-straight}, $\Sht^\loc_w\cong \bB_{\mathrm{profet}} I_{\dot{w}}$. Base change along $\Spec k\to  \bB_{\mathrm{profet}} I_{\dot{w}}$ gives the $I_{\dot{w}}$-torsor $Z^{(\infty)}\to \widetilde{\Hk}_{-\mid w}^{Z}(\Sht^\loc)$. Note that as argued in \Cref{rem: affine of LGw}, since $Z^{(\infty)}\subset LG$ is a pfp closed embedding, it is an affine scheme. It follows that that $\widetilde{\Hk}_{-\mid w}^{Z}(\Sht^\loc)$ is also affine.
\end{proof}

Next we let $LG_{(\leq) v, \leq [w], Z}\subset LG_{(\leq) v}$ be the schematic image of the pfp proper morphism $\tilde{f}_{(\leq) v,\leq w,Z}$ (see \Cref{rem-closed-subscheme}). This is a closed subset of $LG_{(\leq) v}$, and $\tilde{f}_{(\leq) v,\leq w,Z}$ factor as  
\[
\widetilde{\Hk}_{(\leq) v\mid \leq w}^{Z}(\Sht^\loc)\to LG_{(\leq) v, \leq [w], Z}\xrightarrow{i_{(\leq) v, \leq [w], Z}} LG_{(\leq) v},
\] 
with the first map being pfp proper surjective and the second map being pfp closed embedding. This follows from that fact that taking schematic image for quasi-compact morphisms commutes with flat base change (see \Cref{rem-closed-subscheme}) so
such factorization arises as the base along $LG_{(\leq) v}\to \Gr_{(\leq) v}^{(n)}$ of a similar factorization of $f_{(\leq)v,\leq w,Z}^{(n)}: \Hk_{(\leq) v\mid \leq w}^{Z}(\Sht^\loc)^{(n)}\to \Gr_{(\leq) v, \leq [w], Z}^{(n)}\to \Gr_{(\leq) v}^{(n)}$. Note that in particular, $ LG_{(\leq) v, \leq [w], Z}$ is a placid scheme over $k$.

Clearly, for $Z\subset Z'$ and $v\leq v'$, we have
\[
LG_{\leq v, \leq [w], Z}\subset LG_{\leq v, \leq [w], Z'}\subset LG_{\leq v},\quad  LG_{\leq v, \leq [w], Z}= LG_{\leq v', \leq [w], Z}\times_{LG_{\leq v'}} LG_{\leq v}.
\]

Let
\[
LG_{(\leq )v, \leq [w]}:=\colim_{Z} LG_{\leq v, \leq [w], Z}, 
\]
where the colimit is taken over the set pfp closed embeddings $Z\subset \iw\backslash LG/\iw$.
It follows that  $LG_{(\leq )v, \leq [w]}$ is an ind-scheme in $LG_{(\leq) v}$. We will soon show that $LG_{(\leq )v, \leq [w]}$ is in fact a closed subscheme in $LG_{(\leq) v}$. But let us first describe its points.

\begin{lemma}\label{lem: points of Newton bounded by v}
For every $k$-algebra $R$,
\begin{equation*}\label{eq: points of Newton bounded by v}
LG_{(\leq) v, \leq [w]}(R)=\bigl\{g\in LG_{(\leq) v}(R)\mid \forall x\in \Spec R, \ \exists\ h_{\bar x}\in LG(K_{\bar x}), h_{\bar x}^{-1}g_{\bar x}\sigma(h_{\bar x})\in LG_{\leq w}(K_{\bar x})\bigr\},
\end{equation*}
where $\bar x$ denotes a geometric point over $x$, $K_{\bar x}$ denotes the residue field of $\bar x$.
In addition $LG_{(\leq) v, \leq [w],Z}(K)\subset LG_{(\leq) v, \leq [w]}(K)$ consist of those $g$ such that $h$ can be chosen in $Z^{(\infty)}(K)$.
\end{lemma}
\begin{proof}
We first let $R=K$ be an algebraically closed field over $k$.
As $\widetilde{f}_{(\leq v), \leq w, Z}$ is pfp, every $K$-point $h$ of $LG_{(\leq) v, \leq [w]}$ lifts to a $K$-point of $\widetilde{\Hk}_{(\leq v)\mid w}^{Z}(\Sht^\loc)$ for some $Z$, which further lifts to a $K$-point of $(g, g')\in Z^{(\infty)}(K)\times LG_{\leq w}(K)$ such that $h=gg' \sigma(g)^{-1}$ (as $Z^{(\infty)}\times LG_{\leq w}\to \widetilde{\Hk}_{-\mid \leq w}^{Z}(\Sht^\loc)$ is epimorphism in \'etale topology). 

Now the case of general $R$ follows from the field valued description and \Cref{rem-closed-subscheme}.
\end{proof}
\begin{lemma}
If $v$ is minimal length in its $\sigma$-conjugacy class, then either $LG_{v, \leq [w]}=\emptyset$ or $LG_{v, \leq [w]}=LG_v$.
\end{lemma}
\begin{proof}
Note that as already noticed in \Cref{prop: nonempty ADLV}, $LG_v(K)$ is contained in one $\sigma$-conjugacy class of $LG(K)$. 
It then follows from the above description of $K$-points of $LG_{v, \leq [w]}$ that either $LG_{v, \leq [w]}(K)= \emptyset$ or $LG_{v, \leq [w]}(K)= LG_{v}(K)$. 
In the former case, $LG_{v, \leq [w]}=\emptyset$. In the latter case, we have $|LG_v|=\cup_{Z} |LG_{v, \leq[w], Z}|$ at the level of topological spaces, with each $|LG_{v, \leq[w], Z}|$ closed. Now one argues as in \Cref{locally.ind-proper.descent} to conclude that $LG_v=LG_{v,\leq[w],Z}$ for $Z$. Therefore, $LG_{v,\leq[w]}=LG_v$.
\end{proof}

\begin{lemma}
The inclusion
$LG_{(\leq )v, \leq [w]}\to LG_{(\leq) v}$ is a pfp closed embedding. In particular, $LG_{(\leq )v, \leq [w]}$ is a standard placid affine scheme.
\end{lemma}
\begin{proof}
The lemma will follow if we show that for fixed $v,w$, there is a quasi-compact closed substack $Z\subset \iw\backslash LG/\iw$ such that $LG_{v, \leq [w],Z}(K)= LG_{v, \leq [w]}(K)$.  

We prove the last statement by induction on the length of $v$. If $v$ is of minimal length in its $\sigma$-conjugacy class, this has been shown by the previous lemma. 
If $v=xy$ and $v'=y\sigma(x)$ with $\ell(v)=\ell(v')=\ell(x)+\ell(y)$, then the claim holds for $(v,w)$ if and only if it holds for $(v',w)$. Namely, suppose we can find $Z$ for $(v,w)$. Then as $LG_{y}\times^{\iw} LG_{\sigma(x)}\cong LG_{v'}$, we can write $g'\in LG_{v',\leq[w]}(K)$ as $g_1\sigma(g_2)$ for $g_1\in LG_y(K)$ and $g_2\in LG_x(K)$. Then $g:=g_2g_1\in LG_{\leq v, \leq[w], Z}(K)$ and $g'=g_2^{-1}g\sigma(g_2)\in LG_{\leq v',\leq[w], Z'}$ where $Z'$ is the image of $\iw\backslash Z^{(\infty)}\times^{\iw}LG_{x^{-1}}/\iw$ under the convolution $\iw\backslash LG\times^{\iw} LG/\iw\to LG$.

Using the similar argument, one shows that if there is some simple reflection $s$ such that $\ell(v)=\ell(sv\sigma(s))+2$, and if the statement holds for $(sv\sigma(s), w)$ and $(sv, w)$, then the statement holds for $(v,w)$. Now one uses \Cref{thm: reduction to min length elements} \eqref{thm: reduction to min length elements-1} to conclude.
\end{proof}

Now let
\[
LG_{\leq [w]}:=\colim_{v\in \widetilde{W}}LG_{\leq v, \leq [w]}.
\]
By the lemma above, $LG_{\leq [w]}$ is an ind-placid scheme in $LG$ and the morphism $i_{\leq [w]}: LG_{\leq [w]}\subset LG$ is a pfp closed embedding. 
Clearly, if $w'\leq w$, then $LG_{\leq [w]}\subset LG_{\leq [w']}$. We thus can define 
\[
LG_{[w]}:=LG_{\leq [w]}-\cup_{w'\leq w}LG_{\leq [w']}.
\]
As for fixed $w$, $\{w'\leq w\}$  is a finite set, $j_{[w]}:LG_{[w]}\to LG_{\leq [w]}$ is a quasi-compact open embedding. In particular, $LG_{[w]}$ is also a placid ind-scheme.
 
By \Cref{eq: points of Newton bounded by v}, every $k$-algebra $R$,
\begin{equation}\label{eq: points of Newton strata}
LG_{\leq[w]}(R)=\bigl\{g\in LG(R)\mid \forall x\in \Spec R, \ \exists\ h_{\bar x}\in LG(K_{\bar x}), h_{\bar x}g_{\bar x}\sigma(h_{\bar x})^{-1} \in LG_{\leq w}(K_{\bar x})\bigr\}.
\end{equation}
It follows that
\[
LG_{[w]}(R)\subset \bigl\{g\in LG(R)\mid \forall x\in \Spec R, \ \exists\ h_{\bar x}\in LG(K_{\bar x}), h_{\bar x}g_{\bar x}\sigma(h_{\bar x})^{-1} \in LG_{w}(K_{\bar x})\bigr\}.
\]
For general $w$, this inclusion is usually strict. However, we claim that
if $w$ is a $\sigma$-straight element, this inclusion is in fact an equality. In fact, we have

\begin{lemma}\label{lem: identifying two definition of Newton strata}
Let $b\in B(G)$ and $w_b\in\widetilde{W}$ a $\sigma$-straight element corresponding to $b$. Then
\[
LG_{(\leq) b}=LG_{(\leq) [w_b]}.
\]
\end{lemma}
Note that this lemma implies \Cref{newton.stratification-equiv} \eqref{newton.stratification-equiv-1}.
\begin{proof}
To see this, it is enough to check at the level of $K$-points, for $K$ an algebraically closed field over $k$. In other words, we need to show that if $h\in LG(K)$ is $\sigma$-conjugate to an element in $LG_w(K)$, then it cannot be $\sigma$-conjugate to any element in $LG_{w'}(K)$ for $w'\leq w$. By \Cref{prop: Sht-loc-w-straight}, $LG_w(K)$ are contained in a single $\sigma$-conjugacy class.
Let $b\in B(G)$ denote such conjugacy class, coming from a $k$-point of $LG$ (e.g. this point could be a lifting $\dot{w}$ of $w$.)
Then we need to show that $X_{w'}(b)(K)=\emptyset$, which follows from \Cref{prop: nonempty ADLV}.

Now if $g\in LG(K)$ belongs to $b'\leq b$, then $b'\preceq w_b$ so there is a $\sigma$-straight element $w_{b'}$ corresponding to $b'$ such that $w_{b'}\leq w_b$. Therefore, $g\in LG_{\leq [w_b]}$. Conversely, if $g\in LG_{\leq [w_b]}(K)$ so it is $\sigma$-conjugate to an lement in $LG_{w'}(K)$ for some $w'\leq w_b$. Suppose $b'$ is the $\sigma$-conjugacy class given by $g$. Then by \Cref{prop: nonempty ADLV} there is some $w_{b'}\leq w'\leq w_b$, showing that $b'\leq b$.
\end{proof}

The above argument in fact also gives the following statement.

\begin{lemma}\label{lem-ind-finite-etale-chart-of-Newton}
The Newton map $\Nt: \sht^{\loc}\to \kot_G$ restricts to an ind-pfp finite surjective morphism $\Nt_{w_b}: \sht^{\loc}_{w_b}\to \kot_{G,b}$.
\end{lemma}
We refer to \Cref{def:ind-fp.and.proper.morphism} for the notion of ind-pfp finite morphisms between prestacks. 
\begin{proof}
First we notice that $\sht_{\leq w_b}^{\loc}\to \kot_G$ is ind-pfp proper. In addition, \Cref{prop: nonempty ADLV} implies the commutative square in the following diagram is Cartesian
\begin{equation}\label{eq: Cartesian diagram from Sht to Isoc for straight}
\xymatrix{
\sht_{w_b}^{\loc}\ar^-{j_{w_b}}[r]\ar_{\Nt_{w_b}}[d] & \sht_{\leq w_b}^{\loc} \ar[d] \ar^-{i_{\leq w_b}}[r] & \sht^{\loc}\ar^-{\Nt}[dl]\\
\kot_{G,b} \ar^-{i_b}[r]& \kot_G. &
}
\end{equation}

It follows that $\Nt_{w_b}$ is ind-pfp proper. Surjectivity is clear. 
For every $\spec R\to \kot_{G,b}$, we can write $S=\sht^{\loc}_{w_b}\times_{\kot_{G,b}}\spec R$ is an ind-algebraic space $S=\colim_i S_i$ with each $S_i\to \Spec R$ pfp-proper over $\Spec R$. Now by \eqref{eq: ADLV}, the fibers of $S\to \Spec R$ over $K$-points of $\Spec R$ are isomorphic to the affine Deligne-Lusztig variety $X_{w_b}(b)$ which is well-known to be zero dimensional. It follows that each $S_i\to \Spec R$ is quasi-finite, and therefore is perfectly finite.
\end{proof}

Next we prove \Cref{newton.stratification-equiv} \eqref{newton.stratification-equiv-2}.  It is enough to show that for every point $x\in LG_{\leq b}$ admits a generalization $\eta$ in $LG_b$. Let $w_b$ be a $\sigma$-straight element corresponding to $b$.
We lift $x$ to a point $x'\in \widetilde{\Hk}^{Z}_{\leq v\mid \leq w_b }(\Sht^\loc)$ for some $(v,w,Z)$. As mentioned before, the open embedding $\widetilde{\Hk}_{-\mid w_b}^{Z}(\Sht^\loc)\subset \widetilde{\Hk}_{-\mid \leq w_b}^{Z}(\Sht^\loc)$ has dense image. Therefore, after enlarging $v$, we may assume that $x'$ admits a generalization $\eta'$ in $\widetilde{\Hk}_{\leq v\mid w_b}^{Z}(\Sht^\loc)$. Then \Cref{lem: identifying two definition of Newton strata} implies that the image of $\eta'$ is a point $\eta$ in $LG_b$, which is a generalization of $x$, as desired. (Note however, for a fixed $v$, $LG_{\leq v, [w_b]}$ may not be dense in $LG_{\leq v, \leq [w_b]}$.)

Next we prove \Cref{newton.stratification-equiv} \eqref{newton.stratification-equiv-3}.
Fix $b$ and let $w$ be a $\sigma$-straight element corresponding to $b$. Note that $\widetilde{f}_{v,\leq w_b,Z}: \widetilde{\Hk}_{\leq v\mid \leq w_b}^Z(\Sht^\loc)\to LG_{\leq v}\cap LG_{\leq b}$ restricts to $\widetilde{f}_{v,w_b,Z}: \widetilde{\Hk}_{\leq v\mid w_b}^Z(\Sht^\loc)\to LG_{\leq v}\cap LG_b$, which is a surjective pfp proper morphism. 
As $\widetilde{\Hk}_{\leq v\mid w}^Z(\Sht^\loc)$ is affine by \Cref{lem: affineness I}, so is $ LG_{\leq v}\cap LG_b$ by \cite[\href{https://stacks.math.columbia.edu/tag/05YU}{Proposition 05YU}]{stacks-project}. It follows that $LG_b\subset LG_{\leq b}$ is an affine open embedding.

Recall that $\kot_{G,b}(k)$ consists of one point by  \Cref{K-points-of-kot(G)}. By abuse of notation, we also use $b$ to denote such a $k$-point. Let 
\[
   \aut(b)=\Spec k\times_{\kot_{G,b}}\Spec k
\]
be its automorphism group, which is a closed subgroup of $LG$ (after choosing a lift of $b$ to a $k$-point of $LG$ so $\aut(b)$ is the stablizer group of the $\Ad_\sigma$-action) and therefore is a group ind-scheme. Note that $\aut(\dot{w}_b)\subset \aut(b)$.

\begin{lemma}\label{prop-geom-Newton-strata}
Let $b: \Spec k\to LG$ be as above. Recall that we regard the locally pro-finite group $G_b(F)$ as an ind-group over $k$.
Then we have $\aut(b)\cong G_b(F)$.
\end{lemma}

\begin{proof}
Clearly we have $G_b(F)\subset \aut(b)$.
\Cref{lem-ind-finite-etale-chart-of-Newton} implies that that the affine Deligne-Lusztig variety $X_{w_b}(b)$ is zero dimensional, so is an increasing union of $k$-points in the affine flag variety $LG/\iw$, and as a set is a homogeneous space of $G_b(F)$. Now $\aut(b)$ is a $\aut(\dot{w}_b)$-torsor over $X_{w_b}(b)$.
This gives the desired identification $\aut(b)= G_b(F)$. 
\end{proof}

\begin{lemma}\label{lem-ind-finite-etale-chart-of-Newton-2}
The morphism $\sht^{\loc}_{w_b}\to \kot_{G,b}$ is surjective in \'etale topology. 
\end{lemma}
\begin{proof}
First notice that the base change of $\sht^{\loc}_{w_b}\to \kot_{G,b}$  along itself $\Sht^{\loc}_{w_b}\times_{\kot_{G,b}}\Sht^{\loc}_{w_b}\to\Sht^{\loc}_{w_b}$ is \'etale. Indeed, given $\Spec R\to \Sht^{\loc}_{w_b}$, let $\Spec R'\to \Spec R$ be the $\aut(\dot{w}_b)$-torsor.  Then $\Spec R\times_{\Sht^{\loc}_{w_b}} (\Sht^{\loc}_{w_b}\times_{\kot_{G,b}}\Sht^{\loc}_{w_b})\cong \Spec R \times^{\aut(\dot{w}_b)} X_w(b)$ is a disjoint union of schemes finite \'etale over $\Spec R$. The claim follows.

Now let $\Spec R\to \kot_{G,b}$ be a map. We need to show that for every point $x\in \Spec R$, there is an \'etale map $(\Spec R',x')\to (\Spec R, x)$ such that $\Spec R'\to \Spec R$ lifts to $\Spec R'\to \Sht^{\loc}_{w_b}$.

By \Cref{lem-ind-finite-etale-chart-of-Newton}  the base change $S=\sht^{\loc}_{w_b}\times_{\kot_{G,b}}\Spec R$ is an indscheme $S=\colim_i S_i$ with each $S_i\to \Spec R$ pfp-finite (and surjective) over $\Spec R$. Without loss of generality we may assume $S_0$ is the initial member of $\{S_i\}$. Let
\[
T:=S_0\times_{\kot_{G,b}}\Sht^{\loc}_{w_b}=S_0\times_{\Spec R}S=S_0\times_{\Sht^{\loc}_{w_b}}(\Sht^{\loc}_{w_b}\times_{\kot_{G,b}}\Sht^{\loc}_{w_b})
\]
is \'etale over $S_0$. We may write $T=\cup T_j$ as union of open subsets with each $T_j$ qcqs \'etale over $S_0$. Now after a finite \'etale extension, we may assume that $x$ lifts to a point $x'$ in $S$. The preimage of $x'$ under the finite map $T\to S$ is finite. So we may choose $T_j$ such that $T_j$ contains all the preimage of $x'$ in $T$. As $T=S_0\times_{\Spec R}S=\colim_i S_0\times_{\Spec R} S_i$, we see that there is some $S_i$ such that $x'\in S_i$ and $T_j$ is contained in $S_0\times_{\Spec R} S_i$. It follows that there is an affine open neighberhood $U$ of $x'$ in $S_i$ such that $S_0\times_{\Spec R}U\subset T_j$. So $S_0\times_{\Spec R}U$ is \'etale over $S_0$. It follows that $(U,x_1)$ is \'etale over $(\Spec R, x)$, by \cite[\href{https://stacks.math.columbia.edu/tag/0BTP}{Proposition 0BTP}]{stacks-project}.
\end{proof}  
 
 Finally, we prove \Cref{newton.stratification-equiv} \eqref{newton.stratification-equiv-4}.
Given a map $\Spec R\to \kot_{G,b}$, by \Cref{lem-ind-finite-etale-chart-of-Newton-2},  there exists an \'etale cover $\Spec R'\to \Spec R$, and a lifting $\Spec R'\to  \sht^{\loc}_{w_b}$. 
Then there is a pro-\'etale cover $\Spec R''$ of $\Spec R$ such that $\Spec R\to \kot_{G,b}$ lifts to $\Spec R''\to \dot{w}_b$. It follows that we have a map $\kot_{G,b}\to \bB_{\mathrm{proet}} G_b(F)$.

On the other hand, by  \cite[Lemma 11.4]{Anschutz.extending.torsors}, the inclusion $G_b(F)\subset LG$ induces a map, for every ring $R$, the groupoid of $G_b(F)$-torsors on $\spec R$ in $v$-topology, to the groupoid of $G$-torsors on $D_R^*$. By  \cite[Lemma 11.1]{Anschutz.extending.torsors}, if the $G_b(F)$ torsor is pro-\'etale locally trivial, then the corresponding $G$-torsor on $D_R^*$ is trivial after base change along an \'etale covering $R\to R'$. This shows that we also have a morphism $\bB_{\mathrm{proet}} G_b(F)\to \kot_{G,b}$. 

Clearly, the above two maps are inverse to each other, giving the desired isomorphism $\bB_{\mathrm{proet}} G_b(F)\cong \kot_{G,b}$. This finishes the proof of \Cref{newton.stratification} and \Cref{newton.stratification-equiv}.

We extract some corollaries of the proof.

\begin{corollary}\label{cor: ind-integral}
Then map $\dot{w}_b\to \kot_{G,b}$ is ind-integral, and the map
$\bB_{\mathrm{profet}} P_{\dot{w},\breve\bff}\cong \frac{LG_{W_{\breve\bff}w}}{\Ad_\sigma L^+\breve\mP_{\breve\bff}}\to \kot_{G,b}$ from \Cref{prop-local-shtuka-uw} is ind-integral.
\end{corollary}

Let $\breve\mP$ be a standard parahoric.
For $v\in W_{\mP}\backslash \widetilde{W}/W_{\mP}$, we write
\begin{equation*}
    LG_{(\leq) v,(\leq) b}:=LG_{(\leq) v}\cap LG_{(\leq) b}= LG_{(\leq) v, (\leq)[w_b]},
\end{equation*}
where $w_b$ is $\sigma$-straight corresponding to $b$. When $\mP=\mI$, this is equal to $LG_{(\leq v), (\leq) [w_b]}$ introduced before.

\begin{corollary}\label{prop:kot-newton-admissible--fp-closed}
There is some $n\gg 0$ and subschemes $\Gr_{\mP,\leq v,b}^{(n)}\subset \Gr_{\mP,\leq v,\leq b}^{(n)}\subset \Gr_{\mP,\leq v}^{(n)}$ with $\Gr_{\mP, \leq v,\leq b}^{(n)}$ closed in $\Gr_{\mP,\leq v}^{(n)}$ and $\Gr_{\mP,\leq v,b}^{(n)}$ open in $\Gr_{\mP,\leq v, \leq b}^{(n)}$,
such that
\[
LG_{(\leq) v,(\leq) b}=\Gr_{\mP,(\leq) v, (\leq) b}^{(n)}\times_{\Gr_{\mP,(\leq) v}^{(n)}}LG_{(\leq) v}.
\]
In particular, given $v$ and $b$, there is $(m,n)$ large enough and
\[
\Sht^{\loc(m,n)}_{\mP,\leq v,b}\subset \Sht^{\loc(m,n)}_{\mP, \leq v,\leq b}\subset \Sht^{\loc(m,n)}_{\mP, \leq v}
\]
such that $LG_{\leq v,b}\subset LG_{\leq v, \leq b}\subset LG_{\leq v}$ is obtained by pullback along $LG_{\leq v}\to \Sht^{\loc(m,n)}_{\mP, \leq v}$.
\end{corollary}
\begin{proof}
As $LG_{\leq v,\leq b}\to LG_{\leq v}$ is a pfp morphism of perfect qcqs schemes, and $LG_{\leq w}=\varprojlim \Gr_{\mP, \leq w}^{(n)}$ is a placid presentation, it arises as the base change of a morphism $\Gr_{\mP, \leq w,\leq b}^{(n)}\to \Gr_{\mP, \leq w}^{(n)}$ for some $n$ large enough by \Cref{prop:appr-fp-morphism}. As $LG_{\leq w}\to \Gr_{\leq w}^{(n)}$ is an fpqc morphism, the map $\Gr_{\mP, \leq w,\leq b}^{(n)}\to \Gr_{\mP, \leq w}^{(n)}$ is necessarily a closed embedding. The statement for $\Gr_{\mP, \leq w, b}^{(n)}$ is similar.
\end{proof}

\begin{remark}\label{rem: Newton strata in Sht}
It is also an interesting/important problem to determine the closure of $\sht^{\loc}_{\mP, w, b}$. If $\mP$ is hyperspecial, it is known (at least for function fields) that for every dominant coweight $\mu$ we have the closure relation inside $\Sht^{\loc(m,n)}_{\mP,\leq \mu}$
\[
\overline{\Sht^{\loc(m,n)}_{\mP,\leq \mu,b}} = \Sht^{\loc(m,n)}_{\mP,\leq \mu,\leq b}.
\]
In addition, $\dim \Gr^{(n)}_{\mP,\leq \mu,b}=\langle\rho, \mu+\nu_b\rangle+\frac{1}{2}\mathrm{def}_G(b)+n\dim G$. 
In general, the situation is much more complicated, see  \cite{he2016hecke} for a discussion for the Iwahori case.
\end{remark}

\begin{remark}
The $h$-sheafification $\kot_G^h$ of $\kot_G$, as in \Cref{rem: kot(G)-v-topology}, also admits a Newton stratification indexed by $B(G)$. In addition, the stratum $\kot_{G,b}^h$ is then isomorphic to the classifying stack of $G_b(F)$ in $v$-topology.
\end{remark}

In the sequel, by (slightly) abuse of notation, we write 
$I_b$ instead of $I_{\dot{w}_b}$, which is an Iwahori subgroup of $G_b(F)$.

\subsection{Prelude: representations of locally profinite group}\label{SS: rep of loc profinite group}
In this subsection, we identify profinite sets as affine schemes over $k$ as before. 
If $X$ is a locally profinite set, we may write $X=\cup_i X_i$ with each $X_i$ profinite, and the inclusion $X_i\to X_j$ is obtained as the pullback of an inclusion of finite sets. Therefore, we may write $X$ as an ind-affine scheme over $k$.
The goal of this subsection is to relate the category of $\La$-valued smooth representations of a locally profinite group with the category of $\La$-sheaves on its classifying stack, under certain (mild) assumptions. The discussions here serve as a warm-up as well as preparations for the later discussion of $\shv(\kot_G,\La)$. 

\subsubsection{Representations and sheaves}\label{SS: representations and sheaves}
We let $H$ be a locally profinite group in this subsection, which admits an open compact subgroup $K\subset H$. 
We allow $\La$ to be any commutative ring at the beginning.
We recall that there is a Grothendieck abelian category $\rep(H,\La)^{\heartsuit}$ of smooth representations of $H$ on $\La$-modules. We let $\rep(H,\La)$ be the left completion of the derived $\infty$-category $\der(\rep(H,\La)^{\heartsuit})$ with respect to its standard $t$-structure, and call it the ($\infty$-)category of smooth representations of $H$. If $H$ is profinite, we also let $\cRep(H,\La)\subset \rep(H)$ be the full subcategory consisting of those representations whose underlying $\La$-module is perfect. 
By definition, there is a canonical functor
\[
L: \der(\rep(H,\La)^{\heartsuit})\to \rep(H,\La).
\]

\begin{remark}
This functor may not be an equivalence in general. Indeed, if $H$ is profinite, then we can write $H=\varprojlim_i H_i$ with each $H_i$ finite and all the transitioning maps surjective. We may regard each $H_i$ as a constant affine group scheme over $\La$ and then $H$ as a flat affine group scheme over $\La$, denoted by $\underline{H_i}_\La$ and $\underline{H}_\La$ respectively. Then we have the classifying stack $\bB_{\mathrm{fpqc}}\underline{H}_\La$ as in \Cref{ex: fpqc quotient stack}. By \Cref{lem: derived category of qcoh v.s. qcoh}, the functor $L$ is identified with the natural functor $\der(\qcoh(\bB_{\mathrm{fpqc}}\underline{H}_\La)^\heartsuit)\to \qcoh(\bB_{\mathrm{fpqc}}\underline{H}_\La)$, which may not be an equivalence. See \Cref{ex: fpqc quotient stack}.
\end{remark}

However thanks to \Cref{lem: cpt gen. of derived category of Groth ab cat}, we do not need to worry this in the study of local Langlands correspondence.

\begin{example}\label{label:assumption on profinite group}
\Cref{lem: cpt gen. of derived category of Groth ab cat} is applicable in the following situations: (1) $p$ is invertible in $\La$ and $H$ admits a pro-$p$ open compact subgroup;  (2) $\La=\overline{\bF}_p$ and $H$ admits a torsion free pro-$p$ open compact subgroup.
\end{example}

We will mostly work under the following assumption on $H$, which guarantees that \Cref{lem: cpt gen. of derived category of Groth ab cat} is applicable. See \Cref{prop:sheaves-classifying-stack-profinite-group-sheaves-rep}.

\begin{assumption}\label{label:assumption on Haar measure}
We assume that $H$ admits a $\La$-valued left Haar measure, i.e. an $H$-equivariant map $C_c^\infty(H,\La)\to \La$ that sends the characteristic function of some open compact subgroup $K\subset H$ to an invertible element in $\La$, where $C^\infty_c(H,\La)$ is the space of $\La$-valued compactly supported smooth functions on $H$, regarded as a smooth $H\times H$-representation via left and right translations, and the above $H$-equivariance is taken with respect to the $H$-representation structure on $C^\infty_c(H)$ by left translation. 
\end{assumption}

We notice that if a $\La$-valued left Haar measure on $H$ exists, then the set of $\La$-valued left Haar measures form a $\La^\times$-torsor. In this case 
$H$ also admits a $\La$-valued right Haar measure.

We will let $\Delta_H^{-1}$ denote the $\La$-line given by left Haar measures. The right translation action of $H$ on $C_c^\infty(H)$ then induces an $H$-representation on $\Delta_H^{-1}$.
Then $H$ acts on $\Delta_H$ (note the inverse) via the modular character. By abuse of notations, we still use $\Delta_H$ to denote the modular character.
I.e.
Let $\int_H dh'$ be a left Haar measure. Then 
\[
\int_H f(h'h)dh'=\Delta_H(h)\int_H f(h')dh',\quad f\in C_c^\infty(H).
\]

We say an open compact subgroup $K\subset H$ good if its volume with respect one (and therefore any) choice of left Haar measure is invertible in $\La$. In this case the compact induction $\cind_K^H\La$ is a projective object in $\rep(H)^{\heartsuit}$. In addition, any open compact subgroup $K'\subset K$ is also good. Therefore, the collection $\{\cind_K^H\}_K$ with $K$ good form a set of generators of $\rep(H)^{\heartsuit}$ as in \Cref{lem: cpt gen. of derived category of Groth ab cat}.

\begin{remark}\label{rem: adm objects in rep}
It follows from \Cref{lem: char of adm obj in cg cat} that admissible objects of $\rep(H)$ (as defined in \Cref{def:adm vs compact}) consist of those $V$ such that for every $K$ good, $V^K\in \Perf_\La$. In particular, when $\La$ is a field of characteristic zero, $\rep(H)^\adm\cap \rep(H)^{\heartsuit}$ consist of the usual admissible representations of $H$.
\end{remark}

Note that if $K\subset H$ is good, then
\[
\Delta_H(h)=\frac{\mathrm{Vol}(h^{-1}Kh)}{\mathrm{Vol}(K)}.
\]
As a consequence, if $H$ is compact, admitting a left Haar measure, then the modular character is trivial.

Fix a left Haar measure, then for every $V\in \rep(H)^{\heartsuit}$, we have an isomorphism
\[
(V\otimes\Delta_H^{-1})_K\cong V^K,\quad v\mapsto \int_K kvdk,
\]
where as usual $(-)_K$ and $(-)^K$ denote taking coinvariants and invariants.  

We discuss relations between representations and sheaves.
From on now, we regard locally profinite group $H$ as an ind-affine group scheme over $k$ as before. If no confusion will arise, we still write it as $H$ (instead of $\underline{H}_k$).
Let $\La$ be a $\bZ_\ell$-algebra as in \Cref{sec:adic-formalism}. Unless otherwise specified, we will omit $\La$ from the notation. 

First we assume that $H=K$ is profinite satisfying \Cref{label:assumption on Haar measure}. Then by \Cref{ess.pro.p.torsor.hom.descent} and \Cref{prop:sheaves-classifying-stack-profinite-group-sheaves-rep} (and \Cref{ex-constructible-on-BK}),  we see that $\bB_{\mathrm{fpqc}}K$ is placid and there is a canonical equivalence
\begin{equation}\label{eq: sheaf vs rep profinite group}
\shv_{(c)}(\bB_{\mathrm{fpqc}}K)\cong \rep_{(c)}(K).
\end{equation}
By \Cref{lem: cpt gen. of derived category of Groth ab cat}, $\shv(\bB_{\mathrm{fpqc}}K)$ is compactly generated (although $\bB_{\mathrm{fpqc}}K$ is not very placid).
If there is a Haar measure such that the volume of $K$ is invertible in $\La$, then $\shv(\bB_{\mathrm{fpqc}}K)^\cpt=\cshv(\bB_{\mathrm{fpqc}}K)$.
In general, let $K'\subset K$ be an open subgroup such that the volume of $K'$ is invertible in $\La$. Then $*$-pushforward of compact objects along $\bB_{\mathrm{fpqc}}K'\to \bB_{\mathrm{fpqc}}K$ are compact, and they generate $\shv(\bB_{\mathrm{fpqc}}K)$.
It follows that we have the fully faithful embedding 
\[
\Psi_K^L: \shv(\bB_{\mathrm{fpqc}}K)\subset \rshv(\bB_{\mathrm{fpqc}}K),
\] 
left adjoint to the tautological functor 
\begin{equation}\label{eq: renormalized representation}
\Psi_K: \rshv(\bB_{\fpqc}K)\to \shv(\bB_{\fpqc}K)
\end{equation}
obtained by ind-extension of the embedding $\cshv(\bB_{\mathrm{fpqc}}K)\subset \shv(\bB_{\mathrm{fpqc}}K)$.

Next we assume that $H$ is locally profinite satisfying \Cref{label:assumption on Haar measure}.

\begin{lemma}\label{lem: sifted altas of Blocprofingrp}
The morphism $i_K:\bB_{\mathrm{fpqc}} K\to \bB_{\mathrm{fpqc}}H$ is ind-pfp finite morphism. In particular, $\bB_{\mathrm{fpqc}}H$ is sind-placid in the sense of \Cref{def-sifted-placid-stack}.
\end{lemma}
\begin{proof}
Let $\spec R\to  \bB_{\mathrm{fpqc}}H$ be a morphism. Let $P=\Spec R\times_{\bB_{\mathrm{fpqc}}H} \bB_{\mathrm{fpqc}}K$. 
By definition, there is a faithfully flat map $R\to R'$ such that $P':=\Spec R'\times_{\Spec R} P\simeq \Spec R'\times H/K$. As argued in  \cite[Lemma 3.12]{Haines.Richarz.Weil.restriction}, one can write $P'$ as an increasing union of closed  (affine) subschemes $U_i$, such that the descent datum for $P'\to P$ restricts to each $U_i$. Then by fpqc descent of affine schemes, we see that $P$ is  an indscheme ind-pfp finite over $\Spec R$. The lemma is proved.
\end{proof}

Let $\hk_\bullet(\bB_{\mathrm{fpqc}}K)$ denote the \v{C}ech nerve of the map $\bB_{\mathrm{fpqc}} K\to \bB_{\mathrm{fpqc}}H$. Explicitly, 
\[
\hk_n(\bB_{\mathrm{fpqc}}K)=K\backslash H\times^KH\times^K\cdots\times^KH/K\cong K\backslash (H/K)^n,
\] 
where $K$ acts on $(H/K)^n$ diagonally. In particular, it is ind-placid. By writing $(H/K)^n$ as increasing union of $K$-stable finite sets, we see that 
$\shv(\hk_n(\bB_{\mathrm{fpqc}}K))$ is compactly generated. 

\begin{proposition}\label{prop: identifying shv on BG and representation of G}
The category $\shv(\bB_{\mathrm{fpqc}}H)$ is compactly generated.
There is a canonical $t$-exact, symmetric monoidal equivalence 
\[
\shv(\bB_{\mathrm{fpqc}}H)\cong \rep(H)
\] 
such that the $!$-pullback $\shv(\bB_{\mathrm{fpqc}}H)\to \shv(\bB_{\mathrm{fpqc}}K)$ is identified with the forgetful functor $\rep(H)\to \rep(K)$, and such that the $*$-pushfoward $\shv(\bB_{\mathrm{fpqc}}K)\to \shv(\bB_{\mathrm{fpqc}}H)$ is identified with the compact induction functor $\cind: \rep(K)\to \rep(H)$.
\end{proposition}
\begin{proof}
Recall that the adjunction 
\[
\cind^H_K: \rep(K)\leftrightarrows \rep(H): \res_K^H
\] 
identifies $\rep(H)$ as the category of left modules over the monad $\res_K^H\circ \cind_K^H: \rep(K)\to \rep(K)$.

On the other hand, we have
\[
\shv(\bB_{\mathrm{fpqc}}H)=|\hk_\bullet(\shv(\bB_{\mathrm{fpqc}}K))|
\]
(see \eqref{eq-shv-on-sifted-placid}, which in turn follows from \Cref{locally.ind-proper.descent}). Then compact generation of $\shv(\bB_{\mathrm{fpqc}}H)$ follows from  \cite[Proposition 5.5.7.6]{Lurie.higher.topos.theory}.
In addition, by \Cref{lem: sifted altas of Blocprofingrp} and \eqref{eq-shv-on-sifted-placid} again, we may identify $\shv(\bB_{\mathrm{fpqc}}H)$ with left modules of monad $T: \shv(\bB_{\mathrm{fpqc}}K)\to \shv(\bB_{\mathrm{fpqc}}K)$ given by $(p_2)_*(p_1)^!$ where $p_1,p_2:\hk_1(\bB_{\mathrm{fpqc}} K)=K\backslash H/K\to \bB_{\mathrm{fpqc}} K$ are two projections. 
Now under the equivalence $\shv(\bB_{\mathrm{fpqc}} K)\cong \rep(K)$ from \Cref{prop:sheaves-classifying-stack-profinite-group-sheaves-rep}, this monad is nothing but $\res^H_K\circ \cind_K^H$. This shows that $\shv(\bB_{\mathrm{fpqc}}H)\cong \rep(H)$. The rest claims of the proposition are clear.
\end{proof}

We also notice the following.
\begin{lemma}\label{lem: rep(G) tensor product}
Let $X$ be a prestack over $k$ such that $\shv(X)$ is dualizable. Then the natural functor 
\[
\shv(\bB_{\mathrm{fpqc}}H)\otimes_\La \shv(X)\to \shv(\bB_{\mathrm{fpqc}}H\times X)
\] 
is an equivalence. In particular, we have $\Rep(H)\otimes_\La \Rep(H)\cong \Rep(H\times H)$. 
\end{lemma}
Of course, by \Cref{lem: categorical kunneth prestack}, the functor in the lemma is fully faithful for any $X$.
\begin{proof}
Using \eqref{eq-shv-on-sifted-placid}, we reduce to show that $\shv(\bB_{\mathrm{fpqc}}K)\otimes_\La \shv(X)\to \shv(\bB_{\mathrm{fpqc}}K\times X)$ is an equivalence when $H=K$ is profinite (admitting a Haar measure), which follows from \Cref{cor: profinite tensor product}. 
\end{proof}

\subsubsection{Canonical duality}\label{SS: coh. duality of locally profinite gorup}
We continue to assume that $H$ is locally profinite satisfying \Cref{label:assumption on Haar measure}. We first specialize the general discussions 
 from \Cref{ex: duality via Frobenius-structure} and \Cref{SS: cpt gen category} to study
the duality of $\rep(H)$. We refer to \Cref{SS: horizontal trace} for a review of the basic theory of duality for $\La$-linear categories.
We fix a left Haar measure on $H$. 

First recall the notion of Frobenius structure on a symmetric monoidal category (see \Cref{ex: duality via Frobenius-structure}).
\begin{proposition}\label{prop: BZ-duality}
The category $\rep(H)$ admits a Frobenius structure 
\[
\can: \rep(H)\to \Mod_\La, \quad V\mapsto (V\otimes \Delta_H^{-1})_H,
\]
where $(-)_H$ denotes the (derived) functor of $H$-coinvariants. 
\end{proposition}
\begin{proof}
We need to show that the pairing
\[
 \rep(H)\otimes_\La \rep(H)\xrightarrow{\otimes} \rep(H)\xrightarrow{\can}\Mod_\La,
\]
is the co-unit of a duality datum of $\rep(H)$. In fact, we claim that the unit in this duality datum is given by the object 
\[
C_c^\infty(H)\in \rep(H\times H)\cong \rep(H)\otimes_\La\rep(H),
\] 
where as before $C_c^\infty(H)$ is the $H\times H$-representation induced by left and right translation, and where the last isomorphism is from \Cref{lem: rep(G) tensor product}. 

We need to show that 
\[
\rep(H)\xrightarrow{C_c^\infty(H)\otimes -} \rep(H)\otimes_\La\rep(H)\otimes_\La \rep(H)\stackrel{V_1\boxtimes V_2\boxtimes V_3\mapsto V_1\otimes (V_2\otimes V_3\otimes \Delta_H^{-1})_H }{\longrightarrow} \rep(H)
\]
is isomorphic to the identity functor.
Indeed, if $V$ is a representation of $H$, and if we equip $C_c^\infty(H)\otimes V$ with the tensor product $H$-representation structure where $H$ acts on $C_c^\infty(H)$ via right translation, then the map
\begin{equation}\label{eq: integration over H}
C_c^\infty(H)\otimes V\to V, \quad f\otimes v\mapsto \int_H f(h)hv dh
\end{equation}
induces an isomorphism from $(C_c^\infty(H)\otimes V\otimes\Delta_H^{-1})_H$ to $V$. 
\end{proof}

We call the induced duality from the Frobenius structure (as in \Cref{ex: duality via Frobenius-structure})
\begin{equation}\label{eq: coh duality loc pro-finite}
\verd_H^{\can}: \rep(H)^\vee\to \rep(H)
\end{equation}
the canonical self-duality of $\rep(H)$. It is usually also called the cohomological duality or  the Bernstein-Zevelensky duality. 
Let 
\[
(\verd_H^{\can})^\adm\colon (\rep(H)^{\adm})^{\op}\cong \rep(H)^{\adm},\quad (\verd_H^{\can})^\cpt\colon (\rep(H)^\cpt)^{\op}\cong \rep(H)^\cpt
\] 
be its restriction to admissible objects (see \eqref{eq: dual of admissible objects}) and compact objects (see  \eqref{eq:duality on compact objects}) respectively, both of which are anti-involutions (see \eqref{eq: involutive property of duality adm} and \eqref{eq: involutive property of duality cpt}). We describe them more explicitly.

First admissible objects in $\rep(H)$ are admissible representations of $H$ in the usual sense (see \Cref{lem:basic cpt and adm} \eqref{lem:basic cpt and adm-2}), at least when $\La$ is a field of characteristic zero and $V\in \rep(H)^{\heartsuit}$. In addition, the object $\omega^\la$ in \Cref{ex: duality via Frobenius-structure} is 
\[
\omega^\can=\Delta_H,
\] 
and the functor $\bfC^{\op}\xrightarrow{(-)^\vee} \bfC^\vee\xrightarrow{\verd^\la} \bfC$ in  \Cref{ex: duality via Frobenius-structure} is given by 
\[
\rep(H)^{\op}\to \rep(H), \quad V\mapsto \underline\Hom(V,\Delta_H)=:V^{*,\can}.
\] 
In particular, if $H$ is unimodular, i.e. $\Delta_H$ is trivial, and if once a Haar measure of $H$ is chosen, 
then $V^{*,\can}$ is the usual smooth dual of $V$. I.e. when $V$ is a smooth representation of $H$ on a free $\La$-module, then $V^{*,\can}$ is the subspace of smooth vectors in the dual space of $V$, equipped with the subspace representation structure.

Next, for an open compact subgroup $K$ such that $\mathrm{Vol}(K)$ is invertible, the induction $\cind_K^H \La$ is a compact object, and we have
\[
(\verd_H^{\can})^\cpt(\cind_K^{H} \La)\cong \cind_K^{H}\La.
\]
This follows from the fact that \eqref{eq: integration over H} is an $H$-equivariant map if we consider the $H$-action on $C_c^\infty(H)\otimes V$ only through the left translation action on $C_c^\infty(H)$, so taking $K$-invariants (which is exact by our assumption on $K$) gives
\[
(\cind_K^H\La \otimes V\otimes \Delta_H^{-1})_H\cong V^K.
\]

We also recall that there is also the Serre functor $S_{\rep(H)}$ of $\rep(H)$.
For compact object $V$ in $\rep(H)$, we have (see \eqref{eq-Serre-functor-via-duality})
\[
S_{\rep(H)}(V)=(\verd_H^\can)^\cpt(V)^{*,\can}.
\]

\begin{remark}\label{eq-unit-duality-Rep(G)}
\begin{enumerate}
\item\label{eq-unit-duality-Rep(G)-1} It follows that the horizontal trace of $\rep(H)$ (see \eqref{eq: def of hor trace}) is given by
\begin{equation}\label{eq: HH of rep of loc. profinite}
\mathrm{tr}(\rep(H))= (C_c^\infty(H)\otimes\Delta_H^{-1})_H,
\end{equation}
where now $H$ acts on $C_c^\infty(H)$ via the conjugation action of $H$ on itself. When $H$ is unimodular with a Haar measure chosen, its zeroth cohomology $H^0\mathrm{tr}(\rep(H))$ is the cocenter of the Hecke algebra $C_c^\infty(H)$. However, in general $\mathrm{tr}(\rep(H))$ may not concentrate in degree zero. We only have
\[
\mathrm{tr}(\rep(H))\in \Mod_\La^{\leq 0}.
\]

\item If $H=K$ is profinite, with the Haar measure chosen such that the volume of $H$ is one, then $(V\otimes \Delta_H^{-1})_H=V^H$.
\item Of course, $V\mapsto V_H$ is also a Frobenius structure $\la$ on $\rep(H)$, with respect to which the duality $\verd^\la$ sends $\cind_K^H\La$ to $\cind_K^H\La\otimes \Delta_H^{-1}$, but $V^{*,\la}$ then is just the usual smooth dual of $V$.
\end{enumerate}
\end{remark}

Now we explain the above duality in terms of sheaf theory.
First let $H=K$ be a pro-finite group satisfying \Cref{label:assumption on Haar measure}.
The following statement concerning the duality is clear.
\begin{lemma}\label{lem: duality on profinite group via verdier}
The object $\omega_{\bB_{\mathrm{fpqc}} K}$ is a generalized constant sheaf on $\bB_{\mathrm{fpqc}} K$ (in the sense of \Cref{def: generalized constant sheaf}), also denoted by $\La^{\can}=\omega_{\bB_{\mathrm{fpqc}} K}$. The duality 
\[
(\verd^{\can}_K)^c: \cshv(\bB_{\mathrm{fpqc}}K)^{\op}\cong \cshv(\bB_{\mathrm{fpqc}}K)
\] 
induced by $\La^{\can}$ is identified the usual contragredient duality on $\crep(K)$ which sends a representation of $K$ on perfect $\La$-module $V$ to its $\La$-linear dual $V^*$, equipped with the usual dual representation structure. It restricts to an equivalence (still denoted by $\verd^{\can}_K$)
\[
(\verd^{\can}_K)^\cpt: (\shv(\bB_{\mathrm{fpqc}}K)^{\cpt})^{\op}\cong \shv(\bB_{\mathrm{fpqc}}K)^{\cpt}.
\]
\end{lemma}

Unsurprisingly, for $H=K$ being profinite, under the equivalence \eqref{eq: sheaf vs rep profinite group}, the duality in the above lemma coincides with the duality from \eqref{eq: coh duality loc pro-finite} (so our notation is consistent).
We denote the ind-completions of the above equivalences as
\[
\verd^\can_K: \shv(\bB_{\mathrm{fpqc}}K)^\vee\cong \shv(\bB_{\mathrm{fpqc}}K), \quad \verd^\can_{K,\ind\fg}: \rshv(\bB_{\mathrm{fpqc}}K)^\vee\cong \rshv(\bB_{\mathrm{fpqc}}K).
\]

Next we suppose $H$ is locally profinite and $K$ an open compact subgroup whose volume is one with respect to a left Haar measure.
Let $\hk_\bullet(\bB_{\mathrm{fpqc}}K)$ denote the \v{C}ech nerve of the map $\bB_{\mathrm{fpqc}} K\to \bB_{\mathrm{fpqc}}H$ as before.
Then $\shv(\hk_\bullet(\bB_{\mathrm{fpqc}}K))^\cpt=\cshv(\hk_\bullet(\bB_{\mathrm{fpqc}}K))$.
In addition, by applying the construction from \Cref{ex-*-pullback-indfp-generalized-constant} to the map $d_0: \hk_\bullet(\bB_{\mathrm{fpqc}}K)\cong K\backslash (H/K)^n\to \bB_{\mathrm{fpqc}}K$, we obtain a generalized constant sheaf $\La^{\can}_\bullet$ on $\hk_\bullet(\bB_{\mathrm{fpqc}}K)$ from $\La^{\can}=\omega_{ \bB_{\mathrm{fpqc}}K}$. 
It is not difficulty to see that $\La^\can_1$ on $\hk_1(\bB_{\mathrm{fpqc}}K)=K\backslash H/K$ is canonically isomorphic to the $*$-pullback of $\La^\can$ along the face map $d_1$, satisfying a cocycle condition over $\hk_2(\bB_{\mathrm{fpqc}}K)$.
It follows that for each $n$, we have a duality
\[
(\verd_n^{\can})^c: \cshv(\hk_n(\bB_{\mathrm{fpqc}}K))^{\op}\cong \cshv(\hk_n(\bB_{\mathrm{fpqc}}K)),
\]
which commutes with pushfowards along face maps (by \Cref{star.induced.ind.fin.pres.} \eqref{star.induced.ind.fin.pres.-1}). Then by \Cref{lem: duality on profinite group via verdier} together with \eqref{eq-shv-on-sifted-placid}, we obtain the following statement.
\begin{proposition}\label{cohomological.duality.on.locally profinite}
There is a canonical equivalence
\[
(\verd^{\can}_H)^\cpt: (\shv(\bB_{\mathrm{fpqc}}H)^\cpt)^{\op}\to \shv(\bB_{\mathrm{fpqc}}H)^\cpt,
\]
which induces a self duality (by ind-extension), denoted by the same notation
\[
\verd^{\can}_H: \shv(\bB_{\mathrm{fpqc}}H)^\vee\to \shv(\bB_{\mathrm{fpqc}}H),
\]
which under the identification \Cref{prop: identifying shv on BG and representation of G}, corresponds to the above canonicall duality \eqref{eq: coh duality loc pro-finite} of $\rep(H)$. 
\end{proposition}
 \begin{proof}
 The desired duality in question is given by
 \[
 (\shv(\bB_{\mathrm{fpqc}}H)^\cpt)^{\op}\cong |(\shv(\hk_\bullet(\bB_{\mathrm{fpqc}}K))^\cpt)^{\op}|\xrightarrow{|(\verd_\bullet^\can)^\cpt|} |\shv(\hk_\bullet(\bB_{\mathrm{fpqc}}K))^\cpt|\cong \shv(\bB_{\mathrm{fpqc}}H)^\cpt. 
 \]
 By construction, the above duality sends $(\bB_{\mathrm{fpqc}}K\to \bB_{\mathrm{fpqc}}H)_*\mF$ to $(\bB_{\mathrm{fpqc}}K\to \bB_{\mathrm{fpqc}}H)_*(\verd_K^{\can})^\cpt(\mF)$ for $\mF\in \shv(\bB_{\mathrm{fpqc}}K)^\cpt$, and therefore coincides with the duality \eqref{eq: coh duality loc pro-finite}.
 \end{proof}

Alternatively, the generalized constant sheaf $\La^\can_\bullet$ on $\hk_\bullet(\bB_{\mathrm{fpqc}}K)$ as above
gives a simpliciial functor
\begin{equation*}\label{eq.global.sections.hecke.profinite group}
\Hom(\La^{\can}_\bullet,-)\colon \shv(\hk_\bullet(\bB_{\mathrm{fpqc}}K))\to \Mod_\La
\end{equation*}
given by
\[
\rg^{\can}(\Hk_n(\bB_{\mathrm{fpqc}}K),-)=\Hom_{\shv(\Hk_n(\bB_{\mathrm{fpqc}}K))}(\La_n^{\can},-)\colon \shv(\Hk_n(\bB_{\mathrm{fpqc}}K)) \rightarrow \Mod_\La, \quad [n]\in \Delta,
\]
which then induces a Frobenius structure on $\shv(\bB_{\mathrm{fpqc}}H)$
\[
\rg^{\can}(\bB_{\mathrm{fpqc}}H, - )\colon \shv(\bB_{\mathrm{fpqc}}H)= |\shv(\Hk_\bullet(\bB_{\mathrm{fpqc}}K))| \rightarrow \Mod_\La.
\]
This gives a geometric construction of the Frobenius structure in \Cref{prop: BZ-duality}.

\subsubsection{Finitely generated representations}\label{SS: f.g. representations for p-adic groups}
We assume that $\La$ is regular noetherian. Note that in this case, for a profinite group $K$ satisfying \Cref{label:assumption on Haar measure}, the subcategory $\crep(K)\subset \rep(K)$ inherits a standard $t$-structure from $\rep(K)$. The heart $\crep(K)^{\heartsuit}\subset \rep(K)^{\heartsuit}$ consist of smooth $K$-representations with underlying $\La$-modules being finitely generated.
Such $t$-structure extends to an accessible $t$-structure on $\rrep(K)$ such that $\rrep(K)^{\leq 0}$ (resp. $\rrep(K)^{\geq 0}$) is the ind-completion of $\crep(K)^{\leq 0}$ (resp. $\crep(K)^{\geq 0}$). 
The functor \eqref{eq: renormalized representation} is $t$-exact which in addition induces an equivalence $\rrep(K)^+\cong \rep(K)^+$.

To see these facts, we can write $K=\lim_i K_i$ with $K_i$ finite and $\mathrm{Vol}(\ker(K\to K_i))$ invertible. Then 
\[
\crep(K)=\colim_i \crep(K_i)\subset \colim_i \rep(K_i)=\rep(K).
\] 
Each embedding $\crep(K_i)\subset \rep(K_i)$ can be identified with $\Coh(\bB \underline{K_i}_\La)\subset \qcoh(\bB \underline{K_i}_\La)$ so the desired statements follows from discussions in \Cref{SSS: basic theory of coh}.

Our goal is to generalize these statements for $K=G(F)$ being a $p$-adic group.
 
Recall that on a quasi-compact sind-placid stack $X$, there is the category of finitely generated sheaves $\fgshv(X)$ and its ind-completion. There is always a natural functor $\fgshv(X)\to \shv(X)$, which in general may not be fully faithful (see discussions around \Cref{lem: hom of fg sheaves on sifted placid}). 
However, for $X=\bB_{\mathrm{fpqc}} H$, where $H=G(F)$ is a $p$-adic group, the category $\fgshv(\bB_{\mathrm{fpqc}} H)$ is indeed a full subcategory of $\shv(\bB_{\mathrm{fpqc}} H)$ by the following result.

\begin{proposition}\label{prop: f.g. reps as full sub of sm}
The category $\fgshv(\bB_{\mathrm{fpqc}} G(F))$ is generated (as an idempotent complete $\La$-linear category) by objects of the form $(\pi_K)_*V$, where $V\in \cshv(\bB_{\mathrm{fpqc}} K)\cong \crep(K)$ for $K\subset G(F)$ open compact subgroups. In addition, the functor $\fgshv(\bB_{\mathrm{fpqc}} G(F))\to \shv(\bB_{\mathrm{fpqc}} G(F))$ is fully faithful.
\end{proposition}

\begin{corollary}\label{cor: char zero fg=cpt}
If $\La$ is a field of characteristic zero, then $\shv(\bB_{\mathrm{fpqc}} G(F))^\cpt=\fgshv(\bB_{\mathrm{fpqc}} G(F))$.
\end{corollary}

To prove this proposition, we need the following input.
Recall the extended Bruhat-Tits buidling $\scrB^{\mathrm{ext}}(G,F)$ of $G(F)$ (the $\sigma$-fixed points of \eqref{eq: ext BT building}) is a \emph{contractible} simplicial complex acted by $G(F)$ by simplicial automorphisms.
By barycentric subdivision of $\scrB^{\mathrm{ext}}(G,F)$, 
there a finite subcomplex $\Sigma$ (contained in the closure of an alcove) which is the fundamental domain for the $G(F)$-action.
In addition, for each cell $\sigma\subset \Sigma$, the group $K_\sigma=\{g\in G(F)\mid g\sigma=\sigma\}$ is open compact subgroup and  in fact fixes every point of $\sigma$. 
Let $\frakC_{\Sigma}$ be the partially ordered set of simplices in $\Sigma$, regarded as an ordinary category. I.e. for $\sigma, \sigma'\in \frakC_{\Sigma}$, there is a unique arrow from $\sigma$ to $\sigma'$ if $\sigma\subset \bar{\sigma'}$.

The cellular complex computing the homology of $\scrB^{\mathrm{ext}}(G,F)$ gives a resolution of the trivial $G(F)$-module $\La$
\begin{equation}\label{eq:canonical resolution of triv G-module via building-2}
0\to V_l\to \cdots\to V_1\to V_0\to \La\to 0,
\end{equation}
where $V_i$ is the smooth representation of $G(F)$ on the free $\La$-module spanned by cells of dimension $i$.
As $\Sigma$ is a fundamental domain, we have 
\[
V_i=\bigoplus_\sigma \cind_{K_\sigma}^{G(F)}\La,
\]
where the sum is taken over all faces $\sigma\subset \Sigma$ of dimension $i$. We may interpret \eqref{eq:canonical resolution of triv G-module via building-2} as follows.

\begin{lemma}\label{lem:canonical resolution of triv G-module via building}
We have a natural isomorphism
\begin{equation*}\label{eq:canonical resolution of triv G-module via building}
\La\cong \colim_{\frakC^{\op}_{\Sigma}} \cind_{K_{\sigma}}^{G(F)} \La.
\end{equation*}
in $\rep(G(F))$.
\end{lemma}

\begin{proof}[Proof of \Cref{prop: f.g. reps as full sub of sm}]
We consider the full idempotent complete subcategory
\begin{equation}\label{eq: fg sheaf from open compact}
\fgshv(\bB_{\mathrm{fpqc}} G(F))'\subset \fgshv(\bB_{\mathrm{fpqc}} G(F))
\end{equation}
spanned by objects of the form $(i_K)_*V$ for $V\in \cshv(\bB_{\mathrm{fpqc}}K)$ and $K$ open compact.
We first show that the composed functor $\fgshv(\bB_{\mathrm{fpqc}} G(F))'\subset \fgshv(\bB_{\mathrm{fpqc}} G(F))\to \shv(\bB_{\mathrm{fpqc}} G(F))$ is fully faithful.
It is enough to show that
\[
\Hom_{\rshv(\bB_{\mathrm{fpqc}}G(F))}\bigl((\pi_{K_1})_*V_1, (\pi_{K_2})_*V_2\bigr)\cong \Hom_{\shv(\bB_{\mathrm{fpqc}}G(F))}\bigl((\pi_{K_1})_*V_1, (\pi_{K_2})_*V_2\bigr).
\]
That is, we need to show that \eqref{eq: hom of fg sheaves on sifted placid} and \eqref{eq: hom of  sheaves on sifted placid} are isomorphic, for $\mF_i=V_i\in \cshv(K_i)$. We may assume that $V_i\in \crep(K)^{\heartsuit}$. 
Now we notice that the morphism \[
\bB_{\mathrm{fpqc}} K_1\times_{\bB_{\mathrm{fpqc}} G(F)}\bB_{\mathrm{fpqc}} K_2=K_1\backslash G(F)/K_2\to \bB_{\mathrm{fpqc}} K_1
\] 
is ind-finite so using notions there the sheaves $(g_{1j})_*(g_{2j})^!\mF_2\in \cshv(K_1)$ are in the heart of the standard $t$-structure of $\shv(K_1)=\rep(K_1)$. But $\Hom_{\rep(K_1)}(V_1,-)$ does commute with filtered colimits when restricted to $\rep(K_1)^{\geq 0}$. So  \eqref{eq: hom of fg sheaves on sifted placid} and \eqref{eq: hom of  sheaves on sifted placid} are indeed isomorphic in our setting. 

Therefore, it remains to show that \eqref{eq: fg sheaf from open compact} is essentially surjective.  Let $f: X\to \bB_{\mathrm{fpqc}}G(F)$ be an ind-pfp morphism with $X$ being a quasi-compact placid stack and 
let $\mF\in \cshv(X)$. We need to show that $f_*\mF$ belongs to $\fgshv(\bB_{\mathrm{fpqc}}G(F))'$.

Let $\widetilde{X}\to X$ be the $G(F)$-torsor (in fpqc topology) corresponding $f$. Now, for every $\sigma$ as in \Cref{lem:canonical resolution of triv G-module via building}, we have an ind-pfp proper (in fact ind-finite) morphism 
\[
\pi_{X,\sigma}: X_\sigma:=\widetilde{X}/K_\sigma\to X,
\] 
which is the base change of $i_{K_\sigma}: \bB_{\mathrm{fpqc}}K_\sigma\to \bB_{\mathrm{fpqc}} G(F)$. Fibers of this morphism are isomorphic to the discrete set $G(F)/K_\sigma$ (regarded as an ind-affine scheme).
In addition, by \Cref{lem-indfp-to-placid-stacks} we may write $X_\sigma=\colim_j X_{\sigma,j}$ with each $X_{\sigma,j}$ being a placid stack, each map $\pi_{X,\sigma,j}: X_{\sigma,j}\to X$ being representable pfp proper, and each map $X_{\sigma,j}\to \bB_{\mathrm{fpqc}} K_\sigma$ being representable pfp.

\Cref{lem:canonical resolution of triv G-module via building} implies that $\omega_{\bB_{\mathrm{fpqc}} G(F)}$ belongs to $\fgshv(\bB_{\mathrm{fpqc}} G(F))'\subset \shv(\bB_{\mathrm{fpqc}} G(F))$, since $ \colim_{\frakC^{\op}_{\Sigma}} \cind_{K_{\sigma}}^{G(F)} \La$ is a finite colimit. 
It follows that in $\rshv(X)$, we can write
\[
\consdual_X=\colim_\sigma (\pi_{X,\sigma})^{\ind\fg}_*\consdual_{X_\sigma}^{\ind\fg}=\colim_{\sigma,j} (\pi_{X,\sigma,j})_*\consdual_{X_{\sigma,j}}.
\] 
Then by the projection formula, we have the isomorphism in $\rshv(X)$
\[
\mF\cong \mF\otimes^!\omega_X\cong \colim_{\sigma,j} (\pi_{X,\sigma,j})_*((\pi_{X,\sigma,j})^!\mF).
\]
 As $\mF$ is compact in $\rshv(X)$, by a simple fact in category theory given in \Cref{lem: compact object as a colimit} below, we see that $\mF$ is a retract of some $(\pi_{X,\sigma,j})_*((\pi_{X,\sigma,j})^!\mF)$.

Now each object $(\pi_{X,\sigma,j})^!\mF$ belongs to $\cshv(X_{\sigma,j})$, and its $*$-pushforward to $\bB_{\mathrm{fpqc}}K_\sigma$ is constructible. 
It follows that that each $(\pi_{X,\sigma,j})_*((\pi_{X,\sigma,j})^!\mF)$ belongs to $\fgshv(\bB_{\mathrm{fpqc}} G(F))'$. Then the $*$-pushforward of $\mF$ to $\bB_{\mathrm{fpqc}}G(F)$ is a retract of one of these objects and therefore belongs to $\fgshv(\bB_{\mathrm{fpqc}} G(F))'$. This proves the proposition.
\end{proof}

The following lemma is used in the above proof.
\begin{lemma}\label{lem: compact object as a colimit}
Let $\bfC$ be a presentable category and $c\in \bfC$ a compact object. Suppose we can write $c=\colim_i c_i$ in $\bfC$. Then there is some $i$ such that $c$ is a retract of $c_i$.
\end{lemma}

We write 
\[
\fgrep(G(F))\subset \rep(G(F))
\] 
corresponding to $\fgshv(\bB_{\mathrm{fpqc}}G(F))\subset \shv(\bB_{\mathrm{fpqc}}G(F))$, and let 
\[
\rrep(G(F))=\ind \fgrep(G(F)).
\] 
I.e. $\fgrep(G(F))$ is the idempotent complete $\La$-linear subcategory of $\rep(G(F))$ generated by $\cind_K^{G(F)}V$, where $V\in \crep(K)$. By ind-extension of the embedding $\fgrep(G(F))\to \rep(G(F))$, we obtain
\begin{equation}\label{eq: tautological functor for rep of $p$-adic groups}
\Psi_{G(F)}: \rrep(G(F))\to \rep(G(F)),
\end{equation}
which admits a left adjoint $\Psi_{G(F)}^L$.

We give some corollary of (the proof of) \Cref{prop: f.g. reps as full sub of sm} .

\begin{proposition}\label{cor: colim presentation of fgshv}
The natural functor
\[
F:\colim_{\frakC^{\op}_\Sigma} \rrep(K_\sigma)\cong \rrep(G(F)),
\] 
with transitioning functors being (compact) induction, is an equivalence.
\end{proposition}
\begin{proof}
We pass to the right adjoint to prove that 
\begin{equation}\label{eq: lim presentation of fgshv}
F^R: \rrep(G(F))\to \lim_{\frakC_\Sigma} \rrep(K_\sigma)
\end{equation}
is an equivalence.

First, the proof of \Cref{prop: f.g. reps as full sub of sm} shows that for every $\mF\in \rshv(\bB_{\fpqc}G(F))$, we have
\[
F\circ F^R(\mF)=\colim_\sigma (\pi_\sigma)^{\ind\fg}_*((\pi_\sigma)^{\ind\fg,!}\mF)\cong \mF\os \consdual_{\bB_{\fpqc}G(F)}\cong \mF
\]
Therefore $F^R$ is fully faithful. 

\quash{Similarly, for every quasi-compact placid stack $X$ equipped with an ind-pfp map $X\to \bB_{\fpqc}G(F)$, the functor 
\[
F_X^R: \rshv(X)\to \lim_{\frakC_\Sigma} \rshv( X_\sigma)
\] 
is fully faithful. In fact, as $\frakC_\Sigma$ is a finite category, the same fully faithfulness holds if $X$ is replaced by a quasi-compact ind-placid stack. }

It remains to prove that $F^R$ is essential surjective. Recall that for every $\sigma$, there is a $t$-structure on $\rrep(K_\sigma)$, and the $!$-pullbacks preserve bounded from below subcategories (in fact they are $t$-exact).  Therefore, it is enough to show that $\lim_{\frakC_\Sigma} \rrep(K_\sigma)^+$ is contained in the essential image of $F^R$. But the functor $\Psi_{K_\sigma}$ restricts to an equivalence $\rrep(K_\sigma)^+\cong \rep(K_\sigma)^+$,
we can deduce it from the similar version with $\rshv$ replaced by $\shv$. Namely, we claim that
\[
\rep(G(F))\to \lim_{\frakC_\Sigma} \rep(K_\sigma)
\]
is an equivalence. Indeed, by \eqref{eq:Shv-send-colim-to-lim} the right hand side computes $\shv(\colim_{\frakC_{\Sigma}^{\op}}\bB_{\fpqc} K_\sigma)$, where the colimit is taken in the category of prestacks over $k$. However, after $h$-sheafification, $\colim_{\frakC_{\Sigma}^{\op}}\bB_{\fpqc} K_\sigma$ is isomorphic to $\bB_{\fpqc} G(F)$. Therefore, the desired equivalence follows from the $h$-sheaf property of $\shv$ as explained in \Cref{prop-h-descent-shv}.
\end{proof}

\begin{corollary}\label{rem: coh duality on f.g. representations}
When $H=G(F)$ is a $p$-adic group, 
there is a duality 
\[
(\verd_{G(F)}^{\can})^{\mathrm{f.g.}}: \fgshv( \bB_{\mathrm{fpqc}} G(F))^{\op}\cong \fgshv( \bB_{\mathrm{fpqc}} G(F)),
\]
which restricts to the duality in \Cref{cohomological.duality.on.locally profinite}. This functor sends $\cind_K^{G(F)}\La$ to itself, for \emph{any} open compact subgroup $K\subset G(F)$. 
\end{corollary}
\begin{proof}
Notice that the functors in the colimit presentation of $\fgrep(G(F))$ as in \Cref{cor: colim presentation of fgshv} are compatible with the duality $(\verd^\can_{K_\sigma})^c: \crep(K_\sigma)^{\op}\to\crep(K_\sigma)$, and therefore induce the desired duality on $\fgrep(G(F))$.
\end{proof}

\quash{When $\La$ is a regular noetherian ring, the standard $t$-structure on $\rep(G(F))$ restricts to a standard $t$-structure on $\fgrep(G(F))$. It is easy to show (e.g. as argued in \Cref{lem: t-structure renomralized}) that the tautogolical functor $\Psi: \rrep(G(F))\to \rep(G(F))$ restricts to an equivalence
\[
\Psi^{\geq n}: \rrep(G(F))^{\geq n}\cong \rep(G(F))^{\geq n}.
\]
In particular, $\rep(G(F))$ is the left completion of $\rrep(G(F))$.
Note that when $\La$ is a field of characteristic zero, we have $\fgshv(G(F))=\shv(G(F))^{\cpt}$ so
$\rrep(G(F))\cong \rep(G(F))$.
}

\begin{remark}\label{rem: fpqc vs proet}
All the discussions in this subsection remain unchanged if one replaces $\bB_{\mathrm{fpqc}}(-)$ by $\bB_{\mathrm{proet}}(-)$.
\end{remark}

\begin{proposition}
Let 
\[
\ind\fgrep(G(F))^{\leq 0},\quad \mbox{resp.} \ \ \ \ind\fgrep(G(F))^{\geq 0}
\] 
consist of those $V$ such that for every $\sigma\in \frakC_\Sigma$, the image of $V$ under the forgetful functor $\ind\fgrep(G(F))\to \ind\crep(K_\sigma)$ (the right adjoint of the induction functor  $\ind\crep(K_\sigma)\to \ind\fgrep(G(F))$) belongs to $\ind\crep(K_\sigma)^{\leq 0}$ (resp. $\ind\crep(K_\sigma)^{\geq 0}$). 

Then this pair defines a $t$-structure on $\ind\fgrep(G(F))$, which is right complete, compatible with filtered colimits (i.e.  $\ind\fgrep(G(F))^{\geq 0}$ is closed under filtered colimits).  In addition, the functor \eqref{eq: tautological functor for rep of $p$-adic groups} is $t$-exact, which restricts to an equivalence $\rrep(G(F))^{\geq n}\cong \rep(G(F))^{\geq n}$ for each $n$.
\end{proposition}
\begin{proof}
By \Cref{cor: colim presentation of fgshv}, we have $\ind\fgrep(G(F))=\lim_{\frakC_\Sigma}\ind\crep(K_\sigma)$. As mentioned before, each $\ind\crep(K_\sigma)$ has a standard accessible, right complete $t$-structure compatible with filtered colimits, and the forgetful functor $\ind\crep(K_\sigma)\to \ind\crep(K_{\sigma'})$ is $t$-exact (for $\sigma'\subset \overline{\sigma}$). Then the statements of the proposition follow from standard facts about $t$-structures
(e.g. \cite[Lemma 3.1.5.8]{Gaitsgory.Rozenblyum.DAG.vol.I}).
\end{proof}

\subsection{The local Langlands category}\label{SS: local Langlands category}

Now we introduce and study the category $\shv(\kot_G,\La)$ of sheaves on $\kot_G$, which we call the local Langlands category. We will also introduce a variant $\rshv(\kot_G,\La)$. As before, let $\La$ be a $\bZ_\ell$-algebra as in \Cref{sec:adic-formalism}. But unless otherwise specified, we will omit $\La$ from the notation. We will base change all the geometric objects to $k$. But as before, we omit $k$ from the notations.

\subsubsection{The category $\shv(\kot_G)$}\label{SS: definition of LL category} 
We let $\mP$ be a parahoric group scheme of $G$ defined over $\mO_F$.
\begin{proposition}\label{prop: shv on kot via sht}
The category $\shv(\Sht_\mP^{\loc})$ is compactly generated.
The category $\shv(\kot_G)$ is compactly generated.  
There is a pair of adjoint (continuous) functors
\[
(\Nt_\mP)_*:\shv(\Sht_\mP^\loc)  \rightleftharpoons \shv(\kot_G): (\Nt_\mP)^!.
\]
\end{proposition}
\begin{proof}
By the discussion of \Cref{sec.moduli.local.shtukas}, the \v{C}ech nerve the Newton map $\Nt_\mP$ is given by the Hecke groupoid $\Hk_\bullet(\Sht_\mP^{\loc})$, whose $n$th term (for $n\geq 0$) is given by \eqref{eq: nth local Hecke stack for Sht}. Moreover, by \eqref{eq:iterated-local-shtuka} each $\Hk_\bullet(\Sht_\mP^{\loc})$ is an ind-very placid stack written as an increasing union of quotients of standard placid schemes by the group $L^{+}\mP^{n+1}$. By \Cref{compact.generation.admissible.stacks}, for every $n\geq 0$ the category $\shv(\Hk_n(\Sht_\mP^{\loc}))$ is compactly generated. Ind-proper descent gives an equivalence (see \eqref{eq-shv-on-sifted-placid}):
\begin{equation}\label{eq-shv-on-B(G)-as-colim}
\shv(\kot_G)\simeq \tot\big(\shv(\Hk_\bullet(\Sht_\mP^{\loc}))\big)\simeq |\shv(\Hk_\bullet(\Sht_\mP^{\loc}))|.
\end{equation}
Moreover, the face and degeneracy maps of the simplicial object $\Hk_\bullet(\Sht^{\loc})$ are ind-pfp proper so the corresponding $!$-pullback functors are continuous and admit left adjoint given by corresponding $*$-pushfoward (see \Cref{locally.ind-proper.existence.lower-!}). Thus, by \cite[Proposition 5.5.7.6]{Lurie.higher.topos.theory} the category $\shv(\kot_G)$ is compactly generated and the functor $(\Nt_\mP)_*$ is the left adjoint of $(\Nt_\mP)^!$. 
\end{proof}

By \Cref{lem:local-constancy-Kot-map}, there is
a decomposition of the category of sheaves
\begin{equation}\label{eq: decomp shvkotG by connected components}
\shv(\kot_G) = \bigsqcup_{\alpha\in \pi_1(G)_{\Gamma_F}} \shv(\kot_{G,\alpha}), 
\end{equation}
where the coproduct is taken in $\lincat_\La$. For $b\in B(G)$, let
\[
i_{? b }\colon \kot_{G,? b} \rightarrow \kot_G, \quad \mbox{where } ? \in \{\emptyset, \leq, <\}
\] 
be  the corresponding locally closed immersion. From \Cref{newton.stratification}, \Cref{open.closed.recollement} we have: 

\begin{proposition}\label{prop:semi-orthogonal-explicit-fiber-sequences}
For every $b\in B(G)$ and $?\in \{\emptyset, \leq, <\}$ there are pairs of adjunctions
\begin{equation}\label{eq: push-pull-along-ib}
(i_{? b})_! \colon \shv(\kot_{G,? b}) \rightleftharpoons  \shv(\kot_G)\colon (i_{? b})^{!}, \quad 
(i_{? b})^{*} \colon \shv(\kot_G) \rightleftharpoons  \shv(\kot_{G,? b})\colon (i_{? b})_{*}.
\end{equation}
The functors $(i_{? b})_{!}$, $(i_{? b})_{*}$ are fully faithful. Moreover:
\begin{enumerate}
    \item\label{prop:semi-orthogonal-explicit-fiber-sequences-1} $(i_{\leq b})_{*} \simeq (i_{\leq b})_{!}$ and $(i_{< b})_{*} \simeq (i_{< b})_{!}$. 
    \item\label{prop:semi-orthogonal-explicit-fiber-sequences-2} $(i_{\leq b'})^{!}\circ (i_{b})_{*}\simeq 0$ (equivalently, $(i_{b})^{*} \circ (i_{\leq b'})_{!}\simeq 0$) for $b'<b$ or $\kappa(b) \neq \kappa(b')$.
    \item\label{prop:semi-orthogonal-explicit-fiber-sequences-3} \label{fiber.sequence.newton.kot} For $\mathcal{F}\in \shv(\kot_G)$, there are natural fiber sequences
    \begin{align*}
    (i_{< b})_{*}((i_{< b})^{!}\mF)
    \rightarrow (i_{\leq b})_{*}((i_{\leq b})^{!}\mF) &
    \rightarrow
    (i_{b})_{*} ((i_{b})^{!} \mF),\\
    (i_{b})_{!}((i_{b})^{*}\mF)
    \rightarrow (i_{\leq b})_{*}((i_{\leq b})^{*}\mF) &
    \rightarrow
    (i_{< b})_{*} ((i_{< b})^{*} \mF).
    \end{align*}
\end{enumerate}
\end{proposition}

We can identify each category $\shv(\kot_{G,b})$ with category of representations of $G_b(F)$. We choose a $k$-point of $\kot_{G,b}$, still denoted by $b$, to identify $\kot_{G,b}$ with $\bB_{\mathrm{proet}} G_b(F)$ as  \Cref{prop-geom-Newton-strata}. Now the following statement is a consequence of \Cref{prop: identifying shv on BG and representation of G}.
\begin{proposition}\label{shv.on.kot(G).b}
A choice of a $k$-point of $\kot_{G,b}$ induces a natural equivalence
\[
\shv(\kot_{G,b})\cong \Rep(G_b(F)). 
\]
\end{proposition} 

Thus, we may rewrite \eqref{eq: push-pull-along-ib} as pairs of adjoint functors
\[
(i_{b})_{!} \colon \rep(G_b(F)) \rightleftarrows \shv(\kot_G)\colon (i_{b})^{!}, \quad 
(i_{b})^{*} \colon \rep(G_b(F)) \rightleftarrows \shv(\kot_{G,b})\colon (i_{b})_{*}.
\]
If $b$ is basic, the inclusion $i_b$ is a closed embedding and so $(i_{b})_{*}=(i_{b})_{!}$.

\begin{proposition}
\label{compact.objects.of.B(G)}
The functor $(i_b)_!: \shv(\kot_{G,b})\to \shv(\kot_G)$ preserves compact objects.
An object $\mF$ of $\shv(\kot_G)$ is compact if and only if $(i_b)^*\mF = 0$ for all but finitely many $b\in B(G)$, and for every $b\in B(G)$ the object $(i_b)^*\mF$ is compact.
\end{proposition}
\begin{proof}
As $(i_b)_!, (i_b)^*$  are left adjoints of continuous functors, they preserve compact objects. So if $\mF$ is compact then $(i_b)^*(\mF)$ is compact for every $b\in \kot_G$.
As $\kot_{G}=\colim_{B(G)} \kot_{G,\leq b}$, where the colimit is taken over $B(G)$ with the partial order given in \Cref{SS: set BG}, \Cref{colimit.sheaves.pseudo-proper.stacks} implies that
\[
\shv(\kot_G) \simeq \colim_{B(G)} \shv(\kot_{G,\leq b}).
\]
Then every compact object $\mF$ belongs to $\shv(\cup_{i\in I}\kot_{G,\leq b_i})^\cpt$ for a finite set $I$. 
Therefore $(i_b)^*\mF=0$ for all but finitely many $b$.

Now conversely, suppose that $\mF\in\shv(\kot_G)$ satisfies the conditions in the proposition. We may assume that $\mF$ is supported on one connected component of $\kot_G$ and then assume that $\mF=(i_{\leq b_0})_{*}(\mF')$ for some $b_0\in B(G)$.
Now assume $(i_b)^*(\mF)$ is compact for every $b\in B(G)$. As $\mF$ is supported on $\kot_{G,\leq b_0}$, from the fiber sequence
\[
    (i_{b})_{!}((i_{b})^{*}\mF)
    \rightarrow (i_{\leq b})_{!}((i_{\leq b})^{*}\mF) 
    \rightarrow
    (i_{< b})_{!} ((i_{< b})^{*} \mF)
\]
it is enough to show that $(i_{< b_0})^{*} \mathcal{F}$ is compact. Continuing by induction on the finite set $b\leq b_0$ and the corresponding fiber sequences we get that $\mF$ is compact. 
\end{proof}

Later on, using a canonical self-duality on $\shv(\kot_G)$, we will prove the following parallel statement.

\begin{proposition}\label{cor: shrek-restr-reserv-cpt}
The functor $(i_b)_*: \shv(\kot_{G,b})\to \shv(\kot_G)$ preserves compact objects.
In addition, an object $\mF \in \shv(\kot_G)$ is compact if and only if $(i_b)^!\mF = 0$ for all but finitely many $b\in B(G)$, and for every $b\in B(G)$ the object $(i_b)^!\mF$ is compact.
\end{proposition}

Recall we write $j_b: \kot_{G,b}\to \kot_{G,\leq b}$ to be the quasi-compact open embedding.

\begin{corollary}\label{cor:semi-orthogonal decomp for Shv(kot(G))}
For every $b$, both the sequence
\[
\shv(\kot_{G,<b})\xrightarrow{(i_{<b})_*} \shv(\kot_{G,\leq b}) \xrightarrow{(j_b)^!} \shv(\kot_{G,b}) 
\]
and the sequence
\[
\shv(\kot_{G,b})\xrightarrow{(j_b)_!}  \shv(\kot_{G,\leq b})\xrightarrow{ (i_{<b})^*} \shv(\kot_{G,<b})
\]
induce semi-orthogonal decompositions of $\shv_{\kot_{G,\leq b}}$ (in the sense of \Cref{rem:localization sequence}).
\end{corollary}
\begin{proof}
First, the right adjoint of $(j_b)^!$ is $(j_b)_*$, which sends compact objects to compact objects by \Cref{cor: shrek-restr-reserv-cpt}. Indeed, if $\mF\in \shv(\kot_{G,b})^\cpt$, then $(j_b)_*\mF=(i_{\leq b})^*(i_{b})_*\mF$ is compact in $\shv(\kot_{G,\leq b})$. Similarly, $(i_{<b})_*$ sends compact objects to compact objects.
In addition, we see from \Cref{prop:semi-orthogonal-explicit-fiber-sequences} that the above sequences are localization sequences. Therefore, both sequences induce semi-orthogonal decompositions of $\shv(\kot_{G,\leq b})$, as desired.
\end{proof}

\begin{corollary}\label{lem: shvkot(G) tensor product}
For every prestack $X$ over $k$, the natural exterior tensor product functor 
\[
\shv(\kot_G)\otimes_\La \shv(X)\to \shv(\kot_G\times X)\] 
is an equivalence. In particular, $\shv(\kot_G)\otimes_\La \shv(\kot_G)\to \shv(\kot_G\times \kot_G)$ is an equivalence.
\end{corollary}
\begin{proof}
We give an argument and the same argument will be used several times in the sequel.
We can write $\shv(\kot_G\times X)=\colim_{b\in B(G)}\shv(\kot_{G,\leq b}\times X)$ for any $X$. As tensor product in $\lincat_\La$ commutes with colimits, it is enough to show that
$\shv(\kot_{G,\leq b})\otimes_\La \shv(X)\to \shv(\kot_{G,\leq b}\times X)$ is an equivalence, for every $b$. As $\shv(\kot_{G,\leq b})$ is compactly generated, we already know that it is fully faithful by  \Cref{lem: categorical kunneth prestack}. Therefore, it is enough to show that the essential image of the exterior tensor product functor generates $\shv(\kot_{G,\leq b}\times X)$. We have the similar localization sequence  for $\kot_{G,\leq b}\times X$, 
\[
\shv(\kot_{G,<b}\times X)\xrightarrow{(i_{<b})_*} \shv(\kot_{G,\leq b}\times X) \xrightarrow{(j_b)^!} \shv(\kot_{G,b}\times X),
\]
as in \Cref{cor:semi-orthogonal decomp for Shv(kot(G))}. Now for every $\mF\in\shv(\kot_{G,\leq b}\times X)$, we have $(i_{<b})_*(i_{<b})^!\mF\to \mF\to (j_b)_*(j_b)^!\mF$. Recall that $*$-pushforwards and $!$-pullbacks commute with exterior tensor product (as encoded by the sheaf theory, see \eqref{eq:abstract-exterioir-pullback} and \eqref{eq:abstract-Kunneth-formula}).
Therefore, by induction and by \Cref{lem: rep(G) tensor product}, both $(i_{<b})_*(i_{<b})^!\mF$ and $(j_b)_*(j_b)^!\mF$ belong to $\shv(\kot_{G,\leq b})\otimes_\La \shv(X)$. Therefore the essential image of the exterior tensor product functor generates $\shv(\kot_{G,\leq b}\times X)$, as desired.
\end{proof}

\begin{remark}\label{rem: localization sequence vs semiorthogonal decomposition}
For general prestack $X$, the localization sequnce in the proof of \Cref{lem: shvkot(G) tensor product} may not form a semi-orthogonal decomposition in the sense of \Cref{rem:localization sequence}, as $(j_b)_*$ may not admit a continuous right adjoint in general. In addition, we do not have the second localization sequence as in \Cref{cor:semi-orthogonal decomp for Shv(kot(G))}. But this is not a problem if $X$ is a quasi-compact very placid stack.
\end{remark}

We can also consider the horizontal trace of $\shv(\kot_G)$
\[
\mathrm{tr}(\shv(\kot_G))\in \End_{\lincat_\La} \mathbf{1}_{\lincat_\La}=\Mod_\La.
\] 
In addition, for a compact object $\mF\in \shv(\kot_G)$ we have its (Chern) character $\mathrm{ch}(\mF)\in H^0\mathrm{tr}(\shv(\kot_G))$ (see \Cref{rem-Chern-character}).

\begin{corollary}\label{cor: decomposition HH of shvkot}
We have
\[
\bigoplus \mathrm{tr}( (i_b)_*,\id): \bigoplus_{b\in B(G)} \mathrm{tr}(\rep(G_b(F)))\cong \mathrm{tr}(\shv(\kot_G)),
\]
and
\[
\bigoplus \mathrm{tr}( (i_b)_!,\id): \bigoplus_{b\in B(G)} \mathrm{tr}(\rep(G_b(F)))\cong \mathrm{tr}(\shv(\kot_G)).
\]
In particular, $\mathrm{tr}(\shv(\kot_G))\in \Mod_\La^{\leq 0}$. In addition, we have the isomorphism of $K$-groups 
\[
\bigoplus K_0((i_b)_*): \bigoplus_{b\in B(G)}K_0(\rep(G_b(F))^\cpt)\cong K_0(\shv(\kot_G)^\cpt).
\]
\[
\bigoplus K_0((i_b)_!): \bigoplus_{b\in B(G)}K_0(\rep(G_b(F))^\cpt)\cong K_0(\shv(\kot_G)^\cpt).
\]
The Chern character map $\mathrm{ch}: K_0(\shv(\kot_G)^\cpt)\to H^0\mathrm{tr}(\shv(\kot_G))$ is compatible with the above direct sum decompositions.
\end{corollary}
\begin{proof}This follows from \Cref{cor:semi-orthogonal decomp for Shv(kot(G))}, \Cref{cor-hochschild-semi-orthogonal}, and \Cref{rem-Chern-character}. 
More precisely, for each $\al\in \pi_1(G)_{\Ga_F}$, let $B(G)_\al$ denote the corresponding connected component. Then 
we may extend the partial order on $B(G)_\al$ to a total order so identify $B(G)_\al$ with $\bZ_{\geq 0}$ as ordered sets. 
Then we write $\shv(\kot_{G,\al})=\colim_{b\in B(G)_\al} \shv(\kot_{G,\leq b})$ as a direct limit, with the transitioning functors being $*$-pushforwards. Taking horizontal trace commutes with colimits so
\[
\mathrm{tr}(\shv(\kot_{G,\al}))=\colim_{b\in B(G)_\al} \mathrm{tr}(\shv(\kot_{G,\leq b})).
\]
By the first semi-orthogonal decomposition from \Cref{cor:semi-orthogonal decomp for Shv(kot(G))} and \Cref{cor-hochschild-semi-orthogonal}, we have 
\[
\mathrm{tr}(\shv(\kot_{G,\leq b}))=\mathrm{tr}(\shv(\kot_{G,<b}))\oplus \mathrm{tr}(\shv(\kot_{G,b}))
\] 
with inclusions of direct summands induced by $*$-pushforwards. 
Similarly, by the second semi-orthogonal decomposition from \Cref{cor:semi-orthogonal decomp for Shv(kot(G))} and \Cref{cor-hochschild-semi-orthogonal}, we have another decomposition with inclusions of direct summands induced by $!$-pushforwards. 
Therefore, $\tr(\shv(\kot_{G}))=\oplus_b\tr(\shv(\kot_{G,b}))$ as desired.

In addition, by \Cref{eq-unit-duality-Rep(G)} \eqref{eq-unit-duality-Rep(G)-1}, $\mathrm{tr}(\shv(\kot_{G,b}))\cong \mathrm{tr}(\rep(G_b(F)))$ is identified with the derived covariants of $C_c^\infty(G_b(F))$ with respect to the conjugation action of $G_b(F)$ on itself. In particular, $\mathrm{tr}(\rep(G_b(F)))\in \Mod_\La^{\leq 0}$ and so $\mathrm{tr}(\shv(\kot_G))\in \Mod_\La^{\leq 0}$.
\end{proof}

\begin{remark}
It is interesting to know whether above two decompositions of $K$-theory (using $*$-pushforwards and $!$-pushforwards) coincide. Some evidence that this might be the case is provided in \Cref{lem: chern character of ADLV}. 
We also note a closely related conjecture is made by Hansen (see \cite[Conjecture 3.4.3]{Hansen.Beijing}).
\end{remark}

Later on, we will need the following statement which directly follows from \Cref{prop-local-shtuka-uw} and the commutative diagram \eqref{eq: Cartesian diagram from Sht to Isoc for straight}. Let $w\in \widetilde{W}$ be a $\sigma$-straight element corresponding to $b$ and let $\breve\bff$ be a facet as in \Cref{lem: wsigma stable facet}. Let $\breve\mP=\breve\mP_{\breve\bff}$ be the corresponding standard parahoric.
Let $P_b:=P_{\dot{w},\breve\bff}=\{g\in \breve\mP(\breve\mO)\mid g\dot{w}\sigma(g)^{-1}=\dot{w}\}\subset G_b(F)$ as before. Let
\[
i_{\breve\mP, w}: \frac{LG_{W_{\breve\bff}w}}{\Ad_\sigma L^+\breve\mP}\to \sht^\loc_{\breve\mP}
\] 
be the locally closed embedding. 

\begin{lemma}\label{lem:Iwhaori-Hecke-for-Jb}
We have a commutative square
\[
\begin{tikzcd}
\rep(P_{b}) \arrow[r,"(i_{\breve\mP,w})_{*}"]\arrow[d,"\cind_{P_{b}}^{G_b(F)}"] & \shv(\sht_{\breve\mP}^{\loc})\arrow[d,"(\Nt_{\breve\mP})_{*}"]\\
\rep(G_{b}(F)) \arrow[r,"(i_{b})_{*}"] & \shv(\kot_G).
\end{tikzcd}
\]
and similarly with $(i_{\breve\mP,w})_{*}$ and $(i_{b})_{*}$ replaced by $(i_{\breve\mP, w})_{!}$ and $(i_{b})_{!}$.
\end{lemma}
In the sequel, we will denote by 
\begin{equation}\label{eq: vacuum module}
\delta_{P_{b}} = \cind_{P_{b}}^{G_b(F)}\Lambda
\end{equation} 
and denote by $\delta_{P_b,*}$ (resp. $\delta_{P_b,!}$) for its image in $\kot_G$ under $(i_b)_*$ (resp. $(i_{b})_{!}$).

Recall that we also have the notion of $\Mod_\La$-admissible objects as from \Cref{def:adm vs compact}. We simply call them $\La$-admissible, or just admissible objects.  Using \Cref{ex: explicit Serre and admissible}, \Cref{prop:semi-orthogonal-explicit-fiber-sequences}, \Cref{compact.objects.of.B(G)}, 
we obtain the following characterization of admissible objects in $\shv(\kot_G,\La)$. Following the notion from \Cref{sec:right-adjoint-of-upper-shriek}, we will write
$(i_b)_\flat$ for the right adjoint of $(i_b)^!$. We will also write $(i_b)\rstar$ for $((i_b)_*)^R$.

\begin{corollary}\label{cor: characterization of admissible}
An object $\mF\in \shv(\kot_G)$ is admissible if and only if for every $b$, $(i_b)^! \mF\in \rep(G_b(F))^\adm$, if and only if for every $b$, $(i_b)\rstar\mF\in \rep(G_b(F))^\adm$.
In addition, the functors $(i_b)_*, (i_b)_\flat$ preserve admissible objects.
\end{corollary}
\begin{proof}
As all $(i_b)_!, (i_b)^!, (i_b)^*, (i_b)_*$ preserve compact objects, their right adjoints $(i_b)^!, (i_b)_\flat, (i_b)_*, (i_b)\rstar$ preserve admissible objects (see \Cref{lem:basic cpt and adm} \eqref{lem:basic cpt and adm-3}).

In addition, if $(i_b)^!\mF$ is admissible for every $b$, then $\Hom((i_b)_!V, \mF)\in \Perf_\La$ for every $V\in \rep(G_b(F))^\cpt$ by \Cref{ex: explicit Serre and admissible}. The proof of \Cref{compact.objects.of.B(G)} shows that the collection $\{(i_b)_!V\mid b\in B(G), V\in \rep(G_b(F))^\cpt\}$ form a set of compact generators of $\shv(\kot_G)$. Therefore, $\mF$ is admissible, again by \Cref{ex: explicit Serre and admissible}.

The argument for $(i_b)\rstar$ is similar. 
\end{proof}

\begin{example}\label{ex: dualizing isoc adm}
The dualizing sheaf $\omega_{\kot_G}$ is an admissible object. It is not a compact object in $\shv(\kot_G)$. 
\end{example}

\begin{remark}
Note that as in \Cref{cor: characterization of admissible}, all functors $(i_{<b})_*, (i_{<b})\rstar, (i_b)^!, (i_b)_\flat$ preserve admissible objects. It follows that the sequence in \Cref{cor:semi-orthogonal decomp for Shv(kot(G))} restricts to a sequence
\[
\shv(\kot_{G,<b})^\adm\xrightarrow{(i_{<b})_*} \shv(\kot_{G,\leq b})^\adm \xrightarrow{(j_b)^!} \shv(\kot_{G,b})^\adm,
\]
which after ind-completion form a semi-orthogonal decomposition of  $\ind(\shv(\kot_{G,\leq b})^\adm)$.
\end{remark}

For an admissible object $\mF\in\shv(\kot_G)$, let 
\begin{equation}\label{eq:character of admissible in shv(Kot(G))}
\Theta_\mF: H^0\mathrm{tr}(\shv(\kot_G))=\bigoplus_b C_c^\infty(G_b(F),\La)_{G_b(F)}\to \La
\end{equation} 
be its character, as defined in \eqref{eq:abstract character}. Note that for $\mF=(i_b)_*\pi$, the map $\Theta_\mF$ vanishes on direct summand $C_c^\infty(J_{b'}(F),\La)_{J_{b'}(F)}$ for $b'>b$ but may not be trivial for direct summand $C_c^\infty(J_{b'}(F),\La)_{J_{b'}(F)}$ with $b'<b$.

\subsubsection{Canonical duality on $\kot_G$}\label{SS: coh. duality. Kot. stack}
Our next goal is to lift the canonical duality of $\shv(\kot_{G,b})\cong \rep(G_b(F))$ as discussed in \Cref{SS: coh. duality of locally profinite gorup}
to a duality of $\shv(\kot_G)$. We will fix the standard Iwahori $\mI$ and let $\sht^\loc=\sht^\loc_\mI$.

Recall that each $\Hk_n(\Sht^{\loc})$ is an ind-very placid stack, which can be written as
\[
\Hk_n(\Sht^{\loc}) \simeq \colim_{w_0,\dots,w_n} \Sht^\loc_{\leq w_1,\dots,\leq w_n},
\]
where
\[
\Sht^\loc_{\leq w_0,\dots,\leq w_n}=\frac{LG_{\leq w_0}\times^{\iw}\cdots\times^{\iw} LG_{\leq w_n}}{\Ad_\sigma \iw}
\]
is a very placid stack. 
This means that given a compatible system of generalized constant sheaves $\La^\eta$ on $\Hk_n(\sht^\loc)$ in the sense of \Cref{def.dimension.constant.sheaf.ind.placid.stack}, we obtain a Frobenius structure
\[
\rg_{\ind\fg}^\eta(\Hk_n,-): \rshv(\Hk_n(\sht^\loc))\to\Mod_\La,
\]
which induces a canonical equivalence 
\[
(\verd_{\Hk_n}^\eta)^{\mathrm{f.g.}}\colon \fgshv(\Hk_n(\sht^\loc))^{\op} \simeq \fgshv(\Hk_n(\sht^\loc))
\]
such that for every $\mF_1,\mF_2\in \fgshv(\Hk_n(\sht^\loc))$, we have a canonical isomorphism
\[
\Hom_{\fgshv(\Hk_n(\sht^\loc))}\big((\verd_{\Hk_n}^\eta)^{\fg}(\mF_1),\mF_2 \big) \simeq \rg_{\ind\fg}^\eta(\Hk_n(\sht^\loc),\mF_1\os \mF_2).
\]
We denote its ind-extension as
\[
(\verd_{\Hk_n}^\eta)^{\ind\fg}\colon \rshv(\Hk_n(\sht^\loc))^\vee \simeq \rshv(\Hk_n(\sht^\loc)).
\]
In addition, by \Cref{prop: verdier duality on shv for very placid} and \Cref{rem: verdier duality on shv for ind very placid}, the Frobenius structure $\rg_{\ind\fg}^\eta(\Hk_n,-)$ restricts, along the fully faithful embedding $\Psi^L: \shv(\Hk_n(\sht^\loc))\hookrightarrow\rshv(\Hk_n(\sht^\loc))$, to a Frobenius structure 
\[
\rg^\eta(\Hk_n(\sht^\loc),-): \shv(\Hk_n(\sht^\loc))\to \Mod_\La,
\] 
and therefore induces a self duality of $\shv(\Hk_n(\sht^\loc))$, denoted as
\[
\verd_{\Hk_n}^\eta\colon \shv(\Hk_n(\sht^\loc))^\vee \cong \shv(\Hk_n(\sht^\loc)), 
\]
which restricts to
\[
(\verd_{\Hk_n}^\eta)^\cpt\colon (\shv(\Hk_n(\sht^\loc))^\cpt)^\op \cong \shv(\Hk_n(\sht^\loc))^\cpt.
\]

We will fix a particular 
\begin{equation}\label{eq: can gen constant sheaf}
\eta=\can
\end{equation} 
on the Hecke stacks as follows. For $n=0$, we have the stack $\Sht^{\loc}$, equipped with the map $\Sht^{\loc} \rightarrow \iw\backslash LG/\iw$. The perfect ind-scheme $\Gr_\mI=LG/\iw$, which will be denoted by $\Fl$ in the sequel (to be constant with the more standard notations), is ind-finitely presented, and therefore has a canonical compatible system of generalized constant sheaves, whose value at each Schubert variety $\Fl_{\leq w}$ is  
\[
\La^{\can}_{\Fl_{\leq w}}\in \cshv(\Fl_{\leq w}),
\] 
which is defined to be the $*$-pullback of $\consdual_{\Spec k}$ along the pfp morphism $\Fl_{\leq w}\to \Spec k$. 
It uniquely descends to a generalized constant sheaf $\La^{\can}_{\La_{\iw\backslash LG_{\leq w}/\iw}}$ on $\iw\backslash \Fl_{\leq w}$, as equivariance with respect to a connected affine group action is a property rather than a structure of the sheaf. More precisely, one can first choose $n$ large enough such as the action of $\iw$ on $\Fl_{\leq w}$ factors through $\iw^n$. Then we have equivalence (via $!$-pullback) $\shv(\iw^n\backslash \Fl_{\leq w})=\shv(\iw\backslash \Fl_{\leq w})$ as $\iw^{(n)}=\ker(\iw\to \iw^n)$ is pro-unipotent. One can then first descend $\La_{\Fl_{\leq w}}$ to $\La_{\iw^n\backslash \Fl_{\leq w}}$ as usual, e.g. see \cite[Lemma A.1.2]{zhu2016introduction}, which then gives the descent to $\iw\backslash \Fl_{\leq w}$.

Consequently, we can apply the construction as in \Cref{ex-!-pullback-prosm-generalized-constant} to obtain a generalized constant sheaf $\La^\can_{\Sht}$, which is given by the compatible system $\{\La^{\can}_{\Sht_{\leq w}^{\loc}} \}$ with
\[
\La^{\can}_{\Sht_{\leq w}^{\loc}} \simeq \delta^!(\La^{\can}_{\iw\backslash LG_{\leq w^{-1}}/\iw}).
\]
(See \Cref{rem: modification direction of Shtuka} for the appearance of $w^{-1}$.)

We also need to consider
the $*$-restriction of $\La^{\can}_{\sht_{\leq w}^{\loc}}$ to $\sht_w^\loc$, which is a generalized constant sheaf $\La^{\can}_{\sht_{w}^{\loc}}$ on $\sht_w^\loc$ (see \Cref{ex-*-pullback-indfp-generalized-constant}). Alternatively, it can be obtained as the $!$-pullback of $\La^{\can}_{\iw\backslash LG_{w}/\iw}$, which in turn can be obtained from the constant sheaf on $LG_w/\iw$ via descent.

\begin{lemma}\label{duality.representations.shtuka.wb}
   Let $w_b\in\widetilde{W}$ be a $\sigma$-straight element that maps to $b$ under the map $B(\widetilde{W})\to B(G)$. Let $ \verd^{\can}_{\sht^{\loc}_{w_b}}\colon \shv(\sht^{\loc}_{w_b})^\vee\cong \shv(\sht^\loc_{w_b})$ denote the self-duality of $\shv(\sht^{\loc}_{w_b})$ induced by $\La^{\can}_{\sht_{w}^{\loc}}$. Then
   under the identification $\shv(\sht^{\loc}_{w_b}) \simeq \rep(I_b)$ coming from \Cref{prop: Sht-loc-w-straight}, we have a canonical equivalence 
   \[
   \verd^{\can}_{\sht^{\loc}_{w_b}} \simeq \verd^\coh_{I_b}[-2\langle2\rho,\nu_b\rangle](-\langle 2\rho,\nu_b\rangle),
   \]
   where $\verd^\coh_{I_b}$ is the usual contragredient duality of $\rep(I_b,\La)$ as in \Cref{lem: duality on profinite group via verdier}.
\end{lemma}
\begin{proof}
By \Cref{lem: duality on profinite group via verdier}, it will be enough to show there is a canonical equivalence 
\[
\La^\can_{\Sht_{w_b}^{\loc}} \simeq \consdual_{\bB I_b}[-2\langle 2\rho,\nu_b\rangle](-\langle 2\rho,\nu_b\rangle).
\] 
Indeed, this will imply that for every $\mF,\mG\in \shv(\sht^{\loc}_{w_b})^\cpt$, we have
\[
\Hom(\mF,\mG)=\Hom(\La^\can_{\sht_{w_b}^{\loc}},  (\verd^{\can}_{\sht^{\loc}_{w_b}})^\cpt(\mF)\os\mG)=\Hom(\consdual_{\bB I_b}, (\verd^{\can}_{\sht^{\loc}_{w_b}})^\cpt(\mF)[2\langle2\rho,\nu_b\rangle](\langle 2\rho,\nu_b\rangle)\os\mG).
\]

In fact, we claim that for any $w\in \widetilde{W}$, we have
\[
\La^\can_{\sht^\loc_w}\cong \consdual_{\sht^\loc_w}[-2\ell(w)](-\ell(w)).
\]
The lemma then follows from the fact that $\ell(w_b) = \langle2\rho,\nu_b\rangle$.

To prove the claim, note that $\La^\can_{\sht^\loc_w}$ is the $!$-pullback of $\La^\can_{\iw\backslash LG_{w^{-1}}/\iw}$. Therefore, it is enough to supply a canonical isomorphism 
\begin{equation}\label{eq: can constant vs dualizing}
\La^\can_{\iw\backslash LG_w/\iw}\cong \consdual_{\iw\backslash LG_{w^{-1}}/\iw}[-2\ell(w)](-\ell(w)).
\end{equation} 
Note that $\Fl_{w}$ is perfectly smooth of dimension $\ell(w)$, which admits a canonical deperfection as in \cite[Proposition 1.23]{zhu2017affine}. Namely, the left action of $\iw$ on $\Fl_{w}$ is transitive so $\Fl_{w}\simeq \iw/\iw \cap w\iw w^{-1}$. We may then use the canonical deperfection of $\iw$ given by the Greenberg realization to obtain a canonical deperfection of $\iw/\iw \cap w\iw w^{-1}$.
This choice of deperfection then identifies $\La_{\Fl_{w}}=\consdual_{\Fl_{w}}[-2\ell(w)](-\ell(w))$. As mentioned before, equivariance with respect to a connected affine group action is a property rather than a structure of the sheaf. Then \eqref{eq: can constant vs dualizing} follows as desired.
\end{proof}

Similarly, for each $\sht^{\loc}_{\leq w_0,\ldots,\leq w_n}$ there is a generalized constant sheaf $\La^{\can}_{\sht^{\loc}_{\leq w_0,\ldots,\leq w_n}}$, obtained by first descending the constant sheaf on $\Gr_{\leq w_0,\ldots,\leq w_n}$ to $\iw\backslash LG_{\leq w_0}\times^{\iw} LG_{\leq w_1}\times\cdots\times LG_{\leq w_n}/\iw$ and then $!$-pullback to $\sht^{\loc}_{\leq w_0,\ldots,\leq w_n}$.
Alternatively, it can be defined as the $!$-pullback of 
\[
\La^{\can}_{\iw\backslash LG_{\leq w_n^{-1}}/\iw}\boxtimes_{\La}\cdots\boxtimes_{\La} \La^{\can}_{\iw\backslash LG_{\leq w_0^{-1}}/\iw}
\] 
along the map $\sht^{\loc}_{\leq w_0,\ldots,\leq w_n}\to \iw\backslash LG_{\leq w_n^{-1}}/\iw\times \cdots\times  \iw\backslash LG_{\leq w_0^{-1}}/\iw$. The compatible system $\{\La_{\sht^{\loc}_{\leq w_0,\ldots,\leq w_n}}^{\can}\}$ then give a generalized constant sheaf $\La_{\Hk_{n}(\Sht^\loc)}$.

\begin{lemma}\label{compatibility.dimensions.iterated.shtukas}
   For every map $\alpha \colon [0] \rightarrow [n]$, the $*$-pullback of $\La^{\can}_{\Sht^{\loc}}$ along the face map $d_\al:\Hk_n(\sht^{\loc})\to \sht^{\loc}$ is canonically isomorphic to $\La^{\can}_{\Hk_n(\sht^{\loc})}$.
\end{lemma}
\begin{proof}
First note that $d_\al$ is representable pfp morphism between placid stacks, so the $*$-pullback is defined by \Cref{pfp.functors.and.bc.placid.stacks}. 

Now the maps $d_i, i=0,\ldots,n-1: \Hk_{n}(\sht^\loc)\to \Hk_{n-1}(\sht^\loc)$ from \eqref{eq: formula for face maps from iterated Sht to Sht} are the pullback of the corresponding convolution maps of the convolution affine flag varieties, and therefore the desired isomorphism between $(d_i)^*\La^\can_{\Hk_{n-1}(\sht^\loc)}\cong \La^\can_{\Hk_{n}(\sht^\loc)}$ follows from the corresponding statement for affine flag varieties and the base change isomorphism which is also in \Cref{pfp.functors.and.bc.placid.stacks}.

Finally, notice that the partial Frobenius \eqref{eq:partial-Frobenius} is the pullback of the morphism
\[
(\iw\backslash LG/\iw)^{n+1}\xrightarrow{c\circ (\sigma\times \id^n)} (\iw\backslash LG/\iw)^{n+1},
\]
where $c$ denotes the cyclic permutation of sending the first factor to the last. Therefore, again by base change, we have the canonical isomorphism 
\[
(\pFr)^*\La^\can_{\Hk_n(\sht^\loc)}\cong \La^\can_{\Hk_n(\sht^\loc)}.
\]
It follows that $(d_i)^*\La^\can_{\Hk_{n-1}(\sht^\loc)}\cong \La^\can_{\Hk_{n}(\sht^\loc)}$ also holds when $i=n$.
\end{proof}

\begin{remark}\label{rem: generalized constant sheaf on partial affine flag}
For an affine smooth integral model $\breve\mG$ of $G$ such that $L^+\breve\mG\supset \iw$, we will also have a generalized constant sheaf $\La_{\Gr_{\breve\mG}}$ on $\Gr_{\breve\mG}=LG/L^+\breve\mG$ given by the system of \emph{the} constant sheaves on (perfectly) finite type subschemes of $\Gr_{\breve\mG}$, just as in the Iwahori case. It then induces a generalized constant sheaf on $\Sht^{\loc}_{\breve\mG}$ as in the Iwahori case.
Note that the generalized constant sheaf on $\Fl$ is the $*$-pullback of the one on $\Gr_{\breve\mG}$ in the sense of \Cref{ex-*-pullback-indfp-generalized-constant}. Similarly, the generalized constant sheaf on $\Sht^\loc$ is the $*$-pullback of the one on $\Sht^{\loc}_{\breve\mG}$. It follows that the induced corresponding dualities of $\shv(\Sht^\loc)$ and $\Shv(\Sht^\loc_{\breve\mG})$ are compatible under $*$-pushforwards. 
\end{remark}

Now we can use the Verdier duality functors on the stacks $\Hk_n(\sht^\loc)$ to define a duality on $\shv(\kot_G)$.

First, the isomorphisms in
\Cref{compatibility.dimensions.iterated.shtukas} are compatible with each other in an obvious manner (no higher compatibility is needed). Therefore, they together
give rise to a simplicial functor
\begin{equation}\label{eq.global.sections.hecke.shtukas}
\rg^{\can}(\Hk_\bullet(\Sht^\loc),-)\colon \shv(\Hk_\bullet(\sht^{\loc}))\to \Mod_\La
\end{equation}
given by
\[
\rg^{\can}(\Hk_n(\sht^\loc),-)=\Hom_{\shv(\Hk_n(\sht^\loc))}(\La^{\can},-)\colon \shv(\Hk_n(\sht^\loc)) \rightarrow \Mod_\La, \quad [n]\in \Delta,
\]
and for every $\al:[m]\to [n]$ inducing the face map $d_\al: \Hk_n(\sht^\loc)\to \Hk_m(\sht^\loc)$ we have the canonical isomorphism
\[
\rg^{\can}(\Hk_n(\sht^\loc),-)\cong \rg^{\can}(\Hk_m(\sht^\loc), (d_\al)_*(-)).
\]
The system of functors \eqref{eq.global.sections.hecke.shtukas} then induce
\begin{equation}\label{eq:Frob str on shvIsoc}
\rg^{\can}(\kot_G, - )\colon \shv(\kot_G)= |\shv(\Hk_\bullet(\sht^\loc))| \rightarrow \Mod_\La.
\end{equation}

\begin{proposition}\label{cohomological.duality.kottwitz}
The functor \eqref{eq:Frob str on shvIsoc} defines a Frobenius structure on $\shv(\kot_G)$.
It induces a self duality 
\[
\verd^{\can}_{\kot_G}: \shv(\kot_G)^\vee\cong \shv(\kot_G),
\] 
which, when restricted to the anti-involution $(\verd^{\can}_{\kot_G})^\cpt: (\shv(\kot_G)^\cpt)^{\op}\cong \shv(\kot_G)^\cpt$, satifsies
\[
(\verd^{\can}_{\kot_G})^\cpt \circ \Nt_{*}\cong\Nt_*\circ (\verd^{\can}_{\sht^\loc})^\cpt.
\]
\end{proposition}
\begin{proof}
Recall that we have an equivalence in $\lincat_\La$:
\[
\shv(\kot_G) \simeq |\shv(\Hk_\bullet(\sht^\loc))|.
\]
The simplicial functor \eqref{eq.global.sections.hecke.shtukas}  induces a self duality on each of the categories $\shv(\Hk_n)$ and the boundary maps intertwine these dualities, by 
\Cref{star.induced.ind.fin.pres.} \eqref{star.induced.ind.fin.pres.-1} and ind-pfp properness of all boundary maps. 
Therefore, passing to the geometric realization we obtain a self duality $\verd^\can_{\kot_G}$ on $\shv(\kot_G)$, which satisfying $(\verd^{\can}_{\kot_G})^\cpt\circ \Nt_{*}\cong\Nt_*\circ (\verd^{\can}_{\sht^\loc})^\cpt$.

It remains to show that $\verd^\can_{\kot_G}$ can be identified with the one defined by the pairing
\begin{equation}\label{eq.pairing.on.kottwitz.cohomological}
\shv(\kot_{G})\otimes_\La \shv(\kot_{G})\xrightarrow{\os} \shv(\kot_{G})\xrightarrow{\rg^{\can}}\Mod_\La.
\end{equation}
That is, we need to show that for $\mF\in \shv(\kot_G)^\cpt$, we have a canonical isomorphism
\[
\Hom_{\shv(\kot_G)}(\mF, \mG)\cong \rg^{\can}(\kot_{G}, (\verd^{\can}_{\kot_G})^\cpt(\mF) \os \mG).
\]
We may assume that $\mF=\Nt_*\mF'$ for some $\mF'\in \shv(\sht^\loc)^\cpt$. Then by adjunction and the projection formula along the ind-pfp proper morphism $\Nt$, we have canonical isomorphisms
\begin{align*}
&\rg^{\can}(\kot_{G}, (\verd^{\can}_{\kot_G})^\cpt(\Nt_*\mF')\os \mG)=\rg^{\can}(\kot_G, \Nt_*((\verd^{\can}_{\sht^\loc})^\cpt(\mF'))\os \mG)  \\
\cong \ \ &\rg^{\can}(\sht^\loc, (\verd^{\can}_{\sht^\loc})^\cpt(\mF')\os \Nt^!\mG) \cong \Hom(\mF',\Nt^!\mG)\cong \Hom(\Nt_*\mF',\mG),
\end{align*}
as desired. 
This shows that $\eqref{eq:Frob str on shvIsoc}$ is indeed a Frobenius structure on $\shv(\kot_G)$. 
\end{proof}

Next, for every $b\in B(G)$ we define a functor 
\[
\rg^{\can}(\kot_{G,b}, - )\colon \shv(\kot_{G,b}) \rightarrow \Mod_\La, \quad \rg^{\can}(\kot_{G,b}, \mF) = \rg^{\can}(\kot_G, i_{b,*}(\mF))
\]

\begin{lemma}\label{cohomological.duality.on.strata}
For every $b\in B(G)$, under the identification of \Cref{shv.on.kot(G).b}, we have a canonical equivalence of functors
\[
\rg^{\can}(\kot_G, (i_b)_*(-))\cong \rg^\can(\bB_{\proet} G_b(F),-)[2\langle\rho,\nu_b\rangle](\langle\rho,\nu_b\rangle).
\]
In particular, $\rg^{\can}(\kot_G, (i_b)_*(-))$ is a Frobenius structure on $\shv(\kot_{G,b})$ inducing a self duality on $\shv(\kot_{G,b})$, which under the identification of \Cref{shv.on.kot(G).b}, is identified as
\[
\verd^\can_{\kot_{G,b}} \simeq \verd^{\can}_{G_b(F)}[-2\langle2\rho,\nu_b\rangle](-\langle2\rho,\nu_b\rangle).
\]
where $\verd^{\can}_{G_b(F)}$ is the duality on $\rep(G_b(F))$ from \Cref{cohomological.duality.on.locally profinite}.
\end{lemma}
\begin{proof}

We consider the \v{C}ech nerve $(\sht_{w_b}^{\loc}/\kot_{G,b})_\bullet$ of the map $\sht_{w_b}^{\loc}\to \kot_{G,b}$. We have a map of simplicial prestacks
\[
\Hk_\bullet(\bB_{\mathrm{proket}}I_b)\cong (\sht_{w_b}^{\loc}/\kot_{G,b})_\bullet\to \Hk_\bullet(\sht^\loc),
\]
with each $(\sht_{w_b}^{\loc}/\kot_{G,b})_n\to \Hk_n(\sht^\loc)$ being pfp.
By the same argument as in \Cref{duality.representations.shtuka.wb}, we see that the $*$-pullback of $\La^\can_{\Hk_\bullet(\sht^\loc)}$ to $(\sht_{w_b}^{\loc}/\kot_{G,b})_\bullet$ is just $\La^\can_{\Hk_\bullet(\bB_{\mathrm{proket}}I_b)}[-2\langle\rho,\nu_b\rangle](-\langle\rho,\nu_b\rangle)$. This gives the first statement. 
The rest statements follow from the first.
\end{proof}

\begin{proposition}\label{cohomological.duality.kottwitz.pullpush}
The functors $(i_b)_*$ and $(i_b)^!$ preserve compact objects.
We have a canonical equivalences
\[
(\verd_{\kot_G}^{\can})^\cpt\circ (i_{b})_*\cong (i_{b})_!\circ (\verd^{\can}_{G_b(F)})^\cpt[-2\langle2\rho, \nu_b\rangle](-\langle2\rho, \nu_b\rangle),
\] 
\[
(i_{b})^*\circ (\verd_{\kot_G}^{\can})^\cpt\cong (\verd^{\can}_{G_b(F)})^\cpt\circ (i_{b})^! [-2\langle2\rho, \nu_b\rangle](-\langle2\rho, \nu_b\rangle).
\]
In particular, if $b$ is basic, $\bD^{\can}_{\kot_G}$ preserves the full subcategory $\shv(\kot_{G,b})$, and restricts to the canonical duality of $G_b(F)$.
\end{proposition}

\begin{proof}
Let $\Nt_{w_b}:\Sht_{w_b}^\loc\to\kot_{G,b}$ be the restriction of $\Nt$ (see \Cref{lem-ind-finite-etale-chart-of-Newton}). First note that as in \Cref{cohomological.duality.kottwitz}, we have
\[
(\verd_{\kot_{G,b}}^\can)^\cpt\circ (\Nt_{w_b})_*\cong (\Nt_{w_b})_*\circ (\verd_{\sht^{\loc}_{w_b}}^\can)^\cpt.
\]
By  \Cref{cor: shrek pull and star push preserving compact for quotient}, we see that $(i_{w_b})_*: \shv(\sht^\loc_{w_b})\to \shv(\sht^\loc)$ preserves compact objects. Therefore, $(i_b)_*((\Nt_{w_b})_*\mF)=\Nt_*((i_{w_b})_*\mF)$ is compact for any $\mF\in  \shv(\sht^\loc_{w_b})^\cpt$.
As $\shv(\kot_{G,b})^\cpt$ is generated by $(\Nt_{w_b})_*(\sht(\sht^\loc_{w_b})^\cpt)$, we see that $(i_b)_*$ preserve compact objects.

Similarly by \Cref{cor: shrek pull and star push preserving compact for quotient}, the $!$-pullback along $\sht^\loc_b\to \sht^\loc$ preserves compact objects. The map $\sht^\loc_b\to \kot_{G,b}$ obtained by restriction of $\Nt$ is still ind-pfp proper and therefore, the $*$-pushforward along it sends compact objects to compact objects. As $\shv(\kot_G)^\cpt$ is generated by $\Nt_*(\shv(\sht^\loc))$, the base change implies that $(i_b)^!$ also preserves compact objects.

Now for $\mF'\in \shv(\sht^\loc_{w_b})^\cpt$, there are canonical isomorphisms
\begin{align*}
& (\verd_{\kot_G}^{\can})^\cpt((i_{b})_*((\Nt_{w_b})_*\mF'))\cong (\verd_{\kot_G}^{\can})^\cpt(\Nt_*((i_{w_b})_*\mF'))\cong \Nt_*((\verd_{\sht^\loc}^\can)^\cpt((i_{w_b})_*\mF'))\\
\stackrel{(\star)}{\cong} \ \ &\Nt_*((i_{w_b})_!((\verd_{\sht_{w_b}^\loc}^\can)^\cpt(\mF')))\cong (i_b)_!((\Nt_{w_b})_*((\verd_{\sht_{w_b}^\loc}^\can)^\cpt(\mF')))\cong (i_b)_!((\verd_{\kot_{G,b}}^{\can})^\cpt((\Nt_{w_b})_*\mF')),
\end{align*}
where the isomorphism labelled by $(\star)$ follows from \Cref{star.induced.ind.fin.pres.} \eqref{star.induced.ind.fin.pres.-1}.
Together with \Cref{cohomological.duality.on.strata}, this shows the first isomorphism. 

The second isomorphism formally follows from the first and the fact that $(\verd_{\kot_G}^{\can})^\cpt$ and  $(\verd_{\kot_{G,b}}^{\can})^\cpt$ are an anti-involutions (see  \eqref{eq: involutive property of duality cpt}). Namely, let $\mF\in\shv(\kot_G)^\cpt$ and $\mG\in\shv(\kot_{G,b})^\cpt$, we compute
\begin{align*}
&\Hom((i_{b})^*((\verd_{\kot_G}^{\can})^\cpt(\mF)),\mG)\cong \Hom((\verd_{\kot_G}^{\can})^\cpt(\mF), (i_b)_*\mG)\cong \Hom((\verd_{\kot_G}^{\can})^\cpt((i_b)_*\mG), \mF)\\
\cong \ \ & \Hom((i_b)_!((\verd_{\kot_{G,b}}^{\can})^\cpt(\mG)), \mF)\cong \Hom(\verd_{\kot_{G,b}}^{\can})^\cpt(\mG),(i_b)^!\mF)\cong \Hom(\verd_{\kot_{G,b}}^{\can})^\cpt((i_b)^!\mF), \mG).
\end{align*}
The desired isomorphism then follows from this and \Cref{cohomological.duality.on.strata}.
\end{proof}

Now we can give a promised proof of \Cref{cor: shrek-restr-reserv-cpt}.
\begin{proof}[Proof of \Cref{cor: shrek-restr-reserv-cpt}]
That $(i_b)_*$ and $(i_b)^!$ preserve compact objects is contained in \Cref{cohomological.duality.kottwitz.pullpush}.
Now, if $\mF$ is compact, then $\mF=(\verd^{\can}_{\kot_G})^\cpt(\mG)$ for some $\mG\in \shv(\kot_G)^\cpt$. So $(i_b)^!\mF=(\verd^{\can}_{\kot_G})^\cpt((i_b)^*\mG)$ is compact and is zero for all but finitely many $b$s, by \Cref{compact.objects.of.B(G)}. 

Next suppose that $\mF$ satisfies assumptions in the proposition. We argue as in \Cref{compact.objects.of.B(G)}. 
We may assume that that $\mF=(i_{\leq b_0})_{*}(\mF')$ for some $b_0\in B(G)$.
Now assume $(i_b)^*(\mF)$ is compact for every $b\in B(G)$. As $\mF$ is supported on $\kot_{G,\leq b_0}$, from the fiber sequence
\[
    (i_{<b})_{*}((i_{<b})^{!}\mF)
    \rightarrow (i_{\leq b})_{*}((i_{\leq b})^{!}\mF) 
    \rightarrow
    (i_{b})_{*} ((i_{b})^{!} \mF),
\]
it is enough to show that $(i_{< b_0})^{!} \mF$ is compact. Continuing by induction on the finite set $b\leq b_0$ and the corresponding fiber sequences we get that $\mF$ is compact. 
\end{proof}

\begin{remark}\label{rem: right adjoint of rgcan}
One of the corollary of the above discussions is that the functor 
\eqref{eq:Frob str on shvIsoc} sends compact objects to compact objects. Indeed, it is enough to see that $\rg^{\can}(\kot_G,(i_b)_*\cind_K^{G_b(F)}\La)$ is compact, for $K\subset G_b(F)$ pro-$p$ open compact. However, by \Cref{cohomological.duality.on.strata} and \Cref{prop: BZ-duality}, this is nothing but taking (derived) $G_b(F)$-coinvariants of $\cind_K^{G_b(F)}\La$, up to shifts, which then is just $\La$ up to shifts.

Therefore $\rg^{\can}(\kot_G,-)$ admits a continuous right adjoint. In particular,  let $\frobdual^\can$ be the object $\omega^\la$ as in \Cref{ex: duality via Frobenius-structure} associated to the Frobenius structure of $\shv(\kot_G)$ as defined in \eqref{eq:Frob str on shvIsoc}. Then $\frobdual^\can$ is admissible.

Of course, given \Cref{ex: dualizing isoc adm}, these facts also follow from \Cref{rem:unit adm}.

We expect that $\frobdual^\can=\consdual_{\kot_G}$. In fact we expect that $\rg^{\can}(\kot_G,-)$ is the left adjoint of the natural $!$-pullback along $\kot_G\to \Spec k$. However, we cannot prove this yet.
\end{remark}

We recall that $\verd^\can_{\kot_G}$ also restricts to an anti-involution 
\[
(\verd_{\kot_G}^{\can})^\adm: (\shv(\kot_G)^\adm)^{\op}\to \shv(\kot_G)^\adm,
\]
as from \eqref{eq: dual of admissible objects} and  \eqref{eq: involutive property of duality adm}. 
Recall that $(i_b)_*, (i_b)_\flat, (i_b)\rstar, (i_b)^!$ preserve admissible objects (see the proof of \Cref{cor: characterization of admissible}). 

The following statement is dual to \Cref{cohomological.duality.kottwitz.pullpush}.
\begin{proposition}\label{prop: adm.duality.kottwitz.pullpush}
We have
\[
(\verd_{\kot_G}^\can)^\adm\circ (i_b)_* \cong (i_b)_\flat \circ (\verd^\can_{G_b(F)})^{\adm}[-2\langle 2\rho,\nu_b\rangle ](-\langle 2\rho,\nu_b\rangle),
\]
and
\[
 (i_b)\rstar\circ (\verd_{\kot_G}^\can)^\adm[2\langle 2\rho,\nu_b\rangle ](\langle 2\rho,\nu_b\rangle) \cong (\verd^\can_{G_b(F)})^{\adm}\circ (i_b)^!.
\]
\end{proposition}
\begin{proof}
Let $\mF\in \shv(\kot_G)^\adm$. Then 
\[
(\verd_{\kot_G}^\can)^\adm((i_b)_*\mF)=\underline\Hom((i_b)_*\mF, \frobdual^\can)
\]
by \eqref{eq:abstract smooth dual 2}. Now for $\mG\in \shv(\kot_G)$, we have
\begin{align*}
&\Map(\mG, \underline{\Hom}((i_b)_*\mF,\frobdual^\can))\cong \Map(\mG\os (i_b)_*\mF, \frobdual^\can)\\
\cong \ \ & \Map((i_b)_*((i_b)^!\mG\os \mF),\frobdual^\can) \\
\cong \ \ & \Map(\rg^\can(\kot_G, (i_b)_*((i_b)^!\mG\os \mF)),\La)\\
\cong \ \ & \Map(\rg^\can(\bB_{\proet}G_b(F),(i_b)^!\mG\os \mF)[2\langle2\rho,\nu_b\rangle](\langle2\rho,\nu_b\rangle),\La)\\
\cong \ \ & \Map((i_b)^!\mG\os \mF, \consdual_{\bB_{\proet}G_b(F)}[-2\langle2\rho,\nu_b\rangle](-\langle2\rho,\nu_b\rangle))\\
=        \ \ & \Map(\mG, (i_b)_\flat((\verd_{G_b(F)}^{\can})^\adm(\mF)[-2\langle2\rho,\nu_b\rangle](-\langle2\rho,\nu_b\rangle))).
\end{align*}
This gives the desired first isomorphism. The second isomorphism can be proved similarly.
\end{proof}

\begin{remark}
We notice that the composed functor
$(i_{b'})^!\circ (i_b)_\flat\neq 0$ if and only if $b\leq b'$. Informally, this means that $(i_b)_\flat$ sends a sheaf on $\kot_{G,b}$ to a sheaf supported on  $\kot_{G,\geq b}$.
\end{remark}

\begin{remark}\label{rem: canonical triviality of Frobenius endo}
Recall that $\kot_G$ is in fact defined over $k_F$, and therefore admits a $q$-Frobenius endomorphism $\phi=\sigma$. Therefore, we have a functor $\phi_*: \shv(\kot_G)\to \shv(\kot_G)$. We claim that this functor is canonically isomorphic to the identify functor. Indeed, it is enough to show that $(\Nt^u)_*\circ \phi_*\cong (\Nt^u)_*$, which in turn follows from the existence of the following commutative diagram
\[
\xymatrix{
\Sht^\loc\ar^i[r]\ar@/^{15pt}/^{\phi}[rr]\ar_{\id}[dr]& \Hk(\Sht^\loc) \ar_{d_0}[d]\ar^{d_1}[r] & \Sht^{\loc}\ar^{\Nt}[d]\\
& \Sht^{\loc} \ar^{\Nt}[r] & \kot_G, 
}\]
where the map $i$ is given by $g\mapsto (g_0,g_1)=(g,1)$, in terms of notations as in \Cref{example-explicit-1-Hk}. Therefore $d_0\circ i=\id$ via $d_1\circ i=\sigma=\phi$ is the Frobenius endomorphism of $\Sht^\loc$.

Heuristically, $\shv(\kot_G)$ should be identified as the Frobenius twisted categorical trace of appropriately defined category of sheaves on $LG$. Then $\phi_*$ identified with the abstractly defined automorphism of trace category $\tr(\bfA,\phi)$ as in \Cref{rem: imposing commutativity in trace}, which is shown to be canonically isomorphic to the identity functor.
\end{remark}

\subsubsection{The category $\rshv(\kot_G)$}\label{SSS: fgshv for isoc}
Our next goal is to discuss the category of finitely generated sheaves of $\kot_G$. For this purpose, we first discuss $\fgshv(\sht^\loc)$. 

Recall the presentation \eqref{eq:restr local Sht} \eqref{eq: local Sht as lim of restr Sht} of $\Sht^{\loc}_{\leq w}$ as inverse limit of pfp algebraic stacks over $k$. (Here we only need the case $\mP=\mI$.)
Since $\cshv(LG_{\leq w})=\colim_n \cshv(\Gr_{\leq w}^{(n)})$, by descent we see that
\[
\cshv(\Sht^{\loc}_{\leq w})=\colim_{(m,n)} \cshv(\Sht^{\loc(m,n)}_{\leq w}),
\]
with the transitioning map given by $!$-pullback along $\Sht^{\loc(m',n')}_{\leq w}\to \Sht^{\loc(m,n)}_{\leq w}$, and then
\[
\fgshv(\Sht^{\loc})=\colim_{w}\cshv(\Sht^{\loc}_{\leq w}),
\]
with transitioning map given by $*$-pushforwards.

In addition, the map $\Sht^{\loc(m',n')}_{\leq w}\to \Sht^{\loc(m,n)}_{\leq w}$ is weakly coh. pro-smooth (in general not representable) of relative dimension $d=((m'-n')-(m-n))\dim G$, and therefore its $!$-pullback, shifted by $[-d]$ is perverse exact with respect to the usual (dual) perverse $t$-structure for algebraic stacks (as recalled in \Cref{ex: usual perverse t-structure on algebraic stacks}). This implies that
\[
\fgshv(\Sht^{\loc})^{\can, \geq 0}=\colim_{w}\cshv(\Sht^{\loc}_{\leq w})^{\can, \geq 0},
\]
\[
\cshv(\Sht^{\loc}_{\leq w})^{\can, \geq 0}=\colim_{m,n}\cshv(\Sht^{\loc(m,n)}_{\leq w})^{\can, \geq (n-m)\dim G},
\]
where $\cshv(\Sht^{\loc}_{\leq w})^{\can, \geq 0}$ is the coconnective part of the perverse $t$-structure on $\shv(\Sht^\loc)$ as disucssed in \Cref{sec-verdier.duality.perverse.sheaf.ind.placid.stacks}, with respect to the generalized constant sheaf of $\shv(\Sht^\loc)$ as from \eqref{eq: can gen constant sheaf}.

\begin{remark}
The above discussions imply that $\fgshv(\sht^\loc)$ with the perverse $t$-structure associated to $\eta=\can$ can be identified with the category studied in \cite{xiao2017cycles}.
\end{remark}

We also need to make use of the following algebro-geometric version of \Cref{lem:canonical resolution of triv G-module via building}, 
which is a variant of a result of Tao-Trakvin and Varshavsky \cite{TaoTravkin}. 

Let $\scrB^{\mathrm{ext}}(G,\breve F)$ be the extended Bruhat-Tits building of $G$ over $\breve F$, with barycenter subdivision.
Let $\breve \Sigma\subset \overline{\breve{\bfa}}\subset \scrB(G,\breve F)$, where $\breve{\mathbf{a}}$ is the standard alcove and $\breve \Sigma$ is a finite subcomplex of $\overline{\breve\bfa}$ that is a fundamental domain for the $G(\breve F)$-action. 
For every simplex $\breve\sigma\subset \breve \Sigma$. Let $\breve\mG_{\breve\sigma}/\breve\mO$ be the affine smooth integral model of $G_{\breve F}$, such that $\breve\mG_{\breve\sigma}(\breve\mO)$ is the stabilizer group of $\breve\sigma$ for the action of $G(\breve F)$ on $\scrB^{\ext}(G,\breve F)$.
As before, let $\frakC_{\breve \Sigma}$ be the partially ordered set of simplices in $\breve \Sigma$.

\begin{lemma}\label{eq: geometric canonical resolution for p-adic group}
There is a canonical equivalence (in $\shv(\bB LG)$)
\[
\consdual_{\bB LG}=\colim_{\frakC^{\op}_{\breve\Sigma}} \consdual_{\bB L^+\breve\mG_{\breve\sigma}}.
\]
\end{lemma}

Using it, we can prove the analogue of \Cref{prop: f.g. reps as full sub of sm}.

\begin{proposition}\label{prop: fg sheaves on kotG}
The category $\fgshv(\kot_G)$ is generated by $(\Nt_{\breve\mG})_*\mF$, for $\mF\in \fgshv(\sht^\loc_{\breve\mG})$ and $\breve\mG$ affine smooth integral model of $G$ over $\breve \mO$. In addition, the natural functor $ \fgshv(\kot_G)\to \shv(\kot_G)$ is fully faithful.
\end{proposition}

We note that implicitly, it is the version $\fgshv(\kot_G)$ that was used in the work \cite{xiao2017cycles}, i.e. the hom spaces of certain objects in $\fgshv(\kot_G)$ were computed in \emph{loc. cit.} But the proposition says that the result will not change if the hom spaces between these objects are computed in $\shv(\kot_G)$.

\begin{proof}
We follow the same strategy of the proof of  \Cref{prop: f.g. reps as full sub of sm}. 
Let $\fgshv(\kot_G)'\subset\fgshv(\kot_G)$ be the full idempotent complete stable category generated by objects $(\Nt_{\breve\mG})_*\mF$, where $\mF\in \fgshv(\sht^\loc_{\breve\mG})$, and 
$\breve\mG$ is an affine smooth (but not necessarily fiberwise connected) integral model of $G$ over $\breve \mO$ (see \Cref{rem: allowable G-shtuka} for an explanation why $\sht^\loc_{\breve\mG}$ is defined in this generality), and
\[
\Nt_{\breve\mG}:  \sht^\loc_{\breve\mG}=\frac{LG}{\Ad_\sigma L^+\breve\mG}\to \kot_G
\] 
is the corresponding Newton map.

Given \Cref{eq: geometric canonical resolution for p-adic group}, it is enough to show that the composed functor
\[
\fgshv(\kot_G)'\subset\fgshv(\kot_G)\to \shv(\kot_G)
\] 
is fully faithful. Namely, we need to show that for $\mF_i\in \fgshv(\sht^{\loc}_{\breve\mG_{i}})$, for $i=1,2$, 
\begin{equation}\label{ex: fg sheaves hom on kotG}
\Hom_{\fgshv(\kot_G)}((\Nt_{\breve\mG_{1}})_*\mF_1,(\Nt_{\breve\mG_{2}})_*\mF_2)\simeq \Hom_{\shv(\kot_G,\La)}((\Nt_{\breve\mG_{1}})_*\mF_1,(\Nt_{\breve\mG_{2}})_*\mF_2).
\end{equation}
To simplify notations, we assume $\breve\mG_{1}=\breve\mG_{2}=\mI$ but the proof of general cases is the same.

Let $d_0,d_1\colon \Hk(\sht^{\loc})\to \sht^\loc$ be as in \Cref{example-explicit-1-Hk}. In addition, we write $\Hk(\sht^{\loc})=\colim_{w_1,w_2}\sht^{\loc}_{\leq w_1,\leq w_2}$ and let $d_{0,w_1,w_2}$ and $d_{1,w_1,w_2}$ be the restriction of $d_0$ and $d_1$ to $\sht^{\loc}_{\leq w_1,\leq w_2}$. 

As in \eqref{eq: hom of fg sheaves on sifted placid} and  \eqref{eq: hom of  sheaves on sifted placid},
by ind-proper base change, the right hand side of \eqref{ex: fg sheaves hom on kotG} is computed as 
\[
\Hom_{\shv(\sht^{\loc})}(\mF_1, (d_0)_*(d_1)^!\mF_2)=\Hom_{\shv(\sht^\loc)}(\mF_1, \colim_{w_1,w_2}(d_{0,w_1,w_2})_*(d_{1,w_1,w_2})^!\mF_2),
\]
while the left hand side is computed as
\[
\colim_{w_1,w_2}\Hom_{\shv(\sht^\loc)}(\mF_1, (d_{0,w_1,w_2})_*(d_{1,w_1,w_2})^!\mF_2).
\]
As in general $\mF_1$ is not compact in $\shv(\sht^\loc)$, we need to justify why we can pull the colimit out from the hom space. Without loss of generality, we may assume that $\mF_1,\mF_2\in \cshv(\sht^{\loc}_{\leq w})$ for some $w$. 

We consider the above mentioned perverse $t$-structure on $\shv(\Sht^\loc)$.
Each object $\mE\in \fgshv(\Sht^\loc)$ belongs to $\fgshv(\Sht^\loc)^{\can,\geq N}$ for some $N$ negative enough,
and $\Hom_{\shv(\Sht^\loc)^{\can,\geq N}}(\mE,-)$ commutes with filtered colimits in $\shv(\Sht^\loc)^{\can,\geq N}$.

We also recall that affine Deligne-Lusztig varieties are finite dimensional, i.e. each irreducible component of $X_{\leq w}(b)$ is finite dimensional and there is a uniform upper bound (depending on $w,b$) of the dimensions of irreducible components. 

Now note that for each point $x\in \Sht^{\loc}$, the space $d_0^{-1}(x)\cap d_1^{-1}(\Sht^{\loc}_{\leq w})$ is exactly $X_{\leq w}(b_x)$, where $b_x\in B(G)$ is given by $\Nt(x)\in \kot_G$. As $\sht^{\loc}_{\leq w}$ is quasi-compact, the collection $\{b_x\}$ for $x\in \sht^{\loc}_{\leq w}$ is finite.
It follows that the relative dimensions of
\[
d_{0,w_1,w_2}, d_{1,w_1,w_2}\colon \sht^{\loc}_{\leq w_1,\leq w_2}\cap d_0^{-1}(\sht^\loc_{\leq w})\cap d_1^{-1}(\sht^\loc_{\leq w})\to \Sht^\loc_{\leq w}
\] 
are uniformly bounded independent of $w_1$ and $w_2$. 
Therefore, there is some negative integer $N$ such that $\mF_1\in \fgshv(\sht^\loc)^{\can,\geq N}$
and such that $(d_{0,w_1,w_2})_*(d_{1,w_1,w_2})^!\mF_2\in \fgshv(\sht^\loc)^{\geq N}$ for all $w_1,w_2$. It follows that the map
\begin{multline*}
\colim_{w_1,w_2}\Hom_{\shv(\sht^\loc)}(\mF_1, (d_{0,w_1,w_2})_*(d_{1,w_1,w_2})^!\mF_2) \\
\to  \Hom_{\shv(\sht^\loc)}(\mF_1, \colim_{w_1,w_2}(d_{0,w_1,w_2})_*(d_{1,w_1,w_2})^!\mF_2)
\end{multline*}
is an isomorphism, as desired.
\end{proof}

\begin{remark}\label{rem: faihtfully renornalized version of isoc} 
Instead of the $\sigma$-conjugation action of $LG$ on itself, it is also important to consider the usual conjugation action of $LG$ on itself and form the quotient stack $\frac{LG}{\Ad LG}$. 
However, unlike affine Deligne-Lusztig varieties which are always finite dimensional, affine Springer fibers are usually infinite dimensional. Therefore, the tautological functor $\fgshv(\frac{LG}{\Ad LG})\to \shv(\frac{LG}{\Ad LG})$ is not fully faithful.
We refer to \cite{HHZ} for more discussions.
\end{remark}

\begin{corollary}
For every $b$, the natural functor $\fgshv(\kot_{G,\leq b})\to \shv(\kot_{G,\leq b})$ is fully faithful.
\end{corollary}
\begin{proof}
We note that the functor from the category of finitely generated sheaves to all sheaves intertwines $(i_{\leq b})_*$ and $(i_{\leq b})_*^{\ind\fg}$, both of which are and fully faithful embedding. It then follows from \Cref{prop: fg sheaves on kotG} that $\fgshv(\kot_{G,\leq b})\to \shv(\kot_{G,\leq b})$ is fully faithful.
\end{proof}

\begin{proposition}\label{prop: functors between fg sheaves}
For every $b\in B(G)$ and $?\in \{\emptyset, \leq, <\}$, the pairs of adjunctions in \Cref{prop:semi-orthogonal-explicit-fiber-sequences} restrict to pairs of adjunctions
\begin{equation}\label{eq: push-pull-along-ib}
(i_{? b})_! \colon \fgshv(\kot_{G,? b}) \rightleftarrows \fgshv(\kot_G)\colon (i_{? b})^{!}, \quad 
(i_{? b})^{*} \colon \fgshv(\kot_G) \rightleftarrows \fgshv(\kot_{G,? b})\colon (i_{? b})_{*}.
\end{equation}
The functors $(i_{? b})_{!}$, $(i_{? b})_{*}$ are fully faithful.
\end{proposition}
\begin{proof}
First, as $i_{?b}$ is pfp, $(i_{?b})_*:\shv(\kot_{G,{?b}})\to \shv(\kot_G)$  and $(i_{?b})^!:\shv(\kot_G)\to\shv(\kot_{G,?b})$ preserve the subcategory of finitely generated sheaves.

Next we show that $(i_b)_!$ preserves finitely generated sheaves. It is enough to show that for a maximal open compact subgroup $K\subset G_b(F)$ and a representation $V\in \crep(K)$, the object $(i_b)_!\cind_{K}^{G_b(F)} V$ belongs to $\fgshv(\kot_G)$. Recall that every maximal open compact subgroup of $G_b(F)$ is of the for $K_v=\mG_{b,v}(\mO)$, where $v$ is a point in the building $\scrB(G_b, F)$ and $\mG_{b,v}$ is the corresponding stabilizer group scheme (which is an integral model of $G_b$). 

By the classification of $\sigma$-conjugacy classes in $\widetilde{W}$ as discussed at the end of \Cref{sec:sigma-straight-element}, we see that we may find a $\sigma$-straight element $w$ and a standard facet $\breve\bff$ in $\scrA(G_{\breve F},S_{\breve F})$ (i.e. a facet containing $\overline{\breve\bfa}$) so that $w$ is of minimal length in $W_{\breve\bff}w$, $w\sigma(W_{\breve\bff})w^{-1}=W_{\breve\bff}$, such that there is an isomorphism
$\scrB(G_b, \breve F)\cong \scrB(\breve M_w, \breve F)$ and such that under this isomorphism $v$ corresponds to point (still denoted by $v$) on $\breve\bff_{\breve M_w}$, where $\bff_{\breve M_w}$ is the facet determined by $\breve\bff$. We then lift $v$ to a point $v'\in \breve\bff$.

Now by \Cref{eq: non-connected-local-shtuka-uw}, there is an integral model $\breve\mG_{v'}$ of $G$ over $\breve \mO$, and a pfp locally closed embedding $\bB_{\mathrm{profet}} K_v\to \Sht^\loc_{\breve\mG_{v'}}$ such that the following diagram is commutative
\[
\xymatrix{
\bB_{\mathrm{profet}} K_v\ar[r]\ar[d] & \Sht^\loc_{\breve\mG_{v'}}\ar[d]\\
\bB_{\proet} G_b(F) \ar[r] & \kot_G.
}\] 
Therefore, we may regard $V$ as an object in $\fgshv(\bB_{\mathrm{profet}} K_v)$. Its $!$-pushforward to $\sht^\loc_{\breve\mG_{v'}}$ is still a finitely generated sheaf. As $\Nt_*$ preserves finitely generated sheaves, we see that $(i_b)_!\cind_{P_b}^{G_b(F)} V\in \fgshv(\kot_G)$.

Finally, we prove that $(i_{?b})^*$ preserves finitely generated sheaves. By \Cref{prop: fg sheaves on kotG}, it is enough to show that $(i_{?b})^*((\Nt_{\breve\mG})_*\mF)$ is finitely generated for $\mF\in \fgshv(\Sht^{\loc}_{\breve\mG})$, and $\breve\mG$ is maximal. But in this case $\Nt_{\breve\mG}$ is ind-proper so $*$-pushforwards commute with $*$-pullbacks.

The rest of the claims follow easily.
\end{proof}

\begin{remark}
We give a more explicit explanation that why $(i_b)^!$ preserves finitely generated sheaves.
For this, it is enough show that $(i_b)^!(\Nt_{\breve\mG})_*\mF\in \fgshv(\rep(G_b(F))$ for every $\mF\in\fgshv(\sht^\loc_{\breve\mG})$, where $\breve\mG$ is an affine smooth integral model of $G$ over $\breve \mO$ such that $L^+\breve\mG$ contains $\iw_k$. To simplify notations, we assume that $\breve\mG=\mI_{\breve\mO}$.
We may assume that $\mF\in \fgshv(\sht_{\leq w}^\loc)$ for some $w$.
By base change, it is enough to show that $C^*(X_{\leq w}(b),\mF')\in \fgrep(G_b(F))$, where $\mF'$ is the $!$-pullback of $\mF$ along $X_{\leq w}(b)\to \sht^{\loc}_{\breve\mG,\leq w}$.  It is well-known that $X_{\leq w}(b)$ admits a finite partition $X_{\leq w}(b)=\sqcup_\al  X_{\leq w}(b)_\al$ into $G_b(F)$-stable locally closed pieces such that
\begin{itemize}
\item the index set $\{\al\}$ is finite;
\item there is a $G_b(F)$-equivariant isomorphism $X_{\leq w}(b)_\al=G_b(F)\times^K X_{\leq w}(b)_\al^0$, where $K\subset G_b(F)$ is some open compact subgroup, $X_{\leq w}(b)_\al^0$ pfp over $k$ on which $K$ acts through a finite quotient group.
\end{itemize}
It follows that $C^*(X_{\leq w}(b),\mF')\in \fgrep(G_b(F))$. 
\end{remark}

The following result characterizing finitely generated objects is parallel to  \Cref{cor: shrek-restr-reserv-cpt}.

\begin{proposition}\label{cor: shrek-restr-reserv-fg}
An object $\mF$ of $\shv(\kot_G)$ belongs to $\fgshv(\kot_G)$ if and only if $(i_b)^!(\mF) = 0$ for all but finitely many $b\in B(G)$, and for every $b\in B(G)$ the object $(i_b)^!(\mF)$ belongs to $\fgshv(\kot_{G,b})$, if and only if $(i_b)^*(\mF) = 0$ for all but finitely many $b\in B(G)$, and for every $b\in B(G)$ the object $(i_b)^*(\mF)$ belongs to $\fgshv(\kot_{G,b})=\fgrep(G_b(F))$.

The decomposiiton \eqref{eq: decomp shvkotG by connected components} induces a decomposition.
\[
 \fgshv(\kot_G) = \bigoplus_{\alpha\in \pi_1(G)_{\Gamma_F}} \fgshv(\kot_{G,\alpha}).
\]
\end{proposition}
\begin{proof}
Clearly for $X\to \kot_G$ with $X$ quasi-compact placid, the image $|X|\to |\kot_G|$ is the union of finitely many $b$s. This implies that for $\mF\in \fgshv(\kot_G)$, $(i_b)^!\mF=0$ for all but finitely many $b$'s. In addition, we have just explained that $(i_b)^!$ preserves finitely generated sheaves. This gives the ``only if" direction.
The argument as in \Cref{cor: shrek-restr-reserv-cpt} gives the ``if" direction. The statement involves $*$-pullbacks is proved similarly. Finally,
the last statement is also clear.
\end{proof}

By comparing  \Cref{cor: shrek-restr-reserv-cpt} and \Cref{cor: shrek-restr-reserv-fg}, and by \Cref{cor: char zero fg=cpt}, we obtain the following.
\begin{corollary}\label{cor: char zero fg=cpt kot}
If $\La$ is a field of characteristic zero, then $\shv(\kot_G)^\cpt=\fgshv(\kot_G)$.
\end{corollary}

Now let $\mP$ be a standard parahoric of $G$ over $\mO$.
As each $\Hk_n(\Sht_\mP^{\loc})$ is an ind-placid stack, there is the subcategory of finitely generated sheaves
\[
\fgshv(\Hk_n(\Sht_\mP^{\loc}))\subset \shv(\Hk_n(\Sht_\mP^{\loc})),
\]
and the simplicial category $\shv(\Hk_\bullet(\Sht_\mP^{\loc}))$ restricts to a simplicial category $\fgshv(\Hk_\bullet(\Sht_\mP^{\loc}))$. We let
\[
\fgshv(\kot_G)_\mP:=|\fgshv(\Hk_\bullet(\Sht^{\loc}_\mP))|\in \catid_\La.
\]
Tautologically, there is a functor
\[
\fgshv(\kot_G)_\mP\to \fgshv(\kot_G),
\]
which is fully faithful by \Cref{prop: fully faithful for ind-proper surjective fg sheaves}.

\begin{corollary}\label{prop: faihtfully renornalized version of isoc}
The composed functor $\fgshv(\kot_G)_\mP\to\fgshv(\kot_G)\to \shv(\kot_G)$ is fully faithful.  We have $\shv(\kot_G)^\cpt\subset \fgshv(\kot_G)_\mP$.
 If $\La$ is a field of characteristic zero, then we have the equivalence $\shv(\kot_G)^\cpt=\fgshv(\kot_G)_\mP=\fgshv(\kot_G)$
  (in particular $\fgshv(\kot_G)_\mP$ is independent of the choice of $\mP$).
\end{corollary}

\begin{proof}
Fully faithfulness follows from \Cref{prop: fg sheaves on kotG}.
The second statement follows from the fact that we have the inclusions $\shv(\hk_\bullet(\sht_\mP^\loc))^\cpt\subset \fgshv(\hk_\bullet(\sht_\mP^\loc))$, and $|\shv(\hk_\bullet(\sht_\mP^\loc))^\cpt|\cong \shv(\kot_G)^\cpt$ by \Cref{prop: shv on kot via sht}. The last statement follows from \Cref{cor: char zero fg=cpt kot}.
\end{proof}

\begin{remark}\label{rem: Shv on Kot as colocalization}
Recall we for quasi-compact sind-placid stack $X$, we always have a functor $\Psi: \rshv(X)\to \shv(X)$, see \eqref{eq-ind-constr-to-shv}, \eqref{eq-ind-constr-to-shv-sifted}.
By \Cref{compact.generation.admissible.stacks} when $X$ is a quasi-compact very placid stack, $\Psi$ admits a left adjoint $\Psi^L$
realizing $\shv(X,\La)$ as a colocalization of $\rshv(X,\La)$. The same statement extends to quasi-compact ind-very placid stacks, but fails in general for quasi-compact sind-very placid stack. However,  we do have a pair of adjoint functors
\[
\Psi^L: \shv(X)\rightleftharpoons \rshv(X): \Psi,\quad X=\kot_G,\ \kot_{G,\leq b},\ \kot_{G,b}
\]
such that $\Psi\circ \Psi^L\cong \id$.
\end{remark}

We also have the following statement, whose proof is parallel to \Cref{cor: colim presentation of fgshv}.

\begin{proposition}\label{cor: colim presentation of fgshv of kotG}
There is an equivalence 
\[
\colim_{\breve\frakC_\Sigma} \fgshv(\sht^{\loc}_{\breve\mG_\sigma})\cong \fgshv(\kot_G).
\] 
\end{proposition}

Note that \Cref{duality.representations.shtuka.wb} clearly admits a version for finitely generated objects, and \Cref{rem: generalized constant sheaf on partial affine flag} also admits a version for finitely generated sheaves.
Now the following statement is proved as \Cref{rem: coh duality on f.g. representations}.

\begin{corollary}\label{rem: can duality on f.g. sheaves}
There is a canonical duality
\[
(\verd^\can_{\kot_G})^{\mathrm{f.g.}}\colon \fgshv(\kot_G)^{\op}\cong \fgshv(\kot_G)
\]
satisfying
\[
(\verd^\can_{\kot_G})^{\mathrm{f.g.}}\circ (\Nt_\mP)_*\cong (\Nt_\mP)_*\circ (\verd^\can_{\sht^\loc_\mP})^{\mathrm{f.g.}}
\]
and restrict to $(\verd^\can_{\kot_G})^\cpt$.

In addition, we have a canonical equivalences
\[
(\verd_{\kot_G}^{\can})^\fg\circ (i_{b})_*\cong (i_{b})_!\circ (\verd^{\can}_{G_b(F)})^\fg[-2\langle2\rho, \nu_b\rangle](-\langle2\rho, \nu_b\rangle),\] 
\[
(i_{b})^*\circ (\verd_{\kot_G}^{\can})^\fg\cong (\verd^{\can}_{G_b(F)})^\fg\circ (i_{b})^! [-2\langle2\rho, \nu_b\rangle](-\langle2\rho, \nu_b\rangle).
\]
\end{corollary}

\begin{remark}\label{rem: canonical triviality of Frobenius endo-2}
\Cref{rem: canonical triviality of Frobenius endo} continues to hold for $\rshv(\kot_G)$.
\end{remark}

\subsubsection{$t$-structure}\label{SSS: t-structure on shvkot}
Let $\La$ be a $\bZ_\ell$-algebra as in \Cref{sec:adic-formalism}.
We further assume that $\La$ is regular noetherian, and will discuss some natural $t$-structures on $\shv(\kot_G,\La)$. As before, we omit $\La$ for notations. We will also choose, for each $b\in B(G)$, a $k$-point of $\kot_{G,b}$ (still denoted by $b$), to identify $\kot_{G,b}$ with $\bB_{\proet} G_b(F)$ as before.

For each $b\in B(G)$, let 
$(\rep(G_b(F))^{\leq 0}$ be the connective part of the standard $t$-structure on $\rep(G_b(F))$. Note that $(\rep(G_b(F))^{\leq 0}$ is closed under all small colimits and extensions. 

\begin{lemma}
The standard $t$-structure on $\rep(G_b(F))$ restricts to a $t$-structure of $\rep(G_b)^{\adm}$. If $\La$ is a field of characteristic zero, it also restricts to a $t$-structure of $\rep(G_b)^\cpt$.
\end{lemma}
\begin{proof}
Let $\pi\in \rep(G_b)^{\adm}$, which fits into a cofiber sequence $\pi'\to \pi \to \pi''$ with $\pi'\in \rep(G_b)^{\leq 0}$ and $\pi''\in \rep(G_b)^{>0}$. We need to show that $\pi'$ and $\pi''$ are admissible. 
For every pro-$p$-open compact subgroup $K\subset G_b(F)$, we have a cofiber sequence ${\pi'}^K\to \pi^K\to {\pi''}^K$ with $\pi^K\in \Perf_\La$, ${\pi'}^K\in \Mod_\La^{\leq 0}$ and ${\pi''}^K\in \Mod_\La^{>0}$. As $\La$ is regular noetherian, we see that ${\pi'}^K, {\pi''}^K\in \Perf_\La$. Therefore, $\pi',\pi''\in \rep(G_b(F))^{\adm}$.

The case $\rep(G_b(F))^\cpt$ is classical, as in this case $\rep(G_b(F))^{\heartsuit}$ has finite cohomological dimension. 
\end{proof}

\begin{remark}
If $\La$ is a field of characteristic zero, then $\rep(G_b(F))^{\adm,\heartsuit}$ is the usual abelian category of admissible smooth representations of $G_b(F)$ while $\rep(G_b(F))^{\cpt,\heartsuit}$ is the abelian category of finitely generated smooth $G_b(F)$-representations.
\end{remark}

Fix $b\in B(G)$. Recall from \Cref{cor:semi-orthogonal decomp for Shv(kot(G))} we have adjoint functors
    \begin{equation}\label{eq: t-structure-open-closed gluing-1}
    \xymatrix{
    \shv(\kot_{G,b})\ar@/^/[rr]^{(j_b)_!}\ar@/_/[rr]_{(j_b)_*} && \ar[ll]|{(j_b)^!}\shv(\kot_{G,\leq b}) \ar@/^/[rr]^{(i_{<b})^*}\ar@/_/[rr]_{(i_{<b})^!} && \ar[ll]|{(i_{<b})_*}\shv(\kot_{G,<b}),
    }
    \end{equation}
inducing two semi-orthogonal decompositions of $\shv(\kot_{G,\leq b})$. In particular, all the involved functors preserve compact objects. By \Cref{prop: functors between fg sheaves}, all the involved functors also preserve the subcategories of finitely generated sheaves.

Now the standard results on gluing $t$-structures as in \cite[Theorem 1.4.10]{BBDG} give the following.

\begin{proposition}\label{rem: perverse t-structre on fgshvekotG}
Let $\chi$ be a weight of $G$ such that $\langle\chi,\nu_b\rangle\in\bZ$ for every $\nu_b$. For each $\delta\in B(G)$, the pair of subcategories of $\shv(\kot_{G,\leq \delta})$
\[
\shv(\kot_{G,\leq \delta})^{\chi\mbox{-}p,\leq 0}=\bigl\{\mF\in \shv(\kot_{G,\leq  \delta}) \mid (i_{b})^*\mF\in \rep(G_{b}(F))^{\leq \langle\chi,\nu_{b}\rangle}, \ \forall\ b\leq  \delta \bigr\},
\]
\[ 
\shv(\kot_{G,\leq  \delta})^{\chi\mbox{-}p,\geq 0}=\bigl\{\mF\in \shv(\kot_{G,\leq  \delta}) \mid (i_{b})^! \mF\in \rep(G_{b}(F))^{\geq \langle\chi,\nu_{b}\rangle}, \ \forall\ b\leq  \delta \bigr\},
\]
defines a $t$-structure on $\shv(\kot_{G,\leq \delta})$. We similarly have a $t$-structure on $\rshv(\kot_{G,\leq \delta})$. The functor $\Psi: \rshv(\kot_{G,\leq \delta})\to \shv(\kot_{G,\leq \delta})$ is $t$-exact, and restricts to an equivalence 
\[
\rshv(\kot_{G,\leq \delta})^{\chi\mbox{-}p,+}\cong \shv(\kot_{G,\leq \delta})^{\chi\mbox{-}p,+}.
\]
\end{proposition}
\begin{proof}
We only prove the last statement. We first show that for $\mF_1\in\rshv(\kot_{G,\leq \delta}), \mF_2\in \rshv(\kot_{G,\leq \delta})^{\chi\mbox{-}p,+}$, we have
\[
\Hom_{\rshv(\kot_{G,\leq\delta})}(\mF_1,\mF_2)\cong \Hom_{\shv(\kot_{G,\leq\delta})}(\Psi(\mF_1),\Psi(\mF_2)).
\]
Note that $\mF_1$ admits a finite filtrations with associated graded being $(i_b)^{\ind\fg}_!((i_b)^{\ind\fg,*}\mF_1)$ while $\mF_2$ admits finite filtrations with associated graded being $(i_b)^{\ind\fg}_*((i_b)^{\ind\fg,!}\mF_2)$. Therefore we may assume that $\mF_1=(i_{b_1})_!^{\ind\fg}\pi_1$ for $\pi_1 \in \rrep(G_{b_1}(F))$ and $\mF_2=(i_{b_2})_*^{\ind\fg}\pi_2$ for $\pi_2\in \rrep(G_{b_2}(F))^+$. Note that in the case, the hom spaces in question are zero unless $b_1=b_2$. In the later situation, 
 the desired isomorphism follows from
\[
\Hom_{\rrep(G_b(F))}(\pi_1, \pi_2)=\Hom_{\rep(G_b(F))}(\Psi_{G_b(F)}(\pi_1), \Psi_{G_b(F)}(\pi_2)).
\]
This implies that $\rshv(\kot_{G,\leq \delta})^{\chi\mbox{-}p,+}\cong \shv(\kot_{G,\leq \delta})^{\chi\mbox{-}p,+}$ is fully faithful. Now essential surjectivity follows as objects in $\shv(\kot_{G,\leq \delta})^{\chi\mbox{-}p,+}$ is a finite extension of objects of the for $(i_b)_*\pi$ for $\pi\in \rep(G_b(F))^+$, each of which belongs to the essential image of $\Psi(\rshv(\kot_{G,\leq \delta})^{\chi\mbox{-}p,+})$.
\end{proof}

We can pass to the limit to describe a $t$-structure on $\shv$. (We omit the discussion for $\rshv$.)
\begin{proposition}\label{prop: t-structures on llc-1}
Let $\chi$ be a weight of $G$ such that $\langle\chi,\nu_b\rangle\in\bZ$ for every $\nu_b$.
Let $\shv(\kot_G)^{\chi\mbox{-}p,\leq 0}\subset \shv(\kot_G)$ be the full subcategory generated under small colimits and extensions by objects of the form 
\begin{equation}\label{eq: connective part generator}
(i_b)_!\cind_K^{G_b(F)}\La[n-\langle\chi,\nu_b\rangle], \quad b\in B(G),\ n\geq 0, \ K\subset G_b(F) \mbox{ prop-}p \mbox{ open compact}.
\end{equation} 
Then $\shv(\kot_G)^{\chi\mbox{-}p,\leq 0}$ form a connective part of an admissible $t$-structure on $\shv(\kot_G)$. The coconnective part can be described as
\[
\shv(\kot_G)^{\chi\mbox{-}p,\geq 0}=\bigl\{\mF\in\shv(\kot_G)\mid (i_b)^!\mF\in\rep(G_b(F))^{\geq \langle \chi,\nu_b\rangle}\bigr\}.
\]
In addition, if $\La$ is a field of characteristic zero, this $t$-structure restricts to a bounded $t$-structure of $\shv(\kot_G)^\cpt$, whose connective can be described as
\begin{equation}\label{eq: connective part compact}
\shv(\kot_G)^{\chi\mbox{-}p,\leq 0}\cap \shv(\kot_G)^\cpt=\bigl\{\mF\in\shv(\kot_G)^\cpt\mid (i_b)^*\mF\in\rep(G_b(F))^{\leq \langle \chi,\nu_b\rangle}\bigr\}.
\end{equation}
\end{proposition}
\begin{proof}
 That $\shv(\kot_G)^{\chi\mbox{-}p,\leq 0}$ is the connective part of an accessible $t$-structure of $\shv(\kot_G)$ follows directly from \cite[Proposition 1.4.4.11]{Lurie.higher.algebra}. The description of coconnective part follows directly from the fact that $\rep(G_b(F))^{\leq 0}$ is generated by $\cind_K^{G_b(F)}\La[n]$ for $K\subset G_b(F)$ pro-$p$ open compact and $n\geq 0$. 
We also notice that $\shv(\kot_G)^{\chi\mbox{-}p,\leq 0}\cap \shv(\kot_G)^\cpt$ is the full subcategory of $\shv(\kot_G)^\cpt$ generated by objects of the form \eqref{eq: connective part generator} under extensions, finite colimits and idempotent completions. We need to identify it with the one in \eqref{eq: connective part compact} and show that it defines a $t$-structure on $\shv(\kot_G)^\cpt$.

By induction on $ \delta$, it is also easy to see that $\shv(\kot_{G,\leq  \delta})^{\chi\mbox{-}p,\leq 0}$ as in \Cref{rem: perverse t-structre on fgshvekotG} is the full subcategory generated under small colimits and extensions by objects of the form \eqref{eq: connective part generator}, except we only allow those $b\in B(G)$ that is less than or equal to $\delta$. 
As all functors in \eqref{eq: t-structure-open-closed gluing-1} preserve compact objects, this $t$-structure restricts to a $t$-structure on $\shv(\kot_{G,\leq \delta})^\cpt$ when $\La$ is a field of characteristic zero. In addition, $\shv(\kot_{G,\leq \delta})^{\chi\mbox{-}p,\leq 0} \cap  \shv(\kot_{G,\leq \delta})^\cpt$ is the full subcategory of $\shv(\kot_{G,\leq \delta})^\cpt$ generated by objects of the form \eqref{eq: connective part generator} for $b\leq \delta$ under extensions, finite colimits and idempotent completions.

Note that for $\delta\leq \delta'$, the inclusion $\shv(\kot_{G,\leq \delta})\subset \shv(\kot_{G,\leq \delta'})$ induced by $*$-extension is $t$-exact. This implies that 
\[
\shv(\kot_G)^{\chi\mbox{-}p,\leq 0}\cap \shv(\kot_G)^\cpt=\bigcup_{\delta} (\shv(\kot_{G,\leq \delta})^{\chi\mbox{-}p,\leq 0} \cap  \shv(\kot_{G,\leq \delta})^\cpt)
\] 
is the connective part of  a $t$-structure on $\shv(\kot_G)^\cpt=\colim_{B(G)}\shv(\kot_{G,\leq b})^\cpt$, as desired.
\end{proof}

\begin{remark}
It is interesting to know whether $\shv(\kot_{G})^{\chi\mbox{-}p,\leq 0}$ can be identified with the full subcategory of $\shv(\kot_G)$ consisting of $\{\mF\in\shv(\kot_G)\mid (i_b)^*\mF\in\rep(G_b(F))^{\leq \langle \chi,\nu_b\rangle}\}$. This will be the case if the latter category is compactly generated. But we are not able to prove this.

One the other hand, by definition $\shv(\kot_{G})^{\chi\mbox{-}p,\leq 0}$ is compactly generated. In fact, by virtue of \cite[Remark C.6.1.2]{Lurie.SAG} it is a Grothendieck prestable category in the sense of \cite[\textsection{C.1.4}]{Lurie.SAG}, and the heart
$\shv(\kot_G)^{\chi\mbox{-}p,\heartsuit}$ is a Grothendieck abelian category (\cite[Remark 1.3.5.23]{Lurie.higher.algebra}). 
\end{remark}

\begin{definition}
We call the $t$-structure on $\shv(\kot_G)$ defined above the $\chi$-perverse $t$-structure on $\shv(\kot_G)$.  Similarly, we have the perverse $t$-structure of $\shv(\kot_{G,\leq b})$.
\end{definition}

We give a class of examples of perverse sheaves on $\kot_G$ (which play important roles in \cite{xiao2017cycles}).

Suppose $\mP$ is a hyperspecial parahoric group scheme of $G$ (so in particular $G$ is unramified). We consider
\[
L^+\mP\backslash LG/L^+\mP \xleftarrow{\delta_\mP}  \frac{LG}{\Ad_\sigma L^+\mP}=\sht^{\loc}_\mP\xrightarrow{\Nt_\mP} \frac{LG}{\Ad_\sigma LG}=\kot_G.
\]
We endow $\fgshv(L^+\mP\backslash LG/L^+\mP)$ with the perverse $t$-structure induced by the generalized constant sheaf $\La^{\can}$. Its heart $\fgshv(L^+\mP\backslash LG/L^+\mP)^{\heartsuit}$ is the usual Satake category.
\begin{proposition}
Assume that $\La$ is a field of characteristic zero.
Then the functor 
\[
(\Nt_\mP)_*\circ (\delta_\mP)^!: \fgshv(L^+\mP\backslash LG/L^+\mP)\to \shv(\kot_G)
\] 
sends $\shv(L^+\mP\backslash LG/L^+\mP)^{\can,\heartsuit}$ to $\shv(\kot_G)^{2\rho\mbox{-}p,\heartsuit}$.
\end{proposition}
The geometric reason behind this proposition is the dimension formula of affine Deligne-Lusztig varieties in the affine Grassmannians (e.g. see \cite[\textsection{3}]{zhu2017affine}). 
Informally, this dimension formula says that that the Newton map $\Nt_\mP:\sht^\loc_\mP\to \kot_G$
should be a ``stratified semi-small map". 
\begin{proof}
Let $\mF\in \fgshv(L^+\mP\backslash LG/L^+\mP)^{\heartsuit}$. We need to show that 
\[
(i_b)^*((\Nt_\mP)_*((\delta_\mP)^!\mF))\in \rep(G_b(F))^{\leq \langle 2\rho,\nu_b\rangle},\quad (i_b)^!((\Nt_\mP)_*((\delta_\mP)^!\mF))\in \rep(G_b(F))^{\geq \langle2\rho,\nu_b\rangle}.
\]

We may assume that $\mF$ is a constructible perverse sheaf supported on a spherical Schubert variety $\Gr_{\mP,\leq \mu}$, where $\mu\in W_\mP\backslash \widetilde W/W_\mP$ can be represented by a dominant coweight  $\mu$ of $G$ as usual. 
Using \Cref{cohomological.duality.kottwitz.pullpush}, it is enough to prove the second estimate.
We write $\leq $ for the usual Bruhat order on the set of dominant coweights.
Consider
\[
\xymatrix{
X_{\leq \mu}(b)\ar[d]\ar[r]& \sht^\loc_{\mP,\leq\mu}\ar[r]\ar[d] &  L^+\mP\backslash LG_{\leq \mu}/L^+\mP\\
\Spec k \ar^b[r] & \kot_G &
}
\]
By base change, the underlying $\La$-module for $(i_b)^!((\Nt_\mP)_*((\delta_\mP)^!\mF))\in \rep(G_b(F))$ is identified with 
\begin{equation}\label{eq: coh of ADLV}
C^\bullet(X_{\leq \mu}(b),\mF')=\colim_Z C^\bullet(Z, (i_Z)^!\mF'),
\end{equation}
where $\mF'$ is the $!$-pullback of $\mF$ along the first row to $X_{\leq \mu}(b)$, $i_Z: Z\to X_{\leq \mu}(b)$ range over pfp closed subschemes of $X_{\leq \mu}(b)$, and $C^\bullet(Z,-)$ is the usual cohomology of the sheaf $(i_Z)^!\mF'\in \shv(Z)\cong \ind\der_{\ctf}(Z)$ (see \eqref{eq:shv as indcons for pfp} for the last equivalence). We need to show that $C^\bullet(X_{\leq \mu}(b),\mF')\in \Mod_\La^{\geq \langle 2\rho, \nu_b\rangle}$.

We can stratify $X_{\leq \mu}(b)$ as $\sqcup_{\mu'\leq \mu}X_{\mu'}(b)$. Let $\mF'_{\mu'}$ be the $!$-restriction of $\mF'$ to $X_{\mu'}(b)$. This is a finite stratification. 
By the standard spectral sequence for stratified spaces, it is enough to show that 
$C^\bullet(X_{\mu'}(b),\mF'_{\mu'})\in \Mod_\La^{\geq \langle 2\rho, \nu_b\rangle}$. As the $!$-restriction of $\mF$ to $L^+\mP\backslash LG_{\mu'}/L^+\mP$ can be written as extensions of $\consdual_{L^+\mP\backslash LG_{\mu'}/L^+\mP}[-\langle 2\rho, \mu'\rangle-i]$ for $i\geq 0$ (by the definition of perverse $t$-structure on $\fgshv(L^+\mP\backslash LG/L^+\mP)$), we see that $\mF'_{\mu'}$ can be written as extensions of $\consdual_{X_{\mu'}(b)}[-\langle2\rho, \mu'\rangle-i]$, for $i\geq 0$.
As $C^\bullet(X_{\mu'}(b), \consdual_{X_{\mu'}(b)})$ is nothing but the usual Borel-Moore homology of $X_{\mu'}(b)$, it belongs to $\Mod_\La^{\geq -2\dim X_{\mu'}(b)}$.

Also recall that for each $\mu'$ with $X_{\mu'}(b)$ non-empty, we have
\[
\dim X_{\mu'}(b)= \langle\rho, \mu'-\nu_b\rangle-\frac{1}{2}\mathrm{def}_G(b).
\]
It follows that $C^\bullet(X_{\leq \mu}(b),\mF')\in \Mod_\La^{\geq \langle 2\rho, \nu_b\rangle+\mathrm{def}_G(b)}$, as desired.
\end{proof}

Next passing to right adjoints of \eqref{eq: t-structure-open-closed gluing-1}, we also obtain
   \begin{equation}\label{eq: t-structure-open-closed gluing-2}
    \xymatrix{
    \shv(\kot_{G,<b})\ar@/^/[rr]^{(i_{<b})_*} \ar@/_/[rr]_{(i_{<b})_\flat} && \ar[ll]|{(i_{<b})^!} \shv(\kot_{G,\leq b}) \ar@/^/[rr]^{(j_{b})^!} \ar@/_/[rr]_{(j_{b})\rstar} &&  \ar[ll]|{(j_b)_*}\shv(\kot_{G,b}).
    }
    \end{equation}  
 Here all the involved functors preserve admissible objects by \Cref{lem:basic cpt and adm} \eqref{lem:basic cpt and adm-3}.

 \begin{proposition}\label{prop: t-structures on llc-2}
 Let $\chi$ be a character such that $\langle \chi, \nu_b\rangle\in\bZ$ as before. Then
\[
\shv(\kot_G)^{\chi\mbox{-}e,\leq 0}\subset \shv(\kot_G),\quad  \mbox{resp.} \quad \shv(\kot_G)^{\chi\mbox{-}e,\geq 0}\subset \shv(\kot_G)
\]
consisting of those $\mF$ such that 
\[
(i_b)^!\mF\in\rep(G_b(F))^{\leq \langle \chi,\nu_b\rangle },\quad  \mbox{resp.} \quad (i_b)\rstar\mF\in\rep(G_b(F))^{\geq \langle \chi,\nu_b\rangle}
\] 
for each $b\in B(G)$. Then the pair $(\shv(\kot_G)^{\chi\mbox{-}e,\leq 0}, \shv(\kot_G)^{\chi\mbox{-}e,\geq 0})$ defines an accessible $t$-structure on $\shv(\kot_G)$, which restricts to a $t$-structure on $\shv(\kot_G)^\adm$. 
\end{proposition}

\begin{proof}
Using \eqref{eq: t-structure-open-closed gluing-2} and  \cite[Theorem 1.4.10]{BBDG},
we see that for $\delta\in B(G)$, the pair
\[
\shv(\kot_{G,\leq \delta})^{\chi\mbox{-}e,\leq 0}=\bigl\{\mF\mid (i_{b})^!\mF\in \rep(J_{b}(F))^{\leq\langle\chi,\nu_{b}\rangle}, \ \forall\ b\leq \delta \bigr\},
\]
\[ 
\shv(\kot_{G,\leq \delta})^{\chi\mbox{-}e,\geq 0}=\bigl\{\mF\mid (i_{b})\rstar \mF\in \rep(J_{b}(F))^{\geq\langle\chi,\nu_{b}\rangle}, \ \forall\ b\leq \delta \bigr\}.
\]
define a $t$-structure on on $\shv(\kot_{G,\leq \delta})$, which is accessible. In addition, since all functors in \eqref{eq: t-structure-open-closed gluing-2} preserve subcategory of admissible objects, this $t$-structure restricts to a $t$-structure of the subcategory of admissible objects. This means that if we have a cofiber sequence $\mF'\to \mF\to \mF''$ in $\shv(\kot_{G,\leq \delta})$ with $\mF'\in \shv(\kot_{G,\leq \delta})^{\chi\mbox{-}e,\leq 0}, \mF''\in \shv(\kot_{G,\leq \delta})^{\chi\mbox{-}e,\geq 0}$ and $\mF\in \shv(\kot_{G,\leq \delta})^{\adm}$, then $\mF',\mF''\in \shv(\kot_{G,\leq \delta})^{\adm}$.

We note that for $\delta\leq \delta'$, the  $!$-restriction $\shv(\kot_{G,\leq \delta'})\to \shv(\kot_{G,\leq \delta})$ along the pfp closed embedding $\kot_{G,\leq \delta}\hookrightarrow \kot_{G,\leq \delta'}$ is exact with respect to such $t$-structure, while  the inclusion $\shv(\kot_{G,\leq \delta})\subset \shv(\kot_{G,\leq \delta'})$ induced by $*$-extension is right exact.  In addition, if $\mF\in \shv(\kot_{G})^{\chi\mbox{-}e,\leq 0}$ (resp. $\mF\in \shv(\kot_{G})^{e,\geq 0}$), then $(i_{\leq b})^!\mF\in  \shv(\kot_{G,\leq \delta})^{e,\leq 0}$ (resp. $(i_{\leq b})\rstar\mF\in  \shv(\kot_{G,\leq \delta})^{e,\geq 0}$). As for every $\mF\in \shv(\kot_G)$, we have $\mF= \colim_{B(G)} (i_{\leq b})_*((i_{\leq b})^!\mF)$, we see that 
\[
\shv(\kot_{G})^{\chi\mbox{-}e,\leq 0}= \colim_{B(G)}\shv(\kot_{G,\leq \delta})^{\chi\mbox{-}e,\leq 0}.
\]
This in particular implies that $\shv(\kot_{G})^{\chi\mbox{-}e,\leq 0}$ is compactly generated. In fact, the collection of objects
\begin{equation}\label{eq: connective part generator-2}
(i_b)_*\cind_K^{G_b(F)}\La[n-\langle\chi,\nu_b\rangle], \quad b\in B(G),\ n\geq 0, \ K\subset G_b(F) \mbox{ prop-}p \mbox{ open compact}.
\end{equation} 
form a set of compact generators of $\shv(\kot_{G})^{\chi\mbox{-}e,\leq 0}$.

Now \cite[Proposition 1.4.4.11]{Lurie.higher.algebra} implies that  $\shv(\kot_{G})^{\chi\mbox{-}e,\leq 0}$ indeed form the connective part of an accessible $t$-structure of $\shv(\kot_G)$. In addition, the above explicit description of the generators of $\shv(\kot_{G})^{\chi\mbox{-}e,\leq 0}$ implies that coconnective part of this $t$-structure is $\shv(\kot_{G})^{\chi\mbox{-}e,\geq 0}$.

Now let $\mF\in \shv(\kot_G)$, fitting into the cofiber sequence $\mF'\to \mF\to \mF''$ with $\mF'\in \shv(\kot_{G})^{\chi\mbox{-}e,\leq 0}$ and $\mF''\in \shv(\kot_{G})^{\chi\mbox{-}e,\geq 1}$. Then for every $\delta\in B(G)$, we have the following cofiber sequence in $\shv(\kot_{G,\leq \delta})$
\[
(i_{\leq \delta})^!\mF'\to (i_{\leq \delta})^!\mF\to (i_{\leq \delta})^!\mF''
\]  
with $(i_{\leq \delta})^!\mF'\in \shv(\kot_{G,\leq \delta})^{\chi\mbox{-}e,\leq 0}$ and $(i_{\leq \delta})^!\mF''\in \shv(\kot_{G,\leq \delta})^{\chi\mbox{-}e,\geq 1}$. If $\mF$ is admissible, then $(i_{\leq \delta})^!\mF$ is admissible in $\shv(\kot_{G,\leq \delta})$, and so is $(i_{\leq \delta})^!\mF'$ and $(i_{\leq \delta})^!\mF''$. 
It follows from \Cref{cor: characterization of admissible} that $\mF'$ is admissible. This proposition is proved.
\end{proof}

\begin{proposition}\label{lem: adm duality and e-t-structure}
Suppose $\La$ is a field and let $\chi=2\rho$. Then the duality $(\verd^{\can}_{\kot_G})^{\adm}$ interchanges $\shv(\kot_G)^{\chi\mbox{-}e,\leq 0}\cap \shv(\kot_G)^{\adm}$ and $\shv(\kot_G)^{\chi\mbox{-}e,\geq 0}\cap \shv(\kot_G)^{\adm}$.
\end{proposition}
\begin{proof}
This follows from the definition of the $t$-structure from \Cref{prop: t-structures on llc-2}  and \Cref{prop: adm.duality.kottwitz.pullpush}.
\end{proof}

\begin{remark}
We note that the $t$-structure on $\shv(\kot_G)^\adm$ is not bounded. For example, $\consdual_{\kot_G}\in \shv(\kot_G)^{\adm}$ but it has infinite negative cohomological degree with respect to this $t$-structure.
\end{remark}

\begin{remark}\label{rem: two t-str on stratified spaces}
Readers can skip this long remark. The construction of the above two $t$-structures applies in more general setting. Namely, let 
\[
X=\colim_{A} X_{\leq \al}
\] 
be an ind-stack such that for every $\al'<\al$, $X_{\leq \al'}\subset X_{\leq \al}$ is pfp closed, and for each $\al$, $X_{<\al}:=\cup_{\al'<\al} X_{\leq \al'}$ is also pfp closed in $X_{\leq \al}$. 
We write $j_{\al}: X_\al:=X_{\leq \al}-X_{<\al}\hookrightarrow X_{\leq \al}$ for the qcqs open complement of the closed embedding $i_{<\al}: X_{<\al}\subset X_{\leq \al}$. Now suppose 
\[
((j_\al)_!, (i_{<\al})^*), ((j_\al)^!, (i_{<\al})_*), ((j_\al)_*, (i_{<\al})^!), ((j_\al)\rstar, (i_{<\al})_\flat)
\] 
are all defined (each pair are the right adjoints of the previous pair), and suppose on each $\shv(X_\al)$ an accessible $t$-structure is assigned. Then one can define two $t$-structures on $\shv(X)$ by gluing the $t$-structures on various strata in two different ways, just as \Cref{prop: t-structures on llc-1} and \Cref{prop: t-structures on llc-2}. Namely, the first $t$-structure is defined so that
\[
\shv(X)^{p,\geq 0}=\{\mF\in\shv(X)\mid (i_\al)^!\mF\in \shv(X_\al)^{\geq 0}\},
\]
and the second $t$-structure is defined so that
\[
\shv(X)^{e,\geq 0}=\{\mF\in \shv(X)\mid (i_\al)^\sharp\mF\in \shv(X_\al)^{\geq 0}\}.
\]

The first construction is used to defined the usual perverse $t$-structure on stratified spaces. But as far as we know, the second construction is not considered in literature. One of the reasons is that the second construction requires one to work with big (a.k.a. presentable stable) categories of all sheaves, while classically people usually only work with small (a.k.a. idempotent complete) stable categories of constructible sheaves.

However, the second construction sometimes is also interesting.
For example, we consider the space $\iw\backslash LG/\iw$ equipped with the Schubert stratification and with the canonical generalized constant sheaf $\can$ as before.
Then the first gluing gives the usual perverse $t$-structure of $\shv(\iw\backslash LG/\iw)$. The second gluing, on the other hand, also has a nice interpretation when $F=\kappa((\varpi))$ is equal characteristic. 

For simplicity, we assume that $G$ arises as a split reductive group (denoted by the same notation) over $\kappa$. 
Let $\Bun_G(\bP^1)_{(0,\infty)}$ be the moduli space of $G$-bundles on $\bP^1_\kappa$ equipped with Iwahori level structure at $0,\infty$. Recall that geometric points of $\Bun_G(\bP^1)_{(0,\infty)}$ are still parameterized by $\widetilde{W}$. For $w\in \widetilde{W}$, let $\Bun_G(\bP^1)_w$ be the corresponding locally closed substack with underlying point corresponding to $w$. Then $\Bun_G(\bP^1)_e\subset \Bun_G(\bP^1)_{(0,\infty)}$ is open.
 Let $\mathrm{Eis}_e$ be the $!$-extension of the constant sheaf on $\Bun_G(\bP^1)_e$ to $\Bun_G(\bP^1)_{(0,\infty)}$, shifted to degree $\dim T$. Recall that
 $\shv(\iw\backslash LG/\iw)$ acts on $\shv(\Bun_G(\bP^1))_{(0,\infty)}$ as the Hecke operators at $0$ in the usual way (this can be made rigorous in $\infty$-categorical setting by applying the convolution pattern developed in \Cref{sec:trace} to the sheaf theory $\shv$ developed in \Cref{sec:pspl-stacks} ), and the action on $\mathrm{Eis}_e$ induces 
an equivalence
 \[
 \shv(\iw\backslash LG/\iw)\cong \shv(\Bun_G(\bP^1))_{(0,\infty)}, \quad \mF\mapsto \mF\star \mathrm{Eis}_e
 \]
(This functor is also known as the Radon transform.) It is not difficult to show that under this equivalence, the usual perverse $t$-structure on $\shv(\Bun_G(\bP^1))$ corresponds to the exotic $t$-structure on $ \shv(\iw\backslash LG/\iw)$. 
\end{remark}

\begin{remark}\label{rem: comparison between kotG and BunG}
As mentioned at the beginning of this section, one expects that $\shv(\kot_G)$ is equivalent to the appropriately defined category of $\ell$-adic sheaves on the Fargues-Fontaine curve. In addition, one expects that under such equivalence, the $t$-structure of $\shv(\kot_G)$ defined in \Cref{prop: t-structures on llc-2} (for $\chi=2\rho$) corresponds the natural perverse $t$-structure of the category  of $\ell$-adic sheaves on the Fargues-Fontaine curve, analogous to the relation between the exotic $t$-structure on $\shv(\iw\backslash LG/\iw)$ and the perverse $t$-structure on $\Bun_G(\bP^1)_{(0,\infty)}$ as discussed in \Cref{rem: two t-str on stratified spaces}.

On the other hand, one expects for appropriate choice of $\chi$, the $\chi$-perverse $t$-structure on $\shv(\kot_G)$ defined in \Cref{prop: t-structures on llc-1} corresponds the hadal $t$-structure defined by Hansen \cite{Hansen.Beijing}.
\end{remark}

\newpage

\section{Tame and unipotent local Langlands category}\label{sec: unipotent and tame LL category}

As should be clear from the previous discussions, the category $\shv(\kot_G)$ (and the closely related category $\rshv(\kot_G)$) can be regarded as the affine analogue of the category of representations of finite groups of Lie type. On the other hand, a very important approach to representations of finite groups of Lie type is the Deligne-Lusztig theory. In this section, we discuss some natural affine (and categorical) generalization of the Deligne-Lusztig theory, which gives us a way to access the tame part of $\shv(\kot_G)$.

In this section, let $F$ be a non-archimedean local field as in \Cref{sec:kot-stack}, with $k_F$ its residue field and $k$ an algebraic closure of $\kappa_F$.
we assume that $G$ splits over a tamely ramified extension of $F$. We fix a pinning $(B,T,e)$ of $G$ as before. Let 
$A\subset S\subset T$ be subtori of $G$ as before with $\mA\subset \mS\subset \mT$ the corresponding Iwahori group scheme over $\mO_F$. Let $\mI$ be the Iwahori group scheme of $G$ determined by the pinning as before. Let $\iw=L^+\mI$, and $\Sht^\loc=LG/\Ad_\sigma\iw$ as before. We base change everything to $k$. As before, if no confusion will arise, we omit $k$ from the subscript.
We recall that for $w\in \widetilde W$, there is a pfp locally closed embedding $i_w: LG_w\to LG$.

Let $(\hat{G},\hat{B},\hat{T} ,\hat{e})$ be the pinned dual group of $G$ over $\La$ equipped with an action of $\Ga_{\widetilde{F}/F}\subset\Aut(\hat{G},\hat{B},\hat{T} ,\hat{e})$ where $\widetilde{F}/F$ is the finite tame extension. Let ${}^LG=\hat{G}\rtimes \Ga_{\widetilde{F}/F}$ be the $L$-group and ${}^cG=\hat{G}\rtimes(\bG_m\times\Ga_{\widetilde{F}/F})$ be the $C$-group of $G$.

We will let $\La$ be a $\bZ_\ell$-algebra as at the end of \Cref{sec:adic-formalism}.

\subsection{Monodromic and equivariant categories}\label{SS: monodromic sheaves}
In this subsection, we discuss the formalism of monodromic and equivariant categories on a space equipped with an action of an affine algebraic group. We also refer to \cite{DYYZ} for some related discussions. 

\subsubsection{Serre's fundamental group of algebraic groups}
Let $H$ be a connected algebraic group over $k$. The universal cover of $H$ is defined to be the connected affine group scheme
\[
\widetilde{H}=\lim H',
\] 
where the inverse limit is taken over the cofiltered (ordinary) category of finite \'etale homomorphisms $H'\to H$ with $H'$ connected. 
Let 
\[
\pialg(H):=\ker(\widetilde{H}\to H)=\lim \ker(H'\to H),
\] 
regarded as a profinite group over $k$. 
As $\ker(H'\to H)$ must be central in $H'$, $\pialg(H)$ is in fact abelian. As each $H'\to H$ is \'etale, there is a surjective homomorphism $\pi^{\et}_1(H)\to \pialg(H)$. Note that if $H$ is defined over some subfield $k'\subset k$, then $\pialg(H)$ is equipped with a continuous action of $\Ga_{k/k'}$.

\begin{remark}
The supscript $c$ in $\pialg$ stands for ``central", as well as ``character". Namely, this group controls certain central extensions of $H$ as mentioned above, as well as character sheaves on $H$, as we shall see below.
\end{remark}

\begin{example}\label{ex: algpi of ss group}
If $H$ is a semisimple algebraic group over $k$, then its simply-connected cover $H_{\mathrm{sc}}$ is the universal cover of $H$ in the above sense. Therefore, $\pialg(H)=\ker(H_{\mathrm{sc}}(k)\to H(k))$. Note that $\pialg(H)$ is different from the usual algebraic fundamental group of $H$. (E.g. if $H=\mathrm{PGL}_p$ where $p$ is the characteristic of $k$, then $\pialg(H)$ is trivial.)
\end{example}

When $H$ is commutative, the group $\pialg(H)$ was firstly introduced by Serre \cite{Serre} (and was denoted by $\pi_1(H)$). As our base field $k$ is an algebraic closure of $\bF_p$, the group $\pialg(H)$ admits the following rather explicit description. 
\begin{lemma}\label{lem: pi1 of commutative group}
Suppose $H$ is commutative. 
For a  choice of the rational structure of $H$ over a (large enough) finite field $\bF_q\subset k$,
there is a canonical $\Ga_{k/\bF_q}$-equivariant isomorphism
\[
\pialg(H)= \lim_n H(\bF_{q^n}),
\]
where the transitioning maps are given by the norm map.
\end{lemma}
 Note that without considering the $\Ga_{k/\bF_q}$-structures,  the inverse limit on the right hand side is in fact independent of the  choice of the rational structure of $H$ over a finite field in $k$.
\begin{proof}
This is well-known. We sketch a proof for completeness. 
First notice that if $f: H'\to H$ is a finite isogeny, then $H'$ is a central extension of $H$ by $\ker(f)$, and the commutative pairing $H\times H\to \ker(f)$ is necessarily trivial (as $H$ is connected and $\ker(f)$ is finite). It follows that $H'$ is commutative. 
Choose a finite subfield $\bF_q\subset k$ such that the map $f$ and all points of $\ker(f)$ are defined over $\bF_q$. Note that $H'\to H'/H'(\bF_q)\cong H'$, where is the isomorphism is induced by the Lang isogeny of $H'$ (equipped with the $k_1$-rational structure). It follows that the Lang isogeny  $H'\to H'$ factors through $H'\xrightarrow{f} H\to H'$, or equivalently the Lang isogeny $H\to H$ covers $f:H'\to H$. It follows that $\widetilde{H}=\lim H'=\lim_n H$ where the second inverse limit is over $\bF_{q^n}$-Lang isogenies $H\to H$. The lemma follows.
\end{proof}

\begin{remark}\label{ex: algpi of tori}
Suppose $H$ is a torus. Then instead of using the Lang isogeny, one can use the multiplication by $n$ map, where $n$ coprime the characteristic exponent $p'$ of $k$. The same argument as above then shows that
\[
\pialg(H)\cong T^pH:=\lim_{(n,p')=1} H[n]
\]
is isomorphic to the Tate module of $H$. Note that this description of $\pialg(H)$ holds for tori over any algebraically closed field. 
\end{remark}

Now let $f: H_1\to H_2$ be a homomorphism, which induces a homomorphism $\pialg(f): \pialg(H_1)\to \pialg(H_2)$. 
\begin{lemma}\label{lem: map of alg fundamental group between tori}
Suppose $f: H_1\to H_2$ is surjective.
\begin{enumerate}
\item\label{lem: map of alg fundamental group between tori-1} If $\ker(f)$ is finite, then
we have a short exact sequence of profinite groups
\[
1\to \pialg(H_1)\to \pialg(H_2)\to \ker(f)\to 1.
\]
\item\label{lem: map of alg fundamental group between tori-2} If $f$ is surjective with $H_0:=\ker(f)$ connected, we have a right exact sequence of profinite groups
\[
\pialg(H_0)\to \pialg(H_1)\to \pialg(H_2)\to 1. 
\]
\end{enumerate}
\end{lemma}

\begin{proof}
If $f$ is finite, then $H_1\to H_2$ is an isogeny and therefore $\widetilde{H_2}$ maps surjectively to $H_1$. This gives Part \eqref{lem: map of alg fundamental group between tori-1}.

Next assume that $\ker(f)$ is connected. 
First note that any finite isogeny $H'_2\to H_2$ with $H'_2$ connected, its pullback to $H'_1\to H_1$ is still connected. Indeed, let $H_1^{'\circ}$ be the neutral connected component of $H'_1$, and let $\pi$ be the kernel of the map $H_1^{'\circ}\to H_1$. The $\pi$ maps injectively in to $\ker(H'_2\to H_2)$. Quotient out by $\pi$ gives the map $H_1\to H'_2/\pi$ and lifting $H_1\to H_2$. As $H$ is connected and $H'_2/\pi\to H_2$ has finite fibers, we see that $H'_2/\pi\cong H_2$. This shows that $H'_1$ is connected.
It follows that $\widetilde{H_1}\to \widetilde{H_2}$ is surjective and therefore $\pialg(H_1)\to \pialg(H_2)$ is surjective.

Next, let $\widetilde{H_0}'$ be the kernel of $\widetilde{H_1}\to \widetilde{H_2}$. Note that $\widetilde{H_0}'$ is connected. Otherwise, the quotient of $\widetilde{H_1}$ by the neutral connected component of $\widetilde{H_0}'$  would yet a non-trivial pro-finite \'etale cover of $\widetilde{H_2}$, contradicting the universal property of $\widetilde{H_2}$. It follows that there is a surjective map $\widetilde{H_0}\to \widetilde{H_0}'$ and therefore the sequence $\pialg(H_0)\to \pialg(H_1)\to \pialg(H_2)$ is exact in the middle.
\end{proof}

\begin{remark}
If $H_1$ in \Cref{lem: map of alg fundamental group between tori} \eqref{lem: map of alg fundamental group between tori-2} is commutative. Then the map $\pialg(H_0)\to \pialg(H_1)$ is injective, by \cite[\textsection{10.2}]{Serre}.
\end{remark}

\begin{corollary}\label{cor: Serre pi Mazur finiteness}
Let $H$ be an affine algebraic group over $k$. Then for every prime $\ell\neq p$, the pro-$\ell$-quotient of $\pialg(H)$ is topologically finitely generated.
\end{corollary}
\begin{proof}
Let $R_uH$ be the unipotent radical of $H$. 
By write $R_uH$ as successive extensions of $\bG_a$, we see that $\pialg(R_uH)$ is a pro-$p$ group by \Cref{lem: pi1 of commutative group} and \Cref{lem: map of alg fundamental group between tori} \eqref{lem: map of alg fundamental group between tori-2}. 
Then using \Cref{lem: map of alg fundamental group between tori} \eqref{lem: map of alg fundamental group between tori-2} again, we reduce the statement to the case when $H$ is connected reductive. 

Let $H_{\mathrm{der}}$ be the derived group of $H$ and $Z_{H}^\circ$ be the maximal torus in the center of $H$. Let $H_{\mathrm{sc}}\to H_{\mathrm{der}}$ be the simply-connected cover of $H_{\mathrm{der}}$. Let $A=H_{\mathrm{der}}(k)\cap Z_{H}^\circ(k)$ and $B= \ker(H_{\mathrm{sc}}(k)\to H_{\mathrm{der}}(k))$. Both are finite groups of order prime-to-$p$. Then by \Cref{lem: map of alg fundamental group between tori} \eqref{lem: map of alg fundamental group between tori-1}, combining with \Cref{ex: algpi of ss group} and \Cref{ex: algpi of tori},  we have
\[
1\to T^pZ_{H}^\circ\times B\to \pialg(H)\to A\to 1.
\]
The desired statement follows easily.
\end{proof}

Next, we consider the underived moduli space $(R_{\pialg(H),\bG_m})_{\cl}$ of strongly continuous homomorphisms from $\pialg(H)$ to $\bG_m$ over $\bZ_\ell$ as defined in \Cref{SS: Space of continuous representations}. 
\begin{lemma}\label{lem: ind-finite of moduli space R}
The space $(R_{\pialg(H),\bG_m})_\cl$ is represented by an ind-scheme, ind-finite over $\bZ_\ell$.
\end{lemma}
\begin{proof}
Using \Cref{cor: Serre pi Mazur finiteness}, it is enough to notice that (see \Cref{ex: continuous representation of Zhat})
$R_{\bZ_\ell,\bG_m}\subset \bG_m$ is the union of all closed subschemes $i_Z: Z\subset \bG_m$ are finite over $\bZ_\ell$ such that $Z\otimes_{\bZ_\ell}\bF_\ell$ are set theoretically supported at $1\in \bG_m$. 
\end{proof}

We will need to give another description of this space $(R_{\pialg(H),\bG_m})_\cl$ when $H$ is a torus. Let $\hat{H}$ be the dual torus of $H$ over $\bZ_\ell$.
Let $R_{I_F^t,\hat{H}}$ be the moduli space over $\bZ_\ell$ of strongly continuous $\hat{H}$-valued representations of $I_F^t$ (see \Cref{SS: Space of continuous representations}). Again, by \Cref{ex: continuous representation of Zhat}, if we fix a topological generator $\tau$ of $I_F^t$, we may identify
$R_{I_F^t,\hat{H}}\subset \hat{H}$ as the subfunctor $\hat{H}^{\wedge,p}\subset \hat{H}$ which is the union of all closed subschemes $i_Z:Z\subset \hat{H}$ that are finite over $\bZ_\ell$ such that $Z(\overline\bF_\ell)\subset \hat{H}(\overline\bF_\ell)^p$, where $ \hat{H}(\overline\bF_\ell)^p\subset \hat{H}(\overline\bF_\ell)$ consist of points of order prime-to-$p$. 

\begin{remark}\label{rem: formal neighborhood of chi}
For $\La=\overline\bF_\ell,\overline\bQ_\ell$ or $\overline\bZ_\ell$ (the integral closure of $\bZ_\ell$ in $\overline\bQ_\ell$),
we regard $\chi\in\hat{H}(\La)$ as a closed subscheme of $\hat{H}\otimes\La$, and denote by $\hat{\chi}$  the formal completion of $\hat{H}$ along $\chi$. We regard $\hat{\chi}$ as an indscheme. Let $\hat{H}(\overline\bZ_\ell)^p$ be those $\overline\bZ_\ell$-points of $\hat{H}$ whose reduction mod $\ell$ belong to $\hat{H}(\overline\bF_\ell)^p$, and let $\hat{H}(\overline\bQ_\ell)^p$ denote the image of $H(\overline\bZ_\ell)^p$ in $\hat{H}(\overline\bQ_\ell)$. Then for $\La=\overline\bF_\ell$ or $\overline\bQ_\ell$, we have an isomorphism
\[
R_{I_F^t,\hat{H}}\otimes\La\simeq \bigsqcup_{\chi\in \hat{H}(\La)^{p}} \hat{\chi}.
\]
However, $R_{I_F^t,\hat{H}}\otimes\overline\bZ_\ell$ is not the disjoint union of $\hat{\chi}$ over $\chi\in\hat{H}(\overline\bZ_\ell)^p$, as two points $\chi,\chi'\in \hat{H}(\overline\bZ_\ell)^p$ may meet over $\overline\bF_\ell$.
\end{remark}

\begin{lemma}
We have a canonical isomorphism
\[
R_{I_F^t,\hat{H}}\cong (R_{\pialg(H),\bG_m})_{\cl}
\]
where $T^pH=\varprojlim_{(n,p)=1} H[n]$ is the prime-to-$p$ Tate module of $H$.
\end{lemma}
\begin{proof}
As $R_{I_F^t,\hat{H}}$ is classical (i.e. no derived structure) so it is enough to prove $R_{I_F^t,\hat{H}}(A)\cong R_{\pialg(H),\bG_m}(A)$ for any classical $\bZ_\ell$-algebra $A$. In addition, since both spaces are ind-finite over $\bZ_\ell$, it is enough to consider the case $A$ is a finite $\bZ_\ell$-algebra. Then
we have 
\[
R_{I_F^t,\hat{H}}(A)=\Hom_{cts}(I_F^t, \xch(H)\otimes A^\times)=\Hom_{cts}(\xcoch(H)\otimes I_F^t, A^\times)\cong \Hom_{cts}(\pialg(H),A^\times).
\]
Here, the last isomorphism follows from
\[
\pialg(H)=\xcoch(H)\otimes \lim_{(n,p)=1} \mu_n(k)\cong \xcoch(H)\otimes \lim_{(n,p)=1} \mu_n(\breve F)\stackrel{\eqref{eq:iota-vs-tau}}{\cong} \xcoch(H)\otimes I^t_F.
\]
\end{proof}

\subsubsection{Monodromic sheaves on algebraic groups}

Let $H$ be an algebraic group over $k$ as above, and let $m: H\times H\to H$ denote the multiplication of $H$. 
 Let $\La$ be an algebraic extension of $\bF_\ell$, $\bZ_\ell$ or $\bQ_\ell$. Recall a character sheaf  (with coefficient in $\La$) on $H$ is a rank one $\La$-local system $\Ch_\chi$ on $H$ equipped with an isomorphism
\[m^*\Ch_\chi \simeq \Ch_\chi \boxtimes_{\La} \Ch_\chi,\]
satisfying the usual cocycle condition. Note that such an isomorphism necessarily induces a rigidification of $\Ch_\chi$ at the unit $1\in H$. The groupoid of character sheaves $\mathrm{CS}(H,\La)$ on $H$ forms a(n ordinary) Picard groupoid (and so its isomorphism classes form an abelian group). 

It is well-known that when $H$ is connected, the groupoid $\mathrm{CS}(H,\La)$ is discrete and therefore is an abelian group. In fact it is well-known that there is an isomorphism of abelian groups
\begin{equation}\label{eq: character sheaf-0}
(R_{\pialg(H),\bG_m})_\cl(\La)\cong \mathrm{CS}(H,\La),
\end{equation}
sending $\chi\in (R_{\pialg(H),\bG_m})_\cl(\La)$, corresponding to a continuous representation $\pialg(H)\to \La^\times$, to the rank one $\La$-local system  $\Ch_\chi$ on $H$ defined by $\pi_1^{\et}(H)\to \pialg(H)\to \La^\times$.

This isomorphism can be enhanced as follows. Let $\La$ be a Dedekind domain that is an algebraic extension of $\bF_\ell$, $\bQ_\ell$, or $\bZ_\ell$. In particular, $\La$ is regular. We use $(R_{\pialg(H),\bG_m})_\cl$ to denote its  base change to $\La$.
Notice that thanks to \Cref{lem: ind-finite of moduli space R}, the (stable) category of coherent sheaves on $(R_{\pialg(H),\bG_m})_\cl$ make sense, and the abelian category $\Coh((R_{\pialg(H),\bG_m})_\cl)^{\heartsuit}$ of coherent sheaves on $(R_{\pialg(H),\bG_m})_\cl$ is equivalent to the abelian category of continuous representations of $\pialg(H)$ on finite $\La$-modules.
Then we may lift the isomorphism \eqref{eq: character sheaf-0} as a functor
\begin{equation}\label{eq: character sheaf}
\Ch: \Coh((R_{\pialg(H),\bG_m})_\cl)^{\heartsuit}\to \shv(H,\La)^\heartsuit,
\end{equation}
sending $\mO_\chi$ to $\Ch_\chi$. Here $\chi\in (R_{\pialg(H),\bG_m})_\cl(\La')$ for some  finite $\La$-algebra $\La'$, and $\mO_\chi$ is regarded as an ordinary coherent sheaf on $(R_{\pialg(H),\bG_m})_\cl$ via the $*$-pushforward along the finite morphism $\Spec \La'\to (R_{\pialg(H),\bG_m})_\cl$.
This functor is clearly fully faithful, with the essential image denoted by $\shv_\mon(H,\La)^{\cpt,\heartsuit}$. It is 
the thick abelian subcategory generated by character sheaves (with coefficients in possible $\La$-algebras $\La'$). To see the last claim, just notice every continuous representation of $\pialg(H)$ on a finite $\La$-module can be filtered such that the successive quotients are generated over $\pialg(H)$ by one element. But if $M$ is generated over $\pialg(H)$ by one element, then $M$ is free of rank one over some finite $\La$-algebra $\La'$, and the action of $\pialg(H)$ on $M$ factors through $\pialg(H)\to (\La')^\times$. Then $\Ch(M)$ is a character sheaf on $H$ with coefficient in $\La'$.

\begin{lemma}\label{lem: descent of monodromic sheaves}
Let $f: H_1\to H_2$ be a surjective homomorphism. Let $\mF\in \shv(H_2)^{\heartsuit}$ such that $f^*\mF\in \shv_\mon(H_1,\La)^{\cpt,\heartsuit}$, then $\mF\in \shv_{\mon}(H_2,\La)^{\cpt,\heartsuit}$.
\end{lemma}
\begin{proof}
First notice that by descent, $\mF$ is a local system on $H_2$, and therefore corresponds to a representation $\pi_1^{\et}(H_2)$ on a finite $\La$-module $M$. 
By assumption, we know that the induced representation along $\pi_1^{\et}(H_1)\to \pi_1^{\et}(H_2)$ factors through $\pialg(H_1)$.
We need to show that $M$ is in fact a $\pialg(H_2)$-module.
It is enough to consider the case $\ker(f)$ is finite and $\ker(f)$ is connected separately. 

In the first case, we have $\ker(\pi_1^{\et}(H_1)\to \pialg(H_1))\cong \ker(\pi_1^{\et}(H_2)\to \pialg(H_2))$ by \Cref{lem: map of alg fundamental group between tori} \eqref{lem: map of alg fundamental group between tori-1}. Therefore, $M$ is indeed a $\pialg(H_2)$-module. In the second case, let $H_0=\ker(f)$. Then by \Cref{lem: map of alg fundamental group between tori} \eqref{lem: map of alg fundamental group between tori-2}, it is enough to show that the induced representation $\pialg(H_0)\to \pialg(H_1)$ on $M$ is trivial. But $\pi_1^{\et}(H_0)\to \pialg(H_0)$ is surjective and the action of $\pi_1^{\et}(H_0)$ on $M$ is trivial. Therefore, the action of $\pialg(H_0)$ on $M$ is trivial as well.
\end{proof}

\begin{definition}
Suppose $H$ is connected. 
\begin{enumerate}
\item Let $\shv_{\mon}(H,\La)\subset \shv(H,\La)$ denote the full $\La$-linear subcategory generated by $\shv_{\mon}(H,\La)^{\cpt,\heartsuit}$. We call $\shv_{\mon}(H,\La)$ the category of monodromic sheaves on $H$.
\item For a character sheaf $\Ch_\chi$ on $H$ with coefficient in $\La$, let $\shv_{\chi\mbox{-}\mon}(H,\La)$ be the full $\La$-linear subcategory of $\shv(H)$ generated by $\Ch_\chi$. We call $\shv_{\chi\mbox{-}\mon}(H,\La)$ the category of $\chi$-monodromic  sheaves on $H$. 
\end{enumerate}
\end{definition}

To simplify expositions, in the sequel we will use the notation to $(\chi\mbox{-})\mon$ to denote either $\chi$-monodromic or monodromic version.

\begin{remark}
We note that $\shv_{\mon}(H,\La)^\cpt\subset \cshv(H,\La)$ consist of those $\mF$ whose cohomology sheaves belong to $\shv_\mon(H,\La)^{\cpt,\heartsuit}$.
\end{remark}
In the sequel, we will omit $\La$ from the notation if it is clear from the context.

Note that a character sheaf $\Ch_\chi$ on $H$ determines a $\La$-linear functor $\iota_\chi:\Mod_\La\to \shv(H)$, which admits a factorization
\begin{equation}\label{eq: inclusion of character sheaves}
\Mod_\La\to \shv_{\chi\mbox{-}\mon}(H)\subset \shv_{\mon}(H)\subset \shv(H).
\end{equation}
All the above functors admit $\La$-linear right adjoint. In the sequel, we write the inclusion $\shv_{(\chi\mbox{-})\mon}(H)\subset \shv(H)$ as $\iota_{(\chi\mbox{-})\mon}$. 

We let
\begin{equation}\label{eq: monoidal right adjoint of inclusion of character sheaves}
    \av^{(\chi\mbox{-})\mon}:=(\iota_{(\chi\mbox{-})\mon})^R: \shv(H)\to \shv_{(\chi\mbox{-})\mon}(H),\quad \av^{\chi}=(\iota_\chi)^R: \shv(H)\to \Mod_\La
\end{equation}
and let
\begin{equation}\label{eq: monoidal unit character sheave-1}
\Ch_{(\chi\mbox{-})\mon}:=\av^{(\chi\mbox{-})\mon}(\delta_1), 
\end{equation}
where $\delta_1:=(\{1\}\to H)_*\La$ is the delta sheaf at the unit of $H$. Sometimes for simplicity we will also write 
\begin{equation}\label{eq: monoidal unit character sheave-2}
\widetilde{\Ch}=\Ch_\mon,\quad \Ch_{\hchi}=\Ch_{\chi\mbox{-}\mon}.
\end{equation}

\begin{example}\label{ex: monodromic system on unipotent group}
Let $H$ be an unipotent group. Let $\phi: H(\bF_q)\to \La^\times$ be a non-trivial character, giving $\pialg(H)\to H(\bF_q)\xrightarrow{\phi}\La^\times$. In this case, we have $\iota_\phi^{\mon}:\Mod_\La\cong \shv_{\phi-\mon}(\bG_a)$ is an equivalence and $\Ch_{\phi\mbox{-}\mon}=\Ch_\phi$. In the special case $H=\bG_a$, the corresponding character sheaf on $\bG_a$ is usually called the Artin-Schreier sheaf. 
\end{example}

For a description of $\Ch_{\hchi}$ when $H=\bG_m$, we refer to \Cref{ex: char and monodromy local systems}.
\begin{proposition}\label{lem: functoriality depth zero geom Langlands for tori}
 Let $f: H_1\to H_2$ be a homomorphism of connected algebraic groups. 
\begin{enumerate}
\item\label{lem: functoriality depth zero geom Langlands for tori-1} If $\Ch_{\chi_2}\in \mathrm{CS}(H_2)$, then $\Ch_{\chi_1}:=f^*\Ch_{\chi_2}\in \mathrm{CS}(H_1)$. The pullback functor $f^*: \shv(H_2)\to \shv(H_1)$ restricts to a pullback functor $f^*: \shv_{(\chi_2\mbox{-})\mon}(H_2)\to\shv_{(\chi_1\mbox{-})\mon}(H_1)$. The functor $f^*: \shv_{\mon}(H_2)\to \shv_{\mon}(H_1)$  admits a (continuous) right adjoint, denoted by $f_*^{\mon}$. When $f$ is surjective, we have
\[
f_*^{\mon}=f_*|_{\shv_{\mon}(H_1)}.
\] 
In general, we have 
\[
f_*^{\mon}= \av^{\mon}\circ (f_*|_{\shv_{\mon}(H_1)}).
\]

\item\label{lem: functoriality depth zero geom Langlands for tori-3} When $f$ is surjective, the usual compactly supported pushforward functor $f_!$ restricts to a functor between monodromic categories, which is the left adjoint of $f^!=f^*\langle \dim H_2-\dim H_1\rangle$. In addition, we have the following isomorphism of functors
\begin{equation}\label{eq: pseudo-properness for monodromic pushforwards}
f_![d]\cong f_*: \shv_{\mon}(H_1)\to \shv_{\mon}(H_2),
\end{equation}
for some integer $d$ depending on $\ker f$.
\end{enumerate}
\end{proposition}
\begin{proof}
That $f^*$ preserves monodromic categories is clear, and it's clear the right adjoint of $f^*$ is $f_*^{\mon}= \av^{\mon}\circ (f_*|_{\shv_{\mon}(H_1)})$. We show that if $f$ is surjective, both $f_*$ and $f_!$ preserve monodromic subcategories. We deal with the case $f_*$ and the case $f_!$ is similarly. 

It is enough to show that the cohomology sheaf $\mH^i f_*\Ch_\chi \in \shv_{\mon}(H_2)^{\heartsuit}$ for $\Ch_\chi$ a character sheaf on $H_1$ with coefficient in some finite extension $\La'$ of $\La$. By smooth base change, we have
\[
\mH^if^*f_*\Ch_\chi\cong \Ch_\chi\otimes_{\La'} H^i\rg(H_0, \Ch_\chi|_{H_0}).
\]
We apply then \Cref{lem: descent of monodromic sheaves} to conclude.

It remains to show that $f_![d]\cong f_*$ when restricted to the category of monodromic sheaves on $H_1$. We may factors $f$ as a finite isogeny and a homomorphism with connected fibers.
This case of finite isogeny is clear. So we suppose $f$ has connected fibers. Let $K=\ker f$, which is a connected affine group scheme (which is the perfection of an algebraic group). Let $B\subset K$ be a Borel subgroup of $K$. Then $H_1/B$ is proper over $H_2$. As argued in \Cref{lem: monodromic flat shrek pushforward} \Cref{prop: Frob str on coh of group}, it is enough to show that for a connected solvable group $H$, there is some integer $d$, such that $C_c(H,-)[d]\cong C(H,-)$ when restricted the category of monodromic sheaves on $H$. By further writing $H$ as an successive extensions of $\bG_a$ and $\bG_m$, we may assume that $H=\bG_a$ and $\bG_m$. Each case can be treated easily.
\end{proof}

Note that $\shv(H)$ has a natural monoidal structure given by $*$-pushforward along the multiplication map. Formally, it arises via the convolution pattern (see \Cref{rem:segal.objects.morphisms.vert.horiz} \Cref{rem-conv-product}) applied to $H=\pt\times_{\bB H}\pt$. The unit is given by $\delta_1$. \footnote{Note that $\shv(H)$ acquires another monoidal structure given by $!$-pushforward, and as we shall see when restricted to $\shv_{\mon}(H)$, the $!$-monoidal structure differs from the $*$-monoidal structure by a cohomological shift.  We will mainly use the $*$-monoidal structure, as it fits into the sheaf theory formalism for $\shv$. 
}

We need to understand the restriction of the above monoidal structure to $\shv_{(\chi\mbox{-})\mon}(H)$. We start with the following easy but important facts about the category monodromic sheaves.
\begin{lemma}\label{lem: tensor product tensor of monodromic sheaves} 
Let $H_1, H_2$ be two connected algebraic groups over $k$. Then the exterior tensor product functor $\shv(H_1)\otimes_\La \shv(H_2)\to \shv(H_1\times H_2)$ restricts to an equivalence
\begin{equation}\label{eq: monodromic exterior tensor product}
\shv_{\mon}(H_1)\otimes_\La \shv_{\mon}(H_2)\cong \shv_{\mon}(H_1\times H_2),
\end{equation}
which restricts to an equivalence $\shv_{\chi_1\mbox{-}\mon}(H_1)\otimes_\La \shv_{\chi_2\mbox{-}\mon}(H_2)\cong \shv_{(\chi_1\boxtimes \chi_2)\mbox{-}\mon}(H_1\times H_2)$.
\end{lemma}
\begin{proof}
The classical K\"unneth formula (e.g. see \Cref{lem: categorical kunneth} and \Cref{lem: categorical kunneth prestack}) implies that the functor is fully faithful.

On the other hand, we claim that the exterior tensor product induces an equivalence of groupoids 
\begin{equation}\label{eq:character sheaf tensor product}
\mathrm{CS}(H_1,\La')\times \mathrm{CS}(H_2,\La')\xrightarrow{\cong} \mathrm{CS}(H_1\times H_2,\La'),\quad (\Ch_{\chi_1},\Ch_{\chi_2})\mapsto \Ch_{\chi_1}\boxtimes_{\La'}\Ch_{\chi_2}.
\end{equation}
As $\Ch_{\chi_1}\boxtimes_{\La'}\Ch_{\chi_2}$ belongs to the subcategory of $\shv(H_1\times H_2)$ generated by $\Ch_{\chi_1}\boxtimes_{\La}\Ch_{\chi_2}$ under colimits,
this claim clearly implies that the functor is also essential surjective. 

To prove the claim, let $\Ch_{\chi}$ be a character sheaf on $H_1\times H_2$. Let $\Ch_{\chi_1}=\Ch_{\chi}|_{H_1\times \{1\}}$ and $\Ch_{\chi_2}=\Ch_{\chi}|_{\{1\}\times H_2}$. Then using the isomorphism $H_1\times H_2\cong (H_1\times\{1\})\times(\{1\}\times H_2)\xrightarrow{m} H_1\times H_2$ and the character property of $\Ch_{\chi}$, we see that $\Ch_{\chi}\cong \Ch_{\chi_1}\boxtimes_{\La'}\Ch_{\chi_2}$. (Note that the claim holds even without connectedness assumption of $H_1$ and $H_2$.)

The last statement is clear.
\end{proof}

\begin{proposition}\label{lem: monoidal of monodromic}
\begin{enumerate}
\item\label{lem: monoidal of monodromic-3} The functor $\av^{(\chi\mbox{-})\mon}$ is given by  $\Ch_{(\chi\mbox{-})\mon}\star -$.
\item\label{lem: monoidal of monodromic-2}  Both $\shv_{\chi\mbox{-}\mon}(H)$ and $\shv_{\mon}(H)$ have natural monoidal structures. All of the right adjoint of functors in \eqref{eq: inclusion of character sheaves} admit canonical monoidal structures.
\end{enumerate}
\end{proposition}

\begin{proof}
Note that by adjunction, we have a map $\Ch_{(\chi\mbox{-})\mon}\to \delta_1$, which induces for every $\mG\in\shv(H)$ a map 
\begin{equation}\label{eq: counit for avmon}
\Ch_{(\chi\mbox{-})\mon}\star\mG\to \mG.
\end{equation}
For Part \eqref{lem: monoidal of monodromic-3}, we need to show that $\Ch_{(\chi\mbox{-})\mon}\star\mG\in \shv_{(\chi\mbox{-})\mon}(H)$ and 
for every $\mF\in \shv_{(\chi\mbox{-})\mon}(H)$, the map \eqref{eq: counit for avmon} induces an isomorphism
\begin{equation}\label{eq: avmon=conv}
\Hom(\mF, \Ch_{(\chi\mbox{-})\mon}\star\mG)\cong \Hom(\mF,\mG).
\end{equation} 

We first verify $\Ch_{(\chi\mbox{-})\mon}\star\mG\in \shv_{(\chi\mbox{-})\mon}(H)$.
First, notice that if $\mF\in\shv(H)$ such that $m^*\mF\simeq \Ch_\chi\boxtimes \mF$ for some character sheaf $\chi$, then pulling back along $H\xrightarrow{h\mapsto (h,1)} H\times H\xrightarrow{m} H$ shows that 
\[
\mF\simeq \Ch_\chi\otimes^* (\{1\}\to H)^*\mF\in\shv_{\chi\mbox{-}\mon}(H).
\] 
Next let $\chi$ be a character sheaf. Then by smooth base change, we see that $m^*(\Ch_\chi\star \mF)=m^*m_*(\Ch_\chi\boxtimes \mF)\cong \Ch_\chi\boxtimes (\Ch_\chi\star \mF)$. It follows that for $\mF\in \shv_{(\chi\mbox{-})\mon}(H)$ and any $\mG\in \shv(H)$,
we have $\mF\star \mG\in \shv_{(\chi\mbox{-})\mon}(H)$. Similarly, $\mG\star\mF\in \shv_{(\chi\mbox{-})\mon}(H)$. 
It follows that $\Ch_{(\chi\mbox{-})\mon}\star\mG\in \shv_{(\chi\mbox{-})\mon}(H)$.

To show \eqref{eq: avmon=conv},
we may assume that $\mF$ is a character local system on $H$. Then
\[
\Hom(\mF, \widetilde\Ch\star\mG)=\Hom(\mF\boxtimes_{\La}\mF,  \av^{\mon}(\delta_1)\boxtimes_{\La}\mG)=((\{1\}\to H)^*\mF)^\vee\otimes_\La\Hom(\mF,\mG)=\Hom(\mF,\mG),
\]
as desired. 

For Part \eqref{lem: monoidal of monodromic-2},  we shall only prove that $\shv_{\mon}(H)$ has a natural monoidal structure and $\av^{\mon}$ has a natural monoidal structure. All other cases are proved in the same way.

First notice by the above argument and by \Cref{lem: complement on monoidal categories}, $\shv_{(\chi\mbox{-})\mon}(H)$ is a $\shv(H)$-bimodule.
In addition, by \Cref{lem: tensor product tensor of monodromic sheaves}, the natural map 
\[
\av^{\mon}(\mG_1)\boxtimes_{\La}\av^{\mon}(\mG_2)\to \av^{\mon}(\mG_1\boxtimes_{\La}\mG_2)
\] 
is an isomorphism,  for $\mG_i\in \shv(H_i)$ for $i=1,2$.  Now for $\mF\in \shv_{\mon}(H)$, we have 
\begin{multline*}
\Hom(\mF, \av^{\mon}(\mG_1\star\mG_2))=\Hom(\mF,\mG_1\star\mG_2)=\Hom(m^*\mF, \mG_1\boxtimes_{\La}\mG_2)\\
=\Hom(m^*\mF,\av^{\mon}(\mG_1\boxtimes_{\La}\mG_2))=\Hom(\mF,\av^{\mon}(\mG_1)\star\av^{\mon}(\mG_2)).
\end{multline*}

Now by \Cref{lem: complement on monoidal categories}, we see that $\shv_{\mon}(H)$ has a monoidal structure, with $\widetilde\Ch=\av^{\mon}(\delta_1)$ a monoidal unit. In addition, $\av^{\mon}$ is monoidal.
\end{proof}

\begin{remark}\label{rem: of lem monoidal of monodromic} 
One can show that the category  $\shv_{(\chi\mbox{-})\mon}(H)$ can be identified with the category consisting of objects in $\shv(H)$ equipped with an action of $\Ch_{(\chi\mbox{-})\mon}$. We do not need this fact.
\end{remark}

\begin{proposition}\label{prop:semi-rigidity of monodromic category}
Equipped with the above monoidal structure, $\shv_{(\chi\mbox{-})\mon}(H)$ is semi-rigid. 
\end{proposition}
\begin{proof}
We notice that for two character local systems $\Ch_{\chi_1}$ and $\Ch_{\chi_2}$ of $H$ and for $?=*$ or $!$, by base change we have 
\[
m^*m_?(\Ch_{\chi_1}\boxtimes\Ch_{\chi_2})\cong \Ch_{\chi_1}\boxtimes m_?(\Ch_{\chi_1}\boxtimes\Ch_{\chi_2}).
\] 
Then the argument in the proof of \Cref{lem: monoidal of monodromic} implies that $m_?:\shv(H\times H)\to \shv(H)$ sends $\shv_\mon(H)\otimes\shv_\mon(H)\subset \shv(H\times H)$ to $\shv_\mon(H)\subset \shv(H)$. Now we factor the multiplication $m$ as 
\[
H\times H\stackrel{(h_1,h_2)\mapsto (h_1h_2,h_2)}{\cong} H\times H\xrightarrow{\pr_2} H.
\] 
The $*$ and $!$-pushforwards along the first morphism are identified and send monodromic sheaves to monodromic sheaves. Then we apply \Cref{lem: functoriality depth zero geom Langlands for tori} \eqref{lem: functoriality depth zero geom Langlands for tori-3} to $f=\pr_2$ to conclude that $m_*$ and $m_!$ differ by a shift.
It follows that $m_*$ has a continuous right adjoint given by $m^*$ up to shift and from the base change that $m^*$ is $\shv_\mon(H)$-bilinear. On the other hand $\shv_\mon(H)$ is compactly generated by definition. Therefore, $\shv_\mon(H)$ is semi-rigid.

The $\chi$-monodromic case is similar (and in fact simpler).
\end{proof}

Via the monoidal functors in \eqref{eq: monoidal right adjoint of inclusion of character sheaves}, we may regard $\shv_{(\chi\mbox{-})\mon}(H)$ and $\Mod_\La$ as (left) $\shv(H)$-modules. When emphasizing the module structure, $\Mod_\La$ will be denoted as $(\Mod_\La)_\chi$. We note that with the equipped $\shv(H)$-module structures, all the functors in \eqref{eq: inclusion of character sheaves} are $\shv(H)$-linear.

We will let $\lincat_{\shv(H)}$ denote the ($2$-)category of left $\shv(H)$-modules in $\lincat_\La$ (see \Cref{SS: la-linear categories}). Recall that all $\shv(H)$-linear functors between two $\shv(H)$-modules $\bfM$ and $\bfN$ form a $\La$-linear category $\fun^{\mathrm{L}}_{\shv(H)}(\bfM,\bfN)$.

\begin{lemma}\label{lem: self dual of monodromic and equiv sheaves on H}
Let $\bfM$ be any of categories in \eqref{eq: inclusion of character sheaves}. Then $\bfM$ equipped with the right $\shv(H)$-module structure is a left dual (in the sense of \Cref{def: the notation of left dual}) of $\bfM$ as a left $\shv(H)$-module. 
\end{lemma}
\begin{proof}
Let $\bfM= \shv_{(\chi\mbox{-})\mon}(H)$. Notice that we have 
\[
\shv_{(\chi\mbox{-})\mon}(H)\cong \shv_{(\chi\mbox{-})\mon}(H)\otimes_{\shv(H)}\shv_{(\chi\mbox{-})\mon}(H).
\] 
Then the unit $u$ is just given by $\Ch_{(\chi\mbox{-})\mon}$, and the co-unit $e$ is given by 
\[
\shv_\mon(H)\otimes_\La\shv_\mon(H)\to \shv_\mon(H)\to \shv(H),
\]
where the first functor is the tensor product of $\shv_\mon(H)$ and the second  functor is one from \eqref{eq: inclusion of character sheaves}.

Next let $\bfM=(\Mod_\La)_\chi$. As the functor $\av^\chi: \shv(H)\to (\Mod_\La)_\chi$ factors through $\shv_{\chi\mbox{-}\mon}(H)\to (\Mod_\La)_\chi$, it is enough to show the duality as  $\shv_{\chi\mbox{-}\mon}(H)$-modules. But $\shv_{\chi\mbox{-}\mon}(H)$ is semi-rigid (by \Cref{prop:semi-rigidity of monodromic category}), we can apply \Cref{rem-duality-datum-as-plain-cat} \eqref{rem-duality-datum-as-plain-cat-2} to conclude.
\end{proof}

\subsubsection{Monodromic and equivariant categories}

Now let $X$ be a prestack over $k$ acted by an affine algebraic group $H$. Then $\shv(X)$ is an $\shv(H)$-module with the action given by $*$-pushforward. Again formally, it arises via the convolution pattern (see \Cref{rem:segal.objects.morphisms.vert.horiz} \Cref{rem-conv-product}) applied to $H=\pt\times_{\bB H}\pt$ and $X=\pt\times_{\bB H}H\backslash X$.

\begin{definition}\label{def: mono sheaves and equiv sheaves}
Let $X$ be a prestack with an action of an algebraic group $H$ (from the left).
\begin{enumerate}
\item    We define the category of $H$-monodromic sheaves on $X$ as 
    \[
    \shv((H,\mon)\bs X):=\fun^{\mathrm{L}}_{\shv(H)}(\shv_{\mon}(H),\shv(X)).
    \]
\item    For $\Ch_\chi\in\mathrm{CS}(H)$, we define the category of $(H,\chi)$-monodromic sheaves on $X$ as 
    \[
    \shv((H,\chi\mbox{-}\mon) \bs X):=\fun^{\mathrm{L}}_{\shv(H)}(\shv_{\chi\mbox{-}\mon}(H),\shv(X)).
    \]
\end{enumerate}
    In the sequel, if the group $H$ is clearly from the context, we will also just write
\[
\shv_\mon(X)= \shv((H,\mon)\bs X),\quad \shv_{\chi\mbox{-}\mon}(X)= \shv((H,\chi\mbox{-}\mon) \bs X).
\]
For the reason which will be clear later, we will also write $\shv((H,\chi\mbox{-}\mon) \bs X)$ as $\shv((H,\hchi)\bs X)$.
\begin{enumerate}[resume]
\item    We define the category of $(H,\chi)$-equivariant sheaves on $X$ as 
    \[
    \shv((H,\chi)\bs X):=\Hom_{\shv(H)}((\Mod_\La)_\chi,\shv(X)).
    \]
 \end{enumerate}   
\end{definition}

As before, in the sequel we use $(\chi\mbox{-})\mon$ to denote either $\chi$-monodromic or all monodromic version.
Here are a few basic facts about monodromic and equivariant categories of sheaves.

It follows from \Cref{lem: self dual of monodromic and equiv sheaves on H} that for general $X$, the category $\shv_{(\chi\mbox{-})\mon}(X)$ of ($\chi$-)monodromic sheaves on $X$ can be identified with
    \[
    \fun^{\mathrm{L}}_{\shv(H)}(\shv_{(\chi\mbox{-})\mon}(H),\shv(X))\cong \shv_{(\chi\mbox{-})\mon}(H)\otimes_{\shv(H)}\shv(X).
    \]
    In particular, when $X=H$ equipped with the natural left action, the notation is consistent with the previous notation.
As the adjoint pair of functors
    \[
    \iota_{(\chi\mbox{-})\mon}: \shv_{(\chi\mbox{-})\mon}(H) \rightleftharpoons \shv(H): \av^{(\chi\mbox{-})\mon}
    \]
    realize $\shv_{(\chi\mbox{-})\mon}(H)$ as a colocalization of $\shv(H)$ as $\shv(H)$-modules, we see that we have a pair of adjoint functors
    \[
    \iota_{X,(\chi\mbox{-})\mon}:  \shv_{(\chi\mbox{-})\mon}(X)\rightleftharpoons \shv(X): \av_X^{(\chi\mbox{-})\mon}
    \]
    realizing $\shv_{(\chi\mbox{-})\mon}(X)$ as a colocalization of $\shv(X)$. Similarly to \Cref{lem: monoidal of monodromic} \eqref{lem: monoidal of monodromic-3}, we have  
    \begin{equation*}\label{average}
    \av_X^{(\chi\mbox{-})\mon}=\Ch_{(\chi\mbox{-})\mon}\star(-).   
\end{equation*}

Similarly we can identify $\shv((H,\chi)\bs X)$ with
    \begin{equation}\label{eq: monodromic to equivariant}
    (\Mod_\La)_\chi\otimes_{\shv(H)}\shv(X)\cong (\Mod_\La)_\chi\otimes_{\shv_{\mon}(H)}\shv_{\mon}(X)\cong (\Mod_\La)_\chi\otimes_{\shv_{\chi\mbox{-}\mon}(H)}\shv_{\chi\mbox{-}\mon}(X).
    \end{equation}
    
    \begin{remark}\label{rem: equivariant category not subcategory}
    Note that if $H$ is unipotent, then $\shv((H,\chi)\bs X)\cong\shv((H,\hchi)\bs X)\subset \shv(X)$ by virtue of \Cref{ex: monodromic system on unipotent group}.
    
    However, this is not the case in general if $H$ is not unipotent. But as the functor $(\Mod_\La)_\chi\to \shv_{\chi\mbox{-}\mon}(H)$ is $\shv(H)$-linear and the image generates the target, we see that $\shv((H,\hchi)\bs X)$ is generated (as $\La$-linear category) by the essential image of the functor $\shv((H,\chi)\bs X)\to \shv(X)$. On the other hand, using the expression $\shv((H,\hchi)\bs X)\cong \shv_{\chi\mbox{-}\mon}(H)\otimes_{\shv(H)} \shv(X)$, we see that $\shv_{\chi\mbox{-}\mon}(X)$ is generated by $a_*(\Ch_\chi\boxtimes \mF)$ for $\Ch_\chi$ being character sheaves on $H$ and $\mF\in \shv(X)$. Here $a: H\times X\to X$ denotes the action map. Therefore, our definition of the category of $(\chi\mbox{-})$monodromic sheaves on $X$ coincides with other definition used in literature.
\end{remark}

To justify the definition of the category of equivariant sheaves, we notice the following statement.
\begin{lemma}\label{lem: two-definition of equivariant category}
Let $X$ be a prestack with a (left) $H$-action.
Let $u$ be the trivial $\La$-local system on $H$, regarded as a character sheaf. Then 
\[
\shv\bigl((H,u)\backslash X\bigr)\cong \shv(H\backslash X),
\] 
and the natural pair of adjoint functors $\shv((H,u)\backslash X)\rightleftharpoons\shv(X)$ is identified with the natural $*$-pullback and pushforward along
 $X\to H\backslash X$.
\end{lemma}
\begin{proof}
We note that $\shv((H,u)\backslash X)$ can be identified with the geometric realization of the simplicial diagram of categories
\[
\shv_{u\mbox{-}\mon}(H)^{\otimes \bullet}\otimes \shv_{u\mbox{-}\mon}(X)\cong \shv_{u\mbox{-}\mon}(H)^{\otimes \bullet+1}\otimes_{\shv(H)}\shv(X),
\]
with coface maps given by $*$-pushfowards along multiplication maps between adjacent $H$s. However, by virtue of \eqref{eq: pseudo-properness for monodromic pushforwards}, we may pass to the right adjoint to obtain a cosimplicial diagram, and then applying shift $[\dim H]^{\otimes\bullet}$ to the cosimplicial diagram.
The resulting cosimplicial diagram then maps fully faithfully to the cosimplicial diagram $\shv(H^\bullet\times X)$ (with face maps being $!$-pullbacks). By descent, the totalization of later is just $\shv(H\backslash X)$. It follows that we have the fully faithful embedding  $\shv((H,u)\backslash X)\to \shv(H\backslash X)$. On the other hand, its right adjoint is conservative, (as $\shv(H\backslash X)\to \shv(X)$ is conservative). It follows that  $\shv((H,u)\backslash X)\to \shv(H\backslash X)$ is an equivalence. The last identification of functors is also clear.
\end{proof}

\begin{remark}\label{ex: endomorphism of monodromic sheaf}
Let $\Ch_u$ be the trivial local system on $H$. Objects in $\shv_{u\mbox{-}\mon}(X)$ are usually called unipotent monodromic sheaves on $X$\footnote{We caution the readers that in some literature these are simply called monodromic sheaves.}.
It follows from \Cref{rem: equivariant category not subcategory} and \Cref{lem: two-definition of equivariant category} that $\shv_{u\mbox{-}\mon}(X)$ is generated by essential image of the $!$-pullback functor $\shv(H\backslash X)\to \shv(H)$. This coincides with the usual definition of unipotent monodromic categories.
As mentioned at the end of \Cref{rem: equivariant category not subcategory}, the category $\shv_{u\mbox{-}\mon}(X)$ can also be generated by objects $a_*(\La\boxtimes \mF)$ for $\mF\in\shv(X)$, where $a: H\times X\to X$ is the action map.  

Since $\shv_{u\mbox{-}\mon}(X)$ is a module category over $\shv_{u\mbox{-}\mon}(H)$, the algebra $\End(\Ch_{u\mbox{-}\mon})$ 
acts on every object $\mF$ in $\shv_{u\mbox{-}\mon}(X)$.  When $H$ is an algebraic torus, this gives the usual monodromy action. 
\end{remark}

\begin{remark}\label{ex: monodromic sheaf under change of group}
Let $\varphi: H'\to H$ be a homomorphism. Suppose $H'$ acts on $X$ through an action of $H$ on $X$. Using the last statement from \Cref{rem: equivariant category not subcategory}, we see that
$ \shv_{H\mbox{-}\mon}(X)\subset  \shv_{H'\mbox{-}\mon}(X)$, and if $\varphi$ is surjective this inclusion is in fact an equivalence.
\quash{ the left adjoint of the above functor is conservative.
To see this, note that by \Cref{lem: functoriality depth zero geom Langlands for tori}, $\varphi^*$ sends $\shv^{\mon}(H)$ to $\shv^{\mon}(H')$.

it is enough to consider the case $X=H$ equipped with left $H$-action. 
By \Cref{lem: functoriality depth zero geom Langlands for tori} and \Cref{lem: monoidal of monodromic}, we see that $\varphi_*$ sends $\shv^{\mon}(H')\to \shv^{\mon}(H)$, compatible with their module structures (over $\varphi_*:\shv(H')\to \shv(H)$), and therefore induces a functor $\shv^{H'\mbox{-}\mon}(H)=\shv^{\mon}(H')\otimes_{\shv(H')}\shv(H)\to \shv^{\mon}(H)$. This functor is fully faithful, as it is compatible with embeddings into $\shv(H)$. }
\end{remark}

\begin{lemma}\label{lem: functors between monodromic categories}
Let $f: X\to Y$ be an $H$-equivariant morphism of prestacks. 
\begin{enumerate}
\item\label{lem: functors between monodromic categories-1} We have $f^!\circ \av^{\mon}\cong \av^{\mon}\circ f^!$. In particular, the functor $f^!$ restricts to a functor $f^!: \shv_{\mon}(Y)\to \shv_{\mon}(X)$. 
\item\label{lem: functors between monodromic categories-2} If $f$ is in class $\verti$ as in \eqref{eq: pushforward along non-representable morphisms}, then we have $f_*\circ \av^{H\mbox{-}\mon}\cong \av^{H\mbox{-}\mon}\circ f_*$. In particular, $f_*$ restricts to a functor $f_*:   \shv_{\mon}(X)\to   \shv_{\mon}(Y)$. 
\item\label{lem: functors between monodromic categories-4} If $f$ is a representable coh. pro-smooth morphism, then the above statements hold for $f_\flat$ (as defined in \Cref{lem:etale-proper-functoriality-shv-prestack} \eqref{lem:etale-proper-functoriality-shv-prestack-4}) in place of $f_*$.
\end{enumerate}
There are analogous statements for $\chi$-monodromic categories.
\end{lemma}
\begin{proof}
The point is that the functor $f^!: \shv(Y)\to \shv(X)$ is $\shv(H)$-linear, which in turn follows from the base change and projection formula (encoded by the sheaf theory $\shv$ by \Cref{def-dual-sheaves-on-prestacks-vert-indfp}). 
Then we have the following commutative diagram
\[
\xymatrix{
\shv(H)\otimes_{\shv(H)}\shv(Y) \ar@/^/[rr]^-{ \av^{\mon}\otimes \id} \ar_{\id\otimes f^!}[d]&& \ar@/^/[ll] \shv_{\mon}(H)\otimes_{\shv(H)}\shv(Y) \ar^{\id\otimes f^!}[d]\\
\shv(H)\otimes_{\shv(H)}\shv(X)  \ar@/^/[rr]^-{ \av^{\mon}\otimes \id} && \ar@/^/[ll]\shv_{\mon}(H)\otimes_{\shv(H)}\shv(X),
}\]
which implies Part \eqref{lem: functors between monodromic categories-1}. Part \eqref{lem: functors between monodromic categories-2}-\eqref{lem: functors between monodromic categories-4} follow similarly. (For $f_\flat$, the desired base change and projection formula are supplied by \Cref{cor-additional-base-change-for-shv-theory-on-prestacks}.) 
\end{proof}

\Cref{lem: two-definition of equivariant category} also has the following important consequence. For this, we recall that the $!$-pushforwards are defined for representable pfp morphisms (as the left adjoint of $!$-pullbacks) between sind-very placid stacks and satisfy a base change with respect to  weakly coh. pro-smooth pullbacks (see \Cref{pfp.functors.and.bc.ind-placid.stacks}).

\begin{lemma}\label{lem: functors between monodromic categories-5}
Let $X$ be an sind-very placid stack equipped with an $H$-action. Let $f: X\to H\backslash X$ be the quotient morphism. 
Then when restricted to $\shv_\mon(X)$, we have $f_*=f_![\dim H]$. 
\end{lemma}
\begin{proof}
Using the base change \Cref{pfp.functors.and.bc.ind-placid.stacks}, the functor $f_!$ is the geometric realization of the $!$-pushforwards between the following cosimplicial diagrams  
\[
\shv_{u\mbox{-}\mon}(H)^{\bullet+1}\otimes \shv_{u\mbox{-}\mon}(X)\to  \shv_{u\mbox{-}\mon}(H)^{\bullet}\otimes \shv_{u\mbox{-}\mon}(X).
\]
Then the lemma follows from \Cref{lem: functoriality depth zero geom Langlands for tori}.
\end{proof}

In order to apply our general formalism to compute the categorical trace, we upgrade $X\mapsto \shv((H,\mon)\bs X)$ (for a prestack with an $H$-action) as a sheaf theory as follows.
Let $\shv$ be the sheaf theory as in \eqref{eq: pushforward along non-representable morphisms}, and let $\verti$ be the class of morphisms of prestacks as defined there.
We consider the following category $\bfC$ of pairs $(H,X)$ consisting of a prestack $X$ equipped with an action of a torus $H$ over $k$. 
Note that if $H_1\to H, H_2\to H$ are two maps of groups of tori, then the neutral connected component $(H_1\times_HH_2)^\circ$ of fiber product $H_1\times_HH_2$ (in $\prestk^\pf_k$ so we automatically ignore any derived or non-reduced structure) is a torus. Therefore $\bfC$ admits finite products and the forgetful functor $\bfC\to \prestk^\pf_k$ preserves finite products.
We let $\corr(\bfC)_{\verti;\all}$ be the category consisting of those $(H_1, X_1)\leftarrow (H_2,X_2)\rightarrow (H_3,X_3)$ such that $(X_2\to X_3)\in (\prestk^\pf_k)_{\verti}$.
We have a symmetric monoidal functor 
\[
\corr(\bfC)_{\verti;\all}\to \corr(\prestk^\pf_k)_{\verti;\all},\quad (H,X)\mapsto X.
\]

\begin{proposition}\label{prop: sheaf theory for monodromic sheaf}
The assignment $(H,X)\mapsto \shv((H,\mon)\bs X)$ can be upgraded to a sheaf theory
\[
\shv_{\mon}: \corr(\bfC)_{\verti;\all}\to \lincat_\La,
\] 
which sends  $(H_1, X_1)\xleftarrow{g} (H_2,X_2)\xrightarrow{f} (H_3,X_3)$ to $f^{\mon}_*\circ g^{!}$, where $f^{\mon}_*=\av^{H_3\mbox{-}\mon}\circ f_*$.
In addition, the class $\mathrm{HR}$ of morphisms associated to $\shv_\mon$ as defined in \Cref{rem-additional-base.change.sheaf.theory} \eqref{rem-additional-base.change.sheaf.theory-0} (i.e. the class of morphisms satisfying \Cref{assumptions.base.change.sheaf.theory.H}) contain
those morphisms $g:(H_2,X_2)\to (H_1,X_1)$ with $X_2\to X_1$ being representable coh. pro-smooth morphisms.

There is an analogous unipotent version
\[
\shv_{u\mbox{-}\mon}: \corr(\bfC)_{\verti;\all}\to \lincat_\La.
\]
\end{proposition}
\begin{proof}
We consider
\[
 \corr(\prestk_k)_{\verti;\all}\xrightarrow{\shv} \lincat_\La\to \cat\xrightarrow{(-)^{\op}}\cat,
\] 
where $(-)^{\op}$ is the functor sending a category to its opposite category (e.g. see \cite[Remark 2.4.2.7]{Lurie.higher.algebra}).
By symmetric monoidal version of unstraightening (see \cite[Proposition A.2.1]{Hinich}, see also \Cref{rem: Sheaf theory via Grothendieck construction}), 
this functor is classified by a coCartesian fibration $\bfD\to \corr(\bfC)_{\verti;\all}$, where  $\bfD$ consists $(H,X,\mF)$ with $\mF\in \shv(X)^{\op}$, and a morphism $(H_1,X_1,\mF_1)\to (H_2,X_2,\mF_2)$ consists of $(H_3, X_3)\in \bfC$,
a correspondence 
\begin{equation}\label{prop: sheaf theory for monodromic sheaf-correspondence}
X_1\xleftarrow{g} X_3\xrightarrow{f} X_2,
\end{equation} 
where $f\in \verti$ and $g\in \horiz$, both of which are compatible with torus actions, and a morphism $(a: \mF_2\to f_*(g^!\mF_1))\in \Map_{\shv(X_2)^{\op}}(f_*(g^!\mF_1), \mF_2)$.
The category $\bfD$ is endowed with a symmetric monoidal structure $(H,X,\mF)\otimes (H',X',\mF')=(H\times H', X\times X', \mF\boxtimes_{\La} \mF')$ such that the forgetful functor $\bfD\to  \corr(\prestk_k)_{\verti;\all}$ is symmetric monoidal.

The full subcategory $\bfD_{\mon}$ consisting of those $(H,X,\mF)$ with $\mF\in\shv((H,\mon)\bs X)$ is a full symmetric monoidal category and $\bfD_{\mon}\to \corr(\bfC)_{\verti;\all}$ is a coCartesian fibration. Namely, for every $(H_1,X_1,\mF_1)$ and a correspondence as in \eqref{prop: sheaf theory for monodromic sheaf-correspondence}, the coCartesian arrow above it is given by $\mF_2:=f^{\mon}_*(g^!\mF_1)\xrightarrow{\id} f^{\mon}_*(g^!\mF_1)$.
Now straightening gives $\shv^{\mon}: \corr(\bfC)_{\verti;\all}\to \lincat_\La$ as desired. 

That the class of morphisms as defined in the proposition satisfy \Cref{assumptions.base.change.sheaf.theory.H} directly follows from \Cref{cor-additional-base-change-for-shv-theory-on-prestacks}. 

The unipotent version can be treated similarly.
\end{proof}

\begin{remark}\label{rmk: sheaf theory for monodromic sheaf}
\begin{enumerate}
\item Giving $f: (H,X)\to (H',X')$, if the map $H\to H'$ is surjective, then $f_*^{\mon}=f_*|_{\shv^{\mon}(X)}$. 
\item\label{rmk: sheaf theory for monodromic sheaf-2} Let $\bfC'\subset \bfC$ be the full subcategory consisting of those $(H,X)$ such that $X$ is sind-very placid. We restrict $\shv_\mon$ to $\corr(\bfC')_{\verti;\all}$, and let $\mathrm{VR}$ be the class of morphisms associated to $\shv_\mon|_{\corr(\bfC')_{\verti;\all}}$ as defined  \Cref{rem-additional-base.change.sheaf.theory} \eqref{rem-additional-base.change.sheaf.theory-0} (i.e. the class of morphisms satisfying \Cref{assumptions.base.change.sheaf.theory.V}).
Then a morphism $f: (H,X)\to (H',X')$ with $H\to H'$ surjective and $H\backslash X\to H'\backslash X'$ being ind-pfp proper morphisms of sind-very placid stacks belongs to $\mathrm{VR}$. Namely,  by assumption $\ker(H\to H')\backslash X\to X'$ is ind-pfp proper. Therefore, by \Cref{lem: functors between monodromic categories-5}, up to shifts, the right adjoint of $f_*$ is just $f^!$, which then clearly satisfies \Cref{assumptions.base.change.sheaf.theory.V}. 
\end{enumerate}
\end{remark}

Finally let us record the following two statements. The first will be used in the proof of \Cref{lem: convolving central sheaf commutes with killing monodromy}, and the second will be used in the proof of \Cref{prop: categorical property of monodromic Hecke}.

\begin{lemma}\label{lem: convolve with sheaves on torus}  
Let $X$ be a prestack over $k$ equipped with an action $a: H\times X\to X$. Then for every $\mF\in\shv(H)$ and $\mG\in\shv_{\mon}(X)$, we have $a_*(\mF\boxtimes \mG)\cong a_*(\av^{\mon}(\mF)\boxtimes \mG)$.
\end{lemma}
\begin{proof}
We write $\mG\cong a_*(\widetilde{\Ch}\boxtimes \mG)$ so 
\[
a_*(\mF\boxtimes \mG)\cong a_*(\mF\boxtimes a_*(\widetilde{\Ch}\boxtimes \mG))\cong a_*(m_*(\mF\boxtimes\widetilde{\Ch})\boxtimes \mG)\cong a_*(\av^{\mon}(\mF)\boxtimes \mG),
\]
where the last isomorphism follows from \Cref{lem: monoidal of monodromic}. 
\end{proof}

\begin{lemma}\label{lem: monodromic exterior tensor product}
Let $X$ be a prestack over $k$ acted by a torus $H$, and let $H'$ be a torus, acted by itself via left translation. Then
\[
\boxtimes: \shv_{\mon}(H')\times \shv((H,\mon)\bs X)\to \shv((H'\times H,\mon)\bs H'\times X)
\]
is an equivalence.
\end{lemma}
\begin{proof}
The fully faithfulness holds in general, see \Cref{lem: categorical kunneth prestack}. The essential surjectivity follows from the last part of \Cref{ex: endomorphism of monodromic sheaf} and \Cref{lem: tensor product tensor of monodromic sheaves} .
\end{proof}

\subsubsection{Case of tori}
Let $H$ be an algebraic torus over $k$. In this case, we have an appropriate version of the fully faithful embedding \eqref{eq: character sheaf} at the derived level.

In the sequel unless otherwise specified, we will base change $R_{I_F^t,\hat{H}}$ to the coefficient ring $\La$ (which we assume to be either algebraic over $\bF_\ell$ or $\bQ_\ell$, or finite over $\bZ_\ell$) but omit $\La$ from the notation.
If $f: H_1\to H_2$ is a homomorphism of tori, then it induces an ind-finite morphism 
\[
\hat{f}\colon R_{I_F^t,\hat{H}_2}\to R_{I_F^t,\hat{H}_1}.
\]

Let $\ind\Coh (R_{I_F^t,\hat{H}})$ denote the ind-completion of the category of coherent sheaves on $R_{I_F^t,\hat{H}}$. We endow it with a symmetric monoidal structure given by $!$-tensor product.

\begin{proposition}
\label{lem: depth zero geom Langlands for tori}
There is a natural equivalence of $\La$-linear monoidal categories
\[
\Ch\colon \indcoh (R_{I_F^t,\hat{H}})\cong \shv_\mon(H),
\] 
which is $t$-exact with respect to the standard $t$-structures on the source and target. In particular, the monoidal unit $\cohdual_{R_{I_F^t,\hat{H}}}$ of $\indcoh (R_{I_F^t,\hat{H}})$ corresponds to the monoidal unit $\widetilde{\Ch}=\Ch_\mon$ of $\shv_\mon(H)$.

Let $f: H_1\to H_2$ be a homomorphism of tori.
Under the above equivalences, the adjoint functors $(f^*, f_*^{\mon})$ from \Cref{lem: functoriality depth zero geom Langlands for tori} corresponds to the adjoint functor $(\hat{f}^{\indcoh}_*,\hat{f}^{\indcoh,!})$ between ind-coherent sheaves.
\end{proposition}
\begin{proof}
It is enough to notice that for $\mF_1,\mF_2\in \Coh(R_{I_F^t,\hat{H}})^{\heartsuit}$
we have isomorphisms of (complexes of) $\La$-modules 
\[
\Hom_{\Coh(R_{I_F^t,\hat{H}})}(\mF_1, \mF_2)\cong \Hom_{\shv(H)}(\Ch(\mF_1),\Ch(\mF_2)).
\] 
It is enough to show this for $\mF_i=\Ch_{\chi_i}$ for two character sheaves associated to $\chi_i: \pialg(H)\to (\La')^\times$, where $\La'$ is a finite extension of $\La$. But
\[
\Hom_{\Coh(R_{I_F^t,\hat{H}})}(\mF_1, \mF_2)\cong \Hom_{\pialg(H)}(\chi_1,\chi_2).
\]
So it is enough to show that $\Hom_{\pialg(H)}(\chi_1,\chi_2)=\Hom_{\shv(H)}(\Ch_{\chi_1},\Ch_{\chi_2})$. Using \eqref{eq:character sheaf tensor product}, we reduce to the case $H=\bG_m$. Then the claim is clear.

Therefore, $\Ch$ extends to a $t$-exact fully faithful equivalence $\Coh(R_{I_F^t,\hat{H}})\to \shv_{\mon}(H)^\cpt$, and then a $t$-exact equivalence
$\Ch\colon \ind\Coh (R_{I_F^t,\hat{H}})\to \shv_\mon(H)$ as desired.

It is clear that under the equivalence $f^*$ corresponds to $\hat{f}^{\indcoh}_*$. Then $f_*^{\mon}$ corresponds to $\hat{f}^{\indcoh,!}$.

Note that the multiplication map $m: H\times H\to H$ corresponds to the diagonal map $\Delta:R_{I_F^t,\hat{H}}\to R_{I_F^t,\hat{H}}\times R_{I_F^t,\hat{H}}$. If follows that $m^{\mon}_*=m_*: \shv_\mon(H)\otimes \shv_\mon(H)\to \shv_\mon(H)$ corresponds to $\Delta^{\indcoh,!}$. Therefore, the monoidal structure of $\indcoh(R_{I_F^t,\hat{H}})$ given by $!$-tensor product corresponds to the monoidal structure of $\shv_\mon(H)$ given by convolution. 
\end{proof}

\begin{remark}\label{rem: depth zero geom Langlands for tori}
\begin{enumerate}
\item\label{rem: depth zero geom Langlands for tori-1} The above equivalence can be regarded as a version of Mellin transform, and can also be regarded as a version of tame geometric local Langlands correspondence for tori. 
\item\label{rem: depth zero geom Langlands for tori-2} As mentioned before, upon choosing a topological generator of $I_F^t$, we obtain a closed embedding $R_{I_F^t,\hat{H}}\subset \hat{H}$. Then there is a full embedding $\indcoh(R_{I_F^t,\hat{H}})\to \Qcoh(\hat{H})$. However, the equivalence of \Cref{lem: depth zero geom Langlands for tori} does not extend a direct relation between $\shv(H)$ and $\Qcoh(\hat{H})$.
\end{enumerate}
\end{remark}

\begin{example}\label{ex: char and monodromy local systems}
Now let $\chi\in R_{I_F^t,\hat{H}}(\La)$ be a $\La$-point, regarded as a closed subscheme. Let $\hchi$ be the formal completion of $R_{I_F^t,\hat{H}}$ along $\chi$. Let
$\cohdual_\chi=\mO_\chi$ be the dualizing sheaf of $\chi$, and let $\cohdual_{\hchi}$ be the dualizing sheaf of $\hchi$. Then under the equivalence in \Cref{lem: depth zero geom Langlands for tori},  the following sequence of functors correspond to
\[
\Qcoh(\chi)=\indcoh(\chi)\to \indcoh(\hchi)\subset\indcoh(R_{I_F^t,\hat{H}}),
\]
the first three categories of \eqref{eq: inclusion of character sheaves}. In addition, we have 
We have
\[
\Ch(\cohdual_\chi)\cong \Ch_\chi, \quad \Ch(\cohdual_{\hchi})\cong \Ch_{\hchi}.
\]
Using this, we can give a description of $\Ch_{\hchi}$ as an ind-local system on $H$. We let $\{R_{I_F^t,\hat{H},\chi,\al}\}_\al$ be a cofinal system of thickening of $\chi$ in $R_{I_F^t,\hat{H}}$ with each $R_{I_F^t,\hat{H},\chi,\al}\subset R_{I_F^t,\hat{H}}$ a regular embedding, then
\[
\omega_{\hchi}=\colim_\al \omega_{R_{I_F^t,\hat{H},\chi,\al}},\quad  \mbox{with} \quad \omega_{R_{I_F^t,\hat{H},\chi,\al}}=\Hom_{\La[\xcoch(H)]}(\La[R_{I_F^t,\hat{H},\chi,\al}], \omega_{\hat{H}}).
\]
Here $\cohdual_{\hat{H}}$ denotes the dualizing module of $\hat{H}$ (i.e. the sheaf of top differential forms placed in cohomological degree $-\dim H$).
Note that each $\cohdual_{R_{I_F^t,\hat{H},\chi,\al}}$ belongs to $\Coh(R_{I_F^t,\hat{H}})^{\heartsuit}$, and so does $\omega_{\hchi}$. Therefore
\[
\Ch_{\hchi}=\Ch(\cohdual_{\hchi})=\colim_{\al} \Ch( \omega_{R_{I_F^t,\hat{H},\chi,\al}})
\]
is an ind-local system on $H$.
When $\chi=u$ corresponds to the trivial representation of $I_F^t$, the local system $\Ch_u$ is the constant sheaf $\La$ on $H$, and $\Ch_{\hat{u}}$ is an ind-local system with unipotent monodromy. 
\end{example}

\begin{example}\label{ex: chi-monodromic sheaf}
For generally, if $\chi\subset R_{I_F^t,\hat{H}}$ is a closed subscheme. Write $\chi=\Spec \La'$ with $\La'$ finite over $\La$ so $\chi$ gives a homomorphism $\pialg(H)\to {\La'}^\times$, still denoted by $\chi$.
Then $\Ch(\mO_\chi)=\Ch_\chi$ is the character sheaf associated to $\chi$.

Let $\hchi$  denote the formal completion of $R_{I_F^t,\hat{H}}$ along $\chi$. We let $\shv_{\chi\mbox{-}\mon}(H)\subset \shv_\mon(H)$ denote the full subcategory corresponding to $\indcoh(\hchi)$. 
When $\chi$ is given by $\La$-point, the above discussions reduce to the discussions of $\chi$-monodromic sheaves before.  So for general $\chi$, we still call $\shv_{\chi\mbox{-}\mon}(H)$ the category of $\chi$-monodromic sheaves. For a space $X$ acted by $H$, one can similarly define the category of $\chi$-monodromic sheaves on $X$, which will be denoted by the same notions as before.

As $\indcoh(\hchi)$ is a monoidal category with unit $\cohdual_{\hchi}$, we see that $\shv_{\chi\mbox{-}\mon}(H)$ is monoidal with the unit given by $\Ch(\cohdual_{\hchi})$, denoted by $\Ch_{\hchi}$ or $\Ch_{\chi\mbox{-}\mon}$ as before.
\end{example}

\begin{example}\label{ex: kernel of dual hom}
Let $\varphi: H_1\to H_2$ be a homomorphism, inducing $\hat\varphi: R_{I_F^t, \hat{H}_2}\to R_{I_F^t,\hat{H}_1}$. Let
\[
\chi_\varphi:=\ker\hat\varphi=u\times_{R_{I_F^t,\hat{H}_1}}R_{I_F^t,\hat{H}_2},
\]
where $u\in R_{I_F^t,\hat{H}'}$ corresponds to the trivial representation of $I_F^t$. This is a (possibly) derived sub-indscheme in $R_{I_F^t,\hat{H}_2}$. Let $\cohdual_{\chi_\varphi}$ denote its dualizing sheaf, regarded as an ind-coherent sheaf on $R_{I_F^t,\hat{H}_2}$ via $*$-pushforward. Then
\[
\Ch(\cohdual_{\chi_\varphi})\cong \varphi_*\La.
\]
It follows that $\cohdual_{\chi_\varphi}\in \indcoh(R_{I_F^t,\hat{H}_2})^{\geq 0}$ and it belongs to $\indcoh(R_{I_F^t,\hat{H}_2})^{\heartsuit}$ if the map $\hat{H}_2\to \hat{H}_1$ is surjective.

A particular case we need in the sequel is that $\varphi: \bG_m\to H$ is a non-trivial cocharacter. Then $\hat\varphi: \hat{H}\to \bG_m$ is surjective.  We note that we have a short exact sequence in $\shv_\mon(H)^{\heartsuit}$ 
\begin{equation}\label{eq: injective hull general}
1\to \Ch(\cohdual_{\chi_\varphi})\to \widetilde{\Ch}\to \widetilde{\Ch}\to 1.
\end{equation}
\end{example}

\begin{example}\label{ex: Hn-monodromic}
Now suppose $\varphi: H'\to H$ is a finite \'etale homomorphism of tori. Suppose that $H$ acts on $X$, which induces an action of $H'$. 
In this case $\chi_\varphi$ is in fact a finite (classical) closed subscheme of $R_{I_F^t,\hat{H}}$. Indeed, we have
\[
\chi_\varphi\cong (\ker\varphi)^D:=\Spec (\La[\ker\varphi]).
\] 
We have
\[
\shv\bigl((H',u\mbox{-}\mon)\bs X\bigr)= \shv\bigl((H,\chi_\varphi\mbox{-}\mon)\bs X\bigr)
\] 
as subcategories of $\shv\bigl((H',\mon)\bs X\bigr)=\shv\bigl((H,\mon)\bs X\bigr)$. The left hand side is acted by $\shv_{u\mbox{-}\mon}(H')$ while the right hand side is acted by $\shv_{\chi_\varphi\mbox{-}\mon}(H)$. These two actions are compatible via the push-forward $\varphi_*^{\mon}: \shv_{u\mbox{-}\mon}(H')\to \shv_{\chi_\varphi\mbox{-}\mon}(H)$.
\end{example}

\begin{example}\label{ex: Hn-equivariant}
Suppose we are in the situation as in \Cref{ex: Hn-monodromic}. Then
\[
\shv(H'\backslash X)\cong\Mod_\La\otimes_{\shv_{u\mbox{-}\mon}(H')}\shv\bigl((H',u\mbox{-}\mon)\bs X\bigr)\cong \Mod_\La\otimes_{\shv_{\mon}(H')}\shv\bigl((H,\mon)\bs X\bigr).  
\]
A particular case is when $X=H$ equipped with the left $H$-action. Then $H'\backslash H=\bB (\ker\varphi)$. In this case we recover the equivalence
\begin{equation}\label{eq: shv as qcoh on classifying stack of finite group}
\begin{array}{ll}
           & \shv(\bB (\ker\varphi))\cong \Mod_\La\otimes_{\ind\Coh(R_{I_F^t,\hat{H'}})} \ind\Coh(R_{I_F^t,\hat{H}})\\
\cong& \Mod_\La\otimes_{\ind\Coh(R_{I_F^t,\hat{H'}})} \ind\Coh(R_{I_F^t,\hat{H'}})\otimes_{\qcoh(\hat{H}')}\qcoh(\hat{H})\\
\cong& \qcoh(u\times_{R_{I_F^t,\hat{H'}}}R_{I_F^t,\hat{H}})\cong \qcoh(\chi_\varphi).
\end{array}
\end{equation}
For general $X$, we have
\[
\shv(H'\backslash X)\cong   \qcoh(\chi_\varphi)\otimes_{\shv_{\chi_\varphi\mbox{-}\mon}(H)}\shv((H,\chi_\varphi\mbox{-}\mon)\bs X).
\]
It follows that $\shv(H'\backslash X)$ is acted by $\qcoh(\chi_\varphi)$, and for a $\La$-point $\chi\subset\chi_\varphi$ (i.e. for those $\chi$ such that the pullback of $\Ch_\chi$ to $H'$ becomes trivial), we have
\begin{equation}\label{eq: equivariant category via gerbe}
\shv\bigl((H,\chi)\backslash X\bigr)\cong   \Mod_\La\otimes_{\qcoh(\chi_\varphi)}\shv(H'\backslash X).
\end{equation} 
Note that if $\chi$ is a connected component of $\chi_\varphi$, then $\shv\bigl((H,\chi)\backslash X\bigr)$ is canonically a direct summand of $\shv(H'\backslash X)$, induced by the left adjoint of the natural functor $\shv(H'\backslash X)\to \shv\bigl((H,\chi)\backslash X\bigr)$.
\end{example}

We will also make use of the following statement when studying the affine Deligne-Lusztig induction in \Cref{SS: Aff DL induction}. 

\begin{lemma}\label{lem:DL induction for tori}
Let $\varphi: H'\to H$ be a finite \'etale homomorphism (as in  \Cref{ex: Hn-monodromic}). Let $Z\subset R_{I_F^t,\hat{H}}$ be a closed sub-indscheme. 
\begin{enumerate}
\item Under the equivalence \eqref{eq: shv as qcoh on classifying stack of finite group}, the $*$-pushforward of $\Ch(\cohdual_Z)$ along the map $H\to  H'\backslash H= \bB(\ker \varphi)$ corresponds to $\La[\ker \varphi]$-module given by $\omega_{Z\cap \chi_\varphi}$ (which belongs to $\qcoh(\chi_\varphi)\subset \indcoh(\chi_\varphi)$).
\item If $\cohdual_{Z}\in \indcoh(R_{I_F^t,\hat{H}})^{\geq 0}$, so is $\cohdual_{Z\cap \chi_\varphi}$.
\end{enumerate}
\end{lemma}
\begin{proof}
We need to compute the cohomology of the local system $\varphi^!\Ch(\cohdual_Z)=\varphi^*\Ch(\cohdual_Z)$ on $H'$. On the dual side, $\varphi^*\Ch(\cohdual_Z)$ corresponds to the (ind-)coherent sheaf $\hat{\varphi}^{\indcoh}_*\omega_{Z}$. So the cohomology of $\varphi^*\Ch(\cohdual_Z)$ corresponds to $u^{\indcoh,!}(\hat{\varphi}^{\indcoh}_*\omega_{Z})=\omega_{Z\cap \chi_\varphi}$, as desired.

As $\Ch$ is $t$-exact,  the $*$-pushforward of $\Ch(\cohdual_Z)$ sits in cohomological degree $\geq 0$. This gives the second part. (Of course, it can be proved directly in the coherent side.)
\end{proof}

\subsection{Affine Hecke categories}
The goal of this subsection is to review (and generalize) results about affine Hecke categories needed in the sequel. We fix a coefficient ring $\La$ as before, and unless otherwise specified, all geometric spaces are base changed to $k$.
\subsubsection{Convolution pattern}
In order to rigorously define various Hecke categories equipped with a monoidal structure in the $\infty$-categorical setting, we make use of the the convolution pattern.

Let $\breve\mG$ be an affine smooth model of $G$ over $\mO_{\breve F}$. Then $L^+\breve\mG\subset LG$ is a pfp closed embedding and $LG/L^+\breve\mG$ is an ind-scheme ind-pfp over $k$.
Therefore 
\[
X=\bB L^+\breve\mG\to Y=\bB LG
\] 
is an ind-pfp morphism of sind-very placid stacks. Note that 
$X\times_YX$ is identified with $L^+\breve\mG\bs LG/L^+\breve\mG$ and the relative diagonal $\Delta_{X/Y}:X\to X\times_YX$ is identified with $\Delta: \bB L^+\breve\mG=  L^+\breve\mG\bs  L^+\breve\mG/  L^+\breve\mG\to  L^+\breve\mG\bs LG/  L^+\breve\mG$.
By applying the convolution pattern \Cref{rem:segal.objects.morphisms.vert.horiz} and \Cref{rem-conv-product} to the sheaf theory $\shv$ constructed in \Cref{def-dual-sheaves-on-prestacks-vert-indfp} and to $X\to Y$ as above, we see
that the category $\shv(L^+\breve\mG\backslash LG/L^+\breve\mG)$ admits a canonical monoidal structure, with the monoidal unit given by $\Delta_*\consdual_{\bB L^+\breve\mG}$.
 The monoidal product is usually called the convolution product.

\begin{remark}
Informally, the convolution product is induced by the correspondence
\begin{equation}\label{eq: conv mG}
\xymatrix{
\ar_-{m^{\breve\mG}}[d] L^+\breve\mG\backslash LG\times^{L^+\breve\mG}LG/L^+\breve\mG\ar^-{\eta^{\breve\mG}}[r] &  L^+\breve\mG\backslash LG/L^+\breve\mG\times L^+\breve\mG\backslash LG/L^+\breve\mG.\\
L^+\breve\mG\backslash LG/L^+\breve\mG& 
}
\end{equation}
We will denote the convolution product by $\star^{\breve\mG}:= (m^{\breve\mG})_*\circ (\eta^{\breve\mG})^!$.

The the above monoidal structure at the ordinary categorical level (i.e. for the homotopy category of $\shv(L^+\breve\mG\backslash LG/L^+\breve\mG)$) was originally defined by Lusztig and has been considered in literature for a long time. But it is usually constructed in an ad hoc way (e.g. see \cite{Gaitsgorycycles} and \cite{zhu2016introduction}), which cannot be applied in the $\infty$-categorical setting.
\end{remark}

In addition, if $\breve\mG_1,\breve\mG_2$ are two smooth integral models of $G$ as above, then again by the convolution pattern the category $\shv(L^+\breve\mG_1\backslash LG/L^+\breve\mG_2)$ is a $\shv(L^+\breve\mG_1\backslash LG/L^+\breve\mG_1)\mbox{-}\shv(L^+\breve\mG_2\backslash LG/L^+\breve\mG_2)$-bimodule.
Note that $L^+\breve\mG_1\backslash LG/L^+\breve\mG_2$ is an ind-very placid stack, and all the involved convolution products preserve finitely generated subcategories, we have parallel constructions for $\fgshv$. Finally, by passing to the ind-completion (or by applying the convolution pattern to the sheaf theory $\rshv$ constructed in \Cref{thm-rshv-for-ind-quasi-placid-stack}),
we have parallel constructions for $\rshv$.

We need the following variants to deal with monodromic and equivariant affine Hecke categories. 

First, suppose $L^+\breve\mG$ admits a closed normal subgroup $(L^+\breve\mG)^1$ such that $H=L^+\breve\mG/(L^+\breve\mG)^1$ is a connected affine algebraic group over $k$. Then $\bB (L^+\breve\mG)^1$ is equipped with an action of $H$ such that the further quotient of $\bB (L^+\breve\mG)^1$ by $H$ is $\bB L^+\breve\mG$. We equip $\bB LG$ be with trivial group action.
Then we have a morphism in the category $\bfC'$ as in \Cref{rmk: sheaf theory for monodromic sheaf} \eqref{rmk: sheaf theory for monodromic sheaf-2}. 
\[
X^1=\bB (L^+\breve\mG)^1\to Y=\bB LG,
\]
which is still ind-pfp. Then relative diagonal of this map is still a pfp closed embedding. Then we can apply the sheaf theory from \Cref{prop: sheaf theory for monodromic sheaf} and the convolution pattern to obtain a monoidal category $\shv_{(H\times H)\mbox{-}\mon}((L^+\breve\mG)^1\bs LG/ (L^+\breve\mG)^1)$, with the monoidal unit given by $\av^{\mon}\Delta_*\consdual_{\bB (L^+\breve\mG)^1}$. 

Now suppose there are two pairs $(L^+\breve\mG_i)^1\subset L^+\breve\mG_i$ with $H_i=L^+\breve\mG_i/(L^+\breve\mG_i)^1$ as above, for $i=1,2$. Then the category $\shv_{(H_1\times H_2)\mbox{-}\mon}((L^+\breve\mG_1)^1\bs LG/ (L^+\breve\mG_2)^1)$ has a natural structure as a $\shv_{(H_1\times H_1)\mbox{-}\mon}((L^+\breve\mG_1)^1\bs LG/ (L^+\breve\mG_1)^1)\mbox{-}\shv_{(H_2\times H_2)\mbox{-}\mon}((L^+\breve\mG_2)^1\bs LG/ (L^+\breve\mG_2)^1)$-bimodule.

For another variant, suppose that $\widetilde{L^+\breve\mG}\to L^+\breve\mG$ is a surjective homomorphism obtained by pulling back a finite surjective group homomorphism $\widetilde{L^m\breve\mG}\to L^m\breve\mG$ for some $m$. We suppose the kernel $E$ of the homomorphism $\widetilde{L^+\breve\mG}\to L^+\breve\mG$ is an \'etale group over $k$ of order invertible in $\La$.
We consider 
\begin{equation}\label{eq:f in LG}
\widetilde{X}=\bB \widetilde{L^+\breve\mG}\to  \bB L^+\breve\mG\to Y=\bB LG,
\end{equation}
where the first map is an $E$-gerbe (a.k.a. a $\bB E$-torsor) and the second map is ind-pfp.
On the other hand, we consider
\begin{equation}\label{eq:DeltaX/Y in LG}
\bB \widetilde{L^+\breve\mG}=\widetilde{L^+\breve\mG}\backslash  \widetilde{L^+\breve\mG}/ \widetilde{L^+\breve\mG}\to  \widetilde{L^+\breve\mG}\backslash  L^+\breve\mG/ \widetilde{L^+\breve\mG}\to  \widetilde{L^+\breve\mG}\backslash LG/ \widetilde{L^+\breve\mG}=\bB \widetilde{L^+\breve\mG}\times_{\bB LG}\bB \widetilde{L^+\breve\mG},
\end{equation}
where again the first map is an $E$-gerbe and the second map is a pfp closing embedding. 
It follows that both \eqref{eq:f in LG} and \eqref{eq:DeltaX/Y in LG} are in the class $\verti$ as in \Cref{rem: pushforward along non-representable morphisms}.
By applying the convolution pattern to $\widetilde{\bB L^+\breve\mG}\to \bB LG$ we see that
$\shv(\widetilde{L^+\breve\mG}\backslash LG/\widetilde{L^+\breve\mG})$ have canonical monoidal structure, where the monoidal product is induced by the correspondence \eqref{eq: conv mG} with $L^+\breve\mG$ replaced by $\widetilde{L^+\breve\mG}$. In addition, the unit is given by the $*$-pushforward of the dualizing sheaf on $\bB\widetilde{L^+\breve\mG}$ along \eqref{eq:DeltaX/Y in LG}. We denote it by $\mathbf{1}_{\widetilde{L^+\mG}}$.
If $\breve\mG_1,\breve\mG_2$ are two integral models of $G$, with $\widetilde{L^+\breve\mG}_i\to L^+\mG_i$ as above. 
Then $\shv(\widetilde{L^+\breve\mG}_1\backslash LG/\widetilde{L^+\breve\mG}_2)$ is a $\shv(\widetilde{L^+\breve\mG}_1\backslash LG/\widetilde{L^+\breve\mG}_1)\mbox{-}\shv(\widetilde{L^+\breve\mG}_2\backslash LG/\widetilde{L^+\breve\mG}_2)$-bimodule. 
As before, since $\widetilde{L^+\breve\mG}_1\backslash LG/\widetilde{L^+\breve\mG}_2$ is ind-very placid, and all the convolution products preserve $\fgshv$, we have parallel constructions for $\fgshv$ and then passing to ind-completion for $\rshv$.

\subsubsection{Affine Hecke category}\label{SS: Aff Hecke cat}
We recall the usual affine Hecke category and its variants. 
Let $\breve\mG=\mI$ be the standard Iwahori group scheme (defined over $\mO$), the monoidal category 
\[
\fgshv(\iw\backslash LG/\iw)
\] 
is usually called the affine Hecke category in literature. We shall call it the small unipotent affine Hecke category, and call $\shv(\iw\backslash LG/\iw)$ (resp. $\rshv(\iw\backslash LG/\iw)$) the unipotent affine Hecke category (resp. the big unipotent affine Hecke category).
Following traditional notation, for $w\in \widetilde W$, let
\[
\nabla_w\in\fgshv(\iw\backslash LG/\iw), \quad \mbox{resp. } \  \Delta_w\in\fgshv(\iw\backslash LG/\iw)
\] 
denote the $*$-extension (resp. $!$-extension) of the (shifted) dualizing sheaf $\omega_{\iw\backslash\Gr_w}[-\ell(w)]$ on the Schubert cell $\iw\backslash \Gr_w$. It is well-known that $\fgshv(\iw\backslash LG/\iw)$ is the smallest idempotent complete stable subcategory in $\shv(\iw\backslash LG/\iw)$ generated by $\{\Delta_w\}_{w\in \widetilde{W}}$ or by $\{\nabla_w\}_{w\in \widetilde{W}}$. It is also known that $\Delta_w$ is invertible for the monoidal product of $\fgshv(\iw\backslash LG/\iw)$, with an inverse given by $\nabla_{w^{-1}}$.

We need the following categorical properties of affine Hecke categories. They are analogous to  \Cref{compact.objects.of.B(G)} and \Cref{cor: shrek-restr-reserv-cpt}, but are considerably simpler (and are well-known).
\begin{proposition}\label{prop: categorical property of affine Hecke}
\begin{enumerate}
\item\label{prop: categorical property of affine Hecke-1} The category $\shv(\iw\backslash LG/\iw)$ is compactly generated.
An object $\mF\in \shv(\iw\backslash LG/\iw)$ is compact if and only if
$(i_w)^*\mF\in \shv(\iw\bs LG_w/\iw)\cong \shv(\bB \mS_k)$ is compact for every $w$ and and $(i_w)^*\mF=0$ for all but finitely many $w$s, if and only if $(i_w)^!\mF\in \shv(\iw\bs LG_w/\iw)\cong \shv(\bB \mS_k)$ is compact for every $w$ and $(i_w)^!\mF=0$ for all but finitely many $w$s.

The monoidal structure is semi-rigid. For every prestack $X$ over $k$, the functor 
\[
\shv(\iw\backslash LG/\iw)\otimes_\La\shv(X)\to \shv(\iw\backslash LG/\iw\times X)
\] 
is an equivalence.
\item\label{prop: categorical property of affine Hecke-2} The category $\rshv(\iw\backslash LG/\iw)$ is compactly generated. An object $\mF\in \rshv(\iw\backslash LG/\iw)$ is compact if and only if
$(i_w)^*\mF\in \rshv(\iw\bs LG_w/\iw)\cong \rshv(\bB \mS_k)$ is constructible for every $w$ and and $(i_w)^*\mF=0$ for all but finitely many $w$s, if and only if $(i_w)^!\mF\in \rshv(\iw\bs LG_w/\iw)\cong \rshv(\bB \mS_k)$ is constructible for every $w$ and $(i_w)^!\mF=0$ for all but finitely many $w$s.

The monoidal structure is rigid. For every quasi-compact placid stack $X$, the exterior tensor functor 
\[
\rshv(\iw\backslash LG/\iw)\otimes_\La \rshv(X)\to \rshv(\iw\backslash LG/\iw\times X)
\] 
is an equivalence.
\end{enumerate}
\end{proposition}
\begin{proof}
As $\shv(\iw\backslash LG/\iw)=\colim_{w} \shv(\iw\backslash LG_{\leq w}/\iw)$, and each $\iw\backslash LG_{\leq w}/\iw$ is very placid, we see that $\shv(\iw\backslash LG/\iw)$ is compactly generated by \Cref{compact.generation.admissible.stacks}. 
The characterization of compact objects in $\shv(\iw\bs LG/\iw)$ follows from the same arguments as in  \Cref{compact.objects.of.B(G)} and \Cref{cor: shrek-restr-reserv-cpt}.

To show that it is semi-rigid, we apply \Cref{prop-comparison-usual-trace-geo-trac-presemirigid} to $X_1=X\times_YX$, where $X=\bB \iw\to Y=\bB LG$. Note that $X\to Y$ is ind-pfp proper, $X\xrightarrow{\Delta_{X/Y}} X\times_YX$ is a pfp closed embedding, and $X\to X\times X$ is coh. pro-smooth. Thanks to \Cref{cor-additional-base-change-for-shv-theory-on-prestacks}, \Cref{prop-comparison-usual-trace-geo-trac-presemirigid} is applicable, showing that the convolution admits a bilinear right adjoint. Together with compact generation (which in particular implies dualizability), we deduce the semi-rigidity of $\shv(\iw\bs LG/\iw)$.

As argued in \Cref{lem: shvkot(G) tensor product}, the last statement of Part \eqref{prop: categorical property of affine Hecke-1} reduces to show
$\shv(\iw\backslash LG_w/\iw)\otimes\shv(X)\to \shv(\iw\backslash LG_w/\iw\times X)$ is an equivalence. As $\iw\backslash LG_w/\iw\cong \bB H_w$ where $H_w=\iw\cap \dot{w}\iw\dot{w}^{-1}$ (for a lifting $\dot{w}$ of $w$), we can apply \Cref{cor: ess. pro-unip tensor product} to conclude. 

Same arguments apply to $\rshv(\iw\backslash LG/\iw)$. It is in addition a rigid monoidal category, since the unit belongs to $\fgshv(\iw\backslash LG/\iw)$ and therefore is compact. 
Same arguments together with \Cref{cor: fg shv exterior tensor} also implies the equivalence of exterior tensor product functor.
\end{proof}

By  \Cref{cor: rigidity of der(X)} $\rshv(\iw\bs LG/\iw)$, equipped with the convolution product as above, admits a Frobenius structure
\[
\Hom_{\rshv(\iw\bs LG/\iw)}(\Delta_e,-): \rshv(\iw\bs LG/\iw)\to\Mod_\La.
\]
We let $\verd_{\iw\bs LG/\iw}^{\mathrm{sr}}$ denote the self-duality of $\rshv(\iw\bs LG/\iw)$ induced by this Frobenius algebra structure. See \Cref{ex: self-duality of cpt gen rigid monoidal cat}. 

On the other hand, the category $\rshv(\iw\bs LG/\iw)$ is also equipped with a symmetric monoidal product given by the $!$-tensor product $\os$. By See \Cref{prop:verdier.duality.ind.placid.stacks}, upon a choice of a generalized constant sheaf $\La^{\eta}_{\iw\backslash LG/\iw}$, it also admits a  Frobenius structure given by
\[
\rg^\eta_{\ind\fg}(\iw\bs LG/\iw,-): \rshv(\iw\bs LG/\iw)\to \Mod_\La.
\] 
The induced self-duality of $\rshv(\iw\bs LG/\iw)$ is denoted $(\verd^{\eta}_{\iw\bs LG/\iw})^{\ind\fg}$. 
If we let $\eta=\can$, so $\La^{\can}_{\iw\bs LG/\iw}$ is the canonical generalized constant sheaf on $\iw\bs LG/\iw$ given by the compatible system of generalized constant sheaves $\{\La_{\iw\bs LG_{\leq w}/\iw}\}_{w\in \widetilde{W}}$ (see \Cref{SS: coh. duality. Kot. stack}), then $(\verd^{\can}_{\iw\bs LG/\iw})^{\ind\fg}$ gives what people usually call the Verdier duality on $\iw\bs LG/\iw$. Namely, when restricted to the subcategory $\fgshv(\iw\bs LG/\iw)$ of compact objects, the functor $(\verd^{\can}_{\iw\bs LG/\iw})^{\fg}$ is the one interchanging $\Delta_w$ and $\nabla_w$.

Now let $\sw: \iw\bs LG/\iw\to \iw\bs LG/\iw$ be the involution induced by $LG\to LG, \ g\mapsto g^{-1}$. To simplify the notation, we write $\sw$ instead of $\sw^{\ind\fg,!}$, which is an involution of $\rshv(\iw\bs LG/\iw)$.
The following lemma is well-known (e.g. see \cite[\textsection{3.2}]{Zhu2016}).
\begin{lemma}\label{lem: comparison of two duality for affine Hecke}
We have 
\[
(\verd^{\can}_{\iw\bs LG/\iw})^{\ind\fg} \cong \sw\circ \verd_{\iw\bs LG/\iw}^{\mathrm{sr}}.
\] 
In particular, if $\mF\in\fgshv(\iw\bs LG/\iw)$, we have
\[
(\verd^{\can}_{\iw\bs LG/\iw})^{\fg}(\mF)\cong \sw(\mF^\vee).
\]
\end{lemma}
\begin{proof}
We specialize the discussion in \Cref{rem: pivotal structure} to $\rshv(\iw\bs LG/\iw)$. Then $\sw$ there is the automorphism of $\iw\bs LG/\iw$ induced by the morphism $g\mapsto g^{-1}$, and therefore coincides with the map denoted by the same notation here. 
We thus have
\[
\Hom_{\rshv(\iw\bs LG/\iw)}(\Delta_e, \mF\star \mG)=\Hom_{\rshv(\bB \iw)}(\consdual_{\bB \iw}, \pr^{\ind\fg}_*(\mF\os \mathrm{sw}(\mG))).
\]
Notice that
\[
\rg^{\can}_{\ind\fg}(\iw\bs LG/\iw,-)=\Hom(\consdual_{\bB \iw}, \pr^{\ind\fg}_*(-)).
\]
Therefore, 
\[
\Hom_{\rshv(\iw\bs LG/\iw)}(\Delta_e, \mF\star \mG)=\rg^{\can}_{\ind\fg}(\iw\bs LG/\iw,\mF\os \mathrm{sw}(\mG)).
\]
This gives the first isomorphism.

As explained in \Cref{ex: self-duality of cpt gen rigid monoidal cat}, the self-duality $\verd_{\iw\bs LG/\iw}^{\mathrm{sr}}$, when restricted to the subcategory of compact objects, just sends $\mF$ to its right dual $\mF^{\vee}$ (with respect to the convolution monoidal structure). The second isomorphism follows.
\end{proof}

\begin{remark}
We also note that by \Cref{cor: semi rigidity of der(X)}, the restriction of $\Hom_{\rshv(\iw\bs LG/\iw)}(\Delta_e,-)$ to $\shv(\iw\bs LG/\iw)^\cpt$ followed by ind-extension gives the Frobenius structure of $\shv(\iw\bs LG/\iw)$. Similarly, this is the case for $\rg^{\can}_{\ind\fg}$ as well. It follows that the analogous statement of \Cref{lem: comparison of two duality for affine Hecke} holds for $\shv(\iw\bs LG/\iw)$ as well.
\end{remark}

\subsubsection{Monodromic affine Hecke categories}\label{SSS: mono Hecke cat} 
We next turn to monodromic affine Hecke category. We will prove a few basic results for the monodromic affine Hecke category, parallel to those familiar ones for the (small) unipotent affine Hecke categories.

Let $\iw^u=\ker(\iw\to \mS_k)$ be the pro-unipotent radical of $\iw$.
Note that the quotient map $L^+\mS\to \mS_k$ admits a unique splitting given by the maximal torus of $L^+\mS$. Therefore, we have the semi-direct product decomposition $\iw=\mS_k\cdot \iw^u$.

We consider the $(\mS_k\times \mS_k)$-action on $\iw^u\backslash LG/\iw^u$ by left and right multiplication, and form the corresponding monodromic category. Note that as subcategories of $\shv(\iw^u\backslash LG/\iw^u)$,
it coincides with the monodromic category arising from either the left or right $\mS_k$-action on $\iw^u\backslash LG/\iw^u$. Therefore we can use $\shv_{\mon}(\iw^u\backslash LG/\iw^u)$ to denote this category. For each element $w\in\widetilde{W}$, we have similarly defined monodromic categories with respect to the $(\mS_k\times \mS_k)$-action on $\iw^u\backslash LG_{w}/\iw^u$.  

Note that  for a lifting $\dot{w}$  of $w$ in $N_G(S)(\breve F)$, we have the map
\begin{equation}\label{eq: projection w}
\pr_{\dot{w}}: LG_w\cong \iw^u \cdot \mS_k \cdot \dot{w} \times^{\iw^u\cap \Ad_{\dot{w}^{-1}} \iw^u}  \iw^u \to \mS_k.
\end{equation}
which induces a map $\iw^u\backslash LG_w/\iw^u\to \mS_k$ still denoted by $\pr_{\dot{w}}$. Then the functor
\begin{equation}\label{eq: mon sheaf on one strata}
(\pr_{\dot{w}})^![-\ell(w)] \colon \shv_{\mon}(\mS_k)\to \shv_{\mon}(\iw^u\backslash LG_w/\iw^u)
\end{equation}
is a $t$-exact equivalence of categories. Here the $t$-structure on $\shv_{\mon}(\mS_k)$ is the standard one as in \Cref{lem: depth zero geom Langlands for tori}, and the $t$-structure on  $\shv_{\mon}(\iw^u\backslash LG_w/\iw^u)$ is the perverse $t$-structure defined by the following generalized constant sheaf
\begin{equation}\label{eq: canonical constant sheaf on enhanced hecke}
\La^{\can}_{\iw^u\bs LG/\iw^u},
\end{equation} 
whose $!$-pullback to  $LG_{w}/\iw^u$  is the $(-2\dim \mS_k)$-shift of the usual constant sheaf $\La_{LG_{\leq w}/\iw^u}\in \cshv(LG_{w}/\iw^u)$.  
Comparing with \eqref{eq: can gen constant sheaf}, we see that under the natural projection $\iw^u\bs LG/\iw^u\to \iw\bs LG/\iw^u$, $\La^{\can}_{\iw^u\bs LG/\iw^u}$ is isomorphic to the $-(4\dim \mS_k)$-shift of $!$-pullback of $\La^{\can}_{\iw\bs LG/\iw}$. 

On the other hand, the locally closed embedding $i_w: LG_w\to LG$ induces functors
\[
(i_w)_*, (i_w)_!: \shv_{\mon}(\iw^u\backslash LG_w/\iw^u)\to \shv_{\mon}(\iw^u\backslash LG/\iw^u).
\]
Composing with \eqref{eq: mon sheaf on one strata}, we thus obtain two functors
\begin{equation}\label{eq: (co)standard functor}
\Delta^\mon_{\dot{w}},\ \nabla^\mon_{\dot{w}}:   \shv_{\mon}(\mS_k)\to \shv_{\mon}(\iw^u\backslash LG/\iw^u)
\end{equation}
defined as
\[
\Delta^{\mon}_{\dot{w}}(\mL):=(i_w)_!((\pr_{\dot{w}})^!\mL)[-\ell(w)], \quad \nabla^{\mon}_{\dot{w}}(\mL):=(i_w)_*((\pr_{\dot{w}})^!\mL)[-\ell(w)], \quad \mL\in \shv_{\mon}(\mS_k).
\]
In particular, we write
\[
\widetilde{\Delta}^{\mon}_{\dot{w}}=\Delta^{\mon}_{\dot{w}}(\widetilde{\Ch}), \quad \widetilde{\nabla}^{\mon}_{\dot{w}}=\nabla^{\mon}_{\dot{w}}(\widetilde{\Ch}).
\]
They are called cofree monodromic standard and costandard objects.

For every closed subscheme $\chi\subset R_{I_F^t,\hat{S}}$, we let $\hchi$ be the formal completion of $\chi$ in $R_{I_F^t,\hat{S}}$ and write
\[
\Delta_{\dot{w},\hchi}^{\mon}=\Delta_{\dot{w}}^{\mon}(\Ch_{\hchi}),\quad  \nabla_{\dot{w},\hchi}^{\mon}=\nabla_{\dot{w}}^{\mon}(\Ch_{\hchi}),
\] 
where we recall $\Ch_{\hchi}=\Ch(\cohdual_{\hchi})$. They are also called cofree $\hchi$-monodromic standard and costandard objects.

We remark all the above functors depend on the lifting $\dot{w}$ of $w$. 

Given $\hchi, \hchi'\subset R_{I_F^t, \hat{S}}$ being formal completions of $\chi,\chi'$, we have the corresponding monodromic category $\shv_{(\chi,\chi')\mbox{-}\mon}(\iw^u\backslash LG/\iw^u)$, also denoted as $\shv\bigl((\iw,\hchi)\bs LG/(\iw, \hchi')\bigr)$.
Note that 
\begin{equation}\label{eq: non-zero monodromic category}
\shv\bigl((\iw,\hchi)\bs LG_w/(\iw, \hchi')\bigr)=0, \quad \mbox{ if } \chi\cap w\chi'=\emptyset,
\end{equation} 
and $(i_w)_*, (i_w)_!$ preserve $(\chi,\chi')$-monodromic subcategories.
Note that 
\[
\Delta_{\dot{w},\hchi}^{\mon},\quad \nabla_{\dot{w},\hchi}^{\mon}\in \shv\bigl((\iw,\hchi)\bs LG/(\iw, w^{-1}\hchi')\bigr).
\]

We will use 
\[
\star^u:=(m^u)_*\circ (\eta^u)^!
\] 
denote the monoidal structure on $\shv(\iw^u\backslash LG/\iw^u)$ defined as above (via \eqref{eq: conv mG} for $L^+\mG=\iw^u$). 
More generally, we will use $\star^u$ to denote any morphism induced by the multiplication map 
\[
m^{u}: LG\times^{\iw^u}LG\to LG.
\] 
By \Cref{lem: functors between monodromic categories}, $\shv_{\mon}(\iw^u\backslash LG/\iw^u)\subset \shv(\iw^u\backslash LG/\iw^u)$ is closed under the monoidal product. Then by \Cref{lem: complement on monoidal categories} itself is a monoidal category with the unit given by $\widetilde\Delta_e^{\mon}=\widetilde\nabla_e^{\mon}$. We call it the monodromic affine Hecke category. Alternatively, we may apply \Cref{prop: sheaf theory for monodromic sheaf} to obtain the desired monoidal structure of $\shv_{\mon}(\iw^u\backslash LG/\iw^u)$. Namely, $(\iw^u\backslash LG/\iw^u, \mS_k\times\mS_k)$ is an algebra object in the category $\corr(\bfC)_{\verti;\horiz}$ associated to the \v{C}ech nerve of $(\mS_k, \bB \iw^u)\to (\{1\}, \bB LG)$ (via \Cref{thm:gaitsgory-rozen-segal-monad}). Then applying the sheaf theory $\shv_{\mon}$.

Note that for $\chi_i\subset R_{I_F^t,\hat{S}}, \ i=1,\ldots,4$, we have
\[
-\star^u-: \shv\bigl((\iw,\hchi_1)\bs LG_w/(\iw, \hchi_2)\bigr)\otimes_\La \shv\bigl((\iw,\hchi_3)\bs LG_w/(\iw, \hchi_4)\bigr)\to  \shv\bigl((\iw,\hchi_1)\bs LG_w/(\iw, \hchi_4)\bigr).
\]
In addition
\begin{equation}\label{eq: non-zero tensor product monodromic}
-\star^u-=0, \quad \mbox{ if } \chi_2\cap \chi_3= \emptyset.
\end{equation}
It follows that the subcategory 
\[
\shv\bigl((\iw,\hchi)\bs LG_w/(\iw, \hchi)\bigr)
\] 
is closed under monoidal product with the unit $\Delta_{e,\hchi}^{\mon}=\nabla_{e,\hchi}^{\mon}$. We call it the $\chi$-monodromic affine Hecke category. In particular, $\shv\bigl((\iw,\hat{u})\bs LG_w/(\iw, \hat{u})\bigr)$ is called the unipotent monodromic affine Hecke category. Note that the category $\shv\bigl((\iw,\hchi)\bs LG_w/(\iw, \hchi')\bigr))$ has a natural structure as a $\shv\bigl((\iw,\hchi)\bs LG_w/(\iw, \hchi)\bigr)\mbox{-}\shv\bigl((\iw,\hchi')\bs LG_w/(\iw, \hchi')\bigr)$-bimodule.

The following statement is completely parallel to \Cref{prop: categorical property of affine Hecke}.
\begin{proposition}\label{prop: categorical property of monodromic Hecke} 
The category $\shv_{\mon}(\iw^u\backslash LG/\iw^u)$ is compactly generated and the monoidal structure is semi-rigid.
We have
\[
\shv_{\mon}(\iw^u\backslash LG/\iw^u)^{\cpt}=\shv_{\mon}(\iw^u\backslash LG/\iw^u)\cap \shv(\iw^u\backslash LG/\iw^u)^{\cpt}.
\]
For every prestack $X$ equipped with an action by an algebraic group $H$, the natural functor 
\[
\shv_{\mon}(\iw^u\backslash LG/\iw^u)\otimes_\La \shv_\mon(X)\to \shv_{\mon}(\iw^u\backslash LG/\iw^u\times X)
\]
is an equivalence.
\end{proposition}
\begin{proof}
The proof is also completely parallel to \Cref{prop: categorical property of affine Hecke}. For the last statement, we reduce to show that 
\[
\shv_{\mon}(\iw^u\backslash LG_w/\iw^u)\otimes_\La \shv_{\mon}(X)\to \shv_{\mon}(\iw^u\backslash LG_w/\iw^u\times X)
\] 
is an equivalence, which follows from \Cref{lem: monodromic exterior tensor product}.
\end{proof}

Our next goal is to discuss the analogue of \Cref{lem: comparison of two duality for affine Hecke} for monodromic affine Hecke categories. First, as $\shv_{\mon}(\iw^u\bs LG/\iw^u)$ is semi-rigid, we have a self-duality as in  \Cref{ex: self-duality of cpt gen rigid monoidal cat}, which will be denoted as $\verd^{\mathrm{sr}}_{\iw^u\bs LG/\iw^u}$.

On the other hand, recall our choice of  the canonical generalized constant sheaf $\La^{\can}_{\iw^u\bs LG/\iw^u}$ as in \eqref{eq: canonical constant sheaf on enhanced hecke}, which induces a canonical self duality $\verd^{\can}_{\iw^u\bs LG/\iw^u}$ of $\shv(\iw^u\bs LG/\iw^u)$. As usual, let $(\verd^{\can}_{\iw^u\bs LG/\iw^u})^\cpt$ denote its restriction to compact objects. It further restricts to an equivalence
\[
(\verd^{\mon,\can}_{\iw^u\bs LG/\iw^u})^\cpt:  (\shv_{\mon}(\iw^u\bs LG/\iw^u)^\cpt)^{\op}\to \shv_{\mon}(\iw^u\bs LG/\iw^u)^{\cpt}.
\]
Let us denote its ind-completion by
\[ 
\verd^{\mon,\can}_{\iw^u\bs LG/\iw^u}: \shv_{\mon}(\iw^u\bs LG/\iw^u)^\vee \to \shv_{\mon}(\iw^u\bs LG/\iw^u).
\]

We write $\sw$ instead of $\sw^!$ for the involution of $\shv_\mon(\iw^u\bs LG/\iw^u)$ induced by $\sw: LG\to LG, \ g\mapsto g^{-1}$.

\begin{lemma}\label{lem: comparison of two duality for affine Hecke monodromic}
We have
\[
\verd^{\mon,\can}_{\iw^u\bs LG/\iw^u}[\dim \mS_k]\cong \sw\circ \verd^{\mathrm{sr}}_{\iw^u\bs LG/\iw^u}.
\]
Concretely, if $\mF\in \shv_{\mon}(\iw^u\bs LG/\iw^u)^\cpt$ is compact, then there is a canonical isomorphism
\[
(\verd^{\mon,\can}_{\iw^u\bs LG/\iw^u})^\cpt(\mF)[\dim \mS_k]\cong \sw(\mF^\vee).
\]
\end{lemma}
Although the general discussion as in \Cref{rem: pivotal structure}  does not directly apply to the current situation, the basic idea is similar.

\begin{proof}
It is enough to show that for $\mG\in \shv_{\mon}(\iw^u\bs LG/\iw^u)^\cpt$, we have a canonical isomorphism
\[
\Hom(\mG, (\verd^{\mon,\can}_{\iw^u\bs LG/\iw^u})^\cpt(\mF))\cong \Hom(\mG, \sw(\mF)^\vee).
\]
Recall the tensor product of $\fgshv(\iw^u\bs LG/\iw^u)$ associated to $\La^{\can}_{\iw^u\bs LG/\iw^u}$ as in \Cref{rem: equiv between cons-sheaf and cons-cosheaf-2}. Then by \eqref{eq: Verdier duality for quasi-placid-characterization-2}, we have
\[
\Hom(\mG, (\verd^{\mon,\can}_{\iw^u\bs LG/\iw^u})^\cpt(\mF))\cong \Hom(\mF\otimes^{\can}\mG, \consdual_{\iw^u\bs LG/\iw^u}).
\]
On the other hand, we have
\[
\Hom(\mG, \sw(\mF)^\vee)\cong \Hom(\sw(\mF)\star^u \mG, \widetilde\Delta_e^{\mon}).
\]
Therefore, the desired statement is a consequence of the following lemma.
\end{proof}

\begin{lemma}\label{lem: semirigid dual vs verdier dual in monodromic affine Hecke}
For $\mF,\mG\in \shv_\mon(\iw^u\bs LG/\iw^u)$, there is a canonical isomorphism
\[
 \Hom(\sw(\mF)\star^u \mG, \widetilde\Delta_e^{\mon})\cong \Hom(\mF\otimes^{\can}\mG, \consdual_{\iw^u\bs LG/\iw^u})[\dim \mS_k],
\]
functorial in $\mF,\mG$.
\end{lemma}
\begin{proof}
We may assume that $\mF$ and $\mG$ are compact. Note that we have
\[
\Hom_{\shv_{\mon}(\iw^u\bs LG/\iw^u)}(\sw(\mF)\star^u \mG, \widetilde\Delta_e^{\mon})=\Hom_{\shv(\iw^u\bs LG/\iw^u)}(\sw(\mF)\star^u \mG, \consdual_{\iw^u\bs \iw^u/\iw^u}).
\]
Now consider the diagram
\[
\xymatrix{
& \iw^u\bs LG/\iw^u\ar^{\pr}[r]\ar_{g\mapsto (g^{-1},g)}^{i}[d]& \iw^u\bs \iw^u/\iw^u\ar^{i_e}[d] \\
\iw^u\bs LG/\iw^u\times \iw^u\bs LG/\iw^u &\ar_-{\eta^u}[l] \iw^u\bs LG\times^{\iw^u} LG/\iw^u \ar^-{m^u}[r] & \iw^u\bs LG/\iw^u
}\]
with the commutative square Cartesian.
Recall that the inclusion $i_e: \iw^u\bs \iw^u/\iw^u\to \iw^u\bs LG/\iw^u$ is a pfp closed embedding so $(i_e)^*$ exists, and $(m^u)_*=(m^u)_![\dim \mS_k]$ for monodromic sheaves. Then by base change we have
\begin{eqnarray*}
\Hom(\sw(\mF)\star^u \mG, \consdual_{\iw^u\bs \iw^u/\iw^u}) &= &\Hom(\pr_*(i^*((\eta^u)^!( \sw(\mF)\boxtimes \mG))), \consdual_{\iw^u\bs \iw^u/\iw^u})\\
&=& \Hom(i^*((\eta^u)^!(\sw(\mF)\boxtimes \mG)), \consdual_{\iw^u\bs LG/\iw^u}) [-\dim \mS_k].
\end{eqnarray*}
As $\eta^u$ is coh. pro-smooth, as in \Cref{ex-!-pullback-prosm-generalized-constant} we can endow $\iw^u\bs LG\times^{\iw^u} LG/\iw^u$ with a generalized constant sheaf 
\[
\La^{\can}_{\iw^u\bs LG\times^{\iw^u} LG/\iw^u}:=(\eta^u)^!(\La^{\can}_{\iw^u\bs LG/\iw^u}\boxtimes_{\La}\La^{\can}_{\iw^u\bs LG/\iw^u}).
\] 
Recall that by \Cref{compact.generation.admissible.stacks}, we have $\fgshv(\iw^u\bs LG/\iw^u)=\shv(\iw^u\bs LG/\iw^u)^\cpt$, and similarly we have $\fgshv(\iw^u\bs LG\times^{\iw^u} LG/\iw^u)=\shv(\iw^u\bs LG\times^{\iw^u} LG/\iw^u)^\cpt$. 
Then by \Cref{star.induced.ind.fin.pres.} \eqref{star.induced.ind.fin.pres.-2}, we have
\[
(\eta^u)^!((\verd_{\iw^u\bs LG/\iw^u}^{\can})^\cpt\boxtimes (\verd_{\iw^u\bs LG/\iw^u}^{\can})^\cpt)=(\verd^{\can}_{\iw^u\bs LG\times^{\iw^u} LG/\iw^u})^\cpt\circ (\eta^u)^!.
\]

As $i$ is pfp closed embedding, we have the $*$-pullback of $\La^{\can}_{\iw^u\bs LG\times^{\iw^u} LG/\iw^u}$ along $i$, see \Cref{ex-*-pullback-indfp-generalized-constant}. We note that
\[
i^*(\La^{\can}_{\iw^u\bs LG\times^{\iw^u} LG/\iw^u})\cong \La^{\can}_{\iw^u\backslash LG/\iw^u}[-2\dim \mS_k].
\]
Then by \Cref{star.induced.ind.fin.pres.} \eqref{star.induced.ind.fin.pres.-1}, we have 
\[
i^*\cong (\verd^{\can}_{\iw^u\bs LG/\iw^u})^\cpt\circ i^!\circ (\verd^{\can}_{\iw^u\bs LG\times^{\iw^u}LG/\iw^u})^\cpt[-2\dim \mS_k].
\]
Therefore, we have
\begin{eqnarray*}
 \Hom(i^*((\eta^u)^!(\sw(\mF)\boxtimes\mG)), \consdual_{\iw^u\bs LG/\iw^u}) [-\dim \mS_k] & =& \Hom(\mF\otimes^{\can}\mG, \consdual_{\iw^u\bs LG/\iw^u})[\dim \mS_k].
 \end{eqnarray*}
Putting things together gives what we need.
\end{proof}

To state the next result, note that $\widetilde{W}$ acts on $\mS_k$ through $\widetilde{W}\to W_0$ by adjoint action. For $w\in \widetilde{W}$ and $\mL\in \shv_{\mon}(\mS_k)$, we let $w(\mL):= (\Ad_w)_*\mL$.

\begin{proposition}\label{lem: monoidal product computation in affine Hecke category}
For $w\in \widetilde{W}$, let $\dot{w}$ denote a lifting of it to $N_G(S)(\breve F)$. Let $\mL,\mL'\in\shv_{\mon}(\mS_k)$.
\begin{enumerate}
\item\label{lem: monoidal product computation in affine Hecke category-1} For $w,v\in \widetilde{W}$ satisfying $\ell(w)+\ell(v)=\ell(wv)$, we have canonical isomorphisms
\[
\Delta^{\mon}_{\dot{w}}(\mL)\star^u \Delta^{\mon}_{\dot{v}}(\mL')\cong \Delta^{\mon}_{\dot{w}\dot{v}}(\mL\star w(\mL')),\quad  \nabla^{\mon}_{\dot{w}}(\mL)\star^u \nabla^{\mon}_{\dot{v}}(\mL')\cong \nabla^{\mon}_{\dot{w}\dot{v}}(\mL\star w(\mL')).
\]

\item\label{lem: monoidal product computation in affine Hecke category-3} For every $w\in \widetilde{W}$, we have canonical isomorphisms
\[
 \nabla^{\mon}_{\dot{w}}(\mL)\star^u\Delta^{\mon}_{\dot{w}^{-1}}( \mL')\cong \Delta^{\mon}_{\dot{w}}(\mL)\star^u\nabla^{\mon}_{\dot{w}^{-1}}( \mL')\cong \Delta_e^{\mon}(\mL\star w(\mL')).
\]
\end{enumerate}

Now let $s$ be a simple reflection in $\widetilde{W}$. Let $\hat\al_s$ be the vector part of the affine simple coroot corresponding to $s$, regarded as a cocharacter $\bG_m\to \mS_k$. Let $\widetilde{\Ch}_s=\Ch(\cohdual_{\chi_{\hat{\al}_s}})\in \shv_{\mon}(\mS_k)$, where $\chi_{\hat{\al}_s}$ is as in \Cref{ex: kernel of dual hom}.

\begin{enumerate}[resume]
\item\label{lem: monoidal product computation in affine Hecke category-2}
For $s$ a simple reflection in $\widetilde{W}$, we have fiber sequences
\[
 \nabla^{\mon}_{e}(\mL\star s(\mL')) \to  \nabla^{\mon}_{\dot{s}}(\mL)\star^u \nabla^{\mon}_{\dot{s}}(\mL')\to  \nabla^{\mon}_{\dot{s}}(\mL\star \widetilde{\Ch}_{s}\star s(\mL'))[1], 
\]
\[
\Delta^{\mon}_{\dot{s}}(\mL\star \widetilde{\Ch}_{s}\star s(\mL'))\to \Delta^{\mon}_{\dot{s}}(\mL)\star^u \Delta^{\mon}_{\dot{s}}(\mL')\to \Delta^{\mon}_{e}(\mL\star s(\mL'))
\]
\end{enumerate}
\end{proposition}
\begin{proof}
This is a generalization of the well-known corresponding statements for the usual (small unipotent) affine Hecke category (in equal characteristic), e.g. see \cite[Lemma 8]{arkhipov2009perverse}.  
We write down a proof to illustrate where extra cares are needed.
 
We need the following simple observation: for an action $a: H\times X\to X$ of the torus on a prestack inducing an isomorphism $\tilde{a}: H\times X\to X\times H, \ (h,x)\mapsto (hx, h)$,
there is a canonical isomorphism 
\begin{equation}\label{eq: swap H and X}
\tilde{a}_*(\mL\boxtimes \omega_X)\cong \omega_X\boxtimes \mL
\end{equation} 
for any $\mL\in \shv(H)$.

We write $\eta^{u}_{w,v}$ for the base change of $\eta^u$ along $i_w\times i_v$, and $i_{w,v}$ for the base change of $i_w\times i_v$ along $\eta^u$.
Now if $\ell(w)+\ell(v)=\ell(wv)$, there is the natural isomorphism
\begin{eqnarray*}
\iw^u\backslash LG_{w}\times^{\iw^u}LG_v/\iw^u&=&\iw^u\backslash\iw^u\cdot\mS_k\cdot \dot{w}\times^{\iw^u\cap \Ad_{\dot{w}^{-1}} \iw^u} \iw^u\times^{\iw^u}\iw^u\cdot \mS_k\cdot\dot{v}\times^{\iw^u\cap \Ad_{\dot{v}^{-1}} \iw^u} \iw^u/\iw^u\\
&\cong&\iw^u\backslash \iw^u\cdot (\mS_k\times \mS_k) \cdot \dot{w}\dot{v} \times^{\iw^u\cap \Ad_{(\dot{w}\dot{v})^{-1}} \iw^u} \iw^u/\iw^u\\
&\cong& (\mS_k\times \mS_k) \times \bB(\iw^u\cap \Ad_{(\dot{w}\dot{v})^{-1}} \iw^u),
\end{eqnarray*}
Using \eqref{eq: swap H and X}, it follows that under the above isomorphism
\[
(\eta_{w,v}^u)^!((\pr_{\dot{w}})^!(\mL)[-\ell(w)]\boxtimes (\pr_{\dot{v}})^!(\mL')[-\ell(v)])\cong \mL\boxtimes w(\mL')\boxtimes \omega_{\bB(\iw^u\cap \Ad_{(\dot{w}\dot{v})^{-1}} \iw^u)}[-\ell(wv)].
\] 
Now using base change $(\eta^u)^!\circ ((i_{w})_?\boxtimes (i_{v})_?)\cong (i_{w,v})_?(\eta^u_{w,v})^!$ for $?=*,!$ (the $*$-case follows from the formalism of the sheaf theory $\shv$ and the $!$-case follows from \Cref{pfp.functors.and.bc.placid.stacks}), and using the fact that $*$- and $!$- convolutions on monodromic sheaves differ by a shift (which follows by \Cref{lem: functors between monodromic categories-5}), we deduce
Part \eqref{lem: monoidal product computation in affine Hecke category-1}.
Given this, we may assume that $\mL=\mL'=\widetilde{\Ch}$ in Part \eqref{lem: monoidal product computation in affine Hecke category-3} and \eqref{lem: monoidal product computation in affine Hecke category-2} (since $\Delta_{\dot{w}}^{\mon}(\mL)\cong \Delta_e^{\mon}(\mL)\star^u\widetilde{\Delta}_{\dot{w}}^{\mon}$, etc.). 

To proceed, we study convolutions between $\widetilde{\Delta}_{\dot{s}}^{\mon}$ and $\widetilde{\nabla}_{\dot{s}}^{\mon}$ when $s$ is a simple reflection. 
Clearly, any of such convolution is supported on $\iw^u\backslash LG_{\leq s}/\iw^u$. 
Using base change we may restrict the convolution diagram to
\[
\iw^u\backslash LG_{\leq s}/\iw^u\times \iw^u\backslash LG_{\leq s}/\iw^u\leftarrow \iw^u\backslash LG_{\leq s}\times^{\iw^u}LG_{\leq s}/\iw^u\to \iw^u\backslash LG_{\leq s}/\iw^u.
\]
By abuse of notations, the left arrow is still denoted by $\eta^u$ and the right arrow is denoted by $m^u$. For $\mF=\widetilde{\Delta}_{\dot{s}}^{\mon}$ or $\widetilde{\nabla}_{\dot{s}}^{\mon}$,
we need to compute 
\begin{equation}\label{eq: stalk of conv at s}
(\widetilde{\nabla}^{\mon}_{\dot{s}}\star^u\mF)|_{\mS_k\dot{s}},
\end{equation}
where $(-)|_{\mS_k\dot{s}}$ denotes the $!$-pullback along the (representable coh. smooth) morphism $\mS_k \dot{s}\to \iw^u\backslash LG_{\leq s}/\iw^u$. 
We let $\mathring{LG}_s= LG_s- \mS_k\cdot \dot{s}\cdot \iw^u$, which in open in $LG_s$. We have the following commutative diagram with all squares Cartesian
\[
\xymatrix{
\mS_k\times \mathring{LG}_s/\iw^u  \ar@{^{(}->}[d]_j\ar[rr] && \iw^u\bs LG_s\times^{\iw^u}LG_{s}/\iw^u\ar@{^{(}->}[d]\ar[r] &\iw^u\bs LG_s/\iw^u \times \iw^u\bs LG_{s}/\iw^u \ar@{^{(}->}[d]\ar^-{\pr_{\dot{s}}\times\pr_{\dot{s}}}[r] & \mS_k\times\mS_k\\
\mS_k\times LG_s/\iw^u \ar^{(t,g)\mapsto t\dot{s}}_h[d]\ar^-{(t,g)\mapsto (tg, g^{-1}\dot{s})}[rr] && \iw^u\bs LG_s\times^{\iw^u}LG_{\leq s}/\iw^u\ar[r]\ar^{m^u}[d] &\iw^u\bs LG_s/\iw^u \times \iw^u\bs LG_{\leq s}/\iw^u&\\
 \mS_k\dot{s} \ar[rr] &&   \iw^u\bs LG_{\leq s} /\iw^u&&
}\]
All the horizontal maps are (representable) coh. pro-smooth. By base change, the sheaf \eqref{eq: stalk of conv at s} is obtained from $\widetilde{\Ch}[-1]\boxtimes \widetilde{\Ch}[-1]$ on $\mS_k\times\mS_k$ by $!$-pullback along the top maps, followed by $j_!$ or $j_*$, and then followed by $h_*$.

We may write $LG_s/\iw^u\cong \mS_k \times \iw^u\dot{s}\iw^u/\iw^u\cong \mS_k\times\mS_k\times \bA^1$. Then the map $\mathring{LG}_s\cong \mS_k\times \bG_m$ and the map $j$ is induced by the standard inclusion $\bG_m\subset\bA^1$.
Now an $\SL_2$-computation shows that the composition of maps in the top row in the above diagram is identified with
\[
f: \mS_k\times\mS_k\times\bG_m\to \mS_k\times\mS_k,\quad (t, t', x)\mapsto (tt', {t'}^{-1}\hat\al_s(x)^{-1}).
\]
Note that we have the following commutative diagram with Cartesian squares
\[
\xymatrix{
&\ar_{h}[dl]\ar_{\pr_{13}}[d]\mS_k\times\mS_k\times \bA^1&\ar@{_{(}->}_j[l]\mS_k\times\mS_k\times\bG_m\ar^{\pr_{13}}[d]\ar^-f[rr] && \mS_k\times \mS_k\ar^m[d]\\
\mS_k&\mS_k\times\bA^1\ar_{\pr_1}[l] & \ar@{_{(}->}_j[l] \mS_k\times\bG_m \ar^-{(r,x)\mapsto r\hat{\al}_s(x)^{-1}}_{f'}[rr] && \mS_k
}
\]
It follows again by base change and by \Cref{lem: depth zero geom Langlands for tori}
\eqref{eq: stalk of conv at s} is computed as $(\pr_1)_* (j_?({f'}^!\widetilde{\Ch}[-2]))$ for $?=!$ or $*$, depending on whether $\mF$ is standard or costandard.

It is a standard fact that for any character sheaf $\Ch_\chi$ on $\bG_m$, $C_c^\bullet(\bA^1, j_!\Ch_\chi)=0$. It follows that
\begin{equation}\label{eq: restriction to s strata of convolution-1}
(\widetilde{\nabla}^{\mon}_{\dot{s}}\star^u\widetilde{\Delta}^{\mon}_{\dot{s}})|_{\mS_k\dot{s}}=0.
\end{equation}
To compute 
\[
(\pr_1)_* (j_*({f'}^!(\widetilde{\Ch}[-2])))=(\pr_1)_*({f'}^!(\widetilde{\Ch}[-2]))\simeq (\pr_1)_*({f'}^*\widetilde{\Ch}),
\] 
we can pass to the dual group and using the coherent description as in \Cref{lem: depth zero geom Langlands for tori}. So we have
\[
R_{I_F^t,\hat{S}}\xrightarrow{\id\times \hat{\al}_s}R_{I_F^t,\hat{S}}\times R_{I_F^t,\bG_m}\xleftarrow{\id\times \{1\}} R_{I_F^t,\hat{S}}.
\]
Here $\hat{S}$ is the dual torus of $\mS_k$, and $\al_s$ now is regarded as a character $\hat{S}\to \bG_m$.  The fiber product of the above map is nothing but $\ker \hat{\al}_s$.
It follows that
\begin{equation}\label{eq: restriction to s strata of convolution-2}
(\widetilde{\nabla}^{\mon}_{\dot{s}}\star^u\widetilde{\nabla}^{\mon}_{\dot{s}})|_{\mS_k\dot{s}}=\widetilde{\Ch}_s.
\end{equation}

Now we prove Part \eqref{lem: monoidal product computation in affine Hecke category-3} and Part \eqref{lem: monoidal product computation in affine Hecke category-2}. 

We first show that $\widetilde{\nabla}^{\mon}_{\dot{w}}\star^u\widetilde{\Delta}^{\mon}_{\dot{w}^{-1}}$ is supported on $\iw^u\backslash \iw/\iw^u$. Using Part \eqref{lem: monoidal product computation in affine Hecke category-1}, it is enough to prove this when $w=s$ is a simple reflection, and when $w\in \Omega_{\breve\bfa}$ is of length zero (see \eqref{eq-Iwahori-Weyl-semiproduct} for the notation). In fact, the length zero case already follows from Part \eqref{lem: monoidal product computation in affine Hecke category-1} as in this case $\widetilde{\nabla}_{\dot{w}}^{\mon}=\widetilde{\Delta}_{\dot{w}}^{\mon}$. Therefore, we assume that $w=s$ is a simple reflection. But this case follows from \eqref{eq: restriction to s strata of convolution-1}.

Therefore it remains to compute $(\widetilde{\nabla}^{\mon}_{\dot{w}}\star^u\widetilde{\Delta}^{\mon}_{\dot{w}^{-1}})|_{\mS_k e}$. Note that we have
\[
\mS_k\times LG_w/\iw^u\cong (\iw^u\backslash LG_w\times^{\iw^u}LG_{w^{-1}}/\iw^u)\times_{\iw^u\backslash LG/\iw^u} \mS_k e, \quad (t,g)\mapsto (tg,g^{-1}).
\]
Then using similar argument as above, we see that
\[
(\widetilde{\nabla}^{\mon}_{\dot{w}}\star^u\widetilde{\Delta}^{\mon}_{\dot{w}^{-1}})|_{\mS_k e}\simeq \widetilde{\Ch}\otimes C^\bullet(\iw^u\dot{w}\iw^u/\iw^u, \omega_{\iw^u\dot{w}\iw^u/\iw^u}[-2\ell(w)])\cong \widetilde{\Ch},
\]
as desired. Applying the automorphism $LG\to LG, g\mapsto g^{-1}$ gives $\widetilde{\Delta}^{\mon}_{\dot{w}}\star^u\widetilde{\nabla}^{\mon}_{\dot{w}^{-1}}\cong \widetilde{\Delta}_e^{\mon}$ as well.

To prove Part \eqref{lem: monoidal product computation in affine Hecke category-2}, we consider the
the cofiber of $\widetilde{\Delta}_{\dot{s}}^{\mon}\to \widetilde{\nabla}_{\dot{s}}^{\mon}$ is supported on $\iw^u\backslash \iw/\iw^u$ and therefore is of the form $\Delta_e^{\mon}(\mF)$ for some $\mF\in\shv^{\mon}(\mS_k)$. Then by  \eqref{lem: monoidal product computation in affine Hecke category-1} and \eqref{lem: monoidal product computation in affine Hecke category-3}, by convolving $ \widetilde{\nabla}_{\dot{s}}^{\mon} \star^u (-)$, we obtain the following fiber sequence
\[
\widetilde{\Delta}^{\mon}_e\to \widetilde{\nabla}^{\mon}_{\dot{s}}\star^u\widetilde{\nabla}^{\mon}_{\dot{s}}\to  \nabla^{\mon}_{\dot{s}}(\mF)\to.
\]
To compute $\mF$, we may restrict $\widetilde{\nabla}^{\mon}_{\dot{s}}\star^u\widetilde{\nabla}^{\mon}_{\dot{s}}$ to $\iw^u\backslash LG_{s}/\iw^u$. Then by \eqref{eq: restriction to s strata of convolution-2}, $\mF\simeq \widetilde{\Ch}_s[1]$.
The first desired fiber sequence follows. The second fiber sequence can be deduced from the first one using  Part \eqref{lem: monoidal product computation in affine Hecke category-3}.
\end{proof}

It follows that cofree monodromic (co)standard objects are invertible (and in particular dualizable) with respect to the monoidal structure of $\shv\bigl((\iw^u,(\chi\mbox{-})\mon)\backslash LG/(\iw^u,(\chi\mbox{-})\mon)\bigr)$. 
Namely, the inverse of $\Delta_{\dot{w},\chi}^{\mon}$ is given by $\nabla_{\dot{w}^{-1},w^{-1}(\chi)}^{\mon}$. Similarly, in $\shv_{\mon}(\iw^u\backslash LG/\iw^u)$
the inverse of $\widetilde{\Delta}_{\dot{w}}^{\mon}$ is given by $\widetilde{\nabla}_{\dot{w}^{-1}}^{\mon}$. But note that they are not compact objects.

\begin{corollary}\label{cor: hom from standard to costandard}
Let $u_1,u_2\in\widetilde{W}$. Let $\mL_1,\mL_2\in \shv_{\mon}(\mS_k)^{\heartsuit}$.
\begin{enumerate}
\item\label{cor: hom from standard to costandard-1} The object $\Delta^{\mon}_{\dot{u}_1}(\mL_1)\star^u\Delta^{\mon}_{\dot{u}_2}(\mL_2)$ is contained in the subcategory of $\shv_{\mon}(\iw^u\backslash LG/\iw^u)$ generated under extensions by objects of the form $\Delta^{\mon}_w(\mL)[n]$, for $w\in\widetilde{W}, \mL\in \shv_{\mon}(\mS_k)^{\heartsuit}, n\leq 0$.
\item\label{cor: hom from standard to costandard-2} For $w\in \widetilde{W}$, the $\La$-module
\[
\Hom_{\shv_{\mon}(\iw^u\backslash LG/\iw^u)}(\Delta_{\dot{u}_1}^{\mon}(\mL_1)\star^u \Delta^\mon_{\dot{u}_2}(\mL_2), \widetilde\nabla_{\dot{w}}^{\mon})
\]
belongs to $\Mod_\La^{\leq 0}$.
\end{enumerate}
\end{corollary}
\begin{proof}
Notice that convolution of monodromic sheaves are right $t$-exact. In addition, $\widetilde{\Ch}_s$ from \Cref{lem: monoidal product computation in affine Hecke category} belongs to $\shv_\mon(\mS_k)^{\heartsuit}$ (see \Cref{ex: kernel of dual hom}). Now Part \eqref{cor: hom from standard to costandard-1} follows from \Cref{lem: monoidal product computation in affine Hecke category} \eqref{lem: monoidal product computation in affine Hecke category-1}  \eqref{lem: monoidal product computation in affine Hecke category-2}.

By \Cref{lem: monoidal product computation in affine Hecke category} \eqref{lem: monoidal product computation in affine Hecke category-3}, we have 
\begin{eqnarray*}
& &\Hom_{\shv_{\mon}(\iw^u\backslash LG/\iw^u)}(\Delta_{\dot{u}_1}^{\mon}(\mL_1)\star^u \Delta^\mon_{\dot{u}_2}(\mL_2), \widetilde\nabla_{\dot{w}}^{\mon})\\
&\cong &\Hom_{\shv_{\mon}(\iw^u\backslash LG/\iw^u)}(\Delta_{\dot{u}_1}^{\mon}(\mL_1)\star^u \Delta^\mon_{\dot{u}_2}(\mL_2)\star^u\widetilde{\Delta}^\mon_{\dot{w}^{-1}}, \widetilde\Delta^{\mon}_e).
\end{eqnarray*}
Using Part \eqref{cor: hom from standard to costandard-1}
 the sheaf $\Delta_{\dot{u}_1}^{\mon}(\mL_1)\star^u \Delta_{\dot{u}_2}(\mL_2)^{\mon}\star^u\widetilde{\Delta}^\mon_{\dot{w}^{-1}}$ admits a filtration with associated graded being $\Delta^{\mon}_{\dot{v}_i}(\mL_i)[n_i]$ with $v_i\in \widetilde{W}$, $\mL_i\in\shv_\mon(\mS_k)^{\heartsuit}$ and $n_i\leq 0$.
Note that $\Hom(\Delta^{\mon}_{\dot{u}_i}(\mL_i),\widetilde\Delta^{\mon}_e)=0$ unless $u_i=e$.
Therefore, 
the hom space in question admits a filtration with associated graded being $\Hom(\Delta^{\mon}_e(\mL_i),\widetilde\Delta^{\mon}_e)[-n_i]$. Now 
the corollary follows as 
\[
\Hom(\Delta^{\mon}_e(\mL_i),\widetilde\Delta^{\mon}_e)\cong \Hom_{\shv_{\mon}(\mS_k)}(\mL_i,\widetilde{\Ch})=\Hom_{\shv(\mS_k)}(\mL_i,\delta_{1})\in \Mod_\La^{\heartsuit}.
\]
Here we recall $\delta_1:=(\{1\}\to \mS_k)_*\La$ is the delta sheaf at the unit of $\mS_k$. 
\end{proof}

Here is another consequence.

\begin{corollary}\label{cor: standard and costandard same K class}
Let $\shv_{\mon}(\iw^u\bs LG/\iw^u)'\subset \shv_{\mon}(\iw^u\bs LG/\iw^u)$ be the small idempotent complete stable subcategory generated by $\widetilde{\Delta}_{\dot{w}}^{\mon}$ and $\widetilde{\nabla}_{\dot{w}}^{\mon}$. Then in $K_0(\shv_{\mon}(\iw^u\bs LG/\iw^u)')$, we have
\[
[\widetilde\Delta^{\mon}_{\dot{w}}]=[\widetilde{\nabla}^{\mon}_{\dot{w}}].
\]
\end{corollary}
We note that the category $\shv_{\mon}(\iw^u\bs LG/\iw^u)'$ contains, but is strictly larger than, the category $\shv_{\mon}(\iw^u\bs LG/\iw^u)^\cpt$.
\begin{proof}
Given \Cref{lem: monoidal product computation in affine Hecke category}, the proof is similar to the corresponding fact for the usual small unipotent affine Hecke category. Here are the details.
Let $s$ be a simple reflection in $\widetilde{W}$. Then as in the proof of \Cref{lem: monoidal product computation in affine Hecke category}, we have
\begin{equation}\label{eq: fiber of standard to costandard}
\Delta^{\mon}_e(\widetilde{\Ch}_s)\to \widetilde{\Delta}_{\dot{s}}^{\mon}\to \widetilde{\nabla}_{\dot{s}}^{\mon}
\end{equation}
This implies that $\Delta^{\mon}_e(\widetilde{\Ch}_s)\in \shv_{\mon}(\iw^u\bs LG/\iw^u)'$. Then \Cref{lem: monoidal product computation in affine Hecke category} implies that the category $\shv_{\mon}(\iw^u\bs LG/\iw^u)'$ is a monoidal stable subcategory of $\shv_{\mon}(\iw^u\bs LG/\iw^u)$, and therefore $K_0(\shv_{\mon}(\iw^u\bs LG/\iw^u)')$ is equipped with the induced ring structure.

Now we prove the statement by induction on length of $w$. If $\ell(w)=0$, this is clear and if $w$ is a simple reflection, by \eqref{eq: fiber of standard to costandard} above, it is enough to show that $[\Delta^{\mon}_e(\widetilde{\Ch}_s)]=0$. But this follows from the cofiber sequnece
\begin{equation}\label{eq: fiber of standard to costandard-2}
\Delta^{\mon}_e(\widetilde{\Ch}_s)\to \Delta^{\mon}_e(\widetilde\Ch)\to \Delta^{\mon}_e(\widetilde\Ch),
\end{equation}
which in turn follows from  \eqref{eq: injective hull general}.

Next, if we write $w=vs$ with $\ell(w)=\ell(v)+1$ and $s$ is a simple reflection, then we have
\[
[\widetilde\Delta^{\mon}_{\dot{w}}]=[\widetilde\Delta^{\mon}_{\dot{v}}][\widetilde\Delta^{\mon}_{\dot{s}}]=[\widetilde\nabla^{\mon}_{\dot{v}}][\widetilde\nabla^{\mon}_{\dot{s}}]=[\widetilde\nabla^{\mon}_{\dot{w}}].
\]
\end{proof}

We recall the construction of the cofree monodromic tilting sheaves, following \cite{DYYZ}. Namely a cofree monodromic tilting sheaf is an object in $\shv_\mon(\iw^u\bs LG/\iw^u)$ which admits two finite filtrations, one with associated graded being those of the form $\widetilde\Delta_{\dot{w}}^{\mon}, \ w\in \widetilde{W}$  and the other with associated graded being those of the form $\widetilde\nabla_{\dot{v}}^{\mon}, \ v\in\widetilde{W}$. We have the following classification of cofree monodromic tilting sheaves.

\begin{proposition}\label{prop: classification of cofree tilting}
For each $w\in \widetilde{W}$, there is a unique (up to non-unique isomorphism) cofree tilting object $\widetilde{\Til}_{\dot{w}}^{\mon}$ satisfying the following conditions:
\begin{itemize}
\item $\widetilde{\Til}_{\dot{w}}^{\mon}\subset\shv_\mon(\iw^u\bs LG_{\leq w}/\iw^u)$ and $\widetilde{\Til}_{\dot{w}}^{\mon}|_{\iw^u\bs LG_w/\iw^u}\simeq \widetilde{\Ch}$ under the equivalence \eqref{eq: mon sheaf on one strata}.
\item Let $Z\subset R_{I_F^t,\hat{S}}$ be a connected component. Then $\widetilde{\Til}_{\dot{w}}^{\mon}\star^u\Delta_e^{\mon}(\Ch(\cohdual_Z))$ is indecomposable. 
\end{itemize} 
Every cofree monodromic tilting sheaves is a finite direct sum of the above $\widetilde{\Til}_{\dot{w}}^{\mon}$s.
\end{proposition}

This is standard. This type of results have been proved in various settings. In particular a version that is closely related to our situation is proved in \cite[\textsection{5}]{DYYZ}. The same argument applies \emph{mutatis mutandis}. So we only review the main ingredients.

Let $s$ be a simple reflection in $\widetilde{W}$. Then pushing out of \eqref{eq: fiber of standard to costandard}
along the map $\Delta^{\mon}_e(\widetilde{\Ch}_s)\to \Delta^{\mon}_e(\widetilde\Ch)$ in \eqref{eq: fiber of standard to costandard-2} gives the desired object $\widetilde{\Til}_{\dot{s}}^{\mon}$ associated to the simple reflection $s$.

Now if $w\in W_\af$, written as a product of simple reflections $w=s_{i_1}\cdots s_{i_n}$, lifted to $\dot{w}=\dot{s}_{i_1}\cdots\dot{s}_{i_n}$, then for every connected component $Z\subset R_{I_F^t,\hat{S}}$, there is a unique (up to non-unique isomorphism) indecomposable direct summand
\[
\widetilde{\Til}_{\dot{w},Z}^{\mon}\subset \widetilde{\Til}_{\dot{s}_{i_1}}^{\mon}\star^u\cdots\star^u \widetilde{\Til}_{\dot{s}_{i_n}}^{\mon}
\] 
whose restriction to $\iw^u\backslash LG_w/\iw^u$ is $\Ch(\cohdual_Z)$. Then we let
\[
\widetilde{\Til}_{\dot{w}}^{\mon}=\prod_Z \widetilde{\Til}_{\dot{w},Z}^{\mon}
\]
where $Z$ range over all connected components of $R_{I_F^t,\hat{S}}$.

Up to (non-unique) isomorphism, this object $\widetilde{\Til}_{\dot{w}}^{\mon}$ is independent of the choice of the way $w$ written as the product of simple reflections. Finally, if $w\in \widetilde{W}$, written as $w=w_a\tau$ for $w_a\in W_\af$ and $\tau\in \Omega_\mI$ and lifted to $\dot{w}=\dot{w}_a\dot{\tau}$, we have
\[
\widetilde{\Til}_{\dot{w}}^{\mon}=\widetilde{\Til}_{\dot{w}_a}^{\mon}\star^u\widetilde{\Delta}_{\dot{\tau}}^{\mon}.
\]

Now for each $w$, we fix a choice of $\widetilde{\Til}_{\dot{w}}^{\mon}$ together with an isomorphism $\widetilde{\Til}_{\dot{w}}^{\mon}|_{\iw^u\bs LG_w/\iw^u}\simeq \widetilde{\Ch}$ as in \Cref{prop: classification of cofree tilting}. Then 
we can define a tilting extension functor
\[
\Til_{\dot{w}}^{\mon}: \shv_{\mon}(\mS_k)\to \shv_\mon(\iw^u\bs LG/\iw^u),\quad \mL\mapsto \widetilde{\Til}_{\dot{w}}^{\mon}\star^u\Delta_e^{\mon}(\mL)\cong\Delta_e^{\mon}(\mL)\star^u \widetilde{\Til}_{\dot{w}}^{\mon}.
\]
In particular, we have cofree indecomposable $\chi$-monodromic tilting object $\Til_{\dot{w},\hchi}^{\mon}=\Delta^{\mon}_{e,\hchi}\star^u  \widetilde{\Til}_{\dot{w}}^{\mon}$.

\begin{lemma}\label{lem: exactness of convolution with tilting}
The functor $\widetilde{\mathrm{Til}}_{\dot{w}}^{\mon}\star^u(-)$  and $(-)\star^u \widetilde{\mathrm{Til}}_{\dot{w}}^{\mon}$ are perverse exact. The same statement holds for $\hchi$-version.
\end{lemma}
Here we define the perverse $t$-structure on $\iw^u\backslash LG/\iw^u$ using the generalized constant sheaf as in \Cref{SS: coh. duality. Kot. stack}.
\begin{proof}
This is essentially due to I. Mirkovic. Namely, as the multiplication map $m^u: LG_w\times^{\iw^u} LG/\iw^u\to LG/\iw^u$ is an affine morphism, the functor
\[
\shv_{\mon}(\mS_k)\otimes_\La \shv_{\mon}(\iw^u\backslash LG/\iw^u)\to \shv_{\mon}(\iw^u\backslash LG/\iw^u),\quad (\mL, \mF)\mapsto \nabla_{\dot{w}}^{\mon}(\mL)\star^u\mF
\]
is right $t$-exact.  The functor $\mF\mapsto \mF\star^u\widetilde{\Delta}_{\dot{w}}^{\mon}$ on the other hand, is left $t$-exact.
As $\widetilde{\Til}_{\dot{w}}^{\mon}$ admits a filtration with associated graded by $\widetilde{\Delta}_{\dot{w'}}^{\mon}$ as well as a filtration with associated graded by $\widetilde{\nabla}_{\dot{w'}}^{\mon}$, the lemma follows.
\end{proof}

Here is another result we need.  

\begin{proposition}\label{prop: rigid dual of cofree monodromic tilting}
The (right) dual of the cofree monodromic tilting sheaf $\widetilde{\Til}_{\dot{w}}^{\mon}$ with respect to the monoidal structure of $\shv_\mon(\iw^u\bs LG/\iw^u)$ is $\widetilde{\Til}_{\dot{w}^{-1}}^{\mon}$.
\end{proposition}
\begin{proof}
This is a direct consequence of the classification of cofree monodromic tilting sheaves, as $(\widetilde{\Til}_{\dot{w}}^{\mon})^\vee$ clearly satisfies conditions in \Cref{prop: classification of cofree tilting} (with $w$ replaced by $w^{-1}$), and therefore must be isomorphic to $\widetilde{\Til}_{\dot{w}^{-1}}^{\mon}$.
\end{proof}

\subsubsection{Equivariant affine Hecke category}
Let us also discuss $\chi$-equivariant version of the affine Hecke category, for a character $\chi: \pialg(\mS_k)\to\La^\times$. As usual, let $\hchi$ denote the completion of $\chi$ in $R_{I_F^t,\hat{S}}$.
First,  we have the equivariant category $\shv\bigl((\iw,\chi)\backslash LG/(\iw,\chi')\bigr)$ constructed from $\shv_\mon(\iw^u\bs LG/\iw^u)$ as in \eqref{eq: monodromic to equivariant}. Explicitly we have
\begin{eqnarray*}\label{eq: chi chi prime equivariant category}
\shv\bigl((\iw,\chi)\backslash LG/(\iw,\chi')\bigr) & \cong & (\Mod_\La)_\chi\otimes_{\shv_{\mon}(\mS_k)}\shv_\mon(\iw^u\backslash LG/\iw^u)\otimes_{\shv_\mon(\mS_k)} (\Mod_\La)_{\chi'}
\end{eqnarray*}
In particular, 
\[
\shv\bigl((\iw,u)\backslash LG/(\iw,u)\bigr)=\shv(\iw\backslash LG/\iw)
\] 
by \Cref{lem: two-definition of equivariant category}.

We make use of the following lemma to endow $\shv\bigl((\iw,\chi)\backslash LG/(\iw,\chi')\bigr)$ (for $\chi=\chi'$) with a monoidal structure in the $\infty$-categorical setting. Consider
\[
\shv(\iw^u\bs LG/(\iw,\chi))=\shv_\mon(\iw^u\bs LG/(\iw,\chi))=\shv_\mon(\iw^u\bs LG/\iw^u)\otimes_{\shv_\mon(\mS_k)}(\Mod_\La)_{\chi},
\]
which admits a natural left $\shv_\mon(\iw^u\bs LG/\iw^u)$-module structure.
\begin{lemma}\label{lem: monoidal structure of equivariant Hecke cat}
There is a canonical equivalence of $\La$-linear categories
\[
\End_{\shv_{\mon}(\iw^u\backslash LG/\iw^u)} \shv_{\mon}(\iw^u\backslash LG/(\iw,\chi)\bigr)\cong \shv\bigl((\iw,\chi)\backslash LG/(\iw,\chi)\bigr).
\]
\end{lemma}

\begin{proof}
Similar to \Cref{lem: self dual of monodromic and equiv sheaves on H}, $\shv(\iw^u\bs LG/(\iw,\chi))$ as a left $\shv_\mon(\iw^u\bs LG/\iw^u)$-module admits a left dual, given by $\shv((\iw,\chi)\bs LG/\iw^u))$. Then we have
\begin{small}
\begin{eqnarray*}
&  &\End_{\shv_{\mon}(\iw^u\backslash LG/\iw^u)} \shv_{\mon}(\iw^u\backslash LG/(\iw,\chi)\bigr) \\
&=&\shv((\iw,\chi)\bs LG/\iw^u)\otimes_{\shv_\mon(\iw^u\backslash LG/\iw^u)} \shv(\iw^u\backslash LG /(\iw,\chi))\\
&=&(\Mod_\La)_{\chi}\otimes_{\shv_{\mon}(\mS_k)}\shv_\mon(\iw^u\backslash LG/\iw^u)\otimes_{\shv_\mon(\iw^u\backslash LG/\iw^u)} \shv_\mon(\iw^u\backslash LG/\iw^u)\otimes_{\shv_{\mon}(\mS_k)}(\Mod_\La)_{\chi}\\
&=&\shv\bigl((\iw,\chi)\backslash LG/(\iw,\chi)\bigr).
\end{eqnarray*}
\end{small}
\end{proof}

The above lemma in particular endows $\shv\bigl((\iw,\chi)\backslash LG/(\iw,\chi)\bigr)$ with a monoidal structure, namely the one opposite to the natural one on $\End_{\shv_\mon(\iw^u\bs LG/\iw^u)} \shv(\iw^u\bs LG/(\iw,\chi))$.
We shall temporarily call such monoidal structure of $\shv\bigl((\iw,\chi)\backslash LG/(\iw,\chi)\bigr)$ the endomorphism monoidal structure.
When $\chi=u$, it coincides with the natural convolution monoidal structure of $\shv(\iw\backslash LG/\iw)$. In fact, a more general statement is true, as we shall see shortly.

For general $\chi,\chi'$, we can also access the category \eqref{eq: chi chi prime equivariant category} via the equivalence \eqref{eq: equivariant category via gerbe} (under a mild restriction). For simplicity, we will assume that $\La$ is a field in the sequel. Let $p'$ be the product of $p$ and the characteristic exponent of $\La$ (so $p'=p$ if $\La$ is a field of characteristic zero and otherwise $p'=p\cdot\mathrm{char} \La$).  Then every prime-to-$p$ finite order character $\chi: T^p\mS_k\to \La^\times$ has order coprime to $p'$.
For an integer positive $n$ coprime to $p'$, we define $\iw^{[n]}$ via the Cartesian pullback
\begin{equation}\label{eq: iw[n]}
\xymatrix{
\iw^{[n]}\ar^{\varphi^n}[r]\ar[d] & \iw\ar[d]\\
\mS_k\ar^{[n]}[r] & \mS_k.
}
\end{equation}
(Do not confuse $\iw^{[n]}$ with the $n$th congruence subgroup of $\iw$, which we usually denote by $\iw^{(n)}$.)
Sometimes we also write $[n]:\mS_k\to \mS_k$ as $\varphi^n: \mS_k^{[n]}\to\mS_k$. Since $(n,p')=1$, the scheme $\chi_{\varphi^n}$ from \Cref{ex: Hn-monodromic}, denoted by $\chi_n$ for simplicity, is just disjoint union of points, and by \Cref{ex: Hn-equivariant} we have
\begin{equation}\label{eq:decomposition of equivariant cat via characters}
\shv(\iw^{[n]}\backslash LG/\iw^{[n]})=\bigoplus_{\chi,\chi': \mS_k[n]\to\La^\times}\shv\bigl((\iw,\chi)\backslash LG/(\iw,\chi')\bigr).
\end{equation}

For $w\in \widetilde{W}$, the map \eqref{eq: projection w} induces a map
\[
\pr^{[n]}_{\dot{w}}: \iw^{[n]}\backslash LG_w/\iw^{[n]}\to \mS_k^{[n]}\backslash \mS_k  \dot{w} / \mS_k^{[n]}\to  \mS_k^{[n]}\backslash \mS_k=\bB \mS_k[n].
\]
Now given $\chi$ of finite order $n$ (coprime to $p'$), considered as local system on $\bB \mS_k[n]$, we may define
\[
\Delta_{w,\chi}=(i_w)_!(\pr^{[n]}_w)^!\chi[-\ell(w)],\quad \nabla_{w,\chi}=(i_w)_*(\pr^{[n]}_w)^!\chi[-\ell(w)].
\] 
They are standard and costandard objects in $\shv\bigl((\iw,w\chi)\backslash LG/(\iw,\chi)\bigr)$. 

We will use $\star^{[n]}$ to denote the monoidal structure on $\shv(\iw^{[n]}\backslash LG/\iw^{[n]})$ defined as above. More generally, we will use $\star^{[n]}$ to denote any morphism induced by the map $m^{[n]}: LG\times^{\iw^{[n]}}LG\to LG$. Each $\shv\bigl((\iw,\chi)\backslash LG/(\iw,\chi)\bigr)$ is closed under monoidal product, with the unit given by $\Delta_{e,\chi}=\nabla_{e,\chi}$, and therefore acquires a monoidal category structure by \Cref{lem: complement on monoidal categories}. 
We shall temporarily call the corresponding monoidal structure of $\shv\bigl((\iw,\chi)\backslash LG/(\iw,\chi)\bigr)$ the convolution monoidal structure.
Note this
$\shv\bigl((\iw,\chi)\backslash LG/(\iw,\chi')\bigr)$ is a $\shv\bigl((\iw,\chi)\backslash LG/(\iw,\chi)\bigr)\mbox{-}\shv\bigl((\iw,\chi')\backslash LG/(\iw,\chi')\bigr)$-bimodule.

To compare the above two monoidal structures on $\shv\bigl((\iw,\chi)\backslash LG/(\iw,\chi)\bigr)$, we consider the category 
\[
\shv_{\mon}(\iw^{u}\backslash LG/\iw^{[n]})=\shv(\iw^{u}\backslash LG/\iw^{[n]}),
\] 
which is a $\shv_{\mon}(\iw^{u}\backslash LG/\iw^{u})\mbox{-}\shv(\iw^{[n]}\backslash LG/\iw^{[n]})$-bimodule. As before, since $(n,p')=1$, there is the direct sum decomposition
\[
\shv(\iw^{u}\backslash LG/\iw^{[n]})=\bigoplus_{\chi}\shv\bigl(\iw^u\backslash LG/(\iw,\chi)\bigr).
\]
Each $\shv\bigl(\iw^u\backslash LG/(\iw,\chi)\bigr)$ is a $\shv_{\mon}(\iw^u\backslash LG/\iw^u)\mbox{-}\shv\bigl((\iw,\chi)\backslash LG/(\iw,\chi)\bigr)$-bimodule. This shows that the identity functor of $\shv\bigl((\iw,\chi)\backslash LG/(\iw,\chi)\bigr)$ is monoidal, with the resource equipped with the convolution monoidal structure and the target equipped with the endomorphism monoidal structure. Therefore, the two monoidal structures on $\shv\bigl((\iw,\chi)\backslash LG/(\iw,\chi)\bigr)$ coincide.

There are parallel discussions with $\shv(\iw^{u}\backslash LG/\iw^{[n]})$ replaced by $\shv(\iw^{[n]}\backslash LG/\iw^{u})$.

For an object $\mF\in \shv(\iw^{[n]}\backslash LG/\iw^{[n]})$, let $\mF^l$ (resp. $\mF^r$) denote its $!$-pullback to $\shv(\iw^{u}\backslash LG/\iw^{[n]})$ (resp. $\shv(\iw^{[n]}\backslash LG/\iw^{u})$).

\begin{lemma}\label{lem: l and r pullback compatible with conv}
We have $\Delta_{w,\hchi}^{\mon}\star^u \mF^l\cong (\Delta_{w,\chi}\star^{[n]}\mF)^l$ and $\nabla_{w,\chi}^{\mon}\star^u \mF^l\cong (\nabla_{w,\chi}\star^{[n]}\mF)^l$. The similar statements hold with $(-)^l$ replaced by $(-)^r$.
\end{lemma}
\begin{proof}
We have the following diagram with two squares (involving $\eta^u$ and $m^{[n]}$) Cartesian
\[
\xymatrix{
&\ar_-{m^u}[ld]\iw^u\backslash LG_w\times^{\iw^u} LG\ar^-{\eta^u}[r]\ar[d]&\iw^u\backslash LG_w/\iw^u \times \iw^u\backslash LG\ar^{\av^{[n]}_1}[d]& \\
\iw^u\backslash LG\ar[d] & \ar[l] \iw^u\backslash LG_w\times^{\iw^{[n]}} LG\ar[d]\ar[r]& \iw^u\backslash LG_w/\iw^u\times^{\mS_k^{[n]}} \iw^u\backslash LG \ar[r]\ar^{\av^{[n]}_2}[dr]& \iw^u\backslash LG_w/\iw^{[n]}\times \iw^{[n]}\backslash LG\ar[d] \\
\iw^{[n]}\backslash LG & \ar_-{m^{[n]}}[l] \iw^{[n]}\backslash LG_w\times^{\iw^{[n]}} LG\ar^-{\eta^{[n]}}[rr]&& \iw^{[n]}\backslash LG_w/\iw^{[n]}\times \iw^{[n]}\backslash LG.
}\]
Using this diagram and various base change, and the fact that $(m^u)_*$ and $(m^u)_!$ differ by a shift, we can reduce the prove of the first isomorphism 
to proving $(\av^{[n]}_1)_*(\Delta_{\dot{w},\hchi}^{\mon}\boxtimes \mF^l)\cong (\av^{[n]}_2)^!(\Delta_{w,\chi}\boxtimes \mF)$, which in turn follows from the canonical isomorphism $\Ch_{\hchi}\star \Ch_\chi\cong \Ch_\chi$. The proofs of other isomorphisms are similar.
\end{proof}

For later discussion of Whittaker models, we will also take $\widetilde{L^+\mG}\to L^+\mG$ to be $\widetilde{\iw^u}\to \iw^u$ where $\widetilde{\iw^u}$ is defined as follows. Let $\breve\bff\subset\overline{\breve\bfa}$ be a facet contained the in closure of the alcove $\breve\bfa$ (determined by $\iw$). Let
$e_{\breve\bff}: \iw^u\to \bG_a$ be a surjective homomorphism given by 
\[
\iw^u\to \iw^u/[\iw^u,\iw^u]\cong \prod_\al U_\al\to \bG_a,
\] 
where $\al$ ranges over all affine simple roots of $(LG)_k$, such that the restriction of $e_{\breve\bff}$ to $U_\al\to \bG_a$ is an isomorphism for $\al\in \Phi_{\breve\bff}$ and the restriction of $e_{\breve\bff}$ to $U_\al\to \bG_a$ is trivial if $\al\not\in \Phi_{\breve\bff}$. 
Let $\widetilde{\iw^u}\to \iw^u$ be the pullback of the Artin-Scheier isogeny $\bG_a\to \bG_a$. Note that $\widetilde{\iw^u}$ is still coh. pro-unipotent and $\widetilde{\iw^u}\to \iw^u$ is a finite \'etale cover with Galois group $\bG_a(k_F)=k_F$. We write $\underline{k_F}$ for the constant group scheme over $k$ given by $k_F$.

Similar to \eqref{eq:decomposition of equivariant cat via characters}, we have the decomposition
\begin{equation}\label{eq: decomposition of endo whittaker} 
\shv(\widetilde{\iw^u}\backslash LG/\widetilde{\iw^u})\cong \bigoplus_{\psi,\psi': k_F\to \La^\times}\shv\bigl((\iw^u,\psi)\backslash LG/(\iw^u,\psi')\bigr).
\end{equation}
and $\mathbf{1}_{\widetilde{\iw^u}}=\oplus_\psi \mathbf{1}_{(\iw^u,\psi)}$.

We use $\star^{\tilde{u}}$ to denote the monoidal structure on $\shv(\widetilde{\iw^u}\backslash LG/\widetilde{\iw^u})$, and more generally any $*$-pushforward induced by the multiplication $LG\times^{\widetilde{\iw^u}}LG\to LG$.  
Each $\shv\bigl((\iw^u,\psi)\backslash LG/(\iw^u,\psi)\bigr)$ is closed under the convolution, and in fact is monoidal with the unit given by $\mathbf{1}_{(\iw^u,\psi)}$. In addition, each $\shv\bigl((\iw^u,\psi)\backslash LG/(\iw^u,\psi')\bigr)$ is a $\shv\bigl((\iw^u,\psi)\backslash LG/(\iw^u,\psi)\bigr)\mbox{-}\shv\bigl((\iw^u,\psi')\backslash LG/(\iw^u,\psi')\bigr)$-bimodule.
Next  consider
\[
\shv_{\mon}(\iw^u\backslash LG/\widetilde{\iw^u})\cong \bigoplus_{\psi:k_F\to \La^\times}\shv_{\mon}(\iw^u\backslash LG/(\iw^u,\psi)).
\] 
Here the monodromic category is defined using the left action of $\mS_k$ on $\iw^u\backslash LG/\widetilde{\iw^u}$. The category $\shv_{\mon}(\iw^u\backslash LG/(\iw^u,\psi))$ is a $\shv_{\mon}(\iw^u\backslash LG/\iw^u)\mbox{-} \shv((\iw^u,\psi)\backslash LG/(\iw^u,\psi))$-bimodule. 

We will define functors similar to \eqref{eq: (co)standard functor}. Let $i_w: \iw^u\backslash LG_w/\widetilde{\iw^u}\to \iw^u\backslash LG/\widetilde{\iw^u}$ be the locally closed embedding.
Let $w\in\widetilde{W}$ such that $\iw^u\cap \Ad_{\dot{w}^{-1}}\iw^u$ belongs to the kernel of the homomorphism $\iw^u\to \bG_a$. This means that $w$ is the longest element in its coset $w W_{\breve\bff}$, where we recall that $W_{\breve\bff}\subset\widetilde{W}$ is the subgroup generated by affine reflections corresponding to affine simple roots in $\Phi_{\breve\bff}$. 
Then the projection $\pr_{\dot{w}}$ from \eqref{eq: projection w} induces a map
\[
\widetilde{\pr}_{\dot{w}}^r:  \iw^u\backslash LG_w/\widetilde{\iw^u} \cong \mS_k\dot{w}\times (\iw^u\cap \Ad_{\dot{w}^{-1}}\iw^u)\backslash \iw^u/\widetilde{\iw^u}\to \mS_k\times \bB \underline{k_F}.
\]
It follows that $(\widetilde{\pr}_{\dot{w}}^r)^![-\ell(w)]$ induces a $t$-exact equivalence of categories
\[
\shv_{\mon}(\mS_k)\otimes_\La \rep(k_F)\cong \shv_\mon(\iw^u\backslash LG_w/\widetilde{\iw^u}),
\]
where regard representations of $k_F$ as sheaves on $\bB\underline k_F$ as usual.
Similar to \eqref{eq: (co)standard functor}, we can define the following functors
\begin{equation}\label{eq: monodromic whittaker (co)standard}
\Delta_{\dot{w}}^{\mon,\psi}, \nabla_{\dot{w}}^{\mon,\psi}: \shv_{\mon}(\mS_k)\to \shv_{\mon}(\iw^u\backslash LG/(\iw^u,\psi))
\end{equation}
as
\[
\Delta_{\dot{w}}^{\mon,\psi}(\mL)=(i_w)_!((\widetilde{\pr}_{\dot{w}}^r)^!(\mL[-\ell(w)]\boxtimes \psi)),\quad \nabla_{\dot{w}}^{\mon,\psi}(\mL)=(i_w)_*((\widetilde{\pr}_{\dot{w}}^r)^!(\mL[-\ell(w)]\boxtimes \psi)).
\]

Note that when $\psi$ is trivial, we have 
\[
\shv_{\mon}(\iw^u\backslash LG/(\iw^u,\psi))=\shv_{\mon}(\iw^u\backslash LG/\iw^u),
\] 
and under such identification the functor $\Delta_{\dot{w}}^{\mon,\psi}$ (resp. $\nabla_{\dot{w}}^{\mon,\psi}$)
is nothing but the functor $\Delta_{\dot{w}}^{\mon}$ (resp. $\nabla_{\dot{w}}^{\mon}$).  

Thus we can extend the definition of functors \eqref{eq: monodromic whittaker (co)standard} from those $w$ of maximal length in $wW_{\breve\bff}$ to every $w\in \widetilde{W}$. Namely, when $\psi$ is non-trivial, we simply let $\Delta_{\dot{w}}^{\mon,\psi}=\nabla_{\dot{w}}^{\mon,\psi}=0$ if $w$ is not the longest element in its coset $w W_{\breve\bff}$. When $\psi$ is trivial, we let $\Delta_{\dot{w}}^{\mon,\psi}=\Delta_{\dot{w}}^{\mon}$ and $\nabla_{\dot{w}}^{\mon,\psi}=\nabla_{\dot{w}}^{\mon}$ for all $w$.

We write
\[
\widetilde{\Delta}^{\mon,\psi}_{\dot{w}}=\Delta_{\dot{w}}^{\mon,\psi}(\widetilde\Ch),\quad \widetilde{\nabla}_{\dot{w}}^{\mon,\psi} =\nabla_{\dot{w}}^{\mon,\psi}(\widetilde\Ch).
\]

Let $w_0^{\breve\bff}$ be the longest length element in $w_{\breve\bff}$ with $\dot{w}^{\breve\bff}_0$ a lifting. We have the following lemma, which again is the monodromic generalization of well-known facts about Whittaker categories. 
The usual arguments (e.g. see \cite[Lemma 4]{arkhipov2009perverse}) work with appropriate modifications as in the proof of \Cref{lem: monoidal product computation in affine Hecke category}.

\begin{lemma}\label{lem: convolve (co)standard with Whit}
Assume that $\psi$ is non-trivial.
We have a canonical isomorphism of functors
\[
\Delta^{\mon,\psi}_{\dot{w}_0^{\breve\bff}}\cong \nabla^{\mon,\psi}_{\dot{w}_0^{\breve\bff}}.
\]
For $w\in\widetilde{W}$, let $w^{\breve\bff}\in wW_{\breve\bff}$ be the minimal length element in the coset.
We have
\[
\Delta_{\dot{w}}^{\mon}(\mL)\star^u\Delta^{\mon,\psi}_{\dot{w}_0^{\breve\bff}}(\mL')\cong \Delta_{\dot{w}^{\breve\bff}\dot{w}_0^{\breve\bff}}^{\mon,\psi}(\mL\star w(\mL')), \quad \nabla_{\dot{w}}^{\mon}(\mL)\star^u\nabla_{\dot{w}_0^{\breve\bff}}^{\mon,\psi}(\mL')\cong \nabla_{\dot{w}^{\breve\bff}\dot{w}_0^{\breve\bff}}^{\mon,\psi}(\mL\star w(\mL')).
\]
\end{lemma}

There is a parallel story for $\shv^{\mon}(\widetilde{\iw^u}\backslash LG/\iw^u)$
and similarly defined functors as in \eqref{eq: monodromic whittaker (co)standard}, which we denote by ${}^{\psi}\Delta^{\mon}_{\dot{w}}$ and ${}^{\psi}\nabla^{\mon}_{\dot{w}}$. Applying these functors to $\widetilde\Ch$, we obtain the objects ${}^{\psi}\widetilde\Delta^{\mon}_{\dot{w}}$ and ${}^{\psi}\widetilde\nabla^{\mon}_{\dot{w}}$.

\begin{lemma}\label{lem: right adjoint of Whittaker average}
The functor
\[ 
\shv_{\mon}(LG/\iw^u)\to \shv_{\mon}(LG/(\iw^u,\psi)),\quad \mF\mapsto \mF\star^u  \widetilde{\Delta}^{\mon,\psi}_{\dot{w}}  
\]
admits a (continuous) right adjoint given by
\[
\shv_{\mon}(LG/(\iw^u,\psi))\to \shv_{\mon}(LG/\iw^u), \quad \mG\mapsto \mG\star^{\tilde{u}}  {}^{\psi}\widetilde{\nabla}^{\mon}_{\dot{w}^{-1}}. 
\]
\end{lemma}
\begin{proof}
We first deal with the case $w=w_0^{\breve\bff}$. Let $L^+\breve\mP_{\breve\bff}^u=\iw^u\cap \dot{w}^{\breve\bff}_0\iw^u(\dot{w}_0^J)^{-1}$, which is the pro-unipotent radical of the standard parahoric $L^+\breve\mP_{\breve\bff}$ corresponding to ${\breve\bff}$.
Note that $L^+\breve\mP_{\breve\bff}^u$ is a normal subgroup of $\iw^u$ (and of $\widetilde{\iw}^u$), stable under the conjugation by $\dot{w}_0^{\breve\bff}$. Note that we have the natural isomorphisms
\[
\xymatrix{
LG\times^{\iw^u}\iw^u\dot{w}_0^{\breve\bff}\iw^u/\widetilde{\iw^u}& \ar_-{b_l}^-{\cong}[l] LG/\breve\mP_{\breve\bff}^u\times \iw^u/\widetilde{\iw^u} \ar_-\cong^-{b_r}[r] &LG\times^{\widetilde{\iw^u}} \iw^u\dot{w}_0^{\breve\bff}\iw^u/\iw^u,
}\]
where $b_l$ sends $(g, t, h)\in LG/\breve\mP_{\breve\bff}^u\times  \iw^u/\widetilde{\iw^u}$ to $(g, \dot{w}_0^{\breve\bff}h)$ and $b_r$ sends $(g,t,h)$ to $(gh, h^{-1}t\dot{w}_0^{\breve\bff})$. Note that as convolving with $\widetilde{\Ch}$ is an identity functor for monodromic sheaves, we have
\[
 \mF\star^u \widetilde{\Delta}_{\dot{w}}^{\mon,\psi}\cong (a_{l})_!(\mF'\boxtimes \psi[-\ell(w_0^{\breve\bff})]), 
\]
\[
\mG\star^{\tilde{u}}  {}^{\psi}\widetilde{\nabla}_{\dot{w}}^{\mon}\cong (a_r)_*(\mG'\boxtimes \omega_{\iw^u/\widetilde{\iw^u}}[-\ell(w_0^{\breve\bff})]),
\]
where 
\begin{itemize}
\item $a_l: LG/\breve\mP_{\breve\bff}^u\times  \iw^u/\widetilde{\iw^u}\to LG/\widetilde{\iw^u}$ is the map sending $(g, h)$ to $g\dot{w}_0^{\breve\bff}h$ (so $a_l=m^u\circ b_l$);
\item $a_r: LG/\breve\mP_{\breve\bff}^u\times  \iw^u/\widetilde{\iw^u}\to LG/\iw^u$ is the map sending $(g,h)$ to $g\dot{w}_0^{\breve\bff}$ (so $a_r=m^{\tilde{u}}\circ b_r$);
\item $\mF'$ is the $!$-pullback of $\mF$ along $LG/L^+\breve\mP_{\breve\bff}^u\to LG/\iw^u$;
\item $\mG'$ is the $!$-pullback of $\mG$ along $LG/L^+\breve\mP_{\breve\bff}^u\to LG/\widetilde{\iw^u}$.
\end{itemize}

We also notice that $(a_l)^!\mG\cong \mG'\boxtimes  \psi$ and $(a_r)^*\mF\cong \mF'[-2\ell(w_0^{\breve\bff})]\boxtimes \La$.
The lemma for $w=w_0^{\breve\bff}$ then follows from the isomorphism on $LG/\mP_{\breve\bff}^u\times \iw^u/\widetilde{\iw^u}$
\[
\Hom(\mF'\boxtimes \psi, \mG'\boxtimes \psi)=\Hom(\mF'\boxtimes \La, \mG'\boxtimes\La).
\]

To deal with general situation, we note that if $w$ is not the longest length element in $wW_{\breve\bff}$, then the functor is zero. Otherwise, 
we write $w=ww_0^{\breve\bff} w_0^{\breve\bff}$ so $ww_0^{\breve\bff}=w^{\breve\bff}$ as in  \Cref{lem: convolve (co)standard with Whit}.
Then the lemma follows from the case $w=w_0^{\breve\bff}$, together with  \Cref{lem: convolve (co)standard with Whit} and \Cref{lem: monoidal product computation in affine Hecke category} \eqref{lem: monoidal product computation in affine Hecke category-3}.
\end{proof}

\subsection{Affine Deligne-Lusztig theory}
We next generalize some constructions in the Deligne-Lusztig theory to the affine setting.

\subsubsection{Affine Deligne-Lusztig sheaves}\label{SS: affine DL sheaf}
For every $w\in \widetilde{W}$, we define
\begin{equation}\label{eq: unip ADLS}
   R^*_w=\Nt_*((i_{w})_*\omega_{\Sht^\loc_w}[-\ell(w)]),\quad  R_{w}^!= \Nt_* ((i_{w})_!\omega_{\Sht^\loc_w}[-\ell(w)]).
\end{equation}

The costalks of $R^*_w$ and the stalks of $R_w^!$ admit the following interpretations.

\begin{lemma}\label{lem: stalk of ADLS} 
For every $b\in B(G)$, we have
\begin{eqnarray*}
 (i_b)^!R_w^*\cong C_\bullet^{\mathrm{BM}}(X_w(b),\La[-\ell(w)]) \in \fgrep(G_b(F),\La),
\end{eqnarray*}
which is the Borel-Moore homology  of the affine Deligne-Lusztig variety $X_w(b)$. On the other hand,
\[
(i_b)^*R_w^!\cong (\verd^{\can}_{G_b(F)})^{\fg}( (i_b)^!R_w^*)[-2\langle2\rho, \nu_b\rangle].
\]
\end{lemma}
\begin{proof}
Using the base change and the equivalence $\shv(\kot_{G,b})\cong \rep(G_b(F))$, one sees that $(i_b)^!R_w^*$ can be identified with $(\pi_{X_w(b)})_* \omega_{X_w(b)}$ equipped with an action of $G_b(F)$, where we think $b$ as a point $b: \Spec k\to \kot_G$ and $\pi_{X_{w}(b)}\colon X_w(b)\to \Spec k$ is the structural map. As $X_w(b)$ is an ind-scheme, ind pfp over $k$, $(\pi_{X_w(b)})_* \omega_{X_w(b)}$ is nothing but the usual Borel-Moore homology of $X_w(b)$ (see \Cref{ex: Borel-Moore homology}). This gives the first isomorphism.
The second isomorphism then follows from the canonical duality \Cref{rem: can duality on f.g. sheaves}:
\begin{eqnarray*}
 (\verd^{\can}_{G_b(F)})^\fg((i_{b})^!(R_w^*)) & \cong &(i_{b})^*((\verd_{\kot_G}^{\can})^\fg(R_w^*))[2\langle2\rho, \nu_b\rangle]\\
                                                                        & \cong &(i_b)^*(\Nt_*((\verd_{\sht^\loc})^{\fg}((i_w)_*\consdual_{\sht_w^{\loc}}[-\ell(w)])))[2\langle2\rho, \nu_b\rangle]\\
                                                                        & \cong &(i_b)^*R^!_w[2\langle2\rho, \nu_b\rangle].
\end{eqnarray*}
\end{proof}

\begin{remark}
\begin{enumerate}
\item The proposition says that  $R_w^*$ and $R_w^!$ are sheaves on $\kot_G$ obtained by gluing these representations. For this reason, we call $R_w^*$ and $R_w^!$ the \textit{unipotent affine Deligne-Lusztig sheaves}.
\item As we shall see later, for some special $w$, $*$-stalks of $R_w^!$ admit more explicit description.
\end{enumerate}
\end{remark}

As in the usual Deligne-Lusztig theory, each space $\Sht^{\loc}_w=\frac{LG_w}{\Ad_\sigma\iw}$ admits a finite \'etale Galois covering given by
\[
\widetilde{\Sht}^{\loc}_{\dot{w}}:=\frac{\iw^u\dot{w}\iw^u}{\Ad_\sigma\iw^u},
\]
where $\iw^u$ is the pro-unipotent radical of $\iw$ and $\dot{w}$ is a lifting of $w\in\widetilde{W}$ to $N_G(S)(\breve F)$. The corresponding Galois group is the finite abelian group 
\begin{equation}\label{eq:w-twisted torus}
\mS_k^{\bar{w}\sigma}:= \ker (\varphi_{\bar w}: \mS_k\to \mS_k), \quad \varphi_{\bar w}(s)=s^{-1}\bar{w}\sigma(s),
\end{equation} 
where $\bar{w}$ is the image of $w$ in the finite Weyl group $W_0$, which acts on $\mS_k$. Alternatively, we consider the projection \eqref{eq: projection w}, which induces a map
\begin{equation}\label{eq: sigma projection w}
\pr_{\dot{w}}^\sigma: \frac{LG_w}{\Ad_\sigma \iw}\to \frac{\mS_k \cdot \dot{w}}{\Ad_\sigma \mS_k}\cong \bB \mS_k^{\bar{w}\sigma}.
\end{equation}
Then $\widetilde{\Sht}^{\loc}_{\dot{w}}\to \Sht^{\loc}_w$ is the pullback of the $\mS_k^{\bar{w}\sigma}$-torsor $\dot{w}\to \bB \mS_k^{\bar{w}\sigma}$.

Now for a $\La[\mS_k^{\bar{w}\sigma}]$-module $M$, regarded as a local system on $\bB \mS_k^{\bar{w}\sigma}$ in the usual way (e.g. via the equivalence from \Cref{ess.pro.p.torsor.hom.descent}), 
we define two functors
\begin{equation}\label{eq: general ADLI of M}
R^?_{\dot{w}}(-): \rep(\mS_k^{\bar{w}\sigma})\to \shv(\kot_G), \quad R^?_{\dot{w}}(M):=\Nt_*(i_w)_?(\pr_{\dot{w}}^\sigma)^!M[-\ell(w)], \quad ?=*,!.
\end{equation}
When $M=\La[\mS_k^{\bar{w}\sigma}]$ regarded as a left module over itself, we simply write 
\[
\widetilde{R}^*_{\dot{w}}=R^*_{\dot{w}}(\La[\mS_k^{\bar{w}\sigma}]),\quad \widetilde{R}^!_{\dot{w}}=R^!_{\dot{w}}(\La[\mS_k^{\bar{w}\sigma}]).
\]
Then 
\[
R^?_{\dot{w}}(M)=\widetilde{R}^?_{\dot{w}}\otimes_{\La[\mS_k^{\bar{w}\sigma}]} M,\quad \quad ?=*, !.
\]
On the other hand, when $M$ is given by a character $\theta: \mS_k^{\bar{w}\sigma}\to \La^\times$, we write
\[
R^*_{\dot{w},\theta}:=R^*_w(M), \quad R^!_{\dot{w},\theta}:=R^!_w(M).
\]
Note that when $\theta$ is the trivial character, the above two objects reduce to $R_w^*$ and $R_w^!$ above. We call $R^*_{\dot{w},\theta}$ and  $R^!_{\dot{w},\theta}$ affine Deligne-Lusztig sheaves, which are affine analogue of Deligne-Lusztig characters of finite groups of Lie type.

\begin{proposition}\label{lem: compactness of affine DL sheaf}
Both $\widetilde{R}_{\dot{w}}^{*}$ and $\widetilde{R}_{\dot{w}}^{!}$ are compact objects in $\shv(\kot_G)$.
\end{proposition}
The objects $R^*_{\dot{w},\theta}$ and $R^!_{\dot{w},\theta}$, however, may not be compact in $\shv(\kot_G)$. They belong to $\fgshv(\kot_G)$.
\begin{proof}
We have the Cartesian diagram
\[
\xymatrix{
\frac{\iw^u \dot{w} \iw^u}{\Ad_\sigma \iw^u}\ar^{\pi_h}[r] \ar_{\pi_d}[d] & \dot{w} \ar[d]\\
\frac{LG_w}{\Ad_\sigma\iw}\ar^{\pr_{\dot{w}}^\sigma}[r] & \bB \mS_k^{\bar w\sigma}.
}\]
Therefore, we have $(\pr_{\dot{w}}^\sigma)^!(\La[\mS_k^{\bar w\sigma}])\cong (\pi_d)_* (\pi_h)^!\La$. By \Cref{compact.generation.admissible.stacks}, the sheaf $(\pi_h)^!\La$ on $\frac{\iw^u \dot{w} \iw^u}{\Ad_\sigma \iw^u}$ is compact (as it is constructible). Now as $\pi_d$ is proper, $(\pi_d)_*$ admits a continuous right adjoint by $(\pi_d)^!$ and therefore preserves compactness. It follows that $(\pr_{\dot{w}}^\sigma)^!(\La[\mS_k^{\bar w\sigma}])$ is compact.

Next, as $\sht^\loc_{\leq w}\to \kot_G$ is ind-pfp proper, the $*$-pushforward preserves compact objects (as it admits a continuous right adjoint). Therefore, it remains to show that both $*$- and $!$-pushforwards along $\sht^\loc_w\to \sht^\loc_{\leq w}$ preserve compact objects. The case of $!$-pushforward is clear, as it is defined as the left adjoint of $!$-pullback. The $*$-forward case follows from \Cref{verdier.duality.placid.stacks} and \Cref{prop: verdier duality on shv for very placid} that says the Verdier duality on quotient very placid stacks preserves compact objects. 
\end{proof}

Two basic results in the classical Deligne-Lusztig theory are the completeness and disjointness of Deligne-Lusztig characters. There are affine generalizations of such results. 
We formulate a disjointness statement here. The affine analogue of the completeness statement will be discussed in \Cref{prop: second characterization of tame LLCategory}.

We assume that $\La$ is an algebraically closed field of characteristic different from $p$, and
say two pairs $(w,\theta: \mS_k^{\bar{w}\sigma}\to\La^\times)$ and $(w',\theta': \mS_k^{\bar{w}'\sigma}\to \La^\times)$ are geometrically conjugate if the pairs $(\bar{w},\theta)$ and $(\bar{w}',\theta')$ are geometrically conjugate in the sense of Deligne-Lusztig  \cite[\textsection{5}]{Deligne.Lusztig}.
The following statement can be regarded as an affine generalization of \cite[Theorem 6.2]{Deligne.Lusztig}. Probably it can be proved by the similar method as in \emph{loc. cit.} but we will give a more conceptual proof of a more general statement in \Cref{SS: a spectral sequence} (which also works in the finite case).

\begin{proposition}\label{prop: aff DL disjointness}
If $(w_1,\theta_1)$ and $(w_2,\theta_2)$ are not geometrically conjugate, then
\[
\Hom_{\shv(\kot_G)}(R^*_{\dot{w}_1,\theta_1},R^*_{\dot{w}_2,\theta_2})=0.
\] 
\end{proposition}

Therefore, it is important to classify pairs $(w,\theta)$ up to geometric conjugacy.
Deligne-Lusztig interpreted such pairs as semisimple elements in a reductive group defined over a finite field. We need another interpretation in order to connect to the local Langlands correspondence. Recall the notion of tame inertia type from \Cref{def: inertia type}.

\begin{lemma}\label{prop: DL geo conj vs tame inertia type}
  Assume that $\La$ is an algebraically closed field. There is a canonical bijection between the set of geometric conjugacy classes of pairs $(w,\theta)$ and the set of tame inertia types.
\end{lemma}

We note that unlike the interpretation from \cite[\textsection{5}]{Deligne.Lusztig}, the bijection in the proposition is independent of any choice and therefore is completely canonical.
\begin{proof}
Let $\hat{S}=\hat{T}/(1-\tau)\hat{T}$, which is the dual torus of $S$ (or equivalently of $\mS_k$) over $\La$. Indeed, 
\[
\xch(\hat{S})=\xch(\hat{T})^\tau=\xcoch(T)^\tau=\xcoch(S)=\xcoch(\mS_k),
\]
equipped with an action $\bar\sigma$ by the (arithmetic) Frobenius of $k_F$, and an action of $W_0$. By \Cref{lem: classification tame inertia type},
tame inertia types are exactly those $\chi: I_F^t\to \hat{S}$ up to $W_0$-conjugacy, such that there is some $\bar{w}\in W_0$ such that $\bar{w}(\bar\sigma(\chi))=\chi^q$.

 Note that under the isomorphism $\mS_k(k)\cong \xcoch(S)\otimes k^\times$, the homomorphism of $\varphi_w$ from \eqref{eq:w-twisted torus} is given by $q\bar{w}\bar\sigma^{-1}-\id$, so
we have 
\begin{equation}\label{eq:Skwsigma as sub}
0\to \mS_k^{\bar{w}\sigma}\to \xcoch(S)\otimes k^\times\xrightarrow{q\bar{w}\bar\sigma^{-1}-\id} \xcoch(S)\otimes k^\times\to 0.
\end{equation}
Using the canonical isomorphism $k^\times\cong \widehat{\bZ}^p\otimes (\bQ/\bZ)$, and the snake lemma, we get from \eqref{eq:Skwsigma as sub} another short exact sequence
\begin{equation}\label{eq:Skwsigma as quotient}
     0\to \xcoch(S)\otimes \widehat{\bZ}^p(1)\xrightarrow{q\bar{w}\bar\sigma^{-1}-\id}\xcoch(S)\otimes \widehat{\bZ}^p(1)\to \mS_k^{\bar{w}\sigma}\to 0.
\end{equation}

Now by \eqref{eq:Skwsigma as quotient} a character $\theta: \mS_k^{\bar{w}\sigma}\to \La^\times$ gives $\xcoch(S)\otimes\widehat{\bZ}^p(1)\to \La^\times$, or equivalently a homomorphism 
\[
\chi: I_F^t\cong \widehat{\bZ}^p(1)\to \xcoch(\hat{S})\otimes\La^\times\cong \hat{S}(\La)
\] 
by \eqref{eq:iota-vs-tau}. Note that by construction $\bar{\sigma}(\bar{w})^{-1}(\bar{\sigma}(\chi))=\chi^q$ and therefore $\chi$ is an inertia type.
In addition, note that $(w,\theta)$ and $(w',\theta')$ are geometrically conjugate if and only if the corresponding $\chi$ and $\chi'$ are $W_0$-conjugate. 

Clearly, the above construction from geometric conjugacy classes of pairs $(w,\theta)$ to the set of inert types can be reversed. E.g. giving $\chi: \xcoch(S)\otimes \widehat{\bZ}^p(1)\to \La^\times$, \eqref{eq:Skwsigma as quotient} says that if $\chi$ is an inertia type, then there is some $(w,\theta)$ such that $\chi$ factors through a character $\theta: \mS_k^{\bar w\sigma}\to \La^\times$. 
The proposition follows.
\end{proof}

\subsubsection{Affine Deligne-Lusztig induction}\label{SS: Aff DL induction} 
The construction of sheaves on $\kot_G$ from \Cref{SS: affine DL sheaf} can be put in a more general content. Let $\rshv(\iw\backslash LG/\iw)$ be the big unipotent affine Hecke category.
Consider the correspondence
\begin{equation}\label{eq-horocycle-correspondence}
\frac{LG}{\Ad_\sigma LG}=\kot_G \xleftarrow{\Nt} \frac{LG}{\Ad_\sigma \iw}=\sht^\loc_{\mI} \xrightarrow{\delta} \iw\backslash LG/\iw,
\end{equation}
which induces a functor
\begin{equation}\label{eq: unipotent ADLI}
\Ch_{G,\phi}^{\unip}:= \Nt_*^{\ind\fg}\circ \delta^{\ind\fg,!}: \rshv(\iw\backslash LG/\iw,\La)\to \rshv(\kot_G,\La)
\end{equation}
which we call the unipotent affine Deligne-Lusztig induction. Notice that as $\delta$ is representable coh. pro-smooth and $\Nt$ is ind-pfp proper,
$\Ch_{LG,\phi}^{\unip}$ restricts to a functor 
\[
\fgshv(\iw\backslash LG/\iw,\La)\to \fgshv(\kot_G,\La).
\] 
In addition, recall that by \Cref{prop: fg sheaves on kotG} $\fgshv(\kot_G,\La)\to \shv(\kot_G,\La)$ is a fully faithful embedding. Thus we may regard $\Ch_{G,\phi}^{\unip}$ as the ind-completion of the restriction of
the functor 
\[
\Nt_*\delta^!: \shv(\iw\backslash LG/\iw)\to \shv(\kot_G)
\] 
to subcategory of finitely generated objects.
As $LG/\Ad_\sigma \iw\to \iw\backslash LG/\iw$ is representable coh. pro-smooth, the base change gives 
\begin{equation}\label{eq: unip affDL induction vs affDL sheaves}
\Ch_{LG,\phi}^{\unip}(\nabla_w)\cong R^*_{w^{-1}},\quad \Ch_{LG,\phi}^{\unip}(\Delta_w)\cong R^!_{w^{-1}}.
\end{equation}
Note here the inverse sign appears due to \Cref{rem: modification direction of Shtuka}.

\begin{example}\label{eq:DL-induction-IC}
When $\La=\overline\bQ_\ell$, we may consider have
\begin{equation*}
C_w:=\Ch_{LG,\phi}^{\unip}(\mathrm{IC}_{w^{-1}}),
\end{equation*}
where $\mathrm{IC}_{w^{-1}}$ is the perverse sheave on $\iw\backslash LG/\iw$ whose $!$-pullback to $LG/\iw$ is the intersection cohomology sheaf of $LG_{\leq w^{-1}}/\iw$.
\end{example}

We will also need another version of affine Deligne-Lusztig induction. Consider
\begin{equation}\label{eq-horocycle-correspondence-pro-unipotent}
\iw^u\backslash LG/\iw^u   \xleftarrow{\delta^u} \frac{LG}{\Ad_\sigma \iw^u}\xrightarrow{\Nt^u} \frac{LG}{\Ad_\sigma LG}=\kot_G.
\end{equation}
We will call the functor
\begin{equation}\label{eq: monodromic ADLI}
\Ch_{LG,\phi}^{\mon}:= (\Nt^u)_*(\delta^u)^!:  \shv_{\mon}(\iw^u\backslash LG/\iw^u,\La)\to \shv(\kot_G,\La)
\end{equation}
the (monodromic) affine Deligne-Lusztig induction. We note that $(\Nt^u)_*^{\mon}=(\Nt^u)_*$.

For every closed sub-indscheme $Z\subset R_{I_F^t,\hat{S}}$,  let
\begin{equation}\label{eq-affine-DL-induction-chi}
R_{\dot{w},Z}^{\mon,!}:=\Ch_{LG,\phi}^{\mon}(\Delta^{\mon}_{\dot{w}^{-1}}(\Ch(\cohdual_Z))),\quad R_{\dot{w},Z}^{\mon,*}:=\Ch_{LG,\phi}^{\mon}(\nabla^{\mon}_{\dot{w}^{-1}}(\Ch(\cohdual_Z))).
\end{equation}
Recall the isogeny $\varphi_{\bar w}: \mS_k\to \mS_k$ as defined in \eqref{eq:w-twisted torus}. 
Let 
\[
\chi_{\varphi_{\bar w}}=\Spec \La[\ker \varphi_{\bar w}]
\] 
be the Pontryagin dual of $\ker \varphi_{\bar w}$, regarded as a closed subscheme of $R_{I_F^t,\hat{S}}$, as in \Cref{ex: Hn-equivariant}. Also recall $R_w^*(M)$ and $R_w^!(M)$ from \eqref{eq: general ADLI of M}.

\begin{lemma}\label{lem: formula of monodromic ADL induction}
We denote by $Z\cap \chi_{\varphi_{\bar w}}$ to be the intersection of $Z$ and $\chi_{\varphi_{\bar w}}$ in $R_{I_F^t,\hat{S}}$, and regard $\cohdual_{Z\cap \chi_{\varphi_{\bar w}}}$ as a $\La[\ker \varphi_{\bar w}]$-module as in \Cref{lem:DL induction for tori}. Then we have
\[
R_{\dot{w},Z}^{\mon,!}\cong R_{\dot{w}}^!(\cohdual_{Z\cap \chi_{\varphi_{\bar w}}}),\quad R_{\dot{w},Z}^{\mon,*}\cong R_{\dot{w}}^*(\cohdual_{Z\cap \chi_{\varphi_{\bar w}}}).
\]
\end{lemma}
\begin{proof}
We have the following diagram with the square Cartesian
\begin{equation}\label{eq: adiwu-to-adiw}
\xymatrix{
\frac{LG}{\Ad_\sigma \iw^u}\ar^{\av_s}[r] \ar_{\delta^u}[d] & \frac{LG}{\Ad_\sigma \iw}\ar_{\av_u}[d]\ar^{\delta}[dr] \ar^{\Nt}[r]& \kot_G\\
\iw^u\backslash LG/\iw^u \ar^{\av_s}[r] & \frac{\iw^u\backslash LG/\iw^u}{\Ad_\sigma \mS_k} \ar[r] & \iw\backslash LG/\iw.
}
\end{equation}
By base change, it is enough to show that 
\begin{align*}
(\av_u)^!(\av_s)_*\nabla_{\dot{w}}^{\mon}(\Ch(\cohdual_Z))\cong (i_w)_*((\pr_{\dot{w}}^\sigma)^!\omega_{Z\cap \chi_{\varphi_w}}[-\ell(w)]),\\
 (\av_u)^!(\av_s)_*\Delta_{\dot{w}}^{\mon}(\Ch(\cohdual_Z))\cong (i_w)_!((\pr_{\dot{w}}^\sigma)^!\omega_{Z\cap \chi_{\varphi_w}}[-\ell(w)]),
\end{align*}
where we recall $\pr_{\dot{w}}^\sigma$ is from \eqref{eq: sigma projection w}.
By \Cref{lem:DL induction for tori}, the $*$-pushforward of $(\pr_{\dot{w}})^!\Ch(\cohdual_Z)[-\ell(w)]$ along the map $\iw^u\backslash LG_w/\iw^u \to \frac{\iw^u\backslash LG_w/\iw^u}{\Ad_\sigma \mS_k}$ followed by the $!$-pullback along the map $LG_w/\Ad_\sigma\iw\to  \frac{\iw^u\backslash LG_w/\iw^u}{\Ad_\sigma \mS_k}$ is isomorphic to $(\pr_{\dot{w}}^\sigma)^!\omega_{Z\cap \chi_{\varphi_{\bar w}}}[-\ell(w)]$. 

Now the first isomorphism directly follows from this fact. The second isomorphism also follows from this fact, using the isomorphism $(\av_s)_*\Delta_{\dot{w}}^{\mon}(\Ch(\cohdual_Z))\cong (\av_s)_!\Delta^{\mon}_{\dot{w}}(\Ch(\cohdual_Z))[\dim S]$ (see \Cref{lem: functors between monodromic categories-5}), and the base change  (as $\av_u$ is coh. pro-unipotent).
\end{proof}

The above proof also gives the following corollary.
\begin{corollary}\label{cor: av upper star av lower star}
Let $\av_s: \frac{LG}{\Ad_\sigma \iw^u}\to \frac{LG}{\Ad_\sigma\iw}$ be as in \eqref{eq: adiwu-to-adiw}. Then we have
\[
(\av_s)^*(\av_s)_*(\Delta_{\dot{w}}^{\mon}(\cohdual_Z))\cong \Delta^{\mon}_{\dot{w}}(\Ch(\cohdual_{Z\cap \chi_{\varphi_{\bar w}}})),\quad  (\av_s)^*(\av_s)_*(\nabla_{\dot{w}}^{\mon}(\cohdual_Z))\cong \nabla^{\mon}_{\dot{w}}(\Ch(\cohdual_{Z\cap \chi_{\varphi_{\bar w}}})).
\]
\end{corollary}

Now we consider some particular cases of $Z$. First, if $Z= R_{I_F^t, \hat{S}}$, then $Z\cap \chi_{\varphi_{\bar w}}=\chi_{\varphi_{\bar w}}$. We have 
\[
\omega_{\chi_{\varphi_{\bar w}}}\cong \La[\mS_k^{\bar{w}\sigma}],
\]
and
\begin{equation}\label{eq: universal monodromic DL induction formula}
R_{\dot{w}, R_{I_F^t, \hat{S}}}^{\mon,!}\cong R^!_{\dot{w}}(\La[\mS_k^{\bar w\sigma}])=\widetilde{R}^!_{\dot{w}},\quad R_{w, R_{I_F^t, \hat{S}}}^{\mon,!}\cong R^!(\La[\mS_k^{\bar w\sigma}])=\widetilde{R}^*_{\dot{w}}.
\end{equation}

Next let $\chi\in R_{I_F^t,\hat{S}}(\La)$ with $Z=\hchi$ its formal completion in $R_{I_F^t,\hat{S}}$. In this case, $\omega_{\hchi\cap \chi_{\varphi_{\bar w}}}$ belongs to $\indcoh(R_{I_F^t,\hat{S}})^{\heartsuit}$. (Note, however, that $\omega_{\chi\cap \chi_{\varphi_{\bar w}}}$ does not belong to the heart in general). 
If $\La$ is an algebraically closed field and $\chi\in \chi_{\varphi_{\bar w}}$ (otherwise $\hchi\cap \chi_{\varphi_{\bar w}}$ is empty), i.e.
\[
\chi=\bar{w}\bar{\sigma}(\chi)^q  \in \hat{S},
\]
then $\hchi\cap \chi_{\varphi_{\bar w}}$ is the connected component of $\chi_{\varphi_{\bar w}}$ that contains $\chi$, and $\omega_{\hchi\cap \chi_{\varphi_{\bar w}}}$ is isomorphic to its structure sheaf, and is a direct summand of $\La[\mS_k^{\bar w\sigma}]$. It follows that $R^{\mon,!}_{\dot{w},\hchi}$ and $R^{\mon,*}_{\dot{w},\hchi}$ are compact. In addition,
let $\theta: \mS_k^{\bar w\sigma}\to \La^\times$ be the character determined by $\chi$. Then
\[
R_{\dot{w},\theta}^*= \widetilde{R}_{\dot{w}}^{*}\otimes_{\La[\mS_k^{\bar w\sigma}]} \theta= R_{\dot{w},\hchi}^{\mon,*}\otimes_{\La[\mS_k^{\bar w\sigma}]} \theta, \quad R_{\dot{w},\theta}^!= \widetilde{R}_{\dot{w}}^{!}\otimes_{\La[\mS_k^{\bar w\sigma}]} \theta = R_{\dot{w},\hchi}^{\mon,!}\otimes_{\La[\mS_k^{\bar w\sigma}]} \theta.
\]

Note that if $\La$ is of characteristic zero, then $\hat{\varphi}_{\bar w}$ is finite \'etale and $\hchi\cap \chi_{\varphi_{\bar w}}=\chi$. In this case, $\omega_{\hchi\cap \chi_{\varphi_{\bar w}}}$ is isomorphic to the skyscraper sheaf of $R_{I_F^t,\hat{S}}$ supported at $\chi$. So we have
\begin{equation}\label{eq: chimon-equal-to-theta}
R_{\dot{w},\hchi}^{\mon,*}\simeq R_{\dot{w},\theta}^*,\quad R_{\dot{w},\hchi}^{\mon,!}=R_{\dot{w},\theta}^!.
\end{equation}
In this case,  $R_{\dot{w},\theta}^*$ and $R_{\dot{w},\theta}^!$ are compact in $\shv(\kot_G,\La)$.
In general, $\hchi\cap \chi_{\varphi_{\bar w}}$ may be non-reduced, and $R_{\dot{w},\hchi}^{\mon,*}$ and $R_{\dot{w},\theta}^*$ (and similarly $R_{\dot{w},\hchi}^{\mon,!}$ and $R_{\dot{w},\theta}^!$) are different. In addition, the objects $R_{\dot{w},\theta}^*$ and $R_{\dot{w},\theta}^!$ may not be compact in $\shv(\kot_G,\La)$ (as already mentioned before).

\begin{example}
Suppose $\La=\overline{\bF}_\ell$. Let $(\mS_k^{\bar w\sigma})_\ell$ be the Sylow $\ell$-subgroup of $\mS_k^{\bar w\sigma}$. Then
if $\chi=u$ is trivial, we have $\hat{u}\cap \chi_{\varphi_{\bar w}}=\Spec \overline{\bF}_\ell[(\mS_k^{\bar w\sigma})_\ell]$.
\end{example}

The following statement can be regarded as the affine analogue of Deligne-Lusztig reduction method (\cite[Theorem 1.6]{Deligne.Lusztig}). Recall from \Cref{lem: monoidal product computation in affine Hecke category} that for every simple reflection in $\widetilde{W}$, we have a closed sub indscheme $\chi_{\hat{\al}_s}\subset R_{I_F^t, \hat{S}}$. Then we have $R^{\mon,?}_{w\sigma(s),\chi_{\hat{\al}_{\sigma(s)}}}$ as in \eqref{eq-affine-DL-induction-chi}.
\begin{lemma}\label{lem: ADL reduction method} 
For $w,w'\in \widetilde{W}$ and $s$ a simple reflection satisfying $w=sw'\sigma(s)$ and $\ell(w)=\ell(w')+2$, we have cofiber sequencs in $\fgshv(\kot_G)$
\[
R^{*}_{w'}\to R^{*}_{w}\to R^{*}_{w'\sigma(s)}\oplus R^{*}_{w'\sigma(s)}[1],\quad R^{!}_{w'\sigma(s)}\oplus R^{!}_{w'\sigma(s)}[-1]\to R^{!}_{w}\to R^{!}_{w'}.
\]
More generally, we have cofiber sequences in $\shv(\kot_G)^\cpt$
\[
\widetilde{R}^{*}_{w'}\to \widetilde{R}^{*}_{w}\to R^{\mon,*}_{w'\sigma(s),\chi_{\hat{\al}_{\sigma(s)}}}[1],\quad R^{\mon,!}_{w'\sigma(s),\chi_{\hat{\al}_{\sigma(s)}}}\to \widetilde{R}^{!}_{w}\to \widetilde{R}^{!}_{w'},
\]
where $R^{\mon,*}_{w'\sigma(s),\chi_{\hat{\al}_{\sigma(s)}}}[1]$ and $R^{\mon,!}_{w'\sigma(s),\chi_{\hat{\al}_{\sigma(s)}}}$ fits into the cofiber sequences
\[
 \widetilde{R}^{*}_{w'\sigma(s)}\to \widetilde{R}^{*}_{w'\sigma(s)}\to R^{\mon,*}_{w'\sigma(s),\chi_{\hat{\al}_{\sigma(s)}}}[1], \quad R^{\mon,!}_{w'\sigma(s),\chi_{\hat{\al}_{\sigma(s)}}}\to \widetilde{R}^{!}_{w'\sigma(s)}\to \widetilde{R}^{!}_{w'\sigma(s)}.
\]
\end{lemma}
\begin{proof}
First, the last two cofiber sequences follow from \eqref{eq: fiber of standard to costandard-2}. 
We prove the cofiber sequence $\widetilde{R}^{*}_{w'}\to \widetilde{R}^{*}_{w}\to R^{\mon,*}_{w'\sigma(s),\chi_{\hat{\al}_{\sigma(s)}}}[1]$. The rest ones can proved similarly. (In fact the first two can be deduced from the last two.)
We consider the following commutative diagram
\[
\xymatrix{
                                                                              & \frac{LG_s\times^{\iw^u} LG_{w'}\times^{\iw^u} LG_{\sigma(s)}}{\Ad_\sigma\iw^u}\ar[d] \ar^-{\cong}[r]& \frac{ LG_{w'}\times^{\iw^u} LG_{\sigma(s)}\times^{\iw^u}LG_{\sigma(s)}}{\Ad_\sigma\iw^u} \ar[d]\\
\frac{LG_w}{\Ad_\sigma\iw^u}\ar^-{\cong}[r]\ar_{\av_s}[d]& \frac{LG_s\times^{\iw} LG_{w'}\times^{\iw} LG_{\sigma(s)}}{\Ad_\sigma\iw^u}\ar[d] \ar^-{\cong}[r]&  \frac{ LG_{w'}\times^{\iw} LG_{\sigma(s)}\times^{\iw^u}LG_{\sigma(s)}}{\Ad_\sigma\iw} \ar[d] \\
\sht^\loc_{w}\ar^-{\cong}[r]&\sht^{\loc}_{s,w',\sigma(s)} \ar^-{\cong }[r]& \sht^{\loc}_{w',\sigma(s),\sigma(s)}
}\]
and \Cref{lem: monoidal product computation in affine Hecke category} \eqref{lem: monoidal product computation in affine Hecke category-1} (letting $\mL=\mL'=\widetilde{\Ch}$), we see that
\[
(\av_s)_*(\delta^u)^!\widetilde{\nabla}_{\dot{w}}^{\mon}\cong (\av_s)_*(\delta^u)^!(\widetilde{\nabla}_{\dot{w}'}^{\mon}\star^{ u}\widetilde{\nabla}_{\sigma(\dot{s})}^{\mon}\star^{ u}\widetilde{\nabla}_{\sigma(\dot{s})}^{\mon})
\]
Now the claim follows from \Cref{lem: monoidal product computation in affine Hecke category} \eqref{lem: monoidal product computation in affine Hecke category-2}.
\end{proof}

We will also let
\begin{equation}\label{eq:ADLI of tilting sheaf}
\widetilde{R}^{T}_{\dot{w}}:= \Ch_{LG,\phi}^{\mon}(\widetilde{\Til}^{\mon}_{\dot{w}}),\quad R^{\mon,T}_{\dot{w},\hchi}:= \Ch_{LG,\phi}^{\mon}(\Til^{\mon}_{\dot{w},\hchi}).
\end{equation}
It follows by definition that $\widetilde{R}^{T}_{\dot{w}}$ (resp. $R^{\mon,T}_{\dot{w},\hchi}$) admits a filtration with associated graded by $\widetilde{R}_{\dot{w}'}^{*}$ (resp. by $R_{\dot{w}',\hchi}^{\mon,*}$) and another filtration with associated graded by $\widetilde{R}_{\dot{w}'}^{!}$ (resp. by $R_{\dot{w}',\hchi}^{\mon,!}$).

We will need the following computations to understand matching objects under the categorical local Langlands correspondence. Assume that $uw$ is a minimal length element in its $\sigma$-conjugacy class as in \Cref{thm: reduction to min length elements} \eqref{thm: reduction to min length elements-1}. 
Let $\breve\mP_{\breve\bff}, L_{\breve\bff}, B_{L_{\breve\bff}}$ be as in \Cref{prop-local-shtuka-uw}. Let $U_{L_{\breve\bff}}\subset B_{L_{\breve\bff}}$ be the unipotent radical. We let $P_b=P_{\dot{w},\breve\bff}$, which is a parahoric subgroup of $G_b(F)$ with Levi quotient $L_b:=L_{\breve\bff}(k)^{\sigma_{\dot{w}}}$.

\begin{lemma}\label{lem: ADLI for convolution sheaves}
For every $\mF\in \shv_{\mon}(\iw^u\backslash LG_{W_{\breve\bff}}/\iw^u)\cong \shv_{\mon}(U_{L_{\breve\bff}}\backslash L_{\breve\bff} /U_{L_{\breve\bff}})$, we have
\[
\Ch_{LG,\phi}^{\mon}(\mF\star^u \widetilde{\Delta}^{\mon}_{\dot{w}})\cong (i_b)_!(\cind_{P_b}^{G_b(F)}V)[-\langle 2\rho,\nu_b\rangle],
\]
where $V\in \rep(L_b)$ (regarded as a $P_b$-representation by inflation) is obtained as the Deligne-Lusztig induction of $\mF$ along
\[
U_{L_{\breve\bff}}\backslash L_{\breve\bff} /U_{L_{\breve\bff}}\leftarrow \frac{L_{\breve\bff}}{\Ad_{\sigma_{\dot{w}}}U_{L_{\breve\bff}}}\rightarrow \frac{L_{\breve\bff}}{\Ad_{\sigma_{\dot{w}}}L_{\breve\bff}}.
\]
Similarly, we have
\[
\Ch_{G,\phi}^{\mon}(\mF\star^u \widetilde{\nabla}^{\mon}_{\dot{w}})\cong (i_b)_*(\cind_{P_b}^{G_b(F)}V)[-\langle 2\rho,\nu_b\rangle].
\]
\end{lemma}
\begin{proof}
Note that we have the following commutative diagram with squares labeled by $(\mathrm{X})$ Cartesian
\[
\xymatrix{
U_{L_{\breve\bff}}\backslash L_{\breve\bff} /U_{L_{\breve\bff}}\times \mS_k\ar^m[r] \ar@{}[dr] | {(\mathrm{X})}&U_{L_{\breve\bff}}\backslash L_{\breve\bff} /U_{L_{\breve\bff}}&\ar[l]  \frac{L_{\breve\bff}}{\Ad_\sigma  U_{L_{\breve\bff}}}\ar[r]& \frac{L_{\breve\bff}}{L_{\breve\bff}}&\\
\iw^u\backslash L^+\breve\mP_{\breve\bff}\times^{\iw^u}LG_w/\iw^u\ar^-{m^u}[r]\ar[d]\ar^{f_1}[u]& \iw^u\backslash LG_{W_{\breve\bff}w}/\iw^u\ar[d]\ar_{f_2}[u]\ar@{}[dr] | {(\mathrm{X})}& \frac{LG_{W_{\breve\bff}w}}{\Ad_\sigma \iw^u}\ar[l]\ar[d]\ar[r]\ar[u]&\frac{LG_{W_{\breve\bff}w}}{\Ad_\sigma L^+\breve\mP_{\breve\bff}} \ar[r]\ar[u]& \kot_{G,b}\ar^{i_b}[d] \\
\iw^u\backslash LG\times^{\iw^u}LG/\iw^u\ar^-{m_u}[r]&  \iw^u\backslash LG/ \iw^u &\ar_-{\delta^u}[l]  \frac{LG}{\Ad_\sigma\iw^u}\ar^-{\Nt^u}[rr] && \kot_G
}\]
Here the map $f_1$ is given by $\iw^u\backslash L^+\breve\mP_{\breve\bff}\times^{\iw^u}LG_w/\iw^u\cong \iw^u\backslash L^+\breve\mP_{\breve\bff}/\iw^u\times  \mS_k\dot{w} (\iw^u\cap \dot{w}^{-1} \iw^u \dot{w})\backslash \iw^u\to U_{L_{\breve\bff}}\backslash L_{\breve\bff} /U_{L_{\breve\bff}}\times \mS_k$.
Recall that $\ell(w)=\langle 2\rho, \nu_b\rangle$. Now the statement follows from base change.
\end{proof}

\begin{corollary}\label{lem:ADLV-sheaf-minimal-length}
  Assume that $w\in\widetilde{W}$ is a minimal length element in its $\sigma$-conjugacy class. Then there are canonical isomorphisms
  \[
  \widetilde{R}_{\dot{w}}^* \simeq (i_{b})_{*}\cind_{P_b}^{G_b(F)}(\widetilde{R}^{b,*}_{\dot{u}})[-\langle2\rho,\nu_b\rangle], \quad \widetilde{R}^!_{\dot{w}}\simeq (i_{b})_{!}\cind_{P_b}^{G_b(F)}(\widetilde{R}^{b,!}_{\dot{u}})[-\langle2\rho,\nu_b\rangle], 
  \]
where $P_b$ is a parahoric subgroup of $G_b(F)$ and $\widetilde{R}^{b,*}_{\dot{u}} \in \rep(P_b,\La)$ (resp. $\widetilde{R}^{b,!}_{\dot{u}}\in \rep(P_b,\La)$) is a Deligne-Lusztig induction of the Levi quotient of $P_b$.

In particular, when $w$ is $\sigma$-straight, giving $b\in B(G)$, we have
\[
 \widetilde{R}_{\dot{w}}^* \simeq (i_{b})_{*}\cind_{I_b^u}^{G_b(F)}\La, \quad
\]
where $I_b^u$ is the pro-$p$-radical of $I_b$.
We have similarly version for  $R^*_{w}$ and $R^!_{w}$. 
In particular,  when $w$ is a $\sigma$-straight element corresponding to $b$ and $\theta$ is trivial, then
\[
R_w^*\simeq i_{b,*}(\cind_{I_b}^{G_b(F)}\La[-\langle2\rho,\nu_b\rangle]),\quad R^!_w\simeq i_{b,!}(\cind_{I_b}^{G_b(F)}\La[-\langle 2\rho,\nu_b\rangle]).
\]
\end{corollary}
Note that this corollary is consistent with \Cref{lem: stalk of ADLS}.
\begin{proof}
By  \Cref{lem: invariant moduli of Shtuka}, we may assume that $w$ is as in \Cref{thm: reduction to min length elements} \eqref{thm: reduction to min length elements-1}. Now we apply \Cref{lem: ADLI for convolution sheaves} to conclude.
\end{proof}

Recall that for a compactly generated category we have the Chern character \Cref{rem-Chern-character}.
We have the following affine analogue of (\cite[Theorem 1.6]{Deligne.Lusztig}). We thank Xuhua He for drawing our attention to the possibility that such a statement could be true.
\begin{proposition}\label{lem: chern character of ADLV}
Under either decomposition of $\mathrm{tr}(\shv(\kot_G))$ from \Cref{cor: decomposition HH of shvkot}, we have
\[
\mathrm{ch}(\widetilde{R}^*_{\dot{w}})=\mathrm{ch}(\widetilde{R}^!_{\dot{w}})\in C_c(G_b(F),\La)_{G_b(F)},
\] 
where $b\in B(G)$ is the unique element matching the Newton point and the Kottwitz invariant of $w$.
\end{proposition}
\begin{proof}
Let $\shv_{\mon}(\iw^u\bs LG/\iw^u)'\subset \shv_{\mon}(\iw^u\bs LG/\iw^u)$ be as in \Cref{cor: standard and costandard same K class}. Note that 
\[
\Ch_{LG,\phi}^{\mon}: \shv_{\mon}(\iw^u\bs LG/\iw^u)'\to \shv(\kot_G)^\cpt
\]
by \Cref{lem: compactness of affine DL sheaf}. It follows from  \Cref{cor: standard and costandard same K class} that $\widetilde{R}^*_{\dot{w}}$ and $\widetilde{R}^!_{\dot{w}}$ have the same class in $K_0(\shv(\kot_G)^\cpt)$, and therefore
\[
\mathrm{ch}(\widetilde{R}^*_{\dot{w}})=\mathrm{ch}(\widetilde{R}^!_{\dot{w}})\in \mathrm{tr}(\shv(\kot_G)).
\]

Next, we notice that by \Cref{lem: ADL reduction method}, if $w$ and $w'$ are $\sigma$-conjugate, then  $\mathrm{ch}(\widetilde{R}^*_{\dot{w}})= \mathrm{ch}(\widetilde{R}^*_{\dot{w}'})$ and similarly  $\mathrm{ch}(\widetilde{R}^!_{\dot{w}})= \mathrm{ch}(\widetilde{R}^!_{\dot{w}'})$. Therefore, we may assume that $w$ is as in \Cref{thm: reduction to min length elements} \eqref{thm: reduction to min length elements-1}. In this case, the claim follows from \Cref{lem:ADLV-sheaf-minimal-length}. 
\end{proof}

We have analogue of \Cref{lem: Cohdual Ch vs Ch sw dual}.
\begin{proposition}\label{prop: verdual Ch vs Ch sw dual}
Let $\mF\in \shv_\mon(\iw^u\bs LG/\iw)$ which admits a right dual $\mF^\vee$ and suppose $\Ch_{LG,\phi}^{\mon}(\mF)$ is compact in $\shv(\kot_G)$. Then
we have a canonical isomorphism 
\[
(\verd^{\can}_{\kot_G})^\cpt(\Ch_{LG,\phi}^{\mon}(\mF))\cong \Ch_{LG,\phi}^{\mon}(\mathrm{sw}(\mF^\vee)).
\]
\end{proposition}
\begin{proof}
We make use of the commutative diagram \eqref{eq: adiwu-to-adiw}. 
Let $\La^{\can}_{\frac{\iw^u\bs LG/\iw^u}{\Ad_\sigma \mS_k}}$ be the generalized constant sheaf on $\frac{\iw^u\bs LG/\iw^u}{\Ad_\sigma \mS_k}$ obtained by $!$-pullback 
along $\frac{\iw^u\bs LG/\iw^u}{\Ad_\sigma \mS_k}\to \iw\bs LG/\iw$ of the generalized constant sheaf $\La^{\can}_{\iw\bs LG/\iw}$. As explained after \eqref{eq: canonical constant sheaf on enhanced hecke}, we have
$\La^{\can}_{\iw^u\bs LG/\iw^u}=(\av_s)^!\La^{\can}_{\frac{\iw^u\bs LG/\iw^u}{\Ad_\sigma \mS_k}}[-4\dim \mS_k]$. 

Using the base change, it is enough to show that
\begin{equation}\label{eq: verdual Ch vs Ch sw dual-aux}
\verd^{\can}_{\frac{\iw^u\bs LG/\iw^u}{\Ad_\sigma \mS_k}}((\av_s)_*\mF)\cong (\av_s)_*(\sw(\mF^\vee)).
\end{equation}
Let $\mG\in \shv(\frac{\iw^u\bs LG/\iw^u}{\Ad_\sigma \mS_k})^\cpt$.
On the one hand, we have
\begin{eqnarray*}
\Hom(\mG,\verd^{\can}_{\frac{\iw^u\bs LG/\iw^u}{\Ad_\sigma \mS_k}}((\av_s)_*\mF))&=&\Hom((\av_s)_*\mF, \verd^{\can}_{\frac{\iw^u\bs LG/\iw^u}{\Ad_\sigma \mS_k}}(\mG))\\
&=&\Hom(\mF, \verd^{\can}_{\iw^u\bs LG/\iw^u}((\av_s)^*(\mG)))[\dim \mS_k]\\
&=&\Hom( \mF\otimes^{\can}(\av_s)^*\mG, \consdual_{\iw^u\bs LG/\iw^u})[\dim \mS_k]
\end{eqnarray*}
On the other hand, we have
\begin{eqnarray*}
\Hom(\mG,(\av_s)_*(\sw(\mF^\vee)))&=&\Hom((\av_s)^*\mG, \sw(\mF^\vee))\\
&=&\Hom(\sw(\mF)\star^u(\av_s)^*\mG,\widetilde\Delta^{\mon}_e)
\end{eqnarray*}
Now \eqref{eq: verdual Ch vs Ch sw dual-aux} follows form \Cref{lem: semirigid dual vs verdier dual in monodromic affine Hecke}.
\end{proof}

We discuss how to produce $R_{w,\theta}^*$ and $R_{w,\theta}^!$ directly via the affine Deligne-Lusztig induction for equivariant categories. For simplicity, we will assume that $\La$ is a field in the sequel. Let $p'$ be the product of $p$ and the characteristic exponent of $\La$ (so $p'=p$ if $\La$ is a field of characteristic zero and otherwise $p'=p\cdot\mathrm{char} \La$).  Then every prime-to-$p$ finite order character $\chi: T^p\mS_k\to \La^\times$ has order coprime to $p'$.

We can consider analogue of \eqref{eq-horocycle-correspondence}
\begin{equation}\label{eq-horocycle-correspondence-chi}
\iw^{[n]}\backslash LG/\iw^{[n]} \xleftarrow{\delta^{[n]}} \frac{LG}{\Ad_\sigma \iw^{[n]}} \xrightarrow{\Nt^{[n]}} \frac{LG}{\Ad_\sigma LG}=\kot_G.
\end{equation}
and define the affine Deligne-Lusztig induction as
\[
\Ch_{LG,\phi}^{[n]}:=(\Nt^{[n]})^{\ind\fg}_*(\delta^{[n]})^{\ind\fg,!}: \rshv(\iw^{[n]}\backslash LG/\iw^{[n]})\to \rshv(\kot_G).
\] 
Here we note that $(\Nt^{[n]})^{\ind\fg}_*$ is defined thanks to \eqref{eq-rshv-for-placid}. Namely, $\Nt^{[n]}=\Nt\circ \varphi^n$, where $\varphi^n:  \frac{LG}{\Ad_\sigma \iw^{[n]}}\to  \frac{LG}{\Ad_\sigma \iw}$ is a $\bB \mS_k[n]$-gerbe with $n$ in invertible in $\La$.
Note that $(\varphi^n)^{\ind\fg}_*$ is both left and right adjoint of $(\varphi^n)^{\ind\fg,!}$.

\begin{proposition}\label{prop:affine DL induction chi}
We have $\Ch_{LG,\phi}^{[n]} \Delta_{w,\chi}= \Ch_{LG,\phi}^{[n]} \nabla_{w,\chi}=0$ unless $\chi\circ \varphi_w: T^p\mS_k\to \La^\times$ is trivial, in which case 
$\chi$ gives a character $\theta: \mS_{k}^{\bar{w}\sigma}\to \La^\times$ by \eqref{eq:Skwsigma as quotient} and
\[
\Ch_{LG,\phi}^{[n]} \Delta_{w,\chi}\cong R^!_{w,\theta},\quad \Ch_{LG,\phi}^{[n]} \nabla_{w,\chi}\cong R^*_{w,\theta}.
\]
\end{proposition}
\begin{proof}
Let
\[
((\mS_k^{[n]})^{\bar{w}\sigma})'=\left\{g\in \mS_k^{[n]}\mid g^{-1}\bar{w}\sigma(g)\in \mS_k^{[n]}[n]\right\}.
\]
There are exact sequences
\begin{equation}\label{eq:ext of Svarphi to Swsigma}
1\to (\mS_k^{[n]})^{\bar{w}\sigma}\subset  ((\mS_k^{[n]})^{\bar{w}\sigma})'\xrightarrow{g\mapsto g^{-1}\bar{w}\sigma(g)} \mS_k^{[n]}[n]\to 1, \quad 1\to \mS_k^{[n]}[n]\to  ((\mS_k^{[n]})^{\bar{w}\sigma})'\xrightarrow{\varphi^n} (\mS_k^{[n]})^{\bar{w}\sigma}\to 1.
\end{equation}
Note that the composed map $\mS_k^{[n]}[n]\to  ((\mS_k^{[n]})^{\bar{w}\sigma})'\to \mS_k^{[n]}[n]$ from the above two sequences is given by $\varphi_w|_{\mS_k^{[n]}[n]}$. 

Consider the commutative diagram
\[
\xymatrix{
&\iw^{[n]}\backslash LG/\iw^{[n]}&\ar_{\delta^{[n]}}[l] \frac{LG}{\Ad_\sigma \iw^{[n]}}\ar^{\varphi^n}[r]& \frac{LG}{\Ad_\sigma\iw}\\
&\iw^{[n]}\backslash LG_w/ \iw^{[n]}\ar^{i_w}[u]\ar_{\pr_s}[d]\ar_{\pr_w^{[n]}}[dl]&\ar_{\delta^{[n]}}[l] \frac{LG_w}{\Ad_\sigma \iw^{[n]}}\ar^{\varphi^n}[r]\ar_{i_w}[u]\ar^{\pr_s}[d]& \frac{LG_w}{\Ad_\sigma\iw}\ar_{i_w}[u]\ar^{\pr_s}[d]\\
\bB \mS_k[n]&\mS_k^{[n]}\backslash \dot{w}\mS_k/\mS_k^{[n]}\ar[l]&\ar[l] \frac{\mS_k \dot{w}}{\Ad_\sigma \mS_k^{[n]}}\cong \bB ((\mS_k^{[n]})^{\bar{w}\sigma})'\ar^-{\varphi^n}[r] & \frac{\dot{w}\mS_k}{\Ad_\sigma\mS_k}\cong \bB \mS_k^{\bar{w}\sigma} \\
&\mS_k^{[n]}\backslash \dot{w}\mS_k^{[n]}/\mS_k^{[n]}\cong \bB \mS_k^{[n]}\ar[u] & \frac{ \mS_k^{[n]} \dot{w}}{\Ad_\sigma(\mS_k^{[n]})}\cong \bB (\mS_k^{[n]})^{\bar{w}\sigma}.\ar[l] \ar_{\psi_1}[u] &\\
}
\]
Note that all the squares are Cartesian except the left middle one. 

As all the functors below preserves $\fgshv$, we could omit $\ind\fg$ from the superscript when considering pushforwards or pullbacks.
We first compute $(\varphi^n)_*(\delta^{[n]})^!(\pr_w^{[n]})^!\chi[-\ell(w)-2\dim S] \in \fgshv(\frac{LG_w}{\Ad_\sigma \iw})$. Using the base change, it is canonically isomorphic to $(\pr_s)^! M[-\ell(w)]$, where $M$ is the following representation of $\mS_k^{\bar{w}\sigma}$ (regarded as a sheaf on $\bB \mS_k^{\bar{w}\sigma}$): We inflate the character $\chi$ of $\mS_k^{[n]}[n]=\mS_k[n]$ to a representation of $((\mS_k^{[n]})^{\bar{w}\sigma})'$ via the second exact sequence in \eqref{eq:ext of Svarphi to Swsigma} and then taking the (derived) invariants with respect to the subgroup $\mS_k^{[n]}[n]\subset ((\mS_k^{[n]})^{\bar{w}\sigma})'$ from the first exact sequence in \eqref{eq:ext of Svarphi to Swsigma}. Therefore, the resulting representation $M$ of $\mS_k^{\bar{w}\sigma}$ is non-zero if and only if $\chi\circ \varphi_w$ is trivial, in which case it is a character $\theta$ of of $\mS_k^{\bar{w}\sigma}$.

It remains to show that 
\[
(\varphi^n)_*(\delta^{[n]})^!(i_w)_?(\pr_w^{[n]})^!\chi\cong (i_w)_?(\varphi^n)_*(\delta^{[n]})^!(\pr_w^{[n]})^!\chi
\] 
for $?=*$ and $!$. But as $\delta^{[n]}$ is coh. pro-unipotent, $(\delta^{[n]})^!$ commutes with both $*$- and $!$-pushforward along pfp morphisms. In addition, as mentioned before, $(\varphi^n)_*$ is both the left and the right adjoint of $(\varphi^n)^!$, and therefore also commutes with both $*$- and $!$-pushforward along pfp morphisms.
\end{proof}

We finish our general discussion of affine Deligne-Lusztig inductions by the following remark.
\begin{remark}\label{rem: grading on ADLI}
Note that as $\pi_0(\iw^u\backslash LG/\iw^u)=\pi_0(LG)=\pi_1(G)_{I_F}$, there is a decomposition
\[
\shv_{\mon}(\iw^u\backslash LG/\iw^u)=\sqcup_{\al\in \pi_1(G)_{I_F}} \shv_{\mon}(\iw^u\backslash LG_\al/\iw^u),
\]
where $LG_\al$ is the connected component of $LG$ corresponding to $\al\in \pi_1(G)_{I_F}$.
On the other hand, there is a decomposition of $\shv(\kot_G)$ as in \eqref{eq: decomp shvkotG by connected components}.

Since \eqref{eq-horocycle-correspondence-pro-unipotent} induces the following map of connected components
\[
\pi_0(\iw^u\backslash LG/\iw^u)\cong \pi_0(\frac{LG}{\Ad_\sigma\iw^u})\to \pi_0(\kot_G),
\]
which can be identified with the natural quotient map $\pi_1(G)_{I_F}\to \pi_1(G)_{\Ga_F}$, we see that
the functor $\Ch_{G,\phi}^{\mon}$ sends the $\shv_{\mon}(\iw^u\backslash LG_\al/\iw^u)$ to $\shv(\kot_{G,\bar\al})$, where $\bar\al$ is the image of $\al$ in $\pi_1(G)_{\Ga_F}$.
\end{remark}

\subsubsection{A geometric Mackey formula}\label{SS: a spectral sequence}
In representations theory of finite groups, the Mackey formula expresses the composition of an induction functor followed by a restriction functor in terms of a sum of functors which are compositions of restriction functors followed by induction functors.
Our next goal is to discuss an analogue of this result, which will allow us to compute the hom space between certain objects in $\shv(\kot_G)$. 

Although we can directly prove a most general version of the result we need, to benefit readers, we will start with a relatively easy version and explain necessary modifications for variants.

\begin{lemma}\label{lemma-explicit-filtration-unit-affine-flag-unip}
Let $\mF\in \shv(\iw^u\backslash LG/\iw^u)$. There is a filtration of 
\[
(\eta)_\flat((m)^!\mF) \in \shv(\iw^u\backslash LG/\iw\times \iw\backslash LG/\iw^u)
\]
with the associated graded being 
\[
\bigoplus_{w}\nabla_{w}^l\boxtimes_{\La}(\nabla_{w^{-1}}^r\star^{u} \mF),\quad w\in \widetilde{W}.
\]
\end{lemma}
We refer to the paragraph above \Cref{lem: l and r pullback compatible with conv} for the notations $\nabla_{w}^l$ and $\nabla_{w^{-1}}^r$.

\begin{proof}
The lemma in fact follows from \Cref{cor-filtration-unit}. See the remark below. To benefit readers, however, we make the abstract formalism concrete in this case.

We first deal with the case when $\mF=\mathbf{1}_{\iw^u}$ is  the unit of $\shv(\iw^u\backslash LG/\iw^u)$. 
For each $w$, let 
\[
a_{?}\colon \iw^u\backslash LG_{?}\times^{\iw}LG/\iw^u\to \iw^u\backslash LG\times^{\iw}LG/\iw^u,\quad ?=w \mbox{ or } \leq w
\] 
be the pfp (locally) closed embedding. Let $m_{\leq w}= m\circ a_{\leq w}$ and similarly let $m_w=m\circ a_w$.
Then $m^!(\mathbf{1}_{\iw^u})$ admits a filtration with associated graded being $(a_w)_*(m_w)^!\mathbf{1}_{\iw^u}$. We claim that
\begin{equation}\label{eq: stratified unit by Schubert-unip}
\eta_\flat\bigl((a_w)_*((m_w)^!\mathbf{1}_{\iw^u})\bigr)\cong \nabla_{w}^l\boxtimes_{\La}\nabla_{w^{-1}}^r.
\end{equation}
To see this, we can perform the base change along $LG/\iw\times \iw\backslash LG/\iw^u\to \iw^u\backslash LG/\iw\times\iw\backslash LG/\iw^u$, and consider the following sequence of morphisms
\begin{multline}\label{mult: LGw diagonal-unip}
LG_w/\iw\stackrel{g\iw\mapsto (g,g^{-1}\iw^u)}{\xrightarrow{\hspace*{2.5cm}}} LG_w\times^{\iw} LG_w/\iw^u \\ \stackrel{\widetilde{\eta_w}}{\longrightarrow} LG_w/\iw\times\iw\backslash LG_w/\iw^u \hookrightarrow LG/\iw \times \iw\backslash LG/\iw^u.
\end{multline}
Then the base change of $\eta_\flat(a_w)_*(m_w)^!\mathbf{1}_{\iw^u}$ is obtained from the dualizing sheaf on $LG_w/\iw$ by $*$-pushforward along the first map, followed by $\flat$-pushforward along the second map, and then followed by $*$-pushforward along the last map. 

We now consider the following commutative diagram
\[
\xymatrix{
LG_w/\iw\ar[d] \ar[r] & LG_w/\iw^{(r)}\times^{\iw/\iw^{(r)}} LG_w/\iw^u \ar[r] & LG_w/\iw\times (\iw/\iw^{(r)})\backslash LG_w/\iw^u\\
LG_w\times^{\iw} LG_w/\iw^u \ar^-{\widetilde{\eta_1}}[r] \ar@{=}[ur]& LG_w/\iw^{(r)} \times^{\iw/\iw^{(r)}} \iw^{(r)}\backslash LG_w/\iw^u\ar^-{\widetilde{\eta_2}}[r]\ar_{g_1}[u] & LG_w/\iw\times\iw\backslash LG_w/\iw^u\ar_{g_2}[u],
}\]
where $\iw^{(r)}$ is a small enough congruence subgroup of $\iw$ (so $\iw^{(r)}\subset \iw^u \subset \iw$) such that the left action of $\iw^{(r)}$ on $LG_w/\iw^u$ is trivial. Note that the composition of the bottom arrows is the  map $\widetilde{\eta_w}$ in \eqref{mult: LGw diagonal}.
As $\iw^{(r)}$ is coh. pro-unipotent, $(g_i)_\flat$ are equivalences. In addition, $(\widetilde{\eta_2})_\flat$ is isomorphic to $(\widetilde{\eta_2})_*$ up to a shift (and a Tate twist). Therefore, it is enough to compute the $*$-pushforward of the dualizing sheaf of $LG_w/\iw$ along the top arrows, which is $\consdual_{LG_w/\iw}\times \consdual_{(\iw/\iw^{(r)})\backslash LG_w/\iw^u}[2d]$, where $d= \dim (\iw/\iw^{(r)})\backslash LG_w/\iw^u$.
The lemma follows when $\mF=\mathbf{1}_{\iw^u}$.

To deal with general $\mF$, we consider the diagram with both squares Cartesian
\[
\xymatrix{
\ar_{\id\times m^u}[d]\iw^u\backslash LG/\iw \times \iw\backslash LG\times^{\iw^u}LG/\iw^u&\ar_-{\eta\times\id}[l]\iw^u\backslash LG\times^{\iw} LG\times^{\iw^u}LG/\iw^u\ar^-{m\times\id}[r]\ar^{\id\times m^u}[d]&\iw^u\backslash LG\times^{\iw^u}LG/\iw^u\ar^{m^u}[d]\\
\iw^u\backslash LG/\iw \times \iw\backslash LG/\iw^u&\ar_{\eta}[l]\iw^u\backslash LG\times^{\iw} LG/\iw^u\ar^-{m}[r]&\iw^u\backslash LG/\iw^u.
}\]
Now we form the usual ``twisted product" $\mathbf{1}_{\iw^u}\widetilde\boxtimes \mF= (\eta^u)^!(\mathbf{1}_{\iw^u}\boxtimes_{\La} \mF)$ on $\iw^u\backslash LG\times^{\iw^u} LG/\iw^u$. Its $*$-pushfoward along $m^u$ (the rightmost vertical map) is $\mathbf{1}_{\iw^u}\star^{u}\mF=\mF$. By 
stratifying the first $LG$-factor in $\iw^u\backslash LG\times^{\iw} LG\times^{\iw^u}LG/\iw^u$ by $LG_w$, and using the base change between $*$-pushforwards and $!$-pullbacks (as built in the sheaf theory $\shv$) and the base change between $*$-pushforwards and $\flat$-pushforwards (by \Cref{lem-indproper-prosmooth-push-basechange}), we see that we reduce to the case $\mF=\mathbf{1}_{\iw^u}$.
\end{proof}

\begin{remark}
We explain why \Cref{lemma-explicit-filtration-unit-affine-flag} follows from the abstract formalism \Cref{cor-filtration-unit}. As in \Cref{prop: categorical property of affine Hecke}, we have
\[
\shv(\iw^u\backslash LG/\iw)\otimes_\La \shv(\iw\backslash LG/\iw^u)\cong \shv(\iw^u\backslash LG/\iw\times\iw\backslash LG/\iw^u)
\] 
and that $\shv(\iw^u\backslash LG/\iw)$ is dualizable (as $\La$-linear categories) with $\shv(\iw\backslash LG/\iw^u)$ its dual. The unit of the duality datum is nothing but $\eta_\flat(m^!\mathbf{1}_{\iw^u})$. Then one can use \Cref{cor-filtration-unit} to conclude. 
\end{remark}

\begin{proposition}\label{lem-spectral-sequence-hom-unip}
Assumptions are as in \Cref{lemma-explicit-filtration-unit-affine-flag}.
Let $\mF_1$ be an object in $\fgshv(\iw\backslash LG/\iw)$ and $\mF_2$ an object in $\fgshv(\iw^u\backslash LG/\iw^u)$. 
Then there is a filtration on the $\La$-module
\[
\Hom_{\shv(\kot_G)}((\Nt)_*(\delta)^!\mF_1,(\Nt^u)_*(\delta^u)^!\mF_2)
\] 
with the associated graded being $\Hom_{\shv(\iw\backslash LG/\iw^u)}\Bigl(\mF_1\star \Delta_{\sigma(w)}^r, \mF_2\star^u\nabla_{\sigma(w)}^l\Bigr)$ for $w\in \widetilde{W}$. In particular, there is a spectral sequence with
\[
E_1^{p,q}\simeq \bigoplus_{\ell(w)=-p}\mathrm{Ext}^{q+p}_{\shv(\iw\backslash LG/\iw^u)}\Bigl(\mF_1\star \Delta_{\sigma(w)}^r, \mF_2\star^u\nabla_{\sigma(w)}^l\Bigr),
\]
and with abutment $\mathrm{Ext}^*_{\shv(\kot_G)}(\Nt_*(\delta^!\mF_1),(\Nt^u)_*((\delta^u)^!\mF_2))$.
\end{proposition}
\begin{proof}
This proposition also follows from  the abstract formalism \Cref{cor-hom-space-of-compact-object-in-categorical-trace}, \Cref{rem-approx-identity}, and \Cref{cor-filtration-unit}. Again, we make the abstract formalism concrete.

Consider the following commutative diagram
\begin{small}
\[
\xymatrix{ 
\iw\backslash LG/ \iw&\ar_-{m^u}[l]\iw\backslash LG\times^{\iw^u} LG/\iw\ar^-{\eta^u}[r]&\iw\backslash LG/\iw^u\times \iw^u\backslash LG/\iw\ar^{\mathrm{sw}\circ(\sigma\times\id)}[r]&\iw^u\backslash LG/\iw\times \iw\backslash LG/\iw^u\\
 \frac{LG}{\Ad_\sigma \iw}\ar^{\Nt}[rd]\ar_-{\delta}[u] &\ar[u]\ar_{m^u}[l]\frac{LG\times^{\iw^u}LG}{\Ad_\sigma \iw}\ar_-\cong^{\pFr}[r]\ar[dr]&\frac{LG\times^{\iw}LG}{\Ad_\sigma\iw^u}\ar^{m}[d]\ar[r] &\ar_-{\eta}[u]\iw^u\backslash LG\times ^{\iw} LG/\iw^u \ar^{m}[d]\\
 &\kot_G  &\ar_{\Nt^u}[l] \frac{LG}{\Ad_\sigma \iw^u}\ar^-{\delta^u}[r] &\iw^u\backslash LG/\iw^u.
}\]
\end{small}

All the commutative squares and the commutative parallelogram in the diagram are Cartesian. 
By \Cref{prop: fg sheaves on kotG}, we can compute the Hom spaces in $\rshv(\kot_G)$ instead of in $\shv(\kot_G)$. Then using various base change for the sheaf theory $\rshv$ (using the fact that $\delta$ is representable coh. pro-smooth), we see that
\begin{eqnarray*}
&    &\Hom_{\rshv(\kot_G)}(\Nt_*(\delta^!\mF_1),(\Nt^u)_*((\delta^u)^!\mF_2))    \\
& = &\Hom_{\rshv(LG/\Ad_\sigma \iw)}\bigl(\delta^!\mF_1, (\Nt^!((\Nt^u)_*((\delta^u)^!\mF_2)))\bigr)\\
& = &\Hom_{\rshv(\iw\backslash LG/\iw)}\bigl(\mF_1, \delta_\flat(\Nt^!((\Nt^u)_*((\delta^u)^!\mF_2)))\bigr)\\
& = &\Hom_{\rshv(\iw\backslash LG/\iw)}\bigl( \mF_1, (m^u)_*((\eta^u)^!((\mathrm{sw}\circ(\sigma\times\id))^!(\eta_\flat(m^!\mF_2)))) \bigr).
\end{eqnarray*}
Here and below for simplicity we write pull-push functors for the sheaf theory $\rshv$ as  $(-)^!$ and $(-)_*$ instead of $(-)^{\ind\fg,!}$ and $(-)_*^{\ind\fg}$.
All the involved functors in the above sequences of isomorphisms are continuous.
By \Cref{lemma-explicit-filtration-unit-affine-flag}, we see that $\Nt^!((\Nt^u)_*((\delta^u)^!\mF_2))$ admits a filtration, which induces a filtration on
\[
(m^u)_*((\eta^u)^!((\mathrm{sw}\circ(\sigma\times\id))^!(\eta_\flat(m^!\mF_2)))) 
\] 
with associated graded being
$\oplus_{w}\nabla^r_{w^{-1}}\star^u\mF_2\star^u \nabla_{\sigma(w)}^l$. 
It follows that the space 
\[
\Hom_{\rshv(\iw\backslash LG/\iw)}\bigl( \mF_1, (m^u)_*((\eta^u)^!((\mathrm{sw}\circ(\sigma\times\id))^!(\eta_\flat(m^!\mF_2)))) \bigr)
\] 
still admits a filtration with associated graded being
\begin{eqnarray*}
&          &\Hom_{\rshv(\iw\backslash LG/\iw)}\Bigl(\mF_1, \nabla^r_{w}\star^u\mF_2\star^u \nabla_{\sigma(w)^{-1}}^l\Bigr) \\
&\cong &\Hom_{\rshv(\iw\backslash LG/\iw)}\Bigl(\mF_1\star \Delta_{\sigma(w)}, \nabla^r_{w}\star^u\mF_2\star^u \nabla_{\sigma(w)^{-1}}^l\star \Delta_{\sigma(w)}\Bigr)\\
&\cong &\Hom_{\rshv(\iw\backslash LG/\iw)}\Bigl(\mF_1\star \Delta_{\sigma(w)}, \nabla^r_{w}\star^u\mF_2\star^u \Delta_e^l\Bigr)\\
&\cong &\Hom_{\rshv(\iw\backslash LG/\iw^u)}\Bigl(\mF_1\star \Delta_{\sigma(w)}^r, \mF_2\star^u\nabla_{\sigma(w)}^l\Bigr).
\end{eqnarray*}
Here for the last isomorphism, we use the fact that $(-)\star^u\Delta_{e}^l$ is the functor of (shifted) $*$-pushfoward along $\iw\backslash LG/\iw^u\to \iw\backslash LG/\iw$, and therefore its left adjoint is just $!$-pullback along the same map.
\end{proof}

Now we genearlize \Cref{lem-spectral-sequence-hom-unip} to equivariant and monodromic settings, allowing non-trivial monodromy. 
First we need a generalization of \Cref{lemma-explicit-filtration-unit-affine-flag-unip}.

\begin{lemma}\label{lemma-explicit-filtration-unit-affine-flag}
Assume that $\La$ is an algebraically closed field and $(n,p')=1$ where $p'$ is the product of $p$ and the characteristic exponent of $\La$ as before.
Let $\mF\in \shv(\iw^u\backslash LG/\iw^u)$. There is a filtration of 
\[
(\eta^{[n]})_\flat((m^{[n]})^!\mF) \in \shv(\iw^u\backslash LG/\iw^{[n]}\times \iw^{[n]}\backslash LG/\iw^u)
\]
with the associated graded being 
\[
\bigoplus_{w,\chi}\nabla_{w,\chi}^l\boxtimes_{\La}(\nabla_{w^{-1},w(\chi)}^r\star^{u} \mF),\quad w\in \widetilde{W}, \ \chi: \mS_k[n]\to\La^\times.
\]
\end{lemma}
\begin{proof}
The same proof as in \Cref{lemma-explicit-filtration-unit-affine-flag-unip} applies, with a small modification. Again,
we first deal with the case when $\mF=\mathbf{1}_{\iw^u}$. 
For each $w$, let 
\[
a_{?}\colon \iw^u\backslash LG_{?}\times^{\iw^{[n]}}LG/\iw^u\to \iw^u\backslash LG\times^{\iw^{[n]}}LG/\iw^u,\quad ?=w \mbox{ or } \leq w
\] 
be the pfp (locally) closed embedding. Let $m_{\leq w}= m^{[n]}\circ a_{\leq w}$ and similarly let $m_w=m^{[n]}\circ a_w$.
Then $(m^{[n]})^!(\mathbf{1}_{\iw^u})$ admits a filtration with associated graded being $(a_w)_*(m_w)^!\mathbf{1}_{\iw^u}$. The generalization of \eqref{eq: stratified unit by Schubert-unip} now reads as
\begin{equation}\label{eq: stratified unit by Schubert}
(\eta^{[n]})_\flat\bigl((a_w)_*((m_w)^!\mathbf{1}_{\iw^u})\bigr)\cong \bigoplus_{\chi}\nabla_{w,\chi}^l\boxtimes_{\La}\nabla_{w^{-1},w(\chi)}^r.
\end{equation}
Again, by change, it is enough to consider the following sequence of morphisms
\begin{multline}\label{mult: LGw diagonal}
LG_w/\iw^{[n]}\stackrel{g\iw^{[n]}\mapsto (g,g^{-1}\iw^u)}{\xrightarrow{\hspace*{2.5cm}}} LG_w\times^{\iw^{[n]}} LG_w/\iw^u \\ \stackrel{\widetilde{\eta^{[n]}_w}}{\longrightarrow} LG_w/\iw^{[n]}\times\iw^{[n]}\backslash LG_w/\iw^u \hookrightarrow LG/\iw^{[n]} \times \iw^{[n]}\backslash LG/\iw^u.
\end{multline}
Now \eqref{eq: stratified unit by Schubert} would follow if we show that
after the first two pushforwards, we obtain
\begin{equation}\label{eq: stratified unit by Schubert-2}
\bigoplus_{\chi} (\pr^l_w)^!\chi[-\ell(w)]\boxtimes_{\La} (\pr^r_w)^!\chi[-\ell(w)] \in \shv(LG_w/\iw^{[n]}\times \times\iw^{[n]}\backslash LG_w/\iw^u).
\end{equation}
Here, we write $\pr^l_w, \pr^r_w$ for the projections
\[
\pr^l_w\colon LG_w/\iw^{[n]}\cong \iw^u/ (\Ad_{\dot{w}} \iw^u\cap \iw^u)\cdot \dot{w} \cdot \bB \mS_k[n]\longrightarrow \bB \mS_k[n],
\]
and 
\[
\pr^r_w\colon \iw^{[n]}\backslash LG/\iw^u\cong \bB \mS_k[n]\cdot \dot{w} \cdot \bB  (\Ad_{\dot{w}^{-1}} \iw^u\cap \iw^u)\longrightarrow \bB \mS_k[n].
\]
Similar as before,  it is enough to compute the $*$-pushforward of the dualizing sheaf of $LG_w/\iw^{[n]}$ along the maps
\[
LG_w/\iw^{[n]}\to LG_w/\iw^{(r)}\times^{\iw^{[n]}/\iw^{(r)}} LG_w/\iw^u \to LG_w/\iw^{[n]}\times (\iw^{[n]}/\iw^{(r)})\backslash LG_w/\iw^u.
\]
Now using the fact that the $*$-pushforward of the dualizing sheaf of $\bB \mS_k[n]$ along the diagonal map $\bB \mS_k[n]\to \bB \mS_k[n]\times \bB \mS_k[n]$ is $\oplus_{\chi }\chi\boxtimes_{\La} \chi$,
we obtain \eqref{eq: stratified unit by Schubert-2} and therefore  \eqref{eq: stratified unit by Schubert}. 

The case for general $\mF$ follows from the same argument as before.
\end{proof}

Now we have the following generalization of  \Cref{lem-spectral-sequence-hom-unip} and \Cref{rem-spectral-sequence-hom}. Given \Cref{lemma-explicit-filtration-unit-affine-flag}, the proof remains the same.

\begin{proposition}\label{lem-spectral-sequence-hom}
Assumptions are as in \Cref{lemma-explicit-filtration-unit-affine-flag}. Let $\mF_1$ be an object in $\fgshv(\iw^{[n]}\backslash LG/\iw^{[n]})$ and $\mF_2$ an object in $\fgshv(\iw^u\backslash LG/\iw^u)$. 
Then there is a filtration on the $\La$-module
\[
\Hom_{\shv(\kot_G)}((\Nt^{[n]})_*(\delta^{[n]})^!\mF_1,(\Nt^u)_*(\delta^u)^!\mF_2)
\] 
with the associated graded being $\Hom_{\shv(\iw^u\backslash LG/\iw)}(\Delta^l_{w,\chi}\star^{[n]}\mF_1, \mF_2\star^u\nabla^l_{\sigma(w),\sigma(\chi)})$ for $w\in \widetilde{W}$ and $\chi: \mS_k[n]\to \La^\times$. In particular, there is a spectral sequence with
\[
E_1^{p,q}\simeq \bigoplus_{\ell(w)=-p, \chi: \mS_k[n]\to\La^\times}\mathrm{Ext}^{q+p}_{\shv(\iw^u\backslash LG/\iw^{[n]})}(\Delta^l_{w,\chi}\star^{[n]}\mF_1,\mF_2\star^{u}\nabla^l_{\sigma(w),\sigma(\chi)}),
\]
and with abutment $\mathrm{Ext}^*_{\shv(\kot_G)}((\Nt^{[n]})_*(\delta^{[n]})^!\mF_1,(\Nt^u)_*(\delta^u)^!\mF_2)$.
\end{proposition}

\begin{remark}\label{rem-spectral-sequence-hom}
Here is a variant of \Cref{lem-spectral-sequence-hom}. Namely, we factor $\eta^{[n]}$ as 
\[
\iw^{[n]}\backslash LG\times^{\iw^{[n]}} LG/\iw^{[n]}\xrightarrow{\eta_1} \iw^{[n]}\backslash LG/\iw^u\times^{\mS_k^{[n]}}\iw^u\backslash LG/\iw^{[n]}\xrightarrow{\eta_2} \iw^{[n]}\backslash LG/\iw^{[n]}\times \iw^{[n]}\backslash LG/\iw^{[n]}.
\]
Then
if we start with $\mF\in \shv(\iw^{[n]}\backslash LG/\iw^{[n]})$, then 
\[
(\eta_1)_\flat(m^{[n]})^!\mF\in \shv(\iw^{[n]}\backslash LG/\iw^{[n]}\times \iw^{[n]}\backslash LG/\iw^{[n]})
\]
admits a filtration with associated graded being
\[
\bigoplus_{w,\chi}\nabla_{w,\chi} \boxtimes_{\bB \mS_k^{[n]}}(\nabla_{w^{-1},w(\chi)}\star^{[n]} \mF).
\]
Here we recall that the functor
\[
- \boxtimes_{\bB \mS_k^{[n]}} -: \shv(\iw^{[n]}\backslash LG/\iw^{[n]})\otimes_\La \shv(\iw^{[n]}\backslash LG/\iw^{[n]})\to \shv(\iw^{[n]}\backslash LG/\iw^{[n]}\times_{\bB \mS_k^{[n]}} \iw^{[n]}\backslash LG/\iw^{[n]})
\]
is as from \Cref{rem:relative-sheaf-theory}.
Consequently, for $\mF_1\in \fgshv(\iw^{[n]}\backslash LG/\iw^{[n]})$ and $\mF_2\in \shv(\iw^{[n]}\backslash LG/\iw^{[n]})$, the $\La$-module
\[
\Hom_{\shv(\kot_G,\La)}((\Nt^{[n]})_*(\delta^{[n]})^!\mF_1,(\Nt^{[n]})_*(\delta^{[n]})^!\mF_2)
\] 
admits a filtration with associated graded being
\begin{equation}\label{eq:associated graded equiv. version}
\Hom_{\shv(\iw\backslash LG/\iw)}(\mF_1, \nabla_{w^{-1},w(\chi)}\star^{[n]}\mF_2\star^{[n]}\nabla_{\sigma(w),\sigma(\chi)}\otimes_{C^\bullet(\bB \mS_k^{[n]})\otimes C^\bullet(\bB \mS_k^{[n]})} C^\bullet(\bB \mS_k^{[n]})).
\end{equation}
\end{remark}

We will also need a monodromic version of \Cref{lem-spectral-sequence-hom}. For that purpose, we first need an analogue of \Cref{lemma-explicit-filtration-unit-affine-flag}.
We consider the following commutative diagram with Cartesian square
\begin{small}
\begin{equation}\label{eq: diagram for computing unit of monodromic cat}
\xymatrix{
                                         & \ar_-{m^u}[dl]\widetilde{\iw^u}\backslash LG\times^{\iw^u} LG/\widetilde{\iw^u}\ar^-{\eta^u}[r]\ar^{\av_s}[d]& \widetilde{\iw^u}\backslash LG/\iw^u\times\iw^u\backslash LG/\widetilde{\iw^u}\ar^{\av_s}[d]   \\
\widetilde{\iw^u}\backslash LG/\widetilde{\iw^u}&  \ar_-{m}[l]  \widetilde{\iw^u}\backslash LG\times^{\iw} LG/\widetilde{\iw^u}\ar^-{\eta_1}[r]     & \widetilde{\iw^u}\backslash LG/\iw^u\times^{\mS_k}\iw^u\backslash LG/\widetilde{\iw^u}.
}
\end{equation} 
\end{small}

The group $\mS_k$ acts on $\widetilde{\iw^u}\backslash LG/\iw^u\times^{\mS_k}\iw^u\backslash LG/\widetilde{\iw^u}$ through the middle and we can form the corresponding monodromic category $\shv_{\mon}(\widetilde{\iw^u}\backslash LG/\iw^u\times^{\mS_k}\iw^u\backslash LG/\widetilde{\iw^u})$.

\begin{lemma}\label{lem: unit for monodromic category}
The sheaf $\av^{\mon}((\eta_1)_\flat(m^!(\mathbf{1}_{\widetilde{\iw^u}})))\in \shv_{\mon}(\widetilde{\iw^u}\backslash LG/\iw^u\times^{\mS_k}\iw^u\backslash LG/\widetilde{\iw^u})$ admits a filtration with associated graded being
\[
\bigoplus_{w,\psi} (\av_s)_* \bigl({}^\psi\widetilde{\nabla}_{\dot{w}}^{\mon}\boxtimes_{\La}  \widetilde{\nabla}_{\dot{w}^{-1}}^{\mon,\psi}\bigr).
\]
\end{lemma}

\begin{proof}
We still consider the (locally) closed embedding $a_w$ and $a_{\leq w}$ from \Cref{lemma-explicit-filtration-unit-affine-flag}. Then we need to show that
\[
\av^{\mon}((\eta_1)_\flat(((a_w)_*((m_w)^!(\mathbf{1}_{\widetilde{\iw^u}})))))\cong \bigoplus_{\psi}(\av_s)_* \bigl({}^\psi\widetilde{\nabla}_{\dot{w}}^{\mon}\boxtimes_{\La}  \widetilde{\nabla}_{\dot{w}^{-1}}^{\mon,\psi}\bigr).
\]
We follow the same idea of the proof of \eqref{eq: stratified unit by Schubert}. We perform the base change of the Cartesian square in \eqref{eq: diagram for computing unit of monodromic cat} along the map
$LG/\iw^u\times^{\mS_k} \iw^u\backslash LG/\widetilde{\iw^u}\to \widetilde{\iw^u}\backslash LG/\iw^u\times^{\mS_k} \iw^u\backslash LG/\widetilde{\iw^u}$, and consider the following sequence of morphisms
\begin{multline*}
LG_w/\iw\stackrel{g\iw\mapsto (g,g^{-1}\widetilde{\iw^u})}{\xrightarrow{\hspace*{2.5cm}}}    LG_w\times^{\iw}  LG_w/\widetilde{\iw^u} \\ \stackrel{\widetilde{\eta_1}}{\longrightarrow} LG_w/\iw^{u}\times^{\mS_k}\iw^{u}\backslash LG_w/\widetilde{\iw^u} \hookrightarrow LG/\iw^{u} \times^{\mS_k} \iw^{u}\backslash LG/\widetilde{\iw^u}.
\end{multline*}
Then as in  \Cref{lemma-explicit-filtration-unit-affine-flag}, the base change of $\av^{\mon}((\eta_1)_\flat(m^!(\mathbf{1}_{\widetilde{\iw^u}})))$ is obtained from the dualizing sheaf on $LG_w/\iw$ by $*$-pushforward along the first map, followed by $\flat$-pushforward along the second map, and then followed by $*$-pushforward along the last map, and finally followed by the functor $\av^{\mon}$. In addition, using \Cref{lem: functors between monodromic categories}, we see that $\av^{\mon}$ commutes with the last pushforward.

Now, following the proof of \Cref{lemma-explicit-filtration-unit-affine-flag}, we choose $\iw^{(r)}$ sufficiently small congruence subgroup and compute the $*$-pushforward along
\[
LG_w/\iw\to LG_w/\iw^u\times^{\mS_k}  (\iw^u/\iw^{(r)})\backslash LG_w/\widetilde{\iw^u}
\]
followed by $\av^{\mon}$. Using that $\av^{\mon}(\delta_1)=\widetilde\Ch=\widetilde\Ch\star\widetilde\Ch$ (see \Cref{lem: monoidal of monodromic}), the lemma follows.
\end{proof}

\begin{proposition}\label{lem-mon-spectral-sequence-hom}
Let $\mF_1$ be an object in $\shv_{\mon}(\iw^u\backslash LG/\iw^u)$ such that $(\av_s)_*(\delta^u)^!\mF_1$ is compact in $\shv(\sht^{\loc})$,  and let $\mF_2$ be an object in $\shv(\widetilde{\iw^u}\backslash LG/\widetilde{\iw^u})$. 
Then there is a filtration on the $\La$-module 
\[
\Hom_{\shv(\kot_G)}(\Ch_{G,\phi}^{\mon}(\mF_1),(\widetilde{\Nt^u})_*(\widetilde{\delta^u})^!\mF_2)
\] 
indexed by $w\in \widetilde{W}$, with the associated graded being 
\[
\bigoplus_\psi\Hom_{\shv(\iw^u\backslash LG/\widetilde{\iw^u})}((\av_s)^*(\av_s)_*\mF_1\star^u \widetilde{\Delta}_{\sigma(w)}^{\mon,\psi}, \widetilde{\nabla}_{\dot{w}}^{\mon,\psi}\star^{\tilde{u}}\mF_2).
\] 
\end{proposition}
\begin{proof}
We will need to consider a variant of the big commutative diagram in the proof of \Cref{lem-spectral-sequence-hom}. 
\begin{small}
\[
\xymatrix{ 
\iw^{u}\backslash LG/ \iw^{u}\ar_-{\av_s}[d]&\ar_-{\widetilde{m^u}}[l]\iw^{u}\backslash LG\times^{\widetilde{\iw^u}} LG/\iw^u\ar^-{\widetilde{\eta^u}}[r]\ar^-{\av_s}[d]&\iw^u\backslash LG/\widetilde{\iw^u}\times \widetilde{\iw^u}\backslash LG/\iw^u\ar^{\mathrm{sw}\circ(\sigma\times\id)}[r]&\widetilde{\iw^u}\backslash LG/\iw^{u}\times \iw^{u}\backslash LG/\widetilde{\iw^u}\ar^-{\av_s}[d]\\
 \frac{\iw^u\backslash LG/\iw^u}{\Ad_\sigma \mS_k}&\ar_{\widetilde{m^u}}[l]\frac{\iw^u\backslash LG\times^{\widetilde{\iw^u}}LG/\iw^u}{\Ad_\sigma \mS_k}\ar_-\cong^{\pFr}[r]&\frac{LG/\iw^u\times^{\mS_k}\iw^u\backslash LG}{\Ad_\sigma\widetilde{\iw^u}}\ar^-{\widetilde{\delta_u}}[r] &\widetilde{\iw^u}\backslash LG/\iw^u\times^{\mS_k} \iw^u\backslash LG/\widetilde{\iw^u} \\
  \frac{LG}{\Ad_\sigma \iw}\ar^{\Nt}[rd]\ar^-{\av_u=\eta_1}[u] &\ar_{\eta_1}[u]\ar_{\widetilde{m^{u}}}[l]\frac{LG\times^{\widetilde{\iw^u}}LG}{\Ad_\sigma \iw}\ar_-\cong^{\pFr}[r]\ar[dr]&\frac{LG\times^{\iw}LG}{\Ad_\sigma\widetilde{\iw^u}}\ar^{m}[d]\ar^-{\widetilde{\delta^u}}[r]\ar_{\eta_1}[u] &\widetilde{\iw^u}\backslash LG\times ^{\iw} LG/\widetilde{\iw^u} \ar^{m}[d]\ar_{\eta_1}[u]\\
 &\kot_G  &\ar_{\widetilde{\Nt^u}}[l] \frac{LG}{\Ad_\sigma \widetilde{\iw^u}}\ar^-{\widetilde{\delta^u}}[r] &\widetilde{\iw^u}\backslash LG/\widetilde{\iw^u}.
}\]
\end{small}
Using the Cartesian diagram from \eqref{eq: adiwu-to-adiw} and \Cref{lem: functors between monodromic categories} \eqref{lem: functors between monodromic categories-4},
we see that 
\[
(\av_s)_*(\delta^u)^!\mF_1\cong (\av_u)^! (\av_s)_*\mF_1\cong (\av_u)^! (\av_s)_!\mF_1[\dim S].
\]
Therefore,
\begin{eqnarray*}
&  &\Hom_{\shv(\kot_G)}(\Ch_{G,\phi}^{\mon}(\mF_1), (\widetilde{\Nt^u})_*(\widetilde{\eta^u})^!\mF_2)\\
&=&\Hom_{\shv( \frac{LG}{\Ad_\sigma \iw})}((\eta_1)^!((\av_s)_*\mF_1),\Nt^!((\widetilde{\Nt^u})_*(\widetilde{\eta^u})^!\mF_2 ))\\
&=&\Hom_{\shv( \frac{\iw^u\backslash LG/\iw^u}{\Ad_\sigma \mS_k})}((\av_s)_*\mF_1,  (\widetilde{m^u})_*((\widetilde{\delta_u}\circ \pFr)^!((\eta_1)_\flat(m^!\mF_2))))\\
&=&\Hom_{\shv_{\mon}(\iw^u\backslash LG/\iw^u)}(\mF_1,\av^{\mon}((\av_s)^!((\widetilde{m^u})_*((\widetilde{\delta_u}\circ \pFr)^!((\eta_1)_\flat(m^!\mF_2))))))[-\dim S]\\
&=&\Hom_{\shv_{\mon}(\iw^u\backslash LG/\iw^u)}(\mF_1, (\av_s)^!((\widetilde{m^u})_*((\widetilde{\delta_u}\circ \pFr)^!(\av^{\mon}((\eta_1)_\flat(m^!\mF_2)))))).
\end{eqnarray*}

Then using \Cref{lem: unit for monodromic category} and argued as in \Cref{lem-spectral-sequence-hom}, there is a filtration of with associated graded
\[
\bigoplus_\psi\Hom_{\shv_{\mon}(\iw^u\backslash LG/\iw^u)}( (\av_s)^*(\av_s)_*\mF_1, \widetilde{\nabla}_{\dot{w}}^{\mon,\psi}\star^{\tilde{u}} \mF_2\star^{\tilde{u}} {}^\psi\widetilde{\nabla}^{\mon}_{\sigma(\dot{w})^{-1}}),
\]
which by \Cref{lem: right adjoint of Whittaker average} is isomorphic to
\[
\bigoplus_\psi\Hom_{\shv_{\mon}(\iw^u\backslash LG/\widetilde{\iw^u})}((\av_s)^*(\av_s)_*\mF_1\star^u \widetilde{\Delta}_{\sigma(\dot{w})}^{\mon,\psi}, \widetilde{\nabla}_{\dot{w}}^{\mon,\psi}\star^{\tilde{u}}  \mF_2).
\]
The proposition is proved.
\end{proof}

\begin{remark}
Of course, in the above proposition, the case $\widetilde{\iw^u}=\iw^u$ is allowed. 
\end{remark}

\begin{remark}
We suppose $\La=\overline\bQ_\ell$. Let $\mF_1=\nabla_{\dot{w},\hat{u}}^{\mon}$. One checks that \Cref{lem-mon-spectral-sequence-hom} reduces to \Cref{lem-spectral-sequence-hom-unip}.
\end{remark}

As an application of the above discussions, we can now give a proof of \Cref{prop: aff DL disjointness}.

\begin{proof}[Proof of \Cref{prop: aff DL disjointness}]
Let $(w_1,\theta_1)$ and $(w_2,\theta_2)$ be as in the proposition. For $\theta_i$, we let $\chi_i: T^p\mS_k\to \mS_k^{\bar{w}_i\sigma}\to\La^\times$ be the associated character, or equivalent the homomorphism $\chi_i: I_F^t\to \hat{S}(\La)$, as in  the proof of \Cref{prop: DL geo conj vs tame inertia type}. That $(w_1,\theta_1)$ and $(w_2,\theta_2)$ are not geometrically conjugate means that $\chi_1$ and $\chi_2$ are in different $W_0$-orbits.
Recall that on the dual side we have the formal scheme $\hchi_i\subset R_{I_F,\hat{S}}$.

We first show that 
\[
\Hom(R^{\mon,*}_{w_1,\hchi_1}, R^{\mon,*}_{w_2,\hchi_2})=0.
\]
By \Cref{lem-mon-spectral-sequence-hom} and \Cref{cor: av upper star av lower star}, this $\La$-module admits a filtration with associated graded being
\begin{eqnarray*}
&&\Hom_{\shv(\iw^u\backslash LG/\iw^u)}(\nabla^{\mon}_{w_1, \hchi_1\cap \chi_{\varphi_{\bar w_1}}}\star^u \widetilde{\Delta}_{\sigma(w)}^{\mon}, \widetilde{\nabla}_{w}^{\mon}\star^{u}\nabla^{\mon}_{w_2,\hchi_2})\\
&\cong&\Hom_{\shv(\iw^u\backslash LG/\iw^u)}(\nabla^{\mon}_{w_1, \hchi_1\cap \chi_{\varphi_{\bar w_1}}}\star^u \Delta_{\sigma(w), w_1^{-1}(\hchi_1)}^{\mon}, \nabla_{w,ww_2\hchi_2}^{\mon}\star^{u}\nabla^{\mon}_{w_2,\hchi_2})
\end{eqnarray*} 
Here the isomorphism follows from \eqref{eq: non-zero monodromic category} and  \eqref{eq: non-zero tensor product monodromic}. In addition, the above space is non-zero only if there is some $w$ such that 
$\hchi_1\cap ww_2\hchi_2\neq \emptyset$ and 
$\sigma(w)^{-1}w_1^{-1}\hchi_1\cap w_2^{-1}\hchi_2\neq \emptyset$. But this is impossible as $\chi_1$ and $\chi_2$ are in different $W_0$-orbits.

Now note that by \Cref{lem: formula of monodromic ADL induction}, we have $R^{\mon,*}_{w_i,\hchi_i}=R_{\dot{w}}^!(\omega_{\hchi_i\cap \chi_{\varphi_{\bar w_i}}})$. As $\theta_i$ is contained in the full subcategory of $\mS_k^{\bar w_i\sigma}$-modules generated by $\omega_{\hchi_i\cap \chi_{\varphi_{\bar w_i}}}$, the above vanishing result also implies $\Hom_{\shv(\kot_G)}(R^*_{\dot{w}_1,\theta_1},R^*_{\dot{w}_2,\theta_2})=0$, as desired.

Alternatively, the above vanishing result can also be proved using \Cref{rem-spectral-sequence-hom}.
\end{proof}
\quash{As explained before, $\chi_i$ descends to a character $\mS_k[n_i]\to\La^\times$ for some $n_i$ coprime to $p'$. Let $n=n_1n_2$. Now
we let $\mF_i=\nabla_{w_i,\chi_i}$. By applying \Cref{prop:affine DL induction chi}, it is enough to show that unless $(w_1,\chi_1)$ and $(w_2,\chi_2)$ are geometrically conjugate, 
\eqref{eq:associated graded equiv. version} vanishes for for every $w,\chi$.
\[
\Hom_{\shv(\iw^u\backslash LG/\iw)}(\Delta^l_{w,\chi}\star^{[n]}\nabla_{w_1,\chi_1},\nabla_{w_2,\chi_2}\star^{u}\nabla^l_{\sigma(w),\chi})=0.
\]
Note that by \eqref{eq: chi and chi' incompatible}, $\nabla_{w^{-1},w(\chi)}\star^{[n]}\nabla_{w_2,\chi_2}\star^{[n]}\nabla_{\sigma(w),\sigma(\chi)}\neq 0$ only if $\chi_2=\sigma( w(\chi))$ (then $w(\chi)=w_2(\chi_2)$). Then the hom space \eqref{eq:associated graded equiv. version} is non-zero only if $\chi_1=\sigma(\chi)$ (then $w_1(\chi_1)=\chi$).
$\nabla^l_{w,\chi}\star^{[n]}\nabla_{w_1,\chi_1}\neq 0$ only if $\chi=w_1(\chi_1)$ and $\nabla_{w_2,\chi_2}\star^{u}\nabla^l_{\sigma(w),\chi}\neq 0$ only if $\chi_2=w(\chi)$. 
In other others, the above space is non-zero only if $\chi_1$ and $\chi_2$ are conjugated by some element in $W_0$. This means that $\theta_1$ and $\theta_2$ are geometrically conjugate.}

\subsection{Representations of finite group of Lie type: old and new}
We first apply the machinery developed so far to representation theory of finite groups of Lie type, leading to some new results, or alternative (sometimes simpler) proofs of some classical results.

Unlike everywhere else in the article, in this subsection, we let $G$ denote a connected reductive group over a finite field $\kappa$ with a Borel $B$ and its unipotent radical $U$. Let $T=B/U$ be the (abstract) Cartan of $G$. Let $W$ denote the absolute Weyl group of $G$.
We also let $e: U\to \bG_a$ be a homomorphism such that $(G,B,T,e)$ form a pinning of $G$.

As before, we base change everything to a fixed algebraic closure $k$ of $\kappa$. Let $\sigma$ denote the Frobenius endomorphism. 
Our sheaf theory will have coefficient ring $\La$, but otherwise specified $\La$ will be omitted from the notation.

\subsubsection{Deligne-Lusztig inductions}
We have
\[
U\backslash G/U \xleftarrow{\delta^u} \frac{G}{\Ad_\sigma U}\xrightarrow{\Nt^u} \frac{G}{\Ad_\sigma G}\cong \bB G(\kappa), 
\]
and then the Deligne-Lusztig induction functor
\begin{equation}\label{eq: finite mon DL induction}
\Ch_{G,\phi}^{\mon}:=(\Nt^u)_*(\delta^u)^!: \shv_{\mon}(U\backslash G/U)\to \shv(\bB G(\kappa))\cong \rep(G(\kappa)).
\end{equation}

We also have the unipotent version  
\begin{equation}\label{eq: finite unip DL induction}
\Ch_{G,\phi}^{\unip}:=\Nt_*\delta^!: \ind\cshv(B\backslash G/B)\to \ind\cshv(G(\kappa))\cong \ind\crep(G(\kappa)),
\end{equation}
given by the correspondence
\[
B\backslash G/B\xleftarrow{\delta}\frac{G}{\Ad_\sigma B}\xrightarrow{\Nt} \bB G(\kappa).
\]
Here we recall $\crep(G(\kappa))$ is the full subcategory of $\rep(G(\kappa))$ consisting of those objects whose underlying $\La$-module is perfect, which under the equivalence $\shv(\bB G(\kappa))\cong \rep(G(\kappa))$ corresponds to the category of constructible sheaves on $\bB G(\kappa)$. 
We also recall that $ \rep(G(\kappa))^\cpt\subset \crep(G(\kappa))$ but the inclusion might be strict in general.
The functor \eqref{eq: finite unip DL induction} restricts to a functor from $\shv(B\backslash G/B)\to \shv(\bB G(\kappa))=\rep(G(\kappa))$.

For $w\in W$ with a lifting $\dot{w}\in G(k)$, we have $R_w^?, R_{\dot{w},\theta}^?, \widetilde{R}_{\dot{w}}^?, R^{\mon,?}_{\dot{w},\hchi}=R^?_{\dot{w},\hchi\cap \varphi_{\chi_w}}$, for $?=*,!$. We also have $\widetilde{R}^T_{\dot{w}}$ (resp,  $R^{\mon,T}_{\dot{w},\hchi}$), which is the Deligne-Lusztig induction of tilting sheaf $\widetilde{\Til}^\mon_{\dot{w}}$ (resp. $\Til^\mon_{\dot{w},\hchi}$). See \eqref{eq:ADLI of tilting sheaf}. Recall that by (the same argument as in) \Cref{lem: compactness of affine DL sheaf}, $\widetilde{R}_{\dot{w}}^*$ and $\widetilde{R}_{\dot{w}}^!$ are compact objects in $\rep(G(\kappa))$. On the other hand, $R_{\dot{w},\theta}^*$ and $R_{\dot{w},\theta}^!$ (in particular $R_w^*$ and $R_w^!$) are in $\crep(G(\kappa))$. But they may not be compact in general.

The following statement is probably well-known. But in the generality we have not found it in literature.

\begin{lemma}\label{lem: coconnectivity of shrek DL induction}
With respect to the standard $t$-structure on $\shv(\bB G(\kappa))=\rep(G(\kappa))$, we have $\widetilde{R}_{\dot{w}}^!\in \rep(G(\kappa))^{\geq 0}$. Dually, we have $\widetilde{R}_{\dot{w}}^*\in \rep(G(\kappa))^{\leq 0}$.
\end{lemma}
\begin{proof}
Note that
$\widetilde{R}_{\dot{w}}^{!}$ is nothing but the shifted compactly supported cohomology of the Deligne-Lusztig variety $\widetilde X_{\dot{w}}=\{gU\mid g^{-1}\sigma(g)\in U\dot{w}U\}\subset G/U$.
When the cardinality of $\kappa$ is bigger than the Coxeter number $h$ of $G$, it is known (see \cite[Theorem 9.7]{Deligne.Lusztig}) that $\widetilde X_{\dot{w}}$ is affine, and therefore $\widetilde{R}_{\dot{w}}^!\in \rep(G(\kappa))^{\geq 0}$. Affineness of $\widetilde X_{\dot{w}}$ is not known in general, but is known when $w$ is of minimal length in its $\sigma$-conjugacy class (\cite{BR.affineness, He.affineness}). Therefore $\widetilde{R}_{\dot{w}}^!\in \rep(G(\kappa))^{\geq 0}$ for such $w$. Now we apply
the Deligne-Lusztig reduction method (see \cite[Theorem 1.6]{Deligne.Lusztig} and \Cref{lem: ADL reduction method}) and \Cref{thm: reduction to min length elements} \eqref{thm: reduction to min length elements-1}, we can prove the desired estimate by induction on $\ell(w)$.
\end{proof}

We recall the following transitivity of Deligne-Lusztig induction.

\begin{lemma}\label{lem: transitivity DL induction}
Let $P\subset B$ be a standard rational parabolic subgroup of $G$ with $L$ its Levi quotient. Let $B_L$ be the image of $B$ in $L$, which is a rational Borel of $L$.  Let $W_P\subset W$ be the Weyl group of $P$, which is $\sigma$-stable. 
We identify $\shv_{\mon}(U_L\backslash L/U_L)\cong \shv_{\mon}(U\backslash P/U)\subset \shv_{\mon}(U\backslash G/U)$ as a full subcategory. Then for every $\mF\in \shv_{\mon}(U_L\backslash L/U_L)$, we have 
\[
\ind_{P(\kappa)}^{G(\kappa)}\Ch^\mon_{L,\phi}(\mF)\cong \Ch^{\mon}_{G,\phi}(\mF).
\]
A similar statement holds for the unipotent version $\Ch^{\unip}_{G,\phi}$.
\end{lemma}
\begin{proof}
The lemma follows from base change isomorphisms, together with the following commutative diagram with the right square Cartesian
\[
\xymatrix{
U\backslash P/U\ar[d] &\ar[l] \frac{P}{\Ad_\sigma U} \ar[r]\ar[d] &  \frac{P}{\Ad_\sigma P}\ar[r]\ar[d] & \frac{G}{\Ad_\sigma G}\\
U_L\backslash L/U_L & \ar[l]  \frac{L}{\Ad_\sigma U_L} \ar[r] & \frac{L}{\Ad_\sigma L} &
}\]
\end{proof}

Recall the following completeness result of Deligne-Lusztig induction, due to Bonnaf\'e-Rouquier \cite[\textsection{9}, Theorem A]{BR.modular}. The case $\La=\overline\bQ_\ell$ was previously due to Deligne-Lusztig \cite[Corollary 7.7]{Deligne.Lusztig}.
We will also sketch a geometric proof of this fact in \Cref{rem: completeness of DL via Springer sheaf} when the order of the Weyl group $\sharp W$ is invertible in $\La$.

\begin{theorem}\label{thm: DL completeness}
The category $\rep(G(\kappa))$ is generated (as $\La$-linear presentable stable categories) by $\{\widetilde{R}^*_{\dot{w}}\}_{w\in W}$, as well as by $\{\widetilde{R}^!_{\dot{w}}\}_{w\in W}$. 
\end{theorem}

For later purposes, we introduce the following categories, which can be regarded as versions of categories of unipotent representations of $G(\kappa)$.
\begin{definition}\label{def: hunip and unip for finite group}
\begin{enumerate}
\item We let $\rep^{\widehat\unip}(G(\kappa))$ denote the full subcategory of $\rep(G(\kappa))$ generated (as a presentable $\La$-linear category) by $\{R^{\mon,*}_{\dot{w},\hat{u}}\}_{w\in W}$, or equivalently by  $\{R^{\mon,!}_{\dot{w},\hat{u}}\}_{w\in W}$.  
\item We let $\rep^{\unip}_c(G(\kappa))$ be the full idempotent complete stable subcategory of $\crep(G(\kappa))$ generated by $\{R_w^*\}_{w\in W}$, or equivalently by $\{R_w^!\}_{w\in W}$.
\end{enumerate}
\end{definition}

\begin{remark}\label{rem: char zero two versions of unip representation are the same}
Note that when $\La=\overline\bQ_\ell$, then $\rep^{\unip}_c(G(\kappa),\overline\bQ_\ell)=\rep^{\widehat\unip}(G(\kappa),\overline\bQ_\ell)^\cpt$, since $R^{\mon,?}_{\dot{w},\hat{u}}=R_w^?$ for $?=*,!$ (see \eqref{eq: chimon-equal-to-theta}).
Objects in (the heart of the standard $t$-structure of) them are unipotent representations of $G(\kappa)$ in the sense of Deligne-Lusztig. In general, these categories are different. But we always have
\[
\rep^{\widehat\unip}(G(\kappa))^\cpt\subset \crep^{\unip}(G(\kappa))\subset \rep^{\widehat\unip}(G(\kappa)).
\]
In addition, as $R^{\mon,?}_{\dot{w},\hat{u}}\in \crep^{\unip}(G(\kappa))$, we see that $\rep^{\widehat\unip}(G(\kappa))$ is  generated as a presentable $\La$-linear category by $\{R_w^*\}_{w\in W}$, or equivalently by $\{R_w^!\}_{w\in W}$. 

On the other hand, we do not know whether the inclusion 
\[
\crep^{\unip}(G(\kappa))\subset \rep^{\widehat\unip}(G(\kappa))\cap \crep(G(\kappa))
\] 
is an equivalence. For example, the trivial representation of $G(\kappa)$ belongs to $\rep^{\widehat\unip}(G(\kappa))\cap \crep(G(\kappa))$. But we do not know whether it belongs to $\crep^{\unip}(G(\kappa))$ in general.
\end{remark}

Despite the last comment in the above remark, we can show that some induced representations from parabolic subgroups do belong to $\crep^{\unip}(G(\kappa))$, under some restriction of the coefficient ring $\La$.
\begin{lemma}\label{lem: induction of parabolic belonging unipotent}
\begin{enumerate}
\item\label{lem: induction of parabolic belonging unipotent-1} Suppose the image of the $*$-pushforward $\cshv(\bB B,\La)\to \cshv(\bB G,\La)$ generates $\cshv(\bB G,\La)$ (as idempotent complete stable category). Then the trivial representation $\La$ belongs to $\crep^{\unip}(G(\kappa),\La)$.
\item\label{lem: induction of parabolic belonging unipotent-2} If $\sharp W$ is invertible in $\La$, then  the image of the $*$-pushforward $\cshv(\bB B,\La)\to \cshv(\bB G,\La)$ generates the $\cshv(\bB G,\La)$.
\end{enumerate}
\end{lemma}
We note that our assumption on $\La$ in Part \eqref{lem: induction of parabolic belonging unipotent-2} of the lemma is by no means the optimal one.
\begin{proof}
We consider the following Cartesian diagram
\[
\xymatrix{
G(\kappa)\bs G/B\ar[r]\ar_{f'}[d] & \bB B \ar^f[d] \\
\bB G(\kappa) \ar[r] & \bB G.
}\]
It follows that if $\consdual_{\bB G}$ belongs to the idempotent complete subcategory of $\cshv(\bB G,\La)$ generated by the image of the $f_*\consdual_{\bB B}$, then the trivial representation of $G(\kappa)$ is contained in the idempotent complete subcategory generated by  $f'_*\consdual_{G(\kappa)\bs G/B}=\Ch_{G,\phi}^{\unip}(\consdual_{B\bs G/B})$. Part \eqref{lem: induction of parabolic belonging unipotent-1} follows.

For Part \eqref{lem: induction of parabolic belonging unipotent-2}, we just need to observe that under our assumption, $\consdual_{\bB G}=f_*\consdual_{\bB B}\otimes_W\La$ is a direct summand of $f_*\consdual_{\bB B}$.
\end{proof}

\begin{corollary}\label{lem: induction of parabolic belonging unipotent-2}
Suppose $\sharp W$ is invertible in $\La$. Then for every rational parabolic subgroup $P\subset G$, $\ind_{P(\kappa)}^{G(\kappa)}\La\in \crep^{\unip}(G(\kappa))$.
\end{corollary}
\begin{proof}
This is a combination of \Cref{lem: transitivity DL induction} and \Cref{lem: induction of parabolic belonging unipotent}, by noticing that $\sharp W_P\mid \sharp W$.
\end{proof}

When $\La$ is an algebraically closed field, we also have the following disjointness result about Deligne-Lusztig representations \cite[Theorem 6.2]{Deligne.Lusztig} and \cite[Theorem 8.1]{BR.modular}. The argument of  \Cref{prop: aff DL disjointness} (given at the end of \Cref{SS: affine DL sheaf}) applies without change, giving a new proof of this result.

\begin{proposition}\label{prop: DL disjointness}
Let $w_1,w_2\in W_H$ and let $\theta_i: T_H^{w_i\sigma}\to \La^\times$ be two characters. Let $\chi_i: T^pT\to \La^\times$ be associated characters of the prime-to-$p$ Tate module of $T$.
If $(w_1,\theta_1)$ and $(w_2,\theta_2)$ are not geometrically conjugate, then
\[
\Hom_{\rep(G(\kappa))}(R^*_{\dot{w}_1, \hchi_1\cap \chi_{\varphi_{w_1}}}, R^*_{\dot{w}_2, \hchi_2\cap \chi_{\varphi_{w_2}}})=0, \quad \Hom_{\rep(G(\kappa))}(R^*_{\dot{w}_1,\theta_1}, R^*_{\dot{w}_2,\theta_2})=0.
\] 
\end{proposition}

Therefore, we have a decomposition of the category when $\La$ is an algebraically closed field
\begin{equation}\label{eq: Lusztig series finite group Lie type}
\rep(G(\kappa))=\bigoplus_{s} \rep^{\fraks}(G(\kappa)), 
\end{equation}
where $\fraks$ ranges over all geometric conjugacy classes of $(w,\theta)$, and $\rep^{\fraks}(G(\kappa))$ is the full subcategory of $\rep(G(\kappa))$ generated (as presentable $\La$-linear category)
by $R_{\dot{w},\theta}^*$ for $(w,\theta)$ belonging to $\fraks$. Equivalently,  $\rep^{\fraks}(G(\kappa))$ is generated by $R_{\dot{w},\theta}^!$ for $(w,\theta)$ belonging to $\fraks$.
In particular, the block $\rep^{\fraks}(G(\kappa))$ for $\fraks$ containing trivial $\theta$ coincides with $\rep^{\widehat\unip}(G(\kappa))$.

Note that $\rep^{\fraks}(G(\kappa))^\cpt$ is generated (as idempotent complete $\La$-linear category) by
$R^{\mon,*}_{\dot{w},\hchi}$ (or equivalently by $R^{\mon,!}_{\dot{w},\hchi}$), for $\chi$ corresponding to some  $(w,\theta)$ belonging to $\fraks$.

\subsubsection{Projective modules}
The following results seem to be new in the representation theory of finite groups of Lie type. It was also discovered by Arnaud Eteve (see \cite{eteve.tilting}) independently\footnote{Note that the definition of monodromic sheaves in \cite{eteve.monodromic, eteve.tilting, eteve.DL} is different from ours. But probably the resulting monodromic Hecke categories are equivalent.}.

\begin{theorem}\label{prop: exactness of DL induction on tilting} 
We have $\widetilde{R}^T_{\dot{w}}\in \rep(G(\kappa))^\heartsuit$.
I.e. the Deligne-Lusztig induction of the monodromic tilting sheaf is a honest representation of $G(\kappa)$ (rather than a complex). 
In addition, as an object in $\rep(G(\kappa))^{\heartsuit}$, it is projective.
\end{theorem}
\begin{proof}
By \Cref{lem: coconnectivity of shrek DL induction},  $\widetilde{R}_{\dot{w}}^{!}\in  \rep(G(\kappa))^{\geq 0}$. 
In addition as mentioned above, $\widetilde{R}_{w}^{!}$ is compact as an object in $\rep(G(\kappa))$. 
Similarly, $\widetilde{R}_{w}^{*}\in  \rep(G(\kappa))^{\leq 0}$ and is compact as an object in $\rep(G(\kappa))$.

Now, as $\widetilde{R}^T_{\dot{w}}$ admits a finite filtration with associated graded by $\widetilde{R}_{\dot{w}'}^{*}$ as well as a finite filtration with associated graded by $\widetilde{R}_{\dot{w}'}^{!}$, we know that $\widetilde{R}^T_{\dot{w}}\in \rep(G(\kappa))^{\heartsuit}$, and is compact as an object in $\rep(G(\kappa))$.

To prove the second claim, we need to show that $\Hom_{G(\kappa)}(\widetilde{R}^T_{\dot{w}}, \pi)$ concentrates in degree zero for every $\pi\in \rep(G(\kappa))^{\heartsuit}$.

The theory of duality as developed in \Cref{SS: coh. duality of locally profinite gorup} certainly applies to finite groups.
In this case, the Frobenius structure on $\rep(G(\kappa))$ is given by taking the coinvariants. 
Let $\verd^{\can}_{G(\kappa)}$ be the canonical duality induced by such Frobenius structure, and let $(\verd^{\can}_{G(\kappa)})^\cpt: (\rep(G(\kappa))^\cpt)^{\op}\to \rep(G(\kappa))^\cpt$ be the induced anti-involution on compact objects.
Then as $\widetilde{R}^T_{\dot{w}}$ is compact, we have
\[
\Hom(\widetilde{R}^T_{\dot{w}},\pi)=((\verd^{\can}_{G(\kappa)})^\cpt(\widetilde{R}^T_{\dot{w}})\otimes_\La \pi)_{G(\kappa)}.
\]
By \Cref{prop: verdual Ch vs Ch sw dual} and \Cref{prop: rigid dual of cofree monodromic tilting}, applied to the finite setting, we see that
\[
(\verd^{\can}_{G(\kappa)})^\cpt(\widetilde{R}^T_{\dot{w}})\cong \widetilde{R}^T_{\dot{w}}.
\]
Therefore, $H^i\Hom(\widetilde{R}^T_{\dot{w}},\pi)=0$ for $i>0$ as the right hand side concentrates in $\Mod_\La^{\leq 0}$.
\end{proof}

\begin{example}
Let $H=\SL_2$. For $w=e$ being the unit, then $\widetilde{R}^T_{e}=C(H/U_H)$ is the universal principal series representation. 
Let $C(H/U_H)^0$ be the subspace consisting of all functions $f$ such that $\sum_{h\in H(\kappa)/U_H(\kappa)} f(h)=0$.
Let $Y\subset \bA^2$ be the Drinfeld curve.
For $w=s$ being the unique simple reflection, the representation $ \widetilde{R}^T_{\dot{s}}$ fits into the following short exact sequence
\[1\to H^1_c(Y,\La)\to \widetilde{R}^T_{\dot{s}}\to C(H/U_H)^0\to 0.
\]
\end{example}

Here is a direct consequence. We assume that $\La$ is an algebraically closed field.
Recall that a character $\theta: T^{w\sigma}\to \La^\times$ is called non-singular if $W_\chi^\circ$ is trivial, where $\chi: T^pT\to \La^\times$ corresponds to $\theta$, $W_\chi$ is the stabilizer subgroup of $\chi$ under the action of $W$ and $W_\chi^\circ\subset W_\chi$ is the subgroup generated by reflections. When $\theta$ is non-singular, one knows from \eqref{eq: non-zero monodromic category} that for every $v\in W$ such that $v\chi=\sigma(\chi)$, we have
\[
\Delta^{\mon}_{v,\hchi}=\nabla^{\mon}_{v,\hchi}=\Til^{\mon}_{v,\hchi}.
\]

\begin{corollary}
Let $\theta: T^{w\sigma}\to \overline\La^\times$ be a non-singular character. Then $R_{\dot{v},\hchi}^{\mon,!}=R_{\dot{v},\hchi}^{\mon,*}=R_{\dot{v},\hchi}^{\mon,T}$ concentrate in degree zero.
In particular, when $\La=\overline\bQ_\ell$, then $R_{\dot{v},\theta}^!=R_{\dot{v},\theta}^*$  concentrate in degree zero.
\end{corollary}

\begin{lemma}\label{lem: DL induction of tilting generates}
Let $\pi\in \rep(G,\La)^{\heartsuit}$. Then there is some $\widetilde{R}^T_{\dot{w}}$ and a non-zero map $\widetilde{R}^T_{\dot{w}}\to \pi$.
\end{lemma}
\begin{proof}
As mentioned above, the category $\rep(G(\kappa),\La)$ is generated by $\{\widetilde{R}^!_{\dot{w}}\}_{w\in W}$, and therefore is also generated by $\{\widetilde{R}^T_{\dot{w}}\}_{w\in W}$. Therefore, for every $\pi$, there is some  $\widetilde{R}^T_{\dot{w}}$ some such $\Hom(\widetilde{R}^T_{\dot{w}},\pi)\neq 0$. But if $\pi\in \rep(G,\La)^{\heartsuit}$, $\Hom(\widetilde{R}^T_{\dot{w}},\pi)$ concentrates in degree zero. So there is a non-zero map as desired.
\end{proof}

\begin{corollary}\label{lem: rep of H via DL of tilting}
Suppose $\La$ is an algebraically closed field. 
For every irreducible representation $\pi$, there is a minimal length element $w\in W$ (in its $\sigma$-conjugacy class in $W$) such that $\pi$ appears as quotient of $\widetilde{R}^T_{\dot{w}}$. In particular, when  $\La=\overline\bQ_\ell$, $\pi$ appears as a direct summand of $\widetilde{R}^T_{\dot{w}}$.
\end{corollary}

We also recall the following fact for later usage.

\begin{lemma}\label{lem: generation by cuspidal}
Let $\La$ be an algebraically closed field.
The category $\rep(G(\kappa))$ is generated by $\ind_{P(\kappa)}^{G(\kappa)}\pi$ (and their shifts), where $P\subset G$ is a standard (rational) parabolic subgroup and $\pi$ is a cuspidal irreducible representation of $L_P(\kappa)$.
\end{lemma}
\begin{proof}It is enough to prove that for every irreducible representation $V$ of $G(\kappa)$, there is some $(P,\pi)$ as in the lemma such that there is a non-zero map $\ind_{P(\kappa)}^{G(\kappa)}\pi\to V$. But this is standard. Find $P$ with unipotent radical $U_P$ such that $V_{U_P(\kappa)}\neq 0$ but $V_{U_{P'}(\kappa)}=0$ for any $P'\subsetneq P$. Then $V_{U_P(\kappa)}$ contains a cuspidal irreducible representation $\pi$ of $L_P(\kappa)$ as a \emph{subrepresentation}. Then we have a non-zero map $\ind_{P(\kappa)}^{G(\kappa)}\pi\to V$.
\end{proof}

\subsubsection{Deligne-Lusztig induction as categorical trace}\label{SS: Finite Deligne-Lusztig induction}
The following discussions serve as a warm-up for the discussions in the affine setting.

Recall that $\shv_{\mon}(U\bs G/U)$ and $\ind\cshv(B\bs G/B)$ and $\shv(B\bs G/B)$ are monoidal categories. As $G,B,U$ are in fact defined over $\kappa$, the $\sharp \kappa$-Frobenius $\phi$ induces monoidal auto-equivalences of these categories via $*$-pushforwards, still denoted by $\phi$. (See \eqref{eq: phi structure given by pushforward}.)

\begin{theorem}\label{thm: DL theory as categorical trace}
\begin{enumerate}
\item\label{thm: DL theory as categorical trace-1} The (monodromic) Deligne-Lusztig induction \eqref{eq: finite mon DL induction} induces an equivalence
\[
\tr(\shv_{\mon}(U\backslash G/U),\phi)\cong \rep(G(\kappa)).
\]
The unipotent Deligne-Lusztig induction \eqref{eq: finite unip DL induction} induces an equivalence
\[
\tr(\ind\cshv(B\backslash G/B,\La),\phi)\cong \ind\crep^{\unip}(G(\kappa)).
\]
\item\label{thm: DL theory as categorical trace-2} The above equivalences restrict to equivalences
\[
\tr(\shv\bigl((B,\hat{u})\backslash G/(B,\hat{u})\bigr),\phi)\cong \rep^{\widehat{\unip}}(G(\kappa))\cong \tr(\shv(B\backslash G/B),\phi).
\]
\end{enumerate}
\end{theorem}
\begin{proof}
We first prove Part \eqref{thm: DL theory as categorical trace-1}.
For the first case, the fully faithfulness follows from \Cref{cor-fully-faithful-usual-trace-to-loop}, applied to the sheaf theory as in  \Cref{prop: sheaf theory for monodromic sheaf} and \Cref{rmk: sheaf theory for monodromic sheaf}. 

More precisely, we take $\der=\shv_{\mon}$, and let $X=\bB U$ equipped with the natural action of $T$, and let $Y=\bB G$ equipped with the trivial $T$-action. Then as explained in \Cref{rmk: sheaf theory for monodromic sheaf},
the morphism $\bB U\to \bB G$ belongs to the class of morphisms $\mathrm{VR}$ associated to $\shv_{\mon}$, since $\bB U/T=\bB B\to \bB G$ is pfp proper. 
The map $\bB U\to \bB U\times \bB U$ is representable coh. smooth and therefore belongs to  the class of morphisms $\mathrm{HR}$ associated to $\shv_{\mon}$, by \Cref{prop: sheaf theory for monodromic sheaf}.
Thanks to \Cref{prop: categorical property of monodromic Hecke}, \Cref{prop-comparison-usual-trace-geo-trac} and therefore  \Cref{cor-fully-faithful-usual-trace-to-loop} is also applicable. In addition, note that $(\Nt^u)_*^{\mon}=(\Nt^u)_*$.  The fully faithfulness follows.
The essential surjectivity is equivalent to \Cref{thm: DL completeness}. 

For the second case, the fully faithfulness follows from \Cref{thm:trace-constructible} and  \Cref{prop-comparison-usual-trace-geo-trac}. 
Here we let $\der=\ind\cshv$ and $X=\bB B$ and $Y=\bB G$.
Thanks to \Cref{prop: categorical property of affine Hecke}, \Cref{prop-comparison-usual-trace-geo-trac} is applicable so the geometric categorical trace is identified with the usual categorical trace.
Note that for algebraic stacks of finite presentation over $k$, $\fgshv=\cshv$ so $\rshv=\ind\cshv$ is the ind-completion of the category of constructible sheaves.  Finally, the essential surjectivity follows from the definition of $\ind\fgrep^{\unip}(G(\kappa))\subset \ind\crep(G(\kappa))$.

Next we deal with Part \eqref{thm: DL theory as categorical trace-2}.
The fully faithfulness of $\tr(\shv\bigl((B,\hat{u})\backslash G/(B,\hat{u}\bigr),\phi)\to \rep^{\widehat{\unip}}(G(\kappa),\La)$ similarly follows from \Cref{cor-fully-faithful-usual-trace-to-loop}, applied to the sheaf theory as in  \Cref{prop: sheaf theory for monodromic sheaf} and \Cref{rmk: sheaf theory for monodromic sheaf}. The essential surjectivity follows from the definition of  $\rep^{\widehat{\unip}}(G(\kappa))$. 

Next we deal with the second equivalence. Again, the fully faithfulness follows from \Cref{cor-fully-faithful-usual-trace-to-loop}, applied to the sheaf theory $\shv$, and to $X\to Y$ being $\bB B\to \bB G$. 
Again, thanks to \Cref{prop: categorical property of affine Hecke}, \Cref{prop-comparison-usual-trace-geo-trac} is applicable so the geometric categorical trace is identified with the usual categorical trace.

\quash{Note that we cannot directly apply 
\Cref{bimodule.geom.trace.fully.faithfulness.convolution.cat} to compute $\tr(\shv(B\backslash G/B),\phi)$, as the $!$-pullback $\Mod_\La=\shv(\spec k)\to \shv(\bB B)$ along $\bB B\to \spec k$ does not preserve compact objects so its right adjoint is not a continuous functor nor satisfies base change as required in \Cref{bimodule.geom.trace.fully.faithfulness.convolution.cat} Assumption \eqref{bimodule.geom.trace.fully.faithfulness.convolution.cat-3}.  However, by  \Cref{rem: affine hecke vs ind-finite}, $\Psi^L:\shv(B\backslash G/B)\subset \rshv(B\backslash G/B)$ is a colocalization as monoidal categories which induces a fully faithful embedding  $\Psi^L: \tr(\shv(B\backslash G/B),\phi)\to \tr(\ind\cshv(B\backslash G/B),\phi)$.

In addition, we have $\Nt^{\ind\fg}_*\delta^{\ind\fg,!}\Psi^L\cong \Psi^L \Nt_*\delta^!$ (as $\Nt_*$ and $\delta^!$ preserve compact objects). It follows that we have the following commutative diagram
\[
\xymatrix{
\tr(\shv(B\backslash G/B),\phi)\ar_{\Psi^L}[d]\ar[r] & \shv(\frac{G}{\Ad_\sigma G}) \ar^{\Psi^L}[d]\\
\tr(\ind\cshv(B\backslash G/B),\phi)\ar[r] & \ind\cshv(\frac{G}{\Ad_\sigma G}),
}
\]
with both vertical functors, and the bottom functor fully faithful. 
}
We thus obtain a fully faithful embedding
\[
\tr(\shv(B\backslash G/B), \phi)\to \rep(G(\kappa)).
\]
The essential image is generated by $\{R^*_{w}\}_{w\in W}$ and therefore coincides with $\rep^{\widehat\unip}(G(\kappa))$ (see \Cref{rem: char zero two versions of unip representation are the same}).
\end{proof}

\begin{remark}
Let $\La=\overline\bQ_\ell$. We may choose a positive integer $n$ sufficiently large and coprime to $p$, and consider the monoidal category $\ind\cshv(B^{[n]}\backslash G/B^{[n]},\overline\bQ_\ell)$. Then we also have
\[
\tr(\ind\cshv(B^{[n]}\backslash G/B^{[n]},\overline\bQ_\ell),\phi)=\rep(G(\kappa),\overline\bQ_\ell).
\]
Indeed, as argued for the last case in \Cref{thm: DL theory as categorical trace}, the fully faithfulness follows from \Cref{thm:trace-constructible} and  \Cref{prop-comparison-usual-trace-geo-trac}. The essential surjectivity then follows from \Cref{thm: DL completeness}, the isomorphisms \eqref{eq: chimon-equal-to-theta} and the calculation \Cref{prop:affine DL induction chi}.

A version of this construction can be found in \cite{Lusztig2015unipotent, Lusztig2017}. On the other hand, a version more closely to \Cref{thm: DL theory as categorical trace} \eqref{thm: DL theory as categorical trace-1} has also appeared in \cite{eteve.monodromic} recently.
\end{remark}

\begin{remark}\label{rem: completeness of DL via Springer sheaf}
It is well-known to expert that there is also a more geometric argument of essential surjectivity of the functor $\tr(\shv_{\mon}(U\backslash G/U),\phi)\to \rep(G(\kappa))$, at least in the case when $\sharp W$ is invertible in $\La$. Namely, as observed by Mirkovi\'c-Vilonen (\cite{MV.character}), the functor $\Ch_{G,\phi}^{\mon}\circ (\Ch_{G,\phi}^{\mon})^R$ is (essentially) the same as convolution of the Springer sheaf. More precisely, we make use of \Cref{lem: CH HC=conv with Spr}. 
In this case $\mS$ is nothing but the $*$-pushforward of $\consdual_{\frac{U}{\Ad U}}$ along $f:\frac{U}{\Ad U}\to \frac{G}{\Ad G}$. We claim that 
\[
\Delta_1:= (\Delta_{\bB G/\bB G\times \bB G})_*\consdual_{\bB G}\in \shv(\frac{G}{\Ad G}), \quad \mbox{where } \Delta_{\bB G/\bB G\times \bB G}: \bB G\to \mL(\bB G)=\bB G\times_{\bB G\times \bB G}\bB G
\] 
is contained in the presentable stable subcategory generated by $\mS$. This would imply that  the right adjoint of the functor $(\Nt^u)_*(\delta^u)^!$ given by
$(\delta^u)_\flat\circ (\Nt^u)^!: \rep(G(\kappa))\to \shv_{\mon}(U\backslash G/U)$ is conservative, which in turn will imply that $(\Nt^u)_*(\delta^u)^!$ is essentially surjective.

To prove the claim, we factor $f$ as $\frac{U}{\Ad U}\xrightarrow{f_1} \frac{U}{\Ad B}\xrightarrow{f_2}  \frac{G}{\Ad G}$. Clearly, $\consdual_{\frac{U}{\Ad B}}$ is contained in the presentable subcategory generated $(f_1)_*\consdual_{\frac{U}{\Ad U}}$. In fact, by base change, this  follows from the universal situation that $\pi_*\La$ generates $\shv(\bB T)$, where $\pi: \pt \to \bB T$ is the universal $T$-torsor.

On the other hand, $(f_2)_*\consdual_{\frac{U}{\Ad B}}\in\shv(\frac{G}{\Ad G})$ is nothing but the usual Springer sheaf, which is a perverse sheaf equipped with an action of the finite Weyl group $W$. When $\sharp W$ is invertible in $\La$, 
we have $((f_2)_*\consdual_{\frac{U}{\Ad B}})\otimes_{W} \mathrm{triv}=\Delta_1$, which is a direct summand of $(f_2)_*\consdual_{\frac{U}{\Ad B}}$. This proves the claim.
\end{remark}

We also record the following results, which have been a folklore in the geometric representation theory community.

\begin{proposition}\label{prop: fully faithfulness of tensor over finite Hecke category}
Let $H_i\subseteq G, \ i=1,2$ be two closed subgroups, we have fully faithful embedding
  \begin{equation}\label{prop: fully faithfulness of tensor over finite Hecke category-1}
  \shv_{\mon}(H_1\backslash G/U)\otimes_{\shv_{\mon}(U\backslash G/U)}\shv_{\mon}(U\backslash G/H_2)\hookrightarrow \shv(H_1\backslash G/H_2).
  \end{equation}
  If one of $H_i$ is a  standard parabolic subgroup of $G$, then the above fully faithful embedding is an equivalence.

Similarly, we have fully faithful embedding
  \begin{equation}\label{prop: fully faithfulness of tensor over finite Hecke category-2}
  \ind\cshv(H_1\backslash G/B)\otimes_{\ind\cshv(B\backslash G/B)}\ind\cshv(B\backslash G/H_2)\hookrightarrow \ind\cshv(H_1\backslash G/H_2),
  \end{equation}
  which restricts to a fully faithful embedding
  \begin{equation}\label{prop: fully faithfulness of tensor over finite Hecke category-3}
  \shv(H_1\backslash G/B)\otimes_{\shv(B\backslash G/B)}\shv(B\backslash G/H_2)\hookrightarrow \shv(H_1\backslash G/H_2).
  \end{equation}
  If one of $H_i$ is a  standard parabolic subgroup of $G$, then \eqref{prop: fully faithfulness of tensor over finite Hecke category-3} is an equivalence.
\end{proposition}
\begin{proof}
We first deal with the monodromic case. 
As in the proof for \Cref{thm: DL theory as categorical trace}, we let $\der=\shv_{\mon}$ and $X=\bB U$ with the natural $T$-action, and $Y=\bB G$ with the trivial $T$-action.
Then fully faithfulness then follows from \Cref{cor-fully-faithful-usual-trace-to-loop}, by taking $W_i=\bB H_i$. (\Cref{lem: categorical kunneth prestack} guarantees the last assumption of  \Cref{cor-fully-faithful-usual-trace-to-loop} holds.) Next we prove the essential surjectivity. We may assume that $H_1=P_1$ is a standard parabolic subgroup. It is enough to show that the essential image of the $*$-pushforward $\shv_{\mon}(U\backslash G/H_2)\to \shv(P_1\backslash G/H_2)$ generate $\shv(P_1\backslash G/H_2)$, or its right adjoint is conservative. 
Using \Cref{rmk: sheaf theory for monodromic sheaf}, we see that up to shift the right adjoint is given by the $!$-pullback, and is conservative (by descent). The claim follows.

The fully faithfulness of \eqref{prop: fully faithfulness of tensor over finite Hecke category-2} is proved similarly. (\Cref{lem: kunneth for indfg sheaf} guarantees the last assumption of  \Cref{cor-fully-faithful-usual-trace-to-loop} holds.) It also implies the fully faithfulness of \eqref{prop: fully faithfulness of tensor over finite Hecke category-3} by the same argument as in the proof of \Cref{thm: DL theory as categorical trace}. If $H_1=P_1$ is a standard parabolic, the argument above also shows that that \eqref{prop: fully faithfulness of tensor over finite Hecke category-3} is essentially surjective.
\end{proof}

\begin{remark}\label{rem: constructible sheaf BB to BG}
For general coefficient $\La$.
We do not know whether  \eqref{prop: fully faithfulness of tensor over finite Hecke category-2} is an equivalence when $H_1=P_1$ is a standard parabolic, because we do not know whether in general the image of the functor $\ind\cshv(B\backslash G/H_2)\to \ind\cshv(P\backslash G/H_2)$ generates the target category. This is certainly the case in many situations, but for example, we do not know whether this is true when $P=H_2=G$, i.e. we do not know whether the image of the functor $\ind\cshv(\bB B)\to \ind\cshv(\bB G)$ generates $\ind\cshv(\bB G)$  in general. (See \Cref{lem: induction of parabolic belonging unipotent} \eqref{lem: induction of parabolic belonging unipotent-2} for a sufficient condition so that this holds.)
\end{remark}

\begin{remark}
We consider the $\shv(B\bs G/B)\mbox{-}\shv\bigl((U,\hat{u})\bs G/(U,\hat{u})\bigr)$-bimodule 
\[
\shv(B\bs G/U)=\shv_{\mon}(B\bs G/U).
\] 
Using \Cref{prop: fully faithfulness of tensor over finite Hecke category}, we see that
\[
\shv(B\bs G/U)\otimes_{\shv\bigl((U,\hat{u})\bs G/(U,\hat{u})\bigr)}\shv(U\bs G/B)\cong \shv(B\bs G/B),
\]
\[
\shv(U\bs G/B)\otimes_{\shv(B\bs G/B)}\shv(B\bs G/U)\cong \shv\bigl((U,\hat{u})\bs G/(U,\hat{u})\bigr).
\]
In $2$-categorical terms, this says that $\shv(B\bs G/B)$ and $\shv\bigl((U,\hat{u})\bs G/(U,\hat{u})\bigr)$ are Morita equivalent (see \Cref{rem-trace-2-category}). In any case,
it induces an equivalence $\tr(\shv\bigl((U,\hat{u})\bs G/(U,\hat{u})\bigr),\phi)\cong \tr(\shv(B\bs G/B),\phi)$, giving another proof of \Cref{thm: DL theory as categorical trace} \eqref{thm: DL theory as categorical trace-2}. 

On the other hand, if we let
\[
\shv(G)^0:=\shv_{\mon}(G/U)\otimes_{\shv_{\mon}(U\bs G/U)}\shv_{\mon}(U\bs G)\subset \shv(G).
\]
Then $\shv(G)^0$ has a natural monoidal structure such that the inclusion $\shv(G)^0\subset \shv(G)$ is non-unital monoidal. One can show that $\shv(G)^0$ and $\shv_{\mon}(U\backslash G/U)$ are Morita equivalent.
\end{remark}

\begin{remark}\label{rem: replacing Borel by parabolic}
Let $P$ be a standard rational parabolic subgroup of $G$.
In \Cref{thm: DL theory as categorical trace} if we replace $B$ by $P$, we will still have a fully faithful embedding
\[
\tr(\ind\cshv(P\bs G/P),\phi)\hookrightarrow \ind\crep(G(\kappa)),\quad \tr(\shv(P\bs G/P),\phi)\hookrightarrow \rep(G(\kappa)).
\]
The essential image of $\tr(\shv(P\bs G/P),\phi)$ is contained in $\rep^{\widehat\unip}(G(\kappa))$. But we do not know whether the essential image of $\tr(\ind\cshv(P\bs G/P),\phi)$ is contained in $\ind\crep^{\unip}(G(\kappa))$.
Note that tautologically, $\ind_{P(\kappa)}^{G(\kappa)}\La$ is contained in the former but as mentioned in \Cref{rem: char zero two versions of unip representation are the same}, we do not know whether it belongs to the latter.

Similarly, by replacing $B$ by $P$ in \eqref{prop: fully faithfulness of tensor over finite Hecke category-2} and \eqref{prop: fully faithfulness of tensor over finite Hecke category-3}, we obtain fully faithful embeddings
 \begin{equation}\label{prop: fully faithfulness of tensor over finite Hecke category-4}
  \ind\cshv(H_1\backslash G/P)\otimes_{\ind\cshv(P\backslash G/P)}\ind\cshv(P\backslash G/H_2)\hookrightarrow \ind\cshv(H_1\backslash G/H_2),
  \end{equation}
  which restricts to a fully faithful embedding
  \begin{equation}\label{prop: fully faithfulness of tensor over finite Hecke category-5}
  \shv(H_1\backslash G/P)\otimes_{\shv(P\backslash G/P)}\shv(P\backslash G/H_2)\hookrightarrow \shv(H_1\backslash G/H_2).
  \end{equation}
In particular, if we let $P=G$ and $H_1=H_2$ be the trivial group, we obtain fully faithful embeddings
\[
\Mod_\La\otimes_{\ind\cshv(\bB G)}\Mod_\La\subset \shv(G),\quad \Mod_\La\otimes_{\shv(\bB G)}\Mod_\La\subset \shv(G).
\]
The images of both embeddings are generated (as $\La$-linear categories) by the constant sheaf on $G$. 
\end{remark}

\subsubsection{Gelfand-Graev representations}
Finally we review some facts between Gelfand-Graev representations and Deligne-Lusztig representations. 
We fix a non-trivial additive characger $\psi: \kappa\to \La^\times$.
Let
\[
\psi_e: U(\kappa)\xrightarrow{e} \kappa\xrightarrow{\psi} \La^\times.
\]
The Gelfand-Graev representation of $G(\kappa)$ with respect to $\psi_e$ is defined as
\[
\GG_{\psi_e}=\ind_{U(\kappa)}^{G(\kappa)}\psi_e^{-1}.
\]
The following result is originally proved by Dudas (\cite{dudas.DL.restriction.Gelfand-Graev}) using Deodhar decomposition of Richardson variety. 
Using \Cref{lem-mon-spectral-sequence-hom}, we obtain a very short proof\footnote{We notice that Eteve also independently found a proof of Dudas' result similar to ours (see \cite{eteve.DL}).}. We refer to \Cref{prop:Whittaker-exactness} for a generalization of this result in the affine case.

\begin{proposition}\label{prop:multiplicity one GG rep}
We have a natural isomorphism $\Hom_{G(\kappa)}(\widetilde{R}^!_w, \GG_{\psi_e})\cong \La[T^{w\sigma}]$.
\end{proposition}

\begin{proof}
In the finite dimensional case, we have the following equivalence of categories
\[
\Delta_{\dot{w}_0}^{\mon,\psi_e}=\nabla_{\dot{w}_0}^{\mon,\psi_e}\colon \shv_{\mon}(T)\cong \shv_{\mon}\bigl(U\backslash  G/(U,\psi_e)\bigr),
\]
where $w_0$ denotes the longest length element in the Weyl group $W$.
By \Cref{lem-mon-spectral-sequence-hom}, and
using \Cref{cor: av upper star av lower star}, the isomorphisms \eqref{eq: universal monodromic DL induction formula},  we have 
\begin{eqnarray*}
\Hom_{G(\kappa)}(\widetilde{R}^*_{\dot{w}}, \GG_{\overline\psi_e}) & \cong & \Hom_{\shv(U\backslash G/(U,\psi_e))}( (\av_s)^*(\av_s)_*\widetilde\nabla_{\dot{w}}^{\mon}\star^u \widetilde{\Delta}_{\dot{w}_0}^{\mon,\psi_e}, \widetilde{\Delta}_{\dot{w}_0}^{\mon,\psi_e})\\
&\cong &\Hom_{\shv_{\mon}(U\backslash G/(U,\psi_e))}(\nabla_{\dot{w}, \chi_{\varphi_{w}}}^{\mon}\star^u \widetilde{\Delta}_{\dot{w}_0}^{\mon,\psi_e},  \widetilde{\Delta}_{\dot{w}_0}^{\mon,\psi_e})\\
&\cong &\Hom_{\shv_{\mon}(U\backslash G/(U,\psi_e))}(\Delta_{\dot{w}_0}^{\mon,\psi_e}(\Ch_{\chi_{\varphi_w}}),  \widetilde{\Delta}_{\dot{w}_0}^{\mon,\psi_e})\\
&\cong &\Hom_{\shv_{\mon}(T)}(\Ch_{\chi_{\varphi_w}},  \widetilde{\Ch}) \\
&\cong &\Hom_{\qcoh(\chi_{\varphi_w})}(\mO,\mO)=\La[T^{w\sigma}].
\end{eqnarray*}
The proposition is proved.
\end{proof}

\subsection{Tame and unipotent local Langlands category}

Now we generalize the discussions in the previous subsection to the affine case.

\subsubsection{Definition and first properties}\label{SS: 1st approach tame Langlands category}

Recall that there is a notion of ``depth" in the representation theory of $p$-adic groups. We shall not review the general definition here, but only to review the notion of depth zero representations.
As before, we omit the coefficient $\La$ from the notations.

\begin{definition}
Let $H$ be a connected reductive group over $F$ (but we do not assume that it is quasi-split).
Let $\rep^\tame(H(F))\subset \rep(H(F))$ be the presentable $\La$-linear stable subcategory generated by $\cind_{P^u}^{H(F)}\La$, where $P^u$ is the pro-$p$-radical of a parahoric subgroup $P$ of $H(F)$. Objects in $\rep^\tame(H(F))$ are called depth zero representations of $H(F)$. 

On the other hand, we say an object $V\in \rep(H(F))$ has positive depth if $V^{P^u}=0$ for every pro-$p$-radical of a parahoric subgroup $P$ of $H(F)$. The full subcategory of positive depth representations of $H(F))$ is denoted by $\rep^{>0}(H(F))$.
\end{definition}

We note that 
\[
\rep^\tame(H(F))^{\heartsuit}:=\rep^\tame(H(F))\cap \rep(H(F))^{\heartsuit}
\] 
is an abelian subcategory, with a set of projective generators given by $\{\cind_{P^u}^{H(F)}\La\}_P$, where $P$ range over the set of parahoric subgroups of $H(F)$. In addition, we have
\[
\rep^\tame(H(F))=\der(\rep^\tame(H(F))^{\heartsuit}).
\]

Similarly, let $\rep^{>0}(H(F))^{\heartsuit}=\rep^{>0}(H(F))\cap \rep(H(F),\La)^{\heartsuit}$. Then $\rep^{>0}(H(F))=\der(\rep^{>0}(H(F))^{\heartsuit})$.
It is well-known that there is an orthogonal decomposition of categories
\[
\rep(H(F))^{\heartsuit}=\rep^\tame(H(F))^{\heartsuit}\bigoplus \rep^{>0}(H(F))^{\heartsuit},
\]
which then induces an orthogonal decomposition
\begin{equation}\label{eq: depth zero decomposition}
\rep(H(F))=\rep^\tame(H(F))\bigoplus \rep^{>0}(H(F)).
\end{equation}
We let 
\[
\proj^{\tame}: \rep(H(F))\to \rep^\tame(H(F))
\] 
be the continuous right (and the left) adjoint functor of the natural inclusion. 
The decomposition also restricts to the decomposition of the subcategories of compact, admissible and finitely generated objects. It also induces a decomposition 
\[
\mathrm{tr}(\rep(H(F)))=\mathrm{tr}(\rep^\tame(H(F)))\bigoplus \mathrm{tr}(\rep^{>0}(H(F)))
\]
and in particular a decomposition of cocenter of Hecke algebras (once a Haar measure of $H(F)$ is chosen)
\[
C_c^\infty(H(F))_{H(F)}=C_c^{\infty}(H(F))^{\tame}_{H(F)}\bigoplus C_c^{\infty}(H(F))^{>0}_{H(F)}.
\]
If $\pi\in \rep^\tame(H(F))\cap \rep(H(F))^{\adm}$, then its character $\Theta_\pi: C_c^\infty(H(F))_{H(F)}\to \La$ factors as $C_c^\infty(H(F))_{H(F)}\to C_c^{\infty}(H(F))^{\tame}_{H(F)}\to \La$ so we will also regard $\Theta_\pi$ as a functional on $C_c^{\infty}(H(F))^{\tame}_{H(F)}$.

\begin{definition}\label{def:tame-local-Langlands-category}
We let $\shv^{\tame}(\kot_G)\subset \shv(\kot_G)$ be the full subcategory spanned by objects $\mF$ such that for every $b\in B(G)$, we have $(i_b)^!\mF\in \rep^{\tame}(G_b(F))$, and call it the tame local Langlands category.
For $b\in B(G)$ and  for $?$ being either $<$ or $\leq$, we let $\shv^{\tame}(\kot_{G,? b})=\shv^{\tame}(\kot_G)\cap \shv(\kot_{G,? b})$.

We similarly define $\shv^{>0}(\kot_G)\subset \shv(\kot_G)$ be the full subcategory spanned by objects $\mF$ such that $(i_b)^!\mF\in \rep^{>0}(G_b(F))$ for every $b\in B(G)$.
\end{definition}

It follows directly from the definition that
\[
\shv^{\tame}(\kot_G,\La)=\colim_{b\in B(G)} \shv^{\tame}(\kot_{G,\leq b},\La)
\]
is compactly generated, with a set of compact generators given by $\bigl\{(i_b)_*\cind_{P^u}^{G_b(F)}\La\bigr\}$, for $b\in B(G)$ and $P\subset G_b(F)$ parahoric. We still let 
\[
\proj^\tame: \shv(\kot_G)\to \shv^\tame(\kot_G)
\] 
denote the continuous right adjoint of the natural inclusion $\shv^{\tame}(\kot_G)\subset \shv(\kot_G)$. 

Similarly, 
$\shv^{>0}(\kot_G)$ is also compactly generated, with a set of compact generators given by $\bigl\{(i_b)_*\pi\bigr\}$, for $b\in B(G)$ and $\pi\in \rep^{>0}(G_b(F))^\cpt$.

Later on we will prove the following result.
\begin{proposition}\label{prop: can duality of tame category}
The category $\shv^{\tame}(\kot_G)^\cpt$ is stable under the canonical duality $(\verd^{\can}_{\kot_G})^\cpt$ from \Cref{cohomological.duality.kottwitz}. 
\end{proposition}

We will let $(\verd_{\kot_G}^{\tame,\can})^\cpt$ denote the restriction of $(\verd_{\kot_G}^{\can})^\cpt$ to $\shv^{\tame}(\kot_G)^\cpt$, and let
\begin{equation}\label{eq: canonical duality tame part}
\verd_{\kot_G}^{\tame,\can}: \shv^{\tame}(\kot_G)^\vee\cong \shv^{\tame}(\kot_G)
\end{equation}
denote the induced self-duality of $\shv^{\tame}(\kot_G)$.

We thus see that semi-orthogonal decompositions in \Cref{cor:semi-orthogonal decomp for Shv(kot(G))} restrict to tame subcategories.

\begin{corollary}\label{cor:semi-orthogonal decomp for tame Shv(kot(G))}
\begin{enumerate}
\item\label{cor:semi-orthogonal decomp for tame Shv(kot(G))-1} The category $\shv^{\tame}(\kot_G,\La)$ admits a set of compact generators given by $\bigl\{(i_b)_!\cind_{P^u}^{G_b(F)}\La\bigr\}$, for $b\in B(G)$ and $P\subset G_b(F)$ parahoric.
\item\label{cor:semi-orthogonal decomp for tame Shv(kot(G))-2} The pairs of adjoint functors $((i_b)^*,(i_b)_*)$ and $((i_b)_!,(i_b)^!)$ restrict to pairs of adjoint functors between $\rep^{\tame}(G_b(F),\La)$ and $\shv^{\tame}(\kot_G,\La)$. The semi-orthogonal decompositions in  \Cref{cor:semi-orthogonal decomp for Shv(kot(G))} restrict to semi-orthogonal decompositions
\[
\shv^{\tame}(\kot_{G,<b})\xrightarrow{(i_{<b})_*} \shv^{\tame}(\kot_{G,\leq b}) \xrightarrow{(j_b)^!} \shv^{\tame}(\kot_{G,b})
\]
\[
\shv^{\tame}(\kot_{G,b})\xrightarrow{(j_{b})_!} \shv^{\tame}(\kot_{G,\leq b}) \xrightarrow{(i_b)^*} \shv^{\tame}(\kot_{G,<b}).
\]

\item\label{cor:semi-orthogonal decomp for tame Shv(kot(G))-3} 
An object $\mF\in \shv(\kot_G)$ belongs to $\shv^{\tame}(\kot_G)$ if and only if $(i_b)^*\mF\in \rep^{\tame}(G_b(F))$ for every $b\in B(G)$.
\item For every $b\in B(G)$, we have
\[
(i_b)^!\circ \proj^{\tame}\cong \proj^{\tame}\circ (i_b)^!: \shv(\kot_G)\to \rep^{\tame}(G_b(F)), 
\]
\[
(i_b)_*\circ \proj^{\tame}\cong \proj^{\tame}\circ (i_b)_*:  \rep(G_b(F))\to \shv^{\tame}(\kot_G).
\]
In particular, $\proj^{\tame}$ preserves compact objects. On the other hand, being a right adjoint functor, $\proj^{\tame}$ also preserves admissible objects.
\item\label{cor:semi-orthogonal decomp for tame Shv(kot(G))-5} We have $\shv^{>0}(\kot_G)=\ker\proj^\tame$. The inclusion $\shv^{>0}(\kot_G)\subset \shv(\kot_G)$ admits a left adjoint functor, inducing a semi-orthogonal decomposition
\[
\shv^{\tame}(\kot_G)\to \shv(\kot_G)\to \shv^{>0}(\kot_G).
\]
\end{enumerate}
\end{corollary}
\begin{proof}Part \eqref{cor:semi-orthogonal decomp for tame Shv(kot(G))-1} follows directly from \Cref{prop: can duality of tame category} and \Cref{cohomological.duality.kottwitz.pullpush}. The rest parts follow easily.
\end{proof}

\begin{remark}\label{rmk: orthogonal decomposition of shvkot}
Unfortunately, we could not prove that
$\shv(\kot_G)$ admits an orthogonal decomposition by $\shv^{\tame}(\kot_G)$ and $\shv^{>0}(\kot_G)$, although this should be the case, as predicted by the categorical local Langlands conjecture. (See \eqref{eq: depth decomposition spectral side} for the decomposition in the spectral side.) We list a few statements that are equivalent to this orthogonal decomposition.
\begin{enumerate}
\item\label{rmk: orthogonal decomposition of shvkot-1} $\shv(\kot_G)=\shv^{\tame}(\kot_G)\oplus \shv^{>0}(\kot_G)$;
\item\label{rmk: orthogonal decomposition of shvkot-2} $\shv^{>0}(\kot_G)^\cpt$ is stable under the duality $(\verd_{\kot_G}^{\can})^\cpt$;
\item\label{rmk: orthogonal decomposition of shvkot-6} For every $b$, $(i_b)^*$ sends $\shv^{>0}(\kot_G)$ to $\rep^{>0}(G_b(F))$;
\item\label{rmk: orthogonal decomposition of shvkot-5} For every $b$, $(i_b)_!$ sends $\rep^{>0}(G_b(F))$ to $\shv^{>0}(\kot_G)$;
\item\label{rmk: orthogonal decomposition of shvkot-8} For every $b$, $(i_b)\rstar$ sends $\shv^{\tame}(\kot_G)$ to $\rep^{\tame}(G_b(F))$;
\item\label{rmk: orthogonal decomposition of shvkot-7} For every $b$, $(i_b)_\flat$ sends $\rep^{\tame}(G_b(F))$ to $\shv^{\tame}(\kot_G)$;
\item\label{rmk: orthogonal decomposition of shvkot-4} For every $b$, $(i_b)^*\circ \proj^{\tame}\cong \proj^{\tame}\circ (i_b)^*$;
\item\label{rmk: orthogonal decomposition of shvkot-3} For every $b$, $(i_b)_!\circ \proj^{\tame}\cong \proj^{\tame}\circ (i_b)_!$.
\end{enumerate}

We sketch a proof of these equivalences. Let $\mF\in \shv^{>0}(\kot_G)^\cpt$. Then $(\verd^{\can}_{\kot_G})^\cpt(\mF)\in  \shv^{>0}(\kot_G)^\cpt$ if and only if $\Hom(\mG, (\verd^{\can}_{\kot_G})^\cpt(\mF))=\Hom(\mF, (\verd^{\can}_{\kot_G})^\cpt(\mG))=0$ for every $\mG\in \shv^{\tame}(\kot_G)^\cpt$. Using \Cref{prop: can duality of tame category}, we see that this is the case
if and only if $\shv(\kot_G)=\shv^{\tame}(\kot_G)\oplus \shv^{>0}(\kot_G)$. Therefore, \eqref{rmk: orthogonal decomposition of shvkot-1} and \eqref{rmk: orthogonal decomposition of shvkot-2} are equivalent. 

Note that the canonical duality $(\verd_{G_b(F)}^{\can})^\cpt$ of $\rep(G_b(F))^{\cpt}$ preserves $\rep^{>0}(G_b(F))^\cpt$. 
By \Cref{cohomological.duality.kottwitz.pullpush}, 
$(i_b)^!(\verd^{\can}_{\kot_G})^\cpt(\mF)=(\verd^{\can}_{G_b(F)})^\cpt((i_b)^*\mF)[d]$ for some integer $d$ so $(\verd^{\can}_{\kot_G})^\cpt(\mF)\in  \shv^{>0}(\kot_G)^\cpt$ if and only if $(i_b)^*\mF\in  \rep^{>0}(G_b(F))^\cpt$. Thus
\eqref{rmk: orthogonal decomposition of shvkot-2} and \eqref{rmk: orthogonal decomposition of shvkot-6} are equivalent. As $\shv^{>0}(\kot_G)^\cpt$ is generated by compact objects of the form $(i_b)_*\pi$ with $\pi \in \rep^{>0}(G_b(F))^\cpt$, by  \Cref{cohomological.duality.kottwitz.pullpush} again, \eqref{rmk: orthogonal decomposition of shvkot-2} is equivalent to \eqref{rmk: orthogonal decomposition of shvkot-5}.

We also note that \eqref{rmk: orthogonal decomposition of shvkot-6} holds if and only if $\Hom( (i_b)^*((i_{b'})_*\pi), \pi')=\Hom(\pi, (i_{b'})^\sharp((i_b)_*\pi'))=0$ for every $b,b'$, every $\pi\in \rep^{>0}(G_b(F))$ and every $\pi'\in \rep^{\tame}(G_{b'}(F))$, if and only if \eqref{rmk: orthogonal decomposition of shvkot-8} holds.  Similarly, \eqref{rmk: orthogonal decomposition of shvkot-5} holds if and only if $\Hom(\pi',(i_{b'})^!((i_b)_\flat\pi))=\Hom((i_b)^!(i_{b'})_!\pi',\pi)=0$ for every $\pi\in \rep^{\tame}(G_b(F))$ and $\pi'\in \rep^{>0}(G_{b'}(F))$, if and only if \eqref{rmk: orthogonal decomposition of shvkot-7} holds.

Next, for $\mG\in \shv(\kot_G)$, consider the cofiber sequence $\proj^{\tame}\mG\to \mG\to \mG'$, which induces $(i_b)^*(\proj^{\tame}\mG)\to (i_b)^*\mG\to (i_b)^*\mG'$. Note that $\mG'\in \shv^{>0}(\kot_G)$ and $(i_b)^*(\proj^{\tame}\mG)\in \rep^{\tame}(G_b(F))$. Therefore, $(i_b)^*(\proj^{\tame}\mG)\cong \proj^{\tame}((i_b)^*\mG)$ if $(i_b)^*\mG'\in \rep^{>0}(G_b(F))$, showing that \eqref{rmk: orthogonal decomposition of shvkot-6} implies \eqref{rmk: orthogonal decomposition of shvkot-4}. Conversely, if $\eqref{rmk: orthogonal decomposition of shvkot-4}$ holds, then for every $\mG\in \shv^{>0}(\kot_G)$ and every $b\in B(G)$, we have $\proj^{\tame}((i_b)^*\mG)=(i_b)^*(\proj^{\tame}\mG)=0$. Therefore, \eqref{rmk: orthogonal decomposition of shvkot-6} holds.
Similarly, for $\pi\in\rep(G_b(F))$, we have a cofiber sequence $(i_b)_!(\proj^{\tame}\pi)\to (i_b)_!\pi\to (i_b)_!\pi'$ for $\pi'\in \rep^{>0}(G_b(F))$. Then $(i_b)_!(\proj^{\tame}\pi)\cong \proj^{\tame}((i_b)_!\pi)$ if $(i_b)_!\pi'\in \shv^{>0}(\kot_G)$. Therefore, \eqref{rmk: orthogonal decomposition of shvkot-5} implies \eqref{rmk: orthogonal decomposition of shvkot-3}. Conversely, if \eqref{rmk: orthogonal decomposition of shvkot-3} holds, then for every $b\in B(G)$ and $\pi \in \rep^{>0}(G_b(F))$, $\proj^{\tame}((i_b)_!\pi)\cong(i_b)_!(\proj^{\tame}\pi)=0$. Therefore, \eqref{rmk: orthogonal decomposition of shvkot-5} holds.
\end{remark}

\begin{remark}
Despite of \Cref{rmk: orthogonal decomposition of shvkot}, the decomposition 
\[
\mathrm{tr}(\Shv(\kot_G))=\mathrm{tr}(\Shv^{\tame}(\kot_G))\oplus \mathrm{tr}(\Shv^{>0}(\kot_G))
\] 
induced by \Cref{cor:semi-orthogonal decomp for tame Shv(kot(G))} \eqref{cor:semi-orthogonal decomp for tame Shv(kot(G))-5} is compatible with decompositions from
\Cref{cor: decomposition HH of shvkot}. 
\end{remark}

\begin{remark}
More generally, as mentioned before, there is a depth filtration $\rep(H(F))=\cup_r \rep^{\leq r}(H(F))$ of the category of smooth representations of a $p$-adic group $H(F)$. One can then similarly define $\shv^{\leq r}(\kot_G)$ as the full subcategory of $\shv(\kot_G)$ consisting of those $\mF$ such that $(i_b)^!\mF\in \rep^{\leq r}(G_b(F))$ for every $r$. Then
$\shv(\kot_G)$ admits a depth filtration 
\[
\shv(\kot_G)=\cup_{r\geq 0} \shv^{\leq r}(\kot_G),
\] 
and each $\shv^{\leq r}(\kot_G)$ admit a semi-orthogonal decomposition indexed by $\{(i_b)_*(\rep^{\leq r}(G_b(F)))\}_{b\in B(G)}$. However, our later proof of \Cref{prop: can duality of tame category} does not generalize to $\shv^{\leq r}(\kot_G)$, and we do not know whether the above functor $(i_{b'})^!(i_b)_!$ preserves the depth filtration.
\end{remark}

\begin{lemma}\label{lem: tame admissible objects}
The category $\shv^\tame(\kot_G)^{\adm}:=\shv^\tame(\kot_G)\cap \shv(\kot_G)^{\adm}$ consist of admissible objects of $\shv^{\tame}(\kot_G)$.
\end{lemma}
\begin{proof}
Obviously objects in $\shv^\tame(\kot_G)\cap \shv(\kot_G)^{\adm}$ are admissible objects in $\shv^{\tame}(\kot_G)$. Now suppose $\mF$ is an admissible object in $\shv^{\tame}(\kot_G)$. We need to show that  when regarded as an object in $\shv(\kot_G)$, it is still admissible. 

Let $P^u$ be the pro-$p$-radical of a parahoric subgroup of $G_b(F)$. Not that every $V\in \crep(P^u)$ admits a decomposition $V=V_0\oplus V_1$ such that for every field $E$ over $\La$, $V_E=(V_0)_E\oplus (V_1)_E$ is a decomposition of $V_E$ in terms the trivial and non-trivial representations of $P^u$. Now, we have
\[
\Hom((i_b)_!\cind_{P^u}^{G_b(F)}V, \mF)\cong \Hom(\cind_{P^u}^{G_b(F)}V, (i_b)^!\mF)\cong \Hom(\cind_{P^u}^{G_b(F)}V_0, (i_b)^!\mF),
\] 
which is a perfect $\La$-module. Therefore, $\mF\in  \shv(\kot_G)^{\adm}$. 
\end{proof}

Recall that the canonical duality $\verd_{\kot_G}^{\can}: \shv(\kot_G)^\vee\cong \shv(\kot_G)$ also restricts to an equivalence $(\verd_{\kot_G}^{\can})^\adm: (\shv(\kot_G)^\adm)^{\op}\cong \shv(\kot_G)^\adm$. On the other hand, the self-dualty \eqref{eq: canonical duality tame part} induces an equivalence
\begin{equation}\label{eq: admissible duality tame part}
(\verd_{\kot_G}^{\tame,\can})^\adm: (\shv^{\tame}(\kot_G)^\adm)^{\op}\cong \shv^{\tame}(\kot_G)^\adm
\end{equation}

\begin{lemma}
The equivalence \eqref{eq: admissible duality tame part} is identified with the functor
\[
(\shv^{\tame}(\kot_G)^\adm)^{\op} \subset (\shv(\kot_G)^\adm)^{\op} \xrightarrow{(\verd_{\kot_G}^{\can})^\adm}\shv(\kot_G)^\adm\xrightarrow{\proj^{\tame}} \shv^{\tame}(\kot_G)^\adm.
\] 
If equivalent conditions in \Cref{rmk: orthogonal decomposition of shvkot} hold, then $(\verd_{\kot_G}^{\can})^\adm$ restricts to an anti-equivalence of $\shv^{\tame}(\kot_G)^\adm$.
\end{lemma}
\begin{proof}
Since (by definition) the inclusion $\shv^{\tame}(\kot_G)^\cpt\subset \shv(\kot_G)^\cpt$ is compatible with the canonical duality,
we have $(\verd_{\kot_G}^{\tame,\can})^\adm\circ \proj^{\tame}\cong \proj^{\tame}\circ (\verd_{\kot_G}^{\can})^\adm$ by \Cref{lem: adm duality and functors}. The first statement follows.
The second statement follows from \Cref{prop: adm.duality.kottwitz.pullpush} and the fact that $\rep^{\tame}(G_b(F))^{\adm}$ is stable under the usual smooth duality.
\end{proof}

Finally, we also have a $t$-structure on $\shv^{\tame}(\kot_G)$ which restricts to a $t$-structure of $\shv^{\tame}(\kot_G)^\adm$, as in \Cref{prop: t-structures on llc-2}. Namely, by \Cref{cor:semi-orthogonal decomp for tame Shv(kot(G))}, the diagram \eqref{eq: t-structure-open-closed gluing-1} restricts to  a diagram with the subscript $\tame$ added to everywhere.
Then passing to the right adjoints we obtain a tame version of \eqref{eq: t-structure-open-closed gluing-2}
 \begin{equation}\label{eq: t-structure-open-closed gluing-tame}
    \xymatrix{
    \shv^{\tame}(\kot_{G,<b})\ar@/^/[rr]^{(i_{<b})_*} \ar@/_/[rr]_{\proj^{\tame}\circ (i_{<b})_\flat} && \ar[ll]|{(i_{<b})^!} \shv^{\tame}(\kot_{G,\leq b}) \ar@/^/[rr]^{(j_{b})^!} \ar@/_/[rr]_{\proj^{\tame}\circ(j_{b})\rstar} &&  \ar[ll]|{(j_b)_*}\shv^{\tame}(\kot_{G,b}).
    }
\end{equation}  
Note that if equivalent conditions in \Cref{rmk: orthogonal decomposition of shvkot} hold, we can remove $\proj^{\tame}$ in the above diagram. Now we can argue as in \Cref{prop: t-structures on llc-2} to define a $t$-structure on $\shv^{\tame}(\kot_G)$ with
\[
\shv^{\tame}(\kot_G)^{\chi\mbox{-}e,\leq 0}\subset \shv(\kot_G),\quad  \mbox{resp.} \quad \shv^{\tame}(\kot_G)^{\chi\mbox{-}e,\geq 0}\subset \shv(\kot_G)
\]
consisting of those $\mF$ such that 
\[
(i_b)^!\mF\in\rep^{\tame}(G_b(F))^{\leq \langle \chi,\nu_b\rangle },\quad  \mbox{resp.} \quad \proj^{\tame}((i_b)\rstar\mF)\in\rep^{\tame}(G_b(F))^{\geq \langle \chi,\nu_b\rangle}.
\] 
Again, if equivalent conditions in \Cref{rmk: orthogonal decomposition of shvkot} hold, we can remove $\proj^{\tame}$ in the above definition.

\begin{lemma}\label{lem: projection to tame is t-exact}
The functor $\proj^{\tame}: \shv(\kot_G)\to \shv^{\tame}(\kot_G)$ is $t$-exact.
\end{lemma}
\begin{proof}
Suppose $\mF\in \shv(\kot_G)^{\chi\mbox{-}e, \heartsuit}$. We have $(i_b)^!(\proj^{\tame}\mF)=\proj^{\tame}((i_b)^!\mF)\in \rep^{\tame}(G_b(F),\La)^{\leq\langle\chi,\nu_b\rangle}$. On the other hand, we have $\proj^\tame((i_b)\rstar(\proj^{\tame}\mF))\cong \proj^{\tame}((i_b)\rstar\mF)\in  \rep^{\tame}(G_b(F),\La)^{\geq\langle\chi,\nu_b\rangle}$. The lemma follows.
\end{proof}

Similar to \Cref{lem: adm duality and e-t-structure}, we have the following statement.
\begin{proposition}\label{lem: tame adm t-str self-duality}
Suppose $\La$ is a field and let $\chi=2\rho$. Then $(\verd^{\tame,\can}_{\kot_G})^{\adm}$ interchanges $\shv(\kot_G)^{\chi\mbox{-}e,\leq 0}\cap \shv^\tame(\kot_G)^{\adm}$ and $\shv(\kot_G)^{\chi\mbox{-}e,\geq 0}\cap \shv^\tame(\kot_G)^{\adm}$.
\end{proposition}

We also have parallel stories for the unipotent part. Recall \Cref{def: hunip and unip for finite group}.

\begin{definition}
Let $H$ be a connective reductive group over $F$. We let $\rep^{\widehat\unip}(H(F))\subset \rep^{\tame}(H(F))$ be the full subcategory generated by objects of the form $\cind_P^{H(F)}\pi_P$, where $P\subset H(F)$ is a parahoric with $L_P$ its Levi quotient, and $\pi_P\in \rep^{\widehat\unip}(L_P)$, regarded as a representation of $P$ via inflation. We let $\fgrep^{\unip}(H(F))\subset \fgrep(H(F))$ be  the full subcategory generated by objects of the form $\cind_P^{H(F)}\pi_P$, where $P\subset H(F)$ is a parahoric with $L_P$ its Levi quotient, and $\pi_P\in \crep^{\unip}(L_P)$, regarded as a representation of $P$ via inflation.
\end{definition}

\begin{remark}
Clearly we have 
\[
\fgrep^{\unip}(H(F))\subset \rep^{\widehat\unip}(H(F))\cap \fgshv(H(F)),
\] 
but we do not know whether the inclusion is an equality. (Compare to  the last sentence of \Cref{rem: char zero two versions of unip representation are the same}.)
In addition, we have $\rep^{\widehat\unip}(H(F))^\cpt\subset \fgrep^{\unip}(H(F))$ by \Cref{rem: char zero two versions of unip representation are the same}, but the inclusion is strict in general.
\end{remark}

\begin{remark}\label{rem:unipotent Bernstein blocks and Hecke algebra}
When $\La = \overline\bQ_\ell$, we have $\rep^{\widehat\unip}(H(F),\overline\bQ_\ell)^\cpt=\fgrep^{\unip}(H(F),\overline\bQ_\ell)$ by \Cref{rem: char zero two versions of unip representation are the same}.
Irreducible objects in $\fgrep^{\unip}(H(F),\overline\bQ_\ell)^{\heartsuit}$ are just irreducible unipotent representations of $H(F)$ in the sense of \cite{Lusztig.classification.unipotent.IMRN}. Indeed, an irreducible representation $\pi$ of $H(F)$ is called unipotent in \emph{loc. cit.} if it appears as a quotient of some $\cind_P^{H(F)}\pi_P$ for some parahoric subgroup $P$ of $H(F)$ and some irreducible cuspidal unipotent representation $\pi_P$ of $L_P$. But such $(P,\pi_P)$ is a type of $H(F)$ so $\pi$ in fact admits a finite free resolution by $\cind_P^{H(F)}\pi_P$'s. Therefore, $\pi$ indeed belongs to $\fgrep^{\unip}(H(F),\overline\bQ_\ell)$.
In addition, 
\[
\rep^{\widehat\unip}(H(F),\overline\bQ_\ell)=\bigoplus_{(P,\pi_P)/\sim} \lmodu_{H(P,\pi_P)}(\Mod_{\overline\bQ_\ell}),
\]
is a finite union of Bernstein blocks of $\rep(H(F),\overline\bQ_\ell)$, usually called the unipotent blocks. Here $(P,\pi_P)$ range over all pairs as above and $H(P,\pi_P)=\End(\cind_P^{H(F)}\pi_P)$ is the Hecke algebra associated to $(P,\pi_P)$, and
$(P_1,\pi_{P_1})\sim (P_2,\pi_{P_2})$ if there is some $g\in G(F)$ such that $gP_1g^{-1}=P_2$ and $\pi_{P_2}(-)=\pi_{P_1}(g^{-1}(-)g)$.
\end{remark}

\begin{definition}
We let $\shv^{\widehat\unip}(\kot_G)\subset\shv^{\tame}(\kot_G)$ be the full subcategory spanned by those $\mF$ such that $(i_b)^!\mF\in \rep^{\widehat\unip}(G_b(F))$ for every $b\in B(G)$. We let $\fgshv^{\unip}(\kot_G)\subset\fgshv(\kot_G)$ be the full subcategory spanned by those $\mF$ such that $(i_b)^!\mF\in \fgrep^{\unip}(G_b(F))$ for every $b\in B(G)$.
\end{definition}

\begin{example}\label{ex: dualizing sheaf on kot belongs to unip}
The sheaf $\consdual_{\kot_G}$ belongs to $\shv^{\widehat\unip}(\kot_G)$.
\end{example}

\begin{example}\label{ex: deltaP are unipotent}
Suppose $\sharp W_0$ is invertible in $\La$. Let $P_b\subset G_b(F)$ be a parahoric subgroup of $G_b(F)$. Then 
\[
(i_b)_*\cind_{P_b}^{G_b(F)}\La,\quad (i_b)_!\cind_{P_b}^{G_b(F)}\La
\]
belongs to $\fgshv^{\unip}(\kot_G)$. This follows from \Cref{lem: induction of parabolic belonging unipotent-2}.
\end{example}

By definition, the natural inclusion  $\rshv^{\unip}(\kot_G)\subseteq \rshv(\kot_G)$ preserves compact objects and therefore admits a continuous right adjoint 
\begin{equation}\label{eq:projection-to-unipotent-part}
    \proj^{\unip}: \rshv(\kot_G)\to \rshv^{\unip}(\kot_G).
\end{equation}

We similarly have the following statement.

\begin{proposition}\label{prop: can duality of unip category}
The category $\shv^{\widehat\unip}(\kot_G)^\cpt$ is stable under the canonical duality $(\verd^{\can}_{\kot_G})^\cpt$ from \Cref{cohomological.duality.kottwitz}. The category $\fgshv^{\unip}(\kot_G)$ is stable under the duality $(\verd^{\can}_{\kot_G})^\fg$ from \Cref{rem: can duality on f.g. sheaves}. The analogous statements in \Cref{cor:semi-orthogonal decomp for tame Shv(kot(G))} holds with $\shv^{\tame}$ replaced by $\shv^{\widehat\unip}$ and $\rshv^{\unip}$.
\end{proposition}

\subsubsection{Relation with affine Deligne-Lusztig inductions}
Our next goal is to give another characterization of $\shv^{\tame}(\kot_G)$ as well as $\shv^{\widehat\unip}(\kot_G)$ and $\rshv^{\unip}(\kot_G)$, which among other things will imply \Cref{prop: can duality of tame category} and \Cref{prop: can duality of unip category}.

\begin{lemma}\label{lem:generator-of-tame-category}
Let $\bfC\subset\shv(\kot_G)$ (resp. $\bfC^{\hat{u}}\subset\shv(\kot_G)$) be the presentable $\La$-linear category generated by the essential image of $\Ch_{LG,\phi}^{\mon}$ (resp. $\Ch_{LG,\phi}^{\hat{u}\mbox{-}\mon}$).
Then $\bfC$ (resp. $\shv^{\widehat\unip}(\kot_G)$) is generated as presentable $\La$-linear category by either 
\begin{itemize}
\item by objects $\{\widetilde{R}^*_{w}\}$ (resp. $\{R^*_{w,\hat{u}}\}$) for which $w$ is of minimal length in its $\sigma$-conjugacy class in $\widetilde{W}$; or
\item by objects $\{\widetilde{R}^!_{w}\}$ (resp. $\{R^!_{w,\hat{u}}\}$) for which $w$ is as above.
\end{itemize}
Let $\bfC^u\subset \fgshv(\kot_G)$ be the idempotent complete subcategory generated by $\Ch_{LG,\phi}^{\unip}(\fgshv(\iw\backslash LG/\iw))$.
Then $\fgshv^{\unip}(\kot_G)$ is generated as an  idempotent complete $\La$-linear category by either 
\begin{itemize}
\item by objects $\{R^*_w\}$ for which $w$ is of minimal length in its $\sigma$-conjugacy class in $\widetilde{W}$; or
\item by objects $\{R^!_w\}$ for which $w$ is of minimal length in its $\sigma$-conjugacy class in $\widetilde{W}$.
\end{itemize}
\end{lemma}
\begin{proof}
We prove the first statement. The other statements can be proved similarly.
As $\shv_{\mon}(\iw^u\backslash LG/\iw^u)$ is generated by $\{\widetilde\nabla_{\dot{w}}^{\mon}\}_{w\in\widetilde{W}}$, we see that $\bfC$ is generated by $\{\widetilde{R}_w^*\}_{w\in\widetilde{W}}$.
By \Cref{lem: ADL reduction method}, we see that $\bfC$ is generated by those $\widetilde{R}^!_w$ such that for every simple reflection $s$, $\ell(sw\sigma(s))\geq \ell(w)$. By \cite[Theorem A]{He.Nie.minimal.length}, such $w$ are exactly the minimal length elements in the corresponding $\sigma$-conjugacy classes.
\end{proof}

\begin{proposition}\label{prop: second characterization of tame LLCategory}
\begin{enumerate}
\item The category $\shv^{\tame}(\kot_G)$ is generated (as a presentable $\La$-linear category) by the essential image of $\Ch^{\mon}_{LG,\phi}: \shv_{\mon}(\iw^u\backslash LG/\iw^u)\to \shv(\kot_G)$. 
\item The category $\shv^{\widehat\unip}(\kot_G)$ is generated (as a presentable $\La$-linear category) by the essential image of $\Ch^{\hat{u}\mbox{-}\mon}_{LG,\phi}: \shv_{u\mbox{-}\mon}(\iw^u\backslash LG/\iw^u)\to \shv(\kot_G)$. 
\item The category $\rshv^{\unip}(\kot_G)$ is generated (as a presentable $\La$-linear category) by  the essential image of $\Ch^{\unip}_{LG,\phi}: \rshv(\iw\backslash LG/\iw)\to \rshv^{\unip}(\kot_G)$. 
\end{enumerate}
\end{proposition}
\begin{proof}
We only prove the first statement as the others are similar.
Let $\bfC \subseteq \shv(\kot_G)$ be as in \Cref{lem:generator-of-tame-category}. We need to show that $\bfC=\shv^{\tame}(\kot_G)$.
By \Cref{lem:generator-of-tame-category}, we see that $\bfC$ is generated by $\{\widetilde{R}_w^*\}$, with $w$ range over minimal length elements in their $\sigma$-conjugacy class.

We first show that $\bfC\subset \shv^{\tame}(\kot_G,\La)$.
It is enough to show that $\widetilde{R}^{*}_{w}\in \bfC$ with $w$ of minimal length in its $\sigma$-conjugacy class. As $\widetilde{R}^{*}_{w}\cong i_{b,*}(\cind_{P_b}^{G_b(F)}\widetilde{R}^{f,*}_{\dot{u}})[-\langle 2\rho,\nu_b\rangle]$
by \Cref{lem:ADLV-sheaf-minimal-length}, we see that $(i_{b'})^!R_w^* \simeq 0$ for every $b'\neq b$ and
$(i_b)^!R_w^* \simeq \cind_{P_b}^{G_b(F)}\widetilde{R}^{f,*}_{\dot{u}}\in \rep^\tame(G_b(F))$, as desired.

On the other hand, let $\mF\in \shv^{\tame}(\kot_G)$, i.e., such that $i_{b}^{!}\mF\in \rep^\tame(G_b(F),\La)$ for all $b\in B(G)$. We need to show that $\mF\in \bfC$. By \Cref{prop:semi-orthogonal-explicit-fiber-sequences} \eqref{prop:semi-orthogonal-explicit-fiber-sequences-3}, the object $\mF$ is obtained by repeated extensions of $i_{b,*}(i_b^!\mF)$, so we only need to show that for every $b\in B(G)$, every parahoric subgroup $P_b\subset G_b(F)$ with $L_{P_b}$ its Levi quotient, and every $\pi\in \rep(L_{P_b})$,
the object $(i_{b})_*(\cind_{P_b}^{G_{b}(F)}\pi)$ belongs to $\bfC$. We let $P_b=\breve\mP(\breve\mO)^{\sigma_b}$ and $L_{P_b}=L_{\breve\mP}(k)^{\sigma_b}$, where $\breve\mP$ is a parahoric of $G$ over $\breve\mO$ and $\sigma_b$ is the Frobenius structure determined by $b$.
Now by \Cref{lem: transitivity DL induction}, \Cref{thm: DL completeness}, and by applying \Cref{lem:generator-of-tame-category} to the finite case, it is enough to show that
\begin{equation}\label{eq:concrete shape of Rw}
(i_{b})_*(\cind_{P_b}^{G_b(F)} \widetilde{R}_{\dot{u}}^{f,*})\in \bfC.
\end{equation}
Here $\widetilde{R}_{\dot{u}}^{f,*}$ is the Deligne-Lusztig representation of $L_{P_b}$ associated to a minimal length element $u$ in an elliptic $\sigma_b$-conjugacy class of $W_{\breve\mP}$.

Let $uw$ be a minimal length element in its $\sigma$-conjugacy class $C$ as in \Cref{thm: reduction to min length elements} \eqref{thm: reduction to min length elements-1}.
Let $(M, x, \breve\bff_M, \frakc)$ be the standard quadruple constructed from $uw$ (see the end of \Cref{sec:sigma-straight-element}). We recall that $M=Z_G(\nu_w)$, $x=y^{-1}wy$, where $y\in W_0$ is the unique element of minimal length in $yW_M$ such that $y\nu_w=\tilde\nu_w$.
We have $\breve\bff\subset\overline{\breve\bfa}$ as in \Cref{thm: reduction to min length elements} \eqref{thm: reduction to min length elements-1}, and is minimal among such facets. Then $\breve\bff_M=y^{-1}\breve\bff$. Finally, $\frakc$ is the $\Ad_x\sigma$-conjugacy class containing $y^{-1}uy$.

We write $M(\breve F)^{\dot{x}\sigma}\supset \breve\mP_{\breve\bff_M}(\breve\mO)^{\dot{x}\sigma}$ by $G_b(F)\supset P_b(F)$.
Then $\widetilde{R}_{\dot{u}\dot{w}}^*$ is of the form as in \eqref{eq:concrete shape of Rw}. As proved in \cite[\textsection{1.8.3}]{he2016hecke}, and reviewed in \Cref{sec:sigma-straight-element}, every standard quadruple arises in this way. In addition, by \Cref{rem: u min in finite weyl}, every element in $\frakc_{\min}$ is of the form $y^{-1}u'y$ for some $u'\in W_{\breve\bff}$ such that $u'w\in C_{\min}$. It follows that every object of the form \eqref{eq:concrete shape of Rw} is isomorphic to $\widetilde{R}^*_{\dot{v}}$ for some
$v$ minimal length in its $\sigma$-conjugacy class, and therefore belongs to $\bfC$.
\end{proof}

\begin{proof}[Proof of \Cref{prop: can duality of tame category} and \Cref{prop: can duality of unip category}]
As before, we factor $\Nt^u$ as
\[
\frac{LG}{\Ad_\sigma \iw^u}\xrightarrow{\av_s} \frac{LG}{\Ad_\sigma \iw}=\Sht^\loc \xrightarrow{\Nt}\kot_G.
\]
We consider the canonical duality $\verd^{\can}_{\Sht^\loc}$ and $\verd^{\can}_{\kot_G}$. Note that by (the proof of) \Cref{lem: formula of monodromic ADL induction}, we see that
\[
(\verd^{\can}_{\Sht^\loc})^\cpt ((\av_s)_*(\delta^u)^!\widetilde\Delta^{\mon}_{\dot{w}})\cong (\av_s)_*(\delta^u)^!\widetilde\nabla^{\mon}_{\dot{w}}.
\]
It then follows from \Cref{cohomological.duality.kottwitz} that $(\verd^{\can}_{\kot_G})^\cpt(\widetilde{R}^*_{\dot{w}})\cong \widetilde{R}^!_{\dot{w}}$. Therefore, by \Cref{prop: second characterization of tame LLCategory}, we see that $\shv^{\tame}(\kot_G)^\cpt$ is stable under the canonical duality $(\verd_{\kot_G}^\can)^\cpt$.

Clearly the duality further restricts to a duality of $\shv^{\widehat\unip}(\kot_G)$. Finally, arguing similarly using \Cref{rem: can duality on f.g. sheaves} instead of  \Cref{cohomological.duality.kottwitz}, we see that $(\verd_{\kot_G}^{\can})^{\fg}$ preserves $\fgshv^{\unip}(\kot_G)$.
\end{proof}

Now assume that $\La$ is an algebraically closed field. Recall that for finite groups of Lie type, we have a decomposition \eqref{eq: Lusztig series finite group Lie type} of its category of representations.
Here is the analogue in the affine settings. Let $\zeta$ be a tame inertia type, which by  \Cref{prop: DL geo conj vs tame inertia type} is bijective to the set of geometric conjugacy classes of $(w,\theta)$.

We let $\shv^{\hat\zeta}(\kot_G)\subset \shv(\kot_G)$ be the full subcategory generated (as presentable $\La$-linear category) by $R^*_{\dot{w},\theta}$ for $(w,\theta)$ belonging to the geometric conjugacy classes corresponding to $\zeta$. 
Then we have
\begin{equation}\label{eq: decomposition of shv via tame inertia type}
\shv^{\tame}(\kot_G)=\bigoplus_{\zeta} \shv^{\hat\zeta}(\kot_G),
\end{equation}
where the direct sum ranges over all tame inertia types.

\subsubsection{Tame and unipotent Langlands category as a categorical trace}\label{SS: 2nd approach unipotent Langlands category}
Now we turn to another approach to $\shv^{\tame}(\kot_G)$, $\shv^{\widehat\unip}(\kot_G)$ and $\fgshv^{\unip}(\kot_G)$.

\begin{theorem}\label{prop: tame llc as a categorical trace}
The (monodromic) affine Deligne-Lusztig induction $\Ch_{LG,\phi}^{\mon}$ (see \eqref{eq: monodromic ADLI}) induces an equivalence
\begin{equation}\label{eq: tame trace identification}
\tr(\shv_{\mon}(\iw^u\backslash LG/\iw^u), \phi)\cong \shv^{\tame}(\kot_G),
\end{equation}
The unipotent affine Deligne-Lusztig induction $\Ch_{LG,\phi}^{\unip}$ (see \eqref{eq: unipotent ADLI}) induces an equivalence
\begin{equation}\label{eq: unipotent trace identification}
\tr(\rshv(\iw\backslash LG/\iw),\phi)\cong \rshv^{\unip}(\kot_G,\La).
\end{equation}

The above equivalences restrict to equivalences
\begin{equation}\label{eq: unipotent trace identification-2}
\tr(\shv\bigl((\iw,\hat{u})\backslash LG/(\iw, \hat{u}\bigr), \phi)\cong \shv^{\widehat{\unip}}(\kot_G,\La)\cong \tr(\shv(\iw\backslash LG/\iw),\phi).
\end{equation}
\end{theorem}
\begin{proof}
Just as in \Cref{thm: DL theory as categorical trace}, fully faithfulness follows from follows from \Cref{bimodule.geom.trace.fully.faithfulness.convolution.cat} and \Cref{prop-comparison-usual-trace-geo-trac}, applied to the sheaf theory as in  \Cref{prop: sheaf theory for monodromic sheaf} and \Cref{rmk: sheaf theory for monodromic sheaf}. 
Here we let $X=\bB \iw^u$, equipped with the natural action of $\mS_k$, and let $Y=\bB LG$, equipped with the trivial $\mS_k$-action.  
The essential surjectivity follows \Cref{prop: second characterization of tame LLCategory}.

The unipotent case follows by the same argument as in  \Cref{thm: DL theory as categorical trace} as well, applying \Cref{thm:trace-constructible} to $Y=\bB LG$ and $X=\bB \iw$, which is very placid (see \Cref{ex-coh-prosmooth}) and is weakly coh. pro-smooth over $k$ (see \Cref{torsor.descent.gives.placid.stack}).
Here we also use \Cref{prop-comparison-usual-trace-geo-trac} (which is applicable thanks to \Cref{prop: categorical property of affine Hecke}). 

Finally, the equivalence \eqref{eq: unipotent trace identification-2} can be proved by the same argument as in  \Cref{thm: DL theory as categorical trace} as well, taking \Cref{rem: Shv on Kot as colocalization} into account.
\end{proof}

\begin{remark}
We note that we may replace $\phi$ by any other automorphism of $LG$ preserving $\iw$ in \Cref{prop: tame llc as a categorical trace}. The case $\phi=\id$ will be studied in details in \cite{HHZ}.
 However, for the version $\rshv$, we do need $\phi$ to be Frobenius in order to apply \Cref{prop: faihtfully renornalized version of isoc}, see also \Cref{rem: faihtfully renornalized version of isoc}. We also recall that when $\La=\overline\bQ_\ell$, we have $\rshv(\kot_G,\overline\bQ_\ell)=\shv(\kot_G,\overline\bQ_\ell)$ by \Cref{cor: char zero fg=cpt kot}.
\end{remark}

\begin{proposition}\label{prop: automorphic side identifying duality}
Under the canonical equivalence \eqref{eq: tame trace identification}, the self-duality of $\tr(\shv_{\mon}(\iw^u\backslash LG/\iw^u),\phi)$ induced by the one on $\shv_{\mon}(\iw^u\backslash LG/\iw^u)$ is canonically identified with the canonical self-duality of $\shv^{\tame}(\kot_G)$ from \Cref{prop: can duality of tame category}.
\end{proposition}

We have the following affine analogue of \Cref{prop: fully faithfulness of tensor over finite Hecke category}, with the same proof.
\begin{proposition}\label{prop: fully faithfulness of tensor over Hecke category}
Let $\breve\mG_i, \ i=1,2$ be two affine smooth integral model of $G$ over $\breve\mO$. 
Let $\widetilde{L^+\breve\mG_i}\to L^+\breve\mG_i$ be as in \eqref{eq:f in LG}.
Then we have a fully faithful embedding
\begin{equation}\label{prop: fully faithfulness of tensor over Hecke category-1}
\shv_{\mon}(\widetilde{L^+\breve\mG_1}\backslash LG/\iw^u)\otimes_{\shv_{\mon}(\iw^u\backslash LG/\iw^u)}\shv_{\mon}(\iw^u\backslash LG/\widetilde{L^+\breve\mG_2})\to \shv(\widetilde{L^+\breve\mG_1}\backslash LG/\widetilde{L^+\breve\mG_2}),
\end{equation}
If one of $\breve\mG_i$  is a standard parahoric group schemes of $G$ (over $\breve\mO$) and $\widetilde{L^+\breve\mG_i}= L^+\breve\mG_i$,
then the above functor is an equivalence.

Similarly, we have a fully faithful embedding
\begin{equation}\label{prop: fully faithfulness of tensor over Hecke category-2}
\rshv(\widetilde{L^+\breve\mG_1}\backslash LG/\iw)\otimes_{\rshv(\iw\backslash LG/\iw)}\rshv(\iw\backslash LG/\widetilde{L^+\breve\mG_2})\to \rshv(\widetilde{L^+\breve\mG_1}\backslash LG/\widetilde{L^+\breve\mG_2}),
\end{equation}
which restricts to a fully faithful embedding
\begin{equation}\label{prop: fully faithfulness of tensor over Hecke category-3}
\shv(\widetilde{L^+\breve\mG_1}\backslash LG/\iw)\otimes_{\shv(\iw\backslash LG/\iw)}\shv(\iw\backslash LG/\widetilde{L^+\breve\mG_2})\to \shv(\widetilde{L^+\breve\mG_1}\backslash LG/\widetilde{L^+\breve\mG_2}),
\end{equation}
and if one of $\breve\mG_i$  is a standard parahoric group schemes of $G$ (over $\breve\mO$) and $\widetilde{L^+\breve\mG_i}= L^+\breve\mG_i$, then \eqref{prop: fully faithfulness of tensor over Hecke category-3} is an equivalence.
\end{proposition}

\begin{remark}\label{rem: constructible sheaf BB to BG-2}
As in \Cref{rem: constructible sheaf BB to BG}, we do not know whether \eqref{prop: fully faithfulness of tensor over Hecke category-2} is an equivalence when one of $\breve\mG_i$ is a standard parahoric and $\widetilde{L^+\breve\mG_i}= L^+\breve\mG_i$.
\end{remark}

For later applications, we need to understand where certain objects go under the functors. Note that for $\mF\in\shv_{\mon}(\iw^u\backslash LG/\iw^u)$, the object $[\mF]_{\phi}\in \tr(\shv_{\mon}(\iw^u\backslash LG/\iw^u),\phi)$ (see \eqref{eq:universal-trace} and \Cref{E:twisted-trace} for the notation) is identified with $\Ch^{\mon}_{LG,\phi}(\mF)$, so for simplicity we will always use the latter notion if possible.
We refer to \Cref{SS: Aff DL induction}, in particular \Cref{lem: ADLI for convolution sheaves} for descriptions of $\Ch^{\mon}_{LG,\phi}(\mF)$ for certain objects $\mF$.

On the other hand, recall that if $\bfM$ is a (left) dualizable $\shv_{\mon}(\iw^u\backslash LG/\iw^u)$-module, equipped with a left module functor $\phi_{\bfM}: \bfM\to {}^\phi\bfM$, then the map \eqref{eq:class-of-modules}
defines an object 
\[
[\bfM, \phi_\bfM]_{\phi}\in \tr(\shv_{\mon}(\iw^u\backslash LG/\iw^u),\phi)=\shv^{\tame}(\kot_G).
\]
By abuse of notations, we will denote $[\bfM, \phi_\bfM]_{\phi}$ by $\Ch^{\mon}_{LG,\phi}(\bfM,\phi_\bfM)$, although this is not really the monodromic affine Deligne-Lusztig induction of a sheaf.

Similarly, if $\bfM$ is a left $\rshv(\iw\backslash LG/\iw)$-module,  equipped with a left module functor $\phi_{\bfM}: \bfM\to {}^\phi\bfM$, then we write $[\bfM, \phi_\bfM]_{\phi}$ by $\Ch^{\unip}_{LG,\phi}(\bfM,\phi_\bfM)$.

The module category $\bfM$ we will consider arises from the geometry as follows. Let $\breve\mG$ be an affine smooth integral model of $G$ over $\breve\mO$, and let 
\[
W=\bB\widetilde{L^+\breve\mG}\to \bB L^+\breve\mG\to \bB LG
\]
be as in \eqref{eq:f in LG}. We equip $W$ with the trivial $\mS_k$-action. 
Then $\shv_{\mon}(\iw^u\backslash LG/\widetilde{L^+\breve\mG})$ is a left $\shv_{\mon}(\iw^u\backslash LG/\iw^u)$-module by the convolution pattern.

The following statement is proved in a similar way as in \Cref{prop: categorical property of affine Hecke} and  \Cref{prop: categorical property of monodromic Hecke}.
\begin{lemma}\label{lem: exteriori tensor product module category} 
Suppose either $L^+\breve\mG=\iw^u$, or is a standard parahoric subgroup.
\begin{enumerate}
\item For every prestack $X$ with a torus action, the exterior tensor product 
\[
\shv_{\mon}(\iw^u\backslash LG/\widetilde{L^+\breve\mG})\otimes_\La \shv_{\mon}(X)\to \shv_{\mon}(\iw^u\backslash LG/\widetilde{L^+\breve\mG}\times X)
\] 
is an equivalence.
\item For every prestack $X$, the exterior tensor product 
\[
\shv(\iw\backslash LG/\widetilde{L^+\breve\mG})\otimes_\La \shv(X)\to \shv(\iw\backslash LG/\widetilde{L^+\breve\mG}\times X)
\] 
is an equivalence.
\item For every quasi-compact ind-placid stack $X$, the exterior tensor product 
\[
\shv(\iw\backslash LG/\widetilde{L^+\breve\mG})\otimes_\La \rshv(X)\to \rshv(\iw\backslash LG/\widetilde{L^+\breve\mG}\times X)
\] 
is an equivalence.
\end{enumerate}
\end{lemma}

Suppose that $\sharp W_0$ is invertible in $\La$.
Let $w$ be a length zero element in $\widetilde{W}$, determining $b\in B(G)$.
Let $\breve\mP=\breve\mP_{\breve\bff}$ be a standard parahoric group scheme of $G_{\breve F}$ over $\breve\mO$ such that the facet $\breve\bff$ is $w\sigma$-stable. Then $P=\breve\mP(\breve\mO)^{\dot{w}\sigma}\subset G_{b}(F)$ is a parahoric subgroup. 
Let $\delta_P=\cind_P^{G_b(F)}\La\in \fgrep(G_b(F))\subset \rep(G_b(F))$ be as before.

Next, consider the following Frobenius structure on $\iw^u\bs LG/L^+\breve\mP_{\breve\bff}$
\[
 \iw^u\bs LG/L^+\breve\mP_{\breve\bff}\xrightarrow{\sigma} \iw^u\bs LG/L^+\breve\mP_{\sigma(\breve\bff)}\stackrel{g\mapsto g\dot{w}}{\cong}  \iw^u\bs LG/L^+\breve\mP_{\breve\bff}.
\]
We denote this Frobenius structure by $\sigma_{\dot{w}}$. Then the $*$-pushforward along $\sigma_{\dot{w}}$ defines a morphism
\[
\phi: \shv( \iw^u\bs LG/L^+\breve\mP_{\breve\bff}))\to {}^\phi\shv( \iw^u\bs LG/L^+\breve\mP_{\breve\bff})),
\]
 as left $\shv_{\mon}(\iw^u\backslash LG/\iw^u)$-modules.

\begin{proposition}\label{prop-class-of-level} We have
\[
\Ch^{\mon}_{LG,\phi}(\shv_{\mon}(\iw^u\backslash LG/L^+\breve\mP_{\breve\bff}), \phi)\cong \delta_P\in \rep(G_b(F))\subset \shv(\kot_G).
\]
Similarly, regarding $\rshv(\iw\backslash LG/L^+\breve\mP_{\breve\bff})$ as a left $\rshv(\iw\backslash LG/\iw)$-module. Then
\[
\Ch_{LG,\phi}^{\unip}(\rshv(\iw\backslash LG/L^+\breve\mP_{\breve\bff}), \phi)\cong \delta_P\in \fgrep(G_b(F))\subset \fgshv(\kot_G).
\]

\end{proposition}
\begin{proof}
We will apply \Cref{ex: geo phi-trace class}. 
Here as before $\der=\shv_{\mon}$, with $X=\bB \iw^u$ equipped with the natural $\mS_k$-action and $Y=\bB LG$ with the trivial $\mS_k$-action. We let $W=\bB L^+\breve\mP_{\breve\bff}$ equipped with the trivial $\mS_k$-action. Both $X$ and $Y$ are defined over $k_F$, and they admit the $q$-Frobenius endomorphism $\phi_X$ and $\phi_Y$. The space $W$ is equipped with the Frobenius structure $\phi_W$ given by $\Ad_{\dot{w}}\sigma$. Let $h:W\to Y$ be the natural map.
We have a natural isomorphism $\Ad_{\dot{w}}: \phi_Y\circ h\to h\circ \phi_W$. Then
So \Cref{rem-for-assumption-alpha-alphageo}  is applicable.

By \Cref{prop: fully faithfulness of tensor over Hecke category}, \Cref{cor-geom-duality-left-module} \eqref{cor-geom-duality-left-module-1}  is applicable. It remains to notice that
$\mL_\phi(W)=L^+\breve\mP_{\breve\bff}/\Ad_{\dot{w}\sigma} L^+\breve\mP_{\breve\bff}=\bB_{\mathrm{profet}} P$, and $\mL_\phi(h)$ is the natural map 
\[
L^+\breve\mP_{\breve\bff}/\Ad_{\dot{w}\sigma} L^+\breve\mP_{\breve\bff}\to LG/\Ad_\sigma LG=\kot_G,\quad g\mapsto g\cdot \dot{w}.
\]

For the unipotent case, due to \Cref{rem: constructible sheaf BB to BG-2}, we cannot directly apply  \Cref{cor-geom-duality-left-module} \eqref{cor-geom-duality-left-module-1}. But by \Cref{prop: fully faithfulness of tensor over Hecke category} and \Cref{lem: exteriori tensor product module category}, we can apply  \Cref{cor-geom-duality-left-module} \eqref{cor-geom-duality-left-module-2}. 
Note that $\proj_{\trg}=\proj^{\unip}: \rshv(\kot_G)\to \rshv^{\unip}(\kot_G)$. Note that $\sharp W_\mP\mid \sharp W_0$, and therefore $\delta_P\subset \rep^{\widehat\unip}(G(F))\cap \fgrep(G(F))\subset \rshv^{\unip}(G(F))$ by \Cref{lem: induction of parabolic belonging unipotent-2}. We see that
$\proj_{\trg}(\delta_P)=\delta_P$.
\end{proof}


\subsection{Whittaker models}\label{SS: Whittaker model}
In this subsection, we will fix a non-zero additive character $\psi: F\to \La^\times$ such that $\psi(\mO_F\varpi)=1$ and such that $\overline\psi: k_F=\mO_F/\mO_F\varpi\to \La^\times$ is non-trivial. We will discuss Whittaker models of certain objects in $\shv^{\tame}(\kot_G)$.

\subsubsection{Whittaker models of tame and unipotent representations}
We start with Whittaker models of tame and unipotent representations. We refer to the beginning of the section regarding our notations and conventions related to $G$.

For a semi-standard parahoric group scheme $\breve\mP$ of $G$ over $\breve\mO$ (see \Cref{sec:sigma-straight-element} for the meaning), let $L^{++}\breve\mP\subset L^+\breve\mP$ denote the pro-unipotent radical and $L_{\breve\mP}=L^+\breve\mP/L^{++}\breve\mP$ the Levi quotient.  If $\breve\mP$ is defined over $\mO$, we will let $\mP$ denote its model over $\mO$, with $L_\mP$ its Levi quotient over $k_F$. We let $P=\mP(\mO)=L^+\mP(k_F)$, $L_P=L_\mP(k_F)$ and $P^u=L^{++}\mP(k_F)$.
Recall that if $\mP=\mI$, then $L^{++}\mI\subset L^{+}\mI$ are also denoted by $\iw^u\subset \iw$. We also write $\iw^u(k_F)=I^u\subset \iw(k_F)=I$.  

\begin{lemma}\label{lem:unipotent-in-Levi}
Let $\breve\mP$ be a semi-standard parahoric,
  and let $U_{L_{\breve\mP}}:=(LU\cap L^+\breve\mP)/(LU\cap L^{++}\breve\mP)\subset L_{\breve\mP}$. Then $B_{L_{\breve\mP}}:=\mS_k\cdot U_{L_{\breve\mP}}$ is a Borel subgroup of $L_{\breve\mP}$ (so $U_{L_{\breve\mP}}$ is the unipotent radical of $B_{L_{\breve\mP}}$). Similarly statement holds for $U'_{L_{\breve\mP}}:= (\iw^u\cap L^+\breve\mP)/(\iw^u\cap L^{++}\breve\mP)$.

\end{lemma}
\begin{proof}
   This is almost tautological after spreading out definitions.
   Let $\breve\mP'$ be the standard parahoric of the same type as $\breve\mP$.
   Let $w\in W_\af$ such that $\dot{w}L^+\breve\mP\dot{w}^{-1}=L^+\breve\mP'$ for one (and therefore any) lift of $w$ to $N_G(S)(\breve F)$. Then it is enough to show that
   $U_{L_{\breve\mP'},w}:=(\dot{w}LU\dot{w}^{-1}\cap L^+\breve\mP')/(\dot{w}LU\dot{w}^{-1}\cap L^{++}\breve\mP')\subset L_{\breve\mP'}$ is the unipotent radical of the Borel subgroup $\mS_kU_{L_{\breve\mP'},w}$.
   As before, let $\Phi$ be the relative root system of $(G_{\breve F},S_{\breve F})$ and
   $\Phi_\af$ the corresponding affine root system. Let $\Phi_{L_{\breve\mP'}}\subset \Phi_\af$ be the root system corresponding to $L_{\breve\mP'}$.
   Let $\Phi_\af\to\Phi$ be the map sending an affine root $\al$ to its vector part $\dot{\al}$. Then the composition $\pr:\Phi_{L_{\breve\mP'}}\subset \Phi_\af\to \Phi$ is injective and the image can be identified with the root system with respect to $(L_{\breve\mP'}, \mS_k)$.

    It follows that $U_{L_{\breve\mP'},w}$ is generated by the root subgroups of $L_{\breve\mP'}$ corresponding those roots in $\pr(\Phi_{L_{\breve\mP'}})\cap w(\Phi^+)$. But as $\mathrm{pr}(\Phi_{L_{\breve\mP'}})\cap w(\Phi^+)$ is the intersection of $\mathrm{pr}(\Phi_{L_{\breve\mP'}})$ with a half space of $\xcoch(\mS_k)_\bR$, the first claim of the lemma then is clear. For the second claim, notice that $(\iw^u\cap L^+\breve\mP)/(\iw^u\cap L^{++}\breve\mP)$ is generated by the root subgroups of $L_{L_{\breve\mP'}}$ corresponding those roots in $\mathrm{pr}(\Phi_{L_{\breve\mP'}}\cap w(\Phi_\af^+))$. But clearly, $\Phi_{L_{\breve\mP'}}\cap w(\Phi_\af^+)$ form a set of positive roots of $\Phi_{L_{\breve\mP'}}$.
\end{proof}

The additive character $\psi$ and the pinning $(B,T,e)$ together determine a Whittaker datum $(U,\psi_e)$ where
\[
\psi_e: U(F)\xrightarrow{e} F\xrightarrow{\psi} \La^\times. 
\]
Let
\[
\coWhit_{\psi_e}:=\cind_{U(F)}^{G(F)}\psi_e \in \rep(G(F),\La),
\]
be the Whittaker (co)module of $G(F)$. It belongs to the heart of $ \rep(G(F),\La)$.
Recall that it is not in $\fgrep(G(F),\La)^{\heartsuit}$ but can be approximated (i.e. can be written as a filtered colimit of) finitely generated $G(F)$-modules. 
For our purpose, it is enough to consider the first term of the approximation, given by the Iwahori-Whittaker module $\IW_{\psi_1}$ of $G(F)$. Namely, we have the direct product decomposition
\[
I^u=(I^u\cap U^-(F))\cdot (I^u\cap T(F))\cdot (I^u\cap U(F)).
\]
Then there is a unique character, 
\begin{equation}\label{eq-char-psi0}
\psi_1:I^u \to \La^\times
\end{equation}
such that $\psi_1(I^u\cap U^-(F))=\psi_1(I^u\cap T(F))=1$, and that
$\psi_1|_{I^u\cap U(F)}=\psi_e|_{I^u\cap U(F)}$.
Let
\[
\IW_{\psi_1}:=\cind_{I^u}^{G(F)}\psi_1\in \fgrep(G(F),\La)^{\heartsuit},
\]
which is usually called the Iwahori-Whittaker module of $G(F)$. 
Note that $\IW_{\psi_1}\in \rep^{\tame}(G(F),\La)^{\heartsuit}$ as it may be written as
\[
\IW_{\psi_1}=\cind_{K}^{G(F)}\GG_{\overline\psi_1},
\]
where $\GG_{\overline\psi_1}$ is the Gelfand-Graev representation of the finite group $G(k_F)$, defined by unique the character $\overline\psi_1: U(k_F)\to\La^\times$ such that $U(\mO_F)\to U(k_F)\xrightarrow{\overline\psi_1}\La^\times$ is the restriction $\psi_e|_{U(\mO_F)}$.

We have a natural map
\begin{equation}\label{eq-IW-coWhit}
\IW_{\psi_1}\to\coWhit_{\psi_e}
\end{equation}
given by the function $f$ on $G(F)$, supported on $I^u\cdot U(F)$ such that $f(1)=1$.


\begin{lemma}\label{lem-iwahori-whit-apr-cowhit}
  The map \eqref{eq-IW-coWhit} induces an isomorphism $\IW_{\psi_1}\cong \proj^{\tame}(\coWhit_{\psi_e})$.
\end{lemma}
\begin{proof}
Using \Cref{lem: generation by cuspidal}, it is enough to show that for every $(P,\pi)$, where $P=\mP(\mO)$ is a standard parahoric with Levi quotient $L_P=L_{\mP}(k_F)$ and $\pi$ is a \emph{cuspidal} representation of $L_P$, \eqref{eq-IW-coWhit} induces an isomorphism
\[
\Hom(\cind_P^{G(F)}\pi, \IW_{\psi_1})\to \Hom(\cind_P^{G(F)}\pi,\coWhit_{\psi_e}).
\]
We have the Bruhat and the Iwasawa decompositions
\[
G(F)=\bigsqcup_{w\in W_{P}^\sigma\backslash \widetilde{W}^\sigma}PwI^u=\bigsqcup_{w\in W_{P}^\sigma\backslash \widetilde{W}^\sigma}PwU(F),
\]
where $W_P\subset \widetilde W$ is the Weyl group corresponding to $P$, and $(-)^\sigma$ means taking Frobenius invariants. For each $w\in \widetilde W^{\sigma}$, by abuse of notation we write $\mP_w:={}^{w^{-1}}\mP$ (which precisely means $L^+\mP_w={}^{w^{-1}}L^+\mP$), which is a rational semi-standard parahoric. Write $P_w=L^+\mP_w(k_F)$.
Let $\pi_w$ be the representation of $L_{P_w}=L_{\mP_w}(k_F)$ obtained from $\pi$ by transport of structure.

Using notations as in \Cref{lem:unipotent-in-Levi}, we write $U_{L_{P_w}}$ for $U_{L_{\mP_w}}(k_F)$ and similarly $U'_{L_{P_w}}$ for $U'_{L_{\mP_w}}(k_F)$.
It follows from the Frobenius reciprocity law that
\[
\Hom(\cind_{P}^{G(F)}\pi,\IW_{\psi_1})=\bigoplus_{w\in W_P^{\sigma}\backslash\widetilde{W}^\sigma}\Hom_{L_{P_w}} (\pi_w, \ind_{U'_{L_{P_w}}}^{L_{P_w}}\overline{\psi_1}),
\]
where $\overline{\psi_1}= (\psi_1|_{I^u\cap P_w})^{I^u\cap P_w^u}$ is a character of $U'_{L_{P_w}}$. 
Similarly,
\[
\Hom(\cind_{P}^{G(F)}\pi,\coWhit_{\psi_e})=\bigoplus_{w\in W_P^{\sigma}\backslash\widetilde{W}^\sigma}\Hom_{L_{P_w}} (\pi_w, \ind_{U_{L_{P_w}}}^{L_{P_w}}\overline{\psi_e}),
\]
where $\overline{\psi_e}= (\psi_e|_{U(F)\cap P_w})^{U(F)\cap P_w^u}$ is the representation of $U_{L_{P_w}}$. 
We note that both $\Hom$ spaces concentrate in cohomological degree zero.

Note that the restrictions of $\overline{\psi_1}$ and $\overline{\psi_e}$ to $U_{L_{P_w}}\cap U'_{L_{P_w}}$ are the same (by definition of $\psi_1$). As $\pi_w$ is cuspidal, unless $\overline{\psi_1}$ (and therefore $\overline{\psi_e}$) is generic, both Hom space would be zero. But in the generic case, $U_{L_{P_w}}=U'_{L_{P_w}}$ and the map induced by \eqref{eq-IW-coWhit} is just the identity map. The lemma follows.
\end{proof}

We also let
\begin{equation}\label{eq: unipotent IW}
\IW_{\psi_1}^{\unip}:=\proj^{\unip}(\IW_{\psi_1})\in \rshv^{\unip}(G(F)).
\end{equation}

When $\La=\overline\bQ_\ell$, it admits an explicit description as follows.
Let  
\[
   M_{\asp}:=\Hom_{G(F)}(\delta_I,\IW_{\psi_1})\cong \IW_{\psi_1}^I\cong C_c(I\backslash G(F)/(I^u,\psi_1)).
\]
which is an $H_I$-module, usually called the anti-spherical module of $H_I$.
As $H_I$-modules, $M_{\asp}$ is a direct summand of $H_I$. Note that it follows from \Cref{lem-iwahori-whit-apr-cowhit} that 
\[
M_{\asp}=\Hom_{G(F)}(\delta_I,\coWhit_{\psi_e})=\coWhit_{\psi_e}^I.
\]

\begin{corollary}\label{cor-unipotent-part-of-iwahori-whit}
Suppose $\La=\overline\bQ_\ell$. Then we have $\IW_{\psi_1}^{\unip}\cong \delta_I\otimes_{H_I}M_{\asp}$.
\end{corollary}
\begin{proof}
We follow that same argument as above, but now assume that $\pi$ is an irreducible unipotent cuspidal representation of $L_P$. Then it follows from \cite[\textsection{10}]{Deligne.Lusztig} that if $\Hom_{L_{P_w}} (\pi_w, \ind_{U'_{L_{P_w}}}^{L_{P_w}}\overline{\psi_1})\neq 0$ only if $L_P$ is the torus and $\pi$ is trivial. Then $P=I$. The corollary follows. 
\end{proof}

\subsubsection{Iwahori-Whittaker representations as a trace}
Now we apply the discussions of Whittaker categories at the end of \Cref{SS: Aff Hecke cat}. We take $\breve\bff=v_0$ to be absolutely special vertex (determined by the pinning), and let $e_{\breve\bff}=e: \iw^u\to \bG_a$. Let $\widetilde{\iw^u}$ be the pullback of the Artin-Scheier cover of $\bG_a$. We have
\[
\shv_{\mon}(\iw^u\backslash LG/\widetilde{\iw^u})=\bigoplus_{a\in k_F} \shv_{\mon}(\iw^u\backslash LG/(\iw^u,\psi_a)),
\]
where $\psi_a(\cdot)=\overline\psi(a\cdot): k_F\to \La^\times$, inflated as a character of $I^u$ via $I^u\to k_F$. In particular, when $a=1$, $\psi_1$ coincides with the previously defined additive character of $I^u$. In this case,
 $\shv_{\mon}\bigl(\iw^u\backslash LG/(\iw^u,\psi_1)\bigr)$ is called the (monodromic) Iwahori-Whittaker category. On the other hand, note that if $a=0$, then $\psi_a$ is trivial so $\shv_{\mon}(\iw^u\backslash LG/(\iw^u,\psi_a))=\shv_{\mon}(\iw^u\backslash LG/\iw^u)$.

 There is also a unipotent version 
 \[
\rshv(\iw\backslash LG/\widetilde{\iw^u})=\shv(\iw\backslash LG/\widetilde{\iw^u})=\bigoplus_{a\in k_F} \shv(\iw\backslash LG/(\iw^u,\psi_a)),
\]
where the first equality follows from the last statement in \Cref{compact.generation.admissible.stacks}, since $\widetilde{\iw^u}$ is coh. pro-unipotent.

Note that the natural Frobenius endomorphism of $\iw^u\backslash LG/\widetilde{\iw^u}$ and on $\iw\backslash LG/\widetilde{\iw^u}$ preserves the above decompositions. 

\begin{proposition}\label{prop-class-of-iwahori-whit-module}
We have 
\[
\Ch_{LG,\phi}^{\mon}(\shv_{\mon}(\iw^u\backslash LG/(\iw^u,\psi_1)))\cong \IW_{\psi_1}\in \rep(G(F))\stackrel{(i_1)_*}{\hookrightarrow} \shv(\kot_G),
\] 
and similary
\[
\Ch_{LG,\phi}^{\unip}(\shv(\iw\backslash LG/(\iw^u,\psi_1)))\cong \IW_{\psi_1}^{\unip}.
\]
\end{proposition}
\begin{proof}
Thanks to \Cref{lem: exteriori tensor product module category}, assumption of  \Cref{cor-geom-duality-left-module} \eqref{cor-geom-duality-left-module-2} holds, giving a duality datum of $\shv^{\mon}(\iw^u\backslash LG/\widetilde{\iw^u})$ as a left $\shv_{\mon}(\iw^{u}\backslash LG/\iw^u)$-module which in term induces a duality datum of $\shv_{\mon}(\iw^u\backslash LG/(\iw^u,\psi_1))$. Explicitly, the counit is given by
\begin{eqnarray*}
& &\shv_{\mon}(\iw^u\backslash LG/(\iw^{u},\psi_1)) \otimes \shv^{\mon}((\iw^{u},\psi_1)\backslash LG/\iw^u) \\
& \rightarrow & \shv_{\mon}(\iw^u\backslash LG/\widetilde{\iw^{u}}) \otimes \shv_{\mon}(\widetilde{\iw^{u}}\backslash LG/\iw^u) \\
& \xrightarrow{\star^{\tilde{u}}} & \shv_{\mon}(\iw^u\backslash LG/\iw^u),
\end{eqnarray*}
while the unit is given by the image of 
\begin{equation}\label{eq-unit-bi-iwahori-whittaker}
\mathbf{1}_{(\iw^{u},\psi_1)}\in \shv\bigl((\iw^{u},\psi_1)\backslash LG/(\iw^{u},\psi_1)\bigr),
\end{equation}
under the composed functors
\begin{multline*}
\shv\bigl((\iw^{u},\psi_1)\backslash LG/(\iw^{u},\psi_1)\bigr)\to \shv(\widetilde{\iw^{u}}\backslash LG/\widetilde{\iw^{u}}) \\ 
\to  \shv_{\mon}(\widetilde{\iw^{u}}\backslash LG/\iw^u)\otimes_{\shv_{\mon}(\iw^u\backslash LG/\iw^u)}\shv_{\mon}(\iw^u\backslash LG/\widetilde{\iw^{u}}) 
\end{multline*}
where the second functor is the right adjoint of the natural one,
and $\mathbf{1}_{(\iw^{u},\psi_1)}$ is the direct summand $\mathbf{1}_{\widetilde{\iw^u}}$ according to the decomposition \eqref{eq: decomposition of endo whittaker}.

Then by \Cref{prop-geo-class-via-corr} (and \Cref{ex: geo phi-trace class}), $\Ch_{LG,\phi}^{\mon}(\shv_{\mon}(\iw^u\backslash LG/(\iw^u,\psi_1)),\phi)$ is given by the image of $\mathbf{1}_{(\iw^{u},\psi_1)}$ under the functors
\begin{multline*}
\shv((\iw^u,\psi_1)\backslash LG/(\iw^u,\psi_1))\to \shv(\widetilde{\iw^u}\backslash LG/\widetilde{\iw^u})\\
\xrightarrow{(\delta^{\tilde u})^!} \shv( LG/\Ad_\sigma \widetilde{\iw^u})\xrightarrow{(\Nt^{\tilde u})_*} \shv(\kot_G)\xrightarrow{\proj^{\tame}} \shv^{\tame}(\kot_G),
\end{multline*}

The same argument as in \Cref{prop:affine DL induction chi} shows that 
\[
(\Nt^{\tilde u})_*(\delta^{\tilde u})^!(\mathbf{1}_{(\iw^{u},\psi_1)})\cong \IW_{\psi_1}.
\]
Now the first statement follows from \Cref{lem-iwahori-whit-apr-cowhit}.

The second statement can be proved similarly. 
\end{proof}
\subsubsection{Iwahori-Whittaker coefficients}

The next goal of this subsection is to prove the following result, which can be regarded as a vast generalization of \Cref{prop:multiplicity one GG rep} in the affine case.

We recall from \Cref{SSS: mono Hecke cat} that we have the perverse $t$-structure on  $\shv_{\mon}(\iw^u\backslash LG_w/\iw^u)$ defined by the generalized constant sheaf whose $!$-pullback to  $LG_{w}/\iw^u$  is the usual constant sheaf $\La_{LG_{\leq w}/\iw^u}\in \cshv(LG_{w}/\iw^u)$. 

We call $\mZ\in\shv_{\mon}(\iw^u\backslash LG/\iw^u)$ a (monodromic) central sheaf if the following two properties hold:
\begin{itemize}
\item $\mZ$ is a perverse sheaf and $\mZ\star^u(-)$ is convolution exact;
\item there is an isomorphism of functors 
\[
\mZ\star^u(-)\simeq (-)\star^u\mZ: \shv_{\mon}(\iw^u\backslash LG/\iw^u)\to \shv_{\mon}(\iw^u\backslash LG/\iw^u).
\]
\end{itemize}

\begin{theorem}\label{prop:Whittaker-exactness}
Let $\mZ\in \shv_\mon(\iw^u\backslash LG/\iw^u)$ be a monodromic central sheaf. Suppose  $\mZ\star^u \widetilde{\Delta}^{\mon,\psi}_{\dot{w}_0}$
is a cofree tilting object in $\shv_{\mon}(\iw^u\backslash LG/(\iw^u,\psi_1))$. I.e. it
admits a filtration by $\{\widetilde{\Delta}^{\mon,\psi}_{\dot{w}}\}_w$ as well as a filtration by $\{\widetilde{\nabla}^{\mon,\psi}_{\dot{w}}\}_w$.

Let $\mF\in \shv_{\mon}(\iw^u\backslash LG/\iw^u)$ be a cofree monodromic tilting sheaf. Then
\[
\Hom_{\shv(\kot_G)}(\Ch^{\mon}_{LG,\phi}(\mZ\star^u\mF), (i_{1})_*\IW_{\psi_1})\in\Mod_\La^{\heartsuit}.
\] 
\end{theorem}
We expect that $\mZ\star^u \widetilde{\Delta}^{\mon,\psi}_{\dot{w}_0}$
is a always a cofree tilting object in $\shv_{\mon}(\iw^u\backslash LG/(\iw^u,\psi_1))$. But we have not checked this.
\begin{proof}
By \Cref{lem-mon-spectral-sequence-hom}, there is a filtration of $\Hom_{\shv(\kot_G)}(\Ch^{\mon}_{LG,\phi}(\mZ\star^u\mF), (i_{1})_*\IW_{\psi_1})$ with associated graded being
\[
\Hom_{\shv_{\mon}(\iw^u\backslash LG/(\iw^u,\psi_1))}((\av_s)^*(\av_s)_*(\mZ\star^u \mF)\star^u \widetilde{\Delta}^{\mon,\psi}_{\sigma(\dot{w})}, \widetilde{\nabla}^{\mon,\psi}_{\dot{w}}), 
\]
which by \Cref{lem: convolving central sheaf commutes with killing monodromy} below is isomorphic to
\[
\Hom_{\shv_{\mon}(\iw^u\backslash LG/(\iw^u,\psi))}(\mZ\star^u (\av_s)^*(\av_s)_*(\mF)\star^u \widetilde{\Delta}^{\mon,\psi}_{\sigma(\dot{w})}, \widetilde{\nabla}^{\mon,\psi}_{\dot{w}}).
\]

We will show that the $i$th cohomology of the above complex vanishes unless $i=0$.
In fact, we will show that for every $v_1,v_2\in \widetilde{W}$, 
\[
H^i\Hom_{\shv_{\mon}(\iw\backslash LG/(\iw^u,\psi))}(\mZ\star^u (\av_s)^*(\av_s)_*(\mF) \star^u \widetilde{\Delta}^{\mon,\psi}_{\dot{v}_1}, \widetilde{\nabla}^{\mon,\psi}_{\dot{v}_2})=0,\quad \forall i\neq 0.
\] 
First, using \Cref{lem: convolve (co)standard with Whit}, we may write
\[
\mZ\star^u (\av_s)^*(\av_s)_*(\mF) \star^u  \widetilde{\Delta}^{\mon,\psi}_{\dot{v}_1}\cong \mZ\star^u (\av_s)^*(\av_s)_*(\mF) \star^u  \widetilde{\Delta}^{\mon}_{\dot{v}_1\dot{w}_0^{-1}}\star^u \widetilde{\Delta}^{\mon,\psi}_{\dot{w}_0}.
\] 
Since $\mF$ is a cofree monodromic tilting sheaf, by \Cref{cor: av upper star av lower star}
$(\av_s)^*(\av_s)_*(\mF)$ admits a filtration with associated graded being $\{\nabla^{\mon}_{\dot{v}}(\Ch(\cohdual_{\chi_{\varphi_{\bar{v}}}}))\}_{v\in\widetilde{W}}$. Then by a slight variant of \Cref{lem: exactness of convolution with tilting}, we see that $(\av_s)^*(\av_s)_*(\mF) \star^u  \widetilde{\Delta}^{\mon}_{\dot{v}_1\dot{w}_0^{-1}}$ is perverse. Then $ \mZ\star^u (\av_s)^*(\av_s)_*(\mF) \star^u  \widetilde{\Delta}^{\mon}_{\dot{v}_1\dot{w}_0^{-1}}\star^u \widetilde{\Delta}^{\mon,\psi}_{\dot{w}_0}$ is also perverse. Therefore, we have
\[
\Hom_{\shv_{\mon}(\iw\backslash LG/(\iw^u,\psi))}(\mZ\star^u (\av_s)^*(\av_s)_*(\mF) \star^u  \widetilde{\Delta}^{\mon,\psi}_{\dot{v}_1}, \widetilde{\nabla}^{\mon,\psi}_{\dot{v}_2})\in \Mod_\La^{\geq 0}.
\] 
On the other hand, by the centrality of $\mZ$, we have
\[
\mZ\star^u (\av_s)^*(\av_s)_*(\mF) \star^u  \widetilde{\Delta}^{\mon}_{\dot{v}_1\dot{w}_0^{-1}}\star^u \widetilde{\Delta}^{\mon,\psi}_{\dot{w}_0}\cong  (\av_s)^*(\av_s)_*(\mF) \star^u  \widetilde{\Delta}^{\mon}_{\dot{v}_1\dot{w}_0^{-1}}\star^u \mZ\star^u \widetilde{\Delta}^{\mon,\psi}_{\dot{w}_0}.
\] 
By our assumption, $\mZ\star^u \widetilde{\Delta}^{\mon,\psi}_{\dot{w}_0}$ admits a filtration by perverse sheaves with associated graded being as the form $\{\widetilde{\Delta}^{\mon,\psi}_{\dot{v}}\}_{v\in \widetilde{W}}$. On the other hand, $(\av_s)^*(\av_s)_*(\mF)$ admits a filtration with associated graded being $\{\Delta^{\mon}_{\dot{v}}(\Ch(\cohdual_{\chi_{\varphi_{\bar{v}}}}))\}_{v\in\widetilde{W}}$. It follows that
$\mZ\star^u (\av_s)^*(\av_s)_*(\mF) \star^u  \widetilde{\Delta}^{\mon,\psi}_{\dot{v}_1}$ has a filtration with associated graded being 
\[
\Delta^{\mon}_{\dot{v}}(\Ch(\cohdual_{\chi_{\varphi_{\bar{v}}}}))\star^u \widetilde{\Delta}^{\mon}_{\dot{v}_1\dot{w}_0^{-1}}\star^u  \widetilde{\Delta}^{\mon,\psi}_{\dot{v}'}.
\] 
It then follows from \Cref{cor: hom from standard to costandard} that
\[
\Hom_{\shv_{\mon}(\iw\backslash LG/(\iw^u,\psi))}(\mZ\star^u (\av_s)^*(\av_s)_*(\mF) \star^u  \widetilde{\Delta}^{\mon,\psi}_{\dot{v}_1}, \widetilde{\nabla}^{\mon,\psi}_{\dot{v}_2})\in \Mod_\La^{\leq 0}.
\] 
The desired vanishing follows.
\end{proof}

\begin{lemma}\label{lem: convolving central sheaf commutes with killing monodromy}
Let $\mZ\in\shv_{\mon}(\iw^u\backslash LG/\iw^u)$ be a monodromic central sheaf. Then for every $\mF\in\shv_{\mon}(\iw^u\backslash LG/\iw^u)$, we have
\[
(\av_s)^*(\av_s)_*(\mZ\star^u\mF)\cong \mZ\star^u(\av_s)^*(\av_s)_*(\mF).
\]
\end{lemma}
\begin{proof}
By base change, the functor
\[
(\av_s)^*(\av_s)_*= a_*(\La_{\mS_k}\boxtimes -): \shv_{\mon}(\iw^u\backslash LG/\iw^u)\to \shv_{\mon}(\iw^u\backslash LG/\iw^u)
\]
Here, $a: S_k\times \iw^u\backslash LG/\iw^u\to  \iw^u\backslash LG/\iw^u$ is the $\sigma$-conjugation action. We rewrite it as the composition $f: \mS_k\xrightarrow{t\mapsto (t, \sigma(t)^{-1})}\mS_k\times \mS_k$ followed by the action $m_{lr}: \mS_k\times \mS_k\times \iw^u\backslash LG/\iw^u\to \iw^u\backslash LG/\iw^u$ by left and right multiplication. Therefore, 
\[
a_*(\La_{\mS_k}\boxtimes-)\cong(m_{lr})_*(f_*\La_{\mS_k}\boxtimes-)\cong (m_{lr})_*(\av^{\mon}(f_*\La_{\mS_k})\boxtimes-),
\]
where the last isomorphism follows from \Cref{lem: convolve with sheaves on torus}.

On the other hand, recall we have $\Delta_e^{\mon}: \shv_{\mon}(\mS_k)\cong \shv_{\mon}(\iw^u\backslash \iw/\iw^u)\subset \shv_{\mon}(\iw^u\backslash LG/\iw^u)$. Then under this identification, for $\mG_1\boxtimes \mG_2\in \shv_{\mon}(\mS_k)\otimes_\La \shv_{\mon}(\mS_k)$, we have
\[
(m_{lr})_*(\mG_1\boxtimes\mG_2\boxtimes -)\cong \mG_1\star^u(-)\star^u\mG_2: \shv_{\mon}(\iw^u\backslash LG/\iw^u)\to \shv_{\mon}(\iw^u\backslash LG/\iw^u)
\]
In particular, by centrality of $\mZ$, we have
\[
(m_{lr})_*(\mG_1\boxtimes\mG_2\boxtimes (\mZ\star^u\mF))\cong \mZ\star^u((m_{lr})_*(\mG_1\boxtimes\mG_2\boxtimes\mF)).
\] 
As $\shv_{\mon}(\mS_k)\otimes_\La \shv_{\mon}(\mS_k)\cong \shv_{\mon}(\mS_k\times \mS_k)$ (by \Cref{lem: tensor product tensor of monodromic sheaves}), we see that the above isomorphism still holds if we replace $\mG_1\boxtimes\mG_2$ by any $\mG\in \shv_{\mon}(\mS_k\times\mS_k)$, in particular by $\av^{\mon}(f_*\La_{\mS_k})$. Putting everything together gives the lemma.
\end{proof}

\begin{remark}\label{rem:Whittaker-exactness-unip}
There is also a unipotent version of the above theorem. Let $\mZ\in \fgshv(\iw\bs LG/\iw)$ be a central sheaf, by which we a perverse sheaf, which is convolution $t$-exact and such that $\mZ\star \mF\simeq \mF\star \mZ$ for any $\mF\in \rshv(\iw\bs LG/\iw)$. Then 
Note that in this case, it is known that $\mZ\star \Delta^{\psi}_{w_0}$ is a tilting object in $\shv(\iw\bs LG/(\iw^u,\psi_1))$. Then the same argument gives
\[
\Hom_{\rshv(\kot_G)}(\Ch^{\unip}_{LG,\phi}(\mZ), (i_{1})_*\IW^{\unip}_{\psi_1})\in\Mod_\La^{\heartsuit}.
\] 
(Note that there are no tilting object in $\rshv(\iw\bs LG/\iw)$ except $\Delta_e$.)
\end{remark}

\newpage

\section{Tame and unipotent categorical local Langlands correspondence}\label{sec:unip-cat-langlands}
In this section, we put everything together to prove our main theorem. The extra input is the tame local geometric Langlands correspondence as reviewed in \Cref{SS: local geometric Langlands}.

In this section, we will assume that $G$ is an unramified reductive group over $\mO_F$. I.e., we assume that $\bar\tau=1$. We assume that $\La$ is an algebraically closed field over $\bZ_\ell$. We fix 
\begin{itemize}
\item a pinning $(G,B,T,e)$ (over $\mO_F$);
\item an additive character $\psi\colon k_F\to \La^\times$ whose conductor is $\mO_F$ (i.e. $\psi(\mO_F)=1$ but $\psi(\varpi^{-1}\mO_F)$ is non-trivial). 
\end{itemize}
Let $\iw\subset L^+G\subset LG$ be the Iwahori subgroup and the hyperspecial subgroup determined by the pinning as before. Let $\iw^u\subset \iw$ be the pro-unipotent radical of $\iw$ and let $\iw^u\to \bG_a$ be the homomorphism determined by the pinning. 

For a space, a.k.a. a (perfect) prestack $Z$ defined over $k_F$, we use the same notation to denote its base change to $k$, which is equipped with an endomorphism $\phi$ induced by the $\sharp k_F$-Frobenius endomorphism of $Z$ defined over $k_F$. In particular, the category $\shv(Z)$ (and its variants) is equipped with an automorphism $\phi_*$.

On the dual side, we base change everything to $\La$, and omit $\La$ from the notations. Every geometric object in the dual side also admits a $\phi$-action, defined similarly as in \eqref{eq: phi automorphism of loc}. Then all the coherent categories are equipped with an automorphism $\phi_*$.

\subsection{Reminder: unipotent and tame local geometric Langlands correspondence}\label{SS: local geometric Langlands}

We summarize the main results of Arkhipov-Bezrukavnikov's and Bezrukavnikov's works (\cite{bezrukavnikov2016two} \cite{arkhipov2009perverse}), their modular coefficients analogue as established by Bezrukavnikov-Riche  in \cite{Bez.Riche.modular},
 and their monodromic generalizations in \cite{DYYZ2}.  
 
We assume that either $\La=\overline\bQ_\ell$, or $\overline\bF_\ell$. In the latter case, we assume that $\ell$ is large relative to $G$ as in \cite{Bez.Riche.modular}. More precisely, we assume that $\ell$ is bigger than the Coxeter number of any simple factor of (the adjoint group of) $G$ and $\ell\neq 19$ (resp. $\ell\neq 31$) when $G$ has a simple factor of type $E_7$ (resp. $E_8$).

Let $P\subset G$ be a standard parabolic (i.e. a parabolic containing $B$). It determines a parahoric $\mP$ such that $\iw\subset L^+\mP\subset L^+G$.
We let $\hat{P}\supset \hat{B}$ be the corresponding standard parabolic subgroup of $\hat{G}$. Let $\hat{M}$ be the Levi quotient of $\hat{P}$.
Let 
\[
\locsys_{{}^cP,\breve F}^{\tame}\to \locsys_{{}^cM, \breve F}^{\tame}
\]
be as in \eqref{eq:locP and locM} (but with $W_F$ replaced by $I_F^t$).

\begin{theorem}\label{thm:Bez-equivalence}
We fix above choices.
\begin{enumerate}
\item\label{thm:Bez-equivalence-tame-and-unip} There are canonical $\phi$-equivariant equivalence of monoidal categories
\begin{equation}\label{eq:Bez-equivalence-unip}
\bB^{\unip} \colon  \rshv(\iw\backslash LG/\iw)\cong \indcoh(S^{\unip}_{{}^cG,\breve F}),
\end{equation}
and
\begin{equation}\label{eq:Bez-equivalence-mon}
\bB^{\mon} \colon \shv_{\mon}(\iw^u\backslash LG/\iw^u)\cong \indcoh(S^{\tame}_{{}^cG,\breve F}).
\end{equation}
For $\chi,\chi':\pialg(\mS_k)\to\La^\times$, the equivalence \eqref{eq:Bez-equivalence-mon} restricts to an equivalence
\[
\shv\bigl((\iw,\hchi)\backslash LG/(\iw,\hchi')\bigr)\cong \indcoh(S^{\hchi,\hchi'}_{{}^cG,\breve F}).
\]
Under the above equivalences, there is a natural $\phi$-equivariant equivalence of bimodule categories
\begin{equation}\label{eq:Bez-equivalence-bimodule}
 \shv(\iw^u\backslash LG/\iw)\cong \ind\Coh(\locsys^{\tame}_{{}^cB,\breve F}\times_{\locsys^{\tame}_{{}^cG,\breve F}} \locsys_{{}^cB, \breve F}^{\unip}).
\end{equation}

\item\label{thm:Bez-equivalence-central} Under the equivalence \eqref{eq:Bez-equivalence-unip}, the functor $\mZ^{\unip}: \Rep(\hat{G})\to \indcoh(S^{\unip}_{{}^cG,\breve F})$ from \eqref{eq:spec central functor unip}
equipped with the action \eqref{eq: taut action of spectral central sheaf unip}  corresponds to Gaitsgory's central functor $\mZ: \rep(\hat{G})\to \rshv(\iw\backslash LG/\iw)$ equipped with the monodromic action of $I_F^t$ on the nearby cycles.
Similarly, under the equivalence \eqref{eq:Bez-equivalence-mon}, the functor $\mZ^{\tame}: \Rep(\hat{G})\to \indcoh(S^{\tame}_{{}^cG,\breve F})$ from \eqref{eq:spec central functor} equipped with the action \eqref{eq: taut action of spectral central sheaf} 
corresponds to the monodromic central functor $\mZ^{\mon}: \rep(\hat{G})\to \shv_{\mon}(\iw^u\backslash LG/\iw^u)$ equipped with the monodromic action of nearby cycles.

\item\label{thm:Bez-equivalence-grading} Under the above equivalences and under the canonical isomorphism $\pi_1(G)\cong \xch(Z_{\hat{G}})$, the natural $\pi_1(G)$-grading on the left hand side (induced by decomposition of $LG$ into connected components) corresponds to the natural $\xch(Z_{\hat{G}})$-grading on the right hand (induced by the $Z_{\hat{G}}$-gerbe structure on $S_{{}^cG,\breve F}^{\unip}$ and on $S_{{}^cG,\breve F}^{\tame}$).

\item\label{thm:Bez-equivalence-duality} The equivalence $\bB^{\unip}$ intertwines the canonical duality $\verd^{\can}$ of $\rshv(\iw\bs LG/\iw)$ and the twisted Grothendieck-Serre duality $\verd^{\indcoh'}$ of $\indcoh(S^{\unip}_{{}^cG,\breve F})$ (see \eqref{eq: modified GS duality}). Similarly, the equivalence $\bB^{\tame}$ intertwines the canonical duality $\verd^{\can}$ of $\shv_{\mon}(\iw^u\bs LG/\iw^u)$ and the twisted Grothendieck Serre duality $\verd^{\indcoh'}$ of $\Ind\Coh(S^{\tame}_{{}^cG,\breve F})$. 

\item When $\La=\overline\bQ_\ell$, under the equivalences  \eqref{eq:Bez-equivalence-unip} and \eqref{eq:Bez-equivalence-mon},  the following module categories are also $\phi$-equivariantly equivalent
\begin{equation}\label{eq:Bez-equivalence-Iwahori-Whit-mon}
    \shv(\iw\backslash LG/ (\iw^u,\psi))\cong  \indcoh(\locsys^{\unip}_{{}^cB,\breve F}).
\end{equation}
\begin{equation}\label{eq:Bez-equivalence-Iwahori-Whit-mon}
    \shv_{\mon}(\iw^u\backslash LG/ (\iw^u,\psi))\cong \indcoh(\locsys^{\tame}_{{}^cB,\breve F}).
\end{equation}

\item When $\La=\overline\bQ_\ell$, under the equivalence \eqref{eq:Bez-equivalence-mon} and  \eqref{eq:Bez-equivalence-unip},  the following module categories are also $\phi$-equivariantly equivalent
\begin{equation}\label{eq:Bez-equivalence-Iwahori-hyperspecial-mon}
     \shv(\iw^u\backslash LG/ L^+G)\cong \Ind\Coh(\locsys_{{}^cB, \breve F}^{\tame}\times_{\locsys_{{}^cG,\breve F}}  \locsys_{{}^cG,\breve F}^{\unr}).
\end{equation}
Here $\locsys_{{}^cG,\breve F}^{\unr}\cong \bB \hat{G}$ denotes the stack of trivial representations of $I_F$ (see \Cref{ex: unramified component}.) Similarly, we have
\begin{equation}\label{eq:Bez-equivalence-Iwahori-hyperspecial}
     \rshv(\iw\backslash LG/ L^+G)\cong \Ind\Coh(\locsys_{{}^cB, \breve F}^{\unip}\times_{\locsys_{{}^cG,\breve F}}  \locsys_{{}^cG,\breve F}^{\unr}).
\end{equation}
More generally, for the parahoric $L^+\mP$ contained in $L^+G$ determined by a standard parabolic subgroup $P\subset G$. Then there is a canonical equivalence 
\[
\Ind\Coh(\locsys_{{}^cB, \breve F}^{\tame}\times_{\locsys_{{}^cG,\breve F}}  \locsys_{{}^cP,\breve F}^{\tame}\times_{\locsys_{{}^cM,\breve F}^{\tame}}\locsys_{{}^cM,\breve F}^{\mathrm{unr}})\simeq \shv(\iw^u\backslash LG/ L^+\mP).
\]

\item\label{thm:Bez-equivalence-duality-matching object}  The functor $\bB^{\unip}$ sends 
\begin{align}
& J_\la\mapsto \cohdual_{S^{\unip}_1}(\la),\quad \la\in \xcoch(T)   \label{al:Wakimoto-matching}\\
&   \Delta_w\mapsto  \cohdual_{S^{\unip}_w}, \quad \nabla_w\mapsto \mO_{S^{\unip}_w}[-\dim \hat{T}] , \quad w\in W_0. \label{al:finite-Hecke-matching}
\end{align}
Here we recall the Wakimoto sheaf $\{J_\la, \la\in \xcoch(T)\}$ is defined by requiring $J_\la=\nabla_\la$ if $\la$ is anti-dominant and  $J_{\la_1+\la_2}=J_{\la_1}\star J_{\la_2}$. On the coherent side, we use of notations from \Cref{notion: line bundle}.

Similarly, the functor $\bB^{\mon}$ sends
\begin{align}
& \widetilde{J}^{\mon}_\la\mapsto \cohdual_{S^{\tame}_1}(\la),\quad \la\in \xcoch(T)   \label{al:Wakimoto-matching-tame}\\
&   \widetilde\Delta^{\mon}_w\mapsto  \cohdual_{S^{\tame}_w}, \quad w\in W_0. \label{al:finite-Hecke-matching-tame}
\end{align}

\item Suppose $\La=\overline\bQ_\ell$. Let $\underline{c}\leftrightarrow O_{\underline{c}}$ be Lusztig's bijection between two-sided cells of $\widetilde W$ and unipotent conjugacy classes of $\hat{G}$. 
Then the monoidal equivalence $\bB^{\unip}$ induces equivalences of bi-modules
\[
 \fgshv(\iw\backslash LG/\iw, \overline{\bQ}_{\ell})_{\leq \underline c}\cong \Coh_{O_{\leq \underline{c}}}(S^{\unip}_{\hat{G},\overline\bQ_\ell}) 
\]
and left modules
\[
 \fgshv(\iw\backslash LG/ (\iw^u,\psi))_{\leq c}\cong   \Coh_{O_{\leq \underline{c}}}(\locsys^{\unip}_{{}^cB,\breve F}).
\]
\end{enumerate}
\end{theorem}

\begin{remark}
\begin{enumerate}
\item Bezrukavnikov established the equivalence $\bB^{\unip}$ at the level of triangulated category when $F$ is equal characteristic, $G$ is split and $\La=\overline\bQ_\ell$ (see \cite{bezrukavnikov2016two}).  This is usually called the Bezrukavnikov's equivalence. He established, at the same time, various properties of $\bB^{\unip}$ (some of which will be commented below). That such equivalence with its desired properties (in equal characteristic) can be enhanced at the $\infty$-categorical level is well-known to experts. In the Betti setting, such enhancement has been realized in \cite{Dhillon.Taylor}.

\item Note that Bezrukavnikov's original formulation uses the stack $S_{\hat{G}}^{\unip}=\hat{U}/\hat{B}\times_{\hat{G}/\hat{G}}\hat{U}/\hat{B}$ rather than $S_{{}^cG,\breve F}^{\unip}$ and therefore the equivalence of \emph{loc. cit.} depends on a choice of tame generator $\tau$.  Formulated as above, it is canonically independent of any choice. (Of course we still need to fix a pinning of $G$ and an additive character $\psi$.)

In addition, the monoidal structure of $\Coh(S_{\hat{G}}^{\unip})$ used in \cite{bezrukavnikov2016two} is the $*$-convolution. See \Cref{rem: shrek and star convolution of spectral hecke}. Given that remark, our matching of objects in Part \eqref{thm:Bez-equivalence-duality-matching object} differ from the matching of objects in \cite{bezrukavnikov2016two} by a shift. Taking \Cref{rem: Borel Weil} into account, we see that $J_\la$ should be defined such that $J_\la$ is costandard when $\la$ is anti-dominant. We also note the matching of objects in Part \eqref{thm:Bez-equivalence-duality-matching object}  is also consistent with \Cref{lem: comparison of two duality for spectral affine Hecke}.

\item When $G$ is split over $k_F$, that the equivalence $\bB^{\unip}$ is compatible with the $*$-pullback of constructible sheaves along the Frobenius endomorphism of $\iw\bs LG/\iw$, and the $*$-pullback of coherent sheaves along the automorphism of $S^{\unip}_{{}^cG,\breve F}=\hat{U}/\hat{B}\times_{\hat{G}/\hat{G}} \hat{U}/\hat{B}$ given by $(u,g_1\hat{B},g_2\hat{B})\mapsto (u^{\frac{1}{q}}, g_1\hat{B}, g_2\hat{B})$ (see \cite[Proposition 53]{bezrukavnikov2016two}).
This implies the $\phi$-equivariance of $\bB^{\unip}$ in the split case. The general case follows from the fact that $\bB^{\unip}$ is compatible with $\mathrm{Out}(G)=\mathrm{Out}(\hat{G})$-actions on both sides. Unfortunately, this fact has not been documented in literature. Similarly, Part \eqref{thm:Bez-equivalence-duality} has not appeared in literature yet. These compatibilities will be checked in a forthcoming work by Xinyu Li.

\item Bezrukavnikov's equivalence in mixed characteristic has been established in \cite{Ban23} by identifying the affine Hecke category in the mixed characteristic and in equal characteristic. On the other hand, Gaitsgory's central functor in mixed characteristic is constructed in \cite{ALWY}. These two works a priori are not directly related. Verifying unipotent part of Part \eqref{thm:Bez-equivalence-central} of the theorem in mixed characteristic is a subject of a forthcoming work by Bando, Gleason, Louren\c{c}o, and Yu.

\item Extensions of the Bezrukavnikov equivalence from the unipotent case to the tame case, with all the desired properties, will appear in \cite{DYYZ2}. In equal characteristic, it is also be possible to deduce the tame case in the \'etale setting from  \cite{Dhillon.Taylor}.
\end{enumerate}
\end{remark}

\subsection{Categorical equivalences}

\subsubsection{Tame categorical local Langlands correspondence}\label{SS: tame categorical local Langlands}
Now we arrive to our main theorem.
\begin{theorem}\label{eq:main-theorem-categorical-tame-local-Langlands}
Assume that $\La=\overline\bQ_\ell$. 
\begin{enumerate}
\item There is a canonical equivalence of categories
\[
\bL_G^{\tame}\colon \shv^{\tame}(\kot_G)\cong \indcoh(\locsys_{{}^cG,F}^{\tame}),
\]
fitting into the following commutative diagram
\[
\xymatrix{
\shv_{\mon}(\iw^u\backslash LG/\iw^u) \ar^-{\bB^{\mon}}[rr] \ar_{\Ch_{LG,\phi}^{\mon}}[d] &&  \ind\Coh(S^{\tame}_{{}^cG,\breve F}) \ar^{\Ch^{\tame}_{{}^cG,\phi}}[d]\\
\shv^{\tame}(\kot_G) \ar^-{\bL^{\tame}_G}[rr] && \indcoh(\locsys_{{}^cG,F}^{\tame}).
}
\]
In addition, $\bL_G^{\tame}$ restricts to an equivalence
\[
\bL_G^{\widehat\unip} \colon \shv^{\widehat\unip}(\kot_G)\cong \indcoh(\locsys_{{}^cG,F}^{\widehat\unip}).
\]
More generally for every tame inertia type $\zeta$, $\bL_G^{\tame}$ restricts to an equivalence
\[
\bL_G^{\hat\zeta} \colon \shv^{\hat\zeta}(\kot_G)\cong \indcoh(\locsys_{{}^cG,F}^{\hat\zeta}).
\]

\item We have
\begin{eqnarray*}
   & \mathbb{L}^{\tame}_{G}((i_{1})_{*}\delta_{I^u})
    \cong  \cohspr_{^{c}G,F}^{\tame} , & \quad \mathbb{L}^{\tame}_{G}((i_{1})_{*}\delta_{I})
    \cong  \cohspr_{^{c}G,F}^{\unip},\\
  &   \mathbb{L}^{\tame}_{G}((i_{1})_{*}\IW) \cong \mO_{\locsys_{{}^cG,F}^{\tame}}\cong \cohdual_{\locsys_{{}^cG,F}^{\tame}}, & \quad \bL^{\tame}_{G}((i_{1})_{*}\IW^{\unip}) \cong  \mO_{\locsys_{{}^cG,F}^{\widehat\unip}}\cong \cohdual_{\locsys_{{}^cG,F}^{\widehat\unip}},\\
   & \mathbb{L}^{\tame}_G((i_{1})_{*}\delta_{\sph})\cong \mO_{\locsys_{{}^cG,F}^{\mathrm{unr}}}\cong \cohdual_{\locsys_{{}^cG,F}^{\mathrm{unr}}}.  &
\end{eqnarray*}
 \item\label{eq:main-theorem-categorical-tame-local-Langlands-grading} Under the equivalence $\bL^{\tame}_G$, the natural $\pi_1(G)_\sigma$-grading on the left corresponds to the negative of the natural  $\xch(Z(\hat{G}^\sigma))$-grading on the right.
 \item\label{eq:main-theorem-categorical-tame-local-Langlands-duality} The functor $\bL^{\tame}_G$ intertwines the canonical duality $\verd^{\tame,\can}_{\kot_G}$ of $\shv^{\tame}(\kot_G)$ (see \eqref{eq: canonical duality tame part}) 
 and the twisted Grothendieck-Serre duality $\verd^{\indcoh'}$ of $\indcoh(\locsys_{{}^cG,F}^{\tame})$.
\end{enumerate}
\end{theorem}

In \Cref{ss: Matching objects under cat equiv}, we will match more objects under such equivalence.

\begin{proof}
Taking the $\phi$-twisted categorical trace of the equivalence \eqref{eq:Bez-equivalence-mon}, we obtain the following commutative diagram
\[\xymatrix{
\ar@/_8pc/_-{\Ch_{LG,\phi}^{\mon}}[dd] \shv_{\mon}(\iw^u\backslash LG/\iw^u) \ar^-{\cong}[r]\ar[d] & \indcoh(S^{\tame}_{{}^cG,F}) \ar[d] \ar@/^8pc/^-{\Ch_{{}^cG,\phi}^{\tame}}[dd]\\
    \tr(\shv_{\mon}(\iw^u\backslash LG/\iw^u),\phi) \ar^-{\cong}[r]\ar^-{\cong}_{\Cref{prop: tame llc as a categorical trace}}[d] &  \tr(\indcoh(S^{\tame}_{{}^cG,F}),\phi) \ar^-{\Cref{prop:trace-spectral-affine-Hecke-category}}_\cong[d]\\
    \shv^{\tame}(\kot_G) \ar[r] &  \indcoh(\locsys^{\tame}_{{}^cG,F}).
}\]
Then $\bL^{\tame}_G$ is defined to be the functor in the last row that makes the diagram commutative. 

Similarly, taking the $\phi$-twisted categorical trace of the equivalence \eqref{eq:Bez-equivalence-unip}, we obtain the equivalence $\bL_G^{\unip}$ fitting into a commutative diagram as above.

Under the equivalence $\bB^{\mon}$, monoidal units are identified. That is, we have the canonical isomorphism 
\[
\bB^{\mon}(\widetilde{\Delta}^{\mon}_{e})\cong (\Delta_{\locsys_{{}^cB,\breve F}^{\tame}/\locsys_{{}^cG,\breve F}^{\tame}})_*\cohdual_{\locsys_{{}^cB,\breve F}^{\tame}},
\] 
where $\Delta_{\locsys_{{}^cB,\breve F}^{\tame}/\locsys_{{}^cG,\breve F}^{\tame}}:  \locsys_{{}^cB,\breve F}^{\tame}\to \locsys_{{}^cB,\breve F}^{\tame}\times_{\locsys_{{}^cG,\breve F}^{\tame}} \locsys_{{}^cB,\breve F}^{\tame}$ is the diagonal map.
On the representation theory side, $\Ch_{LG,\phi}^{\mon}(\widetilde{\Delta}^{\mon}_{e})\cong (i_1)_*\delta_{I^u}$ by \Cref{lem:ADLV-sheaf-minimal-length}. On the spectral side, $\Ch_{{}^cG,\phi}^\tame((\Delta_{\locsys_{{}^cB,\breve F}^{\tame}/\locsys_{{}^cG,\breve F}^{\tame}})_*\cohdual_{\locsys_{{}^cB,\breve F}^{\tame}})=\cohspr_{{}^cG,F}^{\tame}$ by definition (see \Cref{ex: spectral DL induction w equal 1}). Therefore, 
\[
\bL_G^\tame((i_1)_*\delta_{I^u})\cong  \cohspr_{{}^cG,F}^{\tame}.
\]

Next, under the equivalence $\bB^{\mon}$, the module categories $\indcoh(\locsys_{{}^cB,\breve F}^{\tame})$ and $\shv_{\mon}(\iw^u\backslash LG/(\iw^u,\psi_1))$ gets identified. By \Cref{prop-class-of-iwahori-whit-module} we have 
\[
\Ch_{LG,\phi}^{\mon}(\shv_{\mon}(\iw^u\backslash LG/(\iw^u,\psi_1)),\phi)\cong (i_1)_*\IW_{\psi_1}
\] 
and by \Cref{prop: class of diagonal} we have 
\[
\Ch_{{}^cG,\phi}^{\tame}(\indcoh(\locsys_{{}^cB,\breve F}^{\tame}),\phi)\cong \cohdual_{\locsys_{{}^cG,F}^{\tame}}\cong \mO_{\locsys_{{}^cG,F}^{\tame}}.
\] 
Therefore, we see that
\[
\bL_G^\tame((i_1)_*\IW_{\psi_1})\cong \cohdual_{\locsys_{{}^cG,F}^{\tame}}\cong \mO_{\locsys_{{}^cG,F}^{\tame}}.
\]
By the similar argument, we have
\[\
\bL_G^\tame((i_1)_*\delta_{I})\cong  \cohspr_{^{c}G,F}^{\unip},\quad  \bL^{\tame}_{G}((i_{1})_{*}\IW^{\unip}) \cong  \mO_{\locsys_{{}^cG,F}^{\widehat\unip}}.
\]

Under $\bB^{\mon}$, we have the identification of module categories \eqref{eq:Bez-equivalence-Iwahori-hyperspecial-mon}.
By \Cref{prop-class-of-level}, we have 
\[
\Ch_{LG,\phi}^{\mon}(\shv(\iw^u\backslash LG/L^+G),\phi)\cong (i_1)_*\delta_{\sph},
\] 
and by \Cref{prop: class of parabolic spectral}, we have 
\[
\Ch_{{}^cG,\phi}^{\tame}(\Ind\Coh(\locsys_{{}^cB, \breve F}^{\unip}\times_{\locsys_{{}^cG,\breve F}}  \locsys_{{}^cG,\breve F}^{\unr}))\cong \mO_{\locsys^{\mathrm{unr}}_{{}^cG,F}}.
\] 
Therefore, we have 
\[
\mathbb{L}^{\tame}_G((i_{1})_{*}\delta_{\sph})\cong \mO_{\locsys_{{}^cG,F}^{\mathrm{unr}}}.
\]

Part  \eqref{eq:main-theorem-categorical-tame-local-Langlands-grading} follows directly from \Cref{thm:Bez-equivalence} \eqref{thm:Bez-equivalence-grading}, and the discussions in \Cref{rem: grading of coherent sheaves on loc} and \Cref{rem: grading on ADLI}.

As the self-dualities $\verd^{\indcoh'}$ and $\verd^{\can}$ correspond to each other under $\bB^{\mon}$ and $\bB^{\unip}$ by  \Cref{thm:Bez-equivalence} \eqref{thm:Bez-equivalence-duality}, we see that the induced self-dualities of the $\phi$-twisted categorical trace \Cref{lem-dualizability-of-categorical-trace} match with each other.
Now the claim follows from \Cref{prop: spectral side identifying duality} and \Cref{prop: automorphic side identifying duality}.
\end{proof}

We next consider modular coefficients. We shall only state the unipotent version of the equivalence. The proof is the same as  in \Cref{eq:main-theorem-categorical-tame-local-Langlands}.

\begin{theorem}\label{eq:main-theorem-categorical-tame-local-Langlands-modular}
Suppose $\La=\overline\bF_\ell$. 
\begin{enumerate}
\item\label{eq:main-theorem-categorical-tame-local-Langlands-modular-1} There is a canonical fully faithful embedding of categories
\[
\bL_G^{\unip}\colon \ind\fgshv^{\unip}(\kot_G)\hookrightarrow \indcoh(\locsys_{{}^cG,F}^{\widehat\unip}),
\]
fitting into the following commutative diagram
\[
\xymatrix{
\rshv(\iw\backslash LG/\iw) \ar^-{\bB^{\unip}}[rr] \ar_{\Ch_{LG,\phi}^{\unip}}[d] &&  \ind\Coh(S^{\unip}_{{}^cG,\breve F}) \ar^{\Ch^{\unip}_{{}^cG,\phi}}[d]\\
\rshv^{\unip}(\kot_G) \ar^-{\bL^{\unip}_G}[rr] && \indcoh(\locsys_{{}^cG,F}^{\tame}).
}
\]
In addition, the essential image is stable under the action $\ind\Perf(\locsys_{{}^cG,F}^{\widehat\unip})$ and if $Z_G$ is connected, then the essential image contains $\ind\Perf(\locsys_{{}^cG,F}^{\widehat\unip})$.
\item We have
\begin{align*}
    & \mathbb{L}^{\unip}_{G}((i_{1})_{*}\delta_{I})
    \cong  \cohspr_{^{c}G,F}^{\unip}.
\end{align*}
 \item\label{eq:main-theorem-categorical-unip-local-Langlands-grading} Under the equivalence $\bL^{\tame}_G$, the natural $\xch(Z(\hat{G}^\sigma))$-grading on the left corresponds to the natural $\pi_1(G)_\sigma$-grading on the right.
 \item\label{eq:main-theorem-categorical-unip-local-Langlands-duality} The functor $\bL^{\tame}_G$ intertwines the canonical duality $\verd^{\tame,\can}_{\kot_G}$ of $\shv^{\tame}(\kot_G)$ (see \eqref{eq: canonical duality tame part}) 
 and the modified Grothendieck-Serre duality $\verd^{\indcoh'}$ of $\indcoh(\locsys_{{}^cG,F}^{\tame})$.
\end{enumerate}
\end{theorem}

We have the following corollary about coherent sheaves on the stack of Langlands parameters.
\begin{corollary}
Assume that $\La=\overline\bQ_\ell$.
\begin{enumerate}
\item We have canonical isomorphisms
\begin{align*}
&\End_{\locsys_{{}^c G,F}}(\cohspr^{\tame}_{{}^c G,F})\cong H_{I^u}, \\
&\End_{\locsys_{{}^c G,F}}(\cohspr^{\unip}_{{}^c G,F})\cong H_I.
\end{align*}
The last isomorphism also holds if $\La=\overline\bF_\ell$. 
\item We have canonical isomorphisms
\begin{align*}
&\rg(\locsys^{\tame}_{{}^c G,F}, \cohspr^{\tame}_{{}^c G,F})=C_c(I^u\backslash G(F)/(I^u,\psi_1)),\\
&\rg(\locsys^{\tame}_{{}^c G,F}, \cohspr^{\unip}_{{}^c G,F})=C_c(I \backslash G(F)/(I^u,\psi_1)).
\end{align*}
\item We have canonical isomorphism
\[
\rg(\locsys^{\tame}_{{}^c G,F},\mO)\cong C_c((I^u,\psi_1)\backslash G(F)/(I^u,\psi_1)). 
\]
\end{enumerate}
\end{corollary}
We remind the readers that according to our conventions, both $\End$ and $\rg(\locsys^{\unip}_{{}^cG,F,\iota},-)$ are derived functors.
\begin{proof}
The statements follow from fully faithfulness of $\bL_G$ and explicit matching of objects under $\bL_G$. The first isomorphism follows from
\[
 \End_{\locsys_{{}^c G,F}}(\cohspr^{\tame}_{{}^cG,F})\stackrel{\bL^{\tame}_G}{\cong} \End_{\shv(\kot_G)}(i_{1,*}\delta_{I^u})=\End_{\rep(G(F))}(\delta_{I^u})=H_{I^u}.
\]
The unipotent version is similar.

The second statement follows from
\begin{eqnarray*}
& &\rg(\locsys^{\tame}_{{}^cG,F}, \cohspr^{\tame}_{{}^cG,F})\\
&=&\Hom_{\locsys^{\tame}_{{}^cG,F}}(\mO,\cohspr^{\tame}_{{}^cG,F})\\
&\cong& \Hom_{\shv(\kot_G)}(i_{1,*}(\IW_{\psi_1}), i_{1,*}\delta_{I^u})\\
&=&C_c(I^u\backslash G(F)/(I^u,\psi)).
\end{eqnarray*}
The unipotent version is similar.

The last statement follows from
\[
\rg(\locsys^{\tame}_{{}^c G,F},\mO)=\End(\locsys^{\tame}_{{}^c G,F},\mO)\cong \End_{G(F)}(\IW_{\psi_1})=C_c((I^u,\psi_1)\backslash G(F)/(I^u,\psi_1)).
\]
\end{proof}

We also note the following statement, which gives an explicit formula of the ``spectral action" in certain cases. Recall the ``evaluation bundle" (or called the ``tautological bundle") as from \Cref{ex: evaluation bundle}.

\begin{lemma}\label{lem: spectral action formula}
Let $\mZ^{\mon}: \rep(\hat{G})\to \shv_{\mon}(\iw^u\backslash LG/\iw^u)$ be the monodromic central functor.
Let $\mF\in \shv_{\mon}(\iw^u\backslash LG/\iw^u)$.
Then
\[
\bL_G^{\tame}(\Ch_{LG,\phi}^{\mon}(\mF\star^u \mZ^{\mon}(V)))\cong \bL_G^{\tame}(\Ch_{G,\phi}^{\mon}(\mF))\otimes \widetilde{V}.
\]
Similarly, we have the unipotent central functor $\mZ^{\unip}: \rep(\hat{G})\to \rshv(\iw\bs LG/\iw)$. For $\mF\in \rshv(\iw\backslash LG/\iw)$, we have
\[
\bL_G^{\unip}(\Ch_{LG,\phi}^{\unip}(\mF\star^u \mZ^{\unip}(V)))\cong \bL_G^{\unip}(\Ch_{G,\phi}^{\mon}(\mF))\otimes \widetilde{V}.
\]
\end{lemma}
\begin{proof} 
This follows from  \Cref{lem:image-of-trace-of-some-objects-in-specHecke}  \eqref{lem:image-of-trace-of-some-objects-in-specHecke-2} and the compatibility between $\bL_G^{\tame}$ and $\bB^{\mon}$.  The unipotent case is similar.
\end{proof}

\subsubsection{Matching objects}\label{ss: Matching objects under cat equiv}
We can match more objects under the equivalence in \Cref{eq:main-theorem-categorical-tame-local-Langlands} and \Cref{eq:main-theorem-categorical-tame-local-Langlands-modular}.
We will only state the unipotent case.

To see how to match objects under the functor, we notice that if we write $w=t_\la w_f\in \xcoch\rtimes W_0$, then 
\[
\bB^{\unip}(\cohdual_{S^{\unip}_{{}^cG,\breve F, w_f}}\star \cohdual_{\locsys_{{}^cB,\breve F}^{\unip}}(\la))\cong  \Delta_{w_f}\star J_\la,
\]
and 
\[
\bB^{\unip}( \mO_{S^{\unip}_{{}^cG,\breve F, w_f}}[-\dim\hat{T}]\star \cohdual_{\locsys_{{}^cB,\breve F}^{\unip}}(\la))\cong \nabla_{w_f}\star J_\la.
\]
It follows from \Cref{lem:image-of-trace-of-some-objects-in-specHecke} that 
\begin{equation}\label{eq: Matching objects under cat equiv-1}
\bL_G^{\unip}( \Ch_{{}^cG,\phi}^{\unip}(\nabla_{w_f}\star J_\la))\cong (\widetilde\pi^{\unip}_{w_f})_*\mO_{\widetilde\locsys^{\unip}_{{}^cG,F,w_f}}(\la).
\end{equation}
and
\begin{equation}\label{eq: Matching objects under cat equiv-2}
\bL_G^{\unip}( \Ch_{{}^cG,\phi}^{\unip}(\Delta_{w_f}\star J_\la))\cong (\widetilde\pi^{\unip}_{w_f})_*\cohdual_{\widetilde\locsys^{\unip}_{{}^cG,F,w_f}}(\la).
\end{equation}

We specialize this formula to the following special cases.

\begin{corollary}\label{cor: Wakimoto correspondence}
Let $\la\in \xcoch(T)^+$. Let $b$ be the image of $t_{-\la}\in \widetilde{W}$ under the isomorphism $B(\widetilde{W})_{\mathrm{str}}\cong B(G)$. Then we have
\[
\bL^{\unip}_G(i_{b,!}\cind_{I_b}^{G_b(F)}\La[-\langle 2\rho,\nu_b\rangle])\simeq \pi^{\unip}_*\mO_{\locsys^{\unip}_{{}^cB,F}}(\la),
\]
\[
\bL^{\unip}_G(i_{b,*}\cind_{I_b}^{G_b(F)}\La[-\langle 2\rho,\nu_b\rangle])\simeq \pi^{\unip}_*\mO_{\locsys^{\unip}_{{}^cB,F}}(w_0(\la)).
\]
\end{corollary}
\begin{proof}Notice that $t_{-\la}$ and $t_{-w_0(\la)}$ are in the same $\sigma$-straight conjugacy class, and give the same $b$. We have $J_\la=\Delta_{t_\la}$ and $J_{w_0(\la)}=\nabla_{t_{w_0(\la)}}$. 
Also notice that when $w_f$ is the unit element, then $\widetilde\locsys^{\unip}_{{}^cG,F,w_f}=\locsys^{\unip}_{{}^cB,F}$ whose structure sheaf and dualizing sheaf coincide.
Now the corollary follows from combining \Cref{lem:ADLV-sheaf-minimal-length}, \eqref{eq: unip affDL induction vs affDL sheaves},  \eqref{eq: Matching objects under cat equiv-1} and \eqref{eq: Matching objects under cat equiv-2}, and the above observations.
\end{proof}

Objects appearing in the following corollary can be regarded as generalizations of unipotent coherent springer sheaves. 

\begin{corollary}\label{cor: basic correspondence}
Let $w\in \widetilde{W}$ be a length zero element and let $b\in B(G)$ be the basic element corresponding to $w^{-1}$.
Then
\[
\bL_G^{\unip}((i_b)_*\cind_{I_b}^{G_b(F)}\La)\cong \bL_G^{\unip}((i_b)_!\cind_{I_b}^{G_b(F)}\La) \cong (\widetilde\pi^{\unip}_{w_f})_*\mO_{\widetilde\locsys^{\unip}_{{}^cG,F,w_f}}(\la),
\]
where $\la$ is the unique minuscule coweight such that $t_\la$ and $w^{-1}$ have same image in $\widetilde{W}/W_\af$, and $w_f\in W_0$ such that $w^{-1}=w_ft_\lambda$.
\end{corollary}
\begin{proof}
Let $v_0\in\overline{\breve\bfa}$ be the hyperspecial vertex and the (closure of the) fundamental alcove determined by the pinning as before. Then $w(v_0)\in \overline{\breve\bfa}$ is another hyperspecial vertex. Then there is a minuscule coweight $\la$ such that $w(v_0)=v_0+\la$. It follows that $w=t_{-\la} w^{-1}_f$, with $w_f\in W_0$, or $w^{-1}=w_f t_\la$. 
Then $\Delta_{w^{-1}}=\nabla_{w_f}\star J_\la$.  We then conclude as in \Cref{cor: Wakimoto correspondence}.
\end{proof}

Combining  \Cref{prop: class of parabolic spectral} and \Cref{prop-class-of-level}, we also obtain the following.
\begin{lemma}
Suppose $\La=\overline\bQ_\ell$.
Let $P\subset G$ be a standard parabolic subgroup, determining a parahoric group scheme $\mP$ of $G$ such that $I\subset P\subset \sph$. Let $\hat{P}\subset \hat{G}$ be the corresponding parabolic subgroup with $\hat{M}$ its Levi quotient.
Then
\[
\bL_G^{\tame}((i_1)_*\cind_P^{G(F)}\La)\cong  \pi_* (\cohdual_{\locsys_{{}^cP,F}\times_{\locsys_{{}^cM,F}}\locsys_{{}^cM,F}^{\unr}}).
\]
\end{lemma}

\subsubsection{Some consequences}
Here is an application.

\begin{theorem}\label{thm: cpt induction are coherent sheaf}
Suppose $\La=\overline\bQ_\ell$.
Let $b\in B(G)$ be basic, and let $P_b\subset G_b(F)$ be a parahoric subgroup of $G_b(F)$, with $L_{P_b}$ its Levi quotient. 
Let $\varrho$ be a finite dimensional representation of $L_{P_b}$. Let $\pi=\cind_{P_b}^{G_b(F)}\varrho$ and let
$\frakA_{\pi}:=\bL^{\tame}_G((i_b)_*\pi)$. Then $\frakA_{\pi}\in \Coh(\locsys^{\tame}_{{}^cG,F})^{\heartsuit}$. I.e. it is an honest coherent sheaf rather than a complex. 
\end{theorem}
\begin{proof}
Note that it is enough to show that from every representation $V$ of $\hat{G}$, giving a vector bundle $\widetilde V$ on $\locsys_{{}^cG,F}$ (see \Cref{ex: evaluation bundle}), we have
\[
H^i\rg(\locsys_{{}^cG,F}, \frakA_{\pi}\otimes\widetilde V)=0,\quad \mbox{ for } i\neq 0.
\]
We may assume that $\pi$ is irreducible. We may let $w$ be a length zero element in $\widetilde{W}$ giving $b$. So we may identify $G_b(F)=G(\breve F)^{\dot{w}\sigma}$.

By \Cref{prop: exactness of DL induction on tilting}  and \Cref{lem: rep of H via DL of tilting}, we may assume that there is a minimal length element $w\in W_\mP\subset \widetilde W$ such that $\pi$ appears as a direct summand of $\widetilde{R}_{\dot{w}}^{T}\in \rep(L(\kappa))^{\heartsuit}$, where $T_{w}^{\mon,f}$ is the finite Deligne-Lusztig induction of $\mathrm{Til}_{w}^{\mon}$ with respect to the Levi subgroup $L_\mP$.
Then $\cind_P^{G(F)}\pi$ appears in the affine Deligne-Lusztig induction $T_w^{\mon}$ of $\mathrm{Til}_w^{\mon}$, which by \Cref{prop-local-shtuka-uw} is isomorphic to the compact induction from $P$ to $G(F)$ of $T_w^{\mon,f}$. 

Thus, it is enough to show that 
\[
H^i\Gamma(\locsys^{\tame}_{{}^cG,F}, \bL^{\tame}_G(\cind_P^{G(F)}\mathrm{Til}_w^{\mon,f})\otimes\widetilde V)=0,\quad  i\neq 0.
\]
Under the categorical equivalence $\bL^{\tame}_G$, using \Cref{lem: spectral action formula} this is translated back to the vanishing of
\[
H^i\Hom_{\shv(\kot_G)}(\Ch_{LG,\phi}^{\mon}(\mZ^{\mon}(V)\star^u \mathrm{Til}_w^{\mon}), (i_{1})_*\IW_{\psi_1})=0,\quad i\neq 0,
\]
which follows from \Cref{prop:Whittaker-exactness}.
\end{proof}

\begin{corollary}\label{cor: CohSpr are sheaves}
The coherent Springer sheaf $\cohspr^{\tame}_{{}^c G,F}$ and its unipotent version $\cohspr^{\unip}_{{}^c G,F}$ are honest coherent sheaves. For every $\bar\sigma$-stable standard parabolic subgroup $\hat{P}\subset\hat{G}$,
the sheaf $\pi_* \cohdual_{\locsys_{{}^cP,F}\times_{\locsys_{{}^cM,F}}\locsys_{{}^cM,F}^{\unr}}$ is an honest coherent sheaf.
The coherent complexes $(\widetilde\pi^{\unip}_{w_f})_*\mO_{\widetilde\locsys^{\unip}_{{}^cG,F,w_f}}(\la)$ as in \Cref{cor: basic correspondence} are honest coherent sheaves.
\end{corollary}

That $\cohspr^{\unip}_{{}^c G,F}$ is an honest coherent sheaf was conjectured in \cite{benzvi2020CohSpr, zhu2020coherent} and was also proved in \cite{Propp} by a completely different method.

\begin{corollary}\label{cor: CMness of certain coh sheaf}
Let $\La=\overline\bQ_\ell$. 
The sheaf $\bL_G^{\tame}(\cind_{P_b^u}^{G_b(F)}\La)$ is a maximal Cohen-Macaulay coherent sheaf on $\locsys^{\tame}_{{}^cG,F}$. Here as before, $P_b^u$ denote the pro-$p$-radical of a parahoric subgroup of $G_b(F)$.
\end{corollary}
\begin{proof}
For simplicity, we write $\frakM=\bL_G^{\tame}(\cind_{P_b^u}^{G_b(F)}\La)$. We know that $\frakM$ is an honest coherent sheaf. 
Note that $\cind_{P_b^u}^{G_b(F)}\La$ is self-dual with respect to the cohomological duality. 
It follows that its modified Grothendieck-Serre dual $\verd^{\Coh'}\frakM\cong \frakM$, and therefore is also an honest coherent sheaf. Therefore, the original Grothendieck-Serre dual of $\frakA$ is also an honest coherent sheaf.

Let us write $\locsys^{\tame,\Box}_{{}^cG,F}=\Spec A$. Then the $*$-pullback of $\frakM$ is a finitely generated $A$-module $M$. Note that $A$ is Gorenstein with the dualizing sheaf being $A[\dim G]$. It follows that $\Hom(M, A)=M$. Therefore, $M$ is a maximal Cohen-Macaulay module, as desired.
\end{proof}

\begin{remark}
When $\La=\overline\bF_\ell$,
we expect \Cref{thm: cpt induction are coherent sheaf} holds when $\varrho$ is a projective object in $\rep(L_{P_b})^{\heartsuit}$. In fact, the same argument works, as soon as we know that the categorical equivalence sends $\IW_{\psi_1}$ to the structure sheaf of the stack of Langlands parameters. Similarly, we expect \Cref{cor: CMness of certain coh sheaf} holds in modular coefficient setting.
\end{remark}

If $b$ is not basic, we have the following result.

\begin{proposition}\label{prop: coh amplitude of star extension of depth zero rep}
For every $b\in B(G)$ and every $\pi:=\cind_{P_b^u}^{G_b(F)}\overline\bQ_\ell$, where $P_b^u$ is the pro-$p$-radical of a parahoric subgroup of $G_b(F)$, we have 
\[
\bL_G^{\tame}((i_b)_*\pi[\langle-2\rho,\nu_b\rangle])\in \Coh(\locsys_{{}^cG,F}^{\tame})^{\leq 0},
\]
\[
\bL_G^{\tame}((i_b)_!\pi[\langle-2\rho,\nu_b\rangle])\in \Coh(\locsys_{{}^cG,F}^{\tame})^{\geq 0}.
\]
\end{proposition}
\begin{proof}
It is enough to prove the first statement since the second one follows from the first by taking the duality.

The same argument of  \Cref{thm: cpt induction are coherent sheaf} in fact shows that $\bL_G^{\tame}(\widetilde R^T_{\dot{w}})\in \Coh(\locsys_{{}^cG,F}^{\tame})^{\heartsuit}$ for every $w\in \widetilde{W}$. In addition, it is a maximal Cohen-Macaulay sheaf.

Since $\widetilde\Til^{\mon}_{\dot{w}}$ admits a filtration with associated graded being cofree costandard objects with $\widetilde\nabla^{\mon}_{\dot{w}}$ appears as a quotient, we see that $\widetilde R^T_{\dot{w}}$ admits a filtration, with associated graded being $\widetilde{R}^*_{\dot{v}}$ for $v\leq w$ and such that $\widetilde{R}^*_{\dot{w}}$ appears as the last quotient. By induction, we see that $\bL_G^{\tame}(\widetilde{R}^*_{\dot{w}})\in \Coh(\locsys_{{}^cG,F}^{\tame})^{\leq 0}$ for every $w\in \widetilde{W}$. 

Now we consider the sheaf $\widetilde{\Til}^{\mon}_{\dot{u}}\star^u\widetilde{\nabla}^{\mon}_{\dot{w}}$ with $\dot{w}$ being $\sigma$-straight, of minimal length in $W_{\breve\bff} w$, and $u\in W_{\breve \bff}$, as in \Cref{lem: ADLI for convolution sheaves}. This object then admits a filtration with associated graded being cofree costandard objects. It follows that 
\[
\bL_G^{\tame}((i_b)_*\cind_{P_b}^{G_b(F)}\widetilde{R}^{f,T}_{\dot{u}}[-\langle 2\rho, \nu_b\rangle])\in  \Coh(\locsys_{{}^cG,F}^{\tame})^{\leq 0},
\]
where $\widetilde{R}^{f,T}_{\dot{u}}$ is the finite Deligne-Lusztig induction of monodromic a tilting sheaf, which belongs to $\rep(L_{P_b})^{\heartsuit}$ and is projective by \Cref{prop: exactness of DL induction on tilting}. In addition, by allowing $u\in W_{\breve \bff}$ to vary, these objects form a set of projective generators of $\rep(L_{P_b})^{\heartsuit}$ by \Cref{lem: DL induction of tilting generates}. Therefore, $\bL_G^{\tame}((i_b)_*\pi[\langle-2\rho,\nu_b\rangle])\in \Coh(\locsys_{{}^cG,F}^{\tame})^{\leq 0}$, as desired.
\end{proof}

\begin{remark}
We note that if $b$ is not basic, then in general $\bL_G^{\tame}((i_b)_*\cind_{P_b}^{G_b(F)}\varrho[-\langle2\rho,\nu_b\rangle])$ is not an honest coherent sheaf on $\locsys_{{}^cG,F}^{\tame}$. To see this, we let $w=t_\la$ for $\la$ dominant. It is $\sigma$-straight and $t_{-\la}$ determines an element $b\in B(G)$. Then by \Cref{cor: Wakimoto correspondence} we have
\[
\bL_G^{\tame}((i_b)_*\cind_{I_b}^{G_b(F)}\La[-\langle 2\rho,\nu_b\rangle])\simeq (\pi^{\unip})_*\mO_{\locsys_{{}^cB,F}^{\unip}}(w_0(\la)).
\]
It is known that in general $\locsys_{{}^cB,F}^{\unip}$ has non-trivial derived structure (e.g. see \cite[Remark 2.3.8]{zhu2020coherent}). 
On the other hand, if $\la$ is regular dominant, then $\mO_{{}^{cl}\locsys_{{}^cB,F}^{\unip}}(w_0(\la))$ is relative ample with respect to the proper morphism ${}^{cl}\locsys_{{}^cB,F}^{\unip}\to \locsys_{{}^cG,F}^{\tame}$. See \Cref{rem: Borel Weil}.

Since $\locsys_{{}^cB,F}^{\unip}$ is quasi-smooth, $\mO_{\locsys_{{}^cB,F}^{\unip}}\in \Coh(\locsys_{{}^cB,F}^{\unip})^{[-n,0]}$ for some $n$.
Thus, if we let $\la$ be sufficiently dominant, then for each $-n\leq i\leq 0$, we have 
\[
\mH^i (\pi^{\unip})_*\mO_{\locsys_{{}^cB,F}^{\unip}}(w_0(\la))= \mH^0 (\pi^{\unip})_*\mH^i\mO_{\locsys_{{}^cB,F}^{\unip}}(w_0(\la)),
\]
which in addition is non-zero as soon as $\mH^i\mO_{\locsys_{{}^cB,F}^{\unip}}\neq 0$. This implies that $\bL_G^{\tame}((i_b)_*\cind_{P_b}^{G_b(F)}\varrho[-\langle2\rho,\nu_b\rangle])$ is not an honest coherent sheaf. 
\end{remark}

Here is a corollary of \Cref{thm: cpt induction are coherent sheaf} and \Cref{prop: coh amplitude of star extension of depth zero rep}.
\begin{corollary}\label{cor: honest rep in connective part of coh}
Let $b\in B(G)$ and let $\pi\in \rep(G_b(F),\La)^{\heartsuit}$. Then $\bL^{\tame}_G((i_b)_*\pi[-\langle 2\rho,\nu_b\rangle])\in \indcoh(\locsys_{{}^cG,F}^{\tame})^{\leq 0}$.
\end{corollary}

Similar ideas can be used to prove the following statement. (We do not make use of it in this article.)
\begin{proposition}\label{prop: flatness of universal spherical rep}
Assume $\La=\overline\bQ_\ell$.
Assume that $G$ is unramified and that the GIT quotient map $\hat{G}\sigma/\hat{G}\to \hat{G}\sigma/\!\!/\hat{G}$ is flat. Then $\delta_K$ as a module over the spherical Hecke algebra $C_c(K\backslash G(F)/K)$ is flat.
\end{proposition}
We expect the statement continues to hold when $\La=\overline\bF_\ell$. 
\begin{proof}
As $\delta_I$ is a projective generator of the Iwahori block of $G(F)$, it is enough to show that for any $C_c(K\backslash G(F)/K)$-module $M$,  the following complex
\[
\Hom_{\rep(G(F))}(\delta_I, \delta_K\otimes_{C_c(K\backslash G(F)/K)} M)=\Hom_{\locsys^{\unip}_{{}^cG,F,\iota}}(\cohspr^{\unip}, \mO_{\locsys^{\mathrm{ur}}}\otimes_{\overline{\bQ}_\ell[\hat{G}\sigma]^{\hat{G}}}M)
\]
concentrates in degree zero. Note that the left hand side concentrates in degree $\leq 0$, as $\delta_I$ is projective and $\delta_K\otimes_{C_c(K\backslash G(F)/K} M$ belongs to $\rep(G(F))^{\leq 0}$. On the other hand,
the right hand side concentrates in degree $\geq 0$, as both $\cohspr^{\unip}$ and $\mO_{\locsys^{\mathrm{ur}}}\otimes_{\overline{\bQ}_\ell[\hat{G}\sigma]^{\hat{G}}}M$ are honest coherent sheaves on $\locsys^{\unip}_{{}^cG,F}$. The claim follows.
\end{proof}

\subsection{First applications to the classical Langlands correspondence}
Ideally, one would like to deduce a classical Langlands correspondence from the categorical one.
However, the precise relation between the categorical correspondence and the classical correspondence is not straightforward. In this subsection, we give some first applications of the categorical local Langlands to the classical local Langlands correspondence. In particular, when $\La=\overline\bQ_\ell$, we will be able to attach every depth zero supercuspidal representation $\pi$ of $G$ and its extended inner forms an essential discrete Langlands parameter $\varphi_\pi$ and a representation $r$ of $C_{\hat{G}}(\varphi_\pi)$.

\subsubsection{Semisimple Langlands parameters}
A direct consequence of the categorical local Langlands correspondence is that one can attach a semisimple Langlands parameter to every depth zero irreducible representation of $G(F)$ and its various extended pure inner forms.

We assume that $\La=\overline\bQ_\ell$. The discussion below in fact also applied to the case $\La=\overline\bF_\ell$ but we have to restrict to the unipotent case.

Via the fully faithful embedding
\[
\rep^{\tame}(G_b(F))\xrightarrow{(i_b)_*}\shv^{\tame}(\kot_G)\xrightarrow{\bL_G^{\tame}} \indcoh(\locsys_{{}^cG,F}^{\tame}),
\]
we obtain a map (of $E_1$-algebras,  see \Cref{rem: functoriality between center}) 
\[
Z(\indcoh(\locsys_{{}^cG,F}^{\tame}))\to Z(\shv^{\tame}(\kot_G))\to Z(\rep^{\tame}(G_b(F))), 
\]
such that for every $\pi\in \rep^{\tame}(G_b(F))$, the following diagram is commutative
\[
\xymatrix{
Z(\indcoh(\locsys_{{}^cG,F}^{\tame}))\ar[d]\ar[r] &  \End(\bL_G^{\tame}((i_b)_*\pi))\\
Z(\rep^{\tame}(G_b(F)))\ar[r] & \End(\pi)\ar[u].
}\]

\begin{remark}
Note that for every $\pi\in \rep^{\tame}(G_b(F))$, there is a map $(i_b)_!\pi\to (i_b)_*\pi$ compatible with the $Z(\rep^{\tame}(G_b(F)))$-action. It follows that one can replace $(i_b)_*$ by $(i_b)_!$ in the above construction. The resulting map $Z(\indcoh(\locsys_{{}^cG,F}^{\tame}))\to Z(\rep^{\tame}(G_b(F)))$ does not change. (Of course if $b$ is basic, then $(i_b)_!=(i_b)_*$.)  

Similarly, one can replace $(i_b)_*$ by $(i_b)_\flat$ in the above construction.
\end{remark}

Let 
\[
Z_{G_b,F}^{\tame}=H^0Z(\rep^{\tame}(G_b(F)))
\]
denote the tame Bernstein center of $G_b(F)$. Composed with the map \eqref{eq: stable center}, we obtain a well-defined ring homomorphism
\begin{equation}\label{eq: spec center to Bernstein center}
Z^{\tame}_{{}^cG,F}\to Z^{\tame}_{G_b,F}
\end{equation}

Now, let $\pi$ be a depth zero irreducible representation of $G_b(F)$, or more generally a depth zero representation such that $H^0\End_{G_b(F)}(\pi)=\La$. Then we obtain a homomorphism
\[
Z^{\tame}_{{}^cG,F}\to Z^{\tame}_{G_b,F}\to H^0\End_{G_b(F)}(\pi)=\La,
\]
giving a $\La$-point of $\Spec Z^{\tame}_{{}^cG,F}$.
Such a point gives a semisimple  (or called completely reducible) Langlands parameter
\[
\varphi_\pi^{ss}: W_F\to {}^cG
\]
as desired.

We thus obtain the following theorem.
\begin{theorem}\label{thm: ss parameter for depth zero representations}
There is a map from the isomorphism classes of irreducible depth zero representations of $\rep^{\tame}(G_b(F))$ to the set of tame semisimple Langlands parameters $\pi\mapsto \varphi_\pi^{ss}$.
\end{theorem}

We will discuss the compatibility of the above semisimple Langlands parameters attached to $\pi$ with other parameterizations in another place.

\subsubsection{Coherent sheaves attached supercusipdal representations}\label{SSS: cllc for supercuspidal}
Let $b\in B(G)$, lifted to a $\sigma$-straight element $w_b\in \widetilde{W}$.
Some constructions below also require a choice of a lifting of $w_b$ to $G(\breve F)$. We will fix such a choice, and abuse of notations still denote it by $w_b$. 

Let $G_b$ be the corresponding twisted centralizer group. 
As explained in \Cref{rem: combinatorics of Gb},
$I_b= \iw(k)\cap G_b(F)$ is an Iwahori subgroup of $G_b(F)$. The corresponding (extended) affine Weyl group of $G_b(F)$ is identified with 
\[
\widetilde{W}^{\sigma_b}:=\{w\in \widetilde{W}\mid \sigma_b(w)=w\}.
\]
Here we write $\sigma_b=\Ad_{w_b}\sigma$ for the twisted Frobenius structure. 

In the sequel, we will assume that $b$ is basic, so that $w_b$ is a length zero element.
We can take the $\sigma_b$-invariants of the semi-direct product \eqref{eq-Iwahori-Weyl-semiproduct} gives
\[
 \widetilde{W}^{\sigma_b}=(W_\af)^{\sigma_b}\rtimes \Omega_{\breve\bfa}^\sigma,
\]
where $(W_\af)^{w_b\sigma}$ is the affine Weyl group of $G_b(F)$. 
Here we use the fact that $\Omega_{\breve\bfa}$ is commutative so the action of $\sigma_b$ on $\Omega_{\breve\bfa}$ coincides with the original $\sigma$-action. 
In particular, $ \widetilde{W}^{\sigma_b}$ maps surjectively to $(\widetilde W/W_\af)^\sigma\cong \pi_1(G)_{I_F}^\sigma$.

Let $P\subset G_b(F)$ be a parahoric containing $I_b$, corresponding to a facet $\breve\bff\subset \overline{\breve\bfa}\subset \scrA(G_{\breve F},S_{\breve F})$ stable under the action of $\sigma_b$. Let $W_P\subset \widetilde{W}^{\sigma_b}$ be the corresponding Weyl group. Note that $W_P\subset (W_\af)^{\sigma_b}$. 
Let $N_{\widetilde{W}^{\sigma_b}}(W_P)$ be the normalizer of $W_P$ in $\widetilde{W}^{\sigma_b}$. We have 
\[
N_{\widetilde{W}^{\sigma_b}}(W_P)= W_P \rtimes \Omega_P,
\]
where $\Omega_P$ consist of those $w\in \widetilde{W}^{\sigma_b}$ that fixes each simple reflection (with respect to $\breve\bfa$) of $W_P$. In particular, $N_{\widetilde{W}^{\sigma_b}}(W_P)/W_P\cong \Omega_P$.
Note that we have a natural inclusion $N_{G_b(F)}(P)/P\subset N_{\widetilde{W}^{\sigma_b}}(W_P)/W_P$. 
We have a left exact sequence
\[
1\to N_{(W_\af)^{\sigma_b}}(W_P)/W_P\to  N_{\widetilde{W}^{\sigma_b}}(W_P)/W_P\to (\widetilde W/W_\af)^\sigma.
\]
When $P$ is a maximal parahoric subgroup of $G_b(F)$, we have 
\[
N_{G_b(F)}(P)/P= N_{\widetilde{W}^{\sigma_b}}(W_P)/W_P,\quad N_{(W_\af)^{\sigma_b}}(W_P)/W_P=\{1\}.
\] 
Therefore, in this case we have
\begin{equation}\label{eq: normalizer of WP mod WP}
N_{\widetilde{W}^{\sigma_b}}(W_P)/W_P\cong \Omega_P=\{w\in \Omega_{\breve\bfa}^\sigma \mid  w(v_P)=v_P\}\subset \Omega_{\breve\bfa}^\sigma\cong \pi_1(G)_{I_F}^\sigma,
\end{equation}
where $v_P$ is the vertex in the apartment $\scrA(G_{\breve F},S_{\breve F})$ corresponding to $P$. 

Now we recall some basic facts about Hecke algebras for depth zero Bernstein blocks. We will assume that $b$ is basic and $w_b$ is a length zero element.  
Now let $(V,\varrho)$ be a representation of the Levi quotient $L_P$ of $P$, and let 
\[
\cind_P^{G_b(F)}\varrho=\bigl\{f: G_b(F)\to V\mid f(pg)=\varrho(p)(f(g)), \ \forall p\in P\bigr\}
\] 
be the compact induction.
Recall that the corresponding (underived) Hecke algebra is
\[
H^0H_{P,\varrho}:=H^0\End_{G_b(F)}(\cind_P^{G_b(F)}\varrho)=\bigl\{h: G_b(F)\to \End_{\La}(V)\mid h(p_1gp_2)=\varrho(p_1)h(g)\varrho(p_2)\bigr\},
\] 
with the action given by 
\[
h(f)(g)=\sum_{g'\in P\backslash G_b(F)} h(g{g'}^{-1})(f(g')), \quad g\in G_b(F).
\]
As a vector space, it admits a direct sum decomposition index by 
$W_P\backslash \widetilde{W}^{\sigma_b}/W_P$.
 Namely, for every $w\in W_P\backslash \widetilde{W}^{\sigma_b}/W_P$, with an representative $\dot{w}\in G_b(F)$, we have the corresponding direct summand
\[
H^0H_{P,\varrho,w}:=\bigl\{h: P\dot{w}P\to V\mid h(p_1gp_2)=\varrho(p_1)h(g)\varrho(p_2)\bigr\}\cong H^0\Hom_{P_w}(\varrho|_{P_{w^{-1}}},\varrho|_{P_w}).
\]
Here $P_w=P\cap \dot{w}P\dot{w}^{-1}$ and $P_{w^{-1}}=P\cap \dot{w}^{-1}P\dot{w}$, and we regard $\varrho|_{P_{w^{-1}}}$ as a representation of $P_w$ via the isomorphism $\Ad_{\dot{w}^{-1}}: P_w\cong P_{w^{-1}}$. The isomorphism 
$H^0H_{P,\varrho}\cong \oplus_{w} H^0H_{P,\varrho,w}$  
sends $h$ to the collection $\{h(\dot{w})\in H^0\Hom_{P_w}(\varrho|_{P_{w^{-1}}},\varrho|_{P_w})\}_w$. (Note that we do not claim that $H^0H_{P,\varrho,w}$ is non-zero.)

It is known that if $\varrho$ is an irreducible cuspidal representation of $L_P$, then $H^0H_{P,\varrho,w}=0$ if $w\not\in N_{\widetilde{W}^{\sigma_b}}(W_P)/W_P$. 
On the other hand, when $w\in N_{\widetilde{W}^{\sigma_b}}(W_P)/W_P$, then $\Ad_{w^{-1}}$ induces an outer automorphism of $L_P$, and 
$H^0H_{P,\varrho,w}\neq 0$ if that only if $\varrho$ is isomorphic to its twist by this automorphism, usually denoted by ${}^{w^{-1}}\varrho$. In this case, $H^0H_{P,\varrho,w}\neq 0$  is one dimensional with a basis given by
\begin{equation}\label{eq: Hecke operator h at w}
h_{\dot{w}}(\dot{w}): {}^{w^{-1}}\varrho\simeq \varrho, \quad h_{\dot{w}}(g)=0, \mbox{ if } g\not\in P\dot{w}P.
\end{equation}

Now suppose $P$ is a maximal parahoric subgroup. We consider the subgroup of \eqref{eq: normalizer of WP mod WP}
\[
\Omega_{P,\varrho}=\{w\in \Omega_{\breve\bfa}^{\sigma} \mid w(v_P)=v_P, \ {}^{w^{-1}}\varrho\sim \varrho\}.
\]
The we obtain a projective representation of $\Omega_{P,\varrho}$ on $\varrho$ as usual, giving a central extension $\widetilde\Omega_{P,\varrho}$ of $\Omega_{P,\varrho}$ by $\La^\times$. Then $H_{P,\varrho}$ is isomorphic to the twisted group algebra of $\Omega_{P,\varrho}$ associated to this central extension.
\begin{remark}
We do not know, nor have checked the literature, whether the central extension is always trivial, i.e. whether $H_{P,\varrho}\otimes_{H_{D(\mO_F),\chi}}\La$ is commutative in general. This is known in many cases.
\end{remark}

We write $D=Z_G^\circ$ be the maximal subtorus in $Z_G$. By abuse of notations, we also use it to denote its unique Iwahori group scheme over $\mO_F$.
Note that $D(\mO_F)\subset P$. We let $\chi$ be the restriction of the central character of $\varrho$ to $D(\mO_F)$. Then clearly we have
\begin{equation}\label{eq: center part of the HPrho}
H_{D(\mO_F),\chi}\subset H_{P,\varrho}.
\end{equation}
Namely, if $\la\in\xcoch(D)^\sigma$, giving $t_\la\in \widetilde{W}$. Then the operator $h_{t_\la}\in H_{P,\varrho,t_\la}$ of \eqref{eq: Hecke operator h at w} comes for the corresponding operator of $H_{D(\mO_F),\chi}$.
Here we lift $t_\la$ to $\la(\varpi)\in D(F)$ by chosen a uniformizer $\varpi\in F$. Note that we have inclusions of abelian groups
\[
\xcoch(D)^\sigma\subset \Omega_{P,\varrho}\subset \Omega_{\breve\bfa}^\sigma,
\]
with $\xcoch(D)^\sigma$ finite index in $\Omega_{\breve\bfa}^\sigma$.

We summarize a consequence of the above discussions
\begin{lemma}\label{lem: semisimplicity Hecke algebra}
Suppose $P$ is a maximal parahoric subgroup of $G_b(F)$. Then $H_{P,\varrho}$ is a finite free $H_{D(\mO_F),\chi}$-module. Let $H_{D(\mO_F),\chi}\to \La$ is a homomorphism, then $H_{P,\varrho}\otimes_{H_{D(\mO_F),\chi}}\La$ is a finite dimensional semisimple algebra over $\La$.
\end{lemma}

Having the above quick review of the Hecke algebra associated to some depth zero Bernstein blocks, we can prove the following result.

\begin{proposition}\label{thm: categorical Langlands for tame supercuspidal-1}
Suppose $\La=\overline\bQ_\ell$. Let $b\in B(G)$ be a basic element.
Let $\varrho$ be an irreducible cuspidal representation of $L_{P}$, where $P$ is a \emph{maximal} parahoric subgroup of $G_b(F)$. Let $\pi=\cind_{P_b}^{G_b(F)}\varrho$. Then 
there exist finitely many \emph{disjoint} irreducible components $\frakX_{P,\varrho}:=\frakX_1\sqcup\cdots \sqcup \frakX_r\subset \locsys_{{}^cG,F}^{\tame}$ such that $\bL_G^{\tame}((i_b)_*\pi)$ is a vector bundle on $\frakX_{P,\varrho}$ (regarded as a coherent sheaf on $\locsys^{\tame}_{{}^cG,F}$ via the $*$-pushforward). In addition, each $\frakX_i$ contains an essential discrete parameter.
\end{proposition}

\begin{proof}
Consider the map
\begin{equation}\label{eq: spectral center to endo of supercuspidal}
Z_{{}^cG,F}^{\tame}\to \End(\bL_G((i_b)_*\pi))=H_{P,\varrho}.
\end{equation}
We will let $Z_{{}^cG,F,P,\varrho}^{\tame}$ denote the image of the map.
As $\bL_G^{\tame}((i_b)_*\pi)$ is maximal Cohen-Macaulay, its (set-theoretic) support is the union of several irreducible components $\cup_i \frakX_i$ of $\locsys_{{}^cG,F}^{\tame}$. 
We may write $\Spec Z_{{}^cG,F,P,\varrho}^{\tame}=\cup_i Z_i$ as union of irreducible components so that $\frakX_i$ maps to $Z_i$.

Recall the free action of $C_{{}^cG}$ on $\locsys^{\tame}_{{}^cG,F}$ and on $\Spec Z^{\tame}_{{}^cG, F}$. (See the paragraph before \Cref{prop: geometry of discrete component}.)
It follows from this action that the image of each irreducible component of $\locsys_{{}^cG,F}^{\tame}$ in $\Spec Z_{{}^cG,F}^{\tame}$ at least has dimension $C_{{}^cG}$. Thus each $\dim Z_i\geq  \dim C_{{}^cG}$. Now since $Z_{{}^cG,F,P,\varrho}^{\tame} \subset  H_{P,\varrho}$ and $H_{P,\varrho}$ is finite over $H_{D(\mO_F),\chi}$, which is a commutative algebra of dimension $\dim C_{{}^cG}$, we see that $\dim Z_i\leq \dim C_{{}^cG}$. Therefore $\dim Z_i=\dim C_{{}^cG}$.

As $C_{{}^cG}$ acts free on $Z_i$, we see that $(Z_i)_\red$ is a single $C_{{}^cG}$-orbit.
By \Cref{lem: support of supercusp in one component},  $\frakX_i$ contains an essential discrete parameter $\varphi_i$. Let $z_i\in Z_i\subset \Spec  Z_{{}^cG,F}^{\tame}$ be the image of $\varphi_i$, and let $h_i=\varphi_i^{ss}$ be its semisimplification.
By \Cref{prop: geometry of discrete component}, we have $\frakX_i=C_{{}^cG}\times V_{h_i}$, containing $C_{{}^cG}\times \{\varphi_i\}/C_{\hat{G}}(\varphi_i)$ as an open substack.
If $\frakX_i$ and $\frakX_j$ are two different irreducible components in the support of $\bL_G^{\tame}((i_b)_*\pi)$, then $C_{{}^cG}\times \{\varphi_i\}/C_{\hat{G}}(\varphi_i)$ and $C_{{}^cG}\times \{\varphi_i\}/C_{\hat{G}}(\varphi_j)$ are disjoint. It follows from \Cref{lem: extends ss parameter to discrete parameter} that $z_i$ and $z_j$ are in different $C_{{}^cG}$-orbits. Therefore, $\frakX_i$ and $\frakX_j$ are disjoint.

Therefore, $\bL_G^{\tame}((i_b)_*\pi)$ is set-theoretically supported on several disjoint irreducible components $\frakX_i$ of $\locsys_{{}^cG,F}^{\tame}$, each of which contains an essential discrete parameter.  By \Cref{lem: maximal CM set-supp vs scheme-supp} below, $\bL_G^{\tame}((i_b)_*\pi)$ is scheme-theoretically supported on these irreducible components.
By \Cref{prop: geometry of discrete component}, each of these component is smooth. As $\bL_G^{\tame}((i_b)_*\pi)$ is maximal Cohen-Macaulay, it must be a vector bundle over $\frakX_{P,\varrho}=\sqcup_i \frakX_i$.
\end{proof}

\begin{lemma}\label{lem: maximal CM set-supp vs scheme-supp}
Let $X$ be an equidimensional reduced noetherian scheme of finite Krull dimension. Let $Z$ be an irreducible component of $X$. Suppose $M$ is a maximal Cohen-Macaulay module on $X$ set-theoretically supported on $Z$. Then $M$ is scheme theoretically supported on $Z$.
\end{lemma}
\begin{proof}
We may assume that $X=\Spec R$ is affine.
Write $X=Z\cup Y$ where $Y$ is the Zariski closure of $X-Z$ in $X$. Let $I\subset R$ be the ideal defining $Z$ and $J\subset R$ the ideal defining $Y$. Then $I\cdot J\subset I\cap J=\{0\}$.  Choose an element $0\neq x\in J$. Then $\dim(V(x)\cap Z)=\dim Z-1$. It follows from our assumption and \cite[\href{https://stacks.math.columbia.edu/tag/00N5}{Lemma 00N5}]{stacks-project} that $x$ is a non-zero divisor of $M$. Now for every $am$ with $a\in I$ and $m\in M$, we have $xam=0$. Therefore $am=0$. It follows that $IM=0$ so $M$ is scheme-theoretically supported on $Z$.
\end{proof}

On of the consequences of the above arguments is the following.
\begin{corollary}
The scheme $\Spec Z^{\tame}_{{}^cG,F,P,\varrho}$ is reduced. The map $\frakX_i\to Z_i$ is flat.
\end{corollary}
\begin{proof}Since $\frakX_i$ is reduced, and $H^0\rg(Z_i,\mO)$ is the image of $Z_{{}^cG,F}^{\tame}\to H^0\rg(\frakX_i,\mO)$, we see that $Z_i$ is reduced. 
Then $Z_i$ is a $C_{{}^cG}$-torsor. Since $\frakX_i\to Z_i$ is $C_{{}^cG}$-equivariant, we see that this map is flat.
\end{proof}

We will need an $S=T$ type result in a very special case. Namely, on the one hand, associated to $\la\in \xcoch(D)^{\sigma}\subset \widetilde{W}$ we have the Hecke operator $h_{t_\la}\in H^0H_{P,\varrho}$ supported on $Pt_\la P$, see \eqref{eq: Hecke operator h at w}. On the other hand, via the projection ${}^cG\to \hat{G}_{\mathrm{ab}}\rtimes (\bG_m\times \Ga_{\widetilde F/F})\to (\hat{G}_{\mathrm{ab}})_{\Ga_F}=\widehat{Z_G^s}$, $\la\in \xcoch(Z_G^s)=\xch(\widehat{Z_G^s})$ gives rise to a one dimensional representation $V_\la$ of ${}^cG$. Then we have the $S$-operator $S_{\mZ^{\mon}(V_\la), (\tau,\sigma)}$, given by multiplication by $\chi_{V_\la, (\tau,\sigma)}$, 
see \eqref{eq: tame central S operator} and \Cref{rem: tame central S operator}.

\begin{lemma}\label{lem: S=T for central hecke operator}
We have $h_{t_\la}= S_{V_\la}$.
\end{lemma}

As a consequence of \Cref{lem: S=T for central hecke operator}, we see the map \eqref{eq: center part of the HPrho} fits into the following commutative diagram
\[
\xymatrix{
Z_{{}^cD,F,D(\mO_F),\chi}^{\tame}\ar[r]\ar_{\cong}[d] &  Z_{{}^cG,F,P,\varrho}^{\tame}\ar[d]\\
H_{D(\mO_F),\chi}\ar[r] & H_{P,\varrho}.
}\]
Note that map $\Spec Z_{{}^cG,F,P,\varrho}^{\tame}\to \Spec Z_{{}^cD,F,D(\mO_F),\chi}^{\tame}$ is $C_{{}^cG}$-equivariant, where $C_{{}^cG}$ acts on the target through the homomorphism $C_{{}^cG}\to C_{{}^cD}$, which then acts on $\Spec Z_{{}^cD,F,D(\mO_F),\chi}^{\tame}$. As $C_{{}^cG}\to C_{{}^cD}$ is an isogeny, we see that this map is finite \'etale.

Now if $z: Z_{{}^cD,F,D(\mO_F),\chi}^{\tame}\to \La$ is a homomorphism, by the above discussions and by \Cref{lem: semisimplicity Hecke algebra}, we see that the algebra
\[
H_{P,\varrho,z}:=H_{P,\varrho}\otimes_{Z_{{}^cG,F, P,\varrho}^{\tame},z }\La
\] 
is a finite dimensional semisimple algebra over $\La$.

Now suppose the above homomorphism $z$ gives a point
$z\in Z_i$. Let $h$ be the semisimple parameter associated to $z$. Then the fiber of $z$ in $\frakX_i$ is $V_h$, which is the closure of $\{\varphi\}/C_{\hat{G}}(\varphi)$, where $\varphi$ is an essential discrete parameter such that $\varphi^{ss}=h$.  
We regard $z$ as a character $Z_{{}^cG,F,P,\varrho}^{\tame}$. Then 
\[
\bL_G^{\tame}((i_b)_*(\pi\otimes_{Z_{{}^cG,F,P,\varrho}^{\tame},z }\La))\cong \bL_G^{\tame}((i_b)_*\pi)\otimes_{Z_{{}^cG,F,P,\varrho}^{\tame},z }\La
\]
is a vector bundle on $V_h$, which as usual is regarded as a coherent sheaf on $\locsys^{\tame}_{{}^cG,F}$ via the $*$-pushforward.

\begin{lemma}
The representation $\pi\otimes_{Z_{{}^cG,F,P,\varrho}^{\tame},z }\La\in \rep(G_b(F),\La)^{\heartsuit}$.
\end{lemma}
\begin{proof}
Clearly $\pi\otimes_{Z_{{}^cG,F,P,\varrho}^{\tame},z }\La\in \rep(G_b(F),\La)^{\leq 0}$. As $\pi$ is a projective object in $\rep(G_b(F),\La)^{\heartsuit}$, which is a generator of the Bernstein block it belongs to, we have $\Hom(\pi, \pi\otimes_{Z_{{}^cG,F,P,\varrho}^{\tame},z }\La)\in \Mod^{\leq 0}$. On the other hand, passing to the spectral side, we see that both $\bL_G^{\tame}((i_b)_*\pi)$ and $\bL_G^{\tame}((i_b)_*(\pi\otimes_{Z_{{}^cG,F,P,\varrho}^{\tame},z }\La))$ are in honest coherent sheaves on $\locsys_{{}^cG,F}^{\tame}$. Therefore, their hom space sits in $\Mod_\La^{\geq 0}$. It follows that $\Hom(\pi, \pi\otimes_{Z_{{}^cG,F,P,\varrho}^{\tame},z }\La)\in \Mod_\La^{\heartsuit}$. Therefore, $\pi\otimes_{Z_{{}^cG,F,P,\varrho}^{\tame},z }\La\in \rep(G_b(F),\La)^{\heartsuit}$.
\end{proof}

Note that we have
\[
\pi\otimes_{Z_{{}^cG,F,P,\varrho}^{\tame},z }\La=\pi\otimes_{H_{P,\varrho}}H_{P,\varrho,z}
\]

It follows that if $E$ is a simple $H_{P,\varrho,z}$-module, then $(\cind_{P_b}^{G_b(F)}\varrho)\otimes_{H_{P,\varrho}}E$ is a direct summand of $\pi\otimes_{Z_{{}^cG,F,P,\varrho}^{\tame},z }\La$. Therefore, 
\[
(\cind_{P_b}^{G_b(F)}\varrho)\otimes_{H_{P,\varrho}}E \in \rep(G_b(F),\La)^{\heartsuit}.
\]
In addition, it follows from the classical theory (e.g. \cite[Proposition 1.4]{Morris.tame.supercuspidal}) that it is an irreducible supercuspidal representation of $G_b(F)$. Similarly, 
\[
\bL_G^{\tame}( (i_b)_*((\cind_{P_b}^{G_b(F)}\varrho)\otimes_{H_{P,\varrho}}E))\subset \bL_G^{\tame}((i_b)_*(\pi\otimes_{Z_{{}^cG,F,P,\varrho}^{\tame},z }\La))
\]
is a direct summand.

On the other hand, if $\pi$ is a depth zero irreducible supercupsidal representation of $G_b(F)$, then there is a maximal parahoric subgroup $P\subset G_b(F)$ and an irreducible cuspidal representation $\varrho$ of $L_P$, regarded as a representation of $P$ via inflation, such that $\varrho$ appears as a direct summand of $\pi|_P$ (e.g. see \cite[Proposition 2.2]{Morris.tame.supercuspidal}). In this case, there is a simple $H_{P,\varrho}$-module $E$, such that $\pi\simeq (\cind_P^{G_b(F)}\varrho)\otimes_{H_{P,\varrho}}E$. 
We thus obtain the following theorem.

\begin{theorem}\label{thm: categorical Langlands for tame supercuspidal}
Suppose $\La=\overline\bQ_\ell$. Let $b\in B(G)$ be a basic element.
Let $\pi$ be a depth zero supercuspidal representation of $G_b(F)$. Then $\bL^{\tame}_G((i_b)_*\pi)$ is a vector bundle on $V_h$, where $h=\varphi_\pi^{ss}$ is the semisimple Langlands parameter attached to $\pi$ in \Cref{thm: ss parameter for depth zero representations}, and $V_{h}$ is the stack attached to $h$ in \eqref{eq: vogan stack2}. In addition, $V_h$ contains an open substack $\{\varphi_\pi\}/C_{\hat{G}}(\varphi_\pi)$, where $\varphi_\pi$ is an essential discrete parameter.
\end{theorem}

Again we emphasize that $V_h\subset \varpi_{{}^cG,F}^{-1}(\varpi_{{}^cG,F}(h))$ but this inclusion is usually strict.

\begin{notation}
In the sequel, we shall write $\bar{i}_\varphi: V_h\to \locsys_{{}^cG,F}^{\tame}$ be the closed embedding. Then $i_{\varphi}: \{\varphi\}/C_{\hat{G}}(\varphi)\to \locsys_{{}^cG,F}^{\tame}$ factors as an open embedding $\{\varphi\}/C_{\hat{G}}(\varphi)\subset V_h$ followed by $\bar{i}_\varphi$. We shall write $\bL^{\tame}_G((i_b)_*\pi)$ as $(\bar{i}_\varphi)_*\mE_\pi$.
\end{notation}

Now we recall that
\[
\Vogan_h\cong (\hat{\frakg}\otimes \bZ_\ell(-1))^{h(W_F)}/C_{\hat{G}}(h).
\]
We need the following facts regarding vector bundles on $\Vogan_h$.
\begin{lemma}
Every vector bundle on $\Vogan_h$ is isomorphic to the pullback of a vector bundle on $\bB C_{\hat{G}}(h)$ along the natural map $(\hat{\frakg}\otimes \bZ_\ell(-1))^{h(W_F)}/C_{\hat{G}}(h)\to \bB C_{\hat{G}}(h)$. 
\end{lemma}
\begin{proof}
If $(\hat{\frakg}\otimes \bZ_\ell(-1))^{h(W_F)}=0$, then $\Vogan_h\cong \{\varphi\}/C_{\hat{G}}(\varphi)$. Therefore vector bundles on $\Vogan_h$ correspond to finite dimensional representations $C_{\hat{G}}(h)$. 
If $(\hat{\frakg}\otimes \bZ_\ell(-1))^{h(W_F)}\neq 0$, then $C_{\hat{G}}(h)$ contains a central torus $\bG_m$ that acts on $U=\hat{\frakg}^{h(W_F)}$ by some weight $n\neq 0$. Indeed, let $A$ be the Zariski closure of $\{h(\sigma^n)\}_{n\in\bZ}$ in ${}^cG$. Since $h(\sigma$ acts on $U$ by $q$, we see that $A^\circ$ is a non-trivial torus. It normalizes $h(\tau)$ and therefore centralize $h(\tau)$ (since $h(\tau)$ generates a finite group in ${}^cG$). It follows that $A^\circ$ is a central torus of $C_{\hat{G}(h)}$. In addition, $A^\circ$ acts on $U$ by a non-zero character. 
We one can choose a torus $\bG_m\subset A^\circ$ such that its weight on $U$ is non-zero. Now we apply the following lemma to conclude.
\end{proof}
\begin{lemma}
Let $U$ be a prehomogeneous vector space under the action of a (not necessarily) reductive group $L$. Suppose that there a central torus $\bG_m\subset L$ acting $U$ by some weight $n\neq 0$. Then every vector bundle on $U/L$ is isomorphic the pullback of a vector bundle on $\bB L$.  
\end{lemma}
\begin{proof}
We may assume that $n>0$.
We regard a vector bundle on $U/L$ as a finite free $\mathrm{Sym} U^*$-module $E$ with an $L$-action. Let $E_0\subset E$ be the subspace of highest weight with respect to the central $\bG_m$-action. Then $E_0$ is a subrepresentation of $L$. In addition, the composed map $E_0\subset E\to E/(\mathrm{Sym}^{>0} U^*) E$ is an isomorphism. The graded Nakayama lemma implies that the natural map $E_0\otimes  \mathrm{Sym} U^*\to E$ is an isomorphism.
\end{proof}

As a consequence, we see that given a vector bundle $\mE$ on $\Vogan_h$, $\End(\mE)=\La$ if and only if $\mE$ is the pullback of a vector bundle on $\bB C_{\hat{G}}(h)$, corresponding to an irreducible representation of $C_{\hat{G}}(h)$.

\begin{lemma}
Suppose $\pi$ is a generic supercuspidal representation with respect to our choice of Whittaker datum. I.e. $b=1$ and there is a non-zero map $\IW=\cind_{I^u}^{G(F)}\psi_1\to \pi$. Then the vector bundle $\mE$ on $V_h$ attached to $\pi$ as in \Cref{thm: categorical Langlands for tame supercuspidal} is the trivial bundle.
\end{lemma}
\begin{proof}
Let $\pi'$ be the kernel of the map $\IW\to \pi$. So we have the short exact sequence $0\to\pi'\to \IW\to \pi\to 0$, which via the equivalence $\bL_G^\tame$, gives a fiber sequence of coherent complex on $\locsys_{{}^cG,F}^{\tame}$
\[
\bL_G^{\tame}((i_1)_*\pi')\to \mO_{\locsys_{{}^cG,F}^{\tame}} \to \mE.
\]
By \Cref{cor: honest rep in connective part of coh}, $\bL_G^{\tame}((i_1)_*\pi')\in \Coh(\locsys_{{}^cG,F}^{\tame})^{\leq 0}$ so the map $\mO_{\locsys_{{}^cG,F}^{\tame}} \to \mE$ is non-zero surjective. This forces $\mE=\mO_{V_h}$, as desired.
\end{proof}

Let $\mE$ be the vector bundle on $V_h$ attached to $\pi$. By restriction of $\mE$ to $\{\varphi\}/C_{\hat{G}}(\varphi)$, we obtain a representation $r_\pi$ of $C_{\hat{G}}(\varphi)$. We thus obtain the following result.
\begin{corollary}
For every depth zero supercuspidal representation $\pi$, the semisimple Langlands parameter $\varphi_\pi^{ss}$ attached to $\pi$ in \Cref{thm: ss parameter for depth zero representations} can be lifted to an essentially discrete Langlands parameter $\varphi_\pi$. In addition, the assignment $\pi\mapsto \varphi_\pi$ can be further lifted to an enhanced local langlands parameter $(\varphi_\pi, r_\pi)$. If $\pi$ is generic with respect to the chosen Whittaker datum, then $r_\pi$ is the trivial representation of $C_{\hat{G}}(\varphi_\pi)$.
\end{corollary}

\subsubsection{Representations attached to Langlands parameters}
We assume that $\La=\overline\bQ_\ell$. 
Next we discuss the other direction of the local Langlands correspondence. Name, we discuss how to attach representations (or rather $L$-packets) to local Langlands parameters. We shall mention that some ideas presented below were also observed by David Hansen \cite{{Hansen.Beijing}} in the Fargues-Scholze's setting.

Let $\varphi: W_F^t\to {}^cG$ be a tame Langlands parameter. We suppose the corresponding (finite type) point of $\locsys_{{}^cG, F}^{\tame}$ is a smooth point. Recall that this means that $H^2(W_F^t, \Ad^0_\varphi)=0$, or equivalently $q^{-1}$ is not an eigenvalue of the linear operator $\varphi(\sigma): \hat\frakg^{\varphi(I_F)}\to  \hat\frakg^{\varphi(I_F)}$.
See the proof of \Cref{lem:discrete parameter condition}.
Let $i_\varphi: \{\varphi\}/C_{\hat{G}}(\varphi)\to \locsys_{{}^cG,F}^{\tame}$ be the corresponding locally closed embedding. Note that it is a schematic morphism of finite tor amplitude. We further assume that $C_{\hat{G}}(\varphi)$ is reductive. The discussion in \Cref{rem: non-frob semisimple parameters} implies that $\varphi$ must be Frobenius semisimple. (Using notations there, if $v\neq 0$, then $v$ itself is in the unipotent radical of $C_M(v)=C_{\hat{G}}(\varphi)$.)
Note that as explained at the end of \Cref{rem: SL2 version of WD parameter}, the converse may not be true. However, if $\varphi$ is essentially discrete, then $C_{\hat{G}}(\varphi)$ is reductive.

On the other hand, let $b\in B(G)$ be a basic element.

We will consider the following functor
\begin{multline}\label{eq: Automorphic to Galois pointwise version}
\bL_{G,b,\varphi}\colon \rep^{\tame}(G_b(F))\xrightarrow{(i_b)_*} \shv^{\tame}(\kot_G)\\ \stackrel{\bL_G^{\tame}}{\cong}  \indcoh(\locsys_{{}^cG,F}^{\tame})\xrightarrow{(i_\varphi)^{\indcoh,*}}\indcoh(\{\varphi\}/C_{\hat{G}}(\varphi))=\rep(C_{\hat{G}}(\varphi)).
\end{multline}

It admits a continuous right adjoint functor 
\begin{multline}\label{eq: Galois to automorphic}
\bL_{G,b,\varphi}^R:\rep(C_{\hat{G}}(\varphi))=\indcoh(\{\varphi\}/C_{\hat{G}}(\varphi))\xrightarrow{(i_\varphi)^{\indcoh}_*} \\
\indcoh(\locsys_{{}^cG,F}^{\tame}) \stackrel{(\bL_G^{\tame})^{-1}}{\cong}\shv^{\tame}(\kot_G)\xrightarrow{ (i_b)^!} \rep(G_b(F)).
\end{multline}

Being a right adjoint functor, $\bL_{G,b,\varphi}^R$ sends admissible objects to admissible objects. In particular, if $r$ is a finite dimensional representation of $C_{\hat{G}}(\varphi)$, defining a vector bundle $\mV_r$ on $\{\varphi\}/C_{\hat{G}}(\varphi)$. Then we have
$\bL_{G,b,\varphi}^R(\mV_r)\in \rep(G_b(F))^{\adm}$.

\begin{proposition}
We have 
\[
(\bL_G^{\tame})^{-1}( (i_\varphi)_*^{\indcoh}\mV_r)\in  (\shv(\kot_G)^{\adm})^{2\rho\mbox{-}e,\heartsuit}.
\] 
In particularly, 
\[
\bL_{G,b,\varphi}^R(\mV_r)\in \rep(G_b(F))^{\heartsuit}\cap \rep(G_b(F))^{\adm}
\] 
is an (honest) depth zero admissible representation of $G_b(F)$.
\end{proposition}
\begin{proof}
By \Cref{lem: adm coh duality discrete par} and \Cref{lem: adm duality and e-t-structure}, it is enough to show that  
\[
(\bL_G^{\tame})^{-1}( (i_\varphi)_*^{\indcoh}\mV_r)\in  (\shv(\kot_G)^{\adm})^{2\rho\mbox{-}e,\geq 0}.
\]

Recall that the collection $\{\cind_{P_b}^{G_b(F)}\varrho\}_{(P_b,\varrho)}$ for $P_b$ a parahoric of $G(F)$ and $\varrho$ a representation of $P$ that obtained by inflation from an irreducible representation of $L_P$, form a set of projective generators of $\rep^{\tame}(G_b(F),\overline{\bQ}_\ell)$.
Therefore, it is enough to show that 
\[
\Hom_{\shv(\kot_G)}((i_b)_*\cind_{P_b}^{G(F)}\varrho[-\langle 2\rho,\nu_b\rangle], (\bL^{\tame}_G)^{-1}((j_\varphi)^{\indcoh}_*\mV_r))\in \Mod_{\La}^{\geq 0}.
\] 
This is equivalent to
\[
\Hom_{\rep(C_{\hat{G}}(\rho))}((i_{\varphi})^{\indcoh,*} \bL_G^{\tame}((i_b)_*\cind_P^{G(F)}\varrho),r)\in \Mod_{\La}^{\geq 0}.
\]
But this follows from \Cref{prop: coh amplitude of star extension of depth zero rep}.
\end{proof}

Now, let  $r_0$ be its restriction to $Z_{\hat{G}}^{\Ga_{\widetilde F/F}}$, which corresponds to an element $\al_r\in \pi_1(G)_{\Ga_F}=\xch(Z_{\hat{G}}^{\Ga_{\widetilde F/F}})$. Let $b\in B(G)$ be the unique basic element which maps to $\al_r$ under the Kottwitz map.
We thus attach every enhanced parameter $(\varphi,r)$ a depth zero admissible representation 
\[
\pi_{(\varphi,r)}:= \bL_{G,b,\varphi}^R(\mV_r).
\]

\begin{remark}
Unfortunately, at the moment we can say very little about $\pi_{(\varphi,r)}$. If $r$ is the trivial representation of $C_{\hat{G}}(\varphi)$, so $\mV_r=\mO_{\{\varphi\}/C_{\hat{G}}(\varphi)}$, then $\Hom(\cind_{I^u}^{G(F)}\psi, \bL(j_*\mV_r))=\overline\bQ_\ell$ so the socle $\pi_{(\varphi,r)}$ contains a unique generic representation.
\end{remark}

\begin{remark}
Despite of the above remark, we have the following  formal consequences about the Harish-Chandra characters of the $L$-packets constructed in the above way.

As the representation $\pi_{(\varphi,r)}\in \rep^{\tame}(G_b(F))^{\heartsuit}$ attached to the enhanced Langlands parameter $(\varphi,r)$ is admissible, it admits a character 
\[
\Theta_{\pi_{(\varphi,r)}}: C_c^\infty(G_b(F))^{\tame}_{G_b(F)}\to \La.
\]
The functor $\bL_{G,b,\varphi}$ induces a map of horizontal traces (or Hochschild homology) of categories
\begin{multline*}
\mathrm{tr}(\bL_{G,b,\varphi}): C_c^\infty(G_b(F))^{\tame}_{G_b(F)}=H^0\mathrm{tr}(\rep^{\tame}(G_b(F)))\xrightarrow{\mathrm{tr}((i_b)_*)}  H^0\mathrm{tr}(\shv(\kot_G))\\
\cong H^0\mathrm{tr}(\indcoh(\locsys_{{}^cG,F}^{\tame}))\xrightarrow{\mathrm{tr}((i_\varphi)^{\indcoh,*})} H^0\mathrm{tr}(\indcoh(\{\varphi\}/C_{\hat{G}}(\varphi)))\cong H^0\rg(\frac{C_{\hat{G}}(\varphi)}{C_{\hat{G}}(\varphi)},\mO).
\end{multline*}
Now for $r\in \rep(C_{\hat{G}}(\varphi))$, with $\Theta_r\in H^0\rg(\frac{C_{\hat{G}}(\varphi)}{C_{\hat{G}}(\varphi)},\mO)$ the usual character of $r$. By \eqref{eq-functoriality-abstract-Character}, we have
\[
\Theta_{\pi_{(\varphi,r)}}=\Theta_r\circ \mathrm{tr}(\bL_{G,b,\varphi}).
\]
\end{remark}

\subsubsection{Regular supercuspidal}\label{SSS: regular supercuspidal}
In this subsection, we set $\La=\overline\bQ_\ell$ and assume that $G$ is unramified. We will fix the pinning $(G,B,T,e)$ as before. Let $A\subset S\subset T$, with the Iwahori-Weyl group $\widetilde{W}$ acting on $\scrA(G_{\breve F},S_{\breve F})$, and let $\breve\bfa\subset \scrA(G_{\breve F},S_{\breve F})$ be the alcove as previously before. To simplify our discussions we also assume that $G$ is of adjoint type.

We fix a regular tame inertia type $\zeta$ as in \Cref{ex: regular supercuspidal} and discuss the corresponding categorical local Langlands correspondence.
In this context, all the categorical and geometric complexities associated with the local Langlands correspondence are significantly reduced. The categorical equivalence simplifies to a classical local Langlands correspondence between the set of isomorphism classes of certain supercuspidal representations of $G$ and its inner forms, and the set of equivalence classes of enhanced tame Langlands parameters whose inertia types are $\zeta$. As we shall see, our bijection coincides with the bijection constructed in \cite{Debacker.Reeder}. 
We will work with $L$-group rather than $C$-group in the sequel.

Recall that we let $\Xi(\zeta)$ be the set of homomorphisms $\chi: I_F^t\to \hat{T}$ corresponding to $\zeta$ under \Cref{lem: classification tame inertia type}. Note that $\zeta$ being regular means that $\Xi(\zeta)$ is a $W_0$-torsor. For each $\chi\in\Xi(\zeta)$, let $w_\chi\in W_0$ denote the unique element such that $\chi^q= w_\chi\bar\sigma(\chi)$. We recall the following crucial fact: the map
\begin{equation}\label{eq: wchisigma elliptic}
1-w_\chi\bar\sigma: \xcoch(T_\ad)_\bQ\to \xcoch(T_\ad)_\bQ
\end{equation}
is an isomorphism. It will be convenient to consider the following set $\xcoch(T)\times \Xi(\zeta)$, equipped an action of $\widetilde{W}=\xcoch(T)\rtimes W_0$ given by
\begin{equation}\label{eq: an action of tilW on xcoch times Xi}
w(\la, \chi)=(w(\la), w(\chi)),\ \ \ \mbox{ for } w\in W_0,\quad t_\mu(\la,\chi)=(\la+(1-w_\chi\bar\sigma)(\mu), \chi),\ \ \ \mbox{ for } t_\mu\in \xcoch(T).
\end{equation}
Note that we have a map 
\begin{equation}\label{eq: from xcochXi to tilW}
\xcoch(T)\times \Xi(\zeta)\mapsto \widetilde{W},\quad (\la, \chi)\mapsto t_\la w_\chi,
\end{equation}
which is equivariant with respect to the $\widetilde{W}$ act on the left as defined in \eqref{eq: an action of tilW on xcoch times Xi}, and the $\sigma$-conjugation action of $\widetilde{W}$ on the right. Passing to the quotient induces an \emph{injective} map
\[
\widetilde{W}\bs(\xcoch(T)\times \Xi(\zeta))\hookrightarrow B(\widetilde{W}).
\]
We need the following observation.
\begin{lemma}\label{lem: finiteness of tlawchi}
Let $w=t_\la w_\chi$ be an element in the image of the above map. Then its Newton point is central, i.e. $\nu_{\dot{w}}\in \xcoch(Z_G)\otimes\bQ$. 
\end{lemma}
\begin{proof}
We have $(w\sigma)^n=t_{\sum_{i=0}^{n-1}(w_\chi\sigma)^i(\la)} (w_\chi\sigma)^n$. Since $(\xcoch(T_\ad)\otimes\bQ)^{w_\chi\sigma}=\{0\}$,  we see that for $n$ sufficiently divisible, $(w_\chi\sigma)^n=1$ and $\sum_{i=0}^{n-1}(w_\chi\sigma)^i(\la)\in\xcoch(Z_G)$. Therefore, $\nu_{\dot{w}}\in \xcoch(Z_G)\otimes\bQ$ as desired.
\end{proof}

Now, we consider categorical local Langlands correspondence. We start with the spectral side. 
Recall from \Cref{ex: regular supercuspidal} that if we choose $\chi\in \Xi(\zeta)$, then there is  an isomorphism
\[
\locsys_{{}^cG,F}^{\hat\zeta}\simeq \{\varphi\}/C_{\hat{G}}(\varphi),\quad C_{\hat{G}}(\varphi)=\hat{T}^{w_\chi\sigma}.
\]
Namely, $\varphi$ is a Langlands parameter such that $\varphi|_{I^t_F}=\chi$. Then for a  lifting of the Frobenius $\sigma\in W^t_F$, we have $\varphi(\sigma)=\dot{w}_\chi\bar\sigma\in \hat{G}\bar\sigma$ for some element $\dot{w}_\chi\in N_{\hat{G}}(\hat{T})$ lifting $w_\chi$. It follows that $C_{\hat{G}}(\varphi)=\hat{T}^{w_\chi\sigma}$.

We also recall that in this case, the correspondence defining the spectral Deligne-Lusztig induction can be described explicitly as in \Cref{ex: steinberg in regular case}. 
We can assign every $(\la, \chi)\in \xcoch(T)\times\Xi(\zeta)$ the object 
\[
\mF_{(\la,\chi)}:=\mO_{\locsys_{{}^cB,\breve F}^{\hchi}}(\la)\star \mO_{S^{\hchi, w^{-1}_\chi(\hchi)}_{{}^cG,\breve F, w_\chi}}\in \Coh(\prod_{\chi_1,\chi_2\in \Xi(\zeta)} S_{{}^cG,\breve F}^{\hchi_1,\hchi_2}).
\]
Then
\[
\mV_r:=\Ch^{\tame}_{{}^cG,\phi}(\mF_{(\la,\chi)})
\]
is a vector bundle on $\{\varphi\}/C_{\hat{G}}(\varphi)$ corresponding to the representation $r$ of $C_{\hat{G}}(\varphi)=\hat{T}^{w_\chi\bar\sigma}$ given by the restriction of the character $\la$ along $\hat{T}^{w_\chi\bar\sigma}\subset \hat{T}$. 
If we let $\mathrm{ELP}_\zeta$ denote the set of equivalence classes of enhanced Langlands parameters $(\varphi, r)$ with inertia type $\zeta$. Then the spectral Deligne-Lusztig induction induces a map
\begin{equation}\label{eq: from xcochXi to enhanced paramenters}
\widetilde{W}\bs(\xcoch(T)\times \Xi(\zeta))\cong \mathrm{ELP}_{\zeta},\quad (\la, \chi)\to (\varphi, r=\la|_{\hat{T}^{w_\chi\bar\sigma}}).
\end{equation}

We have the equivalence of monoidal categories
\begin{equation}\label{eq:generic Bez equivalence}
\bB^{\hat\zeta}: \bigoplus_{\chi_1,\chi_2 \in W_0\chi}\shv\bigl((\iw,\hchi_1)\backslash LG/(\iw,\hchi_2)\bigr)\cong \bigoplus_{\chi_1,\chi_2\in \Xi(\zeta)}\indcoh(S_{{}^cG,\breve F}^{\hchi_1,\hchi_2}).
\end{equation}

The equivalence \eqref{eq:generic Bez equivalence} of course is a very special case of \Cref{thm:Bez-equivalence}. But compared with the unipotent case proved by Bezrukavnikov, it is much easier to establish. The key point is that for every $w\in \widetilde{W}$, we have
\begin{equation}\label{eq: cleanness for regular chi}
\Delta^{\mon}_{\dot{w},\hchi}=\nabla^{\mon}_{\dot{w},\hchi},\quad \forall \chi\in \Xi(\zeta), w\in \widetilde{W}. 
\end{equation}
This implies that for every $\chi, w, w'$, we have
\[
\Delta^{\mon}_{\dot{w},\hchi}\star^u\Delta^{\mon}_{\dot{w'}, w^{-1}\hchi}\cong \Delta^{\mon}_{\dot{w}\dot{w}', \hchi},\quad \nabla^{\mon}_{\dot{w},\hchi}\star^u\nabla^{\mon}_{\dot{w'}, w^{-1}\hchi}\cong \nabla^{\mon}_{\dot{w}\dot{w}', \hchi}.
\]
For example, for $\la\in\xcoch(T)$, we have the cofree $\chi$-monodromic Wakimoto sheaf
$J_{\la,\hchi}^{\mon}=\nabla^{\mon}_{\la_1,\hchi}\star^u \Delta^{\mon}_{-\la_2,\hchi}$,
where we write $\la=\la_1-\la_2$ for both $\la_1,\la_2$ dominant. Then the above isomorphisms imply that $J_{\la,\hchi}^{\mon}=\nabla_{\la,\hchi}^{\mon}$. 

We can now associate every $(\la,\chi)\in \xcoch(T)\times \Xi(\zeta)$ an object
\[
\mG_{(\la,\chi)}=J^{\mon}_{\la,\hchi}\star^u \nabla^{\mon}_{w_\chi,\hchi}=\nabla^{\mon}_{t_\la w_\chi, \hchi}.
\]
Note that under the equivalence \eqref{eq:generic Bez equivalence}, the object $\mG_{(\la,\chi)}$ and $\mF_{(\la,\chi)}$ matches to each other.

Now \eqref{eq:generic Bez equivalence} induces the equivalence
\[
\bL^{\hat\zeta}_G: \shv^{\hat\zeta}(\kot_G)\cong \indcoh(\locsys_{{}^cG,F}^{\hat\zeta})\cong \rep(\hat{T}^{w_\chi\bar\sigma}),
\]
which satisfying
\[
\bL_G^{\tame}(R^{\mon,*}_{t_\la w_\chi,\hchi})= \bL_G^{\tame}(R^{\mon,!}_{t_\la w_\chi,\hchi})\simeq \mV_r.
\]
Here $R^{\mon,*}_{w,\hchi}$ and $R^{\mon,!}_{w,\hchi}$ are defined as in \eqref{eq-affine-DL-induction-chi}, and $r$ is the restriction of $\la$ to $\hat{T}^{w_\chi\bar\sigma}$. 
As $\mV_r$ only depends on $(\la, \chi)$ up to the $\widetilde{W}$-action defined in \eqref{eq: an action of tilW on xcoch times Xi} and as the map \eqref{eq: from xcochXi to tilW} is $\widetilde{W}$-equivariant, we see that
\[
R^{\mon,*}_{t_\la w_\chi,\hchi}\simeq R^{\mon,*}_{v(t_\la w_\chi)\sigma(v)^{-1},v(\hchi)}
\]
for every $v\in \widetilde{W}$. It follows that up to $\sigma$-conjugation (and replacing $\chi$ by an element in $\Xi(\zeta)$), we have $w=t_\la w_\chi$ is of minimal length in its $\sigma$-conjugacy class. 
By \Cref{lem: finiteness of tlawchi}, the Newton point of $w$ is central and the Kottwitz invariant $\kappa_G(w)=[\la]$, where $[\la]$ denote the image of $\la$ under the map $\xcoch(T)\to \pi_1(G)\to \pi_1(G)_\sigma$. Let $b=b_\la$ be the corresponding basic element. Then it follows from \Cref{lem:ADLV-sheaf-minimal-length} that, we have
\[
R^{\mon,*}_{t_\la w_\chi,\hchi}\simeq (i_b)_*\cind_{P_b}^{G_b(F)}(R^{\mon,*,f}_{u,\hchi}),
\]
where $P_b$ is a parahoric subgroup of $G_b(F)$, and $R^{\mon,*,f}_{u,\hchi}$ is a finite Deligne-Lusztig character of the Levi subgroup of $L_{P_b}$. 

Since
\[
\End_{\rep(G_b(F)}(\cind_{P_b}^{G_b(F)}(R^{\mon,f}_{u,\hchi}))\cong \End_{\shv(\kot_G)}(R^{\mon,*}_{t_\la w_\chi,\hchi})\cong \End_{\rep(C_{\hat{G}}(\varphi))}(\mV_r)=\La,
\]
we see that $R^{\mon,*,f}_{u,\hchi}$ must concentrate in cohomological degree zero, and is an irreducible representation. In addition, $\cind_{P_b}^{G_b(F)}(R^{\mon,f}_{u,\hchi})$ must be an irreducible supercuspidal representation of $G_b(F)$. 
Together with other properties of the categorical local Langlands correspondence from \Cref{eq:main-theorem-categorical-tame-local-Langlands}, we arrive at the following theorem.

\begin{theorem}
Let $\zeta$ be a tame regular inertia type, and let $\varphi$ be a unique (up to isomorphism) Langlands parameter such that $\varphi|_{I_F}=\zeta$.
For each $\al\in \pi_1(G)_\sigma=\xch(Z_{\hat{G}})$, let $b$ be the corresponding basic element. 
Let $\rep_\al(C_{\hat{G}}(\varphi))\subset \rep(C_{\hat{G}}(\varphi))$ be the full subcategory consisting of those representations $r$ such that $r|_{Z_{\hat{G}}}=\al$. Let $\rep^{\hat\zeta}(G_b(F),\La)=\rep(G_b(F))\cap \shv^{\hat\zeta}(\kot_G)$.
Then there is an equivalence of categories
\[
\rep_\al(C_{\hat{G}}(\varphi))\to \rep^{\hat\zeta}(G_b(F),\La).
\]
The functor sends an irreducible representation of $C_{\hat{G}}(\varphi)$ to a supercuspidal representation of $G_b(F)$. When $\al=1$, the functor sends the trivial representation of $C_{\hat{G}}(\varphi)$ to the supercuspidal representation of $G(F)$ that admits a Whittaker model.
\end{theorem}

Next, we show that the above local Langlands correspondence for the tame inertia type $\zeta$ coincides with the one constructed in \cite{Debacker.Reeder}. This amounts to understanding $R^{\mon,*,f}_{u,\hchi}$ more explicitly.

\begin{lemma}
Let $w=t_\la w_\chi$ be as above. Then there is a unique point $x$ in $\scrA(G_{\breve F}, S_{\breve F})$ such that $w\sigma(x)=x$. In addition, $x$ is a vertex.
\end{lemma}
\begin{proof}
This follows from the fact that \eqref{eq: wchisigma elliptic} is an isomorphism. See \cite[Lemma 4.4.1]{Debacker.Reeder} for details.
\end{proof}

If follows that if $w=t_\la w_\chi$ is a minimal length element in its $\sigma$-conjugacy class as in \Cref{thm: reduction to min length elements}, then the corresponding point $x$ in the above lemma must be the standard facet $\breve\bff$ in \Cref{thm: reduction to min length elements}. In addition, if we write $w=uy$ (here to avoid notation confliction we use $y$ to denote the corresponding $\sigma$-straight element in  \Cref{thm: reduction to min length elements}), then $y$ must be of length zero, and $u$ is elliptic in $W_{\bff}$. It follows that $R^{\mon,*,f}_{u,\hchi}\simeq (i_b)_*\cind_{P_b}^{G_b(F)}(R^*_{\dot{u},\theta})$, where $P_b=\mP_{\breve\bff}^{\dot{y}\sigma}(\breve\mO)$ is a maximal parahoric of $G_b(F)$. The torus $T$ equipped with the Frobenius structure $\sigma_y=\Ad_{\dot{y}}\sigma$ transfers to a maximal torus of the Levi quotient of $\mP_{\breve\bff}$ (equipped with the same Frobenius structure), and
$\theta$ is the character of $T^{w_\chi\sigma}=T^{u \sigma_y}$. It follows that the supercuspidal representation
\[
(i_b)^!R^{\mon,*,f}_{u,\hchi}=\cind_{P_b}^{G_b(F)}R_{\dot{u},\theta}^*
\]
of $G_b(F)$ indeed coincides with the one constructed in \cite[\textsection{4}]{Debacker.Reeder}.

\begin{remark}
In \cite[Lemma 4.5.1]{Debacker.Reeder}, there is an argument showing that $\cind_{P_b}^{G_b(F)}R_{\dot{u},\theta}^*$ is a supercuspidal irreducible representation. But this follows from our categorical equivalence.
\end{remark}

\begin{remark}
For non-singular inertia type $\zeta$ (in the sense of \Cref{ex: regular supercuspidal}), the geometry of $\locsys_{{}^cG,F}^{\hat\zeta}$ is still relatively easy to understand. In addition, the corresponding monodromic affine Hecke category is easy. For example, \eqref{eq:generic Bez equivalence} is still easy to establish directly and \eqref{eq: cleanness for regular chi} continues to hold. It should not be difficult to generalize the above discussions to this case.
\end{remark}

\newpage

\section{Local-global compatibility}\label{sec: local-global-comp}
In this section, we give some first global applications of the unipotent categorical local Langlands correspondence. 

\subsection{Cohomology of Shimura varieties}\label{SS: coh of shimura varieties}
\subsubsection{The categorical local Langlands for non quasi-split group}
We fix once for all a non-zero additive character $\Psi: F\to \La^\times $ of conductor $\mO_F$ as before.

Note that input data of the categorical local Langlands correspondence is a quasi-split reductive group over a non-archimidean local field $F$ equipped with a pinning.
However, in various  applications, one usually starts with a not necessary quasi-split reductive group. Therefore, we need to explain how to extend the correspondence to non quasi-split groups, after making some auxiliary choices. 

Let $G$ be a connected reductive group over a non-archimidean local field $F$. In  \cite[\textsection{4.2}]{zhu2020coherent}, we attach $G$ a groupoid $\mathbf{TS}_G$ of $G$. Choosing $t\in \mathbf{TS}_G$ amounts
to choosing 
\begin{itemize}
\item a pinned quasi-split group $(G^*, B^*, T^*, e^*)$ over $F$;
\item an isomorphism $\eta: G_{\breve F}\cong G^*_{\breve F}$ and an element $b\in G^*(\breve F)$ such that $\eta \sigma \eta^{-1}=\Ad_{b}  \sigma^* $;
\end{itemize}
By $\sigma$-conjugating $b$ by an element in $G^*_{\breve F}$, we may assume that $b=\dot{w}\in N_{G^*}(S^*)$ that also normalizes the Iwahori $\mI^*$ of $G^*$ determined by the pinning. 

Recall that as in \Cref{SS: Iwahori-Weyl-group}, the pinning of $G^*$ determines $A^*\subset S^*\subset T^*$ as well as an alcove $\breve\bfa^*$ of $\scrA(G^*_{\breve F},S^*_{\breve F})$.
Let $\widetilde{W}$ denote the Iwahori-Weyl group of $G^*_{\breve F}$. Then $b=\dot{w}$ for some $w\in \widetilde{W}$ is a length zero element. The tori $S^*\subset T^*$ transfer to $S\subset T\subset G$, as well as the alcove $\breve\bfa\subset \scrA(G_{\breve F},S_{\breve F})$. Note that we have 
\[
\scrA(G,A)=\scrA(G_{\breve F},S_{\breve F})^\sigma\cong \scrA(G^*_{\breve F},S^*_{\breve F})^{w\sigma^*}.
\]

As explained in  \Cref{rem: automorphism of sigma action}, we have an isomorphism
\[
\eta_w: \kot_G:=\frac{LG}{\Ad_\sigma LG}\cong \kot_{G^*}=\frac{LG^*}{\Ad_{\sigma^*}LG^*},\quad g\mapsto \eta(g)\dot{w}.
\]
This map induces a bijection $\eta_w: B(G)=B(G^*)$. Note that this map does not match the Kottwitz invariants nor the Newton points.

In any case, once we fix such $(G^*,B^*,T^*,e^*,\eta, b=\dot{w})$, we thus obtain a fully faithful embedding
\[
\fgshv^{\unip}(\kot_G,\La)\cong \fgshv^{\unip}(\kot_{G^*},\La)\stackrel{\bL_{G^*}^{\unip}}{\hookrightarrow} \Coh(\locsys_{{}^cG,F}^{\tame}\otimes \La).
\] 
We will denote by $\bL_G^{\unip}$ the composed embedding.

\subsubsection{Recollection of mod $p$ geometry of Shimura varieties}
Let $(G,X)$ be a Shimura datum. We fix a prime $p$ and let $K_p$ be a parahoric subgroup of $G(\bQ_p)$. 
Let $K^p\subset G(\bA_f^p)$ be an open compact subgroup away from $p$ that is sufficiently small.
Let $K=K_pK^p$. Let $\bfSh_K(G,X)$ be the associated Shimura variety defined over the reflex field $E=E(G,X)\subset \bC$. 
In the sequel, we will fix an embedding $\iota: E\to \overline\bQ_p$. This determines a place
$v$ of $E$ above $p$. Let $\mO_{E,(v)}$ be the localization of $\mO_E$ at $v$. Let $k$ be the residue field of $\overline\bQ_p$, which is an algebraic closure of $\bF_p$. The map $\iota$ induces a map $\mO_{E,(v)}\to k$.

Now assume that $(G,X)$ is  of abelian type. 
It is by now well-known that there is an (canonical) integral model $\scrS_K(G,X)$ of $\bfSh_K(G,X)$ over $\mO_{E,(v)}$, at least when $p>3$. (See \cite{DY,KPZ}.)
The integral model is canonical in a precise sense, uniquely determined by a list of properties it should satisfy. 
We shall not review all of them, but only mention  in \Cref{thm: integral model of Shimura} some of them that are relevant to our applications.
We let 
\[
\Sh_K(G,X):= (\scrS_K(G,X)\otimes_{\mO_{E,(v)}} k)^{\pf},
\]
be the perfection of the special fiber of the integral model $\scrS_K(G,X)$. If $(G,X)$ is clear from the context, we simply denote $\Sh_K(G,X)$ by $\Sh_K$.

We will let $\mP_{\breve\bff}$ be a standard parahoric of $G_{\breve\bQ_p}$, corresponding to a facet $\breve\bff\subset\overline{\breve\bfa}$. When $\bff=\breve\bff\cap \scrA(G,A)$ is a standard facet, then $\mP_{\breve\bff}$ is a parahoric group scheme of $G$ defined over $\bZ_p$, denoted as $\mP$. We will write $K_{p,\bff}=\mP_{\breve\bff}(\breve\bZ_p)\cap G(\bQ_p)=\mP(\bZ_p)$.

As before, let $(\hat{G},\hat{B},\hat{T},\hat{e})$ be the dual group of $G_{\bQ_p}$ equipped with a pinning, defined over $\bZ$.
Let $\mu\in \xch(\hat{T})^+$ be the minuscule dominant character associated to the Shimura datum $X$. As usual, we let $\mu^*=-w_0(\mu)$, where $w_0$ is the longest length element in the Weyl group of $(\hat{G},\hat{T})$.
We let
\[
\adm(\mu^*)=\{w\in \widetilde{W}\mid w\leq t_{\bar\la} \mbox{ for some } \bar\la \in W_0\bar{\mu^*}\}
\] 
be the admissible set associated to $\mu^*$. Here  for $\la\in \xcoch(T)$, we let $\bar\la$ denote its image in $\xcoch(T)_{I_F}$ and $t_{\bar\la}$ the translation element in $\widetilde{W}$ given by $\xcoch(T)_{I_F}\subset \widetilde{W}$, and
$W_0\bar{\mu^*}\subset \xcoch(T)_{I_F}$ denotes the $W_0$-orbit of $\bar{\mu^*}$ in $\xcoch(T)_{I_F}$.
Let 
\[
\adm^{\breve\bff}(\mu^*)= W_{\breve\bff}\adm(\mu^*)W_{\breve\bff}\subset \widetilde{W}
\]
be the parahoric version. We let
\[
LG_{\mP,\mu^*}:=\cup_{w\in \adm^{\breve\bff}(\mu^*)}LG_w,
\]
which is a closed subset of $LG$. 
We let $\Sht^{\loc}_\mP$ be the moduli of local Shtukas for $\mP$.

Let
\[
\Sht^{\loc}_{\mP,\mu^*}=\frac{LG_{\mP,\mu^*}}{\Ad_\sigma L^+\mP}\subset \Sht^{\loc}_{\mP,\mu^*}.
\]
We again recall that after identification of $G$ and $G^*$ over $\breve\bQ_p$, the Frobenius $\sigma$ becomes $\Ad_{\dot{w}}\sigma^*$.

We summarize the facts about the integral model $\scrS_{K}(G,X)$ we need. We thank Michael Harris for urging us to extract precisely the properties of integral models we need.
\begin{assumption}\label{thm: integral model of Shimura}
The canonical integral model $\scrS_{K}(G,X)$ of $\bfSh_K(G,X)$ satisfies the following properties.
\begin{enumerate}
\item\label{thm: integral model of Shimura-proper} If $\bfSh_K(G,X)$ is proper over $E$, then $\scrS_K(G,X)$ is proper over $\mO_{E,(v)}$.
\item\label{thm: integral model of Shimura-tamelevel} If $K'=K_p(K^p)'\subset K=K_pK^p$ where $(K^p)'\subset K^p$ is an prime-to-$p$ open subgroup, then $\bfSh_{K'}(G,X)\to \bfS_K(G,X)$ is finite \'etale.
\item\label{thm: integral model of Shimura-smooth} If $\mP$ is reductive (so $K_p$ is hyperspecial), then $\scrS_K(G,X)$ is smooth. 
\item\label{thm: integral model of Shimura-nearby} The canonical morphism $\rg(\bfSh_K(G,X)_{\overline\bQ_p},\La)\to \rg(\Sh_K(G,X),R\Psi)$ is an isomorphism, where $R\Psi$ denotes the sheaf of nearby cycles  of $\scrS_K(G,X)$.
\item\label{thm: integral model of Shimura-locSht} There is a morphism 
\begin{equation}\label{eq: from abv to sht}
\loc_p: \Sh_K(G,X)\to \Sht_{\mP,\mu^*}^\loc,
\end{equation}
such that the composed morphism
\[
\loc_p(m,n): \Sh_K(G,X)\to \Sht_{\mP,\mu^*}^{\loc(m,n)},
\]
is coh. smooth. Here $\Sht_{\mP,\mu^*}^{\loc(m,n)}$ is defined as in \eqref{eq:restr local Sht}.
In addition, if $\bff\subset\overline{\bff'}\subset \overline{\bfa}$, the following diagram is Cartesian
\[
\xymatrix{
\Sh_{K'}(G,X)\ar[r]\ar[d] & \Sht_{\mP',\mu^*}^{\loc}\ar[d] \\
\Sh_{K}(G,X) \ar[r] & \Sht_{\mP,\mu^*}^{\loc}.
}\]
\item\label{thm: integral model of Shimura-central} The sheaf $(\loc_p)^!\delta^!\mZ_{\mu}$ is canonically isomorphic to $R\Psi[d]$, where $d=\dim\bfSh_K(G,X)$, where $\mZ_\mu$ is the central sheaf on $L^+\mP\bs LG/L^+\mP$ corresponding to the irreducible representation of $\hat{G}$ of highest weight $\mu$, and $\delta: \Sht^{\loc}_\mP\to L^+\mP\bs LG/L^+\mP$ is the morphism as in \Cref{rem: modification direction of Shtuka}.
\item\label{thm: integral model of Shimura-Hecke} Let $\Hk_\bullet(\Sht_{\mP,\mu^*})$ denote the restriction of the Hecke groupoid $\Hk_\bullet(\Sht_{\mP})$ from \eqref{eq: nth local Hecke stack for Sht} to $\Sht_{\mP,\mu^*}$. Equivalently, $\Hk_\bullet(\Sht_{\mP,\mu^*})$ is the \v{C}ech nerve of $\Sht_{\mP,\mu^*}\to \kot_G$. Then $\Hk_\bullet(\Sht_{\mP,\mu^*})$ pullbacks back to a groupoid over $\Sh_K(G,X)$ under the map $\loc_p$.
\item\label{thm: integral model of Shimura-Igs} The partial minimal compactification of the Igusa variety $\mathrm{Ig}_x$ (as review below) is affine. In particular, if $\bfSh_K(G,X)$ is proper over $E$, then $\mathrm{Ig}_x$ is affine.
\end{enumerate}
\end{assumption}

Note that given \Cref{rem: modification direction of Shtuka}, the appearance of $\Sht_{\mP,\mu^*}^{\loc}$ here is in fact consistent with the appearance of $\Sht^{\loc}_{\mP,\mu}$ in \cite{xiao2017cycles}.

\begin{lemma}\label{lem: smoothness of map from Sh to Sht}
The morphism $\loc_p$ in \eqref{eq: from abv to sht} is representable pseudo coh. pro-smooth in the sense of  \Cref{def-coh-pro-smooth-morphism-space}. In addition, let $\La^{\can}$ be the canonical generalized constant sheaf on $\Sht_{\mP,\mu^*}^{\loc}$ as in \Cref{SS: coh. duality. Kot. stack}, then $(\loc_p)^!\La^{\can}$ is isomorphic to the constant sheaf of $\Sh_K(G,X)$.
\end{lemma}
\begin{proof}
This follows from \Cref{ex: pseudo pro-smooth morphism}. More precisely, although $f:X\to Y$ is assumed to be a morphism of algebraic spaces there, the arguments work without change in the current setting.
\end{proof}

As explained in \cite{zhu2020coherent}, let $x: \Spec k\to \Sht_{\mP,\mu^*}$ be a point. Then
\[
\mathrm{Ig}_x:= \Spec k\times_{\Sht_{\mP,\mu^*}^\loc}\Sh_K(G,X)
\]
is a pro-\'etale cover of the central leaf $C_x\subset \Sh_K(G,X)$. It is known that $C_x$ is perfectly smooth, of dimension $\langle 2\rho, \nu_b\rangle$.

\begin{theorem}
All the above assumptions hold for $\scrS_K(G,X)$, when $(G,X)$ is a Shimura datum of Hodge type (and $p>2$).
\end{theorem}

\subsubsection{The category of sheaves on the perfect Igusa stack}\label{SSS: Igusa stack}

We need the following geometric input.
\begin{proposition}\label{prop: perfect Igusa stack}
There is a stack $\Igs_{K^p}(G,X)\in \prestk_k^\pf$ over $k$ making the following diagram Cartesian
\[
\xymatrix{
\Sh_K(G,X) \ar^-{\loc_p}[rr] \ar_{\Nt_\mP^{\glob}}[d] && \Sht_{\mP,\mu^*}^{\loc} \ar^{\Nt_{\mP,\mu}}[d]\\
\Igs_{K^p}(G,X) \ar^-{\loc_p^0}[rr] &&  \kot_G 
}\]
\end{proposition}
\begin{proof}
Let $\Hk_\bullet(\Sht^{\loc}_\mP)$ be the Hecke groupoid for $\Sht^{\loc}_\mP$ whose $n$th term is given in \eqref{eq: nth local Hecke stack for Sht}. As explained in \Cref{prop: shv on kot via sht}, this can also be regarded as the \v{C}ech nerve of the morphism $\Nt_\mP: \Sht^{\loc}_\mP\to \kot_G$.

Now consider $\Sht^{\loc}_{\mP,\mu^*}\to \kot_G$, and we have the corresponding \v{C}ech nerve $\Hk_\bullet(\Sht^{\loc}_{\mP,\mu^*})$. E.g. 
\[
\Hk(\Sht^{\loc}_{\mP})_{\mu^*\mid\mu^*}:=\Hk_1(\Sht^{\loc}_{\mP,\mu^*})
\] 
can be described as the perfect prestack over $k$ sending a perfect $k$-algebra
$R$ to the groupoid of triples $((\mE_1,\varphi_1), (\mE_2,\varphi_2), \beta)$ where $(\mE_i,\varphi_i)\in \Sht^{\loc}_{\mP,\mu^*}$, and $\beta: \mE_1\dashrightarrow \mE_2$ is a modification compatible with $\varphi_i$.

As \Cref{thm: integral model of Shimura} \eqref{thm: integral model of Shimura-Hecke} holds in our case, this groupoid pullbacks back to a groupoid 
\[
\hk_\bullet(\Sh_K(G,X)):=\Sh_K(G,X)\times_{\Sht_{\mP,\mu^*}^{\loc}}\hk_\bullet(\Sht_{\mP,\mu^*}^{\loc}).
\]
More precisely, there is the following commutative diagram with both squares Cartesian
\[
\xymatrix{
    \Sh_K(G,X)\ar[d]& \ar[l] \hk(\Sh_K(G,X))  \ar[r]\ar[d] & \Sh_K(G,X)\ar[d]\\
    \Sht^{\loc}_{\mP,\mu^*}  & \ar[l] \hk(\Sht^{\loc}_{\mP})_{\mu\mid\mu^*} \ar[r] & \Sht^{\loc}_{\mP,\mu^*}.
}
\]
Then the $n$-term 
\[
\Hk_n(\Sh_K(G,X))=\Hk(\Sh_K(G,X))\times_{\Sh_K(G,X)}\Hk(\Sh_K(G,X))\times_{\Sh_K(G,X)}\cdots\times_{\Sh_K(G,X)} \Hk(\Sh_K(G,X))
\] 
is the $n$-folded product, with the face and boundary maps defined naturally.
Then we define $\Igs_{K^p}(G,X)$ as the \'etale sheafification of the geometric realization of $\hk_\bullet(\Sh_K(G,X))$. 
\end{proof}
\begin{remark}
The stack $\Igs_{K^p}(G,X)$ constructed in \Cref{prop: perfect Igusa stack} is usually called the perfect Igusa stack. Modulo the difference between \'etale sheafification and $h$-sheafification, the above result is in fact a consequence of a much more difficult result on the existence of the Igusa stack as a $v$-stack over $\Spd \bF_p$, as proved in \cite{DHKZ}. In addition, the authors constructed a similar diagram of $v$-stacks with the first row replaced by objects defined over $\Spd \mO_E$, and with $\kot_G$ replaced by the $v$-stack $\Bun_G$ of moduli of $G$-bundles on the Fargues-Fontaine curve. They also explained that the Cartesian diagram in \Cref{prop: perfect Igusa stack} can be obtained by reduction of their Cartesian diagram of $v$-stacks.
\end{remark}

\begin{remark}
It follows by construction that $\Igs_{K^p}(G,X)$ is a quasi-compact sind-very-placid stack, with $\Sh_K(G,X)\to \Igs_{K^p}(G,X)$ a sind-placid atlas. 
The stack $\Igs_{K^p}(G,X)$ is in fact independent of the choice of the level $K_p$.
\end{remark}

As usual, associated to $\mu^*$ there is a finite subset $B(G,\mu^*)\subset B(G)$ consisting of those $b$ such that $b\leq b_{\bar{\mu^*}}$. We let 
\[
\kot_{G,\leq \mu^*}=\cup_{b\in B(G, \mu^*)} \kot_{G,b}.
\]
This is a connected closed substack of $\kot_G$. Clearly, $\loc_p^0$ factors as
\[
\Igs_{K^p}(G,X) \xrightarrow{\loc_p^0} \kot_{G,\leq \mu^*}\subset \kot_G. 
\]

In the sequel, we will omit $(G,X)$ from the notations. E.g. we will write $\Sh_K$ and $\Igs_{K^p}$ instead of $\Sh_K(G,X)$ and $\Igs_{K^p}(G,X)$, etc. 
We first discuss the category of  sheaves on $\Igs_{K^p}$. We fix the coefficient ring $\La$ to be a $\bZ_\ell$-algebra as in \Cref{sec:adic-formalism} as before, but omit it from the notation if it is clear from the context.

\begin{lemma}\label{lem: cat of shv on Igs}
The category $\shv(\Igs_{K^p})$ is compactly generated, with $\shv(\Igs_{K^p})^\cpt$ generated (as idempotent complete category) by objects of the form $\Nt_! \mF$, where $\mF\in \cshv(\Sh_K)$. The dualizing sheaf $\consdual_{\Igs_{K^p}}\in \shv(\Igs_{K^p})^{\adm}$.
\end{lemma}
\begin{proof}
It follows from \eqref{eq-shv-on-sifted-placid} that we have
\begin{equation*}\label{eq: shv on Igusa stack}
\shv(\Igs_{K^p})=|\shv(\hk_\bullet(\Sh_K))|
\end{equation*}
with transitioning functors being $*$-pushforwards. 
The first statement follows from this colimit presentation and the fact that since $\Sh_K$ is a perfect scheme pfp over $k$, we have 
$\shv(\Sh_K)=\ind\cshv(\Sh_K)$ by definition. For the second statement, notice that for every $\mF\in\cshv(\Sh_K)$, we have
\[
\Hom(\Nt_*\mF, \consdual_{\Igs_{K^p}})=\Hom(\mF, \consdual_{\Sh_K(G,X)})\in \Perf_\La.
\]
Then the claim follows from \Cref{lem: char of adm obj in cg cat}.
\end{proof}

We can repeat the construction of the canonical self-duality of $\shv(\kot_G)$ as in  \Cref{SS: coh. duality. Kot. stack},

Let $\La_{\Sh_K}^{\can}\in \cshv(\Sh_K)$ be the constant sheaf on $\Sh_K$, i.e. the $*$-pullback of $\consdual_{\Spec k}$ along the structural map $\Sh_K\to \Spec k$. 
Arguing as in \Cref{SS: coh. duality. Kot. stack}, we have a compatible system of generalized constant sheaves $\La_{\hk_\bullet(\Sh_K)}$. 
Then we have a compatible system of functors
\[
\rg^{\can}(\Hk_\bullet(\Sh_K), -): \shv(\Hk_\bullet(\Sh_K))\to \Mod_\La
\]
defining self dualities 
\[
\verd^{\can}_{\Hk_\bullet(\Sh_K)}: \shv(\Hk_\bullet(\Sh_K))^\vee\cong \shv(\Hk_\bullet(\Sh_K)),
\]
which restricts to anti-involutions $(\verd^{\can}_{\Hk_\bullet(\Sh_K)})^\cpt$ on the subcategories of compact objects. We also note that $(\verd^{\can}_{\Sh_K})^\cpt$ is just the usual Verdier duality on $\Sh_K$. 

As in \Cref{SS: coh. duality. Kot. stack}, the above functors together then induce
\begin{equation}\label{eq: can global section Igs}
\rg^{\can}(\Igs_{K^p},-): \shv(\Igs_{K^p})\to \Mod_\La,
\end{equation}
which then induces the self-duality
\[
\verd^{\can}_{\Igs_{K^p}}:  \shv(\Igs_{K^p})^\vee\cong  \shv(\Igs_{K^p})
\]
which then restricts to an anti-involution on compact objects $(\verd^{\can}_{\Igs_{K^p}})^\cpt$.

We have a canonical isomorphism
\begin{equation}
(\loc_p)^!\La^{\can}_{\Sht^\loc_{\mP,\mu^*}}\cong \La_{\Sh_K}.
\end{equation}

On the other hand, let $\verd_{\kot_G}^{\can}$ be the canonical self-duality of $\shv(\kot_G)$ as constructed in \Cref{cohomological.duality.kottwitz}.
\begin{lemma}\label{lem: compatibility of canonical duality between Igs and kot}
We have
\[
(\loc_p^0)^!\circ (\verd_{\kot_G}^{\can})^{\cpt}\cong (\verd_{\Igs_{K^p}}^{\can})^{\cpt} \circ (\loc_p^0)^!.
\]
\end{lemma}
\begin{proof}
First notice that the $!$-pullback along $\loc_p$ of the canonical generalized constant sheaf $\La_{\Sht^{\loc}_{\mP,\mu^*}}$ is the constant sheaf of $\Sh_K$. 
As $\loc_p$ is pseudo-coh. pro smooth, this implies that 
\[
(\loc_p)^!\circ (\verd_{\Sht^{\loc}_{\mP,\mu^*}}^{\can})^{\cpt}\cong (\verd_{\Sh_K}^{\mathrm{verd}})^{c} \circ (\loc_p)^!.
\]
where $(\verd^{\mathrm{verd}}_{\Sh_K})^c$ denotes the usual Verdier duality functor for $\Sh_K$.
See \Cref{rem:verdier.functioriality.spaces.pseudo.coh.pro.smooth} and \Cref{rem:verdier.functioriality.spaces.pseudo.coh.pro.smooth-2}.
This continues to hold at each level of the \v{C}ech nerve, giving the desired statement.
\end{proof}

\begin{proposition}\label{prop: Igusa sheaf}
The functor $(\loc_p^0)^!: \shv(\kot_G)\to \shv(\Igs_{K^p})$ admits a continuous right adjoint $(\loc_p^0)_\flat$.  The object
\begin{equation*}\label{eq: Igusa sheaf}
\Igss_{K^p}:= (\loc_p^0)_{\flat} \consdual_{\Igs_{K^p}}.
\end{equation*}
belongs in $ \shv(\kot_G)^{\adm}$.
\end{proposition}
We call the sheaf $\Igss_{K^p}$ the Igusa sheaf of $(G,X, K^p)$.
\begin{proof}
As $(\loc_p^0)^!: \shv(\kot_G)\to \shv(\Igs_{K^p})$  sends compact objects to compact objects, the first statement follows.
Since $ \consdual_{\Igs_{K^p}}\in \shv(\Igs_{K^p})^{\adm}$ by \Cref{lem: cat of shv on Igs} and $ (\loc_p^0)_{\flat}$ is continuous admitting left adjoint, we see that $\Igss_{K^p}$ is admissible (by \Cref{lem:basic cpt and adm}). Alternatively, the admissibility of $\Igss_{K^p}$ can also be deduced from \Cref{prop: right adjointability from Sh to Igs} below.
\end{proof}

\begin{remark}
As $\loc_p^0$ factors through $\kot_{G,\leq \mu^*}$ we see that $\Igss_{K^p}$ is in fact the $\flat$-pushforward of an object in $\shv(\kot_{G,\leq \mu^*})$. By abuse of notations, we also use $\Igss_{K^p}$ to denote this object in $\shv(\kot_{G,\leq \mu^*})$.
\end{remark}

\begin{remark}\label{rem: canonical triviality of Frobenius endo-3}
We note that the Igusa stack $\Igs_{K^p}$ is in fact defined over $k_E$, the residual field of $\mO_{E,(v)}$. Therefore, it admits a $q_v$-Frobenius endomorphism $\phi$, which in turn induces an auto-equivalence $\phi_*: \shv(\Igs_{K^p})\to \shv(\Igs_{K^p})$.
The same argument as in \Cref{rem: canonical triviality of Frobenius endo} shows that $\phi_*$ is canonically isomorphic to the identity functor. See also \cite[Proposition 5.2.5]{DHKZ}.
\end{remark}

\subsubsection{Coherent description of cohomology of Shimura varieties}

We need the following strengthening of the first part of \Cref{prop: Igusa sheaf}.
\begin{proposition}\label{prop: right adjointability from Sh to Igs}
The following diagram is right adjointable (in $\lincat_\La$)
\[
\xymatrix{
\shv(\kot_{G,\leq \mu^*})\ar^-{(\loc_p^0)^!}[rr] \ar_{(\Nt_\mP)^!}[d]&& \shv(\mathrm{Igs}_{K^p})\ar^{(\Nt^{\glob})^!}[d] \\
\shv( \Sht_{\mP,\mu^*}^{\loc}) \ar^-{(\loc_p)^!}[rr] && \shv(\Sh_K).
}\]
Consequently, for any $\mF\in \cshv(\Sht_{\mP,\mu^*}^{\loc})$, we have
\[
C(\Sh_K, (\verd^{\mathrm{verd}}_{\Sh_K})^c((\loc_p)^!\mF))\cong \Hom(\Nt_*\mF, \Igss_{K^p}).
\]
\end{proposition}

\begin{proof}
For the first statement, by \Cref{prop:categorical-right-adjointability-colimits}  it is enough to show that for every $n\geq 0$ and $0\leq i\leq n$, the following commutative diagram is right adjointable
\[
\xymatrix{
\shv( \Sht_{\mP,\mu^*}^{\loc})\ar^-{(\loc_p)^!}[rr] \ar_{(d_i)^!}[d]&& \shv(\Sh_K)\ar_{(d_i)^!}[d] \\
\shv(\Hk_n(\Sht_{\mP,\mu^*}^{\loc})) \ar^-{(\loc_{p,n})^!}[rr] && \shv(\Hk_n(\Sh_K)).
}\]
By \Cref{thm: integral model of Shimura} \eqref{thm: integral model of Shimura-locSht}, $\loc_p$  and $\loc_{p,n}$ are representable pseudo coh. pro-smooth and therefore belong to the class of morphisms $\mathrm{HR}$ associated to the sheaf theory $\shv$ by \Cref{cor-additional-base-change-for-shv-theory-on-prestacks}. Therefore, the above diagram is right adjointable.

The last statement follows from
\begin{align*}
\Hom(\Nt_*\mF, \Igss_{K^p})&=\Hom(\mF, \Nt^!\Igss_{K^p}) =\Hom(\mF, (\loc_p)_\flat\consdual_{\Sh_K}) \\
                                                                                         &=\Hom((\loc_p)^!\mF, \consdual_{\Sh_K}) =\Hom((\loc_p)^!\mF, (\pi_{\Sh_K})^!\La)\\                                                                                   
                                                                                         &=C_c(\Sh_K, (\loc_p)^!\mF)^\vee = C(\Sh_K,(\verd^{\mathrm{verd}}_{\Sh_K})^c((\loc_p)^!\mF)).
\end{align*}
\end{proof}

\begin{proposition}\label{prop: coh of shimura var via coh sheaf}
Let $\mZ_\mu:=\mZ(V_{\mu})\in\fgshv(\iw\bs LG/\iw)$ be the central sheaf corresponding to $\mu$. Then
\[
C(\Sh_K,\La[d])\cong C_c(\Sh_K,\La[d])^\vee\cong \Hom(\Nt_*\delta^!\mZ_{\mu}, \Igss).
\]
\end{proposition}

Next, suppose $\mP=\mI$ is the standard Iwahori.
Let $b\in B(G,\mu^*)$ and let $w\in \adm(\mu^*)\subset \widetilde{W}$ be a $\sigma$-straight element corresponding to $b$ under the map \eqref{eq: str-vs-sigma-conj}. Let $\dot{w}$ be a lifting of $w$.
Recall from \Cref{prop: Sht-loc-w-straight} that $\Sht^{\loc}_w\cong \bB_{\mathrm{proket}} I_b$, and $\dot{w}\to \Sht^{\loc}_w$ is the universal $I_b$-torsor. Recall 
\[
\widetilde\Igs_{K^p,\dot{w}}:=\Sh_K\times_{\Sht^\loc_{\mu^*}}\dot{w}
\]
is the Igusa variety at the infinite level, equipped with an action of $G_b(F)$. For an open compact subgroup $H\subset G_b(F)$, let 
\[
\Igs_{HK^p}=\widetilde\Igs_{K^p,\dot{w}}/H,
\] 
which is a perfect scheme of pfp over $k$.
Now if $K_b\subset I_b$, we have a representation $\ind_{K_b}^{I_b}\La=C(I_b/K_b,\La)$ of $I_b$ on finite projective $\La$-module, regarded as a sheaf on $\bB I_b$. 

\begin{proposition}\label{prop: coh of Igus var}
For $\mF=(i_w)_! \ind_{K_b}^{I_b}\La$, we have
\[
C(\Igs_{K_bK^p},\La)\cong \Hom((i_b)_!\cind_{K_b}^{G_b(F)}\La,\Igss_{K^p}).
\]
\end{proposition}

Recall we have a fully faithful embedding $\fgshv(\kot_G)\to \shv(\kot_G)$ extending to a continuous functor $\Psi: \rshv(\kot_G)\to \shv(\kot_G)$ (see \eqref{eq-ind-constr-to-shv-sifted}), which admits a fully faithful left adjoint $\Psi^L$, by \Cref{rem: Shv on Kot as colocalization}. In addition, when restricted to $\kot_{G,\leq \mu^*}$ we have an equivalence by \Cref{rem: perverse t-structre on fgshvekotG}.
\begin{equation}\label{eq: fully faithfulness of Psi in eventually coconn part of Sh section}
\Psi^{L}: \shv(\kot_{G,\leq \mu^*})^{2\rho\mbox{-}p,+}\cong \rshv(\kot_{G,\leq \mu^*})^{2\rho\mbox{-}p,+}: \Psi. 
\end{equation}
The following corollary is observed by Xiangqian Yang.
\begin{corollary}
We regard $\Igss_{K^p}$ as an object in $\Shv(\kot_{G,\leq \mu^*})$. Then $\Igss_{K^p}\in \Shv(\kot_{G,\leq \mu^*})^{2\rho\mbox{-}p, +}$.
In particular, for every $\mF\in \fgshv(\kot_G)$, we have
\[
\Hom_{\shv(\kot_G)}(\mF, \Igss_{K^p})=\Hom_{\rshv(\kot_G)}(\mF, \Psi^L(\Igss_{K^p})).
\]
\end{corollary}
\begin{proof}
By \Cref{prop: coh of Igus var}, $\Hom((i_b)_!\cind_{K_b}^{G_b(F)}\La,\Igss_{K^p})\in \Mod^{\geq 0}_\La$ for every $b$ and every pro-$p$-open compact subgroup. Therefore, $\Igss_{K^p}\in \Shv(\kot_{G,\leq \mu^*})^{2\rho\mbox{-}p, +}$.
We have 
\begin{eqnarray*}
\Hom_{\shv(\kot_G)}(\mF, \Igss_{K^p}) &\cong & \Hom_{\shv(\kot_G,\leq \mu^*)}(\mF, \Igss_{K^p})\\ 
                                                              & \cong & \Hom_{\shv(\kot_G,\leq \mu^*)}(\Psi(\mF), \Igss_{K^p}))\\
                                                              & \cong & \Hom_{\shv(\kot_G,\leq \mu^*)}(\Psi(\mF), \Psi(\Psi^L(\Igss_{K^p})))\\
                                                              & \cong & \Hom_{\rshv(\kot_G,\leq \mu^*)}(\mF, \Psi^L(\Igss_{K^p}))
\end{eqnarray*}
where the first isomorphism is by definition, the second isomorphism follows the fully faithfulness of $\Psi^L$, and the last statement follows from \eqref{eq: fully faithfulness of Psi in eventually coconn part of Sh section}.
\end{proof}

Now we give a formula computing \'etale cohomology of Shimura varieties in terms of coherent sheaves on the stack $\locsys_{{}^cG,\bQ_p}^{\tame}$.
\begin{theorem}\label{thm: etale coh Shimura variety vai coh sheaf}
There is an object $\Igss_{K^p}^{\spec,\unip}\in \indcoh(\locsys_{{}^cG,\bQ_p}^{\tame})$ such that for every $\mF\in \fgshv(\iw\bs LG/\iw)$, such that $\Ch_{LG,\phi}^{\unip}(\mF)$ corresponds to $\frakA$ on $\locsys_{{}^cG,\bQ_p}^{\tame}$, we have
\[
C(\Sh_K, (\verd^{\mathrm{verd}}_{\Sh_K})^c ((\loc_p)^!\mF))\cong \Hom_{\indcoh(\locsys_{{}^cG,\bQ_p}^{\tame})}(\frakA, \Igss^{\spec,\unip}).
\]
\end{theorem}
\begin{proof}
We let 
\[
\Igss_{K^p}^{\unip}:=\proj^{\unip}((i_{\leq \mu^*})_*^{\ind\fg}\Psi^L(\Igss_{K^p}))\in \rshv^{\unip}(\kot_G),
\] 
where $\proj^{\unip}$ is as in \eqref{eq:projection-to-unipotent-part}, and let $\Igss_{K^p}^{\spec,\unip}=\bL_G^{\unip}(\Igss_{K^p}^{\unip})$.
Then we apply \Cref{eq:main-theorem-categorical-tame-local-Langlands-modular}.
\end{proof}

\begin{corollary}\label{cor: etale coh Shimura variety via coh sheaf}
Suppose $G$ is unramified.
If $K_p=I$ is the standard Iwahori, we have an $H_{K^p}\times W_E$-equivariant isomorphism
\begin{equation}\label{cor: etale coh Shimura variety via coh sheaf-1}
C(\bfSh_K,\La[d])\cong C(\Sh_K, R\Psi[d])\cong \Hom_{\indcoh(\locsys_{{}^cG,\bQ_p}^{\tame})}(\cohspr^{\unip}\otimes\widetilde{V}_\mu, \Igss^{\spec,\unip})
\end{equation}
Here $\widetilde{V}_\mu$ is the evaluation bundle associated to $V_\mu$ (see \Cref{ex: evaluation bundle}), which is canonically equipped with an action of $W_E$.

If $K_b=I_b$, we have
\[
C(\Igs_{I_bK^p},\La)=\Hom_{\indcoh(\locsys_{{}^cG,\bQ_p}^{\tame})}((\pi^{\unip})_*\mO(\la_b)), \Igss^{\spec,\unip}).
\]
\end{corollary}

\begin{remark}
In fact, both sides of \eqref{cor: etale coh Shimura variety via coh sheaf-1} also admit the action of the Iwahori-Hecke algebra $H_I$. Namely, it acts on the $C(\bfSh_K(G,X)_{\overline\bQ_p}, \La[d])$ via the usual Hecke algebra, and acts on $\cohspr^{\unip}_{{}^cG,\bQ_p}$ via \Cref{intro: thm: end of coh spr}. It will be shown in \cite{yangzhutorsion} that the above isomorphism is also $H_I$-equivariant.
\end{remark}

\begin{remark}\label{rem: cptly supported coh}
Recall the notion of Serre functor. We may recover the compactly supported cohomology of the Shimura variety as
\[
C_c(\bfSh_K,\La[d])= C(\bfSh_K,\La[d])^\vee=\Hom_{\indcoh(\locsys_{{}^cG,\bQ_p}^{\tame})}(\Igss^{\spec,\unip}, S(\cohspr^{\unip}\otimes\widetilde{V}_\mu)).
\]
\end{remark}

\subsubsection{$t$-structure}
In the sequel, we will assume that $\bfSh_K(G,X)$ is projective. Then $\Sh_K$ is pfp proper over $k$.

\begin{lemma}\label{prop: adm selfduality igusa sheaf}
The sheaf $\Igss_{K^p}$ is self-dual with respect to the duality $(\verd_{\kot_G}^\can)^\adm$.
\end{lemma}
\begin{proof}
By \Cref{lem: compatibility of canonical duality between Igs and kot} and \Cref{lem: adm duality and functors}, we have
\[
 (\verd_{\kot_G}^{\can})^{\adm}\circ (\loc_p^0)_\flat\cong (\loc_p^0)_\flat\circ (\verd_{\Igs_{K^p}}^{\can})^{\adm}.
\]
It remains to show that 
\[
(\verd_{\Igs_{K^p}}^{\can})^{\adm}(\consdual_{\Igs_{K^p}})\cong \consdual_{\Igs_{K^p}}.
\]
Recall that by  \eqref{eq:abstract smooth dual 2}, we have
\[
(\verd_{\Igs_{K^p}}^{\can})^{\adm}(\mF)=\underline\Hom(\mF, \omega_{\Igs_{K^p}}^{\can}),
\]
where the internal hom is with respect to the symmetric monoidal structure on $\shv(\Igs_{K^p})$ given by $\os$, and $\omega_{\Igs_{K^p}}^{\can}\in \shv(\Igs_{K^p})$ is the object given by 
\[
\rg^{\can}(\Igs_{K^p}, \mF)^\vee=\Hom_{\shv(\Igs_{K^p})}(\mF, \omega_{\Igs_{K^p}}^{\can}), \quad \forall \mF\in \shv(\Igs_{K^p}).
\]
As $\consdual_{\Igs_{K^p}}$ is the unit for the symmetric monoidal structure $\os$, we see that 
\[
(\verd_{\Igs_{K^p}}^{\can})^{\adm}(\consdual_{\Igs_{K^p}})=\omega_{\Igs_{K^p}}^{\can}.
\]
It remains to show that $\omega_{\Igs_{K^p}}^{\can}=\consdual_{\Igs_{K^p}}$.
To see this, we compute, for $\mF\in\shv(\Sh_K)$,
\begin{multline*}
\Hom((\Nt^{\glob})_* \mF,  \omega_{\Igs_{K^p}}^{\can})\cong \rg^{\can}(\Igs_{K^p}, (\Nt^{\glob})_* \mF)^\vee \\
 \cong  \rg(\Sh_K,\mF)^\vee \cong \Hom(\mF, \consdual_{\Sh_K})\cong \Hom((\Nt^{\glob})_* \mF,  \consdual_{\Igs_{K^p}}),                                                                    
\end{multline*}
where the third isomorphism follows from the properness of $\Sh_K$. As $\shv(\Igs_{K^p})$ is compactly generated by objects of the form $(\Nt^{\glob})_* \mF$ for $\mF\in\cshv(\Sh_K)=\shv(\Sh_K)^\cpt$, we see that $\omega_{\Igs_{K^p}}^{\can}=\consdual_{\Igs_{K^p}}$ as desired.
\end{proof}

Recall a $t$-structure on $\shv(\kot_G)$ as from \Cref{prop: t-structures on llc-2}. We let $\chi=2\rho$.
\begin{proposition}\label{prop: perversity Igusa sheaf}
We have $\Igss\in \shv(\kot_{G,\leq \mu^*},\La)^{2\rho\mbox{-}e,\geq 0}$. When ($!$-)restricted to $\kot_{G,\leq \mu^*}$, we have $\Igss\in \shv(\kot_{G,\leq \mu^*})^{2\rho\mbox{-}e,\heartsuit}$.
\end{proposition}
\begin{proof}
Given \Cref{prop: adm selfduality igusa sheaf} and \Cref{lem: adm duality and e-t-structure}, it is enough to show that $\Igss\in \shv(\kot_G,\La)^{2\rho\mbox{-}e,\geq 0}$. That is,
for every $b\in B(G)$, and every pro-$p$-open compact subgroup $K_b\subset G_b(F)$, we have
\[
\Hom((i_b)_*\cind_{K_b}^{G_b(F)}\La[-\langle2\rho,\nu_b\rangle], \Igss)\in \Mod_\La^{\leq 0}.
\]
If $b\not\in B(G,\mu^*)$, the above space is simply zero. Suppose $b\in B(G,-\mu)$.
We let $w$ be a $\sigma$-straight element in $\adm(\mu)\subset \widetilde{W}$ corresponding to $b$. Then by \Cref{prop: coh of Igus var}, we have 
\[
\Hom((i_b)_*\cind_{K_b}^{G_b(F)}\La[-\langle2\rho,\nu_b\rangle], \Igss)=C_c(\Igs_{K_bK^p},\La[\langle 2\rho, \nu_b\rangle]).
\]
Now we use the fact that when $\Sh_K$ is projective, $\Igs_{K_bK^p}$ is a (perfect) affine scheme of dimension $\langle 2\rho,\nu_b\rangle$ and the usual Artin vanishing to derive the desired estimate. 
\end{proof}

\begin{remark}
Recall that in \Cref{rem: comparison between kotG and BunG}, we discussed the hope of comparison between $\shv(\kot_G)$ and $D_{\mathrm{lis}}(\Bun_G)$ and comparison of $\shv(\kot_G)^{2\rho\mbox{-}e,\heartsuit}$ and the category of perverse sheaves in $D_{\mathrm{lis}}(\Bun_G)$. Under this comparsion, \Cref{prop: perversity Igusa sheaf} is formally analogue \cite[Theorem 8.6.3]{DHKZ}. However, we note that the actually reasonings are different.
\end{remark}

\begin{remark}
When $\bfSh_K(G,X)$ is not projective, then $\omega^{\can}$ and $\consdual_{\Igs_{K^p}}$ are different objects in general. As will be explained in \cite{yangzhutorsion}, one can define a  different version of the Igusa sheaf as
\[
\Igss_{K^p}^{\can}:= (\loc_p^0)_{\flat} \consdual_{\Igs_{K^p}}^{\can}.
\]
Using it, one can use give a formula computing the compactly supported cohomology of $\bfSh_K$, different from the one in \Cref{rem: cptly supported coh}, and is parallel to \Cref{cor: etale coh Shimura variety via coh sheaf}.

In addition, it will be shown in \cite{yangzhutorsion} that $\Igss_{K^p}^{\can}\in \shv(\kot_G)^{2\rho\mbox{-}e,\leq 0}$.
\end{remark}

\newpage

\part{Toolkits}

This part can be regarded as the appendix of the article. Its purpose is to compile relevant backgrounds on category theory, derived algebraic geometry, and sheaf theory. Compared to the main content, the many materials presented here may be seen as general nonsense. On the other hand, many of the results contained in this part are now standard, particularly within the community of geometric representation theory.

This raises the question: why is such a lengthy part necessary? First, although many of the results in this part are known, they are often scattered across various sources, sometimes appearing in forms that are not convenient for our use. Additionally, some results remain only as folklore theorems within the community. For instance, the six-functor formalism of ind-constructible sheaves has not yet been documented in the literature. Our work requires this formalism to be adapted to the context of perfect algebraic geometry, where certain additional subtleties arise. Therefore,
we believe it will benefit readers if we consolidate all the necessary results in one accessible location, rather than demarcating numerous references throughout the literature.
Moreover, there are indeed some new results proved in this part, as far as we are aware of. We will highlight these new results at the beginning of each section of this part. 

\section{Abstract trace formalism}\label{sec:categorical-preliminaries}

In this section, we review the general categorical trace formalism. Many of results in this section are (essentially) known and have appeared in literature, although we generalize and improve upon some existing results at various places. One possible exception is the notion of admissible objects in general dualizable presentable categories, as introduced and studied in \Cref{SS: admissible objects} and in \Cref{SS: cpt gen category}. This concept generalizes the notion of admissible representations of $p$-adic groups, and can be regarded as a dual notion to that of compact objects.  
Additionally, we take this opportunity to discuss \Cref{thm-second-trace-2}. Although this theorem is likely familiar to experts in the field, it has not yet been thoroughly documented in the literature, as far as we are aware of. 
For further discussions of categorical trace, see also \cite{gaitsgory2019toy}, \cite{hoyois2017higher}. For an elementary account, see \cite{Zhu2016}.

\subsection{Recollections of $\infty$-categories}\label{Subsec: recollection of category}
As the whole work uses theory of $\infty$-categories in a substantial way, we first review the required categorical preliminaries mainly following \cite{Lurie.higher.topos.theory} \cite{Lurie.higher.algebra} and \cite{Gaitsgory.Rozenblyum.DAG.vol.I}. 
We sometimes specialize general discussions of \emph{loc. cit.} to situations that suffice for our purposes, but occasionally will also prove results that we could not find in literature. The main purpose of this subsection is to fix our notations and conventions. Of course we will not be able to review all the necessary background materials and therefore will constantly refer to \emph{loc. cit.} for unexplained concepts and terminologies. 

\subsubsection{Categories of $\infty$-categories}\label{SS-cat-of-cat}
We will do ``higher linear algebras", i.e. to manipulate stable $\infty$-categories as if we manipulate vector spaces. For this purpose we first consider the collection of (certain) $\infty$-categories as a whole.
Unless explicitly saying ``ordinary category", by a category we mean an $(\infty,1)$-category. For two categories, we write $\fun(\bfC,\bfD)$ for the category of functors between them. For a category $\bfC$, let $\mathrm{h}\bfC$ denote its homotopy category, which is an ordinary category. A subcategory $\bfD$ of $\bfC$ is defined the Cartesian pullback of an ordinary (not necessarily full) subcategory $\mathrm{h}\bfD\subset \mathrm{h}\bfC$. (Some authors call such $\bfD$ a $1$-full subcategory of $\bfC$.)

Let $\spc$ be the category of spaces (or nowadays called animas). 
Let $\cat$ be the category of all (not necessarily small) categories. We will mainly use the following subcategories of $\cat$
\begin{equation}\label{eq: working ambient categories}
\catid\cong \cptcat\subset\lincat\subset\prl\subset\cat,
\end{equation}
where
\begin{itemize}
\item $\prl\subset\cat$ is the subcategory of presentable categories with morphisms being continuous (i.e. colimit preserving) functors;
\item $\lincat\subset \prl$ is the full subcategory of presentable stable categories;
\item $\cptcat \subset \lincat$ is the subcategory consisting of compactly generated presentable stable categories and morphisms being continuous functors that preserve compact objects;
\item $\catid\subset \cat$ is the subcategory of idempotent complete small stable categories with functors being exact functors. 
\end{itemize}
There is the ind-completion functor 
\[
\ind: \catid\to \cptcat, \quad \bfC\mapsto \ind(\bfC),
\]
which is an equivalence. (But note that this equivalence is not compatible with the embeddings $\catid\to\cat$ and $\cptcat\to\cat$.) A quasi-inverse is given as follows.
For $\bfC\in\lincat$, we also let $\bfC^\cpt$ denote the full subcategory of compact objects of $\bfC$. So $\bfC\cong (\ind(\bfC))^{\cpt}$ for any $\bfC\in\catid$. The functor $\cptcat\to\catid\colon \bfC\mapsto \bfC^\cpt$ is a quasi-inverse of the ind-completion functor $\ind$. 

\begin{example}
For an $E_\infty$-ring $\La$ with unit, the category $\mathrm{Perf}_\La$ of perfect $\La$-modules belongs to $\catid$, while the $\infty$-category $\Mod_\La$ of all $\La$-modules belongs to $\cptcat$. In addition, $\ind(\mathrm{Perf}_\La)\cong \Mod_\La$.
\end{example}

Recall that $\prl$ has a natural closed symmetric monoidal structure (\cite[Proposition 4.8.1.15, Remark 4.8.1.18]{Lurie.higher.algebra}) such that $\prl\to \cat$ is lax symmetric monoidal, where the latter is equipped with the Cartesian symmetric monoidal structure. 
The unit of $\prl$ with this monoidal structure is $\spc$.
For $\bfC_1,\bfC_2\in \prl$, and $c_i\in\bfC_i$, we write $c_1\boxtimes c_2$ for the image of $(c_1,c_2)$ under the canonical functor $\bfC_1\times\bfC_2\to \bfC_1\otimes\bfC_2$.
We will need the following lemma.

\begin{lemma}\label{lem-tensor-preserve-fully-faithfulness}
Let $\bfC_1\to \bfC_2$ be a continuous fully faithful embedding of presentable categories. Let $\bfD\in \prl$. Then $\bfC_1\otimes \bfD\to \bfC_2\otimes \bfD$ is fully faithful.
\end{lemma}
\begin{proof}
We use \cite[Proposition 4.8.1.17]{Lurie.higher.algebra} to identify $\bfC_i\otimes \bfD$ with the category $\mathrm{R}\Fun(\bfD^{\op},\bfC_i)$ of functors from $\bfD^{\op}$ to $\bfC_i$ that admit left adjoints. As $\bfC_1\to \bfC_2$ is fully faithful, so is $\Fun(\bfD^{\op},\bfC_1)\to\Fun(\bfD^{\op},\bfC_2)$ (see for example \cite[Lemma 5.2]{Gepner.Haugseng.Nikolaus.lax}), the lemma then follows. 
\end{proof}

The category $\lincat$ inherits a symmetric monoidal structure from $\prl$ such that the inclusion $\lincat\subset\prl$ is lax monoidal (\cite[Proposition 4.8.2.18]{Lurie.higher.algebra}).
The inclusions $\cptcat\subset \lincat$ are closed under the monoidal structure.
By transport of structure we also obtain a symmetric monoidal structure on $\catid$. Explicitly, the tensor product in $\catid$ is given by the formula 
\[
\bfC_1\otimes\bfC_2\cong (\ind(\bfC_1)\otimes \ind(\bfC_2))^\cpt.
\]
Note that $\bfC_1\otimes\bfC_2$ is the smallest idempotent complete stable full subcategory of $\ind(\bfC_1)\otimes \ind(\bfC_2)$ containing objects $\{c_1\boxtimes c_2\}_{c_i\in\bfC_i}$. 

We recall that arbitrary (co)limits exist in any of the above categories. The inclusion $\prl\subset \cat$ preserves limits (but not colimits in general). The inclusion $\lincat\subset\prl$ preserves both limits and colimits.
The inclusion $\cptcat\subset \lincat$ preserves colimits (but not limits in general). Finally, the inclusion $\catid\subset \cat$ preserves filtered colimits and limits.

\begin{remark}\label{rem: kappa generated presentable category}
In several places in the article, we will perform certain constructions/arguments to these big categories as if they were small categories. To avoid set-theoretic issues, what we will actually do is the following. We fix a regular cardinal $\kappa$ and let $\prl_\kappa$ denote the $\kappa$-compactly generated (in the sense of \cite[Definition 5.5.7.1]{Lurie.higher.topos.theory}) presentable categories. This is a (non-full) subcategory of $\prl$, with $1$-morphisms being those continuous functors that preserve $\kappa$-compact objects. It is well-known that $\prl_\kappa$ itself is presentable (and $\kappa$-compactly generated by a single object: the arrow category of $\spc$) and is closed under the symmetric monoidal structure on $\prl$. Therefore $\prl_\kappa$ is canonically an object in $\calg(\prl_\kappa)$.
Similarly, we have $\lincat_\kappa=\prl_\kappa\cap \lincat$, which is $\kappa$-compactly generated and closed under the symmetric monoidal structure on $\prl$, and therefore $\lincat_\kappa\in \calg(\prl_\kappa)$. 
For example, when $\kappa=\cpt$ is the countable cardinal, then $\lincat_\kappa=\cptcat$ as mentioned above. 

The cardinal $\kappa$ does not really play any role in the discussion sequel and can be chosen to be large enough in each situation we are considering. Therefore, we will omit it from the notation. That is, when we write $\lincat$ (and similarly other large categories), we implicitly mean $\lincat_\kappa$ for some regular cardinal $\kappa$ large enough. 
\end{remark}

Being categories of categories, all of the categories in \eqref{eq: working ambient categories}
naturally form $(\infty,2)$-categories. We will not seriously make use of such $2$-categorical structure except speaking about functor categories and adjoint functors. For example, for $\bfC, \bfD\in \cat$, we have the usual functor category $\Fun(\bfC,\bfD)$. For $\bfC,\bfD\in\prl$, the corresponding functor category, denoted as $\Fun^{\mathrm{L}}(\bfC,\bfD)$, is the full subcategory of $\Fun(\bfC,\bfD)$ consisting of those functors that commute with arbitrary colimits. For $\bfC,\bfD\in\cptcat$, the corresponding functor category, denoted as $\Fun^{\cpt}(\bfC,\bfD)$, is the full subcategory of $\Fun^{\mathrm{L}}(\bfC,\bfD)$ consisting of those functors that commute with arbitrary colimits and preserve compact objects. More generally, if $\bfC,\bfD\in\lincat_\kappa$, we have $\fun^{\kappa}(\bfC,\bfD)\subset \fun^{\mathrm{L}}(\bfC,\bfD)$ consisting of those functors that commute with arbitrary colimits and preserve $\kappa$-compact objects. 
Finally, for $\bfC,\bfD\in\catid$, the corresponding functor category, denoted as $\Fun^{\mathrm{Ex}}(\bfC,\bfD)$, consist of exact functors. The ind-completion induces an equivalence of categories $\Fun^{\mathrm{Ex}}(\bfC,\bfD)\cong \Fun^{\cpt}(\ind(\bfC),\ind(\bfD))$. 

Then one can talk about adjoint functors. Namely, let $f\in \Fun^?(\bfC,\bfD)$ for $?$ being one of the above supscripts. We say that $f$ admits a right (resp. left) adjoint if it admits a right (resp. left) adjoint $f^R$ (resp. $f^L$) in $\cat$, and $f^R$ (resp. $f^L$) belongs to $\Fun^?(\bfD,\bfC)$.

\begin{definition}\label{def:categories-adjointability}
Consider a commutative square in  one of the categories as above.
\[
\begin{tikzcd}
\bfC \arrow[r,"f"]\arrow[d,"v"] & \bfC' \arrow[d,"u"]\\
\bfD \arrow[r,"g"] & \bfD'.
\end{tikzcd}
\]
That is, we are given a specified isomorphism $u\circ f \simeq g\circ v$. 
Then we say that the square above is \textit{right adjointable} in $?$ if $f$ and $g$ admit right adjoints $f^{R}$ and $g^{R}$ in $?$, and the \textit{Beck-Chevalley map} (or sometimes called the \textit{base change map}) $\beta\colon   v\circ f^{R}\rightarrow g^{R}\circ u $ given by
\begin{equation}\label{eq: Beck-Chevalley map}
v\circ f^{R}  \rightarrow g^{R}\circ g\circ v\circ f^{R} \simeq g^{R}\circ u\circ f \circ f^{R} \rightarrow g^{R} \circ u
\end{equation}
is an isomorphism of functors. Dually, we may say the square is \textit{left adjointable}. 
\end{definition}

We will make use of the following statement.
\begin{lemma}\label{lem: 2-commutative digram for right adjointable}
Given a commutative square as in \Cref{def:categories-adjointability} and suppose it is right adjointable. Then the following diagrams are $2$-commutative.
\[
\xymatrix{
\bfC\ar_v[d] \rrtwocell^{\id_\bfC}_{f^R\circ f}&& \bfC\ar^v[d], & \bfC'\ar_u[d] \rrtwocell_{\id_{\bfC'}}^{f\circ f^R}&& \bfC'\ar^u[d]\\ 
\bfD \rrtwocell^{\id_\bfD}_{g^R\circ g}&& \bfD, & \bfD' \rrtwocell_{\id_{\bfD'}}^{g\circ g^R}&& \bfD' .
}
\]
\end{lemma}

Now let $F\colon  S \rightarrow \lincat$ be a diagram. For an arrow $\varphi\colon  s\rightarrow s'$ the functor $F(\varphi)\colon F(s) \rightarrow F(s')$ preserves colimits and therefore \cite[Corollary 5.5.2.9]{Lurie.higher.topos.theory} admits a right adjoint $F^{R}(\varphi)$ (in $\cat$). By passing to right adjoints we get a diagram $F^{R}\colon  S^{\op} \rightarrow \cat$. By \cite[\textsection {5.5.3}]{Lurie.higher.topos.theory} there is a canonical equivalence
\begin{equation}\label{eq:limit-colimit equivalence}
\colim_{s\in S} F(s) \rightarrow \lim_{s\in S^{op}} F^{R}(s),
\end{equation}
where the morphism is
determined by the maps right adjoint to $\mathrm{ins}_{s}\colon  F(s) \rightarrow \colim_{S} F $, where the left hand side is computed in $\lincat$ and then is mapped to $\cat$, and where the right hand side is computed in $\cat$. In addition, if all $F^{R}(\varphi)$ are continuous, then the right hand side of \eqref{eq:limit-colimit equivalence} can also be computed in $\lincat$ and \eqref{eq:limit-colimit equivalence} is an equivalence in $\lincat$. Denote by $\mathrm{ev}_{s}$ the right adjoint of $\mathrm{ins}_{s}$. It follows from adjunction that for every object $c\in \colim_{S} F$, the natural map
\begin{equation}\label{eq:cat-colimits-objects-adj}
\colim_{s\in S} (\mathrm{ins}_{s} \circ \mathrm{ev}_{s}(c)) \rightarrow c
\end{equation}
is an equivalence in $\colim_{S} F$.

\begin{remark}\label{rem:cat-prelim-limit-colimit}
\begin{enumerate}
    \item Assume that for each $\varphi\colon  s\rightarrow s'$ the functor $F(\varphi)\colon F(s) \rightarrow F(s')$ preserves compact objects. Then the functors $\mathrm{ins}_{s}\colon  F(s) \rightarrow \colim_{S} F$ also preserve compact objects.

    \item If $S$ is filtered and the morphisms in the image of $F$ have continuous right adjoints, then for an object $s$ in $S$ the composition $\mathrm{ev}_{s}\circ \mathrm{ins}_{s}\colon F(s) \rightarrow \colim_{S} F\simeq \lim_{S^{\op}} F^{R} \rightarrow F(s)$ is equivalent to the colimit 
    \[
    \mathrm{ev}_{s}\circ \mathrm{ins}_{s} \simeq \colim_{\varphi\colon s\rightarrow s'} F^{R}(\varphi)\circ F(\varphi).
    \]
\end{enumerate}
\end{remark}

We also review adjointability under taking (co)limits.

\begin{proposition}\label{prop:categorical-right-adjointability-colimits} 
Let $S,T$ be small $\infty$-categories and let $F\colon S\times T \rightarrow \lincat$ be a functor. For $s\rightarrow s'$ in $S$ and $t\rightarrow t'$ in $T$, consider the the square
\begin{equation}\label{eq: right adjointable diagram}
\begin{tikzcd}
F(s,t) \arrow[r] \arrow[d] & F(s',t) \arrow[d]\\
F(s,t') \arrow[r] & F(s',t').
\end{tikzcd}
\end{equation}
If for all $s\rightarrow s'$ in $S$ and $t\rightarrow t'$ in $T$, the square \eqref{eq: right adjointable diagram} is right adjointable (in $\lincat$), then there is an extension $\overline{F}\colon S^{\triangleright}\times T^{\triangleleft} \rightarrow \lincat$ of $F$ such that: 
\begin{enumerate}
    \item\label{prop:categorical-right-adjointability-colimits-1} For each $t\in T$, the diagram $\overline{F}\colon S^{\triangleright}\times \{t\} \rightarrow \lincat$ is a colimit diagram in $\lincat$.
    \item\label{prop:categorical-right-adjointability-colimits-2} For each $s\in S$, the diagram $\overline{F}\colon \{s\}\times T^{\triangleleft} \rightarrow \lincat$ is a limit diagram in $\lincat$.
    \item\label{prop:categorical-right-adjointability-colimits-3} For all $s\rightarrow s'$ in $S^{\triangleright}$ and $t\rightarrow t'$ in $T^{\triangleleft}$ the corresponding square \eqref{eq: right adjointable diagram} is right adjointable (in $\lincat$).    
\end{enumerate}
\end{proposition}
\begin{proof}
As $\lincat\subset \prl$ preserves all limits and colimits, we may replace $\lincat$ by $\prl$. Then this is \cite[Proposition 4.7.4.19]{Lurie.higher.algebra}, except that we need to show that for every $s\rightarrow s'$ in $S^{\triangleright}$ and $t$ in $T^{\triangleleft}$, the right adjoint of $F(s,t)\to F(s',t)$ is continuous.

Indeed, as argued in \cite[Proposition 4.7.4.19]{Lurie.higher.algebra},
by passing to the right adjoint our assumption gives $F^R: S^{\op}\times T\to \prl$. Right Kan extension gives $\overline{F^R}: (S^{\op})^{\triangleleft}\times T^{\triangleleft}\to \prl$. As the inclusion $\prl\to \cat$ commutes with limits, this is also the right Kan extension in $\cat$. It follows from  \cite[Proposition 4.7.4.19]{Lurie.higher.algebra} that for $(s\to s')\in S^{\triangleright}=((S^{\op})^{\triangleleft})^{\op}$ and $t\in T^{\triangleleft}$ the right adjoint of $F(s,t)\to F(s',t)$ is $\overline{F^R}(s',t)\to \overline{F^R}(s,t)$, which is continuous.
\end{proof}

\begin{remark}\label{rmk:categorical-left-adjointability-limits}
Suppose we are given $F\colon S\times T\to \lincat$ as in \Cref{prop:categorical-right-adjointability-colimits} but
now suppose for all $s\rightarrow s'$ in $S$ and $t\rightarrow t'$ in $T$, the square \eqref{eq: right adjointable diagram} is left adjointable (in $\lincat$). Then by passing to the left adjoints and apply \Cref{prop:categorical-right-adjointability-colimits} and \eqref{eq:limit-colimit equivalence}, we obtain $\overline{F}\colon S^{\triangleleft}\times T^{\triangleleft} \rightarrow \lincat$ satisfying conditions parallel \Cref{prop:categorical-right-adjointability-colimits} \eqref{prop:categorical-right-adjointability-colimits-1}-\eqref{prop:categorical-right-adjointability-colimits-3}, with ``colimit" replaced by ``limit" in \eqref{prop:categorical-right-adjointability-colimits-1} and ``right adjointable" replaced by ``left adjointable" in \eqref{prop:categorical-right-adjointability-colimits-3}.
\end{remark}

\subsubsection{Descent}

Recall that for an $\infty$-category $\bfD$, a \textit{monad on} $\bfD$ is an associative algebra object $T$ in the monoidal category $\fun(\bfD,\bfD)$. 
If $G\colon  \bfE \rightarrow \bfD$ is a functor which admits a left adjoint $F$, then the composition $T = G\circ F$ has the structure of a monad on $\bfD$ with identity given by the unit map $\id_{\bfD} \rightarrow G\circ F$ of the adjunction and composition map induced by the co-unit $F\circ G \rightarrow \mathrm{id}_{\bfE}$ via
\[
T\circ T = (G\circ F) \circ (G\circ F) \simeq G\circ (F \circ G) \circ F \rightarrow G\circ F. 
\]

Given a monad $T$ on $\bfD$ one can consider the category $\lmodu_{T}(\bfD)$ of left modules over $T$. The forgetful functor 
$G\colon \lmodu_{T}(\bfD) \rightarrow \bfD$ 
has a left adjoint given by the free construction $A \mapsto T(A)$. An adjunction $F\colon \bfD \rightleftarrows \bfE\colon G$ is called monadic if $\bfE$ is equivalent to $\lmodu_{T}(\bfD)$ for $T = G\circ F$ and $G$ given by the forgetful functor. See \cite[\textsection{4.7.1}]{Lurie.higher.algebra} for detailed discussions.

Now we review (cohomological) descent. For our purpose, we need a slightly stronger version of \cite[Theorem 4.7.5.2, Corollary 4.7.5.3]{Lurie.higher.algebra}. Let $\Delta$ denote the (ordinary) simplex category of non-empty finite linearly ordered sets
and let $\Delta_s\subset \Delta$ denote the subcategory consisting of \emph{injective} maps $[n]\to [m]$.  If one drops the non-emptyness requirement, the resulting categories are denoted by $\Delta_{s,+}\subset \Delta_+$.
Recall that a functor $\Delta\to \cat$ is usually called a cosimplicial category and a functor $\Delta_s\to \cat$ is usually called a semi-cosimplicial category.

\begin{theorem}\label{thm:Beck-Chevalley-descent}
Let $\bfC^{\bullet}\colon  \Delta\rightarrow \cat$ be a cosimplicial category. Assume that for any $\alpha\colon [m] \rightarrow [n]$ in $\Delta_s$, the induced diagram
\begin{equation}\label{eq: BC cosimplicial}
\begin{tikzcd}
\bfC^{m} \arrow[r,"d^{0}"]\arrow[d] &  \bfC^{m+1} \arrow[d]\\
\bfC^{n} \arrow[r,"d^{0}"] & \bfC^{n+1}
\end{tikzcd}
\end{equation}
is left adjointable. We denote the left adjoint of $d^{0}\colon  \bfC^{n} \rightarrow \bfC^{n+1}$ by $F(n)$. Let $\bfC = \tot(\bfC^{\bullet})$. Then the following statements hold.
\begin{enumerate}
    \item The functor $G\colon  \bfC\rightarrow \bfC^{0}$ admits a left adjoint $F$.
    \item The diagram
    \[
    \begin{tikzcd}
    \bfC \arrow[r,"G"]\arrow[d,"G"] & \bfC^{0}\arrow[d,"d^{1}"]\\
    \bfC^{0} \arrow[r,"d_0"] & \bfC^{1}
    \end{tikzcd}
    \]
    is left adjointable. That is, the canonical map $F(0) \circ d^{1} \rightarrow G\circ F$ is an equivalence.  
    \item The adjunction $F\colon \bfC^{0} \rightleftarrows \bfC\colon G$ is monadic. That is, $\bfC$ is equivalent to the category of left modules $\lmodu_{T}(\bfC^{0})$ with $T = F(0) \circ d^{1} \simeq G\circ F$.
\end{enumerate}
    Suppose the above cosimplicial diagram $\bfC^\bullet: \Delta\to \cat$ extends to an augmented cosimplicial diagram $\Delta_+\to \cat$. Let $G'\colon \bfC^{-1} \rightarrow \bfC^0$ denote the augmentation functor.  In addition, assume that  the diagram \eqref{eq: BC cosimplicial} is left adjointable for any $\alpha\colon [m] \rightarrow [n]$ in $\Delta_{s,+}$, and that the category $\bfC^{-1}$ admits geometric realizations of $G'$-split simplicial objects that are preserved by $G'$.  Then 
\begin{enumerate}[resume]
\item the canonical map $\phi\colon  \bfC^{-1} \rightarrow \tot(\bfC^{\bullet})$ admits a fully faithful left adjoint. If, in addition $G':\bfC^{-1} \rightarrow \bfC^{0}$ is conservative, $\phi$ is an equivalence. 
\end{enumerate}
\end{theorem}
\begin{remark}
\begin{enumerate}
\item Comparing with \cite[Theorem 4.7.5.2, Corollary 4.7.5.3]{Lurie.higher.algebra}, we only require left adjointability of \eqref{eq: BC cosimplicial} involving face maps. Such slightly weaker assumption is crucial for our computations of categorical traces.
\item There is also a dual (a.k.a. co-monadic) version by replacing ``left adjoint" with ``right adjoint", and ``realizations of $G$-split simplicial objects"  with ``totalizations of $G$-split cosimplicial objects" in the statement above.
\end{enumerate}
\end{remark}
\begin{proof}[Proof of \Cref{thm:Beck-Chevalley-descent}]
We explain how to modify the argument of \emph{loc. cit.} under this weaker assumption. Namely, we keep the argument of the first paragraph in the proof of  \cite[Theorem 4.7.5.2]{Lurie.higher.algebra} showing that $(\bfC\to \bfC^{0})=\lim_{\Delta}(\bfC^{\bullet}\to \bfC^{\bullet+1})$ in $\fun(\Delta^1, \cat)$. By \cite[Lemma 6.5.3.7]{Lurie.higher.topos.theory}, this is also the limit of the underlying semi-cosimplicial diagram in $\fun(\Delta^1, \cat)$. Now we proceed as in the second paragraph of the proof of  \cite[Theorem 4.7.5.2]{Lurie.higher.algebra}, but with ``cosimplicial" replaced by ``semi-cosimplicial". As in \emph{loc. cit.}, our assumption together with \Cref{prop:categorical-right-adjointability-colimits} and \Cref{rmk:categorical-left-adjointability-limits} then shows that $(\bfC\to \bfC^{0})=\lim_{\Delta_s}(\bfC^{\bullet}\to \bfC^{\bullet+1})$ in $\fun^{\mathrm{LAd}}(\Delta^1, \cat)$. Then one deduces all the desired statements from this fact as in \emph{loc. cit.}
\end{proof}

\subsubsection{Linear categories with $t$-structure}
We use cohomological convention in this article. So
for a stable category $\bfC$, we let $\bfC^{\leq 0}$ denote the connective part of a $t$-structure. Let $\bfC^+=\cup_n \bfC^{\geq n}$ be the bounded from below subcategory.

We will need certain (non-full) subcategory $\lincat^{t,+}$ of $\lincat$ consisting of $(\bfC,\bfC^{\leq 0})$, where $\bfC\in \lincat$ equipped of a $t$-structure which is
\begin{itemize}
\item accessible (\cite[Definition 1.4.4.12]{Lurie.higher.algebra}), compatible with filtered colimits (\cite[Definition 1.3.5.20]{Lurie.higher.algebra}); 
\item $\bfC$ is compactly generated, and $\bfC^\cpt\subset \bfC^b$. 
\end{itemize}
We require morphisms between $(\bfC,\bfC^{\leq 0})$ and $(\bfD,\bfD^{\leq 0})$ in
$\mathrm{h}\lincat^{t,+} $ to be those that are left $t$-exact up to a cohomological shift, i.e. those $F: \bfC\to \bfD$ 
such that $F[d](\bfC^{\geq 0})\subset \bfD^{\geq 0}$ for some integer $d$ (depending on $F$). (Note that we do not require $F$ to preserve compact objects.)

We also need a symmetric monoidal structure on $\lincat^{t,+}$. We start with the following easy statement, whose proof is left to readers. (Note that this slightly generalize \cite[Proposition 2.2.1.1 (1), (2)]{Lurie.higher.algebra}.)

\begin{lemma}\label{lemma: construct functors from homotopy cat}
Let $\bfC^{\otimes}\to \mO^\otimes$ be a map of operads, and let $\mathrm{h}\bfC^{\otimes}\to \mathrm{h}\mO^\otimes$ be the induced map at the homotopy level. Suppose we have a map of (ordinary) operads
$\mathrm{h}\bfD^\otimes\to \mathrm{h}\bfC^\otimes$. Then $\bfD^\otimes:=\mathrm{h}\bfD^\otimes\times_{\mathrm{h}\bfC^\otimes}\bfC^\otimes\to \bfC^\otimes$ is a map of operads. If $\bfC^\otimes\to \mO^\otimes$ is coCartesian and if $\mathrm{h}\bfD\to \mathrm{h}\bfC$ is an $\mathrm{h}\mO$-monoidal functor, then $\bfD^\otimes\to \mO^\otimes$ is coCartesian and $\bfD^\otimes\to \bfC^\otimes$ is $\mO$-monoidal.
\end{lemma}

We apply this lemma to  endow $\lincat^{t,+}$ with a symmetric monoidal structure by endowing $\mathrm{h}\lincat^{t,+}$ with a symmetric monoidal structure such that the natural functor $\mathrm{h}\lincat^{t,+}\to \mathrm{h}\lincat$ is symmetric monoidal. Namely, we define $(\bfC,\bfC^{\leq 0})\otimes (\bfD,\bfD^{\leq 0})=(\bfC\otimes\bfD, \bfC^{\leq 0}\otimes\bfD^{\leq 0})$. As explained in \cite[Remark C.4.2.2]{Lurie.SAG}, the natural functor $\bfC^{\leq 0}\otimes\bfD^{\leq 0}\to \bfC\otimes\bfD$ is indeed fully faithful and defines an accessible $t$-structure of $\bfC\otimes\bfD$ compatible with filtered colimits. In addition, $\bfC\otimes\bfD$ is compactly generated, with $(\bfC\otimes\bfD)^\cpt$ generated as idempotent complete category by objects of the form $c\boxtimes d$ where $c\in \bfC^\cpt$ and $d\in \bfD^\cpt$. Note that there is some $m,n\in \bZ$ such that $c[m]\in \bfC^{\leq 0}$ and $d[n]\in \bfD^{\leq 0}$, and are truncated objects. So $(c\otimes d)[m+n]\in (\bfC\otimes \bfD)^{\leq 0}$, and is truncated. This shows that $(\bfC\otimes\bfD, \bfC^{\leq 0}\otimes\bfD^{\leq 0})$ indeed belongs to $\lincat^{t,+}$. The unit is given by the category of spectra equipped with the natural $t$-structure. In addition, clearly the associativity and commutativity constraints in $\mathrm{h}\lincat$ are $t$-exact with respect to the tensor product $t$-structure. It follows that we have the well-defined symmetric monoidal structure on $\mathrm{h}\lincat^{t,+}\to \mathrm{h}\lincat$. This endows $\lincat^{t,+}$ with a well-defined symmetric monoidal structure.

The following lemma is easy to check.
\begin{lemma}\label{lem: lax symmetric monoidal from lincat with t to cat}
There is a lax symmetric monoidal functor $\lincat^{t,+}\to \cat$ sending $(\bfC,\bfC^{\leq 0})$ to $\bfC^+$. The corresponding operad map $(\lincat^{t,+})^\otimes\to \cat^\otimes$ is a (non-full) subcategory.
\end{lemma}

\subsubsection{Relative tensor product}\label{sec:relative-tensor-product}
Let us review the general formalism of relative tensor products.
Let $\bfR$ be a symmetric monoidal ($\infty$-)category with $\mathbf{1}_\bfR$ its unit.  Let $\alg(\bfR)$ denote the category of associative algebra objects in $\bfR$. Let $\lmodu(\bfR)$ (resp. $\rmodu(\bfR)$) the category of left (resp. right) module objects in $\bfR$. I.e., objects in $\lmodu(\bfR)$ (resp. $\rmodu(\bfR)$) consist of pairs $(A,M)$ with $A\in\alg(\bfR)$ and $M$ a left (resp. right) $A$-module.
For $A\in \alg(\bfR)$ we denote by $\lmodu_{A}(\bfR)=\lmodu(\bfR)\times_{\alg(\bfR)}\{A\}$ (resp. $\rmodu_{A}(\bfR)=\{A\}\times_{\alg(\bfR)}\lmodu(\bfR)$). Recall that if $A$ is a commutative algebra, then $\lmodu_{A}(\bfR)$ inherits a symmetric monoidal structure from $\bfR$, and will be denoted by $\Mod_{A}(\bfR)$.

Similarly, let $\BMod(\bfR)$ denote the category of bimodule objects in $\bfR$.
For $A, B\in \alg(\bfR)$ we denote by ${_{A}}\BMod_{B}=\{A\}\times_{\alg(\bfR)}\BMod(\bfR)\times_{\alg(\bfR)}\{B\}$ the category of $A\mbox{-}B$-bimodules. We identify $A\mbox{-}\mathbf{1}_\bfR$-bimodules with left $A$-modules and $\mathbf{1}_\bfR\mbox{-}A$-bimodules with right $A$-modules.
An $A\mbox{-}A$-bimodule is also called as an $A$-bimodule. For example, $A$ itself can be regarded as $A$-bimodule via the left and the right multiplication. See \cite[\textsection{4.3}]{Lurie.higher.algebra} for detailed discussions. Given associative algebra objects $A,B,C\in \alg(\bfR)$ and bimodules $M\in {_{A}}\BMod_{B}$ and $N \in {_{B}}\BMod_{C}$, the relative tensor product $M\otimes_{B}N$, if exists, is the unique object (up to equivalence) in ${_{A}}\BMod_{C}$, corepresenting the functor sending $X\in {_{A}}\BMod_{C}$ to the space of $B$-bilinear $A\mbox{-}C$-bimodule maps $M\otimes N\to X$ (in appropriate homotopy sense, see \cite[Definition 4.4.2.3]{Lurie.higher.algebra}). On the other hand, there is the two-sided bar construction 
\[
   \lmodu_A(\bfR)\times_{\alg(\bfR)}\rmodu_C(\bfR)\to ({_{A}}\BMod_{C})^{\Delta^\op},\quad (M,N)\mapsto \barcons_B(M,N)_\bullet,
\]
where $\barcons_B(M,N)_\bullet$ is a simplicial object in the category of $A\mbox{-}C$-bimodules, given informally as 
\[
\barcons_B(M,N)_n = M\otimes B^n \otimes N
\]
with face maps induced by the multiplication on $B$ and actions on $M$ and $N$, and degeneracy maps given by insertions of the unit of $B$. See \cite[Notation 4.4.2.4, Construction 4.4.2.7]{Lurie.higher.algebra}. If $A=B=C$ and $M=N=A$, we simply denote $\barcons_A(A,A)_\bullet$ by $\barcons(A)_\bullet$, called the \textit{bar construction} of the bimodule $A$. 

We do not know whether $M\otimes_BN$ is always given by the geometric realization of $\barcons_B(M,N)_\bullet$ as soon as the latter exists. This is the case if $\bfR$ admits geometric realizations and such that the monoidal product $\otimes\colon \bfR \times \bfR \rightarrow \bfR$ preserves geometric realizations in each variable, by
\cite[Theorem 4.4.2.8]{Lurie.higher.algebra}. It is also the case in the following two examples. 

\begin{example}\label{ex:relative-tensor-free}
  Assume that $M=M_0\otimes B$ with $M_0$ a left $A$-module (resp. $N=B\otimes N_0$ with $N_0$ a right $C$-module). Then $M\otimes_BN$ exists and is represented by $M_0\otimes N$ (resp. $M\otimes N_0$). To see this, we follow the argument of \cite[Proposition 5.2.2.6]{Lurie.higher.algebra}. If $\bfR$ admits geometric realizations and such that the monoidal product $\otimes\colon \bfR \times \bfR \rightarrow \bfR$ preserves geometric realizations in each variable, then $M\otimes_BN$ exists and is isomorphic to $M_0\otimes N$ (resp. $M\otimes N_0$) by \cite[Proposition 4.4.3.14, 4.4.3.16]{Lurie.higher.algebra}. The general situation reduces this case via Yoneda embedding.
\end{example}

\begin{example}\label{ex:coCartesian-relative-tensor}
Using a similar argument as above, one can also prove the existence of the relative tensor product in the following situation.
Let $\bfC$ be a (small, $\infty$-)category admitting finite limits. Let $\pt$ denote the final object. Let $\bfC^{\op,\sqcup}$ denote $\bfC^{\op}$ equipped with the coCartesian symmetric monoidal structure: for $X,Y\in \bfC$, the tensor product $X\otimes Y$ in $\bfC^\op$ is the finite product $X\times Y$ in $\bfC$. We note that every object $X$ is a commutative algebra object in $\bfC^{\op}$, with the multiplication given by the diagonal map $\Delta_X: X\to X\times X$ in $\bfC$ and the unit given by the structural map $\pi_X: X\to \pt$. In addition, every morphism $f: X\to Y$ in $\bfC$ gives a commutative algebra homomorphism in $\bfC^\op$. Furthermore, $\lmodu_X(\bfC^{\op,\sqcup})= (\bfC_{/X})^{\op}$.
(Rigorously, these facts follow from \cite[Proposition 2.4.3.9, Corollary 2.4.3.10]{Lurie.higher.algebra}.)
Now given two morphisms $a:M\to X, b:N\to X$ in $\bfC$, we regard $M,N$ as $X$-modules in $\bfC^{\op}$. We claim that 
$M\otimes_X N$ exists and is representable by $M\times_XN$. Namely, the two-sided bar complex $\barcons_X(M,N)$ in $\bfC^{\op}$ is given by
the cosimplicial object in $\bfC$ as
\[
M\times N\substack{\id\times a\times \id \\ \longrightarrow \\ \longrightarrow \\ \id\times b\times \id} M\times X\times N\substack{\longrightarrow\\ \longrightarrow\\ \longrightarrow} M\times X\times X\times N \cdots
\]
We consider the embedding $\bfC^{\op}\to \Ind(\bfC^{\op})$, where $\Ind(\bfC^{\op})$ denotes the ind-completion of $\bfC^{\op}$, equipped with the induced symmetric monoidal structure so that the tensor product preserves filtered colimits in each variable. 
The tensor product then also preserves geometric realizations in each variable. 
Then again by the argument of \cite[Proposition 5.2.2.6]{Lurie.higher.algebra}, it is enough to show that the geometric realization of this simplicial object $\Delta^{\op}\to \Ind(\bfC^{\op})$ is represented (in $\Ind(\bfC^{\op})$) by $M\times_XN$. By \cite[Lemma 6.1.4.7]{Lurie.higher.topos.theory}, its geometric realization can be computed as the colimit of the truncated colimit diagram $(\Delta_{\leq 1})^{\op}\to \Ind(\bfC^{\op})$, which in turn is the limit of $M\times N\substack{\longrightarrow\\ \longrightarrow} M\times X\times N$ in $\bfC$. But this is exactly $M\times_XN$.

Now as the relative tensor products exist in $\bfC^{\op,\sqcup}$, the category ${}_X\BMod_{X}(\bfC^{\op,\sqcup})$ of $X$-bimodules in $\bfC^{\op,\sqcup}$ has a natural monoidal structure (see \cite[Proposition 4.4.3.12]{Lurie.higher.algebra}). 
Its opposite category ${}_X\BMod_{X}(\bfC^{\op,\sqcup})^{\op}$ then also is a monoidal category. We claim there is a canonical functor 
\begin{multline*}
\pi:({}_X\BMod_{X}(\bfC^{\op,\sqcup})^{\op})^{\oast}\to \bfC,\\
(M_1,M_2,\ldots,M_n)\in ({}_X\BMod_{X}(\bfC^{\op,\sqcup})^{\op})^{\oast}_{[n]}\to M_1\times_XM_2\times_X\cdots\times_XM_n.
\end{multline*}
Indeed recall for any symmetric monoidal category $\bfR$ and an associative algebra $A\in\bfR$, the natural forgetful functor ${}_A\BMod_A(\bfR)\to \bfR$ is lax monoidal. It follows that we have 
\[
({}_X\BMod_{X}(\bfC^{\op,\sqcup})^{\op})^{\oast}\to \bfC^{\otimes}\to \bfC,
\]
where the last functor comes from the Cartesian structure of $\bfC$, which sends $(M_1,\ldots,M_n)$ to $M_1\times M_2\times\cdots M_n$. On the other hand, ${}_X\BMod_{X}(\bfC^{\op,\sqcup})^{\op}$ admits final object $X\times X$. So the above functor naturally factors through $({}_X\BMod_{X}(\bfC^{\op,\sqcup})^{\op})^{\oast}\to \bfC_{/X^{2n}}$. Now, the desired functor $\pi$ is the composition of this functor with $\bfC_{/X^{2n}}\to \bfC$ obtained by pullback along the $X\times X^{n-1}\times X\xrightarrow{\id\times\Delta_{23}\times\Delta_{45}\times\cdots\Delta_{2n-2,2n-1}\times \id}X^{2n}$.
\end{example}

We also recall the notion of duality for bimodules. See \cite[\textsection{4.6.2}]{Lurie.higher.algebra} for more details.
\begin{definition}\label{def: the notation of left dual}
 Let $A,B\in \alg(\bfR)$ and let $M\in  {_{A}}\BMod_{B}$. A left dual of $M$ 
is given by an object $N\in {_{B}}\BMod_{A}$ together with a unit (or co-evaluation)
\begin{equation}\label{eq:unit-map-for-dual}
u_M: B \rightarrow N\otimes_{A} M
\end{equation}
which is a morphism in ${_{B}}\BMod_{B}$, and a co-unit (or evaluation) 
\begin{equation}\label{eq:counit-map-for-dual}
e_M:M\otimes_{B} N \rightarrow A
\end{equation}
which is a morphism in ${_{A}}\BMod_{A}$, such that the compositions
\begin{align}
    M\simeq M\otimes_{B} B & \xrightarrow{\id \otimes u_M} M\otimes_{B} N\otimes_{A} M \xrightarrow{e_M\otimes \id} A\otimes_{A} M \simeq M \label{triangle-identity-M}\\
    N \simeq B \otimes_{B} N & \xrightarrow{u_M\otimes \id} N\otimes_{A} M\otimes_{B} N \xrightarrow{\id \otimes e_M} N\otimes_{A} A \simeq N. \label{triangle-identity-N}
\end{align}
are (homotopic to) the identities on $M$ and $N$.

By abuse of notations, we sometimes write ${}^\vee M$ for $N$.
\end{definition}

\begin{remark}\label{rem: the notation of right dual}
Clearly, we also have the notion of right dual. If $N$ is a left dual of $M$, then $M$ is a right dual of $N$.
By abuse of notations we sometimes write it as $N^\vee$. 
\end{remark}

\begin{example}\label{example-free-A-module-dual}
  Let $M=A$, regarded as a left $A$-module. Then $M$ admits a left dual given by $N=A$ regarded as the right $A$-module. The unit and evaluation maps are
  \[
  \mathbf{1}_\bfR\xrightarrow{\mathbf{1}_A}A\cong A\otimes_AA, \quad A\otimes A\xrightarrow{m} A.
  \]
\end{example}

\begin{remark}\label{rem: dualizable object in sym monoidal cat}
\begin{enumerate}
\item\label{rem: dualizable object in sym monoidal cat-1} We recall that given $M$, its left dual $(N, u_M, e_M)$ is unique up to a contractible choice.
\item\label{rem: dualizable object in sym monoidal cat-2} Let $M\in  {_{A}}\BMod_{B}$. If $M$ admits a left dual $N$, then the functor $\lmodu_B(\bfR)\to \lmodu_A(\bfR), \ L\mapsto M\otimes_B(-)$ admits a right adjoint, given by $N\otimes_A(-)$. 
\item\label{rem: dualizable object in sym monoidal cat-3} Let $M\in  {_{A}}\BMod_{B}$ with a left dual $N$. Then for every left $A$-module $L$, the internal hom $\underline\Hom_A(M,L)\in \lmodu_B(\bfR)$ exists and is representable by $N\otimes_AL$. That is, for every $X\in \lmodu_B(\bfR)$, the natural map
\[
\Map_{\lmodu_B(\bfR)}(X, N\otimes_AL)\to \Map_{\lmodu_A(\bfR)}(M\otimes_B X, M\otimes_B N\otimes_AL)\to  \Map_{\lmodu_A(\bfR)}(M\otimes_B X, L)
\]
is an isomorphism.
\item\label{rem: dualizable object in sym monoidal cat-4} Specializing to $B=\mathbf{1}_\bfR$, we obtain the notion of a left dual of a left $A$-module. By \cite[Proposition 4.6.2.13]{Lurie.higher.algebra}, $M$ admits a left dual as an $A\mbox{-}B$-bimodule if and only if it admits a left dual as a left $A$-module. 
\item\label{rem: dualizable object in sym monoidal cat-5} Further specialize to the case $A=B=\mathbf{1}_\bfR$, we arrive to the notion of dualizable objects in the symmetric monoidal category $\bfR$. I.e. $M\in \bfR$ is dualizable in $\bfR$ if there exists $N$ and morphisms $u_M:\mathbf{1}_\bfR\to N\otimes M$ and $e_M:M\otimes N\to \mathbf{1}_\bfR$ such that \eqref{triangle-identity-M} and \eqref{triangle-identity-N} are homotopic to the identities of $M$ and $N$. 
Note that the commutativity constraints also identify $N$ as the right dual of $M$. Following traditional notations, we usual denote $N$ as $M^\vee$.
\end{enumerate}
\end{remark}

\subsubsection{$\bfA$-linear categories}\label{SS: la-linear categories}
We will apply the above discussions to $\bfR=\lincat$. (See \Cref{rem: kappa generated presentable category} for our convention.)

Now let $\bfA\in \alg(\lincat)$, i.e. a monoidal presentable stable category with monoidal product commutes with colimits separately in each variable. Write $\lincat_{\bfA}=\lmodu_\bfA(\lincat)$ for simplicity. Objects in $\lincat_\bfA$ are called presentable $\bfA$-linear stable categories, or sometimes simply called $\bfA$-linear categories. Similarly morphisms in $\lincat_\bfA$ are simply called $\bfA$-linear functors. 

Arbitrary (co)limits exist in $\lincat_{\bfA}$ and
the forgetful functor $\lincat_{\bfA}\to \lincat$ commutes with all (co)limits (using \cite[\textsection{3.4.3}, \textsection{3.4.4}]{Lurie.higher.algebra}). In fact $\lincat_\bfA$ itself is a presentable category, and is compactly generated.
For $\bfC\in \lincat_\bfA$, and $c,d\in\bfC$, we write 
$\Hom_{\bfC/\bfA}(c,d)\in\bfA$ determined (up to equivalence) by 
\begin{equation}\label{eq: internal hom object in A}
\Map_{\bfA}(a,\Hom_{\bfC/\bfA}(c,d))=\Map_{\bfC}(a\otimes c,d),\quad \forall a\in\bfA.
\end{equation} 
On the other hand, for $a\in \bfA$ and $c\in \bfC$, we define $\Hom^{\bfC/\bfA}(a,c)\in\bfC$ (up to equivalence) such that for every $d\in\bfC$,
\begin{equation}\label{eq: internal hom object in C}
\Map_{\bfC}(d, \Hom^{\bfC/\bfA}(a,c))=\Map_{\bfC}(a\otimes d, c).
\end{equation}
Sometimes, we just write $\Hom_{\bfC}(c,d)$ or $\Hom(c,d)$ (and similarly $\Hom^\bfC(a,c)$ or $\Hom(a,c)$) for simplicity if no confusion is likely to arise. In addition, when $\bfC=\bfA$, we write $\Hom_{\bfA/\bfA}=\Hom^{\bfA/\bfA}$ as $\underline{\Hom}$, which is the usual internal hom of $\bfA$.

As $\lincat_\bfA$ is tensored over $\lincat$,
all $\bfA$-linear functors between two $\bfA$-modules $\bfM$ and $\bfN$ form a presentable stable category $\Fun_\bfA^{\mathrm{L}}(\bfM,\bfN)$, equipped with a continuous functor $\Fun_\bfA^{\mathrm{L}}(\bfM,\bfN)\to \Fun^{\mathrm{L}}(\bfM,\bfN)$. 
In particular, giving $(F: \bfM\to \bfN)\in \lincat_\bfA$, it makes sense to ask whether it admits an $\bfA$-linear right or left adjoint. 
To address this question, we suppose the underlying functor $F$ admits a continuous right adjoint $F^R$ (resp. a left adjoint $F^L$). Then $F^R$ (resp. $F^L$) admits a natural lax (resp. oplax) $\bfA$-linear structure, given by the Beck-Chevallay map \eqref{eq: Beck-Chevalley map} associated to the following commutative diagram (in $\lincat$)
\begin{equation}\label{eq: A-linear structure}
\xymatrix{
\bfA\otimes \bfM\ar^{\id\otimes F}[r] \ar_{\mathrm{act}_\bfM}[d] & \bfA\otimes\bfN\ar^{\mathrm{act}_\bfN}[d]\\
\bfM \ar^F[r] &\bfN.
}
\end{equation}
Then $F^R$ (resp. $F^L$) is $\bfA$-linear if this diagram is right (resp. left) adjointable. 

This is not the case in general, but is the case for an important class of algebra objects in $\lincat$. 

\begin{lemma}\label{lem-rigid-lax-is-strict}
Let $\bfA\in \alg(\lincat)$. Suppose the product $m: \bfA\otimes\bfA\to\bfA$ admits an $\bfA\otimes\bfA^{\rev}$-linear right adjoint $m^R$. Then the following statements hold.
\begin{enumerate}
\item\label{lem-rigid-lax-is-strict-1} For every $\bfA$-module $\bfM$, the continuous right adjoint of $\mathrm{act}_\bfM$ exists and is given by 
\[
\bfM\xrightarrow{(m^R\circ \mathbf{1}_\bfA)\otimes\id_\bfM} \bfA\otimes \bfA\otimes \bfM \xrightarrow{\id_\bfA\otimes \mathrm{act}_\bfM}\bfA\otimes \bfM.
\]
In particular, the Beck-Chevalley map $ (\id_\bfA\otimes F)\circ \mathrm{act}_{\bfM}^R\to  \mathrm{act}_{\bfN}^R\circ F$ associated to $F\circ \mathrm{act}_\bfM\cong \mathrm{act}_\bfN\circ (\id_\bfA\otimes F)$ is an isomorphism.
\item\label{lem-rigid-lax-is-strict-2} Every (op)lax $\bfA$-linear functor between $\bfA$-module categories is $\bfA$-linear. 
Consequently, if $F:\bfM\to \bfN$ is an $\bfA$-linear functor between $\bfA$-module categories, with a continuous right adjoint $F^R$ (resp. a left adjoint $F^L$), then $F^R$ (resp. $F^L$) is $\bfA$-linear.
\end{enumerate}
\end{lemma}
\begin{proof}
This is \cite[Lemma 1.9.3.2, Lemma 1.9.3.6]{Gaitsgory.Rozenblyum.DAG.vol.I}. Note that only the above assumption of $\bfA$ is needed in the proof. 
\end{proof}

We note that the relations between adjoints and (co)limits as discussed in \Cref{SS-cat-of-cat} continue to hold in $\lincat_\bfA$, as soon as we require adjoints to be $\bfA$-linear.

Now suppose $\bfA$ is a commutative algebra in $\lincat$, then $\lincat_\bfA$ inherits a closed symmetric monoidal structure from $\lincat$. In this case, for $\bfM,\bfN\in \lincat_\bfA$, we write their tensor product as $\bfM\otimes_\bfA\bfN$. The category $\Fun_\bfA^{\mathrm{L}}(\bfM,\bfN)$ admits a natural $\bfA$-module structure making it the internal hom between $\bfM$ and $\bfN$. In the sequel, for $m\in\bfM$ and $n\in\bfN$, we will let $m\boxtimes_\bfA n$ denote the image of $(m,n)$ under the natural functor $\bfM\times \bfN\to \bfM\otimes \bfN\to \bfM\otimes_\bfA\bfN$.

\begin{example}
The particular important example is the category $\bfA=\Mod_\La$  for an $E_\infty$-ring $\La$, which is a commutative algebra object in $\lincat$. We will write $\lincat_\La$ instead of $\lincat_{\Mod_\La}$, $\bfC\otimes_\La\bfD$ instead of $\bfC\otimes_{\bfA}\bfD$, and $m\boxtimes_\La n$ instead of $m\boxtimes_\bfA n$. We will say $\La$-linearity instead of $\bfA$-linearity in this case. 

Let $\cptcat_{\La}\subset\lincat_\La$ be the subcategory consisting of those $\La$-linear categories $\bfC$ such that the underline category $\bfC$ is compactly generated and those $\La$-linear continuous functors that preserve compact objects.
On the other hand $\mathrm{Perf}_\La$ is a commutative algebra object in $\catid$, and we let $\catid_\La$ denote its module category, usually called the $(\infty,1)$-category of $\La$-linear small idempotent complete stable categories with morphisms being $\La$-linear exact functors. As before, Ind-completion induces an equivalence $\catid_\La\cong \cptcat_{\La}$ of symmetric monoidal categories. As before, $\bfC,\bfD\in \catid_\La$, we use $\bfC\otimes_\La\bfD$ to denote its tensor product. 
\end{example}

The following notion (see  \cite[Definition 1.9.1.2]{Gaitsgory.Rozenblyum.DAG.vol.I}) will play important roles in our discussions.

\begin{definition}\label{def: rigid monoidal category}
An algebra object $\bfA\in \alg(\lincat)$ is called rigid if $\mathbf{1}_\bfA$ is compact and the product $m: \bfA\otimes\bfA\to\bfA$ admits an $(\bfA\otimes\bfA^{\rev})$-linear right adjoint $m^R$ (as in \Cref{lem-rigid-lax-is-strict}).
\end{definition}

We mention that if $\bfA$ is also compactly generated, $\bfA$ being rigid is equivalent to requiring that compact objects of $\bfA$ admit both left and right duals (see \cite[Definition D.7.4.1]{Lurie.SAG} and \cite[Lemma 1.9.1.5]{Gaitsgory.Rozenblyum.DAG.vol.I}).

We record the following statement for applications.
\begin{lemma}\label{lem: complement on monoidal categories}
Let $\bfA$ be a monoidal (resp. symmetric monoidal) presentable stable category with monoidal product commutes with colimits separately in each variable.  Let $\bfI\subset \bfA$ be full subcategory. If for every $m\in \bfA$ and $n\in \bfI$, both $m\otimes n$ and $n\otimes m$ belong to  $\bfI$, then $\bfI$ is has a natural $\bfA$-bimodule structure such that the inclusion $\iota: \bfI\subset \bfA$ is $(\bfA\otimes\bfA^{\rev})$-linear. In addition, if $\iota^R$ is also $(\bfA\otimes\bfA^{\rev})$-linear, then $\bfI$ has a natural monoidal (resp. symmetric monoidal) structure such that $\iota^R$ is monoidal (resp. symmetric monoidal).  
\end{lemma}
\begin{proof}
That $\bfI$ is an $\bfA$-bimodule and $\iota$ is $\bfA$-bilinear follows from  \cite[Proposition 2.2.1.1]{Lurie.higher.algebra} directly. We prove the rest statements.

We notice that $\bfI$ has a natural non-unital (symmetric) monoidal structure, by restriction from $\bfA$. 
Applying \cite[Theorem 5.4.4.5]{Lurie.higher.algebra}, it is then enough to show that at the level of homotopy categories, $\mathrm{h}\bfI$ admits a unit and that $\iota^R: \mathrm{h}\bfA\to \mathrm{h}\bfI$ is (symmetric) monoidal.

Let $\mathbf{1}_\bfA$ be the unit of $\bfA$, and let $\mathbf{1}_\bfI:=\iota^R(\mathbf{1}_\bfA)$. We notice that for every $n\in\bfI$, by assumption we have
\[
\mathbf{1}_\bfI\otimes m\cong \iota^R(\mathbf{1}_\bfA\otimes m)\cong m,
\]
and similarly $m\otimes \mathbf{1}_\bfI$. This gives the desired statement.
\end{proof}

\subsection{Dualizable categories}\label{SS: dualizable categories}
In \Cref{rem: dualizable object in sym monoidal cat}, we have reviewed the notion of dualizable objects in a symmetric monoidal category. We now specialize this notion to the case $\bfR=\lincat_\bfA$, for a fixed commutative algebra $\bfA$ in $\lincat$ (e.g. $\bfA=\Mod_\La$).

\subsubsection{dualizable categories}\label{sec:dual-categories}
For $\bfC\in\lincat_\bfA$, let $\bfC^{\vee,\bfA}=\fun^{\mathrm{L}}_{\bfA}(\bfC,\bfA)$. By definition there is a natural pairing 
\begin{equation}\label{eq-pairing-linear-dual-category}
e_{\bfC/\bfA}\colon \bfC\otimes_\bfA \bfC^{\vee,\bfA}\to \bfA.
\end{equation} 
If $\bfC$ is dualizable in $\lincat_\bfA$, then the above pairing gives the evaluation map in the duality datum and realizes $\bfC^{\vee,\bfA}$ as a dual of $\bfC$. We denote the unit of the duality datum by 
\begin{equation}\label{eq-unit-dual-category}
u_{\bfC/\bfA}:\bfA\to \bfC^{\vee,\bfA}\otimes_\bfA \bfC.
\end{equation} 
Note that by $\bfA$-linearity, $u_\bfC$ is uniquely determined by its value at $\mathbf{1}_\bfA\in\bfA$. Therefore, we usually regard $u_\bfC$ as an object in $\bfC^{\vee,\bfA}\otimes_\bfA\bfC$. Under the canonical equivalence 
\[
\bfC^{\vee,\bfA}\otimes_\bfA\bfC\cong \bfC\otimes_\bfA\bfC^{\vee,\bfA}\cong \fun^L_{\bfA}(\bfC,\bfC),
\] 
$u_\bfC$ corresponds to the identity functor. More generally, an $\bfA$-linear functor $\phi: \bfC\to \bfC$ corresponds to an object (the ``kernel")
\begin{equation}\label{eq: kernel of an endofunctor}
K_\phi=(\id_{\bfC}\otimes \phi)(u_{\bfC})\in \bfC^{\vee,\bfA}\otimes_\bfA\bfC.
\end{equation}

In the sequel, if $\bfA$ is clear from the context, for simplicity we sometimes just write $(\bfC^\vee, u_{\bfC}, e_{\bfC})$ instead of $(\bfC^{\vee,\bfA}, u_{\bfC/\bfA}, e_{\bfC/\bfA})$.

\begin{example}\label{ex: Serre functor}
Let $\bfC$ be a dualizable $\bfA$-module. We define
the Serre functor $S_{\bfC/\bfA}: \bfC\to \bfC$ to be the $\bfA$-linear functor  such that the corresponding object $K_{S_{\bfC/\bfA}}\in \bfC\otimes_\bfA\bfC^\vee$ represents the contravariant functor 
\[
\Map_\bfA(e_{\bfC/\bfA}(-),\mathbf{1}_\bfA): \bfC\otimes_\bfA\bfC^\vee\cong \bfC^\vee\otimes_\bfA\bfC\to \spc.
\] 
(As we shall review in \Cref{ex: explicit Serre and admissible}, when $\bfA=\Mod_\La$ and $\bfC$ is compactly generated this reduces the usual notion of Serre functor.) If $\bfA$ is clear, we also write it as $\bfS_\bfC$ for simplicity.
Recall that $\bfC$ is called $0$-Calabi-Yau if $S_{\bfC}\cong \id_{\bfC}$. 
\end{example}

\begin{remark}\label{rem-duality-dualizable category} 
All dualizable $\bfA$-linear categories can be organized into a (non-full) subcategory $\lincat_\bfA^{\mathrm{dual}}$ with objects being dualizable $\bfA$-linear categories with $1$-morphisms being $\bfA$-linear functors that admit $\bfA$-linear right adjoint.

Let $F: \bfC\to \bfD$ be such a $1$-morphism in $\lincat_\bfA^{\mathrm{dual}}$ with an $\bfA$-linear right adjoint $F^R: \bfD\to \bfC$. Let 
\begin{equation}\label{eq-conjugate.functor}
F^o:=(F^R)^\vee: \bfC^\vee\to \bfD^\vee,
\end{equation}
called the \textit{conjugate} functor to $F$. Note that $F^o$ also admits an $\bfA$-linear right adjoint, namely $F^\vee$.
It follows that there is a symmetric monoidal self-equivalence
\begin{equation}\label{eq: duality functor, cpt}
(-)^\vee: \lincat_\bfA^{\mathrm{dual}}\to \lincat_\bfA^{\mathrm{dual}},\quad  (F\colon \bfC\to\bfD)\mapsto  (F^{o}\colon \bfC^\vee\to \bfD^\vee).
\end{equation}
\end{remark}

\begin{remark}\label{rem-smooth-and-proper}
Note that $\lincat_\bfA^{\mathrm{dual}}$ inherits a symmetric monoidal structure from $\lincat_\bfA$. However, not every object $\bfC$ in $\lincat_\bfA^{\mathrm{dual}}$ is dualizable for the symmetric monoidal structure of $\lincat_\bfA^{\mathrm{dual}}$. Indeed, $\bfC$ is dualizable in $\lincat_\bfA^{\mathrm{dual}}$ if both $u_\bfC: \bfA\to \bfC^\vee\otimes_\bfA\bfC$ and $e_\bfC: \bfC\otimes_\bfA\bfC^\vee\to \bfA$ admit $\bfA$-linear right adjoint. So this is a very restrictive condition on $\bfC$.
Later on in \Cref{SS: cpt gen category}, we will see that when $\bfA=\Mod_\La$, $\cptcat_\La$ is a full subcategory of $\lincat_\La^{\mathrm{dual}}$. Then a compactly generated category $\bfC$ is dualizable in $\cptcat_\La$ is equivalent to $\bfC$ being $2$-dualizable in $\lincat_\La$ (in the sense \Cref{def-2-dualizable-category} below). More explicitly, it means that $u_\bfC$ regarded as an object in $\bfC^\vee\otimes_\La\bfC$ is compact, and $\Hom_{\bfC}(c,d)\in\mathrm{Perf}_\La$ for every $c,d\in\bfC^{\cpt}$. 
\end{remark}

\subsubsection{Localization sequence}\label{SS: localization sequences}
The unit map in the duality datum for a dualizable category is usually hard to write down explicitly. The following result \Cref{prop-gluing-unit-in-semi-orthogonal-decomposition} says that a localization sequence induces a filtration of the unit, which sometimes gives a way to understand it.

\begin{definition}\label{rem:localization sequence}
Let $\bfA$ be an associative algebra in $\lincat$. 
A sequence $\bfM\xrightarrow{F}\bfC\xrightarrow{G}\bfN$ of $\bfA$-linear categories is called a localization sequence if:
\begin{enumerate}
\item both $F$ and $G$ admit $\bfA$-linear right adjoint $F^R$ and $G^R$, and the natural adjunctions $\id_\bfM\to F^R\circ F$ and $G\circ G^R\to \id_{\bfN}$ are equivalences;
\item $G\circ F=0$, and for every $c\in \bfC$ the sequence
\begin{equation}\label{eq:fiber sequence localization}
 F( F^R(c))\to c \to G^R(G(c)).
\end{equation}
is a fiber sequence in $\bfC$.
\end{enumerate}  
If in addition $G^R$ also admits an $\bfA$-linear right adjoint, then we say $(F(\bfM), G^R(\bfN))$ form a semi-orthogonal decomposition of $\bfC$. 
\end{definition}

\begin{remark}
We have not checked whether a localization sequence as defined above is a cofiber sequence in $\lincat_\bfA$. On the other hand, one can define this notion in a more general $(\infty,2)$-categorical setting, see \cite[Definition 3.2]{hoyois2017higher}.
\end{remark}

Now assume that $\bfA$ an a commutative algebra and $\bfM,\bfN,\bfC$ are dualizable.

Let $F^o= (F^R)^\vee\colon \bfM^\vee\to \bfC^\vee$ be the conjugate of $F$, and $G^o= (G^R)^\vee\colon \bfC^\vee\to\bfN^\vee$ be the conjugate of $G$. Then $\bfM^\vee\xrightarrow{F^o}\bfC^\vee\xrightarrow{G^o}\bfN^\vee$ is still a localization sequence.
\quash{
which is still fully faithful preserving compact objects, and $(F^{o})^R=F^\vee$. Similarly we have $G^o=(G^R)^\vee\colon \bfC^\vee\to \bfB^\vee$ with fully faithful right adjoint $(G^o)^R=G^\vee$, which identifies $(G^o)^R(\bfB^\vee)=\ker (F^{o})^R$. Then $(F^o(\bfA^\vee), (G^o)^R(\bfB^\vee))$ form a semi-orthogonal decomposition of $\bfC^\vee$ and there is the fiber sequence in $\bfC^\vee$
\[
 F^o( (F^o)^R(c))\to c \to (G^o)^R(G^o(c)).
\]}

Now we regard $u_{\bfC}$ as an object in $\bfC^\vee\otimes_\bfA\bfC$ and similarly regard $u_\bfM\in \bfM^\vee\otimes_\bfA\bfM$ and $u_\bfN\in  \bfN^\vee\otimes_\bfA\bfN$. 
\begin{lemma}\label{prop-gluing-unit-in-semi-orthogonal-decomposition} 
Under the above situation, there is a fiber sequence 
\[
(F^o\otimes F)u_\bfM\to u_\bfC\to ((G^o)^{R}\otimes G^R)u_\bfN,
\]
where the maps come from \Cref{lem-functoriality-duality-data}.
\end{lemma}
\begin{proof}
First, we have the fiber sequences 
\[
( \id\otimes F)(\id\otimes F^R)u_\bfC\to u_\bfC\to (\id\otimes G^R)(\id\otimes G)u_\bfC,
\] 
\[
(F^{o}\otimes \id)((F^{o})^{R}\otimes\id )(\id\otimes G)u_\bfC\to (\id\otimes G)u_\bfC\to ((G^{o})^{R}\otimes \id) (G^o\otimes \id) (\id\otimes G)u_\bfC,
\] 
\[
(F^{o}\otimes \id)((F^{o})^{R}\otimes\id )(\id\otimes F^R)u_\bfC\to (\id\otimes F^R)u_\bfC\to (\id\otimes ((G^{o})^{R}\otimes \id) (G^o\otimes \id)(\id\otimes F^R)u_\bfC.
\]
Note that $((F^{o})^{R}\otimes\id )(\id\otimes G)u_\bfC=0$ as under duality it corresponds to the functor $G\circ F=0$. Similarly, $(G^o\otimes \id)(\id\otimes F^R)u_\bfC=0$.
In addition, under duality, $((F^o)^{R}\otimes F^R)u_\bfC$ corresponds to the functor $F^R\circ F\cong \id$ and therefore  $u_\bfM\cong ((F^o)^{R}\otimes F^R)u_\bfC$. Similarly, $(G^o\otimes G)u_{\bfC}\cong u_{\bfN}$.
Putting all the considerations together gives the lemma.
\end{proof}

Now, let $S\rightarrow \lincat_\bfA,\ s\mapsto \bfC_s$ be a diagram such that each $\bfC_s$ is dualizable and all transition functors $\bfC_s \rightarrow \bfC_{s'}$ admits an $\bfA$-linear right adjoint. Denote $\bfC = \colim_{s\in S} \bfC_s$. Using \eqref{eq:limit-colimit equivalence} and \eqref{eq:cat-colimits-objects-adj}, it is not difficult to see (e.g. see \cite[Proposition 6.3.4]{Gaitsgory.Rozenblyum.DAG.vol.I}) that the natural map
\begin{equation}\label{eq:colimit-of-dual-category}
\colim_{s\in S} \bfC_s^\vee \rightarrow \bfC^\vee
\end{equation}
obtained by passing to conjugate functors, is an equivalence, and $\bfC$ is dualizable with the unit 
\begin{equation}\label{eq-unit-for-colimit-category}
u_\bfC\cong \colim_s  ((\mathrm{ins}_s)^o\otimes \mathrm{ins}_s)(u_{\bfC_s}).
\end{equation} 

We further assume that  $S=\bN_{\geq 0}$ and every $\bfC_{n-1} \rightarrow \bfC_{n}$ is fully faithful and fits into a localization sequence $\bfC_{n-1}\to \bfC_n\to \bfD_n$. 
Then \Cref{prop-gluing-unit-in-semi-orthogonal-decomposition} and \eqref{eq-unit-for-colimit-category} give the following.

\begin{corollary}\label{cor-filtration-unit}
There is a filtration of the unit $u_{\bfC}$ with associated graded being $((G_n)^\vee\otimes (G_n)^R)u_{\bfD_n}$.
\end{corollary}

Another consequence of \Cref{prop-gluing-unit-in-semi-orthogonal-decomposition} is the well-known localization sequence of Hochschild homology to be discussed in \Cref{cor-hochschild-semi-orthogonal} below. 

\subsubsection{Admissible objects}\label{SS: admissible objects}
Our next goal is to generalize the notion of admissible representations in the representation theory of $p$-adic groups.

\begin{definition}\label{def:adm vs compact}
Let $\bfC\in\lincat_\bfA$. For $c\in \bfC$, we let $F_c: \bfA\to \bfC$ denote the $\bfA$-linear functor determined by $c$ under the equivalence
$\Fun^{\mathrm{L}}_{\bfA}(\bfA, \bfC)\cong \bfC, \ F\mapsto F(\mathbf{1}_\bfA)$. Then $c$ is called 
\begin{itemize}
\item $\bfA$-admissible if $F_c$ admits an $\bfA$-linear left adjoint; and 
\item $\bfA$-compact if $F_c$ admits an $\bfA$-linear right adjoint. 
\end{itemize}
\end{definition}

\begin{example}\label{lem:basic cpt and adm}
\begin{enumerate}
\item If $\bfA$ is rigid, then by \Cref{lem-rigid-lax-is-strict} an object $c\in\bfC$ is $\bfA$-compact if and only if $c\in \bfC^\cpt$. 
\item\label{lem:basic cpt and adm-3} Let $F: \bfC\to \bfD$ be an $\bfA$-linear functor. If $F$ admits an $\bfA$-linear right adjoint (resp. $\bfA$-linear left adjoint), then $F$ sends $\bfA$-compact (resp. $\bfA$-admissible) objects to $\bfA$-compact (resp. $\bfA$-admissible) objects. 
\item Let $\bfA\to \bfA'$ is a map of commutative algebras in $\lincat$, and let $\bfC$ be an $\bfA$-module. For if $c\in \bfC$ is $\bfA$-compact (resp. $\bfA$-admissible), then $c\boxtimes_\bfA \mathbf{1}_{\bfA'} \in \bfC\otimes_{\bfA}\bfA'$ is $\bfA'$-compact (resp. $\bfA'$-admissible).
\item\label{lem:basic cpt and adm-2} Assume that $\bfA=\Mod_\La$, and $\bfC$ is compactly generated. Then later on in \Cref{lem: char of adm obj in cg cat} we will show that $c$ is $\Mod_\La$-admissible if and only if $\Hom(d,c)$ is a perfect $\La$-module for every $d\in \bfC^\cpt$. In this case we simply call $\Mod_\La$-admissible objects being admissible. When $\bfC$ is the category of smooth representations of a $p$-adic group with comp, admissible objects in $\bfC$ specialize to the classical notion of admissible representations (see \Cref{rem: adm objects in rep}). This justifies our terminology. When $\bfC$ is the category of $\ell$-adic sheaves on the classifying stack of an algebraic groups, then admissible objects in $\bfC$ coincide with constructible sheaves (see \Cref{ex-perverse-sheaf-classifying-stack}).
\end{enumerate}
\end{example}

Let us have some more discussions of this notion. 

Let $\bfC$ be a dualizable $\bfA$-module.  
Let 
\begin{equation}\label{eq:dual of adm object}
c^* := (F_c^L)^\vee(\mathbf{1}_\bfA)\in \bfC^\vee.
\end{equation} 
\begin{remark}
We choose $c^*$ rather than $c^\vee$ as the notation, as the latter has been used as the dual of a dualizable object when $\bfC$ has a (symmetric) monoidal structure.
\end{remark}

Then $(F_c^L,F_c)$-adjunction gives
\begin{equation}\label{eq: unit-counit-adjunction-admissible}
u_\bfC\to  c^*\boxtimes_\bfA c ,\quad e_\bfC(c \boxtimes_\bfA c^*)\to \mathbf{1}_\bfA.
\end{equation}

The general facts about adjoint functors give the following lemma.
\begin{lemma}\label{rem: Zorro for admissible objects}
If $c\in \bfC$ is $\bfA$-admissible, so is the object $c^* \in \bfC^\vee$ as in \eqref{eq:dual of adm object}, and $c^{**}\in \bfC^{\vee\vee}$ is canonically isomorphic to $c$.
In addition, the following composed map induced by \eqref{eq: unit-counit-adjunction-admissible}
\begin{equation}\label{eq: Zorro for admissible objects}
\begin{split}
c\cong c\boxtimes_{\bfA}\mathbf{1}_\bfA &\cong (e_\bfC\otimes \id_\bfC)(\id_\bfC\otimes u_\bfC) (c\boxtimes_{\bfA}\mathbf{1}_\bfA)\cong (e_\bfC\otimes \id_\bfC) (c\boxtimes_\bfA u_\bfC) \\
       & \to (e_\bfC\otimes \id_\bfC) (c\boxtimes_\bfA c^* \boxtimes_\bfA c)\cong  e_\bfC(c \boxtimes_{\bfA} c^*) \boxtimes_\bfA c\to    \mathbf{1}_\bfA \boxtimes_{\bfA} c \cong c
\end{split}
\end{equation}
is homotopic to the identity map and so is a similar map for $c^*$. 

Conversely, for $c\in\bfC$, if there is an object $d\in\bfC^\vee$ equipped with $u_\bfC\to d\boxtimes_\bfA c$ and $e_\bfC(c\boxtimes_\bfA d)\to \mathbf{1}_\bfA$ such that the similarly defined maps $c\to c$ and $d\to d$  as above are homotopic to the identity map, then $c$ is $\bfA$-admissible and $c^*\simeq d$.
\end{lemma}

\begin{lemma}\label{lem: dual hom for admissible objects}
If $c$ is $\bfA$-admissible, then we have the canonical isomorphism of functors
\[
\Map_{\bfC^\vee\otimes_\bfA\bfC}(u_\bfC, c^* \boxtimes_\bfA (-))\cong \Map_{\bfC}(c,-): \bfC\to \spc.
\]
\end{lemma}
\begin{proof}
Let $d\in \bfC$. For simplicity, write $-\boxtimes-$ instead of $-\boxtimes_\bfA-$.
The isomorphism in the lemma is given by the following two mutually inverse maps.
\[
\Map(u_\bfC, c^* \boxtimes d)\xrightarrow{c\boxtimes (-)} \Map(c\boxtimes u_\bfC, c \boxtimes c^*\boxtimes d)\xrightarrow{e_\bfC\boxtimes \id_\bfC}\Map(c, e_\bfC(c\boxtimes c^*)\boxtimes d)\rightarrow \Map(c,d).
\]
\[
\Map(c,d)\xrightarrow{c^*\boxtimes (-)}\Map(c^*\boxtimes c, c^*\boxtimes d)\rightarrow \Map(u_\bfC, c^*\boxtimes d).
\]
\end{proof}

We let $\bfC^\adm$ denote the full subcategory of $\bfA$-admissible objects in $\bfC$.

\begin{lemma}\label{lem: characterization of admissible dual}
If $c$ is $\bfA$-admissible, then $c^*$ represents the functor 
\[
(\bfC^\vee)^{\op}\to \spc, \quad d\mapsto \Map_\bfA(e_\bfC(c\boxtimes_\bfA d), \mathbf{1}_\bfA).
\]
The assignment $c\mapsto c^*$ induces an equivalence $(\bfC^{\adm})^{\op}\cong (\bfC^\vee)^{\adm}$.
\end{lemma}
\begin{proof}
We need to show that for $d\in \bfC^\vee$, giving $e_\bfC(c\boxtimes_\bfA d)\to \mathbf{1}_\bfA$ amounts to giving a map $d\to c^*$. Indeed, the desired map is given similar to \eqref{eq: Zorro for admissible objects} as
\begin{multline*} 
d\cong (\mathbf{1}_\bfA \boxtimes_\bfA d )\cong (\id_{\bfC^\vee}\otimes e_\bfC)(u_\bfC \otimes \id_{\bfC^\vee}) (\mathbf{1}_\bfA \boxtimes_\bfA d )\cong (\id_{\bfC^\vee}\otimes e_\bfC) (u_\bfC\boxtimes_{\bfA}d) \\
\to (\id_{\bfC^\vee}\otimes e_\bfC) (c^*\boxtimes_{\bfA}c\boxtimes_\bfA d)\cong c^*\boxtimes_\bfA  e_\bfC(c \boxtimes_{\bfA} d) \to   c^*\boxtimes_\bfA  \mathbf{1}_\bfA \cong c^*.
\end{multline*}
Conversely, a map $d\to c^*$ induces $e_\bfC(c\boxtimes_\bfA d)\to e_\bfC(c\boxtimes_\bfA c^*)\to \mathbf{1}_\bfA$. These two constructions are inverse to each other since
\eqref{eq: Zorro for admissible objects} (for both $c$ and $c^\vee$) is homotopic to the identity map.  

The last statement is clear.
\end{proof}

\begin{lemma}
The category $\bfC^\adm$ is an idempotent complete stable category.
\end{lemma}
\begin{proof}
Let $c_1\to c_2\to c$ be a cofiber sequence, with $c_1, c_2$ admissible. Let $d$ be the fiber of $c_2^*\to c_1^*$ in $\bfC^\vee$. One checks that $d$ gives the desired object needed in \Cref{rem: Zorro for admissible objects} to verify that $c$ is admissible.
\end{proof}

\Cref{lem: characterization of admissible dual} suggests us to extend the assignment $c\mapsto c^*$ for $\bfA$-admissible objects to a functor 
\begin{equation}\label{eq: functor of taking vee}
(-)^*\colon \bfC^{\op}\to \bfC^\vee, \quad \Map_{\bfC^\vee}(d,c^* )=\Map_\bfA(e_\bfC(c\boxtimes_\bfA d), \mathbf{1}_\bfA).
\end{equation}
Note that this is a functor in $\cat$ but not in $\lincat_\bfA$. But iterating it twice gives 
\[
\bfC=(\bfC^{\op})^{\op}\to (\bfC^\vee)^{\op}\to (\bfC^\vee)^\vee\cong \bfC, \quad c\mapsto c^{**}
\]
equipped with a natural functorial transformation $c\to c^{**}$.
We say an object $c\in \bfC$ is $\bfA$-reflexive if this map is an isomorphism. Note that by \Cref{rem: Zorro for admissible objects}, $\bfA$-admissible objects are $\bfA$-reflexive.

\subsubsection{Self-duality}

In many cases, the category $\bfC$ admits a canonical $\bfA$-linear self-duality, i.e. an $\bfA$-linear equivalence
\[
\verd: \bfC^\vee\cong \bfC.
\]
As any such equivalence will preserve admissible objects, by \Cref{lem: characterization of admissible dual} we see that $\verd$ restricts to an equivalence
\begin{equation}\label{eq: dual of admissible objects} 
\verd^{\adm}:=\verd((-)^*): (\bfC^{\adm})^{\op}\cong \bfC^{\adm},\quad c\mapsto \verd(c^*),
\end{equation}

\begin{remark}
When $\bfC$ is compactly generated, then $\verd$ will restrict to an equivalence $\verd^\cpt: (\bfC^\cpt)^{\op}\cong \bfC^\cpt$, as we shall see later. However, \eqref{eq: dual of admissible objects} holds without any compact generation assumption.
\end{remark}

\begin{example}\label{ex: duality via Frobenius-structure}
Suppose that $\bfC$ is an $\bfA$-algebra. Recall a Frobenius structure of $\bfC$ (see \cite[Definition 4.6.5.1]{Lurie.higher.algebra}) is an $\bfA$-module functor $\la: \bfC\to \bfA$ such that the composed functor
\[
\bfC\otimes_\bfA\bfC\xrightarrow{m} \bfC\xrightarrow{\la} \bfA
\]
forms the co-unit map in the duality datum of $\bfC$. 
Therefore, it induces an $\bfA$-module equivalence 
\[
\verd^\la: \bfC^\vee\cong \bfC,
\] 
such that $e_\bfC(c\boxtimes_{\bfA}d)=\la(c\otimes \verd^{\la}(d))$ for every $c\in \bfC$ and every $d\in \bfC^\vee$. 
 Then the  functor 
 \[
 \bfC^{\op}\to \bfC^\vee\xrightarrow{\verd^\la} \bfC, \quad c\mapsto \verd^\la(c^*)
 \] 
 takes a more familiar form as follows. For simplicity, we write $c^{*,\la}$ instead of $\verd^\la(c^{*})$.
Let $\omega^\la\in \bfC$ that represents the functor 
\[
\bfC^{\op}\to \bfA, \quad c\mapsto \Map_\bfA(\la(c),\mathbf{1}_\bfA).
\] 
Then 
\begin{multline*}
\Map_{\bfC}(d,\verd^\la(c^*))\cong\Map_{\bfC^\vee}((\verd^\la)^{-1}(d),c^*)\cong\Map_{\bfA}(e_\bfC(c\boxtimes_{\bfA} (\verd^\la)^{-1}(d)),\mathbf{1}_\bfA)\\
\cong\Map_{\bfA}(\la(c\otimes d),\mathbf{1}_\bfA)= \Map_{\bfA}(\la(d\otimes \sigma_\la(c)),\mathbf{1}_\bfA) =\Map_{\bfC}(d\otimes \sigma_\la(c), \omega^\la).
\end{multline*}
Here
\begin{equation}\label{eq:Serre automorphism of Frob algebra}
\sigma_\la: \bfC\to \bfC,
\end{equation} 
 is the Serre automorphism associated to the Frobenius algebra $(\bfC,\la)$,
which is a monoidal automorphism of $\bfC$ (see \cite[Remark 4.6.5.4, Remark 4.6.5.6]{Lurie.higher.algebra}), characterized such that $\la(a\otimes b)\cong \la(b\otimes \sigma_\la(a))$ for every $a,b\in \bfC$.

Therefore,
\begin{equation}\label{eq:abstract smooth dual 2}
c^{*,\la}=\underline\Hom(\sigma_\la(c),\omega^\la).
\end{equation}
In particular, if $\bfC=\bfA$ with the Frobenius structure given by $\la=\id_\bfA$, then 
\begin{equation}\label{eq:abstract smooth dual 3}
c^{*,\id}=\underline\Hom(c,\mathbf{1}_\bfA).
\end{equation}
\end{example}

\begin{remark}\label{rem: involutive property of duality}
Let $(\bfC,\la)$ be as in \Cref{ex: duality via Frobenius-structure}. We consider $(\verd^\la)^\vee\circ (\verd^\la)^{-1}: \bfC\to \bfC$. Note that by definition, for every $c\in \bfC$ and $d\in\bfC^\vee$ we have
\[
\la( (\verd^\la)^\vee((\verd^\la)^{-1}(c))\otimes \verd^\la(d))=e_\bfC(  (\verd^\la)^\vee((\verd^\la)^{-1}(c)) \boxtimes_{\bfA} d)=e_\bfC(  \verd^\la(d) \boxtimes_\bfA (\verd^\la)^{-1}(c) )=\la( \verd^\la(d)\otimes c).
\]

It follows that
\begin{equation}\label{eq: duality square serre automorphism}
(\verd^\la)^\vee\circ (\verd^\la)^{-1}\cong (\sigma_\la)^{-1}: \bfC\to \bfC.
\end{equation}
In particular, if $\bfC$ is a commutative algebra, then the equivalence $\verd^\la$ as in \Cref{ex: duality via Frobenius-structure} satisfies the following property
\begin{equation}\label{eq: involutive property of duality}
(\verd^\la)^\vee\circ (\verd^\la)^{-1}\cong \id_{\bfC}: \bfC\to \bfC.
\end{equation}

It follows that in this case
\begin{equation}\label{eq: involutive property of duality adm}
((\verd^\la)^{\adm})^2\cong \id_{\bfC^\adm}.
\end{equation}

More generally, a symmetric structure on a Frobenius algebra $(\bfC,\la)$ is an isomorphism $\sigma_\la\cong \id_\bfC$ as algebra automorphisms. (Note that this is stronger than requiring $\sigma_\la\cong \id_\bfC$ as plain functors.) Note that \eqref{eq: involutive property of duality} and \eqref{eq: involutive property of duality adm} continue to hold in this generality.
\end{remark}

\begin{lemma}\label{lem: adm duality and functors} 
Suppose $F: \bfC\to \bfD$ is an $\bfA$-linear functor of dualizable $\bfA$-linear categories with an $\bfA$-linear right adjoint $F^R$. Let $F^o=(F^R)^\vee: \bfC^\vee\to\bfD^\vee$ be the conjugate functor. 
Suppose both $\bfC$ and $\bfD$ admit self-duality $\verd_\bfC: \bfC^\vee\cong \bfC$ and $\verd_\bfD: \bfD^\vee\cong \bfD$, and suppose we are given an isomorphism
$F\circ \verd_\bfC\cong \verd_\bfD\circ F^o$. 
Then there is a natural isomorphism of functors
\[
F^R\circ (\verd_\bfD)^\adm\cong (\verd_\bfC)^\adm\circ  (F^R|_{(\bfD^\adm)^{\op}}).
\]
\end{lemma}
\begin{proof}
For $d\in \bfD^\adm$, let $d^*\in (\bfD^\vee)^\adm$ be the corresponding object. We need to show that $F^R(\verd_\bfD(d^*))\cong \verd_{\bfC}(F^R(d)^*)$. It is enough to show that for every $\phi\in \bfC^\vee$,
\[
\Hom_{\bfC}(\verd_{\bfC}(\phi), F^R(\verd_\bfD(d^*)))\cong \Hom_{\bfC}(\verd_{\bfC}(\phi), \verd_{\bfC}((F^R(d))^*)).
\]
By assumption, the left hand side is isomorphic to $\Hom_{\bfD^\vee}( F^o(\phi), d^*)$, while the right hand side isomorphic to $\Hom_{\bfC^\vee}(\phi, (F^R(d))^*)$. As $(F^o)^R=F^\vee$, it remains to prove that $F^\vee(d^*)\cong (F^R(d))^\vee$. But this follows from $F\circ F_d^L\cong (F_d\circ F^R)^L$.
\end{proof}

\begin{lemma}\label{lem-adm-object-in-Frob-alg}
Assume that $\bfC$ is a commutative $\bfA$-algebra equipped with a Frobenius structure $\la$ as in  \Cref{ex: duality via Frobenius-structure}. 
If $\la$ admits an $\bfA$-linear right adjoint $\la^R$, then $\frobdual^\la$ is $\bfA$-admissible and for every $\bfA$-module $\bfB$ and $b\in \bfB$, we have
\begin{equation}\label{eq:Adm imply linear}
(\id_\bfB\otimes \la)^R (b)\cong b\boxtimes_{\bfA}\frobdual^\la.
\end{equation} 
Without assuming that $\la$ admits an $\bfA$-linear right adjoint, then for $c\in\bfC$ the following are equivalent.
\begin{enumerate}
\item\label{lem-adm-object-in-Frob-alg-1} $c$  is $\bfA$-admissible. 
\item\label{lem-adm-object-in-Frob-alg-2} For every $\bfA$-module $\bfB$ and $b\in\bfB$, there is a natural isomorphism in $\bfB\otimes_\bfA \bfC$
\[
b\boxtimes_\bfA c\cong \Hom^{\bfB\otimes_\bfA\bfC/\bfC}(c^{*,\la},(\id_\bfB\otimes \la)^R (b)).
\]
\item\label{lem-adm-object-in-Frob-alg-3} For every commutative $\bfA$-algebra $\bfB$ and $b\in \bfB$, there is a natural isomorphism 
\[
b\boxtimes_\bfA c\cong\underline{\Hom}(\mathbf{1}_{\bfB}\boxtimes_{\bfA}c^{*,\la}, (\id_\bfB\otimes \la)^R (b)),
\]
where the internal hom is taken in $\bfB\otimes_\bfA\bfC$.
\item\label{lem-adm-object-in-Frob-alg-4} The isomorphism in \eqref{lem-adm-object-in-Frob-alg-3} holds for $\bfB=\bfC$ and $b=c^{*,\la}$.
\end{enumerate}
\end{lemma}
\begin{proof}
For the first statement, note that $\mathbf{1}_\bfA$ is clearly $\bfA$-admissible and if $\la^R$ exists as an $\bfA$-linear functor, then $\frobdual^\la=\la^R(\mathbf{1}_\bfA)$ is $\bfA$-admissible (see \Cref{lem:basic cpt and adm} \eqref{lem:basic cpt and adm-3}). In addition, in this case $(\id_\bfB\otimes \la)^R=\id_\bfB\otimes \la^R$, giving \eqref{eq:Adm imply linear}.

Next, we deduce \eqref{lem-adm-object-in-Frob-alg-2} from \eqref{lem-adm-object-in-Frob-alg-1}. For every $x\in \bfC\otimes_\bfA \bfD$, we need to show that there is a canonical isomorphism
\[
\Map(x, b\boxtimes_\bfA c)\cong  \Map( (\id_\bfB\otimes m)( x\boxtimes_{\bfA} c^{*,\la}), (\id_\bfB\otimes \la)^R(b))\cong \Map((\id_\bfB\otimes \la\circ m)(x \boxtimes_\bfA c^{*,\la}),b).
\]
Given $x\to b\boxtimes_\bfA c$, we obtain 
\[
(\id_\bfB\otimes \la\circ m)(x \boxtimes_\bfA c^{*,\la})\to (\id_\bfB\otimes \la\circ m)(b\boxtimes_\bfA c\boxtimes_\bfA c^{*,\la})\cong b \boxtimes_\bfA  e_\bfC(c\boxtimes_\bfA c^{*})\to b \boxtimes_\bfA \mathbf{1}_\bfA\cong b
\]
and given $(\id_\bfB\otimes \la\circ m)(x \boxtimes_\bfA c^{*,\la})\to b$, or equivalently $(\id_\bfB\otimes e_\bfC)(x\boxtimes_\bfA c^{*})\to b$,
we obtain
\[
x\cong  (\id_\bfB\otimes e_\bfC\otimes \id_\bfC)(\id_\bfB\otimes\id_\bfC\otimes u_\bfC)(x\boxtimes_\bfA \mathbf{1}_\bfA)\to (\id_\bfB\otimes e_\bfC\otimes \id_\bfC)(x\boxtimes_\bfA c^{*}\boxtimes_\bfA c)\to b\boxtimes_\bfA c.
\] 
Again since  \eqref{eq: Zorro for admissible objects} is homotopic to the identity map, the above two constructions give the desired isomorphism.

Clearly, \eqref{lem-adm-object-in-Frob-alg-2} implies \eqref{lem-adm-object-in-Frob-alg-3} and \eqref{lem-adm-object-in-Frob-alg-3} implies \eqref{lem-adm-object-in-Frob-alg-4}. Finally we show that \eqref{lem-adm-object-in-Frob-alg-4} implies \eqref{lem-adm-object-in-Frob-alg-1}. First, there is the tautological map $c\otimes c^{*,\la}\to \omega^\la$, giving $e_\bfC(c\boxtimes_\bfA c^*)\to \mathbf{1}_\bfA$. On the other hand, there is a canonical map $u_\bfC^\la:=(\verd^\la\otimes \id_\bfC)(u_\bfC)\to \underline{\Hom}(\mathbf{1}_\bfC\boxtimes_\bfA c^{*,\la}, (\id_\bfC\otimes \la)^R(c^{*,\la}))$ given by
\[
(\id_\bfC\otimes \la\circ m)( u_\bfC^\la \boxtimes_\bfA c^{*,\la})\cong (\id_\bfC\otimes e_\bfC)(u_\bfC\otimes \id_\bfC)(\mathbf{1}_\bfA \boxtimes_\bfA c^{*,\la})\cong c^{*,\la}.
\]
Now if the natural morphism $c^{*,\la}\boxtimes_\bfA c\to\underline{\Hom}( \mathbf{1}_\bfC\boxtimes c^{*,\la}, (\id_\bfC\otimes \la)^R(c^{*,\la}))$ induced by $(c^{*,\la}\boxtimes_\bfA c)\otimes (\mathbf{1}_\bfC\boxtimes_{\bfA} c^{*,\la})=c^{*,\la} \boxtimes_{\bfA}(c\otimes c^{*,\la})\to(\la\otimes\id_\bfC)^R(c^{*,\la})$ is an isomorphism, then we obtain
\[
e_\bfC(c\boxtimes_\bfA c^*)\to \mathbf{1}_\bfA,\quad  u_\bfC\to c^*\boxtimes_\bfA c.
\]
It is a routine work to check that the maps $c\to c$ and $c^*\to c^*$ induced by the above two maps as in the definition of  \eqref{eq: Zorro for admissible objects} are homotopic to the identity. This shows that $c$ is $\bfA$-admissible.
 
\end{proof}

\begin{remark}\label{rem:unit adm}
Suppose we are in the situation as in \Cref{ex: duality via Frobenius-structure}. If $\mathbf{1}_\bfC$ is $\bfA$-admissible, then so is $\omega^\la$. In this case $\la^R$ is $\bfA$-linear so \eqref{eq:Adm imply linear} holds. In addition, every dualizable object in $\bfC$ is $\bfA$-admissible.
\end{remark}

\subsubsection{Horizontal traces}\label{SS: horizontal trace}
Let $\bfR$ be a symmetric monoidal category. The trace of an endomorphism of a dualizable object in $\bfR$ is a classical notion.
Namely, if $X$ is a dualizable object in $\bfR$ equipped with an endomorphism $f: X\to X$, its trace $\mathrm{tr}(X,f)\in \End(\mathbf{1}_\bfR)$ is given by the composition
\begin{equation}\label{eq-ordinary-trace}
 \mathbf{1}_\bfR\xrightarrow{u} X^\vee\otimes X\xrightarrow{\id_{X^\vee}\otimes f} X^\vee\otimes X\stackrel{\sw}{\cong} X\otimes X^\vee\xrightarrow{e}\mathbf{1}_\bfR,
\end{equation}
where $u$ (resp. $e$) denotes the unit (resp. evaluation) map of $X$. We discuss this notion in the symmetric monoidal category $\lincat_\bfA$.

Let $\phi: \bfC\to\bfC$ be an $\bfA$-linear endofunctor of $\bfC$ and let $K_\phi\subset \bfC^\vee\otimes_\bfA\bfC$ be the ``kernel" representing $\phi$ as from \eqref{eq: kernel of an endofunctor}.
Then the (horizontal) trace of $\phi$ is
an object in $\bfA$ defined as
\begin{equation}\label{eq: def of hor trace}
\mathrm{tr}(\bfC/\bfA, \phi):=e_{\bfC/\bfA}(\sw(K_\phi)).
\end{equation}
This is sometimes also called the Hochschild homology of $\phi$, see \Cref{ex-usual-Hochschild-homology-of-DG-cat} below.
We shall also consider
\begin{equation}\label{eq: def of cat center}
Z(\bfC/\bfA,\phi):=\Hom_{\bfC^{\vee}\otimes_\bfA\bfC/\bfA}(u_\bfC, K_\phi)\in\bfA.
\end{equation}
In particular, if $\phi=\id_{\bfC}$, we write
\begin{equation}\label{eq: hor trace and center for id}
\mathrm{tr}(\bfC/\bfA)=\mathrm{tr}(\bfC/\bfA,\id_{\bfC}),\quad  Z(\bfC/\bfA)=Z(\bfC/\bfA,\id_{\bfC}).
\end{equation}
The object $Z(\bfC/\bfA)\in \bfA$ is sometimes also called the center of $\bfC$. It has a natural $E_2$-algebra structure in $\bfA$ (e.g. see \cite[\textsection{D.1.3.3}]{Lurie.SAG}). 
In addition, $\mathrm{tr}(\bfC/\bfA)$ is naturally a left module over the underlying $E_1$-algebra of $Z(\bfC/\bfA)$.

\begin{remark}\label{rem: action of center of objects}
We note that for every $c\in \bfC$, we have $( e_{\bfC/\bfA}\otimes \id)(c \boxtimes_\bfA u_\bfC)=c$, which induces a canonical morphism of $E_1$-algebras
\begin{equation}\label{eq: center acts on objects}
Z(\bfC/\bfA)=\End_{\bfC^\vee\otimes_\bfA\bfC}(u_\bfC)\xrightarrow{\id_c\boxtimes_\bfA (-)} \End_{\bfC\otimes_\bfA\bfC^\vee\otimes_\bfA\bfC}(c\boxtimes_\bfA u_\bfC)\to \End_{\bfC}(c).
\end{equation}
\end{remark}

\begin{remark}\label{rem: functoriality between center}
Suppose $\bfA=\Mod_\La$. Let $i:\bfC\subset \bfD$ be a fully faithful embedding, both of which are dualizable. Then we note that there is a natural ``restriction" map
\begin{equation}\label{eq: functoriality between center}
Z(\bfD/\bfA)\to Z(\bfC/\bfA),
\end{equation}
defined as follows: We have
\[
\bfD^\vee\otimes_\La \bfD\xrightarrow{i^\vee\otimes \id} \bfC^\vee\otimes \bfD\xleftarrow{\id\otimes i} \bfC^\vee\otimes \bfC.
\]
Then $(i^\vee\otimes \id)(u_\bfD)\cong (\id\otimes i)(u_\bfC)$ is the kernel $K_i$ representing $i$. Then we have
\[
Z(\bfD/\bfA)\to \End(K_i)\leftarrow Z(\bfC/\bfA).
\]
It is known that $\id\otimes i: \bfC^\vee\otimes\bfC\to \bfC^\vee\otimes \bfD$ is fully faithful. (This follows from \Cref{L:hochschild-homology-vs-cohomology} below, relying on \Cref{lem-tensor-preserve-fully-faithfulness}.) Therefore, we can reverse the above left-pointed arrow, giving the desired map.

Clearly, for $c\in \bfC$, \eqref{eq: center acts on objects} and \eqref{eq: functoriality between center} fit into the following commutative diagram
\[
\xymatrix{
Z(\bfD/\bfA)\ar[d]\ar[r] &  \End_{\bfD}(i(c))\\
Z(\bfC/\bfA)\ar[r] & \End_{\bfC}(c)\ar[u].
}\]
\end{remark}

\begin{example}\label{ex-usual-Hochschild-homology-of-DG-cat}
Let $\bfA=\Mod_\La$ for some $E_\infty$-ring $\La$.
If $\bfC=\lmodu_A$ is the category of left $A$-modules for an associative $\La$-algebra $A$, then $\bfC$ is dualizable with $\bfC^\vee\cong \lmodu_{A^{\rev}}=\rmodu_{A}$. The unit $u:\Mod_\La\to \lmodu_{A^{\rev}}\otimes_\La \lmodu_{A}\cong \lmodu_{A\otimes_\La A^{\rev}}$ is given by the $(A\otimes_\La A^{\rev})$-module $A$ and the evaluation map $e: \lmodu_{A}\otimes_\La\lmodu_{A^{\rev}}\to\Mod_\La$ is given by $M\mapsto A\otimes_{A\otimes_\La A^{\rev}}M$.
Then 
\[
\mathrm{tr}(\lmodu_{A})=A\otimes_{A\otimes_\La A^{\rev}} A, \quad Z(\lmodu_{A})=\End_{A\otimes_\La A^{\rev}} (A)
\] 
is the usual Hochschild homology and cohomology of the algebra $A$.
\end{example}

\begin{example}\label{ex: Serre functor-2}
Let $\bfC$ be a dualizable $\bfA$-module. Recall the notion of the Serre functor of $\bfC$ from \Cref{ex: Serre functor}. We have
\[
\Hom_\bfA(\mathrm{tr}(\bfC),\mathbf{1}_\bfA)= Z(\bfC, S_{\bfC}).
\]
In particular, if $\bfC$ is $0$-Calabi-Yau, then $Z(\bfC)=\Hom_\bfA(\mathrm{tr}(\bfC),\mathbf{1}_\bfA)$. 
\end{example}

As before, if $\bfA$ is clear from the context, for simplicity we sometimes just write $(\mathrm{tr}(\bfC,\phi), Z(\bfC,\phi))$ instead of $(\mathrm{tr}(\bfC/\bfA,\phi),Z(\bfC/\bfA,\phi))$.

We review some basic functoriality of the horizontal trace construction. First given $(\bfC,\phi_\bfC)$ and $(\bfD,\phi_\bfD)$, we have a canonical isomorphism in $\bfA$
\begin{equation}\label{eq:trace tensor product}
\mathrm{tr}(\bfC\otimes_\bfA\bfD,\phi_\bfC\otimes\phi_\bfD)\cong \mathrm{tr}(\bfC,\phi_\bfC)\otimes \mathrm{tr}(\bfD,\phi_\bfD).
\end{equation}

\begin{proposition}\label{lem-functoriality-duality-data}
Let $F: \bfC\to \bfD$ be such a $1$-morphism in $\lincat_\bfA^{\mathrm{dual}}$ as in \Cref{rem-duality-dualizable category}. Then there are natural transformations of functors
\[
\al_F: (F^o\otimes F)\circ u_\bfC\Rightarrow u_\bfD,\quad \beta_F: e_{\bfC}\Rightarrow e_{\bfD}\circ (F\otimes F^o).
\] 
Let $\phi_\bfC:\bfC\to\bfC$ and $\phi_\bfD:\bfD\to\bfD$ be $\bfA$-linear endomorphisms. Suppose there is a natural transformation of functors $\eta: F\circ \phi_{\bfC}\Rightarrow \phi_{\bfD}\circ F$. Then
there is a natural morphism in $\bfA$
\begin{align*}
\mathrm{tr}(F,\eta)&\colon \mathrm{tr}(\bfC,\phi_\bfC)=e_\bfC((\id_{\bfC^\vee}\otimes \phi_\bfC)(u_\bfC))\to e_\bfD((F\phi_\bfC\otimes F^o)(u_\bfC))\\
&\xrightarrow{\eta} e_\bfD((\phi_\bfD F\otimes F^o)(u_\bfC))\to e_{\bfD}((\phi_\bfD\otimes\id_{\bfD^\vee})(u_{\bfD}))=\mathrm{tr}(\bfD,\phi_\bfD).
\end{align*}
Suppose we in addition we further have $G: \bfB\to \bfC$ and $\delta: G\circ \phi_\bfB\Rightarrow \phi_\bfC \circ G$. Then 
\[
\mathrm{tr}(G,\delta)\circ \mathrm{tr}(F,\eta)\cong \mathrm{tr}(G\circ F, \delta\circ G(\eta)).
\]
\end{proposition}
\begin{proof}
By definition $(\id_{\bfC^\vee}\otimes F)(u_\bfC)\cong ( F^\vee\otimes \id_\bfD)(u_\bfD)$ in $\bfC^\vee\otimes \bfD$, which gives $(F^o\otimes F)(u_\bfC)\to u_\bfD$ by adjunction. The second natural transformation arises as $e_{\bfC}\to e_{\bfC}\circ ((F^R\circ F)\otimes \id_{\bfC^\vee}) \cong e_{\bfD}\circ (F\otimes F^o)$.
\end{proof}

\begin{remark}\label{trem-trace-as-symmetric-monoidal-functor}
Secretly behind the above discussions, there is a symmetric monoidal $2$-category (or called symmetric monoidal bi-category by some people) structure on $\lincat_\bfA$.
In fact, $\lincat_\bfA$ admits a symmetric monoidal $(\infty,2)$-category structure, and \eqref{eq:trace tensor product} and \Cref{lem-functoriality-duality-data} together can be upgraded as a symmetric monoidal functor from certain symmetric monoidal category $\End(\lincat_\bfA)$ to $\bfA$. 
Informally, $\End(\lincat_\bfA)$ is the symmetric monoidal category with objects being $(\bfC,\phi_\bfC:\bfC\to\bfC)$ and with morphisms from $(\bfC,\phi_\bfC)$ to $(\bfD,\phi_{\bfD})$ being $(F:\bfC\to \bfD, \eta: F\circ \phi_{\bfC}\Rightarrow \phi_{\bfD}\circ F)$ as in  \Cref{lem-functoriality-duality-data}. 

However, we will not systematically explore this approach in this article (but refer to \cite{hoyois2017higher}). On the one hand, we do not want to systematically review the formalism of (symmetric monoidal) $(\infty,2)$-categories (as we are not capable of). On the other hand, when we move to categorical trace, we will implicitly make use some $3$-categorical structures. 
\end{remark}

Now we assume that $\bfC$ is dualizable in $\lincat_\bfA$ as before.
Let $c\in \bfC$ be an $\bfA$-compact object. Then \Cref{lem-functoriality-duality-data} supplies a map in $\bfA$
\begin{equation}\label{eq-Chern-character-1}
\mathrm{tr}(F_c,\id)\colon \mathbf{1}_\bfA=\mathrm{tr}(\bfA)\to  \mathrm{tr}(\bfC),
\end{equation}
which is also denoted as $\mathrm{ch}(c)$ when regarded as a point in $\Map_\bfA(\mathbf{1}_\bfA,\mathrm{tr}(\bfC,\id_\bfC))$, usually called the Chern character of $c$ in literature. 
Note that for $(F: \bfC\to \bfD, \eta=\id)$ as in \Cref{lem-functoriality-duality-data} and $c\in \bfC$ being $\bfA$-compact, we have
\begin{equation}\label{eq-functoriality-Chern-character}
\mathrm{tr}(F,\id)(\mathrm{ch}(c))=\mathrm{ch}(F(c)).
\end{equation}

On the other hand, if $c\in \bfC$ is $\bfA$-admissible, we obtain again by  \Cref{lem-functoriality-duality-data} a map
\begin{equation}\label{eq:abstract character}
\mathrm{tr}(F_c^L, \id): \mathrm{tr}(\bfC)\to \mathbf{1}_\bfA,
\end{equation}
giving a point in $Z(\bfC, S_\bfC)$, which we call the character of $c$ and sometimes denote it by $\Theta_c$ (see \Cref{ex: character of adm repn} below). This map can also be described explicitly as the composition of the two maps in \eqref{eq: unit-counit-adjunction-admissible}. 

Note that for $(F: \bfC\to \bfD, \eta=\id)$ as in \Cref{lem-functoriality-duality-data}, and $d\in \bfD$ being $\bfA$-admissible, we have
\begin{equation}\label{eq-functoriality-abstract-Character}
\Theta_{F^R(d)}=\Theta_d\circ \mathrm{tr}(F,\id).
\end{equation}

\begin{example}\label{ex: character of adm repn}
In the case $\bfA=\Mod_\La$ and $\bfC$ is the category of smooth representations of a $p$-adic group, and $c=\pi$ is an admissible representation of $G$, the above map \eqref{eq:abstract character} is nothing but the usual character $\Theta_\pi$ of $\pi$, which is a conjugate invariant distribution on $G$.
\end{example}

\begin{remark}\label{rem:twisted Chern character}
The Chern character construction admits a twisted generalization.
Let $\phi:\bfC\to\bfC$ be an $\bfA$-linear endomorphism. 
Suppose $c\in \bfC$ is an $\bfA$-compact object equipped a morphism $\phi_c: c\to \phi(c)$. Then again \Cref{lem-functoriality-duality-data} gives
\begin{equation}\label{eq-twisted-Chern-character-1}
\mathrm{tr}(c,\phi_c):\mathbf{1}_\bfA\to  \Hom_{\bfC/\bfA}(c,c)\xrightarrow{\phi_c} \Hom_{\bfC/\bfA}(c,\phi(c))\to\mathrm{tr}(\bfC/\bfA,\phi),
\end{equation}
and we define the twisted Chern character $\mathrm{ch}(c,\phi_c)$ as corresponding point in $\Map_\bfA(\mathbf{1}_\bfA,\mathrm{tr}(\bfC/\bfA,\phi))$. This  implies that $\mathrm{ch}(c,\phi_c)$ is $\End_\bfA (\mathbf{1}_\bfA)$-linear in $\phi_c$. 
Note that when $\phi=\id_{\bfC}$, $\mathrm{ch}(c,\phi_c)$ is the image of $\phi_c\in \Map_{\bfA}(\mathbf{1}_\bfA, \Hom_{\bfC/\bfA}(c,c))$ under the map $\Hom_{\bfC/\bfA}(c,c)\to \mathrm{tr}(\bfC,\id_{\bfC})$. In fact, in this case we consider the full $\bfA$-linear subcategory of $\bfC$ spanned by $c$, which is equivalent to the category of right $B=\Hom_{\bfC/\bfA}(c,c)$-modules. Then $\mathrm{ch}(c,-)$ is nothing but the natural map from $B\otimes_{B\otimes B^{\rev}}B$ to $\mathrm{tr}(\bfC,\id_\bfC)$.

There is the following functoriality of twisted Chern characters.
Let $(F: \bfC\to \bfD, \eta: F\circ\phi_\bfC\Rightarrow\phi_{\bfD}\circ F)$ as in \Cref{lem-functoriality-duality-data}, and let $c\in\bfC^\cpt$ with $\phi_c: c\to\phi_\bfC(c)$. Then $d=F(c)$ is compact in $\bfD$ and we write $\phi_d=\eta\circ F(\phi_c): F(c)\to F(\phi_\bfC(c))\to \phi_\bfD(F(c))$.
It is clear that
\begin{equation}\label{eq-functoriality-twisted-Chern-character}
\mathrm{tr}(F,\eta)(\mathrm{ch}(c,\phi_c))=\mathrm{ch}(d, \phi_d).
\end{equation}
\end{remark}

We recall the well-known localization sequence of Hochschild homology. (See also \cite[Theorem 3.4]{hoyois2017higher}.)
We assume we are In the situation as in \Cref{prop-gluing-unit-in-semi-orthogonal-decomposition}. suppose that there is an $\bfA$-linear functor $\phi_\bfC: \bfC\to \bfC$. Let 
\[
\phi_\bfM:= F^R\circ \phi_\bfC\circ F:\bfM\to\bfM,\quad \phi_\bfN:= G\circ \phi_\bfC\circ G^R:\bfN\to\bfN.
\] 
By adjunction, we obtain
\[
\eta: F\circ \phi_\bfM \Rightarrow \phi_\bfC\circ F,  \quad \delta:  G\circ \phi_\bfC\Rightarrow  \phi_\bfN \circ G.
\]

\begin{proposition}\label{cor-hochschild-semi-orthogonal}
Then there is a canonical fiber sequence in $\bfA$
\[
\mathrm{tr}(\bfM,\phi_\bfM)\xrightarrow{\mathrm{tr}(F,\eta)} \mathrm{tr}(\bfC,\phi_\bfC)\xrightarrow{\mathrm{tr}(G,\delta)} \mathrm{tr}(\bfN,\phi_\bfN).
\]
If in addition, $(F(\bfM), G^R(\bfN))$ form a semi-orthogonal decomposition of $\bfC$,
and the adjunction $\phi_\bfC\circ G^R\Rightarrow G^R\circ \phi_\bfN$ is an isomorphism, then the above sequence canonically splits.
\end{proposition}
\begin{proof}
First note that the natural transformation $e_{\bfM}\Rightarrow e_\bfC\circ (F\otimes F^o)$ from \Cref{lem-functoriality-duality-data} is an isomorphism of functors, so is the natural transformation
$e_\bfC\circ (G^R\otimes (G^o)^R)\Rightarrow e_\bfN \circ (G\otimes G^o)\circ  (G^R\otimes (G^o)^R)\cong e_\bfN$.

Then by our assumption 
\[
e_{\bfM}\circ (\phi_\bfM\otimes\id_{\bfM^\vee})\Rightarrow e_{\bfC}\circ (F\circ\phi_\bfM\otimes F^o)\Rightarrow e_{\bfC}\circ (\phi_\bfC\circ F\otimes F^o)
\] 
is an isomorphism, and so is
\[
e_{\bfC}\circ (\phi_\bfC\circ G^R\otimes (G^o)^R)\Rightarrow e_{\bfN}\circ (G\circ \phi_\bfC\circ G^R\otimes G^o\circ (G^o)^R)\cong e_{\bfN}.
\]

Now we apply $e_\bfC\circ (\phi_\bfC\otimes\id_{\bfC^\vee})$ to the fiber sequence in \Cref{prop-gluing-unit-in-semi-orthogonal-decomposition}. 
Note that
\[
e_{\bfM}((\phi_\bfM\otimes\id_{\bfM^\vee})u_\bfM)\cong e_{\bfC}((F\otimes F^o)(\phi_\bfM\otimes\id_{\bfM^\vee})u_\bfM)\to e_{\bfC}((\phi_\bfC\otimes\id_{\bfC^\vee})(F\otimes F^o)u_\bfM) \to e_{\bfC}((\phi_\bfC\otimes\id_{\bfC^\vee})u_\bfC)
\] 
is identified with $\mathrm{tr}(F,\eta)$, and the middle map is an isomorphism.

On the other hand, one checks (using various adjunctions) that the composed map 
\begin{multline*}
e_{\bfC}((\phi_\bfC\otimes\id_{\bfC^\vee})u_\bfC)\to e_{\bfC}((\phi_\bfC\otimes\id_{\bfC^\vee})(G^{R}\otimes (G^o)^R)u_\bfN)\\
\to e_{\bfN}((G\otimes G^o)(\phi_\bfC\otimes\id_{\bfC^\vee})(G^{R}\otimes (G^o)^R)u_\bfN)\cong e_{\bfN}((\phi_\bfN\otimes\id_{\bfN^\vee})u_{\bfN})
\end{multline*} 
is identified with $\mathrm{tr}(G,\delta)$, and the middle map is an isomorphism. This gives the desired fiber sequence. 

The last statement clearly follows as $\mathrm{tr}(G^R, G^R\circ \phi_\bfN\cong \phi_\bfC\circ G^R)$ gives the desired splitting by \Cref{lem-functoriality-duality-data}.
\end{proof}

\subsubsection{Compactly generated categories}\label{SS: cpt gen category}
We let $\bfA=\Mod_\La$. We will write $\bfC/\bfA$ as $\bfC/\La$.
Recall (e.g. \cite[\textsection{D.7}]{Lurie.SAG} or using \eqref{eq:colimit-of-dual-category} and \eqref{eq-unit-for-colimit-category}) that  every object $\bfC\in \cptcat_{\bfA}$ is dualizable as an object in $\lincat_\bfA$, and some constructions in \Cref{sec:dual-categories}-\Cref{SS: localization sequences} can be made more explicitly.

If $\bfC = \ind(\bfC_0)$ for some $\bfC_0 \in \catid_\La$ we can identify the dual $\bfC^{\vee}$ of $\bfC$ with $\ind(\bfC_0^{\op})$. Explicitly, the evaluation map $\bfC^\vee \otimes_\La \bfC \rightarrow \modu_\La$ is given by the unique continuous extension of the functor given by the unique continuous extension of the functor
\[
   \bfC_0^{\op} \otimes_\La \bfC_0 \rightarrow \modu_\La,\quad   (c,d)  \rightarrow \Hom_{\bfC}(c,d).
\]
Next, let $\bfC= \ind(\bfC_0)$ and $\bfD = \ind(\bfD_0)$ be objects of $\lincat_\La$ which are compactly generated and let $F\colon \bfC \rightarrow \bfD$ be a continuous functor that preserves compact objects. We have a tautological functor $F_0^{\op}\colon \bfC_{0}^{\op} \rightarrow \bfD_0^{\op}$, which after taking its ind-extension, gives the conjugate functor $F^{o}$ as mentioned before.

\begin{remark}\label{rem: cpt object in dual category}
The above discussion says that for a compactly generated $\La$-linear category, there is a canonical equivalence $(\bfC^\cpt)^{\op}\cong (\bfC^\vee)^\cpt$ given by $c\mapsto \Hom(c,-)$. To emphasize the different roles played by $c$ in $\bfC$ and $\bfC^\vee$, we sometimes also write $c^{\op}$ for $\Hom(c,-)$.
Beware that this equivalence $(\bfC^\cpt)^{\op}\cong (\bfC^\vee)^\cpt$ is different from the restriction to $(\bfC^\cpt)^\op$ of the functor $\bfC^\op\to\bfC^\vee$ from \eqref{eq: functor of taking vee}. See \Cref{ex: explicit Serre and admissible} below.
\end{remark}

\begin{lemma}\label{lem: char of adm obj in cg cat}
An object $d\in \bfC$ is $\Mod_\La$-admissible (or simply called admissible in this case) if and only if for every $c\in\bfC^\cpt$, $\Hom_{\bfC/\La}(c,d)$ is a perfect $\La$-module. 
\end{lemma}
\begin{proof}If $F_d^L$ exists, then it sends compact objects in $\bfC$ to compact objects in $\Mod_\La$. So $\Hom_{\bfC/\La}(c,d)=\Hom_{\bfC/\La}(F_d^L(c), \La)$ is perfect for every $c\in\bfC^\cpt$. Conversely, if $\Hom_{\bfC/\La}(c,d)$ is perfect, then we may define $F_d^L$ on compact objects as $F_d^L(c)=\Hom_{\bfC/\La}(c,d)^*$. 
\end{proof}

\begin{remark}\label{ex: explicit Serre and admissible}
For a compactly generated category $\bfC$, the Serre functor (as defined in \Cref{ex: Serre functor}) can be explicitly given as follows: for every compact object $c\in \bfC^{\cpt}$, 
\[
\Hom_{\bfC/\La}(d, S_\bfC(c))=\Hom_{\bfC/\La}(c,d)^*,\quad \forall d\in \bfC.
\]
Here $(-)^*$ is understood as in \eqref{eq:abstract smooth dual 3} (for $\bfA=\Mod_\La$).

It also follows from the proof of \Cref{lem: char of adm obj in cg cat} that for $d$ admissible, $d^*\in \bfC^\vee$ can be given as
\[
d^*=\Hom(d, S_\bfC(-))\colon \bfC^\cpt\to \Mod_\La.
\]
In particular, if $d\in \bfC^\cpt\cap \bfC^\adm$, then as objects in $\bfC^\vee$, we have
\begin{equation*}\label{eq: smooth dual vs tautological dual}
d^*\cong S_{\bfC}^\vee(d^{\op}).
\end{equation*}

Now if $\verd: \bfC^\vee\to \bfC$ is a $\La$-linear equivalence, it restricts to an equivalence
\begin{equation}\label{eq:duality on compact objects}
\verd^\cpt: (\bfC^\cpt)^{\op}\cong (\bfC^\vee)^\cpt\xrightarrow{\verd} \bfC^\cpt.
\end{equation}

In particular,  if $\bfC$ is Calabi-Yau, i.e. $S_\bfC=\id_\bfC$, then for $c\in \bfC^\cpt\cap \bfC^\adm$, we have 
\[
\verd^\adm(c)=\verd^\cpt(c).
\]
Here $\verd^\adm$ is defined as in \eqref{eq: dual of admissible objects}, and for $c\in\bfC^\cpt$, we write $\verd^\cpt(c)$ instead of $\verd^\cpt(c^\op)$ to simplify the notation. 
Without Calabi-Yau assumption, in general $\verd^\adm(c)\neq \verd^\cpt(c)$.
\end{remark}

\begin{remark}\label{ex: explicit Serre and admissible-2}
Let $\bfC$ be as in \Cref{ex: explicit Serre and admissible}.
Now in addition, we assume that $\bfC$ is symmetric monoidal with a Frobenius structure as in \Cref{ex: duality via Frobenius-structure}. Then the functor $(\verd^\la)^\cpt$ satisfies
\[
 \Hom_{\bfC/\La}(c,d)=\la(d\otimes (\verd^\la)^\cpt(c)),\quad \forall c\in \bfC^\cpt, d\in\bfC.
\]
This also gives a relation between $S_\bfC, \verd^\la$ and $(-)^{*,\la}$ on compact objects
\begin{equation}\label{eq-Serre-functor-via-duality}
S_{\bfC}(c)\cong ((\verd^\la)^\cpt(c))^{*,\la}, \quad \forall c\in \bfC^\cpt.
\end{equation}
We also notice that in this case, by \eqref{eq: involutive property of duality}, we have
\begin{equation}\label{eq: involutive property of duality cpt}
((\verd^\la)^\cpt)^2\cong \id_{\bfC^\cpt}.
\end{equation}
\end{remark}

\begin{example}\label{ex: self-duality of cpt gen rigid monoidal cat}
Suppose that $\bfC$ is a compactly generated semi-rigid $\La$-linear category. (See \Cref{def: rigid monoidal category} for the notion of rigid monoidal category and \Cref{def.rigid.from.GR} below for the more general notion of $\La$-linear rigid monoidal category.)  In this case, compact objects are both left and right dualizable objects in $\bfC$, and by \Cref{rem-duality-datum-as-plain-cat} below, the ind-completion of the functor obtained by the restriction of $\Hom(\mathbf{1}_{\bfC},-)$ to $\bfC^\cpt$ defines a Frobenius structure of $\bfC$. We let 
\[
\verd^{\mathrm{sr}}:\bfC^\vee\to\bfC
\]
denote the self duality induced by this Frobenius structure. Note that it is completely determined by the monoidal structure of $\bfC$.

It is easy to see that the induced self-duality $(\verd^{\mathrm{sr}})^{\cpt}: (\bfC^{\cpt})^{\op} \cong \bfC^\cpt$ is given by $c\mapsto c^\vee $. Here $c^\vee $ denotes the right dual of $c$ in $\bfC$ (i.e. the one equipped with $\mathbf{1}_\bfC\to  c\otimes c^\vee$ and $c^\vee\otimes c\to \mathbf{1}_\bfC$ forming a duality datum).

It also follows that the Serre automorphism $\sigma_{\mathrm{sr}}$ of $\bfC$ induced by the above Frobenius structure (see \eqref{eq:Serre automorphism of Frob algebra}) is given by $c\mapsto (c^{\vee})^\vee$ on compact objects. Therefore, to endow $\bfC$ with a symmetric structure amounts to choosing isomorphisms $(c^{\vee})^{\vee}\cong c$ functorial in $c$ and compatible with the monoidal structure. In literature, such additional structure on compactly generated rigid monoidal category is usually called a pivotal structure. See \Cref{def: pivotal semi-rigid monoidal} for a generalization.

We also write  $\omega^{\mathrm{sr}}$ for the object $\omega^{\la}$ as defined in \Cref{ex: duality via Frobenius-structure}. If $\mathbf{1}_\bfC$ is compact, so $\bfC$ is rigid, then
\[
\omega^{\mathrm{sr}}\cong S_\bfC(\mathbf{1}_\bfC).
\]
(Note that the Serre functor $S_\bfC$ of $\bfC$ and the Serre automorphism $\bfC$ are different.) In addition, for $c\in \bfC^\cpt$, 
\[
S_{\bfC}(c)=S_{\bfC}(\mathbf{1}_\bfC)\otimes c=\omega^{\mathrm{sr}}\otimes \sigma^{\mathrm{sr}}(c).
\] 
and for $c\in \bfC^\cpt\cap \bfC^{\adm}$, we have
\[
c^{\mathrm{sr},*}=  \omega^{\mathrm{sr}}\otimes \sigma^{\mathrm{sr}}(c^\vee).
\]
\end{example}

Now, let  $K_0(\bfC^\cpt)$ be the usual Grothendieck K-group of the stable category $\bfC^\cpt$ (the quotient of the free abelian group generated by objects in $\bfC^\cpt$ by the subgroup generated by $[c]-[c']-[c'']$ for any fiber sequence $c'\to c\to c''$ in $\bfC^\cpt$). On the other hand, $\tr(\bfC)$ is a $\La$-module and we let $H^0\mathrm{tr}(\bfC)$ denote its zeroth cohomology (which is the same as $\pi_0\Map_{\Mod_\La}(\La,\tr(\bfC))$).
\begin{proposition}\label{rem-Chern-character}
\begin{enumerate}
\item\label{rem-Chern-character-1} The Chern character construction defines a homomorphism
\begin{equation}
\label{eq-chern-character}
\mathrm{ch}: K_0(\bfC^\cpt)\to H^0\mathrm{tr}(\bfC).
\end{equation}

\item\label{rem-Chern-character-2} For $(F: \bfC\to \bfD, \eta=\id)$ as in \Cref{lem-functoriality-duality-data}, the following diagram is commutative
\begin{equation}\label{eq-functoriality-Chern-character}
\xymatrix{
K_0(\bfC^{\cpt})\ar[r]\ar_{K_0(F)}[d] & H^0\mathrm{tr}(\bfC)\ar^{\mathrm{tr}(F,\id)}[d]\\
K_0(\bfD^{\cpt}) \ar[r] & H^0\mathrm{tr}(\bfD).
}
\end{equation}

\item\label{rem-Chern-character-3} Suppose $\bfM\to\bfC\to\bfN$ is a localization sequence in $\lincat_\La$ which in addition induces a semi-orthogonal decomposition in the sense of \Cref{rem:localization sequence}. Suppose $\bfM,\bfC,\bfN$ are compactly generated. Then $(F,G^R)$ induce $K_0(\bfM^\cpt)\oplus K_0(\bfN^\cpt)\cong K_0(\bfC^\cpt)$ and the Chern character \eqref{eq-chern-character} is compatible with this decomposition and the decomposition from \Cref{cor-hochschild-semi-orthogonal}.
\end{enumerate}
\end{proposition}
Although it is well-known, we sketch a proof for completeness, as the ingredients of the proof are also needed in \Cref{prop-abstract-character additive}.
Part \eqref{rem-Chern-character-2} is sometimes known as the abstract Grothendieck-Riemann-Roch formula.  We also mention that in fact, the Chern character construction can be lifted to a map from the connective $K$-theory spectra of $\bfC$ to $\mathrm{tr}(\bfC)$ (and even to the cyclic homology of $\bfC$). We will not need such refined version. 
\begin{proof}
It is clear that once Part \eqref{rem-Chern-character-1} is established, Part \eqref{rem-Chern-character-2} follows from \Cref{lem-functoriality-duality-data}.

It remains to show that for a fiber sequence $c'\to c\to c''$, we have $\mathrm{ch}(c)=\mathrm{ch}(c')+\mathrm{ch}(c'')$.
Let $S_2\bfC\subset \Fun(\La^2_1,\bfC)$ be the category of fiber sequences in $\bfC$. This is again a compactly generated $\La$-linear category with $(S_2\bfC)^\cpt=S_2\bfC^\cpt$. There is a fully faithful embedding $F: \bfC\to S_2\bfC$ sending $c$ to $c\xrightarrow{\id_c} c\to 0$, with the right adjoint $F^R$ sending $c'\to c\to c''$ to $c'$. The right orthogonal complement of $F(\bfC)$ in $S_2\bfC$ then is still $\bfC$, with $G^R:\bfC\to S_2\bfC$ sending $c$ to $0\to c\xrightarrow{\id_c}c$ (which preserves compact objects). The left adjoint of $G^R$ then is given by $G: S_2\bfC\to \bfC$ sending $c'\to c\to c''$ to $c''$. 
Then by \Cref{cor-hochschild-semi-orthogonal} we have the natural isomorphism $\mathrm{tr}(\bfC)\oplus \mathrm{tr}(\bfC)\xrightarrow{\mathrm{tr}(F)\oplus\mathrm{tr}(G^R)}  \mathrm{tr}(S_2\bfC)$, with inverse map given by
\begin{equation}\label{eq-splitting-HH(S2C)}
\mathrm{tr}(F^R)\oplus \mathrm{tr}(G): \mathrm{tr}(S_2\bfC)\to \mathrm{tr}(\bfC)\oplus \mathrm{tr}(\bfC),
\end{equation} 
which sends  $\mathrm{ch}(c'\to c\to c'')=\mathrm{ch}(c')+\mathrm{ch}(c'')$.

There is another functor $p: S_2\bfC\to \bfC$ sending $c'\to c\to c''$ to $c$, which induces $\mathrm{tr}(p): \mathrm{tr}(S_2\bfC)\to \mathrm{tr}(\bfC)$. As $p\circ F\simeq p\circ G^R\simeq \id_{\bfC}$, we see that under the isomorphism \eqref{eq-splitting-HH(S2C)}, $\mathrm{tr}(p)$ restricts to the identity map of each direct factor. The claim follows.

Finally for Part \eqref{rem-Chern-character-3}, it is enough to notice that if $G^R$ sends compact objects to compact objects, so is $F^R$. It follows that $K_0(\bfM^\cpt)\oplus K_0(\bfN^\cpt)\cong K_0(\bfC^\cpt)$ and by Part \eqref{rem-Chern-character-2} the Chern character is compatible with the direct sum decomposition.
\end{proof}

Dually, we have the following statement for admissible objects and characters.

\begin{proposition}\label{prop-abstract-character additive}
The assignment 
\[
(c\in \bfC^{\adm})\mapsto (\Theta_c: H^0\mathrm{tr}(\bfC)\to \La)
\] 
(where we recall $\Theta_c$ is the map $\mathrm{tr}(F_c^L,\id)$ from \eqref{eq:abstract character}) induces a bilinear map
\[
\langle\cdot,\cdot\rangle_{\bfC}\colon (K_0( \bfC^{\adm})\otimes\La)\otimes_\La H^0\mathrm{tr}(\bfC)\to \La,
\]
such that for $(F: \bfC\to \bfD, \eta=\id)$ as in \Cref{lem-functoriality-duality-data}, then
\[
\langle F^R(d), a \rangle_{\bfC}=\langle d, \mathrm{tr}(F,\id)(a) \rangle_\bfD.
\]
\end{proposition}
\begin{proof}
Given the first statement, the second statement is just a reformulation of \eqref{eq-functoriality-abstract-Character}.

For the first statement, we need to show that for a fiber sequence $c'\to c\to c''$ of admissible objects in $\bfC$, we have $\Theta_c=\Theta_{c'}+\Theta_{c''}$. 
We still make use the constructions from the proof of \Cref{rem-Chern-character}.
Note that a fiber sequence of  admissible objects is an admissible object in $S_2\bfC$. 

We note that in fact $(G^R)^R=p$. Then \eqref{eq-functoriality-abstract-Character} gives
$\Theta_{c'}= \Theta_{c'\to c\to c''}\circ \mathrm{tr}(F)$, $\Theta_{c}=  \Theta_{c'\to c\to c''}\circ \mathrm{tr}(G^R)$. 
Under the isomorphism \eqref{eq-splitting-HH(S2C)}, this means that
\[
\Theta_{c'\to c\to c''}(a_1,a_2)=\Theta_{c'}(a_1)+\Theta_{c}(a_2),\quad a_1,a_2\in H^0\mathrm{tr}(\bfC).
\]
Let $G^L: \bfC\to S_2\bfC$  be the left adjoint of $G$, which sends $c$ to $c[-1]\to 0\to c$. We claim that under the isomorphism \eqref{eq-splitting-HH(S2C)}, we have
$\mathrm{tr}(G^L)(a)=(-a,a)$. Indeed, as $G\circ G^L=\id$, we see that $\mathrm{tr}(G^L)(a)=(b,a)$ for some $b\in \mathrm{tr}(\bfC)$. One the other hand $p\circ G^L=0$ so $b+a=0$. This shows that $b=-a$.

Now by \eqref{eq-functoriality-abstract-Character} again we have $\Theta_{c''}=\Theta_{c'\to c\to c''}\circ \mathrm{tr}(G^L)$, i.e. 
\[
\Theta_{c'\to c\to c''}(-a,a)=\Theta_{c''}(a), \quad a\in H^0\mathrm{tr}(C).
\]
Comparing the above two displayed equations, we see that $\Theta_c=\Theta_{c}+\Theta_{c'}$, as desired.
\end{proof}

Note that in the course of the proof of the above lemma, we also have proved the following statement.
\begin{corollary}\label{lem: shift on trace}
Let $[-1]: \bfC\to\bfC$ be the functor given by looping. Then $\mathrm{tr}([-1],\id): \mathrm{tr}(\bfC)\to\mathrm{tr}(\bfC)$ is given by multiplication by $-1$.
\end{corollary}

\begin{example}\label{ex: pairing between compact and admissible}
Let $c\in \bfC^\cpt$ and $d\in \bfC^{\adm}$. Then $F_d^L\circ F_c=\Hom(c,d)^\vee$ (see \Cref{rem-smooth-and-proper}). It follows that
\[
\Theta_d(\mathrm{ch}(c))= \dim \Hom_{\bfC}(c,d).
\]
Here $\dim \Hom_{\bfC}(c,d)$ is the dimension of $\Hom_\bfC(c,d)$, regarded as a dualizable object in the symmetric monoidal category $\Mod_\La$ (which is nothing but the Euler characteristic of $\Hom_\bfC(c,d)$ if $\La$ is a field of characteristic zero).
\end{example}

\subsubsection{$2$-dualizability and trace formula}
\begin{definition}\label{def-2-dualizable-category}
A dualizable $\bfA$-linear category $\bfC$ is called smooth (over $\bfA$) if $u_\bfC$ admits an $\bfA$-linear right adjoint, and is called proper (over $\bfA$) if $e_\bfC$ admits an $\bfA$-linear right adjoint. A dualizable $\bfA$-linear category $\bfC$ is called $2$-dualizable if it is smooth and proper.
\end{definition}

\begin{remark}\label{rem: smooth/proper vs cpt/adm}
Note that $u_{\bfC}$ (resp.  $e_{\bfC}$) admits an $\bfA$-linear right adjoint if and only if $e_{\bfC}$ (resp. $u_\bfC$) admits an $\bfA$-linear left adjoint. Namely, if $u_\bfC$ admits an $\bfA$-linear right adjoint $u_\bfC^R$, then let
\[
T_\bfC: \bfC=\bfA\otimes_\bfA\bfC\xrightarrow{u_\bfC\boxtimes_\bfA \id_\bfC} \bfC\otimes_\bfA\bfC^\vee\otimes_\bfA\bfC\xrightarrow{\id_\bfC\boxtimes_\bfA u_\bfC^R\circ \mathrm{sw}} \bfC\otimes_\bfA\bfA=\bfC,
\]
which is usually called the dual Serre functor. Then 
\[
e_\bfC^L=(\id_{\bfC^\vee}\otimes T_\bfC)\circ \mathrm{sw}\circ u_\bfC.
\]
On the other hand, if $e_{\bfC}$ admits an $\bfA$-linear right adjoint $e_\bfC^R$,
then 
\[
u_\bfC^L=e_{\bfC}\circ \mathrm{sw}\circ  (S_\bfC\otimes\id_{\bfC^\vee}),
\]
where $S_\bfC$ is the Serre functor of $\bfC$ as before.

It follows that $\bfC$ is smooth (resp. proper) over $\bfA$ if and only if $u_{\bfC}$ is $\bfA$-compact (resp. $\bfA$-admissible) as an object in $\bfC\otimes_\bfA\bfC^\vee$.
\end{remark}

\begin{example}\label{ex: cpt equal to adm in sm and proper cat}
Suppose $\bfC$ is smooth. Then by \Cref{lem: dual hom for admissible objects}, for an $\bfA$-admissible object $c$ we have
$\Map(c,-)= \Map(u_{\bfC}, (c^{\vee} \boxtimes (-))$ preserving colimits. So $c$ is compact.

On the other hand, suppose $\bfA=\Mod_\La$, and $\bfC$ is compactly generated. Then $\bfC$ is proper if and only if for every $c,d\in \bfC^\cpt$, $\Hom_{\bfC}(c,d)\in\Perf_\La$. It follows (by \Cref{ex: explicit Serre and admissible}) that in this case every compact object is admissible.

It follows that for a smooth and proper compactly generated $\La$-linear category $\bfC$, compact objects and admissible objects coincide.
\end{example}

\begin{theorem}\label{thm: classical-trace-formula}
Let $\bfC$ be a $2$-dualizable $\bfA$-linear category, with two right adjointable (in $\lincat_\bfA$) endomorphisms $\phi_i$, and an isomorphism $\eta: \phi_1\circ \phi_2\cong \phi_2\circ \phi_1$ of $\bfA$-linear functors. Then 
$\mathrm{tr}(\bfC,\phi_i)$ is dualizable in $\bfA$, and we have
\[
\mathrm{tr}(\mathrm{tr}(\bfC,\phi_1), \mathrm{tr}(\phi_2, \eta^{-1}))=\mathrm{tr}(\mathrm{tr}(\bfC,\phi_2), \mathrm{tr}(\phi_1,\eta)).
\]
\end{theorem}

We shall not recall its proof here as we will discuss a more general trace formula in \Cref{SS: trace formula}. But given \Cref{trem-trace-as-symmetric-monoidal-functor}, it is a direct consequence of the main result of \cite{Campbell.Ponto}. 

\subsection{Categorical trace}\label{sec:categorical-trace-formalism}
In this article, we will also need a different type of trace construction, known as the vertical trace, or categorical trace. 
Let us first review the general formalism. 

\subsubsection{Vertical trace}\label{sec-hoch-homology}
As before, let $\bfR$ denote a symmetric monoidal category. Let $A$ and $B$ be two associative algebras in $\bfR$. By \cite[Proposition 4.6.3.11]{Lurie.higher.algebra}\footnote{Note that assumption $(\star)$ of \emph{loc. cit.} is not essential, as explained before \cite[Notation 4.6.3.3]{Lurie.higher.algebra}.}, an $A\mbox{-}B$-bimodule can also be regarded as a left $(A\otimes B^{\rev})$-module or a right $(B\otimes A^{\rev})$-module, where $A^{\rev}$ (resp. $B^{\rev}$) is the algebra $A$ (resp. $B$) with the multiplication reversed.
For an associative algebra $A$, and an $A\mbox{-}A$-bimodule $F$,
the Hochschild homology of $F$, if exists, is defined as 
\begin{equation}\label{eq:hochschild-homology-of-an-algebra}
\tr(A,F)= A\otimes_{A\otimes A^{\rev}}F \in \bfR.
\end{equation}
We write
\begin{equation}\label{eq:universal-trace}
[-]_F: F\to \tr(A,F)
\end{equation}
for the natural morphism, sometimes called the universal trace morphism. 

On the other hand, there always exists the Hochschild complex of $F$ defined as
\begin{equation}\label{eq:standard-hochschild-complex}
    \HH(A,F)_\bullet=\barcons(A)_\bullet\otimes_{A\otimes A^{\rev}}F=A^{\otimes\bullet}\otimes F,
\end{equation}
regarded as a simplicial object $\Delta^{\op}\to \bfR$. Explicitly, on the level of simplicies and morphisms, for every $n\geq 0$ we have an equivalence
\begin{equation}\label{eq:nth-term-of-Hochschild-complex}
\HH(A,F)_n \simeq A^{\otimes n}\otimes F.
\end{equation}
Informally under this identification, for $0<i<n$ the face map $d^{\HH}_{i}$ is given by the multiplication map applied to the $i$-th and the $(i+1)$-th factors in $A^{\otimes n}$, and the face map $d^{\HH}_0$ is given by multiplying the first factor in $A^{\otimes n}$ to $F$ from the right and the face map $d^{\HH}_n$ is given by multiplying the $n$-th factor in $A^{\otimes n}$ to $F$ from the left. 
If $\bfR$ admits geometric realizations and the tensor product preserves geometric realizations in each variable, then the Hochschild homology of $F$ exists and can be computed as the geometric realization of the Hochschild complex.

\begin{remark}\label{rem-trace-2-category}
Associated to a symmetric monoidal $(\infty,1)$-category $\bfR$, there is a symmetric monoidal $(\infty,2)$-category $\mathrm{Morita}(\bfR)$ whose objects are associative algebras in $\bfR$ and whose morphism categories are given by categories of bimodules:
\[
\map_{\mathrm{Morita}(\bfR)}(A,B) = {_{B}}\BMod_{A}
\]
and compositions are given by the relative tensor products (assuming relative tensor products exist in $\bfR$). (Do not confuse it with the category $\mathrm{BMod}(\bfR)$ from \Cref{sec:relative-tensor-product}.)
Every $A$-bimodule $F$ gives an endomorphism of $A$ in $\mathrm{Morita}(\bfR)$. For example, when $\bfR=\Mod_\La$, there is a full embedding 
\[
\mathrm{Morita}(\Mod_\La)\subset\lincat_\La
\] 
of symmetric monoidal ($2$-)categories by sending $A$ to $\lmodu_A$ and $M$ to the functor $M\otimes_B(-): \lmodu_B\to\lmodu_A$.

Now for general $\bfR$, every algebra $A$ is a dualizable object in $\mathrm{Morita}(\bfR)$. Under the equivalence ${_{A}}\BMod_{A}\cong {_{\mathbf{1}_\bfR}}\BMod_{A^{\rev}\otimes A}\cong {_{A\otimes A^{\rev}}}\BMod_{\mathbf{1}_\bfR}$, the natural $A$-bimodule structure on $A$ itself gives unit and evaluation maps
\begin{equation}\label{eq-coevaluation-evaluation-module-of-A}
 A^u\in {_{A\otimes A^{\rev}}}\BMod_{\mathbf{1}_\bfR}, \quad A^e\in {_{\mathbf{1}_\bfR}}\BMod_{A^{\rev}\otimes A},
\end{equation}
which identify the dual of $A$ (in $\mathrm{Morita}(\bfR)$) as $A^{\rev}$. (Note that our notations are different from \cite[\textsection{4.6.3}]{Lurie.higher.algebra}).
Then $\tr(A,F)$ is nothing but the trace $F$ in the sense of \eqref{eq-ordinary-trace}, regarded as an endomorphism of $A$ in $\mathrm{Morita}(\bfR)$. This justifies our choice of notations. 

However, as explained in \Cref{trem-trace-as-symmetric-monoidal-functor}, we will not systematically use this approach. 
\end{remark}

\begin{example}\label{ex:relative-tensor-product-as-Hochschild}
  If the $A$-bimodule $F=M\otimes N$, where $M$ is a left $A$-module and $N$ a right $A$-module. Then 
  \[
  \tr(A,F)=A\otimes_{A\otimes A^{\rev}}(M\otimes N) \cong N\otimes_A M.
  \]
  In fact, $\HH(A,F)_\bullet\cong\barcons_A(N,M)_\bullet$.
\end{example}

\begin{example}\label{E:twisted-trace}
Of particular importance in this paper is the following type of bimodules. Let $\phi$ be an endomorphism of the algebra $A$. For an $A$-bimodule $F$ we will denote by ${}^{\phi}F$ the bimodule obtained by the same action on the right but with a pre-composition with $\phi$ for the left action. In this case we will also denote the Hochschild homology of the bimodule ${}^{\phi}A$ by $\tr(A,\phi)$. That is,
\begin{equation}\label{eq:twisted-trace-of-A}
\tr(A,\phi) \simeq A\otimes_{A\otimes A^{\rev}} {}^{\phi}A.
\end{equation}
In this case, we sometimes just write $[-]_{\phi}$ instead of $[-]_{{}^\phi A}$ for simplicity.
\end{example}

\begin{remark}\label{remark-hochschild-cohomology}
For an $A\mbox{-}A$-bimodule $F$, one can also form its Hochschild cohomology 
\[
\Hom_{A\otimes A^{\rev}}(A, F).
\]
We will not make use of this notion.
\end{remark}

\subsubsection{Functoriality of vertical traces}\label{sec-transfer.functors.convolution} 

Now let $F_A\in {_{A}}\BMod_{A}$ and $F_B\in {_{B}}\BMod_{B}$ be two bimodules. Assume that we are given a left dualizable $M\in {_{A}}\BMod_{B}$ together with a morphism of bimodules
\begin{equation}\label{eq:bimodule-and-2-morphism}
    \alpha\colon M\otimes_{B}F_{B}\rightarrow F_{A}\otimes_{A}M.
\end{equation}
Then we can associate to $(M,\al)$ a morphism in $\bfR$
\begin{equation}\label{eq:class-map-of-module}
\tr(M,\alpha)\colon \tr(B,F_B)\rightarrow \tr(A,F_A),
\end{equation}
given by
\begin{align*}
\tr(B,F_B) &= B\otimes_{B\otimes B^{\rev}} F_{B}  \xrightarrow{u_M\otimes \id}
(N\otimes_{A} M) \otimes_{B\otimes B^{\rev}} F_{B}  \simeq A\otimes_{A\otimes A^{\rev}} (M\otimes_{B} F_B\otimes_{B} N)\\
& \xrightarrow{\id\otimes \alpha\otimes \id}
A\otimes_{A\otimes A^{\rev}} (F_{A}\otimes_{A}M\otimes_{B} N)
\xrightarrow{\id\otimes \id\otimes e_M}
A\otimes_{A\otimes A^{\rev}} F_{A} 
=\tr(A,F_A),
\end{align*}
where the isomorphism 
\begin{equation*}\label{eq:cyclic tensor product}
(N\otimes_{A} M) \otimes_{B\otimes B^{\rev}} F_{B}  \simeq A\otimes_{A\otimes A^{\rev}} (M\otimes_{B} F_B\otimes_{B} N)
\end{equation*}
can be established by the same way as in \Cref{ex:relative-tensor-product-as-Hochschild}.

In the particular case when $B=F_B=\mathbf{1}_{\bfR}$ is the unit object of $\bfR$, $\al$ is just a map $M\to F_A\otimes_AM$. Then the above definition of $\tr(M,\al)$ is simplified as
\begin{equation}\label{eq:class-of-modules}
    \mathbf{1}_{\bfR}   \xrightarrow{u_M}
N\otimes_{A} M  \xrightarrow{\id\otimes \alpha}
N\otimes_{A}F_A\otimes_{A}M\cong (M\otimes N)\otimes_{A\otimes A^{\rev}}F_A
\xrightarrow{e_M\otimes \id} \tr(A, F_A).
\end{equation}
In this case, we also denote $\tr(M,\al)$ as $[M,\al]_{F_A}$, thought as a point in the space $\Map(\mathbf{1}_{\bfR}, \tr(A,F_A))$.

\begin{example}\label{ex: vertical trace-simple-functoriality}
Let $\eta: F_1\to F_2$ be an $A$-bimodule homomorphism. Then we obtain a pair $(M,\al)$ where $M=A$ and $\al: M\otimes_AF_1\cong F_1\to F_2\cong F_2\otimes_AM$. It defines a morphism $\tr(M,\al): \tr(A,F_1)\to \tr(A,F_2)$. On the other hand, we may regard $\tr(A,-)$ as a functor from the category of $A$-bimodules in $\bfR$ to $\bfR$. We thus obtain another map $\tr(A,\eta): \tr(A,F_1)\to \tr(A,F_2)$. It is clear that $\tr(M,\al)$ and $\tr(A,\eta)$ are canonically identified. 
\end{example}

\begin{example}\label{ex: class-vs-usual-hochschild}
When $(A,F_A)=(B,F_B)=(\mathbf{1}_\bfR,\mathbf{1}_\bfR)$, an object $M\in\bfR$ regarded as an $A\mbox{-}B$-bimodule admits a left dual if and only if $M$ is dualizable in $\bfR$. In this case, $[M,\al]_{\mathbf{1}_\bfR}=\mathrm{tr}(M,\al)$ from \eqref{eq-ordinary-trace}.
\end{example}

\begin{example}\label{ex: unit as a class}
Let $M$ be a left dualizable $A$-module, with $N$ its left dual. Let $F=M\otimes N$, regarded as an $A$-bimodule. Let $\al: M\to F\otimes_AM$ be the map given by $M\cong M\otimes \mathbf{1}_\bfR\xrightarrow{\id_M\otimes u_M}M\otimes N\otimes_AM$. Then the map $[M,\al]_F: \mathbf{1}_\bfR\to \tr(A,F)=N\otimes_AM$ is nothing but $u_M$.
\end{example}

\begin{example}\label{example-class-module-vs-class-object}
  Suppose we are in the case \Cref{example-free-A-module-dual}. I.e. $M=A$, regarded as a left $A$-module over itself. 
  In this case, giving a left $A$-module morphism $\al: M\to F\otimes_AM$ is equivalent to giving a map $\al_0: \mathbf{1}_\bfR\to F$. Then we have the canonical equivalence of morphisms
  \[
  [A,\al]_{F}\cong [-]_{F}\circ \al_0\colon \mathbf{1}_\bfR\to \tr(A,F).
  \]
\end{example}

We recall the following basic statements. 

\begin{lemma}\label{lem: cyclic invariance of trace}
Let $M$ be an $A$-$B$-bimodule and  $N$ a $B$-$A$-bimodule. Then there is a canonical isomorphism
\begin{equation*}\label{eq: cyclic invariance of trace}
c:  \tr(A, M\otimes_B N)\cong \tr(B, N\otimes_A M),
\end{equation*}
functorial in $M$ and $N$.
\end{lemma}

\begin{lemma}\label{lemma:composition-bimodule-vs-map-trace}
Suppose we have three objects $(A,F_A), (B,F_B), (C,F_C)$ in $\mathrm{BMod}(\bfR)$, an $A\mbox{-}B$-bimodule $S$ with duality datum and $F_S: S\otimes_BF_B\to F_A\otimes_AS$, and a $B\mbox{-}C$-bimodule $T$ with duality datum and $F_T: T\otimes_CF_C\to F_B\otimes_BT$. Let $R=S\otimes_BT$ with induced $F_R=(F_S\otimes 1)\circ (1\otimes F_T)$, then $R$ admits a left dual as an $A\mbox{-}C$-bimodule and
\[
\tr(R,F_R)=\tr(S,F_S)\circ \tr(T,F_T): \tr(C,F_C)\to \tr(A,F_A).
\]
\end{lemma}
Note that in the case as in \Cref{ex: class-vs-usual-hochschild}, the above lemma recovers \eqref{eq:trace tensor product}.

\subsubsection{$2$-dualizability}\label{SS: 2-dual algebra}  
Let $A\in \alg(\bfR)$. We recall the following definitions from \cite[\textsection{4.6.4}]{Lurie.higher.algebra}.
\begin{definition}\label{def-2-dualizable-alg}
\begin{enumerate}
\item We call $A$ a proper algebra in $\bfR$ if $A^e\in {}_{\mathbf{1}_\bfR}\bmodu_{A\otimes A^{\rev}}$ (see \eqref{eq-coevaluation-evaluation-module-of-A}) admits a left dual. In this case, we write $S_A$ for its left dual, and write the unit and evaluation as
\[
\epsilon: A\otimes A^{\rev} \to S_A\otimes A^e,\quad \delta: A^e\otimes_{A\otimes A^{\rev}} S_A\to \mathbf{1}_\bfR.
\]
When regarding $S_A$ as an $A$-bimodule, it is usually called the Serre bimodule of $A$.
\item We call $A$ a smooth algebra in $\bfR$  if $A^u\in {}_{A\otimes A^{\rev}}\bmodu_{\mathbf{1}_\bfR}$ (see \eqref{eq-coevaluation-evaluation-module-of-A}) admits a left dual. In this case, we write $T_A$ for its left dual, and write the unit and evaluation as
\[
\mu:  \mathbf{1}_\bfR \to  T_A\otimes_{A\otimes A^{\rev}} A^u,\quad \nu: A^u\otimes T_A\to A\otimes A^{\rev}.
\]
When regarding $T_A$ as an $A$-bimodule, it is usually called the dual Serre bimodule of $A$. 
\item We call $A$ a $2$-dualizable algebra in $\bfR$, if $A$ is both proper and smooth. 
\end{enumerate}
\end{definition}

\begin{remark}
\begin{enumerate}
\item One should compare the above definition with \Cref{def-2-dualizable-category}. Both definitions are specializations of the notion of proper (resp. smooth, resp. $2$-dualizable) objects in a symmetric monoidal $2$-category. (The symmetric monoidal $2$-category behind \Cref{def-2-dualizable-category} is $\lincat_\bfA$ and behind \Cref{def-2-dualizable-alg} is $\mathrm{Morita}(\bfR)$.) A $\La$-algebra $A$ is a proper (resp. smooth, resp. $2$-dualizable) algebra in $\bfR=\Mod_\La$ the sense of \Cref{def-2-dualizable-alg} if and only if its left module category $\lmodu_A(\Mod_\La)\in\cptcat_\La$ is a proper (resp. smooth, resp. $2$-dualizable) $\La$-linear category in the sense of \Cref{def-2-dualizable-category}. But note that a $\La$-linear monoidal category $\bfA$ (i.e. an algebra in $\lincat_\La$) is a proper (resp. smooth, $2$-dualizable) algebra in $\lincat_\La$ is different from the under $\La$-linear category being proper (resp. smooth, $2$-dualizable). 
\item If $A$ is $2$-dualizable algebra in $\bfR$, then $S_A\otimes_AT_A\cong T_A\otimes_AS_A\cong A$ as $A$-bimodules.
\end{enumerate}
\end{remark}

\begin{notation}\label{not: simplified notation}
In the sequel, to simply notations,
when the algebra $A$ is clear from the context, we simply write $-\odot-$ instead of $-\otimes_A-$, and simply write
$\tr(F)$ or $\langle F\rangle$ instead of $\tr(A,F)$.
\end{notation}

Now suppose $A$ is $2$-dualizable.
Let $F_1$ and $F_2$ be two $A$-bimodules.
We still use $\nu$ to denote the morphism  
\begin{align}\label{eq:map nu TA}
\begin{split}
\langle F_1\odot T_A\odot F_2\rangle&= (A^u\otimes T_A)\otimes_{(A\otimes A^\rev)\otimes (A\otimes A^\rev)^\rev} (F_1\otimes F_2) \\
&\xrightarrow{\nu\otimes \id_{F_1\otimes F_2}}(A\otimes A^{\rev})\otimes_{(A\otimes A^\rev)\otimes (A\otimes A^\rev)^\rev} (F_1\otimes F_2)\cong \langle F_1\rangle \otimes \langle F_2\rangle.
\end{split}
\end{align}

We similarly still use $\epsilon$ to denote the morphism
\begin{align}\label{eq: map epsilon S_A}
\begin{split}
\langle F_1\rangle \otimes \langle F_2\rangle&\cong (A\otimes A^{\rev})\otimes_{(A\otimes A^\rev)\otimes (A\otimes A^\rev)^\rev} (F_1\otimes F_2)\xrightarrow{\epsilon\otimes \id_{F_1\otimes F_2}}\\
&(S_A\otimes A^e)\otimes_{(A\otimes A^\rev)\otimes (A\otimes A^\rev)^\rev} (F_1\otimes F_2)\cong \langle F_2\odot S_A\odot F_1\rangle.
\end{split}
\end{align}
Note that both maps \eqref{eq:map nu TA} and \eqref{eq: map epsilon S_A} are functorial in $F_1$ and $F_2$. 
In addition, there is the following crucial commutative diagram.
\begin{lemma}\label{lem: crucial comm diagram secondary trace}
The following diagram commutative
\begin{small}
\begin{equation}
\xymatrix{
\langle F_1\odot F_2\rangle \otimes\langle F_3\rangle \ar_{c\otimes\id}^\cong[d]&\ar_-{\eqref{eq:map nu TA}}[l]\langle F_1\odot F_2\odot T_A\odot F_3\rangle \cong \langle F_2\odot T_A\odot F_3\odot F_1\rangle\ar^-{\eqref{eq:map nu TA}}[r]&\langle F_2\rangle \otimes\langle F_3\odot F_1\rangle\ar_\cong^{\id\otimes c}[d]\\
\langle F_2\odot F_1\rangle \otimes\langle F_3\rangle \ar^-{\eqref{eq: map epsilon S_A}}[r]&\langle F_3\odot S_A\odot F_2\odot F_1\rangle \cong\langle F_1\odot F_3\odot S_A\odot F_2\rangle &\ar_-{\eqref{eq: map epsilon S_A}}[l]\langle F_2\rangle \otimes\langle F_1\odot F_3\rangle,
}
\end{equation}
\end{small}

\noindent where the isomorphism $c$ comes from the cyclic invariance of trace (see \Cref{lem: cyclic invariance of trace}).
In addition, if there are $A$-bilinear maps $F_i\to F'_i$, then the above diagram maps to the corresponding diagram for $(F'_1,F'_2,F'_3)$, and the resulting cubic diagram is commutative.
\end{lemma}

\begin{lemma}\label{lem-dualizability-of-categorical-trace}
Let $F$ be an $A$-bimodule, with a left dual $G$ (as $A$-bimodules), then $\tr(A,F)$ is dualizable in $\bfR$ with the dual $\tr(A,G)$,
with the unit and evaluation maps given by
\begin{equation}
\mathbf{1}_\bfR\xrightarrow{\mu}T_A\otimes_{A\otimes A^{\rev}}A\xrightarrow{\id_{T_A}\otimes u_F} T_A\otimes_{A\otimes A^{\rev}}(G\otimes_{A}F)\cong \langle F\odot T_A\odot G\rangle\xrightarrow{\nu} \langle F\rangle\otimes \langle G\rangle.
\end{equation}
\begin{align*}
\langle F\rangle \otimes \langle G\rangle \xrightarrow{\epsilon} \langle G\odot S_A\odot F\rangle \cong (F\otimes_AG)\otimes_{A\otimes A^{\rev}}S_A\xrightarrow{e_F\otimes\id_{S_A}} A\otimes_{A\otimes A^{\rev}}S_A\xrightarrow{\delta} \mathbf{1}_\bfR.
\end{align*}
\end{lemma}

\begin{lemma}\label{lem: trace dual morphism vs dual morphism trace}
Now let $f: F_1\to F_2$ be an $A$-bimodule map. Suppose $F_i$ admits a left dual $G_i$. Let 
\[
g: G_2\xrightarrow{u_{F_1}\otimes \id_{G_2}} G_1\otimes_A F_1\otimes_A G_2\xrightarrow{\id_{G_1}\otimes f\otimes \id_{G_2}} G_1\otimes_AF_2\otimes_AG_2\xrightarrow{\id_{G_1}\otimes e_{F_2}}G_1
\]
be the (left) dual of $f$.
Then under the duality from \Cref{lem-dualizability-of-categorical-trace}, the dual of $\tr(A,f): \tr(A,F_1)\to \tr(A,F_2)$ is given by $\tr(A, g): \tr(A,G_2)\to \tr(A,G_1)$. 
\end{lemma}

Now, let $F_1$ and $F_2$ be two $A$-bimodules, both of which admit left duals. Let
\[
\al: F_1\otimes_AF_2\to F_2\otimes_AF_1
\]
be an isomorphism of $A$-bimodules. Then we may form $\tr(A,F_2)\in \bfR$, equipped with an endomorphism
\[
\tr(F_1,\al): \tr(A,F_2)\to \tr(A,F_2).
\] 
By \Cref{lem-dualizability-of-categorical-trace}, $\tr(A,F_2)$ is dualizable in $\bfR$ so one can further form $\mathrm{tr}(\tr(A,F_2),\tr(F_1,\al))\in \End (\mathbf{1}_\bfR)$.  
On the other hand, by switching $F_1$ and $F_2$ one obtains $\mathrm{tr}(\tr(A,F_1),\tr(F_2,\al^{-1}))$.

\begin{theorem}\label{thm-second-trace-1}
Suppose $A$ is $2$-dualizable and $(F_1,F_2,\al)$ are as above. Then there is a canonical isomorphism in $\End(\mathbf{1}_\bfR)$
\[
\mathrm{tr}(\tr(A,F_1), \tr(F_2,\al^{-1}))\cong \mathrm{tr}(\tr(A,F_2), \tr(F_1,\al)).
\]
\end{theorem}
As the case for \Cref{thm: classical-trace-formula}, \Cref{thm-second-trace-1} (and all the lemmas above it)  is a direct consequence of the main result of \cite{Campbell.Ponto}, applied to the symmetric monoidal $2$-category $\mathrm{Morita}(\bfR)$ (as mentioned in \Cref{rem-trace-2-category}). We will discuss it in \Cref{SS: trace formula}.

\subsubsection{Categorical trace}

We fix a \emph{rigid} symmetric monoidal category $\bfR\in \calg(\lincat)$, e.g. $\bfR=\Mod_\La$ for 
an $E_\infty$-ring $\La$. We consider $\lincat_{\bfR}$, equipped with the natural symmetric monoidal structure. It will play the role of the ambient symmetric monoidal category $\bfR$ as in \Cref{sec-hoch-homology}-\Cref{SS: 2-dual algebra}. (We hope this shifting of notations will not cause any confusion.) 
Let $\bfA\in\alg(\lincat_{\bfR})$.
Let $\bfF$ be an $\bfA$-bimodule category. Then $\tr(\bfA,\bfF)$ always exists in $\lincat_{\bfR}$ and is sometimes called the categorical trace of $(\bfA,\bfF)$. We note that if $\bfA$ and $\bfF$ are compactly generated, so is $\tr(\bfA,\bfF)$.

\begin{remark}
Our notation/terminology is slightly abusive as $\tr(\bfA,\bfF)$ depends on the base category $\bfR$ we choose. As in the sequel and in the main body of the work we will alway fix such a base, we omit it from notation/terminology. (So $-\otimes-$ in the sequel will mean $-\otimes_{\bfR}-$, etc.)
\end{remark}

\begin{remark}\label{rem: imposing commutativity in trace}
Let $[-]_{\bfF}: \bfF\to \tr(\bfA,\bfF)$ be the canonical functor \eqref{eq:universal-trace}.  Clearly for $a\in \bfA$ and $f\in\bfF$, we have the canonical isomorphism $[a\otimes f]_\bfF\cong [f\otimes a]_\bfF$ in $\tr(\bfA,\bfF)$. In the case when $\bfF={}^\phi\bfA$ as considered in \Cref{E:twisted-trace}, we see that there is a canonical isomorphism $[\phi(a)\otimes b]_{\phi}\cong [b\otimes a]_{\phi}$. In particular, by setting $b=\bf{1}_{\bfA}$ to be the unit of $\bfA$, we obtain
\[
[\phi(a)]_{\phi}\cong [a]_{\phi}.
\]
It follows that the auto-equivalence of $\tr(\bfA,\phi)$ induced by $\phi: (\bfA,{}^\phi \bfA)\to (\bfA, {}^\phi \bfA)$ is canonically isomorphic to the identity functor.
\end{remark}

As mentioned in \Cref{rem-trace-2-category}, $\tr(\bfA,\bfF)$ can be regarded as the trace of the endomorphism $\bfF$ in the symmetric monoidal ($3$-)category $\mathrm{Morita}(\lincat_{\bfR})$. As such, besides the basic functoriality as discussed in \Cref{sec-transfer.functors.convolution}, there are some further adjointability/functoriality of the categorical trace construction, which we need to discuss.

We fix $(\bfA,\bfF_\bfA)$ and $(\bfB,\bfF_\bfB)$ as before. We first have the following generalization of \Cref{lem-functoriality-duality-data}.
\begin{proposition}\label{lem-functoriality-duality-data-generalization}
Let $\beta: \bfM_1\to \bfM_2$ be a morphism of $\bfA\mbox{-}\bfB$-bimodules. We assume that 
\begin{itemize}
\item $\bfM_i$ admits a left dual $\bfN_i$ as $\bfA\mbox{-}\bfB$-bimodules;
\item $\beta$ admits an $\bfA\mbox{-}\bfB$-linear right adjoint $\beta^R$. 
\end{itemize}
Suppose that there are $\bfA\mbox{-}\bfB$-bimodules maps $\al_i: \bfM_i\otimes_{\bfB} \bfF_\bfB\to \bfF_\bfA\otimes_{\bfA}\bfM_i$ and a natural transformation of functors $\eta: (\id_{\bfF_\bfA}\otimes \beta)\circ\al_1\Rightarrow \al_2\circ (\beta\otimes\id_{\bfF_{\bfB}})$. Then there is a natural transformation of functors
\[
\tr(\beta,\eta)\colon \tr(\bfM_1,\al_1)\Rightarrow \tr(\bfM_2,\al_2): \tr(\bfB,\bfF_{\bfB})\to \tr(\bfA,\bfF_{\bfA}).
\]
\end{proposition}
\begin{proof} We only mention the key point is that since $\beta$ admits an $\bfA\mbox{-}\bfB$-linear right adjoint $\beta^R$, we can define the conjugate functor $\beta^o: \bfN_1\to \bfN_2$, which is a $\bfB\mbox{-}\bfA$-bidmodule map, as the following composition:
\[
\bfN_1\xrightarrow{u_{\bfM_2}\otimes\id_{\bfN_1}} \bfN_2\otimes_{\bfA}\bfM_2\otimes_{\bfB}\bfN_1\xrightarrow{\id_{\bfN_2}\otimes\beta^R\otimes \id_{\bfN_1}} \bfN_2\otimes_{\bfA}\bfM_1\otimes_{\bfB}\bfN_1\xrightarrow{\id_{\bfN_2}\otimes e_{\bfM_1}}\bfN_2.
\] 
Then there are natural transformation of functors 
\[
(\beta^o\otimes \beta)\circ u_{\bfM_1}\Rightarrow u_{\bfM_2},\quad e_{\bfM_1}\circ (\beta\otimes \beta^o)\Rightarrow e_{\bfM_2}.
\] 
The desired natural transformation then follows the construction as in \Cref{lem-functoriality-duality-data}. We leave the details for readers.
\end{proof}

We have the following generalization of \Cref{cor-hochschild-semi-orthogonal}. The proof remains the same.
\begin{proposition}\label{generalization-hochschild-semi-orthogonal}
We let $(\bfA,\bfF)$ as before.
Let $\bfM_1\xrightarrow{F} \bfM_2\xrightarrow{G} \bfM_3$ be a localization sequence of left dualizable $\bfA$-modules in the sense of \Cref{rem:localization sequence}.
Let $\al_2: \bfM_2\to \bfF\otimes_\bfA\bfM_2$ be an $\bfA$-module functor and let $\al_1= (\id_{\bfF}\otimes F^R)\circ\al_2\circ F$, and $\al_3= (\id_{\bfF}\otimes G)\circ \al_2 \circ G^R$.
Then the sequence (from \Cref{lem-functoriality-duality-data-generalization})
\[
[\bfM_1,\al_1]_{\bfF}\to [\bfM_2,\al_2]_{\bfF}\to [\bfM_3,\al_3]_{\bfF}
\] 
is a fiber sequence in $\tr(\bfA,\bfF)$.
\end{proposition}

\begin{definition}\label{def: left proper and smooth module}
Let $\bfA$ and $\bfB$ be two algebras in $\lincat_{\bfR}$, and $\bfM$ is a
left dualizable $\bfA\mbox{-}\bfB$-bimodule. We say
$\bfM$ is left smooth if the functor \eqref{eq:unit-map-for-dual} admits a continuous right adjoint  $u_\bfM^R$ as a $\bfB$-bimodule map, and is left proper if the functor \eqref{eq:counit-map-for-dual} admits a continuous right adjoint  $e_\bfM^R$ as an $\bfA$-bimodule map. When $\bfB=\bfR$, we simply say $\bfM$ is left proper/smooth over $\bfA$.
\end{definition}

\begin{remark}\label{rem-3-category-structure}
\begin{enumerate}
\item\label{rem-3-category-structure-1} If $\bfA=\bfB=\bfR$, then the above notions specialize to the proper and smooth $\bfR$-linear categories as defined in \Cref{def-2-dualizable-category}.
\item\label{rem-3-category-structure-2} Suppose $\bfM$ admits a left dual $\bfN$ as an $\bfA\mbox{-}\bfB$-bimodule. If $\bfM$ is left proper and smooth, then $\bfN$ admits a left dual as an $\bfB\mbox{-}\bfA$-bimodule with the left dual given by $\bfM$. Indeed, the duality datum is just given by
\[
\bfA\xrightarrow{e_\bfM^R} \bfM\otimes_\bfB\bfN,\quad \bfN\otimes_\bfA\bfM\xrightarrow{u_\bfM^R} \bfB.
\]
\item\label{rem-3-category-structure-3} If $\bfM$ is left smooth and proper over $\bfA$, then $\bfM\otimes\bfN$ as an $\bfA$-bimodule  admits a left dual, given by $\bfM\otimes\bfN$ itself with the unit and evaluation maps being
\[
\bfA\xrightarrow{e_\bfM^R} \bfM\otimes_\La\bfN\xrightarrow{\id_\bfM\otimes u_\bfM\otimes\id_\bfN}( \bfM\otimes \bfN)\otimes_\bfA(\bfM\otimes\bfN).
\]
\[
(\bfM\otimes\bfN)\otimes_\bfA(\bfM\otimes\bfN)\xrightarrow{\id_\bfM\otimes u_\bfM^R\otimes\id_\bfN}\bfM\otimes\bfN\xrightarrow{e_\bfM}\bfA.
\]
Note that regarding $e_\bfM:\bfM\otimes\bfN\to\bfA$ as $\bfA$-bimodule map, its conjugate functor (as defined in the proof of \Cref{lem-functoriality-duality-data-generalization}) is $e_\bfM$ itself.

\item\label{rem-3-category-structure-4} The concepts in \Cref{def: left proper and smooth module} should be generalized as $2$-dualizability of morphisms in a symmetric monoidal $3$-category (such as $\mathrm{Morita}(\lincat_\bfR)$ as discussed in \Cref{rem-trace-2-category}).
\item\label{rem-3-category-structure-5} We note that instead of asking $u$ and $e$ to admit (continuous) right adjoints, one could ask them to admit left adjoints. This would lead to another type of $2$-dualizability. Our choice of \Cref{def: left proper and smooth module} is adapted to the applications.
\end{enumerate}
\end{remark}

\begin{example}\label{def.rigid.from.GR} 
Let $\bfA\in \alg(\lincat_\bfR)$. Then $\bfM=\bfA$ as a left $\bfA$-module is 
\begin{enumerate}
    \item\label{item-rigid.cts.right.adjoint} left proper if and only if the monoidal functor $m:\bfA\otimes \bfA \to \bfA$ admits a continuous right adjoint $m^R$ as an $\bfA$-bimodule map;
    \item\label{item-rigid.compact.unit} left smooth if and only if the unit object $\mathbf{1}_{\bfA}$ is compact.
\end{enumerate}
Note that these conditions together look exactly the same as the ones in the definition of rigid monoidal categories as in \Cref{def: rigid monoidal category}, except the tensor product here is taken in $\lincat_\bfR$ rather than in $\lincat$. 
In particular when $\bfR=\Mod_{\La}$, where $\La$ is the sphere spectrum, then $\bfA$ is rigid in the sense of \Cref{def: rigid monoidal category}. On the other hand, it is easy to see that for a symmetric monoidal functor $\bfR'\to \bfR$ between rigid monoidal categories, an algebra $\bfA\in \alg(\lincat_\bfR)$ satisfies the conditions as above if and only if so is its image in $\alg(\lincat_{\bfR'})$. In particular, an object $\bfA\in \alg(\lincat_\bfR)$ satisfying the above conditions is a rigid monoidal category in the sense of \Cref{def: rigid monoidal category}.
Therefore, we will call such $\bfA$ satisfying the above conditions a rigid $\bfR$-linear monoidal category. We also note that \Cref{lem-rigid-lax-is-strict} is applicable to rigid $\bfR$-linear monoidal categories.
\end{example}

We have the following useful observation.

\begin{lemma}\label{lem-adjoint-of-induced-trace-map}
Suppose $\bfM$ is left proper and smooth as an $\bfA\mbox{-}\bfB$-module, and suppose the functor $\al$ in \eqref{eq:bimodule-and-2-morphism} also admits an $\bfA\mbox{-}\bfB$-linear right adjoint $\al^R$, then $\tr(\bfM,\al)$ admits continuous right adjoint $\tr(\bfM,\al)^R$. In particular, if $(\bfB,\bfF_\bfB)=(\bfR,\bfR)$, then $[\bfM,\al]_{\bfF}$ is a compact object in $\tr(\bfA,\bfF_{\bfA})$.
\end{lemma}

\begin{example}\label{ex:compact to compact in trace}
We consider \Cref{example-class-module-vs-class-object} in the current set-up.
Assume that $\bfA$ is a rigid $\bfR$-linear monoidal category, and $\bfF$ an $\bfA$-bimodule. Let $\bfM=\bfA$ regarded as a left $\bfA$-module.
Recall that giving an $\bfR$-linear functor $\al_0:\bfR\to \bfF$ is equivalent to giving an object $X\in \bfF$. We denote the corresponding left $\bfA$-module morphism $\bfM\to \bfF\otimes_\bfA\bfM$ by $\al_X$. 
Then $\al_X$ admits an $\bfA$-linear right adjoint $\al_X^R$ if and only if $X$ is a compact object in $\bfF$ (by \Cref{lem-rigid-lax-is-strict}). It follows that $[X]_{\bfF}=[\bfA,\al_X]_{\bfF}$, regarded as a functor $\bfR\to \tr(\bfA,\bfF)$, admits a continuous right adjoint. That is, $[X]_{\bfF}$ is a compact object in $\tr(\bfA,\bfF)$. This can also be deduced from \Cref{L:hochschild-homology-vs-cohomology} \eqref{cor:right-adjointability-of-module-to-Trace} below. In any case, we see that when $\bfA$ is rigid, the universal trace map \eqref{eq:universal-trace} sends compact objects to compact objects.
\end{example}

\begin{lemma}\label{lem: categorified GRR}
Suppose $\bfM$ is a left proper and smooth $\bfA\mbox{-}\bfB$-bimodule, and $\al$ admits a continuous $\bfA\mbox{-}\bfB$-linear right adjoint $\al^R$ as in \Cref{lem-adjoint-of-induced-trace-map}. Let $\bfN$ be a left dual of $\bfM$. Let
\begin{small}
\begin{multline*}
\hspace*{-0.5cm} 
\delta: \bfN\otimes_\bfA \bfF_\bfA \xrightarrow{\id_{\bfN\otimes_{\bfA}\bfF_\bfA}\otimes e_\bfM^R}  \bfN\otimes_\bfA \bfF_\bfA\otimes_\bfA \bfM\otimes_\bfB\bfN\\
 \xrightarrow{\id_\bfN\otimes \al^R\otimes \id_\bfN}\bfN\otimes_\bfA\bfM\otimes_\bfB\bfF_\bfB\otimes_\bfB\bfN\xrightarrow{u_\bfM^R\otimes\id_{\bfF_\bfB\otimes_{\bfB} \id_{\bfN}}} \bfF_\bfB\otimes_\bfB\bfN.
\end{multline*}
\end{small}
Then we have a canonical isomorphism of functors
\[
 \tr(\bfN, \delta)\cong \tr(\bfM,\al)^R \colon \tr(\bfA,\bfF_\bfA)\to \tr(\bfB,\bfF_\bfB).
\]
\end{lemma}
Here we note that  thanks to \Cref{rem-3-category-structure} \eqref{rem-3-category-structure-2}, $\tr(\bfN, \delta)$ is well-defined.
\begin{proof}The desired isomorphism is a consequence of the following commutative diagram 
\begin{small}
\[
\xymatrix{
\bfF_\bfA\otimes_{\bfA\otimes\bfA^{\rev}}\bfA\ar^-{\id\otimes e_{\bfM}^R}[r] & \bfF_\bfA\otimes_{\bfA\otimes\bfA^{\rev}}(\bfM\otimes_\bfB\bfN)\ar^-\cong[r]& \bfB\otimes_{\bfB\otimes\bfB^{\rev}}(\bfN\otimes_\bfA\bfF_\bfA\otimes_\bfA\bfM)\ar^{\id\otimes \al^R}[d]\ar^-{\id\otimes\delta\otimes\id}[r]& \bfB\otimes_{\bfB\otimes\bfB^{\rev}}(\bfF_\bfB\otimes_\bfB\bfN\otimes_\bfA\bfM)\ar^-{\id\otimes u_\bfM^R}[d] \\
  & &\bfB\otimes_{\bfB\otimes\bfB^{\rev}}(\bfN\otimes_\bfA\bfM\otimes_\bfB\bfF_\bfB)\ar^-{\id\otimes u_\bfM^R\otimes\id}[r] &\bfB\otimes_{\bfB\otimes\bfB^{\rev}}\bfF_\bfB}
\]
\end{small}
\end{proof}

Together with \Cref{lemma:composition-bimodule-vs-map-trace}, we have the following result, which is sometimes referred as the categorified Grothendieck-Riemann-Roch formula.

\begin{corollary}\label{prop: categorified GRR}
Then for a left dualizable $\bfA$-module $\bfL$ equipped with $\bfA$-linear functor $\beta: \bfL\to \bfF_\bfA\otimes_\bfA \bfL$, we have
\[
\tr(\bfM,\al)^R([\bfL,\beta]_{\bfF_\bfA})\cong [\bfN\otimes_\bfA \bfL,\gamma]_{\bfF_\bfB},
\]
where $\gamma$ is the composed functor $\bfN\otimes_\bfA\bfL\xrightarrow{\id_\bfN\otimes \beta} \bfN\otimes_\bfA \bfF_\bfA\otimes_\bfA \bfL\xrightarrow{\delta\otimes\id_\bfL}\bfF_\bfB\otimes_\bfB\bfN\otimes_\bfA\bfL$, with the functor $\delta$ given in \Cref{lem: categorified GRR}.
\end{corollary}

\begin{example}\label{ex: relating horizontal trace and vertical trace}
Suppose we are in the situation as \Cref{ex:compact to compact in trace}. Let $\bfL$ be a left dualizable $\bfA$-module equipped with $\beta: \bfL\to \bfF\otimes_\bfA\bfL$ as in \Cref{prop: categorified GRR}. In this case $\delta$ is the right adjoint of the functor $\bfA\to \bfF$ given by $a\mapsto X\otimes a$. (Note that this is different from the functor $\al_X:\bfA\to\bfF$ which sends $a$ to $a\otimes X$.) 
It follows that we have the following isomorphism in $\bfR$
\begin{equation}\label{eq: endo of X as a HH homology-0}
\Hom_{\tr(\bfA,\bfF)}([X]_\bfF, [\bfL,\beta]_\bfF)\cong \mathrm{tr}(\bfL,\gamma).
\end{equation}
In particular, let $\bfF={}^\phi\bfA$ be as in \Cref{E:twisted-trace}. We write $\phi_{\bfL}: \bfL\to \bfL$ for the underlying $\bfR$-linear functor of $\beta: \bfL\to {}^{\phi}\bfL$. Then $\gamma:\bfL\to\bfL$ is given by $y\mapsto {}^\vee X\otimes\phi_{\bfL}(y)$, where ${}^\vee X$ is the (left) dual of $X$ in $\bfA$, i.e. the one equipped with a duality datum $\mathbf{1}_\bfA\to {}^\vee X\otimes X$ and ${}^\vee X\otimes X\to \mathbf{1}_\bfA$.
We further specialize to the following cases:
\begin{enumerate}
\item When $\bfL=\bfA$ with $\beta: \bfA\to {}^\phi\bfA$ given by $X\in \bfA$ so $[\bfL,\beta]_{{}^\phi\bfA}=[X]_{{}^\phi\bfA}$ (see \Cref{example-class-module-vs-class-object}). Then $\phi_\bfL: \bfA\to\bfA$ is given by $\phi_\bfL(a)=\phi(a)\otimes X$, and $\ga: \bfA\to \bfA$ is given by $\ga(a)={}^\vee X \otimes \phi(a)\otimes X$.
We thus obtain
\begin{equation}\label{eq: endo of X as a HH homology}
\End_{\tr(\bfA,\phi)}([X]_{{}^\phi\bfA})\cong \mathrm{tr}(\bfA,\ga).
\end{equation}
\item When $X=\mathbf{1}_\bfA$, we have
\[
\Hom_{\tr(\bfA,\phi)}([\mathbf{1}_\bfA]_{{}^\phi\bfA}, [\bfL,\beta]_{{}^\phi\bfA})\cong \mathrm{tr}(\bfL,\phi_\bfL).
\]
\item Combining the above two cases so $\bfL=\bfA$ and $X=\mathbf{1}_\bfA$, then we obtain
\[
\End_{\tr(\bfA,\phi)}([\mathbf{1}_\bfA]_{{}^\phi\bfA})\cong \mathrm{tr}(\bfA,\phi).
\]
\end{enumerate}
\end{example}

\begin{remark}
The above corollary in particular endows $\mathrm{tr}(\bfA,\phi)\in\bfR$ with an algebra structure given by $\End_{\tr(\bfA,\phi)}([\mathbf{1}_\bfA]_{{}^\phi\bfA})^{\rev}$ and $\mathrm{tr}(\bfL,\beta)$ with a left $\mathrm{tr}(\bfA,\phi)$-module structure.
On the other hand, \Cref{trem-trace-as-symmetric-monoidal-functor} implies that $\mathrm{tr}(\bfA,\phi)$ acquires another algebra structure and $\mathrm{tr}(\bfL,\phi)$ a module structure over this algebra structure. It turns out these two algebra and module structures coincide. Indeed, both algebra structures can be identified with the monad associated to the $\bfR$-linear functor $\bfR\to \tr(\bfA,\phi),\ \mathbf{1}_\bfR\mapsto [\mathbf{1}_\bfA]_{{}^\phi\bfA}$ (which admits a continuous right adjoint by \Cref{ex:compact to compact in trace}). 
Similarly, the two module structures coincide.
See \cite[Theorem 3.8.5]{gaitsgory2019toy} for more discussions.  
\end{remark}

Now suppose $\bfR=\Mod_\La$ and let $\bfA$ be equipped with $\phi$ as above. Let $\bfF_\bfA={}^\phi\bfA$ as in \Cref{ex: relating horizontal trace and vertical trace}, and let $\bfL=\bfA$ with $\beta: \bfA\to {}^\phi\bfA$ given by $X\in \bfA$ as above. Then $\ga: \bfA\to\bfA$ is given by $a\mapsto {}^\vee X\otimes \phi(a)\otimes X$. 

Suppose that $\bfA$ is compactly generated.  Let $(\bfA^\cpt)^\phi$ be the category consisting of $(Y\in \bfA^\cpt, \eta: Y\to \ga(Y))$. Note that a map $Y\to \ga(Y)$ is equivalent to a map
\begin{equation}\label{eq:input for S operator}
X\otimes Y\to \phi(Y)\otimes X.
\end{equation}
By abuse of notations, we will still denote it by $\eta$.

Recall the twisted Chern characters from \Cref{rem:twisted Chern character} gives a map
\begin{equation}\label{eq: abstract T-operator}
T: K_0((\bfA^\cpt)^\ga)\to H^0\mathrm{tr}(\bfA, \ga),\quad (Y,\eta)\mapsto T_{(Y,\eta)}.
\end{equation}
On the other hand, for every such $(Y,\eta)$, let $Y^\vee$ be the (right) dual of $Y$. Then we have an element in $H^0\End_{\tr(A,\phi)}([X]_\phi)$ defined as
\begin{equation}\label{eq: abstract S-operator}
S_{(Y,\eta)}: [X]_\phi\to [X\otimes Y\otimes Y^\vee]_\phi\to [\phi(Y)\otimes X\otimes Y^\vee]_\phi\cong [X\otimes Y^\vee\otimes Y]_\phi\to  [X]_\phi.
\end{equation}
This is usually called the $S$-operator associated to $(Y,\eta)$.

The following statement can be regarded as an abstract version of $S=T$ theorem \`a la V. Lafforgue. 

\begin{proposition}\label{prop: abstract S equals T}
Suppose $\bfA$ is a compactly generated rigid monoidal category and let $X\in \bfA^\cpt$. Then
under the isomorphism \eqref{eq: endo of X as a HH homology}, we have $S_{(Y,\eta)}=T_{(Y,\eta)}$.
\end{proposition}
\begin{proof}
To avoid confusions, we will write tensor product $a\otimes b$ of two objects in $\bfA$ as $ab$. 

We first make the isomorphism \eqref{eq: endo of X as a HH homology} explicit. Namely, we have the following commutative diagram
\[
\xymatrix{
\Mod_\La\ar^{u_\bfA}[rr]\ar@{=}[d]&& \bfA^\vee\otimes\bfA\ar_{\verd^{\mathrm{sr}}\otimes \id}[d]\ar^-{\id\otimes \ga}[r]\ar[d]& \bfA^\vee\otimes\bfA\ar^-{\mathrm{sw}}[r]\ar_{\verd^{\mathrm{sr}}\otimes\id }[d] & \bfA\otimes\bfA^\vee\ar^-{e_\bfA}[rrr]\ar^{\id\otimes \verd^{\mathrm{sr}}}[d] &&& \Mod_\La\ar@{=}[d] \\
\Mod_\La\ar^-{\mathbf{1}_\bfA}[r]\ar_{X}[dr]& \bfA\ar^-{m^R}[r]\ar^{X\otimes (-)}[d]&  \bfA\otimes\bfA \ar^-{\id\otimes \ga}[r] \ar^{X\otimes (-)\boxtimes \id}[d]& \bfA\otimes\bfA\ar^-{\mathrm{sw}}[r] & \bfA\otimes\bfA\ar^-m[r]& \bfA\ar^-{\Hom(\mathbf{1}_\bfA,-)}[rr] && \Mod_\La\\
&\bfA \ar^-{m^R}[r]\ar_-{[-]_{{}^\phi\bfA}}[drr] & \bfA\otimes \bfA \ar^{\id\otimes\phi}[r] & \bfA\otimes \bfA\ar^-{\mathrm{sw}}[r] &\bfA\otimes \bfA\ar^m[r]&\bfA\ar^-{{}^\vee X\otimes -}[u]{}^\phi\bfA\ar_-{\Hom(X,-)}[urr]&&\\
&&  & \tr(A,\phi) \ar_-{ [-]_{{}^\phi\bfA}^R}[urr] &&   &&
}\]
Here $\verd^{\mathrm{sr}}$ is the self-duality of $\bfA$ as in \Cref{ex: self-duality of cpt gen rigid monoidal cat}. 
Then $\End_{\tr(\bfA,\phi)}([X]_{{}^\phi\bfA})\in\Mod_\La$ is the image of $\La\in \Mod_\La$ under the functor obtained by composing functors along the bottom arrows, and $\mathrm{tr}(\bfA,\ga)$ is the image of $\La\in \Mod_\La$ under the functor obtained by composing functors along the top arrows. 

We recall the construction of $T_{Y,\eta}$ from \Cref{rem:twisted Chern character}. 
By \Cref{ex: self-duality of cpt gen rigid monoidal cat}, under the isomorphism
\[
\End(Y)\cong \Hom(\mathbf{1}_\bfA,  Y\otimes Y^\vee)\cong \Hom( Y^\vee\otimes Y, \mathbf{1}_\bfA)
\]
$\id_Y$ corresponds to the unit $u_Y: \mathbf{1}_\bfA\to Y\otimes Y^\vee$.
The map $\End(Y)\to e_{\bfA}(u_\bfA)$ corresponds to the counit map $e_Y: Y^\vee\otimes Y\to \mathbf{1}_\bfA$ corresponds to $Y^\vee\boxtimes Y\to u_\bfA$.  
\end{proof}

The following description of hom spaces between certain objects in $\tr(\bfA,\phi)$ is useful in practice. 
We refer to \cite[\textsection{3}]{Zhu2016} for more elementary accounts.

\begin{corollary}\label{cor-hom-space-of-compact-object-in-categorical-trace}
Suppose that $\bfR$ is compactly generated (e.g. $\bfR=\Mod_\La$).
 Assume that $\bfA$ is rigid and is compactly generated, with a set of compact generators $\{c_i\}$. Then for $X,Y\in \bfA$ with $X$ compact in $\bfA$,
 \[
 \Hom_{\tr(\bfA,\phi)}([X]_{{}^\phi \bfA},[Y]_{{}^\phi \bfA})\cong \colim_{\bfC\otimes \bfC_{/m^R(Y)}}\Hom_{\bfA}(X, c_j\otimes \phi(c_i)),
 \]
 where $\bfC\otimes \bfC\subset \bfA\otimes \bfA$ denotes the full subcategory spanned by $\{c_i\boxtimes c_j\}_{i,j}$.
\end{corollary}
Note that a morphism $c_i\boxtimes c_j\to m^R(Y)$ in $\bfA\otimes \bfA$ is equivalent to a morphism $c_i\otimes c_j\to Y$ in $\bfA$. So informally, this corollary says that every morphism $[X]_{{}^\phi \bfA}\to [Y]_{{}^\phi \bfA}$ in $\tr(\bfA,\phi)$ can be represented as a pair of morphisms $(X\to c_j\otimes \phi(c_i), c_i\otimes c_j\to Y)$ in $\bfA$ (compare with \cite[\S 3.1]{Zhu2016}).
\begin{proof}
  We have 
  \[
  \Hom_{\tr(\bfA,\phi)}([X]_{{}^\phi \bfA},[Y]_{{}^\phi \bfA})\cong\Hom_{{}^\phi \bfA}(X, m\circ \mathrm{sw}\circ m^R(Y)).
  \]
  As $\bfA$ is compactly generated, so is $\bfA\otimes \bfA$ with a set of compact generators given by $\{c_i\boxtimes c_j\}_{i,j}$. Then 
  \[
  m^R(Y)=\colim_{c_i\boxtimes c_j\to Y} c_i\boxtimes c_j.
  \] 
  As $X$ is compact, for every compact object $r\in \bfR$, $r\otimes X$ is still compact in $\bfA$. Therefore
  \[
  \Map_\bfR(r,\Hom_{{}^\phi \bfA}(X, m\circ \mathrm{sw}\circ m^R(Y)))=\Map_{\bfR}(r, \colim_{i,j} \Hom_{\bfA}(X, \phi(c_j)\otimes c_i)).
  \]
  As $\bfR$ is compactly generated, the above isomorphism implies the lemma.
\end{proof}

\begin{remark}\label{rem-approx-identity}
 In fact, the above corollary admits a more economic form. Namely, suppose we write $m^R(\mathbf{1}_A)\cong \colim_{i} (c_{i,1}\boxtimes c_{i,2})$ as a filtered colimit of (compact) objects in $\bfA\otimes \bfA$. Then as $m^R$ is a right $\bfA$-module homomorphism, we have $m^R(Y)\cong \colim_i ( (Y\otimes c_{i,1})\boxtimes c_{i,2})$. Therefore,
 \[
 \Hom_{\tr(\bfA,\phi)}([X]_{{}^\phi \bfA},[Y]_{{}^\phi \bfA})\cong \colim_i\Hom_\bfA(X, \phi(c_{i,2})\otimes Y\otimes c_{i,1}).
 \]
 Note that $m^R(\mathbf{1}_\bfA)\in \bfA\otimes \bfA$ is in fact isomorphic to the unit of the self-duality of $\bfA$. So in some cases the situation as in \Cref{cor-filtration-unit} is applicable.

Now, let $\bfM$ be an $\bfA$-module with a left dual $\bfN$, and $\al: \bfM\to {}^\phi \bfM$. Then under some similar assumption, one can compute $\Hom_{\tr(\bfA,\phi)}([X]_\phi, [\bfM,\al]_\phi)$. Suppose the image of $\mathbf{1}_\bfR\in \bfR\xrightarrow{u} \bfN\otimes_\bfA\bfM\to \bfN\otimes \bfM$ can be written as $\colim_i (n_i\boxtimes m_i)$, then
\[
\Hom_{\tr(\bfA,\phi)}([X]_\phi, [\bfM,\al]_\phi)\cong \colim_i\Hom_{\bfA}(X, e(\al(m_i)\boxtimes n_i)).
\]
\end{remark}

\subsubsection{Categorical traces of semi-rigid monoidal categories}

First, starting from the Hochschild complex $\HH(\bfA,\bfF)_\bullet$, we obtain a cosimplicial category $\HH(\bfA,\bfF)^R_\bullet$ by passing to the (not necessarily continuous) right adjoint. Then by
 \eqref{eq:limit-colimit equivalence}, we have
\begin{equation}\label{eq:right adjoint Hochschild complex}
\tr(\bfA,\bfF)=|\HH(\bfA,\bfF)_\bullet|\cong \tot(\HH(\bfA,\bfF)^R_\bullet).
\end{equation}

Under further assumption of $\bfA$, the Hochschild complex is monadic, and its categorical traces have nice formal properties.
\begin{lemma}\label{L:hochschild-homology-vs-cohomology}  
Assume that $m:\bfA\otimes \bfA\to \bfA$ admits an $\bfA\otimes \bfA^{\rev}$-linear right adjoint. 
\begin{enumerate}
\item\label{L:hochschild-homology-vs-cohomology-1} Let $\bfF$ be an $\bfA$-bimodule with $a_l: \bfA\otimes \bfF\to \bfF$ and $a_r: \bfF\otimes \bfA\to \bfF$ the left and the right action.
Then
$\tr(\bfA,\bfF)\cong  \lmodu_T(\bfF)$
with $T$ the monad given to $a_r\circ a_l^R$.

\item\label{cor:right-adjointability-of-module-to-Trace} Let $\eta: \bfF_1\to \bfF_2$ be a functor of $\bfA$-bimodules. Then the following diagram is right adjointable 
\[
\xymatrix{
\bfF_1 \ar^-{[-]_{\bfF_1}}[r]\ar_\eta[d] &  \tr(\bfA,\bfF_1) \ar^{\tr(\bfA,\eta)}[d] \\
\bfF_2 \ar^-{[-]_{\bfF_2}}[r]  & \tr(\bfA,\bfF_2).
}
\]
If in addition $\eta$ admits a right adjoint $\eta^R$ (in $\lincat_\bfR$), then the following diagram is right adjointable
\[
\xymatrix{
\bfF_1  \ar_-{[-]_{\bfF_1}}[d]\ar^\eta[rr] &&  \bfF_2 \ar^-{[-]_{\bfF_2}}[d] \\
\tr(\bfA,\bfF_1) \ar^{\tr(\bfA,\eta)}[rr]  && \tr(\bfA,\bfF_2).
}
\]
In fact, the same diagram as above, but with $\eta$ replaced by $\eta^R$, is canonically identified with the right adjoint of the above diagram.

\item\label{cor-fully-faithful-categorical-trace} Let $\bfF_1\to \bfF_2$ be a fully faithful functor of $\bfA$-bimodules. Then the induced functor
\[
\tr(\bfA,\eta): \tr(\bfA,\bfF_1)\to \tr(\bfA,\bfF_2)
\]
is fully faithful.
\end{enumerate}
\end{lemma}
We refer to \Cref{ex: vertical trace-simple-functoriality} for the notation $\tr(\bfA,\eta)$.
\begin{proof}
Using \eqref{eq:right adjoint Hochschild complex} and by \Cref{thm:Beck-Chevalley-descent}, it is enough to show that  for every coface map $\al:[n]\to[m]$, the diagram
\begin{equation}\label{Eq:hochschild-homology-vs-cohomology}
    \xymatrix{
    \HH(\bfA,\bfF)_{m+1}\ar[r] \ar_{d^{\HH}_0}[d]& \HH(\bfA,\bfF)_{n+1} \ar^{d^{\HH}_0}[d]\\
    \HH(\bfA,\bfF)_m\ar[r]  &          \HH(\bfA,\bfF)_n
    }
\end{equation}
is right-adjointable. We may assume that $\al=d_i$ is a coface map. If $i=1$, the desired right adjointability follows from \Cref{lem-rigid-lax-is-strict} \eqref{lem-rigid-lax-is-strict-1}; if  $i\neq 1$, the desired right adjointability follows from \Cref{lem-rigid-lax-is-strict} \eqref{lem-rigid-lax-is-strict-2}. Part \eqref{L:hochschild-homology-vs-cohomology-1} of the lemma follows.

Note that using  \Cref{lem-rigid-lax-is-strict} \eqref{lem-rigid-lax-is-strict-1}, an $\bfA$-bimodule functor $\bfF_1\to\bfF_2$ induces a functor of (semi)cosimplicial categories $\HH(\bfA,\bfF_1)_\bullet^R\to \HH(\bfA,\bfF_1)_\bullet^R$. For this, Part \eqref{cor:right-adjointability-of-module-to-Trace} follows directly from \Cref{prop:categorical-right-adjointability-colimits}.

For Part \eqref{cor-fully-faithful-categorical-trace}, we note that using \eqref{eq:right adjoint Hochschild complex}, it is enough to have level-wise fully faithfulness of the functors $\bfA^{\otimes n}\otimes \bfF_1 \to \bfA^{\otimes n}\otimes \bfF_2$ for all $n\geq 0$. Here we recall that $-\otimes-$ really means $-\otimes_\bfR-$. But as $\bfR$ itself is rigid (as a monoidal category in $\lincat$), applying the same reasoning again, the desired statement follows from \Cref{lem-tensor-preserve-fully-faithfulness}. 
\end{proof}

\begin{corollary}\label{cor: t-structure on categorical trace}
Let $\bfA$ be as in \Cref{L:hochschild-homology-vs-cohomology}, and $\bfF_\bfA$ is an $\bfA$-bimodule.
Suppose that both $\bfA$ and $\bfF_\bfA$ admit $t$-structure that are accessible (i.e. $\bfA^{\leq 0}$ is closed under filtered colimits) and that both action functors $\bfA\otimes\bfF_\bfA\to \bfF_\bfA$ and $\bfF_\bfA\otimes\bfA\to \bfF_\bfA$ are $t$-exact. Then $\tr(\bfA,\bfF_\bfA)$ admits a $t$-structure with $\tr(\bfA,\bfF_\bfA)^{\leq 0}$ generated (under extensions and filtered colimits) by the essential image of $\bfF_\bfA^{\leq 0}$ under the canonical functor $\bfF_\bfA\to \tr(\bfA,\bfF_\bfA)$. In addition, the functor $\bfF_\bfA\to \tr(\bfA,\bfF_\bfA)$ is $t$-exact.
\end{corollary}
\begin{proof}
That $\tr(\bfA,\bfF_\bfA)$  admits a prescribed $t$-structure follows from \cite[Proposition 1.4.4.11]{Lurie.higher.algebra}. To prove the last statement, let $X\in \bfF_\bfA^{\heartsuit}$. Then by definition $[X]_{\bfF_\bfA}\in  \tr(\bfA,\bfF_\bfA)^{\leq 0}$. On the other hand, using \Cref{L:hochschild-homology-vs-cohomology} \eqref{L:hochschild-homology-vs-cohomology-1}, we see that for every $Y\in \bfF_{\bfA}^{\leq -1}$, we have
\[
\Hom_{\tr(\bfA,\bfF_\bfA)}([Y]_{\bfF_\bfA},[X]_{\bfF_\bfA})=\Hom_{\bfF_\bfA}(Y, a_r(a_l^R(X))).
\]
As both $a_l$ and $a_r$ are $t$-exact, we see that $a_r(a_l^R(X))\in \bfF_\bfA^{\geq 0}$. It follows that $[X]_{\bfF_\bfA}\in  \tr(\bfA,\bfF_\bfA)^{\geq 0}$, as desired.
\end{proof}

Therefore, the categorical traces of monoidal categories satisfying the assumption as in  \Cref{L:hochschild-homology-vs-cohomology} have especially good formal properties. 

Here is another consequence.
\begin{lemma}\label{lem: smoothness when mult admits right adjoint}
Let $\bfA$ be as in \Cref{L:hochschild-homology-vs-cohomology}. Then $\bfA$ is a smooth algebra in $\lincat_\bfR$ in the sense of \Cref{def-2-dualizable-alg}.
\end{lemma}
When $\bfA$ is symmetric monoidal, this was proved in \cite[Proposition C.2.3]{AGKRRV.restricted.local.systems}. A slight modification of the argument is needed to deal with general monoidal categories.
\begin{proof}
Let $\bfM$ be a left $(\bfA\otimes\bfA^{\rev})$-module. We regard $\bfA\otimes\bfA^{\rev}$ as a left module over itself.  Then adjoint pair $m:\bfA\otimes \bfA^{\rev} \rightleftharpoons \bfA: m^R$ of $(\bfA\otimes \bfA^{\rev})$-linear functors induce an adjoint pair of $\bfR$-linear functors
\[
\bfM\cong \fun^{\mathrm{L}}_{\bfA\otimes\bfA^{\rev}}(\bfA\otimes\bfA^{\rev},\bfM)\rightleftharpoons \fun^{\mathrm{L}}_{\bfA\otimes\bfA^{\rev}}(\bfA,\bfM)
\]
The functor $\fun^{\mathrm{L}}_{\bfA\otimes\bfA^{\rev}}(\bfA,\bfM)\to \bfM$ is clearly conservative (if $F: \bfA\to \bfM$ is non-zero functor, then $\bfA\otimes \bfA\xrightarrow{m}\bfA\xrightarrow{F}\bfM$ is clearly non-zero). It follows from the Bar-Beck-Lurie theorem that
\[
\fun^{\mathrm{L}}_{\bfA\otimes\bfA^{\rev}}(\bfA,\bfM)\cong \lmodu_{T}\bfM
\]
for some monad $T\in \bfA\otimes\bfA^{\rev}$, which is given by multiplication by $\fraka=m^R(\mathbf{1}_\bfA)\in \bfA\otimes \bfA\cong \bfA\otimes\bfA^{\rev}$. (Note that as an object in $\bfA\otimes\bfA^{\rev}$, $\fraka$ has a natural algebra structure. )

Now let
\[
T_\bfA:=\lmodu_{\fraka}(\bfA\otimes\bfA^{\rev}).
\]
Then $T_\bfA$ has a natural right $(\bfA\otimes\bfA^{\rev})$-module structure. We claim that equipped with this right $(\bfA\otimes\bfA^{\rev})$-module structure, $T_\bfA$ is a left dual of $\bfA^u\in \lmodu_{\bfA\otimes\bfA^{\rev}}(\lincat_\bfR)$ (and therefore is the dual Serre module of $\bfA$ in the sense of  \Cref{def-2-dualizable-alg}). Indeed, for every left $\bfA\otimes\bfA^{\rev}$-module $\bfM$, we have
\[
\fun^{\mathrm{L}}_{\bfA\otimes\bfA^{\rev}}(\bfA,\bfM)=\lmodu_{\fraka}(\bfM)\cong T_\bfA\otimes_{\bfA\otimes\bfA^{\rev}}\bfM,
\]
where the last isomorphism follows form \cite[Theorem 4.8.6.4]{Lurie.higher.algebra}.
\end{proof}

\begin{remark}
When $\bfA$ is symmetric monoidal, $T_\bfA$ is canonically equivalent to $\bfA$ equipped with the right $\bfA\otimes\bfA^{\rev}$-module structure. This recovers \cite[Proposition C.2.3]{AGKRRV.restricted.local.systems}.
\end{remark}

Recall from  \cite[Proposition 4.6.4.4]{Lurie.higher.algebra} that an algebra $\bfA\in\alg(\lincat_\bfR)$ is a proper if and only if the forgetful
$\lmodu_\bfA\to \lincat_\bfR$ sends left dualizable $\bfA$-modules to dualizable $\bfR$-linear categories, and if and only if $\bfA$ is dualizable (in $\lincat_\bfR$). We make the following definition, which is a common generalization of \cite[Definition 3.1]{benzvi2009character} and \cite[\textsection{C.1.1}]{AGKRRV.restricted.local.systems}.

\begin{definition}\label{def: semi-rigid cat}
An algebra $\bfA\in\alg(\lincat_\bfR)$ is called a semi-rigid $\bfR$-linear category if it is proper and if $m:\bfA\otimes \bfA\to \bfA$ admits an $\bfA\otimes \bfA^{\rev}$-linear right adjoint. 
\end{definition}

The following statement is  a direct consequence of the definition and \Cref{lem: smoothness when mult admits right adjoint}.

\begin{proposition}\label{prop: 2-dualizability of semi-rigid monoidal category}
If $\bfA$ is an $\bfR$-linear semi-rigid monoidal category, then $\bfA$ is a $2$-dualizable algebra object in $\lincat_\bfR$ in the sense of \Cref{def-2-dualizable-alg}.
\end{proposition}

\begin{example}\label{ex: duality of rigid monoidal category}
If $\bfA$ is a rigid $\bfR$-linear category, then it is a semi-rigid $\bfR$-linear category. Indeed, in this case $\bfA$ is self-dual as an $\bfR$-linear category,  with unit given by $\bfR\xrightarrow{\mathbf{1}_\bfA} \bfA\xrightarrow{m^R}\bfA\otimes\bfA$ and counit $\bfA\otimes\bfA\xrightarrow{m}\bfA\xrightarrow{\Hom(\mathbf{1}_\bfA,-)}\bfR$. 
\end{example}

We have the following properties of semi-rigid monoidal categories, generalizing some statements from \cite[\textsection{C.2}, \textsection{C.3}]{AGKRRV.restricted.local.systems}. 

\begin{proposition}\label{rem-duality-datum-as-plain-cat}
Let $\bfA$ be a semi-rigid $\bfR$-linear monoidal category, and let $\bfM$ be a left $\bfA$-module. 
\begin{enumerate}
\item\label{rem-duality-datum-as-plain-cat-1}  $\bfA^\vee\cong \bfA$ with the unit of the duality datum given by
$\bfR\xrightarrow{\mathbf{1}_\bfA} \bfA\xrightarrow{m^R}\bfA\otimes\bfA$. In addition, $\bfA\cong \bfA^\vee\xrightarrow{(\mathbf{1}_\bfA)^\vee}\bfR$ is a Frobenius structure on $\bfA$. 
\item\label{rem-duality-datum-as-plain-cat-2} More generally, the category $\bfM$ is left dualizable as a left $\bfA$-module if and only if the underlying category is dualizable in $\lincat_\bfR$.  In addition, if 
\[
\bfR\to \bfN\otimes_\bfA\bfM,\quad \bfM\otimes \bfN\to \bfA
\]
is the duality datum for $\bfM$ as a left $\bfA$-module, then 
\[
\bfR\to \bfN\otimes_\bfA\bfM\xrightarrow{[-]_{\bfM\otimes \bfN}^R}\bfN\otimes \bfM,\quad \bfM\otimes\bfN\to \bfA\xrightarrow{(\mathbf{1}_{\bfA})^\vee}\bfR
\] 
is the duality datum for $\bfM$ as $\bfR$-linear category. 
\item\label{rem-duality-datum-as-plain-cat-3} If $\bfA$ is compactly generated, then the Frobenius structure $\bfA\to \bfR$ as in \eqref{rem-duality-datum-as-plain-cat-1}, when restricted to $\bfA^\cpt$ is given by $\Hom(\mathbf{1}_\bfA,-)$.
\end{enumerate}
\end{proposition}

The discussions in \Cref{def.rigid.from.GR} on the behaviors of rigidity under the change of the base (rigid) symmetric monoidal categories $\bfR'\to\bfR$ also apply to the semi-rigid case. In particular, when $\bfR$ is the category of spectra, then $\bfA$ as above is simply called semi-rigid monoidal category, and every semi-rigid $\bfR$-linear category is semi-rigid.

We make the following definition, generalizing the usual notion of pivotal structure on compactly generated rigid monoidal categories (see \Cref{ex: self-duality of cpt gen rigid monoidal cat}).
\begin{definition}\label{def: pivotal semi-rigid monoidal}
Let $\bfA$ be a semi-rigid monoidal category as above. Let $\sigma_\bfA$ be the Serre automorphism of $\bfA$ associated to the Frobenius structure of $\bfA$ as in \Cref{rem-duality-datum-as-plain-cat} \eqref{rem-duality-datum-as-plain-cat-1}. Then a pivotal structure of $\bfA$ is an isomorphism $\sigma_\bfA\cong \id_\bfA$ as algebra automorphisms of $\bfA$.
\end{definition}

Note that as explained in \cite[Remark 4.6.5.3, 4.6.5.4]{Lurie.higher.algebra} a pivotal structure on $\bfA$ induces isomorphisms 
\[
S_\bfA\cong \bfA\cong T_\bfA
\] 
as $\bfA$-bimodules, where $S_\bfA$ and $T_\bfA$ are the Serre bimodule and the dual Serre bimodule of $\bfA$ defined in \Cref{def-2-dualizable-alg}.

\subsubsection{Trace formula}\label{SS: trace formula}

We assume that $\bfA$ is $\bfR$-linear semi-rigid. By \Cref{prop: 2-dualizability of semi-rigid monoidal category}, it is $2$-dualizable as an algebra object in $\lincat_\bfR$. Let $\bfF_1$ and $\bfF_2$ be $\bfA$-bimodules, both of which admit left duals. Let
\[
\delta: \bfF_1\otimes_\bfA\bfF_2\to \bfF_2\otimes_\bfA\bfF_1
\]
be an isomorphism of $\bfA$-bimodules. 
Recall from \Cref{thm-second-trace-1} that there is a canonical isomorphism of objects in $\bfR$
\begin{equation}\label{eq: second-trace-1}
\mathrm{tr}(\tr(\bfA,\bfF_1), \tr(\bfF_2,\delta^{-1}))\cong \mathrm{tr}(\tr(\bfA,\bfF_2), \tr(\bfF_1,\delta)).
\end{equation}

\begin{theorem}\label{thm-second-trace-2}
Let $(\bfA,\bfF_1,\bfF_2,\delta)$ be as above.
Let $\bfM$ be a left smooth and proper $\bfA$-module. Suppose we are given the following commutative diagram
\begin{equation}\label{eq: pentagon-F1-F2-M}
\xymatrix{
&\ar_-{\beta_1}[ld]\bfM\ar^-{\beta_2}[rd]&\\
\bfF_1\otimes_\bfA\bfM\ar_-{\id\otimes\beta_2}[d]&& \bfF_2\otimes_\bfA\bfM\ar^-{\id\otimes\beta_1}[d]\\ 
\bfF_1\otimes_\bfA\bfF_2\otimes_\bfA\bfM\ar^-{\delta\otimes \id}[rr]&& \bfF_2\otimes_\bfA\bfF_1\otimes_\bfA\bfM,}
\end{equation}
where $\beta_i: \bfM\to \bfF_i\otimes_\bfA\bfM$ are two $\bfA$-linear functors that admit continuous right adjoint. 
\begin{enumerate}
\item\label{thm-second-trace-2-1} The object $[\bfM,\beta_1]_{\bfF_1}\in \tr(\bfA,\bfF_1)$ is compact, and there is a canonical homomorphism 
\[
\eta_1: [\bfM,\beta_1]_{\bfF_1}\to \tr(\bfF_2,\delta^{-1})([\bfM,\beta_1]_{\bfF_1}),
\] 
Similarly, $[\bfM,\beta_2]_{\bfF_2}\in \tr(\bfA,\bfF_2)$ is compact, and there is a canonical homomorphism 
\[
\eta_2: [\bfM,\beta_2]_{\bfF_2}\to \tr(\bfF_1,\delta)([\bfM,\beta_2]_{\bfF_2}).
\] 
 \item\label{thm-second-trace-2-2} Under the isomorphism \eqref{eq: second-trace-1}, there is a canonical isomorphism
\[
\mathrm{ch}([\bfM,\beta_1]_{\bfF_1},\eta_1)=\mathrm{ch}([\bfM,\beta_2]_{\bfF_2},\eta_2),
\]
where the (twisted) Chern character is defined as in \eqref{eq-twisted-Chern-character-1}.
\end{enumerate}
\end{theorem}

\begin{example}
A basic example is when $\bfA=\bfF_1=\bfF_2=\bfR$ with $\al$ being the identity equivalence. Let $\bfM$ be a proper and smooth $\bfR$-linear category equipped with two commuting $\bfR$-linear endomorphisms $\beta_1$ and $\beta_2$ both of which admitting $\bfR$-linear right adjoint. Then via \Cref{ex: class-vs-usual-hochschild}, we recover \Cref{thm: classical-trace-formula}.
\end{example}

\begin{example}
Another special case that is important in representation is as follows. Assume that $\bfA$ is rigid.
 Let $\bfF_1={}^\phi\bfA$ and $\bfF_2=\bfA$, with the equivalence $\al$ being the canonical one. We let $\bfM=\bfA$ regarded as a left $\bfA$-module. We let $\beta_1: \bfM\to \bfF_1\otimes_\bfA\bfM$ be given by a $\phi$-equivariant compact object $Y\in\bfA$ as in \Cref{ex:compact to compact in trace}, and let $\beta_2:  \bfM\to \bfF_2\otimes_\bfA\bfM$ be given by the unit of $\bfA$. 
Then under the isomorphism
\[
\mathrm{tr}(\tr(\bfA,\phi),\id_{\tr(\bfA,\phi)})\cong \mathrm{tr}(\tr(\bfA), \phi)
\]
we have 
\[
\mathrm{ch}([Y]_{{}^\phi\bfA})=\mathrm{ch}([\mathbf{1}_{\bfA}]_{\bfA}, S_Y),
\]
where $S_Y$ is the endomorphism of $[\mathbf{1}_{\bfA}]_{\bfA}$ as constructed in \eqref{eq: abstract S-operator}. 
\end{example}

\begin{proof}[Proof of \Cref{thm-second-trace-2}]
That $[\bfM,\beta_i]_{\bfF_i}\in \tr(\bfA,\bfF_i)$ is compact follows from \Cref{lem-adjoint-of-induced-trace-map}. By \Cref{lemma:composition-bimodule-vs-map-trace},
\[
\tr(\bfF_1,\delta)([\bfM,\beta_2]_{\bfF_2})=[\bfF_1\otimes_\bfA\bfM,  (\delta\otimes\id_{\bfM})\circ(\id_{\bfF_1}\otimes\beta_2)]_{\bfF_2}.
\]
We regard $\beta_1$ is a functor of left $\bfA$-module. The commutative diagram \eqref{eq: pentagon-F1-F2-M} allows us to
apply \Cref{lem-functoriality-duality-data-generalization} and obtain a map
\[
\eta_2: [\bfM,\beta_2]_{\bfF_2}\to [\bfF_1\otimes_\bfA\bfM,  (\delta\otimes\id_{\bfM})\circ(\id_{\bfF_1}\otimes\beta_2)]_{\bfF_2}.
\]
Similarly we have $\eta_1$. This proves Part \eqref{thm-second-trace-2-1}.

To prove Part \eqref{thm-second-trace-2-2}, we first notice that if $\bfF_1=\bfF_2=\bfF$ with $\delta=\id$ and $\beta_1=\beta_2$, then the statement is clear.
In particular, we may apply this observation to the case $\bfF=\bfM\otimes\bfN$, which as an $\bfA$-bimodule admits a left dual given by $\bfM\otimes \bfN$ itself (see  \Cref{rem-3-category-structure} \eqref{rem-3-category-structure-3}).

We also notice that giving a commutative diagram \eqref{eq: pentagon-F1-F2-M} is equivalent to giving a commutative diagram
\[
\xymatrix{
\ar_-{\beta_1^\sharp\otimes \beta_2^\sharp}[d](\bfM\otimes \bfN)\otimes_\bfA (\bfM\otimes \bfN)\ar@{=}[r]&(\bfM\otimes \bfN)\otimes_\bfA (\bfM\otimes \bfN)\ar^-{\beta_2^\sharp\otimes\beta_1^\sharp} [d] \\
\bfF_1\otimes_\bfA \bfF_2\ar^\delta_-\cong[r]& \bfF_2\otimes_\bfA \bfF_1,
}\]
where
\[
\beta_i^\sharp: \bfM\otimes\bfN\xrightarrow{\beta_i\otimes\id_\bfN}\bfF_i\otimes_\bfA\bfM\otimes\bfN\xrightarrow{\id_{\bfF_i}\otimes e_\bfM}\bfF_i,
\]
is an $\bfR$-bilinear functor with an $\bfR$-bilinear right adjoint.
Now the result is a consequence of \Cref{prop: functoriality of second-trace} below.
\end{proof}

To state \Cref{prop: functoriality of second-trace}, let $\bfA$ be as above. Suppose we have left dualizable $\bfA$-bimodules $\bfF'_1,\bfF'_2,\bfF_1,\bfF_2$ together with the following commutative diagram
\[
\xymatrix{
\bfF'_1\otimes_\bfA\bfF'_2\ar^\delta[r]\ar_{\ga_1\otimes\ga_2}[d] & \bfF'_2\otimes_\bfA\bfF'_1\ar^{\ga_2\otimes\ga_1}[d]\\
\bfF_1\otimes_\bfA\bfF_2\ar^{\delta'}[r] & \bfF_2\otimes_\bfA\bfF_1,
}
\]
where $\ga_i: \bfF'_i\to\bfF_i$ is an $\bfA$-bilinear functor that admits an $\bfA$-bilinear right adjoint. 
Then by \Cref{lemma:composition-bimodule-vs-map-trace} we have the following commutative digram
\[
\xymatrix{
\tr(\bfA,\bfF'_1)\ar^-{\tr(\bfF'_2,{\delta'}^{-1})}[rr] \ar_{\tr(\bfA,\ga_1)}[d]&& \tr(\bfA,\bfF'_1) \ar^-{\tr(\bfA,\ga_1)}[d]\\
\tr(\bfA,\bfF_1)\ar^-{\tr(\bfF_2,{\delta}^{-1})}[rr]         &&  \tr(\bfA,\bfF_1),
}
\]
with $\tr(\bfA,\ga_1)$ admitting $\bfR$-linear right adjoint. Let 
\[
\eta_1:  \tr(\bfA,\ga_1)\circ \tr(\bfF'_2,{\delta'}^{-1})\cong \tr(\bfF_2,{\delta'}^{-1})\circ  \tr(\bfA,\ga_1)
\] 
be the isomorphism witnessing the above commutative diagram.
Then by \Cref{lem-functoriality-duality-data}, we have a morphism
\[
\mathrm{tr}(\tr(\bfA,\ga_1),\eta_1): \mathrm{tr}(\tr(\bfA,\bfF'_1), \tr(\bfF'_2, {\delta'}^{-1}))\to \mathrm{tr}(\tr(\bfA,\bfF_1), \tr(\bfF_2, {\delta}^{-1})).
\]
Similarly, we have
\[
\mathrm{tr}(\tr(\bfA,\ga_2),\eta_2): \mathrm{tr}(\tr(\bfA,\bfF'_2), \tr(\bfF'_1, \delta'))\to \mathrm{tr}(\tr(\bfA,\bfF_2), \tr(\bfF_1, \delta)).
\]
\begin{theorem}\label{prop: functoriality of second-trace}
Under the equivalence \eqref{eq: second-trace-1} for $(\bfF'_1,\bfF'_2,\delta')$ and for $(\bfF_1,\bfF_2,\delta)$, the map
$\mathrm{tr}(\tr(\bfA,\ga_1),\eta_1)$ and $\mathrm{tr}(\tr(\bfA,\ga_2),\eta_2)$ are canonically identified.
\end{theorem}

\begin{proof}
We recall the construction of the equivalence \eqref{eq: second-trace-1} following  \cite{Campbell.Ponto}. 
We make use notations as in \Cref{not: simplified notation}. 

We identify $\bfR\cong \fun_\bfR(\bfR,\bfR)$ as before. Then the object
\[
\mathrm{tr}(\tr(\bfA,\bfF_1), \tr(\bfF_2, {\delta}^{-1}))\cong \mathrm{tr}(\tr(\bfA,\bfG_1), \tr(\bfF_2, {\delta}^{-1})^\vee) \in \bfR
\]
is identified with the composition of functors along the left half edges of the above big diagram, while $ \mathrm{tr}(\tr(\bfA,\bfF_2), \tr(\bfF_1, \delta))$ is identified with the composition of functors along the right half edges of the diagram.

\begin{tiny}
\begin{center}
\begin{equation}\label{eq: big diagram for second trace}
\xymatrix{
&&\ar_-{u_{\langle\bfF_1\rangle}}[ddll]\bfR\ar^\mu[d]\ar^-{u_{\langle u_{\bfF_2}\rangle}}[ddrr]&&\\
&&\ar[dl]\langle\bfT_\bfA\rangle\ar[dr]\ar@{}[dd] | {(\mathrm{IV})}&&\\
\langle\bfF_1\rangle\otimes \langle\bfG_1\rangle\ar_{\id\otimes\langle u_{\bfF_2}\odot\id\rangle}[d]\ar@{}[dr] | {(\mathrm{III})}&\ar@{}[u] | {(\mathrm{I})}\ar_-{\eqref{eq:map nu TA}}^-\mu[l]\langle\bfF_1\odot\bfT_\bfA\odot\bfG_1\rangle\ar[d]&&\ar@{}[u] | {(\mathrm{II})}\ar[d]\langle\bfF_2\odot\bfT_\bfA\odot\bfG_2\rangle\ar[r]&\langle\bfF_2\rangle\otimes\langle\bfG_2\rangle\ar[d]\ar@{}[dl] | {(\mathrm{V})}\\
\langle \bfF_1\rangle\otimes \langle\bfG_2\odot\bfF_2\odot \bfG_1\rangle\ar_{\id\otimes c}^{\cong}[d]\ar@{}[drrrr] | {(\mathrm{VI})}&\ar[l] \langle\bfF_1\odot\bfT_\bfA\odot\bfG_2\odot\bfF_2\odot\bfG_1\rangle\ar^-\cong[rr]&&\ar[ll]\langle\bfF_2\odot \bfG_1\odot\bfF_1\odot\bfT_\bfA\odot\bfG_2\rangle\ar[r]& \langle\bfF_2\odot \bfG_1\odot\bfF_1\rangle\otimes\langle\bfG_2\rangle\ar[d]\\
\langle \bfF_1\rangle\otimes \langle\bfF_2\odot \bfG_1\odot\bfG_2\rangle\ar[r]\ar_{\id\otimes\langle\id\odot\delta^\vee\rangle}^{\cong}[d]\ar@{}[dr] | {(\mathrm{VII})}&\langle\bfF_2\odot \bfG_1\odot\bfG_2\odot\bfS_\bfA\odot\bfF_1\rangle\ar^{\langle\id\odot\delta^\vee\odot\id\odot\id\rangle}_\cong[d]\ar^-\cong[rr]&&\ar[ll]\langle\bfG_2\odot\bfS_\bfA\odot \bfF_1\odot\bfF_2\odot \bfG_1\rangle\ar[d]\ar@{}[dr] | {(\mathrm{VIII})} &\ar[l]\langle\bfF_1\odot\bfF_2\odot \bfG_1\rangle\otimes\langle\bfG_2\rangle\ar[d]\\
\langle \bfF_1\rangle\otimes \langle\bfF_2\odot \bfG_2\odot\bfG_1\rangle\ar[r]\ar_{\id\otimes\langle e_{\bfF_2}\odot\id\rangle}[d]\ar@{}[dr] | {(\mathrm{IX})}&\langle\bfF_2\odot \bfG_2\odot\bfG_1\odot\bfS_\bfA\odot\bfF_1\rangle\ar[d]&&\langle\bfG_2\odot\bfS_\bfA\odot \bfF_2\odot\bfF_1\odot \bfG_1\rangle\ar[d]\ar@{}[dr] | {(\mathrm{XI})}&\ar[l]\langle\bfF_2\odot\bfF_1\odot \bfG_1\rangle\otimes\langle\bfG_2\rangle\ar[d]\\
\langle \bfF_1\rangle\otimes \langle\bfG_1\rangle\ar^-{\eqref{eq: map epsilon S_A}}_-\epsilon[r]\ar_-{e_{\langle\bfF_1\rangle}}[ddrr]&\langle\bfG_1\odot\bfS_\bfA\odot\bfF_1\rangle\ar[dr]\ar@{}[d] | {(\mathrm{XII})}&&\ar[dl]\langle\bfG_2\odot\bfS_\bfA\odot \bfF_2\rangle \ar@{}[d] | {(\mathrm{XIII})} &\ar[l]\langle\bfF_2\rangle\otimes\langle\bfG_2\rangle\ar^-{e_{\langle\bfF_2\rangle}}[ddll]\\
&&\langle\bfS_\bfA\rangle\ar^\epsilon[d]\ar@{}[uuu] | {(\mathrm{X})}&&\\
&&\bfR&&
}
\end{equation}
\end{center}
\end{tiny}

We need to explain why this diagram is commutative. Namely,
the commutativity of $(\mathrm{I}), (\mathrm{II}), (\mathrm{XII}), (\mathrm{XIII})$ follows from  \Cref{lem-dualizability-of-categorical-trace}, and the commutativity of $(\mathrm{VI})$ follows from \Cref{lem: crucial comm diagram secondary trace}. The commutativity of  $(\mathrm{III})$,  $(\mathrm{V})$, $(\mathrm{VII})$, $(\mathrm{VIII})$, $(\mathrm{IX})$ and $(\mathrm{XI})$ follows from the fact that \eqref{eq:map nu TA} and \eqref{eq: map epsilon S_A} are functorial in $F_1$ and $F_2$. The commutativity of $(\mathrm{IV})$ comes of the canonical isomorphism of functors
\[
(\id_{\bfG_1\odot\bfF_1\odot\bfT_{\bfA}}\odot u_{\bfF_2})\circ (u_{\bfF_1}\odot \id_{\bfT_{\bfA}})\cong (u_{\bfF_2}\odot \id_{\bfT_{\bfA}\odot \bfG_2\odot\bfF_2})\circ ( \id_{\bfT_{\bfA}}\odot u_{\bfF_2}),
\] 
and the functoriality of cyclic invariance of trace as in \Cref{lem: cyclic invariance of trace}. The commutativity of $(\mathrm{X})$ follows by similar reasoning.

There is also  a corresponding big commutative diagram, which witnesses the equivalence \eqref{eq: second-trace-1} for $(\bfF'_1,\bfF'_2,\delta')$.
Then the morphism $\mathrm{tr}(\tr(\bfA,\ga_1),\eta_1)$ is the composition of $2$-morphisms in the following diagram
\begin{tiny}
\[
\xymatrix{
\bfR\ar[r]\ar@{=}[d]\drtwocell\omit{\quad\quad\; \alpha_{\langle\ga_1\rangle}}&\langle\bfF'_1\rangle\otimes \langle\bfG'_1\rangle\ar[r]\ar[d]\drtwocell\omit{\quad\quad\; \langle\alpha_{\ga_1}\rangle}&\langle \bfF'_1\rangle\otimes \langle\bfG'_2\odot\bfF'_2\odot \bfG'_1\rangle\ar[r]\ar[d]&\langle \bfF'_1\rangle\otimes \langle\bfF'_2\odot \bfG'_2\odot\bfG'_1\rangle\ar[r]\ar[d]\drtwocell\omit{\quad\quad\; \langle\beta_{\ga_1}\rangle}&\langle \bfF'_1\rangle\otimes \langle\bfG'_1\rangle\ar[r]\ar[d]\drtwocell\omit{\quad\quad\; \beta_{\langle\ga_1\rangle}}&\bfR\ar@{=}[d]\\
\bfR\ar[r]&\langle\bfF_1\rangle\otimes \langle\bfG_1\rangle\ar[r]&\langle \bfF_1\rangle\otimes \langle\bfG_2\odot\bfF_2\odot \bfG_1\rangle\ar[r]&\langle \bfF_1\rangle\otimes \langle\bfF_2\odot \bfG_2\odot\bfG_1\rangle\ar[r]&\langle \bfF_1\rangle\otimes \langle\bfG_1\rangle\ar[r]&\bfR,
}\]
\end{tiny}

\noindent where the vertical morphisms are induced by $\ga_i$s and their conjugate functors $\ga_i^o$ (as defined in \Cref{lem-functoriality-duality-data-generalization}), and where the
$2$-morphisms are induced by $2$-morphisms from \Cref{lem-functoriality-duality-data}.

Now each small commutative diagram in \eqref{eq: big diagram for second trace} for $(\bfF'_1,\bfF'_2,\delta')$ maps to the corresponding small commutative diagram for $(\bfF_1,\bfF_2,\delta)$, and we need to show that the resulting diagram is $2$-commutative for our specified $2$-morphisms. 

We start with the observation that for $(\mathrm{VI})$, $(\mathrm{VII})$ and $(\mathrm{VIII})$, the resulting diagrams are strictly commutative. 

Next, we deal with $(\mathrm{IV})$. Using the functoriality of cyclic invariance of vertical trace (see \Cref{lem: cyclic invariance of trace}), it is enough show that the following diagram is $2$-commutative.

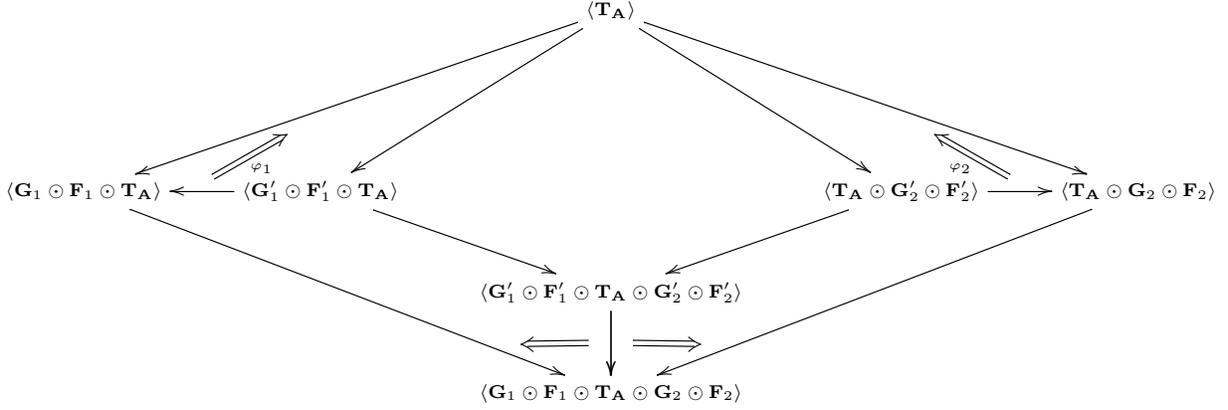
\begin{tiny}
\begin{equation*}
\xymatrix{
&&\ar[ddl]\ar[ddll]^(.6){\;}="b" \langle \bfT_\bfA \rangle\ar[ddr]\ar[ddrr]_(.6){\;}="d" &&\\
&&&&\\
\langle \bfG_1\odot \bfF_1\odot \bfT_\bfA\rangle\ar[ddrr]^(0.8){\;}="i"&\ar[l]_{\;}="a"\langle \bfG'_1\odot\bfF'_1\odot \bfT_\bfA\rangle\ar[dr]&&\langle\bfT_\bfA\odot \bfG'_2\odot\bfF'_2\rangle\ar[dl]\ar[r]^{\;}="c"&\langle \bfT_\bfA\odot \bfG_2\odot \bfF_2\rangle\ar[ddll]_(0.8){\;}="k"\\
&&\ar[d]_{\;}="h"\langle\bfG'_1\odot\bfF'_1\odot\bfT_\bfA\odot\bfG'_2\odot\bfF'_2\rangle\ar[d]^{\;}="j"&&\\
&&\langle\bfG_1\odot\bfF_1\odot\bfT_\bfA\odot\bfG_2\odot\bfF_2\rangle&&
 \ar@{=>}"a";"b"_{\varphi_1}  \ar@{=>}"c";"d"^{\varphi_2} \ar@{=>}"h";"i" \ar@{=>}"j";"k"
 }
\end{equation*}
\end{tiny}

But this is clear. Indeed, both compositions of $2$-morphisms in the left half and in the right half can be identified with the $2$-morphism obtained by taking adjoint of the following isomorphism
\begin{small}
\[
 \langle (\id_{\bfG'_1}\odot \ga_1\odot\id_{\bfT_\bfA}\odot\id_{\bfG'_2}\odot \ga_2)\circ (u_{\bfF'_1}\odot\bfT_\bfA\odot u_{\bfF'_2})\rangle \cong \langle(\ga_1^\vee\odot \id_{\bfF_1}\odot\id_{\bfT_\bfA}\odot\ga_2^\vee\odot \id_{\bfF_2})\circ (u_{\bfF_1}\odot\bfT_\bfA\odot u_{\bfF_2})\rangle.
\]
\end{small}

The proof for $(\mathrm{X})$ is similar.

To deal with the remaining commutative diagrams, the crucial lemma we need is as follows, which follows from \Cref{L:hochschild-homology-vs-cohomology} \eqref{cor:right-adjointability-of-module-to-Trace}.
\begin{lemma}\label{lem: right adjointability of mu and epsilon}
Let $\eta: \bfH_1\to\bfH_2$ be a functor of $\bfA$-bimodules, and let $\bfK$ be an $\bfA$-bimodule. Suppose $\eta^R$ exists (as a $\bfA$-bilinear functor), then the following commutative diagrams (induced by the functoriality of \eqref{eq:map nu TA} and \eqref{lem: trace dual morphism vs dual morphism trace})
\[
\xymatrix{
\langle\bfH_1\odot \bfT_\bfA\odot \bfK\rangle \ar^{\langle\eta\odot \id_{\bfT_\bfA}\odot\id_{\bfK}\rangle}[rr]\ar[d] &&\langle\bfH_2\odot \bfT_\bfA\odot \bfK\rangle \ar[d], & \langle\bfH_1\rangle\otimes\langle\bfK\rangle\ar^{\langle \eta\rangle\otimes \id_{\langle\bfK\rangle}}[rr]\ar[d]&& \langle\bfH_2\rangle\otimes\langle\bfK\rangle\ar[d]\\
\langle\bfH_1\rangle\otimes\langle\bfK\rangle\ar^{\langle \eta\rangle\otimes \id_{\langle\bfK\rangle}}[rr]&& \langle\bfH_2\rangle\otimes\langle\bfK\rangle & \langle\bfH_1\odot \bfS_\bfA\odot \bfK\rangle \ar^{\langle\eta\odot \id_{\bfS_\bfA}\odot\id_{\bfK}\rangle}[rr] && \langle\bfH_2\odot \bfS_\bfA\odot \bfK\rangle
}\]
are right adjointable. In addition, the right adjoint diagrams are induced by $\eta^R$.
\end{lemma}

To see how to apply this lemma, we consider the map of the commutative diagram $(\mathrm{I})$  (as in \eqref{eq: big diagram for second trace}) for $\bfF_i'$ to the corresponding diagram for $\bfF_i$, and show that the resulting diagram is $2$-commutative. That is, we 
claim that the following diagram is $2$-commutative.
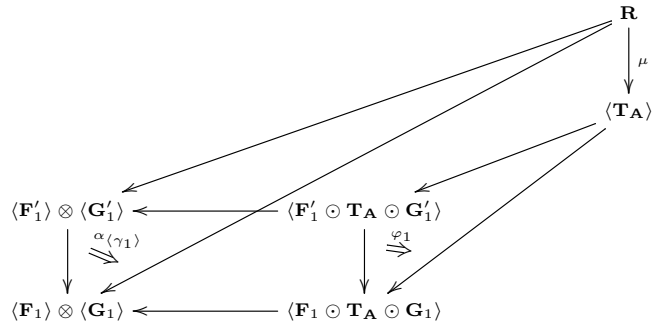
\begin{tiny}
\begin{equation}\label{eq: 2-commutative I}
\xymatrix{
&&&&\ar[ddllll]\bfR\ar^\mu[d]\ar[dddllll]_(.88){\;}="b"\\
&&&&\ar[dll]\langle\bfT_\bfA\rangle\ar[ddll]_(.75){\;}="d"\\
\langle\bfF'_1\rangle\otimes \langle\bfG'_1\rangle\ar[d]^(0.35){\;}="a"&&\ar[ll]\langle\bfF'_1\odot\bfT_\bfA\odot\bfG'_1\rangle\ar[d]^(0.35){\;}="c"&\\
\langle\bfF_1\rangle\otimes \langle\bfG_1\rangle&&\ar[ll]\langle\bfF_1\odot\bfT_\bfA\odot\bfG_1\rangle& 
\ar@{=>}"a";"b"^{\al_{\langle\ga_1\rangle}}  \ar@{=>}"c";"d"^{\varphi_1}}
\end{equation}
\end{tiny}

Unraveling the definition of the $2$-morphisms (as explained in \Cref{lem-functoriality-duality-data}), we see that \eqref{eq: 2-commutative I} can be expanded as the following diagram
\begin{tiny}
\[
\xymatrix{
\bfR\ar[r]&\langle \bfT_\bfA\rangle\ar^-{\langle\id_{\bfT_\bfA}\odot u_{\bfF'_1}\rangle}[r]\ar_-{\langle\id_{\bfT_\bfA}\odot u_{\bfF_1}\rangle}[d]\ar@{}[dr] | {(*)}&\langle \bfT_\bfA\odot\bfG'_1\odot\bfF'_1\rangle\ar[d]\ar[dr]&\\
&\langle \bfT_\bfA\odot\bfG_1\odot\bfF_1\rangle\ar[r]\ar[dr]\ar@{=}[d]^-{\;}="b"&\langle \bfT_\bfA\odot\bfG'_1\odot\bfF_1\rangle\ar[dr]\ar[dl]_(.7){\;}="a" \ar@{}[r] | {(**)}&\langle\bfF'_1\rangle\otimes\langle\bfG'_1\rangle\ar[d]\\
&\langle \bfT_\bfA\odot\bfG_1\odot\bfF_1\rangle\ar[dr]&\langle\bfF_1\rangle\otimes\langle\bfG_1\rangle\ar[r]\ar@{=}[d]^-{\;}="d"&\langle\bfF_1\rangle\otimes\langle\bfG'_1\rangle\ar[dl]_(.7){\;}="c"\\
&&\langle\bfF_1\rangle\otimes\langle\bfG_1\rangle.&
\ar@{=>}"a";"b"  \ar@{=>}"c";"d" 
}\]
\end{tiny}

We explain the unlabelled arrows.
\begin{itemize}
\item All arrows pointing to southeast are given by \eqref{eq:map nu TA}.
\item The two vertical arrows in $(**)$ are induced by $\ga_1: \bfF'_1\to \bfF_1$.
\item All arrows pointing to the southwest are induced by $\ga_1^o: \bfG'_1\to \bfG_1$.
\item Right arrows in the second and the third arrow are induced $\ga_1^\vee: \bfG_1\to \bfG'_1$.
\end{itemize}

Next we explain why this diagram is $2$-commutative. 
\begin{itemize}
\item The commutativity of $(*)$ is due to the canonical isomorphism $(\id_{\bfG'_1}\odot \ga_1)(u_{\bfF'_1})\cong   (\ga_1^\vee\odot \id_{\bfF_1})(u_{\bfF_1})$.
\item The commutativity of $(**)$ is a consequence of the functoriality of $\eqref{eq:map nu TA}$.
\item Since $\ga_1^\vee=(\ga_1^o)^R$ and $\langle\ga_1^\vee\rangle\cong \langle\ga_1\rangle^\vee$ (see \Cref{lem: trace dual morphism vs dual morphism trace}), we see that the right arrows in the second and the third arrow are the right are the right adjoints of the arrows pointing to the southwest. Therefore by  \Cref{lem: 2-commutative digram for right adjointable} and by \Cref{lem: right adjointability of mu and epsilon}, the part below $(*)$ and $(**)$ are $2$-commutative.
\end{itemize}

This shows the $2$-commutativity of \eqref{eq: 2-commutative I}.
The same arguments deal with diagrams involving $(\mathrm{II})$, $(\mathrm{III})$, $(\mathrm{V})$, $(\mathrm{IX})$, $(\mathrm{XI})$, $(\mathrm{XII})$ and $(\mathrm{XIII})$ as in \eqref{eq: big diagram for second trace}. The theorem is proved.
\end{proof}

\newpage

\section{Sheaf theory and traces of convolution categories}\label{sec:trace}

The first goal of this section is to review the abstract formalism of sheaf theory, commonly known as the six-functor formalism, following the works of \cite{liu2012enhanceda}, \cite{liu2012enhanced}, and \cite{Gaitsgory.Rozenblyum.DAG.vol.I} (see also \cite{Lucas.Mann} and \cite{scholze.six.functor}). We aim to formulate the theory in a manner suitable for our applications, and we will extend some results from these sources slightly to construct the sheaf theory we intend to use.
Notably, we have managed to avoid employing a $2$-categorical formalism, as required in \cite{Gaitsgory.Rozenblyum.DAG.vol.I}. As discussed in \Cref{rem-sheaf-theory-for-adjoint-factorization}, $2$-categorical structures in sheaf theory often are not additional structures but are properties inherent to the theory.
\footnote{We initially developed the abstract sheaf theory independently to digest the difficulties associated with constructing sheaf theory via Kan extensions, as outlined in \cite{liu2012enhanced} and \cite{Gaitsgory.Rozenblyum.DAG.vol.I}. Subsequently, \cite{Lucas.Mann} and \cite{scholze.six.functor} were published, streamlining the theory considerably. In particular, \cite{scholze.six.functor} also highlighted that $2$-categorical structures are not essential to the formalism.}

The second goal of this section is to develop a method for calculating the (twisted) categorical trace of monoidal categories arising from convolution patterns in algebraic geometry, building upon ideas from \cite{benzvi2009character} and \cite{ben2017spectral}. In contrast to \emph{loc. cit.}, which typically operates within concrete sheaf theories (primarily the theory of coherent sheaves or the theory of $D$-modules), we will develop the formalism in the context of an abstract sheaf theory. Our aim is to apply this formalism to both the theory of coherent sheaves (which will be developed in \Cref{S: theory of coherent sheaves}) and the theory of $\ell$--adic sheaves (to be explored in \Cref{sec:pspl-stacks}).

Consequently, we will bypass the general integral transform formalism outlined in \cite{benzvi2009character} and \cite{ben2017spectral}. Instead of calculating the categorical trace of a monoidal category directly, we will derive a geometric version of it. In favorable cases (including those considered in this paper), this geometric version coincides with the actual categorical trace. However, this may not hold for future applications, and the geometric version often appears to be the more relevant one. We will also employ similar ideas to compute the categorical trace of a module category arising from a monoidal category, which again originates from a convolution pattern. Additionally, we will investigate the functoriality between categorical traces arising from convolution patterns, which seems to be a novel contribution.

\subsection{Formalism of correspondences}\label{sec:stacks-correspondences-general}
First we review the formalism of correspondences, as first appeared in \cite[\textsection{6.1}]{liu2012enhanced} and \cite[Chapter 7]{Gaitsgory.Rozenblyum.DAG.vol.I}. There are mainly two (closely related) usages of this formalism in the paper. 
First, it provides a convenient framework to discuss convolution pattern arising from algebraic geometry and representation theory, and is useful for our study of (geometric) trace. Second, it encodes various sheaf theories in algebraic geometry in a concise way, as first observed by Lurie. In particular, in \Cref{S: theory of coherent sheaves} and \Cref{sec:pspl-stacks}, we will discuss the theory of coherent sheaves and the theory of $\ell$-adic sheaves using the formalism of correspondences.

\subsubsection{Category of correspondences}
Let $\bfC$ be an $\infty$-category that admits finite limits and finite coproducts. Let $\pt$ denote the final object in $\bfC$. The category $\bfC$ will play the role of the category of geometric objects. 

\begin{definition}\label{def-closure-property-of-morphism}
A class $\mathrm{E}\subset \mathrm{Mor}(\bfC)$ of morphisms in $\bfC$ is called weakly stable if it contains all isomorphisms in $\bfC$, is stable under (homotopy) equivalences of morphisms, and is stable under base change and compositions. It is called strongly stable if it is weakly stable and satisfies the following `$2$ out of $3$' property: for composable morphisms $\al_1,\al_2$ in $\bfC$ with $\al_2\in \mathrm{E}$,  $\al_1$ belongs to $\mathrm{E}$ if and only if $\al_2\circ \al_1$ belongs to $\mathrm{E}$. For a weakly stable class $\mathrm{E}$, we denote by $\bfC_{\mathrm{E}}$ the subcategory of $\bfC$ consisting of morphisms in $\mathrm{E}$. 
\end{definition}

\begin{remark}\label{rem-prop-of-weak-strong-stable-class}
\begin{enumerate}
\item\label{rem-prop-of-weak-strong-stable-class-1} Every weakly stable class is stable under products. That is, if $f_i\colon X_i\rightarrow Y_i, i=1,2$ are in the class, then so is $f_1\times f_2 \colon  X_1 \times X_2 \rightarrow Y_1 \times Y_2$.

\item\label{rem-prop-of-weak-strong-stable-class-2} Let$\bfC_1\subset \bfC_2$ be a fully faithful embedding that preserves finite limits, and let $\mathrm{E}_1$ be a class of morphisms of $\bfC_1$ that stable under base change. Then one can define a class of morphisms $\mathrm{E}_2$ of $\bfC_2$ as those that are representable by morphisms in $\mathrm{E}_1$. That is, we define the class $\mathrm{E}_2$ of consisting of those morphisms $f: X\to Y$ in $\bfC_2$ such that for every $Y'\to Y$ with $Y'\in \bfC_1$, the fiber product $X':=Y'\times_YX$ belongs to $\bfC_1$ and the base change morphism $f': X'\to Y'$ belongs to $\mathrm{E}_1$. If $\mathrm{E}_1$ is  weakly (resp. strongly) stable, so is $\mathrm{E}_2$, and $(\bfC_1)_{\mathrm{E_1}}\to (\bfC_2)_{\mathrm{E_2}}$ is fully faithful.
\end{enumerate}
\end{remark}

Now let $\verti, \horiz$ be two weakly stable classes of morphisms in $\bfC$.
Let $\corr(\bfC)_{\verti;\horiz}$ denote the category of correspondences, as defined in  \cite[\textsection{6.1}]{liu2012enhanced} (as a quasi-category) or \cite[\textsection{7.1}]{Gaitsgory.Rozenblyum.DAG.vol.I} (as a complete Segal space).
Informally, objects of $\corr(\bfC)_{\verti;\horiz}$ are the same as those of $\bfC$ and morphisms from $X$ to $Y$  are given by diagrams
\begin{equation}\label{E: morphism-in-corr}
\begin{tikzcd}
 Z\arrow[r,"g"]\arrow[d,"f"] & X\\
 Y & 
\end{tikzcd}
\end{equation}
with $g\in \horiz$ and $f\in \verti$. We sometimes just write such diagram as $X\xleftarrow{g} Z\xrightarrow{f} Y$ for short, or as $X\stackrel{f\circ g^{-1}}{\dashrightarrow} Y$ or simply as $X\dashrightarrow Y$, to emphasize that such a morphism in $\corr(\bfC)_{\verti;\horiz}$ is a correspondence rather than an actual map.
The composition of the correspondences $X\leftarrow W_1 \rightarrow Y$ and $Y\leftarrow W_2\rightarrow Z$ is given by the correspondence 
\[
X\leftarrow W_1\leftarrow W:=W_1\times_YW_2\rightarrow W_2 \rightarrow Y.
\] 
Given $g\colon  Y\rightarrow X$ in $\horiz$ we will sometimes identify it with the correspondence  $X \xleftarrow{g} Y\xrightarrow{\id} Y$ and we refer to such morphisms as \textit{horizontal}. Similarly, given a morphism $f\colon  X\rightarrow Y$ in $\verti$ we will identify it with the correspondence $X \xleftarrow{\id} X \xrightarrow{f} Y$ and refer to such morphisms of $\corr(\bfC)_{\verti;\horiz}$ as \textit{vertical}. 
We usually write the class of all morphisms (resp. isomorphisms) in $\bfC$ as $\all$ (resp. $\isom$).
We simply write $\corr(\bfC)_{\all;\all}$ by $\corr(\bfC)$. 

\begin{remark}\label{rem-generalization-corr}
\begin{enumerate}
\item\label{rem-generalization-corr-1} The category $\corr(\bfC)_{\verti;\horiz}$ admits an $(\infty,2)$-categorical enhancement $\corr(\bfC)^{\mathrm{T}}_{\verti;\horiz}$ of the category of correspondences, depending on a certain class $\mathrm{T}\subset \verti\cap \horiz$ of morphisms of $\bfC$. A $2$-morphism between $X\xleftarrow{g'} Z'\xrightarrow{f'} Y$ and  $X\xleftarrow{g} Z\xrightarrow{f} Y$ in $\corr(\bfC)^{\mathrm{T}}_{\verti;\horiz}$ is given by a morphism ($r\colon Z'\to Z)\in \mathrm{T}$ with $f'\simeq f\circ r$ and $g'\simeq g\circ r$. See \cite[\textsection{7.1.1.2}]{Gaitsgory.Rozenblyum.DAG.vol.I} for details. We will not make use of such enhancement.
\item\label{rem-generalization-corr-2} In fact, in order to define $\corr(\bfC)_{\verti;\horiz}$ as an $\infty$-category, it is enough to impose weaker conditions on $\bfC$, $\verti$ and $\horiz$. Namely, instead of assuming that finite products exist in $\bfC$ and $\verti$ and $\horiz$ are stable under base change, it is enough to assume that morphisms in $\verti$ are stable under pullbacks by morphisms in $\horiz$ and vice versa (while keeping other assumptions on $\verti$ and $\horiz$ as in the definition of weakly stable class).
\end{enumerate}
\end{remark}

As $\bfC$ admits finite limits it is a symmetric monoidal category under the Cartesian monoidal structure. This induces a symmetric monoidal structure on $\corr(\bfC)_{\verti;\horiz}$, containing subcategories $\bfC_{\verti}$ and $(\bfC_{\horiz})^\op$ as symmetric monoidal subcategories. Informally, the tensor product of objects $X,Y$ in $\corr(\bfC)_{\verti;\horiz}$ is their product $X\times Y$ as objects of $\bfC$. See \cite[\textsection{6.1}]{liu2012enhanced} and  \cite[Chapter 9]{Gaitsgory.Rozenblyum.DAG.vol.I} for details. For our purpose, it is enough to recall the following. We write $\corr(\bfC)_{\verti;\horiz}^\otimes\to \mathrm{Fin}_*$ for the coCartesian fibration encoding the symmetric monoidal structure. Then morphisms over $\al: \langle m\rangle\to \langle n\rangle$ are given by
\begin{equation}\label{E: morphism-in-corr-sym}
\begin{tikzcd}
 (Z_j)_{1\leq j\leq n}\arrow[r,"g"]\arrow[d,"f"] & (X_i)_{1\leq i\leq m}\\
 (Y_j)_{1\leq j\leq n} & 
\end{tikzcd}
\end{equation}
where the vertical map is induced by $(Z_j\to Y_j)\in \verti$ and the horizontal map is given by $(Z_{\al(i)}\to X_i)\in \horiz$ if $\al(i)\in \langle n\rangle^\circ$.
Note that in general $X\times Y$ is not the product of $X$ and $Y$ in $\corr(\bfC)_{\verti;\horiz}$. For this reason, sometimes we write $X\otimes Y$ to emphasize we regard $X\times Y$ as the tensor product of $X$ and $Y$ in $\corr(\bfC)_{\verti;\horiz}$. 
In particular, it makes sense to talk about associative and commutative algebra objects in $\corr(\bfC)_{\verti;\horiz}$.

\begin{example}\label{ex:commutative-algebra-in-corr(C)}
Every object $X\in \bfC$ with the diagonal map $\Delta_X: X\to X\times X$ and the structural map $\pi_X: X\to \pt$ belonging to $\horiz$ has a natural commutative algebra structure in $(\bfC_{\horiz})^{\op}$ with the multiplication given by $\Delta_X$ and the unit given by $\pi_X$. (See \Cref{ex:coCartesian-relative-tensor}.)
Now, if $X$ and $Y$ are two objects satisfying the above properties, then every morphism $f:X\to Y$ belongs to $\horiz$. Namely, we may decompose $f=p_Y\circ (\id\times f)$, where $\id\times f: X\to X\times Y$ is the base change of $\Delta_Y: Y\to Y\times Y$ and therefor belongs to $\horiz$, and $p_Y: X\times Y\to Y$ is the projection which is the base change of $\pi_X: X\to \pt$ and therefore also belongs to $\horiz$.
It follows that we have a commutative algebra homomorphism in $(\bfC_{\horiz})^{\op}$ from $Y$ to $X$ induced by $f$. In particular, $X$ is a (left) $Y$-module.
As $(\bfC_{\horiz})^{\op}\to \corr(\bfC)_{\verti;\horiz}$ is a symmetric monoidal subcategory, we obtain the corresponding (maps between) commutative algebra objects in $\corr(\bfC)_{\verti;\horiz}$. If in addition $f\in \bfC_{\verti}$, then $f: X\to Y$ is naturally a morphism of $Y$-modules from $X$ to $Y$ on $\corr(\bfC)_{\verti;\horiz}$.

Now suppose that both $\pi_X$ and $\Delta_X$ belong to $\verti\cap \horiz$, then $X$ is self dual in $\corr(\bfC)_{\verti;\horiz}$ with unit given by $\Delta_X\circ \pi_X^{-1}$ and evaluation map given by $\pi_X\circ \Delta_X^{-1}$.
In particular, if $\verti=\horiz=\all$ so $\corr(\bfC)_{\verti;\horiz}=\corr(\bfC)$, then every object in $\corr(\bfC)$ is dualizable. This in particular induces a canonical symmetric monoidal equivalence
\[
\corr(\bfC)\cong \corr(\bfC)^{\op},\quad  X\mapsto X^\vee=X.
\]
\end{example}

\subsubsection{Algebras and modules in the category of correspondences}\label{sec:monads-correspondences-review}
We will be interested in a particular class of algebra objects and their bimodules in $\corr(\bfC)$. We review the description of algebras in terms of Segal objects and note how these constructions generalize to describe bimodules. 
The results here reproduce those in \cite[Chapter 9, \textsection{4}]{Gaitsgory.Rozenblyum.DAG.vol.I}. Unlikely \emph{loc. cit.}, our discussions avoid using the $(\infty,2)$-category formalism and stay entirely in the framework as developed in \cite[Chapter 4]{Lurie.higher.algebra}.

First recall the definition of Segal objects (also known as category objects) in an $\infty$-category.
\begin{definition}
A simplicial object $X_{\bullet}\colon  \Delta^{\op}\rightarrow \bfC$ is called a \textit{Segal object} if for every $n\geq 1$, the map
\[
X_n \rightarrow X_1\times_{X_0}X_1\times_{X_0}\cdots \times_{X_0}X_1
\]
induced by the maps $\delta_{i}\colon  [1]\cong \{i,i+1\}\subset [n]$ for $i=0,1,\dots, n-1$, is an equivalence.
\end{definition}

\begin{remark}
If $\bfC$ is an ordinary category, a Segal object is fully determined by the objects $X_0, X_1$, the boundary maps $d_1,d_0\colon  X_1 \rightarrow X_0$, $d_1\colon  X_2 \rightarrow X_1$ and the degeneracy map $s\colon  X_0 \rightarrow X_1$. These define a category object of $\bfC$ in the usual sense. Namely, $X_0$ is the class of objects, $X_1$ the morphism objects, the morphisms $d_1,d_0\colon  X_1 \rightarrow X_0$ as source and target. The composition is given by the morphism $d_1\colon  X_2 \rightarrow X_1$ and unit by $d_1\colon  X_2 \rightarrow X_1$.
\end{remark}

\begin{example}\label{ex:Cech-nerve}
The \v{C}ech nerve $X_{\bullet}\to Y$ of a morphism $f: X\to Y$ (see \cite[\textsection{6.1.2}]{Lurie.higher.topos.theory}), where
\[
X_n= \overbrace{X\times_YX\times\cdots\times_YX}^{n+1},
\]
is easily seen to be a Segal object of $\bfC$. Indeed, it is even a groupoid object, in the sense of \cite[Definition 6.1.2.7]{Lurie.higher.topos.theory}). 
This will be our main example.
\end{example}

\begin{example}\label{ex: monoid object vs algebra object}
A Segal object $X_{\bullet}$ with $X_0=\pt$ is a monoid object (in the sense of \cite[\textsection{4.1.2}]{Lurie.higher.algebra}). Giving such a monoid object is equivalent to giving an associative algebra object in $\bfC$ by \cite[Proposition 4.1.2.10]{Lurie.higher.algebra}.
\end{example}

This last example admits the following generalization. Let $\bfC$ be a category with finite limits. Recall from \Cref{ex:coCartesian-relative-tensor} that every object $X\in\bfC$ is a commutative algebra object in $\bfC^{\op,\sqcup}$, if the category $\bfC^{\op}$ with coCartesian symmetric monoidal structure. In addition, we have the monoidal category ${}_X\BMod_{X}(\bfC^{\op,\sqcup})$ of $X$-bimodules in $\bfC^{\op,\sqcup}$, and its opposite category ${}_X\BMod_{X}(\bfC^{\op,\sqcup})^{\op}$.

\begin{proposition}\label{prop: segal vs algebra object}
Let $\bfC$ be a category with finite limits.
There is natural equivalence from the category of Segal objects in $\bfC$ with $X_0=X$ to the category of algebra objects in ${}_X\BMod_{X}(\bfC^{\op,\sqcup})^{\op}$.
\end{proposition}
See also  \cite[Proposition 9.4.1.5]{Gaitsgory.Rozenblyum.DAG.vol.I}.
Note that if $X=\pt$ is the finite object in $\bfC$, then ${}_X\BMod_{X}(\bfC^{\op,\sqcup})^{\op}$ is nothing but $\bfC$ equipped with the Cartesian symmetric monoidal structure. Therefore, the above statement does generalize \Cref{ex: monoid object vs algebra object}. 

\begin{proof}
As the proof largely follows from the strategy of \cite[Proposition 4.1.2.10]{Lurie.higher.algebra}. We only give a sketch.
We use \cite[Proposition 4.1.3.19]{Lurie.higher.algebra} to identify algebra objects in ${}_X\BMod_{X}(\bfC^{\op,\sqcup})^{\op}$ as functors of planar operads $F:\Delta^{\op}\to ({}_X\BMod_{X}(\bfC^{\op,\sqcup})^{\op})^{\oast}$. 
Let $\pi: ({}_X\BMod_{X}(\bfC^{\op,\sqcup})^{\op})^{\oast}\to \bfC$ be as in \Cref{ex:coCartesian-relative-tensor}.
Then one checks that given a functor of planar operads $F:\Delta^{\op}\to ({}_X\BMod_{X}(\bfC^{\op,\sqcup})^{\op})^{\oast}$ amounts to a Segal object $\pi\circ F$ in $\bfC$ with $X_0=X$.
\end{proof}

\begin{lemma}\label{lem: construction of algebra object in corr}
There is a canonical lax-monoidal functor ${}_X\BMod_{X}(\bfC^{\op,\sqcup})^{\op}\to \corr(\bfC)$.
\end{lemma}
\begin{proof}
The desired functor in the lemma is given by the compositions
\[
{}_X\BMod_{X}(\bfC^{\op,\sqcup})^{\op}\to {}_X\BMod_{X}(\corr(\bfC))^{\op}\to \corr(\bfC)^\op \cong \corr(\bfC).
\]
where the first functor comes from the symmetric monoidal functor $\bfC^{\op,\sqcup}\to \corr(\bfC)$, and the last equivalence comes from the end of \Cref{ex:commutative-algebra-in-corr(C)}.
\end{proof}

We thus recover \cite[Corollary 9.4.4.5]{Gaitsgory.Rozenblyum.DAG.vol.I} as follows. 
We note that unlike \emph{loc. cit.}, our proof of the above statement stays in $(\infty,1)$-categorical formalism.

\begin{corollary}\label{thm:gaitsgory-rozen-segal-monad}
There is a natural functor from the category of Segal objects $X_\bullet$ in $\bfC$ with $X_0=X$ to the category of associative algebra objects in $\corr(\bfC)$.
\end{corollary}

Roughly speaking, the functor sends $X_\bullet$ to $X_1\in \corr(\bfC)$ with multiplication and unit given by the correspondences
\begin{equation}\label{eq-product-Segal-Obj}
\begin{tikzcd}
X_1\times_{X_0} X_1 \simeq X_2 \arrow[d,"m:=d_1"] \arrow[r,"\eta:=d_0\times d_2"]& X_1 \times X_1 \\
X_1 &  &
\end{tikzcd}, \quad
\begin{tikzcd}
X_0 \arrow[d,"u:=s"] \arrow[r, "\pi_{X_0}"]& \pt \\
X_1 & 
\end{tikzcd}.
\end{equation}
In particular, if $X_\bullet$ is the groupoid object arising from the \v{C}ech nerve of a morphism $f:X\to Y$ as in \Cref{ex:Cech-nerve}, then $X\times_YX$ has a natural algebra structure in $\bfC$ with the multiplication and unit maps are given by
\begin{equation}\label{eq-convolution-product}
\begin{tikzcd}
 X\times_{Y} X\times_{Y} X\arrow[rr,"\id\times\Delta_{X}\times\id"]\arrow[d,"\id\times f\times \id"]&& (X\times_{Y}X) \times (X\times_{Y} X)\\
 X\times_{Y} X && 
\end{tikzcd}, \quad
\begin{tikzcd}
 X\arrow[r]\arrow[d,"\Delta_{X/Y}"]& \pt\\
 X\times_{Y} X & 
\end{tikzcd}.
\end{equation}
This multiplication is usually called the convolution product.

\begin{remark}\label{rem:segal.objects.morphisms.vert.horiz}
Assume that the Segal object $X_\bullet$ is such that:
\begin{itemize}
    \item All morphisms of the simplicial object $X_\bullet$ are in $\bfC_{\verti}$;
    \item The diagonal map $\Delta_{X_0}\colon X_0 \rightarrow X_0 \times X_0$ and the structural map $\pi_{X_0}:X_0 \to \pt$ are in $\bfC_{\horiz}$.
\end{itemize}
Then the associative algebra object $X_1$ can be realized as an associative algebra object of the monoidal category $\corr(\bfC)_{\verti;\horiz}$. If $X_\bullet$ is the groupoid object arising from the \v{C}ech nerve of a morphism $f:X\to Y$ as in \Cref{ex:Cech-nerve}, then the above assumptions hold if 
\begin{itemize}
\item $f$ and $\Delta_{X/Y} \colon X\to X\times_Y X$ belong to $\bfC_{\verti}$;
\item $\Delta_X: X\to X\times X$ and $\pi_X: X\to \pt$ belong to $\bfC_{\horiz}$.
\end{itemize}
\end{remark}

As just discussed above, for $X\to Y$ satisfying certain (mild) conditions, the fiber product $X\times_YX$ has a natural algebra structure in $\corr(\bfC)_{\verti;\horiz}$. Our next goal is to produce its left modules in $\corr(\bfC)_{\verti;\horiz}$.

The following definition generalizes \cite[Definition 4.2.2.2]{Lurie.higher.algebra}.
\begin{definition}\label{def: left module over Segal}
A left module over a Segal object in $\bfC$ consists of a map of simplicial object $Q_\bullet\to X_\bullet$ in $\bfC$ with $X_\bullet$ a Segal object
such that for every $n$, the map $Q_n\to X_n$ and $Q_n=Q([n])\to Q(\{n\})\cong Q_0$ exhibits $Q_n$ as the product $X_n\times_{X_0}Q_0$.
\end{definition}

Then analogously to \Cref{thm:gaitsgory-rozen-segal-monad}, we have
\begin{proposition}
There is a natural functor from the category of $Q_\bullet\to X_\bullet$ of left modules over Segal objects in $\bfC$ to the category of algebras and left modules in $\corr(\bfC)$.
\end{proposition}

\begin{example}\label{sec-categorical.action.on.bimodules}
Let $f: X\to Y$ be a morphism and $g: Z\to Y$ another morphism. It follows that $X\times_YZ$ admits a left action of $X\times_YX$. 

Recall that given two algebras $A$ and $B$ in a symmetric monoidal category $\bfR$, the category of $A\mbox{-}B$-bimodules $M$ is equivalent to the category of left $A\otimes B^{\rev}$-modules. It follows that if $Z$ is equipped with a morphism $g: Z\to Y_1\times Y_2$, then $X_1\times_{Y_1}Z\times_{Y_2}X_2$ has a natural $(X_1\times_{Y_1}X_2)\mbox{-}(X_2\times_{Y_2}X_2)$-bimodule structure.
\end{example}

\subsection{Sheaf theories}\label{sec:symmetric-monoidal-and-projection-for-corr}
\subsubsection{The formalism of a sheaf theory}\label{ss: formalism-of-abstract-sheaf-theory}
We fix $\bfC$ as before, and fix a commutative ring $\La$.
\begin{definition}
A(n abstract) sheaf theory with coefficient in $\La$ (or sometimes called a $3$-functor formalism in literature) of $\bfC$ is a lax symmetric monoidal functor 
\begin{equation}\label{eq:abstract-sheaf-theory}
\der\colon \corr(\bfC)_{\verti;\horiz} \rightarrow   \lincat_{\Lambda}.
\end{equation}
For a horizontal morphism
\(
X \xleftarrow{g} Y\xrightarrow{\id} Y \) we will denote the corresponding functor by $g\horizl\colon  \der(X) \rightarrow \der(Y)$. For a vertical morphism
\(
X \xleftarrow{\id} X \xrightarrow{f} Y
\) we denote the corresponding functor by $f\vertl\colon  \der(X) \rightarrow \der(Y)$. Then for a general correspondence 
\(
X \xleftarrow{g} Z \xrightarrow{f} Y
\) 
the associated functor is (isomorphic to) $f\vertl\circ g\horizl$. 
\end{definition}

Let us recall some structures encoded by such a functor. See also  \cite[\textsection{6.2}]{liu2014enhanced} and \cite[Part III, Introduction]{Gaitsgory.Rozenblyum.DAG.vol.I} for some discussions.
\begin{enumerate}
\item The functoriality of $\der$ encodes a ``base change theorem". Namely, let
\begin{equation}\label{eq:base-change-diagram-app}
\begin{tikzcd}
    X' \arrow[r,"g'"]\arrow[d,"f'"] & X\arrow[d,"f"]\\
    Y' \arrow[r,"g"] & Y
\end{tikzcd}
\end{equation}
be a pullback square 
in $\bfC$. If $f\in \verti$ and $g\in \horiz$, then $f'\in \verti$, $g'\in \horiz$ and part of the data of the functor $\der$ is to give an isomorphism of functors
\begin{equation}\label{eq:abstract-base-change}
(f')\vertl\circ (g')\horizl\cong g\horizl\circ f\vertl.
\end{equation}

\item The category $\der(\pt)$ is a commutative algebra in $\lincat_\La$. The lax symmetric monoidal structure of $\der$ provides a functor
\begin{equation}\label{eq:abstract-exterior-tensor-product}
\boxtimes_{\La}\colon   \der(X)\otimes_{\La} \der(Y) \rightarrow \der (X\times Y), \quad X,Y\in \bfC,
\end{equation}
and for $f_i\colon X_i\rightarrow Y_i, \ i=1,2$ in $\bfC_{\horiz}$, a canonical isomorphism
\begin{equation}\label{eq:abstract-exterioir-pullback} 
(f_1)\horizl(\mF_1)\boxtimes_{\La}(f_2)\horizl(\mF_2)\cong (f_1\times f_2)\horizl(\mF_1\boxtimes_{\La} \mF_2),
\end{equation}
and for $f_i\colon X_i\rightarrow Y_i, \ i=1,2$ in $\bfC_{\verti}$, a canonical isomorphism
\begin{equation}\label{eq:abstract-Kunneth-formula}
(f_1)\vertl(\mF_1)\boxtimes_{\La} (f_2)\vertl(\mF_2)\cong (f_1\times f_2)\vertl(\mF_1\boxtimes_{\La} \mF_2),
\end{equation}
together with all necessary higher coherence conditions. When the coefficient $\La$ is clear from the context, we also write $\boxtimes$ instead of $\boxtimes_{\La}$.

Let $X$ be as in Example \ref{ex:commutative-algebra-in-corr(C)}. This induces a symmetric monoidal structure on the category $\der(X)$. Informally, the symmetric monoidal structure is given by the composition
\begin{equation}\label{eq:abstract-interior-tensor-product}
\der(X) \otimes_{\La} \der(X) \xrightarrow{\boxtimes_{\La}} \der(X\times X) \xrightarrow{\Delta_X\horizl} \der(X),\quad \mF,\mG\mapsto \mF \otimes \mG:= \Delta_X\horizl(\mF\boxtimes_{\La} \mG).
\end{equation}
We let $\La_X\in\der(X)$ denote the unit object with respect to this symmetric monoidal structure, which corresponds to the functor
\begin{equation}\label{eq:abstract-unit}
\der(\pt)\xrightarrow{\pi_X\horizl}\der(X).
\end{equation}
In addition, for $f\colon X\rightarrow Y$ in $\bfC_{\horiz}$ as in \Cref{ex:commutative-algebra-in-corr(C)}, the functor $f\horizl: \der(Y)\to \der(X)$ is a symmetric monoidal functor, and therefore endows $\der(X)$ with a structure of a $\der(Y)$-module category. 
If in addition $f\in \bfC_{\verti}$ as well, then $f\vertl \colon \der(X) \rightarrow \der(Y)$
is a morphism of $\der(Y)$-modules. In particular, for $\mF\in \der(X)$, $\mG\in \der(Y)$ we have a canonical equivalence
\begin{equation}\label{eq:abstract-projection-formula}
f\vertl(\mF)\otimes \mG \cong f\vertl(\mF\otimes f\horizl(\mG)),
\end{equation}
which encodes a ``projection formula" for $f\vertl$ and $f\horizl$. 

\item We can pass to (not necessarily continuous) right adjoints. For $(g: X\to Y)\in \bfC_{\horiz}$, let $g\horizr$ be the (not necessarily continuous) right adjoint of $g\horizl$, and for $(f: X\to Y)\in \bfC_{\verti}$, let $f{\vertr}$ be the (not necessarily continuous) right adjoint of $f\vertl$. In addition, for $X$ as in Example \ref{ex:commutative-algebra-in-corr(C)}, the symmetric monoidal structure of $\der(X)$ is closed. That is, for every pair of objects $\mF_1,\mF_2\in \der(X)$ there is an object $\underline\Hom(\mF_1,\mF_2)\in \der(X)$ such that for every $\mG\in \der(X)$ there is a canonical equivalence
\begin{equation}\label{eq:abstract-internal-hom}
\map_{\der(X)}\big(\mG,\underline\Hom(\mF_1,\mF_2)\big) \simeq \map_{\der(X)}\big(\mG\otimes \mF_1,\mF_2\big).
\end{equation}
See \Cref{SS: la-linear categories}.
Note we have
\begin{equation}\label{eq:abstract-internal-hom-tensor}
\rhom(\mF_1\otimes\mF_2,\mF_3)=\rhom(\mF_1,\rhom(\mF_2,\mF_3)),
\end{equation}
and for $(f:X\to Y)\in\bfC_{\horiz}$ and $\mF\in\der(Y)$ and $\mG\in\der(X)$,
\begin{equation}\label{eq:abstract-internal-hom-star-pull-push}
\rhom(\mF, f\horizr\mG)=f\horizr\rhom(f\horizl\mF,\mG).
\end{equation}
In addition, along with the co-unit of the adjunction $(f\vertl,f\vertr)$, the projection formula \eqref{eq:abstract-projection-formula} gives, for every $\mF,\mG\in\der(Y)$, a natural map
\begin{equation}\label{eq:abstract-upper-shriek-projection-formula}
   f\vertr(\mG) \otimes f\horizl(\mF)\to f\vertr(\mG\otimes \mF)
\end{equation}
adjoint to
\[
f\vertl (f\vertr(\mG) \otimes f\horizl(\mF) ) \simeq  
f\vertl f\vertr(\mG) \otimes \mF  \rightarrow 
\mG \otimes \mF, 
\]
In particular, one has the natural transformation of functors
\begin{equation}\label{eq:abstract-!-*-pullback}
    f\vertr(\La_Y)\otimes f\horizl\to f\vertr
\end{equation}
In addition, we have
\begin{equation}\label{eq:abstract-hom-pull-push}
   f\vertr\rhom(\mF,\mG)\simeq \rhom(f\horizl\mF,f\vertr\mG),\quad  f\horizr\rhom(\mF,f\vertr\mG)\simeq \rhom(f\vertl\mF,\mG).
\end{equation}
\end{enumerate}

\begin{remark}\label{rem:relative-sheaf-theory}
Let $S\in \bfC$ be as in Example \ref{ex:commutative-algebra-in-corr(C)}. Then
there is the (non-full) embedding 
\[
  \corr(\bfC_{/S})_{\verti;\horiz}\to \corr(\bfC)_{\verti;\horiz},
\]
which is a lax symmetric monoidal functor. 
Therefore, one can restrict $\der$ along this embedding to obtain a sheaf theory on $\bfC_{/S}$, denoted by $\der_{/S}$. The lax symmetric monoidal structure is provided by 
\[
\der(X)\otimes_\La \der(Y)\xrightarrow{\boxtimes_\La} \der(X\times Y) \xrightarrow{\Delta_S\horizl} \der(X\times_S Y).
\]
\end{remark}

\begin{remark}\label{R:Kunneth-type-formula}
In fact, a sheaf theory $\der$ automatically factors as lax symmetric monoidal functors
\[
\der\colon \corr(\bfC)_{\verti;\horiz}\to \lincat_{\der(\pt)}\to \lincat_\La.
\]
Informally, this means that each $\der(X)$ is a $\der(\pt)$-linear category and
 the functor \eqref{eq:abstract-exterior-tensor-product} factors as
$\der(X)\otimes_\La\der(Y)\to \der(X)\otimes_{\der(\pt)}\der(Y)\to \der(X\times Y).$
In many examples, the functor 
\begin{equation}\label{eq:refined-abstract-exterior-tensor-product}
\boxtimes_{\der(\pt)}: \der(X)\otimes_{\der(\pt)}\der(Y)\to \der(X\times Y)
\end{equation} 
is fully faithful and admits a $\der(\pt)$-linear right adjoint. The fully faithfulness is equivalent to a K\"unneth type formula. 
\end{remark}

\begin{remark}\label{rem-dualizability in corr}
Let $X\in \bfC$ such that both $\pi_X$ and $\Delta_X$ are in $\bfC_{\verti}\cap \bfC_{\horiz}$, then $X$ is self-dual in $\corr(\bfC)_{\verti;\horiz}$, see \Cref{ex:commutative-algebra-in-corr(C)}. It follows that if $\der(X)\otimes_{\der(\pt)}\der(X)\to \der(X\times X)$ is an equivalence (e.g. the sheaf theory $\der$ is symmetric monoidal (rather than just lax symmetric monoidal)), then $\der(X)$ is self-dual as a $\der(\pt)$-linear category. Explicitly,  the unit and the evaluation for the self duality of $\der(X)$ are given by
\begin{equation}\label{eq: unit by diagonal abstract sheaf}
\Delta\vertl \La_X\in \der(X\times X)\cong \der(X)\otimes_{\der(\pt)} \der(X),
\end{equation}
\begin{equation}\label{eq: Frob alg abstract sheaf}
\der(X)\otimes_{\der(\pt)} \der(X)\to \Mod_\La,\quad (\mF_1, \mF_2)\mapsto  (\pi_X)\vertl(\mF_1\otimes\mF_2).
\end{equation}
Note that in this case, $\der(X)$ has a Frobenius algebra structure as in the discussion as in \Cref{ex: duality via Frobenius-structure}. Namely, the functor $\la$ in \Cref{ex: duality via Frobenius-structure} is given by $(\pi_X)\vertl$. In this case, \eqref{eq:abstract-verdier-dual.II} is identified with \eqref{eq:abstract smooth dual 2}. 

Sometimes, even if $\der(X)\otimes_{\der(\pt)}\der(X)\to \der(X\times X)$ is not an equivalence, $\der(X)$ may still has a Frobenius algebra structure given by \eqref{eq: Frob alg abstract sheaf} (see \Cref{rem:Verdier duality via dualizability in corr} for an example).  When this is the case, the self-duality $\verd^\la$ as in \Cref{ex: duality via Frobenius-structure} will be denoted as
\[
\verd_X^{\der}: \der(X)^{\vee}\cong \der(X).
\]
\end{remark}

We also give a very useful criterion to determine when $\der(X)\otimes_{\der(\pt)}\der(Y)\to \der(X\times Y)$ is an equivalence. This is of course well-known, in any concrete sheaf theory.
\begin{lemma}\label{lem: criterion of tensor product equivalence via diagonal}
Suppose \eqref{eq:refined-abstract-exterior-tensor-product} is always fully faithful for every two objects in $\bfC$.
Let $X\in \bfC$ such that both $\pi_X$ and $\Delta_X$ are in $\bfC_{\verti}\cap \bfC_{\horiz}$. Then the following are equivalent.
\begin{enumerate}
\item\label{lem: criterion of tensor product equivalence via diagonal-1} $\der(X)\otimes_{\der(\pt)}\der(Y)\to \der(X\times Y)$ is an equivalence for every $Y$ such that both $\pi_Y$ and $\Delta_Y$ belong to $\bfC_{\verti}\cap \bfC_{\horiz}$;
\item\label{lem: criterion of tensor product equivalence via diagonal-2} $\der(X)\otimes_{\der(\pt)}\der(X)\to \der(X\times X)$ is an equivalence;
\item\label{lem: criterion of tensor product equivalence via diagonal-3} $\der(X)\otimes_{\der(\pt)}\der(X)\xrightarrow{\boxtimes_{\der(\pt)}} \der(X\times X)$ is fully faithful and $(\Delta_X)\vertl\La_X$ belongs to the essential image of $\boxtimes_{\der(\pt)}$.
\end{enumerate}
\end{lemma}
\begin{proof}
Clearly it is enough to show that \eqref{lem: criterion of tensor product equivalence via diagonal-3} implies \eqref{lem: criterion of tensor product equivalence via diagonal-1}.
For simplicity, we write $\otimes$ instead of  $\otimes_{\der(\pt)}$, and write $Z=X\times Y$. Then both $\pi_Z$ and $\Delta_Z$ belong to $\bfC_{\verti}\cap \bfC_{\horiz}$.
We notice that the base change implies that
\[
(p_1)\vertl(\id\times \Delta_Z)\horizl((\Delta_{Z})\vertl\La_{Z}\boxtimes -))\cong \id_{\der(Z)}.
\]

We note that $\La_{Z}=\La_X\boxtimes \La_Y$. It follows that $(\Delta_{Z})\vertl\La_{Z}\in \der(X)\otimes \der(X)\otimes \der(Y\times Y)\subset \der(Z\times Z)$.
On the other hand,
note that for $\mK=\mK_1\boxtimes \mK_2\boxtimes \mK_3\in \der(X)\otimes \der(X)\otimes \der(Y\times Y)\subset \der(Z\times Z)$ and any $\mF\in \der(Z)$, we have
\[
(\id\times \Delta_Z)\horizl(\mK\boxtimes \mF)\cong \mK_1\boxtimes \mF'
\]
for some $\mF'\in \der(Y\times X\times Y)$. It follows that
\[
(p_1)\vertl(\id\times \Delta_Z)\horizl(\mK\boxtimes \mF)=\mK_1\boxtimes (p_1)\vertl\mF'\in \der(X)\otimes \der(Y).
\]
Combining these observations, we see that $\der(X)\otimes\der(Y)\to \der(Z)$ is an equivalence, as desired.
\end{proof}

\begin{remark}\label{rem-conv-product}
Let $X_\bullet$ be a Segal object in $\bfC$ as in \Cref{rem:segal.objects.morphisms.vert.horiz}. Then it follows from \Cref{thm:gaitsgory-rozen-segal-monad} that $\der(X_1)$ has a natural monoidal structure, usually called the convolution monoidal structure. This is different from the natural symmetric monoidal structure on $\der(X_1)$ (assuming $X_1$ is as in \Cref{ex:commutative-algebra-in-corr(C)}). For example, 
the monoidal unit of the former is given by 
\[
\mathbf{1}_{\der(X_1)}=s\vertl(\La_{X_0}),
\] 
while the unit of the latter is $\La_{X_1}$. In addition, the convolution monoidal structure usually is not symmetric.

In the particular case when $X_\bullet$ arises as the \v{C}ech nerve $f:X\to Y$, then $\der(X\times_YX)$ has a monoidal structure given by the convolution product. 
\end{remark}

\begin{remark}\label{rem: codomain of sheaf theory}
\begin{enumerate}
\item In the definition of a sheaf theory, it makes sense to replace $\lincat_\Lambda$ by any other symmetric monoidal $2$-category $\bfR$. For example, one can consider sheaf theory valued in $\bfR=\cat$, or in $\bfR=\cptcat_\La$. Most of the above discussions carry through, except \Cref{lem: criterion of tensor product equivalence via diagonal} and those related to the right adjoint functors $f\vertr$ and $g\horizr$.

\item Recall that natural functor $\lincat_\La\to \cat$ lax symmetric monoidal, realizing $\lincat_\La^\otimes$ as a (non-full) subcategory of $\cat^\otimes$. It follows from \Cref{lemma: construct functors from homotopy cat} below that giving a sheaf theory $\der$ amounts to giving a lax symmetric monoidal functor $\corr(\bfC)_{\verti;\horiz} \rightarrow  \cat$ such that for every $X$ we have $\der(X)\in \lincat_\La$, and such that for every $
X \xleftarrow{g} Z \xrightarrow{f} Y$ the functors $g\horizl$ and $f\vertl$ have $\La$-linear structure.

\item\label{rem: codomain of sheaf theory-3} Suppose the functor \eqref{eq:abstract-exterior-tensor-product} takes values in $\cat$. As explained in \cite{Lu.Zheng.relative.Lefschetz}, via the symmetric monoidal Grothendieck construction (e.g. see \cite[Proposition A.2.1]{Hinich} for an $\infty$-categorical version),
the sheaf theory $\der$ can also be (largely) encoded as a symmetric monoidal $2$-category $\corr^{\der}(\bfC)_{\verti;\horiz}$, usually called the category of cohomological correspondences.
Informally, objects of $\corr^{\der}(\bfC)_{\verti;\horiz}$ consist of pairs $(X,\mF)$ where $X\in \bfC$ and $\mF\in \der(X)$, and morphisms between $(X,\mF)$ and $(Y,\mG)$ consist of pairs $(Z,u)$, where $Y\xleftarrow{g} Z\xrightarrow{f} X$ is a correspondence, and $u: f\vertl(g\horizl\mF)\to \mG$ is a morphism in $\der(Y)$. The symmetric monoidal structure given by $(X,\mF)\otimes (Y,\mG)=(X\times Y, \mF\boxtimes_{\La}\mG)$, with the unit is given by $(\pt,\La_{\pt})$. 

\item\label{rem-sheaf-theory-cpt-and-dualizable} If the functor \eqref{eq:abstract-exterior-tensor-product} takes values in $\cptcat_\La$, by composing $\der$ with the duality functor \eqref{eq: duality functor, cpt}, we may obtain a new sheaf theory still taking value in $\cptcat_\La$.
\end{enumerate}
\end{remark}

\subsubsection{Additional adjunction and base change}\label{sec:additional-base-change}
In practice, a sheaf theory usually satifies additional adjunction and base change properties, besides those already encoded by functoriality as mentioned above. One formulation of these additional structures is via $2$-categorical enhancement of a sheaf theory, as in \cite{Gaitsgory.Rozenblyum.DAG.vol.I}. As we do not make any use of such formalism, we will add certain additional assumptions to a sheaf theory. As will be explained in  \Cref{rem-sheaf-theory-for-adjoint-factorization}, such assumptions are very closely related to the $2$-categorical enhancement. That is to say, the existence $2$-categorical enhancement of a sheaf theory in many cases is a property rather than an additional structure. \footnote{In fact, we find it is more flexible to add these assumptions rather than to adding $2$-categorical enhancement to a sheaf theory. For example,  for the usual sheaf theory of \'etale cohomology, $2$-categorical enhancement as developed in \cite{Gaitsgory.Rozenblyum.DAG.vol.I} can only encode either adjunctions for proper morphisms, or adjunctions for open morphisms, depending on an additional chosen class of morphisms $\mathrm{T}\subset\verti\cap \horiz$, but not both (unless one allows the $2$-morphisms in $\corr(\bfC)$ to be correspondences as well).}

A common setup is as follows. We make use of the Cartesian diagram \eqref{eq:base-change-diagram-app}.

\begin{assumptions}\label{assumptions.base.change.sheaf.theory.H} 
Let $(f_0: X_0\to Y_0)\in \bfC_{\horiz}$. Assume that
\begin{enumerate}
    \item\label{assumptions.base.change.sheaf.theory.H-0} for any of its base change $f:X\to Y$, the functor $f\horizl$ admits a continuous right adjoint $f\horizr=(f\horizl)^R$. 
    \end{enumerate}  
Under this assumption, and given a Cartesian diagram \eqref{eq:base-change-diagram-app}, we further assume that
\begin{enumerate}[resume]    
    \item\label{assumptions.base.change.sheaf.theory.H-1} for $g\in \bfC_{\verti}$, the natural Beck-Chevalley map is an isomorphism $g\vertl\circ (f')\horizr\cong f\horizr\circ (g')\vertl$;
    
     \item\label{assumptions.base.change.sheaf.theory.H-2} for $g\in \bfC_{\horiz}$, the natural Beck-Chevalley map is an isomorphism $g\horizl\circ f\horizr\cong (f')\horizr\circ (g')\horizl$;
               
     \item \label{assumptions.base.change.sheaf.theory.H-3} for $\mF\in \der(X)$, $\mG\in \der(Z)$ we have the natural isomorphism (which is the adjunction of \eqref{eq:abstract-exterioir-pullback})
      $f\horizr(\mF)\boxtimes_{\La} \mG \cong (f\times\id)\horizr(\mF\boxtimes_{\La} \mG)$.
\end{enumerate}  
\end{assumptions}

\begin{assumptions}\label{assumptions.base.change.sheaf.theory.HL}
Let $(f: X_0\to Y_0)\in \bfC_{\horiz}$. Assume that
\begin{enumerate}
    \item\label{assumptions.base.change.sheaf.theory.HL-0} for any of its base change $f:X\to Y$, the functor $f\horizl$ admits a left adjoint $(f\horizl)^L$. 
    \end{enumerate}  
Under this assumption, and given a Cartesian diagram \eqref{eq:base-change-diagram-app}, we further assume that
\begin{enumerate}[resume]
   \item\label{assumptions.base.change.sheaf.theory.H-4} for $g\in \bfC_{\verti}$, the natural Beck-Chevalley map is an isomorphism $(f\horizl)^L\circ (g')\vertl\cong g\vertl\circ ((f')\horizl)^L$;
    
     \item\label{assumptions.base.change.sheaf.theory.H-5} for $g\in \bfC_{\horiz}$, the natural Beck-Chevalley map is an isomorphism $((f')\horizl)^L\circ (g')\horizl\cong g\horizl\circ (f\horizl)^L$;

     \item \label{assumptions.base.change.sheaf.theory.H-6} for $\mF\in \der(X)$, $\mG\in \der(Z)$ we have the natural isomorphism $((f\times\id)\horizl)^L(\mF\boxtimes_{\La}\mG)\simeq (f\horizl)^L(\mF)\boxtimes_{\La}\mG$.
\end{enumerate}  
\end{assumptions}

\begin{assumptions}\label{assumptions.base.change.sheaf.theory.V} 
Let $(f: X_0\to Y_0)\in \bfC_{\verti}$. Assume that
\begin{enumerate}
\item\label{assumptions.base.change.sheaf.theory.V-0} for any of its base change $f:X\to Y$, the functor $f\vertl$ admits a continuous right adjoint $f\vertr=(f\vertl)^R$. 
\end{enumerate} 
Under this assumption, and given a Cartesian diagram \eqref{eq:base-change-diagram-app}, we further assume that
\begin{enumerate}[resume]
   \item\label{assumptions.base.change.sheaf.theory.V-1} for $g\in \bfC_{\verti}$, the natural Beck-Chevalley map is an isomorphism $(g')\vertl\circ (f')\vertr\cong f\vertr\circ g\vertl$;
   
    \item\label{assumptions.base.change.sheaf.theory.V-2} for $g\in \bfC_{\horiz}$, the natural Beck-Chevalley map is an isomorphism $(g')\horizl\circ f\vertr\cong (f')\vertr\circ g\horizl$; 
    
    \item \label{assumptions.base.change.sheaf.theory.V-3} for $\mF\in \der(Y)$, $\mG\in \der(Z)$ we have the natural isomorphism (which is the adjunction of \eqref{eq:abstract-Kunneth-formula}) $f\vertr(\mF)\boxtimes_{\La} \mG\cong (f\times\id)\vertr(\mF\boxtimes_{\La}\mG)$. 
\end{enumerate}   
\end{assumptions}

\begin{assumptions}\label{assumptions.base.change.sheaf.theory.VL} 
Let $(f: X_0\to Y_0)\in \bfC_{\verti}$. Assume that
\begin{enumerate}
\item\label{assumptions.base.change.sheaf.theory.VL-0} for any of its base change $f:X\to Y$, the functor $f\vertl$ admits a left adjoint $(f\vertl)^L$. 
\end{enumerate} 
Under this assumption, and given a Cartesian diagram \eqref{eq:base-change-diagram-app}, we further assume that
\begin{enumerate}[resume]
   \item\label{assumptions.base.change.sheaf.theory.V-4} for $g\in \bfC_{\verti}$, the natural Beck-Chevalley map is an isomorphism $ (f\vertl)^L\circ g\vertl\cong  (g')\vertl\circ ((f')\vertl)^L$;
   
    \item\label{assumptions.base.change.sheaf.theory.V-5} for $g\in \bfC_{\horiz}$, the natural Beck-Chevalley map is an isomorphism $((f')\vertl)^L\circ g\horizl\cong (g')\horizl\circ (f\vertl)^L$;
    
    \item \label{assumptions.base.change.sheaf.theory.V-6} for $\mF\in \der(Y)$, $\mG\in \der(Z)$ we have the natural isomorphism $((f\times\id)\vertl)^L(\mF\boxtimes_{\La}\mG)\cong (f\vertl)^L(\mF)\boxtimes_{\La}\mG$.
     \end{enumerate}   
\end{assumptions}

\begin{remark}\label{rem-additional-base.change.sheaf.theory} 
\begin{enumerate}
\item\label{rem-additional-base.change.sheaf.theory-0}  The class of morphisms satisfying \Cref{assumptions.base.change.sheaf.theory.H} \eqref{assumptions.base.change.sheaf.theory.H-0} - \eqref{assumptions.base.change.sheaf.theory.H-3} is weakly stable. 
We denote this class of morphisms by $\mathrm{HR}$, standing for ``horizontally right adjointable". Similarly, we let $\mathrm{HL}$ denote the class of morphisms satisfying \Cref{assumptions.base.change.sheaf.theory.HL} \eqref{assumptions.base.change.sheaf.theory.HL-0}-\eqref{assumptions.base.change.sheaf.theory.H-6}, let $\mathrm{VR}$ denote the class of morphisms satisfying \Cref{assumptions.base.change.sheaf.theory.V}  \eqref{assumptions.base.change.sheaf.theory.V-0}-\eqref{assumptions.base.change.sheaf.theory.V-3}, and let $\mathrm{VL}$ denote the class of morphisms satisfying \Cref{assumptions.base.change.sheaf.theory.VL}  \eqref{assumptions.base.change.sheaf.theory.VL-0}-\eqref{assumptions.base.change.sheaf.theory.V-6}.  

\item\label{rem-additional-base.change.sheaf.theory-1} Suppose $\bfC_\horiz= \bfC$ (so every object is as in \Cref{ex:commutative-algebra-in-corr(C)}).
Then giving \Cref{assumptions.base.change.sheaf.theory.H} \eqref{assumptions.base.change.sheaf.theory.H-2}, \Cref{assumptions.base.change.sheaf.theory.H} \eqref{assumptions.base.change.sheaf.theory.H-3} is equivalent to the projection formula 
\[
f\horizr(\mF)\otimes \mG\cong f\horizr(\mF\otimes f\horizl(\mG)),\quad \mbox{for }\ \mF\in \der(X), \ \mG\in \der(Y).
\] 
Indeed, by letting $Z=Y$ applying $(\Delta_Y)\horizl$ to \Cref{assumptions.base.change.sheaf.theory.H} \eqref{assumptions.base.change.sheaf.theory.H-3} gives the projection formula. Conversely, Let $p_X: X\times Z\to X, q_Y: Y\times Z\to Y$ and $p_Z: X\times Z\to Z, q_Z: Y\times Z\to Z$ be the projections. Then
\begin{eqnarray*}
f\horizr\mF\boxtimes_{\La}\mG&=&(q_Y)\horizl(f\horizr\mF)\otimes (q_Z)\horizl\mG\cong (f\times \id)\horizr ((p_X)\horizl\mF) \otimes (q_Z)\horizl\mG \\
&\cong& (f\times\id)\horizr( (p_X)\horizl\mF) \otimes (p_Z)\horizl\mG)=(f\times\id)(\mF\boxtimes_{\La}\mG).
\end{eqnarray*}
Similar remarks apply to other situations. E.g. giving \Cref{assumptions.base.change.sheaf.theory.V} \eqref{assumptions.base.change.sheaf.theory.V-2}, then  \Cref{assumptions.base.change.sheaf.theory.V} \eqref{assumptions.base.change.sheaf.theory.V-3} is equivalent to $f\vertr(\mF\otimes\mG)\cong f\vertr\mF\otimes f\horizl\mG$ for $\mF\in \der(Y)$ and $\mG\in\der(Z)$.

\item\label{rem-additional-base.change.sheaf.theory-2} Sometimes the sheaf theory $\shv$ satisfies even stronger assumptions. E.g. it may happen that for some $f\in \mathrm{HL}$, there is a natural isomorphism $(f\horizl)^L\cong f\vertl$ such that the base change isomorphisms and projection formula for $(f\horizl)^L$ as from \eqref{assumptions.base.change.sheaf.theory.H-4}-\eqref{assumptions.base.change.sheaf.theory.H-6} are equivalent to the corresponding base change isomorphisms and projection formula for $f\vertl$ as encoded in the sheaf theory. See \Cref{prop-sheaf-theory-for-adjoint-factorization}. The same remark applies to other cases considered in \Cref{assumptions.base.change.sheaf.theory.H} \Cref{assumptions.base.change.sheaf.theory.V} and \Cref{assumptions.base.change.sheaf.theory.VL}. 
\end{enumerate}
\end{remark}

\subsubsection{Descent}
\begin{definition}\label{def:der-descent}
A morphism $f\colon X\rightarrow Y$ in $\horiz$ (resp. in $\verti$) such that $\Delta_{X/Y}: X\to X\times_YX$ is also in $\horiz$ (resp. in $\verti$)
is called of $\der$-descent (resp. $\der$-codescent) if
the canonical functor
\[
\der(Y) \rightarrow \tot \left(\der(X_\bullet)\right) \quad (\mbox{resp.}   |\der(X_\bullet)|\to \der(Y))
\]
induced by $\horizl$-pullbacks (resp. $\vertl$-pushforward) is an equivalence, where $X_\bullet\to Y$ denotes the \v{Cech} nerve of $f$.
It is said to be of \textit{universal $\der$-descent} (resp. \textit{universal $\der$-codescent}) if its base change along every morphism $Y'\rightarrow Y$  is $\der$-descent (resp. $\der$-codescent).
\end{definition}

Clearly, the collection of morphisms of universal $\der$-descent is closed under base change. We recall the following  result \cite[Lemma 3.1.2]{liu2014enhanced} regarding stability properties of such morphisms. 
\begin{lemma}\label{Liu.Zheng.descent.lemma} 
Assume that $\horiz$ is strongly stable.
Let $f\colon X\rightarrow Y$ and $g\colon Y\rightarrow Z$ be morphisms in $\horiz$.
\begin{enumerate}
    \item If $f$ admits a section, then it is of universal $\der$-descent. 
    \item If $f,g$ are of universal $\der$-descent, then $g\circ f$ is of universal $\der$-descent. 
    \item \label{descent.detection}If $g\circ f$ is of universal  $\der$-descent, then $g$ is of universal $\der$-descent. 
\end{enumerate}
The same statements hold for codescent with $f,g\in \verti$.
\end{lemma}

We will also need the following.
\begin{proposition}\label{prop-descent-codescent-abstract}
\begin{enumerate}
\item\label{prop-descent-codescent-abstract-1} 
 Let $(f: X\to Y)\in \bfC_{\verti}$ such that $(\Delta_{X/Y}: X\to X\times_YX)\in \bfC_\verti$.
Suppose  $f\vertr: \der(Y)\to \der(X)$ is conservative, and suppose $f$ satisfies \Cref{assumptions.base.change.sheaf.theory.V} \eqref{assumptions.base.change.sheaf.theory.V-0} \eqref{assumptions.base.change.sheaf.theory.V-1}. 
Then $f$ is of $\der$-codescent. 
 
\item\label{prop-descent-codescent-abstract-2}  Let $(f: X\to Y)\in \bfC_{\horiz}$ such that $(\Delta_{X/Y}: X\to X\times_YX)\in \bfC_\horiz$. Assume that $(\Delta_Y: Y\to Y\times Y)\in \bfC_\horiz$.
Suppose \Cref{assumptions.base.change.sheaf.theory.H} \eqref{assumptions.base.change.sheaf.theory.H-2} and  \eqref{assumptions.base.change.sheaf.theory.H-3} hold, and $f\in \mathrm{HR}$. If $\La_Y\to f\horizr f\horizl (\La_Y)=f\horizr(\La_X)$ admits a section, then $f$ satisfies $\der$-descent.
\end{enumerate}
\end{proposition}
\begin{proof}Part \eqref{prop-descent-codescent-abstract-1} follows from \Cref{thm:Beck-Chevalley-descent} 
and that  \eqref{eq:limit-colimit equivalence} is an equivalence. We prove Part \eqref{prop-descent-codescent-abstract-2} using the co-monadic version of \Cref{thm:Beck-Chevalley-descent}. 
The right adjointability follows from  \Cref{assumptions.base.change.sheaf.theory.H} \eqref{assumptions.base.change.sheaf.theory.H-2}.
Via the projection formula \Cref{assumptions.base.change.sheaf.theory.H} \eqref{assumptions.base.change.sheaf.theory.H-3} (see \Cref{rem-additional-base.change.sheaf.theory}), we obtain a section $f\horizr f\horizl\rightarrow \id$ of the unit map $\id \rightarrow f\horizr f\horizl$ 
\[
\mF \rightarrow f\horizr (f\horizl \mF)\cong f\horizr(\La_X) \otimes \mF \rightarrow \mF.
\]
This immediately implies that $f\horizl$ is conservative. In addition, if $f\horizl (\mF^\bullet)$ is a split cosimplicial object, the cosimplicial object $f\horizr (f\horizl  (\mF^\bullet))$ is also split and in particular it is a limit diagram in $\der(Y)$. Then the diagram $\mF^{\bullet}\colon \Delta_{+} \rightarrow \der(Y)$ is a retract of a limit diagram, which implies it is a limit diagram as well. This verifies the first assumption in the co-monadic version of \cite[Corollary 4.7.5.3]{Lurie.higher.algebra}.
\end{proof}

\subsubsection{$\der$-Admissibility}\label{SS: der-adm}
Recall the notion of admissible objects \Cref{def:adm vs compact}. We discuss such notion in geometric set-up. 

For $X\in\bfC$ with $\pi_X,\Delta_X\in \verti\cap\horiz$. We regard $\der(X)$ as a $\der(\pt)$-module. If the exterior tensor product \eqref{eq:abstract-exterior-tensor-product} (for $Y=X$, and with $\otimes_\La$ replaced by $\otimes_{\der(\pt)}$) is an equivalence, then $\der(X)$ is self-dual as a $\der(\pt)$-linear category, with unit and counit given in \Cref{rem-dualizability in corr}. In particular, the $\der(\pt)$-linear functor $(\pi_X)\vertl: \der(X)\to \der(\pt)$ induces such self-duality
as in \Cref{ex: duality via Frobenius-structure} so all discussions from \Cref{SS: admissible objects} apply. Although in practice \eqref{eq:abstract-exterior-tensor-product} is often not an equivalence, we can still make sense of such notion in geometric set-up since (the analogue of) the characterization of admissible objects in \Cref{rem: Zorro for admissible objects} always makes sense.

\begin{definition}\label{def:ULA-morphism}
Assume that $\pi_X,\Delta_X\in \verti\cap\horiz$. 
An object $\mF\in \der(X)$ is called $\der$-admissible if there exists another $\mF^\vee\in \der(X)$ equipped with
\[
(\Delta_X)\vertl \La_X\to \mF\boxtimes_{\La} \mF^\vee,\quad (\pi_X)\vertl(\Delta_X)\horizl(\mF^\vee\boxtimes_{\La} \mF)\to \La_\pt
\]
such that both the induced map 
\begin{multline*} 
\mF\cong \La_{\pt}\boxtimes_{\La} \mF\cong (\id_X\times\pi_X)\vertl(\id_X\times\Delta_X)\horizl((\Delta_X)\vertl\La_X\boxtimes_{\La}\mF) \to \\
(\id_X\times\pi_X)\vertl(\id_X\times\Delta_X)\horizl(\mF\boxtimes_{\La} \mF^\vee\boxtimes_{\La}\mF) \cong \mF\boxtimes_{\der(\pt)} (\pi_X)\vertl(\Delta_X)\horizl(\mF^\vee\boxtimes_{\La} \mF)\to   \mF  \boxtimes_{\La}  \La_\pt\cong \mF,
\end{multline*}
and the similarly defined map from $\mF^\vee\to \mF^\vee$ are homotopic the identity map.
We say $X$ is $\der$-admissible if $\La_X=(\pi_X)\horizl\La_\pt$ is $\der$-admissible.

Let $f: X\to Y$ be a morphism in $\verti$ such that $\Delta_{X/Y}\in \verti$. We say $\mF$ is $\der$-admissible with respect to $f$ if $(X,\mF)$ is $\der_{/Y}$-admissible, where the sheaf theory $\der_{/Y}$ is as in Remark \ref{rem:relative-sheaf-theory}. 
We say $f: X\to Y$ is $\der$-admissible if $\La_X$ is $\der$-admissible with respect to $f$.
\end{definition}

For simplicity, from now on we make the following assumption throughout the rest of the section. Recall \Cref{def-closure-property-of-morphism}.
\begin{assumptions}\label{assumptions.horiz.vert}
The class $\horiz=\all$ and the class $\verti$ is strongly stable. 
\end{assumptions}

Note that under the above assumption \Cref{ex:commutative-algebra-in-corr(C)} always applies.
Assume that $\pi_X: X\to \pt$ belongs to $\verti$. Let
\begin{equation}\label{eq:abstract-verdier-dual.II}
   (-)^{\vee,\der}: \der(X)^{\op}\to\der(X),\quad \mF^{\vee,\der}:=\rhom(\mF,\pi_X\vertr\La_{\pt}).
\end{equation}
Then for $f: X\to Y$ in $\verti\cap \horiz$, the isomorphisms in \eqref{eq:abstract-hom-pull-push} specialize to isomorphisms 
\[
(f\horizl(-))^{\vee,\der} =f\vertr ((-)^{\vee,\der}),\quad f\horizr ((-)^{\vee,\der})= (f\vertl(-))^{\vee,\der}.
\]

Now suppose $(\pi_X)\vertl: \der(X)\to\der(\pt)$ is a Frobenius structure as in \Cref{rem-dualizability in corr} (but we do not assume that $\der(X)\otimes_{\der(\pt)}\der(X)\to\der(X\times X)$ is an equivalence).
Then the role of $\verd^\la$ is played by 
$(\pi_X)\vertr\La_\pt$
In particular, if $\mF\in \der(X)$ is $\der$-admissible, then 
\[
\bD_X^\der(\mF)\cong \mF^{\vee,\der}.
\] 
is $\der$-admissible.
We also have the counterpart of \Cref{lem-adm-object-in-Frob-alg} in the geometric setting, as stated below. The same proof is the same.

\begin{lemma}\label{lem-dual-object-in-corr-D} 
    An object $\mF\in\der(X)$ is $\der$-admissible if and only if the natural map 
    \[
     \bD^\der_X(\mF)\boxtimes_{\La}  \mG\to \underline\Hom((p_{X})\horizl\mF, (p_{Y})\vertr\mG)
    \]
    is an isomorphism from every $Y$ and every $\mG\in \der(Y)$, where $p_X: X\times Y\to X$ and $p_Y: X\times Y\to Y$ are two projections; if and only if the above isomorphism holds for $(Y,\mG)=(X,\mF)$.
\end{lemma}

We also have the following geometric counterpart of an observation from \Cref{ex: cpt equal to adm in sm and proper cat}.
\begin{corollary}\label{lem-ula-vs-compact}
  Let $\mF\in\der(X)$ be $\der$-admissible.
  If  $(\Delta_X)\vertl\La_X\in \der(X\times X)$ is a compact object, then $\mF\in \der(X)$ is compact. 
\end{corollary}
\begin{proof}
By \Cref{lem-dual-object-in-corr-D}  and \eqref{eq:abstract-hom-pull-push},
$\Hom_{\der(X)}(\mF, -)=\Hom_{\der(X\times X)}((\Delta_X)\vertl\La_X, \mF^{\vee,\der}\boxtimes_{\La} (-))$ commutes with colimits.
\end{proof}

\begin{remark}\label{rem: Sheaf theory via Grothendieck construction}
As explained in \Cref{rem: codomain of sheaf theory} \eqref{rem: codomain of sheaf theory-3}, 
the sheaf theory $\der$ can also be (largely) encoded as a symmetric monoidal $2$-category $\corr^{\der}(\bfC)_{\verti;\horiz}$. Then $\mF\in \der(X)$ is $\der$-admissible if and only if $(X,\mF)$ is a dualizable object in this category, and \Cref{lem-dual-object-in-corr-D}  also follows from general facts about dualizable objects in a symmetric monoidal category. 
\end{remark}

The importance of this notion lies in the following fact.

\begin{lemma}\label{lem:base-change-for-dualizable} 
  \begin{enumerate}
  \item Assume that $\mF\in \der(X)$ is $\der$-admissible. Then for every $g:Y'\to Y$ and  $\mG\in\der(Y')$, then we have the  natural isomorphism
  \begin{equation}\label{eq:Kunneth-formula-ULA}
     \mF \boxtimes_{\La} g\horizr(\mG)\cong (\id\times g)\horizr(\mF\boxtimes_{\La}\mG).
  \end{equation}
  \item Let $(f: X\to Y)\in \verti$ such that both $\mF\in\der(Y)$ and $f\horizl\mF\in\der(X)$ are $\der$-admissible. Then for every $(Y',\mG)$, we have the natural isomorphism (adjunction of \eqref{eq:abstract-Kunneth-formula})
  \begin{equation}\label{eq:pullback-Kunneth-formula-ULA}
     f\vertr(\bD^\der_Y\mF)\boxtimes_{\La} \mG\to (f\times\id_{Y'})\vertr((\bD^\der_Y\mF)\boxtimes_{\La} \mG).
  \end{equation}
  \end{enumerate}
\end{lemma}
\begin{proof}
This follows from the same proof as in \cite[Lemma 2.11(b)]{Lu.Zheng.relative.Lefschetz}. We sktech a proof for completeness.
First \eqref{eq:Kunneth-formula-ULA} follows from
\[
\mF\boxtimes_{\La} g\horizr\mG=\rhom((p_X)\horizl(\bD^\der_X\mF), (p_Y)\vertr(g\horizr\mG))=(\id_X\times g)\horizr\rhom((p_X)\horizl \bD^\der_X\mF,(p_{Y'})\vertr\mG)=(\id_X\times g)\horizr(\mF\boxtimes_{\La}\mG).
\]

 For \eqref{eq:pullback-Kunneth-formula-ULA},  note that if $f\horizl\mF$ is $\der$-admissible, then so is $f\vertr(
    \bD^\der_Y\mF)$. It follows  \eqref{eq:pullback-Kunneth-formula-ULA} is identified with
   \[
      \rhom((p_{X})\horizl f\horizl\mF,(p_{Y'})\vertr\mG)=\rhom((f\times\id_{Y'})\horizl (p_{Y})\horizl\mF,(p_{Y'})\vertr\mG)=(f\times\id_{Y'})\vertr\rhom((p_Y)\horizl\mF,(p_{Y'})\vertr\mG).
   \]
\end{proof}

Now suppose $f: X\to Y \in \bfC_{\verti}$, and suppose $\mF$ is $\der$-admissible with respect to $f$.
Then for $g: Y'\to Y$, the natural map \eqref{eq:Kunneth-formula-ULA}
specializes to an isomorphism 
\begin{equation*}
\mF\otimes f\horizl(g\horizr\mG)\xrightarrow{\cong} (g')\horizr((g')\horizl\mF\otimes (f')\horizl\mG),
\end{equation*}
where $f',g'$ are as in the Cartesian diagram \eqref{eq:base-change-diagram-app}. In particular, we have the abstract smooth base change isomorphism.
\begin{corollary}\label{cor:abstract-smooth-base-change}
If $f: X\to Y$ is $\der$-admissible, then the Beck-Chevallay map from \eqref{eq:base-change-diagram-app} is an equivalence
  \begin{equation}\label{eq:abstract-smooth-base-change}
        f\horizl\circ g\horizr\simeq (g')\horizr\circ (f')\horizl.
  \end{equation}
\end{corollary}

Similarly, let $(f: X\to Y)\in \bfC_{\verti}$ be $\der$-admissible, and let $g: Y'\to Y$. Then \eqref{eq:pullback-Kunneth-formula-ULA} specializes to an isomorphism
\begin{equation}\label{eq:abstract-uppershrek-as-upperstar}
(g')\horizl(f\vertr\La_Y)\otimes (f')\horizl\mG\xrightarrow{\cong} (f')\vertr\mG,
\end{equation}
where $f',g'$ are as in the Cartesian diagram \eqref{eq:base-change-diagram-app}.

\begin{corollary}\label{cor:abstract-base-change-II}
Assume that $f:X\to Y$ is $\der$-admissible.  \begin{enumerate}
    \item The map \eqref{eq:abstract-!-*-pullback} is an isomorphism. 
    \item The Beck-Chevalley map associated to \eqref{eq:abstract-base-change} is an equivalence
\begin{equation}\label{eq-abstract-comm-uppershrek-upperstar}
(g')\horizl\circ f\vertr\cong (f')\vertr\circ g\horizl.
\end{equation}
    \item In addition if $g$ belongs to $\verti$, then the Beck-Chevallay map from \eqref{eq:base-change-diagram-app}
\begin{equation}\label{eq-abstract-comm-uppershrek-lowershrek}
(g')\vertl\circ (f')\vertr\to f\vertr \circ g\vertl \colon \der(Y') \to \der(X) 
\end{equation}
is an equivalence.  
\end{enumerate}
\end{corollary}
\begin{proof}
 The first statement follows from by letting $g=\id: Y\to Y$ in \eqref{eq:abstract-uppershrek-as-upperstar}. The second statement follows from \eqref{eq:abstract-uppershrek-as-upperstar} by letting $\mG$ be in the essential image of $g\horizl$. The last statement follows from \eqref{eq:abstract-uppershrek-as-upperstar}, the projection formula \eqref{eq:abstract-projection-formula} and the base change isomorphism \eqref{eq:abstract-base-change}.
\end{proof}

\begin{lemma}\label{cor-abstract-composition-ULA}
\begin{enumerate}
    \item\label{cor-abstract-composition-ULA-1}
If $\mF\in\der(X)$ is $\der$-admissible with respect to $f: X\to Y$, then for every $g: Y'\to Y$, $(g')\horizl\mF\in \der(X')$ is $\der$-admissible with respect to $f': X'\to Y'$, where $f',g'$ are as in \eqref{eq:base-change-diagram-app}. In particular, $\der$-admissible morphisms are stable under base change.
   \item\label{cor-abstract-composition-ULA-2} If $\mF\in \der(Y)$ is $\der$-admissible with respect to $g: Y\to Z$ and $f\horizl\mF$ is $\der$-admissible with respect to $f: X\to Y$, then $f\horizl\mF$ is $\der$-admissible with respect to $g\circ f$. In particular, $\der$-admissible morphisms are stable under compositions.
\end{enumerate}
\end{lemma}
\begin{proof}
Part \eqref{cor-abstract-composition-ULA-1} is clear. 
For Part \eqref{cor-abstract-composition-ULA-2}, we may assume that $Z=\pt$.
By \Cref{lem-dual-object-in-corr-D}, $g\circ f$ is $\der$-ULA if and only if $(p_1)\horizl(\pi_X)\vertr\La_{\pt}\cong (p_2)\vertr\La_X$. Using Part \eqref{cor-abstract-composition-ULA-1}  and \eqref{eq-abstract-comm-uppershrek-upperstar}, this follows by the isomorphism $(p_1)\horizl(\pi_Y)\vertr\La_{\pt}\cong (p_2)\vertr\La_Y$ and a similar isomorphism obtain by applying \Cref{lem-dual-object-in-corr-D} to $f$.
\end{proof}

Although we shall not make use of it, let us also explain how Poincar\'e duality fits into the above formalism. 
\begin{proposition}\label{prop: Poincare duality}
Let $(f: X\to Y)\in \horiz\cap \verti$. We suppose $\Delta_X,\pi_X$ and $\Delta_Y,\pi_Y$ belong to $\horiz$ so $f\horizl: \der(Y)\to \der(X)$ is a symmetric monoidal functor. Suppose
$f\vertl$ is the right adjoint of $f\horizl$. (E.g. such situation arises when $\der$ is constructed as \Cref{ex-sheaf-theory-for-adjoint-factorization} below.) In addition, suppose that $f$ is $\der$-admissible.
Then $f\vertl$ sends dualizable objects in $\der(X)$ to dualizable objects in $\der(Y)$. 
\end{proposition}
\begin{proof}
Indeed, suppose $\mF\in\der(X)$ is dualizable with $\mG$ its dual. Then for every $\mG_1,\mG_2\in \der(Y)$, we find
\begin{multline*}
\Hom_{\der(Y)}(f\vertl\mF\otimes \mG_1,\mG_2)=\Hom_{\der(X)}(\mF\otimes f\horizl\mG_1, f\vertr\mG_2)\\
=\Hom_{\der(Y)}(f\horizl\mG_1,\mG\otimes f\horizl\mG_2\otimes f\vertr(\La_Y))=\Hom_{\der(Y)}(\mG_1,f\vertl(\mG\otimes  f\vertr(\La_Y))\otimes \mG_2),
\end{multline*}
showing that the dual of $f\vertl\mF$ is $f\vertl(\mG\otimes f\vertr(\La_Y))$.
\end{proof}

\subsubsection{Extensions of sheaf theories}\label{SS: extension of sheaf theory}
As should be clear from the above discussions, a sheaf theory encodes a huge amount of information. So it must be highly non-trivial to construct a sheaf theory. Now we review (and slightly extend) a few results as in \cite{liu2012enhanced} and \cite{Gaitsgory.Rozenblyum.DAG.vol.I} allowing one to construct sheaf theories from scratch. 

Suppose there are two triples $(\bfC_i,\verti_i,\horiz_i), \ i=1,2$, and a functor $F: \bfC_1\to \bfC_2$ that preserves finite limits and restricts to functors $F_{\verti}: (\bfC_1)_{\verti_1}\to (\bfC_2)_{\verti_2}$ and $F_{\horiz}: (\bfC_1)_{\horiz_1}\to (\bfC_2)_{\horiz_1}$. It then induces a symmetric monoidal functor 
\[
F_{\corr}: \corr(\bfC_1)_{\verti_1;\horiz_1}\to \corr(\bfC_2)_{\verti_2;\horiz_2}.
\]
We suppose $F_{\corr}$ is a (not necessarily full) embedding. Now
giving a sheaf theory $\der_1:  \corr(\bfC_1)_{\verti_1;\horiz_1}\to \lincat_\La$, we would like to ask whether there is a sheaf theory $\der_2: \corr(\bfC_2)_{\verti_2;\horiz_2}\to \lincat_\La$ such that $\der_1\simeq \der_2\circ F_{\corr}$. If so we call such $\der_2$ an extension of $\der_1$.

\begin{remark}
In the discussion below, we will generally ignore uniqueness of extensions. But we expect all the extensions given below should be also unique in appropriate sense.
\end{remark}

We have the following basic result regarding extension of sheaf theories. It is an abstraction of the construction of  \cite[\textsection{3.2}]{liu2012enhanced}. Under slightly different assumptions, it is also proved in  \cite[Theorem 7.5.2.4]{Gaitsgory.Rozenblyum.DAG.vol.I}.

\begin{theorem}\label{prop-sheaf-theory-for-adjoint-factorization}
Let 
\[
\der: \corr(\bfC)_{\verti;\horiz}\to \lincat_\La
\]
be a sheaf theory, and let $\mathrm{HL}$ be the class of morphisms associated $\der$ as defined in \Cref{rem-additional-base.change.sheaf.theory} \eqref{rem-additional-base.change.sheaf.theory-0} (i.e. the class of morphisms satisfying \Cref{assumptions.base.change.sheaf.theory.HL}).

Let $\mathrm{E}\subset \mathrm{HL}$ be a class of morphisms, and let $\verti'$ be another weakly stable class of morphisms in $\bfC$.
Suppose
\begin{enumerate}
\item\label{ex-sheaf-theory-for-adjoint-factorization-condition-1} both $\verti$ and $\mathrm{E}$ are strongly stable; 
\item\label{ex-sheaf-theory-for-adjoint-factorization-condition-2} every $f\in \mathrm{E}\cap \verti$ is $n$-truncated from some $-2\leq n<\infty$ (which may depend on $f$);
\item\label{ex-sheaf-theory-for-adjoint-factorization-condition-3} every $f\in \verti'$ admits a decomposition $f=f_{\verti}\circ f_{\mathrm{E}}$ with $f_{\verti}\in \verti$ and $f_{\mathrm{E}}\in \mathrm{E}$. 
 \end{enumerate}
 Then $\verti'$ is strongly stable and $\der$ admits an extension to a sheaf theory 
 \[
 \der': \corr(\bfC)_{\verti';\horiz}\to \lincat_\La
 \]
 such that $f\vertl=(f\horizl)^L$ for $f\in\mathrm{E}$.
 
 If in addition, $\der$ takes value in $\cptcat_\La$, then $\der'$ also takes value in $\cptcat_\La$.
 \end{theorem}

\begin{proof}
The first statement follows from \cite[Remark 5.5]{liu2012enhanceda}. To prove assertions about extension,
we follow the same arguments of \cite[\textsection{3.2}]{liu2012enhanced}. 
We make use of notations from \emph{loc. cit.} (and therefore use the quasi-category model of $\corr(\bfC)$).
By \cite[Example 4.30]{liu2012enhanceda}, a sheaf theory $\der_1$ is equivalent to a functor $\der: \delta_{2,\{2\}}^*((\bfC^{\op})^{\sqcup,\op})^{\mathrm{cart}}_{\verti,\horiz}\to \lincat_\La$.
By composing with the functor $\delta^*_{3, \{2,3\}}((\bfC^{\op})^{\sqcup,\op})^{\mathrm{cart}}_{\verti,\mathrm{E},\horiz}\to  \delta_{2,\{2\}}^*((\bfC^{\op})^{\sqcup,\op})^{\mathrm{cart}}_{\verti,\horiz}$ obtained by taking the partial diagonal along the 2nd and 3rd factor, we obtain $\delta^*_{3, \{2,3\}}((\bfC^{\op})^{\sqcup,\op})^{\mathrm{cart}}_{\verti,\mathrm{E},\horiz}\to \lincat_\La$. On the other hand, \Cref{assumptions.base.change.sheaf.theory.H} \eqref{assumptions.base.change.sheaf.theory.H-4}-\eqref{assumptions.base.change.sheaf.theory.H-6} allow one to apply \cite[Proposition 1.4.4]{liu2012enhanced} to take the partial adjoint along the second factor, giving $\delta^*_{3, \{3\}}((\bfC^{\op})^{\sqcup,\op})^{\mathrm{cart}}_{\verti_1,\mathrm{E},\horiz}\to \catid_\La$ or $\lincat_\La$. (See \cite[Lemma 3.2.5]{liu2012enhanced} for more details.)  Finally, assumptions of the theorem imply that
 the functor $\delta^*_{3, \{3\}}((\bfC^{\op})^{\sqcup,\op})^{\mathrm{cart}}_{\verti,\mathrm{E},\horiz}\to  \delta_{2,\{2\}}^*((\bfC^{\op})^{\sqcup,\op})^{\mathrm{cart}}_{\verti',\horiz}$ obtained by taking the diagonal along the 1st and the 2nd factor
is a categorical equivalence, by \cite[Theorem 5.4]{liu2012enhanceda}. 
We thus obtained the desired extension. 
\end{proof} 
 
\begin{remark}\label{rem-sheaf-theory-for-adjoint-factorization} 
\begin{enumerate}
\item Being the category of categories, $\lincat_\La$ admits $2$-categorical structures. On the other hand,
as mentioned in \Cref{rem-generalization-corr}, the category of correspondences also admits a $2$-categorical enhancement. 
Sheaf theory constructed in \Cref{prop-sheaf-theory-for-adjoint-factorization} can be enhanced at the $2$-categorical level as a functor 
\[
\corr(\bfC)_{\verti';\horiz}^{\mathrm{E}}\to \lincat_\La,
\] 
at least if all morphisms in $\mathrm{E}$ are $m$-truncated for some $-2\leq m<\infty$, by applying \cite[Theorem 7.4.1.3, Theorem 9.3.1.2]{Gaitsgory.Rozenblyum.DAG.vol.I} (together with \cite[Remark 5.2]{liu2012enhanceda}) inductively to the pair $\mathrm{E}_{i-1}\subset \mathrm{E}_i$, where $\mathrm{E}_i\subset\mathrm{E}$ is the subclass of those morphisms that are $i$-truncated.

\item Suppose we have the sheaf theory constructed in \Cref{prop-sheaf-theory-for-adjoint-factorization}.
Let $X_\bullet$ be a Segal object in $\bfC$ as in \Cref{rem-conv-product}. 
If all morphisms of the simplicial objects $X_\bullet$ are in $\bfC_{\mathrm{HL}}$, then $\La_{X_1}\in \der(X_1)$ is a natural algebra object with respect to the convolution monoidal structure of $\der(X_1)$. Indeed, the multiplication of $\La_{X_1}$ amounts to a morphism
\[
(d_0\times d_2)\horizl(\La_{X_1}\boxtimes_{\La} \La_{X_1})\cong \La_{X_2} \to (d_2)\vertr(\La_{X_1}),
\]
which is given by the adjunction $(d_2)\vertl((d_2)\horizl(\La_{X_1}))=((d_2)\horizl)^L((d_2)\horizl(\La_{X_1}))\to \La_{X_1}$.

\item\label{rem-sheaf-theory-for-adjoint-factorization-3} There are variants of the above theorem. E.g. Instead of assuming that each $f\in \verti'$ admits a decomposition as in the theorem, one could assume that each $f$ admits a decomposition $f=f_{\mathrm{E}}\circ f_{\verti}$.
One could also replace $\mathrm{E}\subset \mathrm{HL}$ by $\mathrm{E}\subset \mathrm{HR}$, $\mathrm{E}\subset \mathrm{VL}$ or $\mathrm{E}\subset \mathrm{VR}$ and the corresponding assumptions. (In the case $\der$ takes value in $\cptcat_\La$, one further requires the right adjoints of morphisms in $\mathrm{E}\subset \mathrm{HR}$ and in $\mathrm{E}\subset \mathrm{VR}$ preserve compact objects.)
The proofs remain the same.

\item\label{rem-sheaf-theory-for-adjoint-factorization-4} Note that in fact in the statement of \Cref{prop-sheaf-theory-for-adjoint-factorization} one may replace $\lincat_\La$ by $\cat$ or other symmetric monoidal $2$-category. The proof does not change.
\end{enumerate}
\end{remark}

Here is the basic example, where things get started.
\begin{corollary}\label{ex-sheaf-theory-for-adjoint-factorization}
Suppose there is a lax symmetric monoidal functor
\[
\der: \bfC^{\op}\to \lincat_\La.
\] 
We regard it as a sheaf theory $\corr(\bfC)_{\iso;\all}\to\lincat_\La$.
\begin{enumerate}
\item\label{ex-sheaf-theory-for-adjoint-factorization-1} Let $\mathrm{L}\subset \mathrm{HL}$ be a weakly stable class. I.e. morphisms in $\mathrm{L}$ satisfy \Cref{assumptions.base.change.sheaf.theory.HL}.
 Then $\der$ extends uniquely to a sheaf theory 
 \[
 \der^{\mathrm{L}}: \corr(\bfC)_{\mathrm{L};\all}\to \lincat_\La,
 \] 
 such that $g\horizl=\der(g)$ and $f\vertl=\der(f)^L$ for $f\in \mathrm{L}$. 
 
\item\label{ex-sheaf-theory-for-adjoint-factorization-2} Dually, let $\mathrm{R}\subset \mathrm{HR}$ be a weakly stable class. I.e. morphisms in $\mathrm{L}$ satisfy \Cref{assumptions.base.change.sheaf.theory.H}. 
Then $\der$ extends to a sheaf theory $\der^{\mathrm{R}}: \corr(\bfC)_{\mathrm{R};\all}\to \lincat_\La$. 

\item\label{ex-sheaf-theory-for-adjoint-factorization-3} Now let $\mathrm{I}\subset \mathrm{L}$ and $\mathrm{P}\subset \mathrm{R}$ be two strongly stable classes of morphisms, and let $\verti$ be a weakly stable class of morphisms containing both $\mathrm{I}$ and $\mathrm{P}$. Suppose
\begin{itemize}
\item for every Cartesian diagram \eqref{eq:base-change-diagram-app} with $f\in \mathrm{I}$ and $g\in \mathrm{P}$, $\der^L(f)\circ \der(g')^R\to \der(g)^R\circ \der(f')^L$ is an isomorphism.
\item every $f\in \mathrm{I}\cap \mathrm{P}$ is $n$-truncated from some $n\geq -2$ (which may depend on $f$);
\item every $f\in \verti$ admits a decomposition $f=f_{\mathrm{P}}\circ f_{\mathrm{I}}$ with $f_{\mathrm{I}}\in \mathrm{I}$ and $f_{\mathrm{P}}\in \mathrm{P}$;
\end{itemize}
Then there is a sheaf theory $\der\colon \corr(\bfC)_{\verti;\all}\to \lincat_\La$ that extends $\der^{\mathrm{L}}|_{\corr(\bfC)_{\mathrm{I};\all}}$ and $\der^{\mathrm{R}}|_{\corr(\bfC)_{\mathrm{P};\all}}$. 
\end{enumerate}
\end{corollary}
\begin{proof}
For Part  \eqref{ex-sheaf-theory-for-adjoint-factorization-1} and \eqref{ex-sheaf-theory-for-adjoint-factorization-2}, we use the same argument as in \Cref{prop-sheaf-theory-for-adjoint-factorization}, except that we do not need the last step (and therefore do not need $\mathrm{L}$ and $\mathrm{R}$ to be strongly stable). The last part follows from \Cref{prop-sheaf-theory-for-adjoint-factorization}.
\end{proof}

We also need another type of extensions of sheaf theory, namely via Kan extensions. 

\begin{proposition}\label{prop-sheaf-theory-right-Kan-extension}
Suppose we have $(\bfC_i, \verti_i, \horiz_i),\ i=1,2$, and a finite limit preserving fully faithful embedding $\bfC_1\subset \bfC_2$ which induces fully faithful embeddings $(\bfC_1)_{\verti_1}\subset (\bfC_2)_{\verti_2}$, $(\bfC_1)_{\horiz_1}\subset (\bfC_2)_{\horiz_2}$. 
We in addition make the following assumptions.
\begin{enumerate}
\item\label{prop-sheaf-theory-right-Kan-extension-1} Let $Y\in\bfC_1$. The main diagonal $\Delta: Y\to Y^m$ belongs to $\horiz_1$ for every $m$, and for every $(Y\to X_1\times\cdots\times X_m)\in (\bfC_2)_{\horiz_2}$, the projection $Y\to X_i$ belongs to $\horiz_2$.
\item\label{prop-sheaf-theory-right-Kan-extension-2} The class $\verti_2$ are representable in $\verti_1$. (See \Cref{rem-prop-of-weak-strong-stable-class} \eqref{rem-prop-of-weak-strong-stable-class-2} for the meaning.) 
\end{enumerate}
Then every sheaf theory $\der_1:\corr(\bfC_1)_{\verti_1, \horiz_1}\to \lincat_\La$ admits an extension
\[
\der_2: \corr(\bfC_2)_{\verti_2;\horiz_2}\to \lincat_\La
\]
such that the restriction $\der_2|_{((\bfC_2)_{\horiz_2})^{\op}}$ is canonically isomorphic to the right Kan extension of $\der_1|_{((\bfC_1)_{\horiz_1})^{\op}}$ along $((\bfC_1)_{\horiz_1})^{\op}\subset ((\bfC_2)_{\horiz_2})^{\op}$ (as plain functors). 
\end{proposition}
\begin{proof}
This is proved in (the easy part of) \cite[Theorem 8.6.1.5, Proposition 9.3.2.4]{Gaitsgory.Rozenblyum.DAG.vol.I}. We include a sketch for completeness. 

We let $\der_2$ be the right Kan extension of $\der_1$ along $\corr(\bfC_1)_{\verti_1;\horiz_1}\to \corr(\bfC_2)_{\verti_2;\horiz_2}$.
Note that Assumption \eqref{prop-sheaf-theory-right-Kan-extension-2} implies that $\corr(\bfC_1)_{\verti_1;\horiz_1}\to \corr(\bfC_2)_{\verti_2;\horiz_2}$ is fully faithful so $\der_2$ is indeed an extension of $\der_1$ (as a functor, but not yet as a sheaf theory). To see that the restriction $\der_2|_{((\bfC_2)_{\horiz_2})^{\op}}$ is the right Kan extension of $\der_1|_{((\bfC_1)_{\horiz_1})^{\op}}$, it is enough to notice 
that for $X\in \bfC_2$, the functor 
\begin{equation}\label{eq: cofinality-horiz-to-corr}
(((\bfC_2)_{\horiz_2})_{/X})^{\op}\times_{((\bfC_2)_{\horiz_2})^{\op}}((\bfC_1)_{\horiz_1})^{\op}\to (\corr(\bfC_2)_{\verti_2;\horiz_2})_{X/}\times_{\corr(\bfC_2)_{\verti_2;\horiz_2}}\corr(\bfC_1)_{\verti_1;\horiz_1}
\end{equation}
is cofinal. Indeed, let we write $\mI$ for the source category and $\mJ$ for the target category of the above functor.
Let $(X\xleftarrow{f} Z\xrightarrow{g} Y)\in \mJ$, then $Y\in \bfC_1$ (so $Z\in \bfC_1$ and $g\in \verti_1$).
Unveiling the definition, $\mI\times_{\mJ}\mJ_{/g\circ f^{-1}}$ is nothing but the category of factorizations of $f$ into $Z\to Z'\to X$ with  $(Z\to Z')\in (\bfC_1)_{\horiz_1}$ and $(Z'\to Z)\in (\bfC_2)_{\horiz_2}$.
It is clear that this category admits a final object given by $Z\xrightarrow{\id_Z}Z\xrightarrow{f} X$. 

It remains to endow $\der_2$ with a lax symmetric monoidal structure. For a symmetric monoidal category $\mE$, let $\mE^\otimes\to \mathrm{Fin}_*$ denote the corresponding coCartesian fibration encoding the symmestric monoidal structure. We may compose $\der_i$ with the lax symmetric monoidal functor  $\lincat_\La\to \cat$ and show the composed functor admits a lax symmetric monoidal structure. As the symmetric monoidal structure on $\cat$ is Cartesian, we may apply \cite[Proposition 2.4.1.7]{Lurie.higher.algebra} to regard $\der_1$ as a lax Cartesian structure from $\corr(\bfC_1)_{\verti_1;\horiz_1}^\otimes$ to $\cat$, sending $(X_j)_{1\leq j\leq m}$ to $\prod_j \der_1(X_j)$.
It is enough to show that its right Kan extension along the (fully faithful) embedding $\corr(\bfC_1)_{\verti_1;\horiz_1}^\otimes\to \corr(\bfC_2)_{\verti_2;\horiz_2}^\otimes$ is a lax Cartesian structure. For this, using $\otimes$-version of \eqref{eq: cofinality-horiz-to-corr}, we reduces to show that for $(X_i)_{1\leq i\leq m}\in \bfC_2^{m}$,
\[
\prod_i\Bigl(((\bfC_2)_{\horiz_2})_{/X_i}\times_{(\bfC_2)_{\horiz_2}}(\bfC_1)_{\horiz_1}\Bigr)\to ((\bfC_2)_{\horiz_2})_{/\prod_i X_i}\times_{(\bfC_2)_{\horiz_2}}(\bfC_1)_{\horiz_1}
\]
is cofinal. Given $(Y\to \prod X_i)\in(\bfC_2)_{\horiz_2}$ with $Y\in \bfC_1$, we need to show that the category of factorizations $Y\to \prod Y_i\to \prod X_i$ with $Y_i\in \bfC_1, (Y_i\to X_i)\in (\bfC_2)_{\horiz_2}$ and $(Y\to \prod Y_i)\in (\bfC_1)_{\horiz_1}$ is contractible. But by Assumption \eqref{prop-sheaf-theory-right-Kan-extension-1}, this category has an initial object given by
factors through $Y\xrightarrow{\Delta_Y} Y^m\to \prod X_i$. Cofinality follows.
\end{proof}

Recall that associated to a sheaf theory there are four classes of morphisms as introduced in \Cref{rem-additional-base.change.sheaf.theory} \eqref{rem-additional-base.change.sheaf.theory-0}. We need to understand how these classes of morphisms behave under the above two types of extensions of sheaf theories.
 
\begin{lemma}\label{prop-sheaf-theory-right-Kan-extension-VR class}
Assumptions are as in \Cref{prop-sheaf-theory-right-Kan-extension}. Let $\der_2$ be the extension of $\der_1$ as constructed in \Cref{prop-sheaf-theory-right-Kan-extension}.
Suppose in addition that $\horiz_2$ is strongly stable, and $\verti_2\subset \horiz_2$.
Let $\mathrm{HR}_i$, $\mathrm{HL}_i$, $\mathrm{VR}_i$, and $\mathrm{VL}_i$ be the classes of morphisms associated to $\der_i$ as in \Cref{rem-additional-base.change.sheaf.theory} \eqref{rem-additional-base.change.sheaf.theory-0}.
Let $f\in (\bfC_2)_{\bfH_2}$ which is representable in $\mathrm{HR}_1$ (resp. $\mathrm{HL}_1$, resp. $\mathrm{VR}_1$, resp. $\mathrm{VL}_1$). Then $f\in \mathrm{HR}_2$ (resp. $\mathrm{HL}_2$, resp. $\mathrm{VR}_2$, resp, $\mathrm{VL}_2$).
\end{lemma}
\begin{proof}
Let $f: X\to Y$ be a morphism in $(\bfC_2)_{\horiz_2}$ that is representable in $\bfC_1$. As we assume that $\horiz_2$ is strongly stable and $(\bfC_1)_{\horiz_1}\subset (\bfC_2)_{\horiz_2}$ is fully faithful,
we have the following natural functor 
\begin{align*}\label{multline}
(((\bfC_2)_{\horiz_2})_{/Y})^{\op}\times_{((\bfC_2)_{\horiz_2})^{\op}}((\bfC_1)_{\horiz_1})^{\op} &\to  (((\bfC_2)_{\horiz_2})_{/X})^{\op}\times_{((\bfC_2)_{\horiz_2})^{\op}}((\bfC_1)_{\horiz_1})^{\op},\\
(Y'\to Y)&\mapsto  (X':=Y'\times_YX\to X),
\end{align*}
which is cofinal. Indeed, let $\mI\to \mJ$ denote the above functor. Then for every $(f: Z\to X)\in \mJ$, the category $\mI\times_{\mJ}{\mJ_{/f}}$ admits a final object, namely $Z\times_YX\to X$. We also notice that by a similar reason, for $X,Y\in \bfC_2$, the opposite of the following functor
\begin{align*}
(\bfC_2)_{\horiz_2})_{/X}\times_{(\bfC_2)_{\horiz_2}}(\bfC_1)_{\horiz_1}\times ((\bfC_2)_{\horiz_2})_{/Y}\times_{(\bfC_2)_{\horiz_2}}(\bfC_1)_{\horiz_1}& \to ((\bfC_2)_{\horiz_2})_{/X\times Y}\times_{(\bfC_2)_{\horiz_2}}(\bfC_1)_{\horiz_1}, \\ ((U\to X), (V\to Y))&\mapsto (U\times V\to X\times Y),
\end{align*}
is cofinal.

Now let $f: X\to Y$ be a morphism, representable in $\mathrm{HR}_1$.
Then for a Cartesian diagram \eqref{eq:base-change-diagram-app} with  $(g:Y'\to Y)\in ((\bfC_2)_{\horiz_2})_{/Y}\times_{(\bfC_2)_{\horiz_2}}(\bfC_1)_{\horiz_1}$, the desired right adjointability of $f\horizl$ and the desired base change isomorphisms with respect to $g\horizl$ and $(g')\horizl$ follow from \Cref{prop:categorical-right-adjointability-colimits} and \Cref{rmk:categorical-left-adjointability-limits}. If we have the Cartesian diagram \eqref{eq:base-change-diagram-app} but with $(g:Y'\to Y)\in (\bfC_2)_{\horiz_2}$, then the corresponding base change isomorphisms  with respect to $g\horizl$ and $(g')\horizl$ can be checked after further $\horizl$-pull backs to $(V\to Y')\in ((\bfC_2)_{\horiz_2})_{/Y'}\times_{(\bfC_2)_{\horiz_2}}(\bfC_1)_{\horiz_1}$, which then follows from the already established cases. On the other hand, if $(g: Y'\to Y)\in (\bfC_2)_{\verti_2}$, the corresponding base change isomorphisms  with respect to $g\vertl$ and $(g')\vertl$ can be similarly checked after further $\horizl$-pull backs to $(V\to Y)\in ((\bfC_2)_{\horiz_2})_{/Y}\times_{(\bfC_2)_{\horiz_2}}(\bfC_1)_{\horiz_1}$, which then also follows from the already established cases. Next, consider $X\times Z\xrightarrow{f\times \id_Z} Y\times Z$. 
Using this and the established base change isomorphisms and cofinality, the corresponding projection formulas can be checked after $\horizl$-pullbacks along $V\times W\to Y\times Z$, with $(V\to Y)\in ((\bfC_2)_{\horiz_2})_{/Y}\times_{(\bfC_2)_{\horiz_2}}(\bfC_1)_{\horiz_1}$ and $(W\to Z)\in ((\bfC_2)_{\horiz_2})_{/Z}\times_{(\bfC_2)_{\horiz_2}}(\bfC_1)_{\horiz_1}$.

This shows that if $f\in \mathrm{HR}_2$. The other three cases can be proved similarly.
\end{proof}

We have a dual version of \Cref{prop-sheaf-theory-right-Kan-extension}.

\begin{proposition}\label{prop-dual-sheaf-theory-right-Kan-extension}
Suppose we have $(\bfC_i, \verti_i, \horiz_i),\ i=1,2$, and a finite limit preserving fully faithful embedding $\bfC_1\subset \bfC_2$ which induces fully faithful embeddings $(\bfC_1)_{\verti_1}\subset (\bfC_2)_{\verti_2}$, $(\bfC_1)_{\horiz_1}\subset (\bfC_2)_{\horiz_2}$. 
Suppose the class $\horiz_2$ are representable in $\horiz_1$.
Then every sheaf theory $\der_1:\corr(\bfC_1)_{\verti_1, \horiz_1}\to \lincat_\La$ admits an extension
\[
\der_2: \corr(\bfC_2)_{\verti_2;\horiz_2}\to \lincat_\La
\]
such that the restriction $\der_2|_{(\bfC_2)_{\verti_2}}$ is canonically isomorphic to the left Kan extension of $\der_1|_{(\bfC_1)_{\verti_1}}$ along $(\bfC_1)_{\verti_1}\subset (\bfC_2)_{\vert_2}$ (as plain functors). If $\der_1$ takes value in $\cptcat_\La$, so is $\der_2$.
\end{proposition}
\begin{proof}
By assumption, $\corr(\bfC_1)_{\verti_1;\horiz_1}\to \corr(\bfC_2)_{\verti_2;\horiz_2}$ is full faithful.
Let 
\[
\der_2 \colon \corr(\bfC_2)_{\verti_2;\horiz_2} \to \lincat_\La 
\]
be a left operadic Kan extension (see \cite[Definition 3.1.2.2]{Lurie.higher.algebra}). 
Recall that tensor product in $\lincat_\La$ preserves colimits separately in each variable.
Then arguing similarly as in \Cref{prop-sheaf-theory-right-Kan-extension}, one shows that the functor analogous to \eqref{eq: cofinality-horiz-to-corr} in the current setting is cofinal. Then using \cite[Proposition 3.1.1.16]{Lurie.higher.algebra}, one sees that the value of $\der_2$ at $X\in \bfC_2$ is 
$
\colim_{X'\in((\bfC_2)_{\verti_2})_{/X}  \times_{(\bfC_2)_{\verti_2}}  (\bfC_1)_{\verti_1}}\der(X')
$, as desired. 
\end{proof}

Next we consider extensions some along non-full embeddings $\corr(\bfC_1)_{\verti_1;\horiz_1}\to \corr(\bfC_2)_{\verti_2;\horiz_2}$. This is in generally difficult, as Kan extensions along non-full embeddings are difficult to compute. However, under certain assumptions, they are still manageable.

\begin{proposition}\label{prop-sheaf-theory-enlarge-H}
Let $\der: \corr(\bfC)_{\verti;\horiz}\to \lincat_\La$ be a sheaf theory. Let $\horiz'$ be a weakly stable class of morphisms. Suppose for every $(f:X\to Y)\in \bfC_{\horiz'}$, there is $(U\to X)\in \bfC_{\horiz}$ that is universally $\der$-descent such that the composed morphism $U\to X\to Y$ belongs to $\horiz$. Then $\der$ admits an extension $\der'\colon  \corr(\bfC)_{\verti;\horiz'} \to \lincat_\La$.
\end{proposition}
\begin{proof}
We take $\der'$ to be the right Kan extension along $\corr(\bfC)_{\verti;\horiz}\to \corr(\bfC)_{\verti;\horiz'}$. Then we need to show that for every $Z\in \bfC$, the functor 
\[
\der(Z)\to \lim_{Y\in (\corr(\bfC)_{\verti;\horiz'})_{Z/}\times_{\corr(\bfC)_{\verti;\horiz'}}\corr(\bfC)_{\verti;\horiz}}\der(Y)=\der'(Z)
\] 
is an equivalence (so $\der'$ is indeed an extension of $\der$). First as argued in the proof of \Cref{prop-sheaf-theory-right-Kan-extension}, $(\bfC_{\horiz'})^{\op}_{Z/}\times_{(\bfC_{\horiz'})^{\op}}(\bfC_{\horiz})^{\op}\to (\corr(\bfC)_{\verti;\horiz'})_{Z/}\times_{\corr(\bfC)_{\verti;\horiz'}}\corr(\bfC)_{\verti;\horiz}$ is cofinal so it is enough to show that
\[
\der(Z)\to \lim_{Y\in (\bfC_{\horiz'})^{\op}_{Z/}\times_{(\bfC_{\horiz'})^{\op}}(\bfC_{\horiz})^{\op}}\der(Y)
\]
is an equivalence. By this will follow if we can show that the functor $(\bfC_{\horiz'})^{\op}_{Z/}\times_{(\bfC_{\horiz'})^{\op}}(\bfC_{\horiz})^{\op}\to \bfC_{\horiz}\xrightarrow{\der}\lincat_\La$ is the right Kan extension of its restriction to $(\bfC_{\horiz})^{\op}_{Z/}$. Let $(g:Y\to Z)\in (\bfC_{\horiz'})^{\op}_{Z/}\times_{(\bfC_{\horiz'})^{\op}}(\bfC_{\horiz})^{\op}$. Then $((\bfC_{\horiz'})^{\op}_{Z/}\times_{(\bfC_{\horiz'})^{\op}}(\bfC_{\horiz})^{\op})_{g/}\times_{(\bfC_{\horiz'})^{\op}_{Z/}\times_{(\bfC_{\horiz'})^{\op}}(\bfC_{\horiz})^{\op}}(\bfC_{\horiz})^{\op}_{Z/}$ can be identified with the category $(\mI_g)^{\op}$, where $\mI_g$ consists of those $(f: Y'\to Y)\in \bfC_{\horiz}$ such that $(gf: Y'\to Z)\in \bfC_{\horiz}$. Therefore, we reduce to show that
\[
\lim_{(Y'\to Y)\in (\mI_g)^{\op}}\der(Y')\cong \der(Y).
\]
By assumption, we can find some $(U\to Y)\in\bfC_\horiz$ which is universal $\der$-descent such that $(U\to Y\to Z)\in \bfC_{\horiz}$. We fix such $U\to Y$.
Let $\mJ_{g,U}\subset \mI_g$ be the subcategory consisting of those $Y'\to Y$ that can be factorized as $Y'\to U\to Y$ with $(Y'\to U)\in \bfC_{\horiz}$. Then it is enough to show that: (a) $\der(Y)\cong  \lim_{(Y'\to Y)\in \mJ_g} \der(Y')$; and (b) $(\mI_g)^{\op}\to \bfC_{\horiz}\xrightarrow{\der} \lincat_\La$ is the right Kan extension of its restriction to $(\mJ_g)^{\op}$.

For (a), let $U^\bullet \to Y$ be the \v{C}ech nerve of $U\to Y$. Then $(U_\bullet)^{\op}\to (\mJ_{g,U})^{\op}$ is cofinal so $\der(Y)\cong \lim \der(U^\bullet)\cong \lim_{\mJ_{g,U}}\der(Y')$.
For (b), let $(f: Y'\to Y)\in \mI_g$. Then $(\mI_g)^{\op}_{f/}\times_{(\mI_g)^{\op}}(\mJ_g)^{\op}$ can be identified with $\mJ_{gf,Y'\times_YU}$, so (b) follows from (a). The proposition if proved.
\end{proof}

Similarly we have the following.

\begin{proposition}\label{lem-sheaf-theory-left-right-Kan-extension}
Let $\der: \corr(\bfC)_{\verti;\horiz}\to \lincat_\La$ be a sheaf theory. Suppose $\verti$ and $\horiz$ are strongly stable. 
Let $\verti'$ be a class of morphisms consisting of those $f: X\to Y$, such that there is some $(U\to Y)\in \bfC_{\horiz}$ that is universally $\der$-descent such that the base change $X\times_YU\to U$ of $f$ belongs to $\verti$. Then the class $\verti'$ is strongly stable and $\der$ admits an extension $\der'\colon  \corr(\bfC)_{\verti';\horiz} \to \lincat_\La$.
\end{proposition}
\begin{proof}
It is clear that $\verti'$ is strongly stable.

We take $\der'$ to be the right Kan extension along $\corr(\bfC)_{\verti;\horiz}\to \corr(\bfC)_{\verti';\horiz}$. Then we need to show that for every $Z\in \bfC$, the functor 
\[
\der(Z)\to \lim_{Y\in (\corr(\bfC)_{\verti';\horiz})_{Z/}\times_{\corr(\bfC)_{\verti';\horiz}}\corr(\bfC)_{\verti;\horiz}}\der(Y)=\der'(Z)
\] 
is an equivalence (so $\der'$ is indeed an extension of $\der$). This will follow if we can show that
the functor $(\corr(\bfC)_{\verti';\horiz})_{Z/}\times_{\corr(\bfC)_{\verti';\horiz}}\corr(\bfC)_{\verti;\horiz}\to \corr(\bfC)_{\verti;\horiz}\xrightarrow{\der}\lincat_\La$ is the right Kan extension of its restriction to $(\corr(\bfC)_{\verti;\horiz})_{Z/}$.

Let $Z\xleftarrow{g} X\xrightarrow{f} Y$ be a correspondence with $f\in \bfC_{\verti'}$. Then the category
\begin{equation}\label{eq: complicated slicing category}
\bigl((\corr(\bfC)_{\verti';\horiz})_{Z/}\times_{\corr(\bfC)_{\verti';\horiz}}\corr(\bfC)_{\verti;\horiz}\bigr)_{f/}\times_{(\corr(\bfC)_{\verti';\horiz})_{Z/}\times_{\corr(\bfC)_{\verti';\horiz}}\corr(\bfC)_{\verti;\horiz}}(\corr(\bfC)_{\verti;\horiz})_{Z/}
\end{equation}
can be identified with $(Y\leftarrow Y'\to W)\in (\corr(\bfC)_{\verti;\horiz})_{Y/}$ such that the composed correspondence $Z\leftarrow X\leftarrow X':=Y'\times_YX\to Y'\to W$ belongs to $\corr(\bfC)_{\verti;\horiz}$. As $\verti$ is strongly stable, this implies that $X'\to Y'$ belongs to $\verti$. Therefore, if we let $\mI_f$ be the full subcategory of $(\bfC_{\horiz})_{/Y}$ consisting of those $Y'\to Y$ such that the base change $f': X'\to Y'$ of $f$ belongs to $\verti$, 
then as argued in the proof of \Cref{prop-sheaf-theory-right-Kan-extension}, $\mI_f$ is initial in \eqref{eq: complicated slicing category}.

Therefore, it is enough to show that $\der(Y)\cong \lim_{\mI_f^{\op}} \der(Y')$. 
Let $U\to Y$ be the universally $\der$-descent morphism associated to $X\to Y$ as in the assumption. We may also consider $\mJ_{f,U}\subset (\bfC_\horiz)_{/Y}$ consisting of those $Y'\to Y$ that can be factorized as $Y'\to U\to Y$, with $(Y'\to U)\in\bfC_\horiz$.  We reduce to show that: (a) $\der(Y)\cong  \lim_{(\mJ_{f,U})^{\op}} \der(Y')$; and (b) $\der: (\mI_f)^{\op}\to \lincat_\La$ is the right Kan extension of its restriction along $(\mJ_{f,U})^{\op}\to (\mI_f)^{\op}$. For (a), let $U^\bullet\to Y$ be the \v{C}ech nerve of $U\to Y$. Then $U^\bullet\to \mJ_{f,U}$ is cofinal and $\der(Y)\cong \lim \der(U^\bullet)$ by assumption. Therefore $\der(Y)\cong  \lim_{(\mJ_{f,U})^{\op}} \der(Y')$. For (b), let $(g:Y'\to Y)\in \mI_f$. Then $(\mI_f)_{/g}\times_{\mI_f}\mJ_{f,U}$ can be identified with $\mJ_{f',U'}$, where $f':X'\to Y'$ is the base change of $f:X\to Y$ along $g$ and $U'$ is the base change of $U$ along $g$. Therefore (b) follows from (a).

Finally one similarly argue as in \Cref{prop-sheaf-theory-right-Kan-extension} to show that $\der'$ is lax symmetric monoidal. The proposition follows. 
\end{proof}

Similar ideas yield the following result, which slightly generalizes \cite[Proposition A.5.14]{Lucas.Mann}. 

\begin{proposition}\label{lem-sheaf-theory-left-Kan-extension}
Let $\der: \corr(\bfC)_{\verti;\horiz}\to \lincat_\La$ be a sheaf theory. Suppose that $\verti$ is strongly stable.
Let $\verti'\supset \verti$ be another weakly stable class of morphisms.  
Suppose that for every $X\in \bfC$, there is a subcategory $\mS_X\subset (\bfC_{\verti})_{/X}$, and for every $(f: X\to Y)\in \bfC_\verti$ there is a full subcategory $\mS_f\subset \mS_X$ satisfying the following properties.
\begin{enumerate}
\item\label{lem-sheaf-theory-left-Kan-extension-1} The inclusions $\mS_f\subset \mS_X\subset (\bfC_{\verti})_{/X}$ respect finite products.
\item\label{lem-sheaf-theory-left-Kan-extension-2} For every $(f: X\to Y)\in \bfC_{\verti}$, the functor 
\[
(\bfC_{\verti})_{/Y}\to (\bfC_{\verti})_{/X}\colon (Y'\to Y)\mapsto (X':=X\times_YY'\to X)
\] 
restricts to a functor $\mS_Y\to \mS_f$, and for every $(X'\to X)\in \mS_f$, the composed map $X'\to X\to Y$ can be factorized as $X'\to Y'\to Y$ with $(Y'\to Y)\in \mS_Y$.
\item\label{lem-sheaf-theory-left-Kan-extension-3} 
The natural functor 
\[
\colim_{X'\in \mS_f} \der(X')\to \der(X)
\] 
is an equivalence.
\item\label{lem-sheaf-theory-left-Kan-extension-4} For every $(f: X\to Y)\in \bfC_{\verti'}$, and $(g:X'\to X)\in \mS_X$, $f\circ g\in \bfC_{\verti}$.
\end{enumerate}
Then $\der$ admits an extension $\der'\colon  \corr(\bfC)_{\verti';\horiz} \to \lincat_\La$. 
\end{proposition}
We introduce a category $\mS_f$ in the proposition to have some extra flexibility to apply this result. In many applications, it is enough to assume that the functor in \eqref{lem-sheaf-theory-left-Kan-extension-2} restricts to a functor $\mS_Y\to \mS_X$ and then take $\mS_f$ to be the span of the essential image of $\mS_Y\to\mS_X$.
\begin{proof}
 Let 
\[
\der' \colon \corr(\bfC)_{\verti';\horiz} \to \lincat_\La 
\]
be a left operadic Kan extension along $\corr(\bfC)_{\verti;\horiz}\to \corr(\bfC)_{\verti';\horiz}$. Then as argued in \Cref{prop-dual-sheaf-theory-right-Kan-extension}, the value of $\der'$ at $X\in \bfC$ is 
$
\colim_{X'\in(\bfC_{\verti'})_{/X}  \times_{\bfC_{\verti'}}  \bfC_{\verti}}\der(X')
$. 
We need to show that it is equivalent to $\der(X)$. For this purpose, it is enough to show that the functor $(\bfC_{\verti'})_{/X}  \times_{\bfC_{\verti'}}  \bfC_{\verti}\to \bfC_{\verti}\xrightarrow{\der} \lincat$ is isomorphic to the left Kan extension of its restriction along $(\bfC_{\verti})_{/X}\to (\bfC_{\verti'})_{/X}  \times_{\bfC_{\verti'}}  \bfC_{\verti}$. For a morphism $(f: X'\to X)\in \bfC_{\verti'}$, we let 
\[
\mI_f= (\bfC_{\verti})_{/X}\times_{\bigl((\bfC_{\verti'})_{/X}  \times_{\bfC_{\verti'}}  \bfC_{\verti}\bigr)}\bigl((\bfC_{\verti'})_{/X}  \times_{\bfC_{\verti'}}  \bfC_{\verti}\bigr)_{/f},
\]
which is nothing but the full subcategory of $(\bfC_{\verti})_{/X'}$ consisting of those $g: X''\to X'$ such that $f\circ g\in \verti$. Then we need to show that 
\begin{equation}\label{eq-sheaf-theory-left-Kan-extension}
\colim_{X''\in \mI_f}\der(X'')\cong \der(X').
\end{equation}

Now for a morphism $(g: X''\to X')\in \bfC_\verti$, let $\mT_g\subset (\bfC_{\verti})_{/X''}$ be the full subcategory consisting of those $Z\to X''$ such that the composed morphism $Z\to X''\to X'$ can be factorized as $Z\to W\to X'$ with $(W\to X')\in \mS_{X'}$. Note that we have a cofinal inclusion $\mS_g\subset \mT_g$. Indeed, for every such $(h: Z\to X'')\in \mT_g$, Assumption \eqref{lem-sheaf-theory-left-Kan-extension-1} \eqref{lem-sheaf-theory-left-Kan-extension-2} (together with the assumption that $\verti$ is strongly stable) implies that  $(\mT_{g})_{h/}\times_{\mT_{g}}\mS_{g}$ is non-empty and admits binary products and therefore is weakly contractible.
It then follows from Assumption \eqref{lem-sheaf-theory-left-Kan-extension-3} that $\colim_{Z\in \mT_g} \der(Z)\to \der(X'')$ is an equivalence.

Now applying the above observation to $g=\id_{X'}: X'\to X'$ (and write $\mT_{X'}$ instead of $\mT_{\id_{X'}}$), we see that $\colim_{Z\in \mT_{X'}}\der(Z)\cong \der(X')$. In addition, by Assumption \eqref{lem-sheaf-theory-left-Kan-extension-4}, we see that $\mT_{X'}\subset \mI_f$. Then \eqref{eq-sheaf-theory-left-Kan-extension} would follow if we show that for every $(g: X''\to X')\in \mI_f$,
\[
\colim_{Z\in  \mT_{X'}\times_{\mI_f}(\mI_f)_{/g}} \der(Z)\cong \der(X'').
\]
But the index category $ \mT_{X'}\times_{\mI_f}(\mI_f)_{/g}$ is nothing but $\mT_g$ as above, so the desired equivalence follows.
\end{proof}

In practice, the category $\mS_X$ could come as certain covering family of $X$ under some Grothendieck topology of $\bfC$. Here is a sample. 
\begin{corollary}\label{lem-sheaf-theory-left-Kan-extension-covering}
Let $\der: \corr(\bfC)_{\verti;\horiz}\to \lincat_\La$ be a sheaf theory. Assume that $\verti$ is strongly stable, and let $\verti'\supset \verti$ be the class of morphisms containing those
$f: X\to Y\in \bfC$ such that there exists universally $\der$-codescent $\varphi: U\to X\in \bfC_{\verti}$ satisfying $f\circ \varphi\in \verti$. Then $\verti'$ is strongly stable and
the $\der$ extends to a sheaf theory $\der': \corr(\bfC)_{\verti';\horiz}\to \lincat_\La$.
\end{corollary}
Note that together with \Cref{prop-descent-codescent-abstract} \eqref{prop-descent-codescent-abstract-1} , this result gives (part of) \cite[\textsection{4}]{liu2012enhanced} and \cite[Proposition A.5.14]{Lucas.Mann}. 
\begin{proof}
We first notice that $\verti'$ is clearly stable under base change. As universal $\der$-codescent morphisms (in $\verti$) are stable under compositions (\Cref{Liu.Zheng.descent.lemma}), one shows that $\verti$ is also stable under compositions and satisfying '$2$ out of $3$' property. 

To prove the extension of $\der$,
the only thing one needs to observe that in \Cref{lem-sheaf-theory-left-Kan-extension}, there is no need to a priori to assign every $(f:X\to Y)\in \bfC_\verti$ the categories $\mS_f\subset \mS_X$. All we need is to prove the equivalence \eqref{eq-sheaf-theory-left-Kan-extension}  for every morphism $(f:X'\to X)\in \bfC_{\verti'}$. Then we just need to assign $\mS_g\subset \mS_{X''}$ for every $(g: X''\to X')\in \bfC_{\verti}$. For this, we choose $\varphi: U\to X'$ as in the assumption, and for every $(g: X''\to X')\in \bfC_{\verti}$ let $\mS_g=\mS_{X''}$ be the base change of the \v{C}ech nerve of $U\to X'$.
\end{proof}

\begin{remark}
Clearly in \Cref{prop-sheaf-theory-right-Kan-extension}-\Cref{lem-sheaf-theory-left-Kan-extension-covering}, we may replace $\lincat_\La$ by $\cat$ as the codomain of the sheaf theory.
\end{remark}

For our purpose, we need another situation such collection $\{\mS_X\}_X$ exists. The following statement is the combination of \Cref{prop-sheaf-theory-right-Kan-extension} and \Cref{lem-sheaf-theory-left-Kan-extension}.

\begin{corollary}\label{prop-sheaf-theory-non-full-right-Kan-extension}
Suppose we have $(\bfC_i, \verti_i, \horiz_i),\ i=1,2$, and a finite limit preserving fully faithful embedding $\bfC_1\subset \bfC_2$ which induces fully faithful embeddings $(\bfC_1)_{\verti_1}\subset (\bfC_2)_{\verti_2}, (\bfC_1)_{\horiz_1}\subset (\bfC_2)_{\horiz_2}$. Suppose Assumption \eqref{prop-sheaf-theory-right-Kan-extension-1} in \Cref{prop-sheaf-theory-right-Kan-extension} holds, and suppose $\verti_1$ is strongly stable. Let $\verti_{2,r}\subset \verti_2$ be the subset of morphisms are that representable in $\verti_1$. Let 
\begin{equation}\label{eq-sheaf-theory-to-be-extended}
\der_1: \corr(\bfC_1)_{\verti_1;\horiz_1}\to \lincat_\La
\end{equation} 
be a sheaf theory.

Suppose that there is a strongly stable class $\mathrm{S}_1\subset\verti_1\cap \horiz_1$, satisfying the following conditions.
Let $\mathrm{S}_2\subset \verti_2\cap \horiz_2$ be the subset of morphisms that are representable in $\mathrm{S}_1$. For every $X\in \bfC_2$, write $\mS_X:=(\bfC_1)_{\mathrm{S}_1}\times_{(\bfC_2)_{\mathrm{S}_2}}((\bfC_2)_{\mathrm{S}_2})_{/X}$. Then
\begin{enumerate}
\item\label{prop-sheaf-theory-non-full-right-Kan-extension-a} For every $(f: X\to Y)\in\bfC_{\verti_2}$ and $(g: X'\to X)\in \mS_X$, the composition $(f\circ g: X'\to Y)\in \bfC_{\verti_{2,r}}$, and for every $(f: X\to Y)\in\bfC_{\verti_{2,r}}$ and $(g: X'\to X)\in \mS_X$, the composition $f\circ g: X'\to Y$ can be factorized as $X'\to Y'\to Y$ with $(Y'\to Y)\in \mS_Y$. 
\item\label{prop-sheaf-theory-non-full-right-Kan-extension-b} The inclusion $\mS_X\to (\bfC_1)_{\horiz_1}\times_{(\bfC_2)_{\horiz_2}}((\bfC_2)_{\horiz_2})_{/X}$ is cofinal.
\item\label{prop-sheaf-theory-non-full-right-Kan-extension-c} The restriction $\der_1|_{\corr((\bfC_1)_{\mathrm{S}_1})}$ is isomorphic to sheaf theory from $\der_1|_{\corr(\bfC_1)_{\isom;\mathrm{S}_1}}$ by applying \Cref{prop-sheaf-theory-for-adjoint-factorization} to $\mathrm{E}=\mathrm{S}_1$.
\end{enumerate}
Then there is an extension of sheaf theory
\[
\der_2: \corr(\bfC_2)_{\verti_2;\horiz_2}\to \lincat_\La
\] 
of \eqref{eq-sheaf-theory-to-be-extended}
along the (non-full) embedding $ \corr(\bfC_1)_{\verti_1;\all}\to  \corr(\bfC_2)_{\verti_2;\all}$, such that
\begin{enumerate}[(a)]
\item\label{prop-sheaf-theory-non-full-right-Kan-extension-1} the restriction $\der_2|_{\corr(\bfC_1)_{\verti_1, \horiz_1}}=\der_1$;
\item\label{prop-sheaf-theory-non-full-right-Kan-extension-2} $\der|_{((\bfC_2)_{\horiz_2})^{\op}}$ is isomorphic to the right Kan extension of $\der|_{((\bfC_1)_{\horiz_1})^{\op}}$ along $(\bfC_1)_{\horiz_1}\subset (\bfC_2)_{\horiz_2}$.
\end{enumerate}
\end{corollary}
\begin{proof}
We factorize the inclusion $\corr(\bfC_1)_{\verti_1;\all}\to\corr(\bfC_2)_{\verti_2;\all}$ as
\[
\corr(\bfC_1)_{\verti_1;\horiz_1}\to \corr(\bfC_2)_{\verti_{2,r};\horiz_2} \to \corr(\bfC_2)_{\verti_2;\horiz_2},
\]
and first apply \Cref{prop-sheaf-theory-right-Kan-extension} to extend $\der_1$ to a sheaf theory $\der_{2,r}: \corr(\bfC_2)_{\verti_{2,r};\horiz_2}\to \lincat_\La$. Then we apply \Cref{lem-sheaf-theory-left-Kan-extension} to extend $\der_{2,r}$ along $\corr(\bfC_2)_{\verti_{2,r};\all}\to \corr(\bfC_2)_{\verti_2;\horiz_2}$ to define $\der_2:  \corr(\bfC_2)_{\verti_2;\horiz_2}\to \lincat_\La$. For this, we need to verify all the assumptions of  \Cref{lem-sheaf-theory-left-Kan-extension}.

Indeed, as $\mathrm{S}_1$ is strongly stable, we see that $\mS_X\subset ((\bfC_2)_{\verti_{2,r}})_{/X}$ is preserved under finite products. We may take $\mS_f=\mS_X$ for any $(f: X\to Y)\in(\bfC_2)_{\verti_{2,r}}$. Then it follows that Assumptions \eqref{lem-sheaf-theory-left-Kan-extension-1} \eqref{lem-sheaf-theory-left-Kan-extension-2} and \eqref{lem-sheaf-theory-left-Kan-extension-4} of \Cref{lem-sheaf-theory-left-Kan-extension} hold. To see \eqref{lem-sheaf-theory-left-Kan-extension-3} also holds, we notice that
as $\der_{2,r}|_{(\bfC_2)^{\op}}$ is canonically isomorphic to the right Kan extension of $\der_1|_{(\bfC_1)^{\op}}$, for $X\in \bfC_2$, we have
\begin{align*}
\der_{2,r}(X)  \cong  \lim_{X'\in (\bfC_1\times_{\bfC_2}(\bfC_2)_{/X})^{\op}}  \der_1(X')  \cong \lim_{X'\in (\mS_X)^{\op}}  \der_1(X') \cong \colim_{X'\in \mS_X}  \der_1(X'),
\end{align*}
where the second equivalence follows from \eqref{prop-sheaf-theory-non-full-right-Kan-extension-b}, and the last equivalence follows by \eqref{prop-sheaf-theory-non-full-right-Kan-extension-c}. The proof is complete.
\end{proof}

\begin{remark}
The above theorem is closely related to \cite[Theorem 8.1.1.9, Proposition 9.3.3.3]{Gaitsgory.Rozenblyum.DAG.vol.I}. In fact, it is easily to deduce \cite[Theorem 8.1.1.9, Proposition 9.3.3.3]{Gaitsgory.Rozenblyum.DAG.vol.I} (under weaker assumptions) by similar reasonings as above. 
\end{remark}

\subsection{Geometric traces in sheaf theory}\label{sec:trace-geometric-trace-main-general}

Now we follow ideas of \cite{benzvi2009character} \cite{ben2017spectral} to develop a method to calculate the (twisted) categorical trace of monoidal categories arising from convolution pattern in  the formalism of category of correspondences and abstract sheaf theory as in \Cref{sec:stacks-correspondences-general} and \Cref{sec:symmetric-monoidal-and-projection-for-corr}. As mentioned before, Compared with the work of \emph{loc. cit.}, we will first calculate a geometric version of categorical trace. Then we will compare the geometric version with the usual version in favorable cases. Our approach allows us to bypass integral transform of sheaf theories, which usually do not hold in the $\ell$-adic setting.

Let $\der: \corr(\bfC)_{\verti;\horiz}\to \lincat_\La$ be a sheaf theory. 
We will make the following assumption on the sheaf theory $\der$.

\begin{assumption}\label{ass: base of sheaf theory}
The symmetric monoidal category $\der(\pt)$ is rigid.
\end{assumption}

In many examples, the canonical functor $\Mod_\La\to \der(\pt)$ is an equivalence so the above assumption is satisfied. 

\subsubsection{Geometric Hochschild homology}\label{sec:trace-set-up-assumptions}

Let $A$ be an associative algebra object in $\corr(\bfC)_{\verti; \horiz}$, and let $M$ be a left $A$-module object in $\corr(\bfC)_{\verti; \horiz}$. As $\der$ is a lax symmetric monoidal functor, $\der(A)$ is an algebra object in $\lincat_{\der(\pt)}$ and $\der(M)$ is a $\der(A)$-module object in $\lincat_{\der(\pt)}$. 
Similarly, if $F$ is an $A$-bimodule, then $\der(F)$ is a $\der(A)$-bimodule. Then one can form its Hochschild homology (a.k.a categorical trace) of $(\der(A),\der(F))$
\[
\tr(\der(A),\der(F))= \der(A)\otimes_{\der(A)\otimes_{\der(\pt)}\der(A)^{\rev}}\der(F)\in \lincat_{\der(\pt)}.
\]
In practice, however, we need to consider a variant $\trg(\der(A),\der(F))$, which we call the geometric trace of $\der(F)$. Namely, we consider the Yoneda embedding 
\[
\corr(\bfC)_{\verti;\horiz}\to \mP(\corr(\bfC)_{\verti;\horiz}),
\]
where $\mP(\corr(\bfC)_{\verti;\horiz})$ is the category of presheaves on $\corr(\bfC)_{\verti;\horiz}$ equipped with the induced symmetric monoidal structure, which by definition preserves colimits in each variable (see \cite[Corollary 4.8.1.12]{Lurie.higher.algebra}). Then we have the Hochschild homology of the $A$-bimodule $F$ in $\mP(\corr(\bfC)_{\verti;\horiz})$
\[
\tr(A,F):=|\HH(A,F)_\bullet|\in \mP(\corr(\bfC)_{\verti;\horiz}).
\]
By the universal property of $\mP(\corr(\bfC)_{\verti;\horiz})$, the functor $\der: \corr(\bfC)_{\verti;\horiz}\to \lincat_\La$ extends to a continuous functor $\der: \mP(\corr(\bfC)_{\verti;\horiz})\to \lincat_\La$. Then we define the \textit{geometric trace} of $\der(F)$ as
\[
\trg(\der(A),\der(F)):=\der(\tr(A,F)).
\]
Explicitly, $\trg(\der(F),\der(A))$ can be computed in the following way.
We first apply the functor $\der$ to the standard Hochschild complex \eqref{eq:standard-hochschild-complex} (which now is a simplicial object in $\corr(\bfC)_{\verti;\horiz}$) to obtain a simplicial object $\der(\HH_\bullet(A,F))$ in $\lincat_{\der(\pt)}$. Then the geometric trace $\trg(\der(A),\der(F))$ is the geometric realization of this simplicial object in $\lincat_{\der(\pt)}$
\begin{equation}\label{eq:geometric-trace}
    \trg(\der(A),\der(F))\cong |\der(\HH(A,F)_\bullet)|.
\end{equation}
We emphasize that $\trg(\der(A),\der(F))$ depends not only on $\der(F)$, but on the $A$-bimodule $F$ itself (and of course the functor $\der$).

In particular, for $A$ equipped with an algebra endomorphism $\phi: A\to A$ we have the $A$-bimodule $F={}^{\phi}A$ in $\corr(\bfC)_{\verti;\horiz}$ as before. We write
\[
\trg(\der(A),\phi) = \trg(\der(A),\der({}^{\phi}A)).
\] 
\begin{remark}\label{rmk:trace-to-geometric-trace}
As $\der$ is equipped with a lax monoidal structure we get a natural comparison functor
\begin{equation}\label{eqn:trace-to-geometric-comparison-map-general}
       \tr(\der(A),\der(F))\simeq |\der(A)^{\otimes_{\der(\pt)} \bullet}\otimes_{\La} \der(F)| \rightarrow |\der(A^{\otimes\bullet}\otimes F)|=\trg(\der(A),\der(F))
\end{equation}
from the usual trace of $\der(F)$ to the geometric trace. This functor is not an equivalence in general. Of course, if for each $n$, the functor $\der(A)^{\otimes_{\der(\pt)} n}\otimes_{\der(\pt)} \der(Q)\to \der(A^n\times Q)$ is an equivalence, then the comparison map \eqref{eqn:trace-to-geometric-comparison-map-general} is an equivalence. We will see later that this functor is an equivalence in many more cases of interest, as shown by \Cref{prop-comparison-usual-trace-geo-trac} below.
\end{remark}


\subsubsection{Fixed point objects and geometric traces of convolution categories}\label{subsec-trace.convolution.categories}
We specialize the previous constructions to the situation appearing in our applications. Let $(f: X\to Y)\in\bfC$ as in \Cref{rem:segal.objects.morphisms.vert.horiz}, that is, $\Delta_X: X\to X\times X$ and $\pi_X: X\to \pt$ belong to $\bfC_{\horiz}$, and $f$ and  the relative diagonal map $\Delta_{X/Y}: X\to X\times_YX$ belongs to $\bfC_{\verti}$. Let $X_\bullet \rightarrow Y$ denote the \v{C}ech nerve of $f$. From \Cref{ex:Cech-nerve} and \Cref{rem:segal.objects.morphisms.vert.horiz}, we see that
\[
X_1 : = X\times_Y X
\]
has a structure of an associative algebra object in $\corr(\bfC)_{\verti;\horiz}$, with the multiplication and unit maps are given in \eqref{eq-convolution-product}.
Let $Z\in \bfC$ equipped with two morphisms $g_{i}\colon Z\to Y,\ i=1,2$ in $\bfC$ and let 
\[
Q=X\times_YZ\times_YX=Z\times_{Y\times Y}(X\times X).
\]
Then  the object $Q$ admits the structure of an $(X\times_Y X)$-bimodule in $\corr(\bfC)_{\verti;\horiz}$ (see \Cref{sec-categorical.action.on.bimodules}). In particular, the left action is given by the diagram
\begin{equation*}
\begin{tikzcd}
X\times_Y X\times_YZ\times_YX \arrow[rr, "\id\times\Delta_X\times\id\times\id"]\arrow[d,"\id\times f\times \id\times\id "]&& (X\times_Y X) \times (X\times_YZ\times_YX) \\
X\times_YZ\times_YX && 
\end{tikzcd}
\end{equation*}
and we have a similar diagram for the right action. Consider the following diagram
\begin{equation}\label{eq:trace-convolution-horocycle-diagram}
\xymatrix{
X\times_{Y\times Y} Z \ar^-{\delta_0=(\Delta_X\times \id_Z)}[rr]\ar_{q=(f\times \id_Z)}[d] && (X\times X)\times_{Y\times Y}Z\\
Y\times_{Y\times Y} Z, &&
}
\end{equation}
which induces a functor $q\vertl\circ (\delta_0)\horizl:\der(X\times_YZ\times_YX)\to\der(Y\times_{Y\times Y}Z)$.

Recall associated to a sheaf theory $\der$, there are classes of morphisms $\mathrm{VR}$ and $\mathrm{HR}$ associated to $\der$, as defined in \Cref{rem-additional-base.change.sheaf.theory} \eqref{rem-additional-base.change.sheaf.theory-0}.
\begin{proposition}\label{bimodule.geom.trace.fully.faithfulness.convolution.cat}
The following diagram is commutative
  \[
  \xymatrix{
  \der(X\times_YZ\times_YX) \ar^{(\delta_0)\horizl}[r]\ar[d] & \der(X\times_{Y\times Y}Z)\ar^-{q\vertl}[dd]\\
  \tr(\der(X\times_YX),\der(X\times_YZ\times_YX))\ar[d]&\\
  \trg(\der(X\times_YX),\der(X\times_YZ\times_YX)) \ar[r] & \der(Y\times_{Y\times Y}Z).
  }
  \]

 Suppose, in addition, 
 \begin{itemize}
\item $f\colon X\to Y \in \bfC_{\mathrm{VR}}$; and 
\item $\Delta_X: X\to X\times X \in \bfC_{\mathrm{HR}}$.
\end{itemize}
  Then the bottom horizontal functor of the above diagram is fully faithful, with the essential image generated (as presentable $\La$-linear categories) by the image of $q\vertl\circ \delta_0\horizl$. The bottom horizontal functor admits a continuous right adjoint, denoted by 
  \[
  \proj_{\trg}: \der(Y\times_{Y\times Y}Z)\to \trg(\der(X\times_YX), \der(X\times_YZ\times_YX)).
  \]
\end{proposition}

The proof of the proposition will be given at the end of \Cref{sec:relative-resolutions-geometric}. Here are some remarks regarding the assumptions of the proposition.

\begin{remark}\label{rem: assumptions of bimodule.geom.trace.fully.faithfulness.convolution.cat}  
\begin{enumerate}
\item We note that there is no assumption on $(g_1,g_2): Z\to Y\times Y$.
\item It will be clear from the proof that \Cref{bimodule.geom.trace.fully.faithfulness.convolution.cat} holds under the weaker assumption that $f\colon X\to Y$ satisfies \Cref{assumptions.base.change.sheaf.theory.V} \eqref{assumptions.base.change.sheaf.theory.V-0}-\eqref{assumptions.base.change.sheaf.theory.V-2} and $\Delta_X: X\to X\times X$ satisfies \Cref{assumptions.base.change.sheaf.theory.H} \eqref{assumptions.base.change.sheaf.theory.H-0}-\eqref{assumptions.base.change.sheaf.theory.H-2}.
\end{enumerate}
\end{remark}

We specialize \Cref{bimodule.geom.trace.fully.faithfulness.convolution.cat} to the following two cases.

First, suppose there are  morphisms $\phi_{X}\colon  X\rightarrow X$ and $\phi_{Y}\colon  Y\rightarrow Y$ in $\bfC$, together with an equivalence 
\begin{equation}\label{eq: meaning of equiv in higher sense}
f\circ \phi_{X} \simeq \phi_{Y}\circ f. 
\end{equation}
We will usually abuse notation and denote both maps by $\phi$ if it is clear from context. 
In this case, if we let $Z=Y$ with the map $g_1=\id$ and $g_2=\phi$, then
 $Z\times_{Y\times Y}Y$ is nothing but the $\phi$-fixed point object $\mathcal{L}_{\phi}(Y)$, defined by the pullback
\begin{equation}\label{eq:tau-fixed-point}
\begin{tikzcd}
\mathcal{L}_{\phi}(Y) \arrow[r]\arrow[d,"p_\phi"] & Y\arrow[d,"\Delta_Y"]\\
Y \arrow[r,"\id \times \phi"] & Y\times Y.
\end{tikzcd}
\end{equation}
We assume in addition that $\phi_X$ is an equivalence. In this case the $(X\times_{Y}X)$-module $X\times_YZ\times_YX$ is isomorphic to the $\phi$-twisted module ${}^{\phi}(X\times_{Y}X)$ (see \Cref{E:twisted-trace} for the notation), with the isomorphism sending $(x,z,x')\in X\times_YZ\times_Y X$ to $(\phi(x),x')\in {}^{\phi}(X\times_YX)$.
Then \eqref{eq:trace-convolution-horocycle-diagram} becomes 
\begin{equation}\label{eq:trace-convolution-horocycle-diagram-special} 
\xymatrix{
X\times_Y\mL_\phi(Y)\ar^-{\delta_0}[r]\ar_{q}[d] & X\times_YX\\
\mL_\phi(Y) & 
}
\end{equation}

\begin{remark}\label{rem: first projection vs second projection}
We note that the composed map $\pr_2\circ \delta_0=\pr_1:X\times_Y\mL_\phi(Y)\to X\times_YX\to X$ is the natural projection to the first factor, while $\pr_1\circ \delta_0=\phi\circ \pr_1: X\times_Y\mL_\phi(Y)\to X\times_YX\to X$.
\end{remark}

It follows that if we let
\begin{equation}\label{eq: phi structure given by pushforward}
\phi=\phi\vertl: \der(X\times_YX)\to \der(X\times_YX),
\end{equation} 
which is a monoidal automorphism, then $\der(X\times_YZ\times_YX)$ as $\der(X\times_YX)$-bimodule is identified with ${}^\phi\der(X\times_YX)$. 

\begin{corollary}\label{trace.phi.convolution.fully.faithful}
Under the same assumption as in \Cref{bimodule.geom.trace.fully.faithfulness.convolution.cat} and given $\phi_X, \phi_Y$, \eqref{eq: meaning of equiv in higher sense} as above with $\phi_X$ an automorphism, there is a canonical factorization
\[
\begin{tikzcd}
\der(X\times_{Y} X) \arrow[d]\arrow[r,"(\delta_{0})\horizl"] & \der(X \times_{X\times X} (X\times_{Y} X))\arrow[d,"q\vertl"] \\
\trg(\der(X\times_{Y} X),\phi) \arrow[r] & \der(\mathcal{L}_{\phi} Y)
\end{tikzcd}
\]
with the lower horizontal arrow is fully faithful. The essential image is generated under colimits by the image of $q\vertl\circ \delta_0\horizl$.
\end{corollary}

We record the observation fact for later purpose.
\begin{lemma}\label{lem: loop of f in verti}
Assume that $f,\Delta_{X/Y}\in \verti$ and assume that $f$ is $\phi$-equivariant with respect to automorphisms $\phi_X$ and $\phi_Y$ as above. Then $\mL_\phi(f): \mL_\phi(X)\to \mL_\phi(Y)$ belongs to $\verti$.
\end{lemma}
\begin{proof}
The map $\mL_\phi(f)$ can be factors as
\[
X\times_{1\times\phi, X\times X}X\xrightarrow{\Delta_{X/Y}\times\Delta_{X/Y}} (X\times_YX)\times_{1\times\phi, X\times X}(X\times_YX) \cong\mL_\phi(Y)\times_{Y\times Y}(X\times X)\xrightarrow{\id\times(f\times f)} \mL_\phi(Y).
\] 
\end{proof}

Another case we need to consider is $Z=W_1\times W_2$ with $g_i:W_i\to Y$ two maps in $\bfC$. In this case, \[
Z\times_{Y\times Y}Y=W_1\times_YW_2, \quad Z\times_{Y\times Y}(X\times X)=(W_1\times_YX)\times(X\times_YW_2),
\]
We denote:
\begin{equation}\label{eq:geometric-relative-tensor-product}
\der(W_1\times_YX)\otimes_{\der(X\times_YX)}^{\geo}\der(X\times_YW_2):=\trg(\der(X\times_YX),\der(Z\times_{Y\times Y}(X\times X)),
\end{equation}
which is the geometric analogue of the relative tensor product.

\begin{corollary}\label{fully.faithfulness.geometric.tensor.product}
   Under the same assumption as in \Cref{bimodule.geom.trace.fully.faithfulness.convolution.cat}, we have a canonical square
   \[
   \xymatrix{
   \der((W_1\times_YX)\times(X\times_Y W_2))\ar^-{(\id_{W_1}\times \Delta_X\times \id_{W_2})\horizl}[rr]\ar[d] && \der(W_1\times_YX\times_YW_2)\ar^{(\id_{W_1}\times f\times \id_{W_2})\vertl}[d]\\
   \der(W_1\times_YX)\otimes_{\der(X\times_YX)}^{\geo}\der(X\times_YW_2)\ar[rr] && \der(W_1\times_YW_2)
   }
   \]
   with the bottom functor fully faithful. The essential image is generated under colimits by the image of $(\id_{W_1}\times f\times\id_{W_2})\vertl\circ ({\id_{W_1}\times\Delta_X\times\id_{W_2}})\horizl$.
\end{corollary}
Again, there is no assumption on $g_1$ and $g_2$.

Let us come back to the set-up of \Cref{bimodule.geom.trace.fully.faithfulness.convolution.cat}. By assumption, we have an adjoint pair in $\lincat_\La$
\[
\mathrm{CH}:=q\vertl\circ(\delta_0)\horizl:  \der(X\times_YZ\times_YX) \rightleftarrows \der(Y\times_{Y\times Y}Z) :  (\delta_0)\horizr\circ q\vertr:= \mathrm{HC}
\]
Then that  $\trg(\der(X\times_YX),\der(X\times_YZ\times_YX)) \to \der(Y\times_{Y\times Y}Z)$ is an equivalence if and only if the image of $\mathrm{CH}$ generates $\der(Y\times_{Y\times Y}Z)$ under colimits, if and only if $\mathrm{HC}$ is conservative. Sometimes, this can be checked by considering the composition $\mathrm{CH}\circ \mathrm{HC}$. To compute this monad,
we make the following further assumptions, in addition to assumptions as in \Cref{bimodule.geom.trace.fully.faithfulness.convolution.cat}. 
\begin{itemize}
\item $\Delta_Y: Y\to Y\times Y$ and $\pi_Y: Y\to \pt$ belong to $\horiz$, and $\Delta_Y$ and $\Delta_{Y/Y\times Y}$ belong to $\verti$;
\item $(f:X\to Y)\in \bfC_{\horiz}$ and there is some integer $m$  for any base change $g:S\to T$ of $f$, we have $g\vertr=g\horizl[m]$;
\item There is some integer $n$ such that for any base change $g: S\to T$ of $\Delta_X: X\to X\times X$, $g\horizr=g\vertl[n]$.
\end{itemize}
Note that by the first assumption, $\mL(Y)=Y\times_{Y\times Y}Y$ is an algebra object in $\corr(\bfC)_{\verti;\horiz}$, and $Y\times_{Y\times Y}Z$ is a left $\mL(Y)$-module. 
Therefore $\der(\mL(Y))$ acts on $\der(Y\times_{Y\times Y}Z)$ by convolution. We use $\star$ to denote the convolution product as usual. Let
\[
\mS:=(\mL f)\vertl \La_{\mL(X)}[m+n] \in \der(\mL(Y)).
\]

\begin{lemma}\label{lem: CH HC=conv with Spr}
Assumptions are as in \Cref{bimodule.geom.trace.fully.faithfulness.convolution.cat}. Then under further assumptions as above, we have
\[
\mathrm{CH}\circ \mathrm{HC}\cong \mS\star(-).
\]
\end{lemma}
\begin{proof}
By our assumption and base change, the functor $\mathrm{CH}\circ \mathrm{HC}$ is equivalent to $[m+n]$-shift of the horizontal pullback followed vertical pushfoward along the following correspondence
\[
\xymatrix{
X\times_{X\times X}X\times_{Y\times Y}Z\ar[d]\ar[r]&X\times_{Y\times Y}Z\ar[d]\ar[r]& Y\times_{Y\times Y}Z \\
X\times_{Y\times Y}Z\ar[d]\ar[r]&X\times_YZ\times_YX &     \\
Y\times_{Y\times Y}Z.  & &   
}
\]
We may factor it as compositions of correspondences 
\[
\xymatrix{
X\times_{X\times X}X\times_{Y\times Y}Z\ar[rr]\ar[d]&& X\times_{X\times X}X\times Y\times_{Y\times Y}Z\ar^-{\pi_{X\times_{X\times X}X}\times \id}[rr]\ar^-{\mL(f)\times\id}[d]&& Y\times_{Y\times Y}Z\\
Y\times_{Y\times Y}Y\times_{Y\times Y}Z\ar^-{\id\times\Delta_Y\times\id}[rr]\ar_{\id\times\Delta_Y\times\id}[d]&& Y\times_{Y\times Y}Y\times Y\times_{Y\times Y}Z &&\\
Y\times_{Y\times Y}Z, &&&&
}\]
from which the lemma follows.
\end{proof}

\subsubsection{The geometric trace and relative resolutions}\label{sec:relative-resolutions-geometric}
Now we prove \Cref{bimodule.geom.trace.fully.faithfulness.convolution.cat}. In fact, (to save notations) we will prove a slightly general statement.
We consider the geometric trace for pair $(X_\bullet, Q_\bullet)$ with $X_\bullet$ a Segal object in $\bfC$ and $Q_\bullet$ a left $(X_\bullet\times X_\bullet)$-module (or equivalently an $X_\bullet$-bimodule) as in \Cref{sec:monads-correspondences-review}. They give  objects in the category $\BMod(\corr(\bfC)_{\verti;\horiz})$, which roughly speaking consist of an algebra $X_1\in \alg(\corr(\bfC)_{\verti;\horiz})$ whose multiplication and unit maps are of the form 
\begin{equation}\label{eq:multiplication-unt-X_1}
\begin{tikzcd}
 X_1\times_{X_0} X_1 \arrow[r,"\eta"]\arrow[d,"m"]& X_1 \times X_1\\
 X_1 & 
\end{tikzcd}, \quad
\begin{tikzcd}
 X_0\arrow[r,"\pi_{X_0}"]\arrow[d,"u"]& \pt\\
 X_1 & 
\end{tikzcd},
\end{equation}
and an $X_1$-bimodule $Q\in {}_{X_1}\BMod_{X_1}(\corr(\bfC)_{\verti;\horiz})$ whose action maps are of the form
\begin{equation}\label{eq:action-X_1-on-Q}
\begin{tikzcd}
 X_1\times_{X_0} Q \arrow[r,"\xi_l"]\arrow[d,"a_l"]& X_1 \times Q\\
 Q & 
\end{tikzcd},\quad
\begin{tikzcd}
 Q\times_{X_0} X_1 \arrow[r,"\xi_r"]\arrow[d,"a_r"]& Q \times X_1\\
 Q & 
\end{tikzcd}.
\end{equation}
Here we require the simplicial object $X_\bullet$ is as in \Cref{rem:segal.objects.morphisms.vert.horiz} so that $m,u,a_l,a_r\in \bfC_{\verti}$ and $\eta, \xi_l, \xi_r\in \bfC_{\horiz}$.

We can then consider the geometric trace
\[
\trg(\der(X_1),\der(Q)) = \der(\tr(X_1,Q))\cong |\der(\HH(X_1,Q)_{\bullet})|
\]
defined in the previous section. On the other hand, the extra structure on the algebra and module allows one to construct a variant of the geometric trace.

In the monoidal category ${}_{X_0}\BMod_{X_0}(\bfC^{\op,\sqcup})^{\op}$  we consider the Bar complex of the algebra object $X_1$, which we denoted by
$\barcons^{X_0}(X_1)_\bullet$. 
Under the lax monoidal functor ${}_{X_0}\BMod_{X_0}(\bfC^{\op,\sqcup})^{\op}\to \corr(\bfC)$, it gives a simplicial object in $\corr(\bfC)_{\verti;\horiz}$ (in fact in $\bfC_{\verti}$), denoted by the same notation. 
The action of $(X_1\times X_1)\times Q\to Q$ by right and left multiplication gives
\[
\barcons^{X_0}(X_1)_\bullet\otimes Q=X_\bullet\times_{X_0\times X_0}(X_1\times X_1)\times Q\to X_\bullet\times_{X_0\times X_0}Q=:\HH^{X_0}(X_1,Q)_\bullet
\]
which is $(X_1\otimes X_1)$-bilinear and therefore induces
\[
\barcons^{X_0}(X_1)_\bullet\otimes_{X_1\otimes X_1} Q\to \HH^{X_0}(X_1,Q)_\bullet.
\]

The lax monoidal functor ${}_{X_0}\BMod_{X_0}(\bfC^{\op,\sqcup})^{\op}\to \corr(\bfC)$ also induces a natural map of simplicial objects
\[
\barcons(X_1)_\bullet\to \barcons^{X_0}(X_1)_\bullet
\]
in $\corr(\bfC)_{\verti;\horiz}$. It follows that we obtain a map of simplicial objects in $\corr(\bfC)_{\verti;\horiz}$
\begin{equation}\label{eq:hochschild-to-relative-hochschild}
\delta_\bullet: \HH(X_1,Q)_\bullet=\barcons(X_1)_\bullet\otimes_{X_1\otimes X_1}Q\to \barcons^{X_0}(X_1)_\bullet\otimes_{X_1\otimes X_1} Q\to \HH^{X_0}(X_1,Q)_\bullet,
\end{equation}
which is given on each level $n\geq 0$ by the horizontal arrow
\[
X_n\times_{X_0\times X_0} Q\xleftarrow{\id}X_n\times_{X_0\times X_0} Q
\xrightarrow{\delta_n}
X_1^{n}\times Q.
\]

Now we define the \textit{$X_0$-relative Hochschild homology} of $Q$ as
\[
\tr^{X_0}(X_1,Q)=|\HH^{X_0}(X_1,Q)_\bullet|\in \mP(\corr(\bfC)_{\verti;\horiz}),
\]
and define the \textit{$X_0$-relative geometric trace} of $\der(Q)$ as the geometric realization in $\lincat_\La$ 
\[
\trg^{X_0}(\der(X_1),\der(Q)) := \der(\tr^{X_0}(X_1,Q))\cong |\der(\HH^{X_0}(X_1,Q)_\bullet)|.
\]
Then \eqref{eq:hochschild-to-relative-hochschild} gives a functor 
\[
\delta\horizl \colon \trg(\der(X_1),\der(Q)) \to \trg^{X_0}(\der(X_1),\der(Q)),
\]
which fits into a commutative diagram
\begin{equation}\label{eq:geo-trace-to-relative-trace}
\begin{tikzcd}
\der(Q) \arrow[d]\arrow[r,"(\delta_0)\horizl"] & \der(X_0 \times_{X_0\times X_0} Q)\arrow[dd] \\
\tr(\der(X_1),\der(Q))\arrow[d] &\\
\trg(\der(X_1),\der(Q)) \arrow[r,"\delta\horizl"] & \trg^{X_0}(\der(X_1),\der(Q)).
\end{tikzcd}
\end{equation}

\begin{proposition}\label{prop:trace.segal.fully.faithful}
We use  notations of \eqref{eq:multiplication-unt-X_1} and \eqref{eq:action-X_1-on-Q}. 
Assume that
\begin{itemize}
    \item\label{prop:trace.segal.fully.faithful-1} $m$  are in $\bfC_{\mathrm{VR}}$, and the diagonal $\Delta_{X_0}\colon X_0 \rightarrow X_0\times X_0$ is in $\bfC_{\mathrm{HR}}$,
    \item\label{prop:trace.segal.fully.faithful-2} $a_l, a_r$ are in $\bfC_{\mathrm{VR}}$. 

\end{itemize}

Then the functor $\delta\horizl$ from \eqref{eq:geo-trace-to-relative-trace}
is fully faithful. The essential image is generated under colimits by the image of $\der(Q)\xrightarrow{(\delta_0)\horizl} \der(X_0 \times_{X_0\times X_0} Q) \rightarrow \trg^{X_0}(\der(X_1),\der(Q))$. 
\end{proposition}
\begin{proof}
Passing to right adjoints gives a natural transformation
\[
(\delta^{\bullet})\horizr\colon   \der(\HH^{X_0}(X_1,Q)^\bullet) \rightarrow \der(\HH(X_1,Q)^\bullet)
\]
of cosimplicial categories. To prove $\delta\horizl$ is fully faithful we will use \Cref{thm:Beck-Chevalley-descent}. 
The first step is to verify each of the underlying semi-cosimplicial categories satisfies the Beck-Chevalley conditions.

\begin{lemma}\label{lem:trace.segal.beck.chevalley}
Under the assumptions of \Cref{prop:trace.segal.fully.faithful}, the underlying semi-cosimplicial categories obtained from the semi-simplicial categories $\der(\HH(X_1,Q)_{\bullet})$ and
$\der(\HH^{X_0}(X_1,Q)_{\bullet})$ by passing to right adjoints satisfy the Beck-Chevalley conditions. 
\end{lemma}
\begin{proof}
As all morphisms of the semi-simplicial object $\HH^{X_0}(X_1,Q)^{\bullet}$ are the base change of $m,a_l,a_r$, the statement for the semi-simplicial categories 
$\der(\HH^{X_0}(X_1,Q)_{\bullet})$ follows directly from our assumption.

It is left to deal with $\der(\HH(X_1,Q)^\bullet)$. For every face map $\alpha\colon [m]\rightarrow [n]$, we have the diagram
\begin{equation}\label{eq:diagram-left-adjoint-beck-chev-corr}
\begin{tikzcd}
X_1^{m}\times Q & \arrow[l,dashed,"d^0_{m}"'] 
X_1^{m+1}\times Q \\
X_1^{n}\times Q 
 \arrow[u,dashed] &
X_1^{n+1}\times Q  \arrow[u,dashed] \arrow[l,dashed,"d^0_{n}"']
\end{tikzcd}
\end{equation}
in $\corr(\bfC)_{\verti;\horiz}$ (see \eqref{E: morphism-in-corr} and after for notations).
We need to show that the induced diagram
\begin{equation}\label{eq:diagram-left-adjoint-beck-chev}
\begin{tikzcd}
\der(X_1^{m}\times Q) \arrow[r,"d^0_{m}"]\arrow[d] &
\der(X_1^{m+1}\times Q) \arrow[d] \\
\der(X_1^{n}\times Q) 
\arrow[r,"d^0_{n}"] &
\der(X_1^{n+1}\times Q)
\end{tikzcd}
\end{equation}
is left adjointable.

We may assume that $\alpha = d_n: [n]\to [n+1]$, that is, $\alpha(i) = i$ for all $i$. (The proof is similar in all other cases.) Then the diagram \eqref{eq:diagram-left-adjoint-beck-chev-corr} is explicitly given by
\[
\begin{tikzcd}
X_1^{n}\times Q &   X_1^{n}\times  Q\times_{X_0}X_1 \arrow[r,"\xi_{r,n}"]\arrow[l,"a_{r,n}"'] & X_1^{n+1}\times Q\\
X_1^{n} \times X_1 \times_{X_0} Q\arrow[u,"a_{l,n}"']\arrow[d,"\xi_{l,n}"]  & X_1^{n}\times X_1\times_{X_0} Q\times_{X_0} X_1\arrow[r,"\tilde\xi_r"]\arrow[l,"\tilde a_r"']\arrow[u,"\tilde a_l"']\arrow[d,"\tilde \xi_l"]  & X_1^{n+1}\times X_1\times_{X_0}Q\arrow[u,"a_{l,n+1}"']\arrow[d,"\xi_{l,n+1}"] \\
X_1^{n+1}\times Q &   X_1^{n+1}\times  Q\times_{X_0}X_1 \arrow[r,"\xi_{r,n+1}"]\arrow[l,"a_{r,n+1}"'] & X_1^{n+2}\times Q,
\end{tikzcd}
\]
as a diagram in $\bfC$, where $a_{l,n}=\id_{X_1^n}\times a_l$, etc. Note that all squares are Cartesian in $\bfC$.

We have $d_{k}^{0}=(\xi_{l,k})\horizr \circ (a_{l,k})\vertr$ (for $k=n,n+1$) with the left adjoint $(a_{l,k})\vertl\circ \xi_{l,k}\horizl$,
and the vertical arrows in \eqref{eq:diagram-left-adjoint-beck-chev} are given by $(\xi_{l,k})\horizr\circ (a_{l,k})\vertr$.
Left adjointability of \ref{eq:diagram-left-adjoint-beck-chev} then means that the natural map
\[
(a_{l,n+1})\vertl\circ (\xi_{r,n+1})\horizl\circ (\xi_{l,n+1})\horizr\circ (a_{l,n+1})\vertr \rightarrow 
(\xi_{l,n})\horizr\circ (a_{l,n})\vertr\circ (a_{r,n})\vertl\circ (\xi_{r,n})\horizl
\]
is an equivalence. It is enough to show that the Beck-Chevalley maps
\begin{align}
    (\tilde a_r)\vertl\circ (\tilde a_l)\vertr &\to (a_{l,n})\vertr\circ (a_{r,n})\vertl, \quad (\xi_{r,n+1})\horizl\circ (\xi_{l,n+1})\horizr \to
    (\tilde\xi_l)\horizr\circ (\tilde\xi_r)\horizl\label{eq-trace.BC.proof.upper}\\
    (\tilde \xi_r)\horizl\circ (a_{l,n+1})\vertr & \to (\tilde a_l)\vertr\circ (\xi_{r,n})\horizl  , \quad
    (a_{r,n+1})\vertl\circ (\tilde \xi_l)\horizr \to  (\xi_{l,n})\horizr\circ (\tilde a_r)\vertl\label{eq-trace.BC.proof.lower}
\end{align}
are equivalences. 
As the maps $a_{l,k}, a_{r,k}, \tilde a_l, \tilde a_r$ are base change of $a_l$ and $\xi_{l,k}, \xi_{r,k},\tilde\xi_l,\tilde\xi_r$ are base change of $\Delta_{X_0}$, the desired equivalences follow from our assumptions.
\end{proof}

We continue to prove \Cref{prop:trace.segal.fully.faithful}. Passing to the right adjoint of \eqref{eq:geo-trace-to-relative-trace} gives
\begin{equation*}
\begin{tikzcd}
\der(Q) & \arrow[l,"(\delta_0)\horizr"'] \der(X_0 \times_{X_0\times X_0} Q) \\
\trg(\der(X_1),\der(Q))\arrow[u] & \arrow[l,"\delta\horizr"']  \trg^{X_0}(\der(X_1),\der(Q)) \arrow[u].
\end{tikzcd}
\end{equation*} 
with vertical arrows monadic (by \Cref{lem:trace.segal.beck.chevalley}). Let $T$ denote the monad corresponding to the cosimplicial category $\der(\HH^{X_0}(X_1,Q)^{\bullet})$ and by $V$ the monad corresponding to $\der(\HH(X_1,Q)^{\bullet})$. Then to show that $\delta\horizl$ is fully faithful it is enough to show that the natural map 
\[
 V \rightarrow (\delta_0)\horizr \circ T \circ (\delta_{0})\horizl,
\]
is an equivalence. The monad $T$ is given by $(d_0)\vertl\circ (d_1)\vertr$ with
\[
d_1,d_0\colon  X_1\times_{X_0\times X_0} Q \rightarrow X_0 \times_{X_0\times X_0} Q.
\]
Recall that the map $d_1$ is induced by the left action of $X_1$ on $Q$ by
\[
X_1\times_{X_0 \times X_0} Q \simeq X_0\times_{X_0\times X_0} (X_1\times_{X_0} Q) \rightarrow X_0\times_{X_0\times X_0} Q
\]
Likewise, the map $d_0$ is induced by the right action via
\[
X_1\times_{X_0 \times X_0} Q \simeq X_0\times_{X_0\times X_0}(Q\times_{X_0} X_1) \rightarrow X_0\times_{X_0\times X_0} Q.
\]
The monad $V$ is given by 
\[
V \simeq (a_r)\vertl\circ (\xi_r)\horizl\circ (\xi_l)\horizr\circ (a_l)\vertr.
\]

These maps fit into a commutative diagram diagram in $\bfC$
\begin{equation}\label{eq:trace-diagram-two-monads}
\begin{tikzcd}
X_1\times Q &
X_1\times_{X_0} Q \arrow[r,"a_l"]\arrow[l,"\xi_l"'] &
 Q\\
 Q\times_{X_0} X_1 \arrow[d,"a_r"] \arrow[u,"\xi_r"'] &
 X_1\times_{X_0\times X_0} Q
 \arrow[r,"d_1"] \arrow[d,"d_0"]
 \arrow[u,"\zeta"'] \arrow[l,"\chi"']
 & 
 X_0\times_{X_0 \times X_0} Q\arrow[u,"\delta_0"']\\
 Q & 
  X_0\times_{X_0 \times X_0} Q
  \arrow[l,"\delta_0"']
 & 
\end{tikzcd}
\end{equation}
with all squares being Cartesian. Then it is enough to show that the natural maps
\begin{align*}
    (\xi_r)\horizl\circ (\xi_l)\horizr& \rightarrow
    \chi\horizr\circ \zeta\horizl,\quad
    \zeta\horizl\circ (a_l)\vertr\to (d_1)\vertr\circ(\delta_0)\horizl
    ,\quad
    (a_r)\vertl \circ \chi\horizr  \rightarrow 
    (\delta_{0})\horizr\circ (d_0)\vertl
\end{align*}
are equivalences, which hold by our assumptions. 
\end{proof}

Now we specialize the above discussions to the case $X_1=X\times_YX$ and $Q=X\times_YZ\times_YX$ as in \Cref{bimodule.geom.trace.fully.faithfulness.convolution.cat}.
In this case, the relative Hochschild complex has a simple interpretation. Consider the fiber product
\[
\begin{tikzcd}
 Y\times_{Y\times Y}Z \arrow[r]\arrow[d] & Y\arrow[d,"\Delta_Y"]\\
 Z \arrow[r,"g=g_1\times g_2"]& Y\times Y
\end{tikzcd}
\]
in $\bfC$ and consider the map $q=f\times\id_Z: X\times_{Y\times Y}Z\rightarrow  Y\times_{Y\times Y}Z$.
\begin{lemma}\label{lemma:cech-nerve-relative-complex}
There is a canonical equivalence of simplicial objects in $\bfC_{\verti}$.
\[
\HH^{X}(X\times_{Y} X,X\times_{Y} Z\times_{Y} X)_\bullet\simeq  X_\bullet \times_{Y\times Y}Z
\]
where the right hand side is the \v{C}ech nerve of $q:X\times_{Y\times Y}Z\rightarrow  Y\times_{Y\times Y}Z$. Under the identification, the map $\delta_0$ from \eqref{eq:hochschild-to-relative-hochschild} is the horizontal map in \eqref{eq:trace-convolution-horocycle-diagram}.
\end{lemma}
\begin{proof}
The construction of the left hand side is natural in $X$ and applying it to the identity map $Y\rightarrow Y$ gives the right hand side. Thus, $f\colon X\rightarrow Y$ induces an augmentation $\HH^{X}(X\times_{Y} X,X\times_{Y} Z\times_{Y} X)_\bullet$ of the corresponding simplicial object. In order to identify this augmented simplicial object with the \v{C}ech nerve of $q$, can use the characterization \cite[Proposition 6.1.2.11]{Lurie.higher.topos.theory} as it is easy to check that the necessary squares are pullbacks. 
\end{proof}

\begin{proof}[Proof of \Cref{bimodule.geom.trace.fully.faithfulness.convolution.cat}]
Only fully faithfulness requires a proof.
Consider the augmented simplicial category associated to the \v{C}ech nerve of $q$. As $f\in \bfC_{\mathrm{VR}}$ so is the map $X\times_{Y\times Y}Y \rightarrow Y\times_{Y\times Y} Z$. Using \Cref{lemma:cech-nerve-relative-complex} we identify the \v{C}ech nerve of this map with the relative Hochschild complex. By passing to right adjoints and using \cite[Corollary 4.7.5.3]{Lurie.higher.algebra} we get a fully faithful functor 
\begin{equation}\label{eq:nerve-convolution-categories}
\left|\der(\HH^{X}(X\times_{Y} X, Z\times_{Y\times Y}(X\times X))_{\bullet})\right| \rightarrow \der(Z\times_{Y\times Y}Y).
\end{equation}
Composition with the fully faithful functor from \Cref{prop:trace.segal.fully.faithful} gives the desired functor. The essential image of \eqref{eq:nerve-convolution-categories} is generated by the image of $q_{\vertl}$ so the description of the essential image follows. That $\proj_{\trg}$ is continuous is also clear.
\end{proof}

We record some observations for later purposes. Let $(f: X\to Y)\in\bfC$ as in \Cref{rem:segal.objects.morphisms.vert.horiz}. Let $Z\leftarrow C\rightarrow Z'$ be a morphism in $\corr(\bfC_{/Y\times Y})_{\verti;\horiz}$, i.e. all $Z,Z',C$ are equipped with morphisms to $Y\times Y$ and $C\to Z$ and $C\to Z'$ are $(Y\times Y)$-morphisms in $\bfC$. 

Let $X_\bullet$ be as above and let $Q'=X\times_YZ'\times_YX$ and $Q=X\times_YZ\times_YX$. Then the following diagram is commutative
\begin{equation}\label{eq-functoriality-trace-to-relative-trace}
    \xymatrix{
    \der(Q')\ar[r]\ar[d]&\tr(\der(X_1),\der(Q'))\ar[r] \ar[d] & \trg(\der(X_1),\der(Q')) \ar[r] \ar[d]  & \der(Y\times_{Y\times Y}Z')\ar[d]\\
    \der(Q)\ar[r]&\tr(\der(X_1),\der(Q))\ar[r]  & \trg(\der(X_1),\der(Q)) \ar[r]   & \der(Y\times_{Y\times Y}Z).
    }
\end{equation}

\subsubsection{Comparison between geometric and ordinary traces}\label{SS:comparison-geom-tr-ord-tr}
In practice, we need to compare the geometric trace defined and studied as above with the ordinary traces reviewed in \Cref{sec-hoch-homology}. The easiest situation has been discussed in \Cref{rmk:trace-to-geometric-trace}. 
On the other hand, the monadicity of the simplicial objects in \Cref{lem:trace.segal.beck.chevalley} can be used to compare the usual trace and the geometric trace in other situations.
Recall notations in \eqref{eq:multiplication-unt-X_1}.

\begin{proposition}\label{prop-comparison-usual-trace-geo-trac-presemirigid}
Let $X_\bullet$ be a Segal object satisfying assumptions of \Cref{prop:trace.segal.fully.faithful}.  In addition, assume that 
\begin{enumerate}
\item\label{item-adjoin.boundary}  The exterior tensor product 
\[
\boxtimes_{\der(\pt)}:\der(X_1)\otimes_{\der(\pt)}\der(X_1)\to \der(X_1\times X_1)
\]
is fully faithful and admits a continuous right adjoint $\boxtimes_{\der(\pt)}^R$ (see \Cref{R:Kunneth-type-formula}). 
\item\label{item-adjoin.boundary-2} The object $(\eta\horizr\circ m\vertr\circ u\vertl)(\La_{X_0})$ belongs to  $\der(X_1)\otimes_{\der(\pt)} \der(X_1)\subset \der(X_1\times X_1)$. 
\end{enumerate}
Then the product  $\der(X_1)\otimes_{\der(\pt)} \der(X_1)\to\der(X_1)$ admits a $\der(X_1)\otimes_{\der(\pt)} \der(X_1)^{\rev}$-linear right adjoint. 
\end{proposition}
Note that the two assumptions automatically hold if $\boxtimes_{\der(\pt)}$ is an equivalence.
\begin{proof}
The multiplication map is given by the composition $m\vertl\circ \eta\horizl \circ \boxtimes$, with the continuous right adjoint given by $\boxtimes^R\circ\eta\horizr\circ m\vertr$. We need to show that that this right adjoint is a $\der(X_1)$-bimodule homomorphism. To see that it is a left $\der(X_1)$-module morphism (the case of right $\der(X_1)$-module structure is similar), consider the following diagram
\[
\xymatrix{
\der(X_1)\otimes_{\der(\pt)} \der(X_1) \ar^-{\id\otimes (\eta\horizr\circ m\vertr)}[rr]\ar_{\boxtimes}[d] && \der(X_1)\otimes_{\der(\pt)}\der(X_1\times X_1) \ar^{\boxtimes}[d] \ar^-{\id\otimes \boxtimes^R}[r] & \der(X_1)\otimes_{\der(\pt)} \der(X_1)\otimes_{\der(\pt)} \der(X_1)\ar[d]\\
\der(X_1\times X_1) \ar^-{(\id\times\eta)\horizr\circ(\id\otimes m)\vertr}[rr]\ar_{m\vertl\circ\eta\horizl}[d] && \der(X_1\times X_1\times X_1)\ar^{(m\times\id)\vertl\circ(\eta\times\id)\horizl}[d] & \der(X_1\times X_1)\otimes_{\der(\pt)} \der(X_1)\ar[d]\\
\der(X_1)\ar^-{\eta\horizr\circ m\vertr}[rr] &&\der(X_1\times X_1)\ar^-{\boxtimes^R}[r]& \der(X_1)\otimes_{\der(\pt)} \der(X_1). 
}
\]
By our assumption on $\der$ and the assumption $m\in\mathrm{VR}$ and $\eta\in \mathrm{HR}$, the left upper square is commutative. By \Cref{lem:trace.segal.beck.chevalley}, the left lower square is commutative. In other words, the functor $\eta\horizr\circ m\vertr: \der(X_1)\to \der(X_1\times X_1)$ is a left $\der(X_1)$-module homomorphism. 
Then together with Assumption \eqref{item-adjoin.boundary} and \eqref{item-adjoin.boundary-2}, we see that the essential images of the functor $\eta\horizr\circ m\vertr$ belong to $\der(X_1)\otimes_{\der(\pt)} \der(X_1)\subset \der(X_1\times X_1)$. Indeed, the unit $\mathbf{1}_{\der(X_1)}$ of $\der(X_1)$ is given by $u\vertl\La_{X_0}$ (see \Cref{rem-conv-product}). Then
$\eta\horizr(m\vertr\mF)=\eta\horizr(m\vertr(\mF\star \mathbf{1}_{\der(X_1)}))=\mF\star\eta\horizr(m\vertr(u\vertr\La_{X_0}))\in \der(X_1)\otimes_{\der(\pt)} \der(X_1)$.
Therefore, the outer square of the above diagram is commutative. (However we do not claim the right square is commutative.) 
This proves the proposition.
\end{proof}

\begin{corollary}\label{cor: semi rigidity of der(X)}
Assumptions are as in \Cref{prop-comparison-usual-trace-geo-trac-presemirigid}. Suppose $\der(X_1)$ is compactly generated. Then $\der(X_1)$ is semi-rigid as a $\der(\pt)$-linear monoidal category (and therefore as a $\La$-linear monoidal category). In addition, it admits a Frobenius structure given by ind-extension of the functor
\[
\Hom(u\vertl\La_{X_0}, -): \der(X_1)^\cpt\to \der(\pt).
\]
\end{corollary}
\begin{proof}This follows from \Cref{prop-comparison-usual-trace-geo-trac-presemirigid} and \Cref{rem-duality-datum-as-plain-cat}.
\end{proof}

Instead of assuming that $\der(X_1)$ is compactly generated, one can impose some other conditions guarantee the semi-rigidity of $\der(X_1)$. We shall not try to give a very general formalism but only discuss a situation that is useful in practice.

\begin{corollary}\label{cor: rigidity of der(X)}
Assumptions are as in \Cref{prop-comparison-usual-trace-geo-trac-presemirigid}. If $(\pi_{X_0})\horizl: \der(\pt)\to \der(X_0)$ and $u\vertl: \der(X_0)\to \der(X_1)$ admit continuous right adjoint, then $\der(X_1)$ is rigid. It admits a Frobenius structure given by
\[
\Hom(u\vertl\La_{X_0},-): \der(X_1)\to \der(\pt).
\]
\end{corollary}
\begin{proof}
By assumption, the unit $(\Delta_{X/Y})\vertl(\pi_{X_0})\horizl : \der(\pt)\to \der(X_1)$ admits continuous right adjoint. Now \Cref{prop-comparison-usual-trace-geo-trac-presemirigid} implies that $\der(X_1)$ is rigid. The last statement again follows from \Cref{rem-duality-datum-as-plain-cat}.
\end{proof}

\begin{remark}\label{rem: pivotal structure}
Suppose we are in the situation as in \Cref{trace.phi.convolution.fully.faithful}. 
Suppose that $\Delta_{X/Y}\in \verti\cap \horiz$ and that for every base change $g$ of $\Delta_{X/Y}$, we have $g\vertr=g\horizl$ (compare with \Cref{rem-additional-base.change.sheaf.theory} \eqref{rem-additional-base.change.sheaf.theory-2}). 
Then by base change, we have the canonical isomorphism
\[
\Hom((\Delta_{X/Y})\vertl\La_X, \mF\star \mG)\cong  \Hom(\La_X, (\pr_1)\vertl
(\mF\otimes \mathrm{sw}\horizl \mG)),\quad \mF,\mG\in \der(X\times_YX).
\]
here $\mathrm{sw}\colon X\times_YX\to X\times_YX$ is the morphism by swapping two factors, and $\pr_1: X\times_YX\to X$ is the first projection, and $\otimes$ is the symmetric monoidal tensor product $\der(X\times_YX)$ as in \eqref{eq:abstract-interior-tensor-product}.

In favorable cases, this will imply that the monoidal structure on $\der(X\times_YX)$ is pivotal (see \Cref{def: pivotal semi-rigid monoidal} for the definition).
\end{remark}

\begin{proposition}\label{prop-comparison-usual-trace-geo-trac}
Let $X_\bullet$, $Q$ be as in the statement of \Cref{prop:trace.segal.fully.faithful}, and keep assumptions as in \Cref{prop-comparison-usual-trace-geo-trac-presemirigid}.
In addition, assume that the exterior tensor product 
\[
\boxtimes_{\der(\pt)}:\der(X_1)\otimes_{\der(\pt)}\der(Q)\to \der(X_1\times Q)
\]
is fully faithful, then the comparison map $\tr(\der(X_1),\der(Q)) \rightarrow \trg(\der(X_1),\der(Q))$ (see \eqref{eqn:trace-to-geometric-comparison-map-general}) is an equivalence. 
\end{proposition}
\begin{proof}
We show that the comparison map \eqref{eqn:trace-to-geometric-comparison-map-general} is an equivalence. Recall that it is induced by the morphism of simplicial objects
\[
\HH(\der(X_1),\der(Q))_\bullet=\der(X_1)^{\otimes \bullet} \otimes_{\der(\pt)} \der(Q) \to \der(X_1^{\bullet}\times Q) =\der(\HH(X_1,Q)_\bullet).
\]
The $0$-th objects of both simplicial objects are given by $\der(Q)$. We use $d_i$ to denote the $i$th face map $d_i: \der(\HH(X_1,Q)_{m+1})\to \der(\HH(X_1,Q)_m)$.

Now we can apply \Cref{L:hochschild-homology-vs-cohomology} to see that the Hochschild complex $\HH(\der(X_1),\der(Q))_\bullet$ is monadic, and that the resulting monad on $\der(Q)$ is given by $d_0\circ\boxtimes\circ \boxtimes^R\circ (d_1)^R$. By \Cref{lem:trace.segal.beck.chevalley}, $\der(\HH(X_1,Q)_\bullet)$  satisfies the Beck-Chevallay conditions, and that the resulting monad on $\der(Q)$ is given by $d_0\circ (d_1)^R$. It remains to show that the two monads are identified which will imply the proposition. Note that $(d_1)^R$ is given by the composition of the bottom functors in the following commutative diagram
\begin{small}
\[
\xymatrix{
&&&&\der(X_1)\otimes_{\der(\pt)} \der(X_1)\otimes_{\der(\pt)} \der(Q) \ar^--{\id\otimes((a_l)\vertl\circ(\xi_l)\horizl\circ\boxtimes_{\La})}[drr] \ar^{\boxtimes_{\La}\otimes\id}[d]&& \\
\der(Q)\ar^-{\mathbf{1}_{X_1}\otimes\id}[rr] \ar@{=}[d]&&\der(X_1)\otimes_{\der(\pt)}\der(Q)\ar^{\boxtimes_{\La}}[d]\ar^-{(\eta\horizr\circ m\vertr)\otimes\id}[rr] \ar^-{(\boxtimes_{\La}^R\circ\eta\horizr\circ m\vertr)\otimes\id}[urr]&& \der(X_1\times X_1)\otimes_{\der(\pt)} \der(Q) \ar^{\boxtimes_{\La}}[d]  && \der(X_1)\otimes_{\der(\pt)} \der(Q)\ar^{\boxtimes_{\La}}[d] \\
\der(Q) \ar^-{(u\times\id)\vertl\circ (\pi_{X_0}\times\id)\horizl}[rr] &&\der(X_1\times Q)\ar^-{(\eta\times\id)\horizr\circ(m\times\id)\vertr}[rr] && \der(X_1\times X_1\times Q) \ar^-{(\id\times a_l)\vertl\circ (\id\times\xi_l)\horizl}[rr] && \der(X_1\times Q).
}
\]
\end{small}

We remark that the commutativity of the triangle follows from the previous discussions in the proof of \Cref{prop-comparison-usual-trace-geo-trac-presemirigid}, and the commutativity of the square below the triangle follows from our assumption on $\der$ and the assumption $m\in\mathrm{VR}$ and $\Delta_{X_0}\in \mathrm{HR}$.  It follows that the essential image of $(d_1)^R$ is contained in $\der(X_1)\otimes_{\der(\pt)}\der(Q)$. Therefore, by our assumption of fully faithfulness of exterior tensor product, we have 
$\boxtimes\circ\boxtimes^R\circ(d_1)^R=(d_1)^R$. In particular, two monads are identified.
\end{proof}

\begin{corollary}\label{cor-fully-faithful-usual-trace-to-loop}
Let $X,Y,Z$ be as in \Cref{bimodule.geom.trace.fully.faithfulness.convolution.cat} and suppose the assumptions in  \Cref{bimodule.geom.trace.fully.faithfulness.convolution.cat} hold. Assume that 
\begin{itemize}
\item for every $W\in\bfC$, the exterior tensor functor 
\[
\der(X\times_YX)\otimes_{\der(\pt)}\der(X\times_YW)\to \der(X\times_YX\times X\times_YW)
\] 
is fully faithful with a continuous right adjoint;
\item $(\id\times\Delta_X\times\id)\horizr\circ(\id\times f\times \id)\vertr\circ(\Delta_{X/Y})\vertl(\La_X)$ belongs to $\der(X\times_YX)\otimes_{\der(\pt)}\der(X\times_YX)$.
\end{itemize}  
Then the canonical map
\[
\tr(\der(X\times_YX), \der(X\times_YZ\times_YX))\to \der(Y\times_{Y\times Y}Z)
\]
is fully faithful. 

If in addition, $Z=W_1\times W_2$ as in \Cref{fully.faithfulness.geometric.tensor.product}, and if the exterior tensor product
\begin{equation}\label{eq-fully-faithful-usual-trace-to-loop}
\der(X\times_Y W_1)\otimes_{\der(\pt)} \der(W_2\times _YX)\xrightarrow{\boxtimes_{\La}}\der(X\times_Y(W_1\times W_2)\times_YX)
\end{equation}
is fully faithful,
then the canonical map
\[
\der(W_1\times_Y X) \otimes_{\der(X\times_Y X)} \der(X\times_Y W_2) \to \der(W_1\times_YW_2)
\]
is fully faithful. 
\end{corollary}
\begin{proof}
The first statement follows directly from \Cref{bimodule.geom.trace.fully.faithfulness.convolution.cat} and  \Cref{prop-comparison-usual-trace-geo-trac}. For the second, note that the functor \eqref{eq-fully-faithful-usual-trace-to-loop} is a $\der(X\times_YX)$-bimodule morphism. By assumption, it is fully faithful.
This implies that
\[
\der(W_2\times _YX)\otimes_{\der(X\times_YX)}\der(X\times_Y W_1)\to \tr(\der(X\times_YX),\der(X\times_Y(W_1\times W_2)\times_YX))
\] 
is fully faithful by \Cref{L:hochschild-homology-vs-cohomology}
\eqref{cor-fully-faithful-categorical-trace} and  \Cref{prop-comparison-usual-trace-geo-trac-presemirigid}. Now the second statement follows from the first.
\end{proof}

\begin{example}\label{ex: relative tensor split case}
Consider the case $Z = W \times X$ for some $g\colon W \rightarrow Y$. Then we have a split augmented simplicial object
\[
\HH(X\times_Y X,(W\times X)\times_{Y\times Y} (X\times X))_\bullet \rightarrow W\times_{Y} X
\]
with the last map given by the action of $X\times_{Y} X$ on $W\times_{Y} X$. In this case, we reduce to the tautological equivalence
\[
\der(W\times_{X} X)\otimes_{\der(X\times_Y X)}\der(X\times_{Y}  X) \xrightarrow{\sim} \der(W\times_Y X).
\]
\end{example}

\begin{example}
When we take $\bfC$ to be the category of (nice) algebraic stacks over $\bC$ and the sheaf theory $\der$ to be the theory of algebraic $\mathrm{D}$-modules, one always has 
\[
\tr(\mathrm{D}(X\times_YX),\mathrm{D}(X\times_YZ\times_YX))=\trg(\mathrm{D}(X\times_YX),\mathrm{D}(X\times_YZ\times_YX))
\] 
as $\mathrm{D}(X)\otimes_\bC\mathrm{D}(Y)\cong \mathrm{D}(X\times Y)$. Therefore, \Cref{trace.phi.convolution.fully.faithful} recovers \cite[Theorem 6.6]{benzvi2009character}. 
In \emph{loc. cit.}, instead of directly considering $\mathrm{D}(\HH^X(X\times_YX, X\times_YX)_\bullet)$, the authors used the relative bar resolution for the monoidal category $\mathrm{D}(X\times_YX)$ and then used integral transforms to embed each level in the resulting simplicial object of this relative resolution fully faithfully into the corresponding level of $\mathrm{D}(\HH^X(X\times_YX, X\times_YX)_\bullet)$. Our method bypasses using the integral transforms, which might fail in other sheaf theoretic content. See \Cref{SS:comparison-geom-tr-ord-tr} for discussions.
\end{example}

For applications to other sheaf theoretic contents, see \Cref{prop: categorical trace for quasi-coherent}, \Cref{prop: control functor image by singular support}, \Cref{thm:trace-constructible}. For more concrete applications, see \Cref{SSS: cat trace spectral side}, \Cref{SS: Finite Deligne-Lusztig induction}, and \Cref{SS: 2nd approach unipotent Langlands category}.

\subsubsection{Functoriality of categorical traces in geometric setting}\label{SS: Functoriality of geometric traces}
Next we discuss functoriality of categorical traces arising from convolution patterns.

We start with the following observation.
Let $f: X\to Y$ be as in \Cref{subsec-trace.convolution.categories}.
Let $Z\dashrightarrow Z'$ be a morphism in $\corr(\bfC_{/Y\times Y})_{\verti;\horiz}$, given as $Z'\leftarrow C \rightarrow Z$.

\begin{lemma}\label{lem-right-adjoint-from-categorical-trace-to-loop-space}
Assumptions are as in \Cref{cor-fully-faithful-usual-trace-to-loop}. 
Then the following diagram is right adjointable
\[
\xymatrix{
\tr(\der(X\times_YX),\der(X\times_YZ\times_YX))\ar[r]\ar[d]&\der(Y\times_{Y\times Y}Z)\ar[d]\\
\tr(\der(X\times_YX),\der(X\times_YZ'\times_YX))\ar[r]&\der(Y\times_{Y\times Y}Z').
}
\]
\end{lemma}
\begin{proof}
Using \Cref{L:hochschild-homology-vs-cohomology} it is enough to prove the right adjointability for the diagram as above but with $\tr(\der(X\times_YX), -)$ replaced by $\HH(\der(X\times_YX),-)_n$. In addition, it is enough to consider the case $n=0$.
I.e., we need to show that the following diagram is right adjointable
\[
\xymatrix{
\der(X\times_YZ\times_YX)\ar[r]\ar[d]&\der(Y\times_{Y\times Y}Z)\ar[d]\\
\der(X\times_YZ'\times_YX)\ar[r]&\der(Y\times_{Y\times Y}Z').
}
\]
But this follows from our assumption of $\der$ and the assumption $\Delta_X\in \bfC_{\mathrm{HR}}$ and $f\in \bfC_{\mathrm{VR}}$. 
\end{proof}

Next we discuss duality of modules arising from the convolution patterns. Consider
\[
(f: X\to Y, g: Z\to Y\times Y),\quad \mbox{and} \quad (f': X'\to Y', g': Z'\to Y'\times Y'),
\] 
where $f: X\to Y$ and $f': X'\to Y'$ are as in \Cref{rem:segal.objects.morphisms.vert.horiz} and $g,g'$ arbitrary.

Let $W\to Y\times Y'$ be a morphism. Let
\[
X_1=X\times_YX, \quad Q=X\times_YZ\times_YX,\quad X'_1=X'\times_{Y'}X',\quad Q'=X'\times_{Y'}Z'\times_{Y'}X',
\]
and let
\[
M=X\times_YW\times_{Y'}X'.
\]
We would like to know when $\der(M)$ is dualizable as a  $\der(X_1)\mbox{-}\der(X'_1)$-bimodule, with the dual of  given by $\der(N)$ where $N=X'\times_{Y'}W\times_YX$. 
We will assume that 
\begin{itemize}
\item $W\to Y'$ and $W\to W\times_{Y'}W$ belong to $\bfC_{\horiz}$; 
\item $W\to Y$ and $W\to W\times_YW$ belong to $\bfC_{\verti}$.
\end{itemize}
Then we have the morphism 
\[
u_{\geo}: Y'\dashrightarrow W\times_YW,\quad e_{\geo}: W\times_{Y'}W\dashrightarrow Y
\] 
in $\corr(\bfC)_{\verti;\horiz}$ given by $ Y' \leftarrow W\to W\times_YW$ and $W\times_{Y'}W  \leftarrow W\to Y$ respectively. They induce
\begin{equation}\label{eq-pre-geo-unit}
\xymatrix{
&  \der(X'_1)\ar^{\der(\id\times u_{\geo}\times \id)}[d] \\
\der(N)\otimes_{\der(X_1)}\der(M)\ar^-{}[r] & \der(X'\times_{Y'}W\times_YW\times_{Y'}X').
}
\end{equation}
and
\begin{equation}\label{eq-pre-geo-evaluation}
\xymatrix{
\der(M)\otimes_{\der(X'_1)}\der(N)\ar^-{}[r]\ar_-{e}[dr]& \der(X\times_YW\times_{Y'}W\times_YX)\ar^{\der(\id\times e_{\geo}\times \id)}[d]\\
& \der(X_1).
}
\end{equation}
Here $e$ is defined to be the composition.

\begin{lemma}\label{lem-geo-duality-I}
Assumptions are as in \Cref{cor-fully-faithful-usual-trace-to-loop} and assume that $f: X\to Y$ and $f':X'\to Y'$ satisfy assumptions of \emph{loc. cit}.
If the vertical morphism in \eqref{eq-pre-geo-unit} factors through a $\der(X'_1)$-bimodule morphism 
\[
u: \der(X'_1)\to \der(N)\otimes_{\der(X_1)}\der(M)
\]
(e.g. if $\der(N)\otimes_{\der(X_1)}\der(M)\to \der(X'\times_{Y'}W\times_YW\times_{Y'}X')$ is an equivalence), then $u$ and $e$ from \eqref{eq-pre-geo-evaluation} give the duality datum of $\der(M)$ as a $\der(X_1)\mbox{-}\der(X'_1)$-module.
\end{lemma}
\begin{proof}
Write $R=X'\times_{Y'}W\times_YW\times_{Y'}X'$ and $S=X\times_{Y}W\times_{Y'}W\times_{Y}X$ for simplicity. Note that (using \eqref{eq-functoriality-trace-to-relative-trace}) we have the following commutative diagram
\begin{equation*}
    \xymatrix{
    \der(X'_1)\otimes_{\der(X'_1)}\der(N)\ar_-{u\otimes\id}[dd] &\ar@{=}[l]\der(X'_1)\otimes_{\der(X'_1)}\der(N)\ar[d] \ar_-{\Cref{ex: relative tensor split case}}^-{\cong}[r]& \der(N)=\der(X'\times_{Y'}Y'\times_{Y'}W\times_YX)\ar^{\der(\id\times u_{\geo}\times \id\times\id)}[dd]\\
    &\der(R)\otimes_{\der(X'_1)}\der(N)\ar[dr]&\\
    \der(N)\otimes_{\der(X_1)}\der(M)\otimes_{\der(X'_1)}\der(N) \ar[ur]\ar_-{\id\otimes e}[dd]\ar[dr]&& \der(X'\times_{Y'}W\times_{Y} W\times_{Y'} W\times_{Y}X)\ar^{\der(\id\times\id\times e_{\geo}\times \id)}[dd]\\
    & \der(N)\otimes_{\der(X_1)}\der(S) \ar[ur]\ar[d]&\\  \der(N)\otimes_{\der(X)}\der(X)&\ar@{=}[l]\der(N)\otimes_{\der(X_1)}\der(X_1)\ar_-{\Cref{ex: relative tensor split case}}^-{\cong}[r]&  \der(N)=\der(X'\times_{Y'}Y'\times_{Y'}W\times_YX)
    }
\end{equation*}

The composition of functors in the right column is isomorphic to the identity functor by the base change isomorphism \eqref{eq:base-change-diagram-app}. It follows that \eqref{triangle-identity-N} in the current setting holds. The same reasoning implies that \eqref{triangle-identity-M} in the current setting also holds.
\end{proof}

In practice, the assumption in \Cref{lem-geo-duality-I} often does not hold. But under some certain technical assumptions, we can still understand the duality datum.

\begin{lemma}\label{lem-geo-duality-II}
Assumptions are as in \Cref{cor-fully-faithful-usual-trace-to-loop} and assume that $f: X\to Y$ and $f':X'\to Y'$ satisfy assumptions of \emph{loc. cit}.
In addition, suppose the following exterior tensor products are fully faithful 
\begin{equation}\label{lem-geo-duality-II-exterior tensor equivalence}
\der(N)\otimes_{\der(\pt)}\der(T)\to \der(N\times T),\quad  T=M, M\times N,R,
\end{equation}
and in addition is an equivalence when $T=M,M\times N$.
Then the functor $e$ from \eqref{eq-pre-geo-evaluation} and the functor
\[
u: \der(X'_1)\to \der(X'\times_{Y'}W\times_YW\times_{Y'}X')\to \der(N)\otimes_{\der(X_1)}\der(M),
\]
where the last functor is the right adjoint of the horizontal morphism in \eqref{eq-pre-geo-unit}, form a duality datum.
\end{lemma}
\begin{proof} 
We will write $T=X'\times_{Y'}W\times_{Y} W\times_{Y'} W\times_{Y}X$.
As in the proof of \Cref{lem-geo-duality-I}, it is enough to establish the following commutative diagram
\quash{\begin{equation*}
    \xymatrix{
    \der(X'_1)\otimes_{\der(X'_1)}\der(N)\ar_-{u\otimes\id}[dd] &\ar@{=}[l]\der(X'_1)\otimes_{\der(X'_1)}\der(N)\ar[d] \ar^-{\cong}[r]& \der(N)\ar[dd]\\
    &\der(R)\otimes_{\der(X'_1)}\der(N)\ar_{(**)}[dl]\ar[dr]&\\
    \der(N)\otimes_{\der(X_1)}\der(M)\otimes_{\der(X'_1)}\der(N) \ar_-{\id\otimes e}[dd]\ar[dr]&& \der(X'\times_{Y'}W\times_{Y} W\times_{Y'} W\times_{Y}X)\ar[dd]\ar_{(**)}[dl]\\
    & \der(N)\otimes_{\der(X_1)}\der(S) \ar[d]&\\  
    \der(N)\otimes_{\der(X_1)}\der(X_1)&\ar@{=}[l]\der(N)\otimes_{\der(X_1)}\der(X_1)\ar^-{\cong}[r]&  \der(N),
    }
\end{equation*}
where the arrows labelled by $(**)$ are right adjoint of the corresponding arrows in the diagram from the proof of \Cref{lem-geo-duality-I}. 
}
\begin{tiny}
\begin{equation*}
    \xymatrix{
    \der(X'_1)\otimes_{\der(X'_1)}\der(N)\ar_-{u\otimes\id}[dd] &\ar@{=}[l]\der(X'_1)\otimes_{\der(X'_1)}\der(N)\ar[d] \ar^-{\cong}[r]& \tr(\der(X'_1),\der(X'_1\times N)) \ar^-{\cong}[r]\ar[d] & \der(N)\ar[dd]\\
    &\der(R)\otimes_{\der(X'_1)}\der(N)\ar_{(**)}[dl]\ar@{^{(}->}[r] \ar@{}[d] | {(\mathrm{I})}&\ar_{(**)}[dl]\tr(\der(X'_1),\der(R\times N)) \ar[dr] &\\
    \der(N)\otimes_{\der(X_1)}\der(M)\otimes_{\der(X'_1)}\der(N) \ar_-{\id\otimes e}[dd]\ar^-{\cong}[r]\ar[dr]&\tr(\der(X_1)\otimes_{\der(\pt)} \der(X'_1), \der(N\times M\times N))\ar[dr] & (\mathrm{II}) & \der(T)\ar[dd]\ar_{(**)}[dl]\\
    & \der(N)\otimes_{\der(X_1)}\der(S) \ar[d]\ar[r]&\tr(\der(X_1), \der(N\times S)) \ar[d] \ar@{}[dr] | {(\mathrm{III})} &\\  
    \der(N)\otimes_{\der(X_1)}\der(X_1)&\ar@{=}[l]\der(N)\otimes_{\der(X_1)}\der(X_1)\ar^-{\cong}[r]&\tr(\der(X_1),\der(N\times X_1))\ar^-{\cong}[r]  & \der(N),
    }
\end{equation*}
\end{tiny}

In the diagram arrows labelled by $(**)$ are right adjoint of natural functors (compare with the diagram from the proof of \Cref{lem-geo-duality-I}). 
Only commutativity of $(\mathrm{I}), (\mathrm{II}), (\mathrm{III})$ requires justification.

We note that both functors in $\der(N)\otimes_{\der(X_1)}\der(X_1)\to \tr(\der(X_1),\der(N\times X_1))\to  \der(N)$ are fully faithful (by \Cref{cor-fully-faithful-usual-trace-to-loop}) and the composition is an equivalence (see \Cref{ex: relative tensor split case}). It follows that the functor $\tr(\der(X_1),\der(N\times X_1))\to  \der(N)$ is an equivalence, as indicated by the diagram. Now the
commutativity of $(\mathrm{III})$ follows from \Cref{lem-right-adjoint-from-categorical-trace-to-loop-space}.

Next we deal with $(\mathrm{II})$.
Consider the commutative diagram
\[
\xymatrix{
\der(N\times M\times N) \ar[r]\ar[d] & \der(R\times N) \ar[d]\\
\der(N\times S) \ar[r] & \der(T).
}
\]
with horizontal morphisms are induced by the correspondence $Y\leftarrow X\to X\times X$ and vertical morphisms induced by $Y'\leftarrow X'\to X'\times X'$. As $f,f'\in \bfC_{\mathrm{VR}}$ and $\Delta_X,\Delta_{X'}\in \bfC_{\mathrm{HR}}$, the above diagram is right adjointable by the same proof as in \Cref{lem:trace.segal.beck.chevalley}. 

Note that the functor $\der(R\times N)\to \der(N\times M\times N)$ obtained by right adjoint is $\der(X'_1)$-bilinear, by \Cref{prop-comparison-usual-trace-geo-trac-presemirigid} and \Cref{lem-rigid-lax-is-strict}. It then follows from \eqref{eq-functoriality-trace-to-relative-trace} that the following diagram is commutative
\[
\xymatrix{
\tr(\der(X'_1),\der(N\times M\times N)) \ar[d] & \tr(\der(X'_1), \der(R\times N)) \ar[d]\ar[l]\\
\der(N\times S)  & \der(T)\ar[l].
}
\]
The left vertical functor is $\der(X_1)$-bilinear. We then use \Cref{L:hochschild-homology-vs-cohomology} \eqref{cor:right-adjointability-of-module-to-Trace} (or rather the proof) to conclude the commutativity of $(\mathrm{II})$.

Finally we consider $(\mathrm{I})$. The hook arrow is fully faithful by \Cref{cor-fully-faithful-usual-trace-to-loop} and by our assumption \eqref{lem-geo-duality-II-exterior tensor equivalence}.  As indicated in the diagram the functor 
\[
\der(N)\otimes_{\der(X_1)}\der(M)\otimes_{\der(X'_1)}\der(N)\to \tr(\der(X_1)\otimes_{\der(\pt)} \der(X_1), \der(N\times M\times N))
\] 
is an equivalence again by \Cref{cor-fully-faithful-usual-trace-to-loop} and by our assumption \eqref{lem-geo-duality-II-exterior tensor equivalence} is an equivalence when $T=M, M\times N$. The diagram $(\mathrm{I})$ is obtained from the following commutative diagram by taking the right adjoint of vertical functors
\[
\xymatrix{
\der(R)\otimes_{\der(X'_1)}\der(N)\ar@{^{(}->}[r]&\tr(\der(X'_1),\der(R\times N))\\
 \der(N)\otimes_{\der(X_1)}\der(M)\otimes_{\der(X'_1)}\der(N) \ar^-{\cong}[r]\ar@{^{(}->}[u] & \tr(\der(X_1)\otimes_{\der(\pt)} \der(X'_1), \der(N\times M\times N))\ar@{^{(}->}[u]
}\]
In this diagram, all functors are fully faithful. It follows that $(\mathrm{I})$ is commutative.
\end{proof}

\begin{corollary}\label{cor-geom-duality-left-module}
Assumptions are as in \Cref{cor-fully-faithful-usual-trace-to-loop}.
Let $W$ be an object such that both morphisms $\Delta_W: W\to W\times W$ and  $\pi_W: W\to \pt$ belonging to $\horiz$. Let $h: W\to Y$ be a morphism in $\verti$ such that $\Delta_{W/Y}:W\to W\times_YW$ also belongs to $\verti$.
\begin{enumerate}
\item\label{cor-geom-duality-left-module-1} If there is an object $u\in \der(W\times_YX)\otimes_{\der(X\times_YX)}\der(X\times_YW)$ whose image in $\der(W\times_YW)$ is $(\Delta_{W/Y})\vertl(\La_W)$, then $\der(X\times_YW)$ is dualizable as a left $\der(X\times_YX)$ with the duality datum given by $(u,e)$.
\item\label{cor-geom-duality-left-module-2} Suppose $\der(X
\times_YW)\otimes_\La\der(T)\to \der(X\times_YW\times T)$ is fully faihtful when $T\in\bfC$ and is an equivalence when $T= X\times_YW$ and $T=X\times_YW\times W\times_YX$. Let $\proj_{\trg}$ be the right adjoint of the functor $\der(W\times_YX)\otimes_{\der(X\times_YX)}\der(X\times_YW)\to \der(W\times_YW)$. Then $\der(X\times_YW)$ is left dualizable as a $\der(X\times_YX)$-module with duality datum given by $(\proj_{\trg}((\Delta_{W/Y})\vertl(\La_W)),e)$.
\end{enumerate}
\end{corollary}

Now suppose we are given a $\der(X_1)$-$\der(X'_1)$-bimodule homomorphism
\[
\al: \der(M)\otimes_{\der(X'_1)}\der(Q')\to \der(Q)\otimes_{\der(X_1)}\der(M).
\]
Then as explained above, under certain dualizability assumption of $\der(M)$, there is a functor
\[
\tr(\der(M),\al): \tr(\der(X'_1),\der(Q'))\to \tr(\der(X_1),\der(Q)).
\]
On the other hand, suppose we are given a correspondence 
\[
\al_{\geo}: W\times_{Y'}Z'\dashrightarrow Z\times_YW
\]
in $\corr(\bfC_{/Y\times Y'})_{\verti;\horiz}$. One can form the correspondence
\[
C(W,\alpha_{\geo}):  Y'\times_{Y'\times Y'} Z'\dashrightarrow
Y\times_{Y\times Y} Z
\]
given by the composition 
\begin{multline}\label{eq: algeo correspondence}
Y'\times_{Y'\times Y'}Z'\stackrel{u_{\geo}\times \id}{\dashrightarrow} (W\times_YW)_{Y'\times Y'}Z'\cong Y\times_{Y\times Y}(W\times_{Y'}Z'\times_{Y'}W)\\
\stackrel{\id\times \al_{\geo}\times\id}{\dashrightarrow} Y\times_{Y\times Y}(Z\times_YW\times_{Y'}W)\stackrel{\id\times\id\times e_{\geo}}{\dashrightarrow} Y\times_{Y\times Y}Z.
\end{multline}

The sheaf theory $\der$ then induces a functor
\[
\der(C(W,\alpha_{\geo})): \der(Y'\times_{Y'\times Y'}Z')\to \der(Y\times_{Y\times Y}Z).
\]
We would like to relate $\tr(\der(M),\al)$ with the above functor under certain assumptions. 

\begin{assumptions}\label{assumption-alpha-alphageo}
\begin{enumerate}[(I)]
\item\label{item-al-algeo} We assume that the following diagram is commutative
\begin{equation}
    \xymatrix{
    \der(M)\otimes_{\der(X'_1)}\der(Q')\ar[r]\ar_{\al}[d] & \der(X\times_YW\times_{Y'}Z'\times_{Y'}X')\ar^{\der(\id\times\al_{\geo}\times\id)}[d]\\
    \der(Q)\otimes_{\der(X_1)}\der(M)\ar[r] & \der(X\times_YZ\times_YW\times_{Y'}X').
    }
\end{equation}
\item\label{item-al-algeo-right} We assume that the following diagram is commutative
\begin{equation}
    \xymatrix{
    \der(M)\otimes_{\der(X'_1)}\der(Q')\ar_{\al}[d] &\ar[l] \der(X\times_YW\times_{Y'}Z'\times_{Y'}X')\ar^{\der(\id\times\al_{\geo}\times\id)}[d]\\
    \der(Q)\otimes_{\der(X_1)}\der(M) &\ar[l] \der(X\times_YZ\times_YW\times_{Y'}X').
    }
\end{equation}
where the horizontal arrows are right adjoints of the natural ones.
\end{enumerate}
\end{assumptions} 

\begin{example}\label{rem-for-assumption-alpha-alphageo}
Suppose that $X,Y,X',Y'$ are equipped with automorphisms $\phi$ and $f$ and $f'$ are $\phi$-equivariant. (Recall that this means we need to supply isomorphisms as in \eqref{eq: meaning of equiv in higher sense}.) In addition, assume that 
 \begin{itemize} 
  \item $Z=Y$ with $g_1=\id$ and $g_2=\phi: Y\to Y$  as in \Cref{trace.phi.convolution.fully.faithful};
 \item $Z'=Y'$ with $g'_1=\id$ and $g'_{2}=\phi:Y'\to Y'$;
 \item there is an automorphism $\phi: W\to W$ and $h:W\to Y\times Y'$ is $\phi$-equivariant. (Again this means that we need to supply an isomorphism as in \eqref{eq: meaning of equiv in higher sense}.) 
 \end{itemize}
 Then $M$ also admits an automorphism, still denoted by $\phi$. Suppose $\al$ is given by
\[
\der(M)\otimes_{\der(X_1')}\der(Q')\cong \der(M)\xrightarrow{\phi\horizl} \der(M)\cong \der(Q)\otimes_{\der(X_1)}\der(M)
\]
and $\al_\geo$ is given by the horizontal map
\[
Z\times_YW\cong W  \xrightarrow{\phi} W\cong W\times_{Y'}Z'.
\]
In this case, \Cref{assumption-alpha-alphageo} holds.
\end{example}

\begin{proposition}\label{prop-geo-class-via-corr}
 Under the assumption in \Cref{lem-geo-duality-I} and \Cref{assumption-alpha-alphageo} \ref{item-al-algeo}, then the following diagram is commutative
 \[
 \xymatrix{
 \tr(\der(X'_1),\der(Q'))\ar^-{\tr(\der(M),\al)}[rr]\ar[d] && \tr(\der(X_1),\der(Q))\ar[d]\\
 \der(Y'\times_{Y'\times Y'}Z') \ar^-{\der(C(W,\al_\geo))}[rr]  && \der(Y\times_{Y\times Y}Z).
 }
 \]
 
 Under the assumption in \Cref{lem-geo-duality-II} and \Cref{assumption-alpha-alphageo} \ref{item-al-algeo-right}, the following diagram is commutative
  \[
 \xymatrix{
 \tr(\der(X'_1),\der(Q'))\ar^-{\tr(\der(M),\al)}[rr]\ar[d] && \tr(\der(X_1),\der(Q))\\
 \der(Y'\times_{Y'\times Y'}Z') \ar^-{\der(C(W,\alpha_{\geo}))}[rr]  && \der(Y\times_{Y\times Y}Z)\ar_{\proj_{\trg}}[u].
 }
 \]
\end{proposition}
\begin{proof}
We write $-\otimes-$ instead of $-\otimes_{\der(\pt)}-$ to simplify notations.
The first case follows from the following commutative diagram
\begin{equation}\label{eq-transfer-categorical-vs-geometric}
\xymatrix{
\der(X'_1)\otimes_{\der(X'_1)\otimes \der(X'_1)^{\rev}}\der(Q')\ar_{u\otimes 1}[d]\ar[r]& \der(Y'\times_{Y'\times Y'}Z')\ar^{\der(u_\geo\times\id)}[d]\\ 
(\der(N)\otimes_{\der(X_1)}\der(M))\otimes_{\der(X'_1)\otimes\der(X'_1)^{\rev}}\der(Q')\ar_\cong[d]\ar[r] &  \der((W\times_Y W)\times_{Y'\times Y'}Z')\ar^\cong[d]\\
\der(X_1)\otimes_{\der(X_1)\otimes\der(X_1)^{\rev}} (\der(M)\otimes_{\der(X'_1)}\der(Q')\otimes_{\der(X'_1)}\der(N))\ar[r] \ar_{1\otimes \al\otimes 1}[d] & \der(Y\times_{Y\times Y}(W\times_{Y'}Z'\times_{Y'}W))\ar^{\der(\id\times \al_\geo\times\id)}[d]\\
\der(X_1)\otimes_{\der(X_1)\otimes\der(X_1)^{\rev}}(\der(Q)\otimes_{\der(X_1)}\der(M)\otimes_{\der(X'_1)}\der(N))\ar[r] \ar_{1\otimes 1\otimes e}[d]& \der(Y\times_{Y\times Y}(Z\times_YW\times_{Y'}W)) \ar^{\der(\id\times\id \times e_{\geo})}[d]\\
\der(X_1)\otimes_{\der(X_1)\otimes\der(X_1)^{\rev}}\der(Q) \ar[r] & \der(Y\times_{Y\times Y}Z)
}
\end{equation}

The second case follows from a similar diagram
\begin{equation}\label{eq-transfer-categorical-vs-geometric-right}
\xymatrix{
\der(X'_1)\otimes_{\der(X'_1)\otimes \der(X'_1)^{\rev}}\der(Q')\ar_{u\otimes 1}[d]\ar[r]& \der(Y'\times_{Y'\times Y'}Z')\ar^{\der(u_\geo\times\id)}[d]\\ 
(\der(N)\otimes_{\der(X_1)}\der(M))\otimes_{\der(X'_1)\otimes\der(X'_1)^{\rev}}\der(Q')\ar_\cong[d] &\ar[l]  \der((W\times_Y W)\times_{Y'\times Y'}Z')\ar^\cong[d]\\
\der(X_1)\otimes_{\der(X_1)\otimes\der(X_1)^{\rev}} (\der(M)\otimes_{\der(X'_1)}\der(Q')\otimes_{\der(X'_1)}\der(N)) \ar_{1\otimes \al\otimes 1}[d] &\ar[l] \der(Y\times_{Y\times Y}(W\times_{Y'}Z'\times_{Y'}W))\ar^{\der(\id\times \al_\geo\times\id)}[d]\\
\der(X_1)\otimes_{\der(X_1)\otimes\der(X_1)^{\rev}}(\der(Q)\otimes_{\der(X_1)}\der(M)\otimes_{\der(X'_1)}\der(N)) \ar_{1\otimes 1\otimes e}[d]&\ar[l] \der(Y\times_{Y\times Y}(Z\times_YW\times_{Y'}W)) \ar^{\der(\id\times\id \times e_{\geo})}[d]\\
\der(X_1)\otimes_{\der(X_1)\otimes\der(X_1)^{\rev}}\der(Q)  &\ar[l] \der(Y\times_{Y\times Y}Z),
}
\end{equation}
where the horizontal left arrows are obtained by the corresponding horizontal right arrows in \eqref{eq-transfer-categorical-vs-geometric} by passing to the right adjoint.
We need to justify the commutativity of this diagram. First we have the commutativity of the following diagram
\begin{small}
\[
\xymatrix{
&\der(X'_1)\otimes_{\der(X'_1)\otimes \der(X'_1)^{\rev}}\der(Q')\ar_{u\otimes 1}[dl]\ar[d]\ar[r]& \der(Y'\times_{Y'\times Y'}Z')\ar^{\der(u_\geo\times\id)}[d]\\ 
(\der(N)\otimes_{\der(X_1)}\der(M))\otimes_{\der(X'_1)\otimes\der(X'_1)^{\rev}}\der(Q') &\ar[l] \der(R)\otimes_{\der(X'_1)\otimes\der(X'_1)^{\rev}}\der(Q')&\ar[l]  \der((W\times_Y W)\times_{Y'\times Y'}Z')
}
\]
\end{small}

Indeed, the left triangle is commutative as we are in the case as in \Cref{lem-geo-duality-II}, and the right square is commutative as the natural functor $\der(R)\otimes_{\der(X'_1)\otimes\der(X'_1)^{\rev}}\der(Q')\to  \der((W\times_Y W)\times_{Y'\times Y'}Z')$ is fully faithful by \Cref{cor-fully-faithful-usual-trace-to-loop}. This justifies the commutativity of the top square in \eqref{eq-transfer-categorical-vs-geometric-right}.

For the commutativity of the third square in \eqref{eq-transfer-categorical-vs-geometric-right}, we can first argue as in the proof of \Cref{lem-geo-duality-II} (more precisely the proof of the commutativity of the square $(\mathrm{II})$ there) to obtain the commutativity of
\[
\xymatrix{
\der(M)\otimes_{\der(X'_1)}\der(Q')\otimes_{\der(X'_1)}\der(N)\ar[d] & \ar[l]\der(X\times_YW\times_{Y'}Z'\times_{Y'}W\times_YX) \ar[d]\\
\der(Q)\otimes_{\der(X'_1)}\der(M)\otimes_{\der(X'_1)}\der(N) & \ar[l] \der(X\times_YZ\times_{Y}W\times_{Y'}W\times_YX),
}\]
using that $\der(M)^{\otimes 2}\otimes \der(T)\to \der(M^2\times T)$ is an equivalence for $T=Q'$.
In addition, this diagram is $\der(X_1)$-bilinear. It follows from \Cref{lem-right-adjoint-from-categorical-trace-to-loop-space} that the third square in \eqref{eq-transfer-categorical-vs-geometric-right} is indeed commutative.
Similar argument also shows that the last square in \eqref{eq-transfer-categorical-vs-geometric-right} is commutative.
\end{proof}

\begin{corollary}\label{ex: geo phi-trace class}
Consider the situation as in \Cref{rem-for-assumption-alpha-alphageo}, with $X'=Y'=\pt$.
Assumptions are as in \Cref{cor-geom-duality-left-module} \eqref{cor-geom-duality-left-module-2}.
Then
\[
[\der(X\times_YW),\al]_{\phi}=\proj_{\trg}(\mL_\phi(h)\vertl(\La_{\mL_\phi(W)}))
\]
as objects in $\tr(\der(X\times_YX),{}^\phi\der(X\times_YX))\subset \der(\mL_\phi(Y))$. In particular, if $W=Y$ with $W\to Y$ being the identity map, then $[\der(X),\alpha]_\phi=\proj_{\trg}(\La_{\mL_\phi(Y)})$.

If assumptions are as in \Cref{cor-geom-duality-left-module} \eqref{cor-geom-duality-left-module-1}, one can remove $\proj_{\trg}$ in the above formulas.
\end{corollary}
\begin{proof}
In this case the correspondence as in \eqref{eq: algeo correspondence} is given by  $\pt \xleftarrow{}  \mL_\phi(W)  \xrightarrow{\mL_\phi(h)} \mL_\phi(Y)$.
\end{proof}

\newpage

\section{Theory of coherent sheaves}
\label{S: theory of coherent sheaves}
Let $\Lambda$ be an (ordinary) regular noetherian ring. In this section, we examine the theory of (ind-)coherent sheaves on (derived ind-)algebraic stacks almost of finite presentation over $\Lambda$. When $\Lambda$ is a field of characteristic zero, such theory was extensively developed in \cite{gaitsgory2013.indcoh, Gaitsgory.Rozenblyum.DAG.vol.I}. The approaches of \emph{loc. cit.} generalize well to schemes (and algebraic spaces) over more general base rings, such as perfect fields of positive characteristic. However, this theory does not seem suitable for addressing certain questions regarding algebraic stacks over fields of positive characteristic, particularly for our intended applications. 

Without delving into details, we note that the category of ind-coherent sheaves on stacks, as developed in the aforementioned references,  is defined via descent, In contrast, in geometric representation theory it appears more natural to consider the ind-completion of the usual category of coherent sheaves on stacks. While these two categories coincide for most algebraic stacks one encounters in practice when $\La$ is a field of characteristic zero, this coincidence breaks down in the case of a field of positive characteristic. In fact, they differ even for the classifying stacks of most algebraic groups.

We will also need the theory of singular support for coherent sheaves. Again, when the base ring $\La$ is a field of characteristic zero, such theory was developed in \cite{arinkin2015singular}. But in positive characteristic, some extra care is needed (even for schemes).

Consequently, we take this opportunity to outline how to establish results for the ind-completion of the category of coherent sheaves, paralleling those proved in the aforementioned works. While we do not aim to develop the theory with maximal generality in this article, we will focus on the aspects that are necessary for our current discussion.

\subsection{Derived algebraic geometry}\label{SS: derived algebraic geometry}
We very briefly review the terminologies and results from derived algebraic geometry we need. As before, one of the purposes of this subsection is to fix the notations.

We allow $\La$ to be a base ordinary commutative ring (not necessarily regular noetherian) in this subsection.
Let $\calg_\La$ be the category of animated $\La$-algebras. Let $\calg_\La^{\heartsuit}\subset\calg_\La$ denote the ordinary category of usual commutative $\La$-algebras.
For an animated $\La$-algebra $A$, let $\Mod_A$ denote the category of $A$-modules. It is equipped with a standard $t$-structure and let $\Mod_A^{\leq 0}$ denote the connective part. 

Recall that a morphism $A\to B$ of animated $\La$-algebras is called flat if $\pi_0(A)\to \pi_0(B)$ is flat and $\pi_0(B)\otimes_{\pi_0(A)}\pi_i(A)\cong \pi_i(A)$. A morphism $A\to B$ is called Zariski open, resp. \'etale, resp. smooth, resp. faithfully flat, if $A\to B$ is flat and the map $\pi_0(A)\to \pi_0(B)$ is Zariski open, resp. \'etale, resp. smooth, resp. faithfully flat. We thus have the usual Zariski, \'etale topology on $\calg_\La$.

\subsubsection{Prestacks and stacks}
\begin{definition}
A prestack is an accessible functor $X: \calg_\La\to \spc$. All prestacks over $\La$ form a full subcategory of $\fun(\calg_\La,\spc)$, denoted by
$\prestk_\La$.  By a(n \'etale) stack, we mean a prestack which is a sheaf with respect to the \'etale topology on $\calg_\La$.
\end{definition}

\begin{remark}\label{rem: why accessibility}
Accessibility is a set theoretic condition that guarantees that for a prestack $X$, the slicing category $(\calg_\La)^{\op}_{/X}=\{(R,x)\mid R\in\calg_\La, x\in X(R)\}^{\op}$ admits a small subcategory that is cofinal. This allows us to take various colimits along $(\calg_\La)^\op_{/X}$. We shall not repeat this remark in the future.
\end{remark}

There is the fully faithful Yoneda embedding $(\calg_\La)^{\op}\subset \prestk_\La$. Essential images are called (derived) affine schemes over $\La$. As usual, the image of $A\in\calg_\La$ in $\prestk_\La$ is denoted as $\spec A$. The essential image of the fully faithful embedding $(\calg_\La)^{\op}\to \prestk_\La$ is also denoted as $\Aff_\La$. Objects are called (derived) affine schemes.  Affine schemes are stacks. A morphism $f:X\to Y$ of prestacks is called affine if for every morphism $Z\to Y$ with $Z$ an affine scheme, the fiber product $Z\times_YX$ is an affine scheme.

Similarly, we let $\prestk_\La^\cl\subset \fun(\calg_\La^{\heartsuit},\spc)$ denote the full subcategory of accessible functors, called the category of classical prestacks.
Restriction along $\calg_\La^{\heartsuit}\subset\calg_\La$ defines a functor 
\[
\prestk_\La\to \prestk_\La^\cl,\quad  X\mapsto X_\cl,
\] 
which admits a fully faithful left adjoint functor given by sending $(F: \calg_\La^{\heartsuit}\to \spc)\in  \prestk_\La^\cl$ to its
left Kan extension along $\calg_\La^{\heartsuit}\subset\calg_\La$. 
We call $X_\cl$ as above the underlying classical prestack associated to $X$, and then regard $X_\cl$ as a prestack. E.g. $(\Spec A)_\cl=\Spec \pi_0(A)$.  Note that there is a canonical morphism $X_\cl\to X$. E.g. if $X=\spec A$, then the map $(\Spec A)_\cl=\Spec \pi_0(A)\to \Spec A$ is given by $A\to \pi_0(A)$. A prestack is called classical if it belongs to $\prestk_\La^\cl$.

\begin{remark}\label{rem: truncated stacks}
Let $\tau_{\leq n}\calg_\La\subset \calg_\La$ be the full subcategory of $n$-truncated objects, i.e. the full subcategory of animated $\La$-algebras $A$ such that $\pi_i(A)=0$ for $i>n$ (so $\tau_{\leq 0}\calg_\La=\calg_\La^{\heartsuit}$). One can then similarly define ${}_{\leq n}\prestk_\La\subset \fun(\tau_{\leq n}\calg_\La, \spc)$ as the full subcategory of accessible functors. Similarly, the restriction along ${}_{\leq n}\calg_\La\subset \calg_\La$ induces $\prestk_\La\to {}_{\leq n}\prestk_\La$ which admits a fully faithful left adjoint by left Kan extensions. We have
\[
\prestk_\La^{\cl}={}_{\leq 0}\prestk_\La\subset {}_{\leq 1}\prestk_{\La}\subset\cdots\subset \prestk_\La.
\]
We will let ${}_{\leq \infty}\prestk_\La=\bigcup_n {}_{\leq n}\prestk_{\La}$. For an object $X\in \prestk_\La$, we let $X_{\leq n}$ be its image in ${}_{\leq n}\prestk_\La$, though as an object in $\prestk_\La$ via the fully faithful embedding ${}_{\leq n}\prestk_\La\subset \prestk_\La$.
\end{remark}

It is convenient to associate to a prestack a topological space. Namely, we consider the left Kan extension  along $(\calg_\La)^{\op}\subset\prestk_\La$ of the usual functor $|\cdot|$ that assigns $R\in\calg_\La$ to the spectrum $|\Spec \pi_0(R)|$ of $\pi_0(R)$. Concretely, this means that
\begin{equation}\label{eq-top-space-prestack}
|X|=\colim_{(\calg_\La)^{\op}_{/X}} |\Spec \pi_0(R)|,
\end{equation} 
where the colimit is taken in the (ordinary) category of topological spaces.  Clearly $|X|=|X_\cl|$.
The topological space $|X|$ could be quite wild in general, i.e. it may not be a spectral space.
However, the underlying set
of points of $|X|$ can be described as in \cite[\href{https://stacks.math.columbia.edu/tag/04XE}{Section 04XE}]{stacks-project}. We say a morphism $f\colon X\to Y$ of prestacks surjective
if the induced map $|X|\to |Y|$ is surjective. Equivalently, for every field $K$ and a point $y \in Y (K)$
there exists a field extension $L/K$ and a point  $x \in X(L)$ lifting $y$. 

\subsubsection{Derived schemes, algebraic spaces, and algebraic stacks}
A morphism $X\to Y$ of prestacks is called an open embedding if there is some open subset $U\subset |Y|$ such that 
\begin{equation}\label{eq: points for open embedding}
X(A)=Y(A)\times_{\Map(|\Spec A|,|Y|)}(\Map(|\Spec A|, U),\quad A\in\calg_\La.
\end{equation}
Clearly, open embeddings form a strongly stable class of morphisms  (in the sense of \Cref{def-closure-property-of-morphism}) in $\prestk_\La$, and are $0$-truncated in the sense of \cite[Definition 5.5.6.8]{Lurie.higher.topos.theory}.
A prestack $X$ is called a (derived) scheme if it is a stack and admits an open covering by (derived) affine schemes, i.e. a collection of open embeddings $\{\Spec A_i\to \mF\}_i$ which is jointly surjective. Let $\Sch_\La\subset\prestk_\La$ be the full subcategory of derived schemes. Just as in the classical algebraic geometry, this subcategory is closed under finite product, and can also be realized as a full subcategory of locally derived ringed spaces. We refer to \cite{Lurie.SAG} for this approach. Note that if $X$ is a derived scheme, then $X_\cl$ is a scheme in the classical sense.

The notion of \'etale, smooth, (faithfully) flat morphisms, etc. between derived schemes make sense, as they are properties local in Zariski topology, and all of them form weakly stable classes of morphisms (the class of \'etale morphisms is strongly stable). Therefore, one can apply \Cref{rem-prop-of-weak-strong-stable-class} \eqref{rem-prop-of-weak-strong-stable-class-2} to make sense of these classes for morphisms between prestacks that are representable (in derived schemes). 
Now a prestack $X$ is called a (derived) algebraic space if it is a stack and admits an \'etale covering by $\{U_i\to X\}$ by derived schemes $U_i$. Let $\algsp_\La$ denote the full subcategory of derived algebraic spaces over $\La$, which again is closed under finite product. It can also be realized as a full subcategory of locally derived ringed topos. Again, we refer to \cite{Lurie.SAG} for this approach. If $X$ is a derived algebraic space, then $X_\cl$ is an algebraic space in the classical sense.
One can then iterate the procedure to make sense of \'etale, smooth, (faithfully) flat morphisms between morphisms between prestacks that are representable (in derived algebraic spaces).

Recall a morphism $f: X\to Y$ of derived algebraic spaces is called proper (closed embedding) if $f_\cl: X_\cl\to Y_\cl$ is a closed embedding in the classical sense. In particular, $X_\cl\to X$ is a closed embedding. We say $f$ to be finite if it is both affine and proper. Closed embeddings are finite morphisms.

All topological notions, such as quasi-compact and quasi-separated (qcqs) make sense in this setting. Therefore in $\Sch_\La$, we have the full subcategory $\Sch_\La^{\qc}$, resp. $\Sch_\La^{\qs}$, resp. $\Sch_\La^{\qcqs}$, resp. $\Sch_\La^{\sep}$, resp. $\Sch_\La^{\qc.\sep}$ of quasi-compact, resp. quasi-separated, resp. qcqs, resp. separated, resp. quasi-compact and separated schemes. We have similarly defined full subcategories of algebraic spaces.

Finally, by an algebraic stack, we mean a stack $X$ over $\La$ such that  the diagonal $X\to X\times_\La X$ is representable by a derived algebraic spaces and there exists a smooth surjective map $U\to X$ with $U$ a derived algebraic space.  An algebraic stack is called quasi-separated if the diagonal is quasi-compact and quasi-separated. It is called quasi-compact if $U$ can be chosen to be quasi-compact, or equivalently the topological space $|X|$ is quasi-compact.
Let $\ArStk_\La\subset\prestk_\La$ denote the full subcategory of algebraic stacks over $\La$.  We have similarly defined full subcategories  $\ArStk_\La^{\qs}\subset \ArStk_\La^{\qcqs}$ of algebraic stacks.

One can inductively define the notion of Artin $n$-stacks, and many discussions below hold for these more general objects. However, we do not need such generalities.

\subsubsection{Almost of finite presentation}
Recall that for a compactly generated category $\bfC$, an object $c$ is called almost compact if for every $n\geq 0$, $\tau_{\leq n}c$ is compact in $\tau_{\leq n}\bfC$ (\cite[Definition 7.2.4.8]{Lurie.higher.algebra}). Almost compact objects in $\calg_\La$ are also called animated rings almost of finite presentation over $\La$.
For an animated $\La$-algebra $A$, almost compact objects in $\Mod_A^{\leq 0}$ are also called connective almost perfect $A$-modules. 
If $\La$ is noetherian, $A$ is almost of finite presentation over $\La$ if and only if $\pi_0(A)$ is a finitely generated $\La$-algebra and each $\pi_i(A)$ is a finitely generated $\pi_0(A)$-module.
In particular, if $\La$ is noetherian, a classical $\La$-algebra of finite type is almost of finite presentation, when regarded as an animated $\La$-algebra.

Let $\calg_{\La}^{\aft}\subset \calg_\La$ be the category $\La$-algebras that are almost of finite presentation over $\La$. We let 
\[
\prestk_\La^{\laft}=\fun(\calg_\La^{\aft}, \spc),
\] 
and call objects in this category prestack locally almost of finite presentation over $\La$. (Note that unlike \cite[2.1.7.2]{Gaitsgory.Rozenblyum.DAG.vol.I}, we do not require prestack locally almost of finite presentation over $\La$ to be nilcomplete/convergent.)
Restriction along $\calg_{\La}^{\aft}\subset \calg_\La$ defines a functor $\prestk_\La\to \prestk_\La^{\laft}$ that admits a fully faithful left adjoint via left Kan extensions. In this way, we regard $\prestk_\La^{\laft}$ as a full subcategory of $\prestk_\La$.
We let 
\[
\Sch_\La^{\aft}=\Sch^{\qcqs}_\La\cap \prestk_\La^{\laft},\quad \algsp_\La^{\aft}=\algsp^{\qcqs}_\La\cap \prestk_\La^{\laft},\quad \ArStk_\La^{\aft}=\ArStk^{\qcqs}_\La\cap \prestk_\La^{\laft}.
\]

For our purpose, we also need ind-objects. 

\begin{definition}\label{def-ind-scheme-derived}
An object $X\in \prestk_\La^{\laft}$ is called an \textit{ind-scheme} (resp. \textit{ind-algebraic space}, resp. \textit{ind-algebraic stack}) if 
\begin{itemize}
\item it is nilcomplete (or sometimes called convergence), i.e. $X(A)=\lim_n X(\tau_{\leq n}A)$ for every $A\in\calg_\La$; 
\item and can be written as a filtered colimit of $X = \colim_{i} X_i$ of $X_i \in \Sch^{\aft}_\La$ (resp. $X_i \in \algsp^{\aft}_\La$ resp. $X_i\in \ArStk_\La^{\aft}$) with transition maps given by closed immersions. 
\end{itemize}
Let $\indsch^{\aft}_\La\subset\indsp^{\aft}_\La\subset \indarstk^{\aft}_\La\subset \prestk^{\laft}_\La$ denote the full subcategory of ind-schemes, ind-algebraic spaces, and ind-algebraic stacks over $\La$.
\end{definition}

\begin{definition}\label{def-ind-aft-morphism}
We will let $\ind\aft$ denote the class of morphisms in $\prestk^{\laft}_\La$ that are representable in $\indsp^{\aft}_\La$. More precisely, a morphism $f: X\to Y$ in $\prestk^{\laft}_\La$ belongs to $\ind\aft$ if for every $Y'\to Y$ with $Y'\in \algsp_\La^{\aft}$, we have $X'=Y'\times_YX\in \indsp^{\aft}_\La$.
\end{definition}

In literature, sometimes ind-schemes are defined as prestacks which can be written as a filtered colimit of $X = \colim_{i} X_i$ with transition maps being closed embeddings as above, but without requiring $X_i$ to be almost of finite presentation.
However, \Cref{def-ind-scheme-derived} is enough to our purpose.

\begin{example}\label{ex: formal completion}
The main example of ind-algebraic stack we need in this article is the formal completion of an algebraic stack along a closed substack. Namely, let $X\in \ArStk_\La^{\aft}$ and let $Z\subset X$ be a closed substack. 
Note that the induced map $|Z|\to |X|$ of topological spaces is a closed inclusion and $|X|\setminus |Z|$ is quasi-compact.
We write $X_Z^\wedge$, or sometimes simply $\widehat{Z}$ if $X$ is clear from the context, for the prestack defined by
\[
X_Z^\wedge(A)=X(A)\times_{\Map(|\Spec A|, |X|)}\Map(|\Spec A|, |Z|), \quad A\in \calg_\La.
\]
This is nilcomplete and can be represented as $X_Z^\wedge=\colim_a Z_a$ where $Z_a$ range over all closed substacks of $X$, of almost of finite presentation over $\La$, with the same underlying topological space of $Z$.
Therefore is an ind-algebraic stack almost of finite presentation over $\La$, called the formal completion of $X$ along $Z$. Clearly, $X_Z^\wedge$ only depends on the underlying topological space $|Z|$ of $Z$.
\end{example}

\subsubsection{Torsors}\label{SSS:torsors convention}
We let $\tau$ be one of the following topology on $\calg_\La$: Zariski, \'etale, fppf, or fpqc.

Let $H$ be a group prestack over $\La$.  An $H$-equivariant morphism $P\to X$ of prestacks is called an $H$-torsor in the $\tau$-topology if the action of $H$ on $X$ is trivial and for every $\Spec R\to X$, there is a cover $R\to R'$ in the $\tau$-topology such that $P\times_X\Spec R'$ is a trivial, i.e. $H$-equivariantly isomorphic to $\Spec R'\times H$. We let $\bB_\tau H$ denote the prestack of $H$-torsors in $\tau$-topology. By definition, this is a $\tau$-stack. Note that it is possible that $\bB_\tau H$ is a $\tau'$-stack for a finer topology $\tau'$. E.g. $\bB_{\mathrm{Zar}}\GL_n$ is a stack in $fpqc$-topology. If $H$ acts on a prestack $X$, by the quotient $(X/H)_\tau$, we mean the $\tau$-sheafification of the prestack quotient of $X$ by $H$. So $(X/H)_\tau$ is the prestack sending $R$ to an $H$-torsor $P$ over $\Spec R$ (in $\tau$-toplogy) and an $H$-equivariant map $P\to X$. 

When $H$ is a  group stack (i.e. group prestack in \'etale topology), and that $\tau=\et$, we simply call $H$-torsors in the \'etale topology by $H$-torsors, and write $\bB H$ for $\bB_{\et} H$, and if $H$ acts on a stack $X$, we write $X/H$ instead of $(X/H)_{\et}$.

\subsection{Quasi-coherent sheaves}\label{SS: quasi-coherent sheaves}
We recall the general theory of quasi-coherent sheaves on prestacks and specialize the general theory to a particular example (see \Cref{ex: fpqc quotient stack}), which is important in this article.
We refer to \cite{Lurie.SAG} for detailed accounts for some of the statements below. (Although the setting of \emph{loc. cit.} is spectral algebraic geometry, many arguments work in derived algebraic geometry as well.)

Namely, there is a lax symmetric monoidal functor
\begin{equation}\label{eq: theory of quasi-coherent sheaf-1}
\qcoh:  (\prestk_\La)^{\op}\to \lincat_\La,
\end{equation}
defined as the right Kan extension along the full embedding $(\calg_\La)^{\op}\to \prestk_\La$ of the symmetric monoidal functor 
\[
\Mod: \calg_\La\to\lincat_\La,\ A\mapsto \Mod_A.
\]  
Recall that inside $\Mod_A$, there is the smallest idempotent complete stable category $\Perf_A$ containing $A\in \Mod_A$, usually called the category of perfect complexes on $\Spec A$. The functor $\Mod: \calg_\La\to\lincat_\La$ restricts to a functor $\mathrm{Perf}: \calg_\La\to\catid_\La$, and therefore via right Kan extension, gives
\begin{equation}\label{eq: theory of perfect complexes-2}
\mathrm{Perf}: (\prestk_\La)^{\op}\to \catid_\La. 
\end{equation}
Explicitly, for a prestack $X$,
\[
\qcoh(X)= \lim_{A\in (\calg_\La)^\op_{/X}} \Mod_A, \quad \mathrm{Perf}(X)=  \lim_{A\in (\calg_\La)^\op_{/X}} \mathrm{Perf}_A.
\]
Note that $\Qcoh(X)$ has a natural symmetric monoidal structure, and $\mathrm{Perf}(X)$ can be identified with the full subcategory $\Qcoh(X)^d$ of dualizable objects in $\Qcoh(X)$.
For $f: X\to Y$, let $f^*: \qcoh(Y)\to\qcoh(X)$ denote the pullback functor, which restricts to a functor $\mathrm{Perf}(Y)\to \mathrm{Perf}(X)$. 

Recall that for $X\in\algsp^{\qcqs}_\La$, the category $\Qcoh(X)$ is compactly generated, and we have 
\[
\Qcoh(X)^\cpt=\Qcoh(X)^d=\mathrm{Perf}(X).
\] 
But for general prestack $X$, $\Qcoh(X)$ may not be compactly generated and compacts objects may not coincide with perfect complexes. See more discussions below (e.g. \Cref{lem: concentrate stacks}).

By \Cref{ex-sheaf-theory-for-adjoint-factorization}, \eqref{eq: theory of quasi-coherent sheaf-1} extends to a sheaf theory
\begin{equation}\label{eq: theory of quasi-coherent sheaf-2}
\qcoh: \corr(\prestk_\La)_{\mathrm{HR};\all}\to \lincat_\La,
\end{equation}
where $\mathrm{HR}$ is the class of morphisms as in \Cref{rem-additional-base.change.sheaf.theory} \eqref{rem-additional-base.change.sheaf.theory-1}.
It is well-known that every $(f: X\to Y)\in \algsp^{\qcqs}_\La$ belongs to $\mathrm{R}$. Then by \Cref{prop-sheaf-theory-right-Kan-extension-VR class}, for morphism $f:X\to Y$ of prestacks that is representable in $\algsp_\La^{\qcqs}$ (i.e. for every morphism $S\to Y$ with $S\in \algsp_\La^{\qcqs}$, the base change $S\times_YX\in \algsp_\La^{\qcqs}$) belongs to $\mathrm{R}$. However, the class $\mathrm{R}$ also contains certain non-representable morphisms, as we shall see shortly.

We also recall that $\qcoh(X)$ admits a standard $t$-structure. It is defined such that 
\[
\qcoh(X)^{\leq 0}=  \lim_{(\calg_\La)^\op_{/X}} \Mod_A^{\leq 0}.
\] 
This $t$-structure is left complete. By definition for a morphism of prestacks $f:X\to Y$, $f^*$ is left exact, and therefore its (not necessarily continuous) right adjoint $f_*$ is right exact. 

Recall that a morphism $f: X\to Y$ of prestacks is called of finite tor amplitude if it is left $t$-exact up to a finite shift, i.e. there is some integer $N$ such that $f^*$ sends $\Qcoh(Y)^{\geq n}$ to $\Qcoh(X)^{\geq n+N}$ for every $n$. 
We recall that $X$ is called eventually coconnective if $X\to \Spec \bZ$ is of finite tor amplitude, or equivalently
\[
\mO_X\in \Qcoh(X)^+:=\bigcup_n \Qcoh(X)^{\geq n}.
\]
Flat morphisms and quasi-smooth morphisms (to be reviewed later) are of finite tor amplitude. If $Y$ is a smooth (and therefore classical) and $X$ is eventually coconnective, then $f$ is of finite tor amplitude. We recall the following facts.

\begin{proposition}\label{prop: proper finite tor morphism}
Suppose $f:X\to Y$ is a morphism of qcqs algebraic spaces almost of finite presentation. 
\begin{enumerate}
\item\label{prop: proper finite tor morphism-1} If $f$ is proper and of finite tor amplitude, then $f_*$ sends $\Perf(X)$ to $\Perf(Y)$.  
\item\label{prop: proper finite tor morphism-2} If $f$ is finite, then $f$ is of finite tor amplitude if and only if $f_*$ sends $\Perf(X)$ to $\Perf(Y)$.
\item\label{prop: proper finite tor morphism-3} If $g: Y\to Z$ is another smooth morphism of qcqs algebraic spaces. Then $g\circ f$ is of finite tor amplitude if and only if $f$ is of finite tor amplitude. 
\end{enumerate}
\end{proposition}

\begin{proof}
The first statement can be proved as in \cite[Theorem 6.1.3.2]{Lurie.SAG}. (Although \emph{loc. cit.} works in the framework of spectral algebraic geometry, the same argument applies in our setting as well.) The second statement then is clear. For the last statement, we immediately reduce to the case $f$ is closed embedding. Then one can use the second statement to conclude.
\end{proof}

If $X$ is an algebraic stack, the heart $\qcoh(X)^{\heartsuit}$ is the usual abelian category of quasi-coherent sheaves on $X_\cl$, which is a Grothendieck abelian category (by  \cite[\href{https://stacks.math.columbia.edu/tag/0781}{Proposition 0781}]{stacks-project}). In fact, the proof of \emph{loc. cit.} applies in a slightly more general situation, giving the first part of the following lemma.

\begin{lemma}\label{lem: derived category of qcoh v.s. qcoh}
Let $X$ be a stack with representable diagonal (by qcqs algebraic spaces) such that there exists an fpqc cover $p:U\to X$ where $U$ is an algebraic space. Then 
\begin{enumerate}
\item\label{lem: derived category of qcoh v.s. qcoh-1} $\qcoh(X)^{\heartsuit}$ is a Grothendieck abelian category.
\item\label{lem: derived category of qcoh v.s. qcoh-2.5} The $t$-structure of $\Qcoh(X)$ is right complete, and is compatible with filtered colimits.
\item\label{lem: derived category of qcoh v.s. qcoh-2} If $X$ is qcqs, then $\tau^{\geq n}\mO_X\in \qcoh(X)^{\geq n}$ is a compact object in $\qcoh(X)^{\geq n}$.
\item\label{lem: derived category of qcoh v.s. qcoh-3} Suppose that in addition the diagonal is affine and  $U$ is quasi-compact and classical,
then the natural $t$-exact functor $\der(\qcoh(X)^{\heartsuit})^+\to \qcoh(X)^+$ (e.g. see \cite[Remark 1.3.5.23]{Lurie.higher.algebra} for the construction of this functor) is an equivalence.
\end{enumerate}
\end{lemma}
\begin{proof}
All these facts have been proved in literature when $X$ is an algebraic stack, but the same proofs go through in this slightly more general setting. The point is that the category $\qcoh(X)^{?}$, for $?=\heartsuit, \geq -n, +$, can be described as objects in $\qcoh(U)^{?}$ equipped with descent datum.
E.g. Part \eqref{lem: derived category of qcoh v.s. qcoh-1} follows from the same proof as in \cite[\href{https://stacks.math.columbia.edu/tag/0781}{Proposition 0781}]{stacks-project}, Part \eqref{lem: derived category of qcoh v.s. qcoh-2.5} follows from the same proof as in \cite[Proposition 3.1.5.7]{Gaitsgory.Rozenblyum.DAG.vol.I}, Part \eqref{lem: derived category of qcoh v.s. qcoh-2} follows from the same proof as in \cite[Corollary 1.3.17]{Drinfeld.Gaitsgory.finiteness.2013} (see \Cref{lem: push-forward of qcoh for coconnective part} \eqref{lem: push-forward of qcoh for coconnective part-1} below for a relative version), and
Part \eqref{lem: derived category of qcoh v.s. qcoh-3} follows from the same proof as in  \cite[Theorem 3.8]{Lurie.Tannaka} (see also \cite[Proposition 3.2.4.3]{Gaitsgory.Rozenblyum.DAG.vol.I}).
\end{proof}

Similarly, arguments of  \cite[\textsection{2.1}]{Drinfeld.Gaitsgory.finiteness.2013} and \cite[Lemma 4.5]{Hall.Rydh} (which relies on the fact that the $t$-structure on $\qcoh(X)$ is left complete for any $X$) give the following statement.

\begin{lemma}\label{lem: concentrate stacks}
Suppose $X$ is a stack such that the diagonal is representable by qcqs algebraic spaces and such that there is an fpqc morphism $p:U\to X$ with $U$ a qcqs algebraic space
(e.g. $X$ is a qcqs algebraic stack), then $\qcoh(X)^\cpt\subset \mathrm{Perf}(X)$.
In addition, the following are equivalent:
\begin{enumerate}
\item $\mO_X$ is compact;
\item $\qcoh(X)^\cpt=\mathrm{Perf}(X)$;
\item There exists $n$ such that for every $\mF\in \qcoh(X)^{\heartsuit}$, $H^i\rg(X,\mF)=0$ for $i>n$.
\end{enumerate}
\end{lemma}
\quash{\begin{proof}
That $\qcoh(X)^\cpt\subset \mathrm{Perf}(X)$ follows from the fact that $p^*$ admits a continuous right adjoint $p_*$, as explained above. 
It is a general fact that in a $\La$-linear symmetric monoidal category, if the unit is compact, then every dualizable object is compact. Therefore the first two conditions are equivalent. For the equivalence of the first and the third condition, we refer to. 
\end{proof}}

Following \cite{Hall.Rydh}, we call a stack $X$ as above satisfying the above equivalent conditions concentrated. 

\begin{lemma}
Let $f: X\to Y$ be a morphism of stacks. Assume that $Y$ is concentrated. Then $X$ is concentrated if and only if for every affine scheme $S\to Y$, the base change $S\times_YX$ is concentrated.
\end{lemma}
\begin{proof}
If $X$ is concentrated, then for every $S\to Y$, then the morphism $S\times_YX\to X$ is representable by a qcqs algebraic space. So the $*$-pullback sends compact objects to compact objects (as it admits a continuous right adjoint given by $*$-pushforwards). Therefore $\mO_{S\times_YX}$ is compact, showing $S\times_YX$ is concentrated. Conversely, suppose $S\times_YX$ is concentrated is for every affine $Y$-scheme $S$. Note that $(f_S)^*: \Qcoh(S)\to \Qcoh(S\times_YX)$ admits a continuous right adjoint in this case. Therefore using the argument as in \Cref{prop-sheaf-theory-right-Kan-extension-VR class}, $f^*$ admits a continuous right adjoint. This implies that
$X$ is concentrated as $\mO_Y$ is compact.
\end{proof}

The above argument also shows that if $f:X\to Y$ is a morphism of concentrated stacks, then $f^*$ admits a continuous right adjoint. In addition, it satisfies the base change and projection formula, as argued in \Cref{prop-sheaf-theory-right-Kan-extension-VR class}. If we apply \Cref{prop-sheaf-theory-right-Kan-extension-VR class} for the right Kan extension from the category of concentrated stacks to all prestacks, we see that the class $\mathrm{R}$ as in \eqref{eq: theory of quasi-coherent sheaf-2} contain morphisms of prestacks that are representable by concentrated stacks.

\begin{lemma}\label{rem: self duality of qcoh for concentrated stacks}
If $X$ is a concentrated stack over $\La$, then the following are equivalent. 
\begin{enumerate}
\item\label{rem: self duality of qcoh for concentrated stacks-1} The pairing 
\[
\qcoh(X)\otimes_\La \qcoh(X)\to \qcoh(X\times_\La X)\xrightarrow{ (\pi_X)_*(\Delta_X)^*}\Mod_\La
\] 
is a co-unit in the duality datum of $\qcoh(X)$.
\item\label{rem: self duality of qcoh for concentrated stacks-2} $\qcoh(X)$ is dualizable.
\item\label{rem: self duality of qcoh for concentrated stacks-3} For every prestack $Y$, the exterior tensor product $\qcoh(X)\otimes_\La\qcoh(Y)\to \qcoh(X\times_\La Y)$ is an equivalence.
\item\label{rem: self duality of qcoh for concentrated stacks-4} $ \qcoh(X)\otimes_\La\qcoh(X)\xrightarrow{\boxtimes_\La} \qcoh(X\times_\La X)$ is fully faithful and $(\Delta_X)_*\mO_X$ belongs to the essential image of $\boxtimes_\La$.
\end{enumerate}
If the above equivalent conditions hold, then the symmetric monoidal structure on $\Qcoh(X)$ given by the usual tensor product is rigid. 
\end{lemma}
\begin{proof}
Clearly \eqref{rem: self duality of qcoh for concentrated stacks-1} implies \eqref{rem: self duality of qcoh for concentrated stacks-2}. That \eqref{rem: self duality of qcoh for concentrated stacks-2} implies \eqref{rem: self duality of qcoh for concentrated stacks-3} is true for $X$ being any prestack (see  \cite[Proposition 3.3.1.7]{Gaitsgory.Rozenblyum.DAG.vol.I}). Finally if \eqref{rem: self duality of qcoh for concentrated stacks-3} holds, then
\Cref{rem-dualizability in corr} is applicable to $X$ giving \eqref{rem: self duality of qcoh for concentrated stacks-1}.
That \eqref{rem: self duality of qcoh for concentrated stacks-3} and \eqref{rem: self duality of qcoh for concentrated stacks-4} are equivalent follows from \Cref{lem: criterion of tensor product equivalence via diagonal}.

Now assume that the above conditions hold. Then the rigidity of $\qcoh(X)$ follows from \Cref{cor: rigidity of der(X)}. In more details,
since we assume that $X$ is concentrated, the unit of $\qcoh(X)$, which is $\mO_X$ is compact. In addition, the base change isomorphism implies that $(\Delta_X)_*: \qcoh(X)\to \qcoh(X\times_\La X)\cong\qcoh(X)\otimes_\La\qcoh(X)$ is a $\qcoh(X)$-bimodule functor. Therefore $\qcoh(X)$ is rigid. 
\end{proof}

\begin{remark}\label{rem: self-duality on qcoh}
Recall that in the situation in \Cref{rem: self duality of qcoh for concentrated stacks}, $\qcoh(X)$ fits into the discussion of \Cref{ex: duality via Frobenius-structure} (as well as \Cref{ex: duality of rigid monoidal category}). We have
\[
\verd_X^{\qcoh}: \qcoh(X)^{\vee}\cong \qcoh(X).
\]
Recall that $\qcoh(X)^\cpt=\Perf(X)$. Suppose that $\qcoh(X)$ is compactly generated (e.g. $X$ is a qcqs algebraic space). Then $\verd_X^{\qcoh}$ restricts to a functor (see \eqref{eq:duality on compact objects})
\[
\verd_X^{\Perf}:=(\verd_X^{\qcoh})^\cpt: \Perf(X)^{\op}\cong \Perf(X).
\]
This is the functor sending $\mE\in \Perf(X)$ to its $\mO_X$-linear dual $\mE^\vee$.
\end{remark}

\begin{example}\label{ex: fpqc quotient stack}
Let $H$ be a classical affine flat group scheme over $\La$, and let $\bB_{\mathrm{fpqc}}H$ denote the stack of $H$-torsors in fpqc topology. (See \Cref{SSS:torsors convention} for our conventions.)
If $H$ is of finite presentation over $\La$, then $\bB_{\mathrm{fpqc}}H=\bB_{\mathrm{fppf}}H$ is the same as the stack of $H$-torsors in fppf topology and therefore is algebraic (e.g. see \cite[\href{https://stacks.math.columbia.edu/tag/06DC}{Theorem 06DC}]{stacks-project}). In general $\bB_{\mathrm{fpqc}}H$ is not algebraic, but still belongs to the class of stacks considered in \Cref{lem: derived category of qcoh v.s. qcoh}. Let $U$ be a qcqs algebraic space over $\La$ equipped with an $H$-action. Similarly, let $(U/H)_{\mathrm{fpqc}}$ be the quotient stack in fpqc topology. 

We call $\qcoh(\bB_{\mathrm{fpqc}}H)$ the ($\infty$-)category of algebraic representations of $H$. This terminology is justified by the fact that $\qcoh(\bB_{\mathrm{fpqc}}H)^{\heartsuit}$ is naturally identified with the abelian category of algebraic representations of $H$ on ordinary $\La$-modules. By \Cref{lem: derived category of qcoh v.s. qcoh} \eqref{lem: derived category of qcoh v.s. qcoh-3}, $\qcoh(\bB_{\mathrm{fpqc}}H)$ is the left completion of $\der(\qcoh(\bB_{\mathrm{fpqc}}H)^{\heartsuit})$. We have natural functors
\begin{equation}\label{eq: derived category of representations}
\ind\Perf(\bB_{\mathrm{fpqc}}H)\to \der(\qcoh(\bB_{\mathrm{fpqc}}H)^{\heartsuit})\to \qcoh(\bB_{\mathrm{fpqc}}H).
\end{equation}
In general, all of these three categories are different. For instance, the trivial representation of $H$ is always a compact object of $\ind\Perf(\bB_{\mathrm{fpqc}}H)$, but it is compact
in $\qcoh(\bB_{\mathrm{fpqc}}H)$ if and only if $\bB_{\mathrm{fpqc}}H$ is concentrate. E.g. if $H=\bZ/p$ is the constant group over $\La=\bF_p$, then $\bB_{\mathrm{fpqc}}H$ is not concentrated, so the composed functor in \eqref{eq: derived category of representations} is not an equivalence. However, in this case the second functor is still an equivalence and $\qcoh(\bB_{\mathrm{fpqc}}H)$ compactly generated. This follows from \Cref{lem: cpt gen. of derived category of Groth ab cat} below, applying to $c_i=\La[H]$, the group algebra of $H$, regarded as a representation of $H$.
On the other hand, if $H$ is a countable product of $\bZ/p$ and $\La=\bF_p$, then $\der(\qcoh(\bB_{\mathrm{fpqc}}H)^{\heartsuit})$ is not left complete so the second functor is not an equivalence. In this case,  $\qcoh(\bB_{\mathrm{fpqc}}H)$ is not compactly generated.

Now suppose $\qcoh(\bB_{\mathrm{fpqc}}H)$ is concentrated (e.g. if $H$ is an affine algebraic group and $\La$ is a field of characteristic zero). In this case, clearly $(U/H)_{\mathrm{fpqc}}$ is also concentrated.
If we suppose in addition the ring of regular functions on $H$, regarded as a representation of $H$ via left translation, can be written as increasing union of $H$-representations on projective $\La$-modules. (e.g. this is always the case when $\La$ is a Dedekind domain.) Note that this condition implies that finite dimensional $H$-equivariant vector bundles on $\spec \La$ form a collection of generators of $\qcoh(\bB_{\mathrm{fpqc}}H)^{\heartsuit}$. Then all categories in \eqref{eq: derived category of representations} are equivalent, and therefore are compactly generated.
If there is some $\mF\in \Perf((U/H)_{\mathrm{fpqc}})$ such that its pullback to $U$ is a generator of $\qcoh(U)$ (e.g. if $U$ is affine, or if $U$ is quasi-projective and $H$ is an algebraic group), then $\qcoh((U/H)_{\mathrm{fpqc}})=\ind\Perf((U/H)_{\mathrm{fpqc}})$.
\end{example}

The following statement is standard.

\begin{lemma}\label{lem: cpt gen. of derived category of Groth ab cat}
For a Grothendieck abelian category $\bfC^{\heartsuit}$ with a set of generators $\{c_i\}_i$ such that $\mathrm{Ext}^\bullet_{\bfC^{\heartsuit}}(c_i,-)$ has finite cohomological dimension, then its derived category $\der(\bfC^{\heartsuit})$ is left complete and is compactly generated by $\{c_i[n]\}_{i,n}$.
\end{lemma}

\begin{remark}
\label{thm: QCA stacks by DG}
Using the above example and the local structures of algebraic stacks, Drinfeld-Gaitsgory proved (see \cite{Drinfeld.Gaitsgory.finiteness.2013}) that every qcqs algebraic stack $X$ over $\bQ$ with affine stabilizers and finitely presented (classical) inertia (such stack is called QCA in \emph{loc. cit.}) is concentrated and $\qcoh(X)$ is dualizable. 
\end{remark}

By combining the above discussions with computations in \Cref{sec:trace-geometric-trace-main-general}, we obtain the following statement.

\begin{proposition}\label{prop: categorical trace for quasi-coherent}
Let $X$ be a concentrated stack satisfying equivalent conditions in \Cref{rem: self duality of qcoh for concentrated stacks}. In addition, we assume that the diagonal $X\to X\times_\La X$ is affine.
Let $Z$ be a prestack equipped with two morphisms $g_i: Z\to X, \ i=1,2$ so $\qcoh(Z)$ is a $\qcoh(X)$-bimodule. Then the $*$-pullback $\qcoh(Z)\to \qcoh(X\times_{X\times X}Z)$ induces an equivalence
\[
\tr(\qcoh(X), \qcoh(Z))\cong \qcoh(X\times_{X\times X}Z).
\]
\end{proposition}
\begin{proof}
By combining the above discussions with  \Cref{trace.phi.convolution.fully.faithful} (for $f=\id_X: X\to X$) and \Cref{prop-comparison-usual-trace-geo-trac}, we obtain a fully faithful embedding 
with essential image generated by the image of $*$-pullback $\qcoh(Z)\to \qcoh(X\times_{X\times X}Z)$. As $X\times_{X\times X}Z\to Z$ is affine, the $*$-pushforward $\qcoh(X\times_{X\times X}Z)\to \qcoh(Z)$ is conservative so the image of $*$-pullback $\qcoh(Z)\to \qcoh(X\times_{X\times X}Z)$ generates the whole category.
\end{proof}

\begin{example}\label{ex: cat trace of qcoh and fixed point}
Let $X$ be as in \Cref{prop: categorical trace for quasi-coherent}.
If we let $Z=Z_1\times_\La Z_2$ such that $Z_1$ satisfies equivalent conditions in \Cref{rem: self duality of qcoh for concentrated stacks}, then we obtain 
\[
\qcoh(Z_1)\otimes_{\qcoh(X)}\qcoh(Z_2)\cong \qcoh(Z_1\times_XZ_2)
\] 
(see \Cref{fully.faithfulness.geometric.tensor.product}), recovering  (and slightly generalizing) \cite[Theorem 4.7]{ben2010integral}. 

On the other hand, if we let $Z=X$ with $g_1=\id$ and $g_2=\phi$ is an automorphism of $X$, then  we obtain 
\[
\tr(\qcoh(X),\phi)=\qcoh(\mL_\phi(X))
\] 
with the canonical functor $[-]_{\phi}: \qcoh(X)\to \tr(\qcoh(X),\phi)$ identified with the $*$-pullback along $\mL_\phi(X)\to X$ (see \Cref{trace.phi.convolution.fully.faithful}).

Now suppose $h:W\to X$ is morphism and suppose $W$ is equipped with an automorphism $\phi=\phi_W$ such that $h$ is $\phi$-equivariant. Then $\qcoh(W)$ is a $\qcoh(X)$-module.
By \Cref{ex: geo phi-trace class}, we have $[\qcoh(W),\phi]_{\phi}=\mL_\phi(h)_*(\mO_{\mL_\phi(W)})$.
\end{example}

Given the above discussions, we see that the notion of concentrated stacks is very useful if $\La$ is a field of characteristic zero. Unfortunately, when $\La$ is of positive characteristic, many important algebraic stacks, including many classifying stacks of algebraic groups in positive characteristic, are often not concentrated. In particular, the $*$-pushforward of all quasi-coherent sheaves does not behave well for general (qcqs) morphisms between algebraic stacks. The situation is much improved if we restrict our attentions to the bounded below subcategories of quasi-coherent sheaves.

\begin{lemma}\label{lem: push-forward of qcoh for coconnective part}
Consider a Cartesian diagram 
\begin{equation*}\label{eq:base-change-diagram-app}
\begin{tikzcd}
    X' \arrow[r,"g'"]\arrow[d,"f'"] & X\arrow[d,"f"]\\
    Y' \arrow[r,"g"] & Y
\end{tikzcd}
\end{equation*}
of stacks that are as in \Cref{lem: derived category of qcoh v.s. qcoh}. Suppose $f$ (and therefore $f'$) is qcqs (but not necessarily representable by algebraic spaces). 
Let $f_*$ (resp. $(f')_*$ ) be the (not necessary continuous) right adjoint of $f^*$ (resp. $(f')^*$). Suppose $g$ is of finite tor amplitude. Then
\begin{enumerate}
\item\label{lem: push-forward of qcoh for coconnective part-1} $f_*|_{\Qcoh(X)^{\geq n}}$ commutes with filtered colimits, for every $n$.
\item\label{lem: push-forward of qcoh for coconnective part-2} The morphism $g'$ is of finite tor amplitude and the Beck-Chevalley map $g^*\circ f_*\to  (f')_*\circ {g'}^*$ is an isomorphism when restricted to $\Qcoh(X)^+$.
\end{enumerate}
\end{lemma}
\begin{proof}
The proof of  \cite[Corollary 1.3.17]{Drinfeld.Gaitsgory.finiteness.2013} (see also \cite[Proposition 3.2.3.2]{Gaitsgory.Rozenblyum.DAG.vol.I}) works in this generality. \quash{We include a sketch of proof for completeness (as the same type of arguments will be used several times in the sequel).

Let $\varphi:U\to X$ be as in \Cref{lem: derived category of qcoh v.s. qcoh} and let $\varphi_\bullet: U_\bullet \to X$ be the \v{C}ech nerve. Let $\mF=\colim_i \mF_i$  be a filtered colimit in $\Qcoh(X)^{\geq 0}$, and let $\mG\in \Qcoh(Y)^{\geq 0}$. By descent, we have
\begin{eqnarray}
\Map(\mG, f_*(\colim_i\mF_i)) &=& \lim_{\Delta_{\leq N}}\Map( \varphi_n^*f^*\mG, \varphi_n^*\colim_i\mF_i)=\Map(f^*\mG, (\varphi_n)_*\colim_i \lim_{\Delta_{\leq N}} (\varphi_n)^*\mF_i)
\end{eqnarray}}
\end{proof}

\subsection{Coherent sheaves}

Now assume that $\Lambda$ is an (ordinary) regular noetherian ring, e.g. $\La$ is a field or more generally a Dedekind domain.

\subsubsection{Basic definitions and properties}\label{SSS: basic theory of coh}
\begin{definition}
Let $X\in \ArStk_\La^{\aft}$ be an algebraic stack almost of finite presentation over $\La$. Let
$\Coh(X)\subset \qcoh(X)$ denote the full subcategory of coherent sheaves, i.e. those $\mF\in\qcoh(X)$ with finitely many cohomological degrees and with each cohomology sheaf being an ordinary coherent sheaf on $X_\cl$. 
Let $\ind\Coh(X)$ be the ind-completion of $\Coh(X)$. 
\end{definition}

We have some immediate remarks concerning the definition.

\begin{remark}\label{rem: first remark of ind-coh}
\begin{enumerate}
\item\label{rem: first remark of ind-coh-1} Our definition of $\Coh(X)$ is consistent with the definition in \cite{gaitsgory2013.indcoh, Gaitsgory.Rozenblyum.DAG.vol.I}, but is different from the definition in \cite{Lurie.SAG} (which works in the setting of spectral algebraic geometry). Note that according to the definition, $\Perf(X)$ may not belong to $\Coh(X)$. In fact, $\Perf(X)\subset \Coh(X)$ if and only if $X$ is eventually coconnective (e.g. $X$ is classical). When $X$ is a scheme, then  $\Perf(X)=\Coh(X)$ if and only if $X$ is classical and is a regular scheme.
Note, however, there is always a monoidal action of 
\begin{equation}\label{eq: monoidal action of perf on coh}
\Perf(X)\otimes_\La \Coh(X)\to \Coh(X),\quad (\mE,\mF)\mapsto \mE\otimes \mF
\end{equation}
obtained by restriction of the monoidal structure of $\Qcoh(X)$.

\item On the other hand, our definition of $\indcoh(X)$ does not coincide with the category of ind-coherent sheaves as defined and studied in \cite{gaitsgory2013.indcoh, Gaitsgory.Rozenblyum.DAG.vol.I}. The category studied in \emph{loc. cit.} will be denoted as $\mathrm{QC}^!(X)$ later on, following the notation from \cite{ben2017spectral}. See more discussions in \Cref{rem: indcoh vs QC} \eqref{rem: indcoh vs QC-2} below.  There will always be a functor  $\indcoh(X)\to \mathrm{QC}^!(X)$, which is an equivalence when $X$ is an algebraic space, or when $\La$ is a field of characteristic zero and $X$ is an algebraic stack over $\La$ with affine diagonal. But this is not the case in general. When $\La$ is of positive characteristic or mixed characteristic, it is $\indcoh(X)$ rather than $\mathrm{QC}^!(X)$ that is more relevant to this work.

\item Clearly $\Coh(X)$ is idempotent complete so $\indcoh(X)^\cpt=\Coh(X)$.
\end{enumerate}
\end{remark}

The category $\Coh(X)$ inherits a standard $t$-structure from $\Qcoh(X)$, with $\Coh(X)^{\heartsuit}$ being the usual abelian category of coherent sheaves on $X_\cl$. Such $t$-structure extends to an accessible $t$-structure on $\indcoh(X)$ such that $\indcoh(X)^{\leq 0}$ (resp. $\indcoh(X)^{\geq 0}$) is the ind-completion of $\Coh(X)^{\leq 0}$ (resp. $\Coh(X)^{\geq 0}$). Let 
\begin{equation}\label{eq: functor Psi from indcoh to qcoh}
\Psi_X: \indcoh(X)\to \qcoh(X)
\end{equation} 
be the ind-completion of the tautological embedding $\Coh(X)\subset\Qcoh(X)$. It is a $t$-exact functor.

The following crucial statement allows one to reduce some questions about  $\indcoh(X)$  to the questions about $\Qcoh(X)$.
\begin{lemma}\label{lem: coconnective part of indcoh}
The functor $\Psi$ restricts to an equivalence 
\[
\Psi_X^{\geq n}\colon \indcoh(X)^{\geq n}\cong \qcoh(X)^{\geq n}
\] 
for every $n$. Consequently, it restricts to an equivalence 
\[
\Psi_X^+: \indcoh(X)^+\cong \Qcoh(X)^+.
\] 
\end{lemma}
\begin{proof}
See \cite[Proposition 4.1.2.2]{Gaitsgory.Rozenblyum.DAG.vol.I} when $X$ is a (derivied) scheme, but all discussions go through for algebraic spaces almost of finite presentation over $\La$ by virtual of \Cref{lem: approximate coh by perf} below. (Or one can deduce the algebraic space case from the scheme case directly using \'etale descent.)

Next we assume that $X$ is an algebraic stack. We need to show that every object in $\Qcoh(X)^{\geq n}$ is a filtered colimit of objects in $\Coh(X)^{\geq n}$ and that for $\mF\in \Coh(X)^{\geq n}$, $\Map(\mF,-)$ commutes with filtered colimits in $\Qcoh(X)^{\geq n}$. One immediately reduces the first statement to the fact that every ordinary quasi-coherent sheaf on $X_{\cl}$ is a filtered limits of ordinary coherent sheaves, which follows from
\cite[\href{https://stacks.math.columbia.edu/tag/0GRF}{Lemma 0GRF}]{stacks-project}. For the second statement, we may assume that $n=0$. we let $\mG=\colim_i \mG_i$ in $\Qcoh(X)^{\geq 0}$. Let $\varphi: U\to X$ be a smooth cover of $X$, and let $\varphi_\bullet: U_\bullet\to X$ be the \v{C}ech nerve of $f$. 
Then $n$th term $\varphi_n: U_n\to X$ is smooth for every $n$, and therefore $(\varphi_n)^*$ is $t$-exact. Therefore by descent $\Qcoh(X)^{\geq 0}\cong \lim\Qcoh(U_\bullet)^{\geq 0}$. For a positive integer $N$, we will let $\Delta_{\leq N}$ denote the finite category of $N$-truncated simplexes. Note that there is some $N$ such that $\mF\in \Coh(X)^{\geq 0}\cap \Coh(X)^{\leq N}$ so $\Map((f_n)^*\mF, \mG)$ is $N$-truncated for all $n$ and all $\mG\in \Qcoh(U_n)^{\geq 0}$. Therefore
\begin{eqnarray*}
\Map(\mF, \colim \mG_i)&\cong&\lim_{\Delta}\Map((\varphi_n)^*\mF, (\varphi_n)^*\colim_i \mG_i)\\
&\cong& \lim_{\Delta_{\leq N}}\Map((\varphi_n)^*\mF, (\varphi_n)^*\colim \mG_i)\\
&\cong&\lim_{\Delta_{\leq N}}\colim_i\Map((\varphi_n)^*\mF, (\varphi_n)^*\mG_i)\\
& \cong&\colim_i\lim_{\Delta_{\leq N}}\Map((\varphi_n)^*\mF, (\varphi_n)^*\mG_i))\\
&\cong& \colim_i\lim_{\Delta}\Map((\varphi_n)^*\mF, (\varphi_n)^*\mG_i)\cong \colim_i\Map(\mF, \mG_i).
\end{eqnarray*}
The lemma is proved.
\end{proof}

\begin{lemma}\label{lem: approximate coh by perf}
Let $X$ be an algebraic space almost of finite presentation over $\La$. Then for every $\mF\in \Coh(X)$, and every $n$, there is some $\mE\in \Perf(X)$ equipped with a map $\mE\to \mF$ in $\Qcoh(X)$ such that the cofiber of this map belongs $\Qcoh(X)^{\leq n}$.
\end{lemma}
\begin{proof}
The case when $X$ is an affine scheme is clear. The reduction of the general case to the affine case  is contained in \cite[\href{https://stacks.math.columbia.edu/tag/08HP}{Theorem 08HP}]{stacks-project}. (The argument was written in classical algebraic geometry but it works in derived algebraic geometry as well.)
\end{proof}

\begin{remark}
The above discussion says that for any $X\in\ArStk_\La^{\aft}$, the category $\indcoh(X)$ is obtained from $\Qcoh(X)$ by regularization in the sense of \cite[\textsection{6}]{benintegralcoh}.
\end{remark}

We also recall the following statement.
\begin{lemma}\label{lem: coh sheaf closed embedded}
Let $\iota: X\to X'$ be a closed embedding such that the induced closed embedding of the underlying classical stack $X_\cl\to X'_{\cl}$ is defined by a nilpotent ideal. Then the essential image $\iota_*: \indcoh(X)\to \indcoh(X')$ generate $\indcoh(X')$ as $\La$-linear category.
\end{lemma}
\begin{proof}
It is enough to show that $\Coh(X')$ is generated by $i_*\Coh(X)$ as idempotent complete stable categories. As $\Coh(X')$ is generated by $\Coh(X')^{\heartsuit}=\Coh(X'_\cl)^{\heartsuit}$, and $i_*$ is $t$-exact, it is enough to notice that every object in $\Coh(X'_{\cl})^{\heartsuit}$ can be written as successive extensions (in the abelian category) by objects in the essential image of $i_*(\Coh(X_\cl)^{\heartsuit})$.
\end{proof}

\begin{lemma}\label{lem: convergence of indcoh}
Let $X\in\ArStk_\La^{\aft}$. For each $n$, let $X_{\leq n}$ denote its $n$-trunction (see \Cref{rem: truncated stacks}). Then $X_{\leq n}\to X$ is a closed embedding (in particular $X_{\leq n}$ is an algebraic stack), and $X_{\leq n}$ is eventually coconnetive and the natural functor
\[
\colim_n\indcoh(X_{\leq n})\to \indcoh(X)
\]
is an equivalence.
\end{lemma}
\begin{proof}
We note by \Cref{lem: coh sheaf closed embedded}, it is enough to show that for every pair $\mF,\mG\in \Coh(X_{\leq n_0})$ for some $n_0$, $\colim_n \Hom( (i_{n,n_0})_*\mF, (i_{n,n_0})_*\mG)\to \Hom( (i_{n_0})_*\mF, (i_{n_0})_*\mG)$ is an isomorphism
Here $i_{n_0,n}: X_{\leq n_0}\to X_{\leq n}$ and $i_{n_0}: X_{\leq n_0}\to X$ are closed embeddings.
We $X$ is a scheme, this is proved in \cite[Proposition 4.3.4]{gaitsgory2013.indcoh} and \cite[Proposition 4.6.4.3]{Gaitsgory.Rozenblyum.DAG.vol.I}. The case of stacks immediately reduces to the scheme case, as in the proof of \Cref{lem: coconnective part of indcoh}. See also \Cref{prop: open-closed gluing coh} below for a very similar type of argument. 
\end{proof}

There is a monoidal action of $\ind\Perf(X)$ on $\indcoh(X)$, obtained as the ind-extension of \eqref{eq: monoidal action of perf on coh}. Fix $\mF\in \indcoh(X)$, the functor
\[
-\otimes \mF: \ind\Perf(X)\to \indcoh(X)
\]
admits a (not necessarily continuous) right adjoint
\begin{equation}\label{eq: internal hom indcoh as quasi}
\underline\Hom(\mF,-): \indcoh(X)\to \ind\Perf(X),
\end{equation}
where $\underline\Hom(-,-)$ is the functor $\Hom_{\bfC/\bfA}(-,-)$ as defined in \eqref{eq: internal hom object in A} applied to $\bfA=\ind\Perf(X)$ and $\bfC=\indcoh(X)$. 
Note that for every $\mE\in \Perf(X)$, $\mE\otimes\underline\Hom(\mF, \mG)\cong \underline\Hom(\mF,\mE\otimes\mG)$. 
Clearly, if $\mF\in \Coh(X)$, then $\underline\Hom(\mF,-)$ is continuous. It follows that if $\mF\in\Coh(X)$, then for every $\mE\in \ind\Perf(X)$, we have
\begin{equation}\label{eq: qcoh coh tensor internel hom}
\mE\otimes \underline\Hom(\mF,-)\cong \underline\Hom(\mF, \mE\otimes -).
\end{equation}

When $X$ is eventually coconnective, i.e. $\mO_X\in \Coh(X)$, then we denote the functor $-\otimes \mO_X$ as 
\begin{equation}\label{eq: XiX functor}
\Xi_X: \ind\Perf(X)\subset \indcoh(X),
\end{equation} 
which is the natural fully faithful embedding.
In addition, if $X\in \algsp_\La^{\aft}$ so $\ind\Perf(X)=\qcoh(X)$, then $\underline\Hom(\mO_X,-)=\Psi_X$, which is the right adjoint of $\Xi_X$.

\subsubsection{The theory $\indcoh^*$}\label{SSS: star coh sheaf theory}
Our next goal is to construct some three functor formalism of (ind-)coherent sheaves. In fact, we will have two such versions. The first one will be denoted as $\indcoh^*$ and the second one will be denoted as $\indcoh^!$ as introduced later. Along the way, we will also review a few other facts related to the category of (ind-)coherent sheaves.

First, by \Cref{lem: push-forward of qcoh for coconnective part} and \Cref{lem: coconnective part of indcoh}, we have the following.
\begin{lemma}\label{eq: indcoh-*-pushforward}
Suppose we have the Cartesian diagram as in \Cref{lem: push-forward of qcoh for coconnective part} with $X,X',Y,Y'\in \ArStk_\La^{\aft}$ and $g$ of finite tor ampltidue.  
\begin{enumerate}
\item There is a unique $\La$-linear functor 
\[
f^{\indcoh}_*:\indcoh(X)\to \indcoh(Y)
\] 
whose restriction to $\indcoh(X)^+\cong \Qcoh(X)^+$ coincides with the functor
$f_*|_{\Qcoh(X)^+}: \Qcoh(X)^+\to \Qcoh(Y)^+$. 

\item\label{eq: indcoh-*-pushforward-2} There is a unique $\La$-linear functor 
\[
g^{\indcoh,*}:\indcoh(Y)\to \indcoh(Y')
\] 
whose restriction to $\indcoh(X)^+\cong \Qcoh(X)^+$ coincides with the functor
$g^*|_{\Qcoh(Y)^+}: \Qcoh(Y)^+\to \Qcoh(Y')^+$. In this case, $g^{\indcoh,*}$ is the left adjoint of $g^{\indcoh}_*$.

\item\label{eq: indcoh-*-pushforward-3} The Beck-Chevalley map
$g^{\indcoh,*}\circ f^{\indcoh}_*\to  (f')^{\indcoh}_*\circ {g'}^{\indcoh,*}$ is an isomorphism.
\end{enumerate}
\end{lemma}

Note that comparing with the theory of quasi-coherent sheaves, we have the $(\indcoh,*)$-pushforward as a $\La$-linear (in particular continuous) functor  for \emph{any} morphism between algebraic stacks almost of finite presentation over $\La$.  
When $Y=\Spec \La$, we write $f_*^{\indcoh}$ as 
\[
\rg^{\indcoh}(X,-)=(\pi_X)_*^{\indcoh}: \indcoh(X)\to \Mod_\La.
\]
When $f$ belongs to the class $\mathrm{R}$ as in \eqref{eq: theory of quasi-coherent sheaf-2} (e.g. $f$ is representable by qcqs algebraic spaces or more generally by concentrated stacks), we have (essentially by definition)
\[
\Psi_Y\circ f_*^{\indcoh}\cong f_*\circ \Psi_X.
\]

If $f: X\to Y$ is a (representable by algebraic spaces) proper morphism between algebraic stacks almost of finite presentation, then $f_*^{\indcoh}$ sends $\Coh(X)$ to $\Coh(Y)$ and therefore admits continuous right adjoint, 
\begin{equation}\label{eq: exceptional pullback}
f^{\indcoh,!}\footnote{Note that our notation is different from \cite{Gaitsgory.Rozenblyum.DAG.vol.I}, where $f^{\indcoh,!}$ is simply denoted by $f^!$.}: \indcoh(Y)\to\indcoh(X),
\end{equation}
 which in  addition sends $\indcoh(Y)^+$ to $\indcoh(X)^+$ (using the fact that $f^{\indcoh}_*$ has finite cohomological amplitude). Given \Cref{lem: coconnective part of indcoh}, the following base change results can be proved exactly as in \cite[Proposition 3.4.2]{gaitsgory2013.indcoh}, \cite[Proposition 4.5.2.2]{Gaitsgory.Rozenblyum.DAG.vol.I} and \cite[Proposition 7.1.6]{gaitsgory2013.indcoh}.

\begin{lemma}
Suppose we have the Cartesian diagram as in \Cref{lem: push-forward of qcoh for coconnective part} with $X,X',Y,Y'\in \ArStk_\La^{\aft}$ and $f$ representable and proper. Then the following Beck-Chevalley map is an isomorphism
\[
(g')^{\indcoh}_*\circ (f')^{\indcoh,!}\cong f^{\indcoh,!}\circ g^{\indcoh}_*.
\]
If $g$ is of finite tor amplitude, then we have the Beck-Chevalley isomorphism
\[
(g')^{\indcoh,*}\circ f^{\indcoh,!}\cong (f')^{\indcoh,!}\circ g^{\indcoh,*}.
\]
\end{lemma}

We recall the following descent properties of $\indcoh$ for morphisms for morphisms between algebraic spaces. They were proved in \cite{Gaitsgory.Rozenblyum.DAG.vol.I} when $\La$ is a field of characteristic zero (but this assumption is not needed in the proof).
\begin{proposition}\label{prop: descent property of coh}
Let $f: X\to Y$ be a morphism in  $\ArStk_\La^{\aft}$.
\begin{enumerate}
\item\label{prop: descent property of coh-1} If $X,Y\in\algsp_\La^{\aft}$ and if $f$ is proper and surjective (for the underlying topological spaces), then the essential image of $f^{\indcoh}_*$ generates $\indcoh(Y)$, or equivalently $f^{\indcoh,!}$ is conservative. It follows from \Cref{prop-descent-codescent-abstract} that proper surjective morphisms between algebraic spaces are of universal $\indcoh$-codescent.
\item\label{prop: descent property of coh-2}
If $X,Y\in\algsp_\La^{\aft}$ and if $f$ is a smooth covering, we have $\indcoh(Y)=\lim \indcoh(X_\bullet)$, where $X_\bullet$ is the \v{C}ech never of $X\to Y$ and the functors are given by $(\indcoh,*)$-pullbacks. 
\item\label{prop: descent property of coh-3} If $f$ is a smooth covering, we have $\indcoh(Y)^+=\lim \indcoh(X_\bullet)^+$.
\end{enumerate}
\end{proposition}

We will also need the following projection formula proved in \cite[Proposition 3.6.11]{gaitsgory2013.indcoh}. Again, the argument as in  does not make use of assumption that $\La$ is a field of characteristic zero.

\begin{lemma}\label{eq: coh perf adjunction}
Let $f: X\to Y$ be a morphism in $\algsp_\La^{\aft}$ of finite tor amplitude. 
Then we have the following projection formula.
\begin{equation*}
f_*\mE\otimes \mF\cong f^{\indcoh}_*(\mE\otimes f^{\indcoh,*}\mF),\quad \forall \mE\in \ind\Perf(X),\quad \mF\in \indcoh(Y).
\end{equation*}
\end{lemma}

\begin{remark}\label{rem: failure coh surjective proper morphism}
We do note know whether \Cref{eq: coh perf adjunction} holds for $X$ and $Y$ being algebraic stacks. On the other hand,
we know that it is crucial to assume that $X$ and $Y$ are algebraic spaces in \Cref{prop: descent property of coh} \eqref{prop: descent property of coh-1} and \eqref{prop: descent property of coh-2}.
The statements generalize to algebraic stacks with affine diagonal and which are almost of finite presentation over characteristic zero field $\La$. But it could fail in positive characteristic. For example, we assume that $\La$ is an algebraically closed field of positive characteristic and we consider $\pt\to \bB H$, where $H$ is a finite group whose order vanishes in $\La$. We regard $H$ as a constant (and therefore smooth) algebraic group over $\La$. Then $\Coh(\bB H)=\Perf(\bB H)=\crep(H,\La)$ and the $(\indcoh,*)$-pushfoward functor $\indcoh(\pt)\to \indcoh(\bB H)$ is not essentially surjective.

The same example shows that smooth descent could fail for algebraic stacks as well. Indeed, $\indcoh(\bB H)=\ind\crep(H,\La)\neq \rep(H,\La)=\Qcoh(\bB H)$. This will cause difficulties study $\indcoh(X)$ for $X$ being general algebraic stacks.
However, it is  easy to see that $\indcoh(X)$ satisfies Zariski descent with respect to $(\indcoh, *)$-pullbacks. (The same proof of \cite[Proposition 4.4.2.2]{Gaitsgory.Rozenblyum.DAG.vol.I} applies.)
\end{remark}

\begin{proposition}\label{lem: coherent Kunneth}
Let $X,Y\in \ArStk^{\aft}_\La$.
\begin{enumerate}
\item\label{lem: coherent Kunneth-1} The exterior product 
\begin{equation*}
\boxtimes\colon  \Qcoh(X) \otimes_\La \Qcoh(Y) \to \Qcoh(X\times_{\La} Y).
\end{equation*}
sends coherent sheaves to coherent sheaves, and
induces a fully faithful embedding
\begin{equation*}
\boxtimes\colon  \indcoh(X) \otimes_\La \indcoh(Y) \to \indcoh(X\times_{\La} Y),
\end{equation*}
which admits a $\La$-linear right adjoint $\boxtimes^R$.
\item\label{lem: coherent Kunneth-2}  Let $f: X\to Z$ be a morphism in $\ArStk^{\aft}_\La$. Then the following diagram is commutative 
 \[
 \xymatrix{
 \indcoh(X)\otimes_\La\indcoh(Y)\ar^-{\boxtimes}[r]\ar_{f^{\indcoh}_*\otimes \id}[d] & \indcoh(X\times_\La Y)\ar^{(f\times \id)^{\indcoh}_*}[d]\\
 \indcoh(Z)\otimes_\La\indcoh(Y)\ar^-{\boxtimes}[r]  & \indcoh(Z\times_\La Y).
 }\]
If $f$ is of finite tor amplitude, then the following diagram is commutative 
 \[
 \xymatrix{
 \indcoh(X)\otimes_\La\indcoh(Y)\ar^-{\boxtimes}[r]& \indcoh(X\times_\La Y)\\
 \indcoh(Z)\otimes_\La\indcoh(Y)\ar^-{\boxtimes}[r] \ar^{f^{\indcoh,*}\otimes \id}[u]  & \indcoh(Z\times_\La Y) \ar_{(f\times \id)^{\indcoh,*}}[u].
 }\]

\end{enumerate}
\end{proposition}
\begin{proof}
We start with Part \eqref{lem: coherent Kunneth-1}.
We need to show that if $\mF\in\Coh(X)$ and $\mF'\in\Coh(Y)$, then $\mF\boxtimes_\La\mG\in \Coh(X\times_\La Y)$. We may assume that $\mF\in\Coh(X)^{\heartsuit}=\Coh(X_\cl)^{\heartsuit}$. Therefore, we may assume that $X$ is classical. Similarly, we may assume that $Y$ is classical. Now both $X\to \Spec \La$ and $Y\to\Spec \La$ are of finite tor dimension. The first statement follows.

For the fully faithfulness statement, we first assume that $X$ and $Y$ are algebraic spaces. In this case, the statement can be proved as in \cite[Proposition 4.6.3.4 (a)]{Gaitsgory.Rozenblyum.DAG.vol.I}, using \Cref{lem: approximate coh by perf}. (A similar type argument is given in \Cref{lem: kunneth for indfg sheaf} below.)

Next we assume that $X$ is an algebraic stack. We need to show that the following map
\[
\Hom(\mF_1,\mF_2)\otimes_\La \Hom(\mG_1,\mG_2)\to \Hom(\mF_1\boxtimes \mG_1, \mF_2\boxtimes \mG_2)
\]
is an isomorphism. Without loss of generality, we may assume that $\mF_1\in \Coh(X)^{\leq 0}$, $\mF_2\in \Coh(X)^{\geq 0}$ and $\mG_1\in \Coh(Y)^{\leq 0}$ and $\mG_2\in \Coh(Y)^{\geq 0}$. It is enough to show that for each $n$, the above map becomes an isomorphism after applying truncation $\tau^{\leq n}$.
We fix such $n$.

Note that as $\Hom(\mG_1,\mG_2)\in \Mod_\La^{\geq 0}$ and $\La$ is regular noetherian, there is some $m$ large enough such that
\[
\tau^{\leq n}(\tau^{\leq m}M\otimes_\La \Hom(\mG_1,\mG_2))\to \tau^{\leq n}(M\otimes_\La \Hom(\mG_1,\mG_2)) 
\]
is an isomorphism for every $M\in \Mod_\La$.

Let $\varphi: U\to X$ be a smooth atlas with $U\in \algsp^{\aft}_\La$. Let $\varphi_\bullet: U_\bullet\to X$ be the \v{C}ech nerve of $\varphi$. By smooth descent, we have $\Hom(\mF_1,\mF_2)=\lim_{\Delta}\Hom(\varphi_j^*\mF_1,\varphi_j^*\mF_2)$.
 Note that for every $m$, there is some $N$ large enough such that 
\[
\tau^{\leq m}\Hom(\mF_1,\mF_2)=\tau^{\leq m}\lim_{\Delta_{\leq N}}\Hom(\varphi_j^*\mF_1,\varphi_j^*\mF_2).
\]
Similarly, for every $m$, there is some $N$ large enough such that
\[
\tau^{\leq m}\Hom(\mF_1\boxtimes\mG_1,\mF_2\boxtimes\mG_2)=\tau^{\leq m}\lim_{\Delta_{\leq N}}\Hom(\varphi_j^*\mF_1\boxtimes\mG_1,\varphi_j^*\mF_2\boxtimes\mG_2).
\]

Now, we choose $m\gg n$ large enough, and $N$ large enough (depending on m). Then we have
\begin{eqnarray*}
\tau^{\leq n}(\Hom(\mF_1,\mF_2)\otimes_\La \Hom(\mG_1,\mG_2)) & \cong & \tau^{\leq n}(\tau^{\leq m}\Hom(\mF_1,\mF_2)\otimes_\La \Hom(\mG_1,\mG_2))\\
                                                                                                            & \cong & \tau^{\leq n}(\tau^{\leq m}\lim_{\Delta_{\leq N}}\Hom(\varphi_j^*\mF_1,\varphi_j^*\mF_2)\otimes_\La \Hom(\mG_1,\mG))\\
                                                                                                            & \cong & \tau^{\leq n}(\lim_{\Delta_{\leq N}}\Hom(\varphi_j^*\mF_1,\varphi_j^*\mF_2)\otimes_\La \Hom(\mG_1,\mG))\\
                                                                                                            & \cong & \tau^{\leq n}\lim_{\Delta_{\leq N}}\Hom(\varphi_j^*\mF_1\boxtimes \mG_1,\varphi_j^*\mF_2\boxtimes \mG_2)\\
                                                                                                            & \cong & \tau^{\leq n}\Hom(\mF_1\boxtimes \mG_1,\mF_2\boxtimes \mG_2),
\end{eqnarray*}   
as desired.                                                                                                         
     
Repeating the argument, we may also allow $Y$ to be an algebraic stack.    

Next we prove Part \eqref{lem: coherent Kunneth-2}.  It is enough to show that $f^{\indcoh}_{*}\mF \boxtimes \mG\cong (f\times \id)^{\indcoh}_{*}(\mF\boxtimes \mG)$ when $\mF\in \Coh(X)$ and $\mG\in\Coh(Y)$. In this case, all involved sheaves are in the bonded from below subcategories and the desired statement follows from the corresponding statement for quasi-coherent sheaves. The case of $(\indcoh, *)$-pullback for morphism of finite tor amplitude is proved similarly.
\end{proof}

Recall the category $\indarstk_\La^{\aft}$ of ind-Artin stacks almost of finite presentation over $\La$.

By  \Cref{lem: push-forward of qcoh for coconnective part}, the class of morphisms that are of finite tor dimension is weakly stable in $\ArStk_\La^{\aft}$. We will denote by $\ftor$ the class morphisms in $\indarstk_\La^{\aft}$ that are representable by algebraic stacks and of finite tor dimension.
Now we are ready to state the first version of $3$-functor formalism for coherent sheaves.
\begin{theorem}\label{thm: indcoh star pullback version} 
There is a sheaf theory
\[
\indcoh^*: \corr(\indarstk_\La^{\aft})_{\all; \ftor}\to \lincat_\La,
\]
which sends $X$ to $\indcoh^*(X)=\indcoh(X)$ and
$Y\xleftarrow{g} Z\xrightarrow{f} X$ to $f_*^{\indcoh}\circ g^{\indcoh,*}$.
The class of morphisms that are representable and proper  satisfy \Cref{assumptions.base.change.sheaf.theory.V}.  
On the other hand the class $\ftor$ satisfy \Cref{assumptions.base.change.sheaf.theory.H}. 
\end{theorem}

\begin{proof}
We start from $\Qcoh: (\ArStk_\La^{\aft})^{\op}\to \lincat_\La\to \cat$, where the first lax symmetric monoidal functor is the restriction of \eqref{eq: theory of quasi-coherent sheaf-2} along $(\ArStk_\La^{\aft})^{\op}\subset \corr(\prestk_\La)_{\mathrm{R};\all}$. Passing to the right adjoint, we obtain a lax symmetric monoidal functor $\Qcoh_*: \corr(\ArStk_\La^{\aft})_{\all;\iso}\to \cat$. 
Via the symmetric monoidal Grothendieck construction (see \Cref{rem: Sheaf theory via Grothendieck construction}), we obtain a coCartesian fibration $\corr^{\Qcoh_*}(\ArStk^{\aft}_\La)_{\all;\iso}\to \corr(\algsp^{\aft}_\La)_{\all;\iso}$ which is symmetric monoidal. 
The full subcategory $\corr^{\Qcoh^+}(\ArStk^{\aft}_\La)_{\all;\iso}$ consisting of $(X,\mF)$ with $\mF\in\Qcoh(X)^+$ is a symmetric monoidal subcategory by \Cref{lem: coherent Kunneth} and $\corr^{\Qcoh^+}(\ArStk^{\aft}_\La)_{\all;\iso}\to \ArStk^{\aft}_\La$ is still coCartesian by \Cref{lem: push-forward of qcoh for coconnective part}. Then we obtain a lax symmetric monoidal functor
\[
\corr(\ArStk^{\aft}_\La)_{\all;\iso}\to \cat,\quad X\mapsto \Qcoh(X)^+=\indcoh(X)^+.
\]
On the other hand, at the level of homotopy categories, we have $\corr(\ArStk^{\aft}_\La)_{\all;\iso}\to \mathrm{h}\corr(\ArStk^{\aft}_\La)_{\all;\iso}\to \mathrm{h}\lincat^{t,+}$ sending $X$ to $\indcoh(X)$ by \Cref{eq: indcoh-*-pushforward} and \Cref{lem: coherent Kunneth}.
Using \Cref{lemma: construct functors from homotopy cat} and \Cref{lem: lax symmetric monoidal from lincat with t to cat}, we can combine the above two constructions into a lax symmetric monoidal functor
\[
\corr(\ArStk^{\aft}_\La)_{\all;\iso}\to \lincat, \quad (f: X\to Y)\mapsto (f_*^{\indcoh}: \indcoh(X)\to \indcoh(Y)).
\]
Taking the operadic left Kan extension along $ \ArStk^{\aft}_\La\to \indarstk^{\aft}_\La$, we obtain (using \Cref{prop-dual-sheaf-theory-right-Kan-extension}) a sheaf theory
\[
\corr(\indarstk^{\aft}_\La)_{\all;\iso}\to \lincat, \quad (f: X\to Y)\mapsto (f_*^{\indcoh}: \indcoh(X)\to \indcoh(Y)).
\]

Next \Cref{eq: indcoh-*-pushforward} \eqref{eq: indcoh-*-pushforward-2} \eqref{eq: indcoh-*-pushforward-3} obviously generalize to the case $g\in \ftor$ (i.e. $g$ is representably by algebraic stacks and is of finite tor amplitude).
Then applying (a variant of) \Cref{ex-sheaf-theory-for-adjoint-factorization} (see \Cref{prop-sheaf-theory-for-adjoint-factorization} and \Cref{rem-sheaf-theory-for-adjoint-factorization} \eqref{rem-sheaf-theory-for-adjoint-factorization-3}), we obtain
\[
\corr(\indarstk^{\aft}_\La)_{\all;\ftor}\to \lincat, 
\]
as desired.
\end{proof}

Let $X\in \ArStk^{\aft}_\La$ and let $\hat\imath: \widehat{Z}\to X$ be the formal completion of $X$ along a closed subset $|Z|\subset |X|$  (see \Cref{ex: formal completion}). Recall that we may write $\widehat{Z}=\colim_a Z_a$ for $\imath_a: Z_a\to X$ closed embedding with $Z_a$ almost of finite presentation over $\La$. Note that by definition
\[
\indcoh(\widehat{Z})=\colim_a \indcoh(Z_a)
\]
with transitioning functors given by $(\indcoh,*)$-pushforwards. 
The functor
\[
(\hat{\imath})^{\indcoh}_*:  \indcoh(\widehat{Z})\to \indcoh(X)
\] 
preserves compact objects. Its continuous right adjoint is denoted as $(\hat{\imath})^{\indcoh,!}$. The following statement is well-known in the classical algebraic geometry. 
\begin{proposition}\label{prop: open-closed gluing coh}
Let $U\subset X$ be the open complement. Then
\[
\indcoh(\widehat{Z})\to \indcoh(X)\to \indcoh(U)
\]
is a localization sequence (in the sense of \Cref{rem:localization sequence}).
\end{proposition}
\begin{proof}
The essential point is to prove that if $\mF_1,\mF_2\in \Coh(\widehat{Z})$, then
\[
\Map(\mF_1,\mF_2)\cong \Map((\hat{\imath})^{\indcoh}_*\mF_1, (\hat{\imath})^{\indcoh}_*\mF_2).
\]
For a proof when $X$ is a derived scheme, see \cite[Proposition 7.4.5]{Gaitsgory.Rozenblyum.indscheme}. (Although \emph{loc. cit.} assumes that ground ring is a field of characteristic zero, such assumption is not needed in the proof.) 
We now assume that $X$ is an algebraic stack.

We suppose $\mF_i=(\imath_a)^{\indcoh}_*\mF'_i$ for some $a$. So $\hat{\imath}^{\indcoh}_{*}\mF_i=(\imath_a)_*\mF'_i$.
For $a'>a$, let $\imath_{a,a'}$ denote the corresponding closed embedding. Note that there is some $N$ such that $\Map((\imath_{a,a'})^{\indcoh}_*\mF'_1, (\imath_{a,a'})^{\indcoh}_*\mF'_2)$ is $N$-truncated for every $a'$. Let $V\to X$ be a smooth atlas, and let $\varphi_{n,a'}:V_{n,a'}\to Z_{a'}$ be the preimage of $Z_a$ in the $n$th term $V_n$ of the \v{C}ech nerve of the cover. Let $\imath_{n,a,a'}$ be the closed embedding from $V_{n,a}$ to $V_{n,a'}$, and let $\imath_{n,a}: V_{n,a}\to V_n$. Finally let $\hat{\imath}_n$ be the formal embedding of the preimage of $\widehat{Z}$ to $V_n$. Then by base change, we have
\begin{eqnarray*}
\Map(\mF_1,\mF_2)&=& \colim_{a'>a}\Map((\imath_{a,a'})_*\mF'_1, (\imath_{a,a'})_*\mF'_2)\\
                                &=& \colim_{a'>a}\lim_{\Delta_{\leq N}}\Map((\imath_{n,a,a'})_* (\varphi_{n,a})^*\mF'_1, (\imath_{n,a,a'})_*(\varphi_{n,a})^*\mF'_2)\\
                                &=& \lim_{\Delta_{\leq N}}\colim_{a'>a}\Map((\imath_{n,a,a'})_* (\varphi_{n,a})^*\mF'_1, (\imath_{n,a,a'})_*(\varphi_{n,a})^*\mF'_2)\\ 
                                &=& \lim_{\Delta_{\leq N}}\Map((\imath_{n,a})_* (\varphi_{n,a})^*\mF'_1, (\imath_{n,a})_*(\varphi_{n,a})^*\mF'_2)\\
                                &=& \Map((\varphi_n)^*(\imath_a)_* \mF'_1,     (\varphi_n)^*(\imath_a)_* \mF'_2)                    
\end{eqnarray*}
as desired.
\end{proof}

One of applications of this result (together with \Cref{lem: coh sheaf closed embedded}) is the following.
\begin{corollary}\label{cor: tensor equiv coh}
The exterior tensor product functor from \Cref{lem: coherent Kunneth} is an equivalence if $X$ and $Y$ admit a finite filtration $X=X_0\supset X_1\supset X_2\supset\cdots$ and  $Y=Y_0\supset Y_1\supset Y_2\supset\cdots$ by closed substacks such that $\indcoh((X_i\setminus X_{i+1})_{\red})\otimes_\La \indcoh((Y_i\setminus Y_{i+1})_{\red})\to \indcoh((X_i\setminus X_{i+1})_{\red}\times_\La (Y_j\setminus Y_{j+1})_{\red})$ is essentially surjective. 
\end{corollary}

Together with \Cref{lem: criterion of tensor product equivalence via diagonal} and the following result, we get the equivalence of $\boxtimes$ in many cases.
\begin{proposition}\label{prop: tensor equiv coh}
The exterior tensor product functor from \Cref{lem: coherent Kunneth} is an equivalence in the following situations.
\begin{enumerate}
\item\label{prop: tensor equiv coh-1} $X,Y\in \algsp_\La^{\aft}$, $X$ is smooth over $\La$ and $Y$ is regular. 
\item\label{prop: tensor equiv coh-2} $X=Y=\bB G$, where $G$ is a smooth affine algebraic group over a field $\La$. 
\end{enumerate}
\end{proposition}

\begin{proof}
In the first case, we note that $X\times_\La Y$ is regular and the statement follows from the theory of quasi-coherent sheaves \Cref{rem: self duality of qcoh for concentrated stacks}.

In the second case, by  \Cref{lem: criterion of tensor product equivalence via diagonal}, it is enough to show that the ring of regular functions $\mO_G$ on $G$, regarded as a $G\times G$-representation, belongs to $\indcoh(\bB G)\otimes \indcoh(\bB G)$. But this is well-known: $\mO_G$ admits an increasing filtration with associated graded being $V_1 \boxtimes V_2$, where $V_1,V_2$ are representations of $G$. 
\end{proof}

\begin{remark}\label{rem: failure of tensor product}
Suppose the base $\La$ is excellent, and $X,Y\in\algsp^{\aft}_\La$. Then the reduced (and therefore classical) subspace $Y_{\red}$ admits an open dense regular subscheme.
Therefore, by \Cref{cor: tensor equiv coh} and \Cref{prop: tensor equiv coh}  \eqref{prop: tensor equiv coh-1}, the exterior tensor product functor from \Cref{lem: coherent Kunneth} is an equivalence if $X$ admits a finite filtration $X=X_0\supset X_1\supset X_2\supset\cdots$ with each $(X_i\setminus X_{i+1})_{\red}$ is smooth over $\La$.
If $\La$ is a perfect field, this assumption always holds, giving \cite[Proposition 4.6.3.4 (b)]{Gaitsgory.Rozenblyum.DAG.vol.I}. However if $\La$ is not perfect, or if $\La$ is not a field, the exterior tensor product functor is in general not an equivalence, even for  $X,Y\in \algsp_\La^{\aft}$.
\end{remark}

We state another result based on ideas of \Cref{lem: criterion of tensor product equivalence via diagonal}.

\begin{lemma}\label{lem: a criterion of generation under star pullback of perf}
Let $f: X\to Y$ be a morphism of finite tor amplitude between algebraic stacks almost of finite presentation over $\La$. Then $\indcoh(X)$ is generated by $f^{\indcoh,*}\indcoh(Y)$ as idempotent complete $\La$-linear category if $(\Delta_{X/Y})_*\mO_X\in\indcoh(X\times_YX)$ is contained in the idempotent complete subcategory generated by the essential images of the $*$-pullback $\indcoh(Y)\to \indcoh(X\times_YX)$.
\end{lemma}
\begin{proof}
The ideal of proof  is similar to \Cref{lem: criterion of tensor product equivalence via diagonal}.  For $\mK\in \indcoh(X\times_YX)$, we consider the functor $F_\mK(-)= \pr_*^{\indcoh}(\pr_1^{\indcoh,*}(-)\otimes \mK): \indcoh(X)\to\indcoh(X)$. If $\mK$ is the $*$-pullback of some object $\mK'\in \indcoh(Y)$, then by projection formula $F_\mK(\mF)\cong f^{\indcoh,*}(f_*^{\indcoh}\mF\otimes \mK)$ belongs to the subcategory of $\indcoh(X)$ generated by $f^{\indcoh,*}(\indcoh(Y))$. On the other hand, $F_{(\Delta_{X/Y})_*\mO_X}$ is the identity functor. The lemma follows by combining these two considerations.
\end{proof}

\subsubsection{The theory $\indcoh^!$}
Our next goal is to construct the exceptional pullback functor. For this purpose, we need to construct another sheaf theory for ind-coherent sheaves.

Recall that the classical Nagata compactification theorem says that  every separated finite type morphism $f: X\to Y$ between classical qcqs algebraic spaces $X$ and $Y$  admits a factorization $X\stackrel{j}{\hookrightarrow} \overline{X}\xrightarrow{\bar{f}} Y$ with $j$ a quasi-compact open embedding and $\bar{f}$ proper (\cite[\href{https://stacks.math.columbia.edu/tag/0F4D}{Theorem 0F4D}]{stacks-project}).  
We have the following derived analogue.
\begin{lemma}\label{lem:derived Nagata}
Suppose $f: X\to Y$ is a separated morphism in $\algsp^{\qcqs}_\La$. Then $f$ factors as $X\stackrel{j}{\hookrightarrow} \overline{X} \xrightarrow{\bar{f}} Y$ with $j$ a quasi-compact open embedding and $\bar{f}$ proper.
\end{lemma}
\begin{proof}We may factors $f_\cl$ as $X_\cl\stackrel{j_\cl}{\hookrightarrow} \overline{X}_\cl\xrightarrow{\bar{f}_\cl} Y_\cl$ with $j_{\cl}$ a quasi-compact open embedding and  $\bar{f}$ proper. Let $\overline{X}:=X\sqcup_{X_\cl} \overline{X}_\cl$. Then $f$ factorizes as claimed.
\end{proof}

We have mentioned in \eqref{eq: exceptional pullback} that for representable proper morphism $f:X\to Y$, there is the exceptional pullback functor $f^{\indcoh,!}$. We now extend it to more general morphisms.

Recall the class of morphisms from \Cref{def-ind-aft-morphism}.

\begin{theorem}\label{lem: indcoh separate scheme} 
There is a sheaf theory
\[
\indcoh^{!}:  \corr(\indarstk_\La^{\aft})_{\ind\aft;\all}\to \lincat_\La,
\]
which  sends $X$ to $\indcoh^!(X)=\indcoh(X)$ and
a correspondence $X\xleftarrow{g} Z\xrightarrow{f} Y$ to the functor 
\[
f_*^{\indcoh}\circ g^{\indcoh,!}: \indcoh(X)\to \indcoh(Y),
\] 
such that if $g$ is an open embedding then $g^{\indcoh,!}$ is the left adjoint of $g^{\indcoh}_{*}$ and when $g$ is ind-proper, $g^{\indcoh,!}$ is the right adjoint of $g^{\indcoh}_*$. 
\end{theorem}
\begin{proof}
We first restrict the sheaf theory from \Cref{thm: indcoh star pullback version} to the category $\algsp_\La^{\aft}\subset \corr(\ArStk_\La^{\aft})_{\all;\ftor}$ of algebraic spaces almost of finite presentation over $\La$.

As mentioned before, open embeddings are $0$-truncated.
Then by (a variant of) \Cref{ex-sheaf-theory-for-adjoint-factorization} (see \Cref{prop-sheaf-theory-for-adjoint-factorization} and \Cref{rem-sheaf-theory-for-adjoint-factorization} \eqref{rem-sheaf-theory-for-adjoint-factorization-3}), we obtain 
\[
\indcoh^{!}: \corr(\algsp_\La^{\aft})_{\all;\sep} \to \lincat_\La,\quad X\mapsto \indcoh(X)
\]
such that if $g$ is an open embedding then $g^{\indcoh,!}=g^{\indcoh,*}$ is the left adjoint of $g^{\indcoh}_{*}$ and when $g$ is proper, $g^{\indcoh,!}$ is the right adjoint of $g^{\indcoh}_*$.  Here $\sep$ denote the class of separated morphisms. 
To continue, one needs some basic properties of this exceptional pullback functor, which are summarized in \Cref{prop: shrek pullback of finite tor amplitude morphism} below. Although they are stated for stacks, at the current stage we only need these properties for separated morphisms between algebraic spaces. In particular, for an \'etale morphism $g$, $g^{\indcoh,*}=g^{\indcoh,!}$.

Then one can further extend the domain of the functor by \Cref{prop-sheaf-theory-enlarge-H}. Namely, for every $f:Y\to Z$, we can find an \'etale cover $U\to Y$ with $U$ affine. Then $U\to X$ is universal $\indcoh^{!}$-descent and $U\to Y$ and $U\to Z$ separated. Then assumptions of  \Cref{prop-sheaf-theory-enlarge-H} hold, and we have an extension
\begin{equation}\label{eq: indcoh for algsp}
\indcoh^{!}: \corr(\algsp_\La^{\aft}) \to \lincat_\La.
\end{equation}
Therefore, $g^{\indcoh,!}$ is defined for any morphism between algebraic spaces almost of finite presentation over $\La$. Again we have \Cref{prop: shrek pullback of finite tor amplitude morphism}, now for any morphisms between algebraic spaces $X$ and $Y$.

Note that both $(\indcoh,*)$-pushforwards and $(\indcoh,!)$-pullbacks preserve the bounded from below subcategories. We thus can consider
\begin{equation*}\label{eq: indcoh for indSch}
\indcoh^{!,+}: \corr(\algsp_\La^{\aft}) \to \cat,
\end{equation*} 
sending $X$ to $\indcoh(X)^+$.
Next we can apply \Cref{prop-sheaf-theory-right-Kan-extension} to obtain
\[
\indcoh^{!,+}:  \corr(\ArStk_\La^{\aft})_{\mathrm{rp};\all}\to \cat,
\]
via right Kan extension. Here $\mathrm{rp}$ denotes the class of morphisms between prestacks almost of finite presentation that are representable in algebraic spaces. For a stack $X$ with a smooth atlas $U\to X$, let $U_\bullet$ be the corresponding \v{C}ech cover. Then we have
\[
\indcoh(X)^{!,+}=\lim\indcoh(U)^{!,+}.
\]
Using \Cref{prop: shrek pullback of finite tor amplitude morphism} below for algebraic spaces, we see that $\indcoh^{!,+}(X)$ is canonically equivalent to $\indcoh^{*,+}(X)=\indcoh(X)^+$. In addition, if $f: X\to Y$ is a representable morphism between stacks, then corresponding $*$-pushforward between bounded from below subcategories is nothing but the restriction of the previously defined functor $f_*^{\indcoh}$.
On the other hand, we have $g^{\indcoh,!}: \indcoh(Y)^+\to \indcoh(X)^+$ for a morphism $g: X\to Y$ between stacks. By restricting it to $\Coh(Y)\subset \indcoh(Y)^+$ and then ind-completion, we see that $g^{\indcoh,!}$ extends to a functor $\indcoh(Y)\to \indcoh(X)$. Therefore at the homotopy category level, we have $\mathrm{h}\corr(\ArStk_\La^{\aft})_{\mathrm{rp};\all}\to \mathrm{h}\lincat_\La^{t,+}$ sending $X$ to $\indcoh(X)$ equipped with the natural $t$-structure.
Then we can argue as in \Cref{thm: indcoh star pullback version} by using \Cref{lemma: construct functors from homotopy cat} and \Cref{lem: lax symmetric monoidal from lincat with t to cat} to obtain
\[
\indcoh^{!}:  \corr(\ArStk_\La^{\aft})_{\mathrm{rp};\all}\to \lincat_\La.
\]
Note that if $g$ is proper, then $g_*^{\indcoh}$ is the left adjoint of $g^{\indcoh,!}$, and when $g$ is an open embedding, then $g_*^{\indcoh}$ is the right adjoint of $g^{\indcoh,!}$.

Finally,  we can then apply \Cref{prop-sheaf-theory-non-full-right-Kan-extension} to further extend the theory to
\[
\indcoh^{!}:  \corr(\indarstk_\La^{\aft})_{\ind\aft;\all}\to \lincat_\La.
\]
Namely, we can let $\bfC_1=\ArStk_\La^{\aft}$, $\verti_1=\mathrm{rp}$, $\horiz_1=\all$, and $\bfC_2=\indarstk_\La^{\aft}$, $\verti_2=\ind\aft$ and $\horiz_2=\all$, and let $\mathrm{S}_1$ be the class of closed embeddings in $\ArStk_\La^{\aft}$.
Note that if $X=\colim X_i$ is a presentation of $X$ as an ind-algebraic stack, then
\[
\indcoh^!(X)=\colim_i \indcoh^!(X_i)=\lim_i\indcoh^!(X_i)
\]
where in the colimit the transitioning functors are given by $(\indcoh,*)$-pushforwards and in the limit the transitioning functors are given by $(\indcoh,!)$-pushbacks. 
In particular, $\indcoh^!(X)=\indcoh(X)$, as desired. In addition, if $g$ is ind-proper, then $g_*^{\indcoh}$ is the left adjoint of $g^{\indcoh,!}$, and when $g$ is an open embedding, then $g_*^{\indcoh}$ is the right adjoint of $g^{\indcoh,!}$.
\end{proof}

We have the following properties of the exceptional pullback functor.
\begin{proposition}\label{prop: shrek pullback of finite tor amplitude morphism}
Let $f:X\to Y$ be a morphism of algebraic stacks almost of finite presentation over $\La$. Then
\begin{enumerate}
\item  $f^{\indcoh,!}$ sends $\indcoh(Y)^+\to \indcoh(X)^+$. If $f:X\to Y$ is of finite tor amplitude, then $f^{\indcoh,!}: \indcoh(Y)\to \indcoh(X)$ restricts to a functor $f^{\indcoh,!}: \Coh(Y)\to \Coh(X)$. 
\item If $f: X\to Y$ is smooth, of relative dimension $d$, then
\begin{equation}\label{eq: shrek pullback star pullback relative cotangent}
f^{\indcoh,!}(-)\cong  \Sym^d(\bL_{X/Y}[1])\otimes f^{\indcoh,*}(-).
\end{equation}
Here $\bL_{X/Y}$ is the relative cotangent complex of $f$, to be reviewed in \Cref{SSS: cotangent complex} below, and $\Sym^d(\bL_{X/Y}[1])$ is its top exterior power shifted to degree $-d$, which is an invertible $\mO_X$-module. 
In particular, if $f: X\to Y$ is \'etale, then $f^{\indcoh, !}\cong f^{\indcoh,*}$.

In fact, the isomorphism \eqref{eq: shrek pullback star pullback relative cotangent} holds if $f$ is a quasi-smooth morphism.
\end{enumerate}
\end{proposition}

\begin{remark}\label{rem: indcoh vs QC}
\begin{enumerate}
\item Note that we have seen that $\indcoh^!(X)=\indcoh^*(X)$ for $X\in \algsp_\La^{\aft}$. The difference is that in $\indcoh^*$ theory, we have $*$-pushforward for arbitrary morphisms (even for non-representable morphisms) but $*$-pullback only for representable morphisms of finite tor amplitude, whereas in $\indcoh^!$ theory, we have $!$-pullback for arbitrary morphisms but $*$-pushforward only for those representable by ind-algebraic spaces. Of course, the $(\indcoh,*)$-pushforwards along representable morphisms coincide in both theories. 
It would be interesting to see whether there is a more general sheaf theory for ind-coherent sheaves combining these two.

\item\label{rem: indcoh vs QC-2} A closely related but different construction was originally given in \cite{Gaitsgory.Rozenblyum.DAG.vol.I}.  We highlight the differences, one small and one large.

For the small one, we use framework from \Cref{SS: extension of sheaf theory} rather than $(\infty,2)$-categorical framework as developed in \emph{loc. cit.} 
Modulo this difference of methods, the restriction of our \eqref{eq: indcoh for algsp} to $\corr(\Sch_\La^{\aft})$ and the sheaf theory constructed in  \cite[Theorem 5.2.1.4]{Gaitsgory.Rozenblyum.DAG.vol.I} are the same. (But note that our result is slightly stronger than \emph{loc. cit.}, even taking \cite[Theorem 5.3.4.3]{Gaitsgory.Rozenblyum.DAG.vol.I} into account, as we can define $(\indcoh,*)$-pushforward along non-separated morphisms.)

However, for algebraic stacks, our theory $\indcoh^!$ and the one constructed in \cite[Theorem 5.3.4.3]{Gaitsgory.Rozenblyum.DAG.vol.I} (then restricted to ind-algebraic stacks) are quite different.
We first construct $\indcoh^{!,+}$ for stacks via right Kan extension from the theory $\indcoh^{!,+}$ for algebraic spaces and then obtain a theory $\indcoh^!$, while \emph{loc. cit.} constructed a theory for stacks via right Kan extension from the theory $\indcoh^!$ for schemes. We denote this latter theory constructed in \cite{Gaitsgory.Rozenblyum.DAG.vol.I} by $\QC$, following the notation from \cite{ben2017spectral}. In general, $\indcoh^!(X)$ and $\QC(X)$ are different.

To summarize, associated to $X\in\ArStk_\La^{\aft}$ we always have the following sequence of functors
\[
\ind\Perf(X)\to \indcoh(X)\to \QC(X)\to \Qcoh(X).
\]
The latter three categories admit $t$-structures and the natural functors are $t$-exact, inducing equivalences
\[
\indcoh(X)^{\geq n}\cong \QC(X)^{\geq n}\cong \Qcoh(X)^{\geq n}
\] 
for every $n$.
The functor $\ind\Perf(X)\to \indcoh(X)$ is an equivalence when $X$ is classical and regular. The functor $\indcoh(X)\to \QC(X)$ is an equivalence when $X$ is an algebraic space, or when $\La$ is a $\bQ$-algebra and the automorphism groups of its geometric points are aﬃne, by \cite[Theorem 3.3.5]{Drinfeld.Gaitsgory.finiteness.2013}.
\end{enumerate}
\end{remark}

\subsubsection{Grothendieck Serre duality}\label{SS: GS duality}
At this point, we know that for an ind-algebraic stack $X$ almost of finite presentation over $\La$, $\indcoh(X)$ admits a symmetric monoidal structure given by 
\[
\mF\boxtimes_\La \mG\mapsto \mF\otimes^!\mG=\Delta_X^{\indcoh,!}(\mF\boxtimes_\La \mG).
\] 
with the monoidal unit given by 
\[
\cohdual_X:=(\pi_X)^{\indcoh,!}\La \in \indcoh(X)^+.
\] 
We call the above tensor product the $!$-tensor product.
By \Cref{prop: shrek pullback of finite tor amplitude morphism}, if $X$ is an algebraic stack eventually coconnective (so $X\to \spec \La$ is of finite tor amplitude), then $\cohdual_X\in \Coh(X)$. In particular, when $X=X_{\cl}$ is a classical algebraic stack, then $\cohdual_X\in\Coh(X)$ is the classical dualizing complex of $X$.

Our goal is to prove the following result.

\begin{theorem}\label{thm: Frob structure for indcoh}
Let $X\in \indarstk_\La^{\aft}$. Then $\rg^{\indcoh}(X,-): \indcoh(X)\to \Mod_\La$ is a Frobenius structure of $\indcoh(X)$. I.e. the functor
\begin{equation}\label{eq: counit coh duality}
e: \indcoh(X)\otimes_\La\indcoh(X)\cong \indcoh(X\times X)\xrightarrow{(\Delta_X)^!} \indcoh(X)\xrightarrow{(\pi_X)_*^{\indcoh}} \Mod_\La
\end{equation}
define a self-duality 
\[
\verd_X^{\indcoh}: \indcoh(X)^\vee\cong \indcoh(X)
\] 
\end{theorem}

During the course of the proof of the theorem, we will see that the restriction $\verd_X^{\indcoh}$ to the subcategory of compact objects gives the usual Grothendieck-Serre duality $\verd_X^{\Coh}: \Coh(X)^{\op}\cong \Coh(X)$ of $X$.

We first deal with a special case.

\begin{lemma}\label{lem: Frob structure for indcoh special case}
Suppose $\indcoh(X)\otimes_\La\indcoh(X)\cong \indcoh(X\times X)$ (e.g. as in \Cref{prop: tensor equiv coh} \eqref{prop: tensor equiv coh-1}). Then
$\rg^{\indcoh}(X,-)$ defines a Frobenius structure of $\indcoh(X)$. In this case, the unit of the self-duality datum are given by
\begin{equation}\label{eq: unit coh duality} 
(\Delta_X)_*^{\indcoh}(\cohdual_X)\in \indcoh(X\times X)\cong \indcoh(X)\otimes_\La\indcoh(X).
\end{equation}
\end{lemma}
\begin{proof}
If $X$ is an algebraic space, we can apply the general consideration  \Cref{rem-dualizability in corr} to the sheaf theory $\indcoh^!$ to conclude. If $X$ is an algebraic stack, currently we cannot put $(\pi_X)_*^{\indcoh}$ and $(\indcoh, !)$-pullbacks into one sheaf theory so \Cref{rem-dualizability in corr} does not apply directly. Nevertheless applying the $\indcoh^!$ theory we see that $(\id\times\Delta_X)^{\indcoh,!}((\Delta_X)^{\indcoh}_*\cohdual_X\boxtimes \mF)\cong (\Delta_X)_*^{\indcoh}(\mF)$. On the other hand, 
\[
(\pi_X)_*^{\indcoh}\otimes \id=(\pr_1)_*^{\indcoh}: \indcoh(X)\otimes\indcoh(X)\cong \indcoh(X\times X)\to \indcoh(X)
\] 
is still defined and we have $(\pr_1)_*^{\indcoh}\circ (\Delta_X)_*^{\indcoh}\cong \id$. This proves the lemma.
\end{proof}

Now we drop the assumption $\indcoh(X)\otimes_\La\indcoh(X)\cong \indcoh(X\times X)$. Our strategy is to first construct the Grothendieck-Serre duality and then use it to prove \Cref{thm: Frob structure for indcoh}.
We start with assuming that $X\in \algsp_\La^{\aft}$.
For $\mF\in\Coh(X)$, we let $\underline{\Hom}(\mF,\cohdual_X)\in\Ind\Perf(X)=\Qcoh(X)$ be as in \eqref{eq: internal hom indcoh as quasi}. 

\begin{lemma}\label{lem: Grothendieck duality commute with proper push}
If $f: Y\to X$ is a proper morphism in $\algsp_\La^{\aft}$, then 
\begin{equation*}
\underline{\Hom}(f^{\indcoh}_*\mF,\cohdual_X)=f_*\underline{\Hom}(\mF,\cohdual_Y).
\end{equation*}
If $f: X\to Y$ is a smooth morphism in $\algsp_\La^{\aft}$, then we have
\[
f^*\underline{\Hom}(\mF,\cohdual_Y)\cong \underline{\Hom}(f^{\indcoh,!}\mF,\cohdual_X).
\]
\end{lemma}

\begin{proof}
For the first isomorphism, let $\mE\in \Qcoh(Y)$, then we have 
\begin{multline*}
\Hom(\mE, \underline{\Hom}(f^{\indcoh}_*\mF,\cohdual_X))\cong \Hom(\mE\otimes f^{\indcoh}_*\mF, \cohdual_X)\\
\cong \Hom(f^{\indcoh}_*(f^*\mE\otimes\mF), \cohdual_X)\cong \Hom(f^*\mE\otimes\mF, \cohdual_Y)\cong \Hom(\mE, f_*\underline{\Hom}(\mF,\cohdual_Y)).
\end{multline*}
For the second isomorphism, we first assume that $f$ is \'etale and $X$ is affine.
We apply \eqref{eq: qcoh coh tensor internel hom} to $\mE=f_*\mO_X$  and use \eqref{eq: coh perf adjunction} to see
\begin{multline*}
f_*f^* \underline{\Hom}(\mF,\cohdual_Y)\cong f_*\mO_X\otimes \underline{\Hom}(\mF,\cohdual_Y)\cong \underline{\Hom}(\mF, f_*\mO_X\otimes \omega_Y) \\
\cong   \underline{\Hom}(\mF, f^{\indcoh}_*\cohdual_X)\cong f_*\underline{\Hom}(f^{\indcoh,*}\mF,\cohdual_X).
\end{multline*}
Next we assume that $f$ is still \'etale but $X$ is separated. Then choosing an \'etale cover of $X$ by an affine scheme, we reduce to the previous case. 

Finally, if $f$ is smooth, again by choose an \'etale cover of $X$ by an affine scheme, we may assume that $X$ is affine. Then we can conclude by the same calculation as above, plus the fact that $f^{\indcoh,*}$ and $f^{\indcoh,!}$ differ by tensoring a line bundle.
\end{proof}

\begin{lemma}\label{lem: Coh GS dual}
If $\mF$ is coherent, then  $\underline{\Hom}(\mF,\cohdual_X)$ is coherent.
\end{lemma}
\begin{proof}
If $X=X_{\cl}$, this is classical. In general, we just need to prove the statement for coherent sheaves of the form $\iota_*\mF$, where $\iota: X_{\cl}\to X$ is a closed embedding. But then this follows from
\Cref{lem: Grothendieck duality commute with proper push}.
\end{proof}

\Cref{lem: Coh GS dual} allows one to define a functor
\begin{equation}\label{eq: GS duality}
\verd_X^{\Coh}=\underline\Hom(-,\cohdual_X): \Coh(X)^{\op}\cong\Coh(X).
\end{equation}
When $X=X_{\cl}$ is classical, this is the classical Grothendieck duality functor, which induces an anti-involution of $\Coh(X)$. This continues to hold in the derived setting.

\begin{proposition}\label{thm: GS duality}
Let $X\in\algsp_\La^{\aft}$. then \eqref{eq: GS duality} is an anti-involution.
\end{proposition}
\begin{proof}
We need to show that the canonical morphism $\mF\to \verd_X^{\Coh}(\verd_X^{\Coh}(\mF))$ is an isomorphism. It is enough to check this for a set of generators. Therefore, we just need to check it for $\iota_*\mF$, where $\iota: X_{\cl}\to X$ is as in the proof of \Cref{lem: Coh GS dual}. Again, using \Cref{lem: Grothendieck duality commute with proper push}, we reduce to the classical Grothendieck duality.
\end{proof}

Via descent, we obtain a duality functor
\begin{equation}\label{eq: GS duality, stack}
\verd_X^{\Coh}: \Coh(X)^{\op}\cong \Coh(X).
\end{equation}
when $X\in \indarstk_\La^{\aft}$ such that \Cref{lem: Grothendieck duality commute with proper push} continue to hold.

\begin{proof}[Proof of \Cref{thm: Frob structure for indcoh}]
We will need to show that for $\mF,\mG\in \Coh(X)$, we have
\[
\Hom(\mF,\mG)\cong \rg^{\indcoh}(X,\verd_X^{\Coh}(\mF)\os \mG).
\]
We can reduce to prove this for $X$ being eventually coconnective using \Cref{lem: convergence of indcoh}.

First, we assume that $X$ is a separated algebraic space. Then 
we have
\[
\rg^{\indcoh}(X,\verd_X^{\Coh}(\mF)\os \mG)=\Hom((\Delta_X)_*\mO_X, \verd_X^{\Coh}(\mF)\boxtimes \mG)=\Hom(\mF\boxtimes \verd_X^{\coh}(\mG), (\Delta_X)_*\cohdual_X).
\]
As explained in \cite[Proposition 9.5.7]{gaitsgory2013.indcoh}, since all the objects are in the bounded from below subcategories, the right hand side can be computed in $\Qcoh(X)$ as 
\[
\Hom_{\Qcoh(X)}(\mF\otimes \verd_X^{\coh}(\mG), \cohdual_X)=\Hom_{\Qcoh(X)}(\mF,\underline\Hom(\verd_X^{\coh}(\mG),\cohdual_X))=\Hom_{\Coh(X)}(\mF,\mG),
\]
as desired. For general coconnective algebraic space $X$, we may choose an \'etale over $f: U\to X$ with $U$ separated. Let $f_\bullet: U_\bullet\to X$ be the correspond \v{C}ech cover. Then $U_n$ is separated and coconnective for each $n$. Let $\mF,\mG\in \Coh(X)$. Then by descent, we have
\begin{eqnarray*}
\Hom(\mF,\mG) & = & \lim_n\Hom((f_n)^*\mF,(f_n)^*\mG)=\lim_n \rg^{\indcoh}(U_n,\verd_{U_n}^{\Coh}((f_n)^*\mF)\os (f_n)^*\mG)\\
                         & = & \lim_n\Hom(\mO_{U_n}, (f_n)^*(\verd_X^{\Coh}(\mF)\os \mG)) =\Hom(\mO_X,\verd_X^{\Coh}(\mF)\os \mG).
\end{eqnarray*}
Next, if $X$ is a coconnective algebraic stack, one can choose a smooth cover $U\to X$ with $U$ an algebraic space, and repeat the argument to conclude. 

Finally, the case of ind-algebraic stacks follows easily as well.
\end{proof}

Recall that in a dualizable category $\bfC$, there is the subcategory $\bfC^{\adm}$ of admissible objects. A self-duality $\verd$ of $\bfC$ induces $\verd^{\adm}: (\bfC^\adm)^{\op}\to \bfC^{\adm}$.
In particular, for $X\in \indarstk_\La^{\aft}$, we have 
\begin{equation}\label{eq: adm duality for indcoh}
(\verd_X^{\indcoh})^{\adm}: (\indcoh(X)^{\adm})^{\op}\to \indcoh(X)^{\adm}.
\end{equation}

We record the following result, which will be used in the main body of the article.
\begin{lemma}\label{lem: adm duality compatible with shrek pullback along closed embedding}
\begin{enumerate}
\item Let $f: X\to Y$ be a representable proper morphism of algebraic stacks almost of finite presentation over $\La$. Then 
\[
f^{\indcoh,!}\circ (\verd_Y^{\indcoh})^\adm\cong (\verd_X^{\indcoh})^\adm\circ f^{\indcoh, !}.
\]
\item Let $f: X\to Y$ be a representable quasi-smooth morphisms of algebraic stacks almost of finite presentation over $\La$. Then
\[
f_*^{\indcoh}((\verd_X^{\indcoh})^\adm(-)\otimes \Sym^d(\bL_{X/Y}[1])^{-1}) \cong  (\verd_Y^{\indcoh})^\adm((f_*^{\indcoh}(-))
\]
\end{enumerate}
\end{lemma}
\begin{proof}
In the first case, we apply \Cref{lem: adm duality and functors}  to $f_*^{\indcoh}: \indcoh(X)\to\indcoh(Y)$, by noticing that $(\indcoh,*)$-pushforwards along proper morphisms commute with Grothendieck-Serre duality.

In the second case, we notice that
\[
f^{*}( \verd_{Y}^{\coh}(-))\cong \verd_{X}^{\coh}(f^{\indcoh,!}(-))\cong  \verd_{X}^{\coh}(f^{*}(-) \otimes \Sym^d(\bL_{X/Y}[1])).
\]
We will again apply \Cref{lem: adm duality and functors} to $f^*$ to conclude.
\end{proof}

\subsubsection{Trace formalism}
We combine the general formalism of (categorical) trace as developed in \Cref{sec:categorical-preliminaries} and \Cref{sec:trace}  with the sheaf theory $\indcoh^*$ and $\indcoh^!$ to deduce the following statements.

\begin{proposition}\label{prop: duality of coh for quotient stack}
Let $X\in\indarstk_\La^{\aft}$, equipped with an automorphism $\phi: X\to X$. Suppose that
$\indcoh(X)\otimes_\La\indcoh(X)\cong \indcoh(X\times X)$ (e.g. $X$ is as in \Cref{prop: tensor equiv coh}). Then
\begin{equation}\label{ordinary Hochschild homology}
\mathrm{tr}(\indcoh(X),\phi)\cong \rg^{\indcoh}(\mL_\phi(X),\omega_{\mL_\phi(X)}).
\end{equation}
\end{proposition}
\begin{proof}
This follows from the fact that \eqref{eq: unit coh duality} and \eqref{eq: counit coh duality} form a duality datum of $\indcoh(X)$ and base change isomorphisms for coherent sheaves.
\end{proof}

\begin{remark}\label{rem: horizontal trace qcoh vs indcoh}
Let $\La$ be a field of characteristic zero. Let $X$ be an algebraic stack over $\La$. We assume that $X$ is a quotient of a scheme $U$ almost of finite presentation over $\La$ by a smooth affine algebraic group $H$. In this case, we know that $\qcoh(X)=\ind\Perf(X)$ is a rigid symmetric monoidal category.
As $\cohdual_X\in\Coh(X)$, we have the functor 
\begin{equation}\label{eq: upsilon from qcoh to coh}
\Upsilon_X: \Qcoh(X)\to \indcoh(X), \quad \mF\mapsto \mF\otimes \omega_X.
\end{equation}
which is $\Qcoh(X)$-linear, with a $\Qcoh(X)$-linear right adjoint. Now suppose $\phi: X\to X$ is an automorphism, inducing autoequivalences of $\Qcoh(X)$ and $\indcoh(X)$ and $\Upsilon_X$ is clearly $\phi$-equivariant.
We may regard \eqref{eq: upsilon from qcoh to coh} as a $\Qcoh(X)$-linear morphism compatible with $\phi$-actions. By \Cref{lem-functoriality-duality-data-generalization}, we have a morphism in $\tr(\Qcoh(X),\phi)=\tr(\Qcoh(X),{}^\phi\Qcoh(X))$
\[
[\Qcoh(X),\phi]_{\phi}\to [\indcoh(X),\phi]_{\phi},
\]
which under the equivalence $\tr(\Qcoh(X),\phi)\cong \Qcoh(\mL_\phi(X))$ from \Cref{ex: cat trace of qcoh and fixed point}, is identified with a morphism 
\begin{equation}\label{eq: canonical volume form on loop space}
\upsilon_X: \mO_{\mL_\phi(X)}\to \cohdual_{\mL_\phi(X)}
\end{equation}
in $\qcoh(\mL_\phi(X))$. Here we use the fact that $\cohdual_{\mL_\phi(X)}\in \indcoh(\mL_\phi(X))^+=\Qcoh(X)^+$.

Taking global section gives
\[
\rg(\mL(X), \mO_{\mL(X)})\to \rg(\mL(X),\cohdual_{\mL(X)})=\rg^{\indcoh}(\mL(X), \cohdual_{\mL(X)}),
\]
which can be identified with
\[
\mathrm{tr}(\Qcoh(X))\to \mathrm{tr}(\indcoh(X)).
\]
from \Cref{lem-functoriality-duality-data}. 
\end{remark}

Next we consider the categorical trace for monoidal categories of (ind-)coherent sheaves arising from the convolution pattern. We follow the notations of \Cref{subsec-trace.convolution.categories}. 

\begin{proposition}\label{prop: trace of convolution category}
Let $f: X\to Y$ be a morphism in $\indarstk_\La^{\aft}$. 
\begin{enumerate}
\item\label{prop: trace of convolution category-1} Suppose $X$ is a smooth algebraic stack over $\La$. Then $\indcoh(X\times_YX)=\indcoh^*(X\times_YX)$ has a natural monoidal structure.
\item\label{prop: trace of convolution category-1.5} On the other hand, if $f:X\to Y$ belongs to $\ind\aft$, then $\indcoh(X\times_YX)=\indcoh^!(X\times_YX)$ has a natural monoidal structure.
\item\label{prop: trace of convolution category-2} Suppose $X\to X\times X$ is of finite tor amplitude and $f$ and the relative diagram $X\to X\times_YX$ are representable proper morphisms and suppose $\indcoh(X\times_YX)\otimes_\La\indcoh(X\times_YX)\to \indcoh(X\times_YX\times X\times_YX)$ is an equivalence. Then $\indcoh(X\times_YX)$ the the monoidal structure from Part \eqref{prop: trace of convolution category-1.5} is a semirigid monoidal category. In addition,
we have a natural fully faithful embedding
 \[
 \tr(\indcoh(X\times_YX),\indcoh(X\times_YZ\times_YX))\hookrightarrow \indcoh(Y\times_{Y\times Y}Z),
 \] 
 with essential image generated (as $\La$-linear categories) by the image of $q_*^{\indcoh}\circ (\delta_0)^{\indcoh,!}: \indcoh(X\times_YZ\times_YX)\to \indcoh(Y\times_{Y\times Y}Z)$.
\end{enumerate}
\end{proposition}
\begin{proof}
As $X$ is smooth, both $\pi_X: X\to\Spec \La$ and $\Delta_X: X\to X\times X$ are of finite tor amplitude. Therefore, \Cref{thm: indcoh star pullback version} together with the general convolution pattern (see \Cref{ex:Cech-nerve} and \Cref{rem:segal.objects.morphisms.vert.horiz}) implies Part \eqref{prop: trace of convolution category-1}, as desired. 
Similarly, we have Part \eqref{prop: trace of convolution category-1.5}.

For Part \eqref{prop: trace of convolution category-2},
we can apply \Cref{bimodule.geom.trace.fully.faithfulness.convolution.cat}. 
In addition, by \Cref{lem: coherent Kunneth} and by our assumption, assumptions of \Cref{prop-comparison-usual-trace-geo-trac-presemirigid} and \Cref{prop-comparison-usual-trace-geo-trac} also hold. 
Part \eqref{prop: trace of convolution category-2} follows from \Cref{cor: semi rigidity of der(X)} and \Cref{prop-comparison-usual-trace-geo-trac}.
\end{proof}

\begin{example}\label{L: Serre functor coh classifying stack}
Let $H$ be an affine smooth algebraic group over a field $\La$. Let $X=\bB H$ be the classifying stack of $H$ with Lie algebra $\frakh$. We write $\pi_X: X\to \Spec \La$ for the structural map. Then $\indcoh(X)=\ind\Perf(X)$ equipped with the $!$-tensor product a rigid symmetric monoidal category (as compact and dualizable objects coincide), with unit
\[
\cohdual_X=(\pi_X)^{\indcoh,!}=(\wedge^{\dim \frakh}\frakh^*)[-\dim \frakh].
\]
It admits two natural Frobenius structures. The first is given by $\rg^{\indcoh}(X,-)$, by \Cref{prop: tensor equiv coh} and \Cref{lem: Frob structure for indcoh special case}. The second is given by
$\Hom(\cohdual_X,-)$, as in \Cref{ex: self-duality of cpt gen rigid monoidal cat}. Associated to these two Frobenius structures, we have corresponding objects $\omega^\la$. The first is given by
\[
\omega^{\can}_X:=((\pi_X)_{*}^{\indcoh})^R(\La),
\]
and the second is given by
\[
\omega^{\mathrm{sr}}_X=\omega^{\can}_X\otimes \omega_X.
\]
Here the tensor product is the $*$-tensor product. 
We note that this is in general different from $\cohdual_X$, even when $\La$ is a field of characteristic zero. Indeed, when $\La$ is a field of characteristic zero, then $\indcoh(X)$ is a proper $\La$-linear category.
Let $\fraku$ be the Lie algebra of its unipotent radical. Then 
\[
\omega^{\can}_X\cong (\wedge^{\dim \fraku}\fraku)[\dim \fraku],\quad \omega^{\mathrm{sr}}_X=\wedge^{\dim \frakh/\fraku} (\frakh/\fraku)^*[\dim(\frakh/\fraku)].
\]

In general, the Serre functor of $\indcoh(X)$ is given by 
\[
S_{\indcoh(X)}(V)= \omega^{\mathrm{sr}}\os V=\omega^{\can}_X\otimes V.
\]

Now, let $\phi: H\to H$ be an automorphism. Then $\tr(\indcoh(X),\phi)\subset \indcoh(H/\Ad_\phi H)$ consisting of those obtained by pullback along $H/\Ad_\phi H\to \bB H$. In addition,
\[
\mathrm{tr}(\indcoh(X),\phi)=\rg(H/\Ad_\phi H,\cohdual_{H/\Ad_\phi H}).
\] 

\end{example}

\subsection{Singular support of coherent sheaves}
We also need to briefly review the theory of singular support of coherent sheaves on quasi-smooth algebraic stacks locally almost of finite presentation over $\La$. Note that the theory of singular support as in \cite{arinkin2015singular} is developed under the assumption that $\La$ is a characteristic zero ground field. We briefly explain why some parts of such theory (with modifications) carry through for general $\La$.

\subsubsection{Cotangent complex}\label{SSS: cotangent complex}
First recall that for an animated $\La$-algebra $A$, the (algebraic) cotangent complex $\bL_{A/\La}$ is a connective $A$-module such that for every $A\to B$ and a connective $B$-module $V$ 
\[
\Map_{\Mod^{\leq 0}_A}(\bL_{A/\La},V)\cong \Map_{{\calg_\La}_{/B}}(A, B\oplus V),
\]
where $B\oplus V\to B$ denotes the trivial square zero extension of $B$  by $V$ in $\calg_\La$, and ${\calg_\La}_{/B}$ denotes the category of animated $\La$-algebras with a $\La$-algebra map to $B$. See \cite[\textsection{25.3.1}, \textsection{25.3.2}]{Lurie.SAG} for a detailed account.
If $A$ is a classical smooth $\La$-algebra, then $\bL_{A/\La}\cong \pi_0(\bL_{A/\La})=\Omega_{A/\La}$ is just the K\"ahler differential of $A$. If $A\to B$ is a morphism in $\calg_\La$, there is a natural morphism $B\otimes_A\bL_{A/\La}\to \bL_{B/\La}$ in $\Mod^{\leq 0}_B$ and the relative cotangent complex $\bL_{B/A}$ (defined as above with $\La$ replaced by $A$) can be identified as its fiber. 

It follows easily from the definition that $\bL_{B/A}=0$ for an \'etale extension. Therefore, the cotangent complex for a derived algebraic space $X$ is well-defined as an object $\bL_{X/\La}\in \qcoh(X)^{\leq 0}$. More generally, if $X\to Y$ is a morphism of prestacks over $\La$ representable by algebraic spaces, then the cotangent complex $\bL_{X/Y}\in \qcoh(X)^{\leq 0}$ is defined such that for every $V\to Y$ with $V$ an algebraic space, the pullback of $\bL_{X/Y}$ to $U:=V\times_YX$ is $\bL_{U/V}$. Then
for an algebraic stack $X$ over $\La$, $\bL_{X/\La}$ can be defined so that its pullback to any smooth cover $U\to X$ is the fiber of $\bL_{U/\La}\to \bL_{U/X}$.

\begin{definition}
A morphism $f: X\to Y$ is called quasi-smooth if it is representable (by algebraic spaces) and is almost of finite presentation, and $\bL_{X/Y}\in \Qcoh(X)^{\leq 0}$ has tor amplitude $\leq 1$.
\end{definition}

Now let $X$ be a quasi-smooth algebraic stack (i.e. there is a smooth cover $U\to X$ with $U$ a quasi-smooth algebraic space over $\La$) almost of finite presentation over $\La$.
Let $\bL_{X/\La}$ be its cotangent complex. 
Let $\bT_{X/\La}$ be the $\mO_X$-linear dual of $\bL_{X/\La}$, which is usually called the tangent complex of $X$. Then $\mH^i\bT_{X/\La}=0$ for $i>1$ and
$\mH^1\bT_{X/\La}$ is a coherent module over $\mO_{X_\cl}$. When $X=\spec A$, we sometimes also write them as $\bL_A, \bT_A, H^1\bT_A$.

We need the relation between (co)tangent complexes and Hochschild (co)homology in order to define
 the singular support of coherent sheaves. Recall that for an animated ring $A$ over $\La$, its Hochschild homology is defined as $A\otimes_{A\otimes_\La A}A$ (see \Cref{ex-usual-Hochschild-homology-of-DG-cat}).  
 Although  $A\otimes_{A\otimes_\La A}A$ has rich algebraic structures, we simply regard it as an $A$-module (via left action). There is a decreasing $\bZ_{\geq 0}$-filtration on $A\otimes_{A\otimes_\La A}A$, sometimes called the HKR filtration, such that
 \begin{equation}\label{eq: HKR filtration on HH}
 \mathrm{gr}_{\mathrm{HKR}}^i (A\otimes_{A\otimes_\La A}A)\cong (\wedge^i_\La \bL_{A/\La})[i].
 \end{equation}
 
We recall one (among various) way to construct such a filtration. 
 We may regard $A\mapsto A\otimes_{A\otimes_\La A}A$ as a functor $\calg_\La\to \lmodu(\Mod_\La)$ (see \Cref{sec:relative-tensor-product} for the notation $\lmodu(\Mod_\La)$), which is isomorphic to the left Kan extension along its restriction to the subcategory of polynomial $\La$-algebras. As a functor from the category of polynomial $\La$-algebras, we may refine it as a functor to the category filtered objects in $\lmodu(\Mod_\La)$, by equipping $A\otimes_{A\otimes_\La A} A$ with a decreasing $\bZ_{\geq 0}$-filtration given by the Postnikov filtration $\mathrm{Fil}^i_{\mathrm{HKR}}(A\otimes_{A\otimes_\La A}A):=\tau^{\leq -i}(A\otimes_{A\otimes_\La A}A)$. Then via the left Kan extension, we thus can refine the Hochschild homology as a functor 
from $\calg_\La$ to the category of filtered objects in $\lmodu(\Mod_\La)$. On the other hand, $A\mapsto \wedge^i\bL_{A/\La}[i]$ can also be regarded as a functor $\calg_\La\to \lmodu(\Mod_\La)$ which is isomorphic to the left Kan extension along its restriction to the category of polynomial $\La$-algebras. When $A$ is a polynomial algebra (or more generally a smooth algebra) over $\La$, the classical Hochschild-Kostant-Rosenberg theorem identifies $\pi_i(A\otimes_{A\otimes_\La A} A)$ with $\Omega_{A/\La}^i=\wedge^i\Omega_{A/\La}$. 
Therefore, \eqref{eq: HKR filtration on HH} holds for polynomial algebras, and therefore holds in general.

Now we suppose $A$ is quasi-smooth over $\La$. It follows that its Hochschild cohomology (see \Cref{ex-usual-Hochschild-homology-of-DG-cat})
 \[
 \Hom_{A\otimes_{\La}A}(A,A)\cong \Hom_{A}(A\otimes_{A\otimes_\La A}A,A)
 \] 
 admits an increasing filtration $\mathrm{Fil}^{\mathrm{HKR}}_\bullet$ with associated graded being $(\wedge^i_\La \bL_{A/\La})^\vee[-i]$. For example, we have a fiber sequence
 \begin{equation}\label{eq:fil one of HH}
 A\to \mathrm{Fil}_1^{\mathrm{HKR}}\Hom_{A\otimes_\La A}(A,A)\to \bT_A[-1].
 \end{equation}
 We thus arrive the following statement.
 \begin{lemma}\label{lem: tangent to HH}
 Let $A$ be a quasi-smooth $\La$-algebra. Then there is a natural injective map
$H^1\bT_A\to H^2\Hom_{A\otimes_\La A}(A,A)=:\mathrm{Ext}^2_{A\otimes_\La A}(A,A)$. 
\end{lemma}
\begin{proof}
Taking $H^2$ of fiber sequence \eqref{eq:fil one of HH} gives
\[
H^1\bT_A\cong H^2 \mathrm{Fil}_1\Hom_{A\otimes_\La A}(A,A)\to H^2 \Hom_{A\otimes_\La A}(A,A),
\]
as desired. Finally, the injectivity follows from the fact that  $\Hom_{A\otimes_\La A}(A,A)/\mathrm{Fil}_1\Hom_{A\otimes_\La A}(A,A)\in \Mod_{\La}^{\geq 2}$.
\end{proof}

\begin{remark}
There is a more concrete description of this map. We suppose $X=\Spec A=\{0\}\times_VU$, where $U=\Spec R_0$ and $V=\Spec R_1$ are smooth over $\La$, and $\{0\}: \Spec \La\to V$ is a point. Then we have the correspondence 
\[
\{0\}\times_V\{0\}\leftarrow X\times_UX\to X\times_\La X.
\] 
Let $\delta_{0}$ be the $*$-pushforward of structure sheaf of $\{0\}$ along the relative diagonal $\{0\}\to \{0\}\times_V\{0\}$. By base change (of quasi-coherent sheaves), it's $*$-pull-push along the above correspondence is canonically $(\Delta_X)_*\mO_X$.
It follows that we obtain a natural morphism
\[
\Sym (\bT_0 V[-2] )\otimes_\La A\cong \End_{\Qcoh(\{0\}\times_V\{0\})}\delta_0\otimes_\La A\to \End_{\Qcoh(X\times_\La X)}((\Delta_X)_*\mO_X)= \Hom_{A\otimes_\La A}(A,A).
\]
This morphism factors through $\bT_0V[-2] \otimes_\La A \to H^1 \bT_A [-2]$, giving the desired morphism as in \Cref{lem: tangent to HH}.
\end{remark}

\subsubsection{Singular support}
\begin{definition}\label{def: Stack of singularieties} 
Let $X$ be a quasi-smooth algebraic stack over $\La$.
The stack of singularities of $X$ is a classical algebraic stack of finite presentation over $\La$ defined as
\[
\sing(X)=\Sym_{\mO_{X_\cl}} (\mH^{1} \bT_{X/\La}).
\]
\end{definition}
By definition, there is a canonical map $\sing(X)\to X_\cl$. The fiber over a (field valued) point $x\in X$ is the vector space $H^{-1}(x^*\bL_{X/\La})$. There is a natural $\bG_m$-action along $\sing(X)\to X$ by dilatation.

Recall that a morphism $f:X\to Y$ of smooth varieties induces a map of cotangent spaces $df: \bT^*Y\times_YX\to \bT^*X$. There is a similar construction for the stack of singularities. Namely,
let $f: X\to Y$ be a morphism of quasi-smooth algebraic stacks almost of finite presentation over $\La$. Then we have a map of coherent sheaves $\mH^1\bT_{X/\La}\to \mH^1 f^*\mH^1\bT_{Y/\La}$ on $X$, inducing a morphism 
\begin{equation}\label{eq: base change singular support}
\Sing(Y)_X:=\Sing(Y)\times_{Y_\cl} X_\cl \rightarrow \Sing(X).
\end{equation} 
Following the notation of \cite{arinkin2015singular}, we denote this map by $\Sing(f)$. We will also use the following notations: if $\mN\subset \Sing(Y)$ is a closed conic subset,  then we let $\Sing(f)(\mN)$ denote the smallest closed conic subset of $\Sing(X)$ containing the image of $\mN\times_{Y_\cl}X_{\cl}$ under the map $\Sing(Y)_X\to \Sing(X)$. If $\mN\subset \Sing(X)$ is a closed conic subset, then we let $\Sing(f)^{-1}(\mN)$ denote the closed conic subset of $\Sing(Y)$ consisting of the image of $\Sing(f)^{-1}(\mN)$ under the map $\Sing(Y)_X\to \Sing(Y)$.

\begin{remark}\label{rem: singular support formal stack}
We give a convenient definition of the stack of singularities for certain formal algebraic stacks.
Let $X$ be a quasi-smooth algebraic stack over $\La$ and let $Z\subset X$ be a closed substack. Let $\widehat{Z}$ be the formal completion of $X$ along $Z$ as in \Cref{ex: formal completion}. We let $\Sing(\widehat{Z})=Z_{\cl}\times_{X_{\cl}}\Sing(X)$. One can show that this only depends on $\widehat{Z}$, i.e. if $\widehat{Z}$ is realized as the formal completion of $X'$ along $Z'$, then $\Sing(\widehat{Z})=\Sing(\widehat{Z}')$. Note that if $Z$ is quasi-smooth, in general $\Sing(\widehat{Z})\neq \Sing(Z)$.
\end{remark}

The goal is to construct, for every coherent complex $\mF\in \Coh(X)$, a conic closed subset $s.s.(\mF)\subset \sing(X)$, called the singular support of $\mF$. Note that the construction of \cite{arinkin2015singular} makes use of some results of Hochschild homology  from Appendix G of \emph{loc. cit.} which are specific to $\bQ$-algebras and therefore not applicable to a general base ring $\La$. Fortunately, to define the singular support of a coherent sheaf, all we need is \Cref{lem: tangent to HH}.

First we assume that $X=\spec A$ is affine. 
As the Hochschild cohomology of $A$ is just the center of the category $\Mod_A$ (see \Cref{ex-usual-Hochschild-homology-of-DG-cat}), we obtain an action of $\Hom_{A\otimes_\La A}(A,A)$ on any $A$-module (see \Cref{rem: action of center of objects}). In particular, if $M$ is a coherent $A$-module, we have a homomorphism of graded (ordinary) commutative algebras
\[
\mathrm{Ext}^{2\bullet}_{A\otimes_\La A}(A,A)\to \mathrm{Ext}_{\Mod_A}^{2\bullet}(M,M),
\]
where 
\[
\mathrm{Ext}_{\Mod_A}^{2\bullet}(M,M)=\oplus_{\bullet} H^{2\bullet} \End_{\Mod_A}(M)
\] 
is a graded algebra under the usual Yoneda product. 
Together with \Cref{lem: tangent to HH}, we obtain a graded commutative algebra map
\begin{equation}\label{eq: HH action on cohomology sheaf}
\Sym^\bullet_{\pi_0(A)} H^1\bT_A  \to \mathrm{Ext}^{2\bullet}(M,M).
\end{equation}

The following technical result is important for our understanding of singular support of coherent sheaves.
Suppose $i:Y=\Spec B\to X=\Spec A$ be a closed embedding, defined by one equation $\La[x]\to A, x\mapsto f$. Suppose $X$ is quasi-smooth over $\La$ (so $Y$ is also quasi-smooth over $\La$).
In this case 
\[
\bL_{B/A}\cong B\otimes_\La \bL_{\La/\La[x]}\cong (B\cdot \eta_f)[1].
\]
Here we use the fact that $\bL_{\La/\La[x]}\cong \La\cdot \eta_x$, where $\eta_x$ is a canonical generator in degree $-1$, and we let $\eta_f=1\otimes \eta_x$, which is a generator of $\bL_{B/A}$. Let $\xi_f$ be the dual basis of $\bT_{B/A}=\bL_{B/A}^\vee$. We have a right exact sequence
\[
H^1\bT_{B/A}\to H^1\bT_B\to \pi_0(B)\otimes_{\pi_0(A)}H^1\bT_A\to 0.
\]
By abuse of notations, we also use the same notation to denote the image of $\xi_f\in H^1\bT_{B/A}$ in $\pi_0(A)\otimes_{\pi_0(B)}H^1\bT_B$.
It follows that 
\[
\pi_0(B)\otimes_{\pi_0(A)}\Sym_{\pi_0(A)} H^1\bT_A\cong \Sym_{\pi_0(B)} H^1\bT_B/(\Sym_{\pi_0(B)} H^1\bT_B\cdot \xi).
\]

Now let $\mF\in \Coh(Y)$. Then $\xi_f\in H^1\bT_B\to \Ext^2(\mF,\mF)$ induces a morphism $\mF\to \mF[2]$, still denoted by $\xi_f$. By tracking of definitions, we have the following statement.
\begin{lemma}\label{lem: cup product with xi}
We have a fiber sequence of coherent $B$-modules.
\[
i^* i_* \mF\to \mF\xrightarrow{\xi_f} \mF[2].
\]
\end{lemma}

The following results were proved in  \cite[Theorem 4.1.8]{arinkin2015singular} under the assumption that $\La$ is a field of characteristic zero. But the proofs go through in the more general base ring we consider. For completeness, we sketch the proof.

\begin{lemma}\label{prop: first property of singular support}
Assume that $X=\Spec A$ is quasi-smooth over $\La$. If $\mF\in \Coh(X)$, then the $\mathrm{Ext}^{2\bullet}(\mF,\mF)$ is a finitely generated graded $\Sym^\bullet_{\pi_0(A)} H^1\bT_A$-module.
\end{lemma}
\begin{proof}
The question is Zariski local on $X$ so we may assume that $X=\{0\}\times_VU$ where $U,V$ are smooth over $\La$. We may choose a regular sequence $(f_1,\ldots, f_m)$ in $\mO_V$ at $0$ and let $V_r=(f_1,\ldots,f_r)$ and $\Spec A_r=X_r:=V_r\times_VU$. It is enough to prove by induction that for $\mF\in \Coh(X_r)$, $\mathrm{Ext}^{2\bullet}_{\Coh(X_r)}(\mF,\mF)$ is a finitely generated graded $\Sym^\bullet_{\pi_0(A_r)} H^1\bT_{A_r}$-module. The case $r=0$ is clear. Suppose this is the case for $r-1$. We let $\imath:X_{r}\to X_{r-1}$ be the closed embedding, defined by the equation $f_r=0$. Now let $\mF\in \Coh(X_r)$. By \Cref{lem: cup product with xi}, there is a cofiber sequence
\[
\mathrm{Ext}^{2\bullet-2}_{\Coh(X_{r})}(\mF,\mF)\xrightarrow{\xi_{f_r}} \mathrm{Ext}^{2\bullet}_{\Coh(X_{r})}(\mF,\mF) \to \mathrm{Ext}^{2\bullet}_{\Coh(X_{r-1})}(\imath_*\mF,\imath_*\mF).
\]
By induction, $\mathrm{Ext}^{2\bullet}_{\Coh(X_{r-1})}(\imath_*\mF,\imath_*\mF)$ is finitely generated over  $\Sym^\bullet_{\pi_0(A_{r})} H^1\bT_{A_{r-1}}$. 
Since the grading of $\mathrm{Ext}^{2\bullet}_{\Coh(X_{r})}(\mF,\mF)$ is bounded from below, a standard argument shows that $ \mathrm{Ext}^{2\bullet}_{\Coh(X_{r})}(\mF,\mF)$ is finitely generated over $\Sym^\bullet_{\pi_0(A_{r})}H^1\bT_{A_{r}}$.
\end{proof}

We thus can regard $\mathrm{Ext}^{2\bullet}(\mF,\mF)$ as a $\bG_m$-equivariant (ordinary) coherent sheaf on $\Sing(X)$. Let $s.s.(\mF)\subset \Sing(X)$ be its support. This is the promised singular support of $\mF$.

Now it is clear that the map \eqref{eq: HH action on cohomology sheaf} $\Sym H^1 \bT_A\to \mathrm{Ext}^{2\bullet}(\mF, \mF)$ is compatible with  smooth morphism. Indeed, if $f: \Spec B=Y\to \Spec A=X$ is smooth, then $\Sing(f): \Sing(X)_Y\to \Sing(Y)$ is an isomorphism and for $\mF\in \Coh(X)$, we have $s.s.(\mF)_Y=s.s.(f^*\mF)$. 
Therefore, $s.s.(\mF)\subset \Sing(X)$ is well-defined for $X$ being a quasi-smooth algebraic stacks over $\La$ and $\mF\in\Coh(X)$.

Now, let $X$ be a quasi-smooth algebraic stack almost of finite presentation over $\La$. Let $\mN\subset \Sing(X)$ be a conic closed subset. We let 
\[
\Coh_\mN(X)\subset \Coh(X)
\]
denote the full subcategory consisting of those $\mF$ such that $s.s.(\mF)\subset \mN$, and let $\indcoh_\mN(X)$ be the ind-completion of $\Coh_\mN(X)$. 

\begin{remark}
Note that $\Coh_\mN(X)$ is clearly idempotent complete so $\indcoh_\mN(X)^\cpt=\Coh_\mN(X)$. In addition, the inclusion $\indcoh_\mN(X)\to \indcoh(X)$ admits a continuous right adjoint. 
We also note that as mentioned in \Cref{rem: first remark of ind-coh}, for general stacks $X$ our definition of $\indcoh_\mN(X)$ is different from the definition given in \cite{arinkin2015singular}. 
\end{remark}

We have the following basic functoriality for singular support of coherent sheaves (compare with \cite[Proposition 4.7.2, Proposition 7.1.3]{arinkin2015singular}).

\begin{proposition}\label{prop: GS duality and Singular support}
Let $\mN\subset \Sing(X)$ be a closed conic subset. Then the Grothendieck-Serre duality $\verd_X^{\coh}$ restricts to an equivalence $\Coh_\mN(X)^{\op}\cong \Coh_\mN(X)$. 
\end{proposition}

\begin{proposition}\label{lem: control functor image by singular support}
Let $f: X\to Y$ be a morphism of quasi-smooth algebraic stacks almost of finite presentation over $\La$.
\begin{enumerate}
\item\label{lem: control functor image by singular support-1} If $f$ is quasi-smooth, then $f^{\indcoh, *}$ sends $\Coh_{\mN}(Y)$ to $\Coh_{\Sing(f)(\mN)}(X)$. 
\item\label{lem: control functor image by singular support-2} If $f: X\to Y$ is a proper morphism. Then $f_{*}^{\indcoh}: \Coh(X)\to \Coh(Y)$ sends $\Coh_{\mN}(X)$ to $\Coh_{\Sing(f)^{-1}(\mN)}(Y)$. 
\end{enumerate}
\end{proposition}

\begin{lemma}\label{lem: pullback quasi-smooth closed embedding}
Let $i: X\to Y$ be a quasi-smooth closed embedding of quasi-smooth algebraic stacks over $\La$. 
Let $\mN_Y\subset \Sing(Y)$ be a closed conic subset and let $\mN_X= \Sing(i)(\mN_Y\times_{Y_\cl}{X_\cl})\subset \Sing(X)$.
Then $i^{\indcoh,*}(\indcoh_{\mN_Y}(Y))$ is contained in $\indcoh_{\mN_X}(X)$ and generates the latter as $\La$-linear presentable stable category. 
\end{lemma}
We note that $\Sing(i): \Sing(Y)_X\to \Sing(X)$ is a closed embedding so $\mN_X$ is automatically a conic closed subset of $\Sing(X)$.
\begin{proof}
As $i_*^{\indcoh}$ sends $\ind(\Coh_{\mN_X}(X))$ to $\ind(\Coh_{\mN_Y}(Y))$,
to show that $i^{\indcoh,*}(\ind(\Coh_{\mN_Y}(Y)))$ generate $\ind(\Coh_{\mN_X}(X))$, it is enough to show that the composed functor $\ind(\Coh_{\mN_X}(X)) \subset \indcoh(X)\xrightarrow{i_*^{\indcoh}}\indcoh(Y)$ is conservative. 

We will show that if $0\neq \mF\in \ind(\Coh_{\mN_X}(X))$, then $i^{\indcoh,!}i_*^{\indcoh}(\mF)\neq 0$.
Namely, by definition, there is some $\mG\in \Coh_{\mN_X}(X)$ such that $\Hom(\mG, \mF)\neq 0$. Then we have
\[
\Hom(\mG, i^{\indcoh,!}i_*^{\indcoh}(\mF))=\Hom(i^*i_*\mG, \mF).
\]
As $s.s.(\mG)\in \mN_X$, we see that $\mG$ is a direct summand of $i^*i_*\mG$ by the following lemma. It follows that  $\Hom(\mG, i^{\indcoh,!}i_*^{\indcoh}(\mF))\neq 0$.
\end{proof}

\begin{lemma}
Assumptions are as in \Cref{lem: pullback quasi-smooth closed embedding}. Let $\mG\in \Coh_{\Sing(i)(\Sing(Y)_X)}(X)$. Consider the cofiber sequence $i^*i_*\mG\to \mG\to \mG'$. Then the map $\mG\to \mG'$ is zero.
\end{lemma}
\begin{proof}
By descent, we may assume that $X\to Y$ is a quasi-smooth closed embedding of quasi-smooth affine schemes. Then working locally on $X$ and $Y$, we may assume that $X$ is defined by one equation $f=0$. Then as $s.s.(\mG)\subset \mN$, we see that the map $\xi_f: \mG\to \mG[2]$ as in \Cref{lem: cup product with xi} is zero, as desired.
\end{proof}

The following result is analogous to \cite[Theorem 4.2.6]{arinkin2015singular}.
\begin{corollary}
Let $\mF\in \Coh(X)$, then $\mF\in \Perf(X)$ if and only if $s.s(\mF)=X_{\cl}\stackrel{0}{\hookrightarrow} \Sing(X)$.
\end{corollary}
\begin{proof}
By descent, we may assume that $X$ is an affine scheme, which can be embedded into a smooth affine scheme $i:X\to Y$. Then by
\Cref{lem: pullback quasi-smooth closed embedding}, $\Coh_{\{0\}}(X)$ is generated by $i^*\Perf(Y)$. The corollary follows.
\end{proof}

\begin{corollary}\label{cor: generation of coh under shrek pullback along closed}
Let $i: X\to Y$ be as in \Cref{lem: pullback quasi-smooth closed embedding}. Let $X^\wedge$ be the formal completion of $Y$ along $X$. Let $\mN_Y\subset \Sing(Y)$ be a closed conic subset. Then 
$i^{\indcoh}_*(i^{\indcoh,!}(\ind(\Coh_{\mN_Y}(Y))))$ generates $\ind(\Coh_{\mN_X}(Y))$ as presentable $\La$-linear category, where $\mN_X=\mN_Y\times_{Y_\cl}X_{\cl}$ is regarded as a conic closed subset of $\Sing(Y)$.
\end{corollary}

The following statements are from \cite[Proposition 7.6.4, Theorem 7.8.2]{arinkin2015singular}. Note that although \emph{loc. cit.} assumes the ground field is of characteristic zero, the proofs actually work for general base ring we consider.
\begin{proposition}\label{prop: control functor image by singular support}
Let $f: X\to Y$ be a morphism of quasi-smooth algebraic spaces almost of finite presentation over $\La$.
\begin{enumerate}
\item\label{prop: control functor image by singular support-1} If $f$ is quasi-smooth, then $f^{\indcoh, *}$ sends $\Coh_{\mN}(Y)$ to $\Coh_{\Sing(f)(\mN)}(X)$ and the essential image generates $\Coh_{\Sing(f)(\mN)}(X)$ as idempotent complete stable categories. 
\item\label{prop: control functor image by singular support-2} If $f: X\to Y$ is a proper morphism. Then $f_{*}^{\indcoh}: \Coh(X)\to \Coh(Y)$ sends $\Coh_{\mN}(X)$ to  $\Coh_{\Sing(f)^{-1}(\mN)}(Y)$ and  the essential image generates $\Coh_{\Sing(f)^{-1}(\mN)}(Y)$ as idempotent complete stable categories.
\end{enumerate}
\end{proposition}

For our applications, it is important to have these results extended to (certain) algebraic stacks. Now situation crucially depends on the base ring $\La$. First, 
if $\La$ is a field of characteristic zero, these statements generalize nicely to a large class of algebraic stacks as shown in \cite[Proposition 8.4.12]{arinkin2015singular} combining with \cite[Corollary 9.2.7, 9.2.8]{arinkin2015singular}.
\begin{proposition}\label{prop: control functor image by singular support stacks}
Let $f: X\to Y$ be a representable morphism of quasi-smooth algebraic stacks of finite presentation over a field $\La$ of characteristic zero. Suppose that $X$ and $Y$ are ``global completion intersection" in the sense of \cite[\textsection{9.2}]{arinkin2015singular}. Then statements of \Cref{prop: control functor image by singular support} hold.
\end{proposition}

\begin{remark}\label{rem: failure of generation by smooth pullback}
Unfortunately, both parts of \Cref{prop: control functor image by singular support} fail for representable morphisms between algebraic stacks in positive characteristic. 
Namely, as mentioned in \Cref{rem: failure coh surjective proper morphism}, \Cref{prop: control functor image by singular support} Part \eqref{prop: control functor image by singular support-2} fails in positive characteristic in general. On the other hand, we consider the affine smooth morphism of smooth algebraic stacks $\mathrm{PGL}_2/\mathrm{PGL}_2\to \bB \mathrm{PGL}_2$, where $\mathrm{PGL}_2$ acts on itself by conjugation action. On shows that if $\La$ is a field of characteristic two, the $*$-pullbacks of $\Perf(\bB \mathrm{PGL}_2)$ does not generate to $\Perf(\mathrm{PGL}_2/\mathrm{PGL}_2)$. 
\end{remark}

\newpage

\section{Theory of $\ell$-adic sheaves}
\label{sec:pspl-stacks}

In this section, we discuss the theory of $\ell$-adic sheaves in algebraic geometry. This theory provides the necessary foundations to work with some exotic algebro-geometric objects, such as the stack of $G$-isocrystals, which is the main focus of this article.

We will begin by reviewing and further developing the theory of $\ell$-adic (ind-constructible) sheaves on quasi-compact and quasi-separated (qcqs) algebraic spaces, and subsequently on general prestacks. This theory was first introduced in \cite{Gaitsgory.Lurie.Weil.I} and further explored in \cite{bouthier2020perverse}, among other works. However, these existing studies are inadequate for our purposes for several reasons. Firstly, we need descent results that are stronger than those proved in \emph{loc. cit.} in order to study the representation theory of (locally) profinite groups\footnote{We caution that pro-\'etale descent fails for ind-constructible sheaves in general. See \Cref{ex:sheaves-example-descent}.} and sheaves of the stack of $G$-isocrystals.
Secondly, neither \cite{Gaitsgory.Lurie.Weil.I} nor \cite{bouthier2020perverse} fully developed the six-functor formalism for such sheaf theory, which is necessary for our work. Finally, we aim to adapt this formalism to the setting of perfect algebraic geometry to study local Langlands correspondence for mixed characteristic local fields. In this context, the usual notions and arguments involving smoothness do not apply. Instead, we employ an appropriate notion of cohomological smoothness, which introduces certain subtleties. For instance, there is no canonically defined trace map in the perfect setting (over a field of characteristic $p>0$), which may present challenges. (See \Cref{lem:existence of generalized constant sheaf} for an example of such subtleties.)

The first goal of this section is to assemble various ingredients from the literature to establish a six-functor formalism for ind-constructible sheaves on prestacks, utilizing the full strength of the extensions of sheaf theories as developed in \Cref{SS: extension of sheaf theory}. We again emphasize the need for a sufficiently general theory capable of addressing profinite groups. Traditionally, the six-functor formalism for $\ell$-adic sheaves only permits pushforward along (ind-)finitely presented morphisms. However, in our formalism, we extend this to allow pushforward along (representable) pro-'etale morphisms. As will be explained in \Cref{def-dual-sheaves-on-prestacks-vert-indfp} and \Cref{rem: pushforward along non-representable morphisms}, we allow pushforward along a much larger class of morphisms, which can sometimes include non-representable cases. The main results regarding this sheaf theory are summarized in \Cref{thm: final theorem of shv}.

In the second part of the section, we specialize the theory to a class of infinite-dimensional stacks known as placid stacks, a concept first introduced in \cite{Raskin.dmodules.infinite.dimensional} and \cite{bouthier2020perverse}. Informally, placid stacks are quotients of algebraic spaces (with finite-type singularities) by pro-smooth relations. In our context, we replace pro-smooth morphisms with cohomologically pro-smooth morphisms. On placid stacks, there are well-defined notions of constructible sheaves, Verdier duality, perverse sheaves, etc., which we will review and further study.
After establishing the foundational theory for constructible sheaves on placid stacks, we will extend this theory to sind-placid stacks, which can be informally described as quotients of (ind-)placid stacks by ind-proper equivalence relations. Examples include classifying stacks of locally profinite groups and the stack of $G$-isocrystals. While this class of prestacks may appear exotic from the classical perspective, the category of $\ell$-adic sheaves on them remains reasonably well-behaved. The second major result of this section is a six-functor formalism for the category of ind-finitely generated ($\ell$-adic) sheaves on sind-placid stacks, as detailed in \Cref{thm-rshv-for-ind-quasi-placid-stack} and \Cref{prop-rshv-for-sifted-placid-stack}.

Finally, we emphasize that although we operate in the perfect setting, all results in this section also apply to the standard algebro-geometric context. To the best of our knowledge, this section presents the first systematic treatment of the theory of $\ell$-adic (co)sheaves that is suitable for applications in geometric representation theory.

\subsection{Perfect algebraic geometry}
We will use the theory of $\ell$-adic sheaves in the setting of perfect algebraic geometry.
For basic definitions and facts regarding perfect schemes and algebraic spaces over $\bF_p$, we refer to \cite[App. A]{zhu2017affine} \cite[App. A]{xiao2017cycles} and \cite[\textsection{3}]{bhatt2017projectivity}. As little extra work is needed, we will work in a slightly more general setting, i.e. we do not require schemes and algebraic spaces are over $\bF_p$. The basic fact we need remains the same as in \emph{loc. cit.}, namely universal homeomorphisms preserve the \'etale topos of schemes (and algebraic spaces).

\subsubsection{Perfect stacks}
We call a commutative ring $R$ \emph{perfect} if the following equivalent conditions are satisfied.
\begin{itemize}
\item The ring $R$ is reduced and every homomorphism $R\to R'$ with $\Spec R'\to \Spec R$ being a universal homeomorphism is an isomorphism. 
\item For all $x,y\in R$ with $x^3=y^2$ there is a unique $r\in R$ with $x=r^2$ and $y=r^3$ (a ring satsifying this condition is called seminormal) and for any prime number $p$ and $x,y\in R$ with $p^px=y^p$ there is a unique $r\in R$ with $x=r^p$ and $y=pr$.
\end{itemize}
Note that such $R$ is called absolutely weak normal in \cite[Appendix B]{rydh2010submersions} (see also \cite[\href{https://stacks.math.columbia.edu/tag/0EUK}{Section 0EUK}]{stacks-project}).
Our choice of terminology is justified as follows: if $pR=0$, then $R$ is perfect in the above sense if and only if it is perfect in the usual sense, i.e. the Frobenius endomorphism $\sigma\colon R\to R, \ r\mapsto r^p$ is an isomorphism. On the other hand, if $R=k$ is a field, then it is perfect in the above sense if and only if it is a perfect field in the usual sense. Another class of perfect commutative rings are Dedekind domains with characteristic zero fractional field.

Now we fix a perfect base commutative ring $k$ and denote by $\calg_k^{\heartsuit}$ the ordinary category of commutative $k$-algebras and $\perf_k$ its full subcategory of perfect $k$-algebras. 
The inclusion $\perf_k\subseteq \calg_k^{\heartsuit}$ admits a left adjoint, called the perfection 
\[
R\mapsto R_{\pf} = \colim_{R\to R'} R',
\]
where the colimit is taken over the (filtered) cateogry of all finitely presented homomorphisms $R\to R'$ with $\Spec R'\to \Spec R$ being universal homeomorphism (see \cite[\href{https://stacks.math.columbia.edu/tag/0EUR}{Lemma 0EUR}]{stacks-project}).
If $pR=0$, one can replace the above filtered colimit by the direct limit of $R$ with transition map being the Frobenius endomorphism. For a $k$-algebra homomorphism $f:R\to R'$, let $f_{\pf}: R_{\pf}\to R'_{\pf}$ denote its perfection. Note that $f\mapsto f_{\pf}$ preserves all topological notions. In addition, if $f$ is \'etale so is $f_{\pf}$ (see \cite[Proposition (B.6)]{rydh2010submersions}). Therefore, it makes sense to talk about Zariski and \'etale topology on $\perf_k$.

We will follow the functor of points approach and identify all our geometric objects with their associated functors from $k$-algebras to the category of sets, the $(2,1)$-category of groupoids, or in general $(\infty,1)$-category $\spc$ of spaces (also called as anima nowadays).\footnote{Although all the perfect prestacks discussed in this section will take values in the $(2,1)$-category of groupoids, it is convenient to allow them to take values in $\spc$ in order to apply higher category formalism.} In the context of perfect algebraic geometry, our test objects will be perfect $k$-algebras instead of all $k$-algebras.

\begin{definition}
A \textit{perfect prestack} $X$ is an accessible functor\footnote{See \Cref{rem: why accessibility} for an explanation.}
\[
X\colon  \perf_{k} \rightarrow \spc.
\]
We write $\prestk^{\pf}_k$ for the category of perfect prestacks (which is a full subcategory of $\fun(\perf_k,\spc)$. We call objects in $(\perf_k)^{\op}\subset\prestk^{\pf}_k$ affine
perfect schemes and write them as $\spec R$ (for $R \in \perf_k$) as usual.
A perfect prestack is called a \textit{perfect stack} if it is a sheaf with respect to the \'{etale} topology on $\perf_k$.
\end{definition}

Restriction along the inclusion $\perf_k\subset \aff_k$ gives the \textit{perfection} functor
\[
    \prestk_k  \rightarrow \prestk_k^\pf,\quad   X  \mapsto X_{\pf}, 
\]
where $\prestk_k$ denotes the category of prestacks over $k$ (as defined in \Cref{SS: derived algebraic geometry}). Note that if $X=\spec R)$ for $k$-algebra $R$, then $X_{\pf}=\spec R_\pf$. In particular, affine morphisms in $\prestk_k^\pf$ make sense and perfection of an affine morphism in $\prestk_k$ is an affine morphism in $\prestk_k^\pf$.

We can associate to every perfect prestack a topological space as in \eqref{eq-top-space-prestack}. Clearly, for a prestack $X$ over $k$ with $X_\pf$ its perfection, we have $|X|=|X_{\pf}|$.

In the rest of this section, we will usually abuse terminology and refer to perfect (pre)stacks simply as (pre)stacks.

\subsubsection{Perfect schemes and algebraic spaces}
We can define perfect schemes (resp. algebraic spaces) as Zariski (resp. \'etale) sheaves $\perf_k\to \spc$ that admit a Zariski cover (resp. \'etale cover) by affine perfect schemes satisfying additional properties as usual (e.g. see \cite[Appendix A]{xiao2017cycles}). As mentioned above, the usual topological notions, such as quasi-compact and quasi-separated (qcqs) make sense in this setting.
Note that if $X$ is a (qcqs) scheme (resp. an algebraic space) over $k$ (in the usual sense), regarded as a functor $\calg^{\heartsuit}_k\to \spc$, then $X_\pf$ as a functor $\perf_k\to\spc$ is a (qcqs) perfect scheme (resp. perfect algebraic space). 
We denote by $\psch_k$ (resp. by $\qcsp_k$) the category of perfect \emph{qcqs} schemes (resp. perfect \emph{qcqs} algebraic spaces) over $k$.

\begin{remark}
We suggest readers to skip this remark.
As in \cite[Appendix B]{rydh2010submersions} (see also \cite[\href{https://stacks.math.columbia.edu/tag/0EUK}{Section 0EUK}]{stacks-project}), there is a notion of absolutely weakly normal schemes and algebraic spaces. They are reduced schemes and algebraic spaces $X$ in the usual sense (in particular are functors $\aff_k\to\spc$) such that every separated universal homeomorphism of schemes (resp. algebraic spaces) $X'\to X$ is an isomorphism. 
If $pk=0$, then they are just schemes (resp. algebraic spaces) whose Frobenius endomorphism is an automorphism, i.e. the category of perfect schemes (resp. algebraic spaces) in usual sense.
On can show (as in \cite[Lemma A.12]{zhu2017affine}) that the restriction of the perfection functor $(-)_\pf: \prestk_k\to\prestk_k^\pf$ to the category of (qcqs) absolutely weakly normal schemes (resp. algebraic spaces) induces an equivalence from it to the category of (qcqs) perfect schemes (resp. algebraic spaces) as defined above.
\end{remark}

\begin{remark}\label{rem-closed-subscheme}
In the classical algebraic geometry (even in the derived algebraic geometry as reviewed in \Cref{SS: derived algebraic geometry}), there is a bijection between open subschemes (open subspaces) of a scheme (algebraic space) $X$ and open subsets of its topological space $|X|$. Indeed, an open subscheme/space $U\subset X$ determines and open subset $|U|\subset |X|$, which in turn determines $U$ as the prestack that represents the functor sending $R$ to those $x\in X(R)$ such that $|\spec R|\to |X|$ factors through $|U|$ (see \eqref{eq: points for open embedding}).
On the other hand, the relation between closed subsets in $|X|$ and closed embeddings $Z\subset X$ is more complicated.

In perfect algebraic geometry, while open embeddings still behave as usual,
the situation for closed embeddings is better.
Suppose $X$ is a qcqs perfect scheme/algebraic space. We say $i:Z\to X$ is a closed embedding if $i$ arises as the perfection of a closed embedding $i':Z'\to X'$ of qcqs schemes/algebraic spaces. Note that if $X'=X=\Spec A$, then $Z'$ is given by $\Spec B$ for some quotient $A\to B$.  Taking the perfection gives $Z=\Spec B_\pf$, which is in fact only depends on the underlying topological space $|Z|\subset |X|$ ( but is independent of the choice of $Z'$).
(Note, however,  that $A\to B_\pf$ may \emph{not} be surjective if $pk\neq 0$!)
In fact, let $|Z|\subset |X|$ be the corresponding closed subset. Then $Z\subset X$ represents the functor sending $R\in\perf_k$ to the subset of $X(R)$ consisting of those $x\in X(R)$ such that $|\Spec R|\to |X|$ factors through $|Z|$.

In particular, if $f:X\to Y$ is a morphism of qcqs algebraic spaces, then we can define its scheme theoretic image $Z$ just as the Zariski closure of $f(|X|)$ in $|Y|$ equipped with the above perfect scheme/algebraic space structure. Then $Z\subset Y$ is a closed embedding. This definition is reasonable (by \cite[\href{https://stacks.math.columbia.edu/tag/01R8}{Lemma 01R8}]{stacks-project}) and the formulation of scheme theoretic image commutes with flat base change (\cite[\href{https://stacks.math.columbia.edu/tag/081I}{Lemma 081I}]{stacks-project}). 
\end{remark}

\begin{definition}\label{def:morphisms-in-qcspk}
A morphism $f: X\to Y$ in $\qcsp_k$ is called perfectly finite type (resp. perfectly finitely presented, resp. perfectly proper, resp. perfectly finite) if it is the perfection of a finite type (resp. finitely presented, resp. proper, resp. finite) morphism $f':X'\to Y'$ of qcqs algebraic spaces. 
We say $f$ is \textit{perfectly smooth at }$x\in X$ if there is an \'{e}tale atlas $U \rightarrow X$ at $x$ and an \'{e}tale atlas $V\rightarrow Y$ at $f(x)$ such that $U \rightarrow Y$ factors through a map $U\rightarrow V$ which has a decomposition of the form $U \xrightarrow{h} V\times (\mathbb{A}^{n})_{\pf} \rightarrow V$ with $h$ \'{e}tale. We say $f$ is perfectly smooth if it is perfectly smooth at every point of $X$ (see \cite[Definition A.25]{zhu2017affine}). We will write pft (resp. pfp) for perfectly finite type (resp. perfectly finitely presented) morphisms for brevity.
\end{definition}

Note that the classes of pft morphisms and pfp morphisms are both strongly stable classes in $\qcsp_k$ (in the sense of \Cref{def-closure-property-of-morphism}).  This follows from the corresponding statements for finite type and finite presented morphisms between qcqs algebraic spaces (in the usual sense)
by some limit and approximation results in the perfect setting, as we now discuss.

\begin{proposition}\label{prop:appr-fp-morphism}
\begin{enumerate}
\item\label{prop:appr-fp-morphism-0} Let $f: X\to Y$ be a  morphism in $\qcsp_k$. Then $f$ is pfp if and only if for every cofiltered limit $Z=\lim_i Z_i$ with $Z_i\to Z_j$ affine, the following natural map is a bijection
\[\colim_i \Map(Z_i, X)\to (\colim_i \Map(Z_i,Y))\times_{\Map(Z,Y)}\Map(Z,X).\]

\item\label{prop:appr-fp-morphism-1}  Let $f: X\to Y$ be a pfp morphism in $\qcsp_k$. Suppose $Y=\lim_{i\in\mI} Y_i$ is a cofiltered limit in $\qcsp_k$ with affine transition maps $Y_i\to Y_j$. Then there exists $i\in \mI$ and a pfp morphism $f_i: X_i\to Y_i$ such that $f$ is the base change of $f_i$ along $Y\to Y_i$. If $f$ is separated, resp. \'etale, resp. perfectly proper, one can choose $i$ such that $f_i$ is also separated, resp. \'etale, resp. perfectly proper.

\item\label{prop:appr-fp-morphism-2} Let $f: X\to Y$ be a morphism in $\qcsp_k$. Then $f$ can be written as a cofiltered limit $X=\lim_i X_i\to Y$ with $X_i\to Y$ pfp and $X_i\to X_j$ affine.
\end{enumerate}
In all statements as above, one can replace $\qcsp_k$ by $\psch_k$.
\end{proposition}
\begin{proof}
   \eqref{prop:appr-fp-morphism-0} follows directly from the classical characterization of finitely presented morphisms.
   For \eqref{prop:appr-fp-morphism-1}, by definition, $f$ is a perfection of a finitely presented morphism $f':X'\to Y'$ (in the usual sense). We may assume that $Y'=Y$. Then $f'$ is the base change of a finitely presented morphism $X'_i\to Y_i$ (e.g. see \cite[\href{https://stacks.math.columbia.edu/tag/07SK}{Lemma 07SK}]{stacks-project}). If $f'$ is \'etale (resp. proper), one can choose $i$ such that so is $X'_i\to Y_i$ by \cite[\href{https://stacks.math.columbia.edu/tag/084V}{Section 084V}]{stacks-project}.
   Taking the perfection gives the desired statement. Similarly, one deduces \eqref{prop:appr-fp-morphism-2} from \cite[\href{https://stacks.math.columbia.edu/tag/09NS}{Lemma 09NS}]{stacks-project}. 
\end{proof}

Parallel to the usual algebraic geometry, we make the following definition.
\begin{definition}\label{def:lpfp morphism prestack}
A morphism $f: X\to Y$ of prestacks is called locally perfectly of finite presentation (lpfp) if for every cofiltered limit $Z=\lim_i Z_i$ with $Z_i\to Z_j$ affine morphisms in $\qcsp_k$, the natural map
$\colim_i \Map(Z_i, X)\to (\colim_i \Map(Z_i,Y))\times_{\Map(Z,Y)}\Map(Z,X)$ of spaces is an equivalence.
\end{definition}

Clearly, the class of lpfp morphisms between prestacks is strongly stable.

\subsubsection{Torsors}\label{SSS: perfect torsors convention}
On $\qcsp_k$, there are several convenient topologies. We have mentioned Zariski and \'etale topology. 
There are also the pro-\'etale, $fpqc$-, and $v$-topology. The $fpqc$ topology can be defined as usual (e.g. see \cite[App. A]{zhu2017affine} \cite[App. A]{xiao2017cycles}) and the $v$-topology was introduced in \cite[\textsection{2}]{bhatt2017projectivity}. We also need 
the $fppf$- and $h$-topology on $\qcsp_k$ by requiring an $fppf$-covering (resp. an $h$-covering)  to be a pfp fpqc-covering (resp. a pfp $v$-covering).
Let $\tau$ be one of these topologies.

We repeat our conventions about torsors (as discussed in \Cref{SSS:torsors convention}) in the perfect setting.
Let $H$ be a perfect group prestack. An $H$-equivariant morphism $P\to X$ of perfect prestacks is called an $H$-torsor in the $\tau$-topology if the action of $H$ on $X$ is trivial and for every $\Spec R\to X$, there is a cover $R\to R'$ in the $\tau$-topology such that $P\times_X\Spec R'$ is a trivial, i.e. $H$-equivariantly isomorphic to $\Spec R'\times H$. We let $\bB_\tau H$ denote the prestack of $H$-torsors in $\tau$-topology. As mentioned in \Cref{SSS:torsors convention}, this is a $\tau$-stack, and sometimes a $\tau'$-stack for a finer topology $\tau'$. E.g. $\bB_{\mathrm{Zar}}\GL_n$ is a stack in $fpqc$-topology, and in fact a also stack in $v$-topology when $k$ is a perfect field of characteristic $p>0$ by \cite{bhatt2017projectivity}. If $H$ acts on a perfect prestack $X$, by the quotient $(X/H)_\tau$, we mean the $\tau$-sheafification of the prestack quotient of $X$ by $H$. So $(X/H)_\tau$ is the prestack sending $R$ to an $H$-torsor $P$ over $\Spec R$ (in $\tau$-toplogy) and an $H$-equivariant map $P\to X$. 
When $H$ is a perfect group stack (i.e. group pre-stack in \'etale topology), and that $\tau=\et$, we simply call $H$-torsors in the \'etale topology by $H$-torsors, and write $\bB H$ for $\bB_{\et} H$, and if $H$ acts on a perfect stack $X$, we write $X/H$ instead of $(X/H)_{\et}$.

\subsection{Ind-constructible sheaves on qcqs algebraic spaces}\label{SS: indconstructible-sheaf-on-qcqs}

For a scheme $X$ and a finite ring $\La$, let $\der(X_{\et},\La)$ denote the derived $\infty$-category of \'etale $\La$-modules on $X$. The assignment $X\leadsto \der(X_{\et},\La)$ can be made in a highly functorial way encoding the usual six functor formalism (e.g. see \cite{liu2012enhanced}). For our applications, however, we need some variants of $\der(X_{\et},\La)$, namely the ind-constructible (co)sheaves on $X$ as first introduced in \cite{Gaitsgory.Lurie.Weil.I}, and further studied in \cite{bouthier2020perverse} (among other works). However, these works are inadequate for our purpose, as explained before. In this section, we assemble various ingredients in literature to write down such formalism for ind-constructible sheaves on arbitrary qcqs schemes (and qcqs algebraic spaces).

\subsubsection{Constructible sheaves}\label{sec:adic-formalism}
For an ordinary topos $\frakX$, and an ordinary ring $\La$ giving a sheaf of rings on $\frakX$, 
let $D(\frakX,\La)$ denote the usual derived $\infty$-category of the abelian category of sheaves of $\La$-modules. As explained in \cite[\textsection{2.2}]{Gaitsgory.Lurie.Weil.I}, this is equivalent to the
$\infty$-category of (hypercomplete) sheaves of $\La$-modules. Let $(D(\frakX,\La)^{\mathrm{std},\leq 0}, D(\frakX,\La)^{\mathrm{std},\geq 0})$ denote its standard t-structure. The heart $D(\frakX,\La)^{\mathrm{std},\heartsuit}$ is identified with the abelian category of abelian sheaves of $\La$-modules.

We fix a perfect base commutative ring $k$.
For $X\in \psch_k$ we denote by $X_{\et}$ the small \'{e}tale site of $X$ consisting of qcqs \'etale $X$-schemes with covers given by \'{e}tale covers of schemes. Let $\mathrm{FinRing}$ denote the category of (ordinary) finite commutative rings.
Then it follows from \cite[\textsection{2}]{liu2012enhanced} that there is a lax symmetric monoidal functor
\begin{equation}\label{eq-*-pullback-etale-sheaf}
\der((-)_{\et},\La): (\psch_k)^{\op}\times \mathrm{FinRing}\to\lincat
\end{equation}
sending $(X,\La)$ to $\der(X_{\et}, \La)$, and sending a morphism from $(X,\La)$ to $(Y,\La')$ given by $(f: X\to Y,\La \rightarrow \La')$ to $\La'\otimes_\La f^*$\footnote{Note that this is equivalent to first extend coefficients to $\La'$ and then apply $f^*$.}, where $f^*$ is the ($\infty$-categorical enhancement) of the natural $*$-pullback functor. We will let $f_*$ denote the (not necessary) continuous right adjoint of $f^*$, usually called $*$-pushforward. Each $\der(X_{\et},\La)$ is a closed symmetric monoidal category (see \Cref{ss: formalism-of-abstract-sheaf-theory} in particular \eqref{eq:abstract-interior-tensor-product} and \eqref{eq:abstract-internal-hom}). We call this monoidal structure the $*$-tensor product, and denote it as $\otimes^*$. We write the internal hom bi-functor as $\underline\Hom(-,-)$.

As mentioned above, for our applications, we need several variants of the functor $X\mapsto \der(X_{\et},\La)$. First, this functor
has a ``small version". Namely, let $\der_{\ctf}(X,\La)\subset \der(X_{\et},\La)$ be the category of constructible sheaves on $X$,
which by definition is the smallest $\La$-linear idempotent complete stable subcategory of $\der(X_{\et},\La)$ spanned by objects of the form $j_!\underline{\La}_U$ for $(j: U\to X)\in X_{\et}$ where $j_!$ denotes the left adjoint of the the functor $j^*$.\footnote{This definition of constructible sheaf is different the traditional one, but is consistent with the one in \cite{Bhatt.Scholze.proetale} and \cite{Gaitsgory.Lurie.Weil.I}. Note that the homotopy category of $\der_{\ctf}(X,\La)$ is $\mathrm{D}_{\ctf}^b(X_{\et},\La)$ in the sense of Deligne.} 
As $*$-pullback and $\otimes^*$ always preserve constructibility, the functor $\der$ induces a lax symmetric monoidal functor 
\begin{equation}\label{eq-*-pullback-ctf}
     \der_{\ctf}(-,\La)\colon (\psch_k)^{\op}\to\catid_\La,\quad X\mapsto \der_{\ctf}(X,\La).
\end{equation}
We also write the functor $\der_{\ctf}$ as $\scshv$ (to be consistent with the notion used later for adic coefficients). 

For the purpose that will be clear in the sequel, we record the following finitary property of $\der_{\ctf}$. Let $\pfpsch_k\subset\psch_k$ be the full subcategory of perfect schemes perfectly finitely presented over $k$. 

\begin{lemma}\label{lem:colim-presentation-of-Dctf} 
Suppose $X = \lim_{i\in\mI} X_i$ is written as a cofiltered limit of qcqs schemes with affine transition maps. Then the natural functor
\begin{equation}\label{eq:sheaves-qcqs-inverse-limit-constructible}
\colim_{i\in \mI^{\op}} \scshv(X_i,\La) \xrightarrow{\sim} \scshv(X,\La)
\end{equation}
is an equivalence. In addition, the functor $\scshv(-,\La)$ is isomorphic to the left Kan extension of its restriction along $\pfpsch_k\subset\psch_k$.
\end{lemma}
\begin{proof}
As every $U\in X_{\et}$ is the pullback of some $U_i\in (X_i)_{\et}$, the functor in question is essentially surjective.
Then we need to show that for every two constructible sheaves $\mF, \mG\in\scshv(X,\La)$ coming as $*$-pullback of a compatible system of constructible sheaves $\mF_i,\mG_i\in \scshv(X_i,\La)$, 
\begin{equation}\label{eq-hom-space-constructible-as-colimit}
\Hom_{\scshv(X,\La)}(\mF,\mG)=\colim_i\Hom_{\scshv(X_i,\La)}(\mF_i,\mG_i).
\end{equation}
One can assume that $\mF_i=j_!\La_{U_i}$ for some $U_i\in (X_i)_{\et}$. Write $U=U_i\times_{X_i}X=\lim_j U_j$, with $U_j=U_i\times_{X_i}X_j$. Then \eqref{eq-hom-space-constructible-as-colimit} reduces to show that $H^*(U, \mG|_U)=\colim_j H^*(U_j, \mG_j|_{U_j})$, which is standard. This proves the equivalence of \eqref{eq:sheaves-qcqs-inverse-limit-constructible}.

For the second statement, we need to show that \eqref{eq:sheaves-qcqs-inverse-limit-constructible} still holds if one replaces $\mI$ by the category $\mJ=(\pfpsch_k)_{X/}$ of maps $X\rightarrow X'$ with $X'$ being pfp over $k$. But we may write $X$ as a direct limit of pfp schemes over $k$ with affine transition maps (see \Cref{prop:appr-fp-morphism} \eqref{prop:appr-fp-morphism-2}) and such system is cofinal in $\mJ$. 
\end{proof}

We have the following statement as a corollary.

\begin{corollary}\label{lem: categorical kunneth}
Assume that $k$ is an algebraically closed field. Then for $X,Y\in \psch_k$,
\[
\scshv(X,\La)\otimes_\La\scshv(Y,\La)\to \scshv(X\times Y,\La)
\]
is fully faithful.
\end{corollary}
\begin{proof}
When $X$ and $Y$ are finite type over $k$, this is well-known. The general case then follows from \Cref{lem:colim-presentation-of-Dctf}.
\end{proof}

Next we consider adic sheaves. 
Let $(\La,\frakm)$ be a pair consisting of a Noetherian ring $\La$ and an ideal $\frakm$ such that $\La$ is complete with respect to the $\frakm$-adic topology and that $\La/\frakm$ is finite.\footnote{For our purpose, it is enough to consider pairs $(\La,\frakm)$ with such assumptions, although it is possible to consider more general pairs $(\La,\frakm)$.} We call such a pair an $\frakm$-adic ring and let $\mathrm{AdicRing}$ denote the corresponding category.
It is natural to define, for a qcqs scheme $X$ over $k$, the category of $\frakm$-adic constructible sheaves\footnote{This is the category of constructible sheaves considered in \cite[\textsection{6.5}]{Bhatt.Scholze.proetale} and in \cite{H.Richarz.Scholbach.Constructible}. It can be embedded into the category $\der(X_{\proet},\La)$ or sometimes even into $\der(X_{\et},\La)$.} as
\begin{equation}\label{eq:traditional-l-adic-constructible-sheaf}
\der_{\ctf}(X,\La) = \lim_{n} \der_{\ctf}(X,\La/\frakm^n)
\end{equation}
with transition maps given by $\mF \mapsto \mF\otimes_{\La/\frakm^n} \La/\frakm^{n-1}$. 
However, for the purposes of this paper, this category is too large in general. For example, if $K$ is a profinite set considered as an affine scheme (via the inverse limit) over an algebraically closed field $k$, we would like to only consider sheaves which are ``locally constant" on $K$, instead of all ``continuous" sheaves. 

\Cref{lem:colim-presentation-of-Dctf} suggests we proceed as follows.
For a pfp scheme over $k$, we still consider the category $\der_{\ctf}(X,\La)$ of $\frakm$-adic constructible sheaves as defined via \eqref{eq:traditional-l-adic-constructible-sheaf}, and denote it by $\scshv(X, \La)$. 
As argued in \cite[\textsection{1.1}]{liu2014enhanced}, via right Kan extension along $\mathrm{FinRing}\subset\mathrm{AdicRing}$, the assignment $X\mapsto \scshv(X,\La)$ upgrades to a lax symmetric monoidal functor 
$\scshv(-,\La): (\pfpsp_k)^{\op}\to\catid_{\La}$
and we then define the functor
\begin{equation}\label{eq:adic-cons-sheaf}
\scshv(-,\La)\colon (\psch_k)^{\op}\to\catid_{\La},\quad (f\colon X\to Y)\mapsto (f^*\colon: \scshv(Y,\La)\to \scshv(X,\La))
\end{equation}
by the left Kan extension along the inclusion $\pfpsch_k\subseteq \psch_k$.

\begin{remark}
For adic ring $\La$, our notation is slightly abusive. Namely, unlike $\der_{\ctf}(X,\La)$ as defined in \eqref{eq:traditional-l-adic-constructible-sheaf}, which only depends on $X$ itself,
the category $\scshv(X,\La)$ depends on $X$ together with a morphism to $\spec k$. 
See Example \ref{ex:l-adic-sheaves-on-SpecK} below.
\end{remark}

\begin{remark}\label{rem:reduction-to-finite-coefficient}
  Note that the functor $\scshv(X,\La)\to \scshv(X,\La/\frakm),\ \mF\mapsto \mF\otimes_\La \La/\frakm$ is conservative and $*$-pullback commutes with reduction mod $\frakm$ (by definition). 
 \end{remark}

Now, let $T\subseteq \La$ be a multiplicatively closed subset and denote by $T^{-1}\La$ the localization. We define
\begin{equation}\label{eq:sheaves-adic-sheaves-localization}
\scshv\colon (\psch_k)^{\op}\to\catid_{T^{-1}\La},\quad X\mapsto \scshv(X,T^{-1}\La)=\scshv(X,\La)\otimes_{\mathrm{Perf}_\La}\mathrm{Perf}_{T^{-1}\La}.
\end{equation}
where the relative tensor product is taken in $\catid_{\La}$. 

Finally, for a filtered colimit $\Lambda = \colim_{i\in \mI} \Lambda_i$ with $\Lambda_i$ a localization of an $\frakm$-adic ring as above, we define a lax symmetric monoidal functor
\begin{equation}\label{eq:sheaves-adic-sheaves-filtered-colimit}
    \scshv\colon (\psch_k)^{\op}\to\catid_{\La},\quad X\mapsto \scshv(X,\La):=\colim_{i} \scshv(X,\La_i),
\end{equation}
by taking the colimit over $i\in\mI$ with transition functors given by extension of scalars. Note that by definition, $*$-pullback and $\otimes^*$ commute with extension of scalars $\La\to \La'$.

\begin{example}\label{ex:l-adic-E-sheaf}
Let $E/\bQ_\ell$ be an algebraic extension with ring of integers $\mO_E$. 
The preceding discussions apply to $\mO_E$ and $E$ and give functors $X\mapsto \scshv(X,\mO_E)$ and $X\mapsto \scshv(X,E)$. Explicitly, by writing $X = \lim_{i} X_i$ as a cofiltered limit of pfp schemes over $k$ with affine transition maps, and by writing $E$ as union of finite extensions $F/\bQ_\ell$ with ring of integers $\mO_F$, we have
\begin{align*}
    \scshv(X,\mO_E)  = \colim_{F\subseteq E, i} \der_{\mathrm{ctf}}(X_i,\mO_F),\quad 
    \scshv(X,E)  = \colim_{F\subseteq E, i} \der_{\mathrm{ctf}}(X_i,\mO_F)[1/\ell].
\end{align*}
\end{example}

\begin{remark}\label{rem-right-Kan-extension-*-constructible}
Note that for any coefficient $\La$ as above, $\scshv(-,\La)\colon \psch_k\to \catid_\La$ is the left Kan extension of its restriction to $\pfpsch_k$, so the analogous equivalence \eqref{eq:sheaves-qcqs-inverse-limit-constructible} still holds for any of such rings.
On the other hand, the functor $\scshv(-,\La)$ is the right Kan extension of its restriction along $\perf_k \subset (\psch_k)^{\op}$. This follows from the fact that $\scshv(-,\La)$ is a hypercomplete \'etale sheaf (in fact a $v$-sheaf, see \Cref{prop:sheaves-cshv-v-descent} below).
\end{remark}

In the rest of this section, we will fix a prime $\ell$ and allow the coefficient ring $\La$ to be any $\bZ_\ell$-algebra of the above form.
When the coefficients are clear from context or when a certain result holds for all such rings, we will occasionally omit $\La$ from the notation.

\subsubsection{Functorialty of $\scshv$}
Now we discuss functoriality of the assignment $X\mapsto \scshv(X)$ (which can be thought as a presheaf of categories on $\psch_k$). First, it is know that $\der_{\ctf}(X)$ is a hypercomplete sheaf with respect to the $v$-topology on $\psch_k$ (\cite[Theorem 2.2]{Hansen.Scholze.relative.perversity}) so in particular $\scshv(-,\La)|_{(\pfpsch_k)^{\op}}=\der_{\ctf}(-,\La)$ is an $h$-sheaf. As $\scshv$ is isomorphic to left Kan extension of $\der_{\ctf}(-,\La)$ along $(\pfpsch_k)^{\op}\subset (\psch_k)^{\op}$, an argument similar to \cite[Theorem 11.2 (2)]{bhatt2017projectivity} also gives $v$-descent of $\scshv$.

\begin{proposition}\label{prop:sheaves-cshv-v-descent}
Assume that $\La$ is a regular neotherian ring.
The functor \eqref{eq:sheaves-adic-sheaves-filtered-colimit} a hypercomplete sheaf of $\infty$-categories for the $v$-topology on $\psch_k$. 
\end{proposition}

We will repeatedly consider the following cartesian diagram in $\psch_k$ (and later on in $\prestk_k$).
\begin{equation}\label{eq:pullback-square-of-algebraic-spaces}
\begin{tikzcd}
    X' \arrow[r,"g'"]\arrow[d,"f'"] & X\arrow[d,"f"]\\
    Y'\arrow[r,"g"] & Y.
\end{tikzcd}
\end{equation}

\begin{proposition}\label{prop-proper-base-change-ind-constructible}
If $f$ is pfp proper, then $f^*$ admits right adjoint $f_*$. In addition,
we have the base change isomorphism
\[
g^*\circ f_*\to (f')_*\circ (g')^* \colon \scshv(X)\to \scshv(Y'),
\]
and  for $\mF\in\scshv(X),\ \mG\in\scshv(Y)$ the projection formula
\begin{equation*}
f_*(\mF)\otimes^* \mG \simeq f_*(\mF\otimes^* f^*(\mG)).
\end{equation*}

If $f$ is \'etale, then $f^*$ admits a left adjoint $f_!$. In addition, we have  the base change isomorphism
  \[
  (f')_! \circ (g')^*\to g^*\circ f_!,
  \]
 and for $\mF\in\scshv(X)$ and $\mG\in \scshv(Y)$, the projection formula
  \[
  f_!(\mF\otimes^* f^*(\mG))\cong f_!(\mF)\otimes^*\mG.
  \]
  
  Finally, if $f$ is pfp proper and $g$ is \'etale, then we have the base change isomorphism
  \[
  g_!\circ (f')_*\to f_*\circ (g')_!.
  \]
\end{proposition}
\begin{proof}
  As we are in the sheaf theory $\sshv$ rather than the usual theory of \'etale sheaves, the statements require justification. We give a detailed proof of the proper base change formula as many statements below can be proved via this type of argument. 
  
  If  $X,Y,X',Y'\in \pfpsch_k$ and $f$ is pfp and proper, this follows from the usual proper change isomorphism $\der_{\ctf}$. We reduce the general case to this case.

  First assume that $X,Y\in \pfpsch_k$. We write $g$ as $Y'=\lim_i Y'_i\to Y$ with each $g_i: Y'_i\to Y$ pfp and $g_{ji}: Y'_j\to Y'_i$ affine (using \Cref{prop:appr-fp-morphism} \eqref{prop:appr-fp-morphism-2}). Let $g'_i: X'_i\to X$, $f_i: X'_i\to Y'_i$ be the corresponding base change. Let $\mG\in \scshv(Y')$ coming from some $\mG_i\in \scshv(Y'_i)$. Let $\mG_j$ be the pullback of $\mG_i$ to $Y'_j$, Then using \eqref{eq-hom-space-constructible-as-colimit} and the proper base change for $\scshv|_{\pfpsch_k}$ we deduce that
  \begin{multline*}
  \Hom(\mG,g^*(f_*\mF))=\colim_j\Hom(\mG_j, g_j^*f_*\mF) =\colim_j\Hom(\mG_j,(f_j)_*(g'_j)^*\mF)\\=\colim_j\Hom((f_j)^*\mG_j,(g'_j)^*\mF)=
  \Hom((f')^*\mG,(g')^*\mF))=\Hom(\mG, (f')_*((g')^*\mF)).
  \end{multline*}
  
  Next we allows $f:X\to Y$ be pfp between arbitrary perfect qcqs schemes over $k$. Then by \Cref{prop:appr-fp-morphism} \eqref{prop:appr-fp-morphism-1}, $f$ is the base change of some pfp proper morphism $f_0: X_0\to Y_0$, with $X_0,Y_0\in\pfpsch_k$ and we may also assume that $\mF$ is a $*$-pullback of $\mF_0\in\scshv(X_0)$ along $X\to X_0$. Then using the case we just established, this situation also follows.

In a similar way, one proves that other statements hold for $\sshv$.
\end{proof}

Now we can upgrade \eqref{eq:sheaves-adic-sheaves-filtered-colimit} as a functor out of category of correspondences. 
Let $\corr(\psch_k)_{\pfp; \all}$ the category of correspondences  (see \S \ref{sec:stacks-correspondences-general}) with objects perfect qcqs schemes over $k$ and morphisms $X\dashrightarrow Y$ given by correspondences $X\xleftarrow{f} Z\xrightarrow{g} Y$ (see \eqref{E: morphism-in-corr})
with $f$ perfectly of finite presentation. As explained in \S \ref{sec:stacks-correspondences-general}, this is a symmetric monoidal category with the tensor product of objects given by the product of perfect qcqs schemes (over $k$). 

\begin{theorem}\label{thm-six-functor-formalism-ind-constructible-sheaf} 
  The theory $\scshv$ extends to a sheaf theory, denoted by the same notation
  \begin{equation}\label{eq-six-operation-*-constructible}
  \scshv(-,\La): \corr(\psch_k)_{\pfp;\all}\to \catid_\La, 
  \end{equation}
which sends a morphism  $X\dashrightarrow Y$  given by $X\xleftarrow{f} Z\xrightarrow{g} Y$ (see \eqref{E: morphism-in-corr}) to the functor denoted as $f_!\circ g^*$.  The functors satisfy the following properties: 
\begin{enumerate}
    \item If $f$ is \'{e}tale, $f_!$ is left adjoint to $f^*$ and if $f$ is pfp proper, $f_!$ is right adjoint to $f^*$. 
    \item In either of the above situation, the base change isomorphism \eqref{eq:abstract-base-change} encoded by the functor $\scshv(-,\La)$ is the Beck-Chevalley map obtained by the adjoint as in \Cref{def:categories-adjointability}.
\end{enumerate}
 \end{theorem}
\begin{proof}

First, it follows from \Cref{prop-proper-base-change-ind-constructible} and \Cref{prop-sheaf-theory-for-adjoint-factorization}, that the restriction of $\scshv|_{\mathrm{Sch}_k^{\mathrm{qc.sep}}}$ extends to a sheaf theory
\[
\scshv(-,\La): \corr(\psch_k)_{\mathrm{Pfp.sep};\all}\to \catid_\La,
\]
where $\mathrm{Pfp.sep}$ denotes the class of pfp and separated morphisms,
by noticing that the class of \'etale and pfp proper morphisms are strongly stable, and pfp separated morphisms admit factorization as a qcqs open embedding followed by a pfp proper morphism. Taking ind-completions, we obtain
\[
\sshv(-,\La): \corr(\psch_k)_{\mathrm{Pfp.sep};\all}\to \lincat_\La.
\]

Next, notice that every qcqs morphism $f: X\to Y$, there is a finite qcqs Zariski cover $\varphi: U\to X$ of $X$ such that $U\to Y$ is separated. Let $\varphi_\bullet: U_\bullet\to Y$ to denote the corresponding \v{C}ech nerve.
 Using the (easy) Zariski (co)descent of $\sshv$, we can apply \Cref{lem-sheaf-theory-left-Kan-extension-covering} to obtain
\begin{equation}\label{eq-*-pullback-ind-constructible}
\sshv(-,\La): \corr(\psch_k)_{\pfp;\all}\to \lincat_\La.
\end{equation}
Explicitly, the functor $f_!: \sshv(X,\La)\to \sshv(Y,\La)$ sends $\mF$ to the geometric realization of the simplicial object $(f\circ \varphi_\bullet)_!(\varphi_\bullet)^*\mF$. As one can replace the \v{C}ech complex by alternating \v{C}ech complex, one sees that $f_!$ preserves constructible sheaves. Therefore, by restriction, we obtain the desired functor \eqref{eq-six-operation-*-constructible}.
\end{proof}

\begin{remark}
The above argument shows that the domain of the sheaf theory $\sshv$ can be extended to $\corr(\psch_k)_{\mathrm{Pft};\all}$, where $\mathrm{Pft}$ denotes the class of perfectly finite type morphisms. 
However, such extended version does not restrict to $\scshv$ as $f_!$ may not preserve constructibility in general. (Consider the case of inclusion of a point in an infinite dimensional affine space.) For our application, it suffices (and is more convenient) to have the domain of the sheaf theory $\sshv$ as in \eqref{eq-*-pullback-ind-constructible}.
\end{remark}

Note that in general, the functor $\scshv$ from \Cref{thm-six-functor-formalism-ind-constructible-sheaf} does not give six functors as we cannot pass to right adjoint (e.g. $*$-pushforward in general does not preserve constructibility). However, under some standard finiteness assumption, right adjoints exist by \cite[Corollaire 1.5]{Deligne.SGA4.5.Finitude}.

\begin{theorem}\label{rem-*-pushforward constructibility}
Assume that $k$ is the perfection of a regular noetherian ring of dimension $\leq 1$ in which $\ell$ is invertible in $k$.\footnote{Another standard assumption is that $k$ is finite dimensional quasi-excellent noetherian ring in which $\ell$ is invertible. We will not work within this setting.} Then when restricted to $\corr(\pfpsch_k)$, the right adjoint of $f_!\circ g^*$ exists, denoted by $g_*\circ f^!$. The internal object $\underline{\Hom}(\mF,\mG)$ for $\otimes^*$-tensor product also exists. For $X\in \psch_k$, let $\consdual_X=\pi_X^!\La_{\spec k}$, where $\pi_X: X\to \spec k$ is the structural morphism. The functor
\begin{equation}\label{eq-ver-dual-*-constructible}
(\verd^{\mathrm{verd}}_X)^c: \scshv(X,\La)^{\op}\cong \scshv(X,\La), \quad \mF\mapsto \underline{\Hom}(\mF,\omega_X),
\end{equation}
is an equivalence satisfying $((\bD^{\mathrm{verd}}_X)^c)^2\cong \id$. 
\end{theorem}
Note that these functors commute with extension of scalars $\La\to \La'$.

\subsubsection{Ind-constructible sheaves}
As mentioned above, one usually cannot pass to right adjoints to obtain six operations for the sheaf theory $\scshv$ in \Cref{thm-six-functor-formalism-ind-constructible-sheaf}, in particular in non finite presentation situation.
For this reason, it is useful to consider its ind-extension as \eqref{eq-*-pullback-ind-constructible}.
Objects in $\sshv(X,\La)$ are usually called ind-constructible sheaves on $X$. In addition, as explained in \Cref{sec:symmetric-monoidal-and-projection-for-corr}, we can always pass to the right adjoint to obtain the usual six functor formalism, with  $g\horizl$ and $f\vertl $ in the abstract setup replaced by $g^*$ and $f_!$. We still write the right adjoint of $f_!$ by $f^!$ and of $g^*$ by $g_*$. Listed properties in \Cref{thm-six-functor-formalism-ind-constructible-sheaf} still hold for $\sshv(-,\La)$, so $f^* \simeq f^!$ if $f$ is \'{e}tale, and $f_!\simeq f_*$ if $f$ is proper.
We still call the monoidal structure on $\sshv(X,\La)$ the $*$-tensor product, and write the internal hom bi-functor as $\underline\Hom(-,-)$. Note that all these right adjoint functors are continuous.

Also recall that under assumptions as in \Cref{rem-*-pushforward constructibility}, there is the Verdier duality \eqref{eq-ver-dual-*-constructible}. Taking its ind-completion gives a self-duality
\begin{equation}\label{eq-ver-dual-*-sheaf}
\verd^{\mathrm{verd}}_X: \sshv(X,\La)^{\vee}\cong \sshv(X,\La),
\end{equation}
which is induced from a pairing (the co-unit in the duality datum)
\begin{equation}\label{eq-verd-pairing} 
\sshv(X,\La)\otimes_\La \sshv(X,\La)\xrightarrow{\boxtimes_{\spec k}} \sshv(X\times_k X,\La)\xrightarrow{(\Delta_X)^!}\sshv(X)\xrightarrow{\Hom(\La_X,-)} \Mod_\La.
\end{equation}

\begin{remark}
Later when we pass from the sheaf theory $\sshv$ to its dual theory, then the Verdier duality fits into the framework as in \Cref{rem-dualizability in corr}. See \Cref{rem:Verdier duality via dualizability in corr}.
\end{remark}

\begin{remark}\label{rem: comparison of ind-const and etale}
When $\La$ is finite, there are tautological functors 
\[
\sshv(X,\La)\xrightarrow{\Psi}\der(X_{\et},\La)\to \der_{\et}(X,\La),
\] 
where $\der_{\et}(X,\La)$ denotes the left-completion of $\der(X_{\et},\La)$ (with respect to the standard $t$-structure), and the first functor
sends an ind-object in $\scshv(X,\La)=\der_{\ctf}(X_{\et},\La)$ to its colimit in $\der(X_{\et},\La)$. They are all equivalences if every $U\in X_{\et}$ has bounded $\La$-cohomological dimension, e.g. $X$ is pft over a finite or an algebraically closed field $k$ (e.g. see \cite[Lemma 6.4.3, Proposition 6.4.8]{Bhatt.Scholze.proetale} or \cite[Proposition 2.2.6.2]{Gaitsgory.Lurie.Weil.I}). 

In general without such assumption, neither functor is an equivalence (see Example \ref{ex:l-adic-sheaves-on-SpecK} below). However, using the fact that for any $U\in X_{\et}$,  $R\Gamma(U_{\et},-): \der(U_{\et},\La)^{\mathrm{std},\geq 0}\to \Mod_\La^{\geq 0}$ commutes with filtered colimits, one see that the above functors restrict to equivalences
\[
\sshv(X,\bF_\ell)^{\mathrm{std},\geq 0}\cong \der(X_{\et},\bF_\ell)^{\mathrm{std},\geq 0}\cong \der_{\et}(X,\bF_\ell)^{\mathrm{std},\geq 0}
\] 
compatible with $*$-pullbacks and $*$-pushforwards, where the t-structure on $\sshv(X,\bF_\ell)$ is defined such that $\sshv(X,\bF_\ell)^{\mathrm{std},\leq 0}$ is the ind-completion of 
\[
\scshv(X,\bF_\ell)^{\mathrm{std},\leq 0}:=\scshv(X,\bF_\ell)\cap \der(X_{\et},\bF_\ell)^{\mathrm{std},\leq 0}.
\] 
It follows that $\der_{\et}(X,\bF_\ell)$ is also the left completion of $\sshv(X,\bF_\ell)$ with respect to the above standard $t$-structure.

We prefer to work with $\sshv(X,\La)$ rather than $\der_{\et}(X,\La)$ is that the former is compactly generated (by definition). On the other hand, it is not clear (to us) that whether $\der_{\et}(X,\La)$ is dualizable. 
\end{remark}

\begin{remark}\label{rmk:stanford t-structure on *-constructible}
Assume that $\La$ is regular noetherian. Recall that there is the standard $t$-structure on $\scshv(X,\La)$. The case $\La=\bF_\ell$ was mentioned in \Cref{rem: comparison of ind-const and etale}. For general $\La$, it is defined as follows: If $X$ pfp over $k$, this is a standard $t$-structure on $\scshv(X,\La)$ whose heart is the usual abelian category constructible $\La$-modules on $X$;  As $*$-pullback is $t$-exact with respect to the standard $t$-structure, we obtain the standard $t$-structure of $\scshv(X,\La)$ for any qcqs $X$ as
\[
\scshv(X,\La)^{\mathrm{std},\leq 0}=\colim_{i\in \mI^{\op}} \scshv(X_i,\La)^{\mathrm{std},\leq 0},
\]
where $X=\lim_{i\in \mI} X_i$ as in \Cref{lem:colim-presentation-of-Dctf}.  Finally, the standard $t$-structure on $\sshv(X,\La)$ is the accessible one such that $\sshv(X,\La)^{\mathrm{std},\leq 0}$ the ind-completion of $\scshv(X,\La)^{\mathrm{std},\leq 0}$.
\end{remark}

\begin{remark}\label{rmk:sheaves-qcqs-inverse-limit-constructible}
Suppose $X=\lim_i X_i$ and $Y=\lim_j Y_j$ are cofiltered limits with affine transition maps, and assume that $f$ is induced from a compatible system of morphisms $f_{ij} \colon X_i\to Y_j$.  Then  $f_*: \sshv(X)\to \sshv(Y)$ in general can be computed as follows. As $f_*$ is continuous, it is enough to compute $\mF\in \scshv(X)$, which comes from some $\mF_i\in\scshv(X_i)$. We write $r_{i'i}: X_{i'}\to X_i$ and $r_i: X\to X_i$, and $s_{j'j}: Y_{j'}\to Y_j$ and $s_j:Y\to Y_j$ to be the natural maps.
Then 
\begin{equation}\label{eq-*-pushforward-computation}
   f_*\mF=\colim_{i',j'} ((s_{j'})^*(f_{i'j'})_*(r_{i'i})^*\mF_{i}). 
\end{equation}
To prove this, we compute $\Hom(\mG,f_*\mF)$ for $\mG\in\scshv(Y)$. We may assume that $\mG=(s_j)^*\mG_j$ for $\mG_j\in \scshv(Y_j,\La)$. By increasing $i$ if necessary we may assume that there is a map $f_{ij}:X_i\to Y_j$.
Then \eqref{eq-*-pushforward-computation} follows from the following isomorphisms
\begin{multline*}
\Hom(f^*\mG, \mF)
=\colim_{i'}\Hom((r_{i'i})^*(f_{ij})^*\mG_{j}, (r_{i'i})^*\mF_{i})\\
=\colim_{i',j'}\Hom((s_{j'j})^*\mG_{j}, (f_{i'j'})_*(r_{i'i})^*\mF_{i})
=\Hom(\mG, \colim_{i',j'} ((s_{j'})^*(f_{i'j'})_*(r_{i'i})^*\mF_{i})).
\end{multline*}

Similarly we can compute $f^!$ explicitly assuming $f:X\to Y$ is pfp. Suppose $f$ is the base change of $f_0: X_0\to Y_0$, and write $f_{j}: X_{j}\to Y_{j}$ the base change of $f_0$ along $Y_j\to Y_0$.
  We suppose $\mG\in \scshv(Y)$ comes from $\mG_0\in \scshv(Y_0)$ and write $\mG_j$ for its pullback to $Y_j$. Then as argued for \eqref{eq-*-pushforward-computation}, we have
\begin{equation}\label{eq-!-pullback-computation}
f^!\mG=\colim_{j} ((r_{j})^*((f_{j})^!\mG_{j}).
\end{equation}
\end{remark}

Now we illustrate the difference between the sheaf theory $\scshv(-,\La)$ (and $\sshv(-,\La)$) and the more traditional construction \eqref{eq:traditional-l-adic-constructible-sheaf} (and \eqref{eq-*-pullback-etale-sheaf}) by two examples, especially when the space is not pfp over $k$. 

\begin{example}\label{ex:sheaves-on-profinite-sets}
Assume $k$ is a separably closed field and let $\La$ be an $\frakm$-adic ring.
For a qcqs algebraic space $X$ over $k$, we write the $*$-pushforward of the ``constant sheaf" $\underline{\La}$ to $k$ in the sheaf theoretic context $\sshv$ as $\rg_*(X,\La)$.
If $X$ is pfp over $k$ or $\La$ is finite, then $\rg_*(X,\La)=\rg(X,\La)$ is the usual \'etale cohomology of $X$ with coefficient $\La$. In general, for adic $\La$,
we have (see \eqref{eq-*-pushforward-computation})
\[
\rg_*(X,\La) \simeq \colim_i \rg(X_i,\La),
\]
where $\lim_{i}X_i$ is a presentation of $X$ as a cofiltered limit of pfp spaces over $k$ with affine transition maps.
That is, we only consider sections which "come from some $X_i$". This is in general quite different from the usual $\frakm$-adic cohomology complex of $X$, as the latter is given by
\[
\lim_{n} \rg(X,\La/\frakm^n) \simeq  \lim_{n} \colim_{i} \rg(X_i,\La/\frakm^n).
\]
In particular, let $S = \lim_{i} S_i$ be a profinite set with each $S_i$ finite, we let $\underline{S_i}=(\spec\bZ)^{S_i}$ be the corresponding constant affine scheme over $\bZ$ and $\underline{S}=\lim_i\underline{S_i}$ as an affine scheme over $\bZ$.
Then 
\[
\rg_*(\underline{S}_k,\La) \simeq \colim_i \rg((\underline{S_i})_k,\La) \simeq \colim_i \Gamma((\underline{S_i})_\La,\mO)  = \Gamma(\underline{S}_\La, \mO).
\]
In addition, for each finite set $S_i$, the category of $\La$-sheaves $\sshv((\underline{S_i})_k,\La)$ simply identifies with the category of quasi-coherent sheaves $\qcoh((\underline{S_i})_\La)$ on $(\underline{S_i})_\La$. So
\[
\sshv(\underline{S}_k,\La) \cong \colim_i\sshv((\underline{S_i})_k,\La) \cong \colim_i \qcoh((\underline{S_i})_\La) \cong \qcoh(\underline{S}_\La).
\]
\end{example}

\begin{example}\label{ex:l-adic-sheaves-on-SpecK}
Let $X=\Spec K$ be a field over $k$, and let $\Ga$ be the Galois group of $K$. 

If $\La$ is finite,  the sequence in \Cref{rem: comparison of ind-const and etale} is identified with the \eqref{eq: derived category of representations}, for $H=\underline{\Ga}_\La$ being the affine group scheme over $\La$ associated to $\Ga$ as in \Cref{ex:sheaves-on-profinite-sets}. 
Indeed, $\der(X_{\et},\La)^{\heartsuit}$ can be identified with the abelian category $\rep(\Ga,\La)^{\heartsuit}$ of smooth representations of $\Ga$, which clearly can be identified with the abelian category $\qcoh(\bB_{\mathrm{fpqc}}\underline{\Ga}_\La)^{\heartsuit}$ of algebraic representations of $\underline{\Ga}_{\Lambda}$.
So $\der(X_{\et},\La)=\der(\qcoh(\bB_{\mathrm{fpqc}}\underline{\Ga}_\La)^{\heartsuit})$ and $\der_{\et}(X,\La)=\qcoh(\bB_{\mathrm{fpqc}}\underline{\Ga}_\La)$. 
Under the equivalence, then $\scshv(X,\La)$ is identified with the idempotent complete stable  full subcategory $\crep(\Ga,\La)\subset \der(\rep(\Ga,\La)^{\heartsuit})$ spanned by the induced representations $\cind_{\Ga'}^{\Ga}\La$, with $\Ga'$ being open subgroups of $\Ga$, which can be further identified with $\Perf(\bB_{\mathrm{fpqc}}\underline{\Ga}_\La)$\footnote{For a general profinite group $\Ga$, under the equivalence $\der(\rep(\Ga,\La)^{\heartsuit})\cong \der(\qcoh(\bB_{\mathrm{fpqc}}\underline{\Ga}_\La)^{\heartsuit})$, every induced representation $\cind_{\Ga'}^{\Ga}\La$ is clearly an object in $\Perf(\bB_{\mathrm{fpqc}}\underline{\Ga}_\La)$. To see that they actually generate $\Perf(\bB_{\mathrm{fpqc}}\underline{\Ga}_\La)$, we may reduce to the case $\Ga$ is finite, and $\La=\overline{\bF}_\ell$, and then an abelian $\ell$-group. In this case, the claim is clear.}. Therefore, neither functor in \Cref{rem: comparison of ind-const and etale} is an equivalence in general, by  \Cref{ex: fpqc quotient stack}.

The situation is more complicated when $\La$ is $\frakm$-adic, as the category $\sshv(\Spec K,\La)$ is computed via an approximation of $\Spec K$ as a cofiltered limit of pfp schemes over $k$. We illustrate it by the case where $K=k(Y)$ is the (perfection of) a function field of a curve over an algebraically closed field $k$, and $\La=\bZ_\ell$. 
As a scheme, $\spec K$ is equivalent to the inverse limit $\lim U$ with $U\subseteq Y$ ranging on affine open subsets. In this case we have an equivalence
\[
\scshv(\Spec K,\bZ_\ell)=\colim_U\cRep^{\cont}(\pi_1(U),\bZ_\ell),\quad \cRep^{\cont}(\pi_1(U),\bZ_\ell):=\varprojlim_n \cRep(\pi_1(U), \bZ/\ell^n).
\]
I.e. $\scshv(\Spec K,\bZ_\ell)$ is identified with the category of (finitely generated) continuous representations of $\Ga$ \emph{unramified} almost everywhere. Indeed, by definition the l.h.s. is equivalent to the colimit of $\scshv(U_i,\bZ_\ell)$. For every affine curve $U$, the functor $\cRep^{\cont}(\pi_1(U),\bZ_\ell) \rightarrow \scshv(U,\bZ_\ell)$ is fully faithful. As any constructible sheaf on some $U$ is lisse on an open subset, we get an equivalence on colimits. Note that on the other hand, $\der_{\ctf}(\Spec K,\bZ_\ell)$ as defined in \eqref{eq:traditional-l-adic-constructible-sheaf} is the category $\cRep^{\cont}(\Ga,\bZ_\ell)$.

One can also show that in the adic case, $\sshv(\Spec K,\La)$ is in general not equivalent to $\Mod_\La$ even when $K$ is algebraically closed. (But $\sshv(\Spec K,\La)^{\heartsuit}\cong \Mod_\La^{\heartsuit}$, where the $t$-structure is defined in \Cref{rmk:stanford t-structure on *-constructible}.)
\end{example}

Now we explain a subtlety when working with $\sshv$. Namely, unlike $\scshv$, which is a $v$-sheaf on $\psch_k$, even \'etale descent can fail for $\sshv$.

\begin{example}\label{ex:sheaves-example-descent}
Assume that $\La$ is finite.
Consider the case $X = \spec K$ where $K$ is a field over $k$. Let $K_s$ be a separable closure of $K$ and $\Ga$ the Galois group of $K$. Let $\underline{\Ga}$ be the affine group scheme over $\bZ$ associated to $\Ga$ (see \Cref{ex:sheaves-on-profinite-sets}). 
Then $\spec K_s \rightarrow \spec K$ is a $\underline{\Gamma}$-torsor and we can identify its \v{Cech} nerve with the simplicial scheme $(\underline{\Gamma}_{K_s})^\bullet$ corresponding to the action of the Galois group on $K_s$. As in \Cref{ex:sheaves-on-profinite-sets}, the cosimplicial category $\sshv((\underline{\Gamma}_{K_s})^\bullet,\La)$ identifies with $\qcoh((\underline{\Ga}_\La)^\bullet)$ and therefore its totalization then identifies with $\qcoh(\bB_{\mathrm{fpqc}}\underline{\Ga}_\La)$, which in general is different from $\sshv(\spec K,\La)\cong \ind\Perf(\bB_{\mathrm{fpqc}}\underline{\Ga}_\La)$ (see \Cref{ex:l-adic-sheaves-on-SpecK} and \Cref{ex: fpqc quotient stack}). Therefore, even \'etale descent can fail for $\sshv$. (Concretely consider $k=\bQ$, $K=\bR$, $K_s=\bC$ and $\La=\bZ/2$.) Under some standard assumptions on $k$, \'etale descent of $\sshv$ can be restored, but pro-\'etale (and therefore $fpqc$-) descent still fails in general. This leads some subtleties for descent along torsors under profinite groups, which plays an important role in our applications.
\end{example}

\begin{proposition}
\label{prop:sheaves-h-descent-shv-*} 
Suppose that $k$ is the perfection of a regular noetherian ring of dimensional $\leq 1$ in which $\ell$ invertible, and suppose $k$ has finite $\ell$-cohomological dimension.
Then the functor $\sshv(-,\La)$ is a sheaf with respect to the \'etale topology on $\psch_k$.
\end{proposition}

\begin{proof}
Since finite products commute with filtered colimits in $\lincat_\La$, the functor $\sshv$ takes finite disjoint unions to products. Let $f\colon X\rightarrow Y$ be a surjective \'etale morphism, and $f_\bullet: X_\bullet\to Y$ the corresponding \v{Cech} nerve. To apply \cite[Corollary 4.7.5.3]{Lurie.higher.algebra} in this situation, it is enough to show that for every $\mF\in \scshv(Y,\La)$, there is an equivalence $|(f_\bullet)_!(f_\bullet)^*\mF|\to \mF$. 
By \Cref{prop-proper-base-change-ind-constructible}, we may reduce to the case $(f: X\to Y)\in \pfpsch_k$. 
Using that $k$  has finite $\bF_\ell$-cohomological dimension, we can further reduce to the case $\La=\bF_\ell$. Such claim then follows from the fact that $\sshv(Y,\bF_\ell)\to \der(Y_{\et},\bF_\ell)$ is an equivalence by \Cref{rem: comparison of ind-const and etale} and \'etale descent tautologically holds for $\der((-)_{\et},\bF_\ell)$.
(See also \cite[Proposition 2.3.5.1]{Gaitsgory.Lurie.Weil.I} and \cite[Corollary 3.35]{H.Richarz.Scholbach.Constructible}.)
\end{proof}

\begin{remark}\label{rem-right-Kan-extension-*-sheaf}
One can identify $\sshv$ with its right Kan extension along the inclusion $\perf_k \subseteq \psch_k$. Indeed, this just means that for any $X \in \psch_k$, the canonical map
\begin{equation}\label{eq-required.limit.affines}
\sshv(X) \rightarrow \lim_{S\in {\perf_k}_{/X}} \sshv(S)
\end{equation}
is an equivalence. Since the \'{e}tale topos of $(\psch_k)_{/X}$ is equivalent to the \'{e}tale topos associated to $(\perf_k)_{/X}$, the identification \eqref{eq-required.limit.affines} follows from the sheaf property as ${\perf_k}_{/X}$ is always a covering sieve. (See \cite[Proposition A.3.3.1]{Lurie.SAG}.)
In addition, using the arguments of \cite[\textsection{3.8}]{H.Richarz.Scholbach.Constructible} one can show that $\sshv$ is a hypercomplete \'etale sheaf.
\end{remark}

It also follows from \Cref{prop:sheaves-h-descent-shv-*} that if $f\colon X\to Y$ is surjective \'etale, then $f^*=f^!\colon \sshv(Y)\to \sshv(X)$ is conservative. For later applications, we also record the following statement.

\begin{lemma}\label{cor-conservativity-!-pullback}
Assumptions are as in \Cref{prop:sheaves-h-descent-shv-*}.
Let $f\colon X\to Y$ be a surjective morphism in $\pfpsch_k$. Then $f^!\colon \sshv(Y)\to \sshv(X)$ is conservative.
\end{lemma}
\begin{proof}
First if $f$ is surjective \'etale, this follows from  \Cref{prop:sheaves-h-descent-shv-*}. Now if $f$ is surjective perfectly smooth, then $f$ admits a section \'etale locally on $Y$, and so $f^!$ is conservative. Note that a surjective morphism of pfp schemes over $k$ is always generically perfectly smooth. (Namely one can assume that $f$ is the perfection of $f': X'\to Y'$ such that the residue field extensions at the generic points are separable). The general case then follows from noetherian induction on the dimension of $X$ and $Y$. 
\end{proof}

\begin{remark}\label{rem-sshv-algebraic-space}
Later on, we will also need the theory of sheaves for algebraic spaces. Most discussions up to now extend from $\psch_k$ to $\qcsp_k$ without change. To wit, the definitions of $\der_{\ctf},\scshv,\sshv$ make sense for qcqs algebraic spaces, and \Cref{lem:colim-presentation-of-Dctf} holds (with $\pfpsch_k$ replaced by the category $\pfpsp_k$ of pfp algebraic spaces over $k$). In addition, it is well-known that for a morphism $f$ of pfp algebraic spaces over $k$, $f_!$ preserves constructibility and satisfies base change with respect to $*$-pullback  (e.g. this can be deduced from the scheme case using noetherian induction and the fact that every qcqs algebraic space has a quasi-compact open subspace that is a qcqs scheme \cite[\href{https://stacks.math.columbia.edu/tag/0A4I}{Section 0A4I}]{stacks-project}). From this \Cref{prop-proper-base-change-ind-constructible} and \Cref{thm-six-functor-formalism-ind-constructible-sheaf} hold for $\qcsp_k$ with the same arguments and \Cref{rem-*-pushforward constructibility} also holds. \'Etale descents for $\scshv(-,\La)$ is clear when $\La$ is finite. From this, \Cref{prop:sheaves-cshv-v-descent} and \Cref{prop:sheaves-h-descent-shv-*} also hold for $\qcsp_k$, and therefore \Cref{cor-conservativity-!-pullback} also holds.
 
So from now on, we will allow the domain of our sheaf theories $\sshv$ and $\scshv$ to be $\corr(\qcsp_k)_{\pfp;\all}$.
\end{remark}

\subsection{Cohomologically (pro-)smooth morphisms}\label{sec:coh-smooth-morphisms}
The notion of a perfectly smooth morphism as in \Cref{def:morphisms-in-qcspk} is not sufficient for the purposes of this paper. Instead we will need to consider a more general class of morphisms that still behaves like smooth morphisms on the categories of sheaves, namely, the class of \textit{cohomologically smooth} morphisms.

Recall that we fix a prime $\ell$ and allow coefficient ring $\La$ to be $\bZ_\ell$-algebras as in \Cref{sec:adic-formalism}.

\subsubsection{Universal local acyclicity}
We will use the notion of universal local acyclicty in a modern formulation following the \cite{Lu.Zheng.relative.Lefschetz} and \cite{Hansen.Scholze.relative.perversity}. It is reviewed in an abstract context of general sheaf theories in \Cref{def:ULA-morphism}.

\begin{definition}\label{def.ula.using.duality}
A map $f\colon X\rightarrow Y$ in $\qcsp_k$ is called $\ell$-\textit{universally locally acyclic}, or $\ell$-ULA for short, if it is $\sshv(-,\bF_\ell)$-admissible in the sense of \Cref{def:ULA-morphism}. More generally, a sheaf $\mF\in \sshv(X,\bF_\ell)$ is $\ell$-ULA with respect to $f$ if it is $\sshv(-,\bF_\ell)$-admissible with respect to a morphism $f: X\to Y$. 
\end{definition}

\begin{remark}\label{rem-ULA-criterion-1}
By definition, $\ell$-ULA morphisms are pfp. Note that if $f:X\to Y$ is $\ell$-ULA, then it is $\sshv(-,\La)$-admissible for every $\bZ_\ell$-algebra $\La$. This follows from the criterion \Cref{lem-dual-object-in-corr-D}, which says that $f$ is $\sshv(-,\La)$-admissible if and only if
\[
(p_1)^*(f^!\La_Y)\cong (p_2)^!\La_X,
\]
where $p_i: X\times_YX\to X$ are two projections. 
\end{remark}

It follows from \Cref{cor-abstract-composition-ULA} that the class of $\ell$-ULA morphisms is weakly stable (in the sense of \Cref{def-closure-property-of-morphism}). They also satisfy several base change properties. Due to the importance, we state them explicitly here.

\begin{proposition}
\label{ula.base.change.algebraic.spaces} 
Consider a pullback square as in \eqref{eq:pullback-square-of-algebraic-spaces} with $f: X\to Y$ being $\ell$-ULA. 
Then the natural transformation (see \eqref{eq:abstract-!-*-pullback})
\begin{equation}\label{eq:natural-transform-*-to-!}
f^{!}\Lambda_Y \otimes^* f^{*} \rightarrow f^! \colon  \sshv(Y,\La)\to \sshv(X,\La) 
\end{equation}
is an isomorphism of functors, and the natural Beck-Chevalley map
\begin{equation}\label{eq:upper-star-upper-shriek-swap}
(g')^{*}\circ f^! \rightarrow (f')^! \circ g^*\colon \sshv(Y,\La)\to \sshv(X',\La),
\end{equation}
and in the case $g$ is pfp, the Beck-Chevalley map
\begin{equation}
(g')_!\circ (f')^!\to f^! \circ g_! \colon \sshv(Y',\La) \to \sshv(X,\La) 
\end{equation}
are isomorphisms of functors. In addition, $f^!$ preserves constructability.
\end{proposition}

\begin{proof}
The isomorphisms follow from \Cref{cor:abstract-base-change-II}. For the last statement, as $f$ is pfp, so is the relative diagonal $\Delta_{X/Y}$. Therefore $(\Delta_{X/Y})_!$ preserves constructibility. It follows from \Cref{lem-ula-vs-compact} that $f^!\La_Y$ is constructible. Then $f^!(-)=f^*(-)\otimes^* f^!\La_Y$ preserves constructibility.
\end{proof}

\begin{proposition}\label{pro.item-smooth.base.change.placid.qcqs.spaces} 
Consider a pullback square as in \eqref{eq:pullback-square-of-algebraic-spaces}. Assume that $f$ can be written as a cofiltered limit of $\ell$-ULA morphisms $f_i\colon X_i\rightarrow Y$ with $X\simeq \lim_i X_i$. Then the Beck-Chevalley map 
\begin{equation}
    f^*\circ g_* \rightarrow (g')_*\circ (f')^*\colon \sshv(Y',\La)\to \sshv(X,\La)
\end{equation} is an isomorphism.
\end{proposition}
\begin{proof}
  If $f$ is $\ell$-ULA, this follows from \Cref{cor:abstract-smooth-base-change}. The general case then follows by similar arguments used in \Cref{prop-proper-base-change-ind-constructible} (using  \eqref{eq-*-pushforward-computation}).
\end{proof}

\begin{lemma}\label{lem-finary-ULA}
In the situation as in \Cref{prop:appr-fp-morphism} \eqref{prop:appr-fp-morphism-1}, if $f$ is $\ell$-ULA, one can choose $f_i$ to be $\ell$-ULA.
\end{lemma}
\begin{proof}
 We may assume that $f$ is the base change of a pfp morphism $f_0: X_0\to Y_0$. Write $f_j:X_j\to Y_j$ the base change, and $p_{1j},p_{2j}: X_{j}\times_{Y_j}X_{j}\to X_j$ two projections. Then by \eqref{eq-!-pullback-computation}, $f^!\La_Y= \colim_j(r_j)^*((f_j)^!\La_{Y_j})$ and $(p_2)^!\La_X=\colim_j(r_j\times r_j)^*((p_{2j})^!\La_{X_j})$. Then isomorphism $(p_1)^*(f^!\La_Y)\cong (p_2)^!\La_X$ then comes from some $(p_{1j})^*((f_j)^!\La_{Y_j})\cong (p_{2j} )^!\La_{X_j}$. Rename $j$ as $0$ gives the claim. 
 \end{proof}

Now we compare the notion of $\ell$-ULA morphisms introduced here and the classical notion of ULA morphisms as in \cite{artin.morphismes.acycliques}.

\begin{proposition}\label{lem-categorical-ULA-vs-classical-ULA}
Suppose that $f: X\to Y$ is pfp. Then $f$ is  $\ell$-ULA if and only if
for every geometric point $\overline{x} \rightarrow X$, and a generalization $\overline{y} \rightarrow f(\overline{x})$ the map 
    \begin{equation}\label{item-classic.ula.condition}
    \La \rightarrow \rg(X_{(\overline{x})}\times_{Y_{(f(\overline{x}))}}\overline{y},\La)
    \end{equation}
    is an equivalence. More generally, $\mF\in \scshv(X,\bF_\ell)$ is $\ell$-ULA with respect to $f$ if and only if $\mF_{(\overline{x})}\to \rg(X_{(\overline{x})}\times_{Y_{(f(\overline{x}))}}\overline{y},\mF)$ is an isomorphism.     
\end{proposition}
\begin{proof}
This is proved in \cite[Theorem 2.16]{Lu.Zheng.relative.Lefschetz} (see also \cite[Theorem 4.4]{Hansen.Scholze.relative.perversity}) for ULA with respect to the usual sheaf theory $\der((-)_{\et},\bF_\ell)$. The same argument works for algebraic spaces. By the second part of \Cref{rem: comparison of ind-const and etale}, it also works for the sheaf theory $\sshv(-,\bF_\ell)$.
\end{proof}

\begin{remark}\label{ULA-v-cover-open}
\begin{enumerate}
\item Morphisms (between schemes) satisfying the condition that \eqref{item-classic.ula.condition} is an isomorphism are called locally acyclic in \cite{artin.morphismes.acycliques}. Note that \Cref{def.ula.using.duality} contains a finiteness condition (i.e. pfp) which is not imposed in the classical formulation of \cite{artin.morphismes.acycliques}. 

\item It follows from \Cref{lem-categorical-ULA-vs-classical-ULA} that $\ell$-ULA morphisms are generalizing (see \cite[\href{https://stacks.math.columbia.edu/tag/0063}{Definition 0063}]{stacks-project}). Therefore, $\ell$-ULA morphisms are universally open by \cite[\href{https://stacks.math.columbia.edu/tag/01U1}{Lemma 01U1}]{stacks-project} and surjective $\ell$-ULA morphisms are $h$-covers by \cite[Proposition (2.1)]{rydh2010submersions}.
\end{enumerate}
\end{remark}

\begin{lemma}\label{lem-ULA-detect}
Assume that $\ell$ is invertible in $Z$.
Let $X\xrightarrow{f} Y \xrightarrow{g} Z$ be a sequence of pfp morphisms with $f$ being surjective $\ell$-ULA and $g\circ f$ being $\ell$-ULA. Then $g$ is $\ell$-ULA.
\end{lemma}
\begin{proof}
We may assume that $Z$ is pfp over some $k$ as in \Cref{rem-*-pushforward constructibility}. 
We need to show that $ (p_1)^*(g^!(\bF_\ell)_Z)\cong (p_2)^!(\bF_\ell)_Y$. As all involved sheaves are constructible and $f$ is an $h$-cover, it is enough to show  such isomorphism after $!$-pullback along $X\times_ZY\to Y\times_ZY$. But this follows from the change base isomorphism \eqref{eq:upper-star-upper-shriek-swap} as both $f$ and $g\circ f$ are $\ell$-ULA.
\end{proof}

\subsubsection{Cohomologically smooth morphisms}
We assume that $\ell$ is invertible in $k$ from now on.

\begin{definition}\label{definition.ell.coh.smooth}
A morphism $f\colon X\rightarrow Y$ in $\qcsp_k$ is called \textit{cohomologically smooth} (or \textit{coh. smooth} for short)  if it is $\ell$-ULA, and the object $f^!(\bF_\ell)_Y \in \scshv(X,\bF_\ell)$ is invertible. 
\end{definition}

\begin{remark}
\begin{enumerate}
\item Clearly this notion depends on $\ell$. But as in this article we fix a prime $\ell$ for most of the time, we suppress $\ell$ from the terminology. Keep in mind that any notion introduced below that is built on \Cref{definition.ell.coh.smooth} will also depend on $\ell$, although $\ell$ may not appear explicitly. 

\item It follows from d\'{e}vissage that for $f^!\La_Y$ is invertible for every $\La$. As we shall see in \Cref{global.iso.coh.smooth}, $f^!\La_Y$ is in fact the constant sheaf $\La_X$ up to a cohomological shift and a Tate twist.

\item The notion of coh. smooth morphisms was firstly introduced in \cite[Definition 23.8]{Scholze.etale.cohomology.diamonds} in the context of $p$-adic analytic geometry. The original definition in \emph{loc. cit.} looks different from what we adapt here, but is shown in \cite[Proposition IV.2.33]{Fargues.Scholze.geometrization}  to be equivalent to a definition involving ULAness (in the corresponding sheaf-theoretic contents). 
In the case when $Y$ is a point, this notion reduces to the classical notion of rational smoothness as below.
\end{enumerate}
\end{remark}

\begin{example}\label{example.coh.smooth.over.point}
Let $K$ be a field over $k$. A morphism $f: X\to \Spec K$ is coh. smooth if and only if it is rationally smooth over $K$. I.e. $f^!\La_{\spec K}|_{X_i}\simeq  \La \langle\dim X_i \rangle$ for each connected component $X_i$ of $X$. Here and below we use $\tateshift{d}:=[2d](d)$, where $(d)$ denotes the usual $d$th Tate twist. Indeed, by choosing a deperfection of $f$, there is always a map $ \La_{X_i} \langle\dim X_i \rangle \rightarrow  f^!\La_{\Spec K}|_{X_i}$  which is non-zero, as it is an isomorphism over a dense open. As it is a map between (shifted) one dimensional local systems it must be an equivalence. 
\end{example}

\begin{example}
For an \'{e}tale morphism $j\colon U\rightarrow X$ we have a canonical identification $j^!\simeq j^*$ so is coh. smooth. In general, if $f\colon  X \rightarrow Y$ be a perfectly smooth morphism between pfp schemes over $k$ (see  \Cref{def:morphisms-in-qcspk}), then $f$ is coh. smooth.  Indeed, by the lemma below, it is enough to show that $f: \bA^1_k\to \Spec k$ is coh. smooth. But this follows as $f^!\La_{\spec k}\simeq \La\tateshift{1}$. 
\end{example}

\begin{lemma}\label{base.change.coh.smooth}
The class of coh. smooth morphisms is weakly stable (in the sense of \Cref{def-closure-property-of-morphism}).
\end{lemma}
\begin{proof}
This follows from  \Cref{ula.base.change.algebraic.spaces} and the corresponding statement for $\ell$-ULA morphisms. 
 \end{proof}

\begin{lemma}\label{coh.smooth.property.continuous}
In the situation as in \Cref{prop:appr-fp-morphism} \eqref{prop:appr-fp-morphism-1}, if $f$ is coh. smooth, one can choose $f_i$ to be coh. smooth.
\end{lemma}
\begin{proof}
We can assume that $f$ is the pullback from an $\ell$-ULA map $f_0\colon X_0 \rightarrow Y_0$ by \Cref{lem-finary-ULA}.
Set $X_i = X_{0}\times_{Y_{0}} Y_i$ and $f_i$ the corresponding map. As $\scshv(X,\bF_\ell) \simeq \colim_i \scshv(X_i,\bF_\ell)$, invertible objects in $\scshv(X,\bF_\ell)$ comes from some $\scshv(X_i,\bF_\ell)$, so $f_i$ would be coh. smooth for some $i\in\mI$.
\end{proof}

\begin{lemma}\label{composition.detection.properties.coh.smooth}
Let $f\colon X\rightarrow Y$ be a surjective $\ell$-ULA and let $g\colon Y \rightarrow Z$ be a pfp morphism. 
If $g\circ f$ is coh. smooth, then both $f$ and $g$ are coh. smooth. 
\end{lemma}
\begin{proof}
By \Cref{lem-ULA-detect}, $g$ is $\ell$-ULA. Since $(g\circ f)^!(\bF_\ell)_Z\cong f^!(\bF_\ell)_Y\otimes f^*(g^!(\bF_\ell)_Z)$, we see that both $f^!(\bF_\ell)_Y$ and $f^*(g^!(\bF_\ell)_Z)$ are invertible. Therefore, $f$ is coh. smooth. In addition, invertibility of $f^*(g^!(\bF_\ell)_Z)$ also implies that $g^!(\bF_\ell)_Z$ is invertible, since $f$ is an $h$-cover, and invertibility of constructible sheaves with respect to the $*$-tensor product can be detected after a $*$-pullback along a $v$-cover (by \Cref{prop:sheaves-cshv-v-descent}).
\end{proof}

For $\ell$-ULA morphisms, coh. smoothness can be checked on geometric fibers. 

\begin{lemma}\label{coh-smooth-ula-smooth-fibers}
Let $f\colon X \rightarrow Y$ be a pfp morphism of qcqs algebraic spaces. If there is a surjective coh. smooth morphism $Y'\to Y$ such that the base change $f': X'\to Y'$ is coh. smooth, so is $f$.  If $f$ is $\ell$-ULA and for every point $y\in |Y|$ the fiber $X_y $ is coh. smooth, then $f$ is coh. smooth.
\end{lemma}
\begin{proof}
The first assertion follows directly from \Cref{composition.detection.properties.coh.smooth}. For the second assertion,
it is enough to show that $\mF:=f^!(\bF_\ell)_Y$ is lisse. As we already know that it is constructible, it would be lisse if and only if for every two geometric points $x,x'\in X$ and a specialization $x' \rightarrow X_{(x)}$, the corresponding specialization map $\mF_{x} \rightarrow \mF_{x'}$ is an equivalence. The case that the points $x,x'$ have the same image $y = f(x) = f(x')$ in $Y$ follows from our assumption that $f$ is fiberwise coh. smooth and \Cref{ula.base.change.algebraic.spaces}. 
The general case reduces to the previous case by local acyclicity. Indeed, since $\mF$ is $\ell$-ULA with respect to $f$, for every specialization $y' \rightarrow Y_{f(x)}$ we have an equivalence $\mF_{x} \rightarrow \mF|_{X_{(x)}\times_{Y_{(f(x))}} y'}$ by \Cref{lem-categorical-ULA-vs-classical-ULA} and so we reduce to the previous case. 
\end{proof}

Given a coh. smooth morphism $f\colon X\rightarrow Y$, 
the sheaf $f^!(\bF_\ell)_Y$ is a shifted local system. We denote by $d_f \colon |X| \to \bZ$ the \textit{cohomological dimension function}, defined by 
\begin{equation}\label{eq-coh-dim-function}
d_f \colon |X|\to \bZ, \quad (f^!(\bF_\ell)_Y)_{\bar{x}} \simeq \bF_\ell \tateshift{d_f(x)}.
\end{equation}
By \Cref{base.change.coh.smooth} and \Cref{example.coh.smooth.over.point} we have
$d_f (x) = \dim_{f(x)}(X_{x})$.
In particular, the function $x \mapsto \dim_{f(x)}(X_{x})$ on $|X|$ is locally constant.

\begin{proposition}\label{global.iso.coh.smooth}
Let $f\colon X\rightarrow Y$ be a coh. smooth morphism. Then there is an isomorphism $\La_X  \tateshift{d_f} \simeq f^!\La_Y$.
\end{proposition}
We do not claim that the above isomorphism is canonical.
\begin{proof}
By d\'{e}visaage, we may assume $\La = \bF_\ell$. The object $\mF:=f^!(\La_Y)\tateshift{-d_f}$ lies in $\scshv(X,\bF_\ell)^{\heartsuit}$ and is a one dimensional \'{e}tale local system on $X$ and we want to show that this invertible local system is constant. It is enough to find a non-zero map $(\bF_\ell)_X \to \mF$. Using \Cref{coh.smooth.property.continuous}, one may assume that both $X$ and $Y$ are pfp over $k$, and in addition we may assume that $f$ arises as the perfection of a morphism $f_0: X_0\to Y_0$ of qcqs space over a regular noetherian ring $k_0$ of dimension $\leq 1$ that is (honestly) finite presented. We drop the subscript $0$ from the notation. We may further assume that $X$ is connected and $f$ is of relative dimension $d$. 

We make use of the following observation: For every open dense subset $U\subset Y$, the restriction map
\begin{equation}\label{eq-ULA-intermediate-extension}
\Hom_{\scshv(X,\bF_\ell)^{\heartsuit}}((\bF_\ell)_X, \mF)\to \Hom_{\scshv(X_U,\bF_\ell)^{\heartsuit}}((\bF_\ell)_{X_U},\mF|_{X_U})
\end{equation}
is injective. Indeed, this follows from \eqref{item-classic.ula.condition}.

Note that to construct a non-zero map $(\bF_\ell)_X\to \mF$, 
one may replace $X$ by an open subset $X'\subset X$ such that $f':X'\to Y$ is fiberwise dense in $f:X\to Y$. Indeed, giving a map $(\bF_\ell)_X\to \mF$ is equivalent to giving a map $H^{2d}f_!(\bF_\ell)_X(d)\to (\bF_\ell)_Y$, and by dimension reasons, $H^{2d}(f')_!(\bF_\ell)_{X'}(d)\cong H^{2d}(f)_!(\bF_\ell)_{X}(d)$. Therefore after replacing $X$ by $X'$ we may assume that there is some dense open subset $U\subset Y$ such that $X_U:=f^{-1}(U)\to U$ is flat. Then the canonical trace map as in \cite[Theorem 2.9]{Deligne.SGA4.Duality} gives a non-zero map $s: (\bF_\ell)_{X_U}\to \mF|_{X_U}$, and we show that it extends to $(\bF_\ell)_X\to \mF$. For this, we choose an \'etale covering $\widetilde{X}\to X$ such that $\mF|_{\widetilde{X}}$ is constant. Then the pullback of $s$ to $\widetilde{X}\times_XX_U$ extends to uniquely to the whole $\widetilde{X}$, and the two pullbacks of such extension to $\widetilde{X}\times_X\widetilde{X}$ must coincide, by the injectivity of the map \eqref{eq-ULA-intermediate-extension}. It follows that $s$ extends over $X$.
\end{proof}
In fact, in the above proof, the existence of trace map for fppf morphisms (as in \cite[Theorem 2.9]{Deligne.SGA4.Duality}) is not needed. It is enough to use the existence of trace maps over generic points of $Y$.

Besides the standard base change results \Cref{ula.base.change.algebraic.spaces} as in \Cref{pro.item-smooth.base.change.placid.qcqs.spaces}, we have the following additional one for coh. smooth morphisms.
\begin{corollary}\label{cor: base.change.coh.smooth}
Consider a pullback square as in \eqref{eq:pullback-square-of-algebraic-spaces} with $f$ coh. smooth and $g$ pfp. Then 
\[
(f')^* \circ g^! \rightarrow (g')^!\circ f^*\colon \sshv(Y)\to \sshv(X')
\] 
is an isomorphism of functors. 
\end{corollary}

Now we discuss a special class of coh. smooth morphisms.
Recall that a topological space is called acyclic if its (co)homology is the same as the (co)homology of a point. The same notion clearly makes sense in the \'etale cohomology.

\begin{definition}
A morphism $f\colon X\rightarrow Y$ in $\qcsp_k$ is called \textit{cohomologically unipotent} if it is coh. smooth and the $\bF_\ell$-cohomology of every fiber over every geometric point of $Y$ is acyclic. 
\end{definition}

Clearly, coh. unipotent morphisms are stable under base change. The following lemma implies that they are stable under compositions and therefore form a weakly stable class of morphisms.
\begin{lemma}\label{lem-equivalence-acyclic-morphism}
Let $f\colon X\rightarrow Y$ be a coh. smooth morphism in $\qcsp_k$. Then the following are equivalent:
\begin{enumerate}
\item\label{lem-equivalence-acyclic-morphism-1} $f$ is coh. unipotent;
\item\label{lem-equivalence-acyclic-morphism-2} $f_!(f^!\bF_\ell)\to \bF_\ell$ is an equivalence;
\item\label{lem-equivalence-acyclic-morphism-3} the pullback functor $f^!\colon \sshv(Y,\bF_\ell) \rightarrow \sshv (X,\bF_\ell)$ (or equivalently the pullback functor $f^*$) is fully faithful.
\item\label{lem-equivalence-acyclic-morphism-4} $\bF_\ell\to f_*(f^*\bF_\ell)$ is an equivalence;
\end{enumerate}
In addition, \eqref{lem-equivalence-acyclic-morphism-2}-\eqref{lem-equivalence-acyclic-morphism-4} hold with $\bF_\ell$ replaced by general $\La$.
\end{lemma}
\begin{proof}
Using base change (including \eqref{eq:upper-star-upper-shriek-swap}), we see that
\eqref{lem-equivalence-acyclic-morphism-1} and  \eqref{lem-equivalence-acyclic-morphism-2} are equivalent by looking at stalks. Fully faithfulness of $f^!$ is equivalent to saying that the map $f_!(f^!\mF)\to \mF$ is an equivalence for every $\mF\in\scshv(Y,\bF_\ell)$. We may assume that $\mF=j_!(\bF_\ell)_U$ for $U\in Y_{\et}$. It follows that  \eqref{lem-equivalence-acyclic-morphism-2} and  \eqref{lem-equivalence-acyclic-morphism-3} are equivalent, again by base change. Similarly, using the projection formula \Cref{prop-proper-base-change-ind-constructible}, we see that  \eqref{lem-equivalence-acyclic-morphism-3} and  \eqref{lem-equivalence-acyclic-morphism-4} are equivalent.
\end{proof}

\subsubsection{Cohomologically pro-smooth morphisms}
We will need various pro-versions of coh. smooth morphisms. 
\begin{definition}\label{def-coh-pro-smooth-morphism-space}
Let $(f\colon X\rightarrow Y)\in\qcsp_k$. 
\begin{enumerate}
\item The morphism $f$ is called \textit{pseudo cohomologically pro-smooth} (or \emph{pseudo coh. pro-smooth} for short) if there exists a presentation $X \simeq \varprojlim_{i} X_i$ as a cofiltered limit of perfect qcqs algebraic spaces with affine transition maps such that every map $X_i \to Y$ is cohomologically smooth.  The morphism $f$ is called  \emph{weakly} pseudo cohomologically pro-smooth if there is a surjective pseudo cohomologically pro-smooth morphism $U\to X$ such that the composed map $U\to Y$ is pseudo cohomogically pro-smooth.
\item The morphism $f$ is called \textit{cohomologically pro-smooth} (or \emph{coh. pro-smooth} for short) if there exists a presentation $X \simeq \varprojlim_{i} X_i$ as a cofiltered limit of perfect qcqs algebraic spaces with cohomologically smooth affine transition maps such that every map $X_i \to Y$ is cohomologically smooth.  The morphism $f$ is called  \emph{weakly} cohomologically pro-smooth if there is a surjective cohomologically pro-smooth morphism $U\to X$ such that the composed map $U\to Y$ is cohomogically pro-smooth. 
\item The morphism $f$ is called \emph{strongly} cohomologically pro-smooth if there exists a presentation $X \simeq \varprojlim_{i} X_i$ as a cofiltered limit of perfect qcqs algebraic spaces with cohomologically smooth affine surjective transition maps such that every map $X_i \to Y$ is cohomologically smooth (but $X_i\to Y$ may not be surjective), . 
\item The morphism $f$ is called \emph{essentially} cohomologically pro-smooth (or \emph{ess. coh. pro-smooth} for short)  if it can be written as $X\to X'\to Y$ with $X\to X'$ cohomologically pro-smooth and $X'\to Y$ perfectly finitely presented. 
\end{enumerate}
\end{definition}

\begin{remark} 
We apologize to introduce several different notions related to cohomological pro-smoothness. The notion of (weakly) coh. pro-smooth introduced as above seems to be a natural notion. But  in our application to Shimura varieties, we could only prove certain map is pseudo coh. pro-smooth in the above sense. That's the reason we introduce this notion. In addition, many properties of coh. pro-smooth morphisms hold for pseudo coh. pro-smooth morphisms. 
\end{remark}

\begin{remark} 
By \Cref{ULA-v-cover-open} and \cite[\href{https://stacks.math.columbia.edu/tag/0EVN}{Lemma 0EVN}]{stacks-project}, surjective weakly pseudo coh. pro-smooth morphisms are $v$-covers. In addition, strongly coh. pro-smooth morphisms are universally open.
\end{remark}

\begin{example}\label{ex: pro-smooth morphism}
We note that any pro-\'etale  (in the sense of \cite{Bhatt.Scholze.proetale}) is coh. pro-smooth. To uniform terminology, we will call a morphism $f: X\to Y$ weakly pro-\'etale if there is a surjective pro-\'etale morphism $U\to X$ such that the composed map $U\to X\to Y$ is pro-\'etale. So weakly pro-\'etale morphisms are weakly coh. pro-smooth.
Note that every  weakly \'etale morphism between affine schemes in the sense of \cite{Bhatt.Scholze.proetale} is weakly pro-\'etale in the above sense. See \cite[Theorem 2.3.4]{Bhatt.Scholze.proetale}.  

Note that a transcendental field extension $K=k(Y)/k$ as in \Cref{ex:l-adic-sheaves-on-SpecK} is coh. pro-smooth (but not strongly coh. pro-smooth).
\end{example}

\begin{example}\label{ex: pseudo pro-smooth morphism}
Here is an example of pseudo coh. pro-smooth morphism we will encounter.
Suppose we have a morphism $f: X\to Y=\lim_i Y_i$ with each $X\to Y_i$ coh. smooth. Let $X_i=X\times_{Y_i}Y$. Then $f_i: X_i\to Y$ is coh. pro-smooth. We have $f=\lim f_i: X=\lim_i X_i\to Y$. Note that for $Y_j\to Y_i$ affine, we have $X\to X\times_{Y_i}Y_j\to X$ with the second map affine and the composed map the identity. Then it is easy to see (e.g. by Serre's criterion of affineness) that $X\to X\times_{Y_i}Y_j$ is affine. Therefore, $X_j\to X_i$ is affine. It follows that $X\to Y$ is pseudo coh. pro-smooth. 
\end{example}

The following claim follows from \Cref{base.change.coh.smooth} and \Cref{coh.smooth.property.continuous} by a standard limit argument. 

\begin{lemma}\label{composition.coh.pro-smooth}
The class of (pseudo/weakly pseudo/weakly/strongly/essentially) coh. pro-smooth morphisms is weakly stable. 
\end{lemma}
\begin{proof}We only prove that ess. coh. pro-smooth morphisms are stable under compositions. It is enough to prove that if we have a pfp morphism $f\colon X\rightarrow Y$ and a coh. pro-smooth morphism $g\colon Y\rightarrow Z$ the composition $g\circ f$ is ess. coh. pro-smooth. Let $g_i \colon Y_i \rightarrow Z$ with $Y \simeq \lim_{i\in \mI} Y_i$ be a presentation of $g$ as a cofiltered limit of coh. smooth morphisms with affine coh. smooth  transition maps. Then for some $i_0\in \mI$ large enough, there exist a pfp map $f_{i_0}\colon X_{i_0} \rightarrow Y_{i_0}$ such that $X = X_{i_0}\times_{Y_{i_0}} Y$. Then, $X\to X_{i_0}\to Z$ give the desired presentation of $g\circ f$.
\end{proof}

\begin{lemma}\label{internal.hom.coh.prosmooth}
Let  $f\colon X\rightarrow Y$ be a weakly pseudo coh. pro-smooth morphism. Then there is a canonical isomorphism
\[
f^*(\rhom(\mF,\mG)) \simeq \rhom(f^*(\mF),f^*(\mG)), \quad \mF,\mG\in \scshv(Y,\La).
\]
\end{lemma}
\begin{proof}
First, if $f$ is coh. smooth, then $f^!$ exists and differs by a shift from $f^*$ so the lemma follows from the canonical isomorphism $f^!(\rhom(\mF,\mG))\simeq \rhom(f^*(\mF),f^!(\mG))$ (see \eqref{eq:abstract-hom-pull-push}).

Next we assume that $f$ is pseudo coh. pro-smooth. We need to show that for every $\mA\in \scshv(X,\La)$, 
\[
\Hom(\mA\otimes^* f^*\mF, f^*\mG)\cong \Hom(\mA, f^*(\rhom(\mF,\mG))).
\] 
We write a presentation $X = \lim_{i\in \mI} X_i$ with maps $f_i\colon X_i \rightarrow Y$ being coh. smooth, and with the transition maps affine. As every object $\scshv(X,\La)$ comes from some $X_i$, we may assume that  $\mA$ is the $*$-pull back of some $\mB_i\in \scshv(X_i,\La)$. For each $j\geq i$, let $\mB_j$ denote the $*$-pullback of $\mB_i$ to $X_j$.
Then the claim follows as
\begin{multline*}
\Hom(\mA\otimes^* f^*\mF, f^*\mG)\cong \colim_j \Hom(\mB_j\otimes^* f_j^*\mF, f_j^*\mG)\\
\cong \colim_j \Hom(\mB_j, f_j^*(\rhom(\mF,\mG)))\cong \Hom(\mA, f^*(\rhom(\mF,\mG))),
\end{multline*}
where the middle equivalence follows as $f_j$ is coh. smooth.

Finally, we assume that $f: X\to Y$ is weakly pseudo coh. pro-smooth. Let $\varphi: U\to X$ be a surjective pseudo coh. pro-smooth morphism such that $f\circ\varphi$ is pseudo coh. pro-smooth, and let $U_\bullet\to X$ be the \v{C}ech nerve of $\varphi$. We write $g_n$ for the composed map $U_n\to X\to Y$, which is pseudo coh. pro-smooth. Then we have canonical isomorphisms $(g_n)^*(\rhom(\mF,\mG)) \simeq \rhom((g_n)^*(\mF),(g_n)^*(\mG))$ of constructible sheaves. By $v$-descent, this gives a canonical isomorphism $f^*(\rhom(\mF,\mG)) \simeq \rhom(f^*(\mF),f^*(\mG)))$, as desired.
\end{proof}

We have the following pro-version of \Cref{cor: base.change.coh.smooth}.
\begin{lemma}\label{base.change.placid.spaces}
Consider a pullback square as in \eqref{eq:pullback-square-of-algebraic-spaces} with $f$  pseudo coh. pro-smooth and $g$ pfp. Then 
\[
(f')^* \circ g^! \rightarrow (g')^!\circ f^*\colon \sshv(Y)\to \sshv(X')
\] 
is an isomorphism of functors. If $g^!$ and $(g')^!$ in addition preserve constructibility, the base change isomorphism holds for $f$ weakly pseudo coh. pro-smooth.
\end{lemma}
\begin{proof}If $f$ is coh. smooth, this is \Cref{cor: base.change.coh.smooth}. Then the pro-version follows as well (using \eqref{eq-!-pullback-computation}). 
For the last statement, one can use $v$-descent for constructible sheaves as in the proof of \Cref{internal.hom.coh.prosmooth}.
\end{proof}

The following ``pro-version" of \Cref{lem-ULA-detect} in particular implies that pfp weakly coh. pro-smooth morphisms are in fact coh. smooth.
\begin{lemma}\label{lem-pro-version-of-lem-ULA-detect}
Let $g\colon Y\to Z$ be a pfp morphism. Suppose both $f\colon X\to Y$ and $g\circ f$ are pseudo coh. pro-smooth. Then $g$ is $\ell$-ULA when restricted to an open subspace of $Y$ containing $f(X)$.
If $g\circ f$ is in addition coh. pro-smooth, then $g$ is coh. smooth  when restricted to an open subspace of $Y$ containing $f(X)$.
\end{lemma}
\begin{proof}
We first use the same strategy for the proof of \Cref{lem-ULA-detect} to prove that $g$ is $\ell$-ULA (after shrinking $Y$). We may assume that $f$ factors as $f: X\to Y'\xrightarrow{g'} Z'\xrightarrow{g''} Y$ such that
\begin{itemize}
\item the composed map $Z'\to Y\to Z$ is coh. smooth;
\item the composed map $Y'\to Z'\to Y$ is coh. smooth and surjective.
\end{itemize}
In addition, we may assume that the chain of morphisms $Y'\to Z'\to Y\to Z$ descend to morphisms in $\pfpsp_k$ satisfying the same properties as above, with $k$ as in \Cref{rem-*-pushforward constructibility}.
Then one checks $ (p_1)^*(g^!(\bF_\ell)_Z)\to (p_2)^!(\bF_\ell)_Y$ is an isomorphism via $*$-pullback along the $h$-cover $Y\times_ZY'\to Y\times_ZY$. Using \Cref{cor: base.change.coh.smooth} twice, we obtain the following commutative diagram
\[
\xymatrix{
\ar_{\cong}[d](\id\times g''g')^* (p_1)^* (g^!(\bF_\ell)_Z)\ar[r]& (\id\times g''g')^*(p_2)^!(\bF_\ell)_Y\ar^\cong[d]\ar[dl] \\
 (\id\times g')^*(Y\times_ZZ'\to Z')^!(\bF_\ell)_{Z'}\ar[r]& (Y\times_ZY'\to Y')^!(\bF_\ell)_{Y'},
}
\]
giving the desired isomorphism. 

If $g\circ f$ is in fact coh. pro-smooth, we may in fact factor $X\to Y$ as $X\to Z''\to Y'\to Z'\to Y$ such that $Z'\to Z, Y'\to Y$ are coh. smooth as before and in addition $Z''\to Z'$ is coh. smooth. In this case, the above argument can be applied to deduce that
$g'': Z'\to Y$ is also a surjective $\ell$-ULA (after possibly shrinking $Z'$ and $Y$). As $Z'\to Z$ is coh. smooth, we conclude that $g$ is coh. smooth by \Cref{composition.detection.properties.coh.smooth}.
\end{proof}

The following statement can be regarded as a cohomological version of \cite[Lemma 1.1.4.(b)]{bouthier2020perverse}. However, the arguments in \emph{loc cit.} are not available in perfect algebraic geometry. This is one of the main reasons we choose to work with coh. smooth morphisms rather than perfectly smooth morphisms. In fact the analogous statement for perfectly smooth morphisms are not known to us.

\begin{lemma}\label{independence.of.placid.presentation}
Let $f_i\colon X\rightarrow Y_i, \ i=1,2$ be coh. pro-smooth morphisms with $Y_i\in \pfpsp_k,\ i=1,2$. Then both $f_i$ factor as $X\to X'\xrightarrow{f'_i} Y_i$ with both $f'_i$ being coh. smooth and $X\to X'$  coh. pro-smooth.
\end{lemma}

\begin{proof}
Let $\{X_\alpha\}_{\alpha\in A}$ be a presentation of $f_1$ as a cofiltered limit of coh. smooth maps $f_\alpha \colon X_\alpha \to Y_1$ with affine coh. smooth transition maps.
Then there is some $\alpha \in A$ such that $f_2$ factors as $X\to  X_\al\to Y_2$. Applying \Cref{lem-pro-version-of-lem-ULA-detect} to this map (and shrink $X_\al$ if necessary), we see that $X_\al\to Y_2$ is coh. smooth. So $X'=X_\al$ does the job.
\end{proof}

We similarly define (essentailly) cohomologically pro-unipotent morphisms. 
\begin{definition}\label{def-ess-coh-pro-unipotent}
A morphism $(f: X\to Y)\in \qcsp_k$ is called \emph{cohomologically pro-unipotent} if $f$ admits a presentation $\lim_i X_i\to Y$ with each $X_i\to Y$ cohomologically unipotent and transition maps affine cohomologically unipotent. The morphism $f$ is called \emph{essentially}  cohomologically pro-unipotent if $f$ admits a decomposition $X\to X'\to Y$ with $X\to X'$ cohomologically pro-unipotent and $X'\to Y$ perfectly finitely presented.
\end{definition}

\begin{remark}\label{rem-fully-faithful-pro-unipotent}
We note that coh. pro-unipotent morphisms are strongly coh. pro-smooth.
The analogue of \Cref{composition.coh.pro-smooth} holds for the class of (ess.) coh. pro-unipotent morphisms.
Since fully faithfulness is preserved under filtered colimits, the functor $f^{*}\colon \sshv(Y,\La) \rightarrow \sshv(X,\La)$ is fully faithful if $f\colon  X\rightarrow Y$ is coh. pro-unipotent. 

We also note that if $\La$ is finite, $X$ is ess. coh. pro-unipotent over an algebraically closed field $k$, then 
\[
\sshv(X,\La)\cong \der(X_{\et},\La)\cong \der_{\et}(X,\La),
\] 
by the reason mentioned in \Cref{rem: comparison of ind-const and etale}.
\end{remark}

\subsubsection{Standard placid spaces}
Placidity in algebraic geometry is meant to capture the property of having singularities of finite type. It was considered and studied in various forms by Drinfeld \cite{drinfeld2006infinite}, Raskin \cite{Raskin.dmodules.infinite.dimensional} (who coined the term), and Bouthier-Kazhdan-Varshavsky \cite{bouthier2020perverse}, among other works.

\begin{definition}\label{def: placid-algsp}
An algebraic space $X\in \qcsp_k$ is called standard placid (over $k$) if the structure morphism $X\rightarrow \spec k$ is essentially cohomologically pro-smooth. 
We denote by $\plsp_k\subset \qcsp_k$ the full subcategory consisting of standard placid algebraic spaces. 
\end{definition}

We caution the readers that although we borrow terminologies from \cite{Raskin.dmodules.infinite.dimensional} and \cite{bouthier2020perverse}, the actual meanings of these terminologies might be different from those in \emph{loc. cit.} (The meanings of the terminologies in  \cite{Raskin.dmodules.infinite.dimensional} and \cite{bouthier2020perverse} are sometimes also different.)

Recall that the lax symmetric monoidal functor \eqref{eq-six-operation-*-constructible}, whose restriction to $\corr(\pfpsp_k)$ extends to a six functor formalism under certain finiteness assumption (as in \Cref{rem-*-pushforward constructibility}).
The following statement essentially says  that six functors for constructible sheaves exist for $\plsp_k$ as well. 

\begin{proposition}\label{constructible.preservation.placid.spaces}
Assume that $k$ is the perfection of a regular noetherian ring of dimension $\leq 1$ in which $\ell$ is invertible. 
Let  $f\colon X\rightarrow Y$ be a pfp morphism between standard placid algebraic spaces. Then both $f_*$ and $f^!$ preserve the constructible subcategories. The internal hom objects between constructible sheaves on standard placid algebraic spaces are constructible. In addition, for an ess. coh. pro-unipotent morphism $f$ between standard placid algebraic spaces, $f_*$ preserves the constructible subcategories.
\end{proposition}
\begin{proof}
The first statement follows from \eqref{eq-*-pushforward-computation}, \eqref{eq-!-pullback-computation}, \Cref{ula.base.change.algebraic.spaces} and the fact that $*$-pushforward and $!$-pullback for morphisms between pfp algebraic spaces over $k$ preserve constructible subcategories. The second statement follows from \Cref{internal.hom.coh.prosmooth}. Using \Cref{lem-equivalence-acyclic-morphism}, the last statement again follows from \eqref{eq-*-pushforward-computation}.
\end{proof}

\subsubsection{Verdier duality and perverse sheaves on standard placid spaces}\label{SS-perverse-sheaf-placid-space}
Assume that $k$ is the perfection of a regular noetherian ring of dimension $\leq 1$ in which $\ell$ is invertible. 
We would like to establish a good notion of Verdier duality and perverse sheaves for $X\in\plsp_k$. In this generality, both notions depend on a choice of ``dualizing sheaf" on $X$ with respect to $k$. This is necessary as a standard placid space could be infinite dimensional, e.g. $\bA^{\infty}_k = \spec(k[x_1,\dots,]) \simeq \lim_{n} \bA^{n}_k$.
For each $n\geq 0$ the dualizing sheaf of $\bA^{n}_k$ is isomorphic to $\La[2n](n)$ and it doesn't really make sense to take $n$ to infinity in that case. Instead, one could take the constant sheaf $\La_{\bA^{\infty}_k }$ as the dualizing sheaf in this case. A slightly more general case is $X = Y\times \bA^{\infty}_k$ for some $Y$ finitely presented over $k$. Then one could take $\omega_Y\boxtimes \La_{\bA^{\infty}_k}$ as a "dualizing sheaf". A similar procedure can be done on a general standard placid space.

\begin{definition}\label{def-renormalized-dualizing-elementary-placid}
Let $X\in \plsp_k$. A \textit{generalized dualizing sheaf} is an object $\eta_X\in \scshv(X,\La)$ isomorphic to $(r^*\omega_{X'})\otimes^*\mL$ for some coh. pro-smooth morphism $r\colon X\rightarrow X'$ with $X'\in\pfpsp_k$ and some invertible object $\mL\in\scshv(X,\La)$.
\end{definition}

By  \Cref{independence.of.placid.presentation} and \Cref{global.iso.coh.smooth}, any two generalized dualizing sheaves on $X$ differ by tensoring an invertible object in $\scshv(X,\La)$. In particular, if $X$ is pfp over $k$, $\eta_X\simeq \omega_X\otimes^*\mL$ for some invertible object $\mL\in\scshv(X,\La)$.

Let $X\in \plsp_k$ equipped a generalized dualizing sheaf $\eta_X$. By \Cref{constructible.preservation.placid.spaces} we can define the corresponding \textit{Verdier duality functor}  by
\begin{equation}\label{eq-verdier.duality.shv.star.placid.space}
(\verd_{X}^{\eta,\mathrm{verd}})^c\colon \scshv(X,\La)\rightarrow \scshv(X,\La)^{\op},\quad (\verd_{X}^{\eta,\mathrm{verd}})^c(\mF) = \rhom(\mF,\eta_X)\in \scshv(X,\La).
\end{equation}
The name is justified by the following. 

\begin{proposition}\label{prop-ver-dual-elementary-placid}
Let $X\in \plsp_k$ equipped a generalized dualizing sheaf $\eta_X$. The functor \eqref{eq-verdier.duality.shv.star.placid.space} defines a bi-duality on $\scshv(X,\La)$. Namely, $((\verd_{X}^{\eta,\mathrm{verd}})^c)^2 \simeq \id_X$. In particular,
\begin{equation*}\label{eq-etaX-dual-to-LambdaX}
(\verd_{X}^{\eta,\mathrm{verd}})^c(\eta_X)\cong \La_X.
\end{equation*}
Moreover, if $\eta_X=r^*\omega_{X'}$ for a coh. pro-smooth morphism $r\colon X\rightarrow X'$ with $X'\in \pfpsp_k$, then we have canonical equivalence
\[
(\verd_{X}^{\eta,\mathrm{verd}})^c(r^*\mF) \simeq r^*((\verd^{\mathrm{verd}}_{X'})^c(\mF)), \quad \mF\in \scshv(X',\La),
\]
where $(\verd^{\mathrm{verd}}_{X'})^c$ is the standard Verdier duality functor for $X'$.
\end{proposition}
\begin{proof}
By \Cref{lem:colim-presentation-of-Dctf}, the first claim follows from the second by Verdier duality on pfp spaces over $k$, and the second statement follows from  \Cref{internal.hom.coh.prosmooth}.
\end{proof}
 
In addition, we have the following functoriality of such duality.

\begin{proposition}\label{pfp.verdier.functioriality.spaces}
Let $f: X\to Y$ is a morphism in $\plsp_k$. Let $\eta_Y$ be a generalized dualizing sheaf on $Y$.
\begin{enumerate}
\item\label{pfp.verdier.functioriality.spaces-1} If $f$ is pfp, then $\phi_X:=f^!\eta_Y$ is a generalized dualizing sheaf on $X$, and we have canonical isomorphisms of contravariant functors between constructible categories
\[
    (\verd_Y^{\eta,\mathrm{verd}})^c\circ f_!  \simeq f_* \circ (\verd_X^{\phi,\mathrm{verd}})^c,\quad
    (\verd_X^{\phi,\mathrm{verd}})^c\circ f^*  \simeq f^! \circ (\verd_Y^{\eta,\mathrm{verd}})^c.
\]

\item\label{pfp.verdier.functioriality.spaces-2} If $f$ is weakly coh. pro-smooth, then $\phi_X:= f^*\eta_Y$ is a generalized dualizing sheaf on $X$, and we have the canonical isomorphism of contravariant functors between constructible categories
\[
(\verd_X^{\phi,\mathrm{verd}})^c\circ f^*  \simeq f^* \circ (\verd_Y^{\eta,\mathrm{verd}})^c.
\]
\end{enumerate}
\end{proposition}
\begin{proof}
For Part \eqref{pfp.verdier.functioriality.spaces-1}, we note that $\phi_X=f^!\eta_Y$ is indeed a generalized dualizing sheaf by
\Cref{base.change.placid.spaces}. The rest follows from \eqref{eq:abstract-hom-pull-push}.

For Part \eqref{pfp.verdier.functioriality.spaces-2}, the case that $f$ is coh. pro-smooth is clear (using \Cref{internal.hom.coh.prosmooth} as before). 
Once this special case of Part \eqref{pfp.verdier.functioriality.spaces-2} is proved, we can use \Cref{lem-gen-dualizing-pro-smooth-local} below to conclude that $\phi_X:= f^*\eta_Y$ is a generalized dualizing sheaf on $X$ even if $f$ is just weakly coh. pro-smooth. Then we can use \Cref{internal.hom.coh.prosmooth} again to conclude that $f^*$ commutes with duality.
\end{proof}

\begin{remark}\label{rem:verdier.functioriality.spaces.pseudo.coh.pro.smooth}
Suppose $f: X\to Y$ is a pseudo coh. pro-smooth morphism between standard placid spaces and suppose for $\eta_Y$ a generalized dualizing sheaf of $Y$. 
We do not know whether $\phi_X=f^*\eta_Y$ is a generalized dualizing sheaf of $X$. 
But if it is the case then we still have $(\verd_X^{\phi,\mathrm{verd}})^c\circ f^*  \simeq f^* \circ (\verd_Y^{\eta,\mathrm{verd}})^c$. This follows from \Cref{internal.hom.coh.prosmooth}.  
\end{remark}

\begin{lemma}\label{lem-gen-dualizing-pro-smooth-local}
Let $f: X\to Y$ be a surjective weakly coh. pro-smooth morphism of standard placid spaces. If $\mF\in \scshv(Y)$ such that $f^*\mF$ is isomorphic to a generalized dualizing sheaf on $X$, then $\mF$ is a generalized dualizing sheaf on $Y$.
\end{lemma}
\begin{proof}
We choose a coh. pro-smooth morphism $r: Y\to Y'$ with $Y'\in \pfpsp_k$ and write $\eta_Y=r^*\omega_{Y'}$, and $\phi_X=f^*\eta_Y$. 
We may write that $f^*\mF\simeq \phi_X\otimes \mL^{-1}$ for some invertible sheaf on $X$. 

We have $f^*((\verd_Y^{\eta,\mathrm{verd}})^c(\mF))\cong (\verd_X^{\phi,\mathrm{verd}})^c(f^*\mF)$ is  isomorphic to $\mL$. As $f$ is an $v$-cover, we see that $(\verd_Y^{\eta,\mathrm{verd}})^c(\mF)$ is invertible.

On the other hand, we have a canonical morphism $(\verd_Y^{\eta,\mathrm{verd}})^c(\mF)\otimes^* \mF\to \eta_{Y}$ of constructible sheaves on $Y$. Taking the $*$-pullback along the $v$-cover $f$ we see that  and
 $(\verd_X^{\phi,\mathrm{verd}})^c(f^*\mF)\otimes^* f^*\mF\to \phi_X$ is an isomorphism. It follows that  and $(\verd_Y^{\eta,\mathrm{verd}})^c(\mF)\otimes^* \mF\to \eta_{Y}$ is an isomorphism. Therefore, $\mF$ is a generalized dualizing sheaf on $Y$.
\end{proof}

\begin{remark}
Our treatment of Verdier duality on standard placid algebraic spaces is inspired by \cite{drinfeld2006infinite} and \cite{Raskin.dmodules.infinite.dimensional}. But unlike \emph{loc. cit.}, we directly choose a generalized dualizing sheaf to define $\bD_X^{\eta,\mathrm{verd}}$ rather than choosing a dimension theory. This is because in perfect algebraic geometry (over a perfect field of characteristic $p>0$), a dimension theory (as defined in \emph{loc. cit.}) only determines a generalized dualizing sheaf up to non-canonical isomorphism. 
\end{remark}

Recall that $\La$ is a $\bZ_\ell$-algebra as in \Cref{sec:adic-formalism}. We now further assume that $\La$ is regular noetherian. 
Recall that in this case, for $X\in\pfpsp_k$ besides the standard $t$-structure (as discussed in \Cref{rmk:stanford t-structure on *-constructible})
there is a perverse $t$-structure $(\scshv(X,\La)^{\leq 0}, \scshv(X,\La)^{\geq 0})$ on $\scshv(X,\La)$. If $\La$ is a field (e.g. $\La=\bF_\ell$ or $\bQ_\ell$), then the perverse $t$-structure is self-dual with respect to the standard Verdier duality $(\verd_X^{\mathrm{verd}})^c$. In addition, if $f: X\to Y$ is a coh. smooth morphism, and $d_f$ is the coh. dimension function of $f$ as defined in \eqref{eq-coh-dim-function}, then $f^*[d_f]: \scshv(Y,\La)\to \scshv(X,\La)$ is perverse exact.
These facts admit a natural generalization for placid algebraic spaces.

Let $X$ be a standard placid space over $k$. For a choice of a generalized dualizing sheaf $\eta_X$, we can define a $t$-structure on $\scshv(X)$, called the $\eta$-perverse $t$-structure. Namely, if $\eta_X=r^*\omega_{X'}$ for some coh. pro-smooth morphism $r:X\to X'$ with $X'\in\pfpsp_k$, then we let
\[
\scshv(X)^{\eta,\leq 0}= \colim_{i\in \mI} \scshv(X_i)^{\leq d_i}\subset \colim_{i\in \mI} \scshv(X_i)\cong \scshv(X),
\]
where $X=\lim_{i\in \mI} X_i\to X'$ is a presentation of $r$, and $d_i$ is the coh. dimension function of the morphism $X_i\to X'$, and transition functors are $*$-pullbacks. Objects in the heart, denoted by $\mathrm{Perv}(X,\La)^\eta$, will be called as $\eta$-perverse sheaves on $X$. 
Clearly, if $\La$ is a field (e.g. $\La=\bF_\ell$ or $\bQ_\ell$), then $\mathrm{Perv}(X,\La)^\eta$ is preserved by $(\verd_X^{\eta,\mathrm{verd}})^c$.

By ind-extension, $\sshv(X)$ is equipped with an accessible $t$-structure with $\sshv(X)^{\eta,\leq 0}$ being the ind-completion of $\scshv(X)^{\eta,\leq 0}$. We still call it the $\eta$-perverse $t$-structure on $\sshv(X)$.

\begin{proposition}\label{lem-descent-perverse-sheaves}
Let $(f: X\to Y)\in \plsp_k$, and let $\eta_Y$ be a generalized dualizing sheaf on $Y$.
\begin{enumerate}
\item If $f$ is a pfp closed embedding and $\phi_X=f^!\eta_Y$, then $f_*=f_!$ is perverse exact (with respect to the $\eta$-perverse $t$-structure on $Y$ and $\phi$-perverse $t$-structure on $X$).
\item If $f$ is weakly coh. pro-smooth and $\phi_X=f^*\eta_Y$, then $f^*$ is perverse exact. If $f$ is in addition surjective, let $X_\bullet\to Y$ be the \v{C}ech nerve and let $\phi_{X_\bullet}$ be the $*$-pullback of $\eta_Y$. Then 
\[
\mathrm{Perv}(Y,\La)^\eta\cong \tot \left(\mathrm{Perv}(X_\bullet,\La)^{\phi_\bullet}\right).
\] 
\end{enumerate}
\end{proposition}

\subsection{Cosheaf theory on prestacks}\label{SS: indconstructible-cosheaf-on-qcqs}
We keep assumptions as in \Cref{rem-*-pushforward constructibility}, i.e., $k$ is the perfection of a regular noetherian ring of dimension $\leq 1$ 
and $\ell$ a prime invertible in $k$.  We allow $\La$ to be $\bZ_\ell$-algebras as in \Cref{sec:adic-formalism}. The functor $\sshv$ does not exactly fit the needs of this paper when considering categories of sheaves on certain ind-objects (see \Cref{shv over sshv}). Instead, we will  consider its dual version, which now we explain.

\subsubsection{Ind-constructible cosheaves on qcqs algebraic spaces}
Note that every $\sshv(X,\La)$ is by definition compactly generated and therefore dualizable. In fact, $\sshv$ takes value in $\cptcat$.
Therefore, as explained in \Cref{rem: codomain of sheaf theory} \eqref{rem-sheaf-theory-cpt-and-dualizable}, we may apply the duality functor $\cptcat_\La \rightarrow \cptcat_\La$ (see \eqref{eq: duality functor, cpt}) to the functor $\sshv$ to obtain a lax symmetric monoidal functor
\begin{equation}\label{category-shv-definition}
\shv(-,\La) \colon \corr(\qcsp_k)_{\pfp; \all} \to \lincat_\La,\quad X\mapsto \shv(X,\La):=\sshv(X,\La)^\vee.
\end{equation}
Explicitly, $\shv(X,\La)$ is compactly generated, with the subcategory of compact objects 
\[
\cshv(X,\La):= \scshv(X,\La)^{\op},
\]
and the functor sends a correspondence $Y\xleftarrow{f} Z\xrightarrow{g} X$ to $(f_!)^{o}\circ (g^*)^{o}$, where the superscript $o$ denotes the conjugate functor, see \eqref{eq-conjugate.functor}.  We will also consider $\cshv(-,\La)$ as a functor:
\begin{equation}\label{category-cshv-definition}
\cshv\colon \corr(\qcsp_k)_{\pfp; \all} \to \catid_\La,\quad X\mapsto \cshv(X,\La),
\end{equation}
and refer to objects in them as constructible cosheaves. 

The functor $\shv$ can be described more concretely in terms of the six functor formalism of $\sshv$. First, when restricted to $\pfpsp_k$, there is the Verdier duality functor \eqref{eq-ver-dual-*-constructible} \eqref{eq-ver-dual-*-sheaf}. 
The canonical isomorphisms of contravariant functors between constructible categories
\[
    (\verd^{\mathrm{verd}}_X)^c\circ f_*  \simeq f_! \circ (\verd^{\mathrm{verd}}_Y)^c, \quad    (\verd^{\mathrm{verd}}_X)^c\circ f^!  \simeq f^* \circ (\verd^{\mathrm{verd}}_Y)^c,
\]
allow us to identify $(f^{*})^{o}$ with $f^!$ and $(f_!)^{o}$ with $f_*$. That is, the restriction of \eqref{category-shv-definition} to
$\corr(\pfpsp_k)$ can be identified with the functor sending $X$ to the category 
\begin{equation}\label{eq:shv as indcons for pfp}
\shv(X,\La)\cong\ind\der_{\ctf}(X,\La),
\end{equation}
and sending a correspondence $X\xleftarrow{g} Z\xrightarrow{f} Y$ to the functor $f_*\circ g^!$, and \eqref{category-shv-definition} itself is isomorphic to  the left Kan extension of $\shv(-,\La)|_{\corr(\qcsp_k)}$ along the full embedding $\corr(\pfpsp_k)\subset \corr(\qcsp_k)_{\pfp;\all}$. In addition,
the restriction of $\shv$ to horizontal morphisms is equivalent to the left Kan extension of $\shv|_{(\pfpsp_k)^{\op}}$. That is, we have an equivalence
\begin{equation*}
    \shv(X,\La) \simeq \colim_{X\rightarrow X'}\shv(X',\La)
\end{equation*}
with $X'\in (\pfpsp_k)_{/X}$ and the transition functors given by $!$-pullbacks. 
Because of the above reasons, from now on, we will always write $(g^*)^{o}$ by $g^!$ and $(f_!)^{o}$ by $f_*$.  This is consistent with the usual notations in sheaf theory.

\begin{remark}
Note that this alternative description of $\shv(X,\La)$ was usually used as the definition, e.g. see \cite{bouthier2020perverse} in the $\ell$-adic sheaf setting (for $k$ an algebraically closed field and $\La=\overline{\bQ}_\ell$), \cite{Raskin.dmodules.infinite.dimensional} in the D-module setting, and \cite{richarz2020intersection} in the motivic sheaf setting. However, as we shall see our definition \eqref{category-shv-definition} allows one to quickly deduce properties of $\shv$ by dualizing the corresponding properties for $\sshv$.
\end{remark} 

\begin{remark}\label{rem:standard t-structure on shv}
Assume that $\La$ is regular noetherian. The standard $t$-structure on $\sshv(X,\La)$ (as discussed in \Cref{rmk:stanford t-structure on *-constructible}) induces
a standard $t$-structure on $\shv(X,\La)$ such that for every $f: X\to Y$, $f^!\colon \shv(Y,\La)\to \shv(X,\La)$ is $t$-exact. 
Namely, the standard $t$-structure of the category $\cshv(X,\La)=\scshv(X,\La)^{\op}$ is defined as 
\[
\cshv(X,\La)^{\mathrm{std},\leq 0}:= (\scshv(X,\La)^{\mathrm{std},\geq 0})^{\op}.
\] 
Finally, the standard $t$-structure on $\shv(X,\La)$ is the accessible one with $\shv(X,\La)^{\mathrm{std},\leq 0}$ is the ind-completion of $\cshv(X,\La)^{\mathrm{std},\leq 0}$.  Note that this $t$-structure on $\cshv(X,\La)$ is bounded, and the $t$-structure on $\shv(X,\La)$ is accessible, compatible with filtered colimits, and right complete. 

Note that if $X$ is pfp over $k$, then under the equivalence \eqref{eq:shv as indcons for pfp}, the standard $t$-structure on $\cshv(X,\La)$ as just described is different from the standard $t$-structure on $\der_\ctf(X,\La)$ as discussed in \Cref{rmk:stanford t-structure on *-constructible}.
\end{remark}

Given the above, we will refer to the symmetric monoidal structure on $\shv(X)$ encoded by the functor $\shv$ as the \textit{$!$-tensor product}. Explicitly, it is given by
\[
\shv(X,\La)\otimes_\La\shv(X,\La)\to\shv(X,\La),\quad  (\mF,\mG)\mapsto \mF\os \mG := \Delta_X^!(\mF\boxtimes_{\La} \mG).
\]
When $X$ is pfp over $k$, under the equivalence \eqref{eq:shv as indcons for pfp}
the unit of the $!$-tensor product in $\shv(X,\La)$ corresponds to the dualizing sheaf $\consdual_X$ in $ \ind\der_\ctf(X,\La)$.
For this reason, we always denote
the unit of $\shv(X,\La)$ (for any $X\in\qcsp_k$) with respect to the $!$-tensor product by $\consdual_X$.

\begin{remark}
\begin{enumerate}
\item As same notions are used in both sheaf theory $\sshv$ and $\shv$, readers should be careful which sheaf-theoretic context we are working with in the sequel. Also note that the notion of $\ell$-ULA and coh. (pro-)smooth morphisms are defined using the sheaf theory $\sshv$. 

\item Recall that the category of cosheaves on a topological space is naturally equivalent to the category of colimit preserving functors from the category of sheaves of $\spc$. For this reason, we may think $\shv(X,\La)$ as the category of ind-$\ell$-adic cosheaves on $X$. The assignment to $X$ the categories $\sshv(X)$ and $\shv(X)$ can be thought as a categorical analogue of assignment to a (nice topological) space its space of functions and its space of measures. 

In general, if $X$ is not pfp over $k$, the categories $\shv(X,\La)$ and $\sshv(X,\La)$ are not equivalent (at least not canonically). However, they are equivalent for standard placid algebraic spaces over $k$, up to a choice of a generalized dualizing sheaf by \Cref{prop-ver-dual-elementary-placid}. 
For example, in the setting as in \Cref{ex:sheaves-on-profinite-sets}, there is a canonical equivalence $\shv(\underline{S}_k)\cong \sshv(\underline{S}_k)$.
\end{enumerate}
\end{remark}

As before, one can pass to right adjoints to obtain additional functoriality encoded by $\shv$. But these right adjoints are exotic. Only some special cases are useful (e.g. see \Cref{sec:right-adjoint-of-upper-shriek} below).
In fact, to be consistent with the usual sheaf theory, we would like to have left adjoints of $g^!$ and of $f_*$, which do not always exist in general. However, we have the following statements, by translating the structures on $\sshv$ as discussed in previous sections.
\begin{proposition}\label{lem:etale-proper-functoriality-shv}
Let $f\colon X\rightarrow Y$ be a morphism, and let $f^!\colon \shv(Y)\to \shv(X)$ the induced $!$-pullback functor. If $f$ is pfp, we also have the $*$-pushforward functor $f_*\colon \shv(X)\to \shv(Y)$.
\begin{enumerate}
    \item\label{lem:etale-proper-functoriality-shv-1} If $f$ is \'{e}tale, $f_*$ is a right adjoint to $f^!$, and if $f$ is pfp proper, $f_*$ is a left adjoint to $f^!$. In either of the above situation, the base change isomorphism \eqref{eq:abstract-base-change} encoded by the functor $\shv$ is the Beck-Chevalley map obtained by the adjoint as in \Cref{def:categories-adjointability}.

     \item\label{lem:etale-proper-functoriality-shv-5} If $f$ is an $\ell$-ULA morphism, then $f_*$ admits a left adjoint $f^*$, which preserves constructibility. In addition, for a pullback square as in \eqref{eq:pullback-square-of-algebraic-spaces} (with $f$ being $\ell$-ULA), there is the natural base change isomorphism of functors from $\shv(Y,\La)$ to $\shv(X',\La)$
     \[
     (f')^*\circ g^!\to (g')^!\circ f^*.
     \] 
     \item\label{lem:etale-proper-functoriality-shv-6} If $f$ is a pfp morphism between standard placid spaces, then $f^!$ and $f_*$ admit \emph{left} adjoints, denoted by $f_!$ and $f^*$ respectively, which preserve constructibility.  In addition, for a pullback square as in \eqref{eq:pullback-square-of-algebraic-spaces} with $g\colon Y'\rightarrow Y$ being weakly pseudo coh. pro-smooth, there are the natural base change isomorphisms of functors
    \begin{equation*}
         (f')_!\circ (g')^!\rightarrow g^!\circ f_!, \quad (f')^*\circ g^! \to  (g')^!\circ f^*.
    \end{equation*}
    
    \item\label{pro-unipotent.are.admissible} Let $f\colon X\rightarrow Y$ be a coh. pro-unipotent morphism of standard placid spaces. Then $f^!$ admits a left adjoint $f_!$, which then automatically preserves constructible subcategories. In addition, for a pullback square as in \eqref{eq:pullback-square-of-algebraic-spaces} with $g\colon Y'\rightarrow Y$ being weakly pseudo coh. pro-smooth, then there is natural base change isomorphism
    \begin{equation*}
         (f')_!\circ (g')^!\rightarrow g^!\circ f_!.
    \end{equation*}
\end{enumerate}
\end{proposition}
\begin{proof} We only discuss Part \eqref{lem:etale-proper-functoriality-shv-5}-\eqref{pro-unipotent.are.admissible}. The existences of left adjoints are based on the following observation: For a morphism $f: X\to Y$, if the $*$-pushforward in the $\sshv$-sheaf theory preserves $\scshv$, then in the $\shv$-sheaf theory, $f^!$ admits a left adjoint $f_!$ preserving $\cshv$. Similarly, if $f$ is a pfp morphism such that the $!$-pullback in the $\sshv$-sheaf theory preserves $\scshv$, then in the $\shv$-sheaf theory, $f_*$ admits a left adjoint $f^*$ preserving $\cshv$. Under our assumptions, the functors in question preserve constructibility by \Cref{ula.base.change.algebraic.spaces}  and \Cref{constructible.preservation.placid.spaces}.

To prove the base change isomorphisms in Part \eqref{lem:etale-proper-functoriality-shv-5}-\eqref{pro-unipotent.are.admissible}, we may restrict our attentions to constructible sheaves, as all involved functors are continuous preserving constructibility. 
Then the base change isomorphism in Part \eqref{lem:etale-proper-functoriality-shv-5} follows from \eqref{eq:upper-star-upper-shriek-swap} by passing to the opposite categories. The base change isomorphisms in Part \eqref{lem:etale-proper-functoriality-shv-6} for $g$ being pseudo coh. pro-smooth follow by restricting \Cref{pro.item-smooth.base.change.placid.qcqs.spaces} and \Cref{base.change.placid.spaces} to constructible subcategories and then passing to the opposite categories. Then for $g$ being weakly pseudo coh. pro-smooth, one can apply $v$-descent to conclude.

To prove the base change isomorphism \eqref{pro-unipotent.are.admissible}, we just notice that every constructible object on $X$ comes from the $!$-pullback of some object on $X_i$ with $X\to X_i$ unipotent. Then the base change follows from Part \eqref{lem:etale-proper-functoriality-shv-6}.
\end{proof}

Descent results for $\sshv$ and for $\scshv$ can also be translated to descent results for $\shv$ and $\cshv$.

\begin{proposition}\label{prop-h-descent-shv}
\begin{enumerate}
\item\label{prop-h-descent-shv-1} The theory $\cshv|_{(\qcsp_k)^{\op}}$ is a hypersheaf for the $v$-topology on $\qcsp_k$.  
\item\label{prop-h-descent-shv-2} Suppose $k$ has finite $\bF_\ell$-cohomological dimension. Then $\shv|_{(\qcsp_k)^{\op}}$ is an $h$-sheaf.  
\end{enumerate}
\end{proposition}
\begin{proof}
Part \eqref{prop-h-descent-shv-1} is obtained from \Cref{prop:sheaves-cshv-v-descent} by passing to opposite categories (taking \Cref{rem-sshv-algebraic-space} into account). 

For  \eqref{prop-h-descent-shv-2}, we first prove descent with respect to surjective pfp proper morphism $f: X\to Y$. In this case $f^!$ admits left adjoint $f_*$. As in the argument of  \Cref{prop:sheaves-h-descent-shv-*}, it is enough to show that $|(f_\bullet)_*(f_\bullet)^!\mF|\to \mF$ is an equivalence for $\mF\in \shv(Y)$, and then one can reduce to the case $X,Y$ are pfp. In this case $\shv\cong \sshv$.
Using \cite[Corollary 4.7.5.3]{Lurie.higher.algebra}, it is enough to show that $f^!\colon \sshv(Y)\to \sshv(X)$ is conservative. But this follows from \Cref{cor-conservativity-!-pullback}. 

Next we note that $\shv|_{(\psch_k)^{\op}}$ satisfies Zariski descent. Indeed, it is enough to check for the cover $X=U\sqcup V\to Y=U\cup V$. In this case $\tot(\shv(X_\bullet))$ can be computed as the finite limit $\shv(U)\times\shv(V)\rightrightarrows \shv(U\cap V)$, which commutes with filtered colimits. Then the desired descent follows from the descent for $\cshv$. This implies that $\shv|_{(\psch_k)^{\op}}$ satisfies $h$-descent, and therefore in particular \'etale descent.

Next, consider the case $X\to Y$ is \'etale with $X$ a perfect qcqs scheme. By \cite[\href{https://stacks.math.columbia.edu/tag/09YC}{Proposition 09YC}]{stacks-project}, there is a pfp proper surjective (in fact the perfection of a finite) morphism $Y'\to Y$ with $Y'$ being a scheme. Base change gives $Y'\to Y$. Now $X\times_YY'\to Y'$ satisfies descent by the scheme case and $Y'\to Y$ satisfies descent by surjective proper case. So $X\times_YY'\to Y$ and then $X\to Y$ satisfy descent by \Cref{Liu.Zheng.descent.lemma}. Finally, the case of general surjective \'etale morphism $X\to Y$ of algebraic spaces also follows from \Cref{Liu.Zheng.descent.lemma} by choosing a surjective \'etale morphism $X'\to X$ with $X'$ being a perfect qcqs scheme.
\end{proof}

\subsubsection{Verdier duality for cosheaves and perverse cosheaves}\label{sec-verdier.duality.dual.setting}
As we have seen in \Cref{prop-ver-dual-elementary-placid}, for $X\in \plsp_k$ equipped with a generalized dualizing sheaf $\eta_X\in \scshv(X)$ there is a self duality on $\sshv(X,\La)$, which can also be interpreted as an equivalence 
\begin{equation}\label{eq: identifying star sheaf with shrek sheaf}
\id^\eta: \shv(X)\simeq \sshv(X),
\end{equation}
which restricts to an equivalence
\begin{equation}\label{eq: identifying star sheaf with shrek sheaf-1}
\id^\eta: \cshv(X)\cong \scshv(X)
\end{equation}
If $r: X\to Y$ is a weakly coh. pro-smooth morphism between standard placid spaces over $k$, and if
 $\phi_X=r^*\eta_Y$, then 
 \[
 \id^\phi\circ r^! \simeq r^*\circ \id^\eta.
 \]

As $\cshv(X)=\scshv(X)^{\op}$ and $\shv(X,\La)=\sshv(X,\La)^\vee$, such duality can also interpreted as the form
\begin{equation}\label{eq-verdier.sshv.shv}
(\verd_X^\eta)^c\colon \cshv(X)^{\op}\simeq \cshv(X),\quad (\verd_X^\eta)^c\colon \shv(X)^\vee \simeq \shv(X).
\end{equation}

We regard $\eta_X$ as an object in $\scshv(X)^{\op}=\cshv(X)$, called the generalized constant sheaf of $X$ and denoted by
$\Lambda_X^{\eta}$. Then we may define an $\La$-linear functor
\begin{equation}\label{eq: eta-cohomology}
\rg^\eta(X,-):= \Hom_{\shv(X)}(\La_X^\eta,-): \shv(X)\to\Mod_\La.
\end{equation}
This is in fact a Frobenius structure of $\shv(X)$ such that  \eqref{eq-verdier.sshv.shv} is the induced self-duality of $\shv(X)$ as in \Cref{ex: duality via Frobenius-structure}. I.e., we have
\begin{equation}\label{eqn:characterization of verdier.pairing.placid.spaces}
\Hom_{\cshv(X)}(\mF,\mG) \simeq \rg^\eta(X, (\verd_{X}^\eta)^c(\mF)\os \mG),\quad \mF, \mG\in \cshv(X).
\end{equation}

\begin{remark}\label{rem:Verdier duality via dualizability in corr}
Our choice of notation is justified by the fact that when $X\in \pfpsp_k$ and when $\eta_X=\consdual_X$ is the usual canonical sheaf of $X$, then $\La_X^\eta=\La_X$
 is the usual constant sheaf on $X$, under the equivalence $\shv(X,\La)\cong \ind\der_{\ctf}(X,\La)$ (see \eqref{eq:shv as indcons for pfp}). In this case the right hand side of \eqref{eqn:characterization of verdier.pairing.placid.spaces} is just \eqref{eq-verd-pairing}. In particular, if $k$ is an algebraically closed field, then $\rg^\eta(X,-)$ given by the $*$-pushforward along $\pi_X: X\to \pt=\spec k$, and therefore fits into the paradigm of \Cref{rem-dualizability in corr}. 
 \end{remark}
 
Note that, if $\eta_X=r^*\consdual_{X'}$ for some coh. pro-smooth morphism $r: X\to X'$ with $X'\in \pfpsp_k$ and $\consdual_{X'}\in \scshv(X')$ is the canonical sheaf of $X'$, then 
\[
\La_X^\eta=r^!\La_X\in \cshv(X).
\]
Also note that for any choice of $\eta_X$, we have
\begin{equation}\label{eq-LaetaX-dual-to-omegaX}
(\verd_X^{\eta})^c(\La^\eta_X)\cong \consdual_X \in \cshv(X,\La).
\end{equation}

We have the dual version of \Cref{pfp.verdier.functioriality.spaces}. More precisely:
\begin{lemma}\label{verdier.functioriality.spaces-shv}
Let $f\colon X\to Y$ be as in  \Cref{pfp.verdier.functioriality.spaces}. 
\begin{enumerate}
\item\label{verdier.functioriality.space-shv-1} If $f$ is perfectly finitely-presented, then $\La_X^{\phi}=f^*\La_Y^{\eta}$ is a generalized constant sheaf, and we have isomorphisms of contravariant functors (for $\cshv$)
\begin{equation*}
    (\verd_Y^{\eta})^c\circ f_*  \simeq f_! \circ (\verd_X^{\phi})^c,\quad
    (\verd_X^{\phi})^c\circ f^!  \simeq f^* \circ (\verd_Y^{\eta})^c.
\end{equation*}
\item\label{verdier.functioriality.space-shv-2} If $f$ is weakly coh. pro-smooth, then $\La^\phi_{X}= f^!\La^\eta_Y$ is a generalized constant sheaf, and we have an isomorphism of contravariant functors:
\begin{equation*}
(\verd_X^{\phi})^c\circ f^!  \simeq f^! \circ (\verd_Y^{\eta})^c.
\end{equation*}
\end{enumerate}
\end{lemma}

\begin{remark}\label{rem-dual-perverse-sheaf-placid-space}
We assume that $\La$ is regular noetherian.
By transport of structure, for a standard placid space $X$, the $\eta$-perverse $t$-structure on $\scshv(X)$ (and on $\sshv(X)$) from \Cref{SS-perverse-sheaf-placid-space} 
corresponds to a $t$-structure on $\cshv(X)$ (and on $\shv(X)$), which we still call the $\eta$-perverse $t$-structure. More precisely, we let
\[
\cshv(X)^{\eta,\leq 0}=(\scshv(X)^{\eta, \geq 0})^{\op},
\]
and let $\shv(X)^{\eta,\leq 0}$ be the ind-completion of $\cshv(X)^{\eta,\leq 0}$. 
Note that when $X$ is pfp over $k$ and $\eta_X=\consdual_X$, then
under the equivalence \eqref{eq:shv as indcons for pfp}, $\eta$-perverse $t$-structure on $\cshv(X)$ corresponds to the Verdier dual of the usual perverse $t$-structure on $\der_{\ctf}(X)$. In particular, when $\La$ is a field, it coincides with the usual perverse $t$-structure.

Note that
\Cref{lem-descent-perverse-sheaves} has a corresponding dual version.
\end{remark}

\subsubsection{The functor $f_\flat$}\label{sec:right-adjoint-of-upper-shriek}
Now we discuss the continuous right adjoint of $f^!\colon \shv(Y)\to \shv(X)$ for  a coh. pro-smooth morphism $f\colon X\rightarrow Y$, denoted by $f_{\flat}$\footnote{
In the context of the standard theory of constructible sheaves (and in the motives literature), it is customary to denote by $f_\sharp$ the left adjoint of $f^*$ for $f$ smooth. For such morphisms $f_\sharp$ identifies with the conjugate $(f_\flat)^\circ$, which is the reason for our notation.
}. Heuristically in this case, the functor $f^!$ should behave like a shifted version of $*$-pullback, so $f_\flat$ should behave like a ``renormalized" version of $*$-pushforward. We summarized some of the important properties of the functor $f_\flat$.

\
First, suppose that $f$ is in fact coh. smooth (so in particular is pfp). 
The isomorphism \eqref{eq:natural-transform-*-to-!} then translates to a natural isomorphism
\[
f^* \xrightarrow{\sim} f^*(\omega_Y) \os f^!,
\]
from which we obtain an expression of $f_\flat$ 
\begin{equation}\label{eq-flat.pushforward.coh.smooth}
f_\flat(\mF) \simeq f_*(f^*(\omega_Y)\os \mF).
\end{equation}
It follows that $f_\flat$ preserves the constructible categories.

\begin{proposition}\label{lem:dual-sheaves-proj-formula-coh-pro-smooth}
Pseudo coh. pro-smooth morphisms in $\qcsp_k$ satisfy \Cref{assumptions.base.change.sheaf.theory.H}. Concretely, this means that if we
consider a pullback square \eqref{eq:pullback-square-of-algebraic-spaces} in $\qcsp_k$ with $f$ pseudo coh. pro-smooth, then for every $\mF\in \shv(X)$, $\mG\in \shv(Y)$, we have the natural isomorphism 
\begin{equation}\label{eq-exotic-projection}
f_\flat (\mF)\os \mG \xrightarrow{\sim} f_\flat(\mF\os f^!\mG).
\end{equation} 
In addition, we have a natural isomorphism of functors
\begin{equation}\label{eq:sheaves-dual-ula-base-change-1}
g^! \circ f_\flat \xrightarrow{\sim} (f')_\flat \circ (g')^!:  \shv(Y')\to \shv(X).
\end{equation}
If in addition $g$ is pfp, then we have a natural isomorphism of functors
\begin{equation}\label{eq:sheaves-dual-ula-base-change-2}
g_*\circ (f')_\flat\to f_\flat\circ (g')_*:  \shv(X')\to \shv(Y).
\end{equation}
\end{proposition}
\begin{proof}
By \Cref{composition.coh.pro-smooth}, the class of pseudo coh. pro-smooth morphisms is weakly stable.

If $f$ is coh. smooth, the isomorphism \eqref{eq-exotic-projection} follows from \eqref{eq-flat.pushforward.coh.smooth} and the usual projection formula for $(f_*,f^!)$ encoded by the sheaf theory $\shv$
\[
f_\flat (\mF)\os \mG \simeq f_*(f^*(\omega_Y)\os \mF)\os \mG \simeq f_*(f^*(\omega_Y)\os \mF\os f^!\mG) \simeq f_\flat(\mF\os f^!\mG).
\]
Now let $f$ be a general pseudo coh. pro-smooth and fix a presentation of $f$ as a cofiltered limit of coh. smooth maps $f_i \colon X_i \rightarrow Y$. 
By \Cref{lem:colim-presentation-of-Dctf} (and \Cref{rem-right-Kan-extension-*-constructible}), we have 
\[
\shv(X,\La)=\colim_i \shv(X_i,\La)
\] 
with transitioning maps in the colimit presentation being $!$-pullbacks and in the limit presentation being $\flat$-pushforwards. Then every $\mF\in\cshv(X,\La)$ arises as the $!$-pullback of some $\mF_i\in\cshv(X_i,\La)$ for some $i$ and 
we have 
\begin{equation}\label{eq-formula-for-fflat-pro-smooth}
f_\flat \mF=\colim_{j\geq i} (f_j)_\flat \mF_j,
\end{equation} 
where $\mF_j$ is the $!$-pullback of $\mF_i$ along $X_j\to X_i$. (See the reasoning in \Cref{rmk:sheaves-qcqs-inverse-limit-constructible}).
Then \eqref{eq-exotic-projection} follows from the projection formula for each $f_j$ via \eqref{eq-formula-for-fflat-pro-smooth}.

Similarly, to prove \eqref{eq:sheaves-dual-ula-base-change-1} and \eqref{eq:sheaves-dual-ula-base-change-2},
we first assume that $f$ is coh. smooth. Then  \eqref{eq:sheaves-dual-ula-base-change-1} follows from the base change encoded by the sheaf theory $\shv$, together with \eqref{eq-flat.pushforward.coh.smooth} and \Cref{lem:etale-proper-functoriality-shv} \eqref{lem:etale-proper-functoriality-shv-5}.  To prove  \eqref{eq:sheaves-dual-ula-base-change-2}, by \eqref{eq-flat.pushforward.coh.smooth} we just need to show 
\[
g_*((f')_*(\mF\otimes^! (f')^*\consdual_X))\cong f_*(  ((g')_*\mF)\otimes^! f^*\consdual_Y).
\] 
Again by \eqref{eq-flat.pushforward.coh.smooth} and \Cref{lem:etale-proper-functoriality-shv} , $(f')^*\consdual_X= (g')^! f^*\consdual_Y$ so the desired isomorphism follows from the usual projection formula encoded by $\shv$ (as in \eqref{eq:abstract-projection-formula}). 
The general case that $f$ is pseudo coh. pro-smooth follows again by using \eqref{eq-formula-for-fflat-pro-smooth}.
\end{proof}

Next we consider $f_\flat$ for weakly coh. pro-smooth morphisms. We do not know whether \Cref{lem:dual-sheaves-proj-formula-coh-pro-smooth} holds in this case. But when we restrict to the category of standard placid spaces, we have similar statements.

\begin{proposition}\label{lem:dual-sheaves-proj-formula-weakly-coh-pro-smooth}
Consider a pullback square of in $\qcsp_k$ as in \eqref{eq:pullback-square-of-algebraic-spaces} with $f: X\to Y$ being a weakly coh. pro-smooth morphism in $\plsp_k$. Then \eqref{eq-exotic-projection}-\eqref{eq:sheaves-dual-ula-base-change-2} hold in this setting.
\end{proposition}
\begin{proof}
We first prove \eqref{eq-exotic-projection} in this setting.
We can assume that $\mF,\mG$ as in \eqref{eq-exotic-projection}  are constructible. 
We need to show that for every $\mA\in \cshv(Y,\La)$,
\begin{equation*}\label{eq-renormalized.projection.formula.proof.plstk}
\Hom_{\shv(Y)}(\mA,f_\flat(\mF)\os \mG) \rightarrow \Hom_{\shv(Y)}(\mA,f_\flat(\mF\os f^!(\mG)))
\end{equation*}
is an isomorphism. Choose a generalized constant sheaf $\La_Y^\eta$ on $Y$ and let $\La_X^\phi=f^!\La_Y^\eta$. Then by \Cref{verdier.functioriality.spaces-shv} \eqref{verdier.functioriality.space-shv-2} and by \eqref{eqn:characterization of verdier.pairing.placid.spaces},
the left hand side can be identified with
\[
\Hom_{\shv(Y)}(\La_Y^\eta, (\verd_Y^{\eta})^c(\mA)\os f_\flat(\mF)\os \mG)=\Hom_{\shv(X)}(f^!((\verd^\eta_Y)^c((\verd_Y^\eta)^c(\mA)\os \mG)),\mF),
\]
while the right hand side can also be identified with
\begin{multline*}
\Hom_{\shv(X)}(f^!\mA,\mF\os f^!\mG)=\Hom_{\shv(X)}(\La_X^\phi,(\verd_X^\phi)^c (f^!\mA)\os \mF\os f^!\mG)\\
=\Hom_{\shv(X)}(f^!((\verd^\eta_Y)^c((\verd_Y^\eta)^c(\mA)\os \mG)),\mF).
\end{multline*}
One checks that the map in \eqref{eq-exotic-projection} is compatible with these two isomorphisms and therefore is an isomorphism.

The isomorphism \eqref{eq:sheaves-dual-ula-base-change-2} follows directly from the second isomorphism in \Cref{lem:etale-proper-functoriality-shv} \eqref{lem:etale-proper-functoriality-shv-6} by passing to the right adjoint. To prove \eqref{eq:sheaves-dual-ula-base-change-1}, we first assume that $g$ is pfp, in which case the desired isomorphism follows from the first isomorphism in \Cref{lem:etale-proper-functoriality-shv} \eqref{lem:etale-proper-functoriality-shv-6} by passing to the right adjoint. For general $g$, we may write $Y'\to Y$ as $Y'=\lim_i Y_i$ with $g_i:Y_i\to Y$ pfp (so $Y_i$ is standard placid). Write $h_i: Y'\to Y_i$.
Let $f_i: X_i\to Y_i$ denote the corresponding base change of $f$ and $g_i': X_i\to X$ the base change of $g_i$. Then as reasoning in \Cref{rmk:sheaves-qcqs-inverse-limit-constructible}, we have
\[
(f')_\flat((g')^! \mF)\cong \colim_i ((h_i)^!((f'_i)_\flat((g'_i)^!\mF)))\cong \colim_i ((h_i)^!((g_i)^! ( f_\flat \mF))) =g^!( f_\flat \mF),
\]
giving the desired isomorphism.
\end{proof}

\begin{remark}\label{rem: duality and tensor product}
In the proof of \Cref{lem:dual-sheaves-proj-formula-weakly-coh-pro-smooth}, we considered the binary operation
\[
\cshv(Y)\otimes \cshv(Y)\to \cshv(Y),\quad (\mF,\mG)\mapsto (\verd_Y^\eta)^c(((\verd_Y^\eta)^c(\mF))\os ((\verd_Y^\eta)^c(\mG))),
\]
which in fact defines a monoidal structure on $\cshv(Y)$, with unit being $\La_Y^\eta$. We shall denote such monoidal structure by $\otimes^\eta$. Note that under the canonical equivalence \eqref{eq: identifying star sheaf with shrek sheaf-1}, this is identified with the usual $*$-monoidal structure $\scshv(X)$.
\end{remark}

Now we let $\mathrm{Eproet}$  denote the class of essentially pro-\'etale morphisms, i.e. those $f: X\to Y$ that can be written as $f: X\to X'\to Y$ with $X\to X'$ pro-\'etale  and $X'\to Y$ pfp. 
\begin{proposition}\label{prop-extension-of-shv-pro-etale} 
The class of ess. pro-\'etale morphisms is strongly stable.
The sheaf theory $\shv$ from \eqref{category-shv-definition} admits an extension
\begin{equation}\label{prop-extension-of-shv-pro-etale-1} 
\shv\colon \corr(\qcsp_k)_{\mathrm{Eproet};\all} \to \lincat_\La
\end{equation}
such that for  pro-\'etale $f: X\to Y$, $f_*=f_\flat$. In addition, if $f: X\to Y$ is an ess. weakly pro-\'etale morphism between standard placid spaces, then $f_*$ admits a left adjoint $f^*$.
\end{proposition}
\begin{proof}We note that the classes of pfp morphisms and pro-\'etale morphisms are strongly stable, and the class of ess. pro-\'etale morphisms is weakly stable (as in the proof of \Cref{composition.coh.pro-smooth}). 
Then \Cref{lem:dual-sheaves-proj-formula-coh-pro-smooth} allow us to apply \Cref{prop-sheaf-theory-for-adjoint-factorization} (and \Cref{rem-sheaf-theory-for-adjoint-factorization})  to obtain the desired extension \eqref{prop-extension-of-shv-pro-etale-1} by letting $\corr(\bfC)_{\verti;\horiz}=\corr(\qcsp_k)_{\pfp;\all}$ and by letting  $\mathrm{E}=\mathrm{HR}$ be the class of pro-\'etale morphisms. 
The last statement is clear.
\end{proof}

\begin{remark}\label{rem:dual-sheaves-proj-formula-coh-pro-smooth}
Clearly, \Cref{lem:dual-sheaves-proj-formula-coh-pro-smooth} continues to hold for the above extended sheaf theory. That is,
pseudo coh. pro-smooth morphisms still satisfy \Cref{assumptions.base.change.sheaf.theory.H} for the sheaf theory \eqref{prop-extension-of-shv-pro-etale-1} .
\end{remark}

\begin{remark}\label{rem: variant of shv}
Note that $\flat$-pushforwards along general coh. pro-smooth morphisms cannot be absorbed into the above sheaf theory. The problem is the class $\mathrm{Prosm}$ of coh. pro-smooth morphisms is not strongly stable so \Cref{prop-sheaf-theory-for-adjoint-factorization} is not applicable. However, this class is still weakly stable. Given \Cref{lem:dual-sheaves-proj-formula-coh-pro-smooth}, we can apply
\Cref{ex-sheaf-theory-for-adjoint-factorization} \eqref{ex-sheaf-theory-for-adjoint-factorization-2} to obtain a variant of \eqref{prop-extension-of-shv-pro-etale-1} 
\[
\shv': \corr(\qcsp_k)_{\mathrm{Prosm};\all} \to \lincat_\La.
\]
which still sends $g: Z\to Y$ to $g^!: \shv(Y)\to \shv(Z)$ but sends $(f:Z\to X)\in \mathrm{Prosm}$ to $f_\flat: \shv(Z)\to \shv(X)$.
When restricted to $\mathrm{Sm}\subset\mathrm{Prosm}$, the two theories $\shv$ and $\shv'$ are essentially equivalent as $f_\flat$ and $f_*$ only deters by a shift (and a twist).
\end{remark}

\subsubsection{Categories of cosheaves on prestacks}

For several purposes, it is convenient to have a definition of the category of (co)sheaves on a general prestack over $k$. 
From now on we assume that $k$ is the perfection of a regular noetherian ring of dimension $\leq 1$ 
and $\ell$ a prime invertible in $k$ such that $k$ has finite $\bF_\ell$-cohomological dimension.  We allow $\La$ to be $\bZ_\ell$-algebras as in \Cref{sec:adic-formalism}.

Let $\mathrm{E}\subset \mathrm{Mor}(\qcsp_k)$ be a class of morphisms in $\qcsp_k$. Recall that if $\mathrm{E}$ is stable under base change, then it extends naturally to a class of morphisms $\mathrm{E}_r$ in $\prestk_k$ consisting of those morphisms $f: X\to Y$ of prestacks that are representable in $\mathrm{E}$. If $\mathrm{E}$ is weakly (resp. strongly) stable, so is $\mathrm{E}_r$. See \Cref{rem-prop-of-weak-strong-stable-class} \eqref{rem-prop-of-weak-strong-stable-class-2}. For example, we may talk representable pfp, pfp proper, coh. smooth, (strongly/weakly/ess.) coh. pro-smooth morphisms, (ess.) pro-\'etale, (ess.) coh. pro-unipotent morphisms between prestacks. 

\begin{remark}\label{rem: qcqs representable}
Note that by definition, for a representable morphism $f: X\to Y$ of prestacks and a morphism $S\to Y$ with $S\in\qcsp_k$ then $S\times_YX$ is qcqs. So our notion of representable morphisms between (pre)stacks is slightly more restrictive than the usual notion of representable morphisms.
\end{remark}

We apply \Cref{prop-sheaf-theory-right-Kan-extension} to define 
\begin{equation}\label{def-dual-sheaves-on-prestacks}
\shv(-,\La): \corr(\prestk_k)_{\mathrm{Eproet}_r;\all}\to \lincat_\La  \quad \mbox{resp. } \cshv(-,\La): \corr(\prestk_k)_{\pfp_r;\all}\to \catid_\La, 
\end{equation}
as the right Kan extension of the functor from \eqref{prop-extension-of-shv-pro-etale-1}, resp. \eqref{category-cshv-definition}, along the full embedding 
\[
\corr(\qcsp_k)_{\mathrm{Eproet};\all}\subset \corr(\prestk_k)_{\mathrm{Eproet}_r;\all}, \quad \mbox{resp.} \quad \corr(\qcsp_k)_{\pfp;\all}\subset \corr(\prestk_k)_{\pfp_r;\all}.
\]

By \Cref{prop-sheaf-theory-right-Kan-extension},  we have
\begin{equation}\label{eq: shv for prestack as limit over qcqs space}
\shv_{(c)}(X,\La)\xrightarrow{\sim} \lim_{S\rightarrow X} \shv_{(c)}(S,\La),
\end{equation}
with $S \in ({\qcsp_k}_{/X})^{\op}$ and transition maps given by $!$-pullbacks.
Informally, this means giving an object $\mF\in \shv(X,\La)$ (resp. $\mF\in \cshv(X,\La)$) amounts to giving for every $S\to X$ with $S\in \qcsp_k$ an object $\mF_S\in\shv(S,\La) $  (resp. $\mF_S\in \cshv(S,\La)$) and to giving for every $g: S'\to S$ an isomorphism $g^!\mF_S\cong \mF_{S'}$ satisfying all (higher) compatibility conditions in a coherent way.  Note that by \cite[Proposition 6.2.1.9]{Lurie.SAG}, for any presentation of a prestack $X$ as a colimit $X\simeq \colim_{\alpha\in A} X_\alpha$ of prestacks, we have an equivalence
\begin{equation}\label{eq:Shv-send-colim-to-lim}
   \shv(X,\La) \xrightarrow{\sim} \lim_{\alpha\in A^{\op}} \shv(X_\alpha,\La).
\end{equation}

Note that $\cshv(X,\La)$ is a full subcategory of $\shv(X,\La)$, but is in general no longer the subcategory of compact objects. In addition, in general $\shv(X)$ may not be compactly generated.
\begin{remark}\label{rem:standard t-structure on prestack shv}
Assume that $\La$ is regular noetherian. It follows from \Cref{rem:standard t-structure on shv} that there is a standard $t$-structure on $\shv_{(c)}(X,\La)$ such that $\shv_{(c)}(-,\La)^{\mathrm{std},\leq 0}$ is the right Kan extension of $\shv_{(c)}(-,\La)^{\mathrm{std},\leq 0}$ from $\qcsp_k$ to $\prestk_k$. Then the $!$-pullback functors are $t$-exact, and the inclusion $\cshv(X,\La)\subset \shv(X,\La)$ is $t$-exact. The standard $t$-structure on $\cshv(X,\La)$ is bounded and on $\shv(X,\La)$ is accessible, compatible with filtered colimits, and right complete.
\end{remark}

\begin{example}\label{ex-dualizing sheaf}
For each prestack $X$, there is an object 
\[
\consdual_X\in\cshv(X)\subset\shv(X),
\] 
whose $!$-pullback to every $S\in\qcsp_k$ is $\consdual_S$. This is in fact the unit of the symmetric monoidal structure on $\shv(X)$.
Note that $\consdual_X$ is a discrete object in $\shv(X)$. I.e. $\map(\consdual_X,\consdual_X)$ is a discrete space, or equivalently $\mathrm{Ext}^i(\consdual_X,\consdual_X)=0$ for $i<0$. Indeed, this is clear if $X$ is a pfp algebraic space over $k$, and then holds for $X\in\qcsp_k$ since if one writes $X=\lim X_i$ as a cofiltered limit of pfp schemes with affine transition maps, then $\map(\consdual_X,\consdual_X)=\colim_i \map(\consdual_{X_i},\consdual_{X_i})$ is discrete (as filtered a colimit of discrete spaces is discrete). Finally, for any prestack $X$, $\map(\consdual_X,\consdual_X)=\lim_{S\to X}\map(\consdual_S,\consdual_S)$ is again discrete, as arbitrary limit of discrete spaces is discrete.
\end{example}

Clearly  \Cref{lem:etale-proper-functoriality-shv} \eqref{lem:etale-proper-functoriality-shv-1}-\eqref{lem:etale-proper-functoriality-shv-5} hold for prestacks. It follows from \Cref{prop-sheaf-theory-right-Kan-extension-VR class} that and \Cref{lem:dual-sheaves-proj-formula-coh-pro-smooth} (\Cref{lem:dual-sheaves-proj-formula-coh-pro-smooth}) also holds for prestacks. We record them in the following statements.

\begin{proposition}\label{lem:etale-proper-functoriality-shv-prestack}
Let $f\colon X\rightarrow Y$ be a morphism of prestacks.
\begin{enumerate}
    \item\label{lem:etale-proper-functoriality-shv-prestack-1} If $f$ is representable and \'{e}tale, $f_*$ is a right adjoint to $f^!$, and if $f$ is representable pfp proper, $f_*$ is a left adjoint to $f^!$. In either of the above situation, the base change isomorphism \eqref{eq:abstract-base-change} encoded by the functor $\shv$ is the Beck-Chevalley map obtained by the adjoint as in \Cref{def:categories-adjointability}.

     \item\label{lem:etale-proper-functoriality-shv-prestack-3} If $f$ is a representable $\ell$-ULA morphism, then $f_*$ admits a left adjoint $f^*$, which preserves constructibility. In addition, if \eqref{eq:pullback-square-of-algebraic-spaces} is a pullback square of prestacks (with $f$ being representable $\ell$-ULA), then we have the base change isomorphism $(f')^*\circ g^!\xrightarrow{\simeq} (g')^!\circ f^*$.

     \item\label{lem:etale-proper-functoriality-shv-prestack-4} Representable pseudo coh. pro-smooth morphisms satisfy  \Cref{assumptions.base.change.sheaf.theory.H}.
\end{enumerate}
\end{proposition}

\Cref{lem:etale-proper-functoriality-shv-prestack} \eqref{lem:etale-proper-functoriality-shv-prestack-1} allows one to apply \eqref{eq:limit-colimit equivalence} to \eqref{eq:Shv-send-colim-to-lim} to give a colimit presentation of $\shv(X)$ is an important case.
\begin{corollary}\label{colimit.sheaves.pseudo-proper.stacks}
Let $X=\colim_{\al} X_\al$ be a (filtered) colimit of prestacks with $X_\al\to X_{\al'}$ representable pfp proper morphisms. Then
\[
\colim_{\al\in A}\shv(X_\al,\La)\to \shv(X,\La)
\]
is an equivalence, where the transition maps are given by $*$-pushforwards.  
\end{corollary}

\begin{remark}\label{shv over sshv}
Of course, one can define another sheaf theory $\sshv$ and its constructible version $\scshv$ for prestacks by right Kan extension of \eqref{eq-*-pullback-ind-constructible} and  \eqref{eq-six-operation-*-constructible} along $(\qcsp_k)^{\op}\subset (\prestk_k)^{\op}$. 
By definition, is a canonical equivalence
\begin{equation}\label{eq: star and shrek constr sheaves on quasi-placid}
\cshv(X)\cong \scshv(X)^{\op},
\end{equation}
but $\shv(X)$ and $\sshv(X)$ are in general unrelated.
The theory $\sshv$ has its own applications.
But \Cref{colimit.sheaves.pseudo-proper.stacks} is the main reason we would like to work with the sheaf theory $\shv$ in this work.
\end{remark}

We record the following statements for further references.
\begin{lemma}\label{lem:sheaves-open-closed-prestacks}
Let $j\colon U\rightarrow X$ be a quasi-compact open embedding of prestacks with a pfp closed complement $i\colon Z\rightarrow X$. Then:
\begin{enumerate}
    \item $i^!\circ j_*\simeq 0$ and $j^!\circ i_{*}\simeq 0$. 
    \item The functors $i_*$ (resp. $j_*$) are fully faithful, with essential image consisting of $\mF\in \shv(X)$ with $j^!\mF\simeq 0$ (resp. $i^!\mF\simeq 0$). 
    \item For every $\mF\in \shv(X,\La)$ we have a canonical fiber sequence 
    \[
    i_*  i^! \mF \rightarrow \mF \rightarrow j_*  j^!\mF,
    \]
    given by the counit of the adjunction $(i_*,i^!)$ and the unit of the adjunction $(j^!,j_*)$.
\end{enumerate}
\end{lemma}
\begin{proof}
Using base change, the statement can be proved after pullback along every $S\to X$ with $S$ varying over perfect qcqs algebraic spaces over $k$. But this is standard.
\end{proof}

We will also need the following K\"unneth type formula.

\begin{proposition}\label{lem: categorical kunneth prestack}
Assume that $k$ is an algebraically closed field.
Let $X$ be a perfect prestack and assume that $\shv(X,\La)$ is dualizable. Then for every $Y$, the exterior tensor product 
\[
\boxtimes: \shv(X,\La)\otimes_\La\shv(Y,\La)\to \shv(X\times Y,\La)
\]
is fully faithful.
\end{proposition}
\begin{proof}
First notice that \Cref{lem: categorical kunneth} continuous holds for any coefficient $\La$ being $\bZ_\ell$-algebras as in \Cref{sec:adic-formalism} and $X,Y\in\qcsp_k$. Then passing to the opposite categories and taking ind-completion, we see that the lemma holds when $X,Y\in\qcsp_k$. 

Now we argue as in \cite[Proposition 3.3.1.7]{Gaitsgory.Rozenblyum.DAG.vol.I}.
Suppose $X\in \qcsp_k$, then as $\shv(X,\La)$ is dualizable, $\shv(X,\La)\otimes_\La-$ commutes with limits. Note that for every prestack $Y$, the functor  $(\qcsp_k)_{/Y}\to (\qcsp_k)_{/X\times Y},\ U\mapsto X\times U$ is cofinal, then follows that
\begin{multline*}
\shv(X,\La)\otimes_\La \shv(Y) =\shv(X,\La)\otimes_\La \lim_{((\qcsp_k)_{/Y})^{\op}}\shv(U,\La)\\
\to \lim_{((\qcsp_k)_{/Y})^{\op}}\shv(X\times U)\cong \shv(X\times Y)
\end{multline*}
is fully faithful. Finally if $\shv(X,\La)$ is dualizable, one can run the above argument again to conclude.
\end{proof}

The above functor is in general not an equivalence. However, it is an equivalence in some important cases. See \Cref{cor: ess. pro-unip tensor product} and \Cref{cor: profinite tensor product}.

\subsubsection{Ind-$\mathrm{E}$ morphisms} 

\begin{definition}\label{def-ind-scheme}
A perfect prestack $X$ is called an \textit{ind-scheme} (resp. \textit{ind-algebraic space}) if it can be written as a filtered colimit of $X = \colim_{i} X_i$ of $X_i \in \psch_k$ (resp. $X_i \in \qcsp_k$) with transition maps given by pfp closed immersions. Let $\ind\psch_k\subset\ind\qcsp_k\subset \prestk^{\pf}_k$ denote the full subcategory of ind-schemes and ind-algebraic spaces over $k$.
\end{definition}

Our definition of ind-schemes/algebraic spaces is not the most general one. In literature, sometimes ind-schemes are defined as a filtered colimit of $X = \colim_{i} X_i$ with transition maps being closed embeddings as above, but without requiring $X_i$ to be qs nor requiring $X_i\to X_j$ to be pfp.
However, the above definition is general enough to our purpose.

\begin{definition}\label{def:ind-fp.and.proper.morphism}
Let $\mathrm{E}$ be a class of morphisms in $\qcsp_k$.
A morphism of prestacks $f\colon X\rightarrow Y$ is called ind-$\mathrm{E}$  if for every map $S\rightarrow Y$ with $S\in\psch_k$ the pullback $f_S\colon X_S \rightarrow S$ admits a presentation as a filtered colimit $X_S = \colim_{i\in \mI} X_i$ with transition maps $X_i \rightarrow X_j$ being (pfp) closed immersions  of algebraic spaces and with each $f_i \colon X_i\rightarrow S$ belonging to $\mathrm{E}$. We let $\mathrm{IndE}$ denote the class of ind-$\mathrm{E}$ morphisms between prestacks.
\end{definition}

Note that this definition in particular applies to $\mathrm{E}$ to be the class of pfp morphisms, pfp proper morphisms, and ess. pro-\'etale morphisms respectively. Sometimes, we will call them ind-pfp, ind-pfp proper, and ind-ess. pro-\'etale morphisms, respectively.

\begin{remark}\label{ind.fp.for.ind.spaces}
\begin{enumerate}
\item\label{ind.fp.for.ind.spaces-1} Note that a map of ind-algebraic spaces $f\colon X\rightarrow Y$ is ind-finitely presented (resp. ind-pfp proper, resp. ind-ess. pro-\'etale) if and only if for every finitely presented closed spaces $X'\subseteq X$, $Y'\subseteq Y$ such that $X'\rightarrow Y$ factors through $Y'$ the map $X'\rightarrow Y'$ is pfp (resp. pfp proper, resp. ess. pro-\'etale). 
It follows that these classes between ind-algebraic spaces are strongly stable.
\item\label{ind.fp.for.ind.spaces-2} In general it is not possible to write an ind-$\mathrm{E}$ morphism $f:X\to Y$ between perfect prestacks as a filtered colimit $\colim_i X_i\to Y$ with each $X_i\to Y$ in $\mathrm{E}_r$ (i.e. representable in $\mathrm{E}$) and transition maps being pfp closed embeddings. See, however, \Cref{lem-indfp-to-placid-stacks}.
\end{enumerate}
\end{remark}

The following lemma says that ind-pfp is a property local for the \'etale topology. But for this being true, it is important to allow algebraic spaces (rather than merely schemes) in our definition.
\begin{lemma}\label{lem-indfp-etale-local} 
Let $f\colon X \to Y$ be a morphism of \'etale stacks with $Y\in\psch_k$. If there is an \'etale cover
$Y'\to Y$ such that the base change $X'\to Y'$ can be written as $X' = \colim_i X'_i$ where $X'_i\in\qcsp_k$ is pfp (resp. pfp proper) over $Y'$, then $f$ is ind-pfp (resp. ind-pfp proper).
\end{lemma}
\begin{proof} This directly follows from \cite[Lemma 3.12]{Haines.Richarz.Weil.restriction}. See also \Cref{lem-indfp-to-placid-stacks} for an argument in a more complicated situation.
\end{proof}

The following lemma follows directly from \Cref{lem:etale-proper-functoriality-shv} \eqref{lem:etale-proper-functoriality-shv-1} \eqref{lem:etale-proper-functoriality-shv-6} and \eqref{pro-unipotent.are.admissible}, together with \Cref{colimit.sheaves.pseudo-proper.stacks}.

\begin{lemma}\label{lem: lower-shrek for ind-ess. coh. pro-unipotent}
Let $f\colon X\rightarrow Y$ be an ind-ess. coh. pro-unipotent morphism with $Y\in \qcsp_k$. Then $f^!$ admits a left adjoint. In addition, for a pullback square as in \eqref{eq:pullback-square-of-algebraic-spaces} with $g\colon Y'\rightarrow Y$ being weakly coh. pro-smooth, then there is natural base change isomorphism $(f')_!\circ (g')^!\rightarrow g^!\circ f_!$.
\end{lemma}

Now we extend the sheaf theory $\shv$ from \eqref{def-dual-sheaves-on-prestacks} to allow $*$-pushforwards along a large class of morphisms. 

\begin{proposition} \label{def-dual-sheaves-on-prestacks-vert-indfp}
The functor from \eqref{def-dual-sheaves-on-prestacks} admits a canonical extension to a functor
\begin{equation}\label{eq-six-operation-prestack}
\shv(-,\La)\colon \corr(\prestk^\pf_k)_{\mathrm{IndEproet};\all}\rightarrow \lincat_\La. 
\end{equation}
\end{proposition}
\begin{proof}

We first apply  \Cref{prop-sheaf-theory-non-full-right-Kan-extension} to $\shv(-,\La)\colon \corr(\qcsp_k)_{\mathrm{Eproet};\all} \rightarrow \lincat_\La$ to obtain an extension
\begin{equation}\label{eq-six-operation-indsp}
\shv(-,\La)\colon \corr(\indsp^\pf_k)_{\mathrm{IndEproet};\all} \rightarrow \lincat_\La.
\end{equation}

To do so, we let $\verti_1$ be the be class $\mathrm{Eproet}$, $\verti_2$ the class $\mathrm{IndEproet}$, $\mathrm{S}_1$ the class of pfp closed embeddings.

Then we let \eqref{eq-six-operation-prestack} to be the right Kan extension of \eqref{eq-six-operation-indsp} along the full embedding 
\[
\corr(\indsp^\pf_k)_{\mathrm{IndEproet};\all}\subset  \corr(\prestk^\pf_k)_{\mathrm{IndEproet};\all}.
\] 
As before, by
\Cref{prop-sheaf-theory-right-Kan-extension},  its restriction to $(\prestk^\pf_k)^{\op}$ is just $\shv$. 
\end{proof}

\begin{remark}\label{rem-*-pushforward-along-indfp-morphisms}
Informally, let $X\xleftarrow{f}Z\xrightarrow{g} Y$ be a correspondence of ind-schemes with $f$ belonging to $\mathrm{IndEproet}$. Suppose we write $Z$ as $Z=\colim_\al Z_\al$ and let $f_\al\colon Z_\al\to X, g_\al\colon Z_\al\to Y$ be composed morphisms so each $f_\al$ is ess. pro-\'etale. Then for $\mF\in\shv(Y)$, 
\[
f_*(g^!\mF)= \colim_i (f_\al)_*((g_\al)^!\mF),
\] 
where the transition maps come from the co-unit adjunction $((\iota_{\al,\beta})_*,(\iota_{\al,\beta})^!)$ for pfp closed immersion $\iota_{\al,\beta}\colon Z_\al\subset Z_\beta$.
\end{remark}

\begin{example}\label{ex: Borel-Moore homology}
Suppose $k$ is an algebraically closed field.
Let $X=\colim X_i$ be an ind-scheme, with each $X_i$ pfp over $k$. We write $\pi_{X_i}: X_i\to \spec k$ and $\pi_X: X\to \spec k$ for the structural maps. Then
\[
(\pi_X)_*\consdual_X=\colim_i (\pi_{X_i})_*\consdual_{X_i}=:\colim_i C^{\mathrm{BM}}_\bullet(X_i,\La)=:C^{\mathrm{BM}}_\bullet(X,\La)
\]
is the usual Borel-Moore homology of $X$.
\end{example}


The following statement follows from the construction (from \Cref{prop-sheaf-theory-non-full-right-Kan-extension}).
\begin{lemma}
\label{locally.ind-proper.existence.lower-!}
Let $f\colon X\rightarrow Y$ be an ind-pfp proper morphism of prestacks. Then $f_*$ is the left adjoint of $f^!$.
\end{lemma}

We mention the following base change result.

\begin{lemma}\label{lem-indproper-prosmooth-push-basechange}
Suppose \eqref{eq:pullback-square-of-algebraic-spaces} is a Cartesian diagram of prestacks. If $f$ is ind-ess. pro-\'etale and $g$ is representable pseudo coh. pro-smooth. Then there is a natural isomorphism of functors $f_*\circ (g')_\flat\to g_\flat\circ (f')_*$.
\end{lemma}
\begin{proof}We may assume that $Y\in\qcsp_k$, and then assume that $f$ is ess. pro-\'etale, which then follows from \Cref{lem:etale-proper-functoriality-shv-prestack} \eqref{lem:etale-proper-functoriality-shv-prestack-4}.
\end{proof}

Recall that associated to a sheaf theory, we have the class of morphisms $\mathrm{HR}$ and $\mathrm{VR}$ as in \Cref{rem-additional-base.change.sheaf.theory} \eqref{rem-additional-base.change.sheaf.theory-1}. 
\begin{corollary}\label{cor-additional-base-change-for-shv-theory-on-prestacks}
The class of representable pseudo coh. pro-smooth morphisms belong to $\mathrm{HR}$, and the class of ind-pfp proper morphisms belong to $\mathrm{VR}$. 
\end{corollary}
\begin{proof}
\Cref{lem-indproper-prosmooth-push-basechange} and  \Cref{lem:etale-proper-functoriality-shv-prestack} \eqref{lem:etale-proper-functoriality-shv-prestack-4} imply that representable pseudo coh. pro-smooth morphisms belong to $\mathrm{HR}$.
As mentioned in \Cref{rem-additional-base.change.sheaf.theory} \eqref{rem-additional-base.change.sheaf.theory-2}, \Cref{locally.ind-proper.existence.lower-!} implies that ind-pfp proper morphisms satisfy \Cref{assumptions.base.change.sheaf.theory.V}. 
\end{proof}

By applying \Cref{lem-sheaf-theory-left-Kan-extension-covering} and \Cref{lem-sheaf-theory-left-right-Kan-extension}, we may further extend \eqref{eq-six-operation-prestack} by allowing pushforward along certain morphisms that are not (ind-)representable. Namely, we inductively define a class $\verti_r$ of morphisms between prestacks as follows. Let $\verti_0=\mathrm{IndEproet}$. Suppose we have $\verti_r$ and an extension of $\shv$ to $\corr(\prestk^\pf_k)_{\verti_r;\all}$. Then let $\verti_{r+1}$ be the class of morphisms constructed from $\verti_r$ as in \Cref{lem-sheaf-theory-left-Kan-extension-covering}. We have an extension of $\shv$ to $\corr(\prestk^\pf_k)_{\verti_{r+1};\all}$. Finally, let $\verti_\infty$ be the union of all $\verti_r$s and let $\verti$ be the class constructed from $\verti_\infty$ as in \Cref{lem-sheaf-theory-left-right-Kan-extension}. Then we have an extension of the sheaf theory 
\begin{equation}\label{eq: pushforward along non-representable morphisms}
\shv\colon\corr(\prestk^\pf_k)_{\verti;\all}\to \lincat_\La.
\end{equation}

\begin{example}\label{rem: pushforward along non-representable morphisms}
For a concrete example of morphisms contained in $\verti$, we note that a morphism $f: X\to Y$ belongs to $\verti$ if there is an \'etale covering $Y'\to Y$ of $Y$ such that $X\times_YY'\to Y'$ is ind-ess. pro-\'etale.

For another example, let $A$ be a finite group, regarded as a constant affine algebraic group over $k$. Then the non-representable morphism $\bB A\to \spec k$ belongs to the class $\verti_1$ as in \Cref{lem-sheaf-theory-left-Kan-extension-covering}. Then if $f: X\to Y$ is an $A$-gerbe over $Y$ (i.e. \'etale locally on $Y$, $X\simeq Y\times \bB A$), then $f\in \verti$. We caution, however, that the pushfoward along $\bB A\to \spec k$ (and therefore along any $A$-gerbe map) encoded in \eqref{eq: pushforward along non-representable morphisms} is the left adjoint of the $!$-pullback along $\bB A\to \spec k$.
Only when the order of $A$ is invertible in $\La$, it is also the right adjoint of the $!$-pullback. We will only use \eqref{eq: pushforward along non-representable morphisms} for $A$-gerbe pushforwards when the order of $A$ is invertible in $\La$. In general, $f_*$ is a ``renormalized" version of the naive right adjoint of the $!$-pullback (which is not continuous).
\end{example}

\subsubsection{Morphisms of universal homological descent}
Now we discuss descent for the sheaf theory $\shv$. Recall that we assume that $k$ has finite $\bF_\ell$-cohomological dimension.

\begin{definition}\label{def:sheaves-prestacks-universal-hom-descent}
A morphism $f\colon X\rightarrow Y$ of perfect prestacks
is said to be of \textit{homological descent} if $f$ is $\shv$-descent in the sense of \Cref{def:der-descent}. I.e.
the  canonical map
\[
\shv(Y,\La) \rightarrow \tot \left(\shv(X_\bullet,\La)\right)
\]
induced by $!$-pullbacks is an equivalence, where $X_\bullet\to Y$ denotes the \v{Cech} nerve of $f$.
It is  said to be of \textit{universal homological descent} if its base change along every morphism $Y'\rightarrow Y$  is of homological descent.
\end{definition}

\begin{remark}\label{lem:sheaves-descent-prestacks-univ-homological-descent}
We use the term "homological" instead of the usual "cohomological" since we are using $!$-pullback functors and the dual category of sheaves. By \eqref{eq: shv for prestack as limit over qcqs space} and \eqref{eq:Shv-send-colim-to-lim}, it is clear that a morphism $f\colon X\rightarrow Y$ of prestacks which is of universal homological descent if only only if its base change along every $S\to Y$ with $S\in\psch_k$ is of homological descent.
\end{remark}

The goal of the rest of this subsection is to exhibit a few classes of morphisms are of universal homological descent.

\begin{proposition}\label{locally.ind-proper.descent}
Let $f\colon X\rightarrow Y$ be an ind-pfp proper and surjective morphism of perfect prestacks. Then $f$ is of universal homological descent.
\end{proposition}
\begin{proof}
\quash{Consider the augmented cosimplicial object $\Delta_{+} \rightarrow \lincat_\La$ given by applying $\shv$ (with horizontal morphisms) $\shv(Y,\La)\to \shv(X_\bullet,\La)$ of $\lincat_\La$. Below, we denote $X_{-1}:= Y$. As all functors in $\lincat_\La$ preserve small colimits, by \Cref{cor:Beck-Chevalley-fully-faitful}, it is enough to check the following conditions:
\begin{enumerate}
    \item For every morphism $\alpha\colon [m]\rightarrow [n]$ in $\Delta_{+}$ we have a pullback diagram
    \[
    \begin{tikzcd}
    X_{n+1} \arrow[r,"f_{\alpha+1}"]\arrow[d,"d^n_0"] &
    X_{m+1} \arrow[d,"d^m_0"]\\
    X_n \arrow[r,"f_\alpha"] & X_m
    \end{tikzcd}
    \]
    the base change map $ (d_0^m)^{*}\circ (f_\alpha)_* \rightarrow (f_{\alpha+1})_{*}\circ (d_0^n)^{*}$ is an equivalence.
    \item The $!$-pullback functor $\shv(Y) \rightarrow \shv(X)$ is conservative. 
\end{enumerate}
Condition (1) follows immediately from locally ind-proper base change, \Cref{ind-proper.bc.prestacks}, as any boundary map in the \v{e}ch nerve is locally ind-proper.  }
By \Cref{prop-descent-codescent-abstract} \eqref{prop-descent-codescent-abstract-1},  it is enough to show that for any qcqs $S$ and a map $S\rightarrow Y$ the functor $f_S^!$ is conservative. We can assume that $X\times_{Y} S$ admits a presentation as a filtered colimit
$X\times_{Y} S = \colim_{i\in \mI} X_{S,i}$
of pfp proper maps $f_{i}\colon X_{S,i} \rightarrow S$ and $X_{S,i}\to X_{S,i'}$ closed immersion. Let $S_i$ denote the image of $f_i(X_{S,i})$ regarded as a (perfectly) finitely presented closed subscheme of $S$ (since we are dealing with perfect schemes there is a unique induced scheme structure on $S_i$). Then the surjectivity of $f$ implies that $S = \cup_{i} S_i$. By \cite[\href{https://stacks.math.columbia.edu/tag/094L}{Lemma 094L}]{stacks-project} the topological space $S$ is spectral, and therefore by \cite[\href{https://stacks.math.columbia.edu/tag/0901}{Lemma 0901}]{stacks-project} we have that $S_{\cons}$ compact, where $S_{\cons}$ is the topological space associated to the scheme $S$ endowed with the constructible topology. The open subsets $(S_i)_{\cons}\subseteq S_\cons$ constitute an open cover of $S_\cons$ and since $\mI$ is filtered $S_i = S$ for some $i$. That is, there exists some $i\in \mI$ such that $f_{i}\colon X_{S,i} \rightarrow S$ is surjective. Since $f_i$ is pfp proper, the functor 
$\shv(S,\La) \rightarrow \shv(X_{S,i},\La)$, and therefore the functor $\shv(S) \rightarrow \shv(X\times_{Y} S)$, is conservative (by \Cref{prop-h-descent-shv} \eqref{prop-h-descent-shv-2}).
\end{proof}

\begin{remark}\label{rem-ind-proper-submersive}
The above argument also shows that if $f\colon X\to Y$ is an ind-pfp proper surjective morphism,
then it is universally submersive. In particular, the map $|X|\to |Y|$ is a quotient map.
\end{remark}

\begin{proposition}\label{pro-unipotent.morphism.descent}   
Let $f\colon X\rightarrow Y$ be a representable ess. coh. pro-unipotent morphism. Then $f$ is of universal homological descent. 
\end{proposition}
\begin{proof}
It is enough to assume that $f$ is coh. pro-unipotent between perfect qcqs algebraic spaces and show that it is of homological descent. By \Cref{rem-fully-faithful-pro-unipotent}, $f^!\colon \shv(Y)\to \shv(X)$ is fully faithful. So the unit map $\id \rightarrow f_\flat f^!$ is an equivalence, and the claim follows from \Cref{prop-descent-codescent-abstract} \eqref{prop-descent-codescent-abstract-2}.
\end{proof}

Next, we need to discuss enough interesting cases of (perfect) affine flat group schemes $H$ over $k$ for which $\spec k\to \bB_{\mathrm{fpqc}} H$ is of universal homological descent, or equivalently
every $H$-torsor (in $fpqc$ topology) $E\rightarrow S$ with $S\in \psch_k$ is homological descent. (Note that $\spec k\to \bB H$ is of universal homological descent, by \Cref{Liu.Zheng.descent.lemma} \eqref{descent.detection} and \Cref{prop-h-descent-shv} \eqref{prop-h-descent-shv-2}. But this is not enough for our purpose.) 

First, we notice that we can reduce this question to any normal subgroup of ``finitely presented index". That is, suppose $H$ admits a short exact sequence of perfect affine flat group schemes
\begin{equation}\label{eq-exact-sequence-flat-group}
1\to H_0\to H\to H'\to 1
\end{equation}
with $H'$ pfp flat over $k$. Then $\spec k\to \bB_{\mathrm{fpqc}} H$ is of universal homological descent if so is $\spec k \to \bB_{\mathrm{fpqc}} H_0$. Indeed, every $H'$-torsor $E\to S$ is an $h$-cover and therefore is of universal homological descent. Then we can utilize \Cref{prop-descent-codescent-abstract} \eqref{prop-descent-codescent-abstract-2}. Now, if $H_0$ in \eqref{eq-exact-sequence-flat-group} is coh. pro-unipotent over $k$, then $\spec k\to \bB_{\mathrm{fpqc}} H$ is of universal homological descent by \Cref{pro-unipotent.morphism.descent}. 
It follows that for such $H$, the classifying stack $\bB H$ in \Cref{cor: ess. pro-unip tensor product} can be replaced by $\bB_{\mathrm{fpqc}} H$.
Of course, instead of considering classifying stack in fpqc topology, the universal $H$-torsor $\spec k \to \bB H$ in \'etale topology is of universal homological descent for any $H$.

\begin{proposition}\label{cor: ess. pro-unip tensor product}
Suppose $k$ is an algebraically closed field.
Let $H$ be an affine group scheme as in \eqref{eq-exact-sequence-flat-group} with $H'$ pfp and $H_0$ coh. pro-unipotent over $k$. Then for every prestack $X$ over $k$, the exterior tensor product
\[
\shv(\bB H,\La)\otimes_\La \shv(X,\La)\to \shv(\bB H\times X,\La)
\]
is an equivalence. The same statement holds with the \'etale quotient replaced by fpqc quotient.
\end{proposition}
\begin{proof}
As usual, we write $\pt$ for $\spec k$, and let $f:\pt\to \bB H$ denote the map of universal $H$-torsor.

By \Cref{compact.generation.admissible.stacks} below, $\shv(\bB H)$ is compactly generated. Then the argument as in \Cref{lem: categorical kunneth prestack} reduces the statement to the case $X\in \qcsp_k$.

Now as $\shv(X)$ is dualizable, we see that $-\otimes_\La\shv(X)=\fun_{\lincat_\La}(\shv(X)^\vee,-)$ commutes with limits. Therefore, $\shv(\bB H,\La)\otimes_\La \shv(X,\La)$ can be computed as the totalization of $\shv(H^\bullet,\La)\otimes_\La\shv(H,\La)$ by descent. Similarly, $\shv(\bB H\times X,\La)$ is computed by $\shv(H^\bullet\times X)$. 
Then by the comonadic version of \cite[Theorem 4.7.3.5]{Lurie.higher.algebra}, it is enough to identify
the two comonads associated to these two cosimplicial diagrams.

We consider the following diagram
\[
\xymatrix{
\shv(\bB H)\otimes_\La \shv(X)\ar^-{f^!\otimes\id }[rr]\ar_{\boxtimes}[d] &&\shv(\pt)\otimes_\La \shv(X)= \shv(X)\ar@{=}[d]\\
 \shv(\bB H\times X)\ar^-{(f\times \id)^!}[rr] && \shv(\pt\times X)=\shv(X).
}
\]
Note that $f^!\otimes \id$ is conservative, as $\boxtimes$ is fully faithful, and $(f\times\id)^!$ is conservative (by descent).

As explained in \Cref{cor-additional-base-change-for-shv-theory-on-prestacks},
\Cref{assumptions.base.change.sheaf.theory.H} holds for $\mathrm{HR}$ being  the class of representable coh. pro-smooth morphisms. Therefore, the above diagram is also right adjointable.
Therefore, $\boxtimes\circ (f_\flat\otimes \id)(f^!\otimes \id)\cong (f\times\id)_\flat (f\times\id)^!\circ \boxtimes$, giving the identification of these two comonads.
\end{proof}

Another case we need is as follows.
We assume that $k$ is an algebraically closed field (and $\ell\neq 0$ in $k$) and identify profinite groups with affine group schemes over $k$ as before.  For a profinite group $K$, let $C^\infty(K,\La)$ denote the space of $\La$-valued smooth functions on $K$ acted by $K$ by right translation.

\begin{proposition}\label{ess.pro.p.torsor.hom.descent} 
If $K$ admits a $\La$-valued Haar measure (i.e. there exists a $K$-equivariant map $C^\infty(K,\La)\to \La$, which sends the characteristic function of some open compact subgroup $K'\subset K$ to an invertible element in $\La$), then $\Spec k\to \bB_{\mathrm{fpqc}} K$ is of universal homological descent.
\end{proposition}

\begin{proof}
Let $K'\subset K$ be an open subgroup such that its volume with respect to one Haar measure is $c\in\La^\times$. By $h$-descent, 
we can assume $K=K'$. Let $\pi\colon E\rightarrow X$ be a $K$-torsor. Again by \Cref{prop-descent-codescent-abstract} \eqref{prop-descent-codescent-abstract-2}, it's enough to construct a section of the natural map $\omega_X\to \pi_*\omega_E$ . Since $E \simeq \lim_{i} E_i$ we have $\shv(E,\La) \simeq \colim_{i} \shv(E_i,\La)$ and under this identification
\[
\pi_* \omega_E \simeq \colim_{i} (\pi_{i})_* \omega_{E_i}.
\]
Each map $\pi_i \colon E_i \rightarrow X$ is a torsor under the finite group $K_i=K/K^i$. In particular, we can identify $(\pi_i)_*$ with $(\pi_i)_!$ and the co-unit gives a natural map $s_i\colon (\pi_i)_* \omega_{E_i}\rightarrow \omega_X$. For each $i$, under any \'{e}tale trivialization $Y\rightarrow X$ of $\pi_i$ the pullback of the natural map $s_i$ identifies with the augmentation map
$C(K_i,\La) \rightarrow \Lambda$.
We can modify each map $s_i$ to a map $t_i \colon (\pi_i)_* \omega_{E_i}\rightarrow \omega_X$ by composing it with multiplication by $\mathrm{Vol}(K^i)=\frac{c}{[K:K^i]}$. The system of maps $\{t_i\}_{i\in \mI}$ is now compatible and $t = \colim_{i} t_i$ gives the desired section.
\end{proof}

Recall that associated to $K$ there is the constant affine group scheme $\underline{K}_\La$ over $\La$ so we have the $\qcoh(\bB_{\mathrm{fpqc}}\underline{K}_\La)$ as in \Cref{ex: fpqc quotient stack}.

\begin{corollary}
\label{prop:sheaves-classifying-stack-profinite-group-sheaves-rep}
If $K$ admits a Haar measure, then there is a canonical $t$-exact equivalence 
\[
\shv(\bB_{\mathrm{fpqc}} K,\La) \simeq \qcoh(\bB_{\mathrm{fpqc}} \underline{K}_\La)
\] 
such that the $!$-pullback functor $\shv(\bB_{\mathrm{fpqc}} K)\to \shv(\Spec k)$ is identified with the $*$-pullback $\qcoh(\bB_{\mathrm{fpqc}} \underline{K}_\La)\to\Mod_\La$. Under the equivalence, the $!$-tensor product on $\shv(\bB_{\mathrm{fpqc}}K)$ is identified with the usual tensor product on $\qcoh(\bB_{\mathrm{fpqc}} \underline{K}_\La)$. In particular, under such equivalence, tensor units are identified, i.e. $\consdual_{\bB_{\mathrm{fpqc}} K}$ corresponds to $\mO_{\bB_{\mathrm{fpqc}} \underline{K}_\La}$. 

In addition, $\shv(\bB_{\mathrm{fpqc}} K,\La)$ is compactly generated.
\end{corollary}
\begin{proof}
Proposition \ref{ess.pro.p.torsor.hom.descent} gives us a comparison between the category of sheaves on $\bB_{\mathrm{fpqc}} K$ with the totalization of the standard cosimplicial object $\shv(K^\bullet,\La)$, which as seen in \Cref{ex:sheaves-example-descent} is equivalent to $\qcoh(\bB_{\mathrm{fpqc}} \underline{K}_\La)$. The identification of tensor structures is clear.

It is enough to prove the compact generation of $\qcoh(\bB_{\mathrm{fpqc}} \underline{K}_\La)$. If $K'\subset K$ is an open compact subgroup such that the volume of $K'$ is invertible in $\La$, then the $*$-pushforward of $\mO_{\bB_{\mathrm{fpqc}}\underline{K'}_\La}$ to $\bB_{\mathrm{fpqc}}\underline{K}_\La$ is a projective object in the abelian category $\qcoh(\bB_{\mathrm{fpqc}} \underline{K}_\La)^{\heartsuit}$. As if $K'$ is such a subgroup, any open compact subgroup of $K'$ is also such a subgroup. Therefore, when $K'$ range over all such open compact subgroups of $K$, these objects form a set of generators of $\qcoh(\bB_{\mathrm{fpqc}} \underline{K}_\La)^{\heartsuit}$. We then apply \Cref{lem: cpt gen. of derived category of Groth ab cat} to conclude that $\qcoh(\bB_{\mathrm{fpqc}} \underline{K}_\La)$ is compactly generated.

The rest claims of the corollary are clear.
\end{proof}

\begin{corollary}\label{cor: profinite tensor product}
Suppose $X$ is a prestack over $k$. Then the exterior tensor product functor $\shv(\bB_{\fpqc} K,\La)\otimes_\La\shv(X,\La)\to \shv(\bB_{\fpqc} K\times X,\La)$ is an equivalence.
\end{corollary}
\begin{proof}
Given \Cref{ess.pro.p.torsor.hom.descent}, the same arguments as in \Cref{cor: ess. pro-unip tensor product} applies.
\end{proof}

\begin{remark}\label{rem: fpqc and proet BK}
One can replace fpqc topology in the above two statements by pro-\'etale or pro-finite \'etale topology. In fact, for $K$ profinite, the natural map $\bB_{\mathrm{profet}}K\to \bB_{\fpqc}K$ is an isomorphism.
\end{remark}

\subsection{Cosheaf theory on placid stacks}\label{sec-placid.stacks}
In this section we introduce a notion of placid stack in the setting of perfect algebraic geometry, following the terminology of \cite{bouthier2020perverse}. (As before, the actual meaning of this notion in this article is different from \emph{loc.cit.}) Recall that in classical algebraic geometry, an Artin stack (locally of finite presentation) over $k$ is a(n \'etale) stack $X$ that admits a smooth atlas $U\to X$ with $U$ an algebraic space (locally of finite presentation) over $k$. 
Roughly speaking, a (quasi-)placid stack generalizes this notion by allowing $U$ to be standard placid algebraic spaces over $k$ and $U\to X$ to be pro-smooth.  The theory of $\ell$-sheaves on this class of stacks are more fruitful than the theory on the general prestacks. For example, there is a good theory of constructible sheaves. Verdier duality and the perverse sheaves behave well on them as well.

We keep assumptions that $k$ is the perfection of a regular noetherian ring of $\leq 1$ and $\ell$ is a prime that is invertible in $k$ such that $k$ has finite $\bF_\ell$-cohomological dimension. We allow $\La$ to be any $\bZ_\ell$-algebra as from \Cref{sec:adic-formalism}.

\subsubsection{Placid stacks}\label{sec-placid.stacks.def}

\begin{definition}\label{def.placid.stack}
\begin{enumerate}
\item An (\'etale) stack $X\colon (\perf_k)^{\op} \rightarrow \spc$ is called \textit{quasi-placid} if 
there exists a family of morphisms  $\{U_i\to X\}_{i\in I}$, where each $U_i \in \plsp_k$ and $U_i\to X$ is representable coh. pro-smooth morphisms, such that for every $S\in \qcsp_k$, there are finite subset $I_S\subset I$ such that $\{U_i\times_XS\to S\}_{i\in I_S}$ is jointly surjective. We call such family $\{U_i\rightarrow X\}_i$ a quasi-placid atlas. We say $X$ is quasi-compact if there is a quasi-placid atlas $U\to X$ with $U\in \plsp_k$.
\item A quasi-placid stack $X$ is called \textit{placid} if there is a quasi-placid atlas $\{U_i\to X\}_i$  such that  
\begin{itemize}
\item each $U_i\to X$ is representable strongly coh. pro-smooth; and 
\item $\sqcup_i U_i\to X$ is of universal homological descent.
\end{itemize}
We call such a quasi-placid atlas as a placid atlas.
\item A quasi-placid stack $X$ is called \textit{very placid} if there is a quasi-placid atlas $\{U_i\to X\}_i$ such that each $U_i\to X$ factors as $U_i\to X_i\to X$ where $U_i\to X_i$ is ess. coh. pro-unipotent and $X_i\to X$ is an open embedding.
\end{enumerate}
Note that by \Cref{pro-unipotent.morphism.descent} (and Zariski descent), very placid stacks are placid. We let $\vplstk_k\subset \plstk_k\subset \qplstk_k\subset \prestk_k$ denote the corresponding full subcategories of very placid, placid and quasi-placid stacks over $k$. 
\end{definition}

\begin{remark}\label{rem-cond-on-placid-stacks}
We make a few remarks. 
\begin{enumerate}
\item Perhaps what we defined should be called quasi-separated (quasi-)placid stacks, but such generality is enough for our purpose. Note that a quasi-placid atlas is an epimorphism in $v$-topology. Therefore,
by $v$-descent of $\cshv$, there is a good theory of constructible sheaves on quasi-placid stacks, as we shall see in \Cref{sec-constructible.sheaves.on.placid.stacks}.  However, \Cref{ex:sheaves-example-descent} shows that $v$-descent fails for $\sshv$ (and for $\shv$) in general so the category of all sheaves on a quasi-placid stack could be wild. We also note that  the topological space associated to a quasi-placid stack could be quite wild. (See \cite[\textsection{11}]{Scholze.etale.cohomology.diamonds} for some discussions in a different but related setup.)
These are the reasons we impose stronger condition in the definition of placid stacks. In particular, if $X$ is placid,  then $|U_i|\to |X|$ is a quasi-compact\footnote{Quasi-compactness follows from our convention of representable morphisms. See \Cref{rem: qcqs representable}.} open map and $|X|$ is a quasi-separated spectral topological space. In particular, a placid stack is quasi-compact if and only if the underlying topological space $|X|$ is quasi-compact.

\item  Note that if the quasi-placid atlas $\{U_i\to X\}_i$ is an effective epimorphism in \'etale topology, then it is a placid atlas. In particular, combined with \cite[\textsection{1.3.3.(a)}]{bouthier2020perverse}, we get that (perfect) $1$-placid stacks in the sense of \cite{bouthier2020perverse} are placid stacks in the terminology of this paper.\footnote{But not conversely. E.g. the prestack $\bB K$ from \Cref{ess.pro.p.torsor.hom.descent} would not be $1$-placid in the sense of \cite{bouthier2020perverse}.}  Moreover, one could develop an analogous theory of $\infty$-stacks and $\infty$-smooth morphisms as in \emph{loc. cit}. The main results of this section generalize to that setting as well.  
\end{enumerate}
\end{remark}

\begin{example}\label{ex-non-standard-placid-space}
Note that $\plsp_k\subset \vplstk_k$. The perfection of a(n quasi-separated) Artin stack locally of finite presentation over $k$ is a very placid stack. On the other hand, $X\in\qcsp_k$ is quasi-placid if there is a surjective coh. pro-smooth morphism $U\to X$ with $U\in\plsp_k$.
\end{example}

\begin{example}\label{torsor.descent.gives.placid.stack} 
Suppose $H$ is an affine group scheme over $k$ that can be written as a cofiltered system $\{H_i\}_{i\in \mI}$ of perfectly smooth group schemes over $k$ with (perfectly) smooth affine transition maps, and suppose $H$ acts on a standard placid space $X$. Then the \'etale quotient stack $X/H$ is placid and the fpqc quotient stack
$(X/H)_{\mathrm{fpqc}}$ is a quasi-placid stack. 
The morphism $X\to (X/H)_{\mathrm{fpqc}}$ is representable strongly coh. pro-smooth, but may not be of universal homological descent in general. However,
suppose in addition we have the short exact sequence \eqref{eq-exact-sequence-flat-group} and suppose $H_0$ is coh. pro-unipotent over $k$. Then $(X/H)_{\mathrm{fpqc}}$ is very placid. 
\end{example}

\begin{remark}\label{rem: representability of diagonal}
We in general do not require the diagonal of a (quasi-)placid stack is representable in algebraic spaces. However, the diagonal of (quasi-)placid stacks from examples in \Cref{torsor.descent.gives.placid.stack} are affine. 
\end{remark}

We also have the following basic representability result, which follows immediately from the definition.  
\begin{lemma}\label{prop:placid.stack.placid.map}
Let $f\colon X \rightarrow Y$ be a representable ess. coh. pro-smooth morphism of \'etale stacks with $Y$ being a quasi-placid stack.  Then $X$ is a quasi-placid stack. If in addition $Y$ is (very) placid, so is $X$.
\end{lemma}
\begin{proof}
Let $V\rightarrow Y$ be a quasi-plaicd atlas. Note that $U= X\times_{Y} V\to V$ is  ess. coh. pro-smooth. Therefore, $U$ is standard placid by \Cref{composition.coh.pro-smooth}. So $U\to X$ is a quasi-placid atlas for $X$. Clearly, if $V\to Y$ is strongly coh. pro-smooth and of universal homological descent, or is ess. coh. pro-unipotent, so is $U\to X$.
\end{proof}

Recall that it makes sense to ask whether a representable morphism between prestacks it is (strongly, weakly) coh. pro-smooth. 
Now we generalize the notion of weakly coh. pro-smooth morphisms to non-representable morphisms.


\begin{definition}\label{def-l-coh-smooth-morphism}
A morphism $f:X\to Y$ of prestacks is called weakly coholomogically pro-smooth (resp. weakly pro-\'etale) if for every map $S\to Y$ with $S\in \qcsp_k$, there is a family of morphisms $\{T_i\to S\times_YX\}_{i\in I}$, where each $T_i\in \qcsp_k$ and 
$T_i\to S\times_YX$ is  representable cohomologically pro-smooth (resp. pro-\'etale)  such that each composed map $T_i\to S\times_YX\to S$ is cohomologically pro-smooth (resp. pro-\'etale). In addition, we require that for every $S'\to S\times_YX$ with $S'\in \qcsp_k$, there is a finite subset $I_{S'}\subset I$ such that $\{T_i\times_{S\times_YX}S'\to S'\}_{i\in I_{S'}}$ is jointly surjective.
\end{definition}

Note that for representable morphisms, this definition coincides with the old definition. 
The follow lemma is easy (using \Cref{composition.coh.pro-smooth}). 

\begin{lemma}\label{lem: atlas for weakly coh prosmooth}
\begin{enumerate}
\item\label{lem: atlas for weakly coh prosmooth-1} The class of weakly coh. pro-smooth morphisms is weakly stable. 
\item\label{lem: atlas for weakly coh prosmooth-2} The class of weakly pro-\'etale morphisms is strongly stable.
\item\label{lem: atlas for weakly coh prosmooth-3} Let $f: X\to Y$ be a weakly coh. pro-smooth morphism with $Y$ quasi-placid. Then $X$ is quasi-placid and there are quasi-placid atlas $\{\varphi_i:U_i\to X\}_i$ and $\{\varphi_j:V_j\to Y\}$, such that for every $i$, there is some $j$ and a coh. pro-smooth morphism $h_{ij}:U_i\to V_j$ such that $\varphi_j\circ h_{ij}=f\circ \varphi_i$.
\end{enumerate}
\end{lemma}

\begin{example}\label{ex-coh-prosmooth}
Let $H$ be as in \Cref{torsor.descent.gives.placid.stack}. Then $\bB H\to \Spec k$ is weakly coh. pro-smooth.
\end{example}

\subsubsection{Constructible sheaves on quasi-placid stacks}\label{sec-constructible.sheaves.on.placid.stacks}

The notion of constructible (co)sheaves can be defined on any prestack via right Kan extension as in \eqref{def-dual-sheaves-on-prestacks} and they form a full subcategory of all (co)sheaves. For a general prestack, this is not a useful notion. For quasi-placid stacks, however, this notion is well-behaved, and plays an important role in this article, as we shall see. 

First, since constructability is local with respect to the $v$-topology (by \Cref{prop-h-descent-shv} \eqref{prop-h-descent-shv-1}), a sheaf $\mF\in \shv(X,\La)$ on a placid stack $X$ is constructible if and only if for some, equivalently any, (quasi-)placid atlas $\{\varphi_i\colon U_i\rightarrow X\}_i$ the pullback $(\varphi_i)^!\mF$ is constructible on $U_i$ for every $i$.

\begin{example}\label{ex-constructible-on-BK}
Consider the situation as in \Cref{torsor.descent.gives.placid.stack}. Then an object $\mF\in \shv((X/H)_{\mathrm{fpqc}},\La)$ is constructible if the its $!$-pullback to $X$ is constructible. In particular, take $H=K$ be as in \Cref{ess.pro.p.torsor.hom.descent}. Then a sheaf $V\in \shv(\bB_{\mathrm{fpqc}} K)$, which identifies with an object of $\qcoh(\underline{K}_\La)$ (by \Cref{prop:sheaves-classifying-stack-profinite-group-sheaves-rep}), is constructible if and only if the underlying object $V\in \modu_\La$ is perfect. That is, $\cshv(\bB_{\mathrm{fpqc}} K)\subset \shv(\bB_{\mathrm{fpqc}} K)$ corresponds to $\Perf(\underline{K}_\La) \subseteq \qcoh(\underline{K}_\La)$. 
\end{example}

By definition, the constructible categories are preserved by $!$-pullback along any morphism. They are also preserved by other functors under usual finiteness assumptions, as we shall see now.

\begin{proposition}\label{functors.preservating.cons.categories.pl.stack} Let $f\colon X\rightarrow Y$ be a morphism of quasi-placid stacks, and let $g: Y'\to Y$ be a weakly coh. pro-smooth morphism. Consider the Cartesian diagram \eqref{eq:pullback-square-of-algebraic-spaces} of prestacks.
\begin{enumerate}
\item\label{functors.preservating.cons.categories.pl.stack-1}  If $f$ is representable ess. coh. pro-unipotent, then $f^!$ admit a left adjoint when restricted to constructible subcategories $f_!: \cshv(X)\to \cshv(Y)$. In addition, there is  the base change isomorphisms $(f')_!\circ(g')^!\xrightarrow{\cong} g^!\circ f_!$ of functors between constructible categories.

\item\label{functors.preservating.cons.categories.pl.stack-3} If $f$ is in addition representable pfp, then $f_*$ preserves constructible objects and admits a left adjoint when restricted to the constructible subcategories $f^*: \cshv(Y)\to \cshv(X)$. In addition, there is the base change isomorphism $(f')^*\circ g^!\xrightarrow{\cong} (g')^!\circ f^*$ of functors between constructible categories.

\item\label{functors.preservating.cons.categories.pl.stack-2} If $f$ is representable coh. smooth, then $f_\flat: \shv(X)\to \shv(Y)$ preserves constructibility, and we have the base change isomorphism  $(f')_\flat\circ(g')^!\xrightarrow{\cong} g^!\circ f_\flat$.
\end{enumerate}
\end{proposition}
\begin{proof}
Note that \Cref{lem:etale-proper-functoriality-shv} \eqref{lem:etale-proper-functoriality-shv-6} \eqref{pro-unipotent.are.admissible} together with descent imply the existence of left adjoints for Part \eqref{functors.preservating.cons.categories.pl.stack-1} and \eqref{functors.preservating.cons.categories.pl.stack-3}. They also imply the base change isomorphism in the special case when $Y'\to Y$ is a quasi-placid atlas.

Next we prove the base change isomorphisms as in Part \eqref{functors.preservating.cons.categories.pl.stack-1} and \eqref{functors.preservating.cons.categories.pl.stack-3} for a general weakly coh. pro-smooth morphism $Y'\to Y$ of quasi-placid stacks.  We can find  quasi-placid atlases $\{\varphi'_i\colon V'_i\rightarrow Y'\}_i$ and $\{\varphi_j\colon V_j\rightarrow Y\}_j$, and coh. pro-smooth morphisms $h_{ij}:V'_i\to V_j$ as in \Cref{lem: atlas for weakly coh prosmooth}. Let $U'_i=X'\times_{Y'}V'_i$ and $U_j=X\times_YV_j$. 
 So  $\psi'_{i}\colon U'_i\to X'$ and $\psi_j\colon U_j\to X$ are quasi-placid atlas of $X'$ and $X$  by the proof of \Cref{prop:placid.stack.placid.map}.
In addition, $h'_{ij}\colon U'_i\to U_j$ is the base change of $h_{ij}\colon V'_i\to V_j$ (e.g. see \cite[Lemma A.2.9]{xiao2017cycles}) and is coh. pro-smooth. By conservativity, it's enough to prove that the natural maps of functors between constructible categories
\[
(\varphi'_i)^!\circ (f')_!\circ (g')^!\rightarrow (\varphi'_i)\circ g^!\circ f_!, \quad (\psi'_i)^!\circ (f')^*\circ g^! \to  (\psi'_i)^!\circ (g')^!\circ f^*
\]
are isomorphisms. But these follow from the above mentioned special case applying to $\varphi'_i$ and $\varphi_j$ and  \Cref{lem:etale-proper-functoriality-shv} \eqref{lem:etale-proper-functoriality-shv-6} applying to $h_{ij}$.

Part \eqref{functors.preservating.cons.categories.pl.stack-2} follows from the above base change isomorphisms and \Cref{lem:etale-proper-functoriality-shv-prestack} \eqref{lem:etale-proper-functoriality-shv-prestack-4}. 
\end{proof}

Recall that in the $\sshv$-sheaf theory, a constructible complex $\mF\in\scshv(X,\La)$ with $X$ pfp over $k$ is ULA.

\begin{lemma}\label{lem: ULA-ish for constructible on quasi-placid}
Assume that $k$ is a field.
Let $X$ be a quasi-placid stack and $\mF\in \cshv(X)$. Let $f: Y'\to Y$ be a representable pfp morphism of quasi-placid stacks. Then for every $\mG\in \cshv(Y')$, we have
$(\id_X\times f)_!(\mF\boxtimes \mG)\cong \mF\boxtimes f_!\mG$.
\end{lemma}
\begin{proof}
By choosing atlas and the base change isomorphisms from \Cref{functors.preservating.cons.categories.pl.stack} \eqref{functors.preservating.cons.categories.pl.stack-1}, we may assume that $X,Y,Y'$ are standard placid spaces. Then we may assume that $X,Y,Y'$ are pfp over $k$. In this case, $\shv=\sshv$. As $k$ is a field, the usual Verdier duality commutes with exterior product and then we reduce to prove that $(\id_X\times f)_*(\mF\boxtimes \mG)\cong \mF\boxtimes f_*\mG$. This follows that $\mF$ is $\ell$-ULA with respect to $\pi_X: X\to \pt=\Spec k$ so \eqref{eq:Kunneth-formula-ULA} from \Cref{lem:base-change-for-dualizable}  is applicable (to the sheaf theory $\sshv$).
\end{proof}

\begin{remark}
We note that unlike the situation as in \Cref{rem-additional-base.change.sheaf.theory} \eqref{rem-additional-base.change.sheaf.theory-1}, the above isomorphism does not imply that $(f_!,f^!)$ satisfies a projection formula. This is because the base change isomorphisms as in  \Cref{functors.preservating.cons.categories.pl.stack} \eqref{functors.preservating.cons.categories.pl.stack-1} only holds for $g$ being weakly coh. pro-smooth.
\end{remark}

\Cref{functors.preservating.cons.categories.pl.stack} in particular says that (under some finiteness assumptions) there is a good six functor formalism for constructible sheaves on quasi-placid stacks, which can be regarded as a generalization of \Cref{constructible.preservation.placid.spaces}. However, unlike the situation of standard placid spaces, $\cshv(X,\La)$ and $\shv(X,\La)^\omega$ usually do not agree. In addition, we do not know whether $\shv(X,\La)$ is compactly generated, or even dualizable. Later on, we will say more about relations between compact objects and constructible objects in $\shv(X,\La)$ when $X$ is very placid stacks.
For quasi-placid stacks, we always have the following statements.

\begin{lemma}\label{lem-cpt-vs-cons-quasi-placid}
Let $X$ be a quasi-placid stack. Then $\shv(X,\La)^\cpt\subset \cshv(X,\La)$. 
\end{lemma}
\begin{proof}
Let $\{\varphi_i\colon U_i\rightarrow X\}_i$ be a quasi-placid atlas with $U_i\in\plsp_k$. By \Cref{lem:etale-proper-functoriality-shv-prestack} \eqref{lem:etale-proper-functoriality-shv-prestack-4}, $(\varphi_i)_\flat$ is continuous. This implies that $(\varphi_i)^!$ preserves compact objects which means $(\varphi_i)^!(\mF)$ is constructible for every compact object $\mF\in \shv(X,\La)$. So $\mF\in\cshv(X,\La)$.
\end{proof}

\subsubsection{Verdier duality and perverse sheaves for quasi-placid stacks}\label{sec-verdier.duality.placid.stacks}

\begin{definition}\label{def: generalized constant sheaf}
Let $X$ be a quasi-placid stack. A generalized constant sheaf of $X$ is an object $\La^\eta_X\in \cshv(X)$ such that for some (and therefore for any by  \Cref{lem-gen-dualizing-pro-smooth-local}) quasi-placid atlas $\{\varphi_i: U_i\to X\}$, $(\varphi_i)^!\La^\eta_X\in \cshv(U_i)$ is a generalized constant sheaf on $U_i$. 
\end{definition}

\begin{lemma}\label{lem:existence of generalized constant sheaf}
Generalized constant sheaves always exist on quasi-compact quasi-placid stacks.
\end{lemma}
The subtlety here lies in the non-canonicity of the isomorphism from \Cref{global.iso.coh.smooth} (as we work in perfect algebraic geometry), so one needs to provide a descent datum to  a generalized constant sheaf $\La^\eta_U$ on a quasi-placid atlas $U\to X$ to get a generalized constant sheaf on $X$.
\begin{proof}
Let $U\to X$ be a quasi-placid atlas and let $U_\bullet$ denote the \v{C}ech nerve. By \Cref{independence.of.placid.presentation}, we may assume that there is the following commutative diagram
\[
\xymatrix{
U_3\ar@<-.8ex>[r]  \ar[r]\ar@<.8ex>[r] \ar[d] & U_2 \ar@<-.4ex>[r] \ar@<.4ex>[r]  \ar[d] & U \ar[d]\\
U'_3\ar@<-.8ex>[r]  \ar[r]\ar@<.8ex>[r]  & U'_2  \ar@<-.4ex>[r] \ar@<.4ex>[r]   & U',
}\]
where the bottom line is a $\Delta_{\leq 2}^{\op}$-object in $\pfpsp_k$ with all morphisms coh. smooth, and where all vertical morphisms coh. pro-smooth. We fix a generalized constant sheaf $\La_{U'}\langle d_{U'}\rangle$ of $U'$, which on a connected component $C\subset U'$ is $\La_C\langle\dim C\rangle$. Note that by \eqref{eq-LaetaX-dual-to-omegaX} and \Cref{ex-dualizing sheaf}, generalized constant sheaves on placid spaces are discrete objects. 
So it is then enough to construct an isomorphism between two $!$-pullbacks
of $\La_{U'}\langle d_{U'}\rangle$ to $U'_2$ that satisfies a cocycle condition when further $!$-pulling back to $U'_3$. 

Choose a deperfection of this $\Delta_{\leq 2}^{\op}$-object  $U''_{3}\substack{\longrightarrow\\[-.9em] \longrightarrow \\[-.9em]  \longrightarrow} U''_{2} \substack{\longrightarrow \\[-.9em] \longrightarrow} U''$. So $U''_3,U''_2,U''$ are finitely presented algebraic spaces over $k$. The argument as in \Cref{global.iso.coh.smooth} gives isomorphisms $(d_i)^!\La_{U''}\cong \La_{U''_2}\langle d_{d_i}\rangle$ where $d_i: U''_{2}\to U'', \ i=0,1$ are two face maps, induced by the the map $(d_i)_!\La_{U''_2}\langle-d_{d_i}\rangle\to \La_{U''}$ (which restricts to the trace map over an open dense subset of $U''$). It follows that we obtain an isomorphism 
\[
\theta: (d_0)^!\La_{U''}\langle d_{U''}\rangle \cong \La_{U_1''}\langle d_{U''_1}\rangle \cong (d_1)^!\La_{U''}\langle d_{U''}\rangle
\] 
by composing (appropriate shift of) these two isomorphisms.
Similarly, the three face maps $d_i: U''_3\to U''_2$ give $(d_i)^!\La_{U''_2}\cong \La_{U''_3}\langle d_{d_i}\rangle$, again induced by trace maps. Since trace maps are compatible with respect to compositions, we see that $\theta$ satisfies the cocycle conditions over $U''_3$.
\end{proof}

Similar to the case in \Cref{sec-verdier.duality.dual.setting}, given a generalized  constant sheaf $\La^\eta_X$ on $X$, we define
\begin{equation}\label{eq: eta-cohomology-stack}
\rg^\eta_{\ind\fg}(X,-): \cshv(X,\La)\to \Mod_\La,\quad \rg^\eta_{\ind\fg}(X,\mF)=\Hom_{\cshv(X,\La)}(\La_X^\eta,\mF).
\end{equation}

With this definition, the following proposition follows from \Cref{verdier.functioriality.spaces-shv}, \Cref{functors.preservating.cons.categories.pl.stack}  and descent.
\begin{proposition}\label{verdier.duality.placid.stacks}
Let $X$ be a quasi-placid stack equipped with a generalized constant sheaf $\La^\eta_X$. There is a canonical equivalence 
\begin{equation}\label{eq: Verdier duality for quasi-placid}
(\verd_X^{\eta})^c \colon \cshv(X,\La)^{\op} \simeq \cshv(X,\La)
\end{equation}
with $((\verd_{X}^{\eta})^c)^2\simeq \id$, uniquely characterized by
\begin{equation}\label{eq: Verdier duality for quasi-placid-characterization}
\Hom_{\cshv(X,\La)}(\mF,\mG) \simeq \rg^\eta_{\ind\fg}(X, (\verd_X^\eta)^c(\mF)\os \mG), \quad \mF,\mG\in \cshv(X,\La).
\end{equation}

Let $f:X\to Y$ be a morphism of quasi-placid stacks, and let $\La_Y^\eta$ be a generalized constant sheaf on $Y$.  
If $f$ is representable pfp, then $\La_X^\phi=f^*\La_Y^{\eta}$ is a generalized constant sheaf on $X$ and we have isomorphisms of contravariant functors between constructible categories
\begin{equation*}
    (\verd_Y^{\eta})^c\circ f_*  \simeq f_! \circ (\verd_X^{\phi})^c,\quad
    (\verd_X^{\phi})^c\circ f^!  \simeq f^* \circ (\verd_Y^{\eta})^c.
\end{equation*}
If $f$ is weakly coh. pro-smooth, then $\La^\phi_{X}: = f^!\La^\eta_Y$ is a generalized constant sheaf on $X$, and we have an isomorphism of contravariant functors between constructible categories
\begin{equation*}
(\verd_X^{\phi})^c\circ f^!  \simeq f^! \circ (\verd_Y^{\eta})^c.
\end{equation*}
\end{proposition}

\begin{remark}\label{rem: equiv between cons-sheaf and cons-cosheaf}
For  a chosen generalized constant sheaf $\La^\eta_X$, the equivalence in \eqref{eq: Verdier duality for quasi-placid} and
the equivalence \eqref{eq: star and shrek constr sheaves on quasi-placid} induces an equivalence 
\begin{equation}\label{eq: identifying star sheaf with shrek sheaf-2}
\id^\eta:\cshv(X,\La)\cong \scshv(X,\La)
\end{equation}
generalizing \eqref{eq: identifying star sheaf with shrek sheaf-1}. Under this equivalence, the usual $*$-tensor product of $\scshv(X,\La)$ becomes the following tensor product on $\cshv(X,\La)$,
\begin{equation}\label{eq: eta tensor product on cshv}
\cshv(X)\otimes \cshv(X)\to \cshv(X),\quad (\mF,\mG)\mapsto (\verd_X^\eta)^c(((\verd_X^\eta)^c(\mF))\os ((\verd_X^\eta)^c(\mG))),
\end{equation}
This generalizes \Cref{rem: duality and tensor product}. Note that \eqref{eq: Verdier duality for quasi-placid-characterization} is equivalent to
\begin{equation}\label{eq: Verdier duality for quasi-placid-characterization-1}
\Hom_{\cshv(X,\La)}(\mF,\mG) \simeq \Hom_{\cshv(X,\La)}((\mF)\otimes^\eta (\verd^{\eta}_X)^c(\mG),\consdual_X), \quad \mF,\mG\in \cshv(X,\La).
\end{equation}
\end{remark}

\begin{remark}\label{rem:verdier.functioriality.spaces.pseudo.coh.pro.smooth-2}
Clearly, \Cref{rem:verdier.functioriality.spaces.pseudo.coh.pro.smooth} also generalizes to pseudo coh. pro-smooth morphisms between quasi-placid stacks.
\end{remark}

We also mention that for quasi-placid $X$ (and regular neotherian $\La$), by \Cref{lem-descent-perverse-sheaves} (and \Cref{rem-dual-perverse-sheaf-placid-space}), after choosing a generalized constant sheaf $\La_X^\eta$ on $X$ one can also define a perverse $t$-structure on $\cshv(X,\La)$. Namely, we let
\[
\cshv(X,\La)^{\eta,\geq 0}=\bigl\{\mF\in\cshv(X,\La)\mid \varphi^!\mF\in \cshv(U,\La)^{\phi,\geq 0}\bigr\}.
\]
Let $\mathrm{Perv}(X,\La)^\eta$ denote the corresponding category of perverse sheaves. Note that when $X$ is quasi-compact, the $\eta$-perverse $t$-structure on $\cshv(X,\La)$ is bounded, as this is the case for standard placid spaces, which in turn follows from the corresponding statement over pfp algebraic spaces over $k$. In addition, if $X$ is quasi-compact and $\La$ is a field, then $\mathrm{Perv}(X,\La)^{\eta}$ is stable under $(\verd_X^\eta)^c$. In addition,

Dualizing \Cref{lem-descent-perverse-sheaves}  gives the following.
\begin{proposition}\label{prop:perverse exact smooth pullback}
Let $f: X\to Y$ be a weakly coh. pro-smooth morphism between quasi-placid stacks, and let $\La_Y^\eta$ be a generalized constant sheaf on $Y$. Then $\La_X^\phi=f^!\La_Y^\eta$ is a generalized constant sheaf on $X$, and $f^!$ sends $\mathrm{Perv}(Y,\La)^\eta$ to $\mathrm{Perv}(X,\La)^\phi$.
\end{proposition}

\begin{example}\label{ex-perverse-sheaf-classifying-stack}
Let $H$ be as in \Cref{torsor.descent.gives.placid.stack} and assume that $H$ is connected. Then $\La^\eta=\consdual_{\bB H}$ is a generalized constant sheaf on $\bB H$. With respect to this choice, $\mathrm{Perv}(X,\La)^\eta\cong \Mod_\La^{\heartsuit}$. 
As the perverse $t$-structure on $\cshv(X,\La)$ is bounded, via perverse truncations, we see that $\cshv(\bB H)$ is generated (as an idempotent complete stable category) by perverse sheaves, and therefore by $\consdual_{\bB H}$. Note that $\consdual_{\bB H}$ is in general not compact (e.g. $H=\bG_m$).
Therefore, this example also shows that constructible sheaves (and therefore perverse sheaves) on a quasi-placid stack $X$ are usually not compact in $\shv(X)$. On the other hand in this case, by \Cref{lem: char of adm obj in cg cat} constructible sheaves are exactly admissible objects (as defined in \Cref{def:adm vs compact}) in $\shv(\bB H,\La)$.
\end{example}

\begin{example}\label{ex: usual perverse t-structure on algebraic stacks}
Let $X$ be the perfection of a fp algebraic stack over $k$. Then there is a generalized constant sheaf $X$ (whose $!$-pullback to a smooth cover $U$ is canonically the constant sheaf of $U$.
Then under the equivalence mentioned in \Cref{rem: equiv between cons-sheaf and cons-cosheaf}, the above perverse $t$-structure on $\cshv(X)$ corresponds to the usual perverse $t$-structure on $\scshv(X)$, as in \Cref{rem-dual-perverse-sheaf-placid-space}.
\end{example}

\subsubsection{The ind-finitely generated sheaves and $\renflat$-pushforward}\label{sec-renormalized.category.placid.stack}
Constructible (co)sheaves play an important role in the article.
However, it is convenient to pass to the compactly generated cocomplete categories to obtain a sheaf theory valued in $\lincat_\La$.

For this purpose, we let $\verti_c$ be a class of morphisms $(f:X\to Y)\in \qplstk_k$ satisfying the following two properties: 
\begin{itemize}
\item For every $Y'\to Y$ in $\qplstk_k$, $X'=Y'\times_YX\in \qplstk_k$. Let $f': X'\to Y'$ be the base change morphism.
\item $f$ belongs to the class $\verti$ from \eqref{eq: pushforward along non-representable morphisms} (so $(f')_*: \shv(X')\to \shv(Y')$ is defined), and $(f')_*$ sends $\cshv(X')$ to $\cshv(Y')$.
\end{itemize}

\begin{example}\label{rem: pushforward along non-representable morphisms constructible}
We note that by \Cref{rem-generalization-corr} \eqref{rem-generalization-corr-2} and by \Cref{functors.preservating.cons.categories.pl.stack}, the class $\pfp_r$ belongs to $\verti_c$. On the other hand, let $f: X\to Y$ be an $A$-gerbe map as in \Cref{rem: pushforward along non-representable morphisms}, with $A$ a finite group of order invertible in $\La$, then $f$ also belongs to $\verti_c$, by \Cref{lem: atlas for weakly coh prosmooth} \eqref{lem: atlas for weakly coh prosmooth-3}.
\end{example}

Note that fiber products in $\qplstk_k$ may not exist. However, thanks to
\Cref{prop:placid.stack.placid.map}, the category $\corr(\qplstk_k)_{\verti_c;\all}$ is still defined. Then we can define

\begin{equation}\label{eq-rshv-for-placid}
\rshv\colon \corr(\qplstk_k)_{\verti_c;\all} \rightarrow \lincat_\La
\end{equation}
obtained by ind-extension of the functor 
$\cshv\colon \corr(\qplstk_k)_{\verti_c;\all} \rightarrow \catid_\La$ sending $Y\xleftarrow{g} Z\xrightarrow{f} X$ to $f^{\ind\fg}_{*}\circ g^{\ind\fg,!}$. 

Here to be consistent with terminology and notations for the later discussions in \Cref{SS: ind and sifted placid}, we use $\rshv$ rather than $\ind\cshv$, and
call $\rshv(X,\La)$ the category of ind-finite generated sheaves on $X$.

All base change isomorphisms and projection formulas of this subsection involving functors which preserve the constructible subcategories as in \Cref{functors.preservating.cons.categories.pl.stack} extend by continuity to the analogous result for the categories of ind-finitely generated sheaves. 
In addition, by taking the colimit, there is a natural symmetric monoidal functor
\begin{equation}\label{eq-ind-constr-to-shv}
\Psi: \rshv(X,\La)\to \shv(X,\La),
\end{equation}
which commutes with the above discussed sheaf operations (such as $!$-pullbacks and $*$-pushforwards, and $!$-pushforwards and $*$-pullbacks when they are defined).

Now let $f\colon X\rightarrow Y$ be a weakly coh. pro-smooth morphism between quasi-placid stacks. Then $f^!$ may not preserve compact objects. As a result, $f_\flat: \shv(X,\La)\to \shv(Y,\La)$ may not be continuous. However, $f^!$ always preserves constructible subcategories. Therefore, the functor $f^{\ind\fg,!}$ admits a \emph{continuous} right adjoint, denoted by
\begin{equation}\label{eq-reflat-push}
f_{\renflat}\colon \rshv(X,\La) \rightarrow \rshv(Y,\La).
\end{equation}
We refer to it as the \textit{$\renflat$-pushforward}. By passing to the right adjoint of the base change isomorphisms from \Cref{functors.preservating.cons.categories.pl.stack}, we also obtain the following.

\begin{lemma}\label{prop-base-change-for-renflat-forward}
Let \eqref{eq:pullback-square-of-algebraic-spaces} be a pullback square of  quasi-placid stacks with $f\colon X\rightarrow Y$ being weakly coh. pro-smooth. Then if $g$ is representable ess. coh. pro-unipotent, there is the natural base change isomorphism of functors \[
g^{\ind\fg,!}\circ f_{\renflat}\xrightarrow{\cong}  (f')_{\renflat}\circ (g')^{\ind\fg,!}.
\]
     If addition $g$ is representable pfp, then there is the natural base change isomorphism 
     \[
     g^{\ind\fg}_{*}\circ (f')_\renflat\xrightarrow{\cong} f_\renflat\circ (g')^{\ind\fg}_{*}.
     \] 
\end{lemma}

It also satisfies the following general projection formula  (compare with \Cref{lem:etale-proper-functoriality-shv-prestack} \eqref{lem:etale-proper-functoriality-shv-prestack-4}) which follows by the same argument as in \Cref{lem:dual-sheaves-proj-formula-weakly-coh-pro-smooth} (using the Verdier duality for $\cshv$, as will be discussed in \Cref{sec-verdier.duality.placid.stacks}). By abuse of notations, the symmetric monoidal structure on $\ind\fgshv$ induced by this sheaf theory is still denoted as $\os$.

\begin{proposition}\label{ren.flat.projection.formula}
Let $f\colon X\rightarrow Y$ be a weakly coh. pro-smooth morphism of quasi-compact quasi-placid stacks. 
 Then the natural map 
 \[
 f_\renflat(\mF)\os \mG \rightarrow f_\renflat(\mF\os f^{\ind\fg,!}(\mG))
 \] 
 is an isomorphism for an $\mF\in \rshv(X)$ and $\mG\in \rshv(Y)$.
\end{proposition}

We also record the following statement.
\begin{lemma}\label{lem-verdier-conjugate-!-push}
If $f: X\to Y$ is a representable ess. coh. pro-unipotent and weakly coh. pro-smooth morphism between quasi-compact quasi-placid stacks, then for a generalized constant sheaf $\La_Y^\eta$ on $Y$ and $\La_X^\phi=f^!\La_Y^\eta$, we have $(\verd_Y^\eta)^c\circ f_!\circ (\verd_X^\phi)^c\cong f_\renflat$. In particular $f_\renflat$ preserves constructibility.
\end{lemma}
\begin{proof}
It is enough to show that 
\[
\Hom_{\cshv(Y)}(\mG,(\verd_Y^\eta)^c(f_!((\verd_X^\phi)^c(\mF))))\cong\Hom_{\cshv(X)}(f^!\mG, \mF),\quad \forall \mF\in \cshv(X),\ \mG\in \cshv(Y).
\] 
This follows easily from \Cref{verdier.duality.placid.stacks}.
\end{proof}

Assume that $\La$ is regular noetherian, and recall the standard $t$-structure on $\shv_{(c)}(-,\La)$ as discussed in \Cref{rem:standard t-structure on prestack shv}. As before, by ind-completion, we obtain the standard $t$-structure on $\rshv(X,\La)$ with $\rshv(X,\La)^{\mathrm{std},\leq 0}=\ind \cshv(X,\La)^{\mathrm{std},\leq 0}$.

\begin{lemma}\label{lem: t-structure renomralized}
The functor $\Psi: \rshv(X,\La)\to \shv(X,\La)$ is exact, which restricts to an equivalence $\rshv(X,\La)^{\mathrm{std},\geq 0}\to \shv(X,\La)^{\mathrm{std},\geq 0}$ when $X$ is quasi-compact placid.
\end{lemma}
\begin{proof}
As $\cshv(X,\La)\to \shv(X,\La)$ is $t$-exact and the $t$-structure on $\shv(X,\La)$ is compatible with filtered colimits, $\Psi$ is $t$-exact. As $X$ is placid, $\shv(X,\La)$ satisfies descent with respect to a placid atlas $U\to X$. Then we can argue as in \Cref{lem: coconnective part of indcoh} to prove the second statement. (We note that the argument as in \cite[\href{https://stacks.math.columbia.edu/tag/0GRF}{Lemma 0GRF}]{stacks-project} is also applicable to the current setting.) 
\end{proof}

We give two applications of \Cref{lem: t-structure renomralized}.
First as \Cref{eq: indcoh-*-pushforward} and \Cref{lem: push-forward of qcoh for coconnective part} , we have the following.
\begin{lemma}\label{prop-base-change-for-renflat-forward-difficult}
Suppose \eqref{eq:pullback-square-of-algebraic-spaces} is a Cartesian diagram of quasi-compact placid stacks, and suppose $f$ is weakly coh. pro-smooth. Then $g^{\ind\fg,!}\circ f_\renflat\cong (f')_\renflat\circ (g')^{\ind\fg,!}$.
\end{lemma}
Next, a similar argument as in \Cref{lem: coherent Kunneth} gives the following.
\begin{proposition}\label{lem: kunneth for indfg sheaf}
Assume that $\La$ is regular noetherian as above.
Let $X$ and $Y$ be two quasi-compact very placid stacks over an algebraically closed field $k$. Then the exterior tensor product functor $\rshv(X)\otimes_\La\rshv(Y)\to \rshv(X\times Y)$ is fully faithful.
\end{proposition}

\begin{corollary}\label{cor: fg shv exterior tensor}
Let $H$ be a connected affine group scheme as in \eqref{eq-exact-sequence-flat-group} with $H_0$ coh. pro-unipotent and $H'$ pfp over $k$. Then for every quasi-compact placid stack $X$,  the exterior tensor product functor $\rshv(\bB H)\otimes_\La\rshv(X)\to \rshv(\bB H\times X)$ is an equivalence.
\end{corollary}
\begin{proof}
By \Cref{lem: kunneth for indfg sheaf}, it remains to prove that the image of the exterior tensor product functor generates $\rshv(\bB H\times X)$. We fix a generalized constant sheaf $\La_X^\eta$ on $X$ (and a standard one on $\bB H$).
As explained in \Cref{ex-perverse-sheaf-classifying-stack}, $\rshv(\bB H\times X)$ is generated by perverse sheaves. As $H$ is connected, the $!$-pullback along $X\to \bB H\times X$ induces an equivalence of categories of perverse sheaves. This shows that the image of the exterior tensor product functor generates $\rshv(\bB H\times X)$.
\end{proof}

\subsubsection{Cosheaves on very placid stacks}
For a very placid stack $X$, we can say more about its category $\shv(X)$ of all cosheaves.  To simplify the exposition, we add the quasi-compactness assumption throughout. But this assumption can be dropped in some statements.

\begin{proposition}\label{compact.generation.admissible.stacks} 
Let $X$ be a quasi-compact very placid stack with a placid atlas $\varphi\colon U\to X$ that is representable ess. coh. pro-unipotent. 
Then the category $\shv(X)$ is compactly generated. The functor $\varphi_!$ preserves compact objects and moreover, the subcategory $\shv(X)^\cpt\subseteq \shv(X)$ is the smallest idempotent complete stable full subcategory containing the objects $\varphi_!(\mF)$ for $\mF\in \cshv(U)$. In particular, \eqref{eq-ind-constr-to-shv} admits a fully faithful left adjoint $\Psi^L: \shv(X,\La)\subseteq \rshv(X,\La)$.

If, in addition $\varphi$ is representable coh. pro-unipotent, then $\shv(X,\La)^\cpt=\cshv(X,\La)$. 
\end{proposition}
Note that \Cref{prop:sheaves-classifying-stack-profinite-group-sheaves-rep} gives an example that $\shv(X)$ is compactly generated when $X$ is not very placid.
\begin{proof}
As $U\to X$ is of universal homological descent, we have $\shv(X)\cong \tot\left( \shv(U_\bullet)\right)$, where $U_\bullet\to X$ is the \v{C}ech nerve of $U\to X$. Using \Cref{lem:etale-proper-functoriality-shv} \eqref{lem:etale-proper-functoriality-shv-6} \eqref{pro-unipotent.are.admissible}, we see that there is the adjunction $\varphi_! \colon \shv(U) \rightleftarrows \shv(X)\colon \varphi^!$.
As $\varphi^!$ is conservative, the category $\shv(X)$ is compactly generated, with a set of generators given by $\varphi_!(\mF)$ with $\mF\in \cshv(U)$. This gives the desired description of $\shv(X)^\cpt$.
The last statement follows as $\shv(X) \rightarrow \shv(U)$ is fully faithful.
\end{proof}

We have a
version of \Cref{functors.preservating.cons.categories.pl.stack} for all sheaves on very placid stacks.

\begin{proposition}\label{pfp.functors.and.bc.placid.stacks}  Let \eqref{eq:pullback-square-of-algebraic-spaces} be a pullback square of quasi-compact very placid stacks with $g$ being weakly coh. pro-smooth. If $f$ is representable ess. coh. pro-unipotent, then $f^!: \shv(Y)\to \shv(X)$ admits a left adjoint $f_!$, and there is the base change isomorphism $(f')_!\circ(g')^!\cong g^!\circ f_!$. 

If in addition $f$ is representable pfp, then $f_*: \shv(X)\to \shv(Y)$ also admits a left adjoint $f^*$, and there is the base change isomorphism 
$(f')^*\circ g^!\cong (g')^!\circ f^*$.
\end{proposition}
Note that from the proof below that the existence of left adjoints (but not base change) only requires that $X$ and $Y$ to be placid. In addition, these left adjoints restricts to the corresponding left adjoints for constructible subcategories from \Cref{functors.preservating.cons.categories.pl.stack}. We also note that \Cref{lem:etale-proper-functoriality-shv-prestack} \eqref{lem:etale-proper-functoriality-shv-prestack-3} says that $f^*$ exists under a different assumption.
\begin{proof}
As placid atlases are of universal homological descent, 
again \Cref{lem:etale-proper-functoriality-shv} \eqref{lem:etale-proper-functoriality-shv-6} \eqref{pro-unipotent.are.admissible} imply the existence of left adjoints. We know the base change isomorphisms hold for constructible sheaves by \Cref{functors.preservating.cons.categories.pl.stack}, and therefore hold for compact objects by \Cref{lem-cpt-vs-cons-quasi-placid}. As the categories are compactly generated by \Cref{compact.generation.admissible.stacks}, the base change isomorphisms hold for all sheaves.
\end{proof}

Using \eqref{eq:limit-colimit equivalence}, we have the following colimit presentation of $\shv(X)$.
\begin{corollary}\label{cor: colim presentation of shv for very placid}
Let $X$ be a quasi-compact very placid stack with a placid altas $U\to X$ that is ess. coh. pro-unipotent, and let $U_\bullet\to X$ be its \v{C}ech nerve as before. Then
\[
\shv(X)=\colim \shv(U_\bullet),
\]
with colimit taken with repsect to $!$-pushforwards along face maps.
\end{corollary}

\begin{example}\label{ex: compactly supported cohomology}
Let $X$ be a quasi-compact very placid stack with a placid altas $U\to X$ that is ess. coh. pro-unipotent and such that $U$ is ess. coh. pro-unipotent over $\pt$. Then the $!$-pushforward $(\pi_X)_!$ along the structural map $\pi_X:X\to \pt$ is defined. In particular, if $\pi=\spec k$ with $k$ an algebraically closed field, we write $C_c(X,-)$ in stead of $(\pi_X)_!$, which should be thought as the compactly supported cohomology of $X$.
\end{example}

Passing to ind-completions,
\Cref{verdier.duality.placid.stacks} in particular implies that if $X$ is a quasi-compact quasi-placid stack, then  $\rshv(X,\La)$ is self-dual, with a Frobenius structure defined by the ind-extension of  \eqref{eq: eta-cohomology-stack}, denoted by the same notation.  In addition, \eqref{eq: Verdier duality for quasi-placid-characterization} holds for $\mF\in\cshv(X)$ and $\mG\in\rshv(X)$.

However, such duality does not give a duality of $\shv(X)$ in general.
In some cases, one can still deduce from \Cref{verdier.duality.placid.stacks} that $\shv(X)$ is also self-dual. One example is $X=\bB_{\mathrm{fpqc}}K$ as from \Cref{prop:sheaves-classifying-stack-profinite-group-sheaves-rep}. Here we shall discuss another situation. Recall that there is the functor $\Psi: \rshv(X)\to \shv(X)$, which admits a fully faithful left adjoint $\Psi^L$ when $X$ is quasi-compact very placid.

\begin{proposition}\label{prop: verdier duality on shv for very placid}
Suppose $k$ is algebraically closed and let $H$ be an affine group scheme over $k$ fitting into a short exact sequence as in \eqref{eq-exact-sequence-flat-group} with $H_0$ coh. pro-unipotent. Suppose $X=U/H$ with $U$ being a standard placid algebraic space (so $X$ is very placid). Then the Frobenius structure \eqref{eq: eta-cohomology-stack} on $\rshv(X)$ restricts to a Frobenius structure along $\shv(X)\stackrel{\Psi^L}{\hookrightarrow}\rshv(X)$, which in turn induces a self-duality
\[
\verd_X^\eta: \shv(X)^\vee\cong \shv(X).
\]
\end{proposition}
\begin{proof}
The proposition is equivalent to saying that the equivalence \eqref{eq: Verdier duality for quasi-placid} restricts to an equivalence
\begin{equation*}\label{eq: Verdier duality compact for very-placid}
(\verd_X^{\eta})^\cpt \colon (\shv(X,\La)^\cpt)^\op \simeq (\cshv(X,\La))^{\cpt}.
\end{equation*}

We write $f: U\to X=U/H$.
By \Cref{compact.generation.admissible.stacks} that $\shv(X)^\cpt$ is spanned by objects of the form $f_!f^!\mF$ for $\mF\in \cshv(X)$. If we choose a generalized constant sheaf $\La^\eta_X$, and let $\La^\phi_{U}=f^!\La_X^{\eta}$ be the generalized constant sheaf on $U$. 

Recall that $f^!$ admits a continuous right adjoint $f_\flat$ by \Cref{lem:etale-proper-functoriality-shv-prestack} \eqref{lem:etale-proper-functoriality-shv-prestack-4}.
In addition, as $f$ is ess. coh. pro-unipotent and coh. pro-smooth, $f_\flat$ preserves constructibility (using \Cref{functors.preservating.cons.categories.pl.stack} \eqref{functors.preservating.cons.categories.pl.stack-2} and the fact that $!$-pullback along a coh. pro-unipotent morphism is fully faithful). Therefore, $f_\flat|_{\cshv(U)}=f_{\renflat}|_{\cshv(U)}$.
Then using \Cref{lem-verdier-conjugate-!-push}, it is enough to show that $f_\flat f^!\mF$ is compact for every $\mF\in \cshv(X)$.  Now \Cref{lem: monodromic flat shrek pushforward} below implies that $f_\flat f^!\mF$ is isomorphic to $\varphi_!\varphi^!\mF$ up to a shift.  The proposition is proved.
\end{proof}

\begin{lemma}\label{lem: monodromic flat shrek pushforward} 
Let $X=U/H$ be as in \Cref{prop: verdier duality on shv for very placid}. Then there exists some integer $d$ and an isomorphism of functors $f_\flat f^!\simeq f_!f^![d]: \cshv(X)\to \cshv(X)$.
\end{lemma}
We refer to \Cref{lem: functors between monodromic categories-5} for a generalization of this result when $H$ is a torus.
\begin{proof}
First we assume that $f=\pr: H\times X\to X$ is the projection to the second factor, and $H$ is pfp, and $X$ is a standard placid space.
We notice that 
\[
\pr_?\pr^!\mF=\pr_?(\omega_{H})\boxtimes \mF
\] 
for $?=\flat$ and $!$ (see \Cref{lem:etale-proper-functoriality-shv-prestack} \eqref{lem:etale-proper-functoriality-shv-prestack-4} and  \Cref{lem: ULA-ish for constructible on quasi-placid}). Therefore \Cref{prop: Frob str on coh of group} below implies that there is an isomorphism $\pr_\flat\pr^!\mF\cong \pr_!\pr^!\mF[d]$.

Next still assume that $f=\pr: H\times X\to X$ with $X$ a standard placid space but with $H$ is general. Let $\pr_0: H\times X\to H'\times X$ be the base change of $H\to H'$.
As $H_0$ is coh. pro-unipotent, we have $(\pr_0)_\flat (\pr_0)^!=(\pr_0)_!(\pr_0)^!=\id$ and we reduce to the previous case.

Note that it follows from the above proof that once we fix an isomorphism $(\pi_H)_\flat\omega_H\cong (\pi_H)_!\omega_H[d]$ (where $\pi_H: H\to \spec k$ is the structural map), then the above isomorphism is functorial in $X$.

Now let $f: U\to X$ be an $H$-torsor and $f^\bullet: U^\bullet\to X$ the associated \v{C}ech nerve. Then we have the Cartesian diagram
\[
\xymatrix{
H\times U^\bullet \ar_{\pr}[d]\ar^-{d^\bullet }[r] & U\ar^f[d]\\
U^\bullet\ar^{f^\bullet}[r] & X.
}
\]
Then for $\mF\in \cshv(X)$, we have
\begin{align*}
f_\flat f^!\mF= | (f^\bullet)_!(f^\bullet)^! f_\flat f^!\mF|= |(f^\bullet)_!\pr_\flat (d^\bullet)^!f^! \mF|= |(f^\bullet)_!\pr_\flat \pr^! (f^\bullet)^! \mF|\\
f_! f^!\mF= | (f^\bullet)_!(f^\bullet)^! f_! f^!\mF|= |(f^\bullet)_!\pr_! (d^\bullet)^!f^! \mF|= |(f^\bullet)_!\pr_! \pr^! (f^\bullet)^! \mF|
\end{align*}
We thus reduce the general case to the special case considered before.
\end{proof}

\begin{lemma}\label{prop: Frob str on coh of group}
Let $H$ be a (perfect) connected affine algebraic group over an algebraically closed field $k$. 
There exists some $d\in \bZ$ and a ``trace map"
$C^\bullet(H,\La)[d]\to \La$ such that the pairing $C^\bullet(H,\La)[d]\otimes C^\bullet(H,\La)\to C^\bullet(H,\La)[d]\to \La$ induces an isomorphism
\[
C^\bullet(H,\La)[d]\cong C^\bullet(H,\La)^\vee\cong C_c^\bullet(H,\omega_H).
\]
\end{lemma}
In fact, one can drop the affineness in the assumption. On the other hand,
the proof presented below is roundabout. It would be good to know what is really going on here.
\begin{proof}
If $k$ is a field of characteristic zero, we may embed $\sigma: k\to \bC$ and
let $K\subset H(\bC)=:\sigma H$ be its maximal compact subgroup, which we recall is homotopic to $\sigma H$. We let $d=\dim_{\bR} K$. 
By the standard \'etale-Betti comparison, we have
\[
C^\bullet(H, \La)\cong C_{\mathrm{Betti}}^{\bullet}(\sigma H,\La)\cong C_{\mathrm{Betti}}^\bullet( K,\La)
\]
The Poincar\'e duality for compact (real) manifolds says that there is a canonical trace map
\[
\tr: C_{\mathrm{Betti}}^\bullet( K,\La)[d]\to \La
\]
that induces a perfect pairing
\[
C_{\mathrm{Betti}}^\bullet( K,\La)[d]\otimes_\La C_{\mathrm{Betti}}^\bullet( K,\La)\to C_{\mathrm{Betti}}^\bullet( K,\La)[d]\xrightarrow{\tr} \La.
\]
This gives the desired ``trace map" on $C^\bullet(H,\La)[d]$ and canonical isomorphism as desired.

Next we assume that $k$ is a field of characteristic $p>0$. We may choose a deperfection $H'$ of $H$ so $H'$ is a smooth affine algebraic group over $k$.
First, if $H$ is reductive,
we may lift $H$ to a split reductive group scheme $\mH$ over $W(k)$. The existence of smooth projective compactification $\overline{\mH}$ of $\mH$ with simple normal crossing boundary divisor implies that 
$C^\bullet(H,\La)\cong C^\bullet(\mH_{W(k)_\bQ},\La)$. If $H$ is affine, then we may write $H$ as the extension of its reductive quotient $H_{\mathrm{red}}$ by its unipotent radical $R_uH$, which is an affine space. Then this case follows from the reductive case. 

Here is a second proof, without relying on Betti cohomology. Let $B\subset H$ be a Borel subgroup. Then $H/B$ is proper. Running the argument of \Cref{lem: monodromic flat shrek pushforward} for $B$, we reduce to prove
the lemma just for $H$ being a connected solvable group. In this case, it is an extension of a torus by a unipotent group. But this case can be dealt with easily.
\end{proof}

\begin{remark}\label{rem:renormalized global section}
We will use 
\[
\rg^\eta(X,-): \shv(X,\La)\to \Mod_\La
\] 
to denote the composed functor $\shv(X,\La)\xrightarrow{\Psi^L}\rshv(X,\La)\xrightarrow{\rg^\eta_{\ind\fg}(X,-)}\Mod_\La$.
Precisely, this functor is obtained by first restricting \eqref{eq: eta-cohomology-stack} to $\shv(X,\La)^\cpt$ followed by ind-extension. Beware that for $\mF\in\cshv(X,\La)$, the $\La$-module
$\rg^\eta(X,\mF)$ is in general different from $\rg^\eta_{\ind\fg}(X,\mF)$.
\end{remark}

\begin{corollary}\label{cor: shrek pull and star push preserving compact for quotient}
Let $f: Y\to X=U/H$ be a representable pfp morphism with $X$ as in \Cref{prop: verdier duality on shv for very placid}. Then both $f^!:\shv(X,\La)\to\shv(Y,\La)$ and $f_*: \shv(Y,\La)\to \shv(X,\La)$ preserve compact objects.
\end{corollary}
\begin{proof}
Let $\varphi_U:U\to X$ be the projection.
Let $\varphi_V: V\to Y$ be the $H$-torsor given by $f$, and let $\tilde{f}:V\to U$ be the base change of $f$. Then it is a very placid atlas of $Y$. Now we know that $\shv(X,\La)^\cpt$ is generated by $(\varphi_U)_\flat \mF$ for $\mF\in\cshv(U,\La)$. By base change, we see that $f^!(\varphi_U)_\flat\mF\cong (\varphi_V)_\flat \tilde{f}^!\mF$ is compact in $\shv(X,\La)$.

Similarly, for $\mG\in \cshv(V,\La)$, we have $f_*((\varphi_V)_\flat \mG)\cong (\varphi_U)_\flat((\tilde{f}_*\mG))$ is compact.
\end{proof}

\subsection{Ind-placid and sind-placid stacks}\label{SS: ind and sifted placid}
Finally, we can study the sheaf theory on a class of geometric objects needed for this work, which we call \textit{sind-placid stacks}. Roughly speaking, the category of sind-placid stacks is the category of \'{e}tale stacks obtained by taking placid stacks over $k$ and adding (certain) filtered colimits and geometric realization with transition maps being pfp-proper.


\subsubsection{Ind-(quasi-)placid and sind-placid stacks}

\begin{definition}\label{def-ind-placid-stack-space}
A prestack $X$ is called an \textit{ind-quasi-placid} stack if it admits a presentation (as a prestack) as a filtered colimit $X = \colim_{i\in \mI} X_i$ of quasi-placid stacks (taking values in $1$-groupoids) with transition maps $X_i \rightarrow X_j$ being pfp closed embeddings. It is called ind-(very) placid (resp. quasi-compact) if each $X_i$ in the above representation is (very) placid (resp. quasi-compact). 
We let $\indvplstk_k\subset \indplstk_k\subset\indqplstk_k\subset \prestk_k$ denote the corresponding full subcategories.
\end{definition}

Again, what we just defined probably should only be called quasi-separated ind-(quasi-)placid stacks, which is general enough for our applications.
As filtered colimits commute with finite limits, we see that ind-quasi-placid stacks are hypercomplete \'{e}tale sheaves on $\perf_k$. Moreover, the following holds.

\begin{lemma}\label{ind-stacks.sheaves.and.colimit.from.algebraic.space} 
Let $X$ be an ind-quasi-placid stack with a presentation $\colim_{i\in \mI} X_{i}$. 
\begin{enumerate}
    \item\label{item-ind.maps.from.space.to.ind-placid.stack} For every $S\in \qcsp_k$, $\colim_i X_i(S) \xrightarrow{\cong} X(S)$. 
    \item \label{item-closed.substack.ind-placid} For every quasi-compact quasi-placid stack $S$,  $\colim_i X_i(S) \xrightarrow{\cong} X(S)$.
\end{enumerate}
\end{lemma}
\begin{proof}
Let $S\in \qcsp_k$. We can use an \'etale cover of it by finitely many affine schemes and the sheaf properties of $X$ (and $X_i$) to get $X(S)=\colim_i X_i(S)$. This gives \eqref{item-ind.maps.from.space.to.ind-placid.stack}.
For \eqref{item-closed.substack.ind-placid}, let $V \rightarrow S$ be a quasi-placid atlas with $V\in\plsp_k$. Then the composites $V\rightarrow S \rightarrow X$ and $V\times_{S} V \rightarrow X$ (for both boundary maps) factor through some $X_i$. It follows that $S$ is contained in $X_i$.  
\end{proof}

We will need the following technical result to define sifted-placid stacks.

\begin{lemma}\label{lem-indfp-to-placid-stacks}
Let $f\colon X\to Y$ be an ind-pfp (resp. ind-pfp proper) morphism (in the sense of \Cref{def:ind-fp.and.proper.morphism}) of \'etale stacks with $Y$ being a quasi-compact (very) placid stack. 
Then $X$ is an ind-(very) placid stack with a presentation $X=\colim_\al X_\al$ such that each $X_\al\to Y$ is representable pfp (resp. representable pfp proper). 
\end{lemma}

To understand the content of this lemma, we refer to \Cref{ind.fp.for.ind.spaces} \eqref{ind.fp.for.ind.spaces-2}. 

\begin{proof}
Let $V\to Y$ be a placid atlas with $V\in \plsp_k$ and let $V_\bullet\to Y$ be its \v{C}ech nerve. Then each $V_n\in \plsp_k$.
Let $U_\bullet \to X$ the base change of $V_\bullet\to Y$ along $f$, which is isomorphic to the \v{C}ech nerve of $U:=U_0\to X$.
We first show that there is a presentation $U=\colim U_\al$ with each $U_\al\in\qcsp_k$ that is  pfp  over $V$ (resp. pfp and proper over $V$ if $f$ is ind-proper) and with transition maps being pfp closed embeddings, such that each $U_\al$ is an invariant subspace of the groupoid $U_2=U\times_XU\rightrightarrows U$ (that is, $U_\al\times_XU=U\times_XU_\al$).

First by definition of ind-finitely presented (resp. ind-proper) morphisms, there exists such a presentation $U=\colim U_\al$ as above except $U_\al$ may not be invariant.
Note that $U_\al\times_XU\cong U_\al\times_VV_2$ is a qcqs algebraic space pfp over $V_2$ (where $\pr_1\colon V_2=V\times_YV\to V$ is the first projection).  As $U_\al\times_XU$ is qcqs,
the second projection $U_\al\times_XU\to U$ factors through $U_\al\times_XU\to U_\beta\subset U$ for some $\beta$.  Let $U'_\al\subset U_\beta$ be the Zariski closure of the image of the map $U_\al\times_XU\to U_\beta$, endowed with the closed subspace structure (see \Cref{rem-closed-subscheme}). It contains $U_\al$. 

By \Cref{lem-base-change-schematic-image} below, $U'_\al\times_XU$ is the closure of the projection $U_\al\times_XU\times_XU\to U\times_XU$. It then follows that $U'_\al$
is an invariant subspace (e.g. see \cite[\href{https://stacks.math.columbia.edu/tag/044G}{Lemma 044G}]{stacks-project} for an argument).

We claim that $U'_\al$ is pfp over $V$. Notice $U_\al\times_XU\to U_\beta$ is over the second projection $\pr_2\colon V_2\to V$, which we recall is surjective \emph{strongly} coh. pro-smooth. Write it as $V_2=\lim_\la W_\la\to V$ with each $W_\la \to V$ surjective coh. smooth and transition maps surjective affine coh. smooth.
By \Cref{prop:appr-fp-morphism} \eqref{prop:appr-fp-morphism-1}, $U_\al\times_XU\to V_2$ is the base change of some pfp morphism $\widetilde{W}_\la\to W_\la$. Write $\widetilde{W}_{\la'}$ for the base change of $\widetilde{W}_\la$ along $W_{\la'}\to W_\la$. Then $U_\al\times_XU=\lim_{\la'}\widetilde{W}_{\la'}$ and each $U_\al\times_XU\to \widetilde{W}_{\la'}$ is surjective. As $U_\beta$ is pfp over $V$, the map $U_\al\times_XU\to U_\beta$ factors through some map $\widetilde{W}_{\la'}\to U_\beta$ over $W_{\la'}\to V$.
So $U'_{\al}$ is the closure of the image of a \emph{pfp} map $\widetilde{W}_{\la'}\to U_\beta$. As $U_\beta$ is pfp over $V$, and therefore is placid, again by \Cref{lem-base-change-schematic-image} below, it arises as a base change of the closure of the image of a morphism in $\pfpsp_k$ and therefore $U'_\al\subset U_\beta$ is a pfp closed embedding. This proves the claim.

So we can write a presentation $U=\colim_\al U_\al$ with each $U_\al$ closed invariant subspace of $U$ and pfp (resp. pfp and proper if $f$ is ind-proper) over $V$. We let $X_\al\subset X$ be the sub-prestack whose $R$-points form the full subgroupoid consisting of those $x\colon \Spec R\to X$ such that the base change map $\Spec R\times_XU\cong \Spec R\times_YV\to U$ factors through $\Spec R\times_XU\to U_\al\subset U$. One checks easily that $X_\al\subset X$ is an \'etale stack.
By invariance of $U_\al$, the natural map $U_\al\to U\to X$ factors through $U_\al\to X_\al$ and the resulting morphism $U_\al\to X_\al\times_XU\cong X_\al\times_YV$ is an isomorphism.
It follows that $X_\al$ is placid with $U_\al\to X_\al$ being a placid atlas (by \Cref{prop:placid.stack.placid.map}). 
In addition, $X_\al\to X$ is a pfp closed embedding. Namely, if $\Spec R\to X$ is a morphism, then $U\times_X \Spec R=V\times_Y\Spec R\to \Spec R$ is surjective strongly coh. pro-smooth and $U_\al\subset U$ is an invariant subspace, we see that the image of  $U_\al\times_X\Spec R\to \Spec R$ is closed, denoted by $Z$, with the complement given by the image of $(U-U_\al)\times_X\Spec R\to\Spec R$. In addition, as $(U-U_\al)\times_X\Spec R$ is quasi-compact open in $U\times_X\Spec R$, we see that $\Spec R-Z$ is quasi-compact open in $\Spec R$. This shows that $Z\to\Spec R$ is pfp closed embedding. Therefore, $X_\al\subset X$ is a pfp closed embedding.
It follows that $X=\colim_\al X_\al$ is a filtered colimit of placid stacks with transition maps pfp closed embeddings. 

Finally, clearly if $Y$ is very placid, so is every $X_\al$.
\end{proof}

\begin{lemma}\label{lem-base-change-schematic-image}
Consider a cartesian diagram \eqref{eq:pullback-square-of-algebraic-spaces} in $\qcsp_k$ with $g: Y'\to Y$ strongly coh. pro-smooth. Then for every closed subspace $Z\subset X$, with $Z':=(g')^{-1}(Z)\subset X'$, we have $g^{-1}(\overline{f(Z)})= \overline{f'(Z')}$.
\end{lemma}
\begin{proof}
The lemma  is purely topological. It is enough to show that for quasi-compact open $U'\subset Y'$ with $f'(Z')\cap U'=\emptyset$, then $g^{-1}(\overline{f(Z)})\cap U'=\emptyset$.

If $Y'\to Y$ is coh. smooth, then it is open. Then $f'(Z')\cap U'=\emptyset\Leftrightarrow (f')^{-1}(U')\cap (g')^{-1}(Z)=\emptyset\Leftrightarrow g(U')\cap f(Z)=\emptyset\Leftrightarrow g(U')\cap \overline{f(Z)}=\emptyset$ (as $g(U')$ in open in $Y$), if and only if $U'\cap g^{-1}(\overline{f(Z)})=\emptyset$.

In general, we may write $Y'=\lim Y_i\to Y$ with $g_i:Y_i\to Y$ coh. smooth and transition maps affine coh. smooth. 
We may assume that $U'$ is the pullback of some quasi-compact open subspace in some $Y_i$. In addition, using the above special case, we may assume that  $Y_i=Y$ so $U'=g^{-1}(U)$ for some quasi-compact open $U\subset Y$. Let $g'_i: X_i\to X$ be the base change of $Y_i\to Y$ along $f$, $Z_i=(g'_i)^{-1}(Z)$, $f_i:X_i\to Y_i$ is the base change of $f$ along $g_i$, and $U_i=g_i^{-1}(U)$. Then $(f')^{-1}(U')\cap Z'=\lim_i (f_i)^{-1}(U_i)\cap Z_i$ with transition maps affine coh. smooth. By \cite[\href{https://stacks.math.columbia.edu/tag/0A2W}{Lemma 0A2W}]{stacks-project} (and \cite[\href{https://stacks.math.columbia.edu/tag/0A4G}{Lemma 0A4G}]{stacks-project}), we see that there is some $i$ such that $Z_i\cap f_i^{-1}(U_i)=\emptyset$. Then $U_i\cap g_i^{-1}(\overline{f(Z)})=\emptyset$ by the above special case. So $U'\cap g^{-1}(\overline{f(Z)})=\emptyset$.
\end{proof}

\begin{definition}\label{def-sifted-placid-stack}
A prestack $X$ is called a sind-(very) placid stack if it is an \'etale stack that admits a surjective ind-pfp proper morphism $V\to X$ from an ind-(very) placid stack $V$. We say such $V\to X$ an ind-atlas of $X$. A sind-placid stack is called quasi-compact if there is an ind-atlas $V\to X$ with $V$ a quasi-compact ind-placid stack.
We let $\sfvplstk_k\subset\sfplstk_k\subset \prestk_k$ denote the full subcategory of sifted-very placid and sifted placid stacks. 
\end{definition}

\begin{remark}
Let $X$ be a quasi-compact sifted-placid stack and $V\to X$ an ind-atlas with $V$ quasi-compact. Let $V_\bullet$ be its \v{C}ech nerve. Then by \Cref{lem-indfp-to-placid-stacks}, all the terms in the \v{C}ech nerve $V_\bullet$ are ind-placid and all boundary maps are ind-pfp proper.
Informally, $X$ can be thought as the combination of filtered colimits and geometric realizations (and therefore sifted-collimits) of placid stacks , thus our choice of terminology.
\end{remark}

\begin{corollary}
Let $Y$ be a quasi-compact sifted (very-)placid stack over $k$, and let $f: X\to Y$ be an ind-pfp morphism. Then $X$ is quasi-compact sind-(very-)placid.
\end{corollary}
\begin{proof}This follows from \Cref{lem-indfp-to-placid-stacks} and \Cref{prop:placid.stack.placid.map}.
\end{proof}

\subsubsection{Ind-finitely generated sheaves on ind-placid stacks}
Recall that constructible sheaves are defined on any prestacks.  Let $X$ be an ind-quasi-placid stack with a presentation $X \simeq \colim_{i\in \mI} X_i $. Then \Cref{colimit.sheaves.pseudo-proper.stacks} gives a functor
$\colim_{i\in \mI} \cshv(X_i) \rightarrow \cshv(X)$
with transition functors being $*$-pushforward.
 As $\mI$ is filtered, and each $\cshv(X_i)\to \cshv(X_j)$  is fully faithful (by \Cref{lem:sheaves-open-closed-prestacks}), we see that the functor is fully faithful. However, it is not essentially surjective. For example $\omega_X$ from  \Cref{ex-dualizing sheaf} belongs to $\cshv(X)$ but does not belong to $\colim_{i\in \mI} \cshv(X_i)$. 
In practice, it is $\colim_{i\in \mI} \cshv(X_i)$ rather than $\cshv(X)$ that is more important for applications. This motivates us to make the following definition.

\begin{definition}\label{def:ind-finite sheaf on ind-placid}
Let $X$ be an ind-quasi-placid stack. An object of $\shv(X,\La)$ is called \textit{finitely generated} if it is of the form $i_{*}(\mF)$ for a pfp closed quasi-placid stack $i\colon Z \rightarrow X$ and $\mF\in \cshv(Z,\La)$. We denote the subcategory of finitely generated sheaves by $\fgshv(X,\La)$.  We denote by $\rshv(X,\La)$ its ind-completion. Let $\Psi: \rshv(X)\to \shv(X)$ be the ind-extension of the natural embedding $\fgshv(X)\subset \shv(X)$. 
\end{definition}

\begin{lemma}\label{cons.sheaves.as.colimit.ind.placid.stack}
Let $X$ be an ind-placid stack with a presentation $X \simeq \colim_{i\in \mI} X_i $. The the $*$-pushforward functors induce a canonical equivalence
\[
\colim_{i\in \mI} \cshv(X_i) \xrightarrow{\sim} \fgshv(X)
\]
\end{lemma}
\begin{proof}
We have seen the fully faithfulness.
Essential surjectivity follows, as any pfp closed qcqs placid $Z\subseteq X$ factors through some $X_i$ by \Cref{ind-stacks.sheaves.and.colimit.from.algebraic.space}\ref{item-closed.substack.ind-placid}.
\end{proof}

\begin{remark}
What we call finitely generated sheaves here on an ind-quasi-placid stack are usually called constructible sheaves in literature (e.g. in \cite{arkhipov2009perverse, bezrukavnikov2016two}). If $\shv(X_i)^\omega=\cshv(X_i)$ for every $i\in \mI$ (e.g. each $X_i$ is a perfect qcqs algebraic space), then $\shv(X)^\omega=\fgshv(X)$ and $\shv(X)=\rshv(X)$.
\end{remark}

\begin{example}\label{ex: indfgdualizing sheaf}
For an ind-quasi-placid stack, we can define $\consdual_X^{\ind\fg}=\colim_i \omega_{X_i}\in \rshv(X)$.  Then $\consdual_X$ is the image of $\consdual_X^{\ind\fg}$ under $\Psi$.
\end{example}

Most of the results of \Cref{sec-constructible.sheaves.on.placid.stacks} and \Cref{sec-renormalized.category.placid.stack} generalize to the ind-placid case as well. We summarize them into the following theorem. For technique reasons, we restrict to quasi-compact ind-quasi-placid stacks. 

\begin{theorem}\label{thm-rshv-for-ind-quasi-placid-stack}
The sheaf theory \eqref{eq-rshv-for-placid} (restricted to $\corr(\qcqplstk_k)_{\verti_c;\all}$) admits a canonical extension
\[
\rshv\colon \corr(\qcindqplstk_k)_{\ind\verti_c;\all}\to \lincat_\La,\quad X\mapsto \rshv(X,\La).
\]
sending $X\xleftarrow{g} Z\xrightarrow{f} Y$ to $f^{\ind\fg}_{*}\circ g^{\ind\fg,!}$. We have the following additional properties.
\begin{enumerate}
\item\label{thm-rshv-for-ind-quasi-placid-stack-0} We have $\Psi\circ f^{\ind\fg}_{*} \circ g^{\ind\fg,!} \cong  f_*\circ g^!\circ \Psi$.

\item\label{thm-rshv-for-ind-quasi-placid-stack-3}  
If  $f:X\to Y$ is representable by quasi-placid stacks (i.e. for every $S\to Y$ with $S$ quasi-compact quasi-placid, $S\times_YX$ is quasi-compact quasi-placid), then $f^{\ind\fg,!}$ preserves finitely generated objects.

If $f$ is representable ess. coh. pro-unipotent or is ind-pfp, then $f^{\ind\fg,!}$ admits a left adjoint (which necessarily preserves finitely generated objects). 
In this case, suppose $f$ fits in to a pullback square of quasi-compact ind-quasi-placid stacks as in \eqref{eq:pullback-square-of-algebraic-spaces} with $g: Y'\to Y$ being weakly coh. pro-smooth. Then there is  the base change isomorphisms 
\[
(f')^{\ind\fg}_!\circ(g')^{\ind\fg,!}\xrightarrow{\cong} g^{\ind\fg,!}\circ f^{\ind\fg}_!.
\]
If $f$ is ind-pfp proper, then $f^{\ind\fg}_{!}=f^{\ind\fg}_{*}$.

\item\label{thm-rshv-for-ind-quasi-placid-stack-2} For $f$ being ind-pfp, the functor $f^{\ind\fg}_{*}$  preserves finitely generated objects. If $f$ is in addition representable pfp, then $f^{\ind\fg}_{*}$ admits a left adjoint  $f^{\ind\fg,*}$ (which necessarily preserve finitely generated objects). In this case, suppose $f$ fits in to a pullback square of quasi-compact ind-quasi-placid stacks as in \eqref{eq:pullback-square-of-algebraic-spaces} with $g: Y'\to Y$ being weakly coh. pro-smooth. Then there is  the base change isomorphisms 
\[
(f')^{\ind\fg,*}\circ g^{\ind\fg,!}\xrightarrow{\cong} (g')^{\ind\fg,!}\circ f^{\ind\fg,*}. 
\]

 \item\label{thm-rshv-for-ind-quasi-placid-stack-4} If $f: X\to Y$ is a weakly coh. pro-smooth morphism between quasi-compact ind-quasi-placid stacks, then $f^{\ind\fg,!}$ admits continuous right adjoint $f_{\renflat}$ satisfies the projection formula
\[
f_\renflat(\mF)\os \mG \xrightarrow{\sim} f_\renflat(\mF\os f^{\ind\fg,!}(\mG)),\quad \mF\in \rshv(X,\La),\mG\in \rshv(Y,\La).
\]

\item\label{thm-rshv-for-ind-quasi-placid-stack-5} 
 Let \eqref{eq:pullback-square-of-algebraic-spaces} be a pullback square of prestacks with $f\colon X\rightarrow Y$ being a weakly coh. pro-smooth morphism between quasi-compact ind-quasi-placid stacks. Then if $g$ is representable ess. coh. pro-unipotent or ind-pfp, there is the natural base change isomorphism of functors 
 \[
  g^{\ind\fg,!}\circ f_{\renflat}\xrightarrow{\cong}  (f')_{\renflat}\circ (g')^{\ind\fg,!}.
  \]

     If addition $g$ is ind-pfp, then there is the natural base change isomorphism 
     \[
     g^{\ind\fg}_{*}\circ (f')_\renflat\xrightarrow{\cong} f_\renflat\circ (g')^{\ind\fg}_{*}.
     \]

\item\label{thm-rshv-for-ind-quasi-placid-stack-6}  Let \eqref{eq:pullback-square-of-algebraic-spaces} be a pullback square of quasi-compact ind-placid stacks with $f\colon X\rightarrow Y$ being weakly coh. pro-smooth. Then $ g^{\ind\fg,!}\circ f_{\renflat}\xrightarrow{\cong}  (f')_{\renflat}\circ (g')^{\ind\fg,!}$.

\end{enumerate}
\end{theorem}
Note that $!$-pullbacks between ind-quasi-placid stacks do not preserve finitely generated sheaves in general.
\begin{proof}
The existence of the extension $\rshv$ follows from the same reasoning as in \Cref{def-dual-sheaves-on-prestacks-vert-indfp}, applying \Cref{prop-sheaf-theory-non-full-right-Kan-extension} to \eqref{eq-rshv-for-placid}. Namely, we still let $\mathrm{S}_1$ be the class of pfp closed embeddings. Let $\verti_1$ be the class of pfp morphisms, and $\verti_2$ be the class of $\mathrm{Indpfp}$ morphisms. Then we can apply \Cref{prop-sheaf-theory-non-full-right-Kan-extension}. 
(We note that the proof of \Cref{prop-sheaf-theory-non-full-right-Kan-extension} does not require the existence of fiber product in $\bfC_1\subset \bfC_2$, the weaker assumption as in \Cref{rem-generalization-corr} \eqref{rem-generalization-corr-2} suffices.)
Note that by construction, for $X$ quasi-compact ind-quasi-placid, we have
\[
\rshv(X)=\lim_{Y\to X} \rshv(Y),
\]
where $Y$ is quasi-compact quasi-placid. By \Cref{ind-stacks.sheaves.and.colimit.from.algebraic.space} \eqref{item-closed.substack.ind-placid} we may replace the index category by those pfp closed embedding $Y\to X$ and $Y$ quasi-compact quasi-placid. By Part \eqref{thm-rshv-for-ind-quasi-placid-stack-3} below, we see that $\lim_{Y\subset X} \rshv(Y)=\colim_{Y\subset X}\rshv(Y)$. So the value of $\rshv(X)$ is indeed the one from \Cref{def:ind-finite sheaf on ind-placid}. 
In addition Part \eqref{thm-rshv-for-ind-quasi-placid-stack-0} follows from definition. 

It follows from \Cref{functors.preservating.cons.categories.pl.stack} \eqref{functors.preservating.cons.categories.pl.stack-1} that when $f$ is representable ess. coh. pro-unipotent then functor $f^{\ind\fg,!}$ preserves finitely generated objects and admits a left adjoint $f^{\ind\fg}_!$. In addition, it satisfies the desired base change isomorphism.
On the other hand, using \Cref{lem-indfp-to-placid-stacks}, the statements for $f$ ind-pfp follows from the case $f$ being representable pfp. 
As in \Cref{locally.ind-proper.existence.lower-!}, $f^{\ind\fg}_!=f^{\ind\fg}_*$  if $f$ is ind-pfp proper. 
This proves Part \eqref{thm-rshv-for-ind-quasi-placid-stack-3},

Part \eqref{thm-rshv-for-ind-quasi-placid-stack-2} for $f$ being representable pfp follows immediately from \Cref{functors.preservating.cons.categories.pl.stack} \eqref{functors.preservating.cons.categories.pl.stack-3}.
Again using \Cref{lem-indfp-to-placid-stacks}, it implies that $f^{\ind\fg}_*$ always preserves finitely generated objects. This proves Part \eqref{thm-rshv-for-ind-quasi-placid-stack-2}. 

Using \Cref{prop-base-change-for-renflat-forward}, we see that $f^{\ind\fg,!}$ admits continuous right adjoint $f_{\renflat}$ if $f$ is weakly coh. pro-smooth. 
In addition,
 \Cref{ren.flat.projection.formula} implies the projection formula in this more general setting, giving Part \eqref{thm-rshv-for-ind-quasi-placid-stack-4}. 
Part \eqref{thm-rshv-for-ind-quasi-placid-stack-5}
also follows from  \Cref{prop-base-change-for-renflat-forward} and  the previous discussions.

Part \eqref{thm-rshv-for-ind-quasi-placid-stack-6} follows from \Cref{prop-base-change-for-renflat-forward-difficult}.
\end{proof}

\Cref{pfp.functors.and.bc.placid.stacks} also admits a generalization. However, we defer the discussion until \Cref{pfp.functors.and.bc.ind-placid.stacks}.

\subsubsection{Verdier duality and perverse sheaves on ind-placid stacks}\label{sec-verdier.duality.perverse.sheaf.ind.placid.stacks}
We now discuss Verdier duality and perverse sheaves  for ind-quasi-placid stacks following the discussion of \Cref{sec-verdier.duality.placid.stacks}. 

\begin{definition}\label{def.dimension.constant.sheaf.ind.placid.stack}
Let $X$ be an ind-placid stack. A generalized constant sheaf $\La_X^\eta$ of $X$ is a rule to assign every pfp closed embedding $Z\subseteq X$ a generalized constant sheaf $\La_Z^\eta$ on $Z$ and for every inclusion $\iota\colon Z \subseteq Z'$ an isomorphism $\iota^*\La_{Z'}^{\eta}\cong \La_Z^{\eta}$ satisfying natural compatibility conditions.
\end{definition}  

\begin{remark}\label{rem: generalized constant sheaf on ind-placid stacks}
Our terminology is abusive. Namely a generalized constant sheaf $\La_X^\eta$ of $X$ is not really an object in $\shv(X)$. Instead, it is an object in $\lim_{Z\to X} \cshv(Z)$ with transition maps given by $*$-pullbacks. I.e. $\La_X^\eta$ can be regarded as an object in $\scshv(X)$ (see \Cref{shv over sshv} for the definition).
Given a presentation $X = \colim_{i\in \mI} X_i$ of an ind-placid stack, to give a generalized constant sheaf of $X$ it is enough to give a collection of generalized constant sheaf $\{\La^{\eta}_{X_i}\}_{i\in \mI}$ satisfying the usual compatibility conditions as in the ordinary category theory. (As discussed in the proof of \Cref{lem:existence of generalized constant sheaf}, higher compatibilities are not needed.)
\end{remark}

\begin{example}\label{dimension.ind.placid.pfp}
Let $X = \colim_i X_i$ be an ind-finitely presented algebraic space. Then $X$ admits a generalized constant sheaf given by $\{\La_{X_i}\}$.
\end{example}

\begin{example}\label{ex-*-pullback-indfp-generalized-constant}
Let $f: X\to Y$ be an ind-pfp morphism between ind-placid stacks. 
If $\La^\eta_Y$ is a generalized constant sheaf of $Y$, then there is a generalized constant sheaf $\La_X^\phi$ of $X$ obtained from $\La_Y^\eta$ as ``$*$-pullback along $f$". Namely, we may write $X=\colim_{i\in \mI} X_i$ and $Y=\colim_{j\in \mJ} Y_j$ with $X_i, Y_j$ placid and transition maps pfp closed embeddings. 
Then for every $i$, there is $j$ such that $f$ restricts to a rfp morphism $f_{ij}: X_i\to Y_j$. Then the collection $\{(f_{ij})^*\La^\eta_{Y_j}\}$ defines $\La_X^\phi$.
\end{example}

\begin{example}\label{ex-!-pullback-prosm-generalized-constant}
Let $f: X\to Y$ be a weakly coh. pro-smooth morphism between ind-placid stacks and let $\La_Y^\eta$ be a generalized constant sheaf of $Y$. Then there is a generalized constant sheaf $\La_X^\phi$ of $X$ obtained from $\La_Y^\eta$ as ``$!$-pullback along $f$". Namely, if $Y=\colim_{i\in \mI}Y_i$ is a presentable of $Y$ with each $Y_i$ placid and transition maps pfp closed embeddings. Then $X=\colim_{i\in \mI}X_i$ with $X_i=X\times_YY_i$ is a presentation of $X$. Then $\La_X^\phi$ is given by the collection $\{(f_i)^!\La_{Y_i}^\eta\}$. 
\end{example}

Given an ind-placid stack $X$ with a generalized sheaf $\La_X^\eta$, by \Cref{verdier.duality.placid.stacks} we get a compatible system of functors 
\[
\rg_{\ind\fg}^{\eta}(Z,-)\colon \cshv(Z,\La) \rightarrow \La, \quad Z\subseteq X
\]
where $Z$ ranges over pfp closed placid stacks of $X$.  Note that for $\iota: Z\subset Z'$, $\Hom(\La_{Z'}^\eta, \iota_*(-))\cong \Hom(\La_Z^\eta,-)$ so these functors together to give a functor which we denote by
\begin{equation}\label{eq: renormalized cohmology, ind-placid}
\rg^\eta_{\ind\fg}(X,-)\colon \fgshv(X,\La) \rightarrow \Mod_\La.
\end{equation}
The presentation of \Cref{cons.sheaves.as.colimit.ind.placid.stack} immediately implies that the pairing induced by the ind-extension of the pairing
\begin{align*}
   \fgshv(X,\La) \otimes_\La \fgshv(X,\La) \xrightarrow{\os} \fgshv(X,\La) \xrightarrow{\rg_{\ind\fg}^\eta(X,-)} \Mod_\La
\end{align*}
induces an equivalence
\[
(\verd_X^\eta)^{\mathrm{f.g.}}: \fgshv(X,\La)^{\op}\cong \fgshv(X,\La)
\] 
such that
\[
\Hom_{X}((\verd_{X}^\eta)^{\mathrm{f.g.}}(\mF),\mG) \simeq \rg_{\ind\fg}^{\eta}(X,\mF\os \mG), \quad \mF,\mG \in \fgshv(X,\La).
\]
In other words, we obtain the following.
\begin{proposition}\label{prop:verdier.duality.ind.placid.stacks}
The ind-completion of \eqref{eq: renormalized cohmology, ind-placid} is a Frobenius structure on $\rshv(X,\La)$, inducing
\[
(\verd_{X}^\eta)^{\ind\fg}: \rshv(X,\La)^\vee\cong \rshv(X,\La).
\]
\end{proposition}

We also get automatically the same results regarding functoriality for ind-pfp morphisms, provided induced generalized constant sheaves are compatible. 
It follows from \Cref{verdier.duality.placid.stacks} that we have the following statement.
\begin{proposition}\label{star.induced.ind.fin.pres.}
\begin{enumerate}
\item\label{star.induced.ind.fin.pres.-1} Let $f\colon X\rightarrow Y$ be an ind-pfp morphism of ind-placid stacks and let $\La_Y^\eta$ be a generalized constant sheaf of $Y$. Let $\La_X^\phi$ be the generalized constant sheaf of $X$ constructed in \Cref{ex-*-pullback-indfp-generalized-constant}, then 
\[
f_*\circ (\verd_X^{\phi})^{\fg}\simeq (\verd_Y^\eta)^{\fg}\circ f_!\colon \fgshv(X)\to \fgshv(Y).
\]
\item\label{star.induced.ind.fin.pres.-2} Let $f\colon X\rightarrow Y$ be a weakly coh. pro-smooth morphism  of ind-placid stacks and let $\La_Y^\eta$ be a generalized constant sheaf of $Y$. Let $\La_X^\phi$ be the generalized constant sheaf of $X$ constructed in \Cref{ex-!-pullback-prosm-generalized-constant}, then
\[
f^!\circ (\verd_Y^{\eta})^{\fg}\simeq (\verd_X^\phi)^{\fg}\circ f^!\colon \fgshv(Y)\to \fgshv(X).
\]
\end{enumerate}
\end{proposition}

\begin{remark}\label{rem: equiv between cons-sheaf and cons-cosheaf-2} 
For quasi-compact ind-placid stacks, we do not have a direct analogue of  \eqref{eq: identifying star sheaf with shrek sheaf-2}. Namely, the category $\cshv(X)$ and $\scshv(X)$ (as defined in \Cref{shv over sshv}) are very different in general, even with a choice of a generalized constant sheaf $\La_X^\eta$ (which we recall is in fact an object in $\scshv(X)$ as explained in \Cref{rem: generalized constant sheaf on ind-placid stacks}). However, we still have the tensor product 
\begin{equation}\label{eq: eta tensor product on fgshv}
\fgshv(X)\otimes \fgshv(X)\to \fgshv(X),\quad (\mF,\mG)\mapsto (\verd_X^\eta)^{\fg}(((\verd_X^\eta)^{\fg}(\mF))\os ((\verd_X^\eta)^{\fg}(\mG))),
\end{equation}
as  \eqref{eq: eta tensor product on cshv}. Note that this tensor product may not have a unit, but we still have
\begin{equation}\label{eq: Verdier duality for quasi-placid-characterization-2}
\Hom_{\rshv(X,\La)}(\mF,\mG) \simeq \Hom_{\rshv(X,\La)}((\mF)\otimes^\eta (\verd^{\eta}_X)^{\fg}(\mG),\consdual^{\ind\fg}_X), \quad \mF,\mG\in \fgshv(X,\La),
\end{equation}
where $\consdual^{\ind\fg}_X$ is as in \Cref{ex: indfgdualizing sheaf}. 

We will still use $\otimes^\eta$ to denote the ind-extension of \eqref{eq: eta tensor product on fgshv}, which now becomes a (non-unital) monoidal product of $\rshv(X,\La)$.
\end{remark}

\begin{remark}\label{rem: verdier duality on shv for ind very placid}
Suppose that $k$ is an algebraically closed field, let $X$ be an ind-placid stack over $k$ of the forma $X=U/H$, where $U=\colim_i U_i$ is an ind-algebraic space with each $U_i$ standard placid and transitioning maps being pfp closed embeddings, and $H$ is affine group scheme as in \Cref{prop: verdier duality on shv for very placid}. Then we have $\Psi^L: \shv(X)\subset \rshv(X)$ and the statements of \Cref{prop: verdier duality on shv for very placid} hold in this generality.
\end{remark}

One can also define the category of perverse sheaves on ind-placid stacks. First, if $f: Z\to X$ is a pfp closed embedding of placid stacks and $\La_Z^\phi=f^*\La_X^\eta$ are generalized constant sheaves. Then $i_*: \cshv(Z)\to \cshv(X)$ is perverse exact sending $\mathrm{Perv}(Z)^\phi$ fully faithfully to $\mathrm{Perv}(X)^\eta$. Now, let $\La_X^\eta$ be a generalized constant sheaf of an ind-placid stack $X=\colim_i X_i$. Then we have
\[
\mathrm{Perv}(X)^\eta:=\colim_i \mathrm{Perv}(X_i)^{\eta}\subset \fgshv(X).
\]

\begin{remark}\label{prop:perverse exact smooth pullback-2}
In the setting as in \Cref{ex-!-pullback-prosm-generalized-constant}, the functor $f^!$ is perverse exact, which follows from \Cref{prop:perverse exact smooth pullback}.
\end{remark}

\subsubsection{Categories of sheaves on sind-placid stacks}
First, it follows from \Cref{locally.ind-proper.descent} and \Cref{colimit.sheaves.pseudo-proper.stacks} that for a sind-placid stack $X$ with an ind-atlas $V\to X$ (see \Cref{def-sifted-placid-stack}), we have 
\begin{equation}\label{eq-shv-on-sifted-placid}
\shv(X)\cong \tot(\shv(V_\bullet))\cong |\shv(V_\bullet)|,
\end{equation}
where for totalization, the transition functors are upper $!$-pullbacks and for geometric realization, the transition functors are lower $*$-pushforward.

We have the following generalization of  \Cref{pfp.functors.and.bc.placid.stacks}.
\begin{proposition}\label{pfp.functors.and.bc.ind-placid.stacks} 
\Cref{pfp.functors.and.bc.placid.stacks} holds with $Y$ (and $Y'$) being sind-very placid. 
In addition, if $f\colon X\to Y$ is ind-pfp morphism, then $f^!$ admits a left adjoint $f_!$ and satisfying the base change isomorphism for weakly coh. pro-smooth morphism $g\colon Y'\to Y$ of sind-very placid stacks. 
\end{proposition}
\begin{proof}
We first assume that $Y=\colim_i Y_i$ is ind-very placid with each $Y_i$ very placid and transition maps pfp closed embedding.
Let $f\colon X\to Y$ be as in \Cref{pfp.functors.and.bc.placid.stacks}. Then $X=\colim_i X_i$ with $X_i=X\times_YY_i$.  Then $*$-pullback (in the case $f$ is rfp) and $!$-pushforward exist for every $f_i:X_i\to Y_i$ by \Cref{pfp.functors.and.bc.placid.stacks}. As
$*$-pushfowards along pfp proper morphisms (in particular closed embeddings) are left adjoint $!$-pullback by \Cref{lem:etale-proper-functoriality-shv-prestack} \eqref{lem:etale-proper-functoriality-shv-prestack-1}, using \Cref{colimit.sheaves.pseudo-proper.stacks} we see that $*$-pullback (in the case $f$ is rfp) and $!$-pushforward exist for $f$, and they satisfy the desired base change. 
It remains to deal with the case when $f$ is ind-finitely presented. But this follows from \Cref{lem-indfp-to-placid-stacks} and the case $f$ is representable pfp.

Next we assume that $Y$ is sind very-placid, and let $V\to Y$ be an ind-atlas with $V$ ind-very placid. Then we can use \eqref{eq-shv-on-sifted-placid} and repeat the above arguments to conclude.
\end{proof}

We summarize main facts we have established for the sheaf theory \eqref{eq: pushforward along non-representable morphisms}.

\begin{theorem}\label{thm: final theorem of shv}
Assume that $k$ is the perfection of a regular noetherian ring of dimension $\leq 1$ 
and $\ell$ a prime invertible in $k$ such that $k$ has finite $\bF_\ell$-cohomological dimension.  Let $\La$ be a $\bZ_\ell$-algebra as in \Cref{sec:adic-formalism}. 
Consider the sheaf theory \eqref{eq: pushforward along non-representable morphisms}. 
\begin{enumerate}
\item If $X$ is a pfp algebraic space over $k$, then $\shv(X)=\ind\der_{\ctf}(X)$. If $X=\lim_i X_i$ is a qcqs algebraic space over $k$ with each $X_i$ pfp over $k$ and affine transitioning maps, then $\shv(X)=\colim_i \shv(X_i)$ with transitioning functors being $!$-pullbacks. If $X$ is a general prestack, then $\shv(X)=\lim_{S\to X}\shv(S)$ with $S\in \qcsp_k$ and with transitioning functors being $!$-pullbacks.
\item If $f:X\to Y$ is a morphism such that there is an \'etale covering $Y'\to Y$ such that $X\times_YY'\to Y'$ is ind-ess. pro-\'etale, then $f\in\verti$. If $f: X\to Y$ is an $A$-gerbe map, with $A$ a finite abelian group of order invertible in $\La$, then $f\in\verti$.
\item If $f$ is ind-pfp proper, then $f_*$ is the left adjoint of $f^!$, so the class of ind-pfp proper morphisms form a class satisfying  \Cref{assumptions.base.change.sheaf.theory.V}. 
\item The class of representable coh. pro-smooth morphisms form a class satisfying \Cref{assumptions.base.change.sheaf.theory.H}.
\item If $f$ is $\ell$-ULA, then $f_*$ admits a left adjoint $f^*$ satisfying base change isomorphisms with respect to arbitrary $!$-pullbacks. 
\item If $f$ is representable pfp between quasi-compact sind very-placid stacks, then $f_*$ admits a left adjoint $f^*$ satisfying base change isomorphisms with respect to $!$-pullbacks along weakly coh. pro-smooth morphisms. 
\item If $f$ is ess. coh. pro-unipotent morphism or ind-pfp morphism between quasi-compact sind very-placid stacks, then $f^!$ admits a left adjoint $f_!$ satisfying base change isomorphisms with respect to $!$-pullbacks along weakly coh. pro-smooth morphisms. 
\end{enumerate}
\end{theorem}

We also record the following result on open-closed gluing, as an extension of \Cref{lem:sheaves-open-closed-prestacks} for sind-placid stacks, for future reference. 
\begin{proposition}\label{open.closed.recollement}
Let $Y$ be a quasi-compact sind-very placid stack and let $j\colon U\rightarrow Y$ be a qcqs open embedding with a closed complement $i\colon Z\rightarrow Y$. Then $U$ and $Z$ are sifted-very placid stacks and $i$ is a pfp closed embedding. In addition:
\begin{enumerate}
   \item  $i_!=i_*$, $j^!\circ i_*\simeq 0$, $i^!\circ j_*\simeq 0$, and $i^*\circ j_! \simeq 0$.
   \item The functor $i_*$ (resp. $j_*$, resp. $j_!$) is fully faithful, with essential image consisting of those $\mF\in \shv(Y)$ satisfying $j^!\mF=0$  (resp. $i^!\mF=0$, resp. $i^*\mF=0$).     
   \item \label{open.shriek.closed} For every $\mF\in \shv(Y,\La)$ we have canonical fiber    sequences 
        \[
        i_*  i^! \mF \rightarrow \mF \rightarrow j_*  j^!\mF,\quad  j_! j^! \mF \rightarrow \mF \rightarrow i_*  i^*\mF.
        \]
\end{enumerate}
\end{proposition}
\begin{proof}
Let $f: X\to Y$ be an ind-atlas of $Y$ and then write $X=\colim_i X_i$ with $X_i$ being quasi-compact very placid stacks and transitioning maps being pfp closed embeddings. Then the $*$-pushforward along $X_i\to Y$ commutes with all functors $(i^*,i_*, i^!, j_!, j^!,j_*)$. Therefore by \eqref{eq-shv-on-sifted-placid}, we reduce to the corresponding statements for $X_i$.  Then we may choose a placid atlas $Z_i\to X_i$ with $Z_i\to X_i$ ess. coh. pro-unipotent. As the $!$-pullback functor along $Z_i\to X_i$ commutes with all functors $(i^*,i_*, i^!, j_!, j^!,j_*)$, we reduce to the case of standard placid algebraic spaces. We may then further reduce to the case that $Y$ is pfp over $k$. This is then standard.
\end{proof}

\subsubsection{Ind-finitely generated sheaves}
Finally, we define the category of ind-finitely generated sheaves on sind-placid stacks. For technical reasons, we restrict to quasi-compact sind-placid stacks.
We let $\horiz$ be the class of morphisms of $\qcsfplstk_k$ that are representable in $\qcindplstk_k$. That is, a morphism $f: X\to Y$ of quasi-compact sifted placid stacks belongs to $\horiz$ if for every map $Z\to Y$ with $Z$ being a quasi-compact ind-placid stack, the fiber product $Z\times_YX$ exists and is represented as a quasi-compact ind-placid stack.
(This is similar to the definition in \Cref{rem-prop-of-weak-strong-stable-class} \eqref{rem-prop-of-weak-strong-stable-class-2}. But as fiber products may not exist in $\qcsfplstk_k$. The definition given there does not directly apply.)
We will let $(\ind\verti_c)_r\subset\horiz$ be the class of those $f$ such that $Z\times_YX\to Z$ belongs to $\ind\verti_c$.
Note that by \Cref{lem-indfp-to-placid-stacks}, the class of ind-pfp morphisms belong to $(\ind\verti_c)_r$. In addition,
the natural inclusion
\begin{equation}\label{eq:ind-placid to sifted-placid}
\corr(\qcindplstk_k)_{\ind\verti_c;\all}\subset \corr(\qcsfplstk_k)_{(\ind\verti_c)_r;\horiz}
\end{equation}
is fully faithful. Therefore we can define 
\begin{equation}\label{eq:renormalized sheaf on sifted}
\rshv:  \corr(\qcsfplstk_k)_{(\ind\verti_c)_r;\horiz}\to \lincat_\La
\end{equation}
as the left operadic Kan extension along \eqref{eq:ind-placid to sifted-placid} of the restriction to $\corr(\qcindplstk_k)_{(\ind\verti_c)_r;\all}$ of the sheaf theory from \Cref{thm-rshv-for-ind-quasi-placid-stack}.
By \Cref{prop-dual-sheaf-theory-right-Kan-extension}, for $X$ quasi-compact sind-placid, we have
\[
\rshv(X,\La)=\colim_{Y\to X} \rshv(Y,\La)
\] 
with $Y$ placid and $Y\to X$ ind-pfp. In particular, it is compactly generated. We let
\[
\fgshv(X,\La):=\rshv(X,\La)^\cpt= \colim_{Y\to X} \cshv(Y,\La),
\]
with $Y$ placid and $Y\to X$ ind-pfp.

Tautologically, there is still the functor 
\begin{equation}\label{eq-ind-constr-to-shv-sifted} 
\Psi\colon \rshv(X,\La)\to \shv(X,\La).
\end{equation} 
Some (but not all) statements \Cref{thm-rshv-for-ind-quasi-placid-stack} continue to hold in this setting without change. We summarize it as follows.

\begin{proposition}\label{prop-rshv-for-sifted-placid-stack} 
Let $X\xleftarrow{g} Y\xrightarrow{f}Z$ be a correspondence between quasi-compact sind-placid stacks, and suppose $g\in \horiz$ and and $f$ is ind-pfp. 
\begin{enumerate}
\item We have $\Psi\circ f^{\ind\fg}_{*} \circ g^{\ind\fg,!} \cong  f_*\circ g^!\circ \Psi$. 
\item The functor $f_*^{\ind\fg}$ preserve finitely generated objects. If $f$ is ind-proper, $f_*^{\ind\fg}$ is the left adjoint of $f^{\ind\fg,!}$. 
\item If $g: Y\to X$ is representable by quasi-compact placid stacks (i.e if $V\to X$ is ind-pfp with $V$ placid, $V\times_XY$ is placid), then $g^{\ind\fg, !}$ preserves finitely generated objects. If $g$ is \'etale, then $g^{\ind\fg}_*$ is the right adjoint of $g^{\ind\fg, !}$. 
\item If $g: Y\to X$ is weakly coh. pro-smooth and representable by quasi-compact placid stacks, then $g^{\ind\fg,!}$ admits a continuous right adjoint $g_\renflat$ satisfying a projection formula, and a base change formula with respect to $(\ind\fg,*)$-pushforward as in \Cref{thm-rshv-for-ind-quasi-placid-stack} \eqref{thm-rshv-for-ind-quasi-placid-stack-5}. If $Y$ and $X$ are sind-very placid, it also satisfies a base change formula with respect to $(\ind\fg,!)$-pullbacks as in  \Cref{thm-rshv-for-ind-quasi-placid-stack} \eqref{thm-rshv-for-ind-quasi-placid-stack-6}.
\item Let $Y$ be a quasi-compact sind-placid stack and let $j\colon U\rightarrow Y$ be a qcqs open embedding with a closed complement $i\colon Z\rightarrow Y$. Then $j^{\ind\fg,!}\circ i^{\ind\fg}_*\simeq 0$, $i^{\ind\fg,!}\circ j^{\ind\fg}_*\simeq 0$.
 In addition, for every $\mF\in \shv(Y,\La)$ we have canonical fiber    sequences 
        \[
        i^{\ind\fg}_*  i^{\ind\fg,!} \mF \rightarrow \mF \rightarrow j^{\ind\fg}_*  j^{\ind\fg,!}\mF.
        \]
 In particular,  the functor $i^{\ind\fg}_*$ (resp. $j^{\ind\fg}_*$) is fully faithful, with essential image consisting of those $\mF\in \shv(Y)$ satisfying $j^{\ind\fg,!}\mF=0$  (resp. $i^{\ind\fg,!}\mF=0$).   
\end{enumerate}
\end{proposition}

\begin{remark}
Unlike \Cref{thm-rshv-for-ind-quasi-placid-stack}, we usually do not have left adjoint functors of $(\ind\fg,*)$-pushforwards and $(\ind\fg,!)$-pullbacks, even in the favorable situations. However, in some special cases, one can prove that such left adjoints exist.
\end{remark}

The functor $\Psi$ in general is far from being equivalence. In fact, unlike the case of ind-placid stacks, even the restriction of $\Psi$ to $\fgshv(X,\La)\to \shv(X,\La)$ may not be fully faithful in general. 

To explain this, consider ind-pfp proper morphisms $f_i:Y_i\to X$ with $Y_i$ quasi-compact placid stacks for $i=1,2$. Let $Z=Y_1\times_XY_2$ with two projections $g_i: Z\to Y_i$. By \Cref{lem-indfp-to-placid-stacks}, $Z=\colim_{j\in J} Z_j$ is quasi-compact ind-placid, with each $g_{ji}: Z_j\to Y_i$ pfp proper. Let $\mF_i\in\cshv(Y_i)$. 
\begin{lemma}\label{lem: hom of fg sheaves on sifted placid}
Notations as above. Then
\begin{equation}\label{eq: hom of fg sheaves on sifted placid} 
\Hom_{\rshv(X)}((f_1)^{\ind\fg}_{*}\mF_1, (f_2)^{\ind\fg}_{*}\mF_2)\cong \colim_j \Hom_{\shv(Y_2)}(\mF_1, (g_{j1})_*(g_{j2})^!\mF_2).
\end{equation}
\end{lemma}
\begin{proof}
Recall that the right adjoint of $(f_1)^{\ind\fg}_{*}$ is $(f_1)^{\ind\fg,!}$. Then by base change, we have
\begin{eqnarray*}
\Hom_{\rshv(X)}((f_1)^{\ind\fg}_{*}\mF_1, (f_2)^{\ind\fg}_{*}\mF_2) &= &\Hom_{\rshv(Y_1)}(\mF_1, (f_1)^{\ind\fg,!}(f_2)^{\ind\fg}_{*}\mF_2)\\
 &=&\Hom_{\rshv(Y_1)}(\mF_1, (g_1)^{\ind\fg}_{*}(g_2)^{\ind\fg,!}\mF_2)\\
 &=& \Hom_{\rshv(Y_1)}(\mF_1, \colim_{j} (g_{j1})^{\ind\fg}_{*}(g_{j2})^{\ind\fg,!}\mF_2).
\end{eqnarray*}
As $\mF_1$ is compact in $\rshv(Y_1)$, we have 
\begin{eqnarray*}
\colim_j \Hom_{\shv(Y_1)}(\mF_1, (g_{j1})_*(g_{j2})^{!}\mF_2) &=& \colim_j \Hom_{\rshv(Y_1)}(\mF_1, (g_{j1})^{\ind\fg}_*(g_{j2})^{\ind\fg,!}\mF_2)\\
 & \cong & \Hom_{\rshv(Y_1)}(\mF_1,  \colim_j (g_{j1})^{\ind\fg}_{*}(g_{j2})^{\ind\fg,!}\mF_2),
\end{eqnarray*}
as desired.
\end{proof}

On the other hand, $\Psi((f_i)^{\ind\fg}_{*}\mF_i)$ is nothing but the one obtained from the usual $*$-pushforward of $\mF_i$ (regarded as an object in $\shv(Y_i)$). Thus the same reasoning implies that
\begin{equation}\label{eq: hom of  sheaves on sifted placid}
\Hom_{\shv(X)}((f_1)_*\mF_1, (f_2)_*\mF_2)=\Hom_{\shv(Y_2)}(\mF_1, \colim_j (g_{j1})_*(g_{j2})^!\mF_2).
\end{equation}
As $\mF_1$ may not be compact in $\shv(Y_2)$, we may not be able to pull the colimit out. Therefore the map 
\[
\Hom_{\rshv(X)}((f_1)^{\ind\fg}_{*}\mF_1, (f_2)^{\ind\fg}_{*}\mF_2)\to \Hom_{\shv(X)}((f_1)_*\mF_1, (f_2)_*\mF_2)
\] 
may not be an isomorphism in general.

On the positive side, using base change, \cite[Corollary 4.7.5.3]{Lurie.higher.algebra} implies the following.
\begin{proposition}\label{prop: fully faithful for ind-proper surjective fg sheaves}
For a sind-placid stack $X$ with an ind-atlas $V\to X$ (see \Cref{def-sifted-placid-stack}), the natural functor
$|\rshv(V^\bullet,\La)|\to \rshv(X,\La)$ is fully faithful.
\end{proposition}

\begin{remark}
Note that unlike \eqref{eq-shv-on-sifted-placid}, the above functor may not be essentially surjective. In fact, this phenomenon already presents for $V=\pt\to \bB G$, where $G$ is a finite group, regarded as a constant group scheme over $k$.  In this case, $\fgshv(\bB G,\La)=\cshv(\bB G,\La)$. Now suppose $\La=\bF_\ell$ with $\ell \mid \sharp G$.
Then the above functor then is the fully faithful (but not essentially surjective) embeddimg $\shv(\bB G,\bF_\ell)\subset \rshv(\bB G,\bF_\ell)$.
\end{remark}

Combining the above discussions with \Cref{bimodule.geom.trace.fully.faithfulness.convolution.cat}, we obtain the following.

\begin{proposition}\label{thm:trace-constructible}
Let $X$ be a quasi-compact very placid stack, weakly coh. pro-smooth over $k$
 such that the diagonal $\Delta_{X}\colon X\rightarrow X\times X$ is representable coh. pro-smooth, and ess. coh. pro-unipotent.
Let $Y$ be a quasi-compact sind-very placid stack and let $f\colon X\rightarrow Y$ be a ind-pfp proper morphism such that the relative diagonal $X\rightarrow X\times_{Y} X$ is also ind-pfp proper. Let $\phi_{X}\colon  X\rightarrow X$ and $\phi_{Y}\colon  Y\rightarrow Y$ be endomorphisms intertwined by $f$. 
Then there is a canonical fully faithful functors
\begin{align*}
    \trg(\rshv(X\times_{Y} X,\La),\phi) & \hookrightarrow \rshv(\mathcal{L}_{\phi}(Y),\La),
\end{align*}
with the essential image generated (as presentable $\La$-linear categories) by the essential of the functor $q_*^{\ind\fg}\circ (\delta_0)^{\ind\fg,!}$, where $\delta_0$ and $q$ are as in  \eqref{eq:trace-convolution-horocycle-diagram-special}.
\end{proposition}

\begin{proof}
We restrict our sheaf theory \eqref{eq:renormalized sheaf on sifted} to sind-very placid stacks.
We verify assumptions as in \Cref{bimodule.geom.trace.fully.faithfulness.convolution.cat} hold. We take the class $\mathrm{VR}$ as in \Cref{assumptions.base.change.sheaf.theory.V} to be the class of ind-pfp proper morphisms. Then \Cref{assumptions.base.change.sheaf.theory.V} \eqref{assumptions.base.change.sheaf.theory.V-1} and \eqref{assumptions.base.change.sheaf.theory.V-2} hold as for ind-pfp proper morphisms $(\ind\fg,*)$-pushforwards are left adjoints of $(\ind\fg,!)$-pullbacks.
We take the class $\mathrm{HR}$ as in \Cref{assumptions.base.change.sheaf.theory.H} to be the class of weakly coh. pro-smooth morphisms. Then \Cref{assumptions.base.change.sheaf.theory.H} \eqref{assumptions.base.change.sheaf.theory.H-1} and \eqref{assumptions.base.change.sheaf.theory.H-2} hold.
\end{proof}

\end{document}